\documentclass{memo-l}

\usepackage{amsmath}
\usepackage{amsfonts}
\usepackage{amssymb}
\usepackage{comment}
\usepackage{graphicx,color}
\usepackage{hyperref}

\newtheorem{theorem}{Theorem}[chapter]
\newtheorem{lemma}[theorem]{Lemma}

\theoremstyle{definition}
\newtheorem{definition}[theorem]{Definition}

\theoremstyle{remark}

\numberwithin{section}{chapter}
\numberwithin{equation}{chapter}

\newcommand{\Q}{{\mathbb Q}}
\newcommand{\Z}{{\mathbb Z}}

\newcommand{\z}{\zeta}

\newcommand{\om}{\omega}
\newcommand{\Om}{\Omega}

\newcommand{\G}{\Gamma}

\renewcommand{\b}{{\mathfrak b}}

\newcommand{\al}{\alpha}
\newcommand{\be}{\beta}
\newcommand{\ga}{\gamma}
\newcommand{\gac}{\gamma_c}

\newcommand{\q}{{\mathfrak q}}

\newcommand{\eps}{\varepsilon}

\newcommand{\ov}{\overline}
\newcommand{\lgs}[2]{\mbox{$\left(\frac{#1}{#2}\right)$}}

\DeclareMathOperator{\sign}{sign}
\DeclareMathOperator{\Ai}{Ai}
\DeclareMathOperator{\Bi}{Bi}
\DeclareMathOperator{\asin}{asin}
\DeclareMathOperator{\acos}{acos}
\DeclareMathOperator{\atan}{atan}
\DeclareMathOperator{\asinh}{asinh}
\DeclareMathOperator{\acosh}{acosh}
\DeclareMathOperator{\atanh}{atanh}
\DeclareMathOperator{\cotan}{cotan}
\DeclareMathOperator{\cotanh}{cotanh}
\DeclareMathOperator{\eintone}{eint1}
\DeclareMathOperator{\Li}{Li}
\DeclareMathOperator{\erf}{erf}
\DeclareMathOperator{\erfc}{erfc}

\DeclareMathOperator{\CS}{CS}
\DeclareMathOperator{\sn}{sn}

\DeclareMathOperator{\cn}{cn}
\DeclareMathOperator{\nc}{nc}
\DeclareMathOperator{\dn}{dn}
\DeclareMathOperator{\nd}{nd}
\DeclareMathOperator{\sd}{sd}

\DeclareMathOperator{\cd}{cd}
\DeclareMathOperator{\dc}{dc}
\DeclareMathOperator{\ssc}{sc}

\DeclareMathOperator{\cm}{cm}
\DeclareMathOperator{\sm}{sm}
\DeclareMathOperator{\cl}{cl}
\DeclareMathOperator{\ssl}{sl}
\newtheorem{proposition}[theorem]{Proposition}
\newtheorem{corollary}[theorem]{Corollary}
\newtheorem{remarks}[theorem]{Remarks}
\newtheorem{cf}{}[section]

\newcommand{\til}{\raise-0.3em\hbox{\~{}}} 

\def\cancelignorespaces{\leavevmode}

\makeindex

\begin{document}

\frontmatter

\title{Continued Fractions of Polynomial Type:\\ Theory and Encyclopedic Dictionary}


\author{Henri Cohen}
\address{Universit\'e de Bordeaux, CNRS, INRIA,
  Bordeaux INP, IMB, UMR 5251,
  F-33400 Talence, FRANCE}
\curraddr{}
\email{henri.cohen2@free.fr}
\thanks{}

\date{}

\subjclass[2020]{Primary }

\keywords{}

\dedicatory{Dedicated to the memory of my friends\\[2pt]
  Michel Mend\`es-France and Alf van der Poorten,\\[2pt]
  great lovers of ``neverending fractions''.}

\begin{abstract}
  After giving a number of properties of continued fractions of polynomial
  type, in particular focusing on convergence properties and
  Bauer--Muir--Ap\'ery acceleration techniques, we give a large list of
  continued fractions, both for specific real numbers,
  and for special functions, some extracted from a number of different
  sources, but most others being probably new. In addition to providing such a
  list, one of our main additions is to include the exact speed of convergence
  of these continued fractions (sometimes only up to a multiplicative
  constant), which is almost always omitted in the literature.
\end{abstract}

\maketitle

\tableofcontents

\chapter*{Preface}

The theory of continued fractions is one of the oldest in mathematics, since
implicitly the Greeks already knew the simple continued fraction of the
Golden Ratio $\phi=(1+\sqrt{5})/2=1+1/(1+1/(1+1/(1+\cdots)))$. Countless
books are devoted to the subject (I have included a few in the bibliography,
but many more exist) which treat different aspects of the subject:
convergence problems, explicit continued fractions for interesting constants
and functions, applications in other fields and for solving certain kinds of
problems, combinatorial aspects, number-theoretic aspects in particular
linked to Diophantine approximation, algorithmic uses, for instance in the
numerical computation of inverse Mellin transforms, etc...

\smallskip

Nevertheless, the need for the present book arose for several reasons.

\begin{enumerate}
\item First, although convergence theorems often form a large and difficult
  part of the literature, the practical applications of the corresponding
  results are rarely given. For instance, in some books giving lists of
  continued fractions, the speed of convergence is mentioned in a
  non-mathematical way, either using vague words like ``slowly'' or ``rapidly''
  convergent, or by giving tables showing numerically how suitable convergents
  approximate the limit.
\item Second, although some books such as \cite{Cuyt} do contain a large
  number of continued fractions, they are far from complete. Inasmuch as
  possible, I have tried to collect everything that I could from the
  abundant literature, and I have also added a large number which seem to be
  new. In addition, as mentioned, I give a precise estimate of their speed of
  convergence.
\item Third, it has apparently not been realized until now that the special
  method used by R.~Ap\'ery to prove the irrationality of $\z(3)$ is in fact
  applicable to a \emph{large} number of continued fractions existing in
  the literature, including continued fractions for \emph{functions}, thus
  providing a large number of completely new ones with
  faster convergence properties. In many cases the corresponding explicit
  formulas are rather unwieldly, so are not given, but are included when they
  are reasonably simple.
\item Finally, since this book contains more than 1900 continued fractions,
  most of them new, it would be wasteful for the reader to have to
  painstakingly copy each of the ones that interests him/her. We thus also
  provide a database of all the CFs and relevant information in a file
  readable by {\tt Pari/GP}, but also by any other software with a suitable
  interface, since the details of the contents and format of this file are
  described in Appendix \ref{chap:file}.
\end{enumerate}

\smallskip

The work described in this book would not have been possible without
a specialized package for working with continued fractions, which will be
described in detail in Chapter \ref{chap:numerics}. This is part of the
{\tt Pari/GP} computer algebra system, and for now (2026) is only available in
the {\tt GIT} branch {\tt henri-ellCF3}, but will hopefully
at some point be available in the main distribution.

\smallskip

The author would evidently appreciate any corrections and additions which
would be included in future editions and on the author's website. Note that
since most of the given continued fractions are new, and almost all the
additional data (exact speed of convergence, asymptotic expansion) do not
exist at all in the literature, it is unavoidable that errors (typos or
mathematical mistakes) have cropped up, which renders user feedback absolutely
essential.

Ideally one should put all this data on a web site which could be
regularly updated by users, but I do not have the knowledge nor the time
to do this.

\medskip

{\bf Acknowledgments} I would first like to thank Wadim Zudilin for numerous conversations, and more
recently for leading me to discover numerous rapidly convergent continued fractions for Chowla--Selberg
gamma quotients, leading to many new irrationality results (see {\tt arXiv:2510.00215v5}). I thank
several people for sending me new interesting CFs, in particular Alex Shvets. Finally, special thanks go to
Francesco Del Re who, using AI, helped me correct a large number of mistakes and typos (perhaps all?).

\medskip

\aufm{Henri Cohen}

\chapter*{Brief Personal Recollections}

\medskip

This book is dedicated to my friends and colleagues Michel Mend\`es-France
and Alf van der Poorten, who passed away in 2018 and 2010 respectively.

I came to settle in Bordeaux in part thanks to the friendly advice of
Michel, and my first serious mathematical paper (my very first one was more
recreational) was directly inspired by his work on continued fractions:
he had studied how the length of a finite CF behaves when multiplied by an
integer. I studied the corresponding problem for the period of a periodic CF,
and for instance proved that, when multiplied by $2$, the length of the
period is at most multiplied by $5$, and $5$ is best possible, see
\cite{Coh1}.

\smallskip

Alf came regularly to Bordeaux, in particular at the invitation of Michel,
and we had numerous mathematical discussions, in particular related to
continued fractions. For instance, I remember telling him that I was surprised
that no textbook on CFs mentions a \emph{necessary and sufficient condition}
for an irreducible fraction $p/q$ to be a partial quotient of the CF of some
irrational number $\al$. These texts all mention that a \emph{necessary}
condition is that $|\al-p/q|\le 1/q^2$, and that a \emph{sufficient} (but
certainly not necessary) condition is that $|\al-p/q|\le 1/(2q^2)$.
But in fact a necessary and sufficient condition is that
$-1/(q(q+p')<\al-p/q<1/(q(2q-p'))$, where $p'$ is the inverse of $p$
modulo $q$ such that $1\le p'\le q$, see for instance Exercise
21 of Chapter 5 of \cite{Coh2}.

Among the many things that Alf taught me was the marvelous subject of
``unlikely intersections'', related simultaneously to torsion points on
Jacobians of curves, to the periodicity of continued fractions of square roots
of polynomials, to Somos sequences, and to the integrability in elementary
terms of hyperelliptic integrals, see for instance Chapter 5 of \cite{BPSZ}.

\aufm{Henri Cohen}

\mainmatter
\chapter{Basic Results on Continued Fractions}\label{chap:basic}

\section{Notation}

Although we assume the reader familiar with the usual terminology and basic
results on continued fractions, we prefer starting from scratch and in
particular explain our specific notation. First, note that to avoid always
writing ``continued fraction'', we will often use the abbreviation ``CF''.

Let $(a(n))_{n\ge0}$ and $(b(n))_{n\ge0}$ be two (finite or infinite)
sequences of complex numbers. The continued fraction $(a,b)$ is the formal
expression
$$(a,b)=a(0)+\dfrac{b(0)}{a(1)+\dfrac{b(1)}{a(2)+\dfrac{b(2)}{a(3)+\ddots}}}$$
which we can write more compactly as
$(a,b)=a(0)+b(0)/(a(1)+b(1)/(a(2)+\cdots))$.
Note that this notation is different from that used by other authors
such as \cite{Cuyt}. I do not claim that my notation
is better but simply that it is the one used in the present book and by
the continued fraction package of {\tt Pari/GP} (currently
the {\tt GIT} branch {\tt origin/henri-ellCF3}).

If $a(n)$ and $b(n)$ are finite sequences, this is simply a rational function
in the coefficients, so not particularly interesting. Similarly if $b(n)=0$
for some $n$. We will therefore always assume that the sequences are
infinite and that $b(n)\ne0$ for all $n$. Although not strictly necessary
we also assume that $a(n)=0$ only for a \emph{finite} number of $n$.

We write as usual formally
$$p(n)/q(n)=a(0)+b(0)/(a(1)+b(1)/(a(2)+\cdots+b(n-1)/a(n)))\;,$$
so that $p(0)=a(0)$, $q(0)=1$, and setting by convention $p(-1)=1$ and
$q(-1)=0$ we have the usual recursions
$$(p(n+1),q(n+1))=a(n+1)(p(n),q(n))+b(n)(p(n-1),q(n-1))\;,$$
or equivalently in matrix form
$$\begin{pmatrix}p(0)&p(-1)\\q(0)&q(-1)\end{pmatrix}=\begin{pmatrix}a(0)&1\\1&0\end{pmatrix}\text{\quad and}$$
$$\begin{pmatrix}p(n+1)&p(n)\\q(n+1)&q(n)\end{pmatrix}=\begin{pmatrix}p(n)&p(n-1)\\q(n)&q(n-1)\end{pmatrix}\begin{pmatrix}a(n+1)&1\\b(n)&0\end{pmatrix}\text{\quad for $n\ge0$}\;,$$
  which implies by induction that
  $$p(n+1)q(n)-p(n)q(n+1)=(-1)^n\prod_{0\le j\le n}b(j)\;.$$
We say that the continued fraction \emph{converges} to a limit $S$ if
$S=\lim_{n\to\infty}p(n)/q(n)$ exists. To estimate the
\emph{speed of convergence} means to give an estimate for $S-p(n)/q(n)$.

As a matter of terminology, $a(n)$ and $b(n)$ are called the numerators
and denominators of the CF respectively, or simply its coefficients.
$(p(n),q(n))$ are called the \emph{partial quotients}, and $p(n)/q(n)$ the
\emph{convergents}.

\smallskip

We will use the following practical notation:

\begin{definition} For any arithmetic function $f(n)$ we define
  $f!(n)=\prod_{1\le j\le n}f(j)$, so in particular $f!(0)=1$.
\end{definition}

Note that we do not include the value $f(0)$, even if it is defined, so
that for instance $\prod_{0\le j\le n-1}f(j)=f(0)f!(n-1)$.

\begin{definition} We say that two continued fractions $(a(n),b(n))$
  and $(a'(n),b'(n))$ are \emph{equivalent} if there exists a nonzero
  sequence $r(n)$ with $r(0)=1$ such that $a'(n)=r(n)a(n)$ and
  $b'(n)=r(n)r(n+1)b(n)$.\end{definition}

The main reason for this notion is the following lemma:

\begin{lemma} If two continued fractions are equivalent they have the
  same convergents $p(n)/q(n)$, and in particular the same limits and
  convergence speed if they converge.\end{lemma}

\begin{proof} Indeed, using the notation in the definition, it is immediate to
  check by induction that if $(p(n),q(n))$ and $(p'(n),q'(n))$ are the
  corresponding partial quotients we have $(p'(n),q'(n))=r!(n)(p(n),q(n))$.
\end{proof}

\begin{remarks}
  \begin{enumerate}\item In particular note that if $a(n)$ and $b(n)$ are
    rational functions, by using a suitable $r(n)$ we can find an equivalent
    continued fraction $(a'(n),b'(n))$ where $a'(n)$ and $b'(n)$ are
    polynomials.
  \item This notion of equivalence is \emph{not} the same as one which is
    often used in the case of simple continued fractions, where it means
    that the CFs have the same \emph{tails} (i.e., there exists $n_0$ such
    that $a'(n)=a(n+n_0)$ for $n$ sufficiently large), or almost equivalently
    that their limits are linked by a M\"obius transformation
    $z'=(Az+B)/(Cz+D)$.
\end{enumerate}\end{remarks}

\smallskip

Finally note the important transformation of continued fractions
called \emph{contractions} (even or odd): this consists in constructing
the natural CF whose partial quotients are $(p(2n),q(2n))$
(or $(p(2n+1),q(2n+1))$ for the odd contraction). For instance, the even contraction is obtained by the continued fraction $(A(n),B(n))$ with
\begin{align*}
  A(0)&=a(0),\ A(1)=a(1)a(2)+b(1)\,\\
  A(n)&=a(2n)a(2n-1)+b(2n-1)+\dfrac{a(2n)b(2n-2)}{a(2n-2)}\text{\quad for $n\ge2$}\\
  B(0)&=b(0)a(2),\ B(n)=-\dfrac{a(2n+2)b(2n)b(2n-1)}{a(2n)}\text{\quad for $n\ge1$\;.}
\end{align*}
The presence of denominators in these formulas means that usually the
corresponding expressions are rational functions,
hence more complicated than the original continued fraction. Thus,
when given continued fractions of period $2$ (see definition below), which
are extremely frequent, although it may be tempting to use contractions to
obtain CFs of period 1, it is often not a good idea.

\section{Types of Continued Fractions}

The simplest type of continued fraction (and historically the first to be
considered) is the \emph{simple} type: a CF is simple if $b(n)=1$ for all
$n$ and $a(n)$ is a strictly positive integer for $n\ge1$ ($a(0)$ can be
arbitrary). This is the type which is used in most elementary questions
of Diophantine approximation because of the following standard properties:
\begin{enumerate}
\item The simple CF expansion of a rational number is finite, and conversely.
\item The simple CF expansion of a real quadratic irrational ($z=a+b\sqrt{d}$
  with $a$, $b$ rational with $b\ne0$ and $d$ a non square positive integer) is periodic
  and conversely.
\item A simple CF always converges (at an exponential rate). In addition,
  if $p(n)/q(n)$ are the convergents of the CF expansion of an irrational
  number $z$, we have $q(n)p(n-1)-p(n)q(n-1)=(-1)^n$, $|z-p(n)/q(n)|<1/q(n)^2$
  (more precise inequalities exist), and in a suitable sense the rational
  numbers $p(n)/q(n)$ are the ``best'' rational approximations to $z$.
\end{enumerate}

\smallskip

Although ``simple'' continued fractions are essential for Diophantine
questions, they are only a tiny fraction of the theory, and it is natural
to consider the more general CFs as presented above. Before coming to
the polynomial type, which is essentially the only type studied in the
present book, we give a few samples of other types, of course far from
being exhaustive.

\begin{enumerate}\item
  Simple continued fractions are intimately linked to the Euclidean
  algorithm and the identity $a=bq+r$ with $0\le r<b$. But we can also
  consider \emph{negative} continued fractions, corresponding to the
  above identity with $-b<r\le0$: this corresponds to choosing $b(n)=-1$
  for all $n$. The first two properties above are preserved (not the third),
  but perhaps what is more interesting is the link with deeper theories,
  for instance in topology and number theory. More generally, we can consider
  $b(n)=e$ a constant, and in that case the corresponding CF of a quadratic
  irrational is \emph{not} necessarily periodic (for instance $e=3$ and
  $z=\sqrt{29}$; the reader can play the game of finding for each small
  positive or negative integer $e$ the list of nonsquare $d\le100$ such
  that $z=\sqrt{d}$ does not apparently have a periodic CF with $b(n)=e$).
\item A very common type of CF that we will not consider is what one
  can call $q$-continued fractions. One of the most famous is the
  Ramanujan CF with $a(0)=0$, $a(n)=1$ for $n\ge1$, $b(0)=q^{1/5}$, and
  $b(n)=q^n$ for $n\ge1$, in other words
  $$R(q)=q^{1/5}/(1+q/(1+q^2/(1+q^3/(1+\cdots))))\;.$$
  This has many modular interpretations and infinite product representations,
  but one of its most striking properties is its \emph{CM-properties},
  in that it takes explicit \emph{algebraic} values when $q=e^{i\pi\tau}$ with
  $\tau$ a quadratic irrational with strictly positive imaginary part:
  for instance $R(e^{-2\pi})$ is a positive root of the polynomial equation
  $x^4+2x^3-6x^2-2x+1=0$.
\item Still other types of CFs are when one chooses $a(n)$ and/or $b(n)$
  to be certain special types of sequences such as \emph{automatic sequences}.
\end{enumerate}

In the present book we will only consider CFs of \emph{polynomial type},
according to the following definition:

\begin{definition} A CF $((a(n))_{n\ge0},(b(n))_{n\ge0})$ is said to be
  of polynomial type if there exists a \emph{period} $T$ such that for
  every $j$ with $0\le j<T$, for $n$ sufficiently large both
  $a(nT+j)$ and $b(nT+j)$ are rational functions of $n$ with
  \emph{rational} coefficients.\end{definition}

The name \emph{polynomial} type is justified since as we have mentioned above,
any CF with rational function coefficients is equivalent to one with polynomial
coefficients, and since we assume rational coefficients, if desired we may
also assume that the polynomials have \emph{integer} coefficients.

In practice $T$ will be equal to $1$ or $2$ (in very exceptional cases $3$,
$4$, or $5$), so for $n$ sufficiently large either $a(n)$ and $b(n)$ are
rational functions, or $a(2n)$, $a(2n+1)$, $b(2n)$, and $b(2n+1)$ are
rational functions.

The only case where $T$ is often large is for continued fractions
associated to real quadratic irrationals, which we will not consider.

\section{Transformations of Series and Products into CFs}

It is well-known at least since Euler that one can transform a series
into a continued fraction in trivial ways, so that the $n$th partial sum
of the series is equal to the $n$th convergent of the continued fraction.
Let us see the different ways in which this is done.

\begin{definition} We will say that a continued fraction is termwise
  equal to an infinite series (resp., to an infinite product) if its $n$th
  convergent $p(n)/q(n)$ is equal to the $n$th partial sum of the series
  (resp., to the $n$th partial product of the infinite product).
\end{definition}

\begin{proposition}[Euler]\label{prop:euler} Let $f$ and $g$ be two arithmetic
  functions, and
  let $S$ be the sum of the series $S=\sum_{m\ge1}g!(m-1)/f!(m)$, assuming
  it converges, using the same factorial notation as above. Then if we set
  \begin{align*}
    a(0)&=0,\ a(1)=f(1),\ a(n)=f(n)+g(n-1)\text{\quad for $n\ge2$}\\
    b(0)&=1,\ b(n)=-f(n)g(n)\text{\quad for $n\ge1$}\;,
  \end{align*}
  the continued fraction $(a(n),b(n))$ is termwise equal to the series $S$.
\end{proposition}

\begin{proof} It is immediate to check by induction that if $(p(n),q(n))$
  are the partial quotients associated to the CF $(a,b)$, we have
  $q(n)=f!(n)$ and $p(n)=f!(n)\sum_{1\le m\le n}g!(m-1)/f!(m)$, so
  $p(n)/q(n)$ is indeed equal to the $n$th partial sum of $S$.\end{proof}

The following variants and special cases can be useful:

\begin{corollary}\label{cor:euler} Keep the above notation.
  The following series $S$ and CFs $(a(n),b(n))$ are termwise equal:
\begin{enumerate}\item $S=\sum_{m\ge0}g!(m)/f!(m)$ and
  \begin{align*}
    a(0)&=1,\ a(1)=f(1),\ a(n)=f(n)+g(n)\text{\quad for $n\ge2$\;,}\\
    b(0)&=g(1),\ b(n)=-f(n)g(n+1)\text{\quad for $n\ge1$}\;.
  \end{align*}
  Equivalently, if $S=\sum_{m\ge0}u(m)$ with $u(0)=1$ and if
  $u(m)/u(m-1)=g(m)/f(m)$ for $m\ge1$, the same formulas hold.
\item $S=\sum_{m\ge1}z^m/f(m)$ and
  \begin{align*}
    a(0)&=0,\ a(1)=f(1),\ a(n)=f(n)+zf(n-1)\text{\quad for $n\ge2$\;,}\\
    b(0)&=z,\ b(n)=-zf(n)^2\text{\quad for $n\ge1$}\;.
  \end{align*}
\item $S=\sum_{m\ge0}f(m)z^m$ and
  \begin{align*}
    a(0)&=f(0),\ a(1)=1,\ a(n)=f(n-1)+zf(n)\text{\quad for $n\ge2$\;,}\\
    b(0)&=zf(1),\ b(1)=-zf(2),\ b(n)=-zf(n-1)f(n+1)\text{\quad for $n\ge2$}\;.
  \end{align*}
\item Let $g$ be any arithmetic function such that $\lim_{n\to\infty}g(n)=0$.
  Then $S=\sum_{m\ge0}f(m)$, and 
  \begin{align*}
    a(0)&=f(0)+g(0),\ a(1)=1,\ a(n)=f(n-1)+f(n)+g(n)-g(n-2)\;,\\
    b(0)&=f(1)+g(1)-g(0),\ b(1)=-(f(2)+g(2)-g(1)),\\
    b(n)&=-(f(n-1)+g(n-1)-g(n-2))(f(n+1)+g(n+1)-g(n))\;,\end{align*}
  the formulas for $a(n)$ and $b(n)$ being for $n\ge2$.
\end{enumerate}\end{corollary}

Because of these (trivial) results, any constant or function which is
given by a \emph{hypergeometric} expression (i.e., the ratio of two
consecutive terms is a rational function of $m$) has automatically
a corresponding CF. Since a very large number of constants and functions are
given in this way, this provides correspondingly a very large number of CFs.
These CFs are not interesting per se (no more and no less than the
corresponding hypergeometric series), but one of the main interests lies
in the fact that these CFs can be \emph{transformed} into completely new
CFs, using transformations which are not possible (at least not directly)
on the series. We will see numerous examples of this.

\medskip

Although less useful in practice, similar results hold for infinite
\emph{products} instead of infinite series: indeed, we have the trivial
identity
$$\prod_{n\ge1}(1+h(n))=1+\sum_{n\ge1}h(n)\prod_{1\le m\le n-1}(1+h(m))\;.$$
For instance, from $\sin(\pi z)/(\pi z)=\prod_{n\ge1}(1-z^2/n^2)$ we deduce
the identity
$$\dfrac{\sin(\pi z)}{\pi z}=1-\dfrac{z^2}{1!^2}-\dfrac{z^2(1^2-z^2)}{2!^2}-\dfrac{z^2(1^2-z^2)(2^2-z^2)}{3!^2}-\cdots$$
We thus obtain the following results:

\begin{corollary} The following infinite products $P$ and CFs $(a(n),b(n))$
  are termwise equal:
  \begin{enumerate}\item $P=\prod_{n\ge1}(1+g(n)/f(n))$ and
\begin{align*}
  a(0)&=1,\ a(1)=f(1),\\
  a(n)&=f(n)g(n-1)+g(n)(f(n-1)+g(n-1))\text{\quad for $n\ge2$\;,}\\
  b(0)&=g(1),\ b(n)=-g(n-1)g(n+1)f(n)(f(n)+g(n))\text{\quad for $n\ge1$\;,}
\end{align*}
where by convention we set $g(0)=1$.
\item $P=\prod_{n\ge1}(1+(-1)^ng(n)/f(n))$ and
\begin{align*}
  a(0)&=1,\ a(1)=f(1),\\
  a(n)&=-f(n)g(n-1)+g(n)(f(n-1)-g(n-1))\text{\quad for $n\ge2$ even},\\
  a(n)&=f(n)g(n-1)-g(n)(f(n-1)+g(n-1))\text{\quad for $n\ge3$ odd}\\
  b(0)&=-g(1),\ b(n)=-g(n-1)g(n+1)f(n)(f(n)+g(n))\text{\quad for $n\ge2$ even}\\
  b(n)&=-g(n-1)g(n+1)f(n)(f(n)-g(n))\text{\quad for $n\ge1$ odd\;.}
\end{align*}
\item $P=\prod_{n\ge1}(1+z/f(n))$ and
\begin{align*}
  a(0)&=1,\ a(1)=f(1),\ a(n)=z+f(n-1)+f(n)\text{\quad for $n\ge2$}\\
  b(0)&=z,\ b(n)=-f(n)(f(n)+z)\text{\quad for $n\ge1$}\;.
\end{align*}
\item Let $g$ be any arithmetic function such that $\lim_{n\to\infty}g(n)=1$.
Then $P=\prod_{n\ge1}h(n)$ and
  \begin{align*}
    a(0)&=g(0),\ a(1)=1,\ a(n)=h(n-1)h(n)g(n)-g(n-2)\text{\quad for $n\ge2$}\\
    b(0)&=f(1),\ b(1)=-h(1)f(2),\ b(n)=-h(n)f(n-1)f(n+1)\text{\quad for $n\ge2$}\;,\end{align*}
where for simplicity we define $f(n)=h(n)g(n)-g(n-1)$.
\end{enumerate}\end{corollary}

\section{Examples of Series and Products}

We give two examples of series. First, $\z(3)=\sum_{m\ge1}1/m^3$.
Using either the first transformation with $f(n)=g(n)=n^3$, or equivalently
the third with $f(n)=n^3$ and $z=1$, we obtain
$a(0)=0$, $a(n)=2n^3-3n^2+3n-1$ for $n\ge1$, $b(0)=1$, and
$b(n)=-n^6$ for $n\ge1$, whence the (a priori uninteresting)
continued fraction
$$\z(3)=\dfrac{1}{1-\dfrac{1}{9-\dfrac{64}{35-\dfrac{729}{91-\dfrac{4096}{189-\dfrac{15625}{341-\ddots}}}}}}$$
with convergence type identical to that of the series, in $1/(2n^2)$.

As mentioned above, the main point of transforming the series into this CF is
that it can now be \emph{transformed}, and in particular be \emph{accelerated}
into new CFs, and this is what led Ap\'ery to his famous proof of irrationality
of $\z(3)$.

Second, we choose $f(n)=n$, and $S=(1/z)\sum_{m\ge1}z^m/m=-\log(1-z)/z$.
We obtain $a(0)=0$, $a(1)=1$, $a(n)=n(1+z)-z$ for $n\ge2$,
$b(0)=1$, and $b(n)=-n^2z$ for $n\ge1$. Changing $z$ into $-z$ and multiplying
by $z$ we obtain the continued fraction
$$\log(1+z)=\dfrac{z}{1+\dfrac{z}{-z+2+\dfrac{4z}{-2z+3+\dfrac{9z}{-3z+4+\dfrac{16z}{-4z+5+\dfrac{25z}{-5z+6+\ddots}}}}}}$$
with convergence type identical to that of the power series, i.e., for
$|z|<1$ and $z=1$. In particular, with $z=1$ and $z=-1/2$ we obtain two
beautiful continued fractions for $\log(2)$:
$$\log(2)=\dfrac{1}{1+\dfrac{1}{1+\dfrac{4}{1+\dfrac{9}{1+\dfrac{16}{1+\dfrac{25}{1+\ddots}}}}}}$$
which converges like the alternating sum slowly in $(-1)^n/n$, and
$$\log(2)=\dfrac{1}{2-\dfrac{2}{5-\dfrac{8}{8-\dfrac{18}{11-\dfrac{32}{14-\dfrac{50}{17-\ddots}}}}}}$$
which converges exponentially in $1/(n2^n)$.

\smallskip

We also give two examples of products, both starting from the identity
$$\dfrac{\sin(\pi z)}{\pi z}=\prod_{n\ge1}\left(1-\dfrac{z^2}{n^2}\right)\;.$$
First, we write $(1-z^2/n^2)=(1-z/n)(1+z/n)$, so we are in the situation
of $P=\prod_{n\ge1}(1+(-1)^ng(n)/f(n))$ with  $g(n)=z$ and $f(n)=n/2$ for
$n$ even and $f(n)=(n+1)/2$ for $n$ odd. After simplification by
$z$ we deduce the continued fraction
$$\dfrac{\sin(\pi z)}{\pi z}=1-\dfrac{z}{1+\dfrac{z-1}{-z-\dfrac{z+1}{-z+1+\dfrac{2z-4}{-z-\dfrac{2z+4}{-z+1+\dfrac{3z-9}{-z-\ddots}}}}}}$$
with very slow convergence type essentially in $1/n$ up to sign. Choosing
$z=1/2$ or $z=1/6$ gives nice but very slowly converging continued fractions
for $\pi$.

\smallskip

We now use directly $\prod_{n\ge1}(1-z^2/n^2)$, so we are in the special
case given above with $z$ replaced by $-z^2$ and $f(n)=n^2$, so that
$a(0)=1$, $a(1)=1$, $a(n)=2n^2-2n+1-z^2$ for $n\ge2$, $b(0)=-z^2$, and
$b(n)=-n^2(n^2-z^2)$ for $n\ge1$, giving the continued fraction
$$\dfrac{\sin(\pi z)}{\pi z}=1-\dfrac{z^2}{1+\dfrac{z^2-1}{-z^2+5+\dfrac{4z^2-16}{-z^2+13+\dfrac{9z^2-81}{-z^2+25+\dfrac{16z^2-256}{-z^2+41+\dfrac{25z^2-625}{-z^2+61+\ddots}}}}}}$$
and the reader can easily check that this is essentially identical to the
contraction of the previous one.

\section{A Generalization}

Let $P(n)$ be for now an arbitrary arithmetic function, and let us
choose $f(n)=\prod_{0\le j\le m}P(n+j)^{e_j}$ and $g(n)=z f(n)$ for some $z$.
We have
\begin{align*}f(n)+g(n-1)&=\prod_{0\le j\le m}P(n+j)^{e_j}+z\prod_{0\le j\le m}P(n+j-1)^{e_j}\\
  &=\prod_{0\le j\le m}P(n+j)^{e_j}+zP(n-1)^{e_0}\prod_{0\le j\le m-1}P(n+j)^{e_{j+1}}\\
  &=\prod_{0\le j\le m-1}P(n+j)^{\min(e_j,e_{j+1})}R(n)\;,\text{\quad with}\\
  R(n)&=P(n+m)^{e_m}\prod_{\substack{0\le j\le m-1\\e_j>e_{j+1}}}P(n+j)^{e_j-e_{j+1}}\\&\phantom{=}+zP(n-1)^{e_0}\prod_{\substack{0\le j\le m-1\\e_j<e_{j+1}}}P(n+j)^{e_{j+1}-e_j}\;,\end{align*}
and of course $$-f(n)g(n)=-z\prod_{0\le j\le m}P(n+j)^{2e_j}\;.$$
Choosing $r(n)=1/\prod_{0\le j\le m-1}P(n+j)^{\min(e_j,e_{j+1})}$ (assuming
of course that the denominator does not vanish), we compute that the CF
is equivalent to one with
$$a(n)=R(n)\text{\quad and\quad}b(n)=-z\prod_{0\le j\le m}P(m+j)^{\max(0,e_j-e_{j+1})+\max(0,e_j-e_{j-1})}\;,$$
where by convention we set $e_{-1}=e_{m+1}=0$.

As a fundamental example, let us take $P(n)=n$, and $1\le e_j\le 2$. By
decomposing into partial fractions, it is clear that when $z=\pm1$, the sum
$\sum_{n\ge1}g!(n-1)/f!(n)$ will be a $\Q$-linear combination of $1$ and
$\pi^2$, and in addition $\log(2)$, so we will obtain CFs for $\pi^2$ for
$z=1$, and some linear combination of $\pi^2$ and $\log(2)$ for $z=-1$.
It is reasonable to restrict to $\deg(a)\le 2$ and $\deg(b)\le 4$, and it then
easy to see that this will be the case if and only if the sequence of $e_j$
is of the form $(1,1,\dotsc,1,2,2\dotsc,2,1,1,\dotsc,1)$ with $m_1$ initial
ones, $m_3$ final ones, and $m_2$ twos, $m_i=0$ being accepted.

We easily check that, for any $P$, we have
\begin{align*}
  a(n)&=R(n)=P(n+m_1+m_2-1)P(n+m-1)+zP(n-1)P(n+m_1-1)\text{\quad and}\\
  b(n)&=-zP(n)P(n+m_1)P(n+m_1+m_2-1)P(n+m-1)\;.
\end{align*}

\smallskip

Similarly, assume that $1\le e_j\le 3$. Here, we will obtain a $\Q$-linear
combination of $1$, $\pi^2$, and $\z(3)$ if $z=1$, and also $\log(2)$ if
$z=-1$. If we restrict to $\deg(a)\le3$ and $\deg(b)\le4$, the sequence
$e_j$ will necessarily be of the form $(1,\dotsc,2,\dotsc,3,\dotsc,2,\dotsc,1,\dotsc)$ with $m_1$, $m_2$, $m_3$, $m_4$, and $m_5$ ones, twos, threes, twos,
and ones respectively, and we easily check that

\begin{align*}
  a(n)&=R(n)=P(n+m_1+m_2+m_3-1)P(n+m_1+m_2+m_3+m_4-1)P(n+m-1)\\
  &\phantom{=}+zP(n-1)P(n+m_1-1)P(n+m_1+m_2-1)\text{\quad and}\\
  b(n)&=-zP(n)P(n+m_1)P(n+m_1+m_2)P(n+m_1+m_2+m_3-1)\cdot\\
  &\phantom{=}\cdot P(n+m_1+m_2+m_3+m_4-1)P(n+m-1)\;.
\end{align*}

\section{Knowing a Particular Solution}

Assume that we know explicitly an arithmetic function $f(n)$ satisfying
the recursion defining a given CF, i.e., such that
$$f(n+1)=a(n+1)f(n)+b(n)f(n-1)\;.$$
It is then easy to find the general solution of this recursion, and in
particular to give explicit formulas for the convergents of the CF:

\begin{proposition} Let $(a(n),b(n))$ be a CF with partial quotients
  $(p(n),q(n))$, and assume given a function $f(n)$ such that $f(n)\ne0$
  for all $n\in\Z_{\ge0}$ and satisfying
  $f(n+1)=a(n+1)f(n)+b(n)f(n-1)$. Then any solution $u(n)$ to this recursion,
  and in particular $p(n)$ and $q(n)$, is given explicitly by
  $$u(n)=f(n)\left(\dfrac{u(0)}{f(0)}+(f(0)u(1)-f(1)u(0))\sum_{1\le m\le n}(-1)^{m-1}\dfrac{b!(m-1)}{f(m)f(m-1)}\right)\;.$$
\end{proposition}

\begin{proof} Set $v(n)=u(n)/f(n)$. A small computation shows that
  $$\dfrac{v(n+1)-v(n)}{v(n)-v(n-1)}=-b(n)\dfrac{f(n-1)}{f(n+1)}\;,$$
  from which by induction we deduce that
  $$\dfrac{v(n+1)-v(n)}{v(1)-v(0)}=(-1)^nb!(n)\dfrac{f(1)f(0)}{f(n+1)f(n)}\;,$$
  and the proposition immediately follows.\end{proof}

The following corollary is clear:

\begin{corollary}\label{cor:expl} Keep the same notation, and set
  $$S(n)=\sum_{1\le m\le n}(-1)^{m-1}\dfrac{b!(m-1)}{f(m)f(m-1)}\;.$$
  We have
  \begin{align*}
    p(n)&=f(n)\left(\dfrac{a(0)}{f(0)}+(f(0)(a(1)a(0)+b(0))-f(1)a(0))S(n)\right)\;,\\
    q(n)&=f(n)\left(\dfrac{1}{f(0)}+(f(0)a(1)-f(1))S(n)\right)\;,\text{\quad and}\\
    \dfrac{p(n)}{q(n)}&=a(0)+f(0)^2b(0)\dfrac{S(n)}{1+f(0)(f(0)a(1)-f(1))S(n)}\;.\end{align*}
  In particular the CF converges to some limit $L$ if and only if either the
  series $S(n)$ converges to some limit $S$ and we have
  $$L=a(0)+f(0)^2b(0)\dfrac{S}{1+f(0)(f(0)a(1)-f(1))S}\;,$$
  or $|S(n)|$ tends to infinity and $f(0)a(1)-f(1)\ne0$, and we have
  $$L=a(0)+\dfrac{f(0)b(0)}{f(0)a(1)-f(1)}\;.$$
\end{corollary}

Note that for simplicity we have ignored cases where denominators vanish.

\smallskip

Let us for example consider the (admittedly trivial) CF with
$a(n)=n+1$ and $b(n)=-(n+1)$. It is clear that $f(n)=1$ satisfies the
recursion, so $S(n)=\sum_{1\le m\le n}m!$ which tends to $\infty$ as
$n\to\infty$. Since $a(0)=1$, $a(1)=2$, and $b(0)=-1$, we deduce that $L=0$.
This is the special case $k=0$, $z=1$ of \ref{6.2.1}.

For a slightly less trivial example, consider $a(n)=n+1$ for $n\ne1$, $a(1)=1$,
$b(n)=-n$ for $n\ge1$, $b(0)=1$. The recursion is
$u(n+1)=(n+2)u(n)-nu(n-1)$ for $n\ne0$, and $u(1)=u(0)-u(-1)$. It
is clear that $f(n)=n!$, $f(-1)=0$ satisfies the recursion. We have
$S(n)=\sum_{1\le m\le n}1/m!$, so $S=\exp(1)-1=e-1$, so
$L=1+S=e$, giving the CF $e=1+1/(1-1/(3-2/(4-3/(5-4/(6-\cdots)))))$,
which is simply the Euler transform of the usual series for $e$.

\chapter{Convergence Results}\label{chap:speed}

A large part of the contents of this chapter can also be found in Chapter 7
of \cite{Bel-Coh}.

\section{Introduction}

In all books devoted to continued fractions, a large part deals with
the difficult problems related to convergence. This can however be considerably
simplified if we impose the following restrictions: first, we only consider
CFs of polynomial type as defined in the previous chapter. And second, we
accept \emph{generic} and \emph{heuristic} results.

Let us give an example to explain these restrictions. The partial
quotients $(p(n),q(n))$ of a CF satisfy a $3$-term recursion, for instance
$q(n+1)=a(n+1)q(n)+b(n)q(n-1)$. As we will see below, studying the speed
of convergence of a CF is essentially equivalent to finding the asymptotic
behavior of $q(n)$. Consider for instance the recursion $u(n+1)=5u(n)-6u(n-1)$.
The general solution is $A3^{n+1}+B2^{n+1}$ for some constants $A$ and $B$
depending on the initial values $u(-1)$ and $u(0)$, more precisely
$A=u(0)-2u(-1)$ and $B=3u(-1)-u(0)$. Thus, if $u(0)\ne 2u(-1)$ we have
$u(n)\sim C\cdot 3^n$ for some nonzero constant $C=3A$, and this is what we
call the \emph{generic} solution. But in the special case $u(0)=2u(-1)$
we have $u(n)\sim C\cdot 2^n$ with $C=2B$, which is completely different
but special, i.e., nongeneric.

We will see the heuristic part in the way the proof of the main result is
obtained.

\section{The Main Result}

We will first assume that $T=1$, so that $a(n)$ and $b(n)$ are
rational functions of $n$ for $n$ sufficiently large. We write
$$a(n)=a_0n^{\al}(1+a_1/n+a_2/n^2+\cdots)\text{\quad and\quad}b(n)=b_0n^{\be}(1+b_1/n+b_2/n^2+\cdots)$$
with $a_0b_0\ne0$ (we only use $a_j$ and $b_j$ in this section, which
of course should not be confused with $a(j)$ and $b(j)$).

We will see below that all of our asymptotics
(excluding the exceptional case of logarithmic
convergence, see below) will be of the form
$$u(n)\sim n!^FE^ne^{(Dn)^{1/2}}n^PC\;,$$
where $C$, $D$, $E$, and $F$ are real constants with $D\ge0$, $E\ne0$ and
$C\ne0$, so to avoid a cumbersome notation we will simply write
$u(n)=[F,E,D,P]$, omitting the constant $C$. Typically this will be
applied to $u(n)=q(n)$ or $u(n)=1/(S-p(n)/q(n))$.

The main result is as follows. Although we call it a ``theorem'', we recall
that it is valid in a generic situation and that the proof is heuristic,
although it can be justified rigorously in many cases. Since it has a large
number of cases it is extremely technical but absolutely fundamental.

\begin{theorem}\label{thm:speed}
  Let $(a(n),b(n))$ be a convergent polynomial type continued
  fraction with limit $S$, $(p(n),q(n))$ the corresponding
  partial quotients, and keep the above notation. We have the following
  convergence types:
  \begin{enumerate}\item The \emph{factorial} type $F>0$ (type $F$): this
    occurs when $\be\le 2\al-1$, and we have two subtypes:
    \begin{enumerate}\item If $\be\le 2\al-2$ we have
      $$(q,r)=([\al,a_0,0,a_1],[2\al-\be,-a_0^2/b_0,0,\al+2a_1-b_1])\;,$$
      where here and below $(q,r)$ is an abbreviation for
      $(q(n),1/(S-p(n)/q(n)))$.
    \item If $\be=2\al-1$ we have
      $$(q,r)=([\al,a_0,0,a_1+b_0/a_0^2],[2\al-\be,-a_0^2/b_0,0,\al+2a_1+2b_0/a_0^2-b_1])\;.$$      
    \end{enumerate}
  \item The \emph{exponential} type $F=0$, $|E|>1$ (type $E$): this occurs when
    $\be=2\al$ and $\Delta=a_0^2+4b_0>0$, and setting $\delta=\sqrt{\Delta}$
    and $c_0=(|a_0|+\delta)/2$ we have
    $$(q,r)=([\al,c_0,0,a_1+(b_1/2-a_1-\al/2)(1-|a_0|/\delta)],[0,-c_0^2/b_0,0,(\al+2a_1-b_1)|a_0|/\delta])\;.$$
    \item The \emph{sub-exponential} type $F=0$, $|E|=1$, $D>0$. We have
      two subtypes:
      \begin{enumerate}\item $E=-1$ (type $D^-$): this occurs when
        $\be=2\al+1$ and $b_0>0$, and we have
        $$(q,r)=([\al+1/2,b_0^{1/2},a_0^2/b_0,(b_1-\al-1/2)/2],[0,-1,4a_0^2/b_0,0])\;.$$
      \item $E=+1$ (type $D^+$): this occurs when $\be=2\al$,
        $\Delta=a_0^2+4b_0=0$, and $\al+2a_1-b_1>0$, and we have
        $$(q,r)=([\al,a_0/2,4(\al+2a_1-b_1),(b_1-\al+1/2)/2],[0,1,16(\al+2a_1-b_1),0])\;.$$
      \end{enumerate}
    \item The \emph{polynomial} convergence type $F=0$, $|E|=1$, $D=0$, $P>0$.
      We have two subtypes:
      \begin{enumerate}\item $E=-1$ (type $P^-$): this occurs when
        $\be=2\al+2$ and $b_0>0$, and we have
        $$(q,r)=([\al+1,b_0^{1/2},0,(b_1+|a_0|b_0^{-1/2}-\al-1)/2],[0,-1,0,|a_0|/b_0^{1/2}])\;.$$
      \item $E=1$ (type $P^+$): this occurs when $\be=2\al$,
        $\Delta=a_0^2+4b_0=0$, $\al+2a_1-b_1=0$, and $B>0$ with
        $B=(2\al+2a_1-1)^2-2\al(\al-1)-4(b_2-2a_2)$, and we have
        $$(q,r)=([\al,|a_0|/2,0,a_1+(1+\sqrt{B})/2],[0,1,0,\sqrt{B}])\;.$$
      \end{enumerate}
    \item The \emph{logarithmic type} (type $L$): this occurs when
      $\be=2\al$, $\Delta=a_0^2+4b_0=0$, $\al+2a_1-b_1=0$, and $B=0$, where
      $B$ is as above, and we have
      $$q(n)=[\al,|a_0|/2,0,a_1+1/2]\log(n)^K\text{\quad and\quad}S-\dfrac{p(n)}{q(n)}\sim\dfrac{C}{\log(n)^{2K-1}}$$
      for some nonzero constant $C$ and $K=0$ or $K=1$, whose value cannot
      be determined only from the asymptotic behavior of $a(n)$ and $b(n)$,
      so convergence only when $K=1$, otherwise logarithmic divergence.
      \item The \emph{non-convergent polynomial} type (type $P^{nc}$):
        the odd and even convergents both converge to usually distinct
        values $S_o$ and $S_e$ respectively. This occurs when $\be\ge2\al+3$,
        and we have
        $$q(n)=[\be/2-\al,(|b_0|/a_0^2)^{1/2},0,b_1/2-a_1-\be/4]$$
        and $$S_e-\dfrac{p(2n)}{q(2n)}\sim\dfrac{C_e}{n^{\be/2-\al-1}}\text{\quad and\quad}S_o-\dfrac{p(2n+1)}{q(2n+1)}\sim\dfrac{C_o}{n^{\be/2-\al-1}}\;.$$
        The contracted continued fraction (taking either even or odd
        convergents) is then of type $P^+$ with $P=\be-2\al-2$.
  \end{enumerate}
\end{theorem}

\begin{remarks}\begin{enumerate}
  \item In the cases not covered by the above theorem, it is reasonable
    to expect that the CF does not converge, but this may not be the case.
  \item As emphasized above, these convergence types are valid generically,
    so in very special cases may be incorrect. Additional errors may come
    from the fact that the proof below is partly heuristic.
  \item An important remark also linked to genericity concerns the square roots
    which must be taken in the following cases: for the exponential type,
    $\delta=\sqrt{\Delta}$, for type $P^-$ we must compute $\sqrt{b_0}$, and
    for type $P^+$ we must compute $\sqrt{B}$.
    In these cases it is understood that we must choose the positive square
    root, which corresponds to the generic situation. However when the CF
    depends on one or more variables, we must be careful to take the correct
    square root, i.e., to add suitable absolute values: $\sqrt{z^2}=|z|$.
  \item A trivial remark: do not confuse the notion of polynomial
    \emph{convergence} type as above, with the notion of CFs of polynomial
    type, which means that $a(n)$ and $b(n)$ are polynomials for large $n$.
    I will try to always include the words ``convergence type'' when necessary.
  \end{enumerate}
\end{remarks}

\section{Sketch of Proof}

The proof is given in \cite{Bel-Coh}, apart from the logarithmic case.
To keep this book self-contained we sketch it here, but advise the reader
to skip the tedious proof. Again emphasizing the above warning, the proof
is correct, but relies on some heuristic assumptions.

Setting $S(n)=p(n)/q(n)$, we have seen that
$$S(n+1)-S(n)=(-1)^n\dfrac{\prod_{0\le j\le n}b(j)}{q(n)q(n+1)}\;.$$
To prove the convergence estimates given in the theorem, we thus need
three results: find an asymptotic estimate for $\prod_{0\le j\le n}b(j)$,
then for $q(n)$, and finally for $S-S(n)$ using the estimate for
$S(n+1)-S(n)$.

The first and third tasks are immediate and are summarized in the following
two lemmas, whose proofs are left to the reader:

\begin{lemma} Using the FEDP notation explained above we have
  $$\prod_{0\le j\le n}b(j)=[\be,b_0,0,b_1]\;.$$\end{lemma}

\begin{lemma}\label{lem:sns} Assume that $S(n+1)-S(n)=1/[F,E,D,P]$, and
  (to insure convergence) that either $F>0$, or $F=0$ and $|E|>1$, or
  $F=0$, $|E|=1$ and $D>0$, or $F=0$, $|E|=1$, $D=0$ and $P>1$, or
  $F=0$, $E=-1$, $D=0$ and $P>0$. If $F=0$, $E=-1$, $D=0$, and $P>0$ we
  assume in addition that $S(2n+2)-S(2n)=1/[F',E',D',P']$ with $E'>0$.
  Then $S-S(n)=1/[F,E,D,P']$, with $P'=P$ unless $F=0$, $E=1$, and either
  $D>0$, in which case $P'=P-1/2$, or $D=0$, in which case $P'=P-1$.\end{lemma}

The more difficult task is to find an asymptotic estimate for $q(n)$, which
we now proceed to achieve.

First an immediate reduction: using $r(n)=1/a(n)$ for $n$ sufficiently
large replaces $(a(n),b(n))$ by the equivalent continued fraction
$(1,b(n)/(a(n)a(n+1)))$, and we have the asymptotic expansion
$$\dfrac{b(n)}{a(n)a(n+1)}=n^{\be-2\al}\dfrac{b_0}{a_0^2}\left(1+\dfrac{B_1}{n}+\dfrac{B_2}{n^2}+\cdots\right)\;,$$
with
\begin{align*}
  B_1&=b_1-2a_1-\al\text{\quad and}\\
  B_2&=b_2-(\al+2a_1)b_1-2a_2+3a_1^2+(2\al+1)a_1+\al(\al+1)/2\;.
\end{align*}

The reduction to $a(n)=1$ being made, the recursion becomes
$q(n+1)=q(n)+b(n)q(n-1)$, so if we set $v(n)=q(n)/q(n-1)$ we have
$v(n+1)=1+b(n)/v(n)$.
We assume heuristically that $v(n)$ has an expansion of the form
$$v(n)=n^{\ga}c_0(1+c_1/n^{1/2}+c_2/n+\cdots)\;.$$
Once these new coefficients found, we have achieved our second goal, and
more precisely we have the following lemma again left to the reader:

\begin{lemma} Assume that $v(n)=q(n)/q(n-1)$ has the above asymptotic
expansion $v(n)=n^{\ga}c_0(1+c_1/n^{1/2}+c_2/n+\cdots)$.
Then generically we have
\begin{align*}
  q(n)&=[\ga,c_0,4c_1^2,c_2-c_1^2/2]\\
  S(n+1)-S(n)&=1/[2\ga-\be,-c_0^2/b_0,16c_1^2,\ga+2c_2-b_1-c_1^2]\\
  S(n+2)-S(n)&=1/[2\ga-\be,c_0^2/b_0,16c_1^2,2\ga+2c_2-b_1-c_1^2]
\end{align*}
\end{lemma}

It remains to find $\ga$ and the $c_i$. We compute that
\begin{align*}
  v(n+1)&=n^{\ga}c_0(1+e_1/n^{1/2}+e_2/n+\cdots)\\
  b(n)/v(n)&=n^{\be-\ga}(b_0/c_0)(1+d_1/n^{1/2}+d_2/n+\cdots)\;,
\end{align*}
with
$$d_1=-c_1,\ d_2=b_1+c_1^2-c_2,\ e_1=c_1,\ e_2=c_2+\ga\;,$$
and easily computed additional formulas for $d_i$ and $e_i$ for $i>2$,
the recursion being $v(n+1)=1+b(n)/v(n)$.

We reason by successive approximations. First we neglect the terms tending
to $0$, so $v(n)$ is close to a root of $X^2-X-n^{\be}b_0=0$, so generically
$n^{\ga}c_0=(1+\sqrt{1+4n^{\be}b_0})/2$. We thus consider three cases:

\smallskip

{\bf Case $\be<0$} The equation gives $n^{\ga}c_0=1$, in other words $\ga=0$
and $c_0=1$. If $\be=-1$ then $e_1=0$, $e_2=b_0$, hence $c_1=0$ and $c_2=b_0$,
so $v(n)=1+b_0/n+\cdots$. If $\be\le-2$ then $e_1=e_2=c_1=c_2=0$,
so $v(n)=1+\cdots$. This leads to case $F$.

\smallskip

{\bf Case $\be>0$} The equation gives $n^{\ga}c_0=(n^{\be}b_0)^{1/2}$,
so $\ga=\be/2$ and $c_0=b_0^{1/2}$ (so we must assume $b_0>0$ in this case).
Simplifying the recursion gives
$$e_1/n^{1/2}+e_2/n+\cdots=b_0^{-1/2}/n^{\be/2}+d_1/n^{1/2}+d_2/n+\cdots$$

\smallskip

{\bf Subcase $\be=1$}, so $\ga=1/2$. Identification easily shows that
$$v(n)=n^{1/2}b_0^{1/2}(1+(b_0^{-1/2}/2)/n^{1/2}+(b_1/2+1/(8b_0)-1/4)/n+\cdots)$$
This leads to case $D^-$.

\smallskip

{\bf Subcase $\be=2$}, so $\ga=1$. Identification easily shows that
$$v(n)=nb_0^{1/2}(1+(b_1/2+b_0^{-1/2}/2-1/2)/n+\cdots)$$
This implies $S(n+1)-S(n)=1/[0,-1,0,b_0^{1/2}]$ which is not sufficient to
prove convergence, but we also compute that
$S(n+2)-S(n)=1/[0,1,0,1+b_0^{1/2}]$, which implies convergence.

This leads to case $P^-$.

\smallskip

{\bf Subcase $\be\ge3$}. Identification easily shows that
$$v(n)=n^{\be/2}b_0^{1/2}(1+(b_1/2-\be/4)/n+\cdots)\;,$$
which implies that $S(n+1)-S(n)=1/[0,-1,0,0]$, so the CF diverges. However
we also compute that $S(n+2)-S(n)=1/[0,1,0,\be/2]$, so separately the
even and odd convergents converge to a limit.

This leads to case $P^{nc}$.

\smallskip

{\bf Case $\be=0$}. From the second degree equation above, we deduce that
$\ga=0$ and $c_0=(1+\sqrt{4b_0+1})/2$, so we assume that $b_0\ge-1/4$.
We again distinguish subcases.

\smallskip

{\bf Subcase $\be=0$, $b_0>-1/4$}. Identification easily shows that
$$v(n)=\dfrac{1+\sqrt{4b_0+1}}{2}(1+b_1(1-1/\sqrt{4b_0+1})/(2n)+\cdots)\;.$$
This leads to case $E$.

\smallskip

{\bf Subcase $\be=0$, $b_0=-1/4$, and $b_1<0$}, so $c_0=1/2$. Identification
(here of terms up to $1/n^{3/2}$) easily shows that
$$v(n)=\dfrac{1}{2}(1+\sqrt{-b_1}/n^{1/2}+1/(4n)+\cdots)$$
This leads to case $D^+$.

\smallskip

{\bf Subcase $\be=0$, $b_0=-1/4$, $b_1=0$, and $b_2<1/4$}, so $c_0=1/2$ and
$c_1=0$. Identification (here of terms up to $1/n^2$) easily shows that if
$b_2\le 1/4$ we have
$$v(n)=\dfrac{1}{2}(1+(1+\sqrt{1-4b_2})/(2n)+\cdots)$$
This implies that $S(n+1)-S(n)=1/[0,1,0,1+\sqrt{1-4b_2}]$, so convergence
at a known speed if $b_2<1/4$.

This leads to case $P^+$.

\smallskip

{\bf Subcase $\be=0$, $b_0=-1/4$, $b_1=0$, and $b_2=1/4$}. This is the
case which was omitted from \cite{Bel-Coh}, more precisely where loc.~cit.~
states that the CF has a logarithmic divergence. This is in fact not always
true, as we will now see. This will lead to case $L$.

By assumption we have $b(n)=-1/4-1/(16n^2)+O(1/n^3)$, and we have just seen
that $v(n)=(1/2)(1+1/(2n)+\cdots)$, where $\cdots$ means as always lower
order terms. We must be more precise about the ``$\cdots$'', so we set
$v(n)=(1/2)(1+1/(2n)+\phi(n))$, with $\phi(n)=o(1/n)$. Identification
(here of terms up to $1/n^3$) shows that
$$\phi(n+1)=\phi(n)-\phi(n)/n-\phi(n)^2+O(1/n^3)\;.$$
Again reasoning heuristically, we assume that $\phi'(n)=o(1/n^2)$ and
$\phi''(n)=o(1/n^3)$, so when $n$ is large we have
$\phi(n+1)=\phi(n)+\phi'(n)+O(1/n^3)$, so $\phi(n)$ should be close to a
solution of the differential equation $y'=-y/n-y^2$. The general solution to
this is $y=1/(n(\log(n)-C))$ for any $C$, plus the special solution $y=0$, and
since we are reasoning up to $O(1/n^2)$, this means that
$\phi(n)=K/(n\log(n))+O(1/n^2)$ with $K=1$ or $K=0$, in other words that
$$v(n)=\dfrac{1}{2}\left(1+\dfrac{1}{2n}+\dfrac{K}{n\log(n)}+\cdots\right)$$
Note the important fact that the numerator $K$ of the term in $1/(n\log(n))$
is either $0$ or $1$, not any constant. Now recall that $v(n)=q(n)/q(n-1)$ and
$q(0)=1$, so $q(n)=\prod_{1\le m\le n}v(m)$, hence taking logarithms we deduce
that $q(n)\sim C2^{-n}n^{1/2}\log(n)^K$ for some nonzero constant $C$. Since
$\prod_{1\le j\le n}b(j)=[0,-1/4,0,0]$, we deduce that
$S(n+1)-S(n)\sim C'/n\log(n)^{2K}$ for some nonzero constant $C'$. We thus
have two different convergence types: either $K=1$, in which case the
continued fraction converges to $S$ and $S-S(n)\sim-C'/\log(n)$, or $K=0$
in which case $S(n)\sim C'\log(n)$, and the continued fraction diverges
logarithmically (to $\pm\infty$ depending on the sign of $C'$). See below
for examples of both cases.

Note that in all other cases, the given formulas is generically correct,
i.e., it is only for special \emph{initial values} of $a(n)$ and $b(n)$
that the formulas must be modified, while in the present logarithmic case,
both $K=0$ and $K=1$ are generically possible.

\smallskip

Putting everything together and coming back to the general case of
$(a(n),b(n))$ proves the result, in a heuristic manner.

\medskip

\section{Additional Remarks}\label{remlog}

\medskip

\subsection{Continued Fractions with Period $T>1$}

\smallskip

In the case of CFs with period $T>1$, the convergence can be given only
for each congruence class modulo $T$. Consider by far the most frequent
case $T=2$. Using the above theorem on the contracted CF we can give
asymptotic estimates for $S-p(2n)/q(2n)$ and $S-p(2n+1)/q(2n+1)$ separately.
Usually, only the constant will change between the two, so we will only
consider $S-p(2n)/q(2n)$. The following lemma is trivial:

\begin{lemma} Assume that $S-p(2n)/q(2n)=[F,E,D,P]$. Then for $n$ even
  we have
  $$S-p(n)/q(n)=[F/2,\sqrt{E}/2^{F/2},D/2,P+F/4]\;.$$
\end{lemma}

Thus, when we give the speed of convergence of a period $2$ CF, we use
the formula given by this lemma.

One important thing to note is that one must compute the square root of $E$:
as mentioned above, this means the positive square root, and care must be
taken when $E$ depends on one or more variables. But in addition, in rare
cases $E$ can be negative: in that case it is possible that the CF does
not converge, but if it does and $E=-1$, we indicate this by $P^i$ or $D^i$
instead of $P^-$ and $D^-$. In more detail: if $E=1$, for $n$ sufficiently
large the signs of successive convergents is $(+,+,+,+,\dotsc)$ or
$(-,-,-,-,\dotsc)$; if $E=-1$ it is $(+,-,+,-,\dotsc)$ or $(-,+,-,+,\dotsc)$;
if $E=i$ it is $(+,+,-,-)$ up to cyclic permutation.

\medskip

\subsection{The Logarithmic Case}

\smallskip

Consider the continued fractions
$$C_1=1-\dfrac{1}{3-\dfrac{4}{5-\dfrac{9}{7-\dfrac{16}{9-\ddots}}}}$$
and
$$C_2=1-\dfrac{1}{2-\dfrac{4}{5-\dfrac{9}{7-\dfrac{16}{9-\ddots}}}}\;.$$
The only difference between the two is that $a(1)=3$ for $C_1$, while
$a(1)=2$ for $C_2$. Both have $a(n)=2n+1$ for $n\ne1$ and $b(n)=-(n+1)^2$ for
all $n$, and we easily check that we are in the logarithmic case. In fact,
as we have seen, the logarithmic case is equivalent to
$$\dfrac{b(n)}{a(n)a(n+1)}=-\dfrac{1}{4}-\dfrac{1}{16n^2}+O(1/n^3)\;.$$
Now it is easily shown that for $C_1$ we have $p(n)/q(n)=1/H(n+1)$, where
$H(m)=\sum_{1\le j\le m}1/j$ is the harmonic sum, so it converges to
$S=0$ and $S-p(n)/q(n)\sim-1/\log(n)$. On the other hand, it is also
immediate to see that for $C_2$ we have $p(n)/q(n)=-H(n+1)+2$, so it
diverges and $p(n)/q(n)\sim-\log(n)$. See \ref{6.2.7} for the general case
$a(1)=z$. In the literature on this case, this specific divergence is
by extension considered as convergent, but I do not agree with this
convention.

\subsection{Asymptotic Expansions}

Apart from the logarithmic case, the above theorem gives results of the type
$q(n)\sim n!^FE^ne^{(Dn)^{1/2}}n^PC$ or $S-p(n)/q(n)\sim C/(n!^FE^ne^{(Dn)^{1/2}}n^P)$, and an algorithm to find $F$, $E$, $D$, and $P$. Finding the
constant $C$ is considerably more difficult and can only be done on a case
by case basis, by first computing it numerically to sufficient accuracy
(using the {\tt cfasymp} command, see the next chapter) and then recognizing
it using a linear dependence program such as {\tt lindep}.

But we may want even more, i.e., if it exists, an expression of
$q(n)/[FEDP]$ or of $(S-p(n)/q(n))FEDP$ as an asymptotic expansion
$C(1+a_1/n+a_2/n^2+\cdots)$, or
$C(1+a_{1/2}/n^{1/2}+a_1/n+a_{3/2}/n^{3/2}+\cdots)$
in the case $D\ne0$. This is in fact much \emph{easier} than computing 
$C$ and can be done without any guesswork, \emph{except} in the case
$P^+$ (i.e., $F=E=D=0$) with $P$ integral. This is explained in detail in
Appendix \ref{chap:asymp}.

\subsection{Irregular Continued Fractions}

Although we only mention CFs with a given period $T$ (almost always $T=1$
or $T=2$), it is worth mentioning that there exist natural CFs whose
coefficients are completely irregular, but nonetheless apparently converge, and
in fact sub-exponentially with type $D^{\pm}$ (although to the author's
knowledge, no proof is known): these are CFs coming from
the transformation of asymptotic expansions of inverse Mellin transforms.
For instance the CF for the inverse Mellin transform of $\G(s)^3$ (which is
easy to create explicitly, for instance using the {\tt Pari/GP} function
{\tt gammamellininvasymp}) seems to converge subexponentially, but as far
as I am aware one does not even know how to prove that it converges at all.
I refer to Chapter 8 of \cite{Bel-Coh} for a detailed description and
discussion of this subject.

\chapter{Numerical Evaluations}\label{chap:numerics}

\section{Introduction}\label{sec:numintro}

Now that we know the convergence behavior of a continued fraction, by far
the most important next step is to numerically compute its limit, and the
method used is of course entirely dependent on its speed of convergence.
We summarize these methods and then give details in the next section.
Write $S-p(n)/q(n)=1/[F,E,D,P]$.

\begin{enumerate}\item Convergence type $F$ and $E$: this is very rapid
  convergence, so we simply compute a sufficiently large $N$
  such that $N!^FE^N>2^B$, where $B$ is the desired number of bits
  (with possibly a small additional correction if $P$ is negative),
  and $p(N)/q(N)$ will be our numerical approximation to the limit.
\item Convergence type $D^{\pm}$: although the convergence is slower
  and we could try to be clever, we simply compute a sufficiently large
  $N$ such that $e^{(DN)^{1/2}}>2^B$, in other words $N>(B\log(2))^2/D$.
  If the desired accuracy is not too large, this is sufficient.
\item Convergence type $P^{\pm}$: in this case we have
  $|S-p(n)/q(n)|\sim C/n^P$, so it is essential to use some kind of
  \emph{extrapolation} method, which we will describe in detail in the
  next section.
\item Logarithmic convergence type $L$: in this case
  $S-p(n)/q(n)\sim C/\log(n)$, so it is absolutely impossible to use any
  kind of extrapolation method, which would require computing $p(n)/q(n)$
  for exponentially large values of $n$. There is, however, a possible
  solution which often works: we use \emph{Bauer--Muir} acceleration
  (see Chapter \ref{chap:bauer}) to transform our CF into one which converges
  to the \emph{same} limit and has convergence type $P^{\pm}$. When this
  is possible we win, otherwise we give up.
\end{enumerate}

\section{The Polynomial Type: Extrapolation Methods}

In this section we assume that we are in convergence type $P^{\pm}$, so
that we want to extrapolate. Note that Chapter 2 of \cite{Bel-Coh} is
entirely devoted to extrapolation methods, however since our present situation
is simpler, it is not difficult to choose a reasonably simple and efficient
method.

Indeed, we know that $S-p(n)/q(n)\sim\eps^n C/n^P$ with $\eps=\pm1$, but
if we refine the proof given in Chapter \ref{chap:speed} it is possible to
show that when $P$ is not integral there is an asymptotic expansion
$$S-\dfrac{p(n)}{q(n)}=\dfrac{\eps^n}{n^P}(C_0+C_1/n+C_2/n^2+\cdots)$$
for suitable constants $C_i$ with $C_0=C\ne0$.

When $P$ is integral, we may also have an asymptotic expansion of the above
shape, but it may also involve $\log(n)$, in which case the extrapolation
methods that we use may only give a few correct decimals. In that case,
it may be useful to first apply one or two Bauer--Muir transformations
before extrapolating.

Note that in the case of convergence type $D^{\eps}$ we would have instead
$$S-p(n)/q(n)=\eps^n/e^{(Dn)^{1/2}}(C_0+C_{1/2}/n^{1/2}+C_1/n+C_{3/2}/n^{3/2}+\cdots)\;,$$
and the asymptotic expansion in powers of $1/n^{1/2}$ would slightly complicate
matters, although it can also be dealt with, see again loc.~cit.

\smallskip

In the case that we have an asymptotic expansion in powers of $1/n$ we can
use a simple method based on Lagrange interpolation. The only thing to be
taken care of is the fact that the exponent $P$ may not be integral. We
extract from \cite{Bel-Coh} the following algorithm adapted to our specific
needs, where in particular is explained the choice of the magic constants
$1/5$, $10$, and $4/3$:

\begin{proposition}\label{prop:accel} Let $B$ be some positive integer, and let
  $p(n)/q(n)$ be the convergents of a continued fraction of polynomial
  convergence type $P^{\pm}$, and assume that $P$ is an \emph{integer}.
  We choose $N=1+\lceil B/5\rceil$, and for $n\le N$ we define
  $w_{N,n}=(-1)^{N-n}\binom{N}{n}n^N$.
  The expression
  $$\dfrac{1}{N!}\sum_{1\le n\le N}(-1)^{N-n}\binom{N}{n}n^N\dfrac{p(10n)}{q(10n)}$$
  is an approximation to $B$ bits to the limit of the continued fraction,
  the computation being done with internal bit accuracy $4B/3$.
\end{proposition}

Of course there are more intelligent ways to compute the above expression
than simply as it is written.

When $P$ is not an integer, there are similar but more complicated expressions
which allow the computation of the limit, see \cite{Bel-Coh} once again.

Once the limit $S$ obtained, one can use similar extrapolation methods
to find the constant $C$ such that
$S-p(n)/q(n)\sim C/(n!^FE^ne^{\sqrt{Dn}}n^P)$. 

\medskip

Let us give an example. Applying Euler's transformation Proposition
\ref{prop:euler} to the series $\pi/4=1-1/3+1/5-1/7+\cdots$ we obtain
the continued fraction
$$\pi=4/(1+1/(2+9/(2+25/(2+49/(2+81/(2+\cdots))))))\;,$$
which is termwise equal to the series hence converges very slowly, like
$(-1)^n/n$, convergence speed which can also be obtained (apart from the
constant) from Theorem \ref{thm:speed}. Thus, to obtain	10 decimals, whether
from the series, or from the CF, would require $10^{10}$ terms.
Using the above method, to obtain 115 decimal digits we
work with 155 digits of internal accuracy and use only 78 terms up to $n=780$,
which shows the amazing efficiency of this method. As an aside, note that
one of the fastest ways to evaluate the \emph{series} (as opposed to the CF)
is to use the {\tt sumalt} function of {\tt Pari/GP} which is a general
method for summing regular alternating series.

\section{The {\tt Pari/GP} Continued Fractions Package}

The rest of this chapter can be skipped on first reading especially if you are
not familiar with {\tt Pari/GP}, but the commands explained in this chapter
have been essential in the construction of the encyclopedic dictionary
which constitutes the rest of this book.

\smallskip

For now, the {\tt Pari/GP} continued fractions package is not part of
the main distribution, but can be downloaded using {\tt GIT} from the
{\tt Pari/GP} branch {\tt henri-ellCF3} (read the explanations for
obtaining such branches on the {\tt Pari/GP} website).

For now I will describe the main available basic commands, and refer to
later chapters for more advanced commands such as {\tt cfbauer}
and {\tt cfapery}. Also, I will not give all the possible options but only
the most useful ones.

\medskip

\subsection{Basic Formats}

\medskip

Let $[A,B]=((a(n))_{n\ge0},(b(n))_{n\ge0})$ be a continued fraction of
polynomial type (not necessarily of polynomial convergence type). If
$A$ (similarly $B$) has period $1$, i.e., if $a(n)$ is a polynomial for
$n\ge n_0$ for some $n_0$, we will represent $A$ as the vector
$$A=[a(0),a(1),\dotsc,a(n_0-1),a(n)]\;,$$ where for $i\le n_0-1$ the $a(i)$
do not depend on $n$, while $a(n)$ is a polynomial (exceptionally also
a rational function) in the reserved variable $n$. Note that in the
package the variable $n$ is \emph{reserved} for this purpose and cannot be
used elsewhere. It would perhaps have been more logical to use the variable
$x$ instead, but the variable $n$ is almost universally used in texts so
I have chosen to use it.

Note also that the representation is not unique, for instance
we also have $A=[a(0),a(1),\dotsc,a(n_0-1),a(n_0),a(n)]$.

For instance, the continued fraction that we have seen above
$$\pi=4/(1+1/(2+9/(2+25/(2+49/(2+81/(2+\cdots))))))$$
is represented by
\begin{verbatim}
[[0,1,2],[4,(2*n-1)^2]]
\end{verbatim}

\smallskip

If $A$ has period $2$, i.e., if $a(2n)$ and $a(2n+1)$ are possibly different
polynomials for $n$ large, we will represent $A$ as a vector of $2$-component
vectors $$A=[[a(0),a(1)],\dotsc,[a(2m_0-2),a(2m_0-1)],[a(2n),a(2n+1)]]\;.$$

For instance, the continued fraction
$$\pi=4-4/(3+9/(5+4/(7+25/(9+16/(11+49/(13+36/(15+\cdots)))))))$$
is represented by
\begin{verbatim}
[[4,2*n+1],[[-4,9],[(2*n)^2,(2*n+3)^2]]]
\end{verbatim}

Note that in this case \emph{all} the components must be two-component
vectors. For instance the above could \emph{not} be written as
\begin{verbatim}
[[4,2*n+1],[-4,[(2*n+1)^2,(2*n)^2]]]
\end{verbatim}
This is a restriction of the package which may be lifted in the future.

\smallskip

Continued fractions of period strictly larger than $2$ are not supported
by the package.

\medskip

A number of special \emph{output} formats (which therefore cannot be used
as inputs to other functions) also exist: see in particular below the output of
{\tt cffromquad} and {\tt cffromser}.

\medskip

\subsection{Creators and Asymptotics}

\medskip

We now list (with examples) the {\tt GP} commands allowing us to \emph{create}
continued fractions, usually infinite, but sometimes also finite. A useful
command to \emph{visualize} continued fractions is {\tt cftochar}, as we will
see below in all the examples.

Note: since the variable {\tt n} plays such a special role, it is useful
(although not strictly necessary) to write ``{\tt n='n;}'' at the beginning of
each {\tt GP} session.

\begin{enumerate}\item As mentioned above, the most common way is using
vectors. Example with period $T=1$:

\begin{verbatim}
? ab=[[1,n-z],[z,z*n]];cftochar(ab,5)
% = "1+z/(-z+1+z/(-z+2+2*z/(-z+3+3*z/(-z+4+4*z/(-z+5)))))"
? log(cflimit(subst(ab,z,2)))
% = 2.0000000000000000000000000000000000000000000000
? cftype(ab)
% = [[1, 1, 0, 0], [1, -1/z, 0, 1], 2]
\end{verbatim}

The {\tt cflimit} command computes the limit, using essentially Proposition
\ref{prop:accel} above (in the present example the limit is equal to
$\exp(z)$). The {\tt cftype} command gives the speed of convergence using
Theorem \ref{thm:speed}: the first component says
that $q(n)\sim C\cdot n!$ for some nonzero constant $C$, the second that
$\exp(z)-p(n)/q(n)\sim C'/(n!(-1/z)^n\cdot n)\sim(-z)^n/(n+1)!$,
so that we are in
the factorial type $F^1$, which for now is indicated by type number $2$.
Note that since equivalent continued fractions can give the same $p(n)/q(n)$,
the really important quantity is $\exp(z)-p(n)/q(n)$, the asymptotics of
$q(n)$ itself being important only for arithmetic applications such as
irrationality measures.

We can even hope to compute the constants $C$ and $C'$ using the command
{\tt cfasymp}: note that {\tt cflimit} and {\tt cfasymp} cannot include
unspecified variables such as $z$, so we must always use {\tt subst}.
For instance, continuing the above example:

\begin{verbatim}
 ? cfasymp(subst(ab,z,2))
% = [7.38..., [1, 1, 0, 0, 0.13533...], 
              [1, -1/2, 0, 1, 109.1963...], 2]
 ? cfasymp(subst(ab,z,3))
% = [20.08..., [1, 1, 0, 0, 0.04978...], 
               [1, -1/3, 0, 1, 1210.28...], 2]
\end{verbatim}
We easily recognize that $0.13533...=\exp(-2)$, that
$109.1963...=2\cdot\exp(4)$, that $0.04978...=\exp(-3)$, and
$1210.28...=3\cdot\exp(6)$. It is thus quite plausible (and indeed true)
that $C=\exp(-z)$ and $C'=z\exp(2z)$.

Remark on {\tt cftype}, {\tt cflimit}, and {\tt cfasymp}: these three functions
play complementary roles. First, {\tt cftype} is the only one which
accepts unspecified variables such as $z$ as above. However, as already
mentioned, it may give an incorrect answer since it often needs to compute
square roots, whose sign may depend on $z$, so care should be
used in interpreting the answer.

Second, {\tt cflimit} is quite robust, while {\tt cfasymp} which also gives
the convergence constants $C$ and $C'$ may fail.

\item Example of period $T=2$:
\begin{verbatim}
? ab=[[0,1],[[1,1/2],[n,n+1]/(4*n+2)]];
? cftochar(ab)
%="1/(1+1/2/(1+1/6/(1+1/3/(1+1/5/(1+3/10/1)))))"
? exp(cflimit(ab))
%=2.0000000000000000000000000000000000000000000
? cftype(ab)
%=[0, 1.207..., 0, 0], [0, -5.828..., 0, 0], 6]
? cfasymp(ab)
%=[0.693..., [0, 1.207..., 0, 0, 1.015...],
             [0, -5.828..., 0, 0, 1.078...], 6]
\end{verbatim}

We recognize $1.207...=(1+\sqrt{2})/2$ and $-5.828...=-(3+2\sqrt{2})$.
The limit is thus $\log(2)$, we are in exponential type $E$, here indicated
by the index $6$, and $q(n)\sim C\cdot ((1+\sqrt{2})/2)^n$ and
$$\log(2)-\dfrac{p(n)}{q(n)}\sim (-1)^n\dfrac{C'}{(3+2\sqrt{2})^n}=(-1)^n\dfrac{C'}{(1+\sqrt{2})^{2n}}\;,$$
and we recognize that $C=1.015...=(1+\sqrt{2})/2^{5/4}$ and
$C'=1.078...=2\pi/(3+2\sqrt{2})$.

Note on constant recognition: this is of course done with the {\tt lindep}
command. In most cases the constants $C$, $C'$, etc... are obtained as
\emph{products}, so we search for linear combinations involving logarithms
of the constants involved, with in addition $\log(\pi)$ and $\log(p)$ for
small primes $p$. A typical command for recognizing the constants in
the above example would be

\centerline{\tt lindep([-log(C),log(1+sqrt(2)),log(Pi),log(2),log(3)])}

It is also sometimes useful to add constants such as $\log(\G(1/3))$ or
$\log(\G(1/4))$.

\item Completely different ways of creating continued fractions in {\tt GP}
are by using the command {\tt cffromser(f,g)}, which can have several different
meanings. The two most important ones are the following:
\begin{enumerate}\item If $f$ and $g$ are both closures, return the continued
  fraction corresponding to Euler's transformation of the series
  $\sum_{j\ge1}g!(j-1)/f!(j)$ as in \ref{prop:euler}. For instance:
\begin{verbatim}
? ab=cffromser(n->n^2,n->n^2)
% = [[0, 2*n^2 - 2*n + 1], [1, -n^4]]
? cftochar(ab)
% = "1/(1-1/(5-16/(13-81/(25-256/(41-625/(61))))))"
? cflimit(ab)-Pi^2/6
% = 0.E-38
? cftype(ab)
% = [[2, 1, 0, 0], [0, 1, 0, 1], 8]
\end{verbatim}
  Of course, the speed of convergence has not been improved since the continued
  fraction is identical to the series.
\item The second important meaning of {\tt cffromser} is when $f$ is not a
  closure, but a formal power series, in which case $g$ must be a small
  nonnegative integer (or omitted, same as $0$) considered as the degree of the
  polynomial $a(n)$. The output is a \emph{finite} continued
  fraction in a special format. We simply give two
  examples, referring to the documentation for details:
\begin{verbatim}
? V2=cffromser(exp(x^2)+O(x^11))
% = [Vecsmall([0, 2]), [1, -1, 1/2, -1/6, 1/6, -1/10]]
? V3=cffromser(1/(x^6+x^2+1))
*** cffromser: impossible inverse in gdiv: 0.
? V3=cffromser(1/(x^6+x^2+1),1)
% = [Vecsmall([0, 2, 4, 4]), [1, 1, -1, 1]]
\end{verbatim}
\item The last way in which one can create continued fractions in {\tt GP}
  is by using the command {\tt cffromquad}, which creates the periodic
  continued fraction corresponding to a real quadratic irrational. Since
  we will not use this in the present book, we refer to the documentation
  for details.
\end{enumerate}
\end{enumerate}

\medskip

\subsection{Other Tools}

\medskip

We have already mentioned the fundamental {\tt cflimit} command which uses
generalizations of Proposition \ref{prop:accel} to give a numerical
approximation to the limit of a continued fraction. Note a limitation of
the present version: the above proposition is useful when the exponent
$P$ of polynomial convergence is integral, and we have mentioned that
generalizations exist, but in fact only when $P$ is a rational number
with denominator $d\le 4$, say. Otherwise the package will try to guess a
numerical limit, but expect it to have only a few if any correct digits.

We have also mentioned {\tt cfasymp} which allows the computation of the
constants entering in the asymptotics of $q(n)$ and of $S-p(n)/q(n)$, and
the same limitations apply, even more strongly.

\smallskip

We have mentioned {\tt cftochar} to visualize a CF, but in addition
all the CFs of this book have been printed using the function {\tt cftotex}
which gives a correct \TeX\ output and attempts to remove some (but not
always all) unnecessary parentheses. Example:

\begin{verbatim}
? A=[[0,(2+z)*(2*n-1)],[2*z,-n^2*z^2]];
? cftochar(A)
% = "2*z/(z+2-z^2/(3*z+6-4*z^2/(5*z+10-9*z^2/
                  (7*z+14-16*z^2/(9*z+18-25*z^2/(11*z+22))))))"
? print(cftotex(A))
\dfrac{2z}{z+2-\dfrac{z^2}{3z+6-\dfrac{4z^2}{5z+10-\dfrac{9z^2}{7z+14
              -\dfrac{16z^2}{9z+18-\dfrac{25z^2}{11z+22-\ddots}}}}}}
\end{verbatim}
(output split for clarity).

\chapter{Bauer--Muir--Ap\'ery Acceleration}\label{chap:bauer}

\section{Introduction}

Many continued fractions occurring in the literature and in the present
book have a polynomial convergence type $P^{\pm}$ as mentioned above.
This is for instance the case for continued fractions obtained by trivial
transformations of slowly convergent series (using Euler's transformation,
i.e., the {\tt cffromser} command). These can be accelerated using the
so-called \emph{Bauer--Muir} method (see for instance \cite{Bel-Coh} for
details), and the crucial contribution of R.~Ap\'ery is to have noticed that
in some cases (for him, the ``trivial'' CFs for $\z(2)$ and $\z(3)$ obtained
by applying Euler's transformation of series) this acceleration method can
be iterated without losing the simplicity of the coefficients, and combined
with a diagonal process, usually leads to an \emph{exponentially convergent}
continued fraction (type $E$), in practice almost always with $|E|$ an integral
power of $(1+\sqrt{2})^2$ or of $((1+\sqrt{5})/2)^5$. This new continued
fraction will have \emph{period 2} in the above sense, meaning that
$(a(2n),b(2n))$ and $(a(2n+1),b(2n+1))$ are given by different formulas.

What is remarkable and may not have been noticed before is that a large
proportion of the CFs which can be found in the literature (and in the
present dictionary) can be accelerated in this way. In every case that it
does and gives reasonably simple coefficients we have given the corresponding
new continued fraction, which usually does not appear in the literature.

In rare cases (such as those used by Ap\'ery in his proof of the
irrationality of $\z(2)$ and $\z(3)$), the even (or odd)
\emph{contraction} of this period $2$ continued fraction gives a period $1$
continued fraction with reasonably simple coefficients.

\section{Bauer--Muir Acceleration}

To keep this book self-contained, we first explain the classical Bauer--Muir
acceleration method. Let $(a(n),b(n))$ be a convergent continued fraction, let
$r(n)$ be (for now) an arbitrary sequence, and set
$$d(n)=r(n)(a(n+1)+r(n+1))-b(n)\;.$$
We will assume that $d(n)\ne0$ for all $n$. We define a new continued fraction
$(a'(n),b'(n))$ by setting
\begin{align*}
  a'(0)&=a(0)+r(0),\quad b'(0)=-d(0),\quad a'(1)=a(1)+r(1)\;,\\
  a'(n)&=a(n)+r(n)-r(n-2)d(n-1)/d(n-2)\quad\text{for $n\ge2$, and}\\
  b'(n)&=b(n-1)d(n)/d(n-1)\quad\text{for $n\ge1$}\;.
\end{align*}
A formal computation shows that if $(p'(n),q'(n))$ are the partial quotients
of this new continued fraction, we have
$(p'(n),q'(n))=(p(n),q(n))+r(n)(p(n-1),q(n-1))$ (this is of course how the
formulas for $a'(n)$ and $b'(n)$ were obtained).

Note that if we denote by $\rho(n)=b(n)/(a(n+1)+b(n+1)/(a(n+2)+...))$ the
$n$th \emph{tail} of the continued fraction, we have trivially
$\rho(n)(a(n+1)+\rho(n+1))-b(n)=0$. Thus to accelerate the CF we want
$r(n)$ be be close to $\rho(n)$, and in view of this recursion this means
that we want $d(n)$ to be small (but of course nonzero for the formulas
of the accelerated CF to make sense).

\smallskip

In our situation, we choose continued fractions such that $a(n)$ and $b(n)$
are rational functions or, up to equivalence, polynomials. We thus would
like to have our accelerated continued fraction to be of a similar type,
and the simplest way to achieve this is to ask for $r(n)$ to be a polynomial,
and $d(n)$ to be a nonzero constant, which is of course a rather severe
constraint (it is possible to have a weaker condition, see Section
\ref{sec:apgen} for details).

An easy examination of the method shows the following, where we set
$S(n)=p(n)/q(n)$ and $s_0=\sign(a_0)$:

\begin{enumerate}\item In the case $P^-$, i.e., $\be=2\al+2$ and $b_0>0$,
  we know that $S-S(n)=1/[0,-1,0,|a_0|/b_0^{1/2}]$. One can show that we
  must choose
  $$r(n)=s_0b_0^{1/2}n^{\be/2}(1+(b_1-\be/2-|a_0|/b_0^{1/2})/(2n)+\cdots)\;,$$
  and then $S-S'(n)=1/[0,-1,0,|a_0|/b_0^{1/2}+2k]$ for some $k\ge1$,
  depending on the further coefficients of the expansion of $r(n)$.
\item In the case $P^+$, i.e., $\be=2\al$, $a_0^2+4b_0=0$, $2a_1=b_1-\al$,
  and $B>0$ with $B=(b_1+\al-1)^2-2\al(\al-1)-4(b_2-2a_2)$, we know
  that $S-S(n)=1/[0,1,0,\sqrt{B}]$. One can show that we must choose
  $$r(n)=-(a_0/2)n^{\al}(1+(b_1+1-\al+\sqrt{B})/(2n)+\cdots)\;,$$
  and as in the previous case we gain a factor of $1/n^{2k}$ for some
  $k\ge1$.
\item In case $E$, i.e., $\be=2\al$ and $d=a_0^2+4b_0>0$, we know that
  $S-S(n)=1/[0,-c_0^2/b_0,0,P]$ with $c_0=(|a_0|+\sqrt{d})/2$ and
  $P=(b_1-2a_1-\al)|a_0|/\sqrt{d}$, and we must choose
  $$r(n)=\dfrac{-a_0+s_0\sqrt{d}}{2}n^{\al}(1+((b_1-\al)(|a_0|+\sqrt{d})/2-|a_0|a_1)/(\sqrt{d}n)+\cdots)\;,$$
  and as in the previous case we gain a factor of $1/n^{2k}$ for some
  $k\ge1$. Of course, this is less important since the initial CF
  already converges exponentially, but we will see that it is essential
  in Ap\'ery's method. Note that since in our applications we want to have
  \emph{rational} coefficients we need $d$ to be the square of a rational
  number, or equivalently $E$ rational.
\item In case $D^-$, i.e., $\be=2\al+1$ and $b_0>0$, we know that
  $S-S(n)=1/[0,-1,4a_0^2/b_0,0]$, and we should choose
  $r(n)\sim r_0n^{\al+1/2}$, but since we only want polynomials, we cannot
  treat this case since $\al\in\Z$.
\item In case $D^+$, i.e., $\be=2\al$, $a_0^2+4b_0=0$, and $2a_1>b_1-\al$,
  we know that $S-S(n)=1/[0,1,16(2a_1-b_1+\al),0]$, and we must choose
  $r(n)\sim(-a_0/2)n^{\al}$, but to obtain any acceleration we would need
  $r(n)=(-a_0/2)n^{\al}(1+r_{1/2}/n^{1/2}+\cdots)$, which is not allowed
  since we want to stay with polynomials. Thus, the best we can do is
  $r(n)=(-a_0/2)n^{\al}(1+r_1/n+\cdots)$, and one can show that this does
  not lead to any acceleration, but it does give a different continued
  fraction converging to the same limit.
\end{enumerate}

Note that by iteration of the Bauer--Muir process we may obtain in this
way an infinite family of new continued fractions, and we will see this
in more detail below in the description of Ap\'ery's method.

We give examples of all this in {\tt Pari/GP}. First, an example
of polynomial type.

\begin{verbatim}
? cffromser(n->2*n-1,n->-(2*n-1))
% = [[0, 1, 2], [1, 4*n^2 - 4*n + 1]] /* Standard CF for Pi/4=atan(1) */
? ab = [[0, 1, 2], [4, 4*n^2 - 4*n + 1]]; /* CF for Pi */
? cftype(ab)[2]
% = [0, -1, 0, 1]
? cftochar(ab)
% = "4/(1+1/(2+9/(2+25/(2+49/(2+81/(2))))))"
/* Convergence P^- with P=1 */
? ab1=cfbauer(ab)
% = [[4, 4, 5, 6], [-4, 4, 4*n^2 - 4*n + 1]]
? cftype(ab1)[2]
% = [0, -1, 0, 3]
? cftochar(ab1)
% = "4-4/(4+4/(5+9/(6+25/(6+49/(6+81/(6))))))"
/* "New" CF for Pi, convergence P^- with P=3 */
? 4-4/(4+4/(5+x))
% = (3*x + 19)/(x + 6)
/* equal to 3+1/(6+x), so we simplify ab1 by setting */
? ab1new=[[3,6],[(2*n+1)^2]];
? cftochar(ab1new)
% = "3+1/(6+9/(6+25/(6+49/(6+81/(6+121/(6))))))"
/* "New" CF for Pi, same P=3 */
? ab2=cfbauer(ab1new);cftype(ab2)[2]
% = [0, -1, 0, 5]
/* We can continue as long as we want */
\end{verbatim}

\smallskip

We now give an example of exponential type. Since we want to keep rational
coefficients, in view of the formulas given above we need $d=a_0^2+4b_0$ to
be a square or, equivalently, $E=-((|a_0|+\sqrt{d})/2)^2/b_0$ to be
a rational number.
\begin{verbatim}
? ab=cffromser(n->2*n,n->n)
% = [[0, 3*n - 1], [1, -2*n^2]]
? cftype(ab)[2]
% = [0, 2, 0, 1]
? cftochar(ab)
% = "1/(2-2/(5-8/(8-18/(11-32/(14-50/(17))))))"
/* Standard CF for log(2)=-log(1-1/2),
   convergence E with E=2 and P=1 */
? ab1=cfbauer(ab)
% = [[1/2, 3*n], [1/2, -2*n^2]]
? cftype(ab1)[2]
% = [0, 2, 0, 3]
? cftochar(ab1)
% = "1/2+1/2/(3-2/(6-8/(9-18/(12-32/(15-50/(18))))))"
/* "New" CF for log(2), still E=2 but P=3 */
? ab2=cfbauer(ab1);cftype(ab2)[2]
% = [0, 2, 0, 5]
/* We can continue as long as we want */
\end{verbatim}

\smallskip

In many cases, a given continued fraction has a convergence type $E$ with
$E=E(z)$ depending on some variable $z$ (or several variables). Thus, if $z$
is chosen to be a rational number such that $E(z)$ is also rational, we
can apply Bauer--Muir acceleration. For instance, in several cases $E(z)$
is a simple expression involving $\sqrt{1+z}$, so we can accelerate by
choosing $z=u^2-1$ with $u$ rational.

\smallskip

Finally here is an example of type $D^+$, where there is no acceleration.
We will see below that if we set
$$S=\int_0^\infty \dfrac{e^{-t}}{t+1}\,dt=0.59634\cdots$$
we have the beautiful continued fraction expansion
$$\int_0^\infty \dfrac{e^{-t}}{t+1}\,dt=\dfrac{1}{2-\dfrac{1}{4-\dfrac{4}{6-\dfrac{9}{8-\dfrac{16}{10-\dfrac{25}{12-\ddots}}}}}}$$
with $a(n)=2n$, $b(0)=1$, and $b(n)=-n^2$ for $n\ge1$,
which converges subexponentially with $D=16$, and more precisely
$S-p(n)/q(n)\sim 2\pi\cdot e/e^{4n^{1/2}}$. By the formulas given above
we must choose $r(n)=-n+r_0$, and we obtain $d(n)=-n+r_0^2+r_0$. Thus $d(n)$
cannot be constant, but all is not lost: if we look at the formulas for
Bauer--Muir acceleration, we only need $b(n)d(n+1)/d(n)$ and $r(n)d(n+1)/d(n)$
to be polynomials. Since $d(n)$ is coprime to $d(n+1)$, we must have
$d(n)\mid b(n)$ and $d(n)\mid r(n)$, which is realized if and only if $r_0=0$.
We thus choose $r(n)=-n$, and we find that $a'(n+1)=2n+1$, $b'(n+1)=-n^2-n$,
and a small computation shows that we obtain the new continued fraction
$$\int_0^\infty \dfrac{e^{-t}}{t+1}\,dt=1-\dfrac{1}{3-\dfrac{2}{5-\dfrac{6}{7-\dfrac{12}{9-\dfrac{20}{11-\dfrac{30}{13-\ddots}}}}}}$$
Thus our work was not for nothing, this is really new (note that if we
iterate, we obtain a continued fraction equivalent to the initial one),
and the reader can easily check that this new continued fraction
converges at exactly the same speed as the previous one, i.e.
$S-p(n)/q(n)\sim-2\pi\cdot e/e^{4n^{1/2}}$.

Note that all this is done automatically by using the {\tt cfbauer} command.

\smallskip

\section{The Special Case $d(n)=0$}

\smallskip

In some cases, it is possible to choose $r(n)$ such that
$$d(n)=r(n)(a(n+1)+r(n+1))-b(n)=0\;.$$
This is apparently useless for Bauer--Muir acceleration since it requires
$d(n)\ne0$, but in fact is much \emph{better} because the continued fraction
can then be computed explicitly. Indeed, if we set $f(n)=(-1)^nr!(n)$,
we check that $f(n)$ satisfies the recursion $f(n+1)=a(n+1)f(n)+b(n)f(n-1)$,
so by Corollary \ref{cor:expl} we have
$$\dfrac{p(n)}{q(n)}=a(0)+b(0)\dfrac{S(n)}{1+(a(1)+r(1))S(n)}\text{\quad with\quad}S(n)=\sum_{1\le m\le n}(-1)^m\dfrac{b!(m-1)}{r!(m)r!(m-1)}\;.$$
Now since the quantity $a(n+1)+r(n+1)$ occurs everywhere in Bauer--Muir,
we set here and in the rest of this book
$$R(n)=a(n+1)+r(n+1)\;,$$
so that $d(n)=0$ is equivalent to $R(n)=b(n)/r(n)$, so
$b!(m-1)/r!(m-1)=R!(m-1)/R(0)$, and since $r(0)R(0)=b(0)$ we deduce that
$$\dfrac{p(n)}{q(n)}=a(0)+r(0)\dfrac{S'(n)}{1+S'(n)}\text{\quad with\quad}
S'(n)=\sum_{1\le m\le n}(-1)^m\dfrac{R!(m-1)}{r!(m)}\;.$$
Thus, computing the limit of the continued fraction is equivalent to
computing the sum of the series $S'=\sum_{m\ge1}(-1)^mR!(m-1)/r!(m)$.

\smallskip

We give a complete example of this which also illustrates the phenomenon of
non genericity.
\begin{verbatim}
? C1=[[1,3*n^3-10*n^2+8*n-2],[1,2*n^5*(2*n+1)]];cftype(A)[2]
% = [0, -4, 0, -5/2]
? cflimit(C1)
% = -8930186..... /* a very large number */
? cftochar(C1)
% = "1+1/(-1+6/(-2+320/(13+3402/(62+18432/(163+68750/(334))))))"
\end{verbatim}
  The {\tt cftype} command indicates that {\tt C1} should converge like
  $1/((-4)^nn^{-5/2})$, so exponentially fast. But the {\tt cflimit}
  command seems to indicate that the CF diverges to $\infty$ (in
  fact, as we will now see, the {\tt cftype} result should be instead
  {\tt [0, -1/4, 0, 5/2]}, wildly divergent). Looking at
  the beginning of the CF, this means that the CF
  $6/(-2+320/(13+\cdots))$ converges to $1$. So let us modify {\tt C1}
  to obtain this new CF:
\begin{verbatim}
? C2=[[0,subst(C1[1][2],n,n+1)],[subst(C1[2][2],n,n+1)]];
? cftochar(C2,5)
% = 6/(-2+320/(13+3402/(62+18432/(163+68750/(334))))))"
? cflimit(C2)
% = 1.0000000000000000000000000000000000291
\end{verbatim}
  Now indeed the speed of convergence (which is the same) is in
  $1/((-4)^nn^{-5/2})$. Let us now \emph{prove}
  that the limit is indeed equal to $1$. Solving linear equations to find
  a suitable $r(n)$, we find that $r(n)=(n+1)^3$ works, and with
  $$a(n)=3n^3-n^2-3n-1\text{\quad and\quad}b(n)=2(n+1)^5(2n+3)$$
  $$d(n)=r(n)(a(n+1)+r(n+1))-b(n)=0\;,$$
  including for $n=0$.
  
  Now for $n\ge1$ we have $R(n-1)=a(n)+r(n)=4n^3+2n^2$, so clearly
  $R!(m-1)/r!(m)$ tends to $+\infty$ exponentially fast. It follows
  from the formula for $p(n)/q(n)$ given above that the CF converges
  to $a(0)+r(0)=1$, as desired, and the speed of convergence is governed
  by $(-1)^nR!(n-1)/r!(n)$ which is indeed of the order of $(-4)^n$.

\medskip

\section{Ap\'ery's Method}

\smallskip

Ap\'ery's method can be summarized in two steps: first, do a suitable
iteration of the Bauer--Muir acceleration method, second use a diagonal
procedure on these iterations.

The first step proceeds as follows. Assume that at iteration $l$
we have a continued fraction $(a(n,l),b(n,l))$ and partial quotients
$(p(n,l),q(n,l))$, so that $(a(n,0),b(n,0))=(a(n),b(n))$, and
$(p(n,0),q(n,0))=(p(n),q(n))$. We choose a suitable sequence $r(n,l)$,
and as above we set $d(n,l)=r(n,l)(r(n+1,l)+a(n+1,l))-b(n,l)$.
Thanks to the Bauer--Muir formulas given above, we can then define
a new continued fraction converging to the same limit with
coefficients $(a'(n),b'(n))$ and partial quotients $(p'(n),q'(n))$.
The first important point is to define the $(l+1)$st iteration by
setting
$$(a(n,l+1),b(n,l+1),p(n,l+1),q(n,l+1))=(a'(n+1),b'(n+1),p'(n+1),q'(n+1))\;,$$
the shift of $n$ by $1$ on the right-hand side being essential.
One can then hope to find a general formula for $(a(n,l),b(n,l))$.

Finding $(a'(n+1),b'(n+1))$ in terms of $(a(n),b(n))$ is done automatically
using the command {\tt cfbauer(,1)} which in fact returns the four sequences
{\tt (a,b,r,d)} as rational functions. Let us apply this to the first
example given above, using the command {\tt cftopol} which gives the
polynomials or rational functions $(a(n),b(n),\dotsc)$ for $n$ sufficiently
large, omitting the preperiod:
\begin{verbatim}
? v0=cffromser(n->2*n-1,n->-(2*n-1))
% = [[0, 1, 2], [1, 4*n^2 - 4*n + 1]]
? cftopol(v0)
% = [2, 4*x^2 - 4*x + 1]
? ab=v0;
? for(l=0,4,abrd=cfbauer(ab,1);print(cftopol(abrd));ab=abrd[1..2]);
[6, 4*x^2 - 4*x + 1, 2*x - 3, -4]
[10, 4*x^2 - 4*x + 1, 2*x - 5, -16]
[14, 4*x^2 - 4*x + 1, 2*x - 7, -36]
[18, 4*x^2 - 4*x + 1, 2*x - 9, -64]
[22, 4*x^2 - 4*x + 1, 2*x - 11, -100]
\end{verbatim}
Note that the $l$th row starting at $l=0$ gives
$(a(n,l+1),b(n,l+1),r(n,l),d(n,l))$. We immediately recognize the
pattern, so we set
\begin{verbatim}
? a(n,l)=4*l+2;
? b(n,l)=(2*n-1)^2;
? r(n,l)=2*n-2*l-3;
? d(n,l)=-4*(l+1)^2;
\end{verbatim}
Two remarks: first, the {\tt cftopol} command used to print the results
only gives the polynomial expression for $n$ sufficiently large, so
the formulas for {\tt a(n,l)} etc... just given are valid only for
$n$ sufficiently large. This will be important for the second stage of
Ap\'ery's method which we will now study.

The second remark is that this ``recognition'' of the pattern for the
coefficients is simply polynomial interpolation (in some cases rational
function interpolation), and this is done
automatically with the {\tt cfapery} command.

\smallskip

The second crucial step in Ap\'ery's method is to use a diagonal process.
Ideally, we would like a continued fraction (in other words a three-term
linear recurrence) whose partial quotients are $(p(n,n),q(n,n))$.
Although possible, this in general leads to complicated formulas.
On the contrary, one can show that if we want a continued fraction
of period 2 whose even partial quotients are $(p(n,n),q(n,n))$,
and odd partial quotients are $(p(n,n+1),q(n,n+1))$, the formulas
become extremely simple (and if desired we can obtain the $(p(n,n),q(n,n))$
by contraction). More precisely, we have the following (see
Appendix \ref{chap:apimpl} for details and the proof of a more general
result):

\begin{proposition} Set $R(n,l)=a(n+1,l)+r(n+1,l)$, and for $n$
  sufficiently large, set
\begin{align*}(A(2n),A(2n+1))&=(R(n,n-1),R(n,n))\text{\quad and\quad}\\
  (B(2n),B(2n+1))&=(b(n,n),-d(n+1,n))\;.\end{align*}
Then if we denote by $S$ the limit of the initial continued
fraction, the one defined by $(A(n),B(n))$ will converge to an
expression of the form $(uS+v)/(wS+x)$ for simple $(u,v,w,x)$ corresponding
to a correct choice of the initial values of $A$ and $B$.
\end{proposition}

Note that it should be possible to find these correct initial values,
but I have been too lazy to do so, since in any case as we will see,
they are easy to determine {\it a posteriori\/} using the {\tt fracdep}
command, see below.

Continuing our example, we thus write:
\begin{verbatim}
? R(n,l)=a(n+1,l)+r(n+1,l);
? AB=[[[R(n,n-1),R(n,n)]],[[b(n,n),-d(n+1,n)]]]
% = [[[4*n - 1, 4*n + 1]], [[4*n^2 - 4*n + 1, 4*n^2 + 8*n + 4]]]
? fracdep(cflimit(AB),cflimit(v0))
% = (-x - 1)/x
? cftype(AB)[2]
% = [0, -5.828..., 0, 0]
\end{verbatim}

Note that the command {\tt fracdep(U, V)} finds if it exists a reasonable
size M\"obius transformation such that $V=(aU+b)/(cU+d)$, and returns $0$
otherwise. Note also the essential double bracket {\tt [[...]]} to indicate
period $T=2$.

\smallskip
  
We have thus indeed obtained an exponentially convergent continued fraction
with $E=-(1+\sqrt{2})^2$, which converges to a quantity trivially linked
to $\pi/4$. All of this is done automatically by the {\tt cfapery} command.

Let us finish this example by cleaning up ``by hand'' (which is
\emph{not} done by {\tt cfapery}):
\begin{verbatim}
? cftochar(AB)
% = "-1+1/(1+4/(3+1/(5+16/(7+9/(9+36/(11))))))"
? u=-1+1/(1+4/z);(-u-1)/u
% = 1/4*z
? CD=[[2*n+3],[[(2*n+1)^2,(2*n+4)^2]]];
? cflimit(CD)
% = 3.1415926...
? cftochar(CD)
% = "3+1/(5+16/(7+9/(9+36/(11+25/(13+64/(15))))))"
\end{verbatim}
We have thus obtained a reasonably simple exponentially convergent
continued fraction for $\pi$.

As mentioned, the {\tt cfapery} command does most of the work. Note,
however the following. In the present version of the program
(but this may of course change), it outputs {\tt AB=[A,B]} as above, but with
some restrictions. First, it is applicable only to polynomials
(as are all of the ones we consider), not rational functions, but only of
period $1$, second, only if the Bauer--Muir accelerations keeps this property,
third since, as we have mentioned, the initial values of $A(n)$ and $B(n)$
are not given by general formulas, the limit of the accelerated continued
fraction {\tt [A,B]} is not necessarily equal to that of the initial one
but a M\"obius transformation of it with small integer coefficients
(so use {\tt fracdep} as above to make the correction). Finally,
the {\tt cftoone} command (see below) is done automatically to $A$ and $B$
if applicable since it can only simplify the result, but on the contrary
the contraction of the result using {\tt cfdoubletosingle} (or in fact
using {\tt cfsimplify}, see below) is not done since in general the result
would be considerably more complicated. 

\smallskip

Two additional commands to know which are especially important in the
context of Ap\'ery's method are the {\tt cfsimplify} and {\tt cftoone}
commands. When given a CF of period 2 {\tt cfsimplify} computes the
contraction of $(A,B)$ and tries to simplify it (in addition, it tries
to find a simpler equivalent CF). The reader can check that applied to our
example it gives a new $(A(n),B(n))$ with $A(n)=2(4n+1)(6n^2+3n-1)$ and
$B(n)=-(n+1)^2(2n+1)^2(4n-1)(4n+7)$, not simple enough to supplant our
period 2 continued fraction.

The {\tt cftoone} command is different: given a period two $A$ (or $B$),
it checks whether (for $n$ large) the polynomial representing $a(2n+1)$
is equal to $P(n+1/2)$, where $P$ is the polynomial representing $a(2n)$,
in which case $A$ has in fact period $1$ and is thus replaced by the
period $1$ vector. As we will see in the examples, note that
a CF which cannot be transformed into a simpler one using the {\tt cfsimplify}
command may be simplified by the {\tt cftoone} command, and conversely.

\smallskip

We now treat our second example, that of an exponentially convergent
continued fraction with $E=2$ for $\log(2)$, and use {\tt cfapery}
directly:
\begin{verbatim}
? ab=cffromser(n->2*n,n->n)
% = [[0, 3*n - 1], [1, -2*n^2]]
? cftype(ab)[2]
% = [0, 2, 0, 1] /* Convergence in 1/2^n */
? AB=cfapery(ab)
% = [[[0,2],[4*n,4*n+2]],[[1,-2],[-2*n^2,-2*n^2-4*n-2]]]
? cflimit(AB)-log(2)
% = 0.E-38
? cftype(AB)[2]
% [0, -5.828..., 0, 0]
? cd=cfsimplify(AB)
% = [[0, 6*n - 3], [2, -n^2]]
? cftype(cd)[2]
% = [0, 33.97..., 0, 0]
? cftochar(cd)
% = "2/(3-1/(9-4/(15-9/(21-16/(27-25/(33))))))"
\end{verbatim}
Here we were doubly lucky: first, Ap\'ery acceleration directly gives
$\log(2)$ and not a M\"obius transformation of it. Note that starting
with an exponentially convergent fraction with $E=2$ we obtain one
with $E=-(1+\sqrt{2})^2=-5.828...$. Second, the contracted fraction
is still extremely simple, and of course the new fraction now has
$E=(1+\sqrt{2})^4$, so we obtain a very simple and very rapidly convergent
continued fraction for $\log(2)$.

\medskip

We now give a particularly spectacular example of Ap\'ery acceleration:

\begin{verbatim}
? ab=[[1/2,7*n-5],[1,-4*n*(3*n-2)]];cftype(ab)[2]
% = [0, 4/3, 0, 5/3]
/* This is a continued fraction for 2^(1/3) with exponential
   convergence in C/(4/3)^nn^(5/3).
   Since 4/3 is rational, we try Apery: */
? AB=cfapery(ab)
% = [[1/2,3,2*n],[[1/2,-5],[-3*n^2+2*n,-3*n^2-8*n-5]]]
? cftype(AB)[2]
% = [0, -3, 0, 0]
/* Apery works, we have faster exponential convergence in C/(-3)^n.
   Note that the first component has automatically been simplified to
   period 1 since cftoone was applicable.
   Since -3 is rational, we can hope to continue, but we first
   need to convert to period 1 using cfsimplify: */
? AB2=cfsimplify(AB)
% = [[1/2, 7, 10*n - 5], [2, -9*n^2 + 4]]
/* We are in luck, the simplification is much simpler, and
   evidently with exponential convergence in C/9^n with 9=(-3)^2.
   Since 9 is rational, we can try Ap\'ery again: */
? AB3=cfapery(AB2)
% = [[1/2,8,18*n-9],[4,-9*n^2+16]];
? cftype(AB3)[2]
% = [0, -33.97..., 0, 0]
/* Apery works again, even faster exponential convergence in
   C/(-(1+sqrt(2))^4)^n. Here, cftoone was in fact applicable to both
   components, done automatically by cfapery. */
/* Limit is not necessarily equal to 2^(1/3) but to a
   Moebius transform of it. We use fracdep to find it: */
? fracdep(cflimit(AB3),2^(1/3))
% = (-10*x - 3)/(12*x - 22)
? %+10/12
% = -32/(18*x - 33)
? cftochar(AB3)
% = "1/2+4/(8+7/27-..."
? -10/12-32/(-33+18*(1/2+4/(8+x)))
% = (x + 13)/(2*x + 10)
? %-1/2
% = 4/(x + 5)
/* Thus, we set */
? AB4=[[1/2,5,9*(2*n-1)],[4,-(9*n^2-16)]]
? cflimit(AB4)-2^(1/3)
% = 0.E-38
\end{verbatim}

This last CF is given below in the parametric family \ref{5.2.16}.

It is interesting to note that the first {\tt cfsimplify} command
is really ``miraculous'' in being so simple, since for instance the
exact same procedure applied to the very analogous CF for $2^{1/3}$
\begin{verbatim}
ab=[[3/2,7*n-1],[-1,-4*n*(3*n+1)]];
\end{verbatim}
would be much more complicated.

\medskip

The main purpose of Ap\'ery's method is of course to obtain exponentially
convergent series from a polynomially convergent one, or to improve the
exponent from an already exponentially convergent one. We can also use it
to obtain polynomially convergent CFs which will usually be completely
different. In particular, we have the following (see again the proof
in Appendix \ref{chap:apimpl}):

\begin{proposition}
  For any fixed $m\ge0$, set
  $$A(n)=R(m+1,n-2)+r(m+1,n-1)\text{\quad and\quad}B(n)=-d(m+1,n-1)\;.$$
  Then if we denote by $S$ the limit of the initial continued
fraction, the one defined by $(A(n),B(n))$ will converge to an
expression of the form $(uS+v)/(wS+x)$ for simple $(u,v,w,x)$ corresponding
to a correct choice of the initial values of $A$ and $B$.
\end{proposition}
Choosing our first example for $\pi/4$, if we choose $m=0$ we find
$A(n)=0$, so this is not applicable, but for $m=1$, we find that
$A(n)=4$ and $B(n)=4n^2$, and an immediate computation using {\tt fracdep}
shows that $\pi/4$ is given by the CF $[[1,4],[-1,4n^2]]$, which can be
simplified using {\tt cfmul} to $[[1,2],[-1/2,n^2]]$, in other words
$$\pi/4=1-1/2/(2+1^2/(2+2^2/(2+3^2/(2+4^2/(2+\cdots)))))$$
(given as \ref{5.3.7}), with polynomial convergence
$\pi/4-p(n)/q(n)\sim(-1)^{n+1}\pi/(32n^2)$.
As mentioned, this is still slow, but is a completely different CF from
the initial one (given as \ref{5.3.3}) which was $[[0,1,2],[1,(2n-1)^2]]$,
in other words
$$\pi/4=1/(1+1^2/(2+3^2/(2+5^2/(2+7^2/(2+\cdots)))))\;.$$

The above is done automatically by the {\tt cfaperydual} command as follows
\begin{verbatim}
? v0=cffromser(n->2*n-1,n->-(2*n-1))
% = [[0, 1, 2], [1, 4*n^2 - 4*n + 1]]
? cfaperydual(v0)
% = [[0, 0], [1, 4*n^2]] /* m = 0, inapplicable */
? A = cfsimplify(cfaperydual(v0,1))
% = [[0, 2], [1/2, n^2]] /* m = 1, result up to Mobius transformation. */
? fracdep(cflimit(A),Pi/4)
% = -x + 1
? A = [[1, 2], [-1/2, n^2]] /* We recover the above result. */
\end{verbatim}

\medskip

The last example below is very instructive and shows another aspect of the
{\tt cfaperydual} command:

\begin{verbatim}
? v0=cffromser(n->(2*n-1)^2,n->-(2*n-1)^2)
% = [[0, 1, 8*n - 8], [1, 16*n^4 - 32*n^3 + 24*n^2 - 8*n + 1]]
/* Standard Euler CF from defining series for Catalan's constant */
? A=cfsimplify(cfaperydual(v0))
% = [[0, 8*n^2 - 8*n + 3], [1/2, -16*n^4]]
? cftype(A)
% = [[2, 4, 0, -1/2], [0, 1, 0, 0], 9]
/* Very slow logarithmic convergence in C/log(n) */
? L = cflimit(A)
% = 0.91596559417721901481115881008038735296
? fracdep(L,Catalan)
% = 0
? L - Catalan
% = -2.4344... E-19
\end{verbatim}

First, thanks to the {\tt cfaperydual} command we obtain a new CF
{\tt A} related to Catalan's constant $G$ (in fact as we will see, whose limit
is already $G$ itself). First, this CF has extremely slow logarithmic
convergence, which is quite exceptional, and second it is possible that
it did not appear in the literature before the work of Y.~Yang \cite{Yan} who
used techniques from modular forms.

Although the convergence is very slow, the {\tt cflimit} command still
gives an answer using the trick mentioned in Section \ref{sec:numintro}:
before computing the limit, the program first tries (and succeeds) using
a Bauer--Muir acceleration {\tt cfbauer}, which here gives
\begin{verbatim}
[[1,8*n^2-8*n+7],[-1/2,-16*n^4]]
\end{verbatim}
which now has polynomial convergence $P^+$ with $P=2$, hence whose
limit {\tt L} can be computed by extrapolation. We now as before try the
{\tt fracdep} command which \emph{fails}. In fact only $19$ decimals are
correct, which we see in the last command. We could have avoided this
apparent failure either by using an additional {\tt cfbauer} command
which would now give $33$ correct decimals, sufficient to satisfy
{\tt fracdep}, or by modifying the {\tt fracdep} command itself by allowing
a last argument giving the tolerable decimal accuracy, exactly as for the
{\tt lindep} command.

This CF as well as its Bauer--Muir accelerated CFs are given as
\ref{5.4.27}, \ref{5.4.28} and the parametric family therein.

\smallskip

Note that exactly the same phenomenon occurs starting from the series
$7\z(3)/8=\sum_{n\ge1}1/(2n-1)^3$, and leads to the CFs \ref{5.5.10} and
\ref{5.5.11}, also discovered by Y.~Yang in loc.~cit.

\section{Basic Ap\'ery Examples}

The initial use (by Ap\'ery and successors) of Ap\'ery's method was to
accelerate sums closely linked to the Riemann zeta function, specifically
$$S(k,\eps)=\sum_{n\ge1}\dfrac{\eps^{n-1}}{n^k}\text{\quad with\quad}\eps=\pm1\;,$$
at least for small values of $k$, so that
$$S(k,1)=\z(k)\text{\quad and\quad}S(k,-1)=(1-1/2^{k-1})\z(k)\;.$$

We will give the explicit two-dimensional arrays for small cases, and
in addition note the following: since the fundamental Bauer--Muir--Ap\'ery
recursions are formal identities between polynomials or rational functions,
they are valid verbatim if we change $n$ into $n+z$ for some real or complex
variable $z$ (assumed not to be a negative integer). Thus, \emph{exactly the
same} formulas are valid for the more general sums
$$S(k,\eps,z)=\sum_{n\ge1}\dfrac{\eps^{n-1}}{(n+z-1)^k}$$
by replacing $n$ by $n+z-1$, so we will not give the arrays explicitly,
but simply the value of the sum of the series.

\smallskip

For each case, we give the initial non-accelerated continued fraction,
the $2$-dimensional arrays satisfying the above recursion, then the Ap\'ery
dual, and finally the Ap\'ery accelerated contracted continued fraction
(both possibly simplified) in the usual format used in this book.

\medskip

{\bf Case} $(k,\eps)=(1,-1)$:

$$\log(2)=\sum_{n\ge1}\dfrac{(-1)^{n-1}}{n}\;,\text{ and more generally }
  \dfrac{\psi((1+z)/2)}{2}-\dfrac{\psi(z/2)}{2}=\sum_{n\ge0}\dfrac{(-1)^n}{n+z}\;.$$

\smallskip

\begin{verbatim}
[()->log(2),[0,1],[1,n^2]]
\end{verbatim}

\begin{align*}
  a(n,l)&=2l+1\;,\\
  b(n,l)&=n^2\;,\\
  r(n,l)&=n-l-1\;,\\
  d(n,l)&=-(l+1)^2\;.\end{align*}

\smallskip

\begin{verbatim}
[()->log(2),[0,1],[1,n^2]]
[()->log(2),[0,6*n-3],[2,-n^2]]
\end{verbatim}

\medskip

{\bf Case} $(k,\eps)=(2,1)$:

$$\z(2)=\sum_{n\ge1}\dfrac{1}{n^2}\;,\text{ and more generally }
  \psi'(z)=\sum_{n\ge0}\dfrac{1}{(n+z)^2}\;.$$

\smallskip

\begin{verbatim}
[()->zeta(2),[0,2*n^2-2*n+1],[1,-n^4]]
\end{verbatim}

\begin{align*}
  a(n,l)&=2n^2-2n+1+l^2+l\;,\\
  b(n,l)&=-n^4\;,\\
  r(n,l)&=-n^2+(l+1)n-(l+1)^2/2\;,\\
  d(n,l)&=-(l+1)^4/4\;.\end{align*}

\smallskip

\begin{verbatim}
[()->zeta(2),[0,2*n-1],[2,n^4]]
[()->zeta(2),[0,11*n^2-11*n+3],[5,n^4]]
\end{verbatim}

\medskip

{\bf Case} $(k,\eps)=(2,-1)$:

$$\dfrac{\z(2)}{2}=\sum_{n\ge1}\dfrac{(-1)^{n-1}}{n^2}\;,\text{ and more generally }
  \dfrac{\psi'(z/2)}{2}-\psi'(z)=\sum_{n\ge0}\dfrac{(-1)^n}{(n+z)^2}\;.$$

\smallskip

\begin{verbatim}
[()->zeta(2)/2,[0,2*n-1],[1,n^4]]
\end{verbatim}

\begin{align*}
  a(n,l)&=(2l+1)(2n-1)\;,\\
  b(n,l)&=n^4\;,\\
  r(n,l)&=n^2-2(l+1)n+2(l+1)^2\;,\\
  d(n,l)&=4(l+1)^4\;.\end{align*}

\smallskip

\begin{verbatim}
[()->zeta(2)/2,[0,2*n^2-2*n+1],[1/2,-n^4]]
[()->zeta(2)/2,[0,11*n^2-11*n+3],[5/2,n^4]]
\end{verbatim}

\medskip

{\bf Case} $(k,\eps)=(3,1)$:

  $$\z(3)=\sum_{n\ge1}\dfrac{1}{n^3}\;,\text{ and more generally }
  -\dfrac{\psi''(z)}{2}=\sum_{n\ge0}\dfrac{1}{(n+z)^3}\;.$$

\smallskip

\begin{verbatim}
[()->zeta(3),[0,(2*n-1)*(n^2-n+1)],[1,-n^6]]
\end{verbatim}

\begin{align*}
  a(n,l)&=2n^3-3n^2+(4l^2+4l+3)n-(2l^2+2l+1)\;,\\
  b(n,l)&=-n^6\;,\\
  r(n,l)&=-n^3+2(l+1)n^2-2(l+1)^2n+(l+1)^3\;,\\
  d(n,l)&=(l+1)^6\;.\end{align*}

\smallskip

\begin{verbatim}
[()->zeta(3),[0,(2*n-1)*(n^2-n+1)],[1,-n^6]]
[()->zeta(3),[0,(2*n-1)*(17*n^2-17*n+5)],[6,-n^6]]
\end{verbatim}

\section{Generalizations of Ap\'ery's Method}\label{sec:zeta4}

There are several useful generalizations of Ap\'ery's method. Although
some are implemented in the package, it is instructive to first give the
corresponding GP scripts.

\smallskip

\subsection{Use of a Simplifier}

\smallskip

The first important generalization is the use of a \emph{simplifier}.
Note that the Bauer--Muir formulas involve denominators. In nice cases
(as all the ones seen above), and in particular when $d(n)$ is constant,
these denominators simplify. But in general applying Bauer--Muir to a
CF with polynomial coefficients leads to one with \emph{rational function}
coefficients. In itself this is not too bad, but it becomes catastrophic
if we want to \emph{iterate}, as in Ap\'ery's method.

Now recall that for any CF $(a(n),b(n))$ and nonzero function $t(n)$ with
$t(0)=1$, the CF $(t(n)a(n),t(n)t(n+1)b(n))$ will have the same convergents
as the initial CF. This allows us to get rid of the denominators introduced
by Bauer--Muir. It is immediate to compute the corresponding formulas
(see Appendix \ref{chap:apimpl} for details), and since this modification
is not available in the {\tt cfapery} command, we give the following
GP script template, where we recall that $tf(n,l)=\prod_{1\le m\le n}t(m,l)$:

\medskip

\begin{verbatim}
a(n,l)=
b(n,l)=
r(n,l)=
R(n,l)=a(n+1,l)+r(n+1,l);
d(n,l)=
t(n,l)=
tf(n,l)=
         
         /* checks: must be 0 */

r(n,l)*R(n,l)-b(n,l)-d(n,l)
t(n,l)*(R(n,l)-r(n-1,l)*d(n,l)/d(n-1,l))-a(n,l+1)
t(n,l)*t(n+1,l)*b(n,l)*d(n+1,l)/d(n,l)-b(n,l+1)
tf(0,l)-1
tf(n,l)/tf(n-1,l)-t(n,l)

VERT(m)=[[0,t(m+1,n-2)*R(m+1,n-2)+r(m+1,n-1)],[1,-t(m+1,n-1)*d(m+1,n-1)]];

DIAG(m)=[[[0,1],[t(n+m,n-1)*R(n+m,n-1),tf(n+m,n)*R(n+m,n)]],
         [[1,1],[tf(n+m,n)*b(n+m,n),-tf(n+m+1,n)*d(n+m+1,n)]]];

\end{verbatim}

This is to be used as follows: $a(n,l)$, $b(n,l)$, $r(n,l)$, and $t(n,l)$
have to be found ``by hand'' using the {\tt cfbauer} command possibly together
with a simplifier $t(n,l)$, followed by polynomial extrapolation.
Then $d(n,l)$ and $tf(n,l)$ are computed from their definition. One then
fills the script with the computed values, and reads it in GP. It must first
print five $0$'s, simply to check that the input data is consistent with the
recursions and definitions. Then {\tt DIAG(0)} gives the period 2 main
staircase walk (use {\tt cfsimplify} to see if it can contract to a nice
period 1 CF). When $t(n,l)=1$, this is essentially identical to the
{\tt cfapery} command. In some cases there is a division by $0$, in which case
try {\tt DIAG(1)}. Similarly {\tt VERT(0)} gives the main vertical walk, and
is essentially identical to the {\tt cfaperydual} command when $t(n,l)=1$.
Note that {\tt VERT(m)} does not make much sense when the arrays have
denominators such as $n-l$ or $2n-l$ since in that case they are not
defined for $l\ge m$ or $l\ge2m$.

\smallskip

To avoid typing the above script, the package provides the function
{\tt cfaperytemplate} whose input is the 3 or 5-component vector of
polynomials or rational functions in the variables $n$ and $l$
{\tt [a(n,l),b(n,l),r(n,l),t(n,l),tf(n,l)]}, where the last two may be
omitted if equal to $1$. Note that it is imperative that the name of the
variables be {\tt n} and {\tt l}.

The output is a two-component vector {\tt [VERT,DIAG]} of polynomials or
rational functions with main variable {\tt x}, and for instance instead of
writing {\tt DIAG(0)} we must write {\tt subst(DIAG,'x,0)}.

\smallskip

\subsection{Use of a Multiplier}

\smallskip

A second ``generalization'' of Ap\'ery's method is the use of a
\emph{multiplier}: an initial CF may have
complicated Bauer--Muir iterations, so that Ap\'ery is not practical.
However, it we use a multiplier on the initial CF, the iterations may in some
cases become simpler. We illustrate this by computing arrays for the same
sums $S(k,\eps)$ for other values of $(k,\eps)$, in the same format as
above.

\medskip

{\bf Case} $(k,\eps)=(4,1)$:

  $$\z(4)=\sum_{n\ge1}\dfrac{1}{n^4}\;,\text{ and more generally }
  \dfrac{\psi'''(z)}{6}=\sum_{n\ge0}\dfrac{1}{(n+z)^4}\;.$$

\smallskip

\begin{verbatim}
[()->zeta(4),[0,2*n^4-4*n^3+6*n^2-4*n+1],[1,-n^8]]
\end{verbatim}

The Bauer--Muir formulas applied to this basic CF are too complicated.
Our salvation comes by using a termwise equivalent CF by multiplying by
$2n-1$ (this factor can be found by using a suitable resultant computation),
so instead of the basic CF above, we use as initial CF:

\begin{verbatim}
[()->zeta(4),[0,(2*n-1)*(2*n^4-4*n^3+6*n^2-4*n+1)],[1,-(4*n^2-1)*n^8]]
\end{verbatim}

\smallskip

We find the following:
  
\begin{align*}
a(n,l)&=(2n-1)\left(2n^4-4n^3+6n^2-4n+1+\dfrac{l(l+1)}{2}(17n^2-17n+5)\right)\;,\\
b(n,l)&=-n^6(4n^2-l^2)(4n^2-(l+1)^2)/4\;,\\
r(n,l)&=(2n+l)\left(-n^4+\dfrac{7}{2}(l+1)n^3-6(l+1)^2n^2+6(l+1)^3n-3(l+1)^4\right)\;,\\
t(n,l)&=(2n-l-2)/(2n-l)\;,\\
d(n,l)&=-9(l+1)^8(4n^2-l^2)\;.\end{align*}

We notice that these arrays are not defined for $l\ge 2n-2$, so
the Ap\'ery dual (i.e., with $l$ fixed) does not exist. On the other hand
the diagonal or the staircases which correspond to $l\approx n$
do exist, and after contraction we find \ref{5.6.7}, in other words

\begin{verbatim}
[()->zeta(4),[0,3*(2*n-1)*(3*n^2-3*n+1)*(15*n^2-15*n+4)],
             [13,3*n^8*(9*n^2-1)]]
\end{verbatim}

\medskip

{\bf Case} $(k,\eps)=(3,-1)$:

  $$\dfrac{3\z(3)}{4}=\sum_{n\ge1}\dfrac{(-1)^{n-1}}{n^3}\;,\text{ and more generally }
  \dfrac{\psi''(z)}{2}-\dfrac{\psi''(z/2)}{8}=\sum_{n\ge0}\dfrac{(-1)^n}{(n+z)^3}\;.$$

\smallskip

Similarly, instead of the basic CF:
\begin{verbatim}
[()->3*zeta(3)/4,[0,3*n^2-3*n+1],[1,n^6]]
\end{verbatim}
which would lead to complicated denominators, we use a termwise equivalent
CF by multiplying by $2n-1$:

\smallskip

\begin{verbatim}
[()->3*zeta(3)/4,[0,(2*n-1)*(3*n^2-3*n+1)],[1,n^6*(4*n^2-1)]]
\end{verbatim}

We find the following:

\begin{align*}
a(n,l)&=(2l+1)(2n-1)(3n^2-3n+1)\;,\\
b(n,l)&=n^6(4n^2-(2l+1)^2)\;,\\
r(n,l)&=(n-2(l+1))(2n+2l+1)(n^2-2(l+1)n+4(l+1)^2)\;,\\
t(n,l)&=(2n-(2l+3))/(2n-(2l+1))\;,\\
d(n,l)&=-64(l+1)^6(4n^2-(2l+1)^2)\;.\end{align*}

\smallskip

After simplification, the Ap\'ery dual is the ordinary CF for $\z(3)$,
so not interesting:

\smallskip

\begin{verbatim}
[()->3*zeta(3)/4,[0,(2*n-1)*(n^2-n+1)],[3/4,-n^6]]
\end{verbatim}

\smallskip

More interesting, after simplification the Ap\'ery accelerated CF
gives a new CF for $\z(3)$ which converges like $64^{-n}$, see
\ref{5.5.21}:

\smallskip

\begin{verbatim}
[()->3*zeta(3)/4,[0,65*n^4-130*n^3+105*n^2-40*n+6],
                 [21/4,-4*n^6*(16*n^2-1)]]
\end{verbatim}

\smallskip

\subsection{Slowing Down Ap\'ery}

\smallskip

A third generalization of Ap\'ery's method consists in \emph{slowing down}
the acceleration. Evidently here the goal is not to find the fastest possible
CF (for instance to prove an irrationality result, as in Ap\'ery's initial
motivation), but simply to find \emph{new} and hopefully interesting CFs.
We illustrate this with $\log(2)$, starting again from the standard CF
already used in our second example above
$\log(2)=((0,3n-1),(1,-2n^2))$. As in the example for $\z(4)$, we first
use a multiplier, here $n$, and write the following:

\begin{verbatim}
? CF=[[0,3*n-1],[1,-2*n^2]]
? CF2=cfmul(CF,n)
% = [[0,3*n^2-n],[1,-2*n^4-2*n^3]]
? cfbauer(CF2,1)[3]
% = [-2,-1,-n^2+n-2]
\end{verbatim}

  This last command tells us that the optimal Bauer--Muir function $r$ is
  $r(n)=-n^2+n-2$. However, let us try something suboptimal and use
  $r(n)=-n^2+n+r_0$ for some unknown $r_0$:

\begin{verbatim}
? a(n)=3*n^2-n;
? b(n)=-2*n^4-2*n^3;
? r(n)=-n^2+n+'r0;
? d(n)=r(n)*(a(n+1)+r(n+1))-b(n);
? a1=a(n+1)+r(n+1)-r(n-1)*d(n)/d(n-1);
? b1=b(n)*d(n+1)/d(n);
? res(P)=polresultant(numerator(P),denominator(P),'n);
? nfroots(,gcd(res(a1),res(b1)))
% = [-2, 0, 2]~
\end{verbatim}

  The rational functions {\tt a1} and {\tt b1} are the coefficients of the new
  CF obtained using Bauer--Muir. We want them to be as close to polynomials
  as possible, so computing the roots in $r_0$ of the resultant with respect
  to the variable {\tt n} gives the values which give some simplification of
  {\tt a1} and {\tt b1} (the command {\tt nfroots(,P)} gives the rational
  roots of {\tt P}).

  The program gives us three possible values of $r_0$: $r_0=-2$ is the optimal
  Bauer--Muir, so we can try $r_0=0$ or $r_0=2$. We can choose any one,
  so arbitrarily let us decide that at each stage we will choose the largest
  value of $r_0$.

  As it happens, this leads to very simple arrays:

\begin{align*}
a(n,l)&=3n^2-n+2l(2n-1)\;,\\
b(n,l)&=-2n^3(n+2l+1)\;,\\
r(n,l)&=-n^2+n+2(l+1)(2l+1)\;,\\
d(n,l)&=4(n+2l+1)(n+2l+2)(l+1)^2\;,\end{align*}
and since no denominator occur, {\tt t(n,l)=tf(n,l)=1}.

We now input this in the template script that we have written above.
The vertical walk {\tt VERT(0)} gives a known CF for $\log(2)$. The
staircase walk {\tt DIAG(0)} gives a reasonable-looking period 2 CF,
but miraculously the {\tt cfsimplify} command gives a very simple
period 1 CF:

\begin{verbatim}
% = [[0,11/2,29*n^2-29*n+8],[5,2,-54*n^4+6*n^2]]
\end{verbatim}

  As usual with Ap\'ery the initial terms need to be determined using
  {\tt fracdep}, and we find finally the following new CF which converges
  like $(27/2)^{-n}$, see \ref{5.3.79}:

\begin{verbatim}
[()->log(2),[0,29*n^2-29*n+8],[5,-6*n^2*(9*n^2-1)]]
\end{verbatim}

A similar procedure leads to the additional new CF which converges
like $(32/27)^{-n}$, see \ref{5.3.62}:

\begin{verbatim}
[()->log(2),[0,59*n^2-59*n+20],[5,-24*n^2*(36*n^2-1)]]
\end{verbatim}

See also the remark after \ref{5.3.79}.

\smallskip

The above procedure can be partly automated thanks to the special flags $-1$
and $-2$ in the {\tt cfbauer} command. For instance, starting again with CF2:

\begin{verbatim}
? CF2=[[0,3*n^2-n],[1,-2*n^4-2*n^3]];
? VR=cfbauer(CF2,-1)
% = [[-n^2+n-2,-n^2+n,-n^2+n+2],\
    [[[0,(3*n^3-n^2)/(n-1)],[1,-2*n^4-4*n^3-2*n^2]],\
     [[0,3*n^2+3*n],[1,-2*n^4-6*n^3-4*n^2]],\
     [[0,3*n^2+3*n-2],[1,-2*n^4-6*n^3]]]]
? R=VR[1][3];CF3=VR[2][3];
\end{verbatim}

{\tt VR[1]} gives the list of possible $r(n)$, and {\tt VR[2]} the
corresponding new CFs.

The flag $-2$ is used when the CF corresponds to a \emph{function}, so which
has extra variables.

\smallskip

\subsection{Irrational Ap\'ery}

\smallskip

We have mentioned that in case of convergence type $E$ with $E$ irrational
the use of Bauer--Muir and a fortiori Ap\'ery would lead to CFs with
irrational coefficients. Nonetheless, if we do apply Bauer--Muir as a first
step for Ap\'ery, in many (all?) cases the successive Bauer--Muir iterations
do not involve any additional irrationalities. I am indebted to
T.~Stachowiak \cite{Sta} for this observation. To take a simple example,
starting from the CF \ref{5.3.83}:
\begin{verbatim}
[()->log(2),[0,6*n-3],[2,-n^2]]
\end{verbatim}
which has $E=(1+\sqrt{2})^4=33.97...$ and applying Ap\'ery leads to the
irrational CF
\begin{verbatim}
[()->log(2),(1+sqrt(2))^2*[0,-4,2*n-1],[16,4,-n^2]]
\end{verbatim}
with convergence type $E=(1+\sqrt{2})^3(9+4\root4\of2-3\sqrt{2})=133.87...$
and $P=0$.

\chapter{Encyclopedic Dictionary: Constants}\label{chap:const}

\section{Introduction}

\medskip

We are now going to give a large list of continued fractions of polynomial type
corresponding either to real numbers, or to functions. Before describing the
format of the entries, two important warnings are necessary.

\smallskip

First, there are literally infinitely many explicit evaluations of infinite
series, and many are discovered everyday. Thanks to Euler's transformation of
series into continued fractions, we thus have an infinite supply of possibly
interesting CFs. Therefore, in addition to the CFs that I could find in the
literature or the new CFs created for instance using Ap\'ery acceleration
and Ap\'ery duality, I have had to choose among the CFs coming ``trivially''
from infinite series, essentially keeping only
``interesting'' constants and/or CFs, this being of course quite subjective.

\smallskip

Second, an important caveat concerning the correctness of the formulas
presented in this dictionary. As already mentioned, there are certainly
many typos or stupid mathematical errors or wrong signs given in the speed
of convergence. But apart from that I must warn the reader about two things:

$\bullet$ As mentioned several times, I believe that the indicated speeds
of convergence are correct, but since they are based on heuristic proofs
and genericity, in some cases they may be wrong, especially when they
depend on one or several variables.

$\bullet$ But much more importantly: I do not claim that I have
\emph{proved} the correctness of the CFs given here. In some cases, they
are either found in the literature, or obtained from those by rigorous
methods such as Bauer--Muir or Ap\'ery. But in other cases, they are
obtained ``\`a la Ramanujan'', by intelligently guessing the general
formulas from a few initial terms. It may of course happen that my guesses
are incorrect, but I am pretty confident (again modulo typos) that very few
will turn out to be wrong. Also, probably in most cases the proofs should
not be difficult, but I have not had the patience of doing them.

\medskip

A given entry in the dictionary contains the following information:

\begin{enumerate}\item The mathematical function or constant, given as a
  {\tt Pari/GP} \emph{closure}, even when the corresponding special function
  or constant does not exist yet.
  There usually are restrictions on the parameters and variables for the
  continued fraction to converge to the value of the function, but in many
  cases these are quite complicated (and in some cases even unknown), so
  we have not given them.
\item The continued fraction itself in the variable {\tt n}, in the format
  given above, except that we add the restriction that all entries have the
  same length as the period.
  
  Both the function and the continued fraction are given as a 3-component
  vector which can be directly input in {\tt Pari/GP}, or if desired,
  shortened to a 2-component vector by omitting the first component.
\item A few terms of the continued fraction, for better human visualization.
\item The type and the speed of convergence in $FEDPC$ format, meaning that
  $$S-\dfrac{p(n)}{q(n)}\sim \dfrac{C}{n!^FE^ne^{\sqrt{Dn}}n^P}\;,$$
  where $S$ is the value of the continued fraction (in some rare cases, the
  convergence is very slow and is instead of the form
  $S-p(n)/q(n)\sim C/\log(n)$). The constant $C$
  does not make much sense for period 2 CFs so is not given, and in
  a few rare cases it is not given for some period 1 CFs because we have not
  found a satisfactory answer. We do not give
  the asymptotic behavior of $q(n)$ since it is not canonical, and can
  anyway be retrieved thanks to the functions {\tt cftype} and {\tt cfasymp}.
\item A few terms of the more precise asymptotic expansion of $S-p(n)/q(n)$
  in the form $A=1+c_{1/2}/n^{1/2}+c_1/n+\cdots$, meaning that
  $$S-\dfrac{p(n)}{q(n)}=\dfrac{C}{n!^FE^ne^{\sqrt{Dn}}n^P}A\;.$$
  To simplify notation, when $E$ is a real quadratic number in $\Q(\sqrt{D})$
  with $D$ squarefree, the $c_i$ are in the same number field, and
  $\sqrt{D}$ is denoted by the variable $d$. When the CF has period 2,
  we give the asymptotic expansion $A$ of the even convergents and change
  $2n$ into $n$ to make it compatible with period 1.
  In the case where the number field depends on one or more variables
  (for example, convergence in $(z+\sqrt{z^2+1})^2$ or similar), determining
  this asymptotic expansion is painful, and I have not bothered to do so.
  Finally, in the special case of convergence of type $P^+$ with $P$
  integral, I do not know of a systematic method to find this expansion
  with terms $1/n^k$ with $k\ge P$, so except in special cases where $p(n)/q(n)$
  is completely explicit and given by a simple formula, the asymptotic
  expansion is only given to $P$ terms (the reader is of course welcome to
  find additional terms). Since I try to give sufficient information on $A$,
  when the entry is $A=1+\cdots$ or $A=1-1/n+\cdots$ for instance, this
  means that I have not been able to find more terms, or possibly that
  I have not recognized the additional terms. Very often, this means that
  there exist terms involving $\log(n)$.

  See Appendix \ref{chap:asymp} for more details.
\item As mentioned above, in many cases the CF seems to correspond to a
  reasonably explicit series, which is then given. Note that I do not
  claim to have proved the correctness of these series. In addition, since the
  inverse of a CF is just as interesting as the CF itself, if it has a
  reasonably explicit series it is also given. In fact, many M\"obius
  transformations of the CF could give additional series, but for simplicity
  we have not tried any other than the CF and its inverse.
  
  Note that in many cases these series
  give nontrivial representations of the constant or function.
\item Warning: the sign given in the asymptotics may be incorrect, especially
  for period $2$ CFs, since it is essentially obtained by taking a suitable
  ``square root'' of the corresponding contracted CF. In many cases, it
  depends on the values of the variables which occur.
\item Evidently a period 2 CF can be contracted (using {\tt cfsimplify},
  which does more) to a period 1 CF. Almost always we will give the latter
  only if either it has an especially simple form, or if the degrees of
  the contracted $a(n)$ and $b(n)$ have not increased.
\item In many cases, Bauer--Muir--Ap\'ery acceleration of a CF leads
  to an infinite number of simple CFs polynomially parameterized by
  an integer constant $k\ge0$. More generally, one can often find
  polynomially parameterized CFs for a given constant.
  It is sometimes quite tedious to find
  the exact rational expressions which are necessary to specify the
  first few terms of the CF, such as $a(0)$ and $b(0)$. In this
  case, we will use a special format, \emph{not} directly readable in
  {\tt GP}, of the form $[f,a(n),b(n)]$, where $f$ is the closure giving the
  constant or function, and $a(n)$ and $b(n)$ are polynomials
  representing the coefficients of the CF for $n$ sufficiently large.
\end{enumerate}

To illustrate this last point, consider one of the simplest CF for
$\pi$ (\ref{5.3.7}): \begin{verbatim}[()->Pi,[4,2],[-2,n^2]]\end{verbatim}
It is immediate to use Ap\'ery acceleration and to obtain the
following (\ref{5.3.14}):
\begin{verbatim}
[()->Pi,[4,2*n+1],[[-4,9],[(2*n)^2,(2*n+3)^2]]]
\end{verbatim}
However in doing so, some information is lost since all the intermediate
Bauer--Muir accelerated CFs are forgotten. In fact, it is easy to
see that the $k$th accelerated CF is
\begin{verbatim}
[()->Pi,[u(k),2*k],[v(k),n^2]]
\end{verbatim}
where {\tt u(k)} and {\tt v(k)} are certain functions of {\tt k}. We
will indicate the parametric family after the CF by
\begin{verbatim}[()->Pi,2*k,n^2]\end{verbatim}
omitting {\tt u(k)} and {\tt v(k)}. Unless mentioned otherwise, it is
understood that $k$ must be an integer such that $k\ge0$: this is especially
important since in many cases the \emph{same} parametric family has
different limits (or does not converge) when $k<0$ or when $k\notin\Z$:
for instance the family \begin{verbatim}[k,n^2]\end{verbatim} converges
(up to a M\"obius transformation) to $\pi$ when $k$ is even, and to $\log(2)$ 
when $k$ is odd.

In this specific example, we have
$$u(k)=4\sum_{j=1}^{k+1}\dfrac{(-1)^j}{2j-1}\text{\quad and\quad}v(k)=(-1)^{k+1}2\;,$$
but in other examples the formulas may be more complicated.

Note also that using the command {\tt cfaperydual} on this example leads
after simplification to the CF 
\begin{verbatim}[z->Pi,[0,1,2],[4,(2*n-1)^2]]\end{verbatim}
which is simply the Euler transformation of the standard series
$\pi/4=1-1/3+1/5-\cdots$, and is given in \ref{5.3.3} below.

\smallskip

Note also the trivial fact that if there exists a family
$[f(z),a(n),b(n)]$, meaning that $f(z)$ can be obtained as the limit of a CF
with coefficients $a(n)$, $b(n)$ for $n$ sufficiently large, we also have
the families $[f(z),a(n+j),b(n+j)]$ for any integer $j\ge0$.

\medskip

For a given constant or function, we usually try to give the continued
 fractions in \emph{increasing} order of convergence speed, so that the
 fastest are at the end.

\smallskip

Recall that continued fractions $(a,b)$ and $(c,d)$ are said to be
\emph{equivalent} if their partial quotients $p(n)/q(n)$, are the same,
or equivalently if there exists a sequence $r(n)$ with $r(n)\ne0$ such
that $c(n)=r(n)a(n)$ and $d(n)=r(n)r(n+1)b(n)$. Among equivalent continued
fractions, we have tried to give the ``simplest'' one, in particular
if possible avoiding unnecessary denominators, except when their use gives
a particularly simple continued fraction. Furthermore, when the even (or odd)
contraction of a continued fraction is sufficiently simple, we usually
only give it and not the initial one.

\smallskip

Note also that many of the CFs given for constants are (hopefully interesting)
specializations of CFs for functions.

\medskip

{\bf Example 1}: one of the entries (same as an example given above) reads
\begin{verbatim}
[z->exp(z),[1,n-z],[z,z*n]]
\end{verbatim}
Calling {\tt v} this vector, it can be used as follows:
\begin{verbatim}
? r=cflimit(cfsubst(v,z,2/5));s=v[1](2/5);[r-s,cftype(v)]
% = [0.E-38, [[1, 1, 0, 0], [1, -1/z, 0, 1], 2]]
\end{verbatim}
This shows numerically that the limit is indeed $\exp(z)$ (at least for
$z=2/5$), and that $\exp(z)-p(n)/q(n)\sim C(z)/(n!\cdot(-1/z)^n\cdot n)$,
and as above we could compute $C(z)$ by using {\tt cfasymp}.

\smallskip

{\bf Example 2}: an entry which initially was in our list was
\begin{verbatim}
[z->exp(z),[1],[z*[1,-1/2],z/(4*n+2)*[1,-1]]]
\end{verbatim}
The presence of denominators is not pleasing, so to get rid of them
we use the function {\tt cfmul} which gives an equivalent continued fraction,
as follows: calling {\tt v} the above vector, we write
\begin{verbatim}
? w=cfmul(v,n->[1,4*n+2])
% = [[[1, 2], [1, 4*n + 2]], [[2*z, -z], [z, -z]]]
\end{verbatim}
This already looks much nicer, but we can also see if the contracted fraction
is better:
\begin{verbatim}
? w2=cfdoubletosingle(w)
% = [[1, -z + 2, 4*n - 2], [2*z, z^2]]
? cftochar(w2)
% = "1+2*z/(-z+2+z^2/(6+z^2/(10+z^2/(14+z^2/(18+z^2/(22))))))"
\end{verbatim}
We recognize the well-known continued fraction for $\exp(z)$, so we do
not include the above entry in our list, but only this last formula.

Most of the above work can be done automatically by using the
{\tt cfsimplify} command. For instance:
\begin{verbatim}
? w2=cfsimplify(v)
% = [[1, -z + 2, 4*n - 2], [2*z, z^2]]
\end{verbatim}

\medskip

{\bf Caveat:} Apart from misprints, the CF given in this dictionary should
be correct, at least for a suitable range of the variables, usually assumed
real. On the other hand, I have usually not given \emph{any} reference for
individual CFs (only global references at the end) since there are so many
sources, and also since most of them, in particular those obtained by Ap\'ery
acceleration and by Ap\'ery duality, are probably new. Also, please take
into consideration the warnings given at the beginning of this chapter.

Finally note that most of the constants $C$ given in the asymptotics
have been found using linear dependence type algorithms (such as the
{\tt lindep} command), but it may happen that (notwithstanding possible
typos) the constant is wrong, although I would really be surprised if it was.
However, in most cases it should not be difficult to \emph{prove} the
correctness of the constant, but I have not done so.

\medskip

\section{Constants: $2^{1/3}$}

\medskip

For many constants, there exist parametric families of (inequivalent)
continued fractions, see \cite{Coh4} for a general discussion. For instance,
when Bauer--Muir acceleration is sufficiently simple to lead to Ap\'ery
acceleration, this automatically gives a one-parameter family of CFs
indexed by the number $k$ of Bauer--Muir accelerations, but there are many
other types. I give the ones that I have had time or been able to find.

\smallskip

Algebraic numbers of degree $2$ can be represented by \emph{simple}
(i.e., with $b(n)=1$ for all $n$) periodic CFs. However some other types
of CFs are less trivial and perhaps more interesting. I give a single
example, where we note that $3^{1/2}$ is a gamma quotient since
$$3^{1/2}=2\dfrac{\G(1/2)^2}{\G(1/3)\G(2/3)}$$

\smallskip

\begin{cf}\label{5.2.1}{\ }
\begin{verbatim}
[()->3^(1/2),[3/2,(2*n+1)^2],[2,-n^2*(n+2)^2]]
\end{verbatim}
$$\sqrt{3}=3/2+\dfrac{2}{9-\dfrac{9}{25-\dfrac{64}{49-\dfrac{225}{81-\dfrac{576}{121-\dfrac{1225}{169-\ddots}}}}}}$$
Convergence type $E$ with $E=(2+\sqrt{3})^2$, $P=0$, and $C=2\sqrt{3}/(2+\sqrt{3})^3$, so that
$$\sqrt{3}-\dfrac{p(n)}{q(n)}\sim\dfrac{2\sqrt{3}}{(2+\sqrt{3})^{2n+3}}\;.$$
$$A=1+d/n+(3/2-3d/2)/n^2+(3d-9/2)/n^3+(45/4-27d/4)/n^4+\cdots$$
\end{cf}

As explained above, in the present case we have $d=\sqrt{3}$.

\medskip

{\bf Remark:} this CF seems to be ``infinitely contractible'', meaning that,
after simplification, its iterated even contractions are all of constant
bidegree $(2,4)$ (i.e., $a(n)$ of degree 2 and $b(n)$ of degree 4). I do not
yet know how to prove this, or how to characterize this property.

\medskip

The following are given as examples of CFs for an algebraic number of degree
at least $3$, other numbers could of course also be studied. Note
that $2^{1/3}$ is a gamma quotient since
$$2^{1/3}=2\dfrac{\G(1/2)\G(1/3)}{\G(1/6)\G(2/3)}$$

\smallskip

\begin{cf}\label{5.2.2}{\ }
\begin{verbatim}
[()->2^(1/3),[0,1],[2,2,3*(n-1)*(3*n-1)]]
\end{verbatim}
$$\root3\of2=\dfrac{2}{1+\dfrac{2}{1+\dfrac{15}{1+\dfrac{48}{1+\dfrac{99}{1+\dfrac{168}{1+\ddots}}}}}}$$
Convergence type $P^-$ with $P=1/3$ and $C=1/\G(2/3)$, so that
$$\root3\of2-\dfrac{p(n)}{q(n)}\sim(-1)^n\dfrac{1/\G(2/3)}{n^{1/3}}\;.$$
$$A=1+(1/18)/n-(5/81)/n^2-(56/2187)/n^3+(1022/19683)/n^4+\cdots$$
Series:
$$\root3\of2=2\sum_{n\ge0}(-1)^n\dfrac{(2/3)_n}{n!}$$
Parametric family for $3\nmid u$:
\begin{verbatim}
[()->2^(1/3),u,3*n*(3*n+3-u)]
\end{verbatim}
Convergence type $P^-$ with $P=u/3$.
\end{cf}

\smallskip

\begin{cf}\label{5.2.3}{\ }
\begin{verbatim}
[()->2^(1/3),[1/2,7*n-5],[1,-4*n*(3*n-2)]]
\end{verbatim}
$$\root3\of2=1/2+\dfrac{1}{2-\dfrac{4}{9-\dfrac{32}{16-\dfrac{84}{23-\dfrac{160}{30-\dfrac{260}{37-\ddots}}}}}}$$
Convergence type $E$ with $E=4/3$, $P=5/3$, and $C=2^{8/3}/\G(1/3)$, so that
$$\root3\of2-\dfrac{p(n)}{q(n)}\sim\dfrac{2^{8/3}/\G(1/3)}{(4/3)^nn^{5/3}}\;.$$
$$A=1-(55/9)/n+(170/3)/n^2-(1579490/2187)/n^3+(230596751/19683)/n^4+\cdots$$
Series:
$$\dfrac{1}{\root3\of2}=2\sum_{n\ge0}\dfrac{(-2/3)_n}{n!}(3/4)^n$$
Parametric families, with $u\ge0$, $v\ge0$ with $v\not\equiv2\pmod3$, and $k\ge0$:
\begin{verbatim}
[()->2^(1/3),7*n-5+k+v+4*u,-4*(n+u)*(3*n+v-2)]
\end{verbatim}
Convergence type $E$ with $E=4/3$, $P=u+(5-v)/3+2k$.

After contraction, Ap\'ery accelerates to the parametric family \ref{5.2.12}.
\end{cf}

\smallskip

\begin{cf}\label{5.2.4}{\ }
\begin{verbatim}
[()->2^(1/3),[2,9,7*n-2],[-10,8*(3*n-1)*(6*n+1)]]
\end{verbatim}
$$\root3\of2=2-\dfrac{10}{9+\dfrac{112}{12+\dfrac{520}{19+\dfrac{1216}{26+\dfrac{2200}{33+\dfrac{3472}{40+\ddots}}}}}}$$
Convergence type $E$ with $E=-16/9$, $P=1/6$, and
$C=-3\sqrt{\pi}/(20^{1/3}\G(1/3))$, so that
$$\root3\of2-\dfrac{p(n)}{q(n)}\sim(-1)^{n+1}\dfrac{3\sqrt{\pi}/(20^{1/3}\G(1/3))}{(4/3)^{2n}n^{1/6}}\;.$$
$$A=1-(223/1800)/n+(28171/2160000)/n^2+\cdots$$
Parametric family:
\begin{verbatim}
[()->2^(1/3),7*n+(25*k-2),8*(3*n-1)*(6*n+1)]
\end{verbatim}
Convergence type $E$ with $E=-16/9$ and $P=2k+1/6$.

There are more general parametric families of the type
\begin{verbatim}
[()->2^(1/3),7*n+u,4*(6*n+v)*(6*n+w)]
\end{verbatim}
with convergence type $E$ with $E=-16/9$ and $P=(2u-7(v+w)/6+7)/25$
but satisfying complicated conditions on {\tt (u,v,w)}.
\end{cf}

\smallskip

\begin{cf}\label{5.2.5}{\ }
\begin{verbatim}
[()->2^(1/3),[1,9*n-4],[1,-6*n*(3*n-1)]]
\end{verbatim}
$$\root3\of2=1+\dfrac{1}{5-\dfrac{12}{14-\dfrac{60}{23-\dfrac{144}{32-\dfrac{264}{41-\dfrac{420}{50-\ddots}}}}}}$$
Convergence type $E$ with $E=2$, $P=4/3$, and $C=2^{2/3}/(3\G(2/3))$, so that
$$\root3\of2-\dfrac{p(n)}{q(n)}\sim\dfrac{2^{2/3}/(3\G(2/3))}{2^nn^{4/3}}\;.$$
$$A=1-(22/9)/n+(679/81)/n^2-(88130/2187)/n^3+(4957862/19683)/n^4+\cdots$$
Series:
$$\dfrac{1}{\root3\of2}=\sum_{n\ge0}\dfrac{(-1/3)_n}{n!}2^{-n}$$
Parametric families, with $u\ge0$, $v\ge0$ with $v\not\equiv2\pmod3$, and $k\ge0$:
\begin{verbatim}
[()->2^(1/3),9*n-5+3*k+v+6*u,-6*(n+u)*(3*n+v-2)]
\end{verbatim}
Convergence type $E$ with $E=2$, $P=u+(5-v)/3+2k$.

After contraction, Ap\'ery accelerates to \ref{5.2.16}
\end{cf}

\smallskip

\begin{cf}\label{5.2.6}{\ }
\begin{verbatim}
[()->2^(1/3),[1,5,2*(4*n-1)],[2,(4*n+1)*(12*n-1)]]
\end{verbatim}
$$\root3\of2=1+\dfrac{2}{5+\dfrac{55}{14+\dfrac{207}{22+\dfrac{455}{30+\dfrac{799}{38+\dfrac{1239}{46+\ddots}}}}}}$$
Convergence type $E$ with $E=-3$, $P=1/6$, and $C=3^{1/4}\G(1/4)/\G(1/12)$,
so that
$$\root3\of2-\dfrac{p(n)}{q(n)}\sim(-1)^n\dfrac{\G(1/4)/\G(1/12)}{3^{n-1/4}n^{1/6}}\;.$$
$$A=1-(61/288)/n+(5159/55296)/n^2+\cdots$$
Parametric family:
\begin{verbatim}
[()->2^(1/3),8*n+(16*k-2),(4*n+1)*(12*n-1)]
\end{verbatim}
Convergence type $E$ with $E=-3$ and $P=2k+1/6$.

There are more general parametric families of the type
\begin{verbatim}
[()->2^(1/3),2*(4*n+2*u-1),(4*n+2*v-1)*(12*n+2*w-1)]
\end{verbatim}
with convergence type $E$ with $E=-3$ and $P=(6u-3v-w+5)/12$
but satisfying complicated conditions on {\tt (u,v,w)}.
\end{cf}
        
\smallskip

\begin{cf}\label{5.2.7}{\ }
\begin{verbatim}
[()->2^(1/3),[2,3,10*n-8],[-2,-(4*n-1)*(4*n-3)]]
\end{verbatim}
$$\root3\of2=2-\dfrac{2}{3-\dfrac{3}{12-\dfrac{35}{22-\dfrac{99}{32-\dfrac{195}{42-\dfrac{323}{52-\ddots}}}}}}$$
Convergence type $E$ with $E=4$, $P=2/3$, and
$C=-\G(1/3)^2/(2^{3/2}3^{2/3}\pi)$, so that
$$\root3\of2-\dfrac{p(n)}{q(n)}\sim-\dfrac{\G(1/3)^2/(2^{3/2}3^{2/3}\pi)}{2^{2n}n^{2/3}}\;.$$
$$A=1-(215/432)/n+(215185/373248)/n^2-(517261175/483729408)/n^3+\cdots$$
Parametric families, with $u\ge0$, $v\ge0$, $t\ge0$ with $t\not\equiv(u+v+2)\pmod{3}$, and $k\ge0$:
\begin{verbatim}
[()->2^(1/3),10*n-8+2*t+6*k,-(4*n-3+4*u)*(4*n-1+4*v)]
\end{verbatim}
Convergence type $E$ with $E=4$ and $P=(2t+2-5(u+v))/3+2k$.

Ap\'ery accelerates to \ref{5.2.15}.
\end{cf}

\smallskip

\begin{cf}\label{5.2.8}{\ }
\begin{verbatim}
[()->2^(1/3),[2,11,10*(3*n-2)],[-6,-(12*n-5)*(12*n+1)]]
\end{verbatim}
$$\root3\of2=2-\dfrac{6}{11-\dfrac{91}{40-\dfrac{475}{70-\dfrac{1147}{100-\dfrac{2107}{130-\dfrac{3355}{160-\ddots}}}}}}$$
Convergence type $E$ with $E=4$, $P=0$, and $C=-3^{1/2}/2^{5/6}$, so that
$$\root3\of2-\dfrac{p(n)}{q(n)}\sim-\dfrac{3^{1/2}}{2^{2n+5/6}}\;.$$
$$A=1-(5/16)/n+(235/1536)/n^2-(25735/73728)/n^3+\cdots$$
Parametric family for $k\ge0$:
\begin{verbatim}
[()->2^(1/3),30*n-20+18*k,-(12*n-5)*(12*n+1)]
\end{verbatim}
Convergence type $E$ with $E=4$ and $P=2k$.
\end{cf}

\smallskip

\begin{cf}\label{5.2.9}{\ }
\begin{verbatim}
[()->2^(1/3),[2,16,21*n-7],[-12,4*(3*n-1)*(6*n-5)]]
\end{verbatim}
$$\root3\of2=2-\dfrac{12}{16+\dfrac{8}{35+\dfrac{140}{56+\dfrac{416}{77+\dfrac{836}{98+\dfrac{1400}{119+\ddots}}}}}}$$
Convergence type $E$ with $E=-8$, $P=7/6$, and $C=-(2^{8/3}/3^{10/3})\sqrt{\pi}/\G(1/3)$, so that
$$\root3\of2-\dfrac{p(n)}{q(n)}\sim(-1)^{n+1}\dfrac{\sqrt{\pi}/(3^{10/3}\G(1/3))}{2^{3n-8/3}n^{7/6}}$$
$$A=1-(203/648)/n+(5257/839808)/n^2+\cdots$$
Parametric family for $k\ge0$:
\begin{verbatim}
[()->2^(1/3),21*n+27*k-7,4*(3*n-1)*(6*n-5)]
\end{verbatim}
Convergence type $E$ with $E=-8$ and $P=2k+7/6$.
\end{cf}

\smallskip

\begin{cf}\label{5.2.10}{\ }
\begin{verbatim}
[()->2^(1/3),[[0,5],[2*n,12*n+6]],
             [[6,-1/3],[-(n+1)*(3*n+2),-n*(3*n+1)]]]
\end{verbatim}
$$\root3\of2=\dfrac{6}{5-\dfrac{1/3}{2-\dfrac{10}{18-\dfrac{4}{4-\dfrac{24}{30-\dfrac{14}{6-\ddots}}}}}}$$
Convergence type $E$ with $E=(1+\sqrt{2})^2$, $P=0$, and $C=2^{1/3}3^{1/2}/(1+\sqrt{2})^2$, so that
$$\root3\of2-\dfrac{p(n)}{q(n)}\sim\dfrac{2^{1/3}3^{1/2}}{(1+\sqrt{2})^{2n+2}}\;.$$
$$A=1+(9d/8)/n+(-9d/8+81/64)/n^2+(67105d/41472-81/32)/n^3+\cdots$$
\end{cf}

\smallskip

\begin{cf}\label{5.2.11}{\ }
\begin{verbatim}
[()->2^(1/3),[2,6,4*(2*n-1)],[[-5,10],[(6*n-2)*(6*n+1),(6*n+2)*(6*n+5)]]]
\end{verbatim}
$$\root3\of2=2-\dfrac{5}{6+\dfrac{10}{12+\dfrac{28}{20+\dfrac{88}{28+\dfrac{130}{36+\dfrac{238}{44+\ddots}}}}}}$$
Convergence type $E$ with $E=-9$, $P=0$, and $C=-2^{1/3}/3^{1/2}$, so that
$$\root3\of2-\dfrac{p(n)}{q(n)}\sim(-1)^{n+1}\dfrac{2^{1/3}}{3^{2n+1/2}}\;.$$
$$A=1-(25/144)/n+(4225/41472)/n^2+\cdots$$
\end{cf}
        
\smallskip

\begin{cf}\label{5.2.12}{\ }
\begin{verbatim}
[()->2^(1/3),[1,10*n-6],[1,-(3*n-2)*(3*n-1)]]
\end{verbatim}
$$\root3\of2=1+\dfrac{1}{4-\dfrac{2}{14-\dfrac{20}{24-\dfrac{56}{34-\dfrac{110}{44-\dfrac{182}{54-\ddots}}}}}}$$
Convergence type $E$ with $E=9$, $P=1$, and $C=3^{1/2}2^{-11/3}$, so that
$$\root3\of2-\dfrac{p(n)}{q(n)}\sim\dfrac{2^{-11/3}}{3^{2n-1/2}n}\;.$$
$$A=1-(85/144)/n+(20995/41472)/n^2-(11698525/17915904)/n^3+\cdots$$
Parametric families with $u\ge0$, $v\ge0$:
\begin{verbatim}
[()->2^(1/3),10*n-6+u+v+8*k,-(3*n-2+3*u)*(3*n-1+3*v)]
\end{verbatim}
Convergence type $E$ with $E=9$ and $P=1-(u+v)+2k$.

Ap\'ery accelerates to \ref{5.2.17}.
\end{cf}

\smallskip

\begin{cf}\label{5.2.13}{\ }
\begin{verbatim}
[()->2^(1/3),[1,9,10*(2*n-1)],[2,-(6*n-1)*(6*n+1)]]
\end{verbatim}
$$\root3\of2=1+\dfrac{2}{9-\dfrac{35}{30-\dfrac{143}{50-\dfrac{323}{70-\dfrac{575}{90-\dfrac{899}{110-\ddots}}}}}}$$
Convergence type $E$ with $E=9$, $P=0$, and $C=2^{1/3}/3$, so that
$$\root3\of2-\dfrac{p(n)}{q(n)}\sim\dfrac{2^{1/3}}{3^{2n+1}}\;.$$
$$A=1-(5/18)/n+(115/648)/n^2-(1595/8748)/n^3+\cdots$$
Parametric family for $k\ge0$:
\begin{verbatim}
[()->2^(1/3),20*n-10+16*k,-(6*n-1)*(6*n+1)]
\end{verbatim}
Convergence type $E$ with $E=9$ and $P=2k$.
\end{cf}

\smallskip

\begin{cf}\label{5.2.14}{\ }
\begin{verbatim}
[()->2^(1/3),[2,5,4*(3*n-2)],[-4,(2*n-1)*(6*n+1)]]
\end{verbatim}
$$\root3\of2=2-\dfrac{4}{5+\dfrac{7}{16+\dfrac{39}{28+\dfrac{95}{40+\dfrac{175}{52+\dfrac{279}{64+\ddots}}}}}}$$
Convergence type $E$ with $E=-(2+\sqrt{3})^2$, $P=0$, and
$C=-2^{1/3}3^{1/2}/(2+\sqrt{3})^{2/3}$, so that
$$\root3\of2-\dfrac{p(n)}{q(n)}\sim(-1)^{n+1}\dfrac{2^{1/3}3^{1/2}}{(2+\sqrt{3})^{2n+2/3}}\;.$$
$$A=1-(5d/72)/n+(5d/216+25/3456)/n^2+(1945d/746496-25/5184)/n^3+\cdots$$
Parametric family with $u\ge0$, $u\not\equiv2\pmod3$
\begin{verbatim}
[()->2^(1/3),4*(3*n+u-2),(2*n-1)*(6*n+4*u+1)]
\end{verbatim}
Convergence type $E$ with $E=-(2+\sqrt{3})^2$ and $P=0$.
\end{cf}

\smallskip

\begin{cf}\label{5.2.15}{\ }
\begin{verbatim}
[()->2^(1/3),[2,11,8*(3*n-2)],
            [[-6,-91],[-9*(4*n-3)*(4*n-1),-(12*n+7)*(12*n+13)]]]
\end{verbatim}
$$\root3\of2=2-\dfrac{6}{11-\dfrac{91}{32-\dfrac{27}{56-\dfrac{475}{80-\dfrac{315}{104-\dfrac{1147}{128-\ddots}}}}}}$$
Convergence type $E$ with $E=(2+\sqrt{3})^2$, $P=0$, and $C=-2^{1/3}3^{1/2}/(2+\sqrt{3})^{2/3}$, so that
$$\root3\of2-\dfrac{p(n)}{q(n)}\sim-\dfrac{2^{1/3}3^{1/2}}{(2+\sqrt{3})^{2n+2/3}}\;.$$
$$A=1-(5d/72)/n+(5d/216+25/3456)/n^2+(1945d/746496-25/5184)/n^3+\cdots$$
Parametric families are of the form
\begin{verbatim}
[()->2^(1/3),8*(3*n+u0),[-9*(4*n+u1)*(4*n+u2),-(12*n+u3)*(12*n+u4)]]
\end{verbatim}
with $3\nmid u0$, $u1$ and $u2$ odd, $u1+u2\equiv0\pmod4$,
$u4-u3\equiv6\pmod{12}$ with $\gcd(u3,6)=\gcd(u4,6)=1$, and
$(u1-u2)=\pm(u3-u4)/3$.

Convergence type $E$ with $E=\pm(2+\sqrt{3})^2$ and $P=0$.
\end{cf}

\smallskip

\begin{cf}\label{5.2.16}{\ }
\begin{verbatim}
[()->2^(1/3),[1,8,9*(2*n-1)],[2,-(9*n^2-1)]]
\end{verbatim}
$$\root3\of2=1+\dfrac{2}{8-\dfrac{8}{27-\dfrac{35}{45-\dfrac{80}{63-\dfrac{143}{81-\dfrac{224}{99-\ddots}}}}}}$$
Convergence type $E$ with $E=(1+\sqrt{2})^4$, $P=0$, and $C=2^{1/3}3^{1/2}/(1+\sqrt{2})^2$, so that
$$\root3\of2-\dfrac{p(n)}{q(n)}\sim\dfrac{2^{1/3}3^{1/2}}{(1+\sqrt{2})^{4n+2}}\;.$$
$$A=1-(5d/48)/n+(5d/96+25/2304)/n^2+(-10175d/331776-25/2304)/n^3+\cdots$$
Parametric family with $3\nmid u$:
\begin{verbatim}
[()->2^(1/3),9*(2*n-1),-(9*n^2-u^2)]
\end{verbatim}
Convergence type $E$ with $E=(1+\sqrt{2})^4$, $P=0$.
\end{cf}

\smallskip

\begin{cf}\label{5.2.17}{\ }
\begin{verbatim}
[()->2^(1/3),[1,9*(2*n-1)],
            [[2,-35],[-(6*n-4)*(6*n-2),-(6*n+5)*(6*n+7)]]]
\end{verbatim}
$$\root3\of2=1+\dfrac{2}{9-\dfrac{35}{27-\dfrac{8}{45-\dfrac{143}{63-\dfrac{80}{81-\dfrac{323}{99-\ddots}}}}}}$$
Convergence type $E$ with $E=(1+\sqrt{2})^4$, $P=0$, and $C=2^{1/3}3^{1/2}/(1+\sqrt{2})^2$, so that
$$\root3\of2-\dfrac{p(n)}{q(n)}\sim\dfrac{2^{1/3}3^{1/2}}{(1+\sqrt{2})^{4n+2}}\;.$$
$$A=1-(5d/48)/n+(5d/96+25/2304)/n^2+(-10175d/331776-25/2304)/n^3+\cdots$$
\end{cf}

\smallskip

\begin{cf}\label{5.2.18}{\ }
\begin{verbatim}
[()->2^(1/3),[5/4,125,122*n+4],[5/4,125*n*(3*n-1)]]
\end{verbatim}
$$\root3\of2=5/4+\dfrac{5/4}{125+\dfrac{250}{248+\dfrac{1250}{370+\dfrac{3000}{492+\dfrac{5500}{614+\dfrac{8750}{736+\ddots}}}}}}$$
Convergence type $E$ with $E=-125/3$, $P=4/3$, and $C=5/(512\G(2/3))$, so that
$$\root3\of2-\dfrac{p(n)}{q(n)}\sim(-1)^n\dfrac{5/(512\G(2/3))}{(125/3)^nn^{4/3}}\;.$$
$$A=1-(311/288)/n+(341047/331776)/n^2-(516559295/573308928)/n^3+\cdots$$
Series:
$$\root3\of2=\dfrac{5}{4}\sum_{n\ge0}\dfrac{(-1/3)_n}{n!}(-3/125)^n$$
Parametric families with $u\ge0$ with $u\not\equiv2\pmod3$ and $k\ge0$:
\begin{verbatim}
[()->2^(1/3),122*n+5-u+128*k,125*n*(3*n+u-2)]
\end{verbatim}
Convergence type $E$ with $E=-125/3$ and $P=(5-u)/3+2k$.
\end{cf}

\smallskip

\begin{cf}\label{5.2.19}{\ }
\begin{verbatim}
[()->2^(1/3),[5/4,131*n-4],[5/4,-128*n*(3*n-1)]]
\end{verbatim}
$$\root3\of2=5/4+\dfrac{5/4}{127-\dfrac{256}{258-\dfrac{1280}{389-\dfrac{3072}{520-\dfrac{5632}{651-\dfrac{8960}{782-\ddots}}}}}}$$
Convergence type $E$ with $E=128/3$, $P=4/3$, and $C=2^{8/3}/(625\G(2/3))$,
so that
$$\root3\of2-\dfrac{p(n)}{q(n)}\sim\dfrac{2^{8/3}/(625\G(2/3))}{(128/3)^nn^{4/3}}\;.$$
$$A=1-(1286/1125)/n+(1550143/1265625)/n^2-(1167321722/854296875)/n^3+\cdots$$
Series:
$$\dfrac{1}{\root3\of2}=\dfrac{4}{5}\sum_{n\ge0}\dfrac{(-1/3)_n}{n!}(3/128)^n$$
Parametric families with $u\ge0$ with $u\not\equiv2\pmod3$ and $k\ge0$:
\begin{verbatim}
[()->2^(1/3),131*n+u-5+125*k,-128*n*(3*n+u-2)]
\end{verbatim}
Convergence type $E$ with $E=128/3$ and $P=(5-u)/3+2k$.
\end{cf}

\smallskip

\begin{cf}\label{5.2.20}{\ }
\begin{verbatim}
[()->2^(1/3),[5/4,252,253*(2*n-1)],[5/2,-(9*n^2-1)]]
\end{verbatim}
$$\root3\of2=\dfrac{5}{4}+\dfrac{5/2}{252-\dfrac{8}{759-\dfrac{35}{1265-\dfrac{80}{1771-\dfrac{143}{2277-\dfrac{224}{2783-\ddots}}}}}}$$
Convergence type $E$ with $E=(16+5\sqrt{10})^4/36$, $P=0$, and
$C=2^{4/3}3^{3/2}/(16+5\sqrt{10})^2$, so that
$$\root3\of2-\dfrac{p(n)}{q(n)}\sim\dfrac{2^{4/3}3^{3/2}}{(16+5\sqrt{10})^{4n+2}6^{-2n}}\;.$$
$$A=1-(253d/5760)/n+(253d/11520+64009/6635520)/n^2+\cdots$$
Parametric family for $3\nmid k$:
\begin{verbatim}
[()->2^(1/3),253*(2*n-1),-(9*n^2-k^2)]
\end{verbatim}
Convergence type $E$ with $E=(16+5\sqrt{10})^4/36$ and $P=0$.
\end{cf}

\smallskip

Note that some of these continued fractions are simply samples coming from the
expansions of the power function.

In particular, the last one comes from the identity $2^7=5^3+3$ applied
to \ref{6.5.1}, and it converges so rapidly that
a study of the arithmetic properties of $p(n)$ and $q(n)$ shows the
highly nontrivial result that there exists an \emph{effective} constant $C>0$ such that
$$\left|\root3\of2-\dfrac{p}{q}\right|>\dfrac{C}{q^{2.9471}}$$
for all integers $p$ and $q$ (much better effective bounds are known, for instance
in \cite{Ben}, it is shown that $|\root3\of2-p/q|>(1/4)/q^{2.47}$).

\smallskip

As mentioned at the beginning of this section, we have given continued
fractions for $\root3\of2$ as examples, but of course
we could do the same for other algebraic numbers. For instance, it is
easy to show that by using the CFs for the power function, any constant
of the form $u^a$ for $u$ and $a$ rational has a CF of the above type.
We give two examples of this, which are continued fractions which converge
extremely fast.

\medskip

\begin{cf}\label{5.2.21}{\ }
\begin{verbatim}
[()->3^(1/3),[13/9,6581,6576*(2*n-1)],[-130/9,-25*(9*n^2-1)]]
\end{verbatim}
$$\root3\of3=\dfrac{13}{9}-\dfrac{130/9}{6581-\dfrac{200}{19728-\dfrac{875}{32880-\dfrac{2000}{46032-\dfrac{3575}{59184-\dfrac{5600}{72336-\ddots}}}}}}$$
Convergence type $E$ with $E=(4\sqrt{39}+25)^4(8-\sqrt{39})^2/25$, $P=0$, and
$C=-3^{5/6}5/((4\sqrt{39}+25)^2(8-\sqrt{39}))$, so that
$$3^{1/3}-\dfrac{p(n)}{q(n)}\sim-\dfrac{3^{5/6}}{((4\sqrt{39}+25)^2(8-\sqrt{39}))^{2n+1}}$$
$$A=1-(2740d/123201)/n+(1370d/123201+3753800/389191959)/n^2+\cdots$$
Parametric family for $3\nmid k$:
\begin{verbatim}
[()->3^(1/3),6576*(2*n-1),-25*(9*n^2-k^2)]
\end{verbatim}
Convergence type $E$ with $E=(4\sqrt{39}+25)^4(8-\sqrt{39})^2/25$ and $P=0$.
\end{cf}

\smallskip

\begin{cf}\label{5.2.22}{\ }
\begin{verbatim}
[()->5^(1/4),[3/2,645,644*(2*n-1)],[-3,-(16*n^2-1)]]
\end{verbatim}
$$\root4\of5=3/2-\dfrac{3}{645-\dfrac{15}{1932-\dfrac{63}{3220-\dfrac{143}{4508-\dfrac{255}{5796-\dfrac{399}{7084-\ddots}}}}}}$$
Convergence type $E$ with $E=((1+\sqrt{5})/2)^{24}$, $P=0$, and
$C=-2^{1/2}5^{1/4}/((1+sqrt(5))/2)^12$, so that
$$\root4\of5-\dfrac{p(n)}{q(n)}\sim-\dfrac{2^{1/2}5^{1/4}}{((1+\sqrt{5})/2)^{24n+12}}$$
$$A=1-(161d/1920)/n+(161d/3840+25921/1474560)/n^2+\cdots$$
Parametric family for $2\nmid k$:
\begin{verbatim}
[()->5^(1/4),644*(2*n-1),-(16*n^2-k^2)]
\end{verbatim}
Convergence type $E$ with $E=((1+\sqrt{5})/2)^{24}$ and $P=0$.
\end{cf}

\medskip

\section{Constants: $\pi$, $\log(2)$, and Periods of Degree $1$}

\medskip

{\bf Preliminaries on CFs for \emph{periods}:} First, we define a
\emph{rational period} of degree $k$ as the sum of a convergent series of the
form $\sum_{n\ge1}\chi(n)f(n)$, where $\chi$ is a periodic arithmetic function
taking \emph{rational} values, and $f\in\Q(x)$ is a rational function with
rational coefficients, whose denominator is of degree $k$, or more
generally a product of coprime polynomials of degree at most $k$, one at
least having degree $k$.

\smallskip

Note that this definition is \emph{not} compatible with the definition
of periods as given by Kontsevitch--Zagier \cite{Kon-Zag}. Note also that
a $\Q$-linear combination of rational periods of degree $\le k$ is again
a rational period of degree $\le k$.

\smallskip

By using Euler's transformation of series seen above, it is clear that
a period has a CF expansion where $a(n)$ and $b(n)$ are polynomials with
rational coefficients for $n$ large. However this is usually uninteresting,
or at least not any more interesting than the series itself. Thus, really
interesting CFs are those for which the degrees of the polynomials are
reasonably small, specifically $\deg(a(n))\le k$ and $\deg(b(n))\le 2k$ for
$n$ sufficiently large (we will say that it has \emph{bidegree} $\le(k,2k)$).
Trivial example: $S=\sum_{n\ge1}1/P(n)$ with
$\deg(P)=k\ge2$ has the CF $S=[[0,P(n)+P(n-1)],[1,-P(n)^2]]$ by Euler, of
the required degrees. Note that applying Euler and simplifying on the sum
$S=\sum_{n\ge1}(n+1)/n^3=\z(3)+\z(2)$ which is a period of degree $3$ leads
to a CF of bidegree $(4,8)$, but there does exist CFs of bidegree $(3,6)$
for $S$, for instance $\z(3)+\z(2)=[[3,2n^3+2n^2+3n+1],[-1,-n(n+1)^5]]$,
coming from the identity $\z(3)+\z(2)=3-\sum_{n\ge1}1/(n(n+1)^3)$.

\smallskip

We will however also often give relatively interesting CFs of larger bidegrees,
and also of period $2$.

\medskip

\begin{cf}\label{5.3.1}{\ }
\begin{verbatim}
[()->Pi,[2,6,8*n^2-2*n+1],[2,-2*n*(2*n+1)^3]]
\end{verbatim}
$$\pi=2+\dfrac{2}{6-\dfrac{54}{29-\dfrac{500}{67-\dfrac{2058}{121-\dfrac{5832}{191-\dfrac{13310}{277-\ddots}}}}}}$$
Convergence type $P^+$ with $P=1/2$ and $C=2/\sqrt{\pi}$, so that
$$\pi-\dfrac{p(n)}{q(n)}\sim\dfrac{2/\sqrt{\pi}}{n^{1/2}}\;.$$
$$A=1-(11/24)/n+(181/640)/n^2-(3599/21504)/n^3+\cdots$$
Series:
$$\pi=2+\sum_{n\ge0}\dfrac{(3/2)_n}{(2n+3)(n+1)!}$$
Parametric family for $k\ge0$:
\begin{verbatim}
[()->Pi,8*n^2-2*n+1+2*k*(2*k+1),-2*n*(2*n+1)^3]
\end{verbatim}
Convergence type $P^+$ with $P=2k+1/2$.
\end{cf}

\smallskip

\begin{cf}\label{5.3.2}{\ }
\begin{verbatim}
[()->Pi,[2,(4*n-1)*(4*n^2-2*n-1)],
        [10,8*n^3*(2*n+1)^3*(4*n-3)*(4*n+5)]]
\end{verbatim}
$$\pi=2+\dfrac{10}{3+\dfrac{1944}{77+\dfrac{520000}{319+\dfrac{11335464}{825+\dfrac{101896704}{1691+\dfrac{565675000}{3013+\ddots}}}}}}$$
Convergence type $P^-$ with $P=1/2$ and $C=\sqrt{\pi}$, so that
$$\pi-\dfrac{p(n)}{q(n)}\sim(-1)^n\dfrac{\sqrt{\pi}}{n^{1/2}}$$
$$A=1-(3/8)/n+(9/128)/n^2+(135/1024)/n^3+\cdots$$
Series:
$$\dfrac{1}{\pi}=\dfrac{1}{2}\sum_{n\ge0}(-1)^n\dfrac{(4n+1)(1/2)_n^3}{n!^3}$$
\end{cf}

\smallskip

\begin{cf}\label{5.3.3}{\ }
\begin{verbatim}
[()->Pi,[0,1,2],[4,(2*n-1)^2]]
\end{verbatim}
$$\pi=\dfrac{4}{1+\dfrac{1}{2+\dfrac{9}{2+\dfrac{25}{2+\dfrac{49}{2+\dfrac{81}{2+\ddots}}}}}}$$
Convergence type $P^-$ with $P=1$ and $C=1$, so that
$$\pi-\dfrac{p(n)}{q(n)}\sim\dfrac{(-1)^n}{n}\;.$$
$$A=1-(1/4)/n^2+(5/16)/n^4-(61/64)/n^6+(1385/256)/n^8-\cdots$$
Series:
$$\pi=4\sum_{n\ge1}\dfrac{(-1)^{n+1}}{2n-1}$$
Parametric families: up to the trivial change $n\mapsto n+j$, with
$k\ge0$ and $0\le u\le k/2$:
\begin{verbatim}
[()->Pi,2*(k+1),(2*n-1)*(2*n+2*k-1-4*u)]
\end{verbatim}
Convergence type $P^-$ with $P=k+1$.
\end{cf}

\smallskip

\begin{cf}\label{5.3.4}{\ }
\begin{verbatim}
[()->Pi,[4,8*n^2+1],[-4,-4*n*(n+1)*(2*n+1)^2]]
\end{verbatim}
$$\pi=4-\dfrac{4}{9-\dfrac{72}{33-\dfrac{600}{73-\dfrac{2352}{129-\dfrac{6480}{201-\dfrac{14520}{289-\ddots}}}}}}$$
Convergence type $P^+$ with $P=1$ and $C=-\pi/4$, so that
$$\pi-\dfrac{p(n)}{q(n)}\sim-\dfrac{\pi/4}{n}\;.$$
$$A=1-(7/8)/n+(23/32)/n^2-(277/512)/n^3+\cdots$$
Series:
\begin{align*}
  \pi&=4-\dfrac{4}{9}\sum_{n\ge0}\dfrac{n!(n+1)!}{(5/2)_n^2}\\
  \dfrac{1}{\pi}&=\dfrac{1}{4}\sum_{n\ge0}\dfrac{(1/2)_n^2}{n!(n+1)!}\end{align*}
Parametric family for $k\ge0$:
\begin{verbatim}
[()->Pi,8*n^2+(2*k+1)^2,-4*n*(n+1)*(2*n+1)^2]
\end{verbatim}
Convergence type $P^+$ with $P=2k+1$.
\end{cf}

\smallskip

\begin{cf}\label{5.3.5}{\ }
\begin{verbatim}
[()->Pi,[2,8*n^2-8*n+3],[2,-4*n^2*(4*n^2-1)]]
\end{verbatim}
$$\pi=2+\dfrac{2}{3-\dfrac{12}{19-\dfrac{240}{51-\dfrac{1260}{99-\dfrac{4032}{163-\dfrac{9900}{243-\ddots}}}}}}$$
Convergence type $P^+$ with $P=1$ and $C=\pi/4$, so that
$$\pi-\dfrac{p(n)}{q(n)}\sim\dfrac{\pi/4}{n}$$
Series:
\begin{align*}
  \pi&=2+2\sum_{n\ge0}\dfrac{n!^2}{(2n+3)(3/2)_n^2}\\
  \dfrac{1}{\pi}&=-\dfrac{1}{2}\sum_{n\ge0}\dfrac{(1/2)_n^2}{(2n-1)n!^2}
\end{align*}
Parametric family for $k\ge0$:
\begin{verbatim}
[()->Pi,8*n^2-8*n+2+(2*k+1)^2,-4*n^2*(4*n^2-1)]
\end{verbatim}
Convergence type $P^+$ with $P=2k+1$.
\end{cf}

\smallskip

\begin{cf}\label{5.3.6}{\ }
\begin{verbatim}
[()->Pi,[3,72*n^2-72*n+35],[3,-36*n^2*(36*n^2-1)]]
\end{verbatim}
$$\pi=3+\dfrac{3}{35-\dfrac{1260}{179-\dfrac{20592}{467-\dfrac{104652}{899-\dfrac{331200}{1475-\dfrac{809100}{2195-\ddots}}}}}}$$
Convergence type $P^+$ with $P=1$ and $C=\pi/36$, so that
$$\pi-\dfrac{p(n)}{q(n)}\sim\dfrac{\pi/36}{n}$$
Series:
\begin{align*}
  \pi&=3+\dfrac{3}{35}\sum_{n\ge0}\dfrac{n!^2}{(11/6)_n(13/6)_n}\\
  \dfrac{1}{\pi}&=\dfrac{1}{3}\sum_{n\ge0}\dfrac{(-1/6)_n(1/6)_n}{n!^2}
\end{align*}
Parametric family for $k\ge0$:
\begin{verbatim}
[()->Pi,72*n^2-72*n+35+36*k*(k+1),-36*n^2*(36*n^2-1)]
\end{verbatim}
Convergence type $P^+$ with $P=2k+1$.
\end{cf}

\smallskip

\begin{cf}\label{5.3.7}{\ }
\begin{verbatim}
[()->Pi,[4,2],[-2,n^2]]
\end{verbatim}
$$\pi=4-\dfrac{2}{2+\dfrac{1}{2+\dfrac{4}{2+\dfrac{9}{2+\dfrac{16}{2+\dfrac{25}{2+\ddots}}}}}}$$
Convergence type $P^-$ with $P=2$ and $C=-\pi/8$, so that
$$\pi-\dfrac{p(n)}{q(n)}\sim(-1)^{n+1}\dfrac{\pi/8}{n^2}\;.$$
$$A=1-1/n+(1/8)/n^2+(3/4)/n^3+(5/256)/n^4-(623/256)/n^5+\cdots$$
Series:
\begin{align*}\pi&=4-4\sum_{n\ge0}\dfrac{n!^2}{(4n+1)(4n+5)(3/2)_n^2}\\
\dfrac{1}{\pi}&=-\sum_{n\ge0}\dfrac{(1/2)_n^2}{(4n-1)(4n+3)n!^2}\end{align*}
Parametric families up to the trivial change $n\mapsto n+j$, with $k\ge0$ and
$k+2u\ge1$:
\begin{verbatim}
[()->Pi,k+1,n*(n+k-1+2*u)]
\end{verbatim}
Convergence type $P^-$ with $P=k+1$.
\end{cf}

\smallskip

\begin{cf}\label{5.3.8}{\ }
\begin{verbatim}
[()->Pi,[2,3,4],[4,4*n^2-1]]
\end{verbatim}
$$\pi=2+\dfrac{4}{3+\dfrac{3}{4+\dfrac{15}{4+\dfrac{35}{4+\dfrac{63}{4+\dfrac{99}{4+\ddots}}}}}}$$
Convergence type $P^-$ with $P=2$ and $C=1/2$, so that
$$\pi-\dfrac{p(n)}{q(n)}\sim(-1)^n\dfrac{1/2}{n^2}\;.$$
$$A=1-1/n+(1/4)/n^2+(1/2)/n^3+(1/16)/n^4-(31/16)/n^5+(1/64)/n^6+\cdots$$
Series:
$$\pi=2+4\sum_{n\ge1}\dfrac{(-1)^{n+1}}{(2n-1)(2n+1)}$$
Parametric family for $k\ge0$:
\begin{verbatim}
[()->Pi,4*(k+1),4*n^2-1]
\end{verbatim}
Convergence type $P^-$ with $P=2k+2$.
\end{cf}

\smallskip

\begin{cf}\label{5.3.9}{\ }
\begin{verbatim}
[()->Pi,[3,35,36*(2*n-1)],[6,(36*n^2-1)^2]]
\end{verbatim}
$$\pi=3+\dfrac{6}{35+\dfrac{1225}{108+\dfrac{20449}{180+\dfrac{104329}{252+\dfrac{330625}{324+\dfrac{808201}{396+\ddots}}}}}}$$
Convergence type $P^-$ with $P=2$ and $C=1/12$, so that
$$\pi-\dfrac{p(n)}{q(n)}\sim(-1)^n\dfrac{1/12}{n^2}$$
$$A=1-1/n+(1/36)/n^2+(17/18)/n^3+(1/1296)/n^4+\cdots$$
Series:
$$\pi=3+6\sum_{n\ge0}\dfrac{(-1)^n}{(6n+5)(6n+7)}$$
Parametric family for $k\ge0$:
\begin{verbatim}
[()->Pi,36*(2*k+1)*(2*n-1),(36*n^2-1)^2]
\end{verbatim}
Convergence type $P^-$ with $P=4k+2$.
\end{cf}

\smallskip

\begin{cf}\label{5.3.10}{\ }
\begin{verbatim}
[()->Pi,[0,3*n-2],[2,-n*(2*n-1)]]
\end{verbatim}
$$\pi=\dfrac{2}{1-\dfrac{1}{4-\dfrac{6}{7-\dfrac{15}{10-\dfrac{28}{13-\dfrac{45}{16-\ddots}}}}}}$$
Convergence type $E$ with $E=2$, $P=1/2$, and $C=2\sqrt{\pi}$, so that
$$\pi-\dfrac{p(n)}{q(n)}\sim\dfrac{\sqrt{\pi}}{2^{n-1}n^{1/2}}\;.$$
$$A=1-(7/8)/n+(241/128)/n^2-(6925/1024)/n^3+(1118859/32768)/n^4+\cdots$$
Series:
$$\pi=2\sum_{n\ge0}\dfrac{n!}{(3/2)_n}2^{-n}$$
Parametric families with $u\ge0$, $v\ge0$, and $k\ge0$:
\begin{verbatim}
[()->Pi,3*n+u+2*v+k-2,-(n+u)*(2*n+2*v-1)]
\end{verbatim}
Convergence type $E$ with $E=2$ and $P=v-u+1/2+2k$.
\end{cf}

\smallskip

\begin{cf}\label{5.3.11}{\ }
\begin{verbatim}
[()->Pi,[3,24,20*n^2+4*n+1],[3,-8*n*(2*n+1)^3]]
\end{verbatim}
$$\pi=3+\dfrac{3}{24-\dfrac{216}{89-\dfrac{2000}{193-\dfrac{8232}{337-\dfrac{23328}{521-\dfrac{53240}{745-\ddots}}}}}}$$
Convergence type $E$ with $E=4$, $P=3/2$, and $C=1/2/\sqrt{\pi}$, so that
$$\pi-\dfrac{p(n)}{q(n)}\sim\dfrac{1/\sqrt{\pi}}{2^{2n+1}n^{3/2}}\;.$$
$$A=1-(21/8)/n+(841/128)/n^2-(18863/1024)/n^3+\cdots$$
Series:
$$\pi=3+\dfrac{3}{8}\sum_{n\ge0}\dfrac{(3/2)_n}{(2n+3)(n+1)!}2^{-2n}$$
Parametric family for $k\ge0$:
\begin{verbatim}
[()->Pi,20*n^2+(8*k+4)*n+1,-4*(2*n-k)*(2*n+1-k)*(2*n+1)^2]
\end{verbatim}
Convergence type $E$ with $E=4$ and $P=3k+3/2$.
\end{cf}

\smallskip

\begin{cf}\label{5.3.12}{\ }
\begin{verbatim}
[()->Pi,[0,2*n-1],[4,n^2]]
[()->Pi,[4,4,2*n+1],[-4,(n+1)^2]]
\end{verbatim}
$$\pi=\dfrac{4}{1+\dfrac{1}{3+\dfrac{4}{5+\dfrac{9}{7+\dfrac{16}{9+\dfrac{25}{11+\ddots}}}}}}=2-\dfrac{2}{4+\dfrac{4}{5+\dfrac{9}{7+\dfrac{16}{9+\dfrac{25}{11+\dfrac{36}{13+\ddots}}}}}}$$
Convergence type $E$ with $E=-(1+\sqrt{2})^2$, $P=0$, and
$C=4\pi/(1+\sqrt{2})$, $C=-4\pi/(1+\sqrt{2})^3$, so that
$$\pi-\dfrac{p(n)}{q(n)}\sim(-1)^n\dfrac{4\pi}{(1+\sqrt{2})^{2n+1}}\text{\quad and\quad}\pi-\dfrac{p(n)}{q(n)}\sim(-1)^{n+1}\dfrac{4\pi}{(1+\sqrt{2})^{2n+3}}\;.$$
\begin{align*}A&=1-(d/8)/n+(d/16+1/64)/n^2+(3d/512-1/64)/n^3+\cdots\\
A&=1-(d/8)/n+(3d/16+1/64)/n^2+(-125d/512-3/64)/n^3+\cdots\end{align*}
Parametric families with $u\ge0$ and $v\ge0$:
\begin{verbatim}
[()->Pi,2*n+2*u+2*v-1,(n+u)*(n+u+2*v)]
\end{verbatim}
Convergence type $E$ with $E=-(1+\sqrt{2})^2$ and $P=0$.
\end{cf}

\smallskip

\begin{cf}\label{5.3.13}{\ }
\begin{verbatim}
[()->Pi,[2,2,2*n+1],[[4,8],[(2*n-1)*(2*n+1),4*(n+1)*(n+2)]]]
\end{verbatim}
$$\pi=2+\dfrac{4}{2+\dfrac{8}{5+\dfrac{3}{7+\dfrac{24}{9+\dfrac{15}{11+\dfrac{48}{13+\ddots}}}}}}$$
Convergence type $E$ with $E=-(1+\sqrt{2})^2$, $P=0$, and
$C=4\pi/(1+\sqrt{2})^3$, so that
$$\pi-\dfrac{p(n)}{q(n)}\sim(-1)^n\dfrac{4\pi}{(1+\sqrt{2})^{2n+3}}\;.$$
$$A=1+(3d/8)/n+(-9d/16+9/64)/n^2+(423d/512-27/64)/n^3+\cdots$$
\end{cf}

\smallskip

\begin{cf}\label{5.3.14}{\ }
\begin{verbatim}
[()->Pi,[4,2*n+1],[[-4,9],[(2*n)^2,(2*n+3)^2]]]
\end{verbatim}
$$\pi=4-\dfrac{4}{3+\dfrac{9}{5+\dfrac{4}{7+\dfrac{25}{9+\dfrac{16}{11+\dfrac{49}{13+\ddots}}}}}}$$
Convergence type $E$ with $E=-(1+\sqrt{2})^2$, $P=0$, and $C=-4\pi/(1+\sqrt{2})^3$, so that
$$\pi-\dfrac{p(n)}{q(n)}\sim(-1)^{n+1}\dfrac{4\pi}{(1+\sqrt{2})^{2n+3}}\;.$$
$$A=1-(d/8)/n+(3d/16+1/64)/n^2+(-125d/512-3/64)/n^3+\cdots$$
Parametric families with $u\ge0$, $0\le v\le 2u+3$:
\begin{verbatim}
[()->Pi,2*n+2*u-1,[(2*n+v-1)^2,(2*n+2*u+2-v)^2]]
\end{verbatim}
Convergence type $E$ with $E=-(1+\sqrt{2})^2$ and $P=0$.

Note that the special case $v=u+1$ gives part of the parametric family
given in \ref{5.3.12}.
\end{cf}

\smallskip

\begin{cf}\label{5.3.15}{\ }
\begin{verbatim}
[()->Pi,[2,2,2*n+1],[[2,-1],-(n+1)*(2*n+1)*[1,1]]]
\end{verbatim}
$$\pi=2+\dfrac{2}{2-\dfrac{1}{5-\dfrac{6}{7-\dfrac{6}{9-\dfrac{15}{11-\dfrac{15}{13-\ddots}}}}}}$$
Convergence type $E$ with $E=(1+\sqrt{2})^2$, $P=0$, and $C=4\pi/(1+\sqrt{2})^3$, so that
$$\pi-\dfrac{p(n)}{q(n)}\sim\dfrac{4\pi}{(1+\sqrt{2})^{2n+3}}\;.$$
$$A=1+(3d/8)/n+(-9d/16+9/64)/n^2+(423d/512-27/64)/n^3+\cdots$$
Parametric families up to the trivial change $n\mapsto n+j$,
with $u\ge0$, $v\ge0$:
\begin{verbatim}
[()->Pi,2*n+2*v+2*u-1,[-n*(2*n+2*u-1),-(n+v+1)*(2*n+2*u+2*v+1)]]
\end{verbatim}
Convergence type $E$ with $E=(1+\sqrt{2})^2$ and $P=0$.
\end{cf}

\smallskip

\begin{cf}\label{5.3.16}{\ }
\begin{verbatim}
[()->Pi,[[3,9],[10*n^2+7*n+1,10*n^2+19*n+9]],
              [[1,-36],[n*(n+1)*(2*n+1)^2,-4*(n+1)^2*(2*n+3)^2]]]
\end{verbatim}
$$\pi=3+\dfrac{1}{9-\dfrac{36}{18+\dfrac{18}{38-\dfrac{400}{55+\dfrac{150}{87-\dfrac{1764}{112+\ddots}}}}}}$$
Convergence type $E$ with $E=-i((1+\sqrt{5})/2)^5$, $P=0$, and $C=4\pi/((1+\sqrt{5})/2)^9$, so that
$$\pi-\dfrac{p(n)}{q(n)}\sim(-1)^n\dfrac{4\pi}{((1+\sqrt{5})/2)^{5n+9}}\;.$$
$$A=1-(d/10)/n+(4d/25+1/40)/n^2+(-2319d/10000-2/25)/n^3+\cdots$$
\end{cf}

\smallskip

\begin{cf}\label{5.3.17}{\ }
\begin{verbatim}
[()->Pi,[[2,3],[20*n^2-1,20*n^2+16*n+3]],
       [[4,9],[-16*n^2*(4*n^2-1),(4*(n+1)^2-1)^2]]]
\end{verbatim}
$$\pi=2+\dfrac{4}{3+\dfrac{9}{19-\dfrac{48}{39+\dfrac{225}{79-\dfrac{960}{115+\dfrac{1225}{179-\ddots}}}}}}$$
Convergence type $E$ with $E=i((1+\sqrt{5})/2)^5$, $P=0$, and
$C=4\pi/((1+\sqrt{5})/2)^5$, so that
$$\pi-\dfrac{p(n)}{q(n)}\sim(-1)^{\lfloor n/2\rfloor}\dfrac{4\pi}{((1+\sqrt{5})/2)^{5n+5}}\;.$$
$$A=1+(d/10)/n+(-d/10+1/40)/n^2+(711d/10000-1/20)/n^3+\cdots$$
\end{cf}

\smallskip

\begin{cf}\label{5.3.18}{\ }
\begin{verbatim}
[()->Pi,[[3,35],[180*n^2-1,180*n^2+144*n+35]],
       [[6,1225],[-144*n^2*(36*n^2-1),(6*n+5)^2*(6*n+7)^2]]]
\end{verbatim}
$$\pi=3+\dfrac{6}{35+\dfrac{1225}{179-\dfrac{5040}{359+\dfrac{20449}{719-\dfrac{82368}{1043+\dfrac{104329}{1619-\ddots}}}}}}$$
Convergence type $E$ with $E=i((1+\sqrt{5})/2)^5$, $P=0$, and
$C=\pi/((1+\sqrt{5})/2)^5$, so that
$$\pi-\dfrac{p(n)}{q(n)}\sim(-1)^{\lfloor n/2\rfloor}\dfrac{\pi}{((1+\sqrt{5})/2)^{5n+5}}\;.$$
$$A=1-(31d/90)/n+(31d/90+961/3240)/n^2+(-6134443d/21870000-961/1620)/n^3+\cdots$$
\end{cf}

It is possible to construct infinitely many continued fractions for $\pi$
of the same type with $E=-((1+\sqrt{5})/2)^5$.

\smallskip

\begin{cf}\label{5.3.19}{\ }
\begin{verbatim}
[()->Pi,[2,290*n^3-213*n^2-29*n+24],
        [78,-3*n*(2*n+1)*(3*n-1)*(3*n+1)*(10*n-7)*(10*n+13)]]
\end{verbatim}
$$\pi=2+\dfrac{78}{72-\dfrac{4968}{1434-\dfrac{450450}{5850-\dfrac{4984560}{15060-\dfrac{27011556}{30804-\ddots}}}}}$$
Convergence type $E$ with $E=27/2$, $P=1/2$, and $C=16\sqrt{3\pi}/45$, so that
$$\pi-\dfrac{p(n)}{q(n)}\sim\dfrac{16\sqrt{3\pi}/45}{(27/2)^nn^{1/2}}\;.$$
$$A=1-(857/1800)/n+(2455249/6480000)/n^2+\cdots$$
Series:
$$\pi=2+2\sum_{n\ge1}\dfrac{(10n+3)(1/2)_n(n-1)!}{(2/3)_n(4/3)_n}(2/27)^n$$
\end{cf}

\smallskip

\begin{cf}\label{5.3.20}{\ }
\begin{verbatim}
[()->Pi,[0,145*n^3-353*n^2+244*n-48],
        [-36,-3*n*(2*n-1)*(3*n-2)*(3*n-1)*(5*n-7)*(5*n+3)]]
\end{verbatim}
$$\pi=\dfrac{36}{12-\dfrac{96}{188-\dfrac{14040}{1422-\dfrac{362880}{4560-\dfrac{2762760}{10472-\dfrac{12383280}{20028-\ddots}}}}}}$$
Convergence type $E$ with $E=27/2$, $P=1/2$, and $C=4\sqrt{3\pi}/5$, so that
$$\pi-\dfrac{p(n)}{q(n)}\sim\dfrac{4\sqrt{3\pi}/5}{(27/2)^nn^{1/2}}\;.$$
$$A=1-(617/1800)/n+(1391569/6480000)/n^2+\cdots$$
Series:
$$\pi=\sum_{n\ge0}\dfrac{(5n+3)(1/2)_nn!}{(4/3)_n(5/3)_n}(2/27)^n$$
\end{cf}

\smallskip

\begin{cf}\label{5.3.21}{\ }
\begin{verbatim}
[()->Pi,[3,1456*n^3-1268*n^2-116*n+33],
        [15,6*n*(2*n-1)*(6*n-1)*(6*n+1)*(14*n-13)*(14*n+15)]]
\end{verbatim}
$$\pi=3+\dfrac{15}{105+\dfrac{6090}{6377+\dfrac{3320460}{27585+\dfrac{48052710}{72465+\dfrac{294919800}{149753+\ddots}}}}}$$
Convergence type $E$ with $E=-27$, $P=3/2$, and $C=\sqrt{\pi}/24$, so that
$$\pi-\dfrac{p(n)}{q(n)}\sim(-1)^n\dfrac{\sqrt{\pi}/24}{3^{3n}n^{3/2}}\;.$$
$$A=1-(37/36)/n+(9265/10368)/n^2-(379505/559872)/n^3+\cdots$$
Series:
$$\pi=3+\dfrac{1}{105}\sum_{n\ge0}(-1)^n\dfrac{(14n+15)(1/2)_nn!}{(11/6)_n(13/6)_n}3^{-3n}$$
\end{cf}

\smallskip

\begin{cf}\label{5.3.22}{\ }
\begin{verbatim}
[()->Pi,[47/15,1456*n^3+1644*n^2+380*n-15],
        [29,6*n*(2*n+1)*(6*n+1)*(6*n+5)*(14*n+1)*(14*n+29)]]
\end{verbatim}
$$\pi=\dfrac{47}{15}+\dfrac{29}{3465+\dfrac{893970}{18969+\dfrac{21918780}{55233+\dfrac{168104286}{120993+\dfrac{758727000}{224985+\ddots}}}}}$$
Convergence type $E$ with $E=-27$, $P=3/2$, and $C=\sqrt{\pi}/216$, so that
$$\pi-\dfrac{p(n)}{q(n)}\sim(-1)^n\dfrac{\sqrt{\pi}/216}{3^{3n}n^{3/2}}\;.$$
$$A=1-(49/36)/n+(13105/10368)/n^2-(242795/559872)/n^3+\cdots$$
Series:
$$\pi=\dfrac{47}{15}-\dfrac{1}{30}\sum_{n\ge1}(-1)^n\dfrac{(14n+15)(1/2)_n(n-1)!}{(7/6)_n(11/6)_n}3^{-3n}$$
\end{cf}

\smallskip

\begin{cf}\label{5.3.23}{\ }
\begin{verbatim}
[()->Pi,[2,3*(4*n-1)*(4*n^2-2*n-1)],
        [10,-n*(n+1)*(4*n^2-1)*(4*n-3)*(4*n+5)]]
\end{verbatim}
$$\pi=2+\dfrac{10}{9-\dfrac{54}{231-\dfrac{5850}{957-\dfrac{64260}{2475-\dfrac{343980}{5073-\dfrac{1262250}{9039-\ddots}}}}}}$$
Convergence type $E$ with $E=(1+\sqrt{2})^4$, $P=0$, and $C=4\pi/(1+\sqrt{2})^3$, so that
$$\pi-\dfrac{p(n)}{q(n)}\sim\dfrac{4\pi}{(1+\sqrt{2})^{4n+3}}\;.$$
$$A=1+(3d/16)/n+(-9d/64+9/256)/n^2+(423d/4096-27/512)/n^3+\cdots$$
\end{cf}

\smallskip

\begin{cf}\label{5.3.24}{\ }
\begin{verbatim}
[()->Pi,[0,4*(4*n-3)*(6*n^2-9*n+2)],
        [-12,-4*n^2*(2*n-1)^2*(4*n+3)*(4*n-5)]]
\end{verbatim}
$$\pi=\dfrac{12}{4-\dfrac{28}{160-\dfrac{4752}{1044-\dfrac{94500}{3224-\dfrac{655424}{7276-\dfrac{2794500}{13776-\ddots}}}}}}$$
Convergence type $E$ with $E=(1+\sqrt{2})^4$, $P=0$, and $C=4\pi/(1+\sqrt{2})$, so that
$$\pi-\dfrac{p(n)}{q(n)}\sim\dfrac{4\pi}{(1+\sqrt{2})^{4n+1}}\;.$$
$$A=1-(d/16)/n+(d/64+1/256)/n^2+(3d/4096-1/512)/n^3+\cdots$$
\end{cf}

\smallskip

\begin{cf}\label{5.3.25}{\ }
\begin{verbatim}
[()->Pi,[4*(4*n-1)*(6*n^2-3*n-1)],
        [-20,-4*n^2*(2*n+1)^2*(4*n-3)*(4*n+5)]]
\end{verbatim}
$$\pi=4-\dfrac{20}{24-\dfrac{324}{476-\dfrac{26000}{1936-\dfrac{269892}{4980-\dfrac{1415232}{10184-\dfrac{5142500}{18124-\ddots}}}}}}$$
Convergence type $E$ with $E=(1+\sqrt{2})^4$, $P=0$, and $C=-4\pi/(1+\sqrt{2})^3$, so that
$$\pi-\dfrac{p(n)}{q(n)}\sim-\dfrac{4\pi}{(1+\sqrt{2})^{4n+3}}\;.$$
$$A=1-(d/16)/n+(3d/64+1/256)/n^2+(-125d/4096-3/512)/n^3+\cdots$$
\end{cf}

\smallskip

\begin{cf}\label{5.3.26}{\ }
\begin{verbatim}
[()->Pi,[16/5,21840*n^4-19408*n^3+192*n^2-12*n-5],
        [-752/5,-512*n^3*(2*n+1)^3*(42*n-37)*(42*n+47)]]
\end{verbatim}
$$\pi=\dfrac{16}{5}-\dfrac{752/5}{2607-\dfrac{6151680}{194915-\dfrac{3152384000}{1246711-\dfrac{73006907904}{4351947-\ddots}}}}$$
Convergence type $E$ with $E=64$, $P=1/2$, and $C=-\sqrt{\pi}/24$, so that
$$\pi-\dfrac{p(n)}{q(n)}\sim-\dfrac{\sqrt{\pi}/3}{2^{6n+3}n^{1/2}}$$
$$A=1-(55/72)/n+(2795/3456)/n^2+\cdots$$
Series:
$$\dfrac{1}{\pi}=\dfrac{1}{16}\sum_{n\ge0}\dfrac{(42n+5)(1/2)_n^3}{n!^3}2^{-6n}$$
\end{cf}

\smallskip

\begin{cf}\label{5.3.27}{\ }
\begin{verbatim}
[()->Pi,[3,220*n^3-176*n^2-7*n+5],
        [6,4*n^2*(2*n+1)^2*(5*n-4)*(5*n+6)]]
\end{verbatim}
$$\pi=3+\dfrac{6}{42+\dfrac{396}{1047+\dfrac{38400}{4340+\dfrac{407484}{11241+\dfrac{2156544}{23070+\dfrac{7877100}{41147+\ddots}}}}}}$$
Convergence type $E$ with $E=-((1+\sqrt{5})/2)^{10}$, $P=0$, and
$C=4\pi/((1+\sqrt{5})/2)^9$, so that
$$\pi-\dfrac{p(n)}{q(n)}\sim(-1)^n\dfrac{4\pi}{((1+\sqrt{5})/2)^{10n+9}}\;.$$
$$A=1-(d/20)/n+(d/25+1/160)/n^2+(-2319d/80000-1/100)/n^3+\cdots$$
\end{cf}

\smallskip

\begin{cf}\label{5.3.28}{\ }
\begin{verbatim}
[()->Pi,[2,880*n^4-1760*n^3+776*n^2+104*n-33],
        [-38,4*n^2*(4*n^2-1)*(20*(n-1)^2-1)*(20*(n+1)^2-1)]]
\end{verbatim}
$$\pi=2+\dfrac{38}{33+\dfrac{948}{3279+\dfrac{816240}{31023+\dfrac{31753260}{125439+\dfrac{360142272}{349887+\dfrac{2270673900}{788847+\ddots}}}}}}$$
Convergence type $E$ with $E=-((1+\sqrt{5})/2)^{10}$, $P=0$, and
$C=4\pi/((1+\sqrt{5})/2)^5$, so that
$$\pi-\dfrac{p(n)}{q(n)}\sim(-1)^n\dfrac{4\pi}{((1+\sqrt{5})/2)^{10n+5}}\;.$$
$$A=1+(d/20)/n+(-d/40+1/160)/n^2+(711d/80000-1/160)/n^3+\cdots$$
\end{cf}

\smallskip

\begin{cf}\label{5.3.29}{\ }
\begin{verbatim}
[()->Pi/sqrt(3),[0,6*n^2-6*n+2],[1,-n^2*(9*n^2-1)]]
\end{verbatim}
$$\dfrac{\pi}{\sqrt{3}}=\dfrac{1}{2-\dfrac{8}{14-\dfrac{140}{38-\dfrac{720}{74-\dfrac{2288}{122-\dfrac{5600}{182-\ddots}}}}}}$$
Convergence type $P^+$ with $P=1/3$ and $C=2\pi\sqrt{3}/\G(1/3)^2$, so that
$$\dfrac{\pi}{\sqrt{3}}-\dfrac{p(n)}{q(n)}\sim\dfrac{2\pi\sqrt{3}/\G(1/3)^2}{n^{1/3}}$$
$$A=1-(1/6)/n+(20/567)/n^2+(1/486)/n^3+\cdots$$
Series:
$$\dfrac{\pi}{\sqrt{3}}=\dfrac{1}{2}\sum_{n\ge0}\dfrac{(4/3)_n}{(n+1)(5/3)_n}$$
Parametric family for $k\ge0$, $u\ge-1$, and $3 \nmid u$:
\begin{verbatim}
[()->Pi/sqrt(3),6*n^2-6*n+(2*k+1)*u+3*(k^2+k+1),-n^2*(9*n^2-u^2)]
\end{verbatim}
Convergence type $P^+$ with $P=2k+1+2u/3$.
\end{cf}

\smallskip

\begin{cf}\label{5.3.30}{\ }
\begin{verbatim}
[()->Pi/sqrt(3),[3/2,18*n^2-18*n+8],[3/2,-9*n^2*(9*n^2-1)]]
\end{verbatim}
$$\dfrac{\pi}{\sqrt{3}}=3/2+\dfrac{3/2}{8-\dfrac{72}{44-\dfrac{1260}{116-\dfrac{6480}{224-\dfrac{20592}{368-\dfrac{50400}{548-\ddots}}}}}}$$
Convergence type $P^+$ with $P=1$ and $C=\pi/3^{5/2}$, so that
$$\dfrac{\pi}{\sqrt{3}}-\dfrac{p(n)}{q(n)}\sim\dfrac{\pi/3^{5/2}}{n}$$
Series:
\begin{align*}
  \dfrac{\pi}{\sqrt{3}}&=\dfrac{3}{2}+\dfrac{3}{16}\sum_{n\ge0}\dfrac{n!^2}{(5/3)_n(7/3)_n}\\
  \dfrac{\sqrt{3}}{\pi}&=\dfrac{2}{3}\sum_{n\ge0}\dfrac{(-1/3)_n(1/3)_n}{n!^2}\end{align*}
Parametric family for $k\ge0$:
\begin{verbatim}
[()->Pi/sqrt(3),18*n^2-18*n+8+9*k*(k+1),-9*n^2*(9*n^2-1)]
\end{verbatim}
Convergence type $P^+$ with $P=2k+1$.
\end{cf}
          
\smallskip

\begin{cf}\label{5.3.31}{\ }
\begin{verbatim}
[()->Pi/sqrt(3),[0,6*n-5],[2,n^2*(3*n-2)^2]]
\end{verbatim}
$$\dfrac{\pi}{\sqrt{3}}=\dfrac{2}{1+\dfrac{1}{7+\dfrac{64}{13+\dfrac{441}{19+\dfrac{1600}{25+\dfrac{4225}{31+\ddots}}}}}}$$
Convergence type $P^-$ with $P=2$ and $C=1/3$, so that
$$\dfrac{\pi}{\sqrt{3}}-\dfrac{p(n)}{q(n)}\sim(-1)^n\dfrac{1/3}{n^2}$$
$$A=1-(1/3)/n-(5/9)/n^2+(11/27)/n^3+(91/81)/n^4+\cdots$$
Series:
$$\dfrac{\pi}{\sqrt{3}}=2\sum_{n\ge0}\dfrac{(-1)^n}{(n+1)(3n+1)}$$
Parametric family for $k\ge0$ and $3\nmid u$:
\begin{verbatim}
[()->Pi/sqrt(3),(2*k+1)*(6*n-3+2*u),n^2*(3*n+2*u)^2]
\end{verbatim}
Convergence type $P^-$ with $P=4k+2$.
\end{cf}

\smallskip

\begin{cf}\label{5.3.32}{\ }
\begin{verbatim}
[()->Pi/sqrt(3),[3/2,8,9*(2*n-1)],[3,(9*n^2-1)^2]]
\end{verbatim}
$$\dfrac{\pi}{\sqrt{3}}=3/2+\dfrac{3}{8+\dfrac{64}{27+\dfrac{1225}{45+\dfrac{6400}{63+\dfrac{20449}{81+\dfrac{50176}{99+\ddots}}}}}}$$
Convergence type $P^-$ with $P=2$ and $C=1/6$, so that
$$\dfrac{\pi}{\sqrt{3}}-\dfrac{p(n)}{q(n)}\sim(-1)^n\dfrac{1/6}{n^2}$$
$$A=1-1/n+(1/9)/n^2+(7/9)/n^3+(1/81)/n^4+\cdots$$
Series:`
$$\dfrac{\pi}{\sqrt{3}}=\dfrac{3}{2}+3\sum_{n\ge0}\dfrac{(-1)^n}{(3n+2)(3n+4)}$$
Parametric family for $k\ge0$:
\begin{verbatim}
[()->Pi/sqrt(3),9*(2*k+1)*(2*n-1),(9*n^2-1)^2]
\end{verbatim}
Convergence type $P^-$ with $P=4k+2$.
\end{cf}
        
\smallskip

\begin{cf}\label{5.3.33}{\ }
\begin{verbatim}
[()->Pi/sqrt(3),[0,7*n-5],[3/2,-6*n*(2*n-1)]]
\end{verbatim}
$$\dfrac{\pi}{\sqrt{3}}=\dfrac{3/2}{2-\dfrac{6}{9-\dfrac{36}{16-\dfrac{90}{23-\dfrac{168}{30-\dfrac{270}{37-\ddots}}}}}}$$
Convergence type $E$ with $E=4/3$, $P=1/2$, and $C=3\sqrt{\pi}/2$, so that
$$\dfrac{\pi}{\sqrt{3}}-\dfrac{p(n)}{q(n)}\sim\dfrac{3\sqrt{\pi}/2}{(4/3)^nn^{1/2}}\;.$$
$$A=1-(15/8)/n+(1249/128)/n^2-(86805/1024)/n^3+(33794859/32768)/n^4+\cdots$$
Series:
$$\dfrac{\pi}{\sqrt{3}}=\dfrac{3}{4}\sum_{n\ge0}\dfrac{n!}{(3/2)_n}(3/4)^n$$
Parametric families up to the trivial change $n\mapsto n+j$, with
$k\ge0$ and $u\ge0$:
\begin{verbatim}
[()->Pi/sqrt(3),7*n+4*u+k-5,-6*n*(2*n+2*u-1)]
\end{verbatim}
Convergence type $E$ with $E=4/3$ and $P=u+2k+1/2$.
\end{cf}

\smallskip

\begin{cf}\label{5.3.34}{\ }
\begin{verbatim}
[()->Pi/sqrt(3),[-4/3,5537*n^3-11307*n^2+6796*n-1188],
                [-240,-144*n*(2*n+1)*(3*n-2)*(3*n-1)*(49*n-58)*(49*n+40)]]
\end{verbatim}
$$\dfrac{\pi}{\sqrt{3}}=-4/3+\dfrac{240}{162-\dfrac{692064}{11472-\dfrac{158976000}{66936-\dfrac{2818392192}{199452-\dfrac{18571576320}{442242-\ddots}}}}}$$
Convergence type $E$ with $E=81/32$, $P=-1/2$, and $C=8\sqrt{3\pi}/9$, so that
$$\dfrac{\pi}{\sqrt{3}}-\dfrac{p(n)}{q(n)}\sim\dfrac{8\sqrt{3\pi}/9}{(81/32)^nn^{-1/2}}$$
$$A=1+(373/504)/n-(1487201/3556224)/n^2+\cdots$$
Series:
$$\dfrac{\pi}{\sqrt{3}}=-\dfrac{4}{3}+\dfrac{1}{24}\sum_{n\ge1}\dfrac{(49n-9)(1/2)_n(n-1)!}{(1/3)_n(2/3)_n}(32/81)^n$$
\end{cf}

\smallskip

\begin{cf}\label{5.3.35}{\ }
\begin{verbatim}
[()->Pi/sqrt(3),[0,4*n-2],[[3/2,-6],-6*[n*(2*n-1),(n+1)*(2*n+1)]]]
\end{verbatim}
$$\dfrac{\pi}{\sqrt{3}}=\dfrac{3/2}{2-\dfrac{6}{6-\dfrac{6}{10-\dfrac{36}{14-\dfrac{36}{18-\dfrac{90}{22-\ddots}}}}}}$$
Convergence type $E$ with $E=3$, $P=0$, and $C=\pi$, so that
$$\dfrac{\pi}{\sqrt{3}}-\dfrac{p(n)}{q(n)}\sim\dfrac{\pi}{3^n}\;.$$
$$A=1-(1/8)/n+(9/128)/n^2+(15/1024)/n^3-(2157/32768)/n^4+\cdots$$
\end{cf}

\smallskip

\begin{cf}\label{5.3.36}{\ }
\begin{verbatim}
[()->Pi/sqrt(3),[0,2*n-1],[3,3*n^2]]
\end{verbatim}
$$\dfrac{\pi}{\sqrt{3}}=\dfrac{3}{1+\dfrac{3}{3+\dfrac{12}{5+\dfrac{27}{7+\dfrac{48}{9+\dfrac{75}{11+\ddots}}}}}}$$
Convergence type $E$ with $E=-3$, $P=0$, and $C=\pi$, so that
$$\dfrac{\pi}{\sqrt{3}}-\dfrac{p(n)}{q(n)}\sim(-1)^n\dfrac{\pi}{3^n}\;.$$
$$A=1-(1/8)/n+(9/128)/n^2+(15/1024)/n^3-(2157/32768)/n^4+\cdots$$
Parametric families up to the trivial change $n\mapsto n+j$, with
$k\ge0$ and $0\le u\le 4k$:
\begin{verbatim}
[()->Pi/sqrt(3),2*n+4*k-1-u,3*n*(n+u)]
\end{verbatim}
Convergence type $E$ with $E=-3$ and $P=2k-u$.

Ap\'ery accelerates to \ref{5.3.40}.
\end{cf}

\smallskip

\begin{cf}\label{5.3.37}{\ }
\begin{verbatim}
[()->Pi/sqrt(3),[0,1,4*n],[2,1,3*(2*n-1)^2]]
\end{verbatim}
$$\dfrac{\pi}{\sqrt{3}}=\dfrac{2}{1+\dfrac{1}{8+\dfrac{27}{12+\dfrac{75}{16+\dfrac{147}{20+\dfrac{243}{24+\ddots}}}}}}$$
Convergence type $E$ with $E=-3$, $P=1$, and $C=3/4$, so that
$$\dfrac{\pi}{\sqrt{3}}-\dfrac{p(n)}{q(n)}\sim(-1)^n\dfrac{3/4}{3^n\cdot n}\;.$$
$$A=1-(1/4)/n-(1/8)/n^2+(7/32)/n^3-(1/32)/n^4-(47/128)/n^5+\cdots$$
Series:
$$\dfrac{\pi}{\sqrt{3}}=2\sum_{n\ge0}(-1)^n\dfrac{3^{-n}}{2n+1}$$
Parametric family up to the trivial change $n\mapsto n+j$
\begin{verbatim}
[()->Pi/sqrt(3),4*n+6*u+8*k,3*(2*n-1)*(2*n-1+2*u)]
\end{verbatim}
Convergence type $E$ with $E=-3$ and $P=2k+u+1$.
\end{cf}

\smallskip

\begin{cf}\label{5.3.38}{\ }
\begin{verbatim}
[()->Pi/sqrt(3),[0,5*n-3],[3,-2*n*(2*n-1)]]
\end{verbatim}
$$\dfrac{\pi}{\sqrt{3}}=\dfrac{3}{2-\dfrac{2}{7-\dfrac{12}{12-\dfrac{30}{17-\dfrac{56}{22-\dfrac{90}{27-\ddots}}}}}}$$
Convergence type $E$ with $E=4$, $P=1/2$, and $C=\sqrt{\pi}$, so that
$$\dfrac{\pi}{\sqrt{3}}-\dfrac{p(n)}{q(n)}\sim\dfrac{\sqrt{\pi}}{2^{2n}n^{1/2}}\;.$$
$$A=1-(13/24)/n+(227/384)/n^2-(9565/9216)/n^3+(2308873/884736)/n^4+\cdots$$
Series:
$$\dfrac{\pi}{\sqrt{3}}=\dfrac{3}{2}\sum_{n\ge0}\dfrac{n!}{(3/2)_n}2^{-2n}$$
Parametric families up to the trivial change $n\mapsto n+j$, with
$k\ge0$ and $u\ge0$:
\begin{verbatim}
[()->Pi/sqrt(3),5*n+4*u-3+3*k,-2*n*(2*n+2*u-1)]
\end{verbatim}
Convergence type $E$ with $E=4$ and $P=2k+u+1/2$.

Ap\'ery accelerates to \ref{5.3.41}.
\end{cf}

\smallskip

\begin{cf}\label{5.3.39}{\ }
\begin{verbatim}
[()->Pi/sqrt(3),[[3/2,8],[45*n^2-1,45*n^2+36*n+8]],
                 [[3,64],[-36*n^2*(9*n^2-1),(3*n+2)^2*(3*n+4)^2]]]
\end{verbatim}
$$\dfrac{\pi}{\sqrt{3}}=3/2+\dfrac{3}{8+\dfrac{64}{44-\dfrac{288}{89+\dfrac{1225}{179-\dfrac{5040}{260+\dfrac{6400}{404-\ddots}}}}}}$$
Convergence type $E$ with $E=i((1+\sqrt{5})/2)^5$, $P=0$, and
$C=\pi\sqrt{3}/((1+\sqrt{5})/2)^5$, so that
$$\dfrac{\pi}{\sqrt{3}}-\dfrac{p(n)}{q(n)}\sim(-1)^{\lfloor n/2\rfloor}\dfrac{\pi\sqrt{3}}{((1+\sqrt{5})/2)^{5n+5}}\;.$$
$$A=1-(8d/45)/n+(8d/45+32/405)/n^2+(-210448d/1366875-64/405)/n^3+\cdots$$
\end{cf}

\smallskip

\begin{cf}\label{5.3.40}{\ }
\begin{verbatim}
[()->Pi/sqrt(3),[0,3*(2*n-1)],[6,3*n^2]]
\end{verbatim}
$$\dfrac{\pi}{\sqrt{3}}=\dfrac{6}{3+\dfrac{3}{9+\dfrac{12}{15+\dfrac{27}{21+\dfrac{48}{27+\dfrac{75}{33+\ddots}}}}}}$$
Convergence type $E$ with $E=-(2+\sqrt{3})^2$, $P=0$, and
$C=2\pi\sqrt{3}/(2+\sqrt{3})$, so that
$$\dfrac{\pi}{\sqrt{3}}-\dfrac{p(n)}{q(n)}\sim(-1)^n\dfrac{2\pi\sqrt{3}}{(2+\sqrt{3})^{2n+1}}\;.$$
$$A=1-(d/8)/n+(d/16+3/128)/n^2+(-11d/1024-3/128)/n^3+\cdots$$
Parametric families up to the trivial change $n\mapsto n+j$, with $u\ge0$:
\begin{verbatim}
[()->Pi/sqrt(3),3*(2*n+2*u-1),3*n*(n+2*u)]
\end{verbatim}
Convergence type $E$ with $E=-(2+\sqrt{3})^2$ and $P=0$.
\end{cf}

\smallskip

\begin{cf}\label{5.3.41}{\ }
\begin{verbatim}
[()->Pi/sqrt(3),[0,4*n-2],[[3,-2],[-2*n*(2*n-1),-2*(n+1)*(2*n+1)]]]
\end{verbatim}
$$\dfrac{\pi}{\sqrt{3}}=\dfrac{3}{2-\dfrac{2}{6-\dfrac{2}{10-\dfrac{12}{14-\dfrac{12}{18-\dfrac{30}{22-\ddots}}}}}}$$
Convergence type $E$ with $E=(2+\sqrt{3})^2$, $P=0$, and $C=2\pi\sqrt{3}/(2+\sqrt{3})$, so that
$$\dfrac{\pi}{\sqrt{3}}-\dfrac{p(n)}{q(n)}\sim\dfrac{2\pi\sqrt{3}}{(2+\sqrt{3})^{2n+1}}\;.$$
$$A=1-(d/8)/n+(d/16+3/128)/n^2+(-11d/1024-3/128)/n^3+\cdots$$
\end{cf}

\smallskip

\begin{cf}\label{5.3.42}{\ }
\begin{verbatim}
[()->Pi/sqrt(3),[0,24,495*n^3-1122*n^2+723*n-120],
                [44,208,n^2*(3*n-2)*(3*n-1)*(15*n-19)*(15*n+11)]]
\end{verbatim}
$$\dfrac{\pi}{\sqrt{3}}=\dfrac{44}{24+\dfrac{208}{798+\dfrac{36080}{5316+\dfrac{733824}{16500+\dfrac{5123360}{37320+\dfrac{21912800}{70746+\ddots}}}}}}$$
Convergence type $E$ with $E=-((1+\sqrt{5})/2)^{10}$, $P=0$, and
$C=2\pi\sqrt{3}/((1+\sqrt{5})/2)^3$, so that
$$\dfrac{\pi}{\sqrt{3}}-\dfrac{p(n)}{q(n)}\sim(-1)^n\dfrac{2\pi\sqrt{3}}{((1+\sqrt{5})/2)^{10n+3}}\;.$$
$$A=1-(2d/45)/n+(4d/225+2/405)/n^2+(-7426d/1366875-8/2025)/n^3+\cdots$$
\end{cf}

\smallskip

\begin{cf}\label{5.3.43}{\ }
\begin{verbatim}
[()->Pi/sqrt(3),[2,(4*n-1)*(28*n^2-14*n-5)],
                [-5,-n^2*(2*n+1)^2*(4*n-3)*(4*n+5)]]
\end{verbatim}
$$\dfrac{\pi}{\sqrt{3}}=2-\dfrac{5}{27-\dfrac{81}{553-\dfrac{6500}{2255-\dfrac{67473}{5805-\dfrac{353808}{11875-\dfrac{1285625}{21137-\ddots}}}}}}$$
Convergence type $E$ with $E=(2+\sqrt{3})^4$, $P=0$, and $C=-2\pi\sqrt{3}/(2+\sqrt{3})^3$, so that
$$\dfrac{\pi}{\sqrt{3}}-\dfrac{p(n)}{q(n)}\sim-\dfrac{2\pi\sqrt{3}}{(2+\sqrt{3})^{4n+3}}\;.$$
$$A=1-(d/16)/n+(3d/64+3/512)/n^2+(-267d/8192-9/1024)/n^3+\cdots$$
\end{cf}

\smallskip

\begin{cf}\label{5.3.44}{\ }
\begin{verbatim}
[()->Pi/sqrt(2),[2,15,16*(2*n-1)],[4,(16*n^2-1)^2]]
\end{verbatim}
$$\dfrac{\pi}{\sqrt{2}}=2+\dfrac{4}{15+\dfrac{225}{48+\dfrac{3969}{80+\dfrac{20449}{112+\dfrac{65025}{144+\dfrac{159201}{176+\ddots}}}}}}$$
Convergence type $P^-$ with $P=2$ and $C=1/8$, so that
$$\dfrac{\pi}{\sqrt{2}}-\dfrac{p(n)}{q(n)}\sim(-1)^n\dfrac{1/8}{n^2}\;.$$
$$A=1-1/n+(1/16)/n^2+(7/8)/n^3+(1/256)/n^4-(691/256)/n^5+\cdots$$
Series:
$$\dfrac{\pi}{\sqrt{2}}=2+4\sum_{n\ge1}\dfrac{(-1)^{n+1}}{(4n-1)(4n+1)}$$
Parametric families for $k\ge0$:
\begin{verbatim}
[()->Pi/sqrt(2),16*(2*k+1)*(2*n-1),(16*n^2-(2*u+1)^2)^2]
\end{verbatim}
Convergence type $P^-$ with $P=4k+2$.

Ap\'ery accelerates to \ref{5.3.48} and \ref{5.3.49}.
\end{cf}

\smallskip

\begin{cf}\label{5.3.45}{\ }
\begin{verbatim}
[()->Pi/sqrt(2),[2,32*n^2-32*n+15],[2,-16*n^2*(16*n^2-1)]]
\end{verbatim}
$$\dfrac{\pi}{\sqrt{2}}=2+\dfrac{2}{15-\dfrac{240}{79-\dfrac{4032}{207-\dfrac{20592}{399-\dfrac{65280}{655-\dfrac{159600}{975-\ddots}}}}}}$$
Convergence type $P^+$ with $P=1$ and $C=\pi/2^{9/2}$, so that
$$\dfrac{\pi}{\sqrt{2}}-\dfrac{p(n)}{q(n)}\sim\dfrac{\pi/2^{9/2}}{n}\;.$$
$$A=1-(17/32)/n+(107/512)/n^2-(1163/32768)/n^3+\cdots$$
Series:
\begin{align*}\dfrac{\pi}{\sqrt{2}}&=2+\dfrac{2}{15}\sum_{n\ge0}\dfrac{n!^2}{(7/4)_n(9/4)_n}\\
\dfrac{\sqrt{2}}{\pi}&=\dfrac{1}{2}\sum_{n\ge0}\dfrac{(-1/4)_n(1/4)_n}{n!^2}\end{align*}
\end{cf}

\smallskip

\begin{cf}\label{5.3.46}{\ }
\begin{verbatim}
[()->Pi/sqrt(2),[2,12,12*n^2+1],[2,-4*n*(2*n+1)^3]]
\end{verbatim}
$$\dfrac{\pi}{\sqrt{2}}=2+\dfrac{2}{12-\dfrac{108}{49-\dfrac{1000}{109-\dfrac{4116}{193-\dfrac{11664}{301-\dfrac{26620}{433-\ddots}}}}}}$$
Convergence type $E$ with $E=2$, $P=3/2$, and $C=1/\sqrt{\pi}$, so that
$$\dfrac{\pi}{\sqrt{2}}-\dfrac{p(n)}{q(n)}\sim\dfrac{1/\sqrt{\pi}}{2^nn^{3/2}}\;.$$
$$A=1-(29/8)/n+(1881/128)/n^2-(77495/1024)/n^3+\cdots$$
Series:
$$\dfrac{\pi}{\sqrt{2}}=2+\dfrac{1}{2}\sum_{n\ge0}\dfrac{(3/2)_n}{(2n+3)(n+1)!}2^{-n}$$
Parametric family for $k\ge0$:
\begin{verbatim}
[()->Pi/sqrt(2),12*n^2+1,-4*(n-k)*(2*n+1)^3]
\end{verbatim}
Convergence type $E$ with $E=2$ and $P=3k+3/2$.
\end{cf}

\smallskip

\begin{cf}\label{5.3.47}{\ }
\begin{verbatim}
[()->Pi/sqrt(2),[13/6,112*n^4+192*n^3+24*n^2-80*n-13],
                [27/2,64*n^3*(2*n-1)*(2*n+3)^4]]
\end{verbatim}
$$\dfrac{\pi}{\sqrt{2}}=\dfrac{13}{6}+\dfrac{27/2}{235+\dfrac{40000}{3251+\dfrac{3687936}{14219+\dfrac{56687040}{41011+\dfrac{419786752}{94187+\ddots}}}}}$$
Convergence type $E$ with $E=-8$, $P=3/2$, and $C=\sqrt{\pi}/27$, so that
$$\dfrac{\pi}{\sqrt{2}}-\dfrac{p(n)}{q(n)}\sim(-1)^n\dfrac{\sqrt{\pi}/27}{2^{3n}n^{3/2}}$$
$$A=1-(41/24)/n+(8755/3456)/n^2+\cdots$$
Series:
$$\dfrac{\sqrt{2}}{\pi}=6\sum_{n\ge0}(-1)^n\dfrac{(2n+1)(3/2)_n^3}{(12n^2+1)(12n^2+24n+13)n!^3}2^{-3n}$$
\end{cf}

\smallskip

\begin{cf}\label{5.3.48}{\ }
\begin{verbatim}
[()->Pi/sqrt(2),[[2,30],[80*n^2-1,80*n^2+96*n+31]],
               [[4,-960],[(16*n^2-1)^2,-64*(n+1)^2*(4*n+3)*(4*n+5)]]]
\end{verbatim}
$$\dfrac{\pi}{\sqrt{2}}=2+\dfrac{4}{30-\dfrac{960}{79+\dfrac{225}{207-\dfrac{16128}{319+\dfrac{3969}{543-\dfrac{82368}{719+\ddots}}}}}}$$
Convergence type $E$ with $E=i((1+\sqrt{5})/2)^5$, $P=0$, and
$C=\pi\sqrt{2}/((1+\sqrt{5})/2)^5$, so that
$$\dfrac{\pi}{\sqrt{2}}-\dfrac{p(n)}{q(n)}\sim(-1)^{\lfloor(n+2)/2\rfloor}\dfrac{\pi\sqrt{2}}{((1+\sqrt{5})/2)^{5n+5}}\;.$$
$$A=1-(11d/40)/n+(11d/40+121/640)/n^2+(-144921d/640000-121/320)/n^3+\cdots$$
\end{cf}

\smallskip

\begin{cf}\label{5.3.49}{\ }
\begin{verbatim}
[()->Pi/sqrt(2),[[2,15],[80*n^2-1,80*n^2+64*n+15]],
               [[4,225],[-64*n^2*(16*n^2-1),(4*n+3)^2*(4*n+5)^2]]]
\end{verbatim}
$$\dfrac{\pi}{\sqrt{2}}=2+\dfrac{4}{15+\dfrac{225}{79-\dfrac{960}{159+\dfrac{3969}{319-\dfrac{16128}{463+\dfrac{20449}{719-\ddots}}}}}}$$
Convergence type $E$ with $E=-i((1+\sqrt{5})/2)^5$, $P=0$, and
$C=\pi\sqrt{2}/((1+\sqrt{5})/2)^5$, so that
$$\dfrac{\pi}{\sqrt{2}}-\dfrac{p(n)}{q(n)}\sim(-1)^{\lfloor(n-1)/2\rfloor}\dfrac{\pi\sqrt{2}}{((1+\sqrt{5})/2)^{5n+5}}\;.$$
$$A=1-(11d/40)/n+(11d/40+121/640)/n^2+(-144921d/640000-121/320)/n^3+\cdots$$
\end{cf}

\smallskip

\begin{cf}\label{5.3.50}{\ }
\begin{verbatim}
[()->Pi/sqrt(2),[0,580*n^3-1296*n^2+803*n-102],
                [-32,-12*n*(2*n-1)*(5*n-6)*(5*n+4)*(6*n-5)*(6*n-1)]]
\end{verbatim}
$$\dfrac{\pi}{\sqrt{2}}=\dfrac{32}{15-\dfrac{540}{960-\dfrac{310464}{6303-\dfrac{6802380}{19494-\dfrac{49335552}{44013-\ddots}}}}}$$
Convergence type $E$ with $E=27/2$, $P=1/2$, and $C=4\sqrt{\pi}/5$, so that
$$\dfrac{\pi}{\sqrt{2}}-\dfrac{p(n)}{q(n)}\sim\dfrac{4\sqrt{\pi}/5}{(27/2)^nn^{1/2}}$$
$$A=1-(407/1800)/n+(748849/6480000)/n^2+\cdots$$
Series:
$$\dfrac{\pi}{\sqrt{2}}=\dfrac{8}{15}\sum_{n\ge0}\dfrac{(5n+4)(1/2)_nn!}{(7/6)_n(11/6)_n}(2/27)^n$$
\end{cf}

\smallskip

\begin{cf}\label{5.3.51}{\ }
\begin{verbatim}
[()->Pi/sqrt(2),[2,1160*n^3-1084*n^2-2*n+31],
                [22,-12*n*(2*n+1)*(6*n-1)*(6*n+1)*(10*n-9)*(10*n+11)]]
\end{verbatim}
$$\dfrac{\pi}{\sqrt{2}}=2+\dfrac{22}{105-\dfrac{26460}{4971-\dfrac{5851560}{21589-\dfrac{70081956}{56919-\dfrac{392720400}{117921-\ddots}}}}}$$
Convergence type $E$ with $E=27/2$, $P=1/2$, and $C=2\sqrt{\pi}/15$, so that
$$\dfrac{\pi}{\sqrt{2}}-\dfrac{p(n)}{q(n)}\sim\dfrac{2\sqrt{\pi}/15}{(27/2)^nn^{1/2}}$$
$$A=1-(1067/1800)/n+(3664969/6480000)/n^2+\cdots$$
Series:
$$\dfrac{\pi}{\sqrt{2}}=2+\dfrac{1}{2}\sum_{n\ge1}\dfrac{(10n+1)(1/2)_n(n-1)!}{(5/6)_n(7/6)_n}(2/27)^n$$
\end{cf}

\smallskip

\begin{cf}\label{5.3.52}{\ }
\begin{verbatim}
[()->Pi/(sqrt(2)-1),[8,63,128*n^2-128*n+62],[-16,-(64*n^2-1)^2]]
\end{verbatim}
$$\dfrac{\pi}{\sqrt{2}-1}=8-\dfrac{16}{63-\dfrac{3969}{318-\dfrac{65025}{830-\dfrac{330625}{1598-\dfrac{1046529}{2622-\dfrac{2556801}{3902-\ddots}}}}}}$$
Convergence type $P^+$ with $P=1$ and $C=-1/4$, so that
$$\dfrac{\pi}{\sqrt{2}-1}-\dfrac{p(n)}{q(n)}\sim-\dfrac{1/4}{n}$$
$$A=1-(3/2)/n+(139/64)/n^2-(387/128)/n^3+\cdots$$
Series:
$$\dfrac{\pi}{\sqrt{2}-1}=-8-16\sum_{n\ge0}\dfrac{1}{64n^2-1}$$
Parametric family for $k\ge0$:
\begin{verbatim}
[()->Pi/(sqrt(2)-1),128*n^2-128*n+62+64*k*(k+1),-(64*n^2-1)^2]
\end{verbatim}
Convergence type $P^+$ with $P=2k+1$.
\end{cf}

\smallskip

\begin{cf}\label{5.3.53}{\ }
\begin{verbatim}
[()->Pi/(sqrt(2)-1),[8,8*n-4],[-2,n^2*(16*n^2-1)]]
\end{verbatim}
$$\dfrac{\pi}{\sqrt{2}-1}=8-\dfrac{2}{4+\dfrac{15}{12+\dfrac{252}{20+\dfrac{1287}{28+\dfrac{4080}{36+\dfrac{9975}{44+\ddots}}}}}}$$
Convergence type $P^-$ with $P=2$ and $C=-\pi/(32(\sqrt{2}-1))$, so that
$$\dfrac{\pi}{\sqrt{2}-1}-\dfrac{p(n)}{q(n)}\sim(-1)^{n+1}\dfrac{\pi/(32(\sqrt{2}-1))}{n^2}$$
$$A=1-1/n+(1/32)/n^2+(15/16)/n^3+(5/4096)/n^4+\cdots$$
Series:
\begin{align*}
  \dfrac{\pi}{\sqrt{2}-1}&=8-\dfrac{8}{63}\sum_{n\ge0}\dfrac{(4n+3)n!^2(3/8)_n(5/8)_n}{(3/2)_n^2(15/8)_n(17/8)_n}\\
  \dfrac{\sqrt{2}-1}{\pi}&=\dfrac{2}{15}\sum_{n\ge0}\dfrac{(4n+1)(1/2)_n^2(-1/8)_n(1/8)_n}{n!^2(11/8)_n(13/8)_n}\end{align*}
Parametric family for $k\ge0$:
\begin{verbatim}
[()->Pi/(sqrt(2)-1),(2*k+1)*(8*n-4),n^2*(16*n^2-1)]
\end{verbatim}
Convergence type $P^-$ with $P=4k+2$.
\end{cf}

\smallskip

\begin{cf}\label{5.3.54}{\ }
\begin{verbatim}
[()->Pi/(sqrt(2)-1),[36/5,81,96*n^2-10],[36,(2*n+1)^2*(4*n-1)^2*(4*n+5)^2]]
\end{verbatim}
$$\dfrac{\pi}{\sqrt{2}-1}=36/5+\dfrac{36}{81+\dfrac{6561}{374+\dfrac{207025}{854+\dfrac{1713481}{1526+\dfrac{8037225}{2390+\dfrac{27300625}{3446+\ddots}}}}}}$$
Convergence type $P^-$ with $P=3$ and $C=9/16$, so that
$$\dfrac{\pi}{\sqrt{2}-1}-\dfrac{p(n)}{q(n)}\sim(-1)^n\dfrac{9/16}{n^3}$$
$$A=1-3/n+(81/16)/n^2-(85/16)/n^3+\cdots$$
Series:
$$\dfrac{\pi}{\sqrt{2}-1}=\dfrac{36}{5}+36\sum_{n\ge0}(-1)^n\dfrac{1}{(2n+3)(4n+3)(4n+9)}$$
\end{cf}

\smallskip

\begin{cf}\label{5.3.55}{\ }
\begin{verbatim}
[()->Pi/(sqrt(2)+1),[8/3,55,128*n^2-128*n+46],[-48,-(64*n^2-9)^2]]
\end{verbatim}
$$\dfrac{\pi}{\sqrt{2}+1}=8/3-\dfrac{48}{55-\dfrac{3025}{302-\dfrac{61009}{814-\dfrac{321489}{1582-\dfrac{1030225}{2606-\dfrac{2531281}{3886-\ddots}}}}}}$$
Convergence type $P^+$ with $P=1$ and $C=-3/4$, so that
$$\dfrac{\pi}{\sqrt{2}+1}-\dfrac{p(n)}{q(n)}\sim-\dfrac{3/4}{n}$$
$$A=1-(3/2)/n+(425/192)/n^2-(411/128)/n^3+\cdots$$
Series:
$$\dfrac{\pi}{\sqrt{2}+1}=-\dfrac{8}{3}-48\sum_{n\ge0}\dfrac{1}{64n^2-9}$$
Parametric family for $k\ge0$:
\begin{verbatim}
[()->Pi/(sqrt(2)+1),128*n^2-128*n+46+64*k*(k+1),-(64*n^2-9)^2]
\end{verbatim}
Convergence type $P^+$ with $P=2k+1$.
\end{cf}

\smallskip

\begin{cf}\label{5.3.56}{\ }
\begin{verbatim}
[()->Pi/(sqrt(2)+1),[8/3,8*n-4],[-6,n^2*(16*n^2-9)]]
\end{verbatim}
$$\dfrac{\pi}{\sqrt{2}+1}=8/3-\dfrac{6}{4+\dfrac{7}{12+\dfrac{220}{20+\dfrac{1215}{28+\dfrac{3952}{36+\dfrac{9775}{44+\ddots}}}}}}$$
Convergence type $P^-$ with $P=2$ and $C=-9\pi/(32(\sqrt{2}+1))$, so that
$$\dfrac{\pi}{\sqrt{2}+1}-\dfrac{p(n)}{q(n)}\sim(-1)^{n+1}\dfrac{9\pi/(32(\sqrt{2}+1))}{n^2}$$
$$A=1-1/n+(9/32)/n^2+(7/16)/n^3+(405/4096)/n^4+\cdots$$
Series:
\begin{align*}
  \dfrac{\pi}{\sqrt{2}+1}&=\dfrac{8}{3}-\dfrac{24}{55}\sum_{n\ge0}\dfrac{(4n+3)n!^2(1/8)_n(7/8)_n}{(3/2)_n^2(13/8)_n(19/8)_n}\\
  \dfrac{\sqrt{2}+1}{\pi}&=\dfrac{6}{7}\sum_{n\ge0}\dfrac{(4n+1)(1/2)_n^2(-3/8)_n(3/8)_n}{n!^2(9/8)_n(15/8)_n}\end{align*}
Parametric family for $k\ge0$:
\begin{verbatim}
[()->Pi/(sqrt(2)+1),(2*k+1)*(8*n-4),n^2*(16*n^2-9)]
\end{verbatim}
Convergence type $P^-$ with $P=4k+2$.
\end{cf}

\smallskip

\begin{cf}\label{5.3.57}{\ }
\begin{verbatim}
[()->Pi/(sqrt(2)+1),[4/3,105,96*n^2+6],[-4,(2*n+1)^2*(4*n+1)^2*(4*n+3)^2]]
\end{verbatim}
$$\dfrac{\pi}{\sqrt{2}+1}=4/3-\dfrac{4}{105+\dfrac{11025}{390+\dfrac{245025}{870+\dfrac{1863225}{1542+\dfrac{8450649}{2406+\dfrac{28227969}{3462+\ddots}}}}}}$$
Convergence type $P^-$ with $P=3$ and $C=-1/16$, so that
$$\dfrac{\pi}{\sqrt{2}+1}-\dfrac{p(n)}{q(n)}\sim(-1)^{n+1}\dfrac{1/16}{n^3}$$
$$A=1-3/n+(73/16)/n^2-(45/16)/n^3+\cdots$$
Series:
$$\dfrac{\pi}{\sqrt{2}+1}=\dfrac{4}{3}-4\sum_{n\ge0}(-1)^n\dfrac{1}{(2n+3)(4n+5)(4n+7)}$$
\end{cf}

\smallskip

In addition to the given parametric families, there exist many CFs similar
to \ref{5.3.52}, \ref{5.3.53}, \ref{5.3.55}, and \ref{5.3.56}.

\smallskip

We have given continued fractions for $\pi$, $\pi/\sqrt{3}$, and
$\pi/\sqrt{2}$ (as well as for $\pi/(\sqrt{2}\pm1)$).
It would be interesting to have other not too complicated
CF for $\pi/\sqrt{d}$ for $d>3$ squarefree. Of course, such CFs can be
obtained by transforming into a CF the series for $L(\chi_{-d},1)$, but this
has little interest.

\smallskip

\begin{cf}\label{5.3.58}{\ }
\begin{verbatim}
[()->log(2),[0,6,36*n^2-54*n+23],[3/2,-3*(2*n-1)*(3*n-1)^2*(6*n+1)]]
\end{verbatim}
$$\log(2)=\dfrac{3/2}{6-\dfrac{84}{59-\dfrac{2925}{185-\dfrac{18240}{383-\dfrac{63525}{653-\dfrac{164052}{995-\ddots}}}}}}$$
Convergence type $P^+$ with $P=1/3$ and $C=\pi\sqrt{3}/(2^{2/3}\G(1/3)^2)$,
so that
$$\log(2)-\dfrac{p(n)}{q(n)}\sim\dfrac{\pi\sqrt{3}/(2^{2/3}\G(1/3)^2)}{n^{1/3}}\;.$$
$$A=1-(5/72)/n-(23/2268)/n^2+(73/8748)/n^3+\cdots$$
Series:
$$\log(2)=\dfrac{1}{2}\sum_{n\ge0}\dfrac{(7/6)_n}{(3n+2)(3/2)_n}$$
Parametric family for $k\ge0$:
\begin{verbatim}
[()->log(2),36*n^2-54*n+23+6*k*(3*k+1),-3*(2*n-1)*(3*n-1)^2*(6*n+1)]
\end{verbatim}
Convergence type $P^+$ with $P=2k+1/3$.
\end{cf}

\smallskip

\begin{cf}\label{5.3.59}{\ }
\begin{verbatim}
[()->log(2),[0,4*n^2-3*n+1],[1/2,-2*n^3*(2*n+1)]]
\end{verbatim}
$$\log(2)=\dfrac{1/2}{2-\dfrac{6}{11-\dfrac{80}{28-\dfrac{378}{53-\dfrac{1152}{86-\dfrac{2750}{127-\ddots}}}}}}$$
Convergence type $P^+$ with $P=1/2$ and $C=1/\sqrt{\pi}$, so that
$$\log(2)-\dfrac{p(n)}{q(n)}\sim\dfrac{1/\sqrt{\pi}}{n^{1/2}}\;.$$
$$A=1-(7/24)/n+(61/640)/n^2-(307/21504)/n^3-\cdots$$
Series:
$$\log(2)=\dfrac{1}{4}\sum_{n\ge0}\dfrac{(3/2)_n}{(n+1)(n+1)!}$$
Parametric family for $k\ge0$:
\begin{verbatim}
[()->log(2),4*n^2-3*n+1+k*(2*k+1),-2*n^3*(2*n+1)]]
\end{verbatim}
Convergence type $P^+$ with $P=2k+1/2$.
\end{cf}

\smallskip

\begin{cf}\label{5.3.60}{\ }
\begin{verbatim}
[()->log(2),[0,1],[1,n^2]]
\end{verbatim}
$$\log(2)=\dfrac{1}{1+\dfrac{1}{1+\dfrac{4}{1+\dfrac{9}{1+\dfrac{16}{1+\dfrac{25}{1+\ddots}}}}}}$$
Convergence type $P^-$ with $P=1$ and $C=1/2$, so that
$$\log(2)-\dfrac{p(n)}{q(n)}\sim\dfrac{(-1)^n}{2n}\;.$$
$$A=1-(1/2)/n+(1/4)/n^3-(1/2)/n^5+(17/8)/n^7+\cdots$$
Series:
$$\log(2)=\sum_{n\ge1}\dfrac{(-1)^{n+1}}{n}$$
Parametric families up to the trivial change $n\mapsto n+j$, with $u\ge0$
and $k\ge0$:
\begin{verbatim}
[()->log(2),u+2*k+1,n*(n+u)]
\end{verbatim}
Convergence type $P^-$ with $P=2k+u+1$.

Ap\'ery accelerates to \ref{5.3.83}.
\end{cf}

\smallskip

\begin{cf}\label{5.3.61}{\ }
\begin{verbatim}
[()->log(2),[1/2,4*n^2-4*n+3],[1/2,-n^2*(4*n^2-1)]]
\end{verbatim}
$$\log(2)=1/2+\dfrac{1/2}{3-\dfrac{3}{11-\dfrac{60}{27-\dfrac{315}{51-\dfrac{1008}{83-\dfrac{2475}{123-\ddots}}}}}}$$
Convergence type $P^+$ with $P=2$ and $C=1/16$, so that
$$\log(2)-\dfrac{p(n)}{q(n)}\sim\dfrac{1/16}{n^2}\;.$$
$$A=1-1/n+(5/8)/n^2-(1/4)/n^3+(1/16)/n^4-(1/16)*n^5+(25/256)/n^6+\cdots$$
Series:
$$\log(2)=\dfrac{1}{2}+\sum_{n\ge1}\dfrac{1}{2n(2n-1)(2n+1)}$$
Parametric families: there exist many families of the type
\begin{verbatim}
[()->log(2),4*n^2+a1*n+a0,-n*(n+c0)*(2*n+b0)*(2*n+b1)]
\end{verbatim}
where the variables satisfy several conditions, including
\begin{verbatim}
a1=2*c0+b0+b1-4, 
8*a0+4*c0^2+(-4*b0+(-4*b1+8))*c0+b0^2+(-2*b1+4)*b0+b1^2+4*b1-12=m^2,
...
\end{verbatim}
all with convergence type $P^+$ with $P=m/2$.
\end{cf}

\smallskip

\begin{cf}\label{5.3.62}{\ }
\begin{verbatim}
[()->log(2),[0,59*n^2-59*n+20],[5,-24*n^2*(36*n^2-1)]]
\end{verbatim}
$$\log(2)=\dfrac{5}{20-\dfrac{840}{138-\dfrac{13728}{374-\dfrac{69768}{728-\dfrac{220800}{1200-\dfrac{539400}{1790-\ddots}}}}}}$$
Convergence type $E$ with $E=32/27$, $P=0$, and $C=3\pi/16$, so that
$$\log(2)-\dfrac{p(n)}{q(n)}\sim\dfrac{3\pi/16}{(32/27)^n}\;.$$
$$A=1-(125/36)/n+(20125/2592)/n^2-(74212325/279936)/n^3+\cdots$$
\end{cf}

\smallskip

\begin{cf}\label{5.3.63}{\ }
\begin{verbatim}
[()->log(2),[3/2,28*n^3-12*n^2-4*n+3],
            [-27/4,-3*n^2*(2*n+1)*(2*n+3)*(4*n-1)*(4*n+1)]]
\end{verbatim}
$$\log(2)=\dfrac{3}{2}-\dfrac{27/4}{15-\dfrac{675}{171-\dfrac{26460}{639-\dfrac{243243}{1587-\dfrac{1211760}{3183-\dfrac{4279275}{5595-\ddots}}}}}}$$
Convergence type $E$ with $E=4/3$, $P=1$, and $C=-9\sqrt{2}/8$, so that
$$\log(2)-\dfrac{p(n)}{q(n)}\sim-\dfrac{3\sqrt{2}/2}{(4/3)^{n+1}n}$$
$$A=1-(61/16)/n+(13529/512)/n^2+\cdots$$
Series:
$$\log(2)=\dfrac{3}{2}-\dfrac{3}{4}\sum_{n\ge1}\dfrac{(1/2)_n(3/2)_n}{n(3/4)_n(5/4)_n}(3/4)^n$$
\end{cf}

\smallskip

\begin{cf}\label{5.3.64}{\ }
\begin{verbatim}
[()->log(2),[0,3*n-1],[1,-2*n^2]]
\end{verbatim}
$$\log(2)=\dfrac{1}{2-\dfrac{2}{5-\dfrac{8}{8-\dfrac{18}{11-\dfrac{32}{14-\dfrac{50}{17-\ddots}}}}}}$$
Convergence type $E$ with $E=2$, $P=1$, and $C=1$, so that
$$\log(2)-\dfrac{p(n)}{q(n)}\sim\dfrac{1}{2^nn}\;.$$
$$A=1-2/n+6/n^2-26/n^3+150/n^4-1082/n^5+9366/n^6-\cdots$$
Series:
$$\log(2)=\dfrac{1}{2}\sum_{n\ge0}\dfrac{2^{-n}}{n+1}$$
Parametric families up to the trivial change $n\mapsto n+j$,
with $u\ge0$ and $k\ge0$:
\begin{verbatim}
[()->log(2),3*n+2*u+k-1,-2*n*(n+u)]
\end{verbatim}
Convergence type $E$ with $E=2$ and $P=2k+u+1$.

Ap\'ery accelerates to \ref{5.3.83}.
\end{cf}

\smallskip

\begin{cf}\label{5.3.65}{\ }
\begin{verbatim}
[()->log(2),[n*(3*n-1)],[1,-2*n^3*(n+1)]]
\end{verbatim}
$$\log(2)=\dfrac{1}{2-\dfrac{4}{10-\dfrac{48}{24-\dfrac{216}{44-\dfrac{640}{70-\dfrac{1500}{102-\ddots}}}}}}$$
Convergence type $E$ with $E=2$, $P=1$, and $C=1$, so that
$$\log(2)-\dfrac{p(n)}{q(n)}\sim\dfrac{1}{2^nn}\;.$$
$$A=1-2/n+6/n^2-26/n^3+150/n^4-1082/n^5+9366/n^6-\cdots$$
Series:
$$\log(2)=\dfrac{1}{2}\sum_{n\ge0}\dfrac{2^{-n}}{n+1}$$
Parametric families up to the trivial change $n\mapsto n+j$, where
some parameters may be negative:
\begin{verbatim}
[()->log(2),3*n^2+(2*(u+v+w)-3)*n+u*v+u*w+v*w-u-v-w+1,
            -2*n*(n+u)*(n+v)*(n+w)]
[()->log(2),3*n^2+(4*(2*u+v+w)+1)*n+4*(u+v)*(u+w),
            -2*n*(n+u)*(n+2*u+2*v+1)*(n+2*u+2*w+1)]
[()->log(2),3*n^2+(4*(u+v+w)-1)*n+2*(2*v-1)*(u+w),
            -2*n*(n+u)*(n+2*v)*(n+2*u+2*w+1)]
\end{verbatim}
where in the first family {\tt u}, {\tt v}, and {\tt w} must not be all of
the same parity.

Convergence type $E$ with $E=2$ and $P=u+v+w$, $u+2v+2w+2$, $-u+2v+2w+1$,
and $-u$ respectively.
\end{cf}

There exist many other parametric families for $\log(2)$ with convergence
type $E=2$, the above are only given as typical examples.

This CF is of course termwise identical to the previous one, but with
completely different parametric families. A special case of the last
family is as follows:

\smallskip

\begin{cf}\label{5.3.66}{\ }
\begin{verbatim}
[()->log(2),[3*n^2-1],[2,-2*n^3*(n+2)]]
\end{verbatim}
$$\log(2)=-1+\dfrac{2}{2-\dfrac{6}{11-\dfrac{64}{26-\dfrac{270}{47-\dfrac{768}{74-\dfrac{1750}{107-\ddots}}}}}}$$
Convergence type $E$ with $E=2$, $P=0$, and $C=1$, so that
$$\log(2)-\dfrac{p(n)}{q(n)}\sim\dfrac{1}{2^n}\;.$$
$$A=1+1/n-2/n^2+6/n^3-26/n^4+150/n^5-1082/n^6+\cdots$$
Series:
$$\log(2)=-1+\dfrac{1}{2}\sum_{n\ge0}\dfrac{n+2}{n+1}2^{-n}$$
Parametric families for $u\ge0$ and $k\ge0$:
\begin{verbatim}
[()->log(2),3*n^2+u*(n-1)-1,-2*n^3*(n+u+2)]
\end{verbatim}
Convergence type $E$ with $E=2$ and $P=-u$.
\end{cf}

\smallskip

\begin{cf}\label{5.3.67}{\ }
\begin{verbatim}
[()->log(2),[0,27*n^2-27*n+8],[3,-18*n^2*(9*n^2-1)]]
\end{verbatim}
$$\log(2)=\dfrac{3}{8-\dfrac{144}{62-\dfrac{2520}{170-\dfrac{12960}{332-\dfrac{41184}{548-\dfrac{100800}{818-\ddots}}}}}}$$
Convergence type $E$ with $E=2$, $P=0$, and $C=\pi/(2\sqrt{3})$, so that
$$\log(2)-\dfrac{p(n)}{q(n)}\sim\dfrac{\pi/\sqrt{3}}{2^{n+1}}\;.$$
$$A=1-(8/9)/n+(68/81)/n^2-(9580/2187)/n^3+(183428/19683)/n^4+\cdots$$
Parametric family for all $k$:
\begin{verbatim}
[()->log(2),27*n^2-27*n+8,-18*n*(n-k)*(9*n^2-1)]
\end{verbatim}
Convergence type $E$ with $E=2$ and $P=3k$.
\end{cf}

\smallskip

\begin{cf}\label{5.3.68}
\begin{verbatim}
[()->log(2),[0,7*n+1],[15/2,72*n*(2*n-1)]]
\end{verbatim}
$$\log(2)=\dfrac{15/2}{8+\dfrac{72}{15+\dfrac{432}{22+\dfrac{1080}{29+\dfrac{2016}{36+\dfrac{3240}{43+\ddots}}}}}}$$
Convergence type $E$ with $E=-16/9$, $P=1/2$, and $C=3\sqrt{\pi}/10$, so that
$$\log(2)-\dfrac{p(n)}{q(n)}\sim(-1)^n\dfrac{3\sqrt{\pi}/10}{(4/3)^{2n}n^{1/2}}$$
$$A=1-(39/200)/n-(3599/80000)/n^2+(323511/3200000)/n^3+\cdots$$
Series:
$$\log(2)=\dfrac{15}{16}\sum_{n\ge0}\dfrac{n!}{(3/2)_n}(-9/16)^n$$
Parametric family for $k\ge0$:
\begin{verbatim}
[()->log(2),7*n+25*k+1,72*n*(2*n-1)]
\end{verbatim}
Convergence type $E$ with $E=-16/9$ and $P=2k+1/2$.
\end{cf}

\smallskip

\begin{cf}\label{5.3.69}
\begin{verbatim}
[()->log(2),[0,17*(2*n-1)],[15/2,-225*n^2]]
\end{verbatim}
$$\log(2)=\dfrac{15/2}{17-\dfrac{225}{51-\dfrac{900}{85-\dfrac{2025}{119-\dfrac{3600}{153-\dfrac{5625}{187-\ddots}}}}}}$$
Convergence type $E$ with $E=25/9$, $P=0$, and $C=3\pi/10$, so that
$$\log(2)-\dfrac{p(n)}{q(n)}\sim\dfrac{3\pi/10}{(5/3)^{2n}}$$
$$A=1-(17/32)/n+(833/2048)/n^2-(81899/65536)/n^3+\cdots$$
Parametric family for $k\ge0$:
\begin{verbatim}
[()->log(2),34*n-17+16*k,-225*n^2]
\end{verbatim}
Convergence type $E$ with $E=25/9$ and $P=2k$.
\end{cf}

\smallskip

\begin{cf}\label{5.3.70}
\begin{verbatim}
[()->log(2),[0,25,4*(17*n-13)],[15,-225*(2*n-1)^2]]
\end{verbatim}
$$\log(2)=\dfrac{15}{25-\dfrac{225}{84-\dfrac{2025}{152-\dfrac{5625}{220-\dfrac{11025}{288-\dfrac{18225}{356-\ddots}}}}}}$$
Convergence type $E$ with $E=25/9$, $P=1$, and $C=15/32$, so that
$$\log(2)-\dfrac{p(n)}{q(n)}\sim\dfrac{15/32}{(5/3)^{2n}n}$$
$$A=1-(17/16)/n+(257/128)/n^2-(12019/2048)/n^3+\cdots$$
Series:
$$\log(2)=\sum_{n\ge0}\dfrac{(3/5)^{2n+1}}{2n+1}$$
Parametric family for $k\ge0$:
\begin{verbatim}
[()->log(2),4*(17*n+8*k-13),-225*(2*n-1)^2]
\end{verbatim}
Convergence type $E$ with $E=25/9$ and $P=2k+1$.
\end{cf}

\smallskip

\begin{cf}\label{5.3.71}{\ }
\begin{verbatim}
[()->log(2),[0,2,10*n^2-10*n+3],[1,-10,-n^2*(16*n^2-1)]]
\end{verbatim}
$$\log(2)=\dfrac{1}{2-\dfrac{10}{23-\dfrac{252}{63-\dfrac{1287}{123-\dfrac{4080}{203-\dfrac{9975}{303-\ddots}}}}}}$$
Convergence type $E$ with $E=4$, $P=0$, and $C=\pi/2^{3/2}$, so that
$$\log(2)-\dfrac{p(n)}{q(n)}\sim\dfrac{\pi}{2^{2n+3/2}}\;.$$
$$A=1-(9/16)/n+(225/512)/n^2-(6579/8192)/n^3+\cdots$$
Parametric family for $k\ge0$:
\begin{verbatim}
[()->log(2),10*n^2-10*n+3+2*k*(2*n-1),-n^2*(4*n-2*k-1)*(4*n-2*k+1)]
\end{verbatim}
Convergence type $E$ with $E=4$ and $P=3k$.
\end{cf}

\smallskip

\begin{cf}\label{5.3.72}{\ }
\begin{verbatim}
[()->log(2),[0,7*n-3],[3,4*n*(2*n-1)]]
\end{verbatim}
$$\log(2)=\dfrac{3}{4+\dfrac{4}{11+\dfrac{24}{18+\dfrac{60}{25+\dfrac{112}{32+\dfrac{180}{39+\ddots}}}}}}$$
Convergence type $E$ with $E=-8$, $P=1/2$, and $C=\sqrt{\pi}/3$, so that
$$\log(2)-\dfrac{p(n)}{q(n)}\sim(-1)^n\dfrac{\sqrt{\pi}/3}{2^{3n}n^{1/2}}\;.$$
$$A=1-(23/72)/n+(347/3456)/n^2+(6625/248832)/n^3-(6868967/71663616)/n^4+\cdots$$
Series:
$$\log(2)=\dfrac{3}{4}\sum_{n\ge0}(-1)^n\dfrac{n!}{(3/2)_n}2^{-3n}$$
Parametric families up to the trivial change $n\mapsto n+j$,
with $u\ge0$, $k\ge0$:
\begin{verbatim}
[()->log(2),7*n-3+8*u+9*k,4*n*(2*n+2*u-1)]
\end{verbatim}
Convergence type $E$ with $E=-8$ and $P=2k+u+1/2$.

Ap\'ery accelerates to \ref{5.3.84}.
\end{cf}

\smallskip

\begin{cf}\label{5.3.73}{\ }
\begin{verbatim}
[()->log(2),[0,n*(7*n-3)],[3,4*n^2*(n+1)*(2*n-1)]]
\end{verbatim}
$$\log(2)=\dfrac{3}{4+\dfrac{8}{22+\dfrac{144}{54+\dfrac{720}{100+\dfrac{2240}{160+\dfrac{5400}{234+\ddots}}}}}}$$
Convergence type $E$ with $E=-8$, $P=1/2$, and $C=\sqrt{\pi}/3$, so that
$$\log(2)-\dfrac{p(n)}{q(n)}\sim(-1)^n\dfrac{\sqrt{\pi}/3}{2^{3n}n^{1/2}}\;.$$
$$A=1-(23/72)/n+(347/3456)/n^2+(6625/248832)/n^3-(6868967/71663616)/n^4+\cdots$$
Series:
$$\log(2)=\dfrac{3}{4}\sum_{n\ge0}(-1)^n\dfrac{n!}{(3/2)_n}2^{-3n}$$
Parametric family for $u\ge0$:
\begin{verbatim}
[()->log(2),7*n^2+(22*u-3)*n+u*(15*u-7),4*(n+u)^2*(n+3*u+1)*(2*n-1)]
\end{verbatim}
Convergence type $E$ with $E=-8$ and $P=u+1/2$.
\end{cf}

This CF is of course termwise identical to the previous one, but with
a completely different parametric family.

\smallskip

\begin{cf}\label{5.3.74}{\ }
\begin{verbatim}
[()->log(2),[1/2,-5,21*n^2+n-12],[-1,-14,4*n*(2*n-1)*(3*n-2)*(3*n+4)]]
\end{verbatim}
$$\log(2)=1/2-\dfrac{1}{-5-\dfrac{14}{74+\dfrac{960}{180+\dfrac{5460}{328+\dfrac{17920}{518+\dfrac{44460}{750+\ddots}}}}}}$$
Convergence type $E$ with $E=-8$, $P=3/2$, and $C=2\sqrt{\pi}/27$, so that
$$\log(2)-\dfrac{p(n)}{q(n)}\sim(-1)^n\dfrac{\sqrt{\pi}/27}{2^{3n-1}n^{3/2}}\;.$$
$$A=1-(29/24)/n+(4315/3456)/n^2-(293855/248832)/n^3+\cdots$$
Series:
$$\log(2)=\dfrac{1}{2}+\dfrac{1}{2}\sum_{n\ge0}(-1)^n\dfrac{(3n+4)n!}{(3n+5)(3n+2)(3/2)_n}2^{-3n}$$
\end{cf}

\smallskip

\begin{cf}\label{5.3.75}{\ }
\begin{verbatim}
[()->log(2),[0,10*n-5],[3,-9*n^2]]
\end{verbatim}
$$\log(2)=\dfrac{3}{5-\dfrac{9}{15-\dfrac{36}{25-\dfrac{81}{35-\dfrac{144}{45-\dfrac{225}{55+\ddots}}}}}}$$
Convergence type $E$ with $E=9$, $P=0$, and $C=\pi/3$, so that
$$\log(2)-\dfrac{p(n)}{q(n)}\sim\dfrac{\pi}{3^{2n+1}}\;.$$
$$A=1-(5/16)/n+(105/512)/n^2-(1725/8192)/n^3+(118515/524288)/n^4+\cdots$$
Parametric families up to the trivial change $n\mapsto n+j$,
with $u\ge0$ and $k\ge0$:
\begin{verbatim}
[()->log(2),10*n-5+u+8*k,-9*n*(n+u)]
\end{verbatim}
Convergence type $E$ with $E=9$ and $P=-u+2k$.

Ap\'ery accelerates to \ref{5.3.83}.
\end{cf}

\smallskip

\begin{cf}\label{5.3.76}{\ }
\begin{verbatim}
[()->log(2),[0,9,20*n-12],[6,-9*(2*n-1)^2]]
\end{verbatim}
$$\log(2)=\dfrac{6}{9-\dfrac{9}{28-\dfrac{81}{48-\dfrac{225}{68-\dfrac{441}{88-\dfrac{729}{108-\ddots}}}}}}$$
Convergence type $E$ with $E=9$, $P=1$, and $C=3/8$, so that
$$\log(2)-\dfrac{p(n)}{q(n)}\sim\dfrac{1/8}{3^{2n-1}n}\;.$$
$$A=1-(5/8)/n+(17/32)/n^2-(175/256)/n^3+(635/512)/n^4-\cdots$$
Series:
$$\log(2)=\dfrac{2}{3}\sum_{n\ge0}\dfrac{3^{-2n}}{2n+1}$$
Parametric families up to the trivial change $n\mapsto n+j$, with
$u\ge0$ and $k\ge0$:
\begin{verbatim}
[()->log(2),20*n-12+18*u+16*k,-9*(2*n-1)*(2*n+2*u-1)]
\end{verbatim}
Convergence type $E$ with $E=9$ and $P=u+1+2k$.

Ap\'ery accelerates to \ref{5.3.85}.
\end{cf}

\smallskip

\begin{cf}\label{5.3.77}{\ }
\begin{verbatim}
[()->log(2),[[1,6],[10*n+3,10*n^2+15*n+6]],
           [[-1,-36],[n*(4*n^2-1),-4*(n+1)*(2*n+3)^2]]]
\end{verbatim}
$$\log(2)=1-\dfrac{1}{6-\dfrac{36}{13+\dfrac{3}{31-\dfrac{200}{23+\dfrac{30}{76-\dfrac{588}{33+\ddots}}}}}}$$
Convergence type $E$ with $E=-i((1+\sqrt{5})/2)^5$, $P=0$, and $C=-2\pi/((1+\sqrt{5})/2)^6$, so that
$$\log(2)-\dfrac{p(n)}{q(n)}\sim(-1)^{\lfloor (n-1)/2\rfloor}\dfrac{2\pi}{((1+\sqrt{5})/2)^{5n+6}}\;.$$
$$A=1-(d/10)/n+(7d/50+1/40)/n^2+(-1719d/10000-7/100)/n^3+\cdots$$
\end{cf}

\smallskip

Note that this is the \emph{unique} CF among those given in this dictionary
for which $a(2n)$ and $a(2n+1)$ do not have the same degrees, although it
is of course possible to construct others (\ref{5.12.17} has different
degrees for $b(2n)$ and $b(2n+1)$).

\smallskip

\begin{cf}\label{5.3.78}{\ }
\begin{verbatim}
[()->log(2),[[0,4],[20*n^2-16*n+3,20*n^2-4*n-3]],
            [[4,12],[-4*n^2*(2*n-1)^2,9*(2*n-1)*(2*n+1)^3]]]
\end{verbatim}
$$\log(2)=\dfrac{4}{4+\dfrac{12}{7-\dfrac{4}{13+\dfrac{243}{51-\dfrac{144}{69+\dfrac{3375}{135-\ddots}}}}}}$$
Convergence type $E$ with $E=-((1+\sqrt{5})/2)^5$, $P=0$, and
$C=2\pi/((1+\sqrt{5})/2)^4$, so that
$$\log(2)-\dfrac{p(n)}{q(n)}\sim(-1)^n\dfrac{2\pi}{((1+\sqrt{5})/2)^{5n+4}}\;.$$
$$A=1-(d/10)/n+(3d/50+1/40)/n^2-(119d/10000+3/100)/n^3+\cdots$$
\end{cf}

\smallskip

\begin{cf}\label{5.3.79}{\ }
\begin{verbatim}
[()->log(2),[0,29*n^2-29*n+8],[5,-6*n^2*(9*n^2-1)]]
\end{verbatim}
$$\log(2)=\dfrac{5}{8-\dfrac{48}{66-\dfrac{840}{182-\dfrac{4320}{356-\dfrac{13728}{588-\dfrac{33600}{878-\ddots}}}}}}$$
Convergence type $E$ with $E=27/2$, $P=0$, and $C=2\pi/3^{3/2}$, so that
$$\log(2)-\dfrac{p(n)}{q(n)}\sim\dfrac{2\pi/3^{3/2}}{(27/2)^n}\;.$$
$$A=1-(88/225)/n+(13772/50625)/n^2-(7861604/34171875)/n^3+\cdots$$
\end{cf}

\smallskip

It is possible that this CF together with several above is part of a family
or pattern. Excluding CFs which can be simplified by $2n-1$ (i.e.,
$a(n)$ multiple of $2n-1$ and $b(n)$ multiple of $4n^2-1$), note the
following:
\begin{verbatim}
[[1/2,4*n^2-4*n+3],[1/2,-n^2*(4*n^2-1)]]         /* P^+=2 */
[[0,2,10*n^2-10*n+3],[1,-10,-n^2*(16*n^2-1)]]    /* E=4 */
[[0,27*n^2-27*n+8],[3,-18*n^2*(9*n^2-1)]]        /* E=2 */
[[0,29*n^2-29*n+8],[5,-6*n^2*(9*n^2-1)]]         /* E=27/2 */
[[0,59*n^2-59*n+20],[5,-24*n^2*(36*n^2-1)]]      /* E=32/27 */
\end{verbatim}
corresponding to \ref{5.3.61}, \ref{5.3.71}, \ref{5.3.67},
\ref{5.3.79}, and \ref{5.3.62} respectively.

\smallskip

\begin{cf}\label{5.3.80}{\ }
\begin{verbatim}
[()->log(2),[[0,5],[17*n-4,17*n+5]],
            [[3,-9],[8*n*(2*n-1),-9*(n+1)*(2*n+1)]]]
\end{verbatim}
$$\log(2)=\dfrac{3}{5-\dfrac{9}{13+\dfrac{8}{22-\dfrac{54}{30+\dfrac{48}{39-\dfrac{135}{47+\ddots}}}}}}$$
Convergence type $E$ with $E=-12\sqrt{2}$, $P=0$, and $C=\pi/3$, so that
$$\log(2)-\dfrac{p(n)}{q(n)}\sim(-1)^n\dfrac{\pi/3}{288^{n/2}}\;.$$
$$A=1-(143/578)/n+(83553/668168)/n^2+\cdots$$
\end{cf}

There exist several CFs of the same type as the one above.

\smallskip

\begin{cf}\label{5.3.81}{\ }
\begin{verbatim}
[()->log(2),[11/16,364*n^3+229*n^2-73*n-40],
            [11/4,6*n*(2*n+1)*(3*n+1)*(3*n+2)*(7*n-3)*(7*n+11)]]
\end{verbatim}
$$\log(2)=\dfrac{11}{16}+\dfrac{11/4}{480+\dfrac{25920}{3642+\dfrac{924000}{11630+\dfrac{7983360}{26628+\dfrac{38329200}{50820+\ddots}}}}}$$
Convergence type $E$ with $E=-27$, $P=3/2$, and $C=\sqrt{3\pi}/432$, so that
$$\log(2)-\dfrac{p(n)}{q(n)}\sim(-1)^n\dfrac{\sqrt{3\pi}/16}{3^{3n+3}n^{3/2}}$$
$$A=1-(16/9)/n+(27445/10368)/n^2+\cdots$$
Series:
$$\log(2)=\dfrac{11}{16}-\dfrac{1}{16}\sum_{n\ge1}(-1)^n\dfrac{(7n+4)(1/2)_n(n-1)!}{(4/3)_n(5/3)_n}3^{-3n}$$
\end{cf}

\smallskip

\begin{cf}\label{5.3.82}{\ }
\begin{verbatim}
[()->log(2),[3/4,364*n^3-135*n^2-43*n+6],
            [-11,6*n*(2*n-1)*(3*n-1)*(3*n+1)*(7*n-3)*(7*n+11)]]
\end{verbatim}
$$\log(2)=\dfrac{3}{4}-\dfrac{11}{192+\dfrac{3456}{2292+\dfrac{346500}{8490+\dfrac{4147200}{20970+\dfrac{23423400}{41916+\ddots}}}}}$$
Convergence type $E$ with $E=-27$, $P=3/2$, and $C=-\sqrt{3\pi}/144$, so that
$$\log(2)-\dfrac{p(n)}{q(n)}\sim(-1)^{n+1}\dfrac{\sqrt{3\pi}/16}{3^{3n+2}n^{3/2}}$$
$$A=1-(11/18)/n+(325/10368)/n^2+\cdots$$
Series:
$$\log(2)=\dfrac{3}{4}-\dfrac{1}{192}\sum_{n\ge0}(-1)^n\dfrac{(7n+11)(1/2)_nn!}{(5/3)_n(7/3)_n}3^{-3n}$$
\end{cf}

\smallskip

\begin{cf}\label{5.3.83}{\ }
\begin{verbatim}
[()->log(2),[0,6*n-3],[2,-n^2]]
\end{verbatim}
$$\log(2)=\dfrac{2}{3-\dfrac{1}{9-\dfrac{4}{15-\dfrac{9}{21-\dfrac{16}{27-\dfrac{25}{33-\ddots}}}}}}$$
Convergence type $E$ with $E=(1+\sqrt{2})^4$, $P=0$, and
$C=2\pi/(3+2\sqrt{2})$, so that
$$\log(2)-\dfrac{p(n)}{q(n)}\sim\dfrac{2\pi}{(1+\sqrt{2})^{4n+2}}\;.$$
$$A=1-(3d/16)/n+(3d/32+9/256)/n^2+(-219d/4096-9/256)/n^3+\cdots$$
Parametric families up to the trivial change $n\mapsto n+j$, with
$u\ge0$:
\begin{verbatim}
[()->log(2),3*(2*n-1+2*u),-n*(n+2*u)]
\end{verbatim}
Convergence type $E$ with $E=(1+\sqrt{2})^4$ and $P=0$.
\end{cf}

\smallskip

\begin{cf}\label{5.3.84}{\ }
\begin{verbatim}
[()->log(2),[0,4*n-2],[[3/2,1],[n*(2*n-1),(n+1)*(2*n+1)]]]
\end{verbatim}
$$\log(2)=\dfrac{3/2}{2+\dfrac{1}{6+\dfrac{1}{10+\dfrac{6}{14+\dfrac{6}{18+\dfrac{15}{22+\ddots}}}}}}$$
Convergence type $E$ with $E=-(1+\sqrt{2})^4$, $P=0$, and $C=2\pi/(1+\sqrt{2})^2$, so that
$$\log(2)-\dfrac{p(n)}{q(n)}\sim(-1)^n\dfrac{2\pi}{(1+\sqrt{2})^{4n+2}}\;.$$
$$A=1-(3d/16)/n+(3d/32+9/256)/n^2+(-219d/4096-9/256)/n^3+\cdots$$
\end{cf}

\smallskip

\begin{cf}\label{5.3.85}{\ }
\begin{verbatim}
[()->log(2),[0,5,6*n-3],[[3,-6],[-(2*n-1)^2,-(2*n+2)^2]]]
\end{verbatim}
$$\log(2)=\dfrac{3}{5-\dfrac{6}{9-\dfrac{1}{15-\dfrac{16}{21-\dfrac{9}{27-\dfrac{36}{33-\ddots}}}}}}$$
Convergence type $E$ with $E=(1+\sqrt{2})^4$, $P=0$, and
$C=2\pi/(1+\sqrt{2})^2$, so that
$$\log(2)-\dfrac{p(n)}{q(n)}\sim\dfrac{2\pi}{(1+\sqrt{2})^{4n+2}}\;.$$
$$A=1-(3d/16)/n+(3d/32+9/256)/n^2+(-219d/4096-9/256)/n^3+\cdots$$
\end{cf}

\smallskip

\begin{cf}\label{5.3.86}{\ }
\begin{verbatim}
[()->log(2),[1,12*n^2-2],[-3,-n*(n+1)*(2*n-1)*(2*n+3)]]
\end{verbatim}
$$\log(2)=1-\dfrac{3}{10-\dfrac{10}{46-\dfrac{126}{106-\dfrac{540}{190-\dfrac{1540}{298-\dfrac{3510}{430-\ddots}}}}}}$$
Convergence type $E$ with $E=(1+\sqrt{2})^4$, $P=0$, and $C=-2\pi/(1+\sqrt{2})^4$, so that
$$\log(2)-\dfrac{p(n)}{q(n)}\sim-\dfrac{2\pi}{(1+\sqrt{2})^{4n+4}}\;.$$
$$A=1+(5d/16)/n+(-5d/16+25/256)/n^2+(1485d/4096-25/128)/n^3+\cdots$$
\end{cf}

\smallskip

One easily finds numerous similar continued fractions for $\log(2)$ with
convergence type $E=\pm(1+\sqrt{2})^4$.

\smallskip

\begin{cf}\label{5.3.87}{\ }
\begin{verbatim}
[()->log(2),[0,220*n^3-506*n^2+334*n-63],
            [-21/2,n^2*(2*n-1)^2*(10*n-13)*(10*n+7)]]
\end{verbatim}
$$\log(2)=-\dfrac{21/2}{-15-\dfrac{51}{341+\dfrac{6804}{2325+\dfrac{141525}{7257+\dfrac{994896}{16457+\dfrac{4270725}{31245+\ddots}}}}}}$$
Convergence type $E$ with $E=-((1+\sqrt{5})/2)^{10}$, $P=0$, and
$C=2\pi/((1+\sqrt{5})/2)^4$, so that
$$\log(2)-\dfrac{p(n)}{q(n)}\sim(-1)^n\dfrac{2\pi}{((1+\sqrt{5})/2)^{10n+4}}\;.$$
$$A=1-(d/20)/n+(3d/200+1/160)/n^2-(119d/80000+3/800)/n^3+\cdots$$
\end{cf}

\smallskip

Many of the above CFs are special cases of the following infinite families:

\begin{cf}\label{5.3.88}{\ }
\begin{verbatim}
[()->log(2),[0,(2^k+1)*(2*n-1)],[2*(2^k-1)/k,-(2^k-1)^2*n^2]]
\end{verbatim}
$$\log(2)=\dfrac{2(2^k-1)/k}{2^k+1-\dfrac{2^{2k}-2\cdot 2^k+1}{3\cdot 2^k+3-\dfrac{4\cdot 2^{2k}-8\cdot 2^k+4}{5\cdot 2^k+5-\dfrac{9\cdot 2^{2k}-18\cdot 2^k+9}{7\cdot 2^k+7-\ddots}}}}$$
Convergence type $E$ with $E=((1+2^{k/2})^4)/(2^k-1)^2$, $P=0$, and $C=...$, so
that
$$\log(2)-\dfrac{p(n)}{q(n)}\sim\dfrac{C}{((1+2^{k/2})^2/(2^k-1))^{2n}}\;.$$
\end{cf}

\smallskip

\begin{cf}\label{5.3.89}{\ }
\begin{verbatim}
[()->log(2),[0,(2^(k+1)-1)*n-(2^k-1)],[(2^k-1)/k,-2^k*(2^k-1)*n^2]]
\end{verbatim}
$$\log(2)=\dfrac{(2^k-1)/k}{2^k-\dfrac{2^{2k}-2^k}{3\cdot 2^k-1-\dfrac{4\cdot 2^{2k}-4\cdot 2^k}{5\cdot 2^k-2-\dfrac{9\cdot 2^{2k}-9\cdot 2^k}{7\cdot 2^k-3-\ddots}}}}$$
Convergence type $E$ with $E=2^k/(2^k-1)$, $P=1$, and $C=...$, so that
$$\log(2)-\dfrac{p(n)}{q(n)}\sim\dfrac{C}{(2^k/(2^k-1))^nn}\;.$$
\end{cf}

\smallskip

\begin{cf}\label{5.3.90}
\begin{verbatim}
[()->log(2),[0,(2^k+1)^2,4*((2^(2*k)+1)*n-2^(2*k)+2^k-1)],
            [2*(2^(2*k)-1)/k,-(2^(2*k)-1)^2*(2*n-1)^2]]
\end{verbatim}
$$\log(2)=\dfrac{(2^{2k+1}-2)/k}{2^{2k}+2^{k+1}+1-\dfrac{2^{4k}-2^{2k+1}+1}{2^{2k+2}+2^{k+2}+4-\dfrac{9\cdot2^{4k}-92^{2k+1}+9}{2^{2k+3}+2^{k+2}+8-\dfrac{25\cdot2^{4k}-25\cdot2^{2k+1}+25}{3\cdot2^{2k+2}+2^{k+2}+12-\ddots}}}}$$
Convergence type $E$ with $E=((2^k+1)/(2^k-1))^2$, $P=1$, and
$C=(2^{2k}-1)/(k\cdot2^{k+2})$, so that
$$\log(2)-\dfrac{p(n)}{q(n)}\sim\dfrac{(2^{2k}-1)/(k\cdot2^{k+2})}{((2^k+1)/(2^k-1))^{2n}n}$$
Series:
$$\log(2)=\dfrac{2}{k}\sum_{n\ge0}\dfrac{((2^k-1)/(2^k+1))^{2n+1}}{2n+1}$$
\end{cf}

\smallskip

We have given continued fractions for $\log(2)$ as examples, but of course
we could also give them for other values of the logarithm.

\smallskip

\begin{cf}\label{5.3.91}{\ }
\begin{verbatim}
[()->log(1+sqrt(2))/sqrt(2),[0,3,32*(n-1)],[2,(4*n-3)^2*(4*n-1)^2]]
\end{verbatim}
$$\dfrac{\log(1+\sqrt{2})}{\sqrt{2}}=\dfrac{2}{3+\dfrac{9}{32+\dfrac{1225}{64+\dfrac{9801}{96+\dfrac{38025}{128+\dfrac{104329}{160+\ddots}}}}}}$$
Convergence type $P^-$ with $P=2$ and $C=1/16$, so that
$$\dfrac{\log(1+\sqrt{2})}{\sqrt{2}}-\dfrac{p(n)}{q(n)}\sim(-1)^n\dfrac{1/16}{n^2}\;.$$
$$A=1-(11/16)/n^2+(361/256)/n^4-(24611/4096)/n^6+\cdots$$
Series:
$$\dfrac{\log(1+\sqrt{2})}{\sqrt{2}}=2\sum_{n\ge1}\dfrac{(-1)^{n+1}}{(4n-1)(4n-3)}$$
Parametric families up to the trivial change $n\mapsto n+j$, with $k\ge0$:
\begin{verbatim}
[()->log(1+sqrt(2))/sqrt(2),32*(2*k+1)*n,(4*n+2*u-1)^2*(4*n-2*u+5)^2]
\end{verbatim}
Convergence type $P^-$ with $P=4k+2$.

Ap\'ery accelerates to \ref{5.3.97}.
\end{cf}

\smallskip

\begin{cf}\label{5.3.92}{\ }
\begin{verbatim}
[()->log(1+sqrt(2))/sqrt(2),[0,32*n^2-32*n+11],[1,-16*n^2*(16*n^2-1)]]
\end{verbatim}
$$\dfrac{\log(1+\sqrt{2})}{\sqrt{2}}=\dfrac{1}{11-\dfrac{240}{75-\dfrac{4032}{203-\dfrac{20592}{395-\dfrac{65280}{651-\dfrac{159600}{971-\ddots}}}}}}$$
Convergence type $L$ with $C=...$, so that
$$\dfrac{\log(1+\sqrt{2})}{\sqrt{2}}-\dfrac{p(n)}{q(n)}\sim\dfrac{C}{\log(n)}\;.$$
\end{cf}

This is only given as an example of a nontrivial CF of type
$L$, and is the case $k=0$ of the parametric family given in the next CF.

\smallskip

\begin{cf}\label{5.3.93}{\ }
\begin{verbatim}
[()->log(1+sqrt(2))/sqrt(2),[2/3,32*n^2-32*n+27],[-1,-16*n^2*(16*n^2-1)]]
\end{verbatim}
$$\dfrac{\log(1+\sqrt{2})}{\sqrt{2}}=2/3-\dfrac{1}{27-\dfrac{240}{91-\dfrac{4032}{219-\dfrac{20592}{411-\dfrac{65280}{667-\dfrac{159600}{987-\ddots}}}}}}$$
Convergence type $P^+$ with $P=2$ and $C=-9\pi^2/2^{25/2}$, so that
$$\dfrac{\log(1+\sqrt{2})}{\sqrt{2}}-\dfrac{p(n)}{q(n)}\sim-\dfrac{9\pi^2/2^{25/2}}{n^2}\;.$$
$$A=1-1/n+\cdots$$
Parametric family for $k\ge0$:
\begin{verbatim}
[()->log(1+sqrt(2))/sqrt(2),32*n^2-32*n+11+16*k^2,-16*n^2*(16*n^2-1)]
\end{verbatim}
Convergence type $P^+$ with $P=2k$.
\end{cf}

\smallskip

\begin{cf}\label{5.3.94}{\ }
\begin{verbatim}
[()->log(1+sqrt(2))/sqrt(2),[0,6*n-3],[1,-8*n^2]]
\end{verbatim}
$$\dfrac{\log(1+\sqrt{2})}{\sqrt{2}}=\dfrac{1}{3-\dfrac{8}{9-\dfrac{32}{15-\dfrac{72}{21-\dfrac{128}{27-\dfrac{200}{33-\ddots}}}}}}$$
Convergence type $E$ with $E=2$, $P=0$, and $C=\pi/4$, so that
$$\dfrac{\log(1+\sqrt{2})}{\sqrt{2}}-\dfrac{p(n)}{q(n)}\sim\dfrac{\pi/4}{2^n}$$
$$A=1-(3/4)/n+(21/32)/n^2-(465/128)/n^3+\cdots$$
Parametric family:
\begin{verbatim}
[()->log(1+sqrt(2))/sqrt(2),6*n-3+2*k,-8*n^2]
\end{verbatim}
Convergence type $E$ with $E=2$ and $P=2k$.
\end{cf}
      
\smallskip

\begin{cf}\label{5.3.95}{\ }
\begin{verbatim}
[()->log(1+sqrt(2))/sqrt(2),[0,4*n-2],[1,-2*n^2]]
\end{verbatim}
$$\dfrac{\log(1+\sqrt{2})}{\sqrt{2}}=\dfrac{1}{2-\dfrac{2}{6-\dfrac{8}{10-\dfrac{18}{14-\dfrac{32}{18-\dfrac{50}{22-\ddots}}}}}}$$
Convergence type $E$ with $E=(1+\sqrt{2})^2$, $P=0$, and $C=\pi/\sqrt{2}/(1+\sqrt{2})$, so that
$$\dfrac{\log(1+\sqrt{2})}{\sqrt{2}}-\dfrac{p(n)}{q(n)}\sim\dfrac{\pi/\sqrt{2}}{(1+\sqrt{2})^{2n+1}}\;.$$
$$A=1-(d/4)/n+(d/8+1/16)/n^2+(-3d/16-1/16)/n^3+(7d/32+55/512)/n^4+\cdots$$
Parametric families up to the trivial change $n\mapsto n+j$, with
$u\ge0$:
\begin{verbatim}
[()->log(1+sqrt(2))/sqrt(2),4*n+4*u-2,-2*n*(n+2*u)]
\end{verbatim}
Convergence type $E$ with $E=(1+\sqrt{2})^2$ and $P=0$.
\end{cf}

\smallskip

\begin{cf}\label{5.3.96}{\ }
\begin{verbatim}
[()->log(1+sqrt(2))/sqrt(2),[0,2*n-1],
                           [[1,2],[2*n*(2*n-1),(2*n+1)*(2*n+2)]]]
\end{verbatim}
$$\dfrac{\log(1+\sqrt{2})}{\sqrt{2}}=\dfrac{1}{1+\dfrac{2}{3+\dfrac{2}{5+\dfrac{12}{7+\dfrac{12}{9+\dfrac{30}{11+\ddots}}}}}}$$
Convergence type $E$ with $E=-(1+\sqrt{2})^2$, $P=0$, and $C=\pi/\sqrt{2}/(1+\sqrt{2})$, so that
$$\dfrac{\log(1+\sqrt{2})}{\sqrt{2}}-\dfrac{p(n)}{q(n)}\sim(-1)^n\dfrac{\pi/\sqrt{2}}{(1+\sqrt{2})^{2n+1}}\;.$$
$$A=1-(d/4)/n+(d/8+1/16)/n^2+(-3d/16-1/16)/n^3+\cdots$$
\end{cf}

\smallskip

\begin{cf}\label{5.3.97}{\ }
\begin{verbatim}
[()->log(1+sqrt(2))/sqrt(2),
[[0,22],[80*n^2-32*n+3,80*n^2+64*n+19]],
[[2,-960],[(4*n-3)^2*(4*n-1)^2,-64*(n+1)^2*(4*n+3)*(4*n+5)]]]
\end{verbatim}
$$\dfrac{\log(1+\sqrt{2})}{\sqrt{2}}=\dfrac{2}{22-\dfrac{960}{51+\dfrac{9}{163-\dfrac{16128}{259+\dfrac{1225}{467-\dfrac{82368}{627+\ddots}}}}}}$$
Convergence type $E$ with $E=-i((1+\sqrt{5})/2)^5$, $P=0$, and
$C=\pi/\sqrt{2}/((1+\sqrt{5})/2)^2$, so that
$$\dfrac{\log(1+\sqrt{2})}{\sqrt{2}}-\dfrac{p(n)}{q(n)}\sim(-1)^{\lfloor (n+1)/2\rfloor}\dfrac{\pi/\sqrt{2}}{((1+\sqrt{5})/2)^{5n+2}}\;.$$
$$A=1-(3d/40)/n+(3d/200+9/640)/n^2+(3759d/640000-9/1600)/n^3+\cdots$$
\end{cf}

\smallskip

\begin{cf}\label{5.3.98}{\ }
\begin{verbatim}
[()->log(1+sqrt(2))/sqrt(2),[3/5,176*n^2-5],
                            [4,16*n*(n+1)*(4*n-1)*(4*n+5)]]
\end{verbatim}
$$\dfrac{\log(1+\sqrt{2})}{\sqrt{2}}=3/5+\dfrac{4}{171+\dfrac{864}{699+\dfrac{8736}{1579+\dfrac{35904}{2811+\dfrac{100800}{4395+\dfrac{228000}{6331+\ddots}}}}}}$$
Convergence type $E$ with $E=-((1+\sqrt{5})/2)^{10}$, $P=0$, and
$C=\pi/\sqrt{2}/(1+\sqrt{5}/2)^{10}$, so that
$$\dfrac{\log(1+\sqrt{2})}{\sqrt{2}}-\dfrac{p(n)}{q(n)}\sim(-1)^n\dfrac{\pi/\sqrt{2}}{((1+\sqrt{5})/2)^{10(n+1)}}$$
$$A=1+(9d/80)/n+(-9d/80+81/2560)/n^2+\cdots$$
\end{cf}

\smallskip

We have given continued fractions for $\log(1+\sqrt{2})/\sqrt{2}$ as examples,
but of course we could give them for many other similar quantities.

\medskip

\section{Constants: $\pi^2$, $G$, and Periods of Degree $2$}

\medskip

\begin{cf}\label{5.4.1}{\ }
\begin{verbatim}
[()->Pi^2,[0,4*n^2-5*n+2],[4,-2*n^3*(2*n-1)]]
\end{verbatim}
$$\pi^2=\dfrac{4}{1-\dfrac{2}{8-\dfrac{48}{23-\dfrac{270}{46-\dfrac{896}{77-\dfrac{2250}{116-\ddots}}}}}}$$
Convergence type $P^+$ with $P=1/2$ and $C=4\sqrt{\pi}$, so that
$$\pi^2-\dfrac{p(n)}{q(n)}\sim\dfrac{4\sqrt{\pi}}{n^{1/2}}\;.$$
$$A=1-(5/24)/n+(21/640)/n^2+(223/21504)/n^3-(671/98304)/n^4+\cdots$$
Series:
$$\pi^2=4\sum_{n\ge0}\dfrac{n!}{(n+1)(3/2)_n}$$
Parametric families up to the trivial change $n\mapsto n+j$, with
all parameters nonnegative integers:
\begin{verbatim}
[()->Pi^2,4*n^2+(6*u+2*v+2*w-5)*n+2*k^2+(2*u+2*v+2*w+1)*k
                                +2*u^2+(2*v+2*w-3)*u+v*(2*w-1)+2,
        -2*n*(n+u)*(n+u+v)*(2*n+2*u+2*w-1)]
\end{verbatim}
Convergence type $P^+$ with $P=u+v+w+2k+1/2$.

Can be Ap\'ery accelerated to a CF with convergence type
$E=-((1+\sqrt{5})/2)^5$, formula too complicated to give here.
\end{cf}

\smallskip

\begin{cf}\label{5.4.2}{\ }
\begin{verbatim}
[()->Pi^2,[0,1,8*n^2-16*n+10],[8,-(2*n-1)^4]]
\end{verbatim}
$$\pi^2=\dfrac{8}{1-\dfrac{1}{10-\dfrac{81}{34-\dfrac{625}{74-\dfrac{2401}{130-\dfrac{6561}{202-\ddots}}}}}}$$
Convergence type $P^+$ with $P=1$ and $C=2$, so that
$$\pi^2-\dfrac{p(n)}{q(n)}\sim\dfrac{2}{n}\;.$$
$$A=1-(1/12)/n^2+(7/240)/n^4-(31/1344)/n^6+(127/3840)/n^8+\cdots$$
Series:
$$\pi^2=8\sum_{n\ge1}\dfrac{1}{(2n-1)^2}$$
Parametric families up to the trivial change $n\mapsto n+j$, with
$v$, $w$, and $k$ nonnegative integers and $-u\le k+\min(v,w)$:
\begin{verbatim}
[()->Pi^2,8*n^2+4*(3*u+v+w-4)*n+4*(k^2+u^2+u*k+v*k+w*k+u*v+u*w+v*w)
         +4*k-10*u-2*v-2*w+10,
         -(2*n-1)*(2*n-1+2*u)*(2*n-1+2*u+2*v)*(2*n-1+2*u+2*w)]
\end{verbatim}
Convergence type $P^+$ with $P=u+v+w+2k+1$.
\end{cf}

\smallskip

\begin{cf}\label{5.4.3}{\ }
\begin{verbatim}
[()->Pi^2,[0,2*n^2-2*n+1],[6,-n^4]]
\end{verbatim}
$$\pi^2=\dfrac{6}{1-\dfrac{1}{5-\dfrac{16}{13-\dfrac{81}{25-\dfrac{256}{41-\dfrac{625}{61-\ddots}}}}}}$$
Convergence type $P^+$ with $P=1$ and $C=6$, so that
$$\pi^2-\dfrac{p(n)}{q(n)}\sim\dfrac{6}{n}\;.$$
$$A=1-(1/2)/n+(1/6)/n^2-(1/30)/n^4+(1/42)/n^6-(1/30)/n^8+\cdots$$
Series:
$$\pi^2=6\sum_{n\ge1}\dfrac{1}{n^2}$$
Parametric families up to the trivial change $n\mapsto n+j$, with
all parameters nonnegative integers:
\begin{verbatim}
[()->Pi^2,2*n^2+(3*u+v+w-2)*n+k^2+(u+v+w+1)*k+u^2+(v+w-1)*u+v*w+1,
         -n*(n+u)*(n+u+v)*(n+u+w)]
\end{verbatim}
Convergence type $P^+$ with $P=u+v+w+2k+1$.

After contraction, Ap\'ery accelerates to Ap\'ery's CF \ref{5.4.14}.
\end{cf}

\smallskip

\begin{cf}\label{5.4.4}{\ }
\begin{verbatim}
[()->Pi^2,[0,2*n-1],[12,n^4]]
\end{verbatim}
$$\pi^2=\dfrac{12}{1+\dfrac{1}{3+\dfrac{16}{5+\dfrac{81}{7+\dfrac{256}{9+\dfrac{625}{11+\ddots}}}}}}$$
Convergence type $P^-$ with $P=2$ and $C=6$, so that
$$\pi^2-\dfrac{p(n)}{q(n)}\sim(-1)^n\dfrac{6}{n^2}\;.$$
$$A=1-1/n+1/n^3-3/n^5+17/n^7-155/n^9+2073/n^{11}-38227/n^{13}+\cdots$$
Series:
$$\pi^2=12\sum_{n\ge1}\dfrac{(-1)^{n+1}}{n^2}$$
Parametric families up to the trivial change $n\mapsto n+j$, with
$u\ge0$, $v\ge0$, $k\ge0$ with $k\ne u+2v-2-2j$ for $j\ge0$:
\begin{verbatim}
[()->Pi^2,(k+1)*(2*n+2*u+2*v-1),n*(n+u)*(n+u+2*v)*(n+2*u+2*v)]
[()->Pi^2,(8*u+7)*n+(24*u^2+24*u+3),n*(n+2*u)*(n+4*u+1)*(n+6*u+4)]
[()->Pi^2,(8*u+7)*n+(24*u^2+42*u+18),n*(n+2*u+3)*(n+4*u+4)*(n+6*u+4)]
[()->Pi^2,(2*u^2+8*u+7)*n+(2*u+3)*(u^2+3*u+1),n*(n+u)*(n+u+1)*(n+2*u+4)]
[()->Pi^2,(2*u^2+8*u+7)*n+(2*u+3)*(u+2)*(u+3),n*(n+u+3)*(n+u+4)*(n+2*u+4)]
\end{verbatim}
Convergence type $P^-$ with $P=2k+2$, $8u+7$, $8u+7$, $2u^2+8u+7$,
and $2u^2+8u+7$ respectively.

After contraction, Ap\'ery accelerates to Ap\'ery's CF \ref{5.4.14}.
\end{cf}

\smallskip

\begin{cf}\label{5.4.5}{\ }
\begin{verbatim}
[()->Pi^2,[0,7*n^2-8*n+3],[27/2,-6*n^3*(2*n-1)]]
\end{verbatim}
$$\pi^2=\dfrac{27/2}{2-\dfrac{6}{15-\dfrac{144}{42-\dfrac{810}{83-\dfrac{2688}{138-\dfrac{6750}{207-\ddots}}}}}}$$
Convergence type $E$ with $E=4/3$, $P=3/2$, and $C=(27/2)\sqrt{\pi}$, so that
$$\pi^2-\dfrac{p(n)}{q(n)}\sim\dfrac{(27/2)\sqrt{\pi}}{(4/3)^nn^{3/2}}\;.$$
$$A=1-(47/8)/n+(6561/128)/n^2-(638517/1024)/n^3+(319608747/32768)/n^4+\cdots$$
Series:
$$\pi^2=\dfrac{27}{4}\sum_{n\ge0}\dfrac{n!}{(n+1)(3/2)_n}(3/4)^n$$
Parametric families, up to the trivial change $n\mapsto n+j$, with $u\ge0$:
\begin{verbatim}
[()->Pi^2,7*n^2+(15*u+3*v-8)*n-3*(2*v+1)*(3*u+3*v-1),
         -6*n*(n+u)*(n+3*u+3*v)*(2*n-1-6*v)]
[()->Pi^2,7*n^2+(15*u+3*v-3)*n-(6*v+8)*(3*u+3*v+2),
         -6*n*(n+u)*(n+3*u+3*v+3)*(2*n-5-6*v)]
[()->Pi^2,7*n^2+(2*u+3*v-6)*n-(5*u^2+(15*v)*u+(18*v^2-2)),
         -6*n*(n+u)*(n+u+3*v)*(2*n+1-4*u-6*v)]
[()->Pi^2,7*n^2+(2*u+3*v-8)*n-(5*u^2+(15*v+2)*u+(18*v^2+3*v-3)),
         -6*n*(n+u)*(n+u+3*v)*(2*n-1-4*u-6*v)]
\end{verbatim}

Apparently sporadic cases:

\begin{verbatim}
[()->Pi^2,7*n^2+17*n,-6*n*(n+3)^2*(2*n+1)]
[()->Pi^2,7*n^2+54*n+5,-6*n*(n+6)*(n+10)*(2*n+1)]
[()->Pi^2,7*n^2+31*n-60,-6*n*(n+1)*(n+10)*(2*n-3)]
[()->Pi^2,7*n^2-16*n+3,-6*n^3*(2*n-5)]
[()->Pi^2,7*n^2+14*n-36,-6*n*(n+3)*(n+5)*(2*n-7)]
[()->Pi^2,7*n^2-19*n-6,-6*n^2*(n+1)*(2*n-11)]
\end{verbatim}
Convergence types all with $E=4/3$ and varying half-integral $P$.
\end{cf}

\smallskip

\begin{cf}\label{5.4.6}{\ }
\begin{verbatim}
[()->Pi^2,[0,3*n^2-3*n+1],[8,-n^3*(2*n-1)]]
\end{verbatim}
$$\pi^2=\dfrac{8}{1-\dfrac{1}{7-\dfrac{24}{19-\dfrac{135}{37-\dfrac{448}{61-\dfrac{1125}{91-\ddots}}}}}}$$
Convergence type $E$ with $E=2$, $P=3/2$, and $C=8\sqrt{\pi}$, so that
$$\pi^2-\dfrac{p(n)}{q(n)}\sim\dfrac{\sqrt{\pi}}{2^{n-3}n^{3/2}}\;.$$
$$A=1-(23/8)/n+(1361/128)/n^2-(54941/1024)/n^3+(11410635/32768)/n^4+\cdots$$
Series:
$$\pi^2=8\sum_{n\ge0}\dfrac{n!}{(n+1)(3/2)_n}2^{-n}$$
Parametric families up to the trivial change $n\mapsto n+j$, with
all parameters nonnegative integers:
\begin{verbatim}
[()->Pi^2,3*n^2+(8*u+4*v-3)*n+(2*u-1)*(2*u+2*v-1),
         -n*(n+2*u)*(n+2*u+2*v)*(2*n+2*u-2*k-1)]
\end{verbatim}
Convergence type $E$ with $E=2$ and $P=u+2v+3/2+3k$.
\end{cf}

\smallskip

\begin{cf}\label{5.4.7}{\ }
\begin{verbatim}
[()->Pi^2,[0,5*n^2-4*n+1],[18,-2*n^3*(2*n-1)]]
\end{verbatim}
$$\pi^2=\dfrac{18}{2-\dfrac{2}{13-\dfrac{48}{34-\dfrac{270}{65-\dfrac{896}{106-\dfrac{2250}{157-\ddots}}}}}}$$
Convergence type $E$ with $E=4$, $P=3/2$, and $C=6\sqrt{\pi}$, so that
$$\pi^2-\dfrac{p(n)}{q(n)}\sim\dfrac{6\sqrt{\pi}}{2^{2n}n^{3/2}}\;.$$
$$A=1-(15/8)/n+(481/128)/n^2-(9749/1024)/n^3+(1008939/32768)/n^4+\cdots$$
Series:
$$\pi^2=9\sum_{n\ge0}\dfrac{n!}{(n+1)(3/2)_n}2^{-2n}$$
Parametric families up to the trivial change $n\mapsto n+j$, with $u\ge0$
and $v\ge0$:
\begin{verbatim}
[()->Pi^2,5*n^2+(14*u+21*v-4)*n+(3*u+3*v-1)*(3*u+6*v-1),
         -2*n*(n+u)*(n+u+3*v)*(2*n-1)]
[()->Pi^2,5*n^2+(14*u+21*v-6)*n+9*u^2+(27*v-8)*u+2*(3*v-1)^2,
         -2*n*(n+u)*(n+u+3*v)*(2*n+1)]
[()->Pi^2,5*n^2+(17*u+21*v+15)*n+2*(2*u+3*v+2)*(3*u+3*v+2),
         -2*n*(n+u)*(n+3*u+3*v+3)*(2*n+2*u+1)]
[()->Pi^2,5*n^2+(17*u+21*v-4)*n+(3*u+3*v-1)*(4*u+6*v-1),
         -2*n*(n+u)*(n+3*u+3*v)*(2*n+2*u-1)]
\end{verbatim}

Apparently sporadic cases:

\begin{verbatim}
[()->Pi^2,5*n^2+15*n+6,-2*n^2*(n+1)*(2*n-1)]
[()->Pi^2,5*n^2+78*n+261,-2*n*(n+3)*(n+7)*(2*n-1)]
[()->Pi^2,5*n^2+38*n+68,-2*n*(n+3)*(n+5)*(2*n+1)]
[()->Pi^2,5*n^2+4*n+1,-2*n^3*(2*n+3)]
[()->Pi^2,5*n^2+19*n+16,-2*n*(n+3)^2*(2*n+3)]
[()->Pi^2,5*n^2+45*n+60,-2*n*(n+1)*(n+10)*(2*n+3)]
\end{verbatim}
Convergence types all with $E=4$ and varying half-integral $P$.
\end{cf}

\smallskip

\begin{cf}\label{5.4.8}{\ }
\begin{verbatim}
[()->Pi^2,[0,(2*n-1)*(3*n^2-3*n+1)],[15,n^4*(16*n^2-1)]]
\end{verbatim}
$$\pi^2=\dfrac{15}{1+\dfrac{15}{21+\dfrac{1008}{95+\dfrac{11583}{259+\dfrac{65280}{549+\dfrac{249375}{1001+\ddots}}}}}}$$
Convergence type $E$ with $E=-4$, $P=0$, and $C=3\pi^2/\sqrt{2}$, so that
$$\pi^2-\dfrac{p(n)}{q(n)}\sim(-1)^n\dfrac{3\pi^2}{2^{2n+1/2}}\;.$$
$$A=1-(5/16)/n+(105/512)/n^2+(65/8192)/n^3-(89125/524288)/n^4+\cdots$$
Parametric family with $k\ge0$:
\begin{verbatim}
[()->Pi^2,(2*n-1)*(3*n^2-3*n+1)+2*k*(11*n^2-11*n+3),n^4*(4*(2*n-k)^2-1)]
\end{verbatim}
Convergence type $E$ with $E=-4$ and $P=5k$.
\end{cf}

\smallskip

\begin{cf}\label{5.4.9}{\ }
\begin{verbatim}
[()->Pi^2,[10,72,24*n^3+44*n^2+2*n+1],[-10,8*n*(2*n+1)^5]]
\end{verbatim}
$$\pi^2=10-\dfrac{10}{72+\dfrac{1944}{373+\dfrac{50000}{1051+\dfrac{403368}{2249+\dfrac{1889568}{4111+\dfrac{6442040}{6781+\ddots}}}}}}$$
Convergence type $E$ with $E=-4$, $P=5/2$, and $C=-1/(2\sqrt{\pi})$, so that
$$\pi^2-\dfrac{p(n)}{q(n)}\sim(-1)^{n+1}\dfrac{1/(2\sqrt{\pi})}{2^{2n}n^{5/2}}$$
$$A=1-(25/8)/n+(785/128)/n^2-(8435/1024)/n^3+\cdots$$
Series:
$$\pi^2=10-\dfrac{5}{4}\sum_{n\ge0}(-1)^n\dfrac{(3/2)_n}{(2n+3)^2(n+1)!}2^{-2n}$$
Parametric family for $k\ge0$:
\begin{verbatim}
[()->Pi^2,24*n^3+44*(2*k+1)*n^2+2*n+2*k+1,4*(2*n-k)*(2*n-k+1)*(2*n+1)^4]
\end{verbatim}
Convergence type $E$ with $E=-4$ and $P=5k+5/2$.
\end{cf}
          
\smallskip

\begin{cf}\label{5.4.10}{\ }
\begin{verbatim}
[()->Pi^2,[0,7*n^2-7*n+2],[24,8*n^4]]
\end{verbatim}
$$\pi^2=\dfrac{24}{2+\dfrac{8}{16+\dfrac{128}{44+\dfrac{648}{86+\dfrac{2048}{142+\dfrac{5000}{212+\ddots}}}}}}$$
Convergence type $E$ with $E=-8$, $P=0$, and $C=2\pi^2$, so that
$$\pi^2-\dfrac{p(n)}{q(n)}\sim(-1)^n\dfrac{\pi^2}{2^{3n-1}}\;.$$
$$A=1-(1/3)/n+(2/9)/n^2-(4/81)/n^3-(20/243)/n^4+(28/243)/n^5+\cdots$$
Parametric families:
\begin{verbatim}
[()->Pi^2,7*n^2+((48-20*v)*u-7)*n+2*(4*u-1)*((8-5*v)*u-1),
         8*n*(n+u*v)*(n+u*(v+2))*(n+4*u)]
[()->Pi^2,7*n^2+(24*v-20*u-7)*n+(4*v-2)*(4*v-5*u-1),
         8*n*(n+u)*(n+2*v)*(n+u+v)]
[()->Pi^2,7*n^2+(12*u+20*v-7)*n+(4*v-2)*(3*u+3*v-1),
         8*n*(n+2*v)*(n+3*u+2*v)*(n+3*u+3*v)] /* u>=-[2*v/3] */
[()->Pi^2,7*n^2+(30*u-7)*n+(3*u-1)*(9*u-2),8*n^2*(n+3*u)^2]
[()->Pi^2,7*n^2+(34*u-7)*n+27*u^2-17*u+2,8*n^2*(n+u)^2]
[()->Pi^2,7*n^2+(22*u-7)*n+(3*u-1)*(5*u-2),8*n*(n+u)^2*(n+3*u)]
[()->Pi^2,7*n^2+(58*u-7)*n+87*u^2-29*u+2,8*n*(n+u)^2*(n+3*u)]
[()->Pi^2,7*n^2+(48*u-23)*n+2*(4*u+1)*(8*u-7),8*n*(n+2)*(n+2*u+3)*(n+4*u+2)]
[()->Pi^2,7*n^2+(34*u+17)*n+3*(3*u+1)*(3*u+2),8*n*(n+1)*(n+u)*(n+u+2)]
\end{verbatim}
Convergence type $E$ with $E=-8$ and varying $P$.

Hundreds of additional sporadic CFs, probably more parametric families.
\end{cf}

\smallskip

\begin{cf}\label{5.4.11}{\ }
\begin{verbatim}
[()->Pi^2,[0,(2*n-1)*(13*n^2-13*n+4)],[42,3*n^4*(9*n^2-1)]]
\end{verbatim}
$$\pi^2=\dfrac{42}{4+\dfrac{24}{90+\dfrac{1680}{410+\dfrac{19440}{1120+\dfrac{109824}{2376+\dfrac{420000}{4334+\ddots}}}}}}$$
Convergence type $E$ with $E=-27$, $P=0$, and $C=4\pi^2/3^{1/2}$, so that
$$\pi^2-\dfrac{p(n)}{q(n)}\sim(-1)^n\dfrac{4\pi^2}{3^{3n+1/2}}\;.$$
$$A=1-(35/72)/n+(3745/10360)/n^2-(410417/2239488)/n^3+\cdots$$
\end{cf}

\smallskip

\begin{cf}\label{5.4.12}{\ }
\begin{verbatim}
[()->Pi^2,
[[12,18],[48*n^4+30*n^3+5*n^2,48*n^4+150*n^3+173*n^2+90*n+18]],
[[-30,-324],[-n^4*(4*n^2-1)*(4*n+1)*(4*n+5),-4*(n+1)^4*(2*n+3)^4]]]
\end{verbatim}
$$\pi^2=12-\dfrac{30}{18-\dfrac{324}{83-\dfrac{135}{479-\dfrac{40000}{1028-\dfrac{28080}{2858-\dfrac{777924}{4743+\ddots}}}}}}$$
Convergence type $E$ with $E=(1+\sqrt{2})^4$, $P=0$, and $C=-24\pi^2/(1+\sqrt{2})^5$, so that
$$\pi^2-\dfrac{p(n)}{q(n)}\sim-\dfrac{24\pi^2}{(1+\sqrt{2})^{4n+5}}\;.$$
$$A=1-(5d/16)/n+(25d/64+25/256)/n^2+(-1793d/4096-125/512)/n^3+\cdots$$
\end{cf}

\smallskip

\begin{cf}\label{5.4.13}{\ }
\begin{verbatim}
[()->Pi^2,[0,1365*n^4-3943*n^3+3879*n^2-1605*n+240],
          [-624,-8*n^3*(2*n-1)^3*(21*n-29)*(21*n+13)]]
\end{verbatim}
$$\pi^2=-\dfrac{624}{-64+\dfrac{2176}{2842-\dfrac{1235520}{34440-\dfrac{69768000}{152972-\dfrac{936911360}{449440-\dfrac{6537672000}{1047606-\ddots}}}}}}$$
Convergence type $E$ with $E=64$, $P=1/2$, and $C=2\pi^{3/2}$, so that
$$\pi^2-\dfrac{p(n)}{q(n)}\sim\dfrac{\pi^{3/2}}{2^{6n-1}n^{1/2}}$$
$$A=1-(37/72)/n+(1139/3456)/n^2+\cdots$$
Series:
$$\pi^2=\dfrac{3}{4}\sum_{n\ge0}\dfrac{(21n+13)n!^3}{(3/2)_n^3}2^{-6n}$$
\end{cf}
          
\smallskip

\begin{cf}\label{5.4.14}{\ }
\begin{verbatim}
[()->Pi^2,[0,11*n^2-11*n+3],[30,n^4]]
\end{verbatim}
$$\pi^2=\dfrac{30}{3+\dfrac{1}{25+\dfrac{16}{69+\dfrac{81}{135+\dfrac{256}{223+\dfrac{625}{333+\ddots}}}}}}$$
Convergence type $E$ with $E=-((1+\sqrt{5})/2)^{10}$, $P=0$, and $C=24\pi^2/((1+\sqrt{5})/2)^5$, so that
$$\pi^2-\dfrac{p(n)}{q(n)}\sim(-1)^n\dfrac{24\pi^2}{((1+\sqrt{5})/2)^{10n+5}}\;.$$
$$A=1-(d/5)/n+(d/10+1/10)/n^2+(-51d/1250-1/10)/n^3+(7d/625+30/625)/n^4+\cdots$$
\end{cf}

This is the famous continued fraction due to R.~Ap\'ery.

\smallskip

Note that there exist continued fractions of period $2$ for $\pi^2$ coming from
those of $\psi'(z)$, with convergence in $C/((1+\sqrt{5})/2)^{5n}$, which have
little interest since their contraction gives Ap\'ery's CF.

\smallskip

\begin{cf}\label{5.4.15}{\ }
\begin{verbatim}
[()->lfun(-3,2),[0,(2*n-1)*(9*n^2-9*n+4)],[2/3,-9*n^4*(9*n^2-1)]]
\end{verbatim}
$$L(\chi_{-3},2)=\dfrac{2/3}{4-\dfrac{72}{66-\dfrac{5040}{290-\dfrac{58320}{784-\dfrac{329472}{1656-\dfrac{1260000}{3014-\ddots}}}}}}$$
Convergence type $L$ with $C=...$, so that
$$L(\chi_{-3},2)-\dfrac{p(n)}{q(n)}\sim\dfrac{C}{\log(n)}\;.$$
$A=1+\cdots$
\end{cf}

This is only given as an example of a nontrivial CF of type
$L$, and is the case $k=0$ of the parametric family given in the next CF.

\smallskip

\begin{cf}\label{5.4.16}{\ }
\begin{verbatim}
[()->lfun(-3,2),[3/4,(2*n-1)*(9*n^2-9*n+22)],[2/3,-9*n^4*(9*n^2-1)]]
\end{verbatim}
$$L(\chi_{-3},2)=\dfrac{3}{4}+\dfrac{2/3}{22-\dfrac{72}{120-\dfrac{5040}{380-\dfrac{58320}{910-\dfrac{329472}{1818-\dfrac{1260000}{3212-\ddots}}}}}}$$
Convergence type $P^+$ with $P=4$ and $C=64\pi^3/3^{23/2}$, so that
$$L(\chi_{-3},2)-\dfrac{p(n)}{q(n)}\sim\dfrac{64\pi^3/3^{23/2}}{n^4}\;.$$
$$A=1-2/n+C_2\log(n)/n^2+\cdots$$

Parametric family:
\begin{verbatim}
[()->lfun(-3,2),(2*n-1)*(9*n^2-9*n+18*k^2+4),-9*n^4*(9*n^2-1)]
\end{verbatim}
Convergence type $P^+$ with $P=4k$.
\end{cf}

These CFs can be Ap\'ery accelerated, formulas too complicated to give here.

\smallskip

\begin{cf}\label{5.4.17}{\ }
\begin{verbatim}
[()->lfun(-3,2),[0,17*n^2-17*n+6],[8/5,-72*n^4]]
\end{verbatim}
$$L(\chi_{-3},2)=\dfrac{8/5}{6-\dfrac{72}{40-\dfrac{1152}{108-\dfrac{5832}{210-\dfrac{18432}{346-\dfrac{45000}{516-\ddots}}}}}}$$
Convergence type $E$ with $E=9/8$, $P=0$, and $C=8\pi^2/135$, so that
$$L(\chi_{-3},2)-\dfrac{p(n)}{q(n)}\sim\dfrac{\pi^2/15}{(9/8)^{n+1}}\;.$$
$$A=1-5/n+15/n^2-789/n^3+4971/n^4-726729/n^5+5648223/n^6-\cdots$$
Parametric families up to the trivial change $n\mapsto n+j$, with $u\ge0$:
\begin{verbatim}
[()->lfun(-3,2),17*n^2+(44*u+52*v-17)*n+2*(2*u+4*v-1)*(5*u+v-3),
               -72*n*(n+u+v)*(n+2*u+v)*(n+2*u+4*v)]
[()->lfun(-3,2),17*n^2+(24*u-17)*n-6*(2*u-1),-72*n^2*(n+u)*(n+2*u)]
[()->lfun(-3,2),17*n^2+(36*u-17)*n-6*(3*u-1),-72*n^2*(n+u)*(n+3*u)]
[()->lfun(-3,2),17*n^2+(42*u-17)*n+3*(3*u-1)*(u-2),-72*n*(n+u)^2*(n+3*u)]
[()->lfun(-3,2),17*n^2+(60*u-17)*n+6*(3*u-1)*(2*u-1),-72*n*(n+2*u)^2*(n+3*u)]
[()->lfun(-3,2),17*n^2+(66*u-17)*n+3*(3*u-1)*(5*u-2),-72*n*(n+2*u)*(n+3*u)^2]
[()->lfun(-3,2),17*n^2+(84*u-17)*n+6*(3*u-1)*(4*u-1),-72*n*(n+3*u)^2*(n+4*u)]
\end{verbatim}
Convergence type $E$ with $E=9/8$ and $P=3u+2v$, $-3u$, $4u$, $-u$, $u$, $-4u$,
and $-2u$ respectively.
\end{cf}

\smallskip

\begin{cf}\label{5.4.18}{\ }
\begin{verbatim}
[()->lfun(-3,2),[0,(2*n-1)*(7*n^2-7*n+3)],[1,-3*n^4*(16*n^2-1)]]
\end{verbatim}
$$L(\chi_{-3},2)=\dfrac{1}{3-\dfrac{45}{51-\dfrac{3024}{225-\dfrac{34749}{609-\dfrac{195840}{1287-\dfrac{748125}{2343-\ddots}}}}}}$$
Convergence type $E$ with $E=4/3$, $P=0$, and $C=\pi^2/(9\sqrt{2})$, so that
$$L(\chi_{-3},2)-\dfrac{p(n)}{q(n)}\sim\dfrac{\pi^2/(9\sqrt{2})}{(4/3)^n}\;.$$
$$A=1-(37/16)/n+(1961/512)/n^2-(518879/8192)/n^3+\cdots$$
\end{cf}

\smallskip

\begin{cf}\label{5.4.19}{\ }
\begin{verbatim}
[()->lfun(-3,2),[0,56*n^3-108*n^2+58*n-3],
                [8/5,-12*n*(n+1)*(2*n-1)^2*(4*n-1)*(4*n-3)]]
\end{verbatim}
$$L(\chi_{-3},2)=\dfrac{8/5}{3-\dfrac{72}{129-\dfrac{22680}{711-\dfrac{356400}{2085-\dfrac{2293200}{4587-\dfrac{9418680}{8553-\ddots}}}}}}$$
Convergence type $E$ with $E=4/3$, $P=1$, and $C=\pi\sqrt{2}/5$, so that
$$L(\chi_{-3},2)-\dfrac{p(n)}{q(n)}\sim\dfrac{\pi\sqrt{2}/5}{(4/3)^nn}\;.$$
$$A=1-(53/16)/n+(11705/512)/n^2-(1953055/8192)/n^3+\cdots$$
Series:
$$L(\chi_{-3},2)=\dfrac{8}{15}\sum_{n\ge0}\dfrac{n!(n+1)!}{(2n+1)(5/4)_n(7/4)_n}(3/4)^n$$
\end{cf}

\smallskip

\begin{cf}\label{5.4.20}{\ }
\begin{verbatim}
[()->lfun(-3,2),[0,10*n^2-10*n+3],[2,-9*n^4]]
\end{verbatim}
$$L(\chi_{-3},2)=\dfrac{2}{3-\dfrac{9}{23-\dfrac{144}{63-\dfrac{729}{123-\dfrac{2304}{203-\dfrac{5625}{303-\ddots}}}}}}$$
Convergence type $E$ with $E=9$, $P=0$, and $C=4\pi^2/27$, so that
$$L(\chi_{-3},2)-\dfrac{p(n)}{q(n)}\sim\dfrac{4\pi^2}{3^{2n+3}}\;.$$
$$A=1-(1/2)/n+(3/8)/n^2-(3/8)/n^3+(51/128)/n^4-(417/512)/n^5+\cdots$$
Parametric families up to the trivial change $n\mapsto n+j$, with $u\ge0$:
\begin{verbatim}
[()->lfun(-3,2),10*n^2-(6*u+10)*n+3*(u+1),-9*n^2*(n+u)^2]
[()->lfun(-3,2),10*n^2+(12*u-10)*n-3*(2*u-1),-9*n^2*(n+2*u)^2]
[()->lfun(-3,2),10*n^2+(26*u-10)*n+16*u^2-13*u+3,-9*n^2*(n+u)^2]
[()->lfun(-3,2),10*n^2+(36*u-10)*n+3*(3*u-1)^2,-9*n^2*(n+u)*(n+3*u)]
\end{verbatim}
Convergence type $E$ with $E=9$ and $P=-4u$, $-2u$, $4u$, and $4u$
respectively.

Many additional sporadic solutions.

\end{cf}

\smallskip

These four continued fractions for $L(\chi_{-3},2)$ apparently cannot be
Ap\'ery accelerated (although their exponential type $E$ is rational).

\smallskip

\begin{cf}\label{5.4.21}{\ }
\begin{verbatim}
[()->lfun(-3,2),[[0,9],[120*n^2-51*n+6,120*n^2+51*n+6]],
                 [[7,-3],[-196*n^4,-(2*n+1)^2*(3*n+1)*(3*n+2)]]]
\end{verbatim}
$$L(\chi_{-3},2)=\dfrac{7}{9-\dfrac{3}{75-\dfrac{196}{177-\dfrac{180}{384-\dfrac{3136}{588-\dfrac{1400}{933-\ddots}}}}}}$$
Convergence type $E$ with $(16+5\sqrt{10})^2/6$, $P=0$, and
$C=8\pi^2(13-4\sqrt{10})/3^{5/2}$, so that
$$L(\chi_{-3},2)-\dfrac{p(n)}{q(n)}\sim\dfrac{8\pi^2(13-4\sqrt{10})/3^{5/2}}{(16+5\sqrt{10})^{2n}6^{-n}}\;.$$
$$A=1-(77d/576)/n+(329d/6400+29645/331776)/n^2+\cdots$$
\end{cf}

\smallskip

\begin{cf}\label{5.4.22}{\ }
\begin{verbatim}
[()->lfun(-3,2)/sqrt(3),[4/9,72,40*n^3+20*n^2+14*n-1],[4/9,-8*n*(2*n+1)^5]]
\end{verbatim}
$$\dfrac{L(\chi_{-3},2)}{\sqrt{3}}=4/9+\dfrac{4/9}{72-\dfrac{1944}{427-\dfrac{50000}{1301-\dfrac{403368}{2935-\dfrac{1889568}{5569-\dfrac{6442040}{9443-\ddots}}}}}}$$
Convergence type $E$ with $E=4$, $P=5/2$, and $C=1/(27\sqrt{\pi})$, so that
$$\dfrac{L(\chi_{-3},2)}{\sqrt{3}}-\dfrac{p(n)}{q(n)}\sim\dfrac{1/(27\sqrt{\pi})}{2^{2n}n^{5/2}}$$
$$A=1-(107/24)/n+(18265/1152)/n^2-\cdots$$
Series:
$$\dfrac{L(\chi_{-3},2)}{\sqrt{3}}=\dfrac{4}{9}+\dfrac{1}{18}\sum_{n\ge0}\dfrac{(3/2)_n}{(2n+3)^2(n+1)!}2^{-2n}$$
\end{cf}

\smallskip

\begin{cf}\label{5.4.23}{\ }
\begin{verbatim}
[()->Catalan,[1/2,9,8*n^2+2*n+1],[1,-2*(n+1)*(2*n+1)^3]]
\end{verbatim}
$$G=1/2+\dfrac{1}{9-\dfrac{108}{37-\dfrac{750}{79-\dfrac{2744}{137-\dfrac{7290}{211-\dfrac{15972}{301-\ddots}}}}}}$$
Convergence type $P^+$ with $P=1/2$ and $C=\sqrt{\pi}/4$, so that
$$G-\dfrac{p(n)}{q(n)}\sim\dfrac{\sqrt{\pi}/4}{n^{1/2}}\;.$$
$$A=1-(13/24)/n+(261/640)/n^2-(6649/21504)/n^3+(21313/98304)/n^4-\cdots$$
Series:
$$G=\dfrac{1}{2}+\dfrac{1}{3}\sum_{n\ge0}\dfrac{(n+1)!}{(2n+3)(5/2)_n}$$
Parametric family for $k\ge0$:
\begin{verbatim}
[()->Catalan,8*n^2+2*n+1+2*k*(2*k+1),-2*(n+1)*(2*n+1)^3]
\end{verbatim}
Convergence type $P^+$ with $P=2k+1/2$.
\end{cf}

Can be Ap\'ery accelerated with convergence type $E=-((1+\sqrt{5})/2)^5$,
formula too complicated to give here.

\smallskip

\begin{cf}\label{5.4.24}{\ }
\begin{verbatim}
[()->Catalan,[0,8*n^2-14*n+7],[1/2,-2*n*(2*n-1)^3]]
\end{verbatim}
$$G=\dfrac{1/2}{1-\dfrac{2}{11-\dfrac{108}{37-\dfrac{750}{79-\dfrac{2744}{137-\dfrac{7290}{211-\ddots}}}}}}$$
Convergence type $P^+$ with $P=1/2$ and $C=\sqrt{\pi}/4$, so that
$$G-\dfrac{p(n)}{q(n)}\sim\dfrac{\sqrt{\pi}/4}{n^{1/2}}\;.$$
$$A=1-(1/24)/n-(19/640)/n^2+(155/21504)/n^3+\cdots$$
Series:
$$G=\dfrac{1}{2}\sum_{n\ge0}\dfrac{n!}{(2n+1)(3/2)_n}$$
Parametric family for $k\ge0$:
\begin{verbatim}
[()->Catalan,8*n^2-14*n+7+2*k*(2*k+1),-2*n*(2*n-1)^3]
\end{verbatim}
Convergence type $P^+$ with $P=2k+1/2$.
\end{cf}

\smallskip

\begin{cf}\label{5.4.25}{\ }
\begin{verbatim}
[()->Catalan,[1,9,8*n],[-1,(2*n+1)^4]]
\end{verbatim}
$$G=1-\dfrac{1}{9+\dfrac{81}{16+\dfrac{625}{24+\dfrac{2401}{32+\dfrac{6561}{40+\dfrac{14641}{48+\ddots}}}}}}$$
Convergence type $P^-$ with $P=2$ and $C=-1/8$, so that
$$G-\dfrac{p(n)}{q(n)}\sim(-1)^{n+1}\dfrac{1/8}{n^2}\;.$$
$$A=1-2/n+(9/4)/n^2-1/n^3-(15/16)/n^4-(3/8)/n^5+(441/64)/n^6+\cdots$$
Series:
$$G=\sum_{n\ge1}\dfrac{(-1)^{n+1}}{(2n-1)^2}$$
Parametric families up to the trivial change $n\mapsto n+j$, with $u\ge0$ and
$v\ge0$:
\begin{verbatim}
[()->Catalan,8*(u+2*v+2*k+1)*(n+u+v),
            (2*n+1)*(2*n+2*u+1)*(2*n+2*u+4*v+1)*(2*n+4*u+4*v+1)]
[()->Catalan,2*((16*u+14)*n+(48*u^2+56*u+13)),
            (2*n+1)*(2*n+4*u+1)*(2*n+8*u+3)*(2*n+12*u+9)]
[()->Catalan,2*((16*u+14)*n+(48*u^2+92*u+43)),
            (2*n+1)*(2*n+4*u+7)*(2*n+8*u+9)*(2*n+12*u+9)]
[()->Catalan,2*((4*u^2+16*u+14)*n+(4*u^3+20*u^2+30*u+13)),
            (2*n+1)*(2*n+2*u+1)*(2*n+2*u+3)*(2*n+4*u+9)]
[()->Catalan,2*((4*u^2+16*u+14)*n+(4*u^3+28*u^2+62*u+43)),
            (2*n+1)*(2*n+2*u+7)*(2*n+2*u+9)*(2*n+4*u+9)]
\end{verbatim}
Convergence type $P^-$ with $P=2u+4v+4k+2$, $8u+7$, $8u+7$,
$2u^2+8u+7$, and $2u^2+8u+7$ respectively.
\end{cf}

\smallskip

\begin{cf}\label{5.4.26}{\ }
\begin{verbatim}
[()->Catalan,[17/18,25,12*n^2+12*n-1],[-8/9,n*(n+2)*(2*n+1)^2*(2*n+3)^2]]
\end{verbatim}
$$G=17/18-\dfrac{8/9}{25+\dfrac{675}{71+\dfrac{9800}{143+\dfrac{59535}{239+\dfrac{235224}{359+\dfrac{715715}{503+\ddots}}}}}}$$
Convergence type $P^-$ with $P=3$ and $C=-1/8$, so that
$$G-\dfrac{p(n)}{q(n)}\sim(-1)^{n+1}\dfrac{1/8}{n^3}$$
$$A=1-(9/2)/n+(25/2)/n^2-(105/4)/n^3+\cdots$$
Series:
$$G=\dfrac{1}{2}+4\sum_{n\ge0}(-1)^n\dfrac{n+1}{(2n+1)^2(2n+3)^2}$$
\end{cf}

\smallskip
  
\begin{cf}\label{5.4.27}{\ }
\begin{verbatim}
[()->Catalan,[0,8*n^2-8*n+3],[1/2,-16*n^4]]
\end{verbatim}
$$G=\dfrac{1/2}{3-\dfrac{16}{19-\dfrac{256}{51-\dfrac{1296}{99-\dfrac{4096}{163-\dfrac{10000}{243-\ddots}}}}}}$$
Convergence type $L$ with $C=\pi^3/8$, so that
$$G-\dfrac{p(n)}{q(n)}\sim\dfrac{\pi^3/8}{\log(n)}\;.$$
$A=1+\cdots$
\end{cf}

This is only given as an example of a nontrivial CF of type
$L$, and is the case $u=k=0$ of the parametric family given in the next CF.

\smallskip

\begin{cf}\label{5.4.28}{\ }
\begin{verbatim}
[()->Catalan,[1,8*n^2-8*n+7],[-1/2,-16*n^4]]
\end{verbatim}
$$G=1-\dfrac{1/2}{7-\dfrac{16}{23-\dfrac{256}{55-\dfrac{1296}{103-\dfrac{4096}{167-\dfrac{10000}{247-\ddots}}}}}}$$
Convergence type $P^+$ with $P=2$ and $C=-\pi^3/2^{10}$, so that
$$G-\dfrac{p(n)}{q(n)}\sim-\dfrac{\pi^3/1024}{n^2}\;.$$
$A=1-1/n+C_2\log(n)/n^2+\cdots$.

Parametric family for $u\ge0$ and $k\ge0$:
\begin{verbatim}
[()->Catalan,8*n^2+(4*u-8)*n+(4*k-2)*u+4*k^2+3,-16*n^3*(n+u)]
\end{verbatim}
Convergence type $P^+$ with $P=2k+u$.
\end{cf}

\smallskip

\begin{cf}\label{5.4.29}{\ }
\begin{verbatim}
[()->Catalan,[0,3*n^2-3*n+1],[1/2,-2*n^4]]
\end{verbatim}
$$G=\dfrac{1/2}{1-\dfrac{2}{7-\dfrac{32}{19-\dfrac{162}{37-\dfrac{512}{61-\dfrac{1250}{91-\ddots}}}}}}$$
Convergence type $E$ with $E=2$, $P=0$, and $C=\pi^2/8$, so that
$$G-\dfrac{p(n)}{q(n)}\sim\dfrac{\pi^2}{2^{n+3}}\;.$$
$$A=1-1/n+1/n^2-5/n^3+11/n^4-130/n^5+(851/2)/n^6-8339/n^7+\cdots$$
Parametric families with all parameters nonnegative integers:
\begin{verbatim}
[()->Catalan,3*n^2+(4*v+4*w-3)*n+(2*v-1)*(2*w-1),
            -2*n*(n+u)*(n+2*v)*(n+2*w)]
[()->Catalan,3*n^2+(2*u+2*v+2*w-3)*n+u*v+u*w+v*w-u-v-w+1,
            -2*n*(n+u)*(n+v)*(n+w)]
\end{verbatim}
where in the second family we must have $u\equiv v\equiv w\pmod{2}$.

Convergence types $E$ with $E=2$ and $P=2v+2w-3u$ or $P=u+v+w$ respectively.
\end{cf}

\smallskip

\begin{cf}\label{5.4.30}{\ }
\begin{verbatim}
[()->Catalan,[[1,6],[6*n+1,6*n+5]],
             [[-1,48],[(2*n+1)^2*(4*n+1),16*(n+1)^2*(4*n+3)]]]
\end{verbatim}
$$G=1-\dfrac{1}{6+\dfrac{48}{7+\dfrac{45}{11+\dfrac{448}{13+\dfrac{225}{17+\dfrac{1584}{19+\ddots}}}}}}$$
Convergence type $E$ with $E=-2$, $P=2$, and $C=-\pi/2^{5/2}$, so that
$$G-\dfrac{p(n)}{q(n)}\sim(-1)^{n+1}\dfrac{\pi}{2^{n+5/2}n^2}\;.$$
$$A=1-(49/8)/n+(4337/128)/n^2-(214043/1024)/n^3+\cdots$$
Series:
$$G=1-\dfrac{1}{30}\sum_{n\ge0}\dfrac{(6n+7)n!^2}{(2n+3)(7/4)_n(9/4)_n}2^{-2n}$$
\end{cf}

\smallskip

\begin{cf}\label{5.4.31}{\ }
\begin{verbatim}
[()->Catalan,[[1/2,2],[6*n,3*n+2]],[[1,2],[-(n+1)^3,(2*n+1)^3]]]
\end{verbatim}
$$G=1/2+\dfrac{1}{2+\dfrac{2}{6-\dfrac{8}{5+\dfrac{27}{12-\dfrac{27}{8+\dfrac{125}{18-\ddots}}}}}}$$
Convergence type $E$ with $E=-2\sqrt{2}$, $P=3/2$, and $C=(2\pi)^{3/2}/27$, so
that
$$G-\dfrac{p(n)}{q(n)}\sim(-1)^n\dfrac{\pi^{3/2}/27}{2^{(3n-3)/2}n^{3/2}}\;.$$
$$A=1-(23/12)/n+(2995/864)/n^2-(197957/31104)/n^3+\cdots$$
Series:
$$G=\dfrac{1}{2}+3\sum_{n\ge0}(-1)^n\dfrac{(n+1)(n+1)!^3}{(3n^2+3n+1)(3n^2+9n+7)(3/2)_n^3}2^{-3n}$$
\end{cf}

\smallskip

\begin{cf}\label{5.4.32}{\ }
\begin{verbatim}
[()->Catalan,[[0,2],[10*n^2-6*n+1,10*n^2+6*n+1]],
             [[1,-4],[-6*n^4,-(3*n+1)*(3*n+2)*(2*n+1)^2]]]
\end{verbatim}
$$G=\dfrac{1}{2-\dfrac{4}{5-\dfrac{6}{17-\dfrac{180}{29-\dfrac{96}{53-\dfrac{1400}{73-\ddots}}}}}}$$
Convergence type $E$ with $E=\sqrt{27/2}$, $P=0$, and $C=\pi^2/(4\sqrt{3})$,
so that
$$G-\dfrac{p(n)}{q(n)}\sim\dfrac{\pi^2/(4\sqrt{3})}{(27/2)^{n/2}}$$
$$A=1-(277/450)/n+(264649/405000)/n^2+\cdots$$
\end{cf}

\smallskip

\begin{cf}\label{5.4.33}{\ }
\begin{verbatim}
[()->Catalan,[[1,14],[20*n^2+8*n+1,20*n^2+32*n+13]],
              [[-1,-64],[(2*n+1)^4,-64*(n+1)^4]]]
\end{verbatim}
$$G=1-\dfrac{1}{14-\dfrac{64}{29+\dfrac{81}{65-\dfrac{1024}{97+\dfrac{625}{157-\dfrac{5184}{205+\ddots}}}}}}$$
Convergence type $E$ with $E=-i((1+\sqrt{5})/2)^5$, $P=0$, and $C=-\pi^2/2/((1+\sqrt{5})/2)^8$, so that
$$G-\dfrac{p(n)}{q(n)}\sim(-1)^{\lfloor (n-1)/2\rfloor}\dfrac{\pi^2/2}{((1+\sqrt{5})/2)^{5n+8}}\;.$$
$$A=1-(d/5)/n+(9d/25+1/10)/n^2+(-384d/625-9/25)/n^3+\cdots$$
\end{cf}

Note: the CF for Catalan's constant given by W.~Zudilin in \cite{Zud} is the
even contracted CF of this one.

\smallskip

\begin{cf}\label{5.4.34}{\ }
\begin{verbatim}
[()->Catalan,[[1,18],[16,53],[20*n^2+5,20*n^2+24*n+9]],
            [[-2,162],[(2*n+3)^4,-64*n^4]]]
\end{verbatim}
$$G=1-\dfrac{2}{18+\dfrac{162}{16+\dfrac{625}{53-\dfrac{64}{85+\dfrac{2401}{137-\dfrac{1024}{185+\ddots}}}}}}$$
Convergence type $E$ with $E=-i((1+\sqrt{5})/2)^5$, $P=0$, and
$C=\pi^2/2/((1+\sqrt{5})/2)^{10}$, so that
$$G-\dfrac{p(n)}{q(n)}\sim(-1)^{\lfloor(n-1)/2\rfloor}\dfrac{\pi^2/2}{((1+\sqrt{5})/2)^{5n+10}}\;.$$
$$A=1+(23d/5)/n+(-21d/5+529/10)/n^2+(53796d/625-483/5)/n^3+\cdots$$
\end{cf}

\smallskip

\begin{cf}\label{5.4.35}{\ }
\begin{verbatim}
[()->log(2)^2,[0,7*n^2+10*n-9],[9/2,72*n^3*(2*n-1)]]
\end{verbatim}
$$\log^2(2)=\dfrac{9/2}{8+\dfrac{72}{39+\dfrac{1728}{84+\dfrac{9720}{143+\dfrac{32256}{216+\dfrac{81000}{303+\ddots}}}}}}$$
Convergence type $E$ with $E=-16/9$, $P=3/2$, and $C=9\sqrt{\pi}/50$, so that
$$\log^2(2)-\dfrac{p(n)}{q(n)}\sim(-1)^n\dfrac{9\sqrt{\pi}/50}{(4/3)^{2n}n^{3/2}}\;.$$
$$A=1-(167/200)/n+(2301/16000)/n^2+(1955127/3200000)/n^3+\cdots$$
Series:
$$\log^2(2)=\dfrac{9}{16}\sum_{n\ge0}(-1)^n\dfrac{n!}{(n+1)(3/2)_n}(3/4)^{2n}$$
\end{cf}

\smallskip

\begin{cf}\label{5.4.36}{\ }
\begin{verbatim}
[()->log(2)^2,[0,7*n^2-2*n-1],[2,4*n^3*(2*n-1)]]
\end{verbatim}
$$\log^2(2)=\dfrac{2}{4+\dfrac{4}{23+\dfrac{96}{56+\dfrac{540}{103+\dfrac{1792}{164+\dfrac{4500}{239+\ddots}}}}}}$$
Convergence type $E$ with $E=-8$, $P=3/2$, and $C=2\sqrt{\pi}/9$, so that
$$\log^2(2)-\dfrac{p(n)}{q(n)}\sim(-1)^n\dfrac{\sqrt{\pi}/9}{2^{3n-1}n^{3/2}}$$
$$A=1-(29/24)/n+(3547/3456)/n^2-(110303/248832)/n^3+\cdots$$
Series:
$$\log^2(2)=\dfrac{1}{2}\sum_{n\ge0}(-1)^n\dfrac{n!}{(n+1)(3/2)_n}2^{-3n}$$
\end{cf}
        
\smallskip

\begin{cf}\label{5.4.37}{\ }
\begin{verbatim}
[()->log(3/2)^2,[0,119*n^2-22*n-25],[25/2,1800*n^3*(2*n-1)]]
\end{verbatim}
$$\log^2(3/2)=\dfrac{25/2}{72+\dfrac{1800}{407+\dfrac{43200}{980+\dfrac{243000}{1791+\dfrac{806400}{2840+\dfrac{2025000}{4127+\ddots}}}}}}$$
Convergence type $E$ with $E=-144/25$, $P=3/2$, and $C=25\sqrt{\pi}/338$, so
that
$$\log^2(3/2)-\dfrac{p(n)}{q(n)}\sim(-1)^n\dfrac{25\sqrt{\pi}/338}{(12/5)^{2n}n^{3/2}}\;.$$
$$A=1-(1559/1352)/n+(3167761/3655808)/n^2-(762881117/4942652416)/n^3+\cdots$$
Series:
$$\log^2(3/2)=\dfrac{25}{144}\sum_{n\ge0}(-1)^n\dfrac{n!}{(n+1)(3/2)_n}(5/12)^{2n}$$
\end{cf}

\smallskip

\begin{cf}\label{5.4.38}{\ }
\begin{verbatim}
[()->log(3/2)^2,[0,23*n^2-10*n-1],[2,12*n^3*(2*n-1)]]
\end{verbatim}
$$\log^2(3/2)=\dfrac{2}{12+\dfrac{12}{71+\dfrac{288}{176+\dfrac{1620}{327+\dfrac{5376}{524+\dfrac{13500}{767+\ddots}}}}}}$$
Convergence type $E$ with $E=-24$, $P=3/2$, and $C=2\sqrt{\pi}/25$, so that
$$\log^2(3/2)-\dfrac{p(n)}{q(n)}\sim(-1)^n\dfrac{2\sqrt{\pi}}{24^nn^{3/2}}$$
$$A=1-(263/200)/n+(21821/16000)/n^2+\cdots$$
Series:
$$\log^2(3/2)=\dfrac{1}{6}\sum_{n\ge0}(-1)^n\dfrac{n!}{(n+1)(3/2)_n}24^{-n}$$
\end{cf}
      
\smallskip

\begin{cf}\label{5.4.39}{\ }
\begin{verbatim}
[()->log((1+sqrt(5))/2)^2,[0,3*n^2-1],[1/2,2*n^3*(2*n-1)]]
\end{verbatim}
$$\log((1+\sqrt{5})/2)^2=\dfrac{1/2}{2+\dfrac{2}{11+\dfrac{48}{26+\dfrac{270}{47+\dfrac{896}{74+\dfrac{2250}{107+\ddots}}}}}}$$
Convergence type $E$ with $E=-4$, $P=3/2$, and $C=\sqrt{\pi}/10$, so that
$$\log((1+\sqrt{5})/2)^2-\dfrac{p(n)}{q(n)}\sim(-1)^n\dfrac{\sqrt{\pi}/5}{2^{2n+1}n^{3/2}}\;.$$
$$A=1-(43/40)/n+(421/640)/n^2+(4243/25600)/n^3-(4368513/4096000)/n^4+\cdots$$
Series:
$$\log((1+\sqrt{5})/2)^2=\dfrac{1}{4}\sum_{n\ge0}(-1)^n\dfrac{n!}{(n+1)(3/2)_n}2^{-2n}$$
\end{cf}

\medskip

Note that we have given continued fractions for periods of degree $2$ only
for $\pi^2$, $L(\chi_{-3},2)$, and $G=L(\chi_{-4},2)$, as well as
samples of $\log^2(z)$ (we will also give linear combinations of periods of
degree $2$ below). Are there any others which are not too complicated ?

\medskip

\section{Constants: $\pi^3$, $\z(3)$, and Periods of Degree $3$}

\smallskip

\begin{cf}\label{5.5.1}{\ }
\begin{verbatim}
[()->Pi^3,[0,1,24*n^2-48*n+26],[32,(2*n-1)^6]]
\end{verbatim}
$$\pi^3=\dfrac{32}{1+\dfrac{1}{26+\dfrac{729}{98+\dfrac{15625}{218+\dfrac{117649}{386+\dfrac{531441}{602+\ddots}}}}}}$$
Convergence type $P^-$ with $P=3$ and $C=2$, so that
$$\pi^3-\dfrac{p(n)}{q(n)}\sim(-1)^n\dfrac{2}{n^3}\;.$$
$$A=1-(3/2)/n^2+(75/16)/n^4-(427/16)/n^6+(62325/256)/n^8+\cdots$$
Series:
$$\pi^3=32\sum_{n\ge1}\dfrac{(-1)^{n+1}}{(2n-1)^3}$$
Parametric family (of an equivalent CF) for $k\ge0$:
\begin{verbatim}
[()->Pi^3,(2*k+1)*(n-1)*(24*n^2-48*n+26),(n-k-1)*(n+k)*(2*n-1)^6]
\end{verbatim}
Convergence type $P^-$ with $P=6k+3$.
\end{cf}

This is simply the Euler transformation into a CF of the series giving
$L(\chi_{-4},3)$.

\smallskip

\begin{cf}\label{5.5.2}{\ }
\begin{verbatim}
[()->Pi^3,[216/7,216,80*n^4+64*n^3+72*n^2+1],[216/7,-8*n*(2*n+1)^7]]
\end{verbatim}
$$\pi^3=216/7+\dfrac{216/7}{216-\dfrac{17496}{2081-\dfrac{1250000}{8857-\dfrac{19765032}{25729-\dfrac{153055008}{59801-\ddots}}}}}$$
Convergence type $E$ with $E=4$, $P=7/2$, $C=9/(7\sqrt{\pi})$, so that
$$\pi^3-\dfrac{p(n)}{q(n)}\sim\dfrac{9/(7\sqrt{\pi})}{2^{2n}n^{7/2}}$$
$$A=1-(151/24)/n+(3705/128)/n^2-(132247/1024)/n^3+\cdots$$
Series:
$$\pi^3=\dfrac{216}{7}+\dfrac{27}{7}\sum_{n\ge0}\dfrac{(3/2)_n}{(2n+3)^3(n+1)!}2^{-2n}$$
\end{cf}

\smallskip

\begin{cf}\label{5.5.3}{\ }
\begin{verbatim}
[()->Pi^3,
[[32,27],[168*n^3+76*n^2+14*n+1,168*n^3+228*n^2+126*n+27]],
[[-32,1458],[4096*n^6*(2*n+1),2*(n+1)*(2*n+3)^6]]]
\end{verbatim}
$$\pi^3=32-\dfrac{32}{27+\dfrac{1458}{259+\dfrac{12288}{549+\dfrac{62500}{1677+\dfrac{1310720}{2535+\dfrac{705894}{5263+\ddots}}}}}}$$
Convergence type $E$ with $E=-8$, $P=1/2$, and $C=-\pi^{5/2}/3/2^{3/2}$,
so that
$$\pi^3-\dfrac{p(n)}{q(n)}\sim(-1)^{n+1}\dfrac{\pi^{5/2}/3}{2^{3n+3/2}n^{1/2}}\;.$$
$$A=1-(67/36)/n+(3827/864)/n^2-(330019/31104)/n^3+\cdots$$
Series:
\begin{align*}
  \dfrac{1}{\pi^3}&=\dfrac{1}{32}\sum_{n\ge0}\dfrac{(6n+1)(28n^2+8n+1)(1/2)_n^7}{n!^7}2^{-6n}\\
  &=\dfrac{1}{32}\sum_{n\ge0}(6n+1)(28n^2+8n+1)\binom{2n}{n}^72^{-20n}\end{align*}
\end{cf}

\smallskip

\begin{cf}\label{5.5.4}{\ }
\begin{verbatim}
[()->Pi^3/sqrt(3),[0,(2*n-1)*(3*n^2-3*n+1)],[81/4,3*n^6]]
\end{verbatim}
$$\dfrac{\pi^3}{\sqrt{3}}=\dfrac{81/4}{1+\dfrac{3}{21+\dfrac{192}{95+\dfrac{2187}{259+\dfrac{12288}{549+\dfrac{46875}{1001+\ddots}}}}}}$$
Convergence type $E$ with $E=-(2+\sqrt{3})^2$, $P=0$, $C=3\pi^3\sqrt{3}/(2+\sqrt{3})$, so that
$$\dfrac{\pi^3}{\sqrt{3}}-\dfrac{p(n)}{q(n)}\sim(-1)^n\dfrac{3\pi^3\sqrt{3}}{(2+\sqrt{3})^{2n+1}}\;.$$
$$A=1-(7d/24)/n+(7d/48+49/384)/n^2+(-469d/27648-49/384)/n^3+\cdots$$
\end{cf}

\smallskip

\begin{cf}\label{5.5.5}{\ }
\begin{verbatim}
[()->zeta(3),[0,1,4*(n-1)*(4*n^2-8*n+7)],[8/7,-(2*n-1)^6]]
\end{verbatim}
$$\z(3)=\dfrac{8/7}{1-\dfrac{1}{28-\dfrac{729}{152-\dfrac{15625}{468-\dfrac{117649}{1072-\dfrac{531441}{2060-\ddots}}}}}}$$
Convergence type $P^+$ with $P=2$ and $C=1/14$, so that
$$\z(3)-\dfrac{p(n)}{q(n)}\sim\dfrac{1/14}{n^2}\;.$$
$$A=1-(1/4)/n^2+(7/48)/n^4-(31/192)/n^6+(381/1280)/n^8+\cdots$$
Series:
$$\z(3)=\dfrac{8}{7}\sum_{n\ge1}\dfrac{1}{(2n-1)^3}$$
Parametric family:
\begin{verbatim}
[()->zeta(3),4*(n+u+v+w-1)*(4*n^2+8*(u+v+w-1)*n+4*u^2+8*v^2+16*w^2
     +12*u*v+16*u*w+24*v*w-4*u+4*w+8*k*(u+2*v+3*w)+8*k^2+8*k+7),
     -(2*n-1)*(2*n-1+2*u)*(2*n-1+2*u+2*v)*
              (2*n-1+2*u+2*v+4*w)*(2*n-1+2*u+4*v+4*w)*(2*n-1+4*u+4*v+4*w)]
\end{verbatim}
Convergence type $P^+$ with $P=2u+4v+6w+4k+2$.

Hundreds of additional sporadic CFs.
\end{cf}

\smallskip

\begin{cf}\label{5.5.6}{\ }
\begin{verbatim}
[()->zeta(3),[0,(2*n-1)*(n^2-n+1)],[1,-n^6]]
\end{verbatim}
$$\z(3)=\dfrac{1}{1-\dfrac{1}{9-\dfrac{64}{35-\dfrac{729}{91-\dfrac{4096}{189-\dfrac{15625}{341-\ddots}}}}}}$$
Convergence type $P^+$ with $P=2$ and $C=1/2$, so that
$$\z(3)-\dfrac{p(n)}{q(n)}\sim\dfrac{1/2}{n^2}\;.$$
$$A=1-1/n+(1/2)/n^2-(1/6)/n^4+(1/6)/n^6-(3/10)/n^8+(5/6)/n^{10}-\cdots$$
Series:
$$\z(3)=\sum_{n\ge1}\dfrac{1}{n^3}$$
Parametric family:
\begin{verbatim}
[()->zeta(3),(2*n+2*(u+v+w)-1)*(n^2+(2*(u+v+w)-1)*n
 +u^2+2*v^2+4*w^2+3*u*v+4*u*w+6*v*w+v+2*w+1+2*k*(u+2*v+3*w)+2*k^2+2*k),
 -n*(n+u)*(n+u+v)*(n+u+v+2*w)*(n+u+2*v+2*w)*(n+2*u+2*v+2*w)]
\end{verbatim}
Convergence type $P^+$ with $P=2u+4v+6w+4k+2$.

Hundreds of additional sporadic CFs.
\end{cf}

\smallskip

\begin{cf}\label{5.5.7}{\ }
\begin{verbatim}
[()->zeta(3),[0,(2*n-1)^2*(n^2-n+1)],[1,-n^6*(4*n^2-1)]]
\end{verbatim}
$$\z(3)=\dfrac{1}{1-\dfrac{3}{27-\dfrac{960}{175-\dfrac{25515}{637-\dfrac{258048}{1701-\dfrac{1546875}{3751-\ddots}}}}}}$$
Convergence type $P^+$ with $P=2$ and $C=1/2$, so that
$$\z(3)-\dfrac{p(n)}{q(n)}\sim\dfrac{1/2}{n^2}\;.$$
$$A=1-1/n+(1/2)/n^2-(1/6)/n^4+(1/6)/n^6-(3/10)/n^8+(5/6)/n^{10}-\cdots$$
Series:
$$\z(3)=\sum_{n\ge1}\dfrac{1}{n^3}$$
Parametric family for $k\ge0$:
\begin{verbatim}
[()->zeta(3),4*n^4-8*n^3-(2*k^2+2*k-9)*n^2+(2*k^2+2*k-5)*n
             -(k^4+2*k^3+3*k^2+2*k-2)/2,-n^6*(4*n^2-(2*k+1)^2)]
\end{verbatim}
Convergence type $P^+$ with $P=2$.
\end{cf}

\smallskip

This CF is termwise identical with the previous one, but the parametric
family is totally different.

\smallskip

\begin{cf}\label{5.5.8}{\ }
\begin{verbatim}
[()->zeta(3),[8/7,27,4*n*(4*n^2+3)],[8/7,-(2*n+1)^6]]
\end{verbatim}
$$\z(3)=8/7+\dfrac{8/7}{27-\dfrac{729}{152-\dfrac{15625}{468-\dfrac{117649}{1072-\dfrac{531441}{2060-\dfrac{1771561}{3528-\ddots}}}}}}$$
Convergence type $P^+$ with $P=2$ and $C=1/14$, so that
$$\z(3)-\dfrac{p(n)}{q(n)}\sim\dfrac{1/14}{n^2}\;.$$
$$A=1-2/n+(11/4)/n^2-3/n^3+(127/48)/n^4+\cdots$$
Series:
$$\z(3)=\dfrac{8}{7}\sum_{n\ge1}\dfrac{1}{(2n-1)^3}$$
Parametric family for $k\ge0$:
\begin{verbatim}
[()->zeta(3),4*n*(4*n^2+3+8*k*(k+1)),-(2*n+1)^6]
\end{verbatim}
Convergence type $P^+$ with $P=4k+2$.
\end{cf}

\smallskip

\begin{cf}\label{5.5.9}{\ }
\begin{verbatim}
[()->zeta(3),[1,(2*n-3)*(2*n+1)*(n^2-n+5)],[-3,-n^6*(4*n^2-9)]]
[()->zeta(3),[1,(2*n+1)*(n^2-n+5)],[3,-n^6*(2*n+3)/(2*n-1)]]
\end{verbatim}
\begin{align*}\z(3)&=1-\dfrac{3}{-15+\dfrac{5}{35-\dfrac{448}{231-\dfrac{19683}{765-\dfrac{225280}{1925-\dfrac{1421875}{4095-\ddots}}}}}}\\&=1+\dfrac{3}{15-\dfrac{5}{35-\dfrac{448/3}{77-\dfrac{6561/5}{153-\dfrac{45056/7}{275-\dfrac{203125/9}{455-\ddots}}}}}}\end{align*}
Convergence type $P^+$ with $P=6$ and $C=1/24$, so that
$$\z(3)-\dfrac{p(n)}{q(n)}\sim\dfrac{1/24}{n^6}\;.$$
$$A=1-3/n+(7/2)/n^2-(87/20)/n^4+(3/4)/n^5+(77/8)/n^6+\cdots$$
Series:
$$\z(3)=1+\sum_{n\ge1}\dfrac{1}{n^3(2n^2-2n+1)(2n^2+2n+1)}$$
\end{cf}

\smallskip

\begin{cf}\label{5.5.10}{\ }
\begin{verbatim}
[()->zeta(3),[0,(2*n-1)*(2*n^2-2*n+1)],[2/7,-4*n^6]]
\end{verbatim}
$$\z(3)=\dfrac{2/7}{1-\dfrac{4}{15-\dfrac{256}{65-\dfrac{2916}{175-\dfrac{16384}{369-\dfrac{62500}{671-\ddots}}}}}}$$
Convergence type $L$ with $C=\pi^4/28$, so that
$$\z(3)-\dfrac{p(n)}{q(n)}\sim\dfrac{\pi^4/28}{\log(n)}\;.$$
$A=1+\cdots$
\end{cf}

This is only given as another example of a nontrivial CF of type
$L$, and is the case $u=0$, $k=0$ of the parameterization given in the
next CF.

\smallskip

\begin{cf}\label{5.5.11}{\ }
\begin{verbatim}
[()->zeta(3),[8/7,(2*n-1)*(2*n^2-2*n+5)],[2/7,-4*n^6]]
\end{verbatim}
$$\z(3)=8/7+\dfrac{2/7}{5-\dfrac{4}{27-\dfrac{256}{85-\dfrac{2916}{203-\dfrac{16384}{405-\dfrac{62500}{715-\ddots}}}}}}$$
Convergence type $P^+$ with $P=4$ and $C=\pi^4/(7\cdot2^{10})$, so that
$$\z(3)-\dfrac{p(n)}{q(n)}\sim\dfrac{\pi^4/7168}{n^4}\;.$$
$$A=1-2/n+(7/4)/n^2-(1/4)/n^3+\cdots$$

Parametric family with $u\ge0$ and $k\ge0$:
\begin{verbatim}
[()->zeta(3),(2*n+2*u-1)*(2*n^2+(4*u-2)*n+(2*k-2)*u+k^2+1),
            -4*n^3*(n+2*u)^3]
\end{verbatim}
Convergence type $P^+$ with $P=2u+2k$.
\end{cf}
  
\smallskip

\begin{cf}\label{5.5.12}{\ }
\begin{verbatim}
[()->zeta(3),[0,3*n^2-3*n+1],[4/3,n^6]]
\end{verbatim}
$$\z(3)=\dfrac{4/3}{1+\dfrac{1}{7+\dfrac{64}{19+\dfrac{729}{37+\dfrac{4096}{61+\dfrac{15625}{91+\ddots}}}}}}$$
Convergence type $P^-$ with $P=3$ and $C=2/3$, so that
$$\z(3)-\dfrac{p(n)}{q(n)}\sim(-1)^n\dfrac{2/3}{n^3}\;.$$
$$A=1-(3/2)/n+(5/2)/n^3-(21/2)/n^5+(153/2)/n^7-(1705/2)/n^9+(26949/2)/n^{11}+\cdots$$
Series:
$$\z(3)=\dfrac{4}{3}\sum_{n\ge1}\dfrac{(-1)^{n+1}}{n^3}$$
\end{cf}

\smallskip

\begin{cf}\label{5.5.13}{\ }
\begin{verbatim}
[()->zeta(3),[0,(2*n-1)*(3*n^2-3*n+1)],[4/3,n^6*(4*n^2-1)]]
\end{verbatim}
$$\z(3)=\dfrac{4/3}{1+\dfrac{3}{21+\dfrac{960}{95+\dfrac{25515}{259+\dfrac{258048}{549+\dfrac{1546875}{1001+\ddots}}}}}}$$
Convergence type $P^-$ with $P=3$ and $C=2/3$, so that
$$\z(3)-\dfrac{p(n)}{q(n)}\sim(-1)^n\dfrac{2/3}{n^3}\;.$$
$$A=1-(3/2)/n+(5/2)/n^3-(21/2)/n^5+(153/2)/n^7-(1705/2)/n^9+(26949/2)/n^{11}+\cdots$$
Series:
$$\z(3)=\dfrac{4}{3}\sum_{n\ge1}\dfrac{(-1)^{n+1}}{n^3}$$
Parametric family for $k\ge0$:
\begin{verbatim}
[()->zeta(3),(2*k+1)*(2*n-1)*(3*n^2-3*n+1),n^6*(4*n^2-(2*k+1)^2)]
\end{verbatim}
Convergence type $P^-$ with $P=6k+3$.
\end{cf}

\smallskip

Note that this CF is termwise equivalent to the previous one, but
contrary to the previous one has a parametric family.

\smallskip
  
\begin{cf}\label{5.5.14}{\ }
\begin{verbatim}
[()->zeta(3),[0,43*n^4-86*n^3+77*n^2-34*n+6],
             [11/3,-12*n^6*(36*n^2-1)]]
\end{verbatim}
$$\z(3)=\dfrac{11/3}{6-\dfrac{420}{246-\dfrac{109824}{1758-\dfrac{2825604}{6606-\dfrac{28262400}{17886-\ddots}}}}}$$
Convergence type $E$ with $E=27/16$, $P=0$, and $C=4\pi^3/81$, so that
$$\z(3)-\dfrac{p(n)}{q(n)}\sim\dfrac{4\pi^3/81}{(27/16)^n}\;.$$
$$A=1-(55/36)/n+(5005/2592)/n^2+\cdots$$
\end{cf}

\smallskip
  
\begin{cf}\label{5.5.15}{\ }
\begin{verbatim}
[()->zeta(3),[[0,-7],[20*n^2-8*n+1,20*n^2+16*n+5]],
             [[2,112],[(2*n-1)^4*(4*n-1),-64*(n+1)^4*(4*n+1)]]]
\end{verbatim}
$$\z(3)=\dfrac{2}{-7+\dfrac{112}{13+\dfrac{3}{41-\dfrac{5120}{65+\dfrac{567}{117-\dfrac{46656}{157+\ddots}}}}}}$$
Convergence type $E$ with $E=2i$, $P=2$, and $C=2^{3/2}\pi/7$, so that
$$\z(3)-\dfrac{p(n)}{q(n)}\sim(-1)^{\lfloor n/2\rfloor}\dfrac{\pi/7}{2^{n-3/2}n^2}\;.$$
$$A=1-(13/8)/n+(41/128)/n^2+(6345/1024)/n^3+\cdots$$
Series:
$$\z(3)=\dfrac{2}{21}\sum_{n\ge0}(-1)^n\dfrac{(20n^2+32n+13)n!^2}{(n+1)(2n+1)^2(5/4)_n(7/4)_n}2^{-2n}$$
\end{cf}

\smallskip
  
\begin{cf}\label{5.5.16}{\ }
\begin{verbatim}
[()->zeta(3),[[1,4],[5*n^2+4*n+1,5*n^2+10*n+5]],
             [[1,8],[-(2*n+1)*(n+1)^4,8*(n+1)^5]]]
\end{verbatim}
$$\z(3)=1+\dfrac{1}{4+\dfrac{8}{10-\dfrac{48}{20+\dfrac{256}{29-\dfrac{405}{45+\dfrac{1944}{58-\ddots}}}}}}$$
Convergence type $E$ with $E=-2$, $P=7/2$, and $C=\sqrt{8\pi}$, so that
$$\z(3)-\dfrac{p(n)}{q(n)}\sim(-1)^n\dfrac{\sqrt{\pi}}{2^{n-3/2}n^{7/2}}\;.$$
$$A=1-(35/4)/n+(1497/32)/n^2-(23177/128)/n^3+\cdots$$
Series:
$$\z(3)=1+\dfrac{1}{12}\sum_{n\ge0}(-1)^n\dfrac{(5n^2+14n+10)n!}{(n+1)^2(n+2)^2(5/2)_n}2^{-2n}$$
\end{cf}

There exist many period 2 CFs for $\z(3)$ with convergence type $E$ with
$|E|=2$, analogous to the above CFs.

\smallskip
  
\begin{cf}\label{5.5.17}{\ }
\begin{verbatim}
[()->zeta(3),[0,(2*n-1)*(5*n^2-5*n+2)],[12/7,-16*n^6]]
\end{verbatim}
$$\z(3)=\dfrac{12/7}{2-\dfrac{16}{36-\dfrac{1024}{160-\dfrac{11664}{434-\dfrac{65536}{918-\dfrac{250000}{1672-\ddots}}}}}}$$
Convergence type $E$ with $E=4$, $P=0$, and $C=\pi^3/14$, so that
$$\z(3)-\dfrac{p(n)}{q(n)}\sim\dfrac{\pi^3/7}{2^{2n+1}}\;.$$
$$A=1-(3/4)/n+(21/32)/n^2-(145/128)/n^3+(3603/2048)/n^4-(56557/8192)/n^5+\cdots$$
\end{cf}

\smallskip

Contrary to the next CFs, I have not been able to find any other CF for
$\z(3)$ with $a(n)=10n^3+\cdots$ and $b(n)=-16n^6+\cdots$.

\smallskip

\begin{cf}\label{5.5.18}{\ }
\begin{verbatim}
[()->zeta(3),[0,3*n^3+n^2-3*n+1],[5/2,2*n^5*(2*n-1)]]
\end{verbatim}
$$\zeta(3)=\dfrac{5/2}{2+\dfrac{2}{23+\dfrac{192}{82+\dfrac{2430}{197+\dfrac{14336}{386+\dfrac{56250}{667+\ddots}}}}}}$$
Convergence type $E$ with $E=-4$, $P=5/2$, and $C=\sqrt{\pi}/2$, so that
$$\z(3)-\dfrac{p(n)}{q(n)}\sim(-1)^n\dfrac{\sqrt{\pi}}{2^{2n+1}n^{5/2}}\;.$$
$$A=1-(15/8)/n+(225/128)/n^2+(235/1024)/n^3-(130261/32768)/n^4+\cdots$$
Series:
$$\z(3)=\dfrac{5}{4}\sum_{n\ge0}(-1)^n\dfrac{n!}{(n+1)^2(3/2)_n}2^{-2n}$$
Parametric family with $k\ge0$:
\begin{verbatim}
[()->zeta(3),3*n^3+(11*k+1)*n^2-(11*k+3)*n+3*k+1,
            n^4*(2*n-k)*(2*n-k-1)]
\end{verbatim}
Convergence type $E$ with $E=-4$ and $P=5(k+1/2)$.
\end{cf}

\smallskip

\begin{cf}\label{5.5.19}{\ }
\begin{verbatim}
[()->zeta(3),[[6/7,21],[56*n^2+4*n+2,56*n^2+76*n+26]],
             [[8,126],[-16*(3*n+1)*(n+1)^4,54*(n+1)*(3*n+2)*(2*n+1)^3]]]
\end{verbatim}
$$\z(3)=6/7+\dfrac{8}{21+\dfrac{126}{62-\dfrac{1024}{158+\dfrac{14580}{234-\dfrac{9072}{402+\dfrac{162000}{518-\ddots}}}}}}$$
Convergence type $E$ with $E=-3\sqrt{3}$, $P=3/2$, and
$C=\pi^{5/2}/(7\sqrt{54})$, so that
$$\z(3)-\dfrac{p(n)}{q(n)}\sim(-1)^n\dfrac{\pi^{5/2}/(7\sqrt{2})}{3^{3n/2+3/2}n^{3/2}}\;.$$
$$A=1-(31/18)/n+(9523/2592)/n^2-(350749/34992)/n^3+\cdots$$
Series:
$$\z(3)=\dfrac{6}{7}+\dfrac{8}{21}\sum_{n\ge0}(-1)^n\dfrac{(28n^2+58n+31)n!(n+1)!^4}{(8n^2+7n+2)(8n^2+23n+17)(3/2)_n^3(5/3)_n(7/3)_n}3^{-3n}$$
\end{cf}

\smallskip

\begin{cf}\label{5.5.20}{\ }
\begin{verbatim}
[()->zeta(3),[0,(2*n-1)*(3*n^2-3*n+1)],[8/7,-n^6]]
\end{verbatim}
$$\z(3)=\dfrac{8/7}{1-\dfrac{1}{21-\dfrac{64}{95-\dfrac{729}{259-\dfrac{4096}{549-\dfrac{15625}{1001-\ddots}}}}}}$$
Convergence type $E$ with $E=(1+\sqrt{2})^4$, $P=0$, and $C=4\pi^3/7/(1+\sqrt{2})^2$, so that
$$\z(3)-\dfrac{p(n)}{q(n)}\sim\dfrac{4\pi^3/7}{(1+\sqrt{2})^{4n+2}}\;.$$
$$A=1-(7d/16)/n+(7d/32+49/256)/n^2+(-511d/4096-49/256)/n^3+\cdots$$
\end{cf}

\smallskip

\begin{cf}\label{5.5.21}{\ }
\begin{verbatim}
[()->zeta(3),[0,65*n^4-130*n^3+105*n^2-40*n+6],[7,-4*n^6*(16*n^2-1)]]
\end{verbatim}
$$\z(3)=\dfrac{7}{6-\dfrac{60}{346-\dfrac{16128}{2586-\dfrac{416988}{9846-\dfrac{4177920}{26806-\dfrac{24937500}{59706-\ddots}}}}}}$$
Convergence type $E$ with $E=64$, $P=0$, and $C=\pi^3/(6\sqrt{2})$, so that
$$\z(3)-\dfrac{p(n)}{q(n)}\sim\dfrac{\pi^3/3}{2^{6n+3/2}}\;.$$
$$A=1-(35/48)/n+(2905/4608)/n^2-(308539/663552)/n^3+\cdots$$
\end{cf}

\smallskip

\begin{cf}\label{5.5.22}{\ }
\begin{verbatim}
[()->zeta(3),[0,(2*n-1)*(17*n^2-17*n+5)],[6,-n^6]]
\end{verbatim}
$$\z(3)=\dfrac{6}{5-\dfrac{1}{117-\dfrac{64}{535-\dfrac{729}{1463-\dfrac{4096}{3105-\dfrac{15625}{5665-\ddots}}}}}}$$
Convergence type $E$ with $E=(1+\sqrt{2})^8$, $P=0$, and
$C=4\pi^3/(1+\sqrt{2})^4$, so that
$$\z(3)-\dfrac{p(n)}{q(n)}\sim\dfrac{4\pi^3}{(1+\sqrt{2})^{8n+4}}\;.$$
$$A=1-(15d/32)/n+(15d/64+225/1024)/n^2+(-3943d/32768-225/1024)/n^3+\cdots$$
\end{cf}

This is the famous continued fraction due to R.~Ap\'ery.

\smallskip

Note that there exist continued fractions of period $2$ for $\z(3)$ coming from
those of $\psi''(z)$, with convergence in $C/(1+\sqrt{2})^{4n}$, which have
little interest since their contraction gives Ap\'ery's CF.

\medskip

\section{Constants: $\pi^4$}

\medskip

\begin{cf}\label{5.6.1}{\ }
\begin{verbatim}
[()->Pi^4/90,[0,2*n^4-4*n^3+6*n^2-4*n+1],[1,-n^8]]
\end{verbatim}
$$\dfrac{\pi^4}{90}=\dfrac{1}{1-\dfrac{1}{17-\dfrac{256}{97-\dfrac{6561}{337-\dfrac{65536}{881-\dfrac{390625}{1921-\ddots}}}}}}$$
Convergence type $P^+$ with $P=3$ and $C=1/3$, so that
$$\dfrac{\pi^4}{90}-\dfrac{p(n)}{q(n)}\sim\dfrac{1/3}{n^3}\;.$$
$$A=1-(3/2)/n+1/n^2-(1/2)/n^4+(2/3)/n^6-(3/2)/n^8+5/n^{10}-(691/30)/n^{12}+\cdots$$
Series:
$$\dfrac{\pi^4}{90}=\sum_{n\ge1}\dfrac{1}{n^4}$$
\end{cf}

\smallskip

\begin{cf}\label{5.6.2}{\ }
\begin{verbatim}
[()->Pi^4/90,[0,1,32*n^4-128*n^3+240*n^2-224*n+82],[16/15,-(2*n-1)^8]]
\end{verbatim}
$$\dfrac{\pi^4}{90}=\dfrac{16/15}{1-\dfrac{1}{82-\dfrac{6561}{706-\dfrac{390625}{3026-\dfrac{5764801}{8962-\dfrac{43046721}{21202-\ddots}}}}}}$$
Convergence type $P^+$ with $P=3$ and $C=1/45$, so that
$$\dfrac{\pi^4}{90}-\dfrac{p(n)}{q(n)}\sim\dfrac{1/45}{n^3}\;.$$
$$A=1-(1/2)/n^2+(7/16)/n^4-(31/48)/n^6+(381/256)/n^8-(2555/512)/n^{10}+\cdots$$
Series:
$$\dfrac{\pi^4}{90}=\dfrac{16}{15}\sum_{n\ge1}\dfrac{1}{(2n-1)^4}$$
\end{cf}

\smallskip

\begin{cf}\label{5.6.3}{\ }
\begin{verbatim}
[()->Pi^4/90,[0,4*n^3-6*n^2+4*n-1],[8/7,n^8]]
\end{verbatim}
$$\dfrac{\pi^4}{90}=\dfrac{8/7}{1+\dfrac{1}{15+\dfrac{256}{65+\dfrac{6561}{175+\dfrac{65536}{369+\dfrac{390625}{671+\ddots}}}}}}$$
Convergence type $P^-$ with $P=4$ and $C=4/7$, so that
$$\dfrac{\pi^4}{90}-\dfrac{p(n)}{q(n)}\sim(-1)^n\dfrac{4/7}{n^4}\;.$$
$$A=1-2/n+5/n^3-28/n^5+255/n^7-3410/n^9+62881/n^{11}-1529080/n^{13}+\cdots$$
Series:
$$\dfrac{\pi^4}{90}=\dfrac{8}{7}\sum_{n\ge1}\dfrac{(-1)^{n+1}}{n^4}$$
\end{cf}

\smallskip

\begin{cf}\label{5.6.4}{\ }
\begin{verbatim}
[()->Pi^4/90,[13/12,2*n^4+4*n^3+23*n^2+21*n+6],
             [-1/18,-n*(n+1)^6*(n+2)]]
\end{verbatim}
$$\dfrac{\pi^4}{90}=13/12-\dfrac{1/18}{56-\dfrac{192}{204-\dfrac{5832}{546-\dfrac{61440}{1226-\dfrac{375000}{2436-\dfrac{1632960}{4416-\ddots}}}}}}\;.$$
Convergence type $P^+$ with $P=9$ and $C=-1/9$, so that
$$\dfrac{\pi^4}{90}-\dfrac{p(n)}{q(n)}\sim-\dfrac{1/9}{n^9}\;.$$
$$A=1-(27/2)/n+(195/2)/n^2-495/n^3+(3927/2)/n^4-6435/n^5+\cdots$$
Series:
$$\dfrac{\pi^4}{90}=\dfrac{13}{12}-\sum_{n\ge1}\dfrac{1}{n(n+1)^4(n+2)(n^2+n+1)(n^2+3n+3)}$$
\end{cf}

\smallskip

\begin{cf}\label{5.6.5}{\ }
\begin{verbatim}
[()->Pi^4/90,[0,5*n^4-6*n^3+6*n^2-4*n+1],[36/17,-2*n^7*(2*n-1)]]
\end{verbatim}
$$\dfrac{\pi^4}{90}=\dfrac{36/17}{2-\dfrac{2}{49-\dfrac{768}{286-\dfrac{21870}{977-\dfrac{229376}{2506-\dfrac{1406250}{5377-\ddots}}}}}}$$
Convergence type $E$ with $E=4$, $P=7/2$, and $C=12\sqrt{\pi}/17$, so that
$$\dfrac{\pi^4}{90}-\dfrac{p(n)}{q(n)}\sim\dfrac{12\sqrt{\pi}/17}{2^{2n}n^{7/2}}\;.$$
$$A=1-(109/24)/n+(2145/128)/n^2-(68821/1024)/n^3+(10294315/32768)/n^4+\cdots$$
Series:
$$\dfrac{\pi^4}{90}=\dfrac{18}{17}\sum_{n\ge0}\dfrac{n!}{(n+1)^3(3/2)_n}2^{-2n}$$
\end{cf}

\smallskip

\begin{cf}\label{5.6.6}
\begin{verbatim}
[()->Pi^4/90,[0,P_{10}(n)],[-16/15,P_{20}(n)]]
\end{verbatim}
with
\begin{align*}P_{10}(n)&=2730n^{10}-20626n^9+68118n^8-129564n^7+157416n^6\\
  &\phantom{=}-127932n^5+70602n^4-26187n^3+6261n^2-873n+54\\
P_{20}(n)&=-2n^7(2n-1)^7(3n-4)(3n+2)(7n^2-19n+13)(7n^2+9n+3)\end{align*}
$$\dfrac{\pi^4}{90}=\dfrac{16/15}{1-\dfrac{190}{13312-\dfrac{1316818944}{7522227-\dfrac{33209936718750}{327303458-\dfrac{22362855885504512}{4962248839-\ddots}}}}}$$
Convergence type $E$ with $E=64$, $P=1/2$, and $C=4\pi^{7/2}/135$, so that
$$\dfrac{\pi^4}{90}-\dfrac{p(n)}{q(n)}\sim\dfrac{\pi^{7/2}/135}{2^{6n-2}n^{1/2}}$$
$$A=1-(85/72)/n+(6563/3456)/n^2+\cdots$$
Series:
$$\dfrac{\pi^4}{90}=\dfrac{8}{45}\sum_{n\ge0}\dfrac{(3n+2)(7n^2+9n+3)n!^7}{(3/2)_n^7}2^{-6n}$$
\end{cf}

\smallskip

\begin{cf}\label{5.6.7}{\ }
\begin{verbatim}
[()->Pi^4/90,[0,3*(2*n-1)*(3*n^2-3*n+1)*(15*n^2-15*n+4)],
             [13,3*n^8*(9*n^2-1)]]
\end{verbatim}
$$\dfrac{\pi^4}{90}=\dfrac{13}{12+\dfrac{24}{2142+\dfrac{26880}{26790+\dfrac{1574640}{142968+\dfrac{28114944}{500688+\dfrac{262500000}{1363362+\ddots}}}}}}$$
Convergence type $E$ with $E=-(2+\sqrt{3})^6$, $P=0$, and $C=(8/3)\pi^4/(2+\sqrt{3})^3$, so that
$$\dfrac{\pi^4}{90}-\dfrac{p(n)}{q(n)}\sim(-1)^n\dfrac{(8/3)\pi^4}{(2+\sqrt{3})^{6n+3}}\;.$$
$$A=1-(13d/24)/n+(13d/48+169/384)/n^2+(-1573d/9216-169/384)/n^3+\cdots$$
\end{cf}

Comment: this continued fraction for $\z(4)=\pi^4/90$ was found in 1980
by G.~Rhin and the author by generalizing Ap\'ery's method, see Section
\ref{sec:zeta4}.

\medskip

\section{Constants: $\z(k)$ and $\pi^k$ for $k\ge5$}

To the author's knowledge, for $k\ge5$ there do not exist any (reasonably
simple) CFs for $\z(k)$ or $\pi^k$ other than the ones trivially obtained
thanks to Euler's transformation from the standard Dirichlet series definitions
and variants, which are not really interesting:

\smallskip

\begin{cf}\label{5.7.1}{\ }
\begin{verbatim}
[k->zeta(k),[0,n^k+(n-1)^k],[1,-n^(2*k)]]
\end{verbatim}
$$\z(k)=\dfrac{1}{1-\dfrac{1^{2k}}{2^k+1^k-\dfrac{2^{2k}}{3^k+2^k-\dfrac{3^{2k}}{4^k+3^k-\dfrac{4^{2k}}{5^k+4^k-\ddots}}}}}\;.$$
Convergence type $P^+$ with $P=k-1$ and $C=1/(k-1)$, so that
$$\z(k)-\dfrac{p(n)}{q(n)}\sim\dfrac{1/(k-1)}{n^{k-1}}\;.$$
$$A=1-((k-1)/2)/n+\sum_{m\ge1}\binom{k+2m-2}{k-2}B_{2m}/n^{2m}\;.$$
Series:
$$\z(k)=\sum_{n\ge1}\dfrac{1}{n^k}$$
\end{cf}

\smallskip

\begin{cf}\label{5.7.2}{\ }
\begin{verbatim}
[k->zeta(k),[0,n^k-(n-1)^k],[2^(k-1)/(2^(k-1)-1),n^(2*k)]]
\end{verbatim}
$$\z(k)=\dfrac{2^{k-1}/(2^{k-1}-1)}{1+\dfrac{1^{2k}}{2^k-1+\dfrac{2^{2k}}{3^k-2^k+\dfrac{3^{2k}}{4^k-3^k+\dfrac{4^{2k}}{5^k-4^k+\ddots}}}}}\;.$$
Convergence type $P^-$ with $P=k$ and $C=2^{k-2}/(2^{k-1}-1)$, so that
$$\z(k)-\dfrac{p(n)}{q(n)}\sim(-1)^n\dfrac{2^{k-2}/(2^{k-1}-1)}{n^k}\;.$$
$$A=1-(2/(k-1))\sum_{m\ge1}(2^{2m}-1)\binom{k+2m-2}{k-2}B_{2m}/n^{2m-1}\;.$$
Series:
$$\z(k)=\dfrac{2^{k-1}}{2^{k-1}-1}\sum_{n\ge1}\dfrac{(-1)^{n+1}}{n^k}$$
\end{cf}

\smallskip

\begin{cf}\label{5.7.3}{\ }
\begin{verbatim}
[k->zeta(k),[0,1,(2*n-1)^k+(2*n-3)^k],[2^k/(2^k-1),-(2*n-1)^(2*k)]]
\end{verbatim}
$$\z(k)=\dfrac{2^k/(2^k-1)}{1-\dfrac{1^{2k}}{3^k+1^k-\dfrac{3^{2k}}{5^k+3^k-\dfrac{5^{2k}}{7^k+5^k-\dfrac{7^{2k}}{9^k+7^k-\ddots}}}}}\;.$$
Convergence type $P^+$ with $P=k-1$ and $C=1/((k-1)(2^k-1))$, so that
$$\z(k)-\dfrac{p(n)}{q(n)}\sim\dfrac{1/((k-1)(2^k-1))}{n^{k-1}}\;.$$
$$A=1-\sum_{m\ge1}(1-1/2^{2m-1})\binom{k+2m-2}{k-2}B_{2m}/n^{2m}\;.$$
Series:
$$\z(k)=\dfrac{2^k}{2^k-1}\sum_{n\ge1}\dfrac{1}{(2n-1)^k}$$
\end{cf}

\smallskip

Of course, for $k$ even the above three CFs also give CFs for $\pi^k$. For
$k$ odd we have the following:

\smallskip

\begin{cf}\label{5.7.4}{\ }
\begin{verbatim}
[k->Pi^(2*k+1),[0,1,(2*n-1)^(2*k+1)-(2*n-3)^(2*k+1)],
               [(-1)^k*(2*k)!*2^(2*k+2)/eulerfrac(2*k),(2*n-1)^(4*k+2)]]
\end{verbatim}
$$\pi^{2k+1}=\dfrac{(-1)^k(2k)!2^{2k+2}/E_{2k}}{1+\dfrac{1^{4k+2}}{3^{2k+1}-1^{2k+1}+\dfrac{3^{4k+2}}{5^{2k+1}-3^{2k+1}+\dfrac{5^{4k+2}}{7^{2k+1}-5^{2k+1}+\ddots}}}}\;.$$
Convergence type $P^-$ with $P=2k+1$ and $C=(-1)^k(2k)!/E_{2k}$, so that
$$\pi^{2k+1}-\dfrac{p(n)}{q(n)}\sim(-1)^{n+k}\dfrac{(2k)!/E_{2k}}{n^{2k+1}}\;.$$
$$A=1+\sum_{m\ge1}\binom{2k+2m}{2k}E_{2m}/(2n)^{2m}\;.$$
Series:
$$\pi^{2k+1}=(-1)^k\dfrac{(2k)!2^{2k+2}}{E_{2k}}\sum_{n\ge1}\dfrac{(-1)^{n+1}}{(2n-1)^{2k+1}}$$
\end{cf}

\medskip

\section{Constants: Linear Combinations of Zeta- and $L$-Values}

\medskip

First note that creating continued fractions for linear combinations of
zeta and $L$-values with rational coefficients is considerably \emph{easier}
than for a single such value:
for instance, if $P(x)$ and $Q(x)$ are polynomials with rational coefficients
such that $P$ is a product of rational linear factors with no nonnegative
integral root, then both $\sum_{n\ge0}Q(n)/P(n)$ and
$\sum_{n\ge0}(-1)^nQ(n)/P(n)$ are such linear combinations when the series
converge ($\deg(P)\ge\deg(Q)+2$ for the former and $\deg(P)\ge\deg(Q)+1$
for the latter), and can then be converted to CFs using Euler's transformation.
Nonetheless I believe it useful to include a large list.

\smallskip

\subsection{Inhomogeneous Periods}

\smallskip

In \cite{Coh3} we have given three systematic ways of constructing such
CFs, and we refer to that paper for the justification of the polynomials
$R_i(k,n)$ which we give below. In the present section we include all the
examples from loc.~cit., and add many more. A notable characteristic of the
present CFs are that their limits are \emph{inhomogenous} periods, i.e.,
linear combinations of periods of different degrees. As mentioned above,
these are much easier to obtain than homogeneous ones, but we give many anyway.

\smallskip

In everything that follows, $k$ stands for a nonnegative integer.

\smallskip

\begin{verbatim}
R1(k,n)=(n^k*(2*n+1)+(n-1)^k*(2*n-3))/(2*n-1);
\end{verbatim}

\smallskip

\begin{cf}\label{5.8.1}{\ }
\begin{verbatim}
[k->sum(j=0,k-1,4^j*zeta(2*(k-j))),[2^(2*k-1),R1(2*k,n)],[-1,-n^(4*k)]]
\end{verbatim}
Convergence type $P^+$ with $P=2k+1$, $C=-1/(4(2k+1))$, so that
$$\sum_{j=0}^{k-1}4^j\z(2(k-j))-\dfrac{p(n)}{q(n)}\sim-\dfrac{1/(4(2k+1))}{n^{2k+1}}\;.$$
Series:
$$\sum_{j=0}^{k-1}4^j\z(2(k-j))=2^{2k-1}-\sum_{n\ge1}\dfrac{1}{n^{2k}(4n^2-1)}$$
\end{cf}

See \ref{5.8.9} for the simplest example.

\smallskip

\begin{verbatim}
R2(k,n)=(n^k*(3*n^2+3*n+1)+(n-1)^k*(3*n^2-9*n+7))/(3*n^2-3*n+1);
\end{verbatim}

\smallskip

\begin{cf}\label{5.8.2}{\ }
\begin{verbatim}
[k->sum(j=0,k,(-3)^(3*j)*(zeta(6*(k-j)+2)
                         +3*zeta(6*(k-j))*(1-(k==j)))),
[-(-3)^(3*k+1)/2,R2(6*k+2,n)],[1,-n^(12*k+4)]]
\end{verbatim}
Convergence type $P^+$ with $P=6k+5$ and $C=1/(9(6k+5))$, so that
$$\sum_{j=0}^k(-3)^{3j}(\z(6(k-j)+2)+3\z(6(k-j))(1-\delta_{k,j}))-\dfrac{p(n)}{q(n)}\sim\dfrac{1/(9(6k+5))}{n^{6k+5}}\;.$$
Series:
$$\sum_{j=0}^k(-3)^{3j}(\z(6(k-j)+2)+3\z(6(k-j))(1-\delta_{k,j}))=(-1)^k\dfrac{3^{3k+1}}{2}+\sum_{n\ge1}\dfrac{1}{n^{6k+2}(9n^4-3n^2+1)}$$
\end{cf}

See \ref{5.8.14} for the simplest example.

\smallskip

\begin{verbatim}
R3(k,n)=(n^k*(2*n^2+2*n+1)+(n-1)^k*(2*n^2-6*n+5))/(2*n^2-2*n+1);
\end{verbatim}

\smallskip

\begin{cf}\label{5.8.3}{\ }
\begin{verbatim}
[k->sum(j=0,k,(-4)^j*zeta(4*(k-j)+3)),
         [(-4)^k,R3(4*k+3,n)],[1,-n^(8*k+6)]]
\end{verbatim}
Convergence type $P^+$ with $P=4k+6$, $C=1/(4(4k+6))$, so that
$$\sum_{j=0}^k(-4)^j\z(4(k-j)+3)-\dfrac{p(n)}{q(n)}\sim\dfrac{1/(4(4k+6))}{n^{4k+6}}\;.$$
Series:
$$\sum_{j=0}^k(-4)^j\z(4(k-j)+3)=(-4)^k+\sum_{n\ge1}\dfrac{1}{n^{4k+3}(4n^4+1)}$$
\end{cf}

See \ref{5.8.13} for the simplest example.

\smallskip

\begin{verbatim}
R4(k,n)=(n^k*(n^2+n+1)+(n-1)^k*(n^2-3*n+3))/(n^2-n+1);
\end{verbatim}

\smallskip

\begin{cf}\label{5.8.4}{\ }
\begin{verbatim}
[k->sum(j=0,k,zeta(6*(k-j)+5)-zeta(6*(k-j)+3)),
        [-1/2,R4(6*k+5,n)],[1,-n^(12*k+10)]]
\end{verbatim}
Convergence type $P^+$ with $P=6k+8$ and $C=1/(6k+8)$, so that
$$\sum_{j=0}^k(\z(6(k-j)+5)-\z(6(k-j)+3))-\dfrac{p(n)}{q(n)}\sim\dfrac{1/(6k+8)}{n^{6k+8}}\;.$$
Series:
$$\sum_{j=0}^k(\z(6(k-j)+5)-\z(6(k-j)+3))=-\dfrac{1}{2}+\sum_{n\ge1}\dfrac{1}{n^{6k+5}(n^4+n^2+1)}$$
\end{cf}

See \ref{5.8.11} for the simplest example.

\smallskip

\begin{verbatim}
R5(k,n)=(n^k*(3*n^2+3*n+1)+(n-1)^k*(3*n^2-9*n+7))/(3*n^2-3*n+1);
\end{verbatim}

\smallskip

\begin{cf}\label{5.8.5}{\ }
\begin{verbatim}
[k->sum(j=0,k,(-3)^(3*j)*(zeta(6*(k-j)+5)+3*zeta(6*(k-j)+3))),
    [(-3)^(3*k+2)/2,R5(6*k+5,n)],[1,-n^(12*k+10)]]
\end{verbatim}
Convergence type $P^+$ with $P=6k+8$ and $C=1/(9(6k+8))$, so that
$$\sum_{j=0}^k(-3)^{3j}(\z(6(k-j)+5)+3\z(6(k-j)+3))-\dfrac{p(n)}{q(n)}\sim\dfrac{1/(9(6k+8))}{n^{6k+8}}\;.$$
Series:
$$\sum_{j=0}^k(-3)^{3j}(\z(6(k-j)+5)+3\z(6(k-j)+3))=(-1)^k\dfrac{3^{3k+2}}{2}+\sum_{n\ge1}\dfrac{1}{n^{6k+5}(9n^4-3n^2+1)}$$
\end{cf}

See \ref{5.8.10} for the simplest example.

\smallskip

\begin{verbatim}
R6(k,n)=(n^k*(2*n+1)-(n-1)^k*(2*n-3))/(2*n-1);
zetastar(k)=if(k==0,0,if(k==1,log(2),(2^(k-1)-1)*zeta(k)));
\end{verbatim}

\smallskip

\begin{cf}\label{5.8.6}{\ }
\begin{verbatim}
[k->sum(j=0,k,2^(4*j)*zetastar(2*(k-j)+1)),
[2^(4*k),R6(2*k+1,n)],[-2^(2*k),n^(4*k+2)]]
\end{verbatim}
Convergence type $P^-$ with $P=2k+3$ and $C=-2^{2k-3}$, so that
$$\sum_{j=0}^k2^{4j}\zeta^*(2(k-j)+1)-\dfrac{p(n)}{q(n)}\sim(-1)^{n+1}\dfrac{2^{2k-3}}{n^{2k+3}}\;.$$
Series:
$$\sum_{j=0}^k2^{4j}\zeta^*(2(k-j)+1)=2^{4k}-2^{2k}\sum_{n\ge1}\dfrac{(-1)^{n+1}}{n^{2k+1}(4n^2-1)}$$
\end{cf}

See \ref{5.8.16} and \ref{5.8.17} for the simplest examples.

\smallskip

\begin{verbatim}
R7(k,n)=(n^k*(2*n^2+2*n+1)-(n-1)^k*(2*n^2-6*n+5))/(2*n^2-2*n+1);
\end{verbatim}

\smallskip

\begin{cf}\label{5.8.7}{\ }
\begin{verbatim}
[k->sum(j=0,k,(-64)^j*zetastar(4*(k-j)+1)),
[(-1)^k*2^(6*k-1),R7(4*k+1,n)],[2^(4*k),n^(8*k+2)]]
\end{verbatim}
Convergence type $P^-$ with $P=4k+5$ and $C=2^{4k-3}$, so that
$$\sum_{j=0}^k(-64)^j\zeta^*(4(k-j)+1)-\dfrac{p(n)}{q(n)}\sim(-1)^n\dfrac{2^{4k-3}}{n^{4k+5}}\;.$$
Series:
$$\sum_{j=0}^k(-64)^j\zeta^*(4(k-j)+1)=(-1)^k2^{6k-1}+2^{4k}\sum_{n\ge1}\dfrac{(-1)^{n+1}}{n^{4k+1}(4n^4+1)}$$
\end{cf}

See \ref{5.8.15} for the simplest example.

\smallskip

\begin{verbatim}
R8(k,n)=(n^k*(n^2+n+1)-(n-1)^k*(n^2-3*n+3))/(n^2-n+1);
\end{verbatim}

\smallskip

\begin{cf}\label{5.8.8}{\ }
\begin{verbatim}
[k->sum(j=0,k,2^(6*j)*(zetastar(6*(k-j)+2)-4*zetastar(6*(k-j)))),
[2^(6*k),R8(6*k+2,n)],[2^(6*k+1),n^(12*k+4)]]
\end{verbatim}
Convergence type $P^-$ with $P=6k+6$ and $C=2^{6k}$, so that
$$\sum_{j=0}^{k}2^{6j}(\zeta^*(6(k-j)+2)-4\zeta^*(6(k-j)))-\dfrac{p(n)}{q(n)}\sim(-1)^n\dfrac{2^{6k}}{n^{6k+6}}\;.$$
Series:
$$\sum_{j=0}^{k}2^{6j}(\zeta^*(6(k-j)+2)-4\zeta^*(6(k-j)))=2^{6k}+2^{6k+1}\sum_{n\ge1}\dfrac{(-1)^{n+1}}{n^{6k+2}(n^4+n^2+1)}$$
\end{cf}

See \ref{5.8.18} for the simplest example.

\smallskip

\begin{cf}\label{5.8.9}{\ }
\begin{verbatim}
[()->zeta(4)+4*zeta(2),[8,n^4+(n-1)^4+2*(n^2+(n-1)^2)],[-1,-n^8]]
\end{verbatim}
$$\z(4)+4\z(2)=8-\dfrac{1}{3-\dfrac{1}{27-\dfrac{256}{123-\dfrac{6561}{387-\dfrac{65536}{963-\dfrac{390625}{2043-\ddots}}}}}}$$
Convergence type $P^+$ with $P=5$ and $C=-1/20$, so that
$$\z(4)+4\z(2)-\dfrac{p(n)}{q(n)}\sim-\dfrac{1/20}{n^5}\;.$$
$$A=1-(5/2)/n+(75/28)/n^2-(5/8)/n^3-(211/144)/n^4-(5/32)/n^5+\cdots$$
Series:
$$\z(4)+4\z(2)=8-\sum_{n\ge1}\dfrac{1}{n^4(2n-1)(2n+1)}$$
\end{cf}

\smallskip

\begin{cf}\label{5.8.10}{\ }
\begin{verbatim}
[()->zeta(5)+3*zeta(3),[9/2,n^5+(n-1)^5+6*(n^3+(n-1)^3)],[1,-n^10]]
\end{verbatim}
$$\z(5)+3\z(3)=9/2+\dfrac{1}{7-\dfrac{1}{87-\dfrac{1024}{485-\dfrac{59049}{1813-\dfrac{1048576}{5283-\dfrac{9765625}{12947-\ddots}}}}}}$$
Convergence type $P^+$ with $P=8$ and $C=1/72$, so that
$$\z(5)+3\z(3)-\dfrac{p(n)}{q(n)}\sim\dfrac{1/72}{n^8}\;.$$
$$A=1-4/n+(94/15)/n^2-(4/3)/n^3-(77/9)/n^4+(32584/945)/n^6+(4/27)/n^7+\cdots$$
Series:
$$\z(5)+3\z(3)=\dfrac{9}{2}+\sum_{n\ge1}\dfrac{1}{n^5(3n^2-3n+1)(3n^2+3n+1)}$$
\end{cf}

\smallskip

\begin{cf}\label{5.8.11}{\ }
\begin{verbatim}
[()->zeta(3)-zeta(5),[1/2,n^5+(n-1)^5+6*(n^3+(n-1)^3)-4*(2*n-1)],
                    [-1,-n^10]]
\end{verbatim}
$$\z(3)-\z(5)=1/2-\dfrac{1}{3-\dfrac{1}{75-\dfrac{1024}{465-\dfrac{59049}{1785-\dfrac{1048576}{5247-\dfrac{9765625}{12903-\ddots}}}}}}$$
Convergence type $P^+$ with $P=8$ and $C=-1/8$, so that
$$\z(3)-\z(5)-\dfrac{p(n)}{q(n)}\sim-\dfrac{1/8}{n^8}\;.$$
$$A=1-4/n+(26/5)/n^2+4/n^3-(55/3)/n^4+(6352/105)/n^6-4/n^7+\cdots$$
Series:
$$\z(3)-\z(5)=\dfrac{1}{2}-\sum_{n\ge1}\dfrac{1}{n^5(n^2-n+1)(n^2+n+1)}$$
\end{cf}

\smallskip

\begin{cf}\label{5.8.12}{\ }
\begin{verbatim}
[()->4*zeta(5)+11*zeta(3),[273/16,n^5+(n-1)^5+
                            16*(n^3+(n-1)^3)-4*(2*n-1)],[4,-n^10]]
\end{verbatim}
$$4\z(5)+11\z(3)=273/16+\dfrac{4}{13-\dfrac{1}{165-\dfrac{1024}{815-\dfrac{59049}{2695-\dfrac{1048576}{7137-\dfrac{9765625}{16313-\ddots}}}}}}$$
Convergence type $P^+$ with $P=12$ and $C=16/75$, so that
$$4\z(5)+11\z(3)-\dfrac{p(n)}{q(n)}\sim\dfrac{16/75}{n^{12}}\;.$$
$$A=1-6/n+(347/35)/n^2+(108/5)/n^3-(9341/100)/n^4-(1218/25)/n^5+\cdots$$
Series:
$$4\z(5)+11\z(3)=\dfrac{273}{16}+64\sum_{n\ge1}\dfrac{1}{n^5(5n^4-10n^3+19n^2-14n+4)(5n^4+10n^3+19n^2+14n+4)}$$
\end{cf}

This is the only CF in this section which is not a special case of the
six infinite families \ref{5.8.1} to \ref{5.8.6}.

\smallskip

\begin{cf}\label{5.8.13}{\ }
\begin{verbatim}
[()->4*zeta(3)-zeta(7),[4,n^7+(n-1)^7+8*(n^5+(n-1)^5)
                       -8*(n^3+(n-1)^3)+4*(2*n-1)],[-1,-n^14]]
\end{verbatim}
$$4\z(3)-\z(7)=4-\dfrac{1}{5-\dfrac{1}{333-\dfrac{16384}{4255-\dfrac{4782969}{28007-\dfrac{268435456}{126225-\dfrac{6103515625}{442585-\ddots}}}}}}$$
Convergence type $P^+$ with $P=10$ and $C=-1/40$, so that
$$4\z(3)-\z(7)-\dfrac{p(n)}{q(n)}\sim-\dfrac{1/40}{n^{10}}\;.$$
$$A=1-5/n+(55/6)/n^2-(2017/84)/n^4+(5/4)/n^5+(2785/24)/n^6+\cdots$$
Series:
$$4\z(3)-\z(7)=4-\sum_{n\ge1}\dfrac{1}{n^7(2n^2-2n+1)(2n^2+2n+1)}$$
\end{cf}

\smallskip

\begin{cf}\label{5.8.14}{\ }
\begin{verbatim}
[()->27*zeta(2)-3*zeta(6)-zeta(8),
[81/2,n^8+(n-1)^8+9*(n^6+(n-1)^6)
    -27/2*(n^4+(n-1)^4)+18*(n^2+(n-1)^2)-15/2],[-1,-n^16]]
\end{verbatim}
$$27\z(2)-3\z(6)-\z(8)=81/2-\dfrac{1}{7-\dfrac{1}{695-\dfrac{65536}{12871-\dfrac{43046721}{111415-\dfrac{4294967296}{622487-\dfrac{152587890625}{2605927-\ddots}}}}}}\;.$$
Convergence type $P^+$ with $P=11$ and $C=-1/99$, so that
$$27\z(2)-3\z(6)-\z(8)-\dfrac{p(n)}{q(n)}\sim-\dfrac{1/99}{n^{11}}\;.$$
$$A=1-(11/2)/n+(440/39)/n^2-(11/6)/n^3-(1309/45)/n^4+(79651/459)/n^6+\cdots$$
Series:
$$27\z(2)-3\z(6)-\z(8)=\dfrac{81}{2}-\sum_{n\ge1}\dfrac{1}{n^8(3n^2-3n+1)(3n^2+3n+1)}$$
\end{cf}

\smallskip

\begin{cf}\label{5.8.15}{\ }
\begin{verbatim}
[()->64*log(2)-15*zeta(5),
[32,9*n^4-18*n^3+30*n^2-21*n+5],[-16,n^10]]
\end{verbatim}
$$64\log(2)-15\z(5)=32-\dfrac{16}{5+\dfrac{1}{83+\dfrac{1024}{455+\dfrac{59049}{1553+\dfrac{1048576}{4025+\dfrac{9765625}{8735+\ddots}}}}}}$$
Convergence type $P^-$ with $P=9$ and $C=-2$, so that
$$64\log(2)-15\z(5)-\dfrac{p(n)}{q(n)}\sim(-1)^{n+1}\dfrac{2}{n^9}\;.$$
$$A=1-(9/2)/n+(165/4)/n^3-(1/4)/n^4-(5135/8)/n^5+(218335/16)/n^7+\cdots$$
Series:
$$64\log(2)-15\z(5)=32-16\sum_{n\ge1}\dfrac{(-1)^{n+1}}{n^5(2n^2-2n+1)(2n^2+2n+1)}$$
\end{cf}

\smallskip

\begin{cf}\label{5.8.16}{\ }
\begin{verbatim}
[()->3*zeta(3)+16*log(2),[16,5*n^2-5*n+3],[-4,n^6]]
\end{verbatim}
$$3\z(3)+16\log(2)=16-\dfrac{4}{3+\dfrac{1}{13+\dfrac{64}{33+\dfrac{729}{63+\dfrac{4096}{103+\dfrac{15625}{153+\ddots}}}}}}$$
Convergence type $P^-$ with $P=5$ and $C=-1/2$, so that
$$3\z(3)+16\log(2)-\dfrac{p(n)}{q(n)}\sim(-1)^{n+1}\dfrac{1/2}{n^5}\;.$$
$$A=1-(5/2)/n+(1/4)/n^2+(63/8)/n^3+(1/16)/n^4-(1857/32)/n^5+(1/64)/n^6+\cdots$$
Series:
$$3\z(3)+16\log(2)=16-4\sum_{n\ge1}\dfrac{(-1)^{n+1}}{n^3(2n-1)(2n+1)}$$
\end{cf}

\smallskip

\begin{cf}\label{5.8.17}{\ }
\begin{verbatim}
[()->15*zeta(5)+48*zeta(3)+256*log(2),
[256,7*n^4-14*n^3+18*n^2-11*n+3],[-16,n^10]]
\end{verbatim}
$$15\z(5)+48\z(3)+256\log(2)=256-\dfrac{16}{3+\dfrac{1}{53+\dfrac{1024}{321+\dfrac{59049}{1143+\dfrac{1048576}{3023+\dfrac{9765625}{6633+\ddots}}}}}}$$
Convergence type $P^-$ with $P=7$ and $C=-2$, so that
$$15\z(5)+48\z(3)+256\log(2)-\dfrac{p(n)}{q(n)}\sim(-1)^{n+1}\dfrac{2}{n^7}\;.$$
$$A=1-(7/2)/n+(1/4)/n^2+(159/8)/n^3+(1/16)/n^4-(7073/32)/n^5+(1/64)/n^6+\cdots$$
Series:
$$15\z(5)+48\z(3)+256\log(2)=256-16\sum_{n\ge1}\dfrac{(-1)^{n+1}}{n^5(2n-1)(2n+1)}$$
\end{cf}

\smallskip

\begin{cf}\label{5.8.18}{\ }
\begin{verbatim}
[()->127*zeta(8)-124*zeta(6)+64*zeta(2),
[64,6*(n^7+(n-1)^7)-8*(n^5+(n-1)^5)+12*(n^3+(n-1)^3)-7*(2*n-1)],
[128,n^16]]
\end{verbatim}
$$127\z(8)-124\z(6)+64\z(2)=64+\dfrac{128}{3+\dfrac{1}{597+\dfrac{65536}{12075+\dfrac{43046721}{102333+\dfrac{4294967296}{536067+\ddots}}}}}$$
Convergence type $P^-$ with $P=12$ and $C=64$, so that
$$127\z(8)-124\z(6)+64\z(2)-\dfrac{p(n)}{q(n)}\sim(-1)^n\dfrac{64}{n^{12}}\;.$$
$$A=1-6/n-1/n^2+98/n^3-2324/n^5+1/n^6+71901/n^7-1/n^8-2767815/n^9+\cdots$$
Series:
$$127\z(8)-124\z(6)+64\z(2)=64+128\sum_{n\ge1}\dfrac{(-1)^{n+1}}{n^8(n^2-n+1)(n^2+n+1)}$$
\end{cf}

\smallskip

\begin{cf}\label{5.8.19}{\ }
\begin{verbatim}
[()->7*zeta(3)-2*Pi*Catalan+3*Pi,[12,16*n^3-8*n^2+22*n-3],
                                 [2,-4*n^2*(2*n+1)^4]]
\end{verbatim}
$$7\z(3)-2\pi G+3\pi=12+\dfrac{2}{27-\dfrac{324}{137-\dfrac{10000}{423-\dfrac{86436}{981-\dfrac{419904}{1907-\dfrac{1464100}{3297-\ddots}}}}}}$$
Convergence type $P^+$ with $P=3$ and $C=\pi/48$, so that
$$7\z(3)-2\pi G+3\pi-\dfrac{p(n)}{q(n)}\sim\dfrac{\pi/48}{n^3}$$
$$A=1-(39/16)/n+\cdots$$
Series:
$$7\z(3)-2\pi G+3\pi=12+\dfrac{2}{9}\sum_{n\ge0}\dfrac{n!^2}{(2n+3)(5/2)_n^2}$$
\end{cf}

\medskip

\subsection{Homogeneous Periods of Degree $1$}

\smallskip

The next CFs come from continued fractions for $\psi(z)$, $\be(z)$, $\psi'(z)$,
$\be_1(z)$, $\psi''(z)$, or similar functions. In view of the trivial
functional equations giving for instance $\psi^{(k)}(z+1)-\psi^{(k)}(z)$,
one can trivially construct infinitely
many analogous continued fractions from the given ones. See the tables
following \ref{6.10.6}, \ref{6.11.6}, \ref{6.12.3}, and \ref{6.13.9} for
additional CFs. Note that the limits of the following CFs are now
\emph{homogeneous} periods, i.e., linear combinations of periods of the
same degree, contrary to those given above. As already mentioned, these are
more difficult to construct than inhomogenous ones.

\smallskip

\begin{cf}\label{5.8.20}{\ }
\begin{verbatim}
[()->Pi-2*log(2),[-2,8*n^2-12*n+5],[6,-n^2*(4*n+1)*(4*n-5)]]
\end{verbatim}
$$\pi-2\log(2)=-2+\dfrac{6}{1+\dfrac{5}{13-\dfrac{108}{41-\dfrac{819}{85-\dfrac{2992}{145-\dfrac{7875}{221-\ddots}}}}}}$$
Convergence type $P^+$ with $P=3/2$ and $C=-2\sqrt{2}\G(3/4)^2/\pi$, so that
$$\pi-2\log(2)-\dfrac{p(n)}{q(n)}\sim-\dfrac{2\sqrt{2}\G(3/4)^2/\pi}{n^{3/2}}$$
$$A=1+(3/10)/n+(11/64)/n^2+(13/128)/n^3+(381/8192)/n^4+\cdots$$
Series:
$$\pi-2\log(2)=-2-6\sum_{n\ge0}\dfrac{(5/4)_n}{(n+1)(4n^2-1)(3/4)_n}$$
Parametric family for $k\ge0$:
\begin{verbatim}
[()->Pi-2*log(2),8*n^2-12*n+5+2*k*(2*k+3),-n^2*(4*n+1)*(4*n-5)]
\end{verbatim}
Convergence type $P^+$ with $P=2k+3/2$.
\end{cf}

\smallskip

\begin{cf}\label{5.8.21}{\ }
\begin{verbatim}
[()->Pi-2*log(2),[-2,2,8*n^2-4*n+3],
                 [10,(n+1)^2*(2*n-1)^2*(4*n+5)*(4*n-3)]]
\end{verbatim}
$$\pi-2\log(2)=-2+\dfrac{10}{2+\dfrac{36}{27+\dfrac{5265}{63+\dfrac{61200}{115+\dfrac{334425}{183+\dfrac{1239300}{267+\ddots}}}}}}$$
Convergence type $P^-$ with $P=1$ and $C=2$, so that
$$\pi-2\log(2)-\dfrac{p(n)}{q(n)}\sim(-1)^n\dfrac{2}{n}$$
$$A=1-(3/4)/n+(7/8)/n^2-(9/8)/n^3+(37/32)/n^4+\cdots$$
Series:
$$\pi-2\log(2)=-2+2\sum_{n\ge0}(-1)^n\dfrac{4n+5}{(n+2)(2n+1)}$$
\end{cf}

\smallskip

\begin{cf}\label{5.8.22}{\ }
\begin{verbatim}
[()->Pi-2*log(2),[0,4*n-3],[2,n^2*(2*n-1)^2]]
\end{verbatim}
$$\pi-2\log(2)=\dfrac{2}{1+\dfrac{1}{5+\dfrac{36}{9+\dfrac{225}{13+\dfrac{784}{17+\dfrac{2025}{21+\ddots}}}}}}$$
Convergence type $P^-$ with $P=2$ and $C=1/2$, so that
$$\pi-2\log(2)-\dfrac{p(n)}{q(n)}\sim(-1)^n\dfrac{1/2}{n^2}$$
$$A=1-(1/2)/n-(1/2)/n^2+(5/8)/n^3+1/n^4-(61/32)/n^5+\cdots$$
Series:
$$\pi-2\log(2)=2\sum_{n\ge0}\dfrac{(-1)^n}{(n+1)(2n+1)}$$
Parametric family for $k\ge0$ and $u$:
\begin{verbatim}
[()->Pi-2*log(2),(2*k+1)*(4*n+4*u+1),n^2*(2*n+4*u+3)^2]
\end{verbatim}
Convergence type $P^-$ with $P=4k+2$.
\end{cf}

\smallskip

\begin{cf}\label{5.8.23}{\ }
\begin{verbatim}
[()->Pi-2*log(2),[0,8*n^2-8*n+3],[2,-n^2*(16*n^2-1)]]
\end{verbatim}
$$\pi-2\log(2)=\dfrac{2}{3-\dfrac{15}{19-\dfrac{252}{51-\dfrac{1287}{99-\dfrac{4080}{163-\dfrac{9975}{243-\ddots}}}}}}$$
Convergence type $P^+$ with $P=1/2$ and $C=2^{5/2}\pi/\G(1/4)^2$, so that
$$\pi-2\log(2)-\dfrac{p(n)}{q(n)}\sim\dfrac{2^{3/2}\pi/\G(1/4)^2}{n^{1/2}}$$
$$A=1-(1/4)/n+(19/320)/n^2+(1/256)/n^3+\cdots$$
Series:
$$\pi-2\log(2)=\dfrac{2}{3}\sum_{n\ge0}\dfrac{(5/4)_n}{(n+1)(7/4)_n}$$
Parametric family for $k\ge0$ and $u$:
\begin{verbatim}
[()->Pi-2*log(2),8*n^2-8*n+3+(4*u+2*k)*(2*k+1),-n^2*(16*n^2-(4*u-1)^2)]
\end{verbatim}
Convergence type $P^+$ with $P=2k+2u+1/2$.
\end{cf}

\smallskip

\begin{cf}\label{5.8.24}{\ }
\begin{verbatim}
[()->Pi-2*log(2),[0,20*n^2-23*n+7],[6,-4*n*(2*n-1)^2*(4*n+1)]]
\end{verbatim}
$$\pi-2\log(2)=\dfrac{6}{4-\dfrac{20}{41-\dfrac{648}{118-\dfrac{3900}{235-\dfrac{13328}{392-\dfrac{34020}{589-\ddots}}}}}}$$
Convergence type $E$ with $E=4$, $P=3/4$, and $C=2^{1/2}\G(3/4)$, so that
$$\pi-2\log(2)-\dfrac{p(n)}{q(n)}\sim\dfrac{\G(3/4)}{2^{2n-1/2}n^{3/4}}$$
$$A=1-(27/32)/n+(2249/2048)/n^2+\cdots$$
Series:
$$\pi-2\log(2)=\dfrac{3}{2}\sum_{n\ge0}\dfrac{(5/4)_nn!}{(3/2)_n^2}2^{-2n}$$
Parametric family for $k\ge0$:
\begin{verbatim}
[()->Pi-2*log(2),20*n^2+(8*k-23)*n+7-5*k,-2*(2*n-k)*(2*n-k-1)*(2*n-1)*(4*n+1)]
\end{verbatim}
Convergence type $E$ with $E=4$ and $P=3k+3/4$.
\end{cf}

\smallskip

\begin{cf}\label{5.8.25}{\ }
\begin{verbatim}
[()->Pi-2*log(2),[0,29*n^2-30*n+7],[10,-3*n*(2*n+1)*(3*n-2)*(3*n-1)]]
\end{verbatim}
$$\pi-2\log(2)=\dfrac{10}{6-\dfrac{18}{63-\dfrac{600}{178-\dfrac{3528}{351-\dfrac{11880}{582-\dfrac{30030}{871-\ddots}}}}}}$$
Convergence type $E$ with $E=27/2$, $P=1/2$, and $C=(8/5)\sqrt{\pi/3}$, so that
$$\pi-2\log(2)-\dfrac{p(n)}{q(n)}\sim\dfrac{(8/5)\sqrt{\pi/3}}{(27/2)^nn^{1/2}}$$
$$A=1-(797/1800)/n+(2054329/6480000)/n^2+\cdots$$
Series:
$$\pi-2\log(2)=\dfrac{5}{3}\sum_{n\ge0}\dfrac{n!(3/2)_n}{(4/3)_n(5/3)_n}(27/2)^{-n}$$
\end{cf}
        
\smallskip

\begin{cf}\label{5.8.26}{\ }
\begin{verbatim}
[()->Pi+2*log(2),[0,8*n^2-12*n+5],[2,-n^2*(4*n-1)*(4*n-3)]]
\end{verbatim}
$$\pi+2\log(2)=\dfrac{2}{1-\dfrac{3}{13-\dfrac{140}{41-\dfrac{891}{85-\dfrac{3120}{145-\dfrac{8075}{221-\ddots}}}}}}$$
Convergence type $P^+$ with $P=1/2$ and $C=\G(1/4)^2/(\pi\sqrt{2})$, so that
$$\pi+2\log(2)-\dfrac{p(n)}{q(n)}\sim\dfrac{\G(1/4)^2/(\pi\sqrt{2})}{n^{1/2}}$$
$$A=1-(1/6)/n+(1/64)/n^2+(25/2688)/n^3-(19/8192)/n^4+\cdots$$
Series:
$$\pi+2\log(2)=2\sum_{n\ge0}\dfrac{(3/4)_n}{(n+1)(5/4)_n}$$
Parametric family for $k\ge0$:
\begin{verbatim}
[()->Pi+2*log(2),8*n^2-12*n+5+2*k*(2*k+1),-n^2*(4*n-1)*(4*n-3)]
\end{verbatim}
Convergence type $P^+$ with $P=2k+1/2$.
\end{cf}

\smallskip

\begin{cf}\label{5.8.27}{\ }
\begin{verbatim}
[()->Pi+2*log(2),[0,8*n^2-12*n+3],[-6,n^2*(2*n-1)^2*(4*n-5)*(4*n+3)]]
\end{verbatim}
$$\pi+2\log(2)=-\dfrac{6}{-1-\dfrac{7}{11+\dfrac{1188}{39+\dfrac{23625}{83+\dfrac{163856}{143+\dfrac{698625}{219+\ddots}}}}}}$$
Convergence type $P^-$ with $P=1$ and $C=2$, so that
$$\pi+2\log(2)-\dfrac{p(n)}{q(n)}\sim(-1)^n\dfrac{2}{n}$$
$$A=1-(1/4)/n-(1/8)/n^2+(1/8)/n^3+(5/32)/n^4-(1/4)/n^5+\cdots$$
Series:
$$\pi+2\log(2)=2\sum_{n\ge0}(-1)^n\dfrac{4n+3}{(n+1)(2n+1)}$$
\end{cf}

\smallskip

\begin{cf}\label{5.8.28}{\ }
\begin{verbatim}
[()->Pi+2*log(2),[4,4*n-1],[2,n^2*(2*n+1)^2]]
\end{verbatim}
$$\pi+2\log(2)=4+\dfrac{2}{3+\dfrac{9}{7+\dfrac{100}{11+\dfrac{441}{15+\dfrac{1296}{19+\dfrac{3025}{23+\ddots}}}}}}$$
Convergence type $P^-$ with $P=2$ and $C=1/2$, so that
$$\pi+2\log(2)-\dfrac{p(n)}{q(n)}\sim(-1)^n\dfrac{1/2}{n^2}$$
$$A=1-(3/2)/n+1/n^2+(3/8)/n^3-(7/8)/n^4+\cdots$$
Series:
$$\pi+2\log(2)=4+2\sum_{n\ge0}\dfrac{(-1)^n}{(n+1)(2n+3)}$$
Parametric family for $k\ge0$ and $u$:
\begin{verbatim}
[()->Pi+2*log(2),(2*k+1)*(4*n+4*u-1),n^2*(2*n+4*u+1)^2]
\end{verbatim}
Convergence type $P^-$ with $P=4k+2$.
\end{cf}

\smallskip

\begin{cf}\label{5.8.29}{\ }
\begin{verbatim}
[()->Pi+2*log(2),[4,8*n^2-8*n+5],[2,-n^2*(16*n^2-1)]]
\end{verbatim}
$$\pi+2\log(2)=4+\dfrac{2}{5-\dfrac{15}{21-\dfrac{252}{53-\dfrac{1287}{101-\dfrac{4080}{165-\dfrac{9975}{245-\ddots}}}}}}$$
Convergence type $P^+$ with $P=3/2$ and $C=\pi/(2^{3/2}\cdot3\G(3/4)^2)$, so
that
$$\pi+2\log(2)-\dfrac{p(n)}{q(n)}\sim\dfrac{\pi/(2^{3/2}\cdot3\G(3/4)^2)}{n^{3/2}}$$
$$A=1-(3/4)/n+(155/448)/n^2-(15/256)/n^3+\cdots$$
Series:
$$\pi+2\log(2)=4+\dfrac{2}{5}\sum_{n\ge0}\dfrac{(3/4)_n}{(n+1)(9/4)_n}$$
Parametric family for $k\ge0$ and $u$:
\begin{verbatim}
[()->Pi+2*log(2),8*n^2-8*n+5+4*u*(2*k+1)+2*k*(2*k+3),-n^2*(16*n^2-(4*u+1)^2)]
\end{verbatim}
Convergence type $P^+$ with $P=2k+2u+3/2$.
\end{cf}

\smallskip

\begin{cf}\label{5.8.30}{\ }
\begin{verbatim}
[()->Pi+2*log(2),[4,(4*n-1)*(5*n-1)],[6,-4*n*(2*n-1)*(2*n+1)*(4*n-1)]]
\end{verbatim}
$$\pi+2\log(2)=4+\dfrac{6}{12-\dfrac{36}{63-\dfrac{840}{154-\dfrac{4620}{285-\dfrac{15120}{456-\dfrac{37620}{667-\ddots}}}}}}$$
Convergence type $E$ with $E=4$, $P=9/4$, and $C=\G(1/4)/2^{5/2}$, so that
$$\pi+2\log(2)-\dfrac{p(n)}{q(n)}\sim\dfrac{\G(1/4)}{2^{2n+5/2}n^{9/4}}$$
$$A=1-(99/32)/n+(17817/2048)/n^2+\cdots$$
Series:
$$\pi+2\log(2)=4+\dfrac{1}{2}\sum_{n\ge0}\dfrac{(3/4)_nn!}{(3/2)_n(5/2)_n}2^{-2n}$$
Parametric family for $k\ge0$:
\begin{verbatim}
[()->Pi+2*log(2),20*n^2+(8*k-9)*n+1-3*k,-2*(2*n-k)*(2*n-k-1)*(2*n+1)*(4*n-1)]
\end{verbatim}
Convergence type $E$ with $E=4$ and $P=3k+9/4$.
\end{cf}

\smallskip

\begin{cf}\label{5.8.31}{\ }
\begin{verbatim}
[()->Pi-4*log(2),[0,10,8*n+3],[4,(n+1)^2*(2*n-1)*(2*n+3)]]
\end{verbatim}
$$\pi-4\log(2)=\dfrac{4}{10+\dfrac{20}{19+\dfrac{189}{27+\dfrac{720}{35+\dfrac{1925}{43+\dfrac{4212}{51+\ddots}}}}}}$$
Convergence type $P^-$ with $P=4$ and $C=3/4$, so that
$$\pi-4\log(2)-\dfrac{p(n)}{q(n)}\sim(-1)^n\dfrac{3/4}{n^4}$$
$$A=1-(9/2)/n+(45/4)/n^2-(75/4)/n^3+(371/16)/n^4+\cdots$$
Series:
$$\pi-4\log(2)=12\sum_{n\ge0}\dfrac{(-1)^n}{(n+2)(2n+1)(2n+3)(2n+5)}$$
\end{cf}

\smallskip

\begin{cf}\label{5.8.32}{\ }
\begin{verbatim}
[()->Pi+4*log(2),[20/3,8*n-3],[-4,n^2*(2*n-1)*(2*n+3)]]
\end{verbatim}
$$\pi+4\log(2)=20/3-\dfrac{4}{5+\dfrac{5}{13+\dfrac{84}{21+\dfrac{405}{29+\dfrac{1232}{37+\dfrac{2925}{45+\ddots}}}}}}$$
Convergence type $P^-$ with $P=4$ and $C=-3/4$, so that
$$\pi+4\log(2)-\dfrac{p(n)}{q(n)}\sim(-1)^{n+1}\dfrac{3/4}{n^4}$$
$$A=1-(7/2)/n+(25/4)/n^2-(25/4)/n^3+(91/16)/n^4+\cdots$$
Series:
$$\pi+4\log(2)=\dfrac{20}{3}-12\sum_{n\ge0}\dfrac{(-1)^n}{(n+1)(2n+1)(2n+3)(2n+5)}$$
\end{cf}

\smallskip

\begin{cf}\label{5.8.33}{\ }
\begin{verbatim}
[()->Pi-6*log(2),[0,3*(2*n-1)],[-4,n^2*(16*n^2-1)]]
\end{verbatim}
$$\pi-6\log(2)=-\dfrac{4}{3+\dfrac{15}{9+\dfrac{252}{15+\dfrac{1287}{21+\dfrac{4080}{27+\dfrac{9975}{33+\ddots}}}}}}$$
Convergence type $P^-$ with $P=3/2$ and $C=-\G(3/4)^2\sqrt{2}/\pi$,
so that
$$\pi-6\log(2)-\dfrac{p(n)}{q(n)}\sim(-1)^{n+1}\dfrac{\G(3/4)^2\sqrt{2}/\pi}{n^{3/2}}\;.$$
$$A=1-(3/4)/n-(1/64)/n^2+(147/256)/n^3+(21/8192)/n^4+\cdots$$
Series:
$$\pi-6\log(2)=\dfrac{4}{3}\sum_{n\ge0}(-1)^{n+1}\dfrac{(5/4)_n}{(n+1)(7/4)_n}$$
Parametric family for $k\ge0$ and $u\ge0$:
\begin{verbatim}
[()->Pi-6*log(2),(8*k+4*u+3)*(2*n-1),n^2*(16*n^2-(4*u-1)^2)]
\end{verbatim}
Convergence type $P^-$ with $P=4k+2u+3/2$.
\end{cf}

\smallskip

\begin{cf}\label{5.8.34}{\ }
\begin{verbatim}
[()->Pi-6*log(2),[0,8*n^2-10*n+5],[-2,-n^2*(4*n-1)^2]]
\end{verbatim}
$$\pi-6\log(2)=-\dfrac{2}{3-\dfrac{9}{17-\dfrac{196}{47-\dfrac{1089}{93-\dfrac{3600}{155-\dfrac{9025}{233-\ddots}}}}}}$$
Convergence type $P^+$ with $P=1$ and $C=-1/2$, so that
$$\pi-6\log(2)-\dfrac{p(n)}{q(n)}\sim-\dfrac{1/2}{n}\;.$$
$$A=1-(3/8)/n+(1/16)/n^2+(9/256)/n^3-(5/256)/n^4+\cdots$$
Series:
$$\pi-6\log(2)=-2\sum_{n\ge1}\dfrac{1}{n(4n-1)}$$
Parametric family for $k\ge0$ and $u\ge0$:
\begin{verbatim}
[()->Pi-6*log(2),8*n^2+(4*u-5)*(2*n-1)+4*k*(k+1),-n^2*(4*n+4*u-1)^2]
\end{verbatim}
Convergence type $P^+$ with $P=2k+1$.
\end{cf}

\smallskip

\begin{cf}\label{5.8.35}{\ }
\begin{verbatim}
[()->Pi-6*log(2),[0,8*n^2-7*n+3],[-2,-4*n^3*(4*n+1)]]
\end{verbatim}
$$\pi-6\log(2)=-\dfrac{2}{4-\dfrac{20}{21-\dfrac{288}{54-\dfrac{1404}{103-\dfrac{4352}{168-\dfrac{10500}{249-\ddots}}}}}}$$
Convergence type $P^+$ with $P=3/4$ and $C=-8/(3\G(1/4))$, so that
$$\pi-6\log(2)-\dfrac{p(n)}{q(n)}\sim-\dfrac{8/(3\G(1/4))}{n^{3/4}}\;.$$
$$A=1-(93/224)/n+(3235/22528)/n^2-(4629/327680)/n^3+\cdots$$
Series:
$$\pi-6\log(2)=-\dfrac{1}{2}\sum_{n\ge0}\dfrac{(5/4)_n}{(n+1)(n+1)!}$$
Parametric family for $k\ge0$ and $u\ge0$:
\begin{verbatim}
[()->Pi-6*log(2),8*n^2-(4*u+7)*n+4*k^2+(4*u+3)*(k+1),
                                       -4*n^3*(4*n-(4*u-1))]
\end{verbatim}
Convergence type $P^+$ with $P=2k+u+3/4$.
\end{cf}

\smallskip

\begin{cf}\label{5.8.36}{\ }
\begin{verbatim}
[()->Pi+6*log(2),[0,2*n-1],[12,n^2*(16*n^2-9)]]
\end{verbatim}
$$\pi+6\log(2)=\dfrac{12}{1+\dfrac{7}{3+\dfrac{220}{5+\dfrac{1215}{7+\dfrac{3952}{9+\dfrac{9775}{11+\ddots}}}}}}$$
Convergence type $P^-$ with $P=1/2$ and $C=\G(1/4)^2\sqrt{2}/\pi$,
so that
$$\pi+6\log(2)-\dfrac{p(n)}{q(n)}\sim(-1)^n\dfrac{\G(1/4)^2\sqrt{2}/\pi}{n^{1/2}}\;.$$
$$A=1-(1/4)/n+(1/64)/n^2+(15/256)/n^3-(19/8192)/n^4+\cdots$$
Series:
$$\pi+6\log(2)=12\sum_{n\ge0}(-1)^n\dfrac{(7/4)_n}{(n+1)(5/4)_n}$$
Parametric family for $k\ge0$ and $u\ge0$:
\begin{verbatim}
[()->Pi+6*log(2),(8*k+4*u+1)*(2*n-1),n^2*(16*n^2-(4*u-3)^2)]
\end{verbatim}
Convergence type $P^-$ with $P=4k+2u+1/2$.
\end{cf}

\smallskip

\begin{cf}\label{5.8.37}{\ }
\begin{verbatim}
[()->Pi+6*log(2),[0,8*n^2-14*n+7],[6,-n^2*(4*n-3)^2]]
\end{verbatim}
$$\pi+6\log(2)=\dfrac{6}{1-\dfrac{1}{11-\dfrac{100}{37-\dfrac{729}{79-\dfrac{2704}{137-\dfrac{7225}{211-\ddots}}}}}}$$
Convergence type $P^+$ with $P=1$ and $C=3/2$, so that
$$\pi+6\log(2)-\dfrac{p(n)}{q(n)}\sim\dfrac{3/2}{n}\;.$$
$$A=1-(9/8)/n+(59/48)/n^2-(333/256)/n^3+(1025/768)/n^4+\cdots$$
Series:
$$\pi+6\log(2)=6\sum_{n\ge1}\dfrac{1}{n(4n-3)}$$
Parametric family for $k\ge0$ and $u\ge0$:
\begin{verbatim}
[()->Pi+6*log(2),8*n^2+(4*u-7)*(2*n-1)+4*k*(k+1),-n^2*(4*n+4*u-3)^2]
\end{verbatim}
Convergence type $P^+$ with $P=2k+1$.
\end{cf}

\smallskip

\begin{cf}\label{5.8.38}{\ }
\begin{verbatim}
[()->Pi+6*log(2),[0,8*n^2-5*n+1],[6,-4*n^3*(4*n+3)]]
\end{verbatim}
$$\pi+6\log(2)=\dfrac{6}{4-\dfrac{28}{23-\dfrac{352}{58-\dfrac{1620}{109-\dfrac{4864}{176-\dfrac{11500}{259-\ddots}}}}}}$$
Convergence type $P^+$ with $P=1/4$ and $C=8/\G(3/4)$, so that
$$\pi+6\log(2)-\dfrac{p(n)}{q(n)}\sim\dfrac{8/\G(3/4)}{n^{1/4}}\;.$$
$$A=1-(23/160)/n+(721/18432)/n^2-(4909/851968)/n^3+\cdots$$
Series:
$$\pi+6\log(2)=\dfrac{3}{2}\sum_{n\ge0}\dfrac{(7/4)_n}{(n+1)(n+1)!}$$
Parametric family for $k\ge0$ and $u\ge0$:
\begin{verbatim}
[()->Pi+6*log(2),8*n^2-(4*u+5)*n+4*k^2+(4*u+1)*(k+1),
                                       -4*n^3*(4*n-(4*u-3))]
\end{verbatim}
Convergence type $P^+$ with $P=2k+u+1/4$.
\end{cf}

\smallskip

Note that there are infinitely many essentially distinct
$\Q$-linear combinations of $\pi$ and $\log(2)$ having a CF of a similar
type, the above are among the simplest.

\smallskip

\begin{cf}\label{5.8.39}{\ }
\begin{verbatim}
[()->Pi/sqrt(3)-log(2),[0,1],[3,9*n^2]]
\end{verbatim}
$$\dfrac{\pi}{\sqrt{3}}-\log(2)=\dfrac{3}{1+\dfrac{9}{1+\dfrac{36}{1+\dfrac{81}{1+\dfrac{144}{1+\dfrac{225}{1+\ddots}}}}}}$$
Convergence type $P^-$ with $P=1/3$ and $C=256\G(1/3)^2/3^{13/2}$, so that
$$\dfrac{\pi}{\sqrt{3}}-\log(2)-\dfrac{p(n)}{q(n)}\sim(-1)^n\dfrac{256\G(1/3)^2/3^{13/2}}{n^{1/3}}\;.$$
$$A=1-(1/6)/n-(7/324)/n^2+(133/1944)/n^3+(713/104976)/n^4+\cdots$$
Parametric family with $k\ge0$:
\begin{verbatim}
[()->Pi/sqrt(3)-log(2),6*k+1,9*n^2]
\end{verbatim}
Convergence type $P^-$ with $P=2k+1/3$.
\end{cf}

\smallskip

\begin{cf}\label{5.8.40}{\ }
\begin{verbatim}
[()->Pi/sqrt(3)-log(2),[0,2,3],[3,(3*n-1)^2]]
\end{verbatim}
$$\dfrac{\pi}{\sqrt{3}}-\log(2)=\dfrac{3}{2+\dfrac{4}{3+\dfrac{25}{3+\dfrac{64}{3+\dfrac{121}{3+\dfrac{196}{3+\ddots}}}}}}$$
Convergence type $P^-$ with $P=1$ and $C=1/2$, so that
$$\dfrac{\pi}{\sqrt{3}}-\log(2)-\dfrac{p(n)}{q(n)}\sim(-1)^n\dfrac{1/2}{n}\;.$$
$$A=1-(1/6)/n-(2/9)/n^2+(13/108)/n^3+(22/81)/n^4-(121/486)/n^5+\cdots$$
Series:
$$\dfrac{\pi}{\sqrt{3}}-\log(2)=3\sum_{n\ge1}\dfrac{(-1)^{n+1}}{3n-1}$$
Parametric family with $k\ge0$:
\begin{verbatim}
[()->Pi/sqrt(3)-log(2),6*k+3,(3*n-1)^2]
\end{verbatim}
Convergence type $P^-$ with $P=2k+1$.
\end{cf}

\smallskip

\begin{cf}\label{5.8.41}{\ }
\begin{verbatim}
[()->Pi/sqrt(3)-log(2),[0,9*n-5],[3,-6*n*(3*n-1)]]
\end{verbatim}
$$\dfrac{\pi}{\sqrt{3}}-\log(2)=\dfrac{3}{4-\dfrac{12}{13-\dfrac{60}{22-\dfrac{144}{31-\dfrac{264}{40-\dfrac{420}{49-\ddots}}}}}}$$
Convergence type $E$ with $E=2$, $P=2/3$, and $C=\G(2/3)$, so that
$$\dfrac{\pi}{\sqrt{3}}-\log(2)-\dfrac{p(n)}{q(n)}\sim\dfrac{\G(2/3)}{2^nn^{2/3}}$$
$$A=1-(11/9)/n+(80/27)/n^2-(24850/2187)/n^3+\cdots$$
Series:
$$\dfrac{\pi}{\sqrt{3}}-\log(2)=\dfrac{3}{4}\sum_{n\ge0}\dfrac{n!}{(5/3)_n}2^{-n}$$
Parametric family for $k\ge0$:
\begin{verbatim}
[()->Pi/sqrt(3)-log(2),9*n+3*k-5,-6*n*(3*n-1)]
\end{verbatim}
Convergence type $E$ with $E=2$ and $P=2k+2/3$.
\end{cf}
      
\smallskip

\begin{cf}\label{5.8.42}{\ }
\begin{verbatim}
[()->Pi/sqrt(3)-log(2),[0,3*n-1],[[3,4],[9*n^2,(3*n+2)^2]]]
\end{verbatim}
$$\dfrac{\pi}{\sqrt{3}}-\log(2)=\dfrac{3}{2+\dfrac{4}{5+\dfrac{9}{8+\dfrac{25}{11+\dfrac{36}{14+\dfrac{64}{17+\ddots}}}}}}$$
Convergence type $E$ with $E=-(1+\sqrt{2})^2$, $P=0$, and
$C=2\pi/(1+\sqrt{2})^{4/3}$, so that
$$\dfrac{\pi}{\sqrt{3}}+\log(2)-\dfrac{p(n)}{q(n)}\sim(-1)^n\dfrac{2\pi}{(1+\sqrt{2})^{2n+4/3}}\;.$$
$$A=1-(19d/72)/n+(19d/108+361/5184)/n^2+(-142705d/1119744-361/3888)/n^3+\cdots$$
\end{cf}

\smallskip

\begin{cf}\label{5.8.43}{\ }
\begin{verbatim}
[()->Pi/sqrt(3)-log(2),[3/2,3*n+2],[[-3,25],[9*n^2,(3*n+5)^2]]]
\end{verbatim}
$$\dfrac{\pi}{\sqrt{3}}-\log(2)=3/2-\dfrac{3}{5+\dfrac{25}{8+\dfrac{9}{11+\dfrac{64}{14+\dfrac{36}{17+\dfrac{121}{20+\ddots}}}}}}$$
Convergence type $E$ with $E=-(1+\sqrt{2})^2$, $P=0$, and
$C=-2\pi/(1+\sqrt{2})^{10/3}$, so that
$$\dfrac{\pi}{\sqrt{3}}+\log(2)-\dfrac{p(n)}{q(n)}\sim(-1)^{n+1}\dfrac{2\pi}{(1+\sqrt{2})^{2n+10/3}}\;.$$
$$A=1+(5d/72)/n+(-25d/216+25/5184)/n^2+(329975d/1119744-125/7776)/n^3+\cdots$$
\end{cf}

\smallskip

\begin{cf}\label{5.8.44}{\ }
\begin{verbatim}
[()->Pi/sqrt(3)+log(2),[3,5],[-3,9*n^2]]
\end{verbatim}
$$\dfrac{\pi}{\sqrt{3}}+\log(2)=3-\dfrac{3}{5+\dfrac{9}{5+\dfrac{36}{5+\dfrac{81}{5+\dfrac{144}{5+\dfrac{225}{5+\ddots}}}}}}$$
Convergence type $P^-$ with $P=5/3$ and $C=-\G(1/3)^2/(2^{5/3}\cdot 9)$,
so that
$$\dfrac{\pi}{\sqrt{3}}+\log(2)-\dfrac{p(n)}{q(n)}\sim(-1)^{n+1}\dfrac{\G(1/3)^2/(2^{5/3}\cdot 9)}{n^{5/3}}\;.$$
$$A=1-(5/6)/n+(25/324)/n^2+(1045/1944)/n^3+(677/104976)/n^4+\cdots$$
Parametric family with $k\ge0$:
\begin{verbatim}
[()->Pi/sqrt(3)+log(2),6*k+5,9*n^2]
\end{verbatim}
Convergence type $P^-$ with $P=2k+5/3$.
\end{cf}

\smallskip

\begin{cf}\label{5.8.45}{\ }
\begin{verbatim}
[()->Pi/sqrt(3)+log(2),[0,1,3],[3,(3*n-2)^2]]
\end{verbatim}
$$\dfrac{\pi}{\sqrt{3}}+\log(2)=\dfrac{3}{1+\dfrac{1}{3+\dfrac{16}{3+\dfrac{49}{3+\dfrac{100}{3+\dfrac{169}{3+\ddots}}}}}}$$
Convergence type $P^-$ with $P=1$ and $C=1/2$, so that
$$\dfrac{\pi}{\sqrt{3}}+\log(2)-\dfrac{p(n)}{q(n)}\sim(-1)^n\dfrac{1/2}{n}\;.$$
$$A=1+(1/6)/n-(2/9)/n^2-(13/108)/n^3+(22/81)/n^4+(121/486)/n^5+\cdots$$
Series:
$$\dfrac{\pi}{\sqrt{3}}+\log(2)=3\sum_{n\ge1}\dfrac{(-1)^{n+1}}{3n-2}$$
Parametric family with $k\ge0$:
\begin{verbatim}
[()->Pi/sqrt(3)+log(2),6*k+3,(3*n-2)^2]
\end{verbatim}
Convergence type $P^-$ with $P=2k+1$.
\end{cf}

\smallskip

\begin{cf}\label{5.8.46}{\ }
\begin{verbatim}
[()->Pi/sqrt(3)+log(2),[3,9*n-1],[-3,-6*n*(3*n+1)]]
\end{verbatim}
$$\dfrac{\pi}{\sqrt{3}}+\log(2)=3-\dfrac{3}{8-\dfrac{24}{17-\dfrac{84}{26-\dfrac{180}{35-\dfrac{312}{44-\dfrac{480}{53-\ddots}}}}}}$$
Convergence type $E$ with $E=2$, $P=4/3$, and $C=-\G(4/3)$, so that
$$\dfrac{\pi}{\sqrt{3}}+\log(2)-\dfrac{p(n)}{q(n)}\sim-\dfrac{\G(4/3)}{2^nn^{4/3}}$$
$$A=1-(26/9)/n+(847/81)/n^2-(110950/2187)/n^3+\cdots$$
Series:
$$\dfrac{\pi}{\sqrt{3}}+\log(2)=3-\dfrac{3}{8}\sum_{n\ge0}\dfrac{n!}{(7/3)_n}2^{-n}$$
Parametric family for $k\ge0$:
\begin{verbatim}
[()->Pi/sqrt(3)+log(2),9*n+3*k-1,-6*n*(3*n+1)]
\end{verbatim}
Convergence type $E$ with $E=2$ and $P=2k+4/3$.
\end{cf}
        
\smallskip

\begin{cf}\label{5.8.47}{\ }
\begin{verbatim}
[()->Pi/sqrt(3)+log(2),[0,3*n-2],[[3,1],[9*n^2,(3*n+1)^2]]]
\end{verbatim}
$$\dfrac{\pi}{\sqrt{3}}+\log(2)=\dfrac{3}{1+\dfrac{1}{4+\dfrac{9}{7+\dfrac{16}{10+\dfrac{36}{13+\dfrac{49}{16+\ddots}}}}}}$$
Convergence type $E$ with $E=-(1+\sqrt{2})^2$, $P=0$, and
$C=2\pi/(1+\sqrt{2})^{2/3}$, so that
$$\dfrac{\pi}{\sqrt{3}}+\log(2)-\dfrac{p(n)}{q(n)}\sim(-1)^n\dfrac{2\pi}{(1+\sqrt{2})^{2n+2/3}}\;.$$
$$A=1+(5d/72)/n+(-5d/216+25/5184)/n^2+(122615d/1119744-25/7776)/n^3+\cdots$$
\end{cf}

\begin{cf}\label{5.8.48}{\ }
\begin{verbatim}
[()->Pi/sqrt(3)+log(2),[3,3*n+1],[[-3,16],[9*n^2,(3*n+4)^2]]]
\end{verbatim}
$$\dfrac{\pi}{\sqrt{3}}+\log(2)=3-\dfrac{3}{4+\dfrac{16}{7+\dfrac{9}{10+\dfrac{49}{13+\dfrac{36}{16+\dfrac{100}{19+\ddots}}}}}}$$
Convergence type $E$ with $E=-(1+\sqrt{2})^2$, $P=0$, and
$C=-2\pi/(1+\sqrt{2})^{8/3}$, so that
$$\dfrac{\pi}{\sqrt{3}}+\log(2)-\dfrac{p(n)}{q(n)}\sim(-1)^{n+1}\dfrac{2\pi}{(1+\sqrt{2})^{2n+8/3}}\;.$$
$$A=1-(19d/72)/n+(19d/54+361/5184)/n^2+(-536689d/1119744-361/1944)/n^3+\cdots$$
\end{cf}

\smallskip

\begin{cf}\label{5.8.49}{\ }
\begin{verbatim}
[()->Pi/sqrt(3)-2*log(2),[-1,12*n^2-18*n+9],[5,-n^2*(6*n-7)*(6*n+1)]]
\end{verbatim}
$$\dfrac{\pi}{\sqrt{3}}-2\log(2)=-1+\dfrac{5}{3+\dfrac{7}{21-\dfrac{260}{63-\dfrac{1881}{129-\dfrac{6800}{219-\dfrac{17825}{333-\ddots}}}}}}$$
Convergence type $P^+$ with $P=5/3$ and $C=-3\pi^2/(2^{4/3}\G(1/3)^4)$, so that
$$\dfrac{\pi}{\sqrt{3}}-2\log(2)-\dfrac{p(n)}{q(n)}\sim-\dfrac{3\pi^2/(2^{4/3}\G(1/3)^4)}{n^{5/3}}$$
$$A=1+(25/1296)/n^2+(15/448)/n^3+\cdots$$
Series:
$$\dfrac{\pi}{\sqrt{3}}-2\log(2)=-1-5\sum_{n\ge0}\dfrac{(7/6)_n}{(n+1)(4n-1)(4n+3)(5/6)_n}$$
Parametric family for $k\ge0$:
\begin{verbatim}
[()->Pi/sqrt(3)-2*log(2),12*n^2-18*n+9+2*k*(3*k+5),-n^2*(6*n-7)*(6*n+1)]
\end{verbatim}
Convergence type $P^+$ with $P=2k+5/3$.
\end{cf}

\smallskip

\begin{cf}\label{5.8.50}{\ }
\begin{verbatim}
[()->Pi/sqrt(3)-2*log(2),[-1,4,9*n^2-3*n+2],
                       [8,(n+1)^2*(3*n-1)^2*(3*n-2)*(3*n+4)]]
\end{verbatim}
$$\dfrac{\pi}{\sqrt{3}}-2\log(2)=-1+\dfrac{8}{4+\dfrac{112}{32+\dfrac{9000}{74+\dfrac{93184}{134+\dfrac{484000}{212+\dfrac{1742832}{308+\ddots}}}}}}$$
Convergence type $P^-$ with $P=1$ and $C=1$, so that
$$\dfrac{\pi}{\sqrt{3}}-2\log(2)-\dfrac{p(n)}{q(n)}\sim\dfrac{(-1)^n}{n}$$
$$A=1-(5/6)/n+(8/9)/n^2-(115/108)/n^3+(92/81)/n^4+\cdots$$
Series:
$$\dfrac{\pi}{\sqrt{3}}-2\log(2)=-1+2\sum_{n\ge0}(-1)^n\dfrac{3n+4}{(n+2)(3n+2)}$$
\end{cf}

\smallskip

\begin{cf}\label{5.8.51}{\ }
\begin{verbatim}
[()->Pi/sqrt(3)-2*log(2),[0,6*n-4],[1,n^2*(3*n-1)^2]]
\end{verbatim}
$$\dfrac{\pi}{\sqrt{3}}-2\log(2)=\dfrac{1}{2+\dfrac{4}{8+\dfrac{100}{14+\dfrac{576}{20+\dfrac{1936}{26+\dfrac{4900}{32+\ddots}}}}}}$$
Convergence type $P^-$ with $P=2$ and $C=1/6$ so that
$$\dfrac{\pi}{\sqrt{3}}-2\log(2)-\dfrac{p(n)}{q(n)}\sim(-1)^n\dfrac{1/6}{n^2}$$
$$A=1-(2/3)/n-(7/18)/n^2+(22/27)/n^3+(61/81)/n^4+\cdots$$
Series:
$$\dfrac{\pi}{\sqrt{3}}-2\log(2)=\sum_{n\ge0}\dfrac{(-1)^n}{(n+1)(3n+2)}$$
Parametric family for $k\ge0$:
\begin{verbatim}
[()->Pi/sqrt(3)-2*log(2),(2*k+1)*(6*n+6*u-4),n^2*(3*n+6*u-1)^2]
\end{verbatim}
Convergence type $P^-$ with $P=4k+2$.
\end{cf}

\smallskip

\begin{cf}\label{5.8.52}{\ }
\begin{verbatim}
[()->Pi/sqrt(3)-2*log(2),[0,12*n^2-12*n+5],[1,-n^2*(36*n^2-1)]]
\end{verbatim}
$$\dfrac{\pi}{\sqrt{3}}-2\log(2)=\dfrac{1}{5-\dfrac{35}{29-\dfrac{572}{77-\dfrac{2907}{149-\dfrac{9200}{245-\dfrac{22475}{365-\ddots}}}}}}$$
Convergence type $P^+$ with $P=2/3$ and $C=2^{2/3}\pi^2/\G(1/3)^4$, so that
$$\dfrac{\pi}{\sqrt{3}}-2\log(2)-\dfrac{p(n)}{q(n)}\sim\dfrac{2^{2/3}\pi^2/\G(1/3)^4}{n^{2/3}}$$
$$A=1-(1/3)/n+(115/1296)/n^2+(5/972)/n^3+\cdots$$
Series:
$$\dfrac{\pi}{\sqrt{3}}-2\log(2)=\dfrac{1}{5}\sum_{n\ge0}\dfrac{(7/6)_n}{(n+1)(11/6)_n}$$
Parametric family for $k\ge0$:
\begin{verbatim}
[()->Pi/sqrt(3)-2*log(2),12*n^2-12*n+2*k*(3*k+2)+5,-36*n^4+n^2]
\end{verbatim}
Convergence type $P^+$ with $P=2k+2/3$.
\end{cf}

\smallskip

\begin{cf}\label{5.8.53}{\ }
\begin{verbatim}
[()->Pi/sqrt(3)+2*log(2),[0,12*n^2-18*n+7],[1,-n^2*(6*n-1)*(6*n-5)]]
\end{verbatim}
$$\dfrac{\pi}{\sqrt{3}}+2\log(2)=\dfrac{1}{1-\dfrac{5}{19-\dfrac{308}{61-\dfrac{1989}{127-\dfrac{6992}{217-\dfrac{18125}{331-\ddots}}}}}}$$
Convergence type $P^+$ with $P=1/3$ and $C=2^{4/3}\pi^2/(3\G(2/3)^4)$, so that
$$\dfrac{\pi}{\sqrt{3}}+2\log(2)-\dfrac{p(n)}{q(n)}\sim\dfrac{2^{4/3}\pi^2/(3\G(2/3)^4)}{n^{1/3}}$$
$$A=1-(1/8)/n+(5/324)/n^2+(7/1620)/n^3+\cdots$$
Series:
$$\dfrac{\pi}{\sqrt{3}}+2\log(2)=\sum_{n\ge0}\dfrac{(5/6)_n}{(n+1)(7/6)_n}$$
Parametric family for $k\ge0$:
\begin{verbatim}
[()->Pi/sqrt(3)+2*log(2),12*n^2-18*n+7+2*k*(3*k+1),-n^2*(6*n-1)*(6*n-5)]
\end{verbatim}
Convergence type $P^+$ with $P=2k+1/3$.
\end{cf}

\smallskip

\begin{cf}\label{5.8.54}{\ }
\begin{verbatim}
[()->Pi/sqrt(3)+2*log(2),[0,9*n^2-15*n+5],[-4,n^2*(3*n-2)^2*(3*n-4)*(3*n+2)]]
\end{verbatim}
$$\dfrac{\pi}{\sqrt{3}}+2\log(2)=-\dfrac{4}{-1-\dfrac{5}{11+\dfrac{1024}{41+\dfrac{24255}{89+\dfrac{179200}{155+\dfrac{790075}{239+\ddots}}}}}}$$
Convergence type $P^-$ with $P=1$ and $C=1$, so that
$$\dfrac{\pi}{\sqrt{3}}+2\log(2)-\dfrac{p(n)}{q(n)}\sim\dfrac{(-1)^n}{n}$$
$$A=1-(1/6)/n-(1/9)/n^2+(7/108)/n^3+(11/81)/n^4+\cdots$$
Series:
$$\dfrac{\pi}{\sqrt{3}}+2\log(2)=2\sum_{n\ge0}(-1)^n\dfrac{3n+2}{(n+1)(3n+1)}$$
\end{cf}

\smallskip

\begin{cf}\label{5.8.55}{\ }
\begin{verbatim}
[()->Pi/sqrt(3)+2*log(2),[3,6*n-2],[1,n^2*(3*n+1)^2]]
\end{verbatim}
$$\dfrac{\pi}{\sqrt{3}}+2\log(2)=3+\dfrac{1}{4+\dfrac{16}{10+\dfrac{196}{16+\dfrac{900}{22+\dfrac{2704}{28+\dfrac{6400}{34+\ddots}}}}}}$$

Convergence type $P^-$ with $P=2$ and $C=1/6$ so that
$$\dfrac{\pi}{\sqrt{3}}+2\log(2)-\dfrac{p(n)}{q(n)}\sim(-1)^n\dfrac{1/6}{n^2}$$
$$A=1-(4/3)/n+(11/18)/n^2+(20/27)/n^3-(59/81)/n^4+\cdots$$
Series:
$$\dfrac{\pi}{\sqrt{3}}+2\log(2)=3+\sum_{n\ge0}\dfrac{(-1)^n}{(n+1)(3n+4)}$$
Parametric family for $k\ge0$:
\begin{verbatim}
[()->Pi/sqrt(3)+2*log(2),(2*k+1)*(6*n+6*u-2),n^2*(3*n+6*u+1)^2]
\end{verbatim}
Convergence type $P^-$ with $P=4k+2$.
\end{cf}

\smallskip

\begin{cf}\label{5.8.56}{\ }
\begin{verbatim}
[()->Pi/sqrt(3)+2*log(2),[3,12*n^2-12*n+7],[1,-n^2*(36*n^2-1)]]
\end{verbatim}
$$\dfrac{\pi}{\sqrt{3}}+2\log(2)=3+\dfrac{1}{7-\dfrac{35}{31-\dfrac{572}{79-\dfrac{2907}{151-\dfrac{9200}{247-\dfrac{22475}{367-\ddots}}}}}}$$
Convergence type $P^+$ with $P=4/3$ and $C=\G(1/3)^4/(2^{17/3}\pi^2)$, so that
$$\dfrac{\pi}{\sqrt{3}}+2\log(2)-\dfrac{p(n)}{q(n)}\sim\dfrac{\G(1/3)^4/(2^{17/3}\pi^2)}{n^{4/3}}$$
$$A=1-(2/3)/n+(112/405)/n^2-(7/243)/n^3+\cdots$$
Series:
$$\dfrac{\pi}{\sqrt{3}}+2\log(2)=3+\dfrac{1}{7}\sum_{n\ge0}\dfrac{(5/6)_n}{(n+1)(13/6)_n}$$
Parametric family for $k\ge0$ and $u$:
\begin{verbatim}
[()->Pi/sqrt(3)+2*log(2),12*n^2-12*n+7+6*u*(2*k+1)+2*k*(3*k+4),
                           -n^2*(36*n^2-(6*u+1)^2)]
\end{verbatim}
Convergence type $P^+$ with $P=2u+2k+4/3$.
\end{cf}

\smallskip

\begin{cf}\label{5.8.57}{\ }
\begin{verbatim}
[()->Pi/sqrt(3)-log(3)/3,[68/45,(2*n-1)*(9*n^2-9*n+29/2)],
                         [-8/9,-n^2*(9*n^2-1)*(9*n^2-4)]]
\end{verbatim}
$$\dfrac{\pi}{\sqrt{3}}-\dfrac{\log(3)}{3}=68/45-\dfrac{8/9}{29/2-\dfrac{40}{195/2-\dfrac{4480}{685/2-\dfrac{55440}{1715/2-\dfrac{320320}{3501/2-\dfrac{1237600}{6259/2-\ddots}}}}}}$$
Convergence type $P^+$ with $P=10/3$ and $C=-(5/3)\G(2/3)^2/2^{25/3}$, so that
$$\dfrac{\pi}{\sqrt{3}}-\dfrac{\log(3)}{3}-\dfrac{p(n)}{q(n)}\sim-\dfrac{(5/3)\G(2/3)^2/2^{25/3}}{n^{10/3}}\;.$$
$$A=1-(5/3)/n+(15655/10368)/n^2-(3175/3888)/n^3+(1930451/6718464)/n^4+\cdots$$
Parametric family for $k\ge0$:
\begin{verbatim}
[()->Pi/sqrt(3)-log(3)/3,(2*n-1)*(9*n^2-9*n+29/2+6*k*(3*k+5)),
                            -n^2*(9*n^2-1)*(9*n^2-4)]
\end{verbatim}
Convergence type $P^+$ with $P=4k+10/3$.
\end{cf}

\smallskip

\begin{cf}\label{5.8.58}{\ }
\begin{verbatim}
[()->Pi/sqrt(3)+log(3)/3,[4/3,(2*n-1)*(9*n^2-9*n+5/2)],
                      [8/9,-n^2*(9*n^2-1)*(9*n^2-4)]]
\end{verbatim}
$$\dfrac{\pi}{\sqrt{3}}+\dfrac{\log(3)}{3}=4/3+\dfrac{8/9}{5/2-\dfrac{40}{123/2-\dfrac{4480}{565/2-\dfrac{55440}{1547/2-\dfrac{320320}{3285/2-\dfrac{1237600}{5995/2-\ddots}}}}}}$$
Convergence type $P^+$ with $P=2/3$ and $C=(1/3)\G(1/3)^2/2^{5/3}$, so that
$$\dfrac{\pi}{\sqrt{3}}+\dfrac{\log(3)}{3}-\dfrac{p(n)}{q(n)}\sim\dfrac{(1/3)\G(1/3)^2/2^{5/3}}{n^{2/3}}\;.$$
$$A=1-(1/3)/n+(565/5184)/n^2-(85/3888)/n^3-(10057/6718464)/n^4+\cdots$$
Parametric family for $k\ge0$:
\begin{verbatim}
[()->Pi/sqrt(3)+log(3)/3,(2*n-1)*(9*n^2-9*n+5/2+6*k*(3*k+1)),
                            -n^2*(9*n^2-1)*(9*n^2-4)]
\end{verbatim}
Convergence type $P^+$ with $P=4k+2/3$.
\end{cf}

\smallskip

Both the above CFs can be Ap\'ery accelerated, formulas too complicated to
give here.

\smallskip

\begin{cf}\label{5.8.59}{\ }
\begin{verbatim}
[()->Pi/sqrt(3)-log(3),[1,(2*n-1)*(9*n^2-9*n+10)],
                       [-8/3,-n^2*(9*n^2-1)*(9*n^2-4)]]
\end{verbatim}
$$\dfrac{\pi}{\sqrt{3}}-\log(3)=1-\dfrac{8/3}{10-\dfrac{40}{84-\dfrac{4480}{320-\dfrac{55440}{826-\dfrac{320320}{1710-\dfrac{1237600}{3080-\ddots}}}}}}$$
Convergence type $P^+$ with $P=8/3$ and $C=-\G(2/3)^4/(6\pi^2)$, so that
$$\dfrac{\pi}{\sqrt{3}}-\log(3)-\dfrac{p(n)}{q(n)}\sim-\dfrac{\G(2/3)^4/(6\pi^2)}{n^{8/3}}\;.$$
$$A=1-(4/3)/n+(598/567)/n^2-(136/243)/n^3+(1613/6561)/n^4+\cdots$$
Series:
$$\dfrac{\pi}{\sqrt{3}}-\log(3)=1-\dfrac{4}{3}\sum_{n\ge0}\dfrac{(3n+1)(1/3)_n^2}{(n+1)(3n+5)(5/3)_n^2}$$
Parametric family for $k\ge0$:
\begin{verbatim}
[()->Pi/sqrt(3)-log(3),(2*n-1)*(9*n^2-9*n+10+6*k*(3*k+4)),
                       -n^2*(9*n^2-1)*(9*n^2-4)]
\end{verbatim}
Convergence type $P^+$ with $P=4k+8/3$.
\end{cf}

\smallskip

\begin{cf}\label{5.8.60}{\ }
\begin{verbatim}
[()->Pi/sqrt(3)-log(3),[0,7*n-2],[4,8*n^2]]
\end{verbatim}
$$\dfrac{\pi}{\sqrt{3}}-\log(3)=\dfrac{4}{5+\dfrac{8}{12+\dfrac{32}{19+\dfrac{72}{26+\dfrac{128}{33+\dfrac{200}{40+\ddots}}}}}}$$
Convergence type $E$ with $E=-8$, $P=1/3$, and $C=3^{-2/3}\G(2/3)^2$, so that
$$\dfrac{\pi}{\sqrt{3}}-\log(3)-\dfrac{p(n)}{q(n)}\sim(-1)^n\dfrac{3^{-2/3}\G(2/3)^2}{2^{3n}n^{1/3}}$$
$$A=1-(31/81)/n+(1379/6561)/n^2-(94451/1594323)/n^3+\cdots$$
Parametric family for $k\ge0$:
\begin{verbatim}
[()->Pi/sqrt(3)-log(3),7*n+9*k-2,8*n^2]
\end{verbatim}
Convergence type $E$ with $E=-8$ and $P=2k+1/3$.
\end{cf}

\smallskip

\begin{cf}\label{5.8.61}{\ }
\begin{verbatim}
[()->Pi/sqrt(3)-log(3),[0,10*n-4],[4,-3*n*(3*n-1)]]
\end{verbatim}
$$\dfrac{\pi}{\sqrt{3}}-\log(3)=\dfrac{4}{6-\dfrac{6}{16-\dfrac{30}{26-\dfrac{72}{36-\dfrac{132}{46-\dfrac{210}{56-\ddots}}}}}}$$
Convergence type $E$ with $E=9$, $P=2/3$, and $C=\G(2/3)/2$, so that
$$\dfrac{\pi}{\sqrt{3}}-\log(3)-\dfrac{p(n)}{q(n)}\sim\dfrac{\G(2/3)/2}{3^{2n}n^{2/3}}$$
$$A=1-(23/36)/n+(55/96)/n^2-(47705/69984)/n^3+\cdots$$
Series:
$$\dfrac{\pi}{\sqrt{3}}-\log(3)=\dfrac{2}{3}\sum_{n\ge0}\dfrac{n!}{(5/3)_n}3^{-2n}$$
Parametric family for $k\ge0$:
\begin{verbatim}
[()->Pi/sqrt(3)-log(3),10*n+8*k-4,-3*n*(3*n-1)]
\end{verbatim}
Convergence type $E$ with $E=9$ and $P=2k+2/3$.
\end{cf}
        
\smallskip

\begin{cf}\label{5.8.62}{\ }
\begin{verbatim}
[()->Pi/sqrt(3)-log(3),[1,3672*n^4-6426*n^3+1980*n^2+650*n-136],
                      [74,-n^2*(3*n-1)*(3*n+2)
                              *(36*n^2+75*n+37)*(36*n^2-69*n+31)]]
\end{verbatim}
$$\dfrac{\pi}{\sqrt{3}}-\log(3)=1+\dfrac{74}{-260+\dfrac{2960}{16428-\dfrac{1959520}{143564-\dfrac{68688576}{562912-\ddots}}}}$$
Convergence type $E$ with $E=(1+\sqrt{2})^8$, $P=0$, and $C=-4\pi/(1+\sqrt{2})^{14/3}$, so that
$$\dfrac{\pi}{\sqrt{3}}-\log(3)-\dfrac{p(n)}{q(n)}\sim\dfrac{4\pi}{(1+\sqrt{2})^{8n+14/3}}\;.$$
$$A=1+(11d/96)/n+(-77d/1152+121/9216)/n^2+(109105d/2654208-847/55296)/n^3+\cdots$$
\end{cf}

\smallskip

\begin{cf}\label{5.8.63}{\ }
\begin{verbatim}
[()->Pi/sqrt(3)+log(3),[2,(2*n-1)*(9*n^2-9*n+4)],
                      [8/3,-n^2*(9*n^2-1)*(9*n^2-4)]]
\end{verbatim}
$$\dfrac{\pi}{\sqrt{3}}+\log(3)=2+\dfrac{8/3}{4-\dfrac{40}{66-\dfrac{4480}{290-\dfrac{55440}{784-\dfrac{320320}{1656-\dfrac{1237600}{3014-\ddots}}}}}}$$
Convergence type $P^+$ with $P=4/3$ and $C=\G(1/3)^4/(12\pi^2)$, so that
$$\dfrac{\pi}{\sqrt{3}}+\log(3)-\dfrac{p(n)}{q(n)}\sim\dfrac{\G(1/3)^4/(12\pi^2)}{n^{4/3}}\;.$$
$$A=1-(2/3)/n+(127/405)/n^2-(22/243)/n^3+(38/6561)/n^4-(148/19683)/n^5+\cdots$$
Series:
$$\dfrac{\pi}{\sqrt{3}}+\log(3)=2+\dfrac{4}{3}\sum_{n\ge0}\dfrac{(3n+2)(2/3)_n^2}{(n+1)(3n+4)(4/3)_n^2}$$
Parametric family for $k\ge0$:
\begin{verbatim}
[()->Pi/sqrt(3)+log(3),(2*n-1)*(9*n^2-9*n+4+6*k*(3*k+2)),
                       -n^2*(9*n^2-1)*(9*n^2-4)]
\end{verbatim}
Convergence type $P^+$ with $P=4k+4/3$.
\end{cf}

\smallskip

\begin{cf}\label{5.8.64}{\ }
\begin{verbatim}
[()->Pi/sqrt(3)+log(3),[0,7*n-5],[8,8*n^2]]
\end{verbatim}
$$\dfrac{\pi}{\sqrt{3}}+\log(3)=\dfrac{8}{2+\dfrac{8}{9+\dfrac{32}{16+\dfrac{72}{23+\dfrac{128}{30+\dfrac{200}{37+\ddots}}}}}}$$
Convergence type $E$ with $E=-8$, $P=-1/3$, and $C=3^{2/3}\G(1/3)^2/2$, so that
$$\dfrac{\pi}{\sqrt{3}}+\log(3)-\dfrac{p(n)}{q(n)}\sim(-1)^n\dfrac{3^{2/3}\G(1/3)^2}{2^{3n+1}n^{-1/3}}$$
$$A=1-(4/81)/n+(527/6561)/n^2+(1390/1594323)/n^3+\cdots$$
Parametric family for $k\ge0$:
\begin{verbatim}
[()->Pi/sqrt(3)+log(3),7*n+9*k-5,8*n^2]
\end{verbatim}
Convergence type $E$ with $E=-8$ and $P=2k-1/3$.
\end{cf}

\smallskip

\begin{cf}\label{5.8.65}{\ }
\begin{verbatim}
[()->Pi/sqrt(3)+log(3),[3,10*n+2],[-1,-3*n*(3*n+1)]]
\end{verbatim}
$$\dfrac{\pi}{\sqrt{3}}+\log(3)=3-\dfrac{1}{12-\dfrac{12}{22-\dfrac{42}{32-\dfrac{90}{42-\dfrac{156}{52-\dfrac{240}{62-\ddots}}}}}}$$
Convergence type $E$ with $E=9$, $P=4/3$, and $C=-\G(1/3)/24$, so that
$$\dfrac{\pi}{\sqrt{3}}+\log(3)-\dfrac{p(n)}{q(n)}\sim-\dfrac{\G(1/3)/8}{3^{2n+1}n^{4/3}}$$
$$A=1-(31/18)/n+(3703/1296)/n^2-(721735/139968)/n^3+\cdots$$
Series:
$$\dfrac{\pi}{\sqrt{3}}+\log(3)=3-\dfrac{1}{12}\sum_{n\ge0}\dfrac{n!}{(7/3)_n}3^{-2n}$$
Parametric family for $k\ge0$:
\begin{verbatim}
[()->Pi/sqrt(3)+log(3),10*n+8*k+2,-3*n*(3*n+1)]
\end{verbatim}
Convergence type $E$ with $E=9$ and $P=2k+4/3$.
\end{cf}
        
\smallskip

\begin{cf}\label{5.8.66}{\ }
\begin{verbatim}
[()->Pi/sqrt(3)+log(3),[2,3672*n^4-8262*n^3+4734*n^2-20*n-260],
                      [-124,-n^2*(3*n-2)*(3*n+1)
                                *(36*n^2-75*n+37)*(36*n^2+69*n+31)]]
\end{verbatim}
$$\dfrac{\pi}{\sqrt{3}}+\log(3)=2-\dfrac{124}{-136+\dfrac{1088}{11292-\dfrac{1086736}{116644-\dfrac{48152160}{486668-\ddots}}}}$$
Convergence type $E$ with $E=(1+\sqrt{2})^8$, $P=0$, and $C=4\pi/(1+\sqrt{2})^{10/3}$, so that
$$\dfrac{\pi}{\sqrt{3}}+\log(3)-\dfrac{p(n)}{q(n)}\sim\dfrac{4\pi}{(1+\sqrt{2})^{8n+10/3}}\;.$$
$$A=1+(11d/96)/n+(-55d/1152+121/9216)/n^2+(58417d/2654208-605/55296)/n^3+\cdots$$
\end{cf}

\smallskip

\begin{cf}\label{5.8.67}{\ }
\begin{verbatim}
[()->Pi/sqrt(3)-3*log(3),[0,2*(2*n-1)],[-4,n^2*(9*n^2-1)]]
\end{verbatim}
$$\dfrac{\pi}{\sqrt{3}}-3\log(3)=-\dfrac{4}{2+\dfrac{8}{6+\dfrac{140}{10+\dfrac{720}{14+\dfrac{2288}{18+\dfrac{5600}{22+\ddots}}}}}}$$
Convergence type $P^-$ with $P=4/3$ and $C=-\G(2/3)^2\sqrt{3}/\pi$,
so that
$$\dfrac{\pi}{\sqrt{3}}-3\log(3)-\dfrac{p(n)}{q(n)}\sim(-1)^{n+1}\dfrac{\G(2/3)^2\sqrt{3}/\pi}{n^{4/3}}\;.$$
$$A=1-(2/3)/n-(1/81)/n^2+(110/243)/n^3+(14/6561)/n^4+\cdots$$
Series:
$$\dfrac{\pi}{\sqrt{3}}-3\log(3)=2\sum_{n\ge0}(-1)^{n+1}\dfrac{(4/3)_n}{(n+1)(5/3)_n}$$
Parametric family for $k\ge0$ and $u\ge0$:
\begin{verbatim}
[()->Pi/sqrt(3)-3*log(3),(6*k+3*u+2)*(2*n-1),n^2*(9*n^2-(3*u-1)^2)]
\end{verbatim}
Convergence type $P^-$ with $P=4k+2u+4/3$.
\end{cf}

\smallskip

\begin{cf}\label{5.8.68}{\ }
\begin{verbatim}
[()->Pi/sqrt(3)-3*log(3),[0,6*n^2-8*n+4],[-2,-n^2*(3*n-1)^2]]
\end{verbatim}
$$\dfrac{\pi}{\sqrt{3}}-3\log(3)=-\dfrac{2}{2-\dfrac{4}{12-\dfrac{100}{34-\dfrac{576}{68-\dfrac{1936}{114-\dfrac{4900}{172-\ddots}}}}}}$$
Convergence type $P^+$ with $P=1$ and $C=-2/3$, so that
$$\dfrac{\pi}{\sqrt{3}}-3\log(3)-\dfrac{p(n)}{q(n)}\sim-\dfrac{2/3}{n}\;.$$
$$A=1-(1/3)/n+(1/27)/n^2+(1/27)/n^3-(1/81)/n^4+\cdots$$
Series:
$$\dfrac{\pi}{\sqrt{3}}-3\log(3)=-2\sum_{n\ge1}\dfrac{1}{n(3n-1)}$$
Parametric family for $k\ge0$ and $u\ge0$:
\begin{verbatim}
[()->Pi/sqrt(3)-3*log(3),6*n^2+(3*u-4)*(2*n-1)+3*k*(k+1),
                                              -n^2*(3*n+3*u-1)^2]
\end{verbatim}
Convergence type $P^+$ with $P=2k+1$.
\end{cf}

\smallskip

\begin{cf}\label{5.8.69}{\ }
\begin{verbatim}
[()->Pi/sqrt(3)-3*log(3),[0,6*n^2-5*n+2],[-2,-3*n^3*(3*n+1)]]
\end{verbatim}
$$\dfrac{\pi}{\sqrt{3}}-3\log(3)=-\dfrac{2}{3-\dfrac{12}{16-\dfrac{168}{41-\dfrac{810}{78-\dfrac{2496}{127-\dfrac{6000}{188-\ddots}}}}}}$$
Convergence type $P^+$ with $P=2/3$ and $C=-3/\G(1/3)$, so that
$$\dfrac{\pi}{\sqrt{3}}-3\log(3)-\dfrac{p(n)}{q(n)}\sim-\dfrac{3/\G(1/3)}{n^{2/3}}\;.$$
$$A=1-(17/45)/n+(7/54)/n^2-(376/24057)/n^3+\cdots$$
Series:
$$\dfrac{\pi}{\sqrt{3}}-3\log(3)=-\dfrac{2}{3}\sum_{n\ge0}\dfrac{(4/3)_n}{(n+1)(n+1)!}$$
Parametric family for $k\ge0$ and $u\ge0$:
\begin{verbatim}
[()->Pi/sqrt(3)-3*log(3),6*n^2-(3*u+5)*n+3*k^2+(3*u+2)*(k+1),
                                        -3*n^3*(3*n-(3*u-1))]
\end{verbatim}
Convergence type $P^+$ with $P=2k+u+2/3$.
\end{cf}

\smallskip

\begin{cf}\label{5.8.70}{\ }
\begin{verbatim}
[()->Pi/sqrt(3)-3*log(3),[-31/2,30,7*n^3+25*n^2+25*n+5],[80,-2*n*(n+2)^3*(2*n+3)*(3*n+8)]]
\end{verbatim}
$$\dfrac{\pi}{\sqrt{3}}-3\log(3)=-31/2+\dfrac{80}{30-\dfrac{2970}{211-\dfrac{25088}{494-\dfrac{114750}{953-\dfrac{380160}{1630-\dfrac{1025570}{2567-\ddots}}}}}}$$
Convergence type $E$ with $E=4/3$, $P=-1/6$, and $C=81\sqrt{\pi}/(8\G(2/3))$,
so that
$$\dfrac{\pi}{\sqrt{3}}-3\log(3)-\dfrac{p(n)}{q(n)}\sim\dfrac{81\sqrt{\pi}/(8\G(2/3))}{(4/3)^nn^{-1/6}}$$
$$A=1+(1/72)/n+(16481/10368)/n^2+\cdots$$
Series:
$$\dfrac{\pi}{\sqrt{3}}-3\log(3)=-18+6\sum_{n\ge1}\dfrac{(n-1)(2/3)_n}{n(1/2)_n}(3/4)^n$$
\end{cf}

\smallskip

\begin{cf}\label{5.8.71}{\ }
\begin{verbatim}
[()->Pi/sqrt(3)+3*log(3),[0,2*n-1],[8,n^2*(9*n^2-4)]]
\end{verbatim}
$$\dfrac{\pi}{\sqrt{3}}+3\log(3)=\dfrac{8}{1+\dfrac{5}{3+\dfrac{128}{5+\dfrac{693}{7+\dfrac{2240}{9+\dfrac{5525}{11+\ddots}}}}}}$$
Convergence type $P^-$ with $P=2/3$ and $C=\G(1/3)^2\sqrt{3}/\pi$,
so that
$$\dfrac{\pi}{\sqrt{3}}+3\log(3)-\dfrac{p(n)}{q(n)}\sim(-1)^n\dfrac{\G(1/3)^2\sqrt{3}/\pi}{n^{2/3}}\;.$$
$$A=1-(1/3)/n+(1/81)/n^2+(26/243)/n^3-(13/6561)/n^4+\cdots$$
Series:
$$\dfrac{\pi}{\sqrt{3}}+3\log(3)=8\sum_{n\ge0}(-1)^n\dfrac{(5/3)_n}{(n+1)(4/3)_n}$$
Parametric family for $k\ge0$ and $u\ge0$:
\begin{verbatim}
[()->Pi/sqrt(3)+3*log(3),(6*k+3*u+1)*(2*n-1),n^2*(9*n^2-(3*u-2)^2)]
\end{verbatim}
Convergence type $P^-$ with $P=4k+2u+2/3$.
\end{cf}

\smallskip

\begin{cf}\label{5.8.72}{\ }
\begin{verbatim}
[()->Pi/sqrt(3)+3*log(3),[0,6*n^2-10*n+5],[4,-n^2*(3*n-2)^2]]
\end{verbatim}
$$\dfrac{\pi}{\sqrt{3}}+3\log(3)=\dfrac{4}{1-\dfrac{1}{9-\dfrac{64}{29-\dfrac{441}{61-\dfrac{1600}{105-\dfrac{4225}{161-\ddots}}}}}}$$
Convergence type $P^+$ with $P=1$ and $C=4/3$, so that
$$\dfrac{\pi}{\sqrt{3}}+3\log(3)-\dfrac{p(n)}{q(n)}\sim\dfrac{4/3}{n}\;.$$
$$A=1-(1/6)/n-(1/54)/n^2+(1/54)/n^3+(1/162)/n^4+\cdots$$
Series:
$$\dfrac{\pi}{\sqrt{3}}+3\log(3)=4\sum_{n\ge1}\dfrac{1}{n(3n-2)}$$
Parametric family for $k\ge0$ and $u\ge0$:
\begin{verbatim}
[()->Pi/sqrt(3)+3*log(3),6*n^2+(3*u-5)*(2*n-1)+3*k*(k+1),
                                       -n^2*(3*n+3*u-2)^2]
\end{verbatim}
Convergence type $P^+$ with $P=2k+1$.
\end{cf}

\smallskip

\begin{cf}\label{5.8.73}{\ }
\begin{verbatim}
[()->Pi/sqrt(3)+3*log(3),[0,6*n^2-4*n+1],[4,-3*n^3*(3*n+2)]]
\end{verbatim}
$$\dfrac{\pi}{\sqrt{3}}+3\log(3)=\dfrac{4}{3-\dfrac{15}{17-\dfrac{192}{43-\dfrac{891}{81-\dfrac{2688}{131-\dfrac{6375}{193-\ddots}}}}}}$$
Convergence type $P^+$ with $P=1/3$ and $C=6/\G(2/3)$, so that
$$\dfrac{\pi}{\sqrt{3}}+3\log(3)-\dfrac{p(n)}{q(n)}\sim\dfrac{6/\G(2/3)}{n^{1/3}}\;.$$
$$A=1-(7/36)/n+(65/1134)/n^2-(391/43740)/n^3+\cdots$$
Series:
$$\dfrac{\pi}{\sqrt{3}}+3\log(3)=\dfrac{4}{3}\sum_{n\ge0}\dfrac{(5/3)_n}{(n+1)(n+1)!}$$
Parametric family for $k\ge0$ and $u\ge0$:
\begin{verbatim}
[()->Pi/sqrt(3)+3*log(3),6*n^2-(3*u+4)*n+3*k^2+(3*u+1)*(k+1),
                                        -3*n^3*(3*n-(3*u-2))]
\end{verbatim}
Convergence type $P^+$ with $P=2k+u+1/3$.
\end{cf}

\smallskip

\begin{cf}\label{5.8.74}{\ }
\begin{verbatim}
[()->Pi/sqrt(3)+3*log(3),[0,7*n^2-7*n+2],[6,-2*n^2*(2*n-1)*(3*n+1)]]
\end{verbatim}
$$\dfrac{\pi}{\sqrt{3}}+3\log(3)=\dfrac{6}{2-\dfrac{8}{16-\dfrac{168}{44-\dfrac{900}{86-\dfrac{2912}{142-\dfrac{7200}{212-\ddots}}}}}}$$
Convergence type $E$ with $E=4/3$, $P=7/6$, and $C=18\sqrt{\pi}/\G(1/3)$, so that
$$\dfrac{\pi}{\sqrt{3}}+3\log(3)-\dfrac{p(n)}{q(n)}\sim\dfrac{18\sqrt{\pi}/\G(1/3)}{(4/3)^nn^{7/6}}$$
$$A=1-(335/72)/n+(121867/3456)/n^2+\cdots$$
Series:
$$\dfrac{\pi}{\sqrt{3}}+3\log(3)=6\sum_{n\ge1}\dfrac{(1/3)_n}{n(1/2)_n}(3/4)^n$$
\end{cf}

\smallskip

\begin{cf}\label{5.8.75}{\ }
\begin{verbatim}
[()->Pi/sqrt(3)-log(432)/3,[0,5*(2*n-1)],[-4/3,n^2*(36*n^2-1)]]
\end{verbatim}
$$\dfrac{\pi}{\sqrt{3}}-\dfrac{\log(432)}{3}=-\dfrac{4/3}{5+\dfrac{35}{15+\dfrac{572}{25+\dfrac{2907}{35+\dfrac{9200}{45+\dfrac{22475}{55+\ddots}}}}}}$$
Convergence type $P^-$ with $P=5/3$ and $C=-\G(5/6)^2/(3\pi)$,
so that
$$\dfrac{\pi}{\sqrt{3}}-\dfrac{\log(432)}{3}-\dfrac{p(n)}{q(n)}\sim(-1)^{n+1}\dfrac{\G(5/6)^2/(3\pi)}{n^{5/3}}\;.$$
$$A=1-(5/6)/n-(5/324)/n^2+(1375/1944)/n^3+(121/52488)/n^4+\cdots$$
Series:
$$\dfrac{\pi}{\sqrt{3}}-\dfrac{\log(432)}{3}=\dfrac{4}{15}\sum_{n\ge0}(-1)^{n+1}\dfrac{(7/6)_n}{(n+1)(11/6)_n}$$
Parametric family for $k\ge0$ and $u\ge0$:
\begin{verbatim}
[()->Pi/sqrt(3)-log(432)/3,(12*k+6*u+5)*(2*n-1),n^2*(36*n^2-(6*u-1)^2)]
\end{verbatim}
Convergence type $P^-$ with $P=4k+2u+5/3$.
\end{cf}

\smallskip

\begin{cf}\label{5.8.76}{\ }
\begin{verbatim}
[()->Pi/sqrt(3)-log(432)/3,[0,12*n^2-14*n+7],[-2/3,-n^2*(6*n-1)^2]]
\end{verbatim}
$$\dfrac{\pi}{\sqrt{3}}-\dfrac{\log(432)}{3}=-\dfrac{2/3}{5-\dfrac{25}{27-\dfrac{484}{73-\dfrac{2601}{143-\dfrac{8464}{237-\dfrac{21025}{355-\ddots}}}}}}$$
Convergence type $P^+$ with $P=1$ and $C=-1/9$, so that
$$\dfrac{\pi}{\sqrt{3}}-\dfrac{\log(432)}{3}-\dfrac{p(n)}{q(n)}\sim-\dfrac{1/9}{n}\;.$$
$$A=1-(5/12)/n+(5/54)/n^2+(25/864)/n^3-(17/648)/n^4+\cdots$$
Series:
$$\dfrac{\pi}{\sqrt{3}}-\dfrac{\log(432)}{3}=-\dfrac{2}{3}\sum_{n\ge1}\dfrac{1}{n(6n-1)}$$
Parametric family for $k\ge0$ and $u\ge0$:
\begin{verbatim}
[()->Pi/sqrt(3)-log(432)/3,12*n^2+(6*u-7)*(2*n-1)+6*k*(k+1),
                                          -n^2*(6*n+6*u-1)^2]
\end{verbatim}
Convergence type $P^+$ with $P=2k+1$.
\end{cf}

\smallskip

\begin{cf}\label{5.8.77}{\ }
\begin{verbatim}
[()->Pi/sqrt(3)-log(432)/3,[0,12*n^2-11*n+5],[-2/3,-6*n^3*(6*n+1)]]
\end{verbatim}
$$\dfrac{\pi}{\sqrt{3}}-\dfrac{\log(432)}{3}=-\dfrac{2/3}{6-\dfrac{42}{31-\dfrac{624}{80-\dfrac{3078}{153-\dfrac{9600}{250-\dfrac{23250}{371-\ddots}}}}}}$$
Convergence type $P^+$ with $P=5/6$ and $C=-4/(5\G(1/6))$, so that
$$\dfrac{\pi}{\sqrt{3}}-\dfrac{\log(432)}{3}-\dfrac{p(n)}{q(n)}\sim-\dfrac{4/(5\G(1/6))}{n^{5/6}}\;.$$
$$A=1-(355/792)/n+(27265/176256)/n^2+\cdots$$
Series:
$$\dfrac{\pi}{\sqrt{3}}-\dfrac{\log(432)}{3}=-\dfrac{1}{9}\sum_{n\ge0}\dfrac{(7/6)_n}{(n+1)(n+1)!}$$
Parametric family for $k\ge0$ and $u\ge0$:
\begin{verbatim}
[()->Pi/sqrt(3)-log(432)/3,12*n^2-(6*u+11)*n+6*k^2+(6*u+5)*(k+1),
                                            -6*n^3*(6*n-(6*u-1))]
\end{verbatim}
Convergence type $P^+$ with $P=2k+u+5/6$.
\end{cf}

\smallskip

\begin{cf}\label{5.8.78}{\ }
\begin{verbatim}
[()->Pi/sqrt(3)+log(432)/3,[0,2*n-1],[20/3,n^2*(36*n^2-25)]]
\end{verbatim}
$$\dfrac{\pi}{\sqrt{3}}+\dfrac{\log(432)}{3}=\dfrac{20/3}{1+\dfrac{11}{3+\dfrac{476}{5+\dfrac{2691}{7+\dfrac{8816}{9+\dfrac{21875}{11+\ddots}}}}}}$$
Convergence type $P^-$ with $P=1/3$ and $C=\G(1/6)^2/(3\pi)$,
so that
$$\dfrac{\pi}{\sqrt{3}}+\dfrac{\log(432)}{3}-\dfrac{p(n)}{q(n)}\sim(-1)^n\dfrac{\G(1/6)^2/(3\pi)}{n^{1/3}}\;.$$
$$A=1-(1/6)/n+(5/324)/n^2+(49/1944)/n^3-(217/104976)/n^4+\cdots$$
Series:
$$\dfrac{\pi}{\sqrt{3}}+\dfrac{\log(432)}{3}=\dfrac{20}{3}\sum_{n\ge0}(-1)^n\dfrac{(11/6)_n}{(n+1)(7/6)_n}$$
Parametric family for $k\ge0$ and $u\ge0$:
\begin{verbatim}
[()->Pi/sqrt(3)+log(432)/3,(12*k+6*u+1)*(2*n-1),n^2*(36*n^2-(6*u-5)^2)]
\end{verbatim}
Convergence type $P^-$ with $P=4k+2u+1/3$.
\end{cf}

\smallskip

\begin{cf}\label{5.8.79}{\ }
\begin{verbatim}
[()->Pi/sqrt(3)+log(432)/3,[0,12*n^2-22*n+11],[10/3,-n^2*(6*n-5)^2]]
\end{verbatim}
$$\dfrac{\pi}{\sqrt{3}}+\dfrac{\log(432)}{3}=\dfrac{10/3}{1-\dfrac{1}{15-\dfrac{196}{53-\dfrac{1521}{115-\dfrac{5776}{201-\dfrac{15625}{311-\ddots}}}}}}$$
Convergence type $P^+$ with $P=1$ and $C=5/9$, so that
$$\dfrac{\pi}{\sqrt{3}}+\dfrac{\log(432)}{3}-\dfrac{p(n)}{q(n)}\sim\dfrac{5/9}{n}\;.$$
$$A=1-(1/12)/n-(1/54)/n^2+(5/864)/n^3+(17/3240)/n^4+\cdots$$
Series:
$$\dfrac{\pi}{\sqrt{3}}+\dfrac{\log(432)}{3}=\dfrac{10}{3}\sum_{n\ge1}\dfrac{1}{n(6n-5)}$$
Parametric family for $k\ge0$ and $u\ge0$:
\begin{verbatim}
[()->Pi/sqrt(3)+log(432)/3,12*n^2+(6*u-11)*(2*n-1)+6*k*(k+1),
                                           -n^2*(6*n+6*u-5)^2]
\end{verbatim}
Convergence type $P^+$ with $P=2k+1$.
\end{cf}

\smallskip

\begin{cf}\label{5.8.80}{\ }
\begin{verbatim}
[()->Pi/sqrt(3)+log(432)/3,[0,12*n^2-7*n+1],[10/3,-6*n^3*(6*n+5)]]
\end{verbatim}
$$\dfrac{\pi}{\sqrt{3}}+\dfrac{\log(432)}{3}=\dfrac{10/3}{6-\dfrac{66}{35-\dfrac{816}{88-\dfrac{3726}{165-\dfrac{11136}{266-\dfrac{26250}{391-\ddots}}}}}}$$
Convergence type $P^+$ with $P=1/6$ and $C=4/\G(5/6)$, so that
$$\dfrac{\pi}{\sqrt{3}}+\dfrac{\log(432)}{3}-\dfrac{p(n)}{q(n)}\sim\dfrac{4/\G(5/6)}{n^{1/6}}\;.$$
$$A=1-(47/504)/n+(341/14976)/n^2+\cdots$$
Series:
$$\dfrac{\pi}{\sqrt{3}}+\dfrac{\log(432)}{3}=\dfrac{5}{9}\sum_{n\ge0}\dfrac{(11/6)_n}{(n+1)(n+1)!}$$
Parametric family for $k\ge0$ and $u\ge0$:
\begin{verbatim}
[()->Pi/sqrt(3)+log(432)/3,12*n^2-(6*u+7)*n+6*k^2+(6*u+1)*(k+1),
                                           -6*n^3*(6*n-(6*u-5))]
\end{verbatim}
Convergence type $P^+$ with $P=2k+u+1/6$.
\end{cf}

\smallskip

\begin{cf}\label{5.8.81}{\ }
\begin{verbatim}
[()->Pi/sqrt(2)+sqrt(2)*log(1+sqrt(2)),[2,1],[4,4*n*(n+1)]]
\end{verbatim}
$$\dfrac{\pi}{\sqrt{2}}+\sqrt{2}\log(1+\sqrt{2})=2+\dfrac{4}{1+\dfrac{8}{1+\dfrac{24}{1+\dfrac{48}{1+\dfrac{80}{1+\dfrac{120}{1+\ddots}}}}}}$$
Convergence type $P^-$ with $P=1/2$ and $C=\G(1/4)^2/2^{5/2}$, so that
$$\dfrac{\pi}{\sqrt{2}}+\sqrt{2}\log(1+\sqrt{2})-\dfrac{p(n)}{q(n)}\sim(-1)^n\dfrac{\G(1/4)^2/2^{5/2}}{n^{1/2}}\;.$$
$$A=1-(1/2)/n+(37/128)/n^2-(25/256)/n^3-(1017/32768)/n^4+\cdots$$
Parametric family for $k\ge0$:
\begin{verbatim}
[()->Pi/sqrt(2)+sqrt(2)*log(1+sqrt(2)),4*k+1,4*n*(n+1)]
\end{verbatim}
Convergence type $P^-$ with $P=2k+1/2$.
\end{cf}

\smallskip

\begin{cf}\label{5.8.82}{\ }
\begin{verbatim}
[()->Pi/sqrt(2)+sqrt(2)*log(1+sqrt(2)),[4,3],[-2,4*n^2]]
\end{verbatim}
$$\dfrac{\pi}{\sqrt{2}}+\sqrt{2}\log(1+\sqrt{2})=4-\dfrac{2}{3+\dfrac{4}{3+\dfrac{16}{3+\dfrac{36}{3+\dfrac{64}{3+\dfrac{100}{3+\ddots}}}}}}$$
Convergence type $P^-$ with $P=3/2$ and $C=\G(1/4)^2/2^{11/2}$, so that
$$\dfrac{\pi}{\sqrt{2}}+\sqrt{2}\log(1+\sqrt{2})-\dfrac{p(n)}{q(n)}\sim(-1)^n\dfrac{\G(1/4)^2/2^{11/2}}{n^{3/2}}\;.$$
$$A=1-(3/4)/n+(7/128)/n^2+(231/512)/n^3+(75/32768)/n^4+\cdots$$
Parametric families with $k\ge0$:
\begin{verbatim}
[()->Pi/sqrt(2)+sqrt(2)*log(1+sqrt(2)),4*k+3,4*n^2]
[()->Pi/sqrt(2)+sqrt(2)*log(1+sqrt(2)),4*k+1,4*n*(n+1)]
\end{verbatim}
Convergence type $P^-$ with $P=2k+3/2$ and $P=2k+1/2$ respectively.
\end{cf}

\smallskip

\begin{cf}\label{5.8.83}{\ }
\begin{verbatim}
[()->Pi/sqrt(2)+sqrt(2)*log(1+sqrt(2)),[0,1,4],[4,(4*n-3)^2]]
\end{verbatim}
$$\dfrac{\pi}{\sqrt{2}}+\sqrt{2}\log(1+\sqrt{2})=\dfrac{4}{1+\dfrac{1}{4+\dfrac{25}{4+\dfrac{81}{4+\dfrac{169}{4+\dfrac{289}{4+\ddots}}}}}}$$
Convergence type $P^-$ with $P=1$ and $C=1/2$, so that
$$\dfrac{\pi}{\sqrt{2}}+\sqrt{2}\log(1+\sqrt{2})-\dfrac{p(n)}{q(n)}\sim(-1)^n\dfrac{1/2}{n}\;.$$
$$A=1+(1/4)/n-(3/16)/n^2-(11/64)/n^3+(57/256)/n^4+(361/1024)/n^5+\cdots$$
Series:
$$\dfrac{\pi}{\sqrt{2}}+\sqrt{2}\log(1+\sqrt{2})=4\sum_{n\ge1}\dfrac{(-1)^{n+1}}{4n-3}$$
Parametric family with $k\ge0$:
\begin{verbatim}
[()->Pi/sqrt(2)+sqrt(2)*log(1+sqrt(2)),8*k+4,(4*n-3)^2]
\end{verbatim}
Convergence type $P^-$ with $P=2k+1$.
\end{cf}

\smallskip

\begin{cf}\label{5.8.84}{\ }
\begin{verbatim}
[()->Pi/sqrt(2)+sqrt(2)*log(1+sqrt(2)),[2,5,8],[8,(4*n+1)*(4*n-3)]]
\end{verbatim}
$$\dfrac{\pi}{\sqrt{2}}+\sqrt{2}\log(1+\sqrt{2})=2+\dfrac{8}{5+\dfrac{5}{8+\dfrac{45}{8+\dfrac{117}{8+\dfrac{221}{8+\dfrac{357}{8+\ddots}}}}}}$$
Convergence type $P^-$ with $P=2$ and $C=1/4$, so that
$$\dfrac{\pi}{\sqrt{2}}+\sqrt{2}\log(1+\sqrt{2})-\dfrac{p(n)}{q(n)}\sim(-1)^n\dfrac{1/4}{n^2}\;.$$
$$A=1-(1/2)/n-(5/16)/n^2+(7/16)/n^3+(181/256)/n^4-(691/512)/n^5+\cdots$$
Series:
$$\dfrac{\pi}{\sqrt{2}}+\sqrt{2}\log(1+\sqrt{2})=2+8\sum_{n\ge1}\dfrac{(-1)^{n+1}}{(4n-3)(4n+1)}$$
Parametric family with $k\ge0$:
\begin{verbatim}
[()->Pi/sqrt(2)+sqrt(2)*log(1+sqrt(2)),8*k+8,(4*n+1)*(4*n-3)]
\end{verbatim}
Convergence type $P^-$ with $P=2k+2$.
\end{cf}

\smallskip

\begin{cf}\label{5.8.85}{\ }
\begin{verbatim}
[()->Pi/sqrt(2)+sqrt(2)*log(1+sqrt(2)),[4,6*n-1],[-2,-2*n*(4*n+1)]]
\end{verbatim}
$$\dfrac{\pi}{\sqrt{2}}+\sqrt{2}\log(1+\sqrt{2})=4-\dfrac{2}{5-\dfrac{10}{11-\dfrac{36}{17-\dfrac{78}{23-\dfrac{136}{29-\dfrac{210}{35-\ddots}}}}}}$$
Convergence type $E$ with $E=2$, $P=5/4$, and $C=-\G(1/4)/4$, so that
$$\dfrac{\pi}{\sqrt{2}}+\sqrt{2}\log(1+\sqrt{2})-\dfrac{p(n)}{q(n)}\sim-\dfrac{\G(1/4)}{2^{n+2}n^{5/4}}$$
$$A=1-(85/32)/n+(18825/2048)/n^2+\cdots$$
Series:
$$\dfrac{\pi}{\sqrt{2}}+\sqrt{2}\log(1+\sqrt{2})=4-\dfrac{2}{5}\sum_{n\ge0}\dfrac{n!}{(9/4)_n}2^{-n}$$
Parametric family for $k\ge0$:
\begin{verbatim}
[()->Pi/sqrt(2)+sqrt(2)*log(1+sqrt(2)),6*n+2*k-1,-2*n*(4*n+1)]
\end{verbatim}
Convergence type $E$ with $E=2$ and $P=2k+5/4$.
\end{cf}
    
\smallskip

\begin{cf}\label{5.8.86}{\ }
\begin{verbatim}
[()->Pi/sqrt(2)+sqrt(2)*log(1+sqrt(2)),[0,-2,4*n-3],
                                      [[4,16],[(4*n-3)^2,(4*n+4)^2]]]
\end{verbatim}
$$\dfrac{\pi}{\sqrt{2}}+\sqrt{2}\log(1+\sqrt{2})=\dfrac{4}{-2+\dfrac{16}{5+\dfrac{1}{9+\dfrac{64}{13+\dfrac{25}{17+\dfrac{144}{21+\ddots}}}}}}$$
Convergence type $E$ with $E=-(1+\sqrt{2})^2$, $P=0$, and
$C=2\pi/(1+\sqrt{2})^{1/2}$, so that
$$\dfrac{\pi}{\sqrt{2}}+\sqrt{2}\log(1+\sqrt{2})-\dfrac{p(n)}{q(n)}\sim(-1)^n\dfrac{2\pi}{(1+\sqrt{2})^{2n+1/2}}\;.$$
$$A=1+(3d/16)/n+(-3d/64+9/256)/n^2+(645d/4096-9/512)/n^3+\cdots$$
\end{cf}

\smallskip

\begin{cf}\label{5.8.87}{\ }
\begin{verbatim}
[()->Pi/sqrt(2)+sqrt(2)*log(1+sqrt(2)),[4,4*n+1],
                                      [[-4,25],[16*n^2,(4*n+5)^2]]]
\end{verbatim}
$$\dfrac{\pi}{\sqrt{2}}+\sqrt{2}\log(1+\sqrt{2})=4-\dfrac{4}{5+\dfrac{25}{9+\dfrac{16}{13+\dfrac{81}{17+\dfrac{64}{21+\dfrac{169}{25+\ddots}}}}}}$$
Convergence type $E$ with $E=-(1+\sqrt{2})^2$, $P=0$, and
$C=-2\pi/(1+\sqrt{2})^{5/2}$, so that
$$\dfrac{\pi}{\sqrt{2}}+\sqrt{2}\log(1+\sqrt{2})-\dfrac{p(n)}{q(n)}\sim(-1)^{n+1}\dfrac{2\pi}{(1+\sqrt{2})^{2n+5/2}}\;.$$
$$A=1-(5d/16)/n+(25d/64+25/256)/n^2+(-2115d/4096-125/512)/n^3+\cdots$$
\end{cf}

\smallskip

\begin{cf}\label{5.8.88}{\ }
\begin{verbatim}
[()->Pi/sqrt(2)+sqrt(2)*log(1+sqrt(2)),
[2,4*n+1],[[8,5],[16*n*(n+1),(4*n+1)*(4*n+5)]]]
\end{verbatim}
$$\dfrac{\pi}{\sqrt{2}}+\sqrt{2}\log(1+\sqrt{2})=2+\dfrac{8}{5+\dfrac{5}{9+\dfrac{32}{13+\dfrac{45}{17+\dfrac{96}{21+\dfrac{117}{25+\ddots}}}}}}$$
Convergence type $E$ with $E=-(1+\sqrt{2})^2$, $P=0$, and
$C=2\pi/(1+\sqrt{2})^{5/2}$, so that
$$\dfrac{\pi}{\sqrt{2}}+\sqrt{2}\log(1+\sqrt{2})-\dfrac{p(n)}{q(n)}\sim(-1)^n\dfrac{2\pi}{(1+\sqrt{2})^{2n+5/2}}\;.$$
$$A=1+(11d/16)/n+(-55d/64+121/256)/n^2+(6093d/4096-605/512)/n^3+\cdots$$
\end{cf}

\smallskip

\begin{cf}\label{5.8.89}{\ }
\begin{verbatim}
[()->Pi/sqrt(2)+sqrt(2)*log(1+sqrt(2)),[0,116*n^2-236*n+117],
                                       [-10,-12*n*(2*n-1)*(6*n-7)*(6*n-5)]]
\end{verbatim}
$$\dfrac{\pi}{\sqrt{2}}+\sqrt{2}\log(1+\sqrt{2})=\dfrac{10}{3-\dfrac{12}{109-\dfrac{2520}{453-\dfrac{25740}{1029-\dfrac{108528}{1837-\dfrac{310500}{2877-\ddots}}}}}}$$
Convergence type $E$ with $E=27/2$, $P=1/2$, and $C=6\sqrt{\pi}/5$, so that
$$\dfrac{\pi}{\sqrt{2}}+\sqrt{2}\log(1+\sqrt{2})-\dfrac{p(n)}{q(n)}\sim\dfrac{6\sqrt{\pi}/5}{(27/2)^nn^{1/2}}$$
$$A=1-(347/1800)/n+(509929/6480000)/n^2+\cdots$$
Series:
$$\dfrac{\pi}{\sqrt{2}}+\sqrt{2}\log(1+\sqrt{2})=\dfrac{10}{3}\sum_{n\ge0}\dfrac{(1/2)_nn!}{(5/6)_n(7/6)_n}(2/27)^n$$
\end{cf}

\smallskip

\begin{cf}\label{5.8.90}{\ }
\begin{verbatim}
[()->Pi/sqrt(2)-sqrt(2)*log(1+sqrt(2)),[0,1],[2,4*n^2]]
\end{verbatim}
$$\dfrac{\pi}{\sqrt{2}}-\sqrt{2}\log(1+\sqrt{2})=\dfrac{2}{1+\dfrac{4}{1+\dfrac{16}{1+\dfrac{36}{1+\dfrac{64}{1+\dfrac{100}{1+\ddots}}}}}}$$
Convergence type $P^-$ with $P=1/2$ and $C=\G(3/4)^2/2^{1/2}$, so that
$$\dfrac{\pi}{\sqrt{2}}-\sqrt{2}\log(1+\sqrt{2})-\dfrac{p(n)}{q(n)}\sim(-1)^n\dfrac{\G(3/4)^2/2^{1/2}}{n^{1/2}}\;.$$
$$A=1-(1/4)/n-(3/128)/n^2+(55/512)/n^3+(199/32768)/n^4+\cdots$$
Parametric family for $k\ge0$:
\begin{verbatim}
[()->Pi/sqrt(2)-sqrt(2)*log(1+sqrt(2)),4*k+1,4*n^2]
\end{verbatim}
Convergence type $P^-$ with $P=2k+1/2$.
\end{cf}

\smallskip

\begin{cf}\label{5.8.91}{\ }
\begin{verbatim}
[()->Pi/sqrt(2)-sqrt(2)*log(1+sqrt(2)),[2/3,3],[4/3,4*n*(n+1)]]
\end{verbatim}
$$\dfrac{\pi}{\sqrt{2}}-\sqrt{2}\log(1+\sqrt{2})=2/3+\dfrac{4/3}{3+\dfrac{8}{3+\dfrac{24}{3+\dfrac{48}{3+\dfrac{80}{3+\dfrac{120}{3+\ddots}}}}}}$$
Convergence type $P^-$ with $P=3/2$ and $C=3\G(3/4)^2/2^{7/2}$, so that
$$\dfrac{\pi}{\sqrt{2}}-\sqrt{2}\log(1+\sqrt{2})-\dfrac{p(n)}{q(n)}\sim(-1)^n\dfrac{3\G(3/4)^2/2^{7/2}}{n^{3/2}}\;.$$
$$A=1-(3/2)/n+(199/128)/n^2-(273/256)/n^3+(14043/32768)/n^4+\cdots$$
Parametric families with $k\ge0$:
\begin{verbatim}
[()->Pi/sqrt(2)-sqrt(2)*log(1+sqrt(2)),4*k+1,4*n^2]
[()->Pi/sqrt(2)-sqrt(2)*log(1+sqrt(2)),4*k+3,4*n*(n+1)]
\end{verbatim}
Convergence type $P^-$ with $P=2k+1/2$ and $P=2k+3/2$ respectively.
\end{cf}

\smallskip

\begin{cf}\label{5.8.92}{\ }
\begin{verbatim}
[()->Pi/sqrt(2)-sqrt(2)*log(1+sqrt(2)),[0,3,4],[4,(4*n-1)^2]]
\end{verbatim}
$$\dfrac{\pi}{\sqrt{2}}-\sqrt{2}\log(1+\sqrt{2})=\dfrac{4}{3+\dfrac{9}{4+\dfrac{49}{4+\dfrac{121}{4+\dfrac{225}{4+\dfrac{361}{4+\ddots}}}}}}$$
Convergence type $P^-$ with $P=1$ and $C=1/2$, so that
$$\dfrac{\pi}{\sqrt{2}}-\sqrt{2}\log(1+\sqrt{2})-\dfrac{p(n)}{q(n)}\sim(-1)^n\dfrac{1/2}{n}\;.$$
$$A=1-(1/4)/n-(3/16)/n^2+(11/64)/n^3+(57/256)/n^4-(361/1024)/n^5+\cdots$$
Series:
$$\dfrac{\pi}{\sqrt{2}}-\sqrt{2}\log(1+\sqrt{2})=4\sum_{n\ge1}\dfrac{(-1)^{n+1}}{4n-1}$$
Parametric family with $k\ge0$:
\begin{verbatim}
[()->Pi/sqrt(2)-sqrt(2)*log(1+sqrt(2)),8*k+4,(4*n-1)^2]
\end{verbatim}
Convergence type $P^-$ with $P=2k+1$.
\end{cf}

\smallskip

\begin{cf}\label{5.8.93}{\ }
\begin{verbatim}
[()->Pi/sqrt(2)-sqrt(2)*log(1+sqrt(2)),[2/3,7,8],[8/3,(4*n-1)*(4*n+3)]]
\end{verbatim}
$$\dfrac{\pi}{\sqrt{2}}-\sqrt{2}\log(1+\sqrt{2})=2/3+\dfrac{8/3}{7+\dfrac{21}{8+\dfrac{77}{8+\dfrac{165}{8+\dfrac{285}{8+\dfrac{437}{8+\ddots}}}}}}$$
Convergence type $P^-$ with $P=2$ and $C=1/4$, so that
$$\dfrac{\pi}{\sqrt{2}}-\sqrt{2}\log(1+\sqrt{2})-\dfrac{p(n)}{q(n)}\sim(-1)^n\dfrac{1/4}{n^2}\;.$$
$$A=1-(3/2)/n+(19/16)/n^2-(3/16)/n^3-(59/256)/n^4-(873/512)/n^5+\cdots$$
Series:
$$\dfrac{\pi}{\sqrt{2}}-\sqrt{2}\log(1+\sqrt{2})=\dfrac{2}{3}+8\sum_{n\ge1}\dfrac{(-1)^{n+1}}{(4n-1)(4n+3)}$$
Parametric family with $k\ge0$:
\begin{verbatim}
[()->Pi/sqrt(2)-sqrt(2)*log(1+sqrt(2)),8*k+8,(4*n-1)*(4*n+3)]
\end{verbatim}
Convergence type $P^-$ with $P=2k+2$.
\end{cf}

\smallskip

\begin{cf}\label{5.8.94}{\ }
\begin{verbatim}
[()->Pi/sqrt(2)-sqrt(2)*log(1+sqrt(2)),[0,6*n-3],[2,-2*n*(4*n-1)]]
\end{verbatim}
$$\dfrac{\pi}{\sqrt{2}}-\sqrt{2}\log(1+\sqrt{2})=\dfrac{2}{3-\dfrac{6}{9-\dfrac{28}{15-\dfrac{66}{21-\dfrac{120}{27-\dfrac{190}{33-\ddots}}}}}}$$
Convergence type $E$ with $E=2$, $P=3/4$, and $C=\G(3/4)$, so that
$$\dfrac{\pi}{\sqrt{2}}-\sqrt{2}\log(1+\sqrt{2})-\dfrac{p(n)}{q(n)}\sim\dfrac{\G(3/4)}{2^nn^{3/4}}$$
$$A=1-(45/32)/n+(7385/2048)/n^2+\cdots$$
Series:
$$\dfrac{\pi}{\sqrt{2}}-\sqrt{2}\log(1+\sqrt{2})=\dfrac{2}{3}\sum_{n\ge0}\dfrac{n!}{(7/4)_n}2^{-n}$$
Parametric family for $k\ge0$:
\begin{verbatim}
[()->Pi/sqrt(2)-sqrt(2)*log(1+sqrt(2)),6*n+2*k-3,-2*n*(4*n-1)]
\end{verbatim}
Convergence type $E$ with $P=2k+3/4$.
\end{cf}

\smallskip

\begin{cf}\label{5.8.95}{\ }
\begin{verbatim}
[()->Pi/sqrt(2)-sqrt(2)*log(1+sqrt(2)),[0,2,4*n-1],
                                      [[4,16],[(4*n-1)^2,(4*n+4)^2]]]
\end{verbatim}
$$\dfrac{\pi}{\sqrt{2}}-\sqrt{2}\log(1+\sqrt{2})=\dfrac{4}{2+\dfrac{16}{7+\dfrac{9}{11+\dfrac{64}{15+\dfrac{49}{19+\dfrac{144}{23+\ddots}}}}}}$$
Convergence type $E$ with $E=-(1+\sqrt{2})^2$, $P=0$, and
$C=2\pi/(1+\sqrt{2})^{3/2}$, so that
$$\dfrac{\pi}{\sqrt{2}}-\sqrt{2}\log(1+\sqrt{2})-\dfrac{p(n)}{q(n)}\sim(-1)^n\dfrac{2\pi}{(1+\sqrt{2})^{2n+3/2}}\;.$$
$$A=1-(5d/16)/n+(15d/64+25/256)/n^2+(-835d/4096-75/512)/n^3+\cdots$$
\end{cf}

\smallskip

\begin{cf}\label{5.8.96}{\ }
\begin{verbatim}
[()->Pi/sqrt(2)-sqrt(2)*log(1+sqrt(2)),[0,4*n-1],
                                      [[4,9],[16*n^2,(4*n+3)^2]]]
\end{verbatim}
$$\dfrac{\pi}{\sqrt{2}}-\sqrt{2}\log(1+\sqrt{2})=\dfrac{4}{3+\dfrac{9}{7+\dfrac{16}{11+\dfrac{49}{15+\dfrac{64}{19+\dfrac{121}{23+\ddots}}}}}}$$
Convergence type $E$ with $E=-(1+\sqrt{2})^2$, $P=0$, and
$C=2\pi/(1+\sqrt{2})^{3/2}$, so that
$$\dfrac{\pi}{\sqrt{2}}-\sqrt{2}\log(1+\sqrt{2})-\dfrac{p(n)}{q(n)}\sim(-1)^n\dfrac{2\pi}{(1+\sqrt{2})^{2n+3/2}}\;.$$
$$A=1-(5d/16)/n+(15d/64+25/256)/n^2+(-835d/4096-75/512)/n^3+\cdots$$
\end{cf}

\smallskip

\begin{cf}\label{5.8.97}{\ }
\begin{verbatim}
[()->Pi/sqrt(2)-sqrt(2)*log(1+sqrt(2)),
[2/3,4*n+3],[[8/3,21],[16*n*(n+1),(4*n+3)*(4*n+7)]]]
\end{verbatim}
$$\dfrac{\pi}{\sqrt{2}}-\sqrt{2}\log(1+\sqrt{2})=2/3+\dfrac{8/3}{7+\dfrac{21}{11+\dfrac{32}{15+\dfrac{77}{19+\dfrac{96}{23+\dfrac{165}{27+\ddots}}}}}}$$
Convergence type $E$ with $E=-(1+\sqrt{2})^2$, $P=0$, and
$C=2\pi/(1+\sqrt{2})^{7/2}$, so that
$$\dfrac{\pi}{\sqrt{2}}-\sqrt{2}\log(1+\sqrt{2})-\dfrac{p(n)}{q(n)}\sim(-1)^n\dfrac{2\pi}{(1+\sqrt{2})^{2n+7/2}}\;.$$
$$A=1+(3d/16)/n+(-21d/64+9/256)/n^2+(1413d/4096-63/512)/n^3+\cdots$$
\end{cf}

\smallskip

Using \ref{6.10.5} and \ref{6.10.11}, exactly similar CFs can be given for
$\pi\pm\sqrt{3}\log(2+\sqrt{3})$ for instance.

\smallskip

\subsection{Homogeneous Periods of Degree $2$}

\smallskip

\begin{cf}\label{5.8.98}{\ }
\begin{verbatim}
[()->Pi^2-9*lfun(-3,2),[0,(2*n-1)*(3*n^2-3*n+2)],[4,-n^4*(9*n^2-1)]]
\end{verbatim}
$$\pi^2-9L(\chi_{-3},2)=\dfrac{4}{2-\dfrac{8}{24-\dfrac{560}{100-\dfrac{6480}{266-\dfrac{36608}{558-\dfrac{140000}{1012-\ddots}}}}}}$$
Convergence type $P^+$ with $P=4/3$ and $C=\G(2/3)/\G(4/3)$, so that
$$\pi^2-9L(\chi_{-3},2)-\dfrac{p(n)}{q(n)}\sim\dfrac{\G(2/3)/\G(4/3)}{n^{4/3}}$$
$$A=1-(2/3)/n+(103/405)/n^2+(2/243)/n^3+\cdots$$
Series:
$$\pi^2-9L(\chi_{-3},2)=2\sum_{n\ge0}\dfrac{(4/3)_n}{(n+1)^2(5/3)_n}$$
Parametric family for $k\ge0$:
\begin{verbatim}
[()->Pi^2-9*lfun(-3,2),(2*n-1)*(3*n^2-3*n+2+2*k*(3*k+2)),-n^4*(9*n^2-1)]
\end{verbatim}
Convergence type $P^+$ with $P=4k+4/3$.
\end{cf}

\smallskip

\begin{cf}\label{5.8.99}{\ }
\begin{verbatim}
[()->Pi^2+9*lfun(-3,2),[0,(2*n-1)*(3*n^2-3*n+1)],[8,-n^4*(9*n^2-4)]]
\end{verbatim}
$$\pi^2+9L(\chi_{-3},2)=\dfrac{8}{1-\dfrac{5}{21-\dfrac{512}{95-\dfrac{6237}{259-\dfrac{35840}{549-\dfrac{138125}{1001-\ddots}}}}}}$$
Convergence type $P^+$ with $P=2/3$ and $C=4\G(1/3)/\G(5/3)$, so that
$$\pi^2+9L(\chi_{-3},2)-\dfrac{p(n)}{q(n)}\sim\dfrac{4\G(1/3)/\G(5/3)}{n^{2/3}}$$
$$A=1-(1/3)/n+(31/324)/n^2-(1/243)/n^3+\cdots$$
Series:
$$\pi^2+9L(\chi_{-3},2)=8\sum_{n\ge0}\dfrac{(5/3)_n}{(n+1)^2(4/3)_n}$$
Parametric family for $k\ge0$:
\begin{verbatim}
[()->Pi^2+9*lfun(-3,2),(2*n-1)*(3*n^2-3*n+1+2*k*(3*k+1)),-n^4*(9*n^2-4)]
\end{verbatim}
Convergence type $P^+$ with $P=4k+2/3$.
\end{cf}

\smallskip

\begin{cf}\label{5.8.100}{\ }
\begin{verbatim}
[()->4*Pi^2-27*lfun(-3,2),[0,2*n-1],[36,9*n^4]]
\end{verbatim}
$$4\pi^2-27L(\chi_{-3},2)=\dfrac{36}{1+\dfrac{9}{3+\dfrac{144}{5+\dfrac{729}{7+\dfrac{2304}{9+\dfrac{5625}{11+\ddots}}}}}}$$
Convergence type $P^-$ with $P=2/3$ and $C=6\G(2/3)^4$, so that
$$4\pi^2-27L(\chi_{-3},2)-\dfrac{p(n)}{q(n)}\sim(-1)^n\dfrac{6\G(2/3)^4}{n^{2/3}}\;.$$
$$A=1-(1/3)/n-(5/81)/n^2+(50/243)/n^3+(134/6561)/n^4-(7406/19683)/n^5+\cdots$$
Parametric family with $k\ge0$:
\begin{verbatim}
[()->4*Pi^2-27*lfun(-3,2),(6*k+1)*(2*n-1),9*n^4]
\end{verbatim}
Convergence type $P^-$ with $P=4k+2/3$.
\end{cf}

\smallskip

\begin{cf}\label{5.8.101}{\ }
\begin{verbatim}
[()->4*Pi^2-27*lfun(-3,2),[0,4,18*n^2-30*n+17],[54,-(3*n-1)^4]]
\end{verbatim}
$$4\pi^2-27L(\chi_{-3},2)=\dfrac{54}{4-\dfrac{16}{29-\dfrac{625}{89-\dfrac{4096}{185-\dfrac{14641}{317-\dfrac{38416}{485-\ddots}}}}}}$$
Convergence type $P^+$ with $P=1$ and $C=6$, so that
$$4\pi^2-27L(\chi_{-3},2)-\dfrac{p(n)}{q(n)}\sim\dfrac{6}{n}\;.$$
$$A=1-(1/6)/n-(1/18)/n^2+(1/27)/n^3+(13/810)/n^4+\cdots$$
Series:
$$4\pi^2-27L(\chi_{-3},2)=54\sum_{n\ge1}\dfrac{1}{(3n-1)^2}$$
Parametric family with $k\ge0$:
\begin{verbatim}
[()->4*Pi^2-27*lfun(-3,2),18*n^2-30*n+17+9*k*(k+1),-(3*n-1)^4]
\end{verbatim}
Convergence type $P^+$ with $P=2k+1$.
\end{cf}

\smallskip

\begin{cf}\label{5.8.102}{\ }
\begin{verbatim}
[()->4*Pi^2-27*lfun(-3,2),[[0,2],[45*n^2-18*n+2,45*n^2+36*n+8]],
                         [[27,-16],[81*n^4,-4*(3*n+2)^4]]]
\end{verbatim}
$$4\pi^2-27L(\chi_{-3},2)=\dfrac{27}{2-\dfrac{16}{29+\dfrac{81}{89-\dfrac{2500}{146+\dfrac{1296}{260-\dfrac{16384}{353+\ddots}}}}}}$$
Convergence type $E$ with $E=-i((1+\sqrt{5})/2)^5$, $P=0$, and $C=24\pi^2/((1+\sqrt{5})/2)^{11/3}$, so that
$$4\pi^2-27L(\chi_{-3},2)-\dfrac{p(n)}{q(n)}\sim(-1)^n\dfrac{24\pi^2}{((1+\sqrt{5})/2)^{5n+11/3}}\;.$$
$$A=1-(14d/45)/n+(6d/25+98/405)/n^2+(-175708d/1366875-28/75)/n^3+\cdots$$
\end{cf}

\smallskip

\begin{cf}\label{5.8.103}{\ }
\begin{verbatim}
[()->4*Pi^2-27*lfun(-3,2),[[0,1],[45*n^2-18*n+2,45*n^2+18*n+5]],
                         [[36,27],[-4*(3*n-1)^4,81*(n+1)^4]]]
\end{verbatim}
$$4\pi^2-27L(\chi_{-3},2)=\dfrac{36}{1+\dfrac{27}{29-\dfrac{64}{68+\dfrac{1296}{146-\dfrac{2500}{221+\dfrac{6561}{353-\ddots}}}}}}$$
Convergence type $E$ with $E=i((1+\sqrt{5})/2)^5$, $P=0$, and $C=24\pi^2/((1+\sqrt{5})/2)^{11/3}$, so that
$$4\pi^2-27L(\chi_{-3},2)-\dfrac{p(n)}{q(n)}\sim(-1)^{\lfloor n/2\rfloor}\dfrac{24\pi^2}{((1+\sqrt{5})/2)^{5n+11/3}}\;.$$
$$A=1-(14d/45)/n+(6d/25+98/405)/n^2+(-175708d/1366875-28/75)/n^3+\cdots$$
\end{cf}

\smallskip

\begin{cf}\label{5.8.104}{\ }
\begin{verbatim}
[()->4*Pi^2+27*lfun(-3,2),[54,5*(2*n-1)],[36,9*n^4]]
\end{verbatim}
$$4\pi^2+27L(\chi_{-3},2)=54+\dfrac{36}{5+\dfrac{9}{15+\dfrac{144}{25+\dfrac{729}{35+\dfrac{2304}{45+\dfrac{5625}{55+\ddots}}}}}}$$
Convergence type $P^-$ with $P=10/3$ and $C=2\G(1/3)^4/27$, so that
$$4\pi^2+27L(\chi_{-3},2)-\dfrac{p(n)}{q(n)}\sim(-1)^n\dfrac{2\G(1/3)^4/27}{n^{10/3}}\;.$$
$$A=1-(5/3)/n+(5/81)/n^2+(740/243)/n^3-(271/6561)/n^4+\cdots$$
Parametric family with $k\ge0$:
\begin{verbatim}
[()->4*Pi^2+27*lfun(-3,2),(6*k+5)*(2*n-1),9*n^4]
\end{verbatim}
Convergence type $P^-$ with $P=4k+10/3$.
\end{cf}

\smallskip

\begin{cf}\label{5.8.105}{\ }
\begin{verbatim}
[()->4*Pi^2+27*lfun(-3,2),[0,1,18*n^2-42*n+29],[54,-(3*n-2)^4]]
\end{verbatim}
$$4\pi^2+27L(\chi_{-3},2)=\dfrac{54}{1-\dfrac{1}{17-\dfrac{256}{65-\dfrac{2401}{149-\dfrac{10000}{269-\dfrac{28561}{425-\ddots}}}}}}$$
Convergence type $P^+$ with $P=1$ and $C=6$, so that
$$4\pi^2+27L(\chi_{-3},2)-\dfrac{p(n)}{q(n)}\sim\dfrac{6}{n}\;.$$
$$A=1+(1/6)/n-(1/18)/n^2-(1/27)/n^3+(13/810)/n^4+\cdots$$
Series:
$$4\pi^2+27L(\chi_{-3},2)=54\sum_{n\ge1}\dfrac{1}{(3n-2)^2}$$
Parametric family with $k\ge0$:
\begin{verbatim}
[()->4*Pi^2+27*lfun(-3,2),18*n^2-42*n+29+9*k*(k+1),-(3*n-2)^4]
\end{verbatim}
Convergence type $P^+$ with $P=2k+1$.
\end{cf}

\smallskip

\begin{cf}\label{5.8.106}{\ }
\begin{verbatim}
[()->4*Pi^2+27*lfun(-3,2),[[54,8],[45*n^2+18*n+2,45*n^2+72*n+32]],
                         [[27,-256],[81*n^4,-4*(3*n+4)^4]]]
\end{verbatim}
$$4\pi^2+27L(\chi_{-3},2)=54+\dfrac{27}{8-\dfrac{256}{65+\dfrac{81}{149-\dfrac{9604}{218+\dfrac{1296}{356-\dfrac{40000}{461+\ddots}}}}}}$$
Convergence type $E$ with $E=-i((1+\sqrt{5})/2)^5$, $P=0$, and $C=24\pi^2/((1+\sqrt{5})/2)^{19/3}$,
so that
$$4\pi^2+27L(\chi_{-3},2)-\dfrac{p(n)}{q(n)}\sim(-1)^{\lfloor (n+1)/2\rfloor}\dfrac{4\pi^2}{((1+\sqrt{5})/2)^{5n+19/3}}\;.$$
$$A=1-(14d/45)/n+(86d/225+98/405)/n^2+(-564508d/1366875-1204/2025)/n^3+\cdots$$
\end{cf}

\smallskip

\begin{cf}\label{5.8.107}{\ }
\begin{verbatim}
[()->4*Pi^2+27*lfun(-3,2),[[54,5],[45*n^2+18*n+2,45*n^2+54*n+17]],
                         [[36,27],[-4*(3*n+1)^4,81*(n+1)^4]]]
\end{verbatim}
$$4\pi^2+27L(\chi_{-3},2)=54+\dfrac{36}{5+\dfrac{27}{65-\dfrac{1024}{116+\dfrac{1296}{218-\dfrac{9604}{305+\dfrac{6561}{461-\ddots}}}}}}$$
Convergence type $E$ with $E=i((1+\sqrt{5})/2)^5$, $P=0$, and  $C=24\pi^2/((1+\sqrt{5})/2)^{19/3}$,
so that
$$4\pi^2+27L(\chi_{-3},2)-\dfrac{p(n)}{q(n)}\sim(-1)^{\lfloor n/2\rfloor}\dfrac{24\pi^2}{((1+\sqrt{5})/2)^{5n+19/3}}\;.$$
$$A=1-(14d/45)/n+(86d/225+98/405)/n^2+(-564508d/1366875-1204/2025)/n^3+\cdots$$
\end{cf}

\smallskip

\begin{cf}\label{5.8.108}{\ }
\begin{verbatim}
[()->4*Pi^2-45*lfun(-3,2),[0,25,72*n^2-96*n+50],[72,-(6*n-1)^4]
\end{verbatim}
$$4\pi^2-45L(\chi_{-3},2)=\dfrac{72}{25-\dfrac{625}{146-\dfrac{14641}{410-\dfrac{83521}{818-\dfrac{279841}{1370-\dfrac{707281}{2066-\ddots}}}}}}$$
Convergence type $P^+$ with $P=1$ and $C=2$, so that
$$4\pi^2-45L(\chi_{-3},2)-\dfrac{p(n)}{q(n)}\sim\dfrac{2}{n}\;.$$
$$A=1-(1/3)/n+(1/36)/n^2+(5/108)/n^3+\cdots$$
Series:
$$4\pi^2-45L(\chi_{-3},2)=72\sum_{n\ge1}\dfrac{1}{(6n-1)^2}$$
Parametric family for $k\ge0$:
\begin{verbatim}
[()->4*Pi^2-45*lfun(-3,2),72*n^2-96*n+50+36*k*(k+1),-(6*n-1)^4]
\end{verbatim}
Convergence type $P^+$ with $P=2k+1$.
\end{cf}

\smallskip

\begin{cf}\label{5.8.109}{\ }
\begin{verbatim}
[()->4*Pi^2-45*lfun(-3,2),[0,2*(2*n-1)],[12,9*n^4]]
\end{verbatim}
$$4\pi^2-45L(\chi_{-3},2)=\dfrac{12}{2+\dfrac{9}{6+\dfrac{144}{10+\dfrac{729}{14+\dfrac{2304}{18+\dfrac{5625}{22+\ddots}}}}}}$$
Convergence type $P^-$ with $P=4/3$ and $C=2\G(5/6)^4$, so that
$$4\pi^2-45L(\chi_{-3},2)-\dfrac{p(n)}{q(n)}\sim(-1)^n\dfrac{2\G(5/6)^4}{n^{4/3}}\;.$$
$$A=1-(2/3)/n-(4/81)/n^2+(125/243)/n^3+(469/52488)/n^4+\cdots$$
Parametric family with $k\ge0$:
\begin{verbatim}
[()->4*Pi^2-45*lfun(-3,2),(6*k+2)*(2*n-1),9*n^4]
\end{verbatim}
Convergence type $P^-$ with $P=4k+4/3$.
\end{cf}

\smallskip

\begin{cf}\label{5.8.110}{\ }
\begin{verbatim}
[()->4*Pi^2-45*lfun(-3,2),[[0,25],[90*n^2-18*n+1,90*n^2+90*n+25]],
                         [[72,-625],[324*n^4,-(6*n+5)^4]]]
\end{verbatim}
$$4\pi^2-45L(\chi_{-3},2)=\dfrac{72}{25-\dfrac{625}{73+\dfrac{324}{205-\dfrac{14641}{325+\dfrac{5184}{565-\dfrac{83521}{757+\ddots}}}}}}$$
Convergence type $E$ with $E=-i((1+\sqrt{5})/2)^5$, $P=0$, and $C=8\pi^2/((1+\sqrt{5})/2)^{13/3}$, so that
$$4\pi^2-45L(\chi_{-3},2)-\dfrac{p(n)}{q(n)}\sim(-1)^{\lfloor(n+1)/2\rfloor}\dfrac{8\pi^2}{((1+\sqrt{5})/2)^{5n+13/3}}\;.$$
$$A=1-(17d/45)/n+(74d/225+289/810)/n^2+(-295633d/1366875-1258/2025)/n^3+\cdots$$
\end{cf}

\smallskip

\begin{cf}\label{5.8.111}{\ }
\begin{verbatim}
[()->4*Pi^2-45*lfun(-3,2),[[0,1],[90*n^2-18*n+1,90*n^2+54*n+13]],
                         [[6,27],[-(6*n-1)^4,324*(n+1)^4]]]
\end{verbatim}
$$4\pi^2-45L(\chi_{-3},2)=\dfrac{6}{1+\dfrac{27}{73-\dfrac{625}{157+\dfrac{5184}{325-\dfrac{14641}{481+\dfrac{26244}{757-\ddots}}}}}}$$
Convergence type $E$ with $E=i((1+\sqrt{5})/2)^5$, $P=0$, and  $C=8\pi^2/((1+\sqrt{5})/2)^{13/3}$, so that
$$4\pi^2-45L(\chi_{-3},2)-\dfrac{p(n)}{q(n)}\sim(-1)^{\lfloor n/2\rfloor}\dfrac{8\pi^2}{((1+\sqrt{5})/2)^{5n+13/3}}\;.$$
$$A=1-(17d/45)/n+(74d/225+289/810)/n^2+(-295633d/1366875-1258/2025)/n^3+\cdots$$
\end{cf}

\smallskip

\begin{cf}\label{5.8.112}{\ }
\begin{verbatim}
[()->4*Pi^2+45*lfun(-3,2),[72,49,72*n^2-48*n+26],[72,-(6*n+1)^4]]
\end{verbatim}
$$4\pi^2+45L(\chi_{-3},2)=72+\dfrac{72}{49-\dfrac{2401}{218-\dfrac{28561}{530-\dfrac{130321}{986-\dfrac{390625}{1586-\dfrac{923521}{2330-\ddots}}}}}}$$
Convergence type $P^+$ with $P=1$ and $C=2$, so that
$$4\pi^2+45L(\chi_{-3},2)-\dfrac{p(n)}{q(n)}\sim\dfrac{2}{n}\;.$$
$$A=1-(2/3)/n+(13/36)/n^2-(7/54)/n^3+\cdots$$
Series:
$$4\pi^2+45L(\chi_{-3},2)=72\sum_{n\ge1}\dfrac{1}{(6n-5)^2}$$
Parametric family for $k\ge0$:
\begin{verbatim}
[()->4*Pi^2+45*lfun(-3,2),72*n^2-48*n+26+36*k*(k+1),-(6*n+1)^4]
\end{verbatim}
Convergence type $P^+$ with $P=2k+1$.
\end{cf}

\smallskip

\begin{cf}\label{5.8.113}{\ }
\begin{verbatim}
[()->4*Pi^2+45*lfun(-3,2),[72,4*(2*n-1)],[12,9*n^4]]
\end{verbatim}
$$4\pi^2+45L(\chi_{-3},2)=72+\dfrac{12}{4+\dfrac{9}{12+\dfrac{144}{20+\dfrac{729}{28+\dfrac{2304}{36+\dfrac{5625}{44+\ddots}}}}}}$$
Convergence type $P^-$ with $P=8/3$ and $C=\G(1/6)^4/648$, so that
$$4\pi^2+45L(\chi_{-3},2)-\dfrac{p(n)}{q(n)}\sim(-1)^n\dfrac{\G(1/6)^4/648}{n^{8/3}}\;.$$
$$A=1-(4/3)/n+(4/81)/n^2+(434/243)/n^3-(211/26244)/n^4+\cdots$$
Parametric family with $k\ge0$:
\begin{verbatim}
[()->4*Pi^2+45*lfun(-3,2),(6*k+4)*(2*n-1),9*n^4]
\end{verbatim}
Convergence type $P^-$ with $P=4k+8/3$.
\end{cf}

\smallskip

\begin{cf}\label{5.8.114}{\ }
\begin{verbatim}
[()->4*Pi^2+45*lfun(-3,2),[[72,49],[90*n^2+18*n+1,90*n^2+126*n+49]],
                         [[72,-2401],[324*n^4,-(6*n+7)^4]]]
\end{verbatim}
$$4\pi^2+45L(\chi_{-3},2)=72+\dfrac{72}{49-\dfrac{2401}{109+\dfrac{324}{265-\dfrac{28561}{397+\dfrac{5184}{661-\dfrac{130321}{865+\ddots}}}}}}$$
Convergence type $E$ with $E=-i((1+\sqrt{5})/2)^5$, $P=0$, and $C=8\pi^2/((1+\sqrt{5})/2)^{17/3}$, so that
$$4\pi^2+45L(\chi_{-3},2)-\dfrac{p(n)}{q(n)}\sim(-1)^{\lfloor(n+1)/2\rfloor}\dfrac{8\pi^2}{((1+\sqrt{5})/2)^{5n+17/3}}\;.$$
$$A=1-(17d/45)/n+(32d/75+289/810)/n^2+(-562933d/1366875-544/675)/n^3+\cdots$$
\end{cf}

\smallskip

\begin{cf}\label{5.8.115}{\ }
\begin{verbatim}
[()->4*Pi^2+45*lfun(-3,2),[[72,2],[90*n^2+18*n+1,90*n^2+90*n+25]],
                         [[6,27],[-(6*n+1)^4,324*(n+1)^4]]]
\end{verbatim}
$$4\pi^2+45L(\chi_{-3},2)=72+\dfrac{6}{2+\dfrac{27}{109-\dfrac{2401}{205+\dfrac{5184}{397-\dfrac{28561}{565+\dfrac{26244}{865-\ddots}}}}}}$$
Convergence type $E$ with $E=i((1+\sqrt{5})/2)^5$, $P=0$, and $C=8\pi^2/((1+\sqrt{5})/2)^{17/3}$, so that
$$4\pi^2+45L(\chi_{-3},2)-\dfrac{p(n)}{q(n)}\sim(-1)^{\lfloor n/2\rfloor}\dfrac{8\pi^2}{((1+\sqrt{5})/2)^{5n+17/3}}\;.$$
$$A=1-(17d/45)/n+(32d/75+289/810)/n^2+(-562933d/1366875-544/675)/n^3+\cdots$$
\end{cf}

\smallskip

\begin{cf}\label{5.8.116}{\ }
\begin{verbatim}
[()->4*Pi^2-81*lfun(-3,2),[0,18*n^2-18*n+7],[-54,-81*n^4]]
\end{verbatim}
$$4\pi^2-81L(\chi_{-3},2)=-\dfrac{54}{7-\dfrac{81}{43-\dfrac{1296}{115-\dfrac{6561}{223-\dfrac{20736}{367-\dfrac{50625}{547-\ddots}}}}}}$$
Convergence type $P^+$ with $P=1/3$ and $C=-3\cdot2^{5/3}\G(5/6)^4$, so that
$$4\pi^2-81L(\chi_{-3},2)-\dfrac{p(n)}{q(n)}\sim-\dfrac{3\cdot2^{5/3}\G(5/6)^4}{n^{1/3}}\;.$$
$$A=1-(1/6)/n+(337/11340)/n^2+(83/9720)/n^3-(3327707/375289200)/n^4+\cdots$$
Parametric family with $k\ge0$:
\begin{verbatim}
[()->4*Pi^2-81*lfun(-3,2),18*n^2-18*n+9*k^2+3*k+7,-81*n^4]
\end{verbatim}
Convergence type $P^+$ with $P=2k+1/3$.
\end{cf}

\smallskip

\begin{cf}\label{5.8.117}{\ }
\begin{verbatim}
[()->4*Pi^2-81*lfun(-3,2),[-27/2,18*n^2-9*n+4],[-81,-81*n^3*(n+1)]]
\end{verbatim}
$$4\pi^2-81L(\chi_{-3},2)=-\dfrac{27}{2}-\dfrac{81}{13-\dfrac{162}{58-\dfrac{1944}{139-\dfrac{8748}{256-\dfrac{25920}{409-\dfrac{60750}{598-\ddots}}}}}}$$
Convergence type $P^+$ with $P=4/3$ and $C=-2^{8/3}\G(5/6)^4$, so that
$$4\pi^2-81L(\chi_{-3},2)-\dfrac{p(n)}{q(n)}\sim-\dfrac{2^{8/3}\G(5/6)^4}{n^{4/3}}\;.$$
$$A=1-(10/7)/n+(143/81)/n^2-(10942/5265)/n^3+(15854/6561)/n^4+\cdots$$
Parametric family with $k\ge0$:
\begin{verbatim}
[()->4*Pi^2-81*lfun(-3,2),18*n^2-9*n+(3*k+2)^2,-81*n^3*(n+1)]
\end{verbatim}
Convergence type $P^+$ with $P=2k+4/3$.
\end{cf}

\smallskip

\begin{cf}\label{5.8.118}{\ }
\begin{verbatim}
[()->4*Pi^2-81*lfun(-3,2),[0,4,3*(6*n-5)],[-108,(3*n-1)^4]]
\end{verbatim}
$$4\pi^2-81L(\chi_{-3},2)=-\dfrac{108}{4+\dfrac{16}{21+\dfrac{625}{39+\dfrac{4096}{57+\dfrac{14641}{75+\dfrac{38416}{93+\ddots}}}}}}$$
Convergence type $P^-$ with $P=2$ and $C=-6$, so that
$$4\pi^2-81L(\chi_{-3},2)-\dfrac{p(n)}{q(n)}\sim-\dfrac{6}{n^2}\;.$$
$$A=1-(1/3)/n-(2/3)/n^2+(13/27)/n^3+(110/81)/n^4+\cdots$$
Series:
$$4\pi^2-81L(\chi_{-3},2)=-108\sum_{n\ge1}\dfrac{(-1)^{n+1}}{(3n-1)^2}$$
Parametric family with $k\ge0$:
\begin{verbatim}
[()->4*Pi^2-81*lfun(-3,2),(6*k+3)*(6*n-5),(3*n-1)^4]
\end{verbatim}
Convergence type $P^-$ with $P=4k+2$.
\end{cf}
      
\smallskip

\begin{cf}\label{5.8.119}{\ }
\begin{verbatim}
[()->4*Pi^2-81*lfun(-3,2),[[0,4],[45*n^2-12*n+1,45*n^2+24*n+4]],
                         [[-108,16],[-324*n^4,(3*n+2)^4]]]
\end{verbatim}
$$4\pi^2-81L(\chi_{-3},2)=-\dfrac{108}{4+\dfrac{16}{34-\dfrac{324}{73+\dfrac{625}{157-\dfrac{5184}{232+\dfrac{4096}{370-\ddots}}}}}}$$
Convergence type $E$ with $E=i((1+\sqrt{5})/2)^5$, $P=0$, and
$C=-24\pi^2/((1+\sqrt{5})/2)^3$, so that
$$4\pi^2-81L(\chi_{-3},2)-\dfrac{p(n)}{q(n)}\sim(-1)^{n+1}\dfrac{24\pi^2}{((1+\sqrt{5})/2)^{5n+3}}\;.$$
$$A=1-(14d/45)/n+(118d/675+98/405)/n^2+(-56908d/1366875-1652/6075)/n^3+\cdots$$
\end{cf}

\smallskip

\begin{cf}\label{5.8.120}{\ }
\begin{verbatim}
[()->4*Pi^2+81*lfun(-3,2),[108,18*n^2-18*n+13],[-54,-81*n^4]]
\end{verbatim}
$$4\pi^2+81L(\chi_{-3},2)=108-\dfrac{54}{13-\dfrac{81}{49-\dfrac{1296}{121-\dfrac{6561}{229-\dfrac{20736}{373-\dfrac{50625}{553-\ddots}}}}}}$$
Convergence type $P^+$ with $P=5/3$ and $C=-\G(1/6)^4/(135\cdot 2^{5/3})$,
so that
$$4\pi^2+81L(\chi_{-3},2)-\dfrac{p(n)}{q(n)}\sim-\dfrac{\G(1/6)^4/(135\cdot 2^{5/3})}{n^{5/3}}\;.$$
$$A=1-(5/6)/n+(1385/3564)/n^2-(65/1944)/n^3-(1042067/12492144)/n^4+\cdots$$
Parametric family with $k\ge0$:
\begin{verbatim}
[()->4*Pi^2+81*lfun(-3,2),18*n^2-18*n+9*k^2+15*k+13,-81*n^4]
\end{verbatim}
Convergence type $P^+$ with $P=2k+5/3$.
\end{cf}

\smallskip

\begin{cf}\label{5.8.121}{\ }
\begin{verbatim}
[()->4*Pi^2+81*lfun(-3,2),[54,18*n^2-9*n+1],[162,-81*n^3*(n+1)]]
\end{verbatim}
$$4\pi^2+81L(\chi_{-3},2)=54+\dfrac{162}{10-\dfrac{162}{55-\dfrac{1944}{136-\dfrac{8748}{253-\dfrac{25920}{406-\dfrac{60750}{595-\ddots}}}}}}$$
Convergence type $P^+$ with $P=2/3$ and $C=\G(1/6)^4/(9\cdot 2^{2/3})$, so that
$$4\pi^2+81L(\chi_{-3},2)-\dfrac{p(n)}{q(n)}\sim\dfrac{\G(1/6)^4/(9\cdot 2^{2/3})}{n^{2/3}}\;.$$
$$A=1-(1/5)/n+(1/81)/n^2+(158/6237)/n^3-(766/45927)/n^4+\cdots$$
Parametric family with $k\ge0$:
\begin{verbatim}
[()->4*Pi^2+81*lfun(-3,2),18*n^2-9*n+(3*k+1)^2,-81*n^3*(n+1)]
\end{verbatim}
Convergence type $P^+$ with $P=2k+2/3$.
\end{cf}

\smallskip

\begin{cf}\label{5.8.122}{\ }
\begin{verbatim}
[()->4*Pi^2+81*lfun(-3,2),[0,1,3*(6*n-7)],[108,(3*n-2)^4]]
\end{verbatim}
$$4\pi^2+81L(\chi_{-3},2)=\dfrac{108}{1+\dfrac{1}{15+\dfrac{256}{33+\dfrac{2401}{51+\dfrac{10000}{69+\dfrac{28561}{87+\ddots}}}}}}$$
Convergence type $P^-$ with $P=2$ and $C=6$, so that
$$4\pi^2+81L(\chi_{-3},2)-\dfrac{p(n)}{q(n)}\sim\dfrac{6}{n^2}\;.$$
$$A=1+(1/3)/n-(2/3)/n^2-(13/27)/n^3+(110/81)/n^4+\cdots$$
Series:
$$4\pi^2+81L(\chi_{-3},2)=108\sum_{n\ge1}\dfrac{(-1)^{n+1}}{(3n-2)^2}$$
Parametric family with $k\ge0$:
\begin{verbatim}
[()->4*Pi^2+81*lfun(-3,2),(6*k+3)*(6*n-7),(3*n-2)^4]
\end{verbatim}
Convergence type $P^-$ with $P=4k+2$.
\end{cf}

\smallskip

\begin{cf}\label{5.8.123}{\ }
\begin{verbatim}
[()->4*Pi^2+81*lfun(-3,2),[[108,16],[45*n^2+12*n+1,45*n^2+48*n+16]],
                          [[-108,256],[-324*n^4,(3*n+4)^4]]]
\end{verbatim}
$$4\pi^2+81L(\chi_{-3},2)=108-\dfrac{108}{16+\dfrac{256}{58-\dfrac{324}{109+\dfrac{2401}{205-\dfrac{5184}{292+\dfrac{10000}{442-\ddots}}}}}}$$
Convergence type $E$ with $E=i((1+\sqrt{5})/2)^5$, $P=0$, and
$C=-24\pi^2/((1+\sqrt{5})/2)^7$, so that
$$4\pi^2+81L(\chi_{-3},2)-\dfrac{p(n)}{q(n)}\sim(-1)^{\lfloor(n+1)/2\rfloor}\dfrac{24\pi^2}{((1+\sqrt{5})/2)^{5n+7}}\;.$$
$$A=1-(14d/45)/n+(302d/675+98/405)/n^2+(-802108d/1366875-4228/6075)/n^3+\cdots$$
\end{cf}

\smallskip

\begin{cf}\label{5.8.124}{\ }
\begin{verbatim}
[()->11*Pi^2-135*lfun(-3,2),[0,(2*n-1)*(6*n^2-6*n+5)],[12,-n^4*(36*n^2-1)]]
\end{verbatim}
$$11\pi^2-135L(\chi_{-3},2)=\dfrac{12}{5-\dfrac{35}{51-\dfrac{2288}{205-\dfrac{26163}{539-\dfrac{147200}{1125-\dfrac{561875}{2035-\ddots}}}}}}$$
Convergence type $P^+$ with $P=5/3$ and $C=(36/5)(\G(5/6)/\G(1/6))$, so that
$$11\pi^2-135L(\chi_{-3},2)-\dfrac{p(n)}{q(n)}\sim\dfrac{(36/5)(\G(5/6)/\G(1/6))}{n^{5/3}}$$
$$A=1-(5/6)/n+(1295/3564)/n^2+(25/1944)/n^3+\cdots$$
Series:
$$11\pi^2-135L(\chi_{-3},2)=\dfrac{12}{5}\sum_{n\ge0}\dfrac{(7/6)_n}{(n+1)^2(11/6)_n}$$
Parametric family for $k\ge0$:
\begin{verbatim}
[()->11*Pi^2-135*lfun(-3,2),(2*n-1)*(6*n^2-6*n+5+2*k*(6*k+5)),-n^4*(36*n^2-1)]
\end{verbatim}
Convergence type $P^+$ with $P=4k+5/3$.
\end{cf}

\smallskip

\begin{cf}\label{5.8.125}{\ }
\begin{verbatim}
[()->11*Pi^2+135*lfun(-3,2),[0,(2*n-1)*(6*n^2-6*n+1)],[60,-n^4*(36*n^2-25)]]
\end{verbatim}
$$11\pi^2+135L(\chi_{-3},2)=\dfrac{60}{1-\dfrac{11}{39-\dfrac{1904}{185-\dfrac{24219}{511-\dfrac{141056}{1089-\dfrac{546875}{1991-\ddots}}}}}}$$
Convergence type $P^+$ with $P=1/3$ and $C=36(\G(1/6)/\G(5/6))$, so that
$$11\pi^2+135L(\chi_{-3},2)-\dfrac{p(n)}{q(n)}\sim\dfrac{36(\G(1/6)/\G(5/6))}{n^{1/3}}$$
$$A=1-(1/6)/n+(89/2268)/n^2-(5/1944)/n^3+\cdots$$
Series:
$$11\pi^2+135L(\chi_{-3},2)=60\sum_{n\ge0}\dfrac{(11/6)_n}{(n+1)^2(7/6)_n}$$
Parametric family for $k\ge0$:
\begin{verbatim}
[()->11*Pi^2+135*lfun(-3,2),(2*n-1)*(6*n^2-6*n+1+2*k*(6*k+1)),-n^4*(36*n^2-25)]
\end{verbatim}
Convergence type $P^+$ with $P=4k+1/3$.
\end{cf}

\smallskip

\begin{cf}\label{5.8.126}{\ }
\begin{verbatim}
[()->Pi^2-8*Catalan,[0,2*n-1],[4,4*n^4]]
\end{verbatim}
$$\pi^2-8G=\dfrac{4}{1+\dfrac{4}{3+\dfrac{64}{5+\dfrac{324}{7+\dfrac{1024}{9+\dfrac{2500}{11+\ddots}}}}}}$$
Convergence type $P^-$ with $P=1$ and $C=\G(3/4)^4$, so that
$$\pi^2-8G-\dfrac{p(n)}{q(n)}\sim(-1)^n\dfrac{\G(3/4)^4}{n}\;.$$
$$A=1-(1/2)/n-(1/16)/n^2+(11/32)/n^3+(1/64)/n^4-(89/128)/n^5+\cdots$$
Parametric family with $k\ge0$:
\begin{verbatim}
[()->Pi^2-8*Catalan,(4*k+1)*(2*n-1),4*n^4]
\end{verbatim}
Convergence type $P^-$ with $P=4k+1$.
\end{cf}

\smallskip

\begin{cf}\label{5.8.127}{\ }
\begin{verbatim}
[()->Pi^2-8*Catalan,[0,9,32*n^2-48*n+26],[16,-(4*n-1)^4]]
\end{verbatim}
$$\pi^2-8G=\dfrac{16}{9-\dfrac{81}{58-\dfrac{2401}{170-\dfrac{14641}{346-\dfrac{50625}{586-\dfrac{130321}{890-\ddots}}}}}}$$
Convergence type $P^+$ with $P=1$ and $C=1$, so that
$$\pi^2-8G-\dfrac{p(n)}{q(n)}\sim\dfrac{1}{n}\;.$$
$$A=1-(1/4)/n-(1/48)/n^2+(3/64)/n^3+(7/3840)/n^4+\cdots$$
Series:
$$\pi^2-8G=16\sum_{n\ge1}\dfrac{1}{(4n-1)^2}$$
Parametric family with $k\ge0$:
\begin{verbatim}
[()->Pi^2-8*Catalan,32*n^2-48*n+26+16*k*(k+1),-(4*n-1)^4]
\end{verbatim}
Convergence type $P^+$ with $P=2k+1$.
\end{cf}

\smallskip

\begin{cf}\label{5.8.128}{\ }
\begin{verbatim}
[()->Pi^2-8*Catalan,[[0,9],[40*n^2-12*n+1,40*n^2+36*n+9]],
                   [[16,-81],[64*n^4,-(4*n+3)^4]]]
\end{verbatim}
$$\pi^2-8G=\dfrac{16}{9-\dfrac{81}{29+\dfrac{64}{85-\dfrac{2401}{137+\dfrac{1024}{241-\dfrac{14641}{325+\ddots}}}}}}$$
Convergence type $E$ with $E=-i((1+\sqrt{5})/2)^5$, $P=0$, and $C=4\pi^2/((1+\sqrt{5})/2)^4$, so that
$$\pi^2-8G-\dfrac{p(n)}{q(n)}\sim(-1)^{\lfloor(n+1)/2\rfloor}\dfrac{4\pi^2}{((1+\sqrt{5})/2)^{5n+4}}\;.$$
$$A=1-(7d/20)/n+(57d/200+49/160)/n^2+(-13377d/80000-399/800)/n^3+\cdots$$
\end{cf}

\smallskip

\begin{cf}\label{5.8.129}{\ }
\begin{verbatim}
[()->Pi^2-8*Catalan,[[0,1],[40*n^2-12*n+1,40*n^2+20*n+5]],
                   [[4,16],[-(4*n-1)^4,64*(n+1)^4]]]
\end{verbatim}
$$\pi^2-8G=\dfrac{4}{1+\dfrac{16}{29-\dfrac{81}{65+\dfrac{1024}{137-\dfrac{2401}{205+\dfrac{5184}{325-\ddots}}}}}}$$
Convergence type $E$ with $E=i((1+\sqrt{5})/2)^5$, $P=0$, and $C=4\pi^2/((1+\sqrt{5})/2)^4$, so that
$$\pi^2-8G-\dfrac{p(n)}{q(n)}\sim(-1)^{\lfloor n/2\rfloor}\dfrac{4\pi^2}{((1+\sqrt{5})/2)^{5n+4}}\;.$$
$$A=1-(7d/20)/n+(57d/200+49/160)/n^2+(-13377d/80000-399/800)/n^3+\cdots$$
\end{cf}

\smallskip

\begin{cf}\label{5.8.130}{\ }
\begin{verbatim}
[()->Pi^2+8*Catalan,[16,3*(2*n-1)],[4,4*n^4]]
\end{verbatim}
$$\pi^2+8G=16+\dfrac{4}{3+\dfrac{4}{9+\dfrac{64}{15+\dfrac{324}{21+\dfrac{1024}{27+\dfrac{2500}{33+\ddots}}}}}}$$
Convergence type $P^-$ with $P=3$ and $C=\G(1/4)^4/256$, so that
$$\pi^2+8G-\dfrac{p(n)}{q(n)}\sim(-1)^n\dfrac{\G(1/4)^4/256}{n^3}\;.$$
$$A=1-(3/2)/n+(1/16)/n^2+(75/32)/n^3-(5/256)/n^4+\cdots$$
Parametric family with $k\ge0$:
\begin{verbatim}
[()->Pi^2+8*Catalan,(4*k+3)*(2*n-1),4*n^4]
\end{verbatim}
Convergence type $P^-$ with $P=4k+3$.
\end{cf}

\smallskip

\begin{cf}\label{5.8.131}{\ }
\begin{verbatim}
[()->Pi^2+8*Catalan,[0,1,32*n^2-80*n+58],[16,-(4*n-3)^4]]
\end{verbatim}
$$\pi^2+8G=\dfrac{16}{1-\dfrac{1}{26-\dfrac{625}{106-\dfrac{6561}{250-\dfrac{28561}{458-\dfrac{83521}{730-\ddots}}}}}}$$
Convergence type $P^+$ with $P=1$ and $C=1$, so that
$$\pi^2+8G-\dfrac{p(n)}{q(n)}\sim\dfrac{1}{n}\;.$$
$$A=1+(1/4)/n-(1/48)/n^2-(3/64)/n^3+(7/3840)/n^4+\cdots$$
Series:
$$\pi^2+8G=16\sum_{n\ge1}\dfrac{1}{(4n-3)^2}$$
Parametric family with $k\ge0$:
\begin{verbatim}
[()->Pi^2+8*Catalan,32*n^2-80*n+58+16*k*(k+1),-(4*n-3)^4]
\end{verbatim}
Convergence type $P^+$ with $P=2k+1$.
\end{cf}

\smallskip

\begin{cf}\label{5.8.132}{\ }
\begin{verbatim}
[()->Pi^2+8*Catalan,[[16,25],[40*n^2+12*n+1,40*n^2+60*n+25]],
                   [[16,-625],[64*n^4,-(4*n+5)^4]]]
\end{verbatim}
$$\pi^2+8G=16+\dfrac{16}{25-\dfrac{625}{53+\dfrac{64}{125-\dfrac{6561}{185+\dfrac{1024}{305-\dfrac{28561}{397+\ddots}}}}}}$$
Convergence type $E$ with $E=-i((1+\sqrt{5})/2)^5$, $P=0$, and $C=4\pi^2/((1+\sqrt{5})/2)^6$, so that
$$\pi^2+8G-\dfrac{p(n)}{q(n)}\sim(-1)^{\lfloor(n+1)/2\rfloor}\dfrac{4\pi^2}{((1+\sqrt{5})/2)^{5n+6}}\;.$$
$$A=1-(7d/20)/n+(83d/200+49/160)/n^2+(-34177d/80000-581/800)/n^3+\cdots$$
\end{cf}

\smallskip

\begin{cf}\label{5.8.133}{\ }
\begin{verbatim}
[()->Pi^2+8*Catalan,[[16,3],[40*n^2+12*n+1,40*n^2+44*n+13]],
                   [[4,16],[-(4*n+1)^4,64*(n+1)^4]]]
\end{verbatim}
$$\pi^2+8G=16+\dfrac{4}{3+\dfrac{16}{53-\dfrac{625}{97+\dfrac{1024}{185-\dfrac{6561}{261+\dfrac{5184}{397-\ddots}}}}}}$$
Convergence type $E$ with $E=i((1+\sqrt{5})/2)^5$, $P=0$, and $C=4\pi^2/((1+\sqrt{5})/2)^6$, so that
$$\pi^2+8G-\dfrac{p(n)}{q(n)}\sim(-1)^{\lfloor n/2\rfloor}\dfrac{4\pi^2}{((1+\sqrt{5})/2)^{5n+6}}\;.$$
$$A=1-(7d/20)/n+(83d/200+49/160)/n^2+(-34177d/80000-581/800)/n^3+\cdots$$
\end{cf}

\smallskip

\begin{cf}\label{5.8.134}{\ }
\begin{verbatim}
[()->Pi^2-16*Catalan,[0,(4*n-3)*(2*n^2-3*n+2)],
                     [-4,-2*n^3*(2*n-1)^3]]
\end{verbatim}
$$\pi^2-16G=-\dfrac{4}{1-\dfrac{2}{20-\dfrac{432}{99-\dfrac{6750}{286-\dfrac{43904}{629-\dfrac{182250}{1176-\ddots}}}}}}$$
Convergence type $P^+$ with $P=3/2$ and $C=-2\sqrt{\pi}/3$, so that
$$\pi^2-16G-\dfrac{p(n)}{q(n)}\sim-\dfrac{2\sqrt{\pi}/3}{n^{3/2}}$$
$$A=1-(3/8)/n-(17/896)/n^2+(87/1024)/n^3+\cdots$$
Series:
$$\pi^2-16G=-4\sum_{n\ge0}\dfrac{n!}{(n+1)(2n+1)(3/2)_n}$$
Parametric family for $k\ge0$:
\begin{verbatim}
[()->Pi^2-16*Catalan,(4*n-3)*(2*n^2-3*n+2+k*(4*k+3)),
                     -2*n^3*(2*n-1)^3]
\end{verbatim}
Convergence type $P^+$ with $P=4k+3/2$.
\end{cf}

\smallskip

\begin{cf}\label{5.8.135}{\ }
\begin{verbatim}
[()->Pi^2+16*Catalan,[24,(4*n-1)*(2*n^2-n+2)],
                     [4,-2*n^3*(2*n+1)^3]]
\end{verbatim}
$$\pi^2+16G=24+\dfrac{4}{9-\dfrac{54}{56-\dfrac{2000}{187-\dfrac{18522}{450-\dfrac{93312}{893-\dfrac{332750}{1564-\ddots}}}}}}$$
Convergence type $P^+$ with $P=5/2$ and $C=\sqrt{\pi}/5$, so that
$$\pi^2+16G-\dfrac{p(n)}{q(n)}\sim\dfrac{\sqrt{\pi}/5}{n^{5/2}}$$
$$A=1-(15/8)/n+(835/384)/n^2-(1845/1024)/n^3+\cdots$$
Series:
$$\pi^2+16G=24+\dfrac{4}{3}\sum_{n\ge0}\dfrac{n!}{(n+1)(2n+3)(5/2)_n}$$
Parametric family for $k\ge0$:
\begin{verbatim}
[()->Pi^2+16*Catalan,(4*n-1)*(2*n^2-n+2+k*(4*k+5)),
                     -2*n^3*(2*n+1)^3]
\end{verbatim}
Convergence type $P^+$ with $P=4k+5/2$.
\end{cf}

\smallskip

\begin{cf}\label{5.8.136}{\ }
\begin{verbatim}
[()->5*Pi^2-48*Catalan,[0,(2*n-1)*(4*n^2-4*n+3)],[12,-n^4*(16*n^2-1)]]
\end{verbatim}
$$5\pi^2-48G=\dfrac{12}{3-\dfrac{15}{33-\dfrac{1008}{135-\dfrac{11583}{357-\dfrac{65280}{747-\dfrac{249375}{1353-\ddots}}}}}}$$
Convergence type $P^+$ with $P=3/2$ and $C=8\G(3/4)/\G(1/4)$, so that
$$5\pi^2-48G-\dfrac{p(n)}{q(n)}\sim\dfrac{8\G(3/4)/\G(1/4)}{n^{3/2}}$$
$$A=1-(3/4)/n+(137/448)/n^2+(3/256)/n^3+\cdots$$
Series:
$$5\pi^2-48G=4\sum_{n\ge0}\dfrac{(5/4)_n}{(n+1)^2(7/4)_n}$$
Parametric family for $k\ge0$:
\begin{verbatim}
[()->5*Pi^2-48*Catalan,(2*n-1)*(4*n^2-4*n+3+2*k*(4*k+3)),-n^4*(16*n^2-1)]
\end{verbatim}
Convergence type $P^+$ with $P=4k+3/2$.
\end{cf}

\smallskip

\begin{cf}\label{5.8.137}{\ }
\begin{verbatim}
[()->5*Pi^2+48*Catalan,[0,(2*n-1)^3],[36,-n^4*(16*n^2-9)]]
\end{verbatim}
$$5\pi^2+48G=\dfrac{36}{1-\dfrac{7}{27-\dfrac{880}{125-\dfrac{10935}{343-\dfrac{63232}{729-\dfrac{244375}{1331-\ddots}}}}}}$$
Convergence type $P^+$ with $P=1/2$ and $C=24\G(1/4)/\G(3/4)$, so that
$$5\pi^2+48G-\dfrac{p(n)}{q(n)}\sim\dfrac{24\G(1/4)/\G(3/4)}{n^{1/2}}$$
$$A=1-(1/4)/n+(21/320)/n^2-(1/256)/n^3+\cdots$$
Series:
$$5\pi^2+48G=36\sum_{n\ge0}\dfrac{(7/4)_n}{(n+1)^2(5/4)_n}$$
Parametric family for $k\ge0$:
\begin{verbatim}
[()->5*Pi^2+48*Catalan,(2*n-1)*(4*n^2-4*n+1+2*k*(4*k+1)),-n^4*(16*n^2-9)]
\end{verbatim}
Convergence type $P^+$ with $P=4k+1/2$.
\end{cf}

\smallskip

\begin{cf}\label{5.8.138}{\ }
\begin{verbatim}
[()->Pi^2-6*log(2)^2,[0,3*n^2-2*n+1],[12,-2*n^4]]
\end{verbatim}
$$\pi^2-6\log^2(2)=\dfrac{12}{2-\dfrac{2}{9-\dfrac{32}{22-\dfrac{162}{41-\dfrac{512}{66-\dfrac{1250}{97-\ddots}}}}}}$$
Convergence type $E$ with $E=2$, $P=2$, and $C=12$, so that
$$\pi^2-6\log^2(2)-\dfrac{p(n)}{q(n)}\sim\dfrac{3}{2^{n-2}n^2}$$
$$A=1-4/n+18/n^2-104/n^3+750/n^4+\cdots$$
Series:
$$\pi^2-6\log(2)^2=6\sum_{n\ge0}\dfrac{2^{-n}}{(n+1)^2}$$
\end{cf}

\smallskip

\begin{cf}\label{5.8.139}{\ }
\begin{verbatim}
[()->Pi^2-12*log(2)^2,[0,29*n^3-32*n^2+10*n-1],
                      [24,-3*n^3*(2*n+1)*(3*n-2)*(3*n-1)]]
\end{verbatim}
$$\pi^2-12\log^2(2)=\dfrac{24}{6-\dfrac{18}{123-\dfrac{2400}{524-\dfrac{31752}{1383-\dfrac{190080}{2874-\dfrac{750750}{5171-\ddots}}}}}}$$
Convergence type $E$ with $E=27/2$, $P=3/2$, and $C=32\sqrt{3\pi}/25$, so that
$$\pi^2-12\log^2(2)-\dfrac{p(n)}{q(n)}\sim\dfrac{32\sqrt{3\pi}/25}{(27/2)^nn^{3/2}}$$
$$A=1-(2741/1800)/n+(2710229/1296000)/n^2+\cdots$$
Series:
$$\pi^2-12\log^2(2)=4\sum_{n\ge0}\dfrac{n!(3/2)_n}{(n+1)(4/3)_n(5/3)_n}(27/2)^{-n}$$
\end{cf}

\smallskip

\begin{cf}\label{5.8.140}{\ }
\begin{verbatim}
[()->Pi^2-3*log(3)^2,[0,113*n^3-161*n^2+82*n-16],
                     [96,-144*n^3*(2*n+1)*(3*n-2)*(3*n-1)]]
\end{verbatim}
$$\pi^2-3\log^2(3)=\dfrac{96}{18-\dfrac{864}{408-\dfrac{115200}{1832-\dfrac{1524096}{4968-\dfrac{9123840}{10494-\dfrac{36036000}{19088-\ddots}}}}}}$$
Convergence type $E$ with $E=81/32$, $P=3/2$, and $C=128\sqrt{3\pi}/49$,
so that
$$\pi^2-3\log^2(3)-\dfrac{p(n)}{q(n)}\sim\dfrac{128\sqrt{3\pi}/49}{(81/32)^nn^{3/2}}$$
$$A=1-(8405/3528)/n+\cdots$$
Series:
$$\pi^2-3\log^2(3)=\dfrac{16}{3}\sum_{n\ge0}\dfrac{n!(3/2)_n}{(n+1)(4/3)_n(5/3)_n}(81/32)^{-n}$$
\end{cf}

\smallskip

\begin{cf}\label{5.8.141}{\ }
\begin{verbatim}
[()->Pi^2-18*log((1+sqrt(5))/2)^2,[6,18,6*n^2+9*n+2],[-6,2*(n+1)*(2*n+1)^3]]
\end{verbatim}
$$\pi^2-18\log^2((1+\sqrt{5})/2)=6-\dfrac{6}{18+\dfrac{108}{44+\dfrac{750}{83+\dfrac{2744}{134+\dfrac{7290}{197+\dfrac{15972}{272+\ddots}}}}}}$$
Convergence type $E$ with $E=-4$, $P=3/2$, and $C=-3\sqrt{\pi}/10$, so that
$$\pi^2-18\log^2((1+\sqrt{5})/2)-\dfrac{p(n)}{q(n)}\sim(-1)^{n+1}\dfrac{3\sqrt{\pi}/5}{2^{2n+1}n^{3/2}}$$
$$A=1-(83/40)/n+(2101/640)/n^2-(104917/25600)/n^3+\cdots$$
Series:
$$\pi^2-18\log^2((1+\sqrt{5})/2)=6\sum_{n\ge0}(-1)^n\dfrac{n!}{(2n+1)(3/2)_n}2^{-2n}$$
\end{cf}

\smallskip

In the next CFs involving $L(\chi_8,2)$, note that of course
$L(\chi_8,2)=\pi^2/(8\sqrt{2})$.

\smallskip

\begin{cf}\label{5.8.142}{\ }
\begin{verbatim}
[()->lfun(-8,2)-lfun(8,2),[0,9,8*(4*n-3)],[2,(4*n-1)^4]]
\end{verbatim}
$$L(\chi_{-8},2)-L(\chi_8,2)=\dfrac{2}{9+\dfrac{81}{40+\dfrac{2401}{72+\dfrac{14641}{104+\dfrac{50625}{136+\dfrac{130321}{168+\ddots}}}}}}$$
Convergence type $P^-$ with $P=2$ and $C=1/16$, so that
$$L(\chi_{-8},2)-L(\chi_8,2)-\dfrac{p(n)}{q(n)}\sim(-1)^n\dfrac{1/16}{n^2}\;.$$
$$A=1-(1/2)/n-(9/16)/n^2+(11/16)/n^3+(285/256)/n^4-(1083/512)/n^5+\cdots$$
Series:
$$L(\chi_{-8},2)-L(\chi_8,2)=2\sum_{n\ge1}\dfrac{(-1)^{n+1}}{(4n-1)^2}$$
Parametric family with $k\ge0$:
\begin{verbatim}
[()->lfun(-8,2)-lfun(8,2),8*(2*k+1)*(4*n-3),(4*n-1)^4]
\end{verbatim}
Convergence type $P^-$ with $P=4k+2$.
\end{cf}

\smallskip

\begin{cf}\label{5.8.143}{\ }
\begin{verbatim}
[()->lfun(-8,2)-lfun(8,2),[0,32*n^2-32*n+13],[1,-256*n^4]]
\end{verbatim}
$$L(\chi_{-8},2)-L(\chi_8,2)=\dfrac{1}{13-\dfrac{256}{77-\dfrac{4096}{205-\dfrac{20736}{397-\dfrac{65536}{653-\dfrac{160000}{973-\ddots}}}}}}$$
Convergence type $P^+$ with $P=1/2$ and $C=\pi^5/\G(1/4)^6$, so that
$$L(\chi_{-8},2)-L(\chi_8,2)-\dfrac{p(n)}{q(n)}\sim\dfrac{\pi^5/\G(1/4)^6}{n^{1/2}}\;.$$
$$A=1-(1/4)/n+(419/7680)/n^2+(61/6144)/n^3-(470839/33030144)/n^4+\cdots$$
Parametric family with $k\ge0$:
\begin{verbatim}
[()->lfun(-8,2)-lfun(8,2),32*n^2-32*n+16*k^2+8*k+13,-256*n^4]
\end{verbatim}
Convergence type $P^+$ with $P=2k+1/2$.
\end{cf}

\smallskip

\begin{cf}\label{5.8.144}{\ }
\begin{verbatim}
[()->lfun(-8,2)-lfun(8,2),[[0,26],[80*n^2-16*n+1,80*n^2+80*n+25]],
                         [[2,-1024],[(4*n-1)^4,-1024*(n+1)^4]]]
\end{verbatim}
$$L(\chi_{-8},2)-L(\chi_8,2)=\dfrac{2}{26-\dfrac{1024}{65+\dfrac{81}{185-\dfrac{16384}{289+\dfrac{2401}{505-\dfrac{82944}{673+\ddots}}}}}}$$
Convergence type $E$ with $E=-i((1+\sqrt{5})/2)^5$, $P=0$, and
$C=(\pi^2/4)/((1+\sqrt{5})/2)^{7/2}$, so that
$$L(\chi_{-8},2)-L(\chi_8,2)-\dfrac{p(n)}{q(n)}\sim(-1)^{\lfloor(n-1)/2\rfloor}\dfrac{\pi^2/4}{((1+\sqrt{5})/2)^{5n+7/2}}\;.$$
$$A=1-(7d/20)/n+(6d/25+49/160)/n^2+(-7977d/80000-21/50)/n^3+\cdots$$
\end{cf}

\smallskip

\begin{cf}\label{5.8.145}{\ }
\begin{verbatim}
[()->lfun(-8,2)-lfun(8,2),[[0,9],[80*n^2-16*n+1,80*n^2+48*n+9]],
                         [[2,81],[-1024*n^4,(4*n+3)^4]]]
\end{verbatim}
$$L(\chi_{-8},2)-L(\chi_8,2)=\dfrac{2}{9+\dfrac{81}{65-\dfrac{1024}{137+\dfrac{2401}{289-\dfrac{16384}{425+\dfrac{14641}{673-\ddots}}}}}}$$
Convergence type $E$ with $E=i((1+\sqrt{5})/2)^5$, $P=0$, and
$C=(\pi^2/4)/((1+\sqrt{5})/2)^{7/2}$, so that
$$L(\chi_{-8},2)-L(\chi_8,2)-\dfrac{p(n)}{q(n)}\sim(-1)^{\lfloor(n+1)/2\rfloor}\dfrac{\pi^2/4}{((1+\sqrt{5})/2)^{5n+7/2}}\;.$$
$$A=1-(7d/20)/n+(6d/25+49/160)/n^2+(-7977d/80000-21/50)/n^3+\cdots$$
\end{cf}

\smallskip

\begin{cf}\label{5.8.146}{\ }
\begin{verbatim}
[()->lfun(-8,2)+lfun(8,2),[2,25,8*(4*n-1)],[-2,(4*n+1)^4]]
\end{verbatim}
$$L(\chi_{-8},2)+L(\chi_8,2)=2-\dfrac{2}{25+\dfrac{625}{56+\dfrac{6561}{88+\dfrac{28561}{120+\dfrac{83521}{152+\dfrac{194481}{184+\ddots}}}}}}$$
Convergence type $P^-$ with $P=2$ and $C=-1/16$, so that
$$L(\chi_{-8},2)+L(\chi_8,2)-\dfrac{p(n)}{q(n)}\sim(-1)^{n+1}\dfrac{1/16}{n^2}\;.$$
$$A=1-(3/2)/n+(15/16)/n^2+(9/16)/n^3-(275/256)/n^4+\cdots$$
Series:
$$L(\chi_{-8},2)+L(\chi_8,2)=2\sum_{n\ge1}\dfrac{(-1)^{n+1}}{(4n-3)^2}$$
Parametric family with $k\ge0$:
\begin{verbatim}
[()->lfun(-8,2)+lfun(8,2),8*(2*k+1)*(4*n-1),(4*n+1)^4]
\end{verbatim}
Convergence type $P^-$ with $P=4k+2$.
\end{cf}

\smallskip

\begin{cf}\label{5.8.147}{\ }
\begin{verbatim}
[()->lfun(-8,2)+lfun(8,2),[2,32*n^2-32*n+21],[-1,-256*n^4]]
\end{verbatim}
$$L(\chi_{-8},2)+L(\chi_8,2)=2-\dfrac{1}{21-\dfrac{256}{85-\dfrac{4096}{213-\dfrac{20736}{405-\dfrac{65536}{661-\dfrac{160000}{981-\ddots}}}}}}$$
Convergence type $P^+$ with $P=3/2$ and $C=-\G(1/4)^6/(24576\pi)$, so that
$$L(\chi_{-8},2)+L(\chi_8,2)-\dfrac{p(n)}{q(n)}\sim-\dfrac{\G(1/4)^6/(24576\pi)}{n^{3/2}}\;.$$
$$A=1-(3/4)/n+(1177/3584)/n^2-(57/2048)/n^3-(1809699/28835840)/n^4+\cdots$$
Parametric family with $k\ge0$:
\begin{verbatim}
[()->lfun(-8,2)+lfun(8,2),32*n^2-32*n+16*k^2+24*k+21,-256*n^4]
\end{verbatim}
Convergence type $P^+$ with $P=2k+3/2$.
\end{cf}

\smallskip

\begin{cf}\label{5.8.148}{\ }
\begin{verbatim}
[()->lfun(-8,2)+lfun(8,2),[[2,42],[80*n^2+16*n+1,80*n^2+112*n+41]],
                         [[-2,-1024],[(4*n+1)^4,-1024*(n+1)^4]]]
\end{verbatim}
$$L(\chi_{-8},2)+L(\chi_8,2)=2-\dfrac{2}{42-\dfrac{1024}{97+\dfrac{625}{233-\dfrac{16384}{353+\dfrac{6561}{585-\dfrac{82944}{769+\ddots}}}}}}$$
Convergence type $E$ with $E=-i((1+\sqrt{5})/2)^5$, $P=0$, and
$C=-(\pi^2/4)/((1+\sqrt{5})/2)^{13/2}$, so that
$$L(\chi_{-8},2)+L(\chi_8,2)-\dfrac{p(n)}{q(n)}\sim(-1)^{\lfloor(n+1)/2\rfloor}\dfrac{\pi^2/4}{((1+\sqrt{5})/2)^{5n+13/2}}\;.$$
$$A=1-(7d/20)/n+(23d/50+49/160)/n^2+(-43177d/80000-161/200)/n^3+\cdots$$
\end{cf}

\smallskip

\begin{cf}\label{5.8.149}{\ }
\begin{verbatim}
[()->lfun(-8,2)+lfun(8,2),[[2,25],[80*n^2+16*n+1,80*n^2+80*n+25]],
                         [[-2,625],[-1024*n^4,(4*n+5)^4]]]
\end{verbatim}
$$L(\chi_{-8},2)+L(\chi_8,2)=2-\dfrac{2}{25+\dfrac{625}{97-\dfrac{1024}{185+\dfrac{6561}{353-\dfrac{16384}{505+\dfrac{28561}{769-\ddots}}}}}}$$
Convergence type $E$ with $E=i((1+\sqrt{5})/2)^5$, $P=0$, and
$C=-(\pi^2/4)/((1+\sqrt{5})/2)^{13/2}$, so that
$$L(\chi_{-8},2)+L(\chi_8,2)-\dfrac{p(n)}{q(n)}\sim(-1)^{n+1}\dfrac{\pi^2/4}{((1+\sqrt{5})/2)^{5n+13/2}}\;.$$
$$A=1-(7d/20)/n+(23d/50+49/160)/n^2+(-43177d/80000-161/200)/n^3+\cdots$$
\end{cf}

\smallskip

\begin{cf}\label{5.8.150}{\ }
\begin{verbatim}
[()->8*Catalan-Pi*log(2+sqrt(3)),[3,18,10*n^2+7*n+2],
                                 [3,-2*(n+1)*(2*n+1)^3]]
\end{verbatim}
$$8G-\pi\log(2+\sqrt{3})=3+\dfrac{3}{18-\dfrac{108}{56-\dfrac{750}{113-\dfrac{2744}{190-\dfrac{7290}{287-\dfrac{15972}{404-\ddots}}}}}}$$
Convergence type $E$ with $E=4$, $P=3/2$, and $C=\sqrt{\pi}/4$, so that
$$8G-\pi\log(2+\sqrt{3})-\dfrac{p(n)}{q(n)}\sim\dfrac{\sqrt{\pi}}{2^{2n+2}n^{3/2}}$$
$$A=1-(23/8)/n+(2963/384)/n^2-(208037/9216)/n^3+\cdots$$
Series:
$$8G-\pi\log(2+\sqrt{3})=3\sum_{n\ge0}\dfrac{n!}{(2n+1)(3/2)_n}2^{-2n}$$
\end{cf}

\smallskip

\begin{cf}\label{5.8.151}{\ }
\begin{verbatim}
[()->(4*Catalan+Pi*log(2))/sqrt(2),[4,36,24*n^3+4*n^2+10*n-1],
                                   [4,-4*n*(2*n+1)^5]]
\end{verbatim}
$$\dfrac{4G+\pi\log(2)}{\sqrt{2}}=4+\dfrac{4}{36-\dfrac{972}{227-\dfrac{25000}{713-\dfrac{201684}{1639-\dfrac{944784}{3149-\dfrac{3221020}{5387-\ddots}}}}}}$$
Convergence type $E$ with $E=2$, $P=5/2$, and $C=1/\sqrt{\pi}$, so that
$$\dfrac{4G+\pi\log(2)}{\sqrt{2}}-\dfrac{p(n)}{q(n)}\sim\dfrac{1/\sqrt{\pi}}{2^nn^{5/2}}\;.$$
$$A=1-(49/8)/n+(4481/128)/n^2-(237899/1024)/n^3+\cdots$$
Series:
$$\dfrac{4G+\pi\log(2)}{\sqrt{2}}=4+\sum_{n\ge0}\dfrac{(3/2)_n}{(2n+3)^2(n+1)!}2^{-n}$$
\end{cf}

\smallskip

\subsection{Homogeneous Periods of Degree $3$}

\smallskip

\begin{cf}\label{5.8.152}{\ }
\begin{verbatim}
[()->28*zeta(3)-Pi^3,[0,(2*n-1)*(8*n^2-8*n+5)],[8,-64*n^6]]
\end{verbatim}
$$28\z(3)-\pi^3=\dfrac{8}{5-\dfrac{64}{63-\dfrac{4096}{265-\dfrac{46656}{707-\dfrac{262144}{1485-\dfrac{1000000}{2695-\ddots}}}}}}$$
Convergence type $P^+$ with $P=1$ and $C=16\pi^7/\G(1/4)^8$, so that
$$28\z(3)-\pi^3-\dfrac{p(n)}{q(n)}\sim\dfrac{16\pi^7/\G(1/4)^8}{n}\;.$$
$$A=1-c/n+(c^2-7/64)/n^2+(c^3-(7/32)c+7/128)/n^3+\cdots\;,$$
with $c=0.4738999777589193\cdots$, but I have not been able to recognize $c$.
Parametric family for $k\ge0$:
\begin{verbatim}
[()->28*zeta(3)-Pi^3,(2*n-1)*(8*n^2-8*n+4+(4*k+1)^2),-64*n^6]
\end{verbatim}
Convergence type $P^+$ with $P=4k+1$.
\end{cf}

\smallskip

\begin{cf}\label{5.8.153}{\ }
\begin{verbatim}
[()->28*zeta(3)-Pi^3,[0,27,2*(4*n-3)*(16*n^2-24*n+21)],[64,-(4*n-1)^6]]
\end{verbatim}
$$28\z(3)-\pi^3=\dfrac{64}{27-\dfrac{729}{370-\dfrac{117649}{1674-\dfrac{1771561}{4706-\dfrac{11390625}{10234-\dfrac{47045881}{19026-\ddots}}}}}}$$
Convergence type $P^+$ with $P=2$ and $C=1/2$, so that
$$28\z(3)-\pi^3-\dfrac{p(n)}{q(n)}\sim\dfrac{1/2}{n^2}\;.$$
$$A=1-(1/2)/n-(1/16)/n^2+(3/16)/n^3+(7/768)/n^4+\cdots$$
Series:
$$28\z(3)-\pi^3=64\sum_{n\ge1}\dfrac{1}{(4n-1)^3}$$
Parametric family with $k\ge0$:
\begin{verbatim}
[()->28*zeta(3)-Pi^3,2*(4*n-3)*(16*n^2-24*n+21+32*k*(k+1)),-(4*n-1)^6]
\end{verbatim}
Convergence type $P^+$ with $P=4k+2$.
\end{cf}
        
\smallskip

\begin{cf}\label{5.8.154}{\ }
\begin{verbatim}
[()->28*zeta(3)-Pi^3,
[[0,27],[(8*n-1)*(48*n^2-12*n+1),(8*n+3)*(48*n^2+36*n+9)]],
[[64,-729],[-4096*n^6,-(4*n+3)^6]]]
\end{verbatim}
$$28\z(3)-\pi^3=\dfrac{64}{27-\dfrac{729}{259-\dfrac{4096}{1023-\dfrac{117649}{2535-\dfrac{262144}{5187-\dfrac{1771561}{9131-\ddots}}}}}}$$
Convergence type $E$ with $E=(1+\sqrt{2})^4$, $P=0$, and $C=4\pi^3/(1+\sqrt{2})^3$, so that
$$28\z(3)-\pi^3-\dfrac{p(n)}{q(n)}\sim\dfrac{4\pi^3}{(1+\sqrt{2})^{4n+3}}\;.$$
$$A=1-(13d/16)/n+(39d/64+169/256)/n^2+(-1833d/4096-507/512)/n^3+\cdots$$
\end{cf}

\smallskip

\begin{cf}\label{5.8.155}{\ }
\begin{verbatim}
[()->28*zeta(3)+Pi^3,[64,(2*n-1)*(8*n^2-8*n+13)],[8,-64*n^6]]
\end{verbatim}
$$28\z(3)+\pi^3=64+\dfrac{8}{13-\dfrac{64}{87-\dfrac{4096}{305-\dfrac{46656}{763-\dfrac{262144}{1557-\dfrac{1000000}{2783-\ddots}}}}}}$$
Convergence type $P^+$ with $P=3$ and $C=\G(1/4)^8/(49152\pi)$, so that
$$28\z(3)+\pi^3-\dfrac{p(n)}{q(n)}\sim\dfrac{\G(1/4)^8/(49152\pi)}{n^3}\;.$$
$$A=1-(3/2)/n+(67/64)/n^2+\cdots$$
Note: this asymptotic expansion seems to continue, but I have not been
able to recognize the next term $-0.12030551\cdots/n^3$, situation identical
to \ref{5.8.152}.
Parametric family for $k\ge0$:
\begin{verbatim}
[()->28*zeta(3)+Pi^3,(2*n-1)*(8*n^2-8*n+4+(4*k+3)^2),-64*n^6]
\end{verbatim}
Convergence type $P^+$ with $P=4k+3$.
\end{cf}

\smallskip

\begin{cf}\label{5.8.156}{\ }
\begin{verbatim}
[()->28*zeta(3)+Pi^3,[0,1,2*(4*n-5)*(16*n^2-40*n+37)],[64,-(4*n-3)^6]]
\end{verbatim}
$$28\z(3)+\pi^3=\dfrac{64}{1-\dfrac{1}{126-\dfrac{15625}{854-\dfrac{531441}{2926-\dfrac{4826809}{7110-\dfrac{24137569}{14174-\ddots}}}}}}$$
Convergence type $P^+$ with $P=2$ and $C=1/2$, so that
$$28\z(3)+\pi^3-\dfrac{p(n)}{q(n)}\sim\dfrac{1/2}{n^2}\;.$$
$$A=1+(1/2)/n-(1/16)/n^2-(3/16)/n^3+(7/768)/n^4+\cdots$$
Series:
$$28\z(3)+\pi^3=64\sum_{n\ge1}\dfrac{1}{(4n-3)^3}$$
Parametric family with $k\ge0$:
\begin{verbatim}
[()->28*zeta(3)+Pi^3,2*(4*n-5)*(16*n^2-40*n+37+32*k*(k+1)),-(4*n-3)^6]
\end{verbatim}
Convergence type $P^+$ with $P=4k+2$.
\end{cf}

\smallskip

\begin{cf}\label{5.8.157}{\ }
\begin{verbatim}
[()->28*zeta(3)+Pi^3,
[[64,125],[(8*n+1)*(48*n^2+12*n+1),(8*n+5)*(48*n^2+60*n+25)]],
[[64,-15625],[-4096*n^6,-(4*n+5)^6]]]
\end{verbatim}
$$28\z(3)+\pi^3=64+\dfrac{64}{125-\dfrac{15625}{549-\dfrac{4096}{1729-\dfrac{531441}{3689-\dfrac{262144}{7077-\dfrac{4826809}{11725-\ddots}}}}}}$$
Convergence type $E$ with $E=(1+\sqrt{2})^4$, $P=0$, and $C=4\pi^3/(1+\sqrt{2})^5$, so that
$$28\z(3)+\pi^3-\dfrac{p(n)}{q(n)}\sim\dfrac{4\pi^3}{(1+\sqrt{2})^{4n+5}}\;.$$
$$A=1-(13d/16)/n+(65d/64+169/256)/n^2+(-5161d/4096-845/512)/n^3+\cdots$$
\end{cf}

\smallskip

\begin{cf}\label{5.8.158}{\ }
\begin{verbatim}
[()->91*zeta(3)-2*Pi^3*sqrt(3),
[0,(2*n-1)*(18*n^2-18*n+13)],[18,-324*n^6]]
\end{verbatim}
$$91\z(3)-2\pi^3\sqrt{3}=\dfrac{18}{13-\dfrac{324}{147-\dfrac{20736}{605-\dfrac{236196}{1603-\dfrac{1327104}{3357-\dfrac{5062500}{6083-\ddots}}}}}}$$
Convergence type $P^+$ with $P=4/3$ and
$C=3^{3/2}\G(5/6)^7/(2^{8/3}\pi^{1/2})$, so that
$$91\z(3)-2\pi^3\sqrt{3}-\dfrac{p(n)}{q(n)}\sim\dfrac{3^{3/2}\G(5/6)^7/(2^{8/3}\pi^{1/2})}{n^{4/3}}\;.$$
$$A=1-(2/3)/n+(77/324)/n^2+(35/972)/n^3-(2324029/26873856)/n^4+\cdots$$
Parametric family for $k\ge0$:
\begin{verbatim}
[()->91*zeta(3)-2*Pi^3*sqrt(3),
     (2*n-1)*(18*n^2-18*n+13+36*k^2+24*k),-324*n^6]
\end{verbatim}
Convergence type $P^+$ with $P=4k+4/3$.
\end{cf}

\smallskip

\begin{cf}\label{5.8.159}{\ }
\begin{verbatim}
[()->91*zeta(3)-2*Pi^3*sqrt(3),
[0,125,4*(3*n-2)*(36*n^2-48*n+43)],[216,-(6*n-1)^6]]
\end{verbatim}
$$91\z(3)-2\pi^3\sqrt{3}=\dfrac{216}{125-\dfrac{15625}{1456-\dfrac{1771561}{6244-\dfrac{24137569}{17080-\dfrac{148035889}{36556-\ddots}}}}}$$
Convergence type $P^+$ with $P=2$ and $C=1/2$, so that
$$91\z(3)-2\pi^3\sqrt{3}-\dfrac{p(n)}{q(n)}\sim\dfrac{1/2}{n^2}\;.$$
$$A=1-(2/3)/n+(1/12)/n^2+(5/27)/n^3+\cdots$$
Series:
$$91\z(3)-2\pi^3\sqrt{3}=216\sum_{n\ge1}\dfrac{1}{(6n-1)^3}$$
Parametric family for $k\ge0$:
\begin{verbatim}
[()->91*zeta(3)-2*Pi^3*sqrt(3),
4*(3*n-2)*(36*n^2-48*n+43+72*k*(k+1)),-(6*n-1)^6]
\end{verbatim}
Convergence type $P^+$ with $P=4k+2$.
\end{cf}

\smallskip

\begin{cf}\label{5.8.160}{\ }
\begin{verbatim}
[()->91*zeta(3)-2*Pi^3*sqrt(3),
[[0,125],[(12*n-1)*(108*n^2-18*n+1),(12*n+5)*(108*n^2+90*n+25)]],
[[216,-15625],[-46656*n^6,-(6*n+5)^6]]]
\end{verbatim}
$$91\z(3)-2\pi^3\sqrt{3}=\dfrac{216}{125-\dfrac{15625}{1001-\dfrac{46656}{3791-\dfrac{1771561}{9131-\dfrac{2985984}{18473-\dfrac{24137569}{32165-\ddots}}}}}}$$
Convergence type $E$ with $E=(1+\sqrt{2})^4$, $P=0$, and $C=4\pi^3/(1+\sqrt{2})^{10/3}$, so that
$$91\z(3)-2\pi^3\sqrt{3}-\dfrac{p(n)}{q(n)}\sim\dfrac{4\pi^3}{(1+\sqrt{2})^{4n+10/3}}\;.$$
$$A=1-(127d/144)/n+(635d/864+16129/20736)/n^2+\cdots$$
\end{cf}

\smallskip

\begin{cf}\label{5.8.161}{\ }
\begin{verbatim}
[()->91*zeta(3)+2*Pi^3*sqrt(3),
[216,(2*n-1)*(18*n^2-18*n+25)],[18,-324*n^6]]
\end{verbatim}
$$91\z(3)+2\pi^3\sqrt{3}=216+\dfrac{18}{25-\dfrac{324}{183-\dfrac{20736}{665-\dfrac{236196}{1687-\dfrac{1327104}{3465-\dfrac{5062500}{6215-\ddots}}}}}}$$
Convergence type $P^+$ with $P=8/3$ and
$C=\G(1/6)^7/(2^{34/3}3^{9/2}\pi^{1/2})$, so that
$$91\z(3)+2\pi^3\sqrt{3}-\dfrac{p(n)}{q(n)}\sim\dfrac{\G(1/6)^7/(2^{34/3}3^{9/2}\pi^{1/2})}{n^{8/3}}\;.$$
$$A=1-(4/3)/n+(137/162)/n^2-(35/486)/n^3-(2460977/8398080)/n^4+\cdots$$
Parametric family for $k\ge0$:
\begin{verbatim}
[()->91*zeta(3)+2*Pi^3*sqrt(3),
     (2*n-1)*(18*n^2-18*n+25+36*k^2+48*k),-324*n^6]
\end{verbatim}
Convergence type $P^+$ with $P=4k+8/3$.
\end{cf}

\smallskip

\begin{cf}\label{5.8.162}{\ }
\begin{verbatim}
[()->91*zeta(3)+2*Pi^3*sqrt(3),
[216,343,4*(3*n-1)*(36*n^2-24*n+31)],[216,-(6*n+1)^6]]
\end{verbatim}
$$91\z(3)+2\pi^3\sqrt{3}=216+\dfrac{216}{343-\dfrac{117649}{2540-\dfrac{4826809}{9056-\dfrac{47045881}{22484-\dfrac{244140625}{45416-\ddots}}}}}$$
Convergence type $P^+$ with $P=2$ and $C=1/2$, so that
$$91\z(3)+2\pi^3\sqrt{3}-\dfrac{p(n)}{q(n)}\sim\dfrac{1/2}{n^2}\;.$$
$$A=1-(4/3)/n+(13/12)/n^2-(14/27)/n^3+\cdots$$
Series:
$$91\z(3)+2\pi^3\sqrt{3}=216\sum_{n\ge1}\dfrac{1}{(6n-5)^3}$$
Parametric family for $k\ge0$:
\begin{verbatim}
[()->91*zeta(3)+2*Pi^3*sqrt(3),
4*(3*n-1)*(36*n^2-24*n+31+72*k*(k+1)),-(6*n+1)^6]
\end{verbatim}
Convergence type $P^+$ with $P=4k+2$.
\end{cf}

\smallskip

\begin{cf}\label{5.8.163}{\ }
\begin{verbatim}
[()->91*zeta(3)+2*Pi^3*sqrt(3),
[[216,343],[(12*n+1)*(108*n^2+18*n+1),(12*n+7)*(108*n^2+126*n+49)]],
[[216,-117649],[-46656*n^6,-(6*n+7)^6]]]
\end{verbatim}
$$91\z(3)+2\pi^3\sqrt{3}=216+\dfrac{216}{343-\dfrac{117649}{1651-\dfrac{46656}{5377-\dfrac{4826809}{11725-\dfrac{2985984}{22723-\dfrac{47045881}{37999-\ddots}}}}}}$$
Convergence type $E$ with $E=(1+\sqrt{2})^4$, $P=0$, and $C=4\pi^3/(1+\sqrt{2})^{14/3}$, so that
$$91\z(3)+2\pi^3\sqrt{3}-\dfrac{p(n)}{q(n)}\sim\dfrac{4\pi^3}{(1+\sqrt{2})^{4n+14/3}}\;.$$
$$A=1-(127d/144)/n+(889d/864+16129/20736)/n^2+\cdots$$
\end{cf}

\smallskip

\begin{cf}\label{5.8.164}{\ }
\begin{verbatim}
[()->117*zeta(3)-2*Pi^3*sqrt(3),[0,(2*n-1)*(9*n^2-9*n+5)],[81,-81*n^6]]
\end{verbatim}
$$117\z(3)-2\pi^3\sqrt{3}=\dfrac{81}{5-\dfrac{81}{69-\dfrac{5184}{295-\dfrac{59049}{791-\dfrac{331776}{1665-\dfrac{1265625}{3025-\ddots}}}}}}$$
Convergence type $P^+$ with $P=2/3$ and $C=81\G(2/3)^{10}/(8\pi^2)$, so that
$$117\z(3)-2\pi^3\sqrt{3}-\dfrac{p(n)}{q(n)}\sim\dfrac{81\G(2/3)^{10}/(8\pi^2)}{n^{2/3}}\;.$$
$$A=1-(1/3)/n+(11/162)/n^2+(8/243)/n^3-(29513/918540)/n^4+\cdots$$
Parametric family for $k\ge0$:
\begin{verbatim}
[()->117*zeta(3)-2*Pi^3*sqrt(3),
     (2*n-1)*(9*n^2-9*n+5+18*k^2+6*k),-81*n^6]
\end{verbatim}
Convergence type $P^+$ with $P=4k+2/3$.
\end{cf}

\smallskip

\begin{cf}\label{5.8.165}{\ }
\begin{verbatim}
[()->117*zeta(3)-2*Pi^3*sqrt(3),
[0,8,(6*n-5)*(9*n^2-15*n+13)],[243,-(3*n-1)^6]]
\end{verbatim}
$$117\z(3)-2\pi^3\sqrt{3}=\dfrac{243}{8-\dfrac{64}{133-\dfrac{15625}{637-\dfrac{262144}{1843-\dfrac{1771561}{4075-\dfrac{7529536}{7657-\ddots}}}}}}$$
Convergence type $P^+$ with $P=2$ and $C=9/2$, so that
$$117\z(3)-2\pi^3\sqrt{3}-\dfrac{p(n)}{q(n)}\sim\dfrac{9/2}{n^2}\;.$$
$$A=1-(1/3)/n-(1/6)/n^2+(4/27)/n^3+(13/162)/n^4+\cdots$$
Series:
$$117\z(3)-2\pi^3\sqrt{3}=243\sum_{n\ge1}\dfrac{1}{(3n-1)^3}$$
Parametric family with $k\ge0$:
\begin{verbatim}
[()->117*zeta(3)-2*Pi^3*sqrt(3),
(6*n-5)*(9*n^2-15*n+13+18*k*(k+1)),-(3*n-1)^6]
\end{verbatim}
Convergence type $P^+$ with $P=4k+2$.
\end{cf}

\smallskip

\begin{cf}\label{5.8.166}{\ }
\begin{verbatim}
[()->117*zeta(3)-2*Pi^3*sqrt(3),
[[0,8],[(6*n-1)*(27*n^2-9*n+1),(6*n+2)*(27*n^2+18*n+4)]],
[[243,-64],[-729*n^6,-(3*n+2)^6]]]
\end{verbatim}
$$117\z(3)-2\pi^3\sqrt{3}=\dfrac{243}{8-\dfrac{64}{95-\dfrac{729}{392-\dfrac{15625}{1001-\dfrac{46656}{2072-\dfrac{262144}{3689-\ddots}}}}}}$$
Convergence type $E$ with $E=(1+\sqrt{2})^4$, $P=0$, and
$C=36\pi^3/(1+\sqrt{2})^{8/3}$, so that
$$117\z(3)-2\pi^3\sqrt{3}-\dfrac{p(n)}{q(n)}\sim\dfrac{36\pi^3}{(1+\sqrt{2})^{4n+8/3}}\;.$$
$$A=1-(103d/144)/n+(103d/216+10609/20736)/n^2+\cdots$$
\end{cf}

\smallskip

\begin{cf}\label{5.8.167}{\ }
\begin{verbatim}
[()->117*zeta(3)+2*Pi^3*sqrt(3),[243,(2*n-1)*(9*n^2-9*n+17)],[81,-81*n^6]]
\end{verbatim}
$$117\z(3)+2\pi^3\sqrt{3}=243+\dfrac{81}{17-\dfrac{81}{105-\dfrac{5184}{355-\dfrac{59049}{875-\dfrac{331776}{1773-\dfrac{1265625}{3157-\ddots}}}}}}$$
Convergence type $P^+$ with $P=10/3$ and $C=\G(1/3)^{10}/(1440\pi^2)$, so that
$$117\z(3)+2\pi^3\sqrt{3}-\dfrac{p(n)}{q(n)}\sim\dfrac{\G(1/3)^{10}/(1440\pi^2)}{n^{10/3}}\;.$$
$$A=1-(5/3)/n+(205/162)/n^2-(40/243)/n^3-(146849/288684)/n^4+\cdots$$
Parametric family for $k\ge0$:
\begin{verbatim}
[()->117*zeta(3)+2*Pi^3*sqrt(3),
     (2*n-1)*(9*n^2-9*n+17+18*k^2+30*k),-81*n^6]
\end{verbatim}
Convergence type $P^+$ with $P=4k+10/3$.
\end{cf}

\smallskip

\begin{cf}\label{5.8.168}{\ }
\begin{verbatim}
[()->117*zeta(3)+2*Pi^3*sqrt(3),
[0,1,(6*n-7)*(9*n^2-21*n+19)],[243,-(3*n-2)^6]]
\end{verbatim}
$$117\z(3)+2\pi^3\sqrt{3}=\dfrac{243}{1-\dfrac{1}{65-\dfrac{4096}{407-\dfrac{117649}{1343-\dfrac{1000000}{3197-\dfrac{4826809}{6293-\ddots}}}}}}$$
Convergence type $P^+$ with $P=2$ and $C=9/2$, so that
$$117\z(3)+2\pi^3\sqrt{3}-\dfrac{p(n)}{q(n)}\sim\dfrac{9/2}{n^2}\;.$$
$$A=1+(1/3)/n-(1/6)/n^2-(4/27)/n^3+(13/162)/n^4+\cdots$$
Series:
$$117\z(3)+2\pi^3\sqrt{3}=243\sum_{n\ge1}\dfrac{1}{(3n-2)^3}$$
Parametric family with $k\ge0$:
\begin{verbatim}
[()->117*zeta(3)+2*Pi^3*sqrt(3),
(6*n-7)*(9*n^2-21*n+19+18*k*(k+1)),-(3*n-2)^6]
\end{verbatim}
Convergence type $P^+$ with $P=4k+2$.
\end{cf}
      
\smallskip

\begin{cf}\label{5.8.169}{\ }
\begin{verbatim}
[()->117*zeta(3)+2*Pi^3*sqrt(3),
[[243,64],[(6*n+1)*(27*n^2+9*n+1),(6*n+4)*(27*n^2+36*n+16)]],
[[243,-4096],[-729*n^6,-(3*n+4)^6]]]
\end{verbatim}
$$117\z(3)+2\pi^3\sqrt{3}=243+\dfrac{243}{64-\dfrac{4096}{259-\dfrac{729}{790-\dfrac{117649}{1651-\dfrac{46656}{3136-\dfrac{1000000}{5149-\ddots}}}}}}$$
Convergence type $E$ with $E=(1+\sqrt{2})^4$, $P=0$, and $C=36\pi^3/(1+\sqrt{2})^{16/3}$, so that
$$117\z(3)+2\pi^3\sqrt{3}-\dfrac{p(n)}{q(n)}\sim\dfrac{36\pi^3}{(1+\sqrt{2})^{4n+16/3}}\;.$$
$$A=1-(103d/144)/n+(103d/108+10609/20736)/n^2+\cdots$$
\end{cf}

\smallskip

\begin{cf}\label{5.8.170}{\ }
\begin{verbatim}
[()->351*zeta(3)-10*Pi^3*sqrt(3),[0,8,9*(9*n^2-15*n+7)],[-972,(3*n-1)^6]]
\end{verbatim}
$$351\z(3)-10\pi^3\sqrt{3}=-\dfrac{972}{8+\dfrac{64}{117+\dfrac{15625}{387+\dfrac{262144}{819+\dfrac{1771561}{1413+\dfrac{7529536}{2169+\ddots}}}}}}$$
Convergence type $P^-$ with $P=3$ and $C=-18$, so that
$$351\z(3)-10\pi^3\sqrt{3}-\dfrac{p(n)}{q(n)}\sim-\dfrac{18}{n^3}\;.$$
$$A=1-(1/2)/n-(4/3)/n^2+(65/54)/n^3+(110/27)/n^4+\cdots$$
Series:
$$351\z(3)-10\pi^3\sqrt{3}=-972\sum_{n\ge1}\dfrac{(-1)^{n+1}}{(3n-1)^3}$$
Parametric family (for a termwise equivalent CF) for $k\ge0$:
\begin{verbatim}
[()->351*zeta(3)-10*Pi^3*sqrt(3),
 9*(2*k+1)*(6*n-5)*(9*n^2-15*n+7),(6*n-6*k-5)*(6*n+6*k+1)*(3*n-1)^6]
\end{verbatim}
Convergence type $P^-$ with $P=6k+3$.
\end{cf}

\smallskip

\begin{cf}\label{5.8.171}{\ }
\begin{verbatim}
[()->351*zeta(3)+10*Pi^3*sqrt(3),[0,1,9*(9*n^2-21*n+13)],[972,(3*n-2)^6]]
\end{verbatim}
$$351\z(3)+10\pi^3\sqrt{3}=\dfrac{972}{1+\dfrac{1}{63+\dfrac{4096}{279+\dfrac{117649}{657+\dfrac{1000000}{1197+\dfrac{4826809}{1899+\ddots}}}}}}$$
Convergence type $P^-$ with $P=3$ and $C=18$, so that
$$351\z(3)+10\pi^3\sqrt{3}-\dfrac{p(n)}{q(n)}\sim\dfrac{18}{n^3}\;.$$
$$A=1+(1/2)/n-(4/3)/n^2-(65/54)/n^3+(110/27)/n^4+\cdots$$
Series:
$$351\z(3)+10\pi^3\sqrt{3}=972\sum_{n\ge1}\dfrac{(-1)^{n+1}}{(3n-2)^3}$$
Parametric family (for a termwise equivalent CF) for $k\ge0$:
\begin{verbatim}
[()->351*zeta(3)+10*Pi^3*sqrt(3),
 9*(2*k+1)*(6*n-7)*(9*n^2-21*n+13),(6*n-6*k-7)*(6*n+6*k-1)*(3*n-2)^6]
\end{verbatim}
Convergence type $P^-$ with $P=6k+3$.
\end{cf}

\smallskip

For the CFs that follow, note that $L(\chi_{-8},3)=3\pi^3/(64\sqrt{2})$.

\smallskip

\begin{cf}\label{5.8.172}{\ }
\begin{verbatim}
[()->lfun(-8,3)-lfun(8,3),[0,27,192*n^2-288*n+124],[2,(4*n-1)^6]]
\end{verbatim}
$$L(\chi_{-8},3)-L(\chi_8,3)=\dfrac{2}{27+\dfrac{729}{316+\dfrac{117649}{988+\dfrac{1771561}{2044+\dfrac{11390625}{3484+\dfrac{47045881}{5308+\ddots}}}}}}$$
Convergence type $P^-$ with $P=3$ and $C=1/64$, so that
$$L(\chi_{-8},3)-L(\chi_8,3)-\dfrac{p(n)}{q(n)}\sim\dfrac{1/64}{n^3}\;.$$
$$A=1-(3/4)/n-(9/8)/n^2+(55/32)/n^3+(855/256)/n^4+\cdots$$
Series:
$$L(\chi_{-8},3)-L(\chi_8,3)=2\sum_{n\ge1}\dfrac{(-1)^{n+1}}{(4n-1)^3}$$
Parametric family (for a termwise equivalent CF) for $k\ge0$:
\begin{verbatim}
[()->lfun(-8,3)-lfun(8,3),
(2*k+1)*(4*n-3)*(192*n^2-288*n+124),(4*n-4*k-3)*(4*n+4*k+1)*(4*n-1)^6]
\end{verbatim}
Convergence type $P^-$ with $P=6k+3$.
\end{cf}

\smallskip

\begin{cf}\label{5.8.173}{\ }
\begin{verbatim}
[()->lfun(-8,3)+lfun(8,3),[0,1,192*n^2-480*n+316],[2,(4*n-3)^6]]
\end{verbatim}
$$L(\chi_{-8},3)+L(\chi_8,3)=\dfrac{2}{1+\dfrac{1}{124+\dfrac{15625}{604+\dfrac{531441}{1468+\dfrac{4826809}{2716+\dfrac{24137569}{4348+\ddots}}}}}}$$
Convergence type $P^-$ with $P=3$ and $C=1/64$, so that
$$L(\chi_{-8},3)+L(\chi_8,3)-\dfrac{p(n)}{q(n)}\sim\dfrac{1/64}{n^3}\;.$$
$$A=1+(3/4)/n-(9/8)/n^2-(55/32)/n^3+(855/256)/n^4+\cdots$$
Series:
$$L(\chi_{-8},3)+L(\chi_8,3)=2\sum_{n\ge1}\dfrac{(-1)^{n+1}}{(4n-3)^3}$$
Parametric family (for a termwise equivalent CF) for $k\ge0$:
\begin{verbatim}
[()->lfun(-8,3)+lfun(8,3),
(2*k+1)*(4*n-5)*(192*n^2-480*n+316),(4*n-4*k-5)*(4*n+4*k-1)*(4*n-3)^6]
\end{verbatim}
Convergence type $P^-$ with $P=6k+3$.
\end{cf}

\smallskip

\begin{cf}\label{5.8.174}{\ }
\begin{verbatim}
[()->3*zeta(3)-2*log(2)^3,[0,7*n^3-n^2-3*n+1],[12,4*n^5*(2*n-1)]]
\end{verbatim}
$$3\z(3)-2\log(2)^3=\dfrac{12}{4+\dfrac{4}{47+\dfrac{384}{172+\dfrac{4860}{421+\dfrac{28672}{836+\dfrac{112500}{1459+\ddots}}}}}}$$
Convergence type $E$ with $E=-8$, $P=5/2$, and $C=4\sqrt{\pi}/3$, so that
$$3\z(3)-2\log(2)^3-\dfrac{p(n)}{q(n)}\sim(-1)^n\dfrac{4\sqrt{\pi}/3}{2^{3n}n^{5/2}}\;.$$
$$A=1-(151/72)/n+(27409/10368)/n^2-\cdots$$
Series:
$$3\z(3)-2\log(2)^3=3\sum_{n\ge0}(-1)^n\dfrac{n!}{(n+1)^2(3/2)_n}2^{-3n}$$
\end{cf}

\smallskip

\begin{cf}\label{5.8.175}{\ }
\begin{verbatim}
[()->Pi^3+12*Pi*log(2)^2,[48,54,32*n^4-8*n^3+36*n^2-6*n+1],
                         [48,-2*n*(2*n+1)^7]]
\end{verbatim}
$$\pi^3+12\pi\log(2)^2=48+\dfrac{48}{54-\dfrac{4374}{581-\dfrac{312500}{2683-\dfrac{4941258}{8233-\dfrac{38263752}{19871-\ddots}}}}}$$
Convergence type $P^+$ with $P=5/2$ and $C=12/(5\sqrt{\pi})$, so that
$$\pi^3+12\pi\log(2)^2-\dfrac{p(n)}{q(n)}\sim\dfrac{12/(5\sqrt{\pi})}{n^{5/2}}$$
$$A=1-(135/56)/n+(4265/1152)/n^2-(48425/11264)/n^3+\cdots$$
Series:
$$\pi^3+12\pi\log(2)^2=48+24\sum_{n\ge0}\dfrac{(3/2)_n}{(2n+3)^3(n+1)!}$$
\end{cf}

\smallskip

\begin{cf}\label{5.8.176}{\ }
\begin{verbatim}
[()->(21/4)*zeta(3)+log(2)^3-Pi^2*log(2)/2,
                            [0,3*n^3-3*n^2+3*n-1],[6,-2*n^6]]
\end{verbatim}
$$\dfrac{21}{4}\z(3)+\log^3(2)-\dfrac{\pi^2}{2}\log(2)=\dfrac{6}{2-\dfrac{2}{17-\dfrac{128}{62-\dfrac{1458}{155-\dfrac{8192}{314-\dfrac{31250}{557-\ddots}}}}}}$$
Convergence type $E$ with $E=2$, $P=3$, and $C=6$, so that
$$\dfrac{21}{4}\z(3)+\log^3(2)-\dfrac{\pi^2}{2}\log(2)-\dfrac{p(n)}{q(n)}\sim\dfrac{3}{2^{n-1}n^3}\;.$$
$$A=1-6/n+36/n^2-260/n^3+2250/n^4-22722/n^5+262248/n^6+\cdots$$
Series:
$$\dfrac{21}{4}\z(3)+\log^3(2)-\dfrac{\pi^2}{2}\log(2)=3\sum_{n\ge0}\dfrac{2^{-n}}{(n+1)^3}$$
\end{cf}

\medskip

\section{Miscellaneous Constants Related to $L$-Values}

\smallskip

\begin{cf}\label{5.9.1}{\ }
\begin{verbatim}
[()->Catalan/Pi,[1/4,12,16*n^3-8*n^2+6*n-1],[1/4,-4*n^2*(2*n+1)^4]]
\end{verbatim}
$$\dfrac{G}{\pi}=1/4+\dfrac{1/4}{12-\dfrac{324}{107-\dfrac{10000}{377-\dfrac{86436}{919-\dfrac{419904}{1829-\dfrac{1464100}{3203-\ddots}}}}}}$$
Convergence type $P^+$ with $P=1$ and $C=1/(8\pi)$, so that
$$\dfrac{G}{\pi}-\dfrac{p(n)}{q(n)}\sim\dfrac{1/(8\pi)}{n}\;.$$
$$A=1-(7/8)/n+(65/96)/n^2-(225/512)/n^3+(6721/30720)/n^4-\cdots$$
Series:
$$\dfrac{G}{\pi}=\dfrac{1}{4}+\dfrac{1}{16}\sum_{n\ge0}\dfrac{(3/2)_n^2}{(2n+3)(n+1)!^2}$$
\end{cf}

\smallskip

\begin{cf}\label{5.9.2}{\ }
\begin{verbatim}
[()->(1+2*Catalan)/Pi,[1,12,16*n^3-16*n^2+14*n-3],
                      [-1,-4*n^2*(2*n-1)*(2*n+1)^3]]
\end{verbatim}
$$\dfrac{1+2G}{\pi}=1-\dfrac{1}{12-\dfrac{108}{89-\dfrac{6000}{327-\dfrac{61740}{821-\dfrac{326592}{1667-\dfrac{1197900}{2961-\ddots}}}}}}$$
Convergence type $P^+$ with $P=2$ and $C=-1/(8\pi)$, so that
$$\dfrac{1+2G}{\pi}-\dfrac{p(n)}{q(n)}\sim-\dfrac{1/(8\pi)}{n^2}\;.$$
$$A=1-(7/6)/n+(57/64)/n^2-(451/960)/n^3+(297/2048)/n^4-\cdots$$
Series:
$$\dfrac{1+2G}{\pi}=1-\dfrac{1}{4}\sum_{n\ge0}\dfrac{(2n+1)(1/2)_n^2}{(2n+3)(n+1)!^2}$$
\end{cf}

Similar CFs exist for infinitely many linear combinations of $1/\pi$ and
$G/\pi$.

\smallskip

\begin{cf}\label{5.9.3}{\ }
\begin{verbatim}
[()->log(2)-2*Catalan/Pi,[1/16,32,8*n^3+16*n^2+13*n+4],
                             [9/16,-4*(n+1)^4*(2*n+3)^2]]
\end{verbatim}
$$\log(2)-\dfrac{2G}{\pi}=\dfrac{1}{16}+\dfrac{9/16}{32-\dfrac{1600}{158-\dfrac{15876}{403-\dfrac{82944}{824-\dfrac{302500}{1469-\dfrac{876096}{2386-\ddots}}}}}}$$
Convergence type $P^+$ with $P=1$ and $C=1/(4\pi)$, so that
$$\log(2)-\dfrac{2G}{\pi}-\dfrac{p(n)}{q(n)}\sim\dfrac{1/(4\pi)}{n}\;.$$
$$A=1-(13/8)/n+(245/96)/n^2-(1975/512)/n^3+\cdots$$
Series:
$$\log(2)-\dfrac{2G}{\pi}=\dfrac{1}{16}+\dfrac{9}{64}\sum_{n\ge0}\dfrac{(5/2)_n^2}{(n+2)(n+2)!^2}$$
\end{cf}

\medskip

\section{Constants Related to Powers of $e=\exp(1)$}

\smallskip

In what follows, the variables $k$, $l$, and $m$ are implicitly assumed to be
strictly positive integers.

\smallskip

\begin{cf}\label{5.10.1}{\ }
\begin{verbatim}
[()->exp(1),[1,1,n+1],[1,-n]]
\end{verbatim}
$$e=1+\dfrac{1}{1-\dfrac{1}{3-\dfrac{2}{4-\dfrac{3}{5-\dfrac{4}{6-\dfrac{5}{7-\ddots}}}}}}$$
Convergence type $F^1$ with $E=1$, $P=1$, and $C=1$, so that
$$e-\dfrac{p(n)}{q(n)}\sim\dfrac{1}{n!n}\;.$$
$$A=1-1/n^2+1/n^3+2/n^4-9/n^5+9/n^6+50/n^7+\cdots$$
Series:
$$e=1+\sum_{n\ge0}\dfrac{1}{(n+1)!}$$
\end{cf}

\smallskip

\begin{cf}\label{5.10.2}{\ }
\begin{verbatim}
[()->exp(1),[5/2,n+1],[1/2,n]]
\end{verbatim}
$$e=5/2+\dfrac{1/2}{2+\dfrac{1}{3+\dfrac{2}{4+\dfrac{3}{5+\dfrac{4}{6+\dfrac{5}{7+\ddots}}}}}}$$
Convergence type $F^1$ with $E=-1$, $P=5$, and $C=2e^2$, so that
$$e-\dfrac{p(n)}{q(n)}\sim(-1)^n\dfrac{2e^2}{n!n^5}\;.$$
$$A=1-14/n+126/n^2-945/n^3+6511/n^4-43300/n^5+\cdots$$
Series:
$$e^{-1}=2\sum_{n\ge0}\dfrac{(-1)^n}{(n^2+3n+1)(n^2+5n+5)n!}$$
Parametric family up to the trivial change $n\mapsto n+j$ with $u\ge0$:
\begin{verbatim}
[()->exp(1),n+u,n]
\end{verbatim}
Convergence type $F^1$ with $E=-1$ and $P=2u+3$.
\end{cf}

\smallskip

\begin{cf}\label{5.10.3}{\ }
\begin{verbatim}
[()->exp(1),[2,n+1],[n+2]]
\end{verbatim}
$$e=2+\dfrac{2}{2+\dfrac{3}{3+\dfrac{4}{4+\dfrac{5}{5+\dfrac{6}{6+\dfrac{7}{7+\ddots}}}}}}$$
Convergence type $F^1$ with $E=-1$, $P=3$, and $C=e^2$, so that
$$e-\dfrac{p(n)}{q(n)}\sim(-1)^n\dfrac{e^2}{n!\cdot n^3}$$
$$A=1-7/n+36/n^2-171/n^3+813/n^4-4012/n^5+20891/n^6+\cdots$$
Series:
$$e^{-1}=\sum_{n\ge0}\dfrac{(-1)^n}{(n+2)!}$$
\end{cf}

\smallskip

\begin{cf}\label{5.10.4}{\ }
\begin{verbatim}
[()->exp(1),[n+3],[-(n+1)]]
\end{verbatim}
$$e=3-\dfrac{1}{4-\dfrac{2}{5-\dfrac{3}{6-\dfrac{4}{7-\dfrac{5}{8-\dfrac{6}{9-\ddots}}}}}}$$
Convergence type $F^1$ with $E=1$, $P=4$, and $C=-1$, so that
$$e-\dfrac{p(n)}{q(n)}\sim-\dfrac{1}{n!\cdot n^4}\;.$$
$$A=1-5/n+13/n^2-14/n^3-44/n^4+260/n^5-405/n^6-2189/n^7+\cdots$$
Series:
$$e=3-\sum_{n\ge0}\dfrac{1}{(n+1)(n+2)(n+2)!}$$
Parametric family up to the trivial change $n\mapsto n+j$ with $u\ge0$:
\begin{verbatim}
[()->exp(1),n+u,-n]
\end{verbatim}
Convergence type $F^1$ with $E=1$ and $P=2u-1$.
\end{cf}

\smallskip

\begin{cf}\label{5.10.5}{\ }
\begin{verbatim}
[()->exp(1),[5/2,n^2+2*n],[3/4,n^2*(n+3)]]
\end{verbatim}
$$e=5/2+\dfrac{3/4}{3+\dfrac{4}{8+\dfrac{20}{15+\dfrac{54}{24+\dfrac{112}{35+\dfrac{200}{48+\ddots}}}}}}$$
Convergence type $F^1$ with $E=-1$, $P=5$, and $C=2e^2$, so that
$$e-\dfrac{p(n)}{q(n)}\sim(-1)^n\dfrac{2e^2}{n!n^5}\;.$$
$$A=1-14/n+126/n^2-945/n^3+6511/n^4-43300/n^5+286160/n^6+\cdots$$
Series:
$$e^{-1}=2\sum_{n\ge0}\dfrac{(-1)^n}{(n^2+3n+1)(n^2+5n+5)n!}$$
\end{cf}

\smallskip

\begin{cf}\label{5.10.6}{\ }
\begin{verbatim}
[()->exp(1),[11/4,n^2+3*n+1],[-1/6,n^2*(n+2)]]
\end{verbatim}
$$e=11/4-\dfrac{1/6}{5+\dfrac{3}{11+\dfrac{16}{19+\dfrac{45}{29+\dfrac{96}{41+\dfrac{175}{55+\ddots}}}}}}$$
Convergence type $F^1$ with $E=-1$, $P=8$, and $C=-12e^2$, so that
$$e-\dfrac{p(n)}{q(n)}\sim(-1)^{n+1}\dfrac{12e^2}{n!n^8}\;.$$
$$A=1-28/n+464/n^2-5982/n^3+66729/n^4-680622/n^5+6562330/n^6+\cdots$$
Series:
$$e^{-1}=-12\sum_{n\ge0}(-1)^n\dfrac{n+1}{(n^4+8n^3+17n^2+8n-1)(n^4+12n^3+47n^2+70n+33)n!}$$
\end{cf}

\smallskip

\begin{verbatim}
S1(k)=sum(j=0,k-2,(-1)^(j+1)*binomial(k,j+2)/j!);
\end{verbatim}

\smallskip

\begin{cf}\label{5.10.7}{\ }
\begin{verbatim}
[k->exp(1)*k/S1(k),[1,n^2-n-k],[k+1,n^2*(n+k+1)]]
\end{verbatim}
$$e\cdot k/S_1(k)=1+\dfrac{k+1}{-k+\dfrac{k+2}{-k+2+\dfrac{4k+12}{-k+6+\dfrac{9k+36}{-k+12+\dfrac{16k+80}{-k+20+\dfrac{25k+150}{-k+30+\ddots}}}}}}$$
Convergence type $F^1$ with $E=-1$, $P=1-k$, and $C=e^2\cdot k/((k-1)!S1(k)^2)$, so that
$$e\cdot k/S_1(k)-\dfrac{p(n)}{q(n)}\sim(-1)^n\dfrac{e^2\cdot k/((k-1)!S(k)^2)}{n!n^{1-k}}\;.$$
$$A=1+((k-1)(k+4)/2)/n+((k-1)(3k^3+17k^2+2k-120)/24)/n^2+\cdots$$
\end{cf}

\smallskip

\begin{verbatim}
S2(k)=sum(j=0,k,binomial(k,j)/j!);
\end{verbatim}

\smallskip

\begin{cf}\label{5.10.8}{\ }
\begin{verbatim}
[k->exp(1)*S2(k),[2+k,4,(n+1)*(n+2)+k],[(k+1)*(k+2),-(n+1)^2*(n+2+k)]]
\end{verbatim}
$$eS_2(k)=k+2+\dfrac{k^2+3k+2}{4-\dfrac{4k+12}{k+12-\dfrac{9k+36}{k+20-\dfrac{16k+80}{k+30-\dfrac{25k+150}{k+42-\dfrac{36k+252}{k+56-\ddots}}}}}}$$
Convergence type $F^1$ with $E=1$, $P=2-k$, and $C=1/k!$, so that
$$eS_2(k)-\dfrac{p(n)}{q(n)}\sim\dfrac{1/k!}{n!n^{2-k}}\;.$$
$$A=1+((k^2+5k-4)/2)/n+((3k^4+26k^3+21k^2-170k+48)/24)/n^2+\cdots$$
Series:
$$eS_2(k)=k+2+\dfrac{1}{k!}\sum_{n\ge0}\dfrac{(n+k+2)!}{(n+2)!^2}$$
\end{cf}

More generally, there are infinitely many CF for $e=\exp(1)$ with
$a(n)=a_2n^2+a_1n+a_0$ and $b(n)=\pm a_2(n+u)(n+v)(n+w)$ for $n$ sufficiently
large.

\smallskip

\begin{cf}\label{5.10.9}{\ }
\begin{verbatim}
[()->exp(1),[1,1,4*n-2],[2,1]]
\end{verbatim}
$$e=1+\dfrac{2}{1+\dfrac{1}{6+\dfrac{1}{10+\dfrac{1}{14+\dfrac{1}{18+\dfrac{1}{22+\ddots}}}}}}$$
Convergence type $F^2$ with $E=-16$, $P=0$, and $C=\pi e/2$, so that
$$e-\dfrac{p(n)}{q(n)}\sim(-1)^n\dfrac{\pi e}{n!^22^{4n+1}}\;.$$
$$A=1-(1/8)/n+(9/128)/n^2+(283/3072)/n^3-(18847/98304)/n^4+\cdots$$
\end{cf}

\smallskip

\begin{cf}\label{5.10.10}{\ }
\begin{verbatim}
[()->exp(1),[2,3,4*n^2+2*n-1],[-(2*n-1)*(2*n+2)]]
\end{verbatim}
$$e=2+\dfrac{2}{3-\dfrac{4}{19-\dfrac{18}{41-\dfrac{40}{71-\dfrac{70}{109-\dfrac{108}{155-\ddots}}}}}}$$
Convergence type $F^2$ with $E=4$, $P=5/2$, and $C=e^2\sqrt{\pi}/8$, so
that
$$e-\dfrac{p(n)}{q(n)}\sim\dfrac{e^2\sqrt{\pi}}{n!^22^{2n+3}n^{5/2}}\;.$$
$$A=1-(27/8)/n+(1097/128)/n^2-(20769/1024)/n^3+(1580315/32768)/n^4+\cdots$$
Series:
\begin{align*}e^{-1}&=1-\dfrac{1}{2}\sum_{n\ge0}\dfrac{1}{(1/2)_n(n+1)!}2^{-2n}\\
&=1-\dfrac{1}{2}\sum_{n\ge0}\dfrac{1}{(n+1)(2n)!}\end{align*}
\end{cf}

\smallskip

\begin{cf}\label{5.10.11}{\ }
\begin{verbatim}
[()->exp(1),[3,4*n^2+6*n+1],[-3,-2*n*(2*n+3)]]
\end{verbatim}
$$e=3-\dfrac{3}{11-\dfrac{10}{29-\dfrac{28}{55-\dfrac{54}{89-\dfrac{88}{131-\dfrac{130}{181-\ddots}}}}}}$$
Convergence type $F^2$ with $E=4$, $P=7/2$, and $C=-e^2\sqrt{\pi}/16$, so
that
$$e-\dfrac{p(n)}{q(n)}\sim-\dfrac{e^2\sqrt{\pi}}{n!^22^{2n+4}n^{7/2}}\;.$$
$$A=1-(43/8)/n+(2505/128)/n^2-(62993/1024)/n^3+(5982107/32768)/n^4+\cdots$$
Series:
\begin{align*}e^{-1}&=\dfrac{1}{3}\sum_{n\ge0}\dfrac{1}{(5/2)_nn!}2^{-2n}\\
&=\sum_{n\ge0}\dfrac{1}{(2n+3)(2n+1)!}\end{align*}
\end{cf}

\smallskip

\begin{verbatim}
S3(k)=k!*sum(j=0,k,(-1)^j*(k-j+1)/j!);
\end{verbatim}

\smallskip

\begin{cf}\label{5.10.12}{\ }
\begin{verbatim}
[k->((S3(k-2)+(-1)^k)*exp(1)-k*(k-2)!)/(S3(k-2)*exp(1)-k*(k-2)!),
[4*n^2+(4*k+2)*n+k^2+k-1],[-(4*n^2+(4*k+2)*n+k^2+k-2)]]
\end{verbatim}
$$\dfrac{(S_3(k-2)+(-1)^k)e-k(k-2)!}{S_3(k-2)e-k(k-2)!}=k^2+k-1-\dfrac{k^2+k-2}{k^2+5k+5-\dfrac{k^2+5k+4}{k^2+9k+19-\ddots}}$$
Convergence type $F^2$ with $E=4$, $P=k+5/2$, and
$C=-e^2\sqrt{\pi}k(k-2)!/(2^{k+3}(S3(k-2)e-k(k-2)!)^2)$, so that
$$\dfrac{(S_3(k-2)+(-1)^k)e-k(k-2)!}{S_3(k-2)e-k(k-2)!}-\dfrac{p(n)}{q(n)}\sim-\dfrac{e^2\sqrt{\pi}k(k-2)!/(2^{k+3}(S3(k-2)e-k(k-2)!)^2)}{n!^22^{2n}n^{k+5/2}}\;.$$
$$A=1-((2k^2+14k+27)/8)/n+((12k^4+184k^3+1080k^2+2948k+3291)/384)/n^2+\cdots$$
\end{cf}

\smallskip

\begin{cf}\label{5.10.13}{\ }
\begin{verbatim}
[()->(exp(1)-1)/(exp(1)+1),[0,4*n-2],[1]]
\end{verbatim}
$$\dfrac{e-1}{e+1}=\dfrac{1}{2+\dfrac{1}{6+\dfrac{1}{10+\dfrac{1}{14+\dfrac{1}{18+\dfrac{1}{22+\ddots}}}}}}$$
Convergence type $F^2$ with $E=-16$, $P=0$, and $C=\pi/(e+2+1/e)$,
so that
$$\dfrac{e-1}{e+1}-\dfrac{p(n)}{q(n)}\sim(-1)^n\dfrac{\pi/(e+2+1/e)}{n!^22^{4n}}\;.$$
$$A=1-(1/8)/n+(9/128)/n^2+(283/3072)/n^3-(18847/98304)/n^4+\cdots$$
\end{cf}

This is of course the same continued fraction as \ref{5.10.9}, presented
differently.

\smallskip

\begin{cf}\label{5.10.14}{\ }
\begin{verbatim}
[()->exp(2),[7,[3*n-1,1,1,3*n,12*n+6]],[1]]
\end{verbatim}
$$e^2=7+\dfrac{1}{2+\dfrac{1}{1+\dfrac{1}{1+\dfrac{1}{3+\dfrac{1}{18+\dfrac{1}{5+\dfrac{1}{1+\dfrac{1}{1+\dfrac{1}{6+\dfrac{1}{30+\ddots}}}}}}}}}}\;.$$
\end{cf}


\smallskip

\begin{cf}\label{5.10.15}{\ }
\begin{verbatim}
[()->exp(2),[1,1,n+2],[2,-2*n]]
\end{verbatim}
$$e^2=1+\dfrac{2}{1-\dfrac{2}{4-\dfrac{4}{5-\dfrac{6}{6-\dfrac{8}{7-\dfrac{10}{8-\ddots}}}}}}$$
Convergence type $F^1$ with $E=1/2$, $P=1$, and $C=2$, so that
$$e^2-\dfrac{p(n)}{q(n)}\sim\dfrac{1}{n!(1/2)^{n+1}n}\;.$$
$$A=1+1/n-1/n^2-3/n^3+7/n^4+13/n^5-89/n^6+\cdots$$
Series:
$$e^2=1+2\sum_{n\ge0}\dfrac{2^n}{(n+1)!}$$
\end{cf}
            
\smallskip

\begin{cf}\label{5.10.16}{\ }
\begin{verbatim}
[()->exp(2),[5,n+1],[8,2*n+4]]
[()->exp(2),[5,1],[4,2/(n+1)]]
\end{verbatim}
$$e^2=5+\dfrac{8}{2+\dfrac{6}{3+\dfrac{8}{4+\dfrac{10}{5+\dfrac{12}{6+\dfrac{14}{7+\ddots}}}}}}=5+\dfrac{4}{1+\dfrac{1}{1+\dfrac{2/3}{1+\dfrac{1/2}{1+\dfrac{2/5}{1+\dfrac{1/3}{1+\ddots}}}}}}$$
Convergence type $F^1$ with $E=-1/2$, $P=5$, and $C=16e^4$, so that
$$e^2-\dfrac{p(n)}{q(n)}\sim(-1)^n\dfrac{e^4}{n!2^{-n-4}n^{5/2}}\;.$$
$$A=1-19/n+232/n^2-2354/n^3+21839/n^4-194377/n^5+\cdots$$
Series:
$$e^{-2}=8\sum_{n\ge0}\dfrac{n+3}{(n+5)!}(-2)^n$$
\end{cf}

\smallskip

\begin{cf}\label{5.10.17}{\ }
\begin{verbatim}
[()->exp(2),[7,2*n+3],[2,1]]
\end{verbatim}
$$e^2=7+\dfrac{2}{5+\dfrac{1}{7+\dfrac{1}{9+\dfrac{1}{11+\dfrac{1}{13+\dfrac{1}{15+\ddots}}}}}}$$
Convergence type $F^2$ with $E=-4$, $P=4$, and $C=\pi e^2/16$, so that
$$e^2-\dfrac{p(n)}{q(n)}\sim(-1)^n\dfrac{\pi e^2}{n!^22^{2n+4}n^4}\;.$$
$$A=1-(23/4)/n+(669/32)/n^2-(23359/384)/n^3+(931529/6144)/n^4+\cdots$$
\end{cf}

\smallskip

\begin{cf}\label{5.10.18}{\ }
\begin{verbatim}
[()->exp(3),[1,1,n+3],[3,-3*n]]
\end{verbatim}
$$e^3=1+\dfrac{3}{1-\dfrac{3}{5-\dfrac{6}{6-\dfrac{9}{7-\dfrac{12}{8-\dfrac{15}{9-\ddots}}}}}}$$
Convergence type $F^1$ with $E=1/3$, $P=1$, and $C=3$, so that
$$e^3-\dfrac{p(n)}{q(n)}\sim\dfrac{1}{n!(1/3)^{n+1}n}$$
$$A=1+2/n+1/n^2-7/n^3-8/n^4+65/n^5+37/n^6+\cdots$$
Series:
$$e^3=1+3\sum_{n\ge0}\dfrac{3^n}{(n+1)!}$$
\end{cf}
          
\smallskip

\begin{cf}\label{5.10.19}{\ }
\begin{verbatim}
[()->exp(3),[13,n+1],[27,3*n+6]]
[()->exp(3),[13,1],[27/2,3/(n+1)]]
\end{verbatim}
$$e^3=13+\dfrac{27}{2+\dfrac{9}{3+\dfrac{12}{4+\dfrac{15}{5+\dfrac{18}{6+\dfrac{21}{7+\ddots}}}}}}=13+\dfrac{27/2}{1+\dfrac{3/2}{1+\dfrac{1}{1+\dfrac{3/4}{1+\dfrac{3/5}{1+\dfrac{1/2}{1+\ddots}}}}}}$$
Convergence type $F^1$ with $E=-1/3$, $P=7$, and $C=486e^6$, so that
$$e^3-\dfrac{p(n)}{q(n)}\sim(-1)^n\dfrac{2e^6}{n!(1/3)^{n+5}n^7}\;.$$
$$A=1-37/n+825/n^2-14508/n^3+223077/n^4-3161955/n^5+\cdots$$
Series:
$$e^{-3}=162\sum_{n\ge0}\dfrac{1}{(n^2+11n+27)(n^2+13n+39)(n+2)!}(-3)^n$$
\end{cf}

\smallskip

\begin{cf}\label{5.10.20}{\ }
\begin{verbatim}
[()->exp(3),[13,7,4*n+6],[54,9]]
\end{verbatim}
$$e^3=13+\dfrac{54}{7+\dfrac{9}{14+\dfrac{9}{18+\dfrac{9}{22+\dfrac{9}{26+\dfrac{9}{30+\ddots}}}}}}$$
Convergence type $F^2$ with $E=-16/9$, $P=4$, and $C=243\pi e^3/512$, so that
$$e^3-\dfrac{p(n)}{q(n)}\sim(-1)^n\dfrac{2\pi e^3}{n!^2(4/3)^{2n+5}n^4}\;.$$
$$A=1-(41/8)/n+(2041/128)/n^2-(37439/1024)/n^3+(1862763/32768)/n^4+\cdots$$
\end{cf}

Equivalently, we have
\begin{verbatim}
[()->9*(exp(3)-1)/(exp(3)+5),[4*n+6],[9]]
\end{verbatim}

\smallskip

\begin{cf}\label{5.10.21}{\ }
\begin{verbatim}
[(k,l)->(exp(2*k/l)-1)/(exp(2*k/l)+1),[0,(2*n-1)*l],[k,k^2]]
\end{verbatim}
$$\dfrac{e^{2k/l}-1}{e^{2k/l}+1}=\dfrac{k}{l+\dfrac{k^2}{3 l+\dfrac{k^2}{5 l+\dfrac{k^2}{7 l+\dfrac{k^2}{9 l+\dfrac{k^2}{11 l+\ddots}}}}}}$$
Convergence type $F^2$ with $E=-4l^2/k^2$, $P=0$, and $C=2\pi k/l/(e^{k/l}+e^{-k/l})^2$, so that
$$\dfrac{e^{2k/l}-1}{e^{2k/l}+1}-\dfrac{p(n)}{q(n)}\sim(-1)^n\dfrac{4\pi/(e^{k/l}+e^{-k/l})^2}{n!^2(2l/k)^{2n+1}}\;.$$
\end{cf}

\smallskip

\begin{cf}\label{5.10.22}{\ }
\begin{verbatim}
[k->exp(1/k),[[1,(2*n+1)*k-1,1]],[1]]
\end{verbatim}
$$e^{1/k}=1+\dfrac{1}{k-1+\dfrac{1}{1+\dfrac{1}{1+\dfrac{1}{3k-1+\dfrac{1}{1+\dfrac{1}{1+\dfrac{1}{5k-1+\ddots}}}}}}}\;.$$
Convergence type $F^{2/3}$ with $E=-2^{2/3}$ and $P=0$. Since it has period
$3$, $C$ has no meaning, and
$$e^{1/k}-\dfrac{p(n)}{q(n)}\approx(-1)^n\dfrac{C}{n!^{2/3}2^{2n/3}}\;.$$
\end{cf}

\smallskip

\begin{cf}\label{5.10.23}{\ }
\begin{verbatim}
[k->exp(2/k),[1,k-1,k*(2*n-1)],[2,1]]
\end{verbatim}
$$e^{2/k}=1+\dfrac{2}{k-1+\dfrac{1}{3k+\dfrac{1}{5k+\dfrac{1}{7k+\dfrac{1}{9k+\dfrac{1}{11k+\ddots}}}}}}$$
Convergence type $F^2$ with $E=-4k^2$, $P=0$, and $C=(\pi/k)e^{2/k}$, so that
$$e^{2/k}-\dfrac{p(n)}{q(n)}\sim(-1)^n\dfrac{\pi e^{2/k}/k}{n!^2(2k)^{2n}}\;.$$
$$A=1-((k^2-2)/(4k^2))/n+((k^2-2)(5k^2-2)/(32k^4))/n^2+\cdots$$
\end{cf}

\smallskip

\begin{cf}\label{5.10.24}{\ }
\begin{verbatim}
[k->k*exp(1/k),[k+1,[2*k-1,2*n,1]],[1]]
\end{verbatim}
$$ke^{1/k}=k+1+\dfrac{1}{2k-1+\dfrac{1}{2+\dfrac{1}{1+\dfrac{1}{2k-1+\dfrac{1}{4+\dfrac{1}{1+\dfrac{1}{2k-1+\ddots}}}}}}}\;.$$
\end{cf}


\smallskip

\begin{cf}\label{5.10.25}{\ }
\begin{verbatim}
[k->exp(k),[sum(j=0,k,k^j/j!),n+1],[2*k^k/(k-1)!,k*(n+2)]]
[k->exp(k),[sum(j=0,k,k^j/j!),1],[k^k/(k-1)!,k/(n+1)]]
\end{verbatim}
$$e^k=\sum_{j=0}^k\dfrac{k^j}{j!}+\dfrac{2k^k/(k-1)!}{2+\dfrac{3k}{3+\dfrac{4k}{4+\dfrac{5k}{5+\dfrac{6k}{6+\dfrac{7k}{7+\ddots}}}}}}=\sum_{j=0}^k\dfrac{k^j}{j!}+\dfrac{k^k/(k-1)!}{1+\dfrac{k/2}{1+\dfrac{k/3}{1+\dfrac{k/4}{1+\dfrac{k/5}{1+\dfrac{k/6}{1+\ddots}}}}}}$$
Convergence type $F^1$ with $E=-1/k$, $P=2k+1$, and $C=k^{k+2}(k-1)!e^{2k}$,
so that
$$e^k-\dfrac{p(n)}{q(n)}\sim(-1)^n\dfrac{(k-1)!e^{2k}}{n!(1/k)^{n+k+2}n^{2k+1}}\;.$$
$$A=1-(3k^2+3k+1)/n+(27k^4+74k^3+72k^2+37k+6)/n^2+\cdots$$
\end{cf}
For example, for $k=1$ this gives \ref{5.10.3}, for $k=2$ this gives
\ref{5.10.16}, and for $k=3$ this gives \ref{5.10.19}.

\smallskip

\begin{cf}\label{5.10.26}{\ }
\begin{verbatim}
[()->exp(1)*eint1(1),[2*n],[1,-n^2]]
\end{verbatim}
$$\int_0^\infty\dfrac{e^{-t}}{1+t}\,dt=e\cdot\eintone(1)=\dfrac{1}{2-\dfrac{1}{4-\dfrac{4}{6-\dfrac{9}{8-\dfrac{16}{10-\dfrac{25}{12-\ddots}}}}}}$$
Convergence type $D^+$ with $D=16$ and $C=2\pi\cdot e$, so that
$$e\cdot\eintone(1)-\dfrac{p(n)}{q(n)}\sim\dfrac{2\pi\cdot e}{e^{4n^{1/2}}}\;.$$
$$A=1-(31/24)/n^{1/2}+(961/1152)/n-(103433/414720)/n^{3/2}-\cdots$$
\end{cf}

Note that this is the special case $b=0$ and $z=-1$ of \ref{6.23.21}.

\smallskip

More generally:

\smallskip

\begin{verbatim}
S4(k)=sum(j=0,k-1,(-1)^j*j!);
\end{verbatim}

\begin{cf}\label{5.10.27}{\ }
\begin{verbatim}
[k->exp(1)*eint1(1),[S4(k),2*n+k],[(-1)^k*k!,-n*(n+k)]]
\end{verbatim}
$$\int_0^\infty\dfrac{e^{-t}}{1+t}\,dt=e\cdot\eintone(1)=S_4(k)+\dfrac{(-1)^kk!}{k+2-\dfrac{k+1}{k+4-\dfrac{2k+4}{k+6-\dfrac{3k+9}{k+8-\dfrac{4k+16}{k+10-\ddots}}}}}$$
Convergence type $D^+$ with $D=16$ and $C=(-1)^k2\pi\cdot e$, so that
$$e\cdot\eintone(1)-\dfrac{p(n)}{q(n)}\sim(-1)^k\dfrac{2\pi\cdot e}{e^{4n^{1/2}}}\;.$$
$$A=1+(k^2/2-k-31/24)/n^{1/2}+(k^4/8-k^3/2-7k^2/48+31k/24+961/1152)/n+\cdots$$
\end{cf}

\smallskip

\begin{cf}\label{5.10.28}{\ }
\begin{verbatim}
[()->exp(Pi/2),[3,2],[4*n^2+12*n+10]]
\end{verbatim}
$$e^{\pi/2}=3+\dfrac{10}{2+\dfrac{26}{2+\dfrac{50}{2+\dfrac{82}{2+\dfrac{122}{2+\dfrac{170}{2+\ddots}}}}}}$$
Convergence type $P^-$ with $P=1$ and $C=e^{\pi/2}\cosh(\pi/2)/2$, so that
$$e^{\pi/2}-\dfrac{p(n)}{q(n)}\sim(-1)^n\dfrac{e^{\pi/2}\cosh(\pi/2)/2}{n}\;.$$
$$A=1-2/n+(59/16)/n^2-(49/8)/n^3+(1143/128)/n^4+\cdots$$
Parametric family up to the trivial change $n\mapsto n+j$:
\begin{verbatim}
[()->exp(Pi/2),2+4*k,(2*n+1)^2+1]
\end{verbatim}
Convergence type $P^-$ with $P=2k+1$.
\end{cf}

\smallskip

\begin{cf}\label{5.10.29}{\ }
\begin{verbatim}
[()->exp(Pi/2),[4],[4*n^2+8*n+5]]
\end{verbatim}
$$e^{\pi/2}=4+\dfrac{5}{4+\dfrac{17}{4+\dfrac{37}{4+\dfrac{65}{4+\dfrac{101}{4+\dfrac{145}{4+\ddots}}}}}}$$
Convergence type $P^-$ with $P=2$ and $C=e^{\pi/2}\sinh(\pi/2)/4$, so that
$$e^{\pi/2}-\dfrac{p(n)}{q(n)}\sim(-1)^n\dfrac{e^{\pi/2}\sinh(\pi/2)/4}{n^2}\;.$$
$$A=1-3/n+6/n^2-9/n^3+(641/64)/n^4+\cdots$$
Parametric family up to the trivial change $n\mapsto n+j$:
\begin{verbatim}
[()->exp(Pi/2),4*k,4*n^2+1]
\end{verbatim}
Convergence type $P^-$ with $P=2k$.
\end{cf}

\smallskip

\begin{cf}\label{5.10.30}{\ }
\begin{verbatim}
[()->exp(Pi/2),[4,2*n+3],[n^2+4*n+5]]
\end{verbatim}
$$e^{\pi/2}=4+\dfrac{5}{5+\dfrac{10}{7+\dfrac{17}{9+\dfrac{26}{11+\dfrac{37}{13+\dfrac{50}{15+\ddots}}}}}}$$
Convergence type $E$ with $E=-(1+\sqrt{2})^2$, $P=0$, and $C=2e^{\pi/2}\sinh(\pi)/(1+\sqrt{2})^5$, so that
$$e^{\pi/2}-\dfrac{p(n)}{q(n)}\sim(-1)^n\dfrac{2e^{\pi/2}\sinh(\pi)}{(1+\sqrt{2})^{2n+5}}\;.$$
$$A=1-(5d/8)/n+(25d/16+25/64)/n^2+(-5635d/1536-125/64)/n^3+\cdots$$
\end{cf}

\smallskip

\begin{cf}\label{5.10.31}{\ }
\begin{verbatim}
[()->exp(Pi),[1,-1,2],[4,4*n^2-4*n+5]]
\end{verbatim}
$$e^{\pi}=1+\dfrac{4}{-1+\dfrac{5}{2+\dfrac{13}{2+\dfrac{29}{2+\dfrac{53}{2+\dfrac{85}{2+\ddots}}}}}}$$
Convergence type $P^-$ with $P=1$ and $C=e^{\pi}\cosh(\pi)$, so that
$$e^{\pi}-\dfrac{p(n)}{q(n)}\sim(-1)^n\dfrac{e^{\pi}\cosh(\pi)}{n}\;.$$
$$A=1-(1/2)/n^2+(7/8)/n^4-(99/32)/n^6+\cdots$$
Parametric family up to the trivial change $n\mapsto n+j$:
\begin{verbatim}
[()->exp(Pi),4*k+2,(2*n-1)^2+4]
\end{verbatim}
Convergence type $P^-$ with $P=2k+1$.
\end{cf}

\smallskip

\begin{cf}\label{5.10.32}{\ }
\begin{verbatim}
[()->exp(Pi),[15,2],[40,n^2+6*n+10]]
\end{verbatim}
$$e^{\pi}=15+\dfrac{40}{2+\dfrac{17}{2+\dfrac{26}{2+\dfrac{37}{2+\dfrac{50}{2+\dfrac{65}{2+\ddots}}}}}}$$
Convergence type $P^-$ with $P=2$ and $C=(5/8)e^{\pi}\sinh(\pi)$, so that
$$e^{\pi}-\dfrac{p(n)}{q(n)}\sim(-1)^n\dfrac{(5/8)e^{\pi}\sinh(\pi)}{n^2}\;.$$
$$A=1-7/n+(285/8)/n^2-(623/4)/n^3+(157501/256)/n^4+\cdots$$
Parametric family up to the trivial change $n\mapsto n+j$:
\begin{verbatim}
[()->exp(Pi),2*k,n^2+1]
\end{verbatim}
Convergence type $P^-$ with $P=2k$.
\end{cf}

\smallskip

\begin{cf}\label{5.10.33}{\ }
\begin{verbatim}
[()->exp(Pi),[1,-1,2*n-1],[n^2+4]]
\end{verbatim}
$$e^{\pi}=1+\dfrac{4}{-1+\dfrac{5}{3+\dfrac{8}{5+\dfrac{13}{7+\dfrac{20}{9+\dfrac{29}{11+\ddots}}}}}}$$
Convergence type $E$ with $E=-(1+\sqrt{2})^2$, $P=0$, and $C=2e^{\pi}\sinh(2\pi)/(1+\sqrt{2})$, so that
$$e^{\pi}-\dfrac{p(n)}{q(n)}\sim(-1)^n\dfrac{2e^{\pi}\sinh(2\pi)}{(1+\sqrt{2})^{2n+1}}\;.$$
$$A=1-(17d/8)/n+(17d/16+289/64)/n^2+(-2023d/1536-289/64)/n^3+\cdots$$
\end{cf}

These CFs for $e^{\pi/2}$ and $e^{\pi}$ are simply specializations of
corresponding CFs for $e^{\pi z/2}$, see \ref{6.3.8} and following.
See \cite{Riv} for other more complicated CFs for these quantities.

\smallskip

\begin{cf}\label{5.10.34}{\ }
\begin{verbatim}
[()->2*exp(Pi/sqrt(3)),[11,6*n+9],[3*n^2+12*n+21]]
\end{verbatim}
$$2e^{\pi/\sqrt{3}}=11+\dfrac{21}{15+\dfrac{36}{21+\dfrac{57}{27+\dfrac{84}{33+\dfrac{117}{39+\dfrac{156}{45+\ddots}}}}}}$$
Convergence type $E$ with $E=-(2+\sqrt{3})^2$, $P=0$, and $C=4e^{\pi/\sqrt{3}}\sinh(\pi\sqrt{3})/(2+\sqrt{3})^5$, so that
$$2e^{\pi/\sqrt{3}}-\dfrac{p(n)}{q(n)}\sim(-1)^n\dfrac{4e^{\pi/\sqrt{3}}\sinh(\pi\sqrt{3})}{(2+\sqrt{3})^{2n+5}}\;.$$
$$A=1-(13d/8)/n+(65d/16+507/128)/n^2+(-11375d/1024-2535/128)/n^3+\cdots$$
\end{cf}

\smallskip

\begin{cf}\label{5.10.35}{\ }
\begin{verbatim}
[()->exp(Pi/(3*sqrt(3))),[1,2,6*n-3],[2,3*n^2+1]]
\end{verbatim}
$$e^{\pi/(3\sqrt{3})}=1+\dfrac{2}{2+\dfrac{4}{9+\dfrac{13}{15+\dfrac{28}{21+\dfrac{49}{27+\dfrac{76}{33+\ddots}}}}}}$$
Convergence type $E$ with $E=-(2+\sqrt{3})^2$, $P=0$, and $C=2e^{\pi/(3\sqrt{3})}\sinh(\pi/\sqrt{3})/(2+\sqrt{3})$, so that
$$e^{\pi/(3\sqrt{3})}-\dfrac{p(n)}{q(n)}\sim(-1)^n\dfrac{2e^{\pi/(3\sqrt{3})}\sinh(\pi/\sqrt{3})}{(2+\sqrt{3})^{2n+1}}\;.$$
$$A=1-(7d/24)/n+(7d/48+49/384)/n^2+(-469d/27648-49/384)/n^3+\cdots$$
\end{cf}

\smallskip

These are simply specializations to $z=3$ and $z=1$ of the CF \ref{6.3.14} for
$e^{\pi z/(3\sqrt{3})}$.
  
\medskip

\section{Constants Involving Gamma Quotients}\label{sec:gam0}

\medskip

Preliminary remark: if we define the \emph{weight} $w$ of a gamma quotient
by $w(\G(a))=a$, hence $w(\pi^{1/2})=w(\G(1/2))=1/2$, $w(u)=0$ for all
algebraic numbers $u$, and extended additively, then all the gamma quotients
occurring in this book have \emph{integral} weight. For example, you will
find CFs for $\G(1/4)^2/\sqrt{\pi}$, but not for $\G(1/4)^2$ or
$\G(1/4)^2/\pi$.

In addition, please note that all the irrational auxiliary constants which are used are
themselves pure gamma quotients, in particular

\begin{align*}
  \sqrt{\pi}&=\G(1/2)\;,\ \sqrt{2}=\dfrac{\G(1/4)\G(3/4)}{\G(1/2)^2}\;,\ \sqrt{2\pi}=\dfrac{\G(1/4)\G(3/4)}{\G(1/2)}\;,\\
  3^{1/2}&=2\dfrac{\G(1/2)^2}{\G(1/3)\G(2/3)}\;,\ 2^{1/3}=2\dfrac{\G(1/2)\G(1/3)}{\G(1/6)\G(2/3)}\;,\ 2^{2/3}=\dfrac{\G(1/6)\G(2/3)}{\G(1/2)\G(1/3)}\;.\end{align*}

Thus, if you find it more aesthetic, in the following CFs you can replace these constants by
the corresponding gamma quotients.

\smallskip

\subsection{Around the Lemniscatic Constants}

\medskip

There is a wide variety of sometimes conflicting notation related to the
so-called lemniscatic constants. We have essentially four closely related
constants: the lemniscatic constant $\varpi$, Gauss' constant $G$ (warning:
not to be confused with Catalan's constant, also denoted $G$), and
the first and second lemniscatic constants $L_1$ and $L_2$. Here are
their expressions in terms of gamma products, recalling that of course
$\G(1/4)\G(3/4)=\pi\sqrt{2}$:
\begin{align*}
  \varpi&=\dfrac{\G(1/4)^2}{2\sqrt{2\pi}}=\dfrac{\sqrt{\pi}}{2}\dfrac{\G(1/4)}{\G(3/4)}=2.622057...\\
    G&=\dfrac{\G(1/4)^2}{(2\pi)^{3/2}}=\dfrac{1}{2\sqrt{\pi}}\dfrac{\G(1/4)}{\G(3/4)}=0.834626...\\
    L_1&=\dfrac{\G(1/4)^2}{4\sqrt{2\pi}}=\dfrac{\sqrt{\pi}}{4}\dfrac{\G(1/4)}{\G(3/4)}=1.311028...\\
    L_2&=\dfrac{\G(3/4)^2}{\sqrt{2\pi}}=\sqrt{\pi}\dfrac{\G(3/4)}{\G(1/4)}=0.599070...\end{align*}
Note the evident relations $L_1=\varpi/2$, $G=\varpi/\pi$, $L_2=1/(2G)$,
$L_1L_2=\pi/4$, and
$$\sqrt{\pi}\dfrac{\G(3/4)}{\G(1/4)}=L_2\text{\quad and\quad}\sqrt{\pi}\dfrac{\G(1/4)}{\G(3/4)}=4L_1=2\varpi\;.$$

To simplify the classification of the following CFs, note (up to a rational
multiple) all these constants are of the form
$2^{j/2}\pi^{1/2-k}\G((2m+1)/4)^2$ with $j=0,1$, $k=0,1,2,3$, and $m=0,1$, so
we will always use this uniform representation, recalling other representations
as we go along.

Note that
$$(2^{j_1/2}\pi^{1/2-k_1}\G(1/4)^2)(2^{j_2/2}\pi^{1/2-k_2}\G(3/4)^2)=2^{1+(j_1+j_2)/2}\pi^{3-k_1-k_2}\;,$$
so the CF for $(j_2,k_2,1)$ can be trivially deduced from the CF for
$(1-j_2,3-k_2,0)$ and conversely, so we will only mention them.

\medskip

{\bf $k=0$: Continued Fractions for $2^{j/2}\pi^{1/2}\G((2m+1)/4)^2$}

\smallskip

I have not found any with $j=0$, so I give only those with $j=1$.

\smallskip

\begin{cf}\label{5.11.1}{\ }
\begin{verbatim}
[()->sqrt(2*Pi)*gamma(1/4)^2,
[32,50,(4*n-1)*(16*n^2-8*n+9)],[32,-2*n*(2*n+1)*(4*n+1)^4]]
\end{verbatim}
$$\sqrt{2\pi}\G(1/4)^2=32+\dfrac{32}{50-\dfrac{3750}{399-\dfrac{131220}{1419-\dfrac{1199562}{3495-\dfrac{6013512}{7011-\ddots}}}}}$$
Convergence type $P^+$ with $P=3/2$ and $C=4/(3\sqrt{\pi})$, so that
$$\sqrt{2\pi}\G(1/4)^2-\dfrac{p(n)}{q(n)}\sim\dfrac{4/(3\sqrt{\pi})}{n^{3/2}}\;.$$
$$A=1-(9/8)/n+(799/896)/n^2-(507/1024)/n^3+\cdots$$
Series:
$$\sqrt{2\pi}\G(1/4)^2=32+16\sum_{n\ge0}\dfrac{(3/2)_n}{(4n+5)^2(n+1)!}$$
Parametric family for $k\ge0$:
\begin{verbatim}
[()->sqrt(2*Pi)*gamma(1/4)^2,(4*n-1)*(16*n^2-8*n+32*k^2+24*k+9),
                             -2*n*(2*n+1)*(4*n+1)^4]
\end{verbatim}
Convergence type $P^+$ with $P=4k+3/2$.
\end{cf}

\smallskip

\begin{cf}\label{5.11.2}{\ }
\begin{verbatim}
[()->sqrt(2*Pi)*gamma(3/4)^2,
[4,18,(4*n-3)*(16*n^2-24*n+25)],[-4,-2*n*(2*n-1)*(4*n-1)^4]]
\end{verbatim}
$$\sqrt{2\pi}\G(3/4)^2=4-\dfrac{4}{18-\dfrac{162}{205-\dfrac{28812}{873-\dfrac{439230}{2405-\dfrac{2835000}{5185-\dfrac{11728890}{9597-\ddots}}}}}}$$
Convergence type $P^+$ with $P=5/2$ and $C=-1/(20\sqrt{\pi})$, so that
$$\sqrt{2\pi}\G(3/4)^2-\dfrac{p(n)}{q(n)}\sim-\dfrac{1/(20\sqrt{\pi})}{n^{5/2}}\;.$$
$$A=1-(5/8)/n-(55/1152)/n^2+(265/1024)/n^3+\cdots$$
Series:
$$\sqrt{2\pi}\G(3/4)^2=4-2\sum_{n\ge0}\dfrac{(1/2)_n}{(4n+3)^2(n+1)!}$$
Parametric family for $k\ge0$:
\begin{verbatim}
[()->sqrt(2*Pi)*gamma(3/4)^2,(4*n-3)*(16*n^2-24*n+32*k^2+40*k+25),
                             -2*n*(2*n-1)*(4*n-1)^4]
\end{verbatim}
Convergence type $P^+$ with $P=4k+5/2$.
\end{cf}

\medskip

{\bf $k=1,j=0$: Continued Fractions for $\pi^{-1/2}\G((2m+1)/4)^2$}

\medskip

For $m=0$, note that

$$\dfrac{\G(1/4)^2}{\sqrt{\pi}}=\sqrt{2\pi}\dfrac{\G(1/4)}{\G(3/4)}=4\sqrt{2}L_1=2\sqrt{2}\om\;.$$

\smallskip

\begin{cf}\label{5.11.3}{\ }
\begin{verbatim}
[()->gamma(1/4)^2/sqrt(Pi),
[4,20,32*n^2-12*n+3],[12,-4*n*(4*n+1)^2*(4*n+3)]]
\end{verbatim}
$$\dfrac{\G(1/4)^2}{\sqrt{\pi}}=4+\dfrac{12}{20-\dfrac{700}{107-\dfrac{7128}{255-\dfrac{30420}{467-\dfrac{87856}{743-\dfrac{202860}{1083-\ddots}}}}}}$$
Convergence type $P^+$ with $P=1/4$ and $C=4/\G(3/4)$, so that
$$\dfrac{\G(1/4)^2}{\sqrt{\pi}}-\dfrac{p(n)}{q(n)}\sim\dfrac{4/\G(3/4)}{n^{1/4}}\;.$$
$$A=1-(31/160)/n+(491/6144)/n^2-(25653/851968)/n^3+\cdots$$
Series:
$$\dfrac{\G(1/4)^2}{\sqrt{\pi}}=4+3\sum_{n\ge0}\dfrac{(7/4)_n}{(4n+5)(n+1)!}$$
Parametric family for $k\ge0$:
\begin{verbatim}
[()->gamma(1/4)^2/sqrt(Pi),
32*n^2-12*n+3+4*k*(4*k+1),-4*n*(4*n+1)^2*(4*n+3)]
\end{verbatim}
Convergence type $P^+$ with $P=2k+1/4$.
\end{cf}

\smallskip

\begin{cf}\label{5.11.4}{\ }
\begin{verbatim}
[()->gamma(1/4)^2/sqrt(Pi),
[16,32*n^2+8*n+1],[-144,-16*n*(n+1)*(4*n+3)^2]]
\end{verbatim}
$$\dfrac{\G(1/4)^2}{\sqrt{\pi}}=16-\dfrac{144}{41-\dfrac{1568}{145-\dfrac{11616}{313-\dfrac{43200}{545-\dfrac{115520}{841-\dfrac{253920}{1201-\ddots}}}}}}$$
Convergence type $P^+$ with $P=1/2$ and $C=-\pi^3/(2\G(3/4)^6)$, so that
$$\dfrac{\G(1/4)^2}{\sqrt{\pi}}-\dfrac{p(n)}{q(n)}\sim-\dfrac{\pi^3/(2\G(3/4)^6)}{n^{1/2}}$$
$$A=1-(31/48)/n+(309/512)/n^2-(103735/172032)/n^3+\cdots$$
Series:
$$\dfrac{\sqrt{\pi}}{\G(1/4)^2}=\dfrac{1}{16}\sum_{n\ge0}\dfrac{(3/4)_n^2}{(n+1)n!^2}$$
Parametric family for $k\ge0$:
\begin{verbatim}
[()->gamma(1/4)^2/sqrt(Pi),
32*n^2+8*n+(4*k+1)^2,-16*n*(n+1)*(4*n+3)^2]
\end{verbatim}
Convergence type $P^+$ with $P=2k+1/2$.
\end{cf}

\smallskip

\begin{cf}\label{5.11.5}{\ }
\begin{verbatim}
[()->gamma(1/4)^2/sqrt(Pi),[16,-50,144*n^2-216*n+49],
[432,2*n*(2*n-1)*(4*n-5)*(4*n-3)*(4*n+1)*(4*n+3)]]
\end{verbatim}
$$\dfrac{\G(1/4)^2}{\sqrt{\pi}}=16+\dfrac{432}{-50-\dfrac{70}{193+\dfrac{17820}{697+\dfrac{368550}{1489+\dfrac{2586584}{2569+\dfrac{11084850}{3937+\ddots}}}}}}$$
Convergence type $P^-$ with $P=9/2$ and $C=-9/(16\sqrt{\pi})$, so that
$$\dfrac{\G(1/4)^2}{\sqrt{\pi}}-\dfrac{p(n)}{q(n)}\sim(-1)^{n+1}\dfrac{9/(16\sqrt{\pi})}{n^{9/2}}$$
$$A=1-(9/8)/n-(223/128)/n^2+(3757/1024)/n^3+\cdots$$
Series:
$$\dfrac{\G(1/4)^2}{\sqrt{\pi}}=16-72\sum_{n\ge0}(-1)^n\dfrac{(4n+3)(1/2)_n}{(4n+1)^2(4n+5)^2(n+1)!}$$
\end{cf}

\smallskip

\begin{cf}\label{5.11.6}{\ }
\begin{verbatim}
[()->gamma(1/4)^2/sqrt(Pi),
[16,25,32*n^2-16*n+1],[-144,-8*(n+1)*(2*n-1)*(4*n+1)^2]]
\end{verbatim}
$$\dfrac{\G(1/4)^2}{\sqrt{\pi}}=16-\dfrac{144}{25-\dfrac{400}{97-\dfrac{5832}{241-\dfrac{27040}{449-\dfrac{80920}{721-\dfrac{190512}{1057-\ddots}}}}}}$$
Convergence type $P^+$ with $P=1$ and $C=-9\pi^{3/2}/(8\G(3/4)^2)$, so that
$$\dfrac{\G(1/4)^2}{\sqrt{\pi}}-\dfrac{p(n)}{q(n)}\sim-\dfrac{9\pi^{3/2}/(8\G(3/4)^2)}{n}$$
$$A=1-(15/32)/n+(313/1536)/n^2+\cdots$$
Series:
\begin{align*}\dfrac{\G(1/4)^2}{\sqrt{\pi}}&=16-\dfrac{144}{25}\sum\dfrac{(n+1)!(1/2)_n}{(9/4)_n^2}\\
  \dfrac{\sqrt{\pi}}{\G(1/4)^2}&=-\dfrac{1}{2}+\dfrac{9}{16}\sum_{n\ge0}\dfrac{(1/4)_n^2}{(n+1)!(1/2)_n}\end{align*}
Parametric family for $k\ge0$:
\begin{verbatim}
[()->gamma(1/4)^2/sqrt(Pi),
32*n^2-16*n+16*k*(k+1)+1,-8*(n+1)*(2*n-1)*(4*n+1)^2]
\end{verbatim}
Convergence type $P^+$ with $P=2k+1$.
\end{cf}

\smallskip

\begin{cf}\label{5.11.7}{\ }
\begin{verbatim}
[()->gamma(1/4)^2/sqrt(Pi),
[16,33,24*(2*n-1)],[-288,(4*n-3)^2*(4*n+3)^2]]
\end{verbatim}
$$\dfrac{\G(1/4)^2}{\sqrt{\pi}}=16-\dfrac{288}{33+\dfrac{49}{72+\dfrac{3025}{120+\dfrac{18225}{168+\dfrac{61009}{216+\dfrac{152881}{264+\ddots}}}}}}$$
Convergence type $P^-$ with $P=3$ and $C=-\pi^2/(8\G(3/4)^4)$, so that
$$\dfrac{\G(1/4)^2}{\sqrt{\pi}}-\dfrac{p(n)}{q(n)}\sim(-1)^{n+1}\dfrac{\pi^2/(8\G(3/4)^4)}{n^3}$$
$$A=1-(3/2)/n+(29/32)/n^2+(15/64)/n^3+(1607/2048)/n^4+\cdots$$
Series:
$$\dfrac{\G(1/4)^2}{\sqrt{\pi}}=16-288\sum_{n\ge0}(-1)^n\dfrac{(7/4)_n^2}{(16n^2+16n+1)(16n^2+48n+33)(5/4)_n^2}$$
Parametric family for $k\ge0$:
\begin{verbatim}
[()->gamma(1/4)^2/sqrt(Pi),
8*(4*k+3)*(2*n-1),(4*n-3)^2*(4*n+3)^2]
\end{verbatim}
Convergence type $P^-$ with $P=4k+3$.
\end{cf}

\smallskip

\begin{cf}\label{5.11.8}{\ }
\begin{verbatim}
[()->gamma(1/4)^2/sqrt(Pi),[8,10,12*n-3],[-8,2*n*(2*n+1)*(4*n+1)^2]]
\end{verbatim}
$$\dfrac{\G(1/4)^2}{\sqrt{\pi}}=8-\dfrac{8}{10+\dfrac{150}{21+\dfrac{1620}{33+\dfrac{7098}{45+\dfrac{20808}{57+\dfrac{48510}{69+\ddots}}}}}}$$
Convergence type $P^-$ with $P=3/2$ and $C=-1/\sqrt{\pi}$, so that
$$\dfrac{\G(1/4)^2}{\sqrt{\pi}}-\dfrac{p(n)}{q(n)}\sim(-1)^{n+1}\dfrac{1/\sqrt{\pi}}{n^{3/2}}$$
$$A=1-(9/8)/n+(73/128)/n^2+(357/1024)/n^3+\cdots$$
Series:
$$\dfrac{\G(1/4)^2}{\sqrt{\pi}}=8\sum_{n\ge0}(-1)^n\dfrac{(1/2)_n}{(4n+1)n!}$$
Parametric family for $k\ge0$:
\begin{verbatim}
[()->gamma(1/4)^2/sqrt(Pi),(8*k+3)*(4*n-1),2*n*(2*n+1)*(4*n+1)^2]
\end{verbatim}
Convergence type $P^-$ with $P=4k+3/2$.
\end{cf}

\smallskip

\begin{cf}\label{5.11.9}{\ }
\begin{verbatim}
[()->gamma(1/4)^2/sqrt(Pi),[8,32*n^2-40*n+25],[-8,-16*n^2*(4*n-1)^2]]
\end{verbatim}
$$\dfrac{\G(1/4)^2}{\sqrt{\pi}}=8-\dfrac{8}{17-\dfrac{144}{73-\dfrac{3136}{193-\dfrac{17424}{377-\dfrac{57600}{625-\dfrac{144400}{937-\ddots}}}}}}$$
Convergence type $P^+$ with $P=3/2$ and $C=-\pi^3/(48\G(3/4)^6)$, so that
$$\dfrac{\G(1/4)^2}{\sqrt{\pi}}-\dfrac{p(n)}{q(n)}\sim-\dfrac{\pi^3/(48\G(3/4)^6)}{n^{3/2}}$$
$$A=1-(9/16)/n+(475/3584)/n^2+(465/8192)/n^3+\cdots$$
Series:
$$\dfrac{\sqrt{\pi}}{\G(1/4)^2}=\dfrac{1}{8}\sum_{n\ge0}\dfrac{(-1/4)_n^2}{n!^2}$$
Parametric family for $k\ge0$:
\begin{verbatim}
[()->gamma(1/4)^2/sqrt(Pi),32*n^2-40*n+25+8*k*(2*k+3),-16*n^2*(4*n-1)^2]
\end{verbatim}
Convergence type $P^+$ with $P=2k+3/2$.
\end{cf}

\smallskip

\begin{cf}\label{5.11.10}{\ }
\begin{verbatim}
[()->gamma(1/4)^2/sqrt(Pi),
[4,5,12*n^2-4*n-3],[12,-8*n^2*(2*n-1)*(2*n+3)]]
\end{verbatim}
$$\dfrac{\G(1/4)^2}{\sqrt{\pi}}=4+\dfrac{12}{5-\dfrac{40}{37-\dfrac{672}{93-\dfrac{3240}{173-\dfrac{9856}{277-\dfrac{23400}{405-\ddots}}}}}}$$
Convergence type $E$ with $E=2$, $P=1$, and $C=\pi^2/\G(3/4)^4$, so that
$$\dfrac{\G(1/4)^2}{\sqrt{\pi}}-\dfrac{p(n)}{q(n)}\sim\dfrac{\pi^2/\G(3/4)^4}{2^nn}\;.$$
$$A=1-(5/4)/n+(105/32)/n^2-(1795/128)/n^3+(165675/2048)/n^4+\cdots$$
Series:
$$\dfrac{\sqrt{\pi}}{\G(1/4)^2}=\dfrac{1}{4}\sum_{n\ge0}\dfrac{(-1/2)_n(3/2)_n}{n!^2}2^{-n}$$
Parametric family for $k\ge0$:
\begin{verbatim}
[()->gamma(1/4)^2/sqrt(Pi),
12*n^2-4*n-3,-8*(n-k)*n*(2*n-1)*(2*n+3)]
\end{verbatim}
Convergence type $E$ with $E=2$ and $P=3k+1$.
\end{cf}

\smallskip

\begin{cf}\label{5.11.11}{\ }
\begin{verbatim}
[()->gamma(1/4)^2/sqrt(Pi),[[32/5,16],[12*n+4,6]],
                           [[96/5,40],[-(2*n+1)*(2*n+5),2*(4*n+1)*(4*n+5)]]]
\end{verbatim}
$$\dfrac{\G(1/4)^2}{\sqrt{\pi}}=32/5+\dfrac{96/5}{16+\dfrac{40}{16-\dfrac{21}{6+\dfrac{90}{28-\dfrac{45}{6+\dfrac{234}{40-\ddots}}}}}}$$
Convergence type $E$ with $E=-2\sqrt{2}$, $P=3/2$, and
$C=2^{7/2}\pi/(27\G(3/4)^2)$, so that
$$\dfrac{\G(1/4)^2}{\sqrt{\pi}}-\dfrac{p(n)}{q(n)}\sim(-1)^n\dfrac{\pi/(27\G(3/4)^2)}{2^{(3n-7)/2}n^{3/2}}\;.$$
$$A=1-(31/24)/n+(331/3456)/n^2+(1961195/248832)/n^3+\cdots$$
Series:
$$\dfrac{\G(1/4)^2}{\sqrt{\pi}}=-32\sum_{n\ge0}(-1)^n\dfrac{(3n+1)(2n+1)(2n+3)(1/2)_n^2}{(12n^2-4n-3)(12n^2+20n+5)(1/4)_n(5/4)_n}2^{-3n}$$
\end{cf}

\medskip

For $m=1$, note that

$$\dfrac{\G(3/4)^2}{\sqrt{\pi}}=\sqrt{2\pi}\dfrac{\G(3/4)}{\G(1/4)}=\sqrt{2}L_2=\dfrac{1}{\sqrt{2}G}\;.$$

\smallskip

\begin{cf}\label{5.11.12}{\ }
\begin{verbatim}
[()->gamma(3/4)^2/sqrt(Pi),
[2/3,28,32*n^2-4*n+3],[2,-4*n*(4*n+1)*(4*n+3)^2]]
\end{verbatim}
$$\dfrac{\G(3/4)^2}{\sqrt{\pi}}=2/3+\dfrac{2}{28-\dfrac{980}{123-\dfrac{8712}{279-\dfrac{35100}{499-\dfrac{98192}{783-\dfrac{222180}{1131-\ddots}}}}}}$$
Convergence type $P^+$ with $P=3/4$ and $C=2/(3\G(1/4))$, so that
$$\dfrac{\G(3/4)^2}{\sqrt{\pi}}-\dfrac{p(n)}{q(n)}\sim\dfrac{2/(3\G(1/4))}{n^{3/4}}\;.$$
$$A=1-(165/224)/n+(13459/22528)/n^2-(31113/65536)/n^3+\cdots$$
Series:
$$\dfrac{\G(3/4)^2}{\sqrt{\pi}}=\dfrac{2}{3}+\dfrac{1}{2}\sum_{n\ge0}\dfrac{(5/4)_n}{(4n+7)(n+1)!}$$
Parametric family for $k\ge0$:
\begin{verbatim}
[()->gamma(3/4)^2/sqrt(Pi),
32*n^2-4*n+3+4*k*(4*k+3),-4*n*(4*n+1)*(4*n+3)^2]
\end{verbatim}
Convergence type $P^+$ with $P=2k+3/4$.
\end{cf}

\smallskip

\begin{cf}\label{5.11.13}{\ }
\begin{verbatim}
[()->gamma(3/4)^2/sqrt(Pi),
[8/9,32*n^2-8*n+9],[-8/9,-16*n*(n+1)*(4*n+1)^2]]
\end{verbatim}
$$\dfrac{\G(3/4)^2}{\sqrt{\pi}}=8/9-\dfrac{8/9}{33-\dfrac{800}{121-\dfrac{7776}{273-\dfrac{32448}{489-\dfrac{92480}{769-\dfrac{211680}{1113-\ddots}}}}}}$$
Convergence type $P^+$ with $P=3/2$ and $C=-3\pi^3/\G(1/4)^6$, so that
$$\dfrac{\G(3/4)^2}{\sqrt{\pi}}-\dfrac{p(n)}{q(n)}\sim-\dfrac{3\pi^3/\G(1/4)^6}{n^{3/2}}$$
$$A=1-(117/80)/n+(877/512)/n^2-(14791/8192)/n^3+\cdots$$
Series:
$$\dfrac{\sqrt{\pi}}{\G(3/4)^2}=\dfrac{9}{8}\sum_{n\ge0}\dfrac{(1/4)_n^2}{(n+1)n!^2}$$
Parametric family for $k\ge0$:
\begin{verbatim}
[()->gamma(3/4)^2/sqrt(Pi),
32*n^2-8*n+(4*k+3)^2,-16*n*(n+1)*(4*n+1)^2]
\end{verbatim}
Convergence type $P^+$ with $P=2k+3/2$.
\end{cf}

\smallskip

\begin{cf}\label{5.11.14}{\ }
\begin{verbatim}
[()->gamma(3/4)^2/sqrt(Pi),[8/9,98/3,112*n^2-56*n-25],
[-40/27,2*n*(2*n+1)*(4*n-3)*(4*n-1)*(4*n+3)*(4*n+5)]]
\end{verbatim}
$$\dfrac{\G(3/4)^2}{\sqrt{\pi}}=8/9-\dfrac{40/27}{98/3+\dfrac{1134}{311+\dfrac{100100}{815+\dfrac{1060290}{1543+\dfrac{5601960}{2495+\dfrac{20429750}{3671+\ddots}}}}}}$$
Convergence type $P^-$ with $P=7/2$ and $C=-1/(16\sqrt{\pi})$, so that
$$\dfrac{\G(3/4)^2}{\sqrt{\pi}}-\dfrac{p(n)}{q(n)}\sim(-1)^{n+1}\dfrac{1/(16\sqrt{\pi})}{n^{7/2}}$$
$$A=1-(21/8)/n+(377/128)/n^2+(33/1024)/n^3+\cdots$$
Series:
$$\dfrac{\G(3/4)^2}{\sqrt{\pi}}=\dfrac{8}{9}-4\sum_{n\ge0}(-1)^n\dfrac{(4n+5)(3/2)_n}{(4n+3)^2(4n+7)^2(n+1)!}$$
\end{cf}

\smallskip

\begin{cf}\label{5.11.15}{\ }
\begin{verbatim}
[()->gamma(3/4)^2/sqrt(Pi),
[8/9,49,32*n^2+16*n+9],[-8/9,-8*(n+1)*(2*n+1)*(4*n+3)^2]]
\end{verbatim}
$$\dfrac{\G(3/4)^2}{\sqrt{\pi}}=8/9-\dfrac{8/9}{49-\dfrac{2352}{169-\dfrac{14520}{345-\dfrac{50400}{585-\dfrac{129960}{889-\dfrac{279312}{1257-\ddots}}}}}}$$
Convergence type $P^+$ with $P=1$ and $C=-\pi^{3/2}/(8\G(1/4)^2)$, so that
$$\dfrac{\G(3/4)^2}{\sqrt{\pi}}-\dfrac{p(n)}{q(n)}\sim-\dfrac{\pi^{3/2}/(8\G(1/4)^2)}{n}$$
$$A=1-(39/32)/n+(723/512)/n^2+\cdots$$
Series:
\begin{align*}\dfrac{\G(3/4)^2}{\sqrt{\pi}}&=\dfrac{8}{9}-\dfrac{8}{441}\sum\dfrac{(n+1)!(3/2)_n}{(11/4)_n^2}\\
  \dfrac{\sqrt{\pi}}{\G(3/4)^2}&=1+\dfrac{1}{8}\sum_{n\ge0}\dfrac{(3/4)_n^2}{(n+1)!(3/2)_n}\end{align*}
Parametric family for $k\ge0$:
\begin{verbatim}
[()->gamma(3/4)^2/sqrt(Pi),
32*n^2+16*n+16*k*(k+1)+9,-8*(n+1)*(2*n+1)*(4*n+3)^2]
\end{verbatim}
Convergence type $P^+$ with $P=2k+1$.
\end{cf}

\smallskip

\begin{cf}\label{5.11.16}{\ }
\begin{verbatim}
[()->gamma(3/4)^2/sqrt(Pi),
[8/9,41,40*(2*n-1)],[-16/9,(4*n-1)^2*(4*n+1)^2]]
\end{verbatim}
$$\dfrac{\G(3/4)^2}{\sqrt{\pi}}=8/9-\dfrac{16/9}{41+\dfrac{225}{120+\dfrac{3969}{200+\dfrac{20449}{280+\dfrac{65025}{360+\dfrac{159201}{440+\ddots}}}}}}$$
Convergence type $P^-$ with $P=5$ and $C=-9\pi^2/(16\G(1/4)^4)$, so that
$$\dfrac{\G(3/4)^2}{\sqrt{\pi}}-\dfrac{p(n)}{q(n)}\sim(-1)^{n+1}\dfrac{9\pi^2/(16\G(1/4)^4)}{n^3}$$
$$A=1-(5/2)/n-(5/32)/n^2+(595/64)/n^3-(597/2048)/n^4+\cdots$$
Series:
$$\dfrac{\G(3/4)^2}{\sqrt{\pi}}=\dfrac{8}{9}-16\sum_{n\ge0}(-1)^n\dfrac{(5/4)_n^2}{(16n^2+16n+9)(16n^2+48n+41)(7/4)_n^2}$$
Parametric family for $k\ge0$:
\begin{verbatim}
[()->gamma(3/4)^2/sqrt(Pi),
8*(4*k+5)*(2*n-1),(4*n-1)^2*(4*n+1)^2]
\end{verbatim}
Convergence type $P^-$ with $P=4k+5$.
\end{cf}

\smallskip

\begin{cf}\label{5.11.17}{\ }
\begin{verbatim}
[()->gamma(3/4)^2/sqrt(Pi),[4/3,14,4*n+1],[-12,2*n*(2*n+3)*(4*n+3)^2]]
\end{verbatim}
$$\dfrac{\G(3/4)^2}{\sqrt{\pi}}=4/3-\dfrac{12}{14+\dfrac{490}{9+\dfrac{3388}{13+\dfrac{12150}{17+\dfrac{31768}{21+\dfrac{68770}{25+\ddots}}}}}}$$
Convergence type $P^-$ with $P=1/2$ and $C=-1/\sqrt(\pi)$, so that
$$\dfrac{\G(3/4)^2}{\sqrt{\pi}}-\dfrac{p(n)}{q(n)}\sim(-1)^{n+1}\dfrac{1/\sqrt{\pi}}{n^{1/2}}$$
$$A=1-(5/8)/n+(65/128)/n^2-(375/1024)/n^3+\cdots$$
Series:
$$\dfrac{\G(3/4)^2}{\sqrt{\pi}}=4\sum_{n\ge0}(-1)^n\dfrac{(3/2)_n}{(4n+3)n!}$$
Parametric family for $k\ge0$:
\begin{verbatim}
[()->gamma(3/4)^2/sqrt(Pi),(8*k+1)*(4*n+1),2*n*(2*n+3)*(4*n+3)^2]
\end{verbatim}
Convergence type $P^-$ with $P=4k+1/2$.
\end{cf}

\smallskip

\begin{cf}\label{5.11.18}{\ }
\begin{verbatim}
[()->gamma(3/4)^2/sqrt(Pi),[4/3,32*n^2-56*n+49],[-12,-16*n^2*(4*n-3)^2]]
\end{verbatim}
$$\dfrac{\G(3/4)^2}{\sqrt{\pi}}=4/3-\dfrac{12}{25-\dfrac{16}{65-\dfrac{1600}{169-\dfrac{11664}{337-\dfrac{43264}{569-\dfrac{115600}{865-\ddots}}}}}}$$
Convergence type $P^+$ with $P=5/2$ and $C=-27\pi^3/(40\G(1/4)^6)$, so that
$$\dfrac{\G(3/4)^2}{\sqrt{\pi}}-\dfrac{p(n)}{q(n)}\sim-\dfrac{27\pi^3/(40\G(1/4)^6)}{n^{5/2}}$$
$$A=1-(5/16)/n-(595/4608)/n^2+(805/8192)/n^3+\cdots$$
Series:
$$\dfrac{\sqrt{\pi}}{\G(3/4)^2}=\dfrac{3}{4}\sum_{n\ge0}\dfrac{(-3/4)_n^2}{n!^2}$$
Parametric family for $k\ge0$:
\begin{verbatim}
[()->gamma(3/4)^2/sqrt(Pi),32*n^2-56*n+49+8*k*(2*k+5),-16*n^2*(4*n-3)^2]
\end{verbatim}
Convergence type $P^+$ with $P=2k+5/2$.
\end{cf}

\smallskip

\begin{cf}\label{5.11.19}{\ }
\begin{verbatim}
[()->gamma(3/4)^2/sqrt(Pi),
[1,9,12*n^2-4*n+1],[-1,-8*n^2*(2*n+1)^2]]
\end{verbatim}
$$\dfrac{\G(3/4)^2}{\sqrt{\pi}}=1-\dfrac{1}{9-\dfrac{72}{41-\dfrac{800}{97-\dfrac{3528}{177-\dfrac{10368}{281-\dfrac{24200}{409-\ddots}}}}}}$$
Convergence type $E$ with $E=2$, $P=1$, and $C=-4\pi^2/\G(1/4)^4$, so that
$$\dfrac{1}{\sqrt{2\pi}}\dfrac{\G(1/4)}{\G(3/4)}-\dfrac{p(n)}{q(n)}\sim-\dfrac{\pi^2/\G(1/4)^4}{2^{n-2}n}\;.$$
$$A=1-(9/4)/n+(225/32)/n^2-(3927/128)/n^3+(362619/2048)/n^4+\cdots$$
Series:
$$\dfrac{\sqrt{\pi}}{\G(3/4)^2}=\sum_{n\ge0}\dfrac{(1/2)_n^2}{n!^2}2^{-n}$$
Parametric family for $k\ge0$:
\begin{verbatim}
[()->gamma(3/4)^2/sqrt(Pi),
12*n^2-4*n+1,-8*(n-k)*n*(2*n+1)^2]
\end{verbatim}
Convergence type $E$ with $E=2$ and $P=3k+1$.
\end{cf}

\smallskip

\begin{cf}\label{5.11.20}{\ }
\begin{verbatim}
[()->gamma(3/4)^2/sqrt(Pi),
[[[8/9,-32],[12*n+4,6]],[[16/9,-144],[-(2*n+3)^2,2*(4*n+3)^2]]]
\end{verbatim}
$$\dfrac{\G(3/4)^2}{\sqrt{\pi}}=8/9+\dfrac{16/9}{-32-\dfrac{144}{16-\dfrac{25}{6+\dfrac{98}{28-\dfrac{49}{6+\dfrac{242}{40-\ddots}}}}}}$$
Convergence type $E$ with $E=-2\sqrt{2}$, $P=3/2$, and
$C=-2^{7/2}\pi/(27\G(1/4)^2)$, so that
$$\dfrac{\G(3/4)^2}{\sqrt{\pi}}-\dfrac{p(n)}{q(n)}\sim(-1)^{n+1}\dfrac{\pi/(27\G(1/4)^2)}{2^{(3n-7)/2}n^{3/2}}\;.$$
$$A=1-(67/24)/n+(23875/3456)/n^2-(4211977/248832)/n^3+\cdots$$
Series:
$$\dfrac{\G(3/4)^2}{\sqrt{\pi}}=8\sum_{n\ge0}(-1)^n\dfrac{(3n+1)(3/2)_n^2}{(12n^2-4n+1)(12n^2+20n+9)(3/4)_n^2}2^{-3n}$$
\end{cf}

\medskip

{\bf $k=1,j=1$: Continued Fractions for $(2\pi)^{-1/2}\G((2m+1)/4)^2$}

\medskip

For $m=0$, note that

$$\dfrac{\G(1/4)^2}{\sqrt{2\pi}}=\sqrt{\pi}\dfrac{\G(1/4)}{\G(3/4)}=4L_1=2\varpi\;.$$

\smallskip

\begin{cf}\label{5.11.21}{\ }
\begin{verbatim}
[()->gamma(1/4)^2/sqrt(2*Pi),[0,1,2],[8,(4*n-1)*(4*n-3)]]
\end{verbatim}
$$\dfrac{\G(1/4)^2}{\sqrt{2\pi}}=\dfrac{8}{1+\dfrac{3}{2+\dfrac{35}{2+\dfrac{99}{2+\dfrac{195}{2+\dfrac{323}{2+\ddots}}}}}}$$
Convergence type $P^-$ with $P=1/2$ and $C=\G(1/4)/\G(3/4)$, so that
$$\dfrac{\G(1/4)^2}{\sqrt{2\pi}}-\dfrac{p(n)}{q(n)}\sim(-1)^n\dfrac{\G(1/4)/\G(3/4)}{n^{1/2}}\;.$$
$$A=1-(7/64)/n^2+(861/8192)/n^4-\cdots$$
Series:
$$\dfrac{\G(1/4)^2}{\sqrt{2\pi}}=8\sum_{n\ge0}(-1)^n\dfrac{(3/4)_n}{(5/4)_n}$$
Parametric families for $v$, $w$ and $k$ nonnegative:
\begin{verbatim}
[()->gamma(1/4)^2/sqrt(2*Pi),
8*k+4*w+2,(4*n+4*u-1)*(4*n+4*u-8*v+4*w-3)]
[()->gamma(1/4)^2/sqrt(2*Pi),
8*k+4*w+2,(4*n+4*u-3)*(4*n+4*u+8*v-4*w-1)]
\end{verbatim}
Convergence types $P^-$ with $P=2k+w+1/2$.
\end{cf}

\smallskip

\begin{cf}\label{5.11.22}{\ }
\begin{verbatim}
[()->gamma(1/4)^2/sqrt(2*Pi),[4,2],[4,2*n*(2*n+1)]]
\end{verbatim}
$$\dfrac{\G(1/4)^2}{\sqrt{2\pi}}=4+\dfrac{4}{2+\dfrac{6}{2+\dfrac{20}{2+\dfrac{42}{2+\dfrac{72}{2+\dfrac{110}{2+\ddots}}}}}}$$
Convergence type $P^-$ with $P=1$ and $C=\sqrt{\pi}\G(1/4)/\G(3/4)/4$,
so that
$$\dfrac{\G(1/4)^2}{\sqrt{2\pi}}-\dfrac{p(n)}{q(n)}\sim(-1)^n\dfrac{\sqrt{\pi}\G(1/4)/\G(3/4)/4}{n}\;.$$
$$A=1-(3/4)/n+(21/64)/n^2+(27/256)/n^3-(387/2048)/n^4+\cdots$$
Series:
\begin{align*}\dfrac{\G(1/4)^2}{\sqrt{2\pi}}&=4+\dfrac{4}{5}\sum_{n\ge0}\dfrac{(3/4)_nn!}{(9/4)_n(3/2)_n}\\
\dfrac{\sqrt{2\pi}}{\G(1/4)^2}&=\dfrac{1}{6}\sum_{n\ge0}\dfrac{(1/4)_n(1/2)_n}{(7/4)_nn!}\end{align*}
Parametric families for $v$ and $k$ nonnegative:
\begin{verbatim}
[()->gamma(1/4)^2/sqrt(2*Pi),4*(k-u)+2*v,2*(n+u)*(2*n-2*u+2*v-1)]
[()->gamma(1/4)^2/sqrt(2*Pi),4*(k+u)-2*v,2*(n+u)*(2*n-2*u+2*v-1)]
\end{verbatim}
Convergence type $P^-$ with $P=|2k+v-2u|$ and $P=|2k-v+2u|$.
\end{cf}

\smallskip

\begin{cf}\label{5.11.23}{\ }
\begin{verbatim}
[()->gamma(1/4)^2/sqrt(2*Pi),[4,4*n+1],[12,4*n*(n+1)*(2*n+1)*(2*n+3)]]
\end{verbatim}
$$\dfrac{\G(1/4)^2}{\sqrt{2\pi}}=4+\dfrac{12}{5+\dfrac{120}{9+\dfrac{840}{13+\dfrac{3024}{17+\dfrac{7920}{21+\dfrac{17160}{25+\ddots}}}}}}$$
Convergence type $P^-$ with $P=1$ and $C=\G(1/4)^4/(8\pi^2)$, so that
$$\dfrac{\G(1/4)^2}{\sqrt{2\pi}}-\dfrac{p(n)}{q(n)}\sim(-1)^n\dfrac{\G(1/4)^4/(8\pi^2)}{n}$$
$$A=1-(5/4)/n+(45/32)/n^2-(175/128)/n^3+\cdots$$
Series:
$$\dfrac{\sqrt{2\pi}}{\G(1/4)^2}=\dfrac{1}{4}\sum_{n\ge0}(-1)^n\dfrac{(2n+1)(1/2)_n^2}{(n+1)n!^2}$$
Parametric family for $k\ge0$:
\begin{verbatim}
[()->gamma(1/4)^2/sqrt(2*Pi),(4*k+1)*(4*n+1),4*n*(n+1)*(2*n+1)*(2*n+3)]
\end{verbatim}
Convergence type $P^-$ with $P=4k+1$.
\end{cf}

\smallskip

\begin{cf}\label{5.11.24}{\ }
\begin{verbatim}
[()->gamma(1/4)^2/sqrt(2*Pi),
[4,10,16*n^2-24*n+21],[12,-2*n*(2*n-1)*(4*n-3)*(4*n+1)]]
\end{verbatim}
$$\dfrac{\G(1/4)^2}{\sqrt{2\pi}}=4+\dfrac{12}{10-\dfrac{10}{37-\dfrac{540}{93-\dfrac{3510}{181-\dfrac{12376}{301-\dfrac{32130}{453-\ddots}}}}}}$$
Convergence type $P^+$ with $P=5/2$ and $C=3/(20\sqrt{\pi})$, so that
$$\dfrac{\G(1/4)^2}{\sqrt{2\pi}}-\dfrac{p(n)}{q(n)}\sim\dfrac{3/(20\sqrt{\pi})}{n^{5/2}}$$
$$A=1-(5/8)/n+(35/384)/n^2+(105/1024)/n^3+\cdots$$
Series:
$$\dfrac{\G(1/4)^2}{\sqrt{2\pi}}=4+6\sum_{n\ge0}\dfrac{(1/2)_n}{(4n+1)(4n+5)(n+1)!}$$
Parametric family for $k\ge0$:
\begin{verbatim}
[()->gamma(1/4)^2/sqrt(2*Pi),
16*n^2-24*n+21+4*k*(2*k+5),-2*n*(2*n-1)*(4*n-3)*(4*n+1)]
\end{verbatim}
Convergence type $P^+$ with $P=2k+5/2$.
\end{cf}
      
\smallskip

\begin{cf}\label{5.11.25}{\ }
\begin{verbatim}
[()->gamma(1/4)^2/sqrt(2*Pi),[0,6*n-4],[4,-2*n*(4*n-1)]]
\end{verbatim}
$$\dfrac{\G(1/4)^2}{\sqrt{2\pi}}=\dfrac{4}{2-\dfrac{6}{8-\dfrac{28}{14-\dfrac{66}{20-\dfrac{120}{26-\dfrac{190}{32-\ddots}}}}}}$$
Convergence type $E$ with $E=2$, $P=-1/4$, and $C=2\pi^2/\G(3/4)^3$, so that
$$\dfrac{\G(1/4)^2}{\sqrt{2\pi}}-\dfrac{p(n)}{q(n)}\sim\dfrac{2\pi^2/\G(3/4)^3}{2^nn^{-1/4}}\;.$$
$$A=1-(21/32)/n-(343/2048)/n^2-(189567/65536)/n^3+\cdots$$
Parametric family with $k\ge0$:
\begin{verbatim}
[()->gamma(1/4)^2/sqrt(2*Pi),6*n-4+2*k,-2*n*(4*n-1)]
\end{verbatim}
Convergence type $E$ with $E=2$ and $P=2k-1/4$.
\end{cf}

\smallskip

\begin{cf}\label{5.11.26}{\ }
\begin{verbatim}
[()->gamma(1/4)^2/sqrt(2*Pi),[4,5,6*n],[4,-(2*n+1)*(4*n+1)]]
\end{verbatim}
$$\dfrac{\G(1/4)^2}{\sqrt{2\pi}}=4+\dfrac{4}{5-\dfrac{15}{12-\dfrac{45}{18-\dfrac{91}{24-\dfrac{153}{30-\dfrac{231}{36-\ddots}}}}}}$$
Convergence type $E$ with $E=2$, $P=3/4$, and $C=\G(1/4)/\sqrt{\pi}$, so that
$$\dfrac{\G(1/4)^2}{\sqrt{2\pi}}-\dfrac{p(n)}{q(n)}\sim\dfrac{\G(1/4)/\sqrt{\pi}}{2^nn^{3/4}}\;.$$
$$A=1-(57/32)/n+(10241/2048)/n^2-(1320627/65536)/n^3+\cdots$$
Series:
$$\dfrac{\G(1/4)^2}{\sqrt{2\pi}}=4+\dfrac{4}{5}\sum_{n\ge0}\dfrac{(3/2)_n}{(9/4)_n}2^{-n}$$
Parametric family with $k\ge0$:
\begin{verbatim}
[()->gamma(1/4)^2/sqrt(2*Pi),6*n+2*k,-(2*n+1)*(4*n+1)]
\end{verbatim}
Convergence type $E$ with $E=2$ and $P=2k+3/4$.
\end{cf}

\smallskip

\begin{cf}\label{5.11.27}{\ }
\begin{verbatim}
[()->gamma(1/4)^2/sqrt(2*Pi),[16/3,40*n^2+22*n+3],
                             [-16/3,-32*n*(n+1)*(2*n+1)*(4*n+1)]]
\end{verbatim}
$$\dfrac{\G(1/4)^2}{\sqrt{2\pi}}=16/3-\dfrac{16/3}{65-\dfrac{960}{207-\dfrac{8640}{429-\dfrac{34944}{731-\dfrac{97920}{1113-\dfrac{221760}{1575-\ddots}}}}}}$$
Convergence type $E$ with $E=4$, $P=9/4$, and $C=-\G(1/4)^3/(32\pi^{3/2})$,
so that
$$\dfrac{\G(1/4)^2}{\sqrt{2\pi}}-\dfrac{p(n)}{q(n)}\sim-\dfrac{\G(1/4)^3/\pi^{3/2}}{2^{2n+5}n^{9/4}}$$
$$A=1-(135/32)/n+(29985/2048)/n^2+\cdots$$
Series:
$$\dfrac{\sqrt{2\pi}}{\G(1/4)^2}=\dfrac{3}{16}\sum_{n\ge0}\dfrac{(1/2)_n(1/4)_n}{n!(n+1)!}2^{-2n}$$
Parametric family for $k\ge0$:
\begin{verbatim}
[()->gamma(1/4)^2/sqrt(2*Pi),40*n^2+(16*k+22)*n+2*k+3,
                            -16*(n+1)*(2*n-k)*(2*n+1-k)*(4*n+1)]
\end{verbatim}
Convergence type $E$ with $E=4$ and $P=3k+9/4$.
\end{cf}

\smallskip

\begin{cf}\label{5.11.28}{\ }
\begin{verbatim}
[()->gamma(1/4)^2/sqrt(2*Pi),[4,10,(2*n-1)*(20*n-9)],
                             [12,-2*(2*n-1)*(4*n-3)*(4*n-1)*(4*n+1)]]
\end{verbatim}
$$\dfrac{\G(1/4)^2}{\sqrt{2\pi}}=4+\dfrac{12}{10-\dfrac{30}{93-\dfrac{1890}{255-\dfrac{12870}{497-\dfrac{46410}{819-\dfrac{122094}{1221-\ddots}}}}}}$$
Convergence type $E$ with $E=4$, $P=9/4$, and $C=\pi^{3/2}/(4\G(3/4)^3$, so that
$$\dfrac{\G(1/4)^2}{\sqrt{2\pi}}-\dfrac{p(n)}{q(n)}\sim\dfrac{\pi^{3/2}/\G(3/4)^3}{2^{2n+2}n^{9/4}}$$
$$A=1-(81/32)/n+(13137/2048)/n^2+\cdots$$
Series:
$$\dfrac{\G(1/4)^2}{\sqrt{2\pi}}=4+\dfrac{6}{5}\sum_{n\ge0}\dfrac{(3/4)_n(1/2)_n}{(5/4)_n(9/4)_n}2^{-2n}$$
Parametric family for $k\ge0$:
\begin{verbatim}
[()->gamma(1/4)^2/sqrt(2*Pi),40*n^2+(16*k-38)*n+9-10*k,
                             -2*(2*n-1)*(4*n+1)*(4*n-1-2*k)*(4*n-3-2*k)]
\end{verbatim}
Convergence type $E$ with $E=4$ and $P=3k+9/4$.
\end{cf}

\smallskip

\begin{cf}\label{5.11.29}{\ }
\begin{verbatim}
[()->gamma(1/4)^2/sqrt(2*Pi),[0,4*n-3],[4,-n*(2*n-1)]]
\end{verbatim}
$$\dfrac{\G(1/4)^2}{\sqrt{2\pi}}=\dfrac{4}{1-\dfrac{1}{5-\dfrac{6}{9-\dfrac{15}{13-\dfrac{28}{17-\dfrac{45}{21-\ddots}}}}}}$$
Convergence type $E$ with $E=(1+\sqrt{2})^2$, $P=0$, and $C=2\G(1/4)^2/(\sqrt{\pi}(1+\sqrt{2})^{1/2})$, so that
$$\dfrac{\G(1/4)^2}{\sqrt{2\pi}}-\dfrac{p(n)}{q(n)}\sim\dfrac{2\G(1/4)^2/\sqrt{\pi}}{(1+\sqrt{2})^{2n+1/2}}\;.$$
$$A=1-(3d/16)/n+(3d/64+9/256)/n^2-(429d/4096+9/512)/n^3+\cdots$$
\end{cf}

\smallskip

\begin{cf}\label{5.11.30}{\ }
\begin{verbatim}
[()->gamma(1/4)^2/sqrt(2*Pi),
[4*n+5],[[(4*n+1)*(4*n+3),(4*n+6)*(4*n+8)]]]
\end{verbatim}
$$\dfrac{\G(1/4)^2}{\sqrt{2\pi}}=5+\dfrac{3}{9+\dfrac{48}{13+\dfrac{35}{17+\dfrac{120}{21+\dfrac{99}{25+\dfrac{224}{29+\ddots}}}}}}$$
Convergence type $E$ with $E=-(1+\sqrt{2})^2$, $P=0$, and
$C=\sqrt{8\pi}(\G(1/4)/\G(3/4))/(1+\sqrt{2})^{9/2}$, so that
$$\dfrac{\G(1/4)^2}{\sqrt{2\pi}}-\dfrac{p(n)}{q(n)}\sim(-1)^n\dfrac{\sqrt{8\pi}\G(1/4)/\G(3/4)}{(1+\sqrt{2})^{2n+9/2}}\;.$$
$$A=1-(3d/16)/n+(27d/64+9/256)/n^2-(4269d/4096+81/512)/n^3+\cdots$$
\end{cf}

\smallskip

\begin{cf}\label{5.11.31}{\ }
\begin{verbatim}
[()->gamma(1/4)^2/sqrt(2*Pi),[26/5,88*n^2+44*n+3],[6,2*n*(2*n+3)*(4*n+1)*(4*n+5)]]
\end{verbatim}
$$\dfrac{\G(1/4)^2}{\sqrt{2\pi}}=\dfrac{26}{5}+\dfrac{6}{135+\dfrac{450}{443+\dfrac{3276}{927+\dfrac{11934}{1587+\dfrac{31416}{2423+\dfrac{68250}{3435+\ddots}}}}}}$$
Convergence type $E$ with $E=-((1+\sqrt{5})/2)^{10}$, $P=0$, and $C=2\G(1/4)^2((1+\sqrt{5})/2)^{-25/2}/\sqrt{\pi}$, so that
$$\dfrac{\G(1/4)^2}{\sqrt{2\pi}}-\dfrac{p(n)}{q(n)}\sim(-1)^n\dfrac{2\G(1/4)^2/\sqrt{\pi}}{((1+\sqrt{5})/2)^{10n+25/2}}$$
$$A=1+(9d/80)/n+(-9d/64+81/2560)/n^2+\cdots$$
\end{cf}

\medskip

For $m=1$, note that

$$\dfrac{\G(3/4)^2}{\sqrt{2\pi}}=\sqrt{\pi}\dfrac{\G(3/4)}{\G(1/4)}=L_2\;,$$
and since $G=\varpi/\pi=1/(2L_2)$ the following CFs also trivially give CFs for
$G$ and $\varpi/\pi$ which we do not include.

\smallskip

\begin{cf}\label{5.11.32}{\ }
\begin{verbatim}
[()->gamma(3/4)^2/sqrt(2*Pi),[1,3,2],[-2,(4*n-1)*(4*n+1)]]
\end{verbatim}
$$\dfrac{\G(3/4)^2}{\sqrt{2\pi}}=1-\dfrac{2}{3+\dfrac{15}{2+\dfrac{63}{2+\dfrac{143}{2+\dfrac{255}{2+\dfrac{399}{2+\ddots}}}}}}$$
Convergence type $P^-$ with $P=1/2$ and $C=-\pi\sqrt{2}/\G(1/4)^2$, so that
$$\dfrac{\G(3/4)^2}{\sqrt{2\pi}}-\dfrac{p(n)}{q(n)}\sim(-1)^{n+1}\dfrac{\pi\sqrt{2}/\G(1/4)^2}{n^{1/2}}$$
$$A=1-(1/4)/n-(1/64)/n^2+(25/256)/n^3+21/8192/n^4+\cdots$$
Series:
$$\dfrac{\G(3/4)^2}{\sqrt{2\pi}}=-1+2\sum_{n\ge0}(-1)^n\dfrac{(1/4)_n}{(3/4)_n}$$
Parametric families for $v$, $w$ and $k$ nonnegative:
\begin{verbatim}
[()->gamma(3/4)^2/sqrt(2*Pi),
8*k+4*w+2,(4*n+4*u-3)*(4*n+4*u-8*v+4*w-5)]
[()->gamma(3/4)^2/sqrt(2*Pi),
8*k+4*w+2,(4*n+4*u-1)*(4*n+4*u+8*v-4*w+1)]
\end{verbatim}
Convergence types $P^-$ with $P=2k+w+1/2$.
\end{cf}

\smallskip

\begin{cf}\label{5.11.33}{\ }
\begin{verbatim}
[()->gamma(3/4)^2/sqrt(2*Pi),[1,2],[-1,2*n*(2*n-1)]]
\end{verbatim}
$$\dfrac{\G(3/4)^2}{\sqrt{2\pi}}=1-\dfrac{1}{2+\dfrac{2}{2+\dfrac{12}{2+\dfrac{30}{2+\dfrac{56}{2+\dfrac{90}{2+\ddots}}}}}}$$
Convergence type $P^-$ with $P=1$ and $C=-\G(3/4)^2/\sqrt{32\pi}$, so that
$$\dfrac{\G(3/4)^2}{\sqrt{2\pi}}-\dfrac{p(n)}{q(n)}\sim(-1)^{n+1}\dfrac{\G(3/4)^2/\sqrt{32\pi}}{n}$$
$$A=1-(1/4)/n-(11/64)/n^2+(41/256)/n^3+\cdots$$
Series:
\begin{align*}\dfrac{\G(3/4)^2}{\sqrt{2\pi}}&=1-\dfrac{1}{3}\sum_{n\ge0}\dfrac{n!(1/4)_n}{(3/2)_n(7/4)_n}\\
\dfrac{\sqrt{2\pi}}{\G(3/4)^2}&=2\sum_{n\ge0}\dfrac{(1/2)_n(-1/4)_n}{n!(5/4)_n}\end{align*}
Parametric families for all parameters nonnegative:
\begin{verbatim}
[()->gamma(3/4)^2/sqrt(2*Pi),4*(k+u)-2*v+2,2*(n+u)*(2*n-2*u+2*v-1)]
[()->gamma(3/4)^2/sqrt(2*Pi),4*(k-u)+2*v-2,2*(n+u)*(2*n-2*u+2*v-1)]
\end{verbatim}
Convergence type $P^-$ with $P=|2k+2u-v+1|$ and $P=|2k-2u+v-1|$.
\end{cf}

\smallskip

\begin{cf}\label{5.11.34}{\ }
\begin{verbatim}
[()->gamma(3/4)^2/sqrt(2*Pi),[1/2,4*n-1],[1/2,4*n^2*(2*n+1)^2]]
\end{verbatim}
$$\dfrac{\G(3/4)^2}{\sqrt{2\pi}}=1/2+\dfrac{1/2}{3+\dfrac{36}{7+\dfrac{400}{11+\dfrac{1764}{15+\dfrac{5184}{19+\dfrac{12100}{23+\ddots}}}}}}$$
Convergence type $P^-$ with $P=1$ and $C=(\G(3/4)/\G(1/4))^2$, so that
$$\dfrac{\G(3/4)^2}{\sqrt{2\pi}}-\dfrac{p(n)}{q(n)}\sim(-1)^n\dfrac{(\G(3/4)/\G(1/4))^2}{n}\;.$$
$$A=1-(3/4)/n+(9/32)/n^2+(27/128)/n^3-(549/2048)/n^4+\cdots$$
Series:
$$\dfrac{\sqrt{2\pi}}{\G(3/4)^2}=2\sum_{n\ge0}(-1)^n\dfrac{(1/2)_n^2}{n!^2}=\sum_{n\ge0}(-1)^n\dfrac{(2n)!^2}{2^{4n-1}n!^4}$$
Parametric family for $k\ge0$:
\begin{verbatim}
[()->gamma(3/4)^2/sqrt(2*Pi),(4*k+1)*(4*n-1),4*n^2*(2*n+1)^2]
\end{verbatim}
Convergence type $P^-$ with $P=4k+1$.
\end{cf}

\smallskip

\begin{cf}\label{5.11.35}{\ }
\begin{verbatim}
[()->gamma(3/4)^2/sqrt(2*Pi),
[1/2,6,16*n^2-24*n+15],[1/2,-2*n*(2*n-1)*(4*n-1)^2]]
\end{verbatim}
$$\dfrac{\G(3/4)^2}{\sqrt{2\pi}}=1/2+\dfrac{1/2}{6-\dfrac{18}{31-\dfrac{588}{87-\dfrac{3630}{175-\dfrac{12600}{295-\dfrac{32490}{447-\ddots}}}}}}$$
Convergence type $P^+$ with $P=3/2$ and $C=1/(24\sqrt{\pi})$, so that
$$\dfrac{\G(3/4)^2}{\sqrt{2\pi}}-\dfrac{p(n)}{q(n)}\sim\dfrac{1/(24\sqrt{\pi})}{n^{3/2}}\;.$$
$$A=1-(3/8)/n-(5/896)/n^2+(75/1024)/n^3+\cdots$$
Series:
$$\dfrac{\G(3/4)^2}{\sqrt{2\pi}}=\dfrac{1}{2}+\dfrac{1}{4}\sum_{n\ge0}\dfrac{(1/2)_n}{(4n+3)(n+1)!}$$
Parametric family for $k\ge0$:
\begin{verbatim}
[()->gamma(3/4)^2/sqrt(2*Pi),
16*n^2-24*n+15+8*k^2+12*k,-2*n*(2*n-1)*(4*n-1)^2]
\end{verbatim}
Convergence type $P^+$ with $P=2k+3/2$.
\end{cf}

\smallskip

\begin{cf}\label{5.11.36}{\ }
\begin{verbatim}
[()->gamma(3/4)^2/sqrt(2*Pi),[0,6*n-4],[-1,-2*n*(4*n+1)]]
\end{verbatim}
$$\dfrac{\G(3/4)^2}{\sqrt{2\pi}}=-\dfrac{1}{2-\dfrac{10}{8-\dfrac{36}{14-\dfrac{78}{20-\dfrac{136}{26-\dfrac{210}{32-\ddots}}}}}}$$
Convergence type $E$ with $E=2$, $P=-7/4$, and $C=-8\pi^2/\G(1/4)^3$, so that
$$\dfrac{\G(3/4)^2}{\sqrt{2\pi}}-\dfrac{p(n)}{q(n)}\sim-\dfrac{\pi^2/\G(1/4)^3}{2^{n-3}n^{-7/4}}\;.$$
$$A=1-(61/32)/n-(8231/2048)/n^2-(1275655/65536)/n^3+\cdots$$
Parametric family with $k\ge0$:
\begin{verbatim}
[()->gamma(3/4)^2/sqrt(2*Pi),6*n-4+2*k,-2*n*(4*n+1)]
\end{verbatim}
Convergence type $E$ with $E=2$ and $P=2k-7/4$.
\end{cf}

\smallskip

\begin{cf}\label{5.11.37}{\ }
\begin{verbatim}
[()->gamma(3/4)^2/sqrt(2*Pi),[1,3,6*n-4],[-(2*n-1)*(4*n-1)]]
\end{verbatim}
$$\dfrac{\G(3/4)^2}{\sqrt{2\pi}}=1-\dfrac{1}{3-\dfrac{3}{8-\dfrac{21}{14-\dfrac{55}{20-\dfrac{105}{26-\dfrac{171}{32-\ddots}}}}}}$$
Convergence type $E$ with $E=2$, $P=5/4$ and $C=-\G(3/4)/(2\sqrt{\pi})$, so
that
$$\dfrac{\G(3/4)^2}{\sqrt{2\pi}}-\dfrac{p(n)}{q(n)}\sim-\dfrac{\G(3/4)/(2\sqrt{\pi})}{2^nn^{5/4}}\;.$$
$$A=1-(65/32)/n+(13425/2048)/n^2-(2011035/65536)/n^3+\cdots$$
Series:
$$\dfrac{\G(3/4)^2}{\sqrt{2\pi}}=1-\dfrac{1}{3}\sum_{n\ge0}\dfrac{(1/2)_n}{(7/4)_n}2^{-n}$$
Parametric family with $k\ge0$:
\begin{verbatim}
[()->gamma(3/4)^2/sqrt(2*Pi),6*n-4+2*k,-(2*n-1)*(4*n-1)]
\end{verbatim}
Convergence type $E$ with $E=2$ and $P=2k+5/4$.
\end{cf}

\smallskip

\begin{cf}\label{5.11.38}{\ }
\begin{verbatim}
[()->gamma(3/4)^2/sqrt(2*Pi),[2/3,40*n^2-6*n+1],
                             [-2,-32*n^2*(2*n+1)*(4*n+3)]]
\end{verbatim}
$$\dfrac{\G(3/4)^2}{\sqrt{2\pi}}=2/3-\dfrac{2}{35-\dfrac{672}{149-\dfrac{7040}{343-\dfrac{30240}{617-\dfrac{87552}{971-\dfrac{202400}{1405-\ddots}}}}}}$$
Convergence type $E$ with $E=4$, $P=3/4$, and $C=-\G(3/4)^3/(4\pi^{3/2})$,
so that
$$\dfrac{\G(3/4)^2}{\sqrt{2\pi}}-\dfrac{p(n)}{q(n)}\sim-\dfrac{\G(3/4)^3/\pi^{3/2}}{2^{2n+2}n^{3/4}}$$
$$A=1-(39/32)/n+(4097/2048)/n^2+\cdots$$
Series:
$$\dfrac{\sqrt{2\pi}}{\G(3/4)^2}=\dfrac{3}{2}\sum_{n\ge0}\dfrac{(1/2)_n(3/4)_n}{n!^2}2^{-2n}$$
Parametric family for $k\ge0$:
\begin{verbatim}
[()->gamma(3/4)^2/sqrt(2*Pi),40*n^2+(16*k-6)*n+1-2*k,
                            -16*n*(2*n-k)*(2*n+1-k)*(4*n+3)]
\end{verbatim}
Convergence type $E$ with $E=4$ and $P=3k+3/4$.
\end{cf}

\smallskip

\begin{cf}\label{5.11.39}{\ }
\begin{verbatim}
[()->gamma(3/4)^2/sqrt(2*Pi),[1/2,18,40*n^2-26*n+5],
                             [3/2,-2*(2*n+1)*(4*n-1)^2*(4*n+1)]]
\end{verbatim}
$$\dfrac{\G(3/4)^2}{\sqrt{2\pi}}=\dfrac{1}{2}+\dfrac{3/2}{18-\dfrac{270}{113-\dfrac{4410}{287-\dfrac{22022}{541-\dfrac{68850}{875-\dfrac{166782}{1289-\ddots}}}}}}$$
Convergence type $E$ with $E=4$, $P=3/4$, and $C=\pi^{3/2}/\G(1/4)^3$, so that
$$\dfrac{\G(3/4)^2}{\sqrt{2\pi}}-\dfrac{p(n)}{q(n)}\sim\dfrac{\pi^{3/2}/\G(1/4)^3}{2^{2n}n^{3/4}}$$
$$A=1-(33/32)/n+(3089/2048)/n^2+\cdots$$
Series:
$$\dfrac{\G(3/4)^2}{\sqrt{2\pi}}=-1+\dfrac{3}{2}\sum_{n\ge0}\dfrac{(1/4)_n(1/2)_n}{(3/4)_n^2}2^{-2n}$$
Parametric family for $k\ge0$:
\begin{verbatim}
[()->gamma(3/4)^2/sqrt(2*Pi),40*n^2+(16*k-26)*n+5-6*k,
                             -2*(2*n+1)*(4*n-1)*(4*n-1-2*k)*(4*n+1-2*k)]
\end{verbatim}
Convergence type $E$ with $E=4$ and $P=3k+3/4$.
\end{cf}

\smallskip

\begin{cf}\label{5.11.40}{\ }
\begin{verbatim}
[()->gamma(3/4)^2/sqrt(2*Pi),[1,4*n-1],[-1,-n*(2*n+1)]]
\end{verbatim}
$$\dfrac{\G(3/4)^2}{\sqrt{2\pi}}=1-\dfrac{1}{3-\dfrac{3}{7-\dfrac{10}{11-\dfrac{21}{15-\dfrac{36}{19-\dfrac{55}{23-\ddots}}}}}}$$
Convergence type $E$ with $E=(1+\sqrt{2})^2$, $P=0$, and $C=-2\G(3/4)^2/(\sqrt{\pi}(1+\sqrt{2})^{3/2})$, so that
$$\dfrac{\G(3/4)^2}{\sqrt{2\pi}}-\dfrac{p(n)}{q(n)}\sim(-1)^{n+1}\dfrac{2\G(3/4)^2/\sqrt{\pi}}{(1+\sqrt{2})^{2n+3/2}}$$
$$A=1-(3d/16)/n+(9d/64+9/256)/n^2+\cdots$$
\end{cf}
      
\smallskip

\begin{cf}\label{5.11.41}{\ }
\begin{verbatim}
[()->gamma(3/4)^2/sqrt(2*Pi),
[1,4,4*n-1],[[-2,8],[(4*n-1)*(4*n+1),8*(n+1)*(2*n+1)]]]
\end{verbatim}
$$\dfrac{\G(3/4)^2}{\sqrt{2\pi}}=1-\dfrac{2}{4+\dfrac{8}{7+\dfrac{15}{11+\dfrac{48}{15+\dfrac{63}{19+\dfrac{120}{23+\ddots}}}}}}$$
Convergence type $E$ with $E=-(1+\sqrt{2})^2$, $P=0$, and
$C=-2\G(3/4)^2/(\sqrt{\pi}(1+\sqrt{2})^{3/2})$, so that
$$\dfrac{\G(3/4)^2}{\sqrt{2\pi}}-\dfrac{p(n)}{q(n)}\sim(-1)^{n+1}\dfrac{2\G(3/4)^2/\sqrt{\pi}}{(1+\sqrt{2})^{2n+3/2}}$$
$$A=1-(3d/16)/n+(9d/64+9/256)/n^2+\cdots$$
\end{cf}

\smallskip

\begin{cf}\label{5.11.42}{\ }
\begin{verbatim}
[()->gamma(3/4)^2/sqrt(2*Pi),
             [0,4*n-5],[[-1,-1],[-2*n*(4*n+1),-(2*n-1)*(4*n-1)]]]
\end{verbatim}
$$\dfrac{\G(3/4)^2}{\sqrt{2\pi}}=-\dfrac{1}{-1-\dfrac{1}{3-\dfrac{10}{7-\dfrac{3}{11-\dfrac{36}{15-\dfrac{21}{19-\ddots}}}}}}$$
Convergence type $E$ with $E=(1+\sqrt{2})^2$, $P=0$, and $C=-2\G(3/4)^2(1+\sqrt{2})^{1/2}/\pi^{1/2}$, so that
$$\dfrac{\G(3/4)^2}{\sqrt{2\pi}}-\dfrac{p(n)}{q(n)}\sim-\dfrac{2\G(3/4)^2/\pi^{1/2}}{(1+\sqrt{2})^{2n-1/2}}\;.$$
$$A=1+(13d/16)/n+(13d/64+169/256)/n^2+(99d/4096+169/512)/n^3+\cdots$$
\end{cf}

\smallskip

\begin{cf}\label{5.11.43}{\ }
\begin{verbatim}
[()->gamma(3/4)^2/sqrt(2*Pi),[1/2,880*n^3-1408*n^2+518*n-35],
[-9/2,2*n*(2*n+1)*(4*n-1)^2*(10*n-11)*(10*n+9)]]
\end{verbatim}
$$\dfrac{\G(3/4)^2}{\sqrt{2\pi}}=1/2-\dfrac{9/2}{-45-\dfrac{1026}{2409+\dfrac{255780}{12607+\dfrac{3765762}{35829+\dfrac{23020200}{77355+\ddots}}}}}$$
Convergence type $E$ with $E=-((1+\sqrt{5})/2)^{10}$, $P=0$, and
$C=\sqrt{8\pi}(\G(3/4)/\G(1/4))/((1+\sqrt{5})/2)^{11/2}$, so that
$$\dfrac{\G(3/4)^2}{\sqrt{2\pi}}-\dfrac{p(n)}{q(n)}\sim(-1)^n\dfrac{\sqrt{8\pi}\G(3/4)/\G(1/4)}{((1+\sqrt{5})/2)^{10n+11/2}}\;,$$
$$A=1-(3d/80)/n+(39d/1600+9/2560)/n^2+\cdots$$
\end{cf}

\smallskip

\begin{cf}\label{5.11.44}{\ }
\begin{verbatim}
[()->gamma(3/4)^2/sqrt(2*Pi),[-1,88*n^2-220*n+135],[5,2*n*(2*n-3)*(4*n-5)*(4*n-1)]]
\end{verbatim}
$$\dfrac{\G(3/4)^2}{\sqrt{2\pi}}=-1+\dfrac{5}{3+\dfrac{6}{47+\dfrac{84}{267+\dfrac{1386}{663+\dfrac{6600}{1235+\dfrac{19950}{1983+\ddots}}}}}}$$
Convergence type $E$ with $E=-((1+\sqrt{5})/2)^{10}$, $P=0$, and $C=2\G(3/4)^2((1+\sqrt{5})/2)^{5/2}/\sqrt{\pi}$, so that
$$\dfrac{\G(3/4)^2}{\sqrt{2\pi}}-\dfrac{p(n)}{q(n)}\sim(-1)^n\dfrac{2\G(3/4)^2/\sqrt{\pi}}{((1+\sqrt{5})/2)^{10n-5/2}}$$
$$A=1+(9d/80)/n+(9d/320+81/2560)/n^2+\cdots$$
\end{cf}

\medskip

{\bf $k\ge2$:}

\medskip

As mentioned above, CFs for $2^{j/2}\pi^{1/2-k}\G((2m+1)/4)^2$ can be
trivially deduced from those for $2^{j/2}\pi^{1/2-(3-k)}\G((3-2m)/4)^2$,
so we do not give them explicitly but simply refer to the corresponding CFs.

\begin{itemize}\item We have
  $$\pi^{-3/2}\G(1/4)^2=\dfrac{2}{\pi^{-1/2}\G(3/4)^2}\text{\quad and\quad}\pi^{-3/2}\G(3/4)^2=\dfrac{2}{\pi^{-1/2}\G(1/4)^2}\;,$$
  so invert the CFs beginning at \ref{5.11.12} and \ref{5.11.3} respectively.
\item We have
  $$(2\pi)^{-3/2}\G(1/4)^2=\dfrac{1/2}{(2\pi)^{-1/2}\G(3/4)^2}\text{\quad and\quad}
  (2\pi)^{-3/2}\G(3/4)^2=\dfrac{1/2}{(2\pi)^{-1/2}\G(1/4)^2}\;,$$
  so invert the CFs beginning at \ref{5.11.32} and \ref{5.11.21} respectively.
\item We have
  $$(2\pi)^{-5/2}\G(1/4)^2=\dfrac{1/2}{(2\pi)^{1/2}\G(3/4)^2}\text{\quad and\quad}
  (2\pi)^{-5/2}\G(3/4)^2=\dfrac{1/2}{(2\pi)^{1/2}\G(1/4)^2}\;,$$
  so invert \ref{5.11.2} and \ref{5.11.1} respectively.
\end{itemize}

\medskip

\subsection{CFs for $\pi^{-j}\G((2m+1)/4)^4$}

\medskip

We now give CFs for $\G(1/4)^4/\pi^j$ and $\G(3/4)^4/\pi^j$ for $0\le j\le 4$.
Note that $\G(3/4)^4/\pi^j=4/(\G(1/4)^4/\pi^{4-j})$, so CFs for
$\G(3/4)^4/\pi^j$ can trivially be deduced from those for $\G(1/4)^4/\pi^j$
and conversely, so we will only mention the latter. Also, note that
$(\G(1/4)/\G(3/4))^2=(1/2)(\G(1/4)^4/\pi^2)$
so that CFs for the LHS can trivially be deduced from CFs for the RHS
and conversely, but since the LHS is more symmetrical we will in fact give the
CFs for $(\G(1/4)/\G(3/4))^2$.

\medskip

{\bf $j=0$: Continued Fractions for $\G((2m+1)/4)^4$}

\medskip

\begin{cf}\label{5.11.45}{\ }
\begin{verbatim}
[()->gamma(1/4)^4,[512/3,3*(4*n+1)*(48*n^4+48*n^3+8*n^2-2*n-1)],
                  [3584,32*n*(n+1)^4*(2*n+1)^4*(2*n+3)*(4*n-1)*(4*n+7)]]
\end{verbatim}
$$\G(1/4)^4=512/3+\dfrac{3584}{1515+\dfrac{6842880}{31833+\dfrac{2381400000}{204711+\dfrac{110992121856}{789429+\dfrac{1991919600000}{2279907+\ddots}}}}}$$
Convergence type $P^-$ with $P=9/2$ and $C=\G(1/4)^8/(128\pi^{5/2})$, so that
$$\G(1/4)^4-\dfrac{p(n)}{q(n)}\sim(-1)^n\dfrac{\G(1/4)^8/(128\pi^{5/2})}{n^{9/2}}$$
$$A=1-(45/8)/n+(2121/128)/n^2-(30615/1024)/n^3+\cdots$$
Series:
$$\dfrac{1}{\G(1/4)^4}=\dfrac{1}{512}\sum_{n\ge0}(-1)^n\dfrac{(2n+1)(4n+3)(1/2)_n^5}{n!(n+1)!^4}$$
\end{cf}

\smallskip

\begin{cf}\label{5.11.46}{\ }
\begin{verbatim}
[()->gamma(3/4)^4,[2,(4*n-1)*(48*n^4-48*n^3+8*n^2+2*n-1)],
                  [10,32*n^5*(2*n+1)^5*(4*n-3)*(4*n+5)]]
\end{verbatim}
$$\G(3/4)^4=2+\dfrac{10}{27+\dfrac{69984}{2933+\dfrac{208000000}{29359+\dfrac{19995758496}{140265+\dfrac{528232513536}{459971+\ddots}}}}}$$
Convergence type $P^-$ with $P=3/2$ and $C=\G(3/4)^8/\pi^{5/2}$, so that
$$\G(3/4)^4-\dfrac{p(n)}{q(n)}\sim(-1)^n\dfrac{\G(3/4)^8/\pi^{5/2}}{n^{3/2}}$$
$$A=1-(9/8)/n+(65/128)/n^2+(525/1024)/n^3+\cdots$$
Series:
$$\dfrac{1}{\G(3/4)^4}=\dfrac{1}{2}\sum_{n\ge0}(-1)^n\dfrac{(4n+1)(1/2)_n^5}{n!^5}$$
\end{cf}

\medskip

{\bf $j=1$: Continued Fractions for $\pi^{-1}\G((2m+1)/4)^4$}

\medskip

For $m=0$:

\smallskip

\begin{cf}\label{5.11.47}{\ }
\begin{verbatim}
[()->gamma(1/4)^4/Pi,[32,15,16*n^2+2],[96,-(2*n+1)^2*(4*n+1)*(4*n+3)]]
\end{verbatim}
$$\dfrac{\G(1/4)^4}{\pi}=32+\dfrac{96}{15-\dfrac{315}{66-\dfrac{2475}{146-\dfrac{9555}{258-\dfrac{26163}{402-\dfrac{58443}{578-\ddots}}}}}}$$
Convergence type $P^+$ with $P=1/2$ and $C=2^{5/2}\G(1/4)^2/\pi$, so that
$$\dfrac{\G(1/4)^4}{\pi}-\dfrac{p(n)}{q(n)}\sim\dfrac{2^{5/2}\G(1/4)^2/\pi}{n^{1/2}}$$
$$A=1-(1/2)/n+(109/320)/n^2-(29/128)/n^3+\cdots$$
Series:
$$\dfrac{\G(1/4)^4}{\pi}=32+\dfrac{96}{5}\sum_{n\ge0}\dfrac{(7/4)_n}{(2n+3)(9/4)_n}$$
Parametric family for $k\ge0$:
\begin{verbatim}
[()->gamma(1/4)^4/Pi,16*n^2+2+4*k*(2*k+1),-(2*n+1)^2*(4*n+1)*(4*n+3)]
\end{verbatim}
Convergence type $P^+$ with $P=2k+1/2$.
\end{cf}

\medskip

\begin{cf}\label{5.11.48}{\ }
\begin{verbatim}
[()->gamma(1/4)^4/Pi,[32,8*n-6],[48,n^2*(2*n-1)^2]]
\end{verbatim}
$$\dfrac{\G(1/4)^4}{\pi}=32+\dfrac{48}{2+\dfrac{1}{10+\dfrac{36}{18+\dfrac{225}{26+\dfrac{784}{34+\dfrac{2025}{42+\ddots}}}}}}$$
Convergence type $P^-$ with $P=4$ and $C=3\G(1/4)^4/(32\pi)$, so that
$$\dfrac{\G(1/4)^4}{\pi}-\dfrac{p(n)}{q(n)}\sim(-1)^n\dfrac{3\G(1/4)^4/(32\pi)}{n^4}$$
$$A=1-1/n-(7/4)/n^2+(13/4)/n^3+(231/32)/n^4+\cdots$$
Series:
\begin{align*}
  \dfrac{\G(1/4)^4}{\pi}&=32+96\sum_{n\ge0}\dfrac{(8n+5)n!^2(3/4)_n^2}{(16n^2+4n+1)(16n^2+36n+21)(3/2)_n^2(5/4)_n^2}\\
  \dfrac{\pi}{\G(1/4)^4}&=\dfrac{3}{8}\sum_{n\ge0}\dfrac{(8n+1)(1/2)_n^2(1/4)_n^2}{(16n^2-12n+3)(16n^2+20n+7)n!^2(3/4)_n^2}
\end{align*}
Parametric family for $k\ge0$:
\begin{verbatim}
[()->gamma(1/4)^4/Pi,(k+1)*(8*n-6),n^2*(2*n-1)^2]
\end{verbatim}
Convergence type $P^-$ with $P=4(k+1)$.
\end{cf}

\smallskip

\begin{cf}\label{5.11.49}
\begin{verbatim}
[()->gamma(1/4)^4/Pi,[48,(4*n-3)*(4*n^2-6*n+9)],
                                [48,-8*n^3*(2*n-1)^3]]
\end{verbatim}
$$\dfrac{\G(1/4)^4}{\pi}=48+\dfrac{48}{7-\dfrac{8}{65-\dfrac{1728}{243-\dfrac{27000}{637-\dfrac{175616}{1343-\dfrac{729000}{2457-\ddots}}}}}}$$
Convergence type $P^+$ with $P=7/2$ and $C=\pi^{9/2}/(84\G(3/4)^8)$, so that
$$\dfrac{\G(1/4)^4}{\pi}-\dfrac{p(n)}{q(n)}\sim\dfrac{\pi^{9/2}/(84\G(3/4)^8)}{n^{7/2}}$$
$$A=1-(7/8)/n-(21/1408)/n^2+(483/1024)/n^3+\cdots$$
Series:
$$\dfrac{\pi}{\G(1/4)^4}=\dfrac{1}{48}\sum_{n\ge0}\dfrac{(-1/2)_n^3}{n!^3}$$
Parametric family for $k\ge0$:
\begin{verbatim}
[()->gamma(1/4)^4/Pi,(4*n-3)*(4*n^2-6*n+9+2*k*(4*k+7)),
                     -8*n^3*(2*n-1)^3]
\end{verbatim}
Convergence type $P^+$ with $P=4k+7/2$.
\end{cf}

\smallskip

\begin{cf}\label{5.11.50}{\ }
\begin{verbatim}
[()->gamma(1/4)^4/Pi,[48,25,20*n^2+7*n+1],
                     [144,-(n+1)*(4*n+1)^2*(4*n+3)]]
\end{verbatim}
$$\dfrac{\G(1/4)^4}{\pi}=48+\dfrac{144}{25-\dfrac{350}{95-\dfrac{2673}{202-\dfrac{10140}{349-\dfrac{27455}{536-\dfrac{60858}{763-\ddots}}}}}}$$
Convergence type $E$ with $E=4$, $P=3/4$, and $C=\G(1/4)^2/\G(3/4)$, so that
$$\dfrac{\G(1/4)^4}{\pi}-\dfrac{p(n)}{q(n)}\sim\dfrac{\G(1/4)^2/\G(3/4)}{2^{2n}n^{3/4}}$$
$$A=1-(45/32)/n+(5273/2048)/n^2+\cdots$$
Series:
$$\dfrac{\G(1/4)^4}{\pi}=48\sum_{n\ge0}\dfrac{n!(3/4)_n}{(5/4)_n^2}2^{-2n}$$
Parametric family for $k\ge0$:
\begin{verbatim}
[()->gamma(1/4)^4/Pi,20*n^2+(8*k+7)*n+k+1,
                     -(n+1)*(4*n+1)*(4*n+1-2*k)*(4*n+3-2*k)]
\end{verbatim}
Convergence type $E$ with $E=4$ and $P=3k+3/4$.
\end{cf}

\smallskip

For $m=1$:

\smallskip

\begin{cf}\label{5.11.51}
\begin{verbatim}
[()->gamma(3/4)^4/Pi,[1,(4*n-1)*(4*n^2-2*n+1)],
                               [-1,-8*n^3*(2*n+1)^3]]
\end{verbatim}
$$\dfrac{\G(3/4)^4}{\pi}=1-\dfrac{1}{9-\dfrac{216}{91-\dfrac{8000}{341-\dfrac{74088}{855-\dfrac{373248}{1729-\dfrac{1331000}{3059-\ddots}}}}}}$$
Convergence type $P^+$ with $P=1/2$ and $C=-2\G(3/4)^8/\pi^{7/2}$, so that
$$\dfrac{\G(3/4)^4}{\pi}-\dfrac{p(n)}{q(n)}\sim-\dfrac{2\G(3/4)^8/\pi^{7/2}}{n^{1/2}}$$
$$A=1-(3/8)/n+(109/640)/n^2-(57/1024)/n^3+\cdots$$
Series:
$$\dfrac{\pi}{\G(3/4)^4}=\sum_{n\ge0}\dfrac{(1/2)_n^3}{n!^3}$$
Parametric family for $k\ge0$:
\begin{verbatim}
[()->gamma(3/4)^4/Pi,(4*n-1)*(4*n^2-2*n+1+2*k*(4*k+1)),
                     -8*n^3*(2*n+1)^3]
\end{verbatim}
Convergence type $P^+$ with $P=4k+1/2$.
\end{cf}

\medskip

\begin{cf}\label{5.11.52}{\ }
\begin{verbatim}
[()->gamma(3/4)^4/Pi,[2/3,21,16*n^2+6],[2/3,-(2*n+1)^2*(4*n+1)*(4*n+3)]]
\end{verbatim}
$$\dfrac{\G(3/4)^4}{\pi}=2/3+\dfrac{2/3}{21-\dfrac{315}{70-\dfrac{2475}{150-\dfrac{9555}{262-\dfrac{26163}{406-\dfrac{58443}{582-\ddots}}}}}}$$
Convergence type $P^+$ with $P=3/2$ and $C=\G(3/4)^2/(2^{3/2}3\pi)$, so that
$$\dfrac{\G(3/4)^4}{\pi}-\dfrac{p(n)}{q(n)}\sim\dfrac{\G(3/4)^2/(2^{3/2}3\pi)}{n^{3/2}}$$
$$A=1-(3/2)/n+(785/448)/n^2-(225/128)/n^3+\cdots$$
Series:
$$\dfrac{\G(3/4)^4}{\pi}=\dfrac{2}{3}+\dfrac{2}{21}\sum_{n\ge0}\dfrac{(5/4)_n}{(2n+3)(11/4)_n}$$
Parametric family for $k\ge0$:
\begin{verbatim}
[()->gamma(3/4)^4/Pi,16*n^2+6+4*k*(2*k+3),-(2*n+1)^2*(4*n+1)*(4*n+3)]
\end{verbatim}
Convergence type $P^+$ with $P=2k+3/2$.
\end{cf}

\medskip

\begin{cf}\label{5.11.53}{\ }
\begin{verbatim}
[()->gamma(3/4)^4/Pi,[2/3,8*n-2],[1/3,n^2*(2*n+1)^2]]
\end{verbatim}
$$\dfrac{\G(3/4)^4}{\pi}=2/3+\dfrac{1/3}{6+\dfrac{9}{14+\dfrac{100}{22+\dfrac{441}{30+\dfrac{1296}{38+\dfrac{3025}{46+\ddots}}}}}}$$
Convergence type $P^-$ with $P=4$ and $C=3\G(3/4)^4/(32\pi)$, so that
$$\dfrac{\G(3/4)^4}{\pi}-\dfrac{p(n)}{q(n)}\sim(-1)^n\dfrac{3\G(3/4)^4/(32\pi)}{n^4}$$
$$A=1-3/n+(13/4)/n^2+(9/4)/n^3-(217/32)/n^4+\cdots$$
Series:
\begin{align*}
  \dfrac{\G(3/4)^4}{\pi}&=\dfrac{2}{3}+\dfrac{2}{3}\sum_{n\ge0}\dfrac{(8n+7)n!^2(5/4)_n^2}{(16n^2+12n+3)(16n^2+44n+31)(3/2)_n^2(7/4)_n^2}\\
  \dfrac{\pi}{\G(3/4)^4}&=6\sum_{n\ge0}\dfrac{(8n+3)(1/2)_n^2(3/4)_n^2}{(16n^2-4n+1)(16n^2+28n+13)n!^2(5/4)_n^2}
\end{align*}
Parametric family for $k\ge0$:
\begin{verbatim}
[()->gamma(3/4)^4/Pi,(k+1)*(8*n-2),n^2*(2*n+1)^2]
\end{verbatim}
Convergence type $P^-$ with $P=4(k+1)$.
\end{cf}

\smallskip

\begin{cf}\label{5.11.54}{\ }
\begin{verbatim}
[()->gamma(3/4)^4/Pi,[2/3,20*n^2+n],[1,-n*(4*n-1)*(4*n+1)*(4*n+3)]]
\end{verbatim}
$$\dfrac{\G(3/4)^4}{\pi}=2/3+\dfrac{1}{21-\dfrac{105}{82-\dfrac{1386}{183-\dfrac{6435}{324-\dfrac{19380}{505-\dfrac{45885}{726-\ddots}}}}}}$$
Convergence type $E$ with $E=4$, $P=9/4$, and $C=\G(3/4)^2/(4\G(1/4))$, so that
$$\dfrac{\G(3/4)^4}{\pi}-\dfrac{p(n)}{q(n)}\sim\dfrac{\G(3/4)^2/\G(1/4)}{2^{2n+2}n^{9/4}}$$
$$A=1-(117/32)/n+(23433/2048)/n^2+\cdots$$
Series:
$$\dfrac{\G(3/4)^4}{\pi}=\dfrac{2}{3}+\dfrac{1}{21}\sum_{n\ge0}\dfrac{n!(5/4)_n}{(7/4)_n(11/4)_n}2^{-2n}$$
Parametric family for $k\ge0$:
\begin{verbatim}
[()->gamma(3/4)^4/Pi,20*n^2+(8*k+1)*n-k,
                     -n*(4*n+3)*(4*n-1-2*k)*(4*n+1-2*k)]
\end{verbatim}
Convergence type $E$ with $E=4$ and $P=3k+9/4$.
\end{cf}

\medskip

{\bf $j=2$: Continued Fractions for $\pi^{-2}\G((2m+1)/4)^4$}

\medskip

Since
$$\left(\dfrac{\G(1/4)}{\G(3/4)}\right)^2=\dfrac{1}{2}\dfrac{\G(1/4)^4}{\pi^2}\text{\quad and\quad}\left(\dfrac{\G(3/4)}{\G(1/4)}\right)^2=\dfrac{1}{2}\dfrac{\G(3/4)^4}{\pi^2}\;,$$
we will give CFs for the more symmetrical $(\G(1/4)/\G(3/4))^2$, so that CFs
for $\G(1/4)^4/\pi^2$ are trivially obtained by multiplying by $2$, and
CFs for $\G(3/4)^4/\pi^2$ by inverting and multiplying by $2$.

\smallskip

In addition, note that $(\G(1/4)/\G(3/4))^2=\CS(-4)$, see Section \ref{sec:CS}
for CFs related to general Chowla--Selberg gamma quotients.

\smallskip

\begin{cf}\label{5.11.55}{\ }
\begin{verbatim}
[()->(gamma(1/4)/gamma(3/4))^2,[4,1,2],[8,4*n^2-1]]
\end{verbatim}
$$\left(\dfrac{\G(1/4)}{\G(3/4)}\right)^2=4+\dfrac{8}{1+\dfrac{3}{2+\dfrac{15}{2+\dfrac{35}{2+\dfrac{63}{2+\dfrac{99}{2+\ddots}}}}}}$$
Convergence type $P^-$ with $P=1$ and $C=(\G(1/4)/\G(3/4))^2/2$, so that
$$\left(\dfrac{\G(1/4)}{\G(3/4)}\right)^2-\dfrac{p(n)}{q(n)}\sim(-1)^n\dfrac{(\G(1/4)/\G(3/4))^2/2}{n}\;.$$
$$A=1-(1/2)/n+(1/16)/n^2+(5/32)/n^3-(1/128)/n^4-(83/256)/n^5+\cdots$$
Series:
\begin{align*}\left(\dfrac{\G(1/4)}{\G(3/4)}\right)^2&=4+16\sum_{n\ge0}\dfrac{(3/4)_n^2}{(4n+5)(5/4)_n^2}\\
\left(\dfrac{\G(3/4)}{\G(1/4)}\right)^2&=-\dfrac{1}{4}+\sum_{n\ge0}\dfrac{(1/4)_n^2}{(4n+3)(3/4)_n^2}\end{align*}
Parametric family with $k\ge0$:
\begin{verbatim}
[()->(gamma(1/4)/gamma(3/4))^2,4*k+2,4*n^2-1]
\end{verbatim}
Convergence type $P^-$ with $P=2k+1$.
\end{cf}

\smallskip

\begin{cf}\label{5.11.56}{\ }
\begin{verbatim}
[()->(gamma(1/4)/gamma(3/4))^2,[8,4],[4,(2*n+1)^2]]
\end{verbatim}
$$\left(\dfrac{\G(1/4)}{\G(3/4)}\right)^2=8+\dfrac{4}{4+\dfrac{9}{4+\dfrac{25}{4+\dfrac{49}{4+\dfrac{81}{4+\dfrac{121}{4+\ddots}}}}}}$$
Convergence type $P^-$ with $P=2$ and $C=(\G(1/4)/\G(3/4))^2/8$, so that
$$\left(\dfrac{\G(1/4)}{\G(3/4)}\right)^2-\dfrac{p(n)}{q(n)}\sim(-1)^n\dfrac{(\G(1/4)/\G(3/4))^2/8}{n^2}\;.$$
$$A=1-2/n+(19/8)/n^2-(3/2)/n^3+(5/256)/n^4-(143/128)/n^5+\cdots$$
Series:
\begin{align*}\left(\dfrac{\G(1/4)}{\G(3/4)}\right)^2&=8+\dfrac{16}{25}\sum_{n\ge0}\dfrac{(3/4)_n^2}{(9/4)_n^2}\\
\left(\dfrac{\G(3/4)}{\G(1/4)}\right)^2&=\dfrac{1}{9}\sum_{n\ge0}\dfrac{(1/4)_n^2}{(7/4)_n^2}\end{align*}
Parametric family with $k\ge0$:
\begin{verbatim}
[()->(gamma(1/4)/gamma(3/4))^2,4*(k+1),(2*n+1)^2]
\end{verbatim}
Convergence type $P^-$ with $P=2k+2$.
\end{cf}

\smallskip

\begin{cf}\label{5.11.57}
\begin{verbatim}
[()->(gamma(1/4)/gamma(3/4))^2,[8,(4*n-3)*(8*n^2-12*n+15)],
                       [8,-4*n^2*(2*n-1)^2*(4*n-1)^2]]
\end{verbatim}
$$\left(\dfrac{\G(1/4)}{\G(3/4)}\right)^2=8+\dfrac{8}{11-\dfrac{36}{115-\dfrac{7056}{459-\dfrac{108900}{1235-\dfrac{705600}{2635-\dfrac{2924100}{4851-\ddots}}}}}}$$
Convergence type $P^+$ with $P=3$ and $C=\pi^3/(96\G(3/4)^8)$, so that
$$\left(\dfrac{\G(1/4)}{\G(3/4)}\right)^2-\dfrac{p(n)}{q(n)}\sim\dfrac{\pi^3/(96\G(3/4)^8)}{n^3}$$
$$A=1-(3/4)/n-(1/32)/n^2+\cdots$$
Series:
$$\left(\dfrac{\G(3/4)}{\G(1/4)}\right)^2=-\dfrac{1}{8}\sum_{n\ge0}\dfrac{(-1/2)_n^2}{(4n-1)n!^2}$$
Parametric family for $k\ge0$:
\begin{verbatim}
[()->(gamma(1/4)/gamma(3/4))^2,(4*n-3)*(8*n^2-12*n+15+8*k*(2*k+3)),
                       -4*n^2*(2*n-1)^2*(4*n-1)^2]
\end{verbatim}
Convergence type $P^+$ with $P=4k+3$.
\end{cf}

\smallskip

\begin{cf}\label{5.11.58}
\begin{verbatim}
[()->(gamma(1/4)/gamma(3/4))^2,[8,20,(4*n-1)*(8*n^2-4*n+3)],
                       [8,-4*n^2*(2*n+1)^2*(4*n+1)^2]]
\end{verbatim}
$$\left(\dfrac{\G(1/4)}{\G(3/4)}\right)^2=8+\dfrac{8}{20-\dfrac{900}{189-\dfrac{32400}{693-\dfrac{298116}{1725-\dfrac{1498176}{3477-\dfrac{5336100}{6141-\ddots}}}}}}$$
Convergence type $P^+$ with $P=1$ and $C=2/\pi$, so that
$$\left(\dfrac{\G(1/4)}{\G(3/4)}\right)^2-\dfrac{p(n)}{q(n)}\sim\dfrac{2/\pi}{n}$$
$$A=1-(3/4)/n+(15/32)/n^2-(27/128)/n^3+\cdots$$
Series:
$$\left(\dfrac{\G(1/4)}{\G(3/4)}\right)^2=8+2\sum_{n\ge0}\dfrac{(3/2)_n^2}{(4n+5)(n+1)!^2}$$
Parametric family for $k\ge0$:
\begin{verbatim}
[()->(gamma(1/4)/gamma(3/4))^2,(4*n-1)*(8*n^2-4*n+3+8*k*(2*k+1)),
                       -4*n^2*(2*n+1)^2*(4*n+1)^2]
\end{verbatim}
Convergence type $P^+$ with $P=4k+1$.
\end{cf}

\smallskip

\begin{cf}\label{5.11.59}{\ }
\begin{verbatim}
[()->(gamma(1/4)/gamma(3/4))^2,[0,1,6*(n-1)],[4,-2*(2*n-1)^2]]
\end{verbatim}
$$\left(\dfrac{\G(1/4)}{\G(3/4)}\right)^2=\dfrac{4}{1-\dfrac{2}{6-\dfrac{18}{12-\dfrac{50}{18-\dfrac{98}{24-\dfrac{162}{30-\ddots}}}}}}$$
Convergence type $E$ with $E=2$, $P=0$, and $C=2(\G(1/4)/\G(3/4))^2$, so that
$$\left(\dfrac{\G(1/4)}{\G(3/4)}\right)^2-\dfrac{p(n)}{q(n)}\sim\dfrac{(\G(1/4)/\G(3/4))^2}{2^{n-1}}$$
$$A=1-(3/4)/n+(9/32)/n^2-(405/128)/n^3+\cdots$$
Parametric family for $k\ge0$:
\begin{verbatim}
[()->(gamma(1/4)/gamma(3/4))^2,6*n-6+2*k,-2*(2*n-1)^2]
\end{verbatim}
Convergence type $E$ with $E=2$ and $P=2k$.
\end{cf}

\smallskip

\begin{cf}\label{5.11.60}{\ }
\begin{verbatim}
[()->(gamma(1/4)/gamma(3/4))^2,[0,1,14*(n-1)],[12,2*(4*n-1)*(4*n-3)]]
\end{verbatim}
$$\left(\dfrac{\G(1/4)}{\G(3/4)}\right)^2=\dfrac{12}{1+\dfrac{6}{14+\dfrac{70}{28+\dfrac{198}{42+\dfrac{390}{56+\dfrac{646}{70+\ddots}}}}}}$$
Convergence type $E$ with $E=-8$, $P=0$, and $C=2^{3/2}(\G(1/4)/\G(3/4))^2$,
so that
$$\left(\dfrac{\G(1/4)}{\G(3/4)}\right)^2-\dfrac{p(n)}{q(n)}\sim(-1)^n\dfrac{(\G(1/4)/\G(3/4))^2}{2^{3n-3/2}}$$
$$A=1-(7/48)/n+(49/4608)/n^2+\cdots$$
Parametric family for $k\ge0$:
\begin{verbatim}
[()->(gamma(1/4)/gamma(3/4))^2,14*n-14+18*k,2*(4*n-1)*(4*n-3)]
\end{verbatim}
Convergence type $E$ with $E=-8$ and $P=2k$.
\end{cf}

\smallskip

\begin{cf}\label{5.11.61}{\ }
\begin{verbatim}
[()->(gamma(1/4)/gamma(3/4))^2,[0,3,40*n-42],[24,-9*(4*n-3)^2]]
\end{verbatim}
$$\left(\dfrac{\G(1/4)}{\G(3/4)}\right)^2=\dfrac{24}{3-\dfrac{9}{38-\dfrac{225}{78-\dfrac{729}{118-\dfrac{1521}{158-\dfrac{2601}{198-\ddots}}}}}}$$
Convergence type $E$ with $E=9$, $P=1/2$ and $C=3\G(1/4)/\G(3/4)$, so that
$$\left(\dfrac{\G(1/4)}{\G(3/4)}\right)^2-\dfrac{p(n)}{q(n)}\sim\dfrac{\G(1/4)/\G(3/4)}{3^{2n-1}n^{1/2}}$$
$$A=1-(17/64)/n+(757/8192)/n^2-(69095/524288)/n^3+\cdots$$
Parametric family for $k\ge0$:
\begin{verbatim}
[()->(gamma(1/4)/gamma(3/4))^2,40*n-42+32*k,-9*(4*n-3)^2]
\end{verbatim}
Convergence type $E$ with $E=9$ and $P=2k+1/2$.
\end{cf}

\smallskip

\begin{cf}\label{5.11.62}{\ }
\begin{verbatim}
[()->(gamma(1/4)/gamma(3/4))^2,[0,5,40*n-38],[24,-9*(4*n-1)^2]]
\end{verbatim}
$$\left(\dfrac{\G(1/4)}{\G(3/4)}\right)^2=\dfrac{24}{5-\dfrac{81}{42-\dfrac{441}{82-\dfrac{1089}{122-\dfrac{2025}{162-\dfrac{3249}{202-\ddots}}}}}}$$
Convergence type $E$ with $E=9$, $P=-1/2$ and
$C=(8/3)(\G(1/4)/G(3/4))^3$, so that
$$\left(\dfrac{\G(1/4)}{\G(3/4)}\right)^2-\dfrac{p(n)}{q(n)}\sim\dfrac{8(\G(1/4)/\G(3/4))^3}{3^{2n+1}n^{-1/2}}$$
$$A=1-(17/64)/n-(179/8192)/n^2-(53183/524288)/n^3+\cdots$$
Parametric family for $k\ge0$:
\begin{verbatim}
[()->(gamma(1/4)/gamma(3/4))^2,40*n-38+32*k,-9*(4*n-1)^2]
\end{verbatim}
Convergence type $E$ with $E=9$ and $P=2k-1/2$.
\end{cf}

\smallskip

\begin{cf}\label{5.11.63}{\ }
\begin{verbatim}
[()->(gamma(1/4)/gamma(3/4))^2,[8,11,12*n],[8,-(2*n+1)^2]]
\end{verbatim}
$$\left(\dfrac{\G(1/4)}{\G(3/4)}\right)^2=8+\dfrac{8}{11-\dfrac{9}{24-\dfrac{25}{36-\dfrac{49}{48-\dfrac{81}{60-\dfrac{121}{72-\ddots}}}}}}$$
Convergence type $E$ with $E=(1+\sqrt{2})^4$, $P=0$, and $C=4(\G(1/4)/\G(3/4))^2/(1+\sqrt{2})^4$, so that
$$\left(\dfrac{\G(1/4)}{\G(3/4)}\right)^2-\dfrac{p(n)}{q(n)}\sim\dfrac{4(\G(1/4)/\G(3/4))^2}{(1+\sqrt{2})^{4n+4}}\;.$$
$$A=1-(3d/16)/n+(3d/16+9/256)/n^2-(795d/4096+9/128)/n^3+\cdots$$
\end{cf}

\smallskip

\begin{cf}\label{5.11.64}{\ }
\begin{verbatim}
[()->(gamma(1/4)/gamma(3/4))^2,[8,255,(2*n-1)*(272*n^2-272*n-1)],
                               [192,-450,-(n^2-1)*(16*n^2-1)^2]]
\end{verbatim}
$$\left(\dfrac{\G(1/4)}{\G(3/4)}\right)^2=8+\dfrac{192}{255-\dfrac{450}{1629-\dfrac{11907}{8155-\dfrac{163592}{22841-\dfrac{975375}{48951-\dfrac{3820824}{89749-\ddots}}}}}}$$
Convergence type $E$ with $E=(1+\sqrt{2})^8$, $P=0$, and $C=4(\G(1/4)/\G(3/4))^2/(1+\sqrt{2})^4$, so that
$$\left(\dfrac{\G(1/4)}{\G(3/4)}\right)^2-\dfrac{p(n)}{q(n)}\sim\dfrac{4(\G(1/4)/\G(3/4))^2}{(1+\sqrt{2})^{8n+4}}\;.$$
$$A=1-(3d/32)/n+(3d/64+9/1024)/n^2-(795d/32768+9/1024)/n^3+\cdots$$
\end{cf}

This is simply the contraction of the preceding CF.

\smallskip

\begin{cf}\label{5.11.65}{\ }
\begin{verbatim}
[()->(gamma(1/4)/gamma(3/4))^2,[4,7,12*n^2-12*n-1],
                               [32,-6,-(n^2-1)*(4*n^2-1)]]
\end{verbatim}
$$\left(\dfrac{\G(1/4)}{\G(3/4)}\right)^2=4+\dfrac{32}{7-\dfrac{6}{23-\dfrac{45}{71-\dfrac{280}{143-\dfrac{945}{239-\dfrac{2376}{359-\ddots}}}}}}$$
Convergence type $E$ with $E=(1+\sqrt{2})^4$, $P=0$, and $C=4(\G(1/4)/\G(3/4))^2/(1+\sqrt{2})^2$, so that
$$\left(\dfrac{\G(1/4)}{\G(3/4)}\right)^2-\dfrac{p(n)}{q(n)}\sim\dfrac{4(\G(1/4)/\G(3/4))^2}{(1+\sqrt{2})^{4n+2}}\;.$$
$$A=1+(d/16)/n+(-d/32+1/256)/n^2-(87d/4096+1/256)/n^3+\cdots$$
\end{cf}

\smallskip

\begin{cf}\label{5.11.66}{\ }
\begin{verbatim}
[()->(gamma(1/4)/gamma(3/4))^2,[0,15,378*(n-1)],[132,6*(6*n-1)*(6*n-5)]]
\end{verbatim}
$$\left(\dfrac{\G(1/4)}{\G(3/4)}\right)^2=\dfrac{132}{15+\dfrac{30}{378+\dfrac{462}{756+\dfrac{1326}{1134+\dfrac{2622}{1512+\dfrac{4350}{1890+\ddots}}}}}}$$
Convergence type $E$ with $E=-(1327+231\sqrt{33})/4$, $P=0$, and
$C=4(\G(1/4)/\G(3/4))^2$, so that
$$\left(\dfrac{\G(1/4)}{\G(3/4)}\right)^2-\dfrac{p(n)}{q(n)}\sim(-1)^n\dfrac{4(\G(1/4)/\G(3/4))^2}{((1327+231\sqrt{33})/4)^n}$$
$$A=1-(35d/1452)/n+(1225/127776)/n^2+\cdots$$
\end{cf}

\smallskip

\begin{cf}\label{5.11.67}{\ }
\begin{verbatim}
[()->12^(1/4)*(gamma(1/4)/gamma(3/4))^2,[0,3,56*(n-1)],[48,-(4*n-1)*(4*n-3)]]
\end{verbatim}
$$12^{1/4}\left(\dfrac{\G(1/4)}{\G(3/4)}\right)^2=\dfrac{48}{3-\dfrac{3}{56-\dfrac{35}{112-\dfrac{99}{168-\dfrac{195}{224-\dfrac{323}{280-\ddots}}}}}}$$
Convergence type $E$ with $E=(2+\sqrt{3})^4$, $P=0$, and
$C=2\cdot3^{5/4}(\G(1/4)/\G(3/4))^2$, so that
$$12^{1/4}\left(\dfrac{\G(1/4)}{\G(3/4)}\right)^2-\dfrac{p(n)}{q(n)}\sim\dfrac{2\cdot3^{5/4}(\G(1/4)/\G(3/4))^2}{(2+\sqrt{3})^{4n}}$$
$$A=1-(7d/64)/n+(147/8192)/n^2+\cdots$$
\end{cf}

\smallskip

\begin{cf}\label{5.11.68}{\ }
\begin{verbatim}
[()->5^(-1/4)*(gamma(1/4)/gamma(3/4))^2,
[0,41,1288*(n-1)],[240,-(4*n-1)*(4*n-3)]]
\end{verbatim}
$$5^{-1/4}\left(\dfrac{\G(1/4)}{\G(3/4)}\right)^2=\dfrac{240}{41-\dfrac{3}{1288-\dfrac{35}{2576-\dfrac{99}{3864-\dfrac{195}{5152-\dfrac{323}{6440-\ddots}}}}}}$$
Convergence type $E$ with $E=((1+\sqrt{5})/2)^{24}$, $P=0$, and
$C=2^{1/2}5^{3/4}(\G(1/4)/\G(3/4))^2$, so that
$$5^{-1/4}\left(\dfrac{\G(1/4)}{\G(3/4)}\right)^2-\dfrac{p(n)}{q(n)}\sim\dfrac{2^{1/2}5^{3/4}(\G(1/4)/\G(3/4))^2}{((1+\sqrt{5})/2)^{24n}}$$
$$A=1-(161d/1920)/n+(25921/1474560)/n^2+\cdots$$
\end{cf}

\smallskip

{\bf $j\ge3$:}

\medskip

As mentioned above, CFs for $\pi^{-j}\G((2m+1)/4)^4$ can be trivially
deduced from those for $\pi^{4-j}\G((3-2m)/4)^4$, so we do not give them
but simply refer to the corresponding CFs.

\begin{itemize}\item We have
  $$\pi^{-3}\G(1/4)^4=\dfrac{4}{\pi^{-1}\G(3/4)^4}\text{\quad and\quad}
  \pi^{-3}\G(3/4)^4=\dfrac{4}{\pi^{-1}\G(1/4)^4}\;,$$
  so invert the CFs beginning with \ref{5.11.51} and \ref{5.11.47}
  respectively.
\item We have
  $$\pi^{-4}\G(1/4)^4=\dfrac{4}{\G(3/4)^4}\text{\quad and\quad}
  \pi^{-4}\G(3/4)^4=\dfrac{4}{\G(1/4)^4}\;,$$
  so invert \ref{5.11.46} and \ref{5.11.45} respectively.
\end{itemize}

\medskip

\subsection{CFs for Higher Powers of $\G((2m+1)/4)$}

\medskip

The only CFs of this type that I know are for $(\G(1/4)/\G(3/4))^4$
and one for $(\G(1/4)/\G(3/4))^8$.

\begin{cf}\label{5.11.69}{\ }
\begin{verbatim}
[()->(gamma(1/4)/gamma(3/4))^4,[80,8*n^2+5],[-32,-(2*n+1)^4]]
\end{verbatim}
$$\left(\dfrac{\G(1/4)}{\G(3/4)}\right)^4=80-\dfrac{32}{13-\dfrac{81}{37-\dfrac{625}{77-\dfrac{2401}{133-\dfrac{6561}{205-\dfrac{14641}{293-\ddots}}}}}}$$
Convergence type $P^+$ with $P=2$ and $C=-\pi(\G(1/4)/\G(3/4))^4/64$, so that
$$\left(\dfrac{\G(1/4)}{\G(3/4)}\right)^4-\dfrac{p(n)}{q(n)}\sim-\dfrac{\pi(\G(1/4)/\G(3/4))^4/64}{n^2}\;.$$
$$A=1-2/n+\cdots$$
Parametric family for $k\ge0$ and $u\ge0$:
\begin{verbatim}
[()->(gamma(1/4)/gamma(3/4))^4,8*n^2+1-4*u^2+4*k^2,
                               -((2*n+1)^2-4*u^2)*(2*n+1)^2]
\end{verbatim}
Convergence type $P^+$ with $P=2k$.
\end{cf}

\smallskip
        
\begin{cf}\label{5.11.70}{\ }
\begin{verbatim}
[()->(gamma(1/4)/gamma(3/4))^4,[48,4,8*n^2-8*n+5],[96,-(4*n^2-1)^2]]
\end{verbatim}
$$\left(\dfrac{\G(1/4)}{\G(3/4)}\right)^4=48+\dfrac{96}{4-\dfrac{9}{21-\dfrac{225}{53-\dfrac{1225}{101-\dfrac{3969}{165-\dfrac{9801}{245-\ddots}}}}}}$$
Convergence type $P^+$ with $P=2$ and $C=3\pi(\G(1/4)/\G(3/4))^4/64$, so that
$$\left(\dfrac{\G(1/4)}{\G(3/4)}\right)^4-\dfrac{p(n)}{q(n)}\sim\dfrac{3\pi(\G(1/4)/\G(3/4))^4/64}{n^2}\;.$$
$$A=1-1/n+\cdots$$
Parametric family for $k\ge0$:
\begin{verbatim}
[()->(gamma(1/4)/gamma(3/4))^4,8*n^2-8*n+1+4*(k+1)^2,-(4*n^2-1)^2]
\end{verbatim}
Convergence type $P^+$ with $P=2k+2$.
\end{cf}

\smallskip
        
\begin{cf}\label{5.11.71}{\ }
\begin{verbatim}
[()->(gamma(1/4)/gamma(3/4))^4,[48,5,8*n],[192,(2*n-1)*(2*n+1)^2*(2*n+3)]]
\end{verbatim}
$$\left(\dfrac{\G(1/4)}{\G(3/4)}\right)^4=48+\dfrac{192}{5+\dfrac{45}{16+\dfrac{525}{24+\dfrac{2205}{32+\dfrac{6237}{40+\dfrac{14157}{48+\ddots}}}}}}$$
Convergence type $P^-$ with $P=2$ and $C=(\G(1/4)/\G(3/4))^4/2$, so that
$$\left(\dfrac{\G(1/4)}{\G(3/4)}\right)^4-\dfrac{p(n)}{q(n)}\sim(-1)^n\dfrac{(\G(1/4)/\G(3/4))^4/2}{n^2}\;.$$
$$A=1-2/n+(11/4)/n^2-3/n^3+(25/8)/n^4-(19/4)/n^5+\cdots$$
Series:
\begin{align*}\left(\dfrac{\G(1/4)}{\G(3/4)}\right)^4&=48+1024\sum_{n\ge0}\dfrac{(n+1)(4n+3)(3/4)_n^4}{(4n+5)^3(5/4)_n^4}\\
\left(\dfrac{\G(3/4)}{\G(1/4)}\right)^4&=-\dfrac{1}{16}+\dfrac{2}{81}\sum_{n\ge0}\dfrac{(2n+1)(4n+1)(4n+3)(1/4)_n^4}{(7/4)_n^4}\end{align*}
Parametric family for $k\ge0$:
\begin{verbatim}
[()->(gamma(1/4)/gamma(3/4))^4,8*(2*k+1)*n,(2*n-1)*(2*n+1)^2*(2*n+3)]
\end{verbatim}
Convergence type $P^-$ with $P=4k+2$.
\end{cf}

\smallskip
        
\begin{cf}\label{5.11.72}{\ }
\begin{verbatim}
[()->(gamma(1/4)/gamma(3/4))^4,[80,17,16*n],[-64,(2*n+1)^4]]
\end{verbatim}
$$\left(\dfrac{\G(1/4)}{\G(3/4)}\right)^4=80-\dfrac{64}{17+\dfrac{81}{32+\dfrac{625}{48+\dfrac{2401}{64+\dfrac{6561}{80+\dfrac{14641}{96+\ddots}}}}}}$$
Convergence type $P^-$ with $P=4$ and $C=-(\G(1/4)/\G(3/4))^4/8$, so that
$$\left(\dfrac{\G(1/4)}{\G(3/4)}\right)^4-\dfrac{p(n)}{q(n)}\sim-\dfrac{(\G(1/4)/\G(3/4))^4/8}{n^4}\;.$$
$$A=1-4/n+(15/2)/n^2-5/n^3-(107/16)/n^4+\cdots$$
Series:
\begin{align*}\left(\dfrac{\G(1/4)}{\G(3/4)}\right)^4&=80-\dfrac{2048}{625}\sum_{n\ge0}\dfrac{(n+1)(3/4)_n^4}{(9/4)_n^4}\\
\left(\dfrac{\G(3/4)}{\G(1/4)}\right)^4&=\dfrac{1}{16}-\dfrac{4}{81}\sum_{n\ge0}\dfrac{(2n+1)(1/4)_n^4}{(7/4)_n^4}\end{align*}
Parametric family for $k\ge0$ and $u\ge0$:
\begin{verbatim}
[()->(gamma(1/4)/gamma(3/4))^4,16*(k+1)*n,((2*n+1)^2-4*u^2)^2]
\end{verbatim}
Convergence type $P^-$ with $P=4(k+1)$.
\end{cf}

\smallskip
        
\begin{cf}\label{5.11.73}{\ }
\begin{verbatim}
[()->(gamma(1/4)/gamma(3/4))^4,[48,15,20*n^2+1],
                               [288,-(2*n+1)^2*(4*n+1)*(4*n+3)]]
\end{verbatim}
$$\left(\dfrac{\G(1/4)}{\G(3/4)}\right)^4=48+\dfrac{288}{15-\dfrac{315}{81-\dfrac{2475}{181-\dfrac{9555}{321-\dfrac{26163}{501-\dfrac{58443}{721-\ddots}}}}}}$$
Convergence type $E$ with $E=4$, $P=0$, and $C=2^{-1/2}(\G(1/4)/\G(3/4))^4$,
so that
$$\left(\dfrac{\G(1/4)}{\G(3/4)}\right)^4-\dfrac{p(n)}{q(n)}\sim\dfrac{(\G(1/4)/\G(3/4))^4}{2^{2n+1/2}}$$
$$A=1-(9/16)/n+(369/512)/n^2-(11331/8192)/n^3+\cdots$$
Parametric family for $k\ge0$:
\begin{verbatim}
[()->(gamma(1/4)/gamma(3/4))^4,20*n^2+8*k*n+1,
                               -(2*n+1)^2*(4*n-2*k+1)*(4*n-2*k+3)]
\end{verbatim}
Convergence type $E$ with $E=4$ and $P=3k$.
\end{cf}

\smallskip

\begin{cf}\label{5.11.74}{\ }
\begin{verbatim}
[()->(gamma(1/4)/gamma(3/4))^4,[96,125,96*n^2+2],
                               [-2592,3*(4*n+1)^2*(4*n+3)^2]]
\end{verbatim}
$$\left(\dfrac{\G(1/4)}{\G(3/4)}\right)^4=96-\dfrac{2592}{125+\dfrac{3675}{386+\dfrac{29403}{866+\dfrac{114075}{1538+\dfrac{312987}{2402+\dfrac{699867}{3458+\ddots}}}}}}$$
Convergence type $E$ with $E=-(2+\sqrt{3})^2$, $P=0$, and
$C=-4(\G(1/4)/\G(3/4))^4/(2+\sqrt{3})^2$, so that
$$\left(\dfrac{\G(1/4)}{\G(3/4)}\right)^4-\dfrac{p(n)}{q(n)}\sim(-1)^{n+1}\dfrac{4(\G(1/4)/\G(3/4))^4}{(2+\sqrt{3})^{2n+2}}$$
$$A=1-(d/6)/n+(d/6+1/24)/n^2-(119d/864+1/12)/n^3+\cdots$$
\end{cf}
    
\smallskip

\begin{cf}\label{5.11.75}{\ }
\begin{verbatim}
[()->(gamma(1/4)/gamma(3/4))^4,[80,47,44*n^2+1],[-160,(2*n+1)^4]]
\end{verbatim}
$$\left(\dfrac{\G(1/4)}{\G(3/4)}\right)^4=80-\dfrac{160}{47+\dfrac{81}{177+\dfrac{625}{397+\dfrac{2401}{705+\dfrac{6561}{1101+\dfrac{14641}{1585+\ddots}}}}}}$$
Convergence type $E$ with $E=-((1+\sqrt{5})/2)^{10}$, $P=0$, and
$C=-2\G(1/4)^8/(((1+\sqrt{5})/2)^{10}\pi^4)$, so that
$$\left(\dfrac{\G(1/4)}{\G(3/4)}\right)^4-\dfrac{p(n)}{q(n)}\sim(-1)^{n+1}\dfrac{2\G(1/4)^8/\pi^4}{((1+\sqrt{5})/2)^{10n+10}}\;.$$
$$A=1-(d/5)/n+(d/5+1/10)/n^2-(477d/2500+1/5)/n^3+\cdots$$
\end{cf}

Note that this CF converges quite fast.

\smallskip

\begin{cf}\label{5.11.76}{\ }
\begin{verbatim}
[()->3^(1/2)*(gamma(1/4)/gamma(3/4))^4,[144,100,90*n^2+2],
                                       [-1152,2*(3*n+1)^2*(3*n+2)^2]]
\end{verbatim}
$$3^{1/2}\left(\dfrac{\G(1/4)}{\G(3/4)}\right)^4=144-\dfrac{1152}{100+\dfrac{800}{362+\dfrac{6272}{812+\dfrac{24200}{1442+\dfrac{66248}{2252+\dfrac{147968}{3242+\ddots}}}}}}$$
Convergence type $E$ with $E=-(2+\sqrt{3})^3$, $P=0$, and
$C=-2\cdot3^{3/2}(\G(1/4)/\G(3/4))^4/(2+\sqrt{3})^3$, so that
$$3^{1/2}\left(\dfrac{\G(1/4)}{\G(3/4)}\right)^4-\dfrac{p(n)}{q(n)}\sim(-1)^{n+1}\dfrac{2\cdot3^{3/2}\CS(-4)^2}{(2+\sqrt{3})^{3n+3}}$$
$$A=1-(2d/9)/n+(2d/9+2/27)/n^2-(50d/243+4/27)/n^3+\cdots$$
\end{cf}

\smallskip

\begin{cf}\label{5.11.77}{\ }
\begin{verbatim}
[()->(gamma(1/4)/gamma(3/4))^8,[2816,635,688*n^4+200*n^2+3],
                   [720896,-192*(2*n+1)^6*(3*n+1)*(3*n+2)]]
\end{verbatim}
$$\left(\dfrac{\G(1/4)}{\G(3/4)}\right)^8=2816+\dfrac{720896}{635-\dfrac{2799360}{11811-\dfrac{168000000}{57531-\dfrac{2484746880}{179331-\ddots}}}}$$
Convergence type $E$ with $E=27/16$, $P=0$, and $C=16\cdot3^{-5/2}(\G(1/4)/\G(3/4))^8$, so that
$$\left(\dfrac{\G(1/4)}{\G(3/4)}\right)^8-\dfrac{p(n)}{q(n)}\sim\dfrac{3^{1/2}(\G(1/4)/\G(3/4))^8}{(27/16)^{n+1}}$$
$$A=1-(55/36)/n+(6985/2592)/n^2+\cdots$$
\end{cf}

\smallskip

\begin{cf}\label{5.11.78}{\ }
\begin{verbatim}
[()->(gamma(1/4)/gamma(3/4))^8,
[5376,1065,1040*n^4+120*n^2+1],[516096,-16*(2*n+1)^6*(4*n+1)*(4*n+3)]]
\end{verbatim}
$$\left(\dfrac{\G(1/4)}{\G(3/4)}\right)^8=5376+\dfrac{516096}{1065-\dfrac{408240}{17121-\dfrac{24750000}{85321-\dfrac{367064880}{268161-\dfrac{2746487088}{653001-\dfrac{13690623408}{1352161-\ddots}}}}}}$$
Convergence type $E$ with $E=64$, $P=0$, and $C=2^{-5/2}(\G(1/4)/\G(3/4))^8$,
so that
$$\left(\dfrac{\G(1/4)}{\G(3/4)}\right)^8-\dfrac{p(n)}{q(n)}\sim\dfrac{(\G(1/4)/\G(3/4))^8}{2^{6n+5/2}}$$
$$A=1-(35/48)/n+(4585/4608)/n^2+\cdots$$
\end{cf}

\medskip

\subsection{CFs Involving $\G(1/3)$ and $\G(2/3)$}

{\ }
\medskip

In the CFs that follow involving $2^{1/3}$ or $2^{2/3}$, note that
$$\dfrac{2\G(1/2)\G(1/3)}{\G(1/6)\G(2/3)}=2^{1/3}\text{\quad and\quad}
  \dfrac{\G(1/6)\G(2/3)}{\G(1/2)\G(1/3)}=2^{2/3}$$

\smallskip
  
\begin{cf}\label{5.11.79}{\ }
\begin{verbatim}
[()->gamma(1/3)^3,[0,-2,4*(27*n^4-108*n^3+159*n^2-102*n+23)],
                 [-27,-(2*n-3)*(2*n+1)*(3*n-2)^3*(3*n-1)^3]]
\end{verbatim}
$$\G(1/3)^3=-\dfrac{27}{-2+\dfrac{24}{92-\dfrac{40000}{1676-\dfrac{3687936}{8636-\dfrac{59895000}{27452-\ddots}}}}}$$
Convergence type $P^+$ with $P=2/3$ and $C=2\cdot3^{1/2}\pi/\Gamma(2/3)^2$,
so that
$$\G(1/3)^3-\dfrac{p(n)}{q(n)}\sim\dfrac{2\cdot3^{1/2}\pi/\Gamma(2/3)^2}{n^{2/3}}\;.$$
$$A=1-(7/162)/n^2+(566/45927)/n^4-(27917/3188646)/n^6+\cdots$$
Series:
$$\G(1/3)^3=27\sum_{n\ge0}\dfrac{(2n+1)(2/3)_n}{(3n+1)(3n+2)(4/3)_n}$$
\end{cf}

\smallskip

\begin{cf}\label{5.11.80}{\ }
\begin{verbatim}
[()->gamma(1/3)^3,[81/4,80,3*(2*n-1)*(9*n^2-9*n+26)],
                 [-81,-(9*n^2-4)*(9*n^2-1)^2]]
\end{verbatim}
$$\G(1/3)^3=81/4-\dfrac{81}{80-\dfrac{320}{396-\dfrac{39200}{1200-\dfrac{492800}{2814-\dfrac{2862860}{5562-\dfrac{11088896}{9768-\ddots}}}}}}$$
Convergence type $P^+$ with $P=14/3$ and $C=-\G(1/3)^2/(14\cdot3^{1/2}\pi)$,
so that
$$\G(1/3)^3-\dfrac{p(n)}{q(n)}\sim-\dfrac{\G(1/3)^2/(14\cdot3^{1/2}\pi)}{n^{14/3}}\;.$$
$$A=1-(7/3)/n+(217/81)/n^2-(385/243)/n^3+(28154/85293)/n^4+\cdots$$
Series:
$$\G(1/3)^3=\dfrac{81}{4}-\dfrac{81}{2}\sum_{n\ge0}\dfrac{(2/3)_n}{(3n+1)(3n+4)(3n+2)(3n+5)(7/3)_n}$$
Parametric family for $k\ge0$:
\begin{verbatim}
[()->gamma(1/3)^3,3*(2*n-1)*(9*n^2-9*n+26+3*k*(6*k+14)),
                  -(9*n^2-4)*(9*n^2-1)^2]
\end{verbatim}
Convergence type $P^+$ with $P=4k+14/3$.
\end{cf}

\smallskip

\begin{cf}\label{5.11.81}{\ }
\begin{verbatim}
[()->gamma(1/3)^3,[27/2,(6*n-5)*(9*n^2-15*n+25)],[108,-27*n^3*(3*n-2)^3]]
\end{verbatim}
$$\G(1/3)^3=\dfrac{27}{2}+\dfrac{108}{19-\dfrac{27}{217-\dfrac{13824}{793-\dfrac{250047}{2071-\dfrac{1728000}{4375-\dfrac{7414875}{8029-\ddots}}}}}}$$
Convergence type $P^+$ with $P=4$ and $C=4\G(1/3)^3/729$, so that
$$\G(1/3)^3-\dfrac{p(n)}{q(n)}\sim\dfrac{4\G(1/3)^3/729}{n^4}\;.$$
$$A=1-(2/3)/n-(25/81)/n^2+(40/81)/n^3+\cdots$$
Series:
\begin{align*}\G(1/3)^3&=\dfrac{27}{2}+432\sum_{n\ge0}\dfrac{n!^3}{(27n^2+9n+2)(27n^2+63n+38)(4/3)_n^3}\\
  \dfrac{1}{\G(1/3)^3}&=\dfrac{2}{27}\sum_{n\ge0}\dfrac{(-2/3)_n^3}{n!^3}\end{align*}
Parametric family for $k\ge0$:
\begin{verbatim}
[()->gamma(1/3)^3,(6*n-5)*(9*n^2-15*n+18*(k+1)^2+7),-27*n^3*(3*n-2)^3]
\end{verbatim}
Convergence type $P^+$ with $P=4k+4$.
\end{cf}
        
\smallskip

\begin{cf}\label{5.11.82}{\ }
\begin{verbatim}
[()->gamma(1/3)^3,[81/4,(6*n+1)*(9*n^2+3*n+4)],
                  [-81,-27*n*(n+1)^2*(3*n+1)^2*(3*n+4)]]
\end{verbatim}
$$\G(1/3)^3=81/4-\dfrac{81}{112-\dfrac{12096}{598-\dfrac{238140}{1786-\dfrac{1684800}{4000-\dfrac{7300800}{7564-\ddots}}}}}$$
Convergence type $P^+$ with $P=2$ and $C=-2\G(1/3)^3/27$, so that
$$\G(1/3)^3-\dfrac{p(n)}{q(n)}\sim-\dfrac{2\G(1/3)^3/27}{n^2}\;.$$
$$A=1-(7/3)/n+(109/27)/n^2+\cdots$$
Series:
\begin{align*}\G(1/3)^3&=\dfrac{81}{4}-\dfrac{81}{16}\sum_{n\ge0}\dfrac{(n+1)^2n!^3}{(3n+7)(7/3)_n^3}\\
\dfrac{1}{\G(1/3)^3}&=\dfrac{4}{81}\sum_{n\ge0}\dfrac{(3n+1)(n+1)(1/3)_n^3}{(n+1)!^3}\end{align*}
Parametric family for $k\ge0$:
\begin{verbatim}
[()->gamma(1/3)^3,(6*n+1)*(9*n^2+3*n+4+18*k*(k+1)),
                   -27*n*(n+1)^2*(3*n+1)^2*(3*n+4)]
\end{verbatim}
Convergence type $P^+$ with $P=4k+2$.
\end{cf}

Can both be Ap\'ery accelerated with type $E=-(1+\sqrt{2})^4$.

\smallskip

\begin{cf}\label{5.11.83}{\ }
\begin{verbatim}
[()->gamma(1/3)^3,[18,48,(6*n-1)*(9*n^2-3*n+4)],
                  [36,-3*n*(3*n+1)^4*(3*n+2)]]
\end{verbatim}
$$\G(1/3)^3=18+\dfrac{36}{48-\dfrac{3840}{374-\dfrac{115248}{1292-\dfrac{990000}{3128-\dfrac{4798248}{6206-\dfrac{16711680}{10850-\ddots}}}}}}$$
Convergence type $P^+$ with $P=4/3$ and $C=3/(2\G(2/3))$, so that
$$\G(1/3)^3-\dfrac{p(n)}{q(n)}\sim\dfrac{3/(2\G(2/3))}{n^{4/3}}$$
$$A=1-(10/9)/n+(383/405)/n^2-(1370/2187)/n^3+\cdots$$
Series:
$$\G(1/3)^3=18+12\sum_{n\ge0}\dfrac{(5/3)_n}{(3n+4)^2(n+1)!}$$
Parametric family for $k\ge0$:
\begin{verbatim}
[()->gamma(1/3)^3,(6*n-1)*(9*n^2-3*n+4+6*k*(3*k+2)),
                  -3*n*(3*n+1)^4*(3*n+2)]
\end{verbatim}
Convergence type $P^+$ with $P=4k+4/3$.
\end{cf}

\smallskip

\begin{cf}\label{5.11.84}{\ }
\begin{verbatim}
[()->gamma(1/3)^3,[77/4,(6*n+1)*(9*n^2+3*n+37)],[-8,-27*n^3*(3*n+4)^3]]
\end{verbatim}
$$\G(1/3)^3=77/4-\dfrac{8}{343-\dfrac{9261}{1027-\dfrac{216000}{2413-\dfrac{1601613}{4825-\dfrac{7077888}{8587-\dfrac{23149125}{14023-\ddots}}}}}}$$
Convergence type $P^+$ with $P=6$ and $C=-256\G(1/3)^3/3^{10}$, so that
$$\G(1/3)^3-\dfrac{p(n)}{q(n)}\sim-\dfrac{256\G(1/3)^3/3^{10}}{n^6}$$
$$A=1-7/n+(1001/36)/n^2-(245/3)/n^3+\cdots$$
Series:
$$\G(1/3)^3=\dfrac{77}{4}-\dfrac{8}{343}\sum_{n\ge0}\dfrac{n!^3}{(10/3)_n^3}$$
There also exists a hypergeometric series for $1/\G(1/3)^3$ with polynomials
of degree $4$, not given.

Parametric family for $k\ge0$:
\begin{verbatim}
[()->gamma(1/3)^3,(6*n+1)*(9*n^2+3*n+37+18*k*(k+3)),-27*n^3*(3*n+4)^3]
\end{verbatim}
Convergence type $P^+$ with $P=4k+6$.
\end{cf}
        
\smallskip

\begin{cf}\label{5.11.85}{\ }
\begin{verbatim}
[()->gamma(2/3)^3,[0,-16,4*(27*n^4-108*n^3+186*n^2-156*n+43)],
[-81/2,-(2*n-3)*(2*n+1)*(3*n-4)*(3*n+1)*(3*n-2)^2*(3*n-1)^2]]
\end{verbatim}
$$\G(2/3)^3=-\dfrac{81/2}{-16-\dfrac{48}{172-\dfrac{28000}{2080-\dfrac{3292800}{9580-\dfrac{56628000}{29152-\ddots}}}}}$$
Convergence type $P^+$ with $P=10/3$ and $C=-3^{1/2}\G(2/3)^2/(20\pi)$, so that
$$\G(2/3)^3-\dfrac{p(n)}{q(n)}\sim-\dfrac{3^{1/2}\G(2/3)^2/(20\pi)}{n^{10/3}}\;.$$
$$A=1-(10/81)/n^2+(11674/72171)/n^4-(381200/1594323)/n^6+\cdots$$
Series:
$$\G(2/3)^3=-\dfrac{81}{4}\sum_{n\ge0}\dfrac{(2n+1)(1/3)_n}{(3n-1)(3n+2)(3n+1)(3n+4)(5/3)_n}$$
\end{cf}

\smallskip

\begin{cf}\label{5.11.86}{\ }
\begin{verbatim}
[()->gamma(2/3)^3,[9/4,40,3*(2*n-1)*(9*n^2-9*n+14)],
                 [9,-(9*n^2-4)*(9*n^2-1)^2]]
\end{verbatim}
$$\G(2/3)^3=9/4+\dfrac{9}{40-\dfrac{320}{288-\dfrac{39200}{1020-\dfrac{492800}{2562-\dfrac{2862860}{5238-\dfrac{11088896}{9372-\ddots}}}}}}$$
Convergence type $P^+$ with $P=10/3$ and $C=\G(2/3)^2/(10\cdot3^{1/2}\pi)$,
so that
$$\G(2/3)^3-\dfrac{p(n)}{q(n)}\sim\dfrac{\G(2/3)^2/(10\cdot3^{1/2}\pi)}{n^{10/3}}\;.$$
$$A=1-(5/3)/n+(125/81)/n^2-(220/243)/n^3+(28414/72171)/n^4+\cdots$$
Series:
$$\G(2/3)^3=\dfrac{9}{4}+\dfrac{9}{5}\sum_{n\ge0}\dfrac{(1/3)_n}{(3n+2)(3n+4)(8/3)_n}$$
Parametric family for $k\ge0$:
\begin{verbatim}
[()->gamma(2/3)^3,3*(2*n-1)*(9*n^2-9*n+14+6*k*(3*k+5)),
                  -(9*n^2-4)*(9*n^2-1)^2]
\end{verbatim}
Convergence type $P^+$ with $P=4k+10/3$.
\end{cf}

\smallskip

\begin{cf}\label{5.11.87}{\ }
\begin{verbatim}
[()->gamma(2/3)^3,[9/4,(6*n-1)*(9*n^2-3*n+4)],
                  [9,-27*n^2*(n+1)*(3*n-1)*(3*n+2)^2]]
\end{verbatim}
$$\G(2/3)^3=9/4+\dfrac{9}{50-\dfrac{2700}{374-\dfrac{103680}{1292-\dfrac{940896}{3128-\dfrac{4656960}{6206-\dfrac{16386300}{10850-\ddots}}}}}}$$
Convergence type $P^+$ with $P=2$ and $C=2\G(2/3)^3/27$, so that
$$\G(2/3)^3-\dfrac{p(n)}{q(n)}\sim\dfrac{2\G(2/3)^3/27}{n^2}\;.$$
$$A=1-(5/3)/n+(53/27)/n^2+\cdots$$
Series:
\begin{align*}\G(2/3)^3&=\dfrac{9}{4}+\dfrac{9}{2}\sum_{n\ge0}\dfrac{(n+1)n!^3}{(3n+5)^2(5/3)_n^3}\\
\dfrac{1}{\G(2/3)^3}&=-\dfrac{4}{9}\sum_{n\ge0}\dfrac{(2/3)_n^3}{(3n-1)(n+1)n!^3}\end{align*}
Parametric family for $k\ge0$:
\begin{verbatim}
[()->gamma(2/3)^3,(6*n-1)*(9*n^2-3*n+4+18*k*(k+1)),
                  -27*n^2*(n+1)*(3*n-1)*(3*n+2)^2]
\end{verbatim}
Convergence type $P^+$ with $P=4k+2$.
\end{cf}

\smallskip

\begin{cf}\label{5.11.88}{\ }
\begin{verbatim}
[()->gamma(2/3)^3,[5/2,(6*n-1)*(9*n^2-3*n+19)],
                  [-2,-27*n^3*(3*n+2)^3]]
\end{verbatim}
$$\G(2/3)^3=5/2-\dfrac{2}{125-\dfrac{3375}{539-\dfrac{110592}{1547-\dfrac{970299}{3473-\dfrac{4741632}{6641-\dfrac{16581375}{11375-\ddots}}}}}}$$
Convergence type $P^+$ with $P=4$ and $C=-4\G(2/3)^3/729$, so that
$$\G(2/3)^3-\dfrac{p(n)}{q(n)}\sim-\dfrac{4\G(2/3)^3/729}{n^4}\;.$$
$$A=1-(10/3)/n+(515/81)/n^2-(700/81)/n^3+\cdots$$
Series:
\begin{align*}\G(2/3)^3&=\dfrac{5}{2}-\dfrac{2}{125}\sum_{n\ge0}\dfrac{n!^3}{(8/3)_n^3}\\
  \dfrac{1}{\G(2/3)^3}&=16\sum_{n\ge0}\dfrac{(2/3)_n^3}{(27n^2-9n+2)(27n^2+45n+20)n!^3}\end{align*}
Parametric family for $k\ge0$:
\begin{verbatim}
[()->gamma(2/3)^3,(6*n-1)*(9*n^2-3*n+19+18*k*(k+2)),-27*n^3*(3*n+2)^3]
\end{verbatim}
Convergence type $P^+$ with $P=4k+4$.
\end{cf}

\smallskip

\begin{cf}\label{5.11.89}{\ }
\begin{verbatim}
[()->gamma(2/3)^3,[-27/8,(6*n-7)*(9*n^2-21*n+49)],
                  [-216,-27*n^3*(3*n-4)^3]]
\end{verbatim}
$$\G(2/3)^3=-27/8-\dfrac{216}{-37+\dfrac{27}{215-\dfrac{1728}{737-\dfrac{91125}{1853-\dfrac{884736}{3887-\dfrac{4492125}{7163-\ddots}}}}}}$$
Convergence type $P^+$ with $P=6$ and $C=256\G(2/3)^3/3^{10}$, so that
$$\G(2/3)^3-\dfrac{p(n)}{q(n)}\sim\dfrac{256\G(2/3)^3/3^{10}}{n^6}\;.$$
$$A=1+1/n-(7/36)/n^2-(7/9)/n^3+(371/486)/n^4+\cdots$$
Series:
$$\dfrac{1}{\G(2/3)^3}=-\dfrac{8}{27}\sum_{n\ge0}\dfrac{(-4/3)_n^3}{n!^3}$$
There also exists a hypergeometric series for $\G(2/3)^3$ with polynomials
of degree $4$, not given.

Parametric family for $k\ge0$:
\begin{verbatim}
[()->gamma(2/3)^3,(6*n-7)*(9*n^2-21*n+18*(k+1)*(k+2)+13),
                  -27*n^3*(3*n-4)^3]
\end{verbatim}
Convergence type $P^+$ with $P=4k+6$.
\end{cf}
        
The preceding three CFs can all be Ap\'ery accelerated with type $E=-(1+\sqrt{2})^4$.

\smallskip

\begin{cf}\label{5.11.90}{\ }
\begin{verbatim}
[()->gamma(2/3)^3,[3,12,(6*n-5)*(9*n^2-15*n+16)],
                  [-6,-3*n*(3*n-1)^4*(3*n-2)]]
\end{verbatim}
$$\G(2/3)^3=3-\dfrac{6}{12-\dfrac{48}{154-\dfrac{15000}{676-\dfrac{258048}{1900-\dfrac{1756920}{4150-\dfrac{7491120}{7750-\ddots}}}}}}$$
Convergence type $P^+$ with $P=8/3$ and $C=-1/(12\G(1/3))$, so that
$$\G(2/3)^3-\dfrac{p(n)}{q(n)}\sim-\dfrac{1/(12\G(2/3))}{n^{8/3}}$$
$$A=1-(4/9)/n-(38/189)/n^2+(496/2187)/n^3+\cdots$$
Series:
$$\G(2/3)^3=3-2\sum_{n\ge0}\dfrac{(1/3)_n}{(3n+2)^2(n+1)!}$$
Parametric family for $k\ge0$:
\begin{verbatim}
[()->gamma(2/3)^3,(6*n-5)*(9*n^2-15*n+16+6*k*(3*k+4)),
                  -3*n*(3*n-1)^4*(3*n-2)]
\end{verbatim}
Convergence type $P^+$ with $P=4k+8/3$.
\end{cf}

\medskip

For the next CFs, note that
$$\dfrac{\G(1/3)^3}{2\pi/\sqrt{3}}=\dfrac{\G(1/3)^2}{\G(2/3)}=\dfrac{4\pi^2}{9\G(2/3)^3}\text{\quad and\quad}\dfrac{\G(2/3)^3}{2\pi/\sqrt{3}}=\dfrac{\G(2/3)^2}{\G(1/3)}=\dfrac{4\pi^2}{9\G(1/3)^3}\;.$$

\smallskip

\begin{cf}\label{5.11.91}{\ }
\begin{verbatim}
[()->gamma(1/3)^2/gamma(2/3),
[3,12,18*n^2-6*n+2],[6,-3*n*(3*n+1)^2*(3*n+2)]]
\end{verbatim}
$$\dfrac{\G(1/3)^2}{\G(2/3)}=3+\dfrac{6}{12-\dfrac{240}{62-\dfrac{2352}{146-\dfrac{9900}{266-\dfrac{28392}{422-\dfrac{65280}{614-\ddots}}}}}}$$
Convergence type $P^+$ with $P=1/3$ and $C=3/\G(2/3)$, so that
$$\dfrac{\G(1/3)^2}{\G(2/3)}-\dfrac{p(n)}{q(n)}\sim\dfrac{3/\G(2/3)}{n^{1/3}}\;.$$
$$A=1-(5/18)/n+(76/567)/n^2-(265/4374)/n^3+(5021/255879)/n^4-\cdots$$
Series:
$$\dfrac{\G(1/3)^2}{\G(2/3)}=3+2\sum_{n\ge0}\dfrac{(5/3)_n}{(3n+4)(n+1)!}$$
Parametric family for $k\ge0$:
\begin{verbatim}
[()->gamma(1/3)^2/gamma(2/3),
18*n^2-6*n+2+3*k*(3*k+1),-3*n*(3*n+1)^2*(3*n+2)]
\end{verbatim}
Convergence type $P^+$ with $P=2k+1/3$.
\end{cf}

\smallskip

\begin{cf}\label{5.11.92}{\ }
\begin{verbatim}
[()->gamma(1/3)^2/gamma(2/3),[9,18*n^2+3*n+1],[-36,-9*n*(n+1)*(3*n+2)^2]]
\end{verbatim}
$$\dfrac{\G(1/3)^2}{\G(2/3)}=9-\dfrac{36}{22-\dfrac{450}{79-\dfrac{3456}{172-\dfrac{13068}{301-\dfrac{35280}{466-\dfrac{78030}{667-\ddots}}}}}}$$
Convergence type $P^+$ with $P=2/3$ and $C=-8\pi^4/(27\G(2/3)^8)$,
so that
$$\dfrac{\G(1/3)^2}{\G(2/3)}-\dfrac{p(n)}{q(n)}\sim-\dfrac{8\pi^4/(27\G(2/3)^8)}{n^{2/3}}$$
$$A=1-(37/45)/n+(22/27)/n^2-(19946/24057)/n^3+\cdots$$
Series:
$$\dfrac{\G(2/3)}{\G(1/3)^2}=\dfrac{1}{9}\sum_{n\ge0}\dfrac{(2/3)_n^2}{(n+1)n!^2}$$
Parametric family for $k\ge0$:
\begin{verbatim}
[()->gamma(1/3)^2/gamma(2/3),18*n^2+3*n+(3*k+1)^2,-9*n*(n+1)*(3*n+2)^2]
\end{verbatim}
Convergence type $P^+$ with $P=2k+2/3$.
\end{cf}

\smallskip

\begin{cf}\label{5.11.93}{\ }
\begin{verbatim}
[()->gamma(1/3)^2/gamma(2/3),
[5,24,18*n^2+3*n+5],[4,-3*(n+1)*(3*n+1)^2*(3*n+2)]]
\end{verbatim}
$$\dfrac{\G(1/3)^2}{\G(2/3)}=5+\dfrac{4}{24-\dfrac{480}{83-\dfrac{3528}{176-\dfrac{13200}{305-\dfrac{35490}{470-\dfrac{78336}{671-\ddots}}}}}}$$
Convergence type $P^+$ with $P=4/3$ and $C=1/(2\G(2/3))$, so that
$$\dfrac{\G(1/3)^2}{\G(2/3)}-\dfrac{p(n)}{q(n)}\sim\dfrac{1/(2\G(2/3))}{n^{4/3}}$$
$$A=1-(94/63)/n+(151/81)/n^2-(60110/28431)/n^3+\cdots$$
Series:
$$\dfrac{\G(1/3)^2}{\G(2/3)}=5+\dfrac{4}{3}\sum_{n\ge0}\dfrac{(5/3)_n}{(3n+4)(n+2)!}$$
Parametric family for $k\ge0$:
\begin{verbatim}
[()->gamma(1/3)^2/gamma(2/3),
18*n^2+3*n+5+3*k*(3*k+4),-3*(n+1)*(3*n+1)^2*(3*n+2)]
\end{verbatim}
Convergence type $P^+$ with $P=2k+4/3$.
\end{cf}
      
\smallskip

\begin{cf}\label{5.11.94}{\ }
\begin{verbatim}
[()->gamma(1/3)^2/gamma(2/3),
[9,16,18*n^2-6*n+1],[-36,-3*(n+1)*(3*n-1)*(3*n+1)^2]]
\end{verbatim}
$$\dfrac{\G(1/3)^2}{\G(2/3)}=9-\dfrac{36}{16-\dfrac{192}{61-\dfrac{2205}{145-\dfrac{9600}{265-\dfrac{27885}{421-\dfrac{64512}{613-\ddots}}}}}}$$
Convergence type $P^+$ with $P=1$ and $C=-16\pi^2/(27\G(2/3)^3)$, so that
$$\dfrac{\G(1/3)^2}{\G(2/3)}-\dfrac{p(n)}{q(n)}\sim-\dfrac{16\pi^2/(27\G(2/3)^3)}{n}$$
$$A=1-(11/18)/n+(169/486)/n^2-(767/4374)/n^3+\cdots$$
Series:
\begin{align*}\dfrac{\G(1/3)^2}{\G(2/3)}&=9-\dfrac{9}{4}\sum_{n\ge0}\dfrac{(n+1)!(2/3)_n}{(7/3)_n^2}\\
\dfrac{\G(2/3)}{\G(1/3)^2}&=-\dfrac{1}{3}+\dfrac{4}{9}\sum_{n\ge0}\dfrac{(1/3)_n^2}{(n+1)!(2/3)_n}\end{align*}
Parametric family for $k\ge0$:
\begin{verbatim}
[()->gamma(1/3)^2/gamma(2/3),18*n^2-6*n+1+9*k*(k+1),
                                 -3*(n+1)*(3*n-1)*(3*n+1)^2]
\end{verbatim}
Convergence type $P^+$ with $P=2k+1$.
\end{cf}

\smallskip

\begin{cf}\label{5.11.95}{\ }
\begin{verbatim}
[()->gamma(1/3)^2/gamma(2/3),[9,19,15*(2*n-1)],
                                 [-72,(3*n-2)^2*(3*n+2)^2]]
\end{verbatim}
$$\dfrac{\G(1/3)^2}{\G(2/3)}=9-\dfrac{72}{19+\dfrac{25}{45+\dfrac{1024}{75+\dfrac{5929}{105+\dfrac{19600}{135+\dfrac{48841}{165+\ddots}}}}}}$$
Convergence type $P^-$ with $P=10/3$ and $C=-\G(1/3)^4/(12\pi^2)$, so that
$$\dfrac{\G(1/3)^2}{\G(2/3)}-\dfrac{p(n)}{q(n)}\sim(-1)^{n+1}\dfrac{\G(1/3)^4/(12\pi^2)}{n^{10/3}}$$
$$A=1-(5/3)/n+(65/81)/n^2+(260/243)/n^3+\cdots$$
Series:
$$\dfrac{\G(1/3)^2}{\G(2/3)}=9-72\sum_{n\ge0}(-1)^n\dfrac{(5/3)_n^2}{(9n^2+9n+1)(9n^2+27n+19)(4/3)_n^2}$$
Parametric family for $k\ge0$:
\begin{verbatim}
[()->gamma(1/3)^2/gamma(2/3),3*(6*k+5)*(2*n-1),(3*n-2)^2*(3*n+2)^2]
\end{verbatim}
Convergence type $P^-$ with $P=4k+10/3$.
\end{cf}

\smallskip

\begin{cf}\label{5.11.96}{\ }
\begin{verbatim}
[()->gamma(1/3)^2/gamma(2/3),[3,5*(6*n-5)],[12,9*n^2*(3*n-2)^2]]
\end{verbatim}
$$\dfrac{\G(1/3)^2}{\G(2/3)}=3+\dfrac{12}{5+\dfrac{9}{35+\dfrac{576}{65+\dfrac{3969}{95+\dfrac{14400}{125+\dfrac{38025}{155+\ddots}}}}}}$$
Convergence type $P^-$ with $P=10/3$ and $C=\G(1/3)^4/(18\pi^2)$, so that
$$\dfrac{\G(1/3)^2}{\G(2/3)}-\dfrac{p(n)}{q(n)}\sim(-1)^n\dfrac{\G(1/3)^4/(18\pi^2)}{n^{10/3}}\;.$$
$$A=1-(5/9)/n-(110/81)/n^2+(2900/2187)/n^3+(91202/19683)/n^4+\cdots$$
Series:
$$\dfrac{\G(2/3)}{\G(1/3)^2}=\dfrac{1}{3}\sum_{n\ge0}(-1)^n\dfrac{(-2/3)_n^2}{n!^2}$$
Parametric family for $k\ge0$:
\begin{verbatim}
[()->gamma(1/3)^2/gamma(2/3),(6*k+5)*(6*n-5),9*n^2*(3*n-2)^2]
\end{verbatim}
Convergence type $P^-$ with $P=4k+10/3$.
\end{cf}

\smallskip

\begin{cf}\label{5.11.97}{\ }
\begin{verbatim}
[()->gamma(1/3)^2/gamma(2/3),
[5,234*n^2-156*n-34],[14,9*n^2*(3*n+1)^2*(6*n-5)*(6*n+7)]]
\end{verbatim}
$$\dfrac{\G(1/3)^2}{\G(2/3)}=5+\dfrac{14}{44+\dfrac{1872}{590+\dfrac{234612}{1604+\dfrac{2632500}{3086+\dfrac{14333904}{5036+\dfrac{53280000}{7454+\ddots}}}}}}$$
Convergence type $P^-$ with $P=13/3$ and $C=(2/27)(\G(1/3)/\G(2/3))^2$, so that
$$\dfrac{\G(1/3)^2}{\G(2/3)}-\dfrac{p(n)}{q(n)}\sim(-1)^n\dfrac{(2/27)(\G(1/3)/\G(2/3))^2}{n^{13/3}}$$
$$A=1-(26/9)/n+(202/81)/n^2+(8588/2187)/n^3+\cdots$$
Series:
$$\dfrac{\G(2/3)}{\G(1/3)^2}=2\sum_{n\ge0}(-1)^n\dfrac{(6n+1)(1/3)_n^2}{(9n^2-6n+2)(9n^2+12n+5)n!^2}$$
\end{cf}
      
\smallskip

\begin{cf}\label{5.11.98}{\ }
\begin{verbatim}
[()->gamma(1/3)^2/gamma(2/3),[12,7*n],[-30,-2*n*(6*n+5)]]
\end{verbatim}
$$\dfrac{\G(1/3)^2}{\G(2/3)}=12-\dfrac{30}{7-\dfrac{22}{14-\dfrac{68}{21-\dfrac{138}{28-\dfrac{232}{35-\dfrac{350}{42-\ddots}}}}}}$$
Convergence type $E$ with $E=4/3$, $P=7/6$, and
$C=-(3\pi)^{3/2}/(2^{2/3}\G(2/3)^2)$, so that
$$\dfrac{\G(1/3)^2}{\G(2/3)}-\dfrac{p(n)}{q(n)}\sim-\dfrac{(3\pi)^{3/2}/(2^{2/3}\G(2/3)^2)}{(4/3)^nn^{7/6}}$$
$$A=1-(287/72)/n+(35833/1152)/n^2+\cdots$$
Parametric families for $k\ge0$ and $u\ge0$:
\begin{verbatim}
[()->gamma(1/3)^2/gamma(2/3),7*n+k,-2*n*(6*n+11-6*u)]
[()->gamma(1/3)^2/gamma(2/3),7*n+k,-2*(2*n-1)*(3*n+5-3*u)]
\end{verbatim}
Convergence type $E$ with $E=4/3$ and $P=2k+7u-35/6$ or $P=2k+7u-7/6$
respectively.
\end{cf}

\smallskip

\begin{cf}\label{5.11.99}{\ }
\begin{verbatim}
[()->gamma(1/3)^2/gamma(2/3),[3,2*n],[15/2,n*(3*n+5)]]
\end{verbatim}
$$\dfrac{\G(1/3)^2}{\G(2/3)}=3+\dfrac{15/2}{2+\dfrac{8}{4+\dfrac{22}{6+\dfrac{42}{8+\dfrac{68}{10+\dfrac{100}{12+\ddots}}}}}}$$
Convergence type $E$ with $E=-3$, $P=-1/3$, and $C=2^{10/3}\pi^3/(27\G(2/3)^5)$, so that
$$\dfrac{\G(1/3)^2}{\G(2/3)}-\dfrac{p(n)}{q(n)}\sim(-1)^n\dfrac{2^{10/3}\pi^3/\G(2/3)^5}{3^{n+3}n^{-1/3}}$$
$$A=1+(47/72)/n-(1223/3456)/n^2+\cdots$$
Parametric families for $k\ge0$ and $u\ge0$:
\begin{verbatim}
[()->gamma(1/3)^2/gamma(2/3),2*n+4*k+2+u,n*(3*n-3*u-1)]
[()->gamma(1/3)^2/gamma(2/3),4*n+8*k+4+2*u,(2*n-1)*(6*n-6*u+1)]
\end{verbatim}
Convergence type $E$ with $E=-3$ and $P=2k+u+5/3$.
\end{cf}

\smallskip

\begin{cf}\label{5.11.100}{\ }
\begin{verbatim}
[()->gamma(1/3)^2/gamma(2/3),[6,5*n],[-3,-2*n*(2*n+1)]]
\end{verbatim}
$$\dfrac{\G(1/3)^2}{\G(2/3)}=6-\dfrac{3}{5-\dfrac{6}{10-\dfrac{20}{15-\dfrac{42}{20-\dfrac{72}{25-\dfrac{110}{30-\ddots}}}}}}$$
Convergence type $E$ with $E=4$, $P=5/6$, and $C=-\G(1/3)^4/(8\pi^{3/2}3^{1/3})$, so that
$$\dfrac{\G(1/3)^2}{\G(2/3)}-\dfrac{p(n)}{q(n)}\sim-\dfrac{\G(1/3)^4/(\pi^{3/2}3^{1/3})}{2^{2n+3}n^{5/6}}$$
$$A=1-(265/216)/n+(189985/93312)/n^2+\cdots$$
Parametric family for $k$ and $u\ge0$:
\begin{verbatim}
[()->gamma(1/3)^2/gamma(2/3),5*n-u+3*k,-2*n*(2*n+1-2*u)]
\end{verbatim}
Convergence type $E$ with $E=4$ and $P=u+2k+5/6$.
\end{cf}

\smallskip

\begin{cf}\label{5.11.101}{\ }
\begin{verbatim}
[()->gamma(1/3)^2/gamma(2/3),[9/2,14,15*n+1],[9,-2*(3*n+2)*(6*n+1)]]
\end{verbatim}
$$\dfrac{\G(1/3)^2}{\G(2/3)}=9/2+\dfrac{9}{14-\dfrac{70}{31-\dfrac{208}{46-\dfrac{418}{61-\dfrac{700}{76-\dfrac{1054}{91-\ddots}}}}}}$$
Convergence type $E$ with $E=4$, $P=1/2$, and $C=\pi^{3/2}/(2^{1/3}3^{1/2}\G(2/3)^3)$, so that
$$\dfrac{\G(1/3)^2}{\G(2/3)}-\dfrac{p(n)}{q(n)}\sim\dfrac{\pi^{3/2}/(3^{1/2}\G(2/3)^3)}{2^{2n+1/3}n^{1/2}}$$
$$A=1-(7/8)/n+(499/384)/n^2-(70975/27648)/n^3+\cdots$$
Series:
$$\dfrac{\G(1/3)^2}{\G(2/3)}=\dfrac{9}{2}+\dfrac{9}{14}\sum_{n\ge0}\dfrac{(5/3)_n}{(13/6)_n}2^{-2n}$$
Parametric family for $u$, $v$, $w$ with $v+2w\le0$:
\begin{verbatim}
[()->gamma(1/3)^2/gamma(2/3),15*n+3*u+1,
                             -2*(3*n+9*v+3*u+3*w-4)*(6*n+12*w-5)]
\end{verbatim}
Convergence type $E$ with $E=4$ and $P=-(u+5v+5w)+11/2$.
\end{cf}

\smallskip

\begin{cf}\label{5.11.102}{\ }
\begin{verbatim}
[()->gamma(1/3)^2/gamma(2/3),[4,49,270*n^3-216*n^2-54*n+59],
[60,-2*(3*n-1)^2*(3*n+5)*(6*n-5)*(6*n+1)^2]]
\end{verbatim}
$$\dfrac{\G(1/3)^2}{\G(2/3)}=4+\dfrac{60}{49-\dfrac{3136}{1247-\dfrac{650650}{5243-\dfrac{8409856}{13667-\dfrac{48853750}{28139-\ddots}}}}}$$
Convergence type $E$ with $E=4$, $P=3/2$, and
$C=8\pi^{9/2}/(3^{9/2}\G(2/3)^9)$, so that
$$\dfrac{\G(1/3)^2}{\G(2/3)}-\dfrac{p(n)}{q(n)}\sim\dfrac{\pi^{9/2}/(3^{9/2}\G(2/3)^9)}{2^{2n-3}n^{3/2}}$$
$$A=1-(9/8)/n+(611/384)/n^2+\cdots$$
Series:
$$\dfrac{\G(1/3)^2}{\G(2/3)}=-2+6\sum_{n\ge0}\dfrac{(-1/3)_n^2(5/3)_n}{(1/6)_n(7/6)_n^2}2^{-2n}$$
\end{cf}

\smallskip

\begin{cf}\label{5.11.103}{\ }
\begin{verbatim}
[()->gamma(1/3)^2/gamma(2/3),
[12,2*(4*n-1)*(10*n^2-5*n-1)],
[-150,-n*(2*n+1)*(3*n-1)*(4*n-3)*(4*n+5)*(6*n+5)]]
\end{verbatim}
$$\dfrac{\G(1/3)^2}{\G(2/3)}=12-\dfrac{150}{24-\dfrac{594}{406-\dfrac{55250}{1628-\dfrac{591192}{4170-\dfrac{3135132}{8512-\dfrac{11453750}{15134-\ddots}}}}}}$$
Convergence type $E$ with $E=9$, $P=0$, and $C=-\G(1/3)^2/\G(2/3)/\sqrt{3}$,
so that
$$\dfrac{\G(1/3)^2}{\G(2/3)}-\dfrac{p(n)}{q(n)}\sim-\dfrac{\G(1/3)^2/\G(2/3)}{3^{2n+1/2}}$$
$$A=1+(59/144)/n-(9263/41472)/n^2+\cdots$$
\end{cf}
        
\smallskip

\begin{cf}\label{5.11.104}{\ }
\begin{verbatim}
[()->gamma(1/3)^2/gamma(2/3),[0,6*n-5],[6,n*(3*n-2)]]
\end{verbatim}
$$\dfrac{\G(1/3)^2}{\G(2/3)}=\dfrac{6}{1+\dfrac{1}{7+\dfrac{8}{13+\dfrac{21}{19+\dfrac{40}{25+\dfrac{65}{31+\ddots}}}}}}$$
Convergence type $E$ with $E=-(2+\sqrt{3})^2$, $P=0$, and
$C=3(\G(1/3)^2/\G(2/3))/(2+\sqrt{3})^{1/3}$, so that
$$\dfrac{\G(1/3)^2}{\G(2/3)}-\dfrac{p(n)}{q(n)}\sim(-1)^n\dfrac{3(\G(1/3)^2/\G(2/3))}{(2+\sqrt{3})^{2n+1/3}}$$
$$A=1-(5d/72)/n+(5d/432+25/3456)/n^2+\cdots$$
\end{cf}

\smallskip

\begin{cf}\label{5.11.105}{\ }
\begin{verbatim}
[()->gamma(1/3)^2/gamma(2/3),
[6,1326/31,1485*n^3-1881*n^2+420*n+44],
[-936/31,1872,3*n*(3*n-1)*(3*n+1)^2*(15*n-14)*(15*n+16)]]
\end{verbatim}
$$\dfrac{\G(1/3)^2}{\G(2/3)}=6-\dfrac{936/31}{1326/31+\dfrac{1872}{5240+\dfrac{1081920}{24470+\dfrac{13615200}{66668+\dfrac{77988768}{140744+\dfrac{298421760}{255608+\ddots}}}}}}$$
Convergence type $E$ with $E=-((1+\sqrt{5})/2)^{10}$, $P=0$, and $C=-3(\G(1/3)^2/\G(2/3))/((1+\sqrt{5})/2)^{19/3}$, so that
$$\dfrac{\G(1/3)^2}{\G(2/3)}-\dfrac{p(n)}{q(n)}\sim(-1)^{n+1}\dfrac{3\G(1/3)^2/\G(2/3)}{((1+\sqrt{5})/2)^{10n+19/3}}$$
$$A=1-(2d/45)/n+(22d/675+2/405)/n^2+\cdots$$
\end{cf}

\smallskip

For the next CFs, note that
$$2^{-1/3}\dfrac{\G(1/3)^2}{\G(2/3)}=\dfrac{\G(1/6)\G(2/3)}{2\G(1/2)}$$

\smallskip

\begin{cf}\label{5.11.106}{\ }
\begin{verbatim}
[()->2^(-1/3)*gamma(1/3)^2/gamma(2/3),[6,4,2],[-12,(3*n+1)*(3*n+2)]]
\end{verbatim}
$$2^{-1/3}\dfrac{\G(1/3)^2}{\G(2/3)}=6-\dfrac{12}{4+\dfrac{20}{2+\dfrac{56}{2+\dfrac{110}{2+\dfrac{182}{2+\dfrac{272}{2+\ddots}}}}}}$$
Convergence type $P^-$ with $P=2/3$ and $C=-\G(1/3)/\G(2/3)$, so that
$$2^{-1/3}\dfrac{\G(1/3)^2}{\G(2/3)}-\dfrac{p(n)}{q(n)}\sim(-1)^{n+1}\dfrac{\G(1/3)/\G(2/3)}{n^{2/3}}$$
$$A=1-(2/3)/n+(65/162)/n^2-(20/243)/n^3-(913/6561)/n^4+\cdots$$
Series:
$$2^{-1/3}\dfrac{\G(1/3)^2}{\G(2/3)}=6-3\sum_{n\ge0}(-1)^n\dfrac{(5/3)_n}{(7/3)_n}$$
Parametric family for $k\ge0$:
\begin{verbatim}
[()->2^(-1/3)*gamma(1/3)^2/gamma(2/3),6*k+2,(3*n+1)*(3*n+2)]
\end{verbatim}
Convergence type $P^-$ with $P=2k+2/3$.
\end{cf}

\smallskip

\begin{cf}\label{5.11.107}{\ }
\begin{verbatim}
[()->2^(-1/3)*gamma(1/3)^2/gamma(2/3),
[9/2,36*n^2-18*n+10],[-9/2,-9*n*(2*n+1)*(3*n+1)*(6*n+1)]]
\end{verbatim}
$$2^{-1/3}\dfrac{\G(1/3)^2}{\G(2/3)}=9/2-\dfrac{9/2}{28-\dfrac{756}{118-\dfrac{8190}{280-\dfrac{35910}{514-\dfrac{105300}{820-\dfrac{245520}{1198-\ddots}}}}}}$$
Convergence type $P^+$ with $P=1$ and $C=-(2^{-4/3}/9)\G(1/3)^2/\G(2/3)$, so
that
$$2^{-1/3}\dfrac{\G(1/3)^2}{\G(2/3)}-\dfrac{p(n)}{q(n)}\sim-\dfrac{(2^{-4/3}/9)\G(1/3)^2/\G(2/3)}{n}$$
Series:
\begin{align*}
  2^{-1/3}\dfrac{\G(1/3)^2}{\G(2/3)}&=\dfrac{9}{2}-\dfrac{9}{56}\sum_{n\ge0}\dfrac{n!(3/2)_n}{(7/3)_n(13/6)_n}\\
  2^{1/3}\dfrac{\G(2/3)}{\G(1/3)^2}&=\dfrac{2}{9}\sum_{n\ge0}\dfrac{(1/6)_n(1/3)_n}{n!(3/2)_n}\end{align*}
Parametric family for $k\ge0$:
\begin{verbatim}
[()->2^(-1/3)*gamma(1/3)^2/gamma(2/3),
36*n^2-18*n+10+18*k*(k+1),-9*n*(2*n+1)*(3*n+1)*(6*n+1)]
\end{verbatim}
Convergence type $P^+$ with $P=2k+1$.
\end{cf}

\smallskip

\begin{cf}\label{5.11.108}{\ }
\begin{verbatim}
[()->2^(-1/3)*gamma(1/3)^2/gamma(2/3),
[3,36*n^2-54*n+25],[6,-9*n*(2*n-1)*(3*n-2)*(6*n+1)]]
\end{verbatim}
$$2^{-1/3}\dfrac{\G(1/3)^2}{\G(2/3)}=3+\dfrac{6}{7-\dfrac{63}{61-\dfrac{2808}{187-\dfrac{17955}{385-\dfrac{63000}{655-\dfrac{163215}{997-\ddots}}}}}}$$
Convergence type $P^+$ with $P=1$ and $C=(2^{-1/3}/9)\G(1/3)^2/\G(2/3)$, so
that
$$2^{-1/3}\dfrac{\G(1/3)^2}{\G(2/3)}-\dfrac{p(n)}{q(n)}\sim\dfrac{(2^{-1/3}/9)\G(1/3)^2/\G(2/3)}{n}$$
Series:
\begin{align*}
  2^{-1/3}\dfrac{\G(1/3)^2}{\G(2/3)}&=3+\dfrac{6}{7}\sum_{n\ge0}\dfrac{n!(1/2)_n}{(4/3)_n(13/6)_n}\\
  2^{1/3}\dfrac{\G(2/3)}{\G(1/3)^2}&=\dfrac{1}{3}\sum_{n\ge0}\dfrac{(1/6)_n(-2/3)_n}{n!(1/2)_n}\end{align*}
Parametric family for $k\ge0$:
\begin{verbatim}
[()->2^(-1/3)*gamma(1/3)^2/gamma(2/3),
36*n^2-54*n+25+18*k*(k+1),-9*n*(2*n-1)*(3*n-2)*(6*n+1)]
\end{verbatim}
Convergence type $P^+$ with $P=2k+1$.
\end{cf}

\smallskip

\begin{cf}\label{5.11.109}{\ }
\begin{verbatim}
[()->2^(-1/3)*gamma(1/3)^2/gamma(2/3),[3],[6,3*n*(3*n+2)]]
\end{verbatim}
$$2^{-1/3}\dfrac{\G(1/3)^2}{\G(2/3)}=3+\dfrac{6}{3+\dfrac{15}{3+\dfrac{48}{3+\dfrac{99}{3+\dfrac{168}{3+\dfrac{255}{3+\ddots}}}}}}$$
Convergence type $P^-$ with $P=1$ and $C=(2^{-1/3}/3)\G(1/3)^2/\G(2/3)$, so that
$$2^{-1/3}\dfrac{\G(1/3)^2}{\G(2/3)}-\dfrac{p(n)}{q(n)}\sim(-1)^n\dfrac{(2^{-1/3}/3)\G(1/3)^2/\G(2/3)}{n}$$
$$A=1-(5/6)/n+(17/36)/n^2-(5/216)/n^3-(77/432)/n^4+\cdots$$
Series:
\begin{align*}
2^{-1/3}\dfrac{\G(1/3)^2}{\G(2/3)}&=3+\dfrac{3}{4}\sum_{n\ge0}\dfrac{n!(5/6)_n}{(3/2)_n(7/3)_n}\\
2^{1/3}\dfrac{\G(2/3)}{\G(1/3)^2}&=\dfrac{1}{5}\sum_{n\ge0}\dfrac{(1/3)_n(1/2)_n}{n!(11/6)_n}\end{align*}
Parametric family for $k\ge0$:
\begin{verbatim}
[()->2^(-1/3)*gamma(1/3)^2/gamma(2/3),6*k+3,3*n*(3*n+2)]
\end{verbatim}
Convergence type $P^-$ with $P=2k+1$.
\end{cf}

\smallskip

\begin{cf}\label{5.11.110}{\ }
\begin{verbatim}
[()->2^(-1/3)*gamma(1/3)^2/gamma(2/3),
[3,7,9*(2*n-1)],[12,(9*n^2-4)*(36*n^2-1)]]
\end{verbatim}
$$2^{-1/3}\dfrac{\G(1/3)^2}{\G(2/3)}=3+\dfrac{12}{7+\dfrac{175}{27+\dfrac{4576}{45+\dfrac{24871}{63+\dfrac{80500}{81+\dfrac{198679}{99+\ddots}}}}}}$$
Convergence type $P^-$ with $P=1$ and $C=(2^{-5/3}/3)(\G(1/3)/\G(2/3))^3$,
so that
$$2^{-1/3}\dfrac{\G(1/3)^2}{\G(2/3)}-\dfrac{p(n)}{q(n)}\sim(-1)^n\dfrac{(2^{-5/3}/3)(\G(1/3)/\G(2/3))^3}{n}$$
$$A=1-(1/2)/n+(1/18)/n^2+(1/6)/n^3-(1/432)/n^4+\cdots$$
Series:
$$2^{-1/3}\dfrac{\G(1/3)^2}{\G(2/3)}=3+\dfrac{12}{7}\sum_{n\ge0}(-1)^n\dfrac{(5/3)_n(5/6)_n}{(4/3)_n(13/6)_n}$$
Parametric family for $k\ge0$:
\begin{verbatim}
[()->2^(-1/3)*gamma(1/3)^2/gamma(2/3),
9*(4*k+1)*(2*n-1),(9*n^2-4)*(36*n^2-1)]
\end{verbatim}
Convergence type $P^-$ with $P=4k+1$.
\end{cf}

\smallskip

\begin{cf}\label{5.11.111}{\ }
\begin{verbatim}
[()->2^(-1/3)*gamma(1/3)^2/gamma(2/3),
[9/2,28,27*(2*n-1)],[-9,(9*n^2-1)*(36*n^2-1)]]
\end{verbatim}
$$2^{-1/3}\dfrac{\G(1/3)^2}{\G(2/3)}=9/2-\dfrac{9}{28+\dfrac{280}{81+\dfrac{5005}{135+\dfrac{25840}{189+\dfrac{82225}{243+\dfrac{201376}{297+\ddots}}}}}}$$
Convergence type $P^-$ with $P=3$ and $C=-(2^{-8/3}/9)(\G(1/3)/\G(2/3))^3$,
so that
$$2^{-1/3}\dfrac{\G(1/3)^2}{\G(2/3)}-\dfrac{p(n)}{q(n)}\sim(-1)^{n+1}\dfrac{(2^{-8/3}/9)(\G(1/3)/\G(2/3))^3}{n^3}$$
$$A=1-(3/2)/n+(1/6)/n^2+(25/12)/n^3+(7/432)/n^4+\cdots$$
Series:
$$2^{-1/3}\dfrac{\G(1/3)^2}{\G(2/3)}=\dfrac{9}{2}-\dfrac{9}{28}\sum_{n\ge0}(-1)^n\dfrac{(2/3)_n(5/6)_n}{(7/3)_n(13/6)_n}$$
Parametric family for $k\ge0$:
\begin{verbatim}
[()->2^(-1/3)*gamma(1/3)^2/gamma(2/3),
9*(4*k+3)*(2*n-1),(9*n^2-1)*(36*n^2-1)]
\end{verbatim}
Convergence type $P^-$ with $P=4k+3$.
\end{cf}

\smallskip

\begin{cf}\label{5.11.112}{\ }
\begin{verbatim}
[()->2^(-1/3)*gamma(1/3)^2/gamma(2/3),[4,49,72*n^2-72*n+50],
                                      [8,-(36*n^2-1)^2]]
\end{verbatim}
$$2^{-1/3}\dfrac{\G(1/3)^2}{\G(2/3)}=4+\dfrac{8}{49-\dfrac{1225}{194-\dfrac{20449}{482-\dfrac{104329}{914-\dfrac{330625}{1490-\dfrac{808201}{2210-\ddots}}}}}}$$
Convergence type $P^+$ with $P=5/3$ and $C=(1/135)2^{-1/3}(\G(1/3)/\G(2/3))^4$,
so that
$$2^{-1/3}\dfrac{\G(1/3)^2}{\G(2/3)}-\dfrac{p(n)}{q(n)}\sim\dfrac{(1/135)2^{-1/3}(\G(1/3)/\G(2/3))^4}{n^{5/3}}$$
$$A=1-(5/6)/n+(365/891)/n^2-(35/486)/n^3+\cdots$$
Series:
$$2^{-1/3}\dfrac{\G(1/3)^2}{\G(2/3)}=4+\dfrac{8}{49}\sum_{n\ge0}\dfrac{(5/6)_n^2}{(13/6)_n^2}$$
Parametric family for $k\ge0$:
\begin{verbatim}
[()->2^(-1/3)*gamma(1/3)^2/gamma(2/3),
72*n^2-72*n+50+12*k*(3*k+5),-(36*n^2-1)^2]
\end{verbatim}
Convergence type $P^+$ with $P=2k+5/3$.
\end{cf}

\smallskip

\begin{cf}\label{5.11.113}{\ }
\begin{verbatim}
[()->2^(-1/3)*gamma(1/3)^2/gamma(2/3),
[4,3*(6*n-1)],[4,3*n*(3*n+1)^2*(3*n+2)]]
\end{verbatim}
$$2^{-1/3}\dfrac{\G(1/3)^2}{\G(2/3)}=4+\dfrac{4}{15+\dfrac{240}{33+\dfrac{2352}{51+\dfrac{9900}{69+\dfrac{28392}{87+\dfrac{65280}{105+\ddots}}}}}}$$
Convergence type $P^-$ with $P=2$ and $C=(1/18)2^{-1/3}\G(1/3)^2/\G(2/3)$,
so that
$$2^{-1/3}\dfrac{\G(1/3)^2}{\G(2/3)}-\dfrac{p(n)}{q(n)}\sim(-1)^n\dfrac{(1/18)2^{-1/3}\G(1/3)^2/\G(2/3)}{n^2}$$
$$A=1-(5/3)/n+(25/18)/n^2-(425/432)/n^4+\cdots$$
Series:
\begin{align*}
  2^{-1/3}\dfrac{\G(1/3)^2}{\G(2/3)}&=4+\dfrac{4}{245}\sum_{n\ge0}\dfrac{(12n+11)n!(2/3)_n^2(4/3)_n}{(3/2)_n(11/6)_n(13/6)_n^2}\\
  2^{1/3}\dfrac{\G(2/3)}{\G(1/3)^2}&=\dfrac{3}{64}\sum_{n\ge0}\dfrac{(12n+5)(1/2)_n(1/6)_n^2(5/6)_n}{n!(4/3)_n(5/3)_n^2}\end{align*}
Parametric family for $k\ge0$:
\begin{verbatim}
[()->2^(-1/3)*gamma(1/3)^2/gamma(2/3),
3*(2*k+1)*(6*n-1),3*n*(3*n+1)^2*(3*n+2)]
\end{verbatim}
Convergence type $P^-$ with $P=4k+2$.
\end{cf}

\smallskip

\begin{cf}\label{5.11.114}{\ }
\begin{verbatim}
[()->2^(-1/3)*gamma(1/3)^2/gamma(2/3),[8,17*n],[-40,-12*n*(6*n+5)]]
\end{verbatim}
$$2^{-1/3}\dfrac{\G(1/3)^2}{\G(2/3)}=8-\dfrac{40}{17-\dfrac{132}{34-\dfrac{408}{51-\dfrac{828}{68-\dfrac{1392}{85-\dfrac{2100}{102-\ddots}}}}}}$$
Convergence type $E$ with $E=9/8$, $P=17/6$, and $C=-2^{20/3}\G(1/3)^4/(\sqrt{3}\pi^{3/2})$, so that
$$2^{-1/3}\dfrac{\G(1/3)^2}{\G(2/3)}-\dfrac{p(n)}{q(n)}\sim-\dfrac{2^{20/3}\G(1/3)^4/(\sqrt{3}\pi^{3/2})}{(9/8)^nn^{17/6}}$$
$$A=1-(2737/72)/n+(14373857/10368)/n^2+\cdots$$
Parametric family for $k\ge0$:
\begin{verbatim}
[()->2^(-1/3)*gamma(1/3)^2/gamma(2/3),17*n+k,-12*n*(6*n+5)]
\end{verbatim}
Convergence type $E$ with $E=9/8$ and $P=2k+17/6$.

There exist two-variable families as above which I have been too lazy to
write down explicitly.
\end{cf}

\smallskip

\begin{cf}\label{5.11.115}{\ }
\begin{verbatim}
[()->2^(-1/3)*gamma(1/3)^2/gamma(2/3),[16,14*n^2+13*n+1],
                             [-192,-8*n*(n+1)*(2*n+3)*(3*n+4)]]
\end{verbatim}
$$2^{-1/3}\dfrac{\G(1/3)^2}{\G(2/3)}=16-\dfrac{192}{28-\dfrac{560}{83-\dfrac{3360}{166-\dfrac{11232}{277-\dfrac{28160}{416-\dfrac{59280}{583-\ddots}}}}}}$$
Convergence type $E$ with $E=4/3$, $P=1/6$, and $C=-3^{1/2}\pi^{5/2}/(2^{2/3}\G(2/3)^5)$, so that
$$2^{-1/3}\dfrac{\G(1/3)^2}{\G(2/3)}-\dfrac{p(n)}{q(n)}\sim-\dfrac{3^{1/2}\pi^{5/2}/(2^{2/3}\G(2/3)^5)}{(4/3)^nn^{1/6}}$$
$$A=1-(77/72)/n+(17267/3456)/n^2+\cdots$$
Series:
$$2^{1/3}\dfrac{\G(2/3)}{\G(1/3)^2}=\dfrac{1}{16}\sum_{n\ge0}\dfrac{(3/2)_n(4/3)_n}{n!(n+1)!}(3/4)^n$$
\end{cf}

\smallskip

\begin{cf}\label{5.11.116}{\ }
\begin{verbatim}
[()->2^(-1/3)*gamma(1/3)^2/gamma(2/3),[3,24,28*n^2+1],
                                [15,-8*n*(2*n+1)^2*(6*n+5)]]
\end{verbatim}
$$2^{-1/3}\dfrac{\G(1/3)^2}{\G(2/3)}=3+\dfrac{15}{24-\dfrac{792}{113-\dfrac{6800}{253-\dfrac{27048}{449-\dfrac{75168}{701-\dfrac{169400}{1009-\ddots}}}}}}$$
Convergence type $E$ with $E=4/3$, $P=7/6$, and $C=3^{3/2}\pi^{1/2}/(2^{1/3}\G(2/3)^2)$, so that
$$2^{-1/3}\dfrac{\G(1/3)^2}{\G(2/3)}-\dfrac{p(n)}{q(n)}\sim\dfrac{3^{3/2}\pi^{1/2}/(2^{1/3}\G(2/3)^2)}{(4/3)^nn^{7/6}}$$
$$A=1-(377/72)/n+(46793/1152)/n^2+\cdots$$
Series:
$$2^{-1/3}\dfrac{\G(1/3)^2}{\G(2/3)}=3\sum_{n\ge0}\dfrac{(5/6)_n}{(2n+1)n!}(3/4)^n$$
\end{cf}

\smallskip

\begin{cf}\label{5.11.117}{\ }
\begin{verbatim}
[()->2^(-1/3)*gamma(1/3)^2/gamma(2/3),[6,9*n],[-12,-6*n*(3*n+2)]]
\end{verbatim}
$$2^{-1/3}\dfrac{\G(1/3)^2}{\G(2/3)}=6-\dfrac{12}{9-\dfrac{30}{18-\dfrac{96}{27-\dfrac{198}{36-\dfrac{336}{45-\dfrac{510}{54-\ddots}}}}}}$$
Convergence type $E$ with $E=2$, $P=1$, and $C=-(2^{1/3}/3)\G(1/3)^2/\G(2/3)$,
so that
$$2^{-1/3}\dfrac{\G(1/3)^2}{\G(2/3)}-\dfrac{p(n)}{q(n)}\sim-\dfrac{(1/3)\G(1/3)^2/\G(2/3)}{2^{n-1/3}n}$$
$$A=1-2/n+(56/9)/n^2-(2188/81)/n^3+\cdots$$
Parametric family for $k\ge0$:
\begin{verbatim}
[()->2^(-1/3)*gamma(1/3)^2/gamma(2/3),9*n+3*k,-6*n*(3*n+2)]
\end{verbatim}
Convergence type $E$ with $E=2$ and $P=2k+1$.

There exist two-variable families as above which I have been too lazy to
write down explicitly.
\end{cf}

\smallskip

\begin{cf}\label{5.11.118}{\ }
\begin{verbatim}
[()->2^(-1/3)*gamma(1/3)^2/gamma(2/3),[8/3,7*n],[40/3,4*n*(2*n+5)]]
\end{verbatim}
$$2^{-1/3}\dfrac{\G(1/3)^2}{\G(2/3)}=8/3+\dfrac{40/3}{7+\dfrac{28}{14+\dfrac{72}{21+\dfrac{132}{28+\dfrac{208}{35+\dfrac{300}{42+\ddots}}}}}}$$
Convergence type $E$ with $E=-8$, $P=-7/6$, and
$C=3^{7/3}\pi^{5/2}/(2^{20/3}\G(2/3)^4)$, so that
$$2^{-1/3}\dfrac{\G(1/3)^2}{\G(2/3)}-\dfrac{p(n)}{q(n)}\sim(-1)^n\dfrac{3^{7/3}\pi^{5/2}/\G(2/3)^4}{2^{3n+20/3}n^{-7/6}}$$
$$A=1+(1813/648)/n+(540841/839808)/n^2+\cdots$$
Parametric family for $k\ge0$:
\begin{verbatim}
[()->2^(-1/3)*gamma(1/3)^2/gamma(2/3),7*n+9*k,4*n*(2*n+5)]
\end{verbatim}
Convergence type $E$ with $E=-8$ and $P=2k-7/6$.

There exist two-variable families as above which I have been too lazy to
write down explicitly.
\end{cf}

\smallskip

\begin{cf}\label{5.11.119}{\ }
\begin{verbatim}
[()->2^(-1/3)*gamma(1/3)^2/gamma(2/3),[9/2,28,21*n+5],[-9,4*(3*n+2)*(6*n+1)]]
\end{verbatim}
$$2^{-1/3}\dfrac{\G(1/3)^2}{\G(2/3)}=9/2-\dfrac{9}{28+\dfrac{140}{47+\dfrac{416}{68+\dfrac{836}{89+\dfrac{1400}{110+\dfrac{2108}{131+\ddots}}}}}}$$
Convergence type $E$ with $E=-8$, $P=1/2$, and $C=-\G(1/3)^3/(2^{10/3}\pi^{3/2})$, so that
$$2^{-1/3}\dfrac{\G(1/3)^2}{\G(2/3)}-\dfrac{p(n)}{q(n)}\sim(-1)^{n+1}\dfrac{\G(1/3)^3/\pi^{3/2}}{2^{3n+10/3}}$$
$$A=1-(47/72)/n+(2027/3456)/n^2+\cdots$$
Series:
$$2^{-1/3}\dfrac{\G(1/3)^2}{\G(2/3)}=\dfrac{9}{2}-\dfrac{9}{28}\sum_{n\ge0}(-1)^n\dfrac{(5/3)_n}{(13/6)_n}2^{-3n}$$
Parametric family for $k\ge0$:
\begin{verbatim}
[()->2^(-1/3)*gamma(1/3)^2/gamma(2/3),21*n+5+27*k,4*(3*n+2)*(6*n+1)]
\end{verbatim}
Convergence type $E$ with $E=-8$ and $P=2k+1/2$.
\end{cf}

\smallskip

\begin{cf}\label{5.11.120}{\ }
\begin{verbatim}
[()->2^(-1/3)*gamma(1/3)^2/gamma(2/3),[4,10*n],[2,-3*n*(3*n-1)]]
\end{verbatim}
$$2^{-1/3}\dfrac{\G(1/3)^2}{\G(2/3)}=4+\dfrac{2}{10-\dfrac{6}{20-\dfrac{30}{30-\dfrac{72}{40-\dfrac{132}{50-\dfrac{210}{60-\ddots}}}}}}$$
Convergence type $E$ with $E=9$, $P=5/3$, and $C=\pi^3/(2^{11/3}\cdot3\G(2/3)^5)$, so that
$$2^{-1/3}\dfrac{\G(1/3)^2}{\G(2/3)}-\dfrac{p(n)}{q(n)}\sim\dfrac{\pi^3/(2^{11/3}\G(2/3)^5)}{3^{2n+1}n^{5/3}}$$
$$A=1-(245/144)/n+(12455/4608)/n^2+\cdots$$
Parametric family for $k\ge0$:
\begin{verbatim}
[()->2^(-1/3)*gamma(1/3)^2/gamma(2/3),10*n+8*k,-3*n*(3*n-1)]
\end{verbatim}
Convergence type $E$ with $E=9$ and $P=2k+5/3$.

There exist two-variable families as above which I have been too lazy to
write down explicitly.
\end{cf}

\smallskip

\begin{cf}\label{5.11.121}{\ }
\begin{verbatim}
[()->2^(-1/3)*gamma(1/3)^2/gamma(2/3),[4,21,20*n+2],[4,-3*(2*n+1)*(6*n+1)]]
\end{verbatim}
$$2^{-1/3}\dfrac{\G(1/3)^2}{\G(2/3)}=4+\dfrac{4}{21-\dfrac{63}{42-\dfrac{195}{62-\dfrac{399}{82-\dfrac{675}{102-\dfrac{1023}{122-\ddots}}}}}}$$
Convergence type $E$ with $E=9$, $P=2/3$, and $C=2^{-1/3}\G(1/3)/(6\G(2/3))$, so that
$$2^{-1/3}\dfrac{\G(1/3)^2}{\G(2/3)}-\dfrac{p(n)}{q(n)}\sim\dfrac{2^{-4/3}\G(1/3)/\G(2/3)}{3^{2n+1}n^{2/3}}$$
$$A=1-(35/36)/n+(1075/864)/n^2+\cdots$$
Series:
$$2^{-1/3}\dfrac{\G(1/3)^2}{\G(2/3)}=4\sum_{n\ge0}\dfrac{(1/2)_n}{(7/6)_n}3^{-2n}$$
Parametric family for $k\ge0$:
\begin{verbatim}
[()->2^(-1/3)*gamma(1/3)^2/gamma(2/3),20*n+2+16*k,-3*(2*n+1)*(6*n+1)]
\end{verbatim}
Convergence type $E$ with $E=9$ and $P=2k+2/3$.
\end{cf}

\smallskip

\begin{cf}\label{5.11.122}{\ }
\begin{verbatim}
[()->2^(-1/3)*gamma(1/3)^2/gamma(2/3),[9/2,18*n+3],[-6,-3*n*(3*n+4)]]
\end{verbatim}
$$2^{-1/3}\dfrac{\G(1/3)^2}{\G(2/3)}=9/2-\dfrac{6}{21-\dfrac{21}{39-\dfrac{60}{57-\dfrac{117}{75-\dfrac{192}{93-\dfrac{285}{111-\ddots}}}}}}$$
Convergence type $E$ with $E=(1+\sqrt{2})^4$, $P=0$, and
$C=-3\G(1/3)^3/(2^{1/3}\pi(1+\sqrt{2})^{14/3})$, so that
$$2^{-1/3}\dfrac{\G(1/3)^2}{\G(2/3)}-\dfrac{p(n)}{q(n)}\sim-\dfrac{3\G(1/3)^3/(2^{1/3}\pi)}{(1+\sqrt{2})^{4n+14/3}}$$
$$A=1+(7d/48)/n+(-49d/288+49/2304)/n^2+\cdots$$
\end{cf}

\smallskip

\begin{cf}\label{5.11.123}{\ }
\begin{verbatim}
[()->2^(-1/3)*gamma(1/3)^2/gamma(2/3),
[6,12*(9*n^2-6*n-1)],[-42,-3*n*(3*n+1)*(6*n-5)*(6*n+7)]]
\end{verbatim}
$$2^{-1/3}\dfrac{\G(1/3)^2}{\G(2/3)}=6-\dfrac{42}{24-\dfrac{156}{276-\dfrac{5586}{744-\dfrac{29250}{1428-\dfrac{91884}{2328-\dfrac{222000}{3444-\ddots}}}}}}$$
Convergence type $E$ with $E=(1+\sqrt{2})^4$, $P=0$, and $C=-2^{2/3}3^{1/2}\G(1/3)^2/(\G(2/3)(1+\sqrt{2})^{8/3})$, so that
$$2^{-1/3}\dfrac{\G(1/3)^2}{\G(2/3)}-\dfrac{p(n)}{q(n)}\sim-\dfrac{2^{2/3}3^{1/2}\G(1/3)^2/\G(2/3)}{(1+\sqrt{2})^{4n+8/3}}$$
$$A=1+(d/16)/n+(-d/24+1/256)/n^2+\cdots$$
\end{cf}

\smallskip

For the next CFs, note that
$$2^{-2/3}\dfrac{\G(1/3)^2}{\G(2/3)}=\dfrac{\G(1/2)\G(1/3)^3}{\G(1/6)\G(2/3)^2}$$

\smallskip

\begin{cf}\label{5.11.124}{\ }
\begin{verbatim}
[()->2^(-2/3)*gamma(1/3)^2/gamma(2/3),
[9/4,6*n+2],[45/2,9*n*(n+1)*(3*n+2)*(3*n+5)]]
\end{verbatim}
$$2^{-2/3}\dfrac{\G(1/3)^2}{\G(2/3)}=9/4+\dfrac{45/2}{8+\dfrac{720}{14+\dfrac{4752}{20+\dfrac{16632}{26+\dfrac{42840}{32+\dfrac{91800}{38+\ddots}}}}}}$$
Convergence type $P^-$ with $P=2/3$ and $C=(2^{-4/3}/3)(\G(1/3)/\G(2/3))^4$,
so that
$$2^{-2/3}\dfrac{\G(1/3)^2}{\G(2/3)}-\dfrac{p(n)}{q(n)}\sim\dfrac{(2^{-4/3}/3)(\G(1/3)/\G(2/3))^4}{n^{2/3}}$$
$$A=1-(8/9)/n+(49/54)/n^2-(1936/2187)/n^3+\cdots$$
Series:
$$2^{2/3}\dfrac{\G(2/3)}{\G(1/3)^2}=\dfrac{4}{9}\sum_{n\ge0}(-1)^n\dfrac{(2/3)_n(5/3)_n}{n!(n+1)!}$$
Parametric family for $k\ge0$ and $u\ge0$:
\begin{verbatim}
[()->2^(-2/3)*gamma(1/3)^2/gamma(2/3),
2*(6*k+1)*(3*n+1+3*u),9*n*(n+1)*(3*n+2+6*u)*(3*n+5+6*u)]
\end{verbatim}
Convergence type $P^-$ with $P=4k+2/3$.
\end{cf}

\smallskip

\begin{cf}\label{5.11.125}{\ }
\begin{verbatim}
[()->2^(-2/3)*gamma(1/3)^2/gamma(2/3),
[3,8*(3*n-2)],[3,9*n^2*(3*n-1)^2]]
\end{verbatim}
$$2^{-2/3}\dfrac{\G(1/3)^2}{\G(2/3)}=3+\dfrac{3}{8+\dfrac{36}{32+\dfrac{900}{56+\dfrac{5184}{80+\dfrac{17424}{104+\dfrac{44100}{128+\ddots}}}}}}$$
Convergence type $P^-$ with $P=8/3$ and $C=(2^{-7/3}/27)(\G(1/3)/\G(2/3))^4$,
so that
$$2^{-2/3}\dfrac{\G(1/3)^2}{\G(2/3)}-\dfrac{p(n)}{q(n)}\sim\dfrac{(2^{-7/3}/27)(\G(1/3)/\G(2/3))^4}{n^{8/3}}$$
$$A=1-(8/9)/n-(16/27)/n^2+(3248/2187)/n^3+\cdots$$
Series:
$$2^{2/3}\dfrac{\G(2/3)}{\G(1/3)^2}=\dfrac{1}{3}\sum_{n\ge0}(-1)^n\dfrac{(-1/3)_n^2}{n!^2}$$
Parametric family for $k\ge0$ and $u\ge0$:
\begin{verbatim}
[()->2^(-2/3)*gamma(1/3)^2/gamma(2/3),
4*(3*k+2)*(3*n-2+3*u),9*n^2*(3*n-1+6*u)^2]
\end{verbatim}
Convergence type $P^-$ with $P=4k+8/3$.
\end{cf}

\smallskip

\begin{cf}\label{5.11.126}{\ }
\begin{verbatim}
[()->2^(-2/3)*gamma(1/3)^2/gamma(2/3),
[3,21,36*n^2-24*n+10],[3,-3*n*(3*n+1)*(6*n+1)^2]]
\end{verbatim}
$$2^{-2/3}\dfrac{\G(1/3)^2}{\G(2/3)}=3+\dfrac{3}{21-\dfrac{588}{106-\dfrac{7098}{262-\dfrac{32490}{490-\dfrac{97500}{790-\dfrac{230640}{1162-\ddots}}}}}}$$
Convergence type $P^+$ with $P=2/3$ and $C=3/(4\G(1/3))$, so that
$$2^{-2/3}\dfrac{\G(1/3)^2}{\G(2/3)}-\dfrac{p(n)}{q(n)}\sim\dfrac{3/(4\G(1/3))}{n^{2/3}}$$
$$A=1-(4/9)/n+(85/432)/n^2-(125/2187)/n^3+\cdots$$
Series:
$$2^{-2/3}\dfrac{\G(1/3)^2}{\G(2/3)}=3\sum_{n\ge0}\dfrac{(1/3)_n}{(6n+1)n!}$$
Parametric family for $k\ge0$:
\begin{verbatim}
[()->2^(-2/3)*gamma(1/3)^2/gamma(2/3),
36*n^2-24*n+10+6*k*(3*k+2),-3*n*(3*n+1)*(6*n+1)^2]
\end{verbatim}
Convergence type $P^+$ with $P=2k+2/3$.
\end{cf}

\smallskip

\begin{cf}\label{5.11.127}{\ }
\begin{verbatim}
[()->2^(-2/3)*gamma(1/3)^2/gamma(2/3),
[9/4,21,36*n^2-60*n+55],[45/2,-3*n*(3*n-2)*(6*n-5)*(6*n+1)]]
\end{verbatim}
$$2^{-2/3}\dfrac{\G(1/3)^2}{\G(2/3)}=9/4+\dfrac{45/2}{21-\dfrac{21}{79-\dfrac{2184}{199-\dfrac{15561}{391-\dfrac{57000}{655-\dfrac{151125}{991-\ddots}}}}}}$$
Convergence type $P^+$ with $P=8/3$ and $C=5/(64\G(1/3))$, so that
$$2^{-2/3}\dfrac{\G(1/3)^2}{\G(2/3)}-\dfrac{p(n)}{q(n)}\sim\dfrac{5/(64\G(1/3))}{n^{8/3}}$$
$$A=1-(4/9)/n-(11/189)/n^2+(253/2187)/n^3+\cdots$$
Series:
$$2^{-2/3}\dfrac{\G(1/3)^2}{\G(2/3)}=-\dfrac{15}{4}\sum_{n\ge0}\dfrac{(-2/3)_n}{(6n-5)n!}$$
Parametric family for $k\ge0$:
\begin{verbatim}
[()->2^(-2/3)*gamma(1/3)^2/gamma(2/3),
36*n^2-60*n+55+6*k*(3*k+8),-3*n*(3*n-2)*(6*n-5)*(6*n+1)]
\end{verbatim}
Convergence type $P^+$ with $P=2k+8/3$.
\end{cf}

\smallskip

\begin{cf}\label{5.11.128}{\ }
\begin{verbatim}
[()->2^(-2/3)*gamma(1/3)^2/gamma(2/3),[4,14*n^2+3*n+1],
                             [-8,-8*n*(n+1)*(2*n+1)*(3*n+2)]]
\end{verbatim}
$$2^{-2/3}\dfrac{\G(1/3)^2}{\G(2/3)}=4-\dfrac{8}{18-\dfrac{240}{63-\dfrac{1920}{136-\dfrac{7392}{237-\dfrac{20160}{366-\dfrac{44880}{523-\ddots}}}}}}$$
Convergence type $E$ with $E=4/3$, $P=11/6$, and
$C=-2^{2/3}\pi^{7/2}/(3\G(2/3)^7)$, so that
$$2^{-2/3}\dfrac{\G(1/3)^2}{\G(2/3)}-\dfrac{p(n)}{q(n)}\sim-\dfrac{2^{2/3}\pi^{7/2}/(3\G(2/3)^7)}{(4/3)^nn^{11/6}}$$
$$A=1-(617/72)/n+(912401/10368)/n^2+\cdots$$
Series:
$$2^{2/3}\dfrac{\G(2/3)}{\G(1/3)^2}=\dfrac{1}{4}\sum_{n\ge0}\dfrac{(1/2)_n(2/3)_n}{n!(n+1)!}(3/4)^n$$
\end{cf}

\smallskip

\begin{cf}\label{5.11.129}{\ }
\begin{verbatim}
[()->gamma(2/3)^2/gamma(1/3),
[1/2,15,18*n^2-3*n+2],[1,-3*n*(3*n+1)*(3*n+2)^2]]
\end{verbatim}
$$\dfrac{\G(2/3)^2}{\G(1/3)}=1/2+\dfrac{1}{15-\dfrac{300}{68-\dfrac{2688}{155-\dfrac{10890}{278-\dfrac{30576}{437-\dfrac{69360}{632-\ddots}}}}}}$$
Convergence type $P^+$ with $P=2/3$ and $C=1/(2\G(1/3))$, so that
$$\dfrac{\G(2/3)^2}{\G(1/3)}-\dfrac{p(n)}{q(n)}\sim\dfrac{1/(2\G(1/3))}{n^{2/3}}\;.$$
$$A=1-(29/45)/n+(13/27)/n^2-(8422/24057)/n^3+(4661/19683)/n^4+\cdots$$
Series:
$$\dfrac{\G(2/3)^2}{\G(1/3)}=\dfrac{1}{2}+\dfrac{1}{3}\sum_{n\ge0}\dfrac{(4/3)_n}{(3n+5)(n+1)!}$$
Parametric family for $k\ge0$:
\begin{verbatim}
[()->gamma(2/3)^2/gamma(1/3),
18*n^2-3*n+2+3*k*(3*k+2),-3*n*(3*n+1)*(3*n+2)^2]
\end{verbatim}
Convergence type $P^+$ with $P=2k+2/3$.
\end{cf}

\smallskip

\begin{cf}\label{5.11.130}{\ }
\begin{verbatim}
[()->gamma(2/3)^2/gamma(1/3),
[3/4,18*n^2-3*n+4],[-3/4,-9*n*(n+1)*(3*n+1)^2]]
\end{verbatim}
$$\dfrac{\G(2/3)^2}{\G(1/3)}=3/4-\dfrac{3/4}{19-\dfrac{288}{70-\dfrac{2646}{157-\dfrac{10800}{280-\dfrac{30420}{439-\dfrac{69120}{634-\ddots}}}}}}$$
Convergence type $P^+$ with $P=4/3$ and $C=-16\pi^4/(9\G(1/3)^8)$,
so that
$$\dfrac{\G(2/3)^2}{\G(1/3)}-\dfrac{p(n)}{q(n)}\sim-\dfrac{16\pi^4/(9\G(1/3)^8)}{n^{4/3}}$$
$$A=1-(86/63)/n+(127/81)/n^2-(46990/28431)/n^3+\cdots$$
Series:
$$\dfrac{\G(1/3)}{\G(2/3)^2}=\dfrac{4}{3}\sum_{n\ge0}\dfrac{(1/3)_n^2}{(n+1)n!^2}$$
Parametric family for $k\ge0$:
\begin{verbatim}
[()->gamma(2/3)^2/gamma(1/3),18*n^2-3*n+(3*k+2)^2,-9*n*(n+1)*(3*n+1)^2]
\end{verbatim}
Convergence type $P^+$ with $P=2k+4/3$.
\end{cf}

\smallskip

\begin{cf}\label{5.11.131}{\ }
\begin{verbatim}
[()->gamma(2/3)^2/gamma(1/3),
[3/4,25,18*n^2+6*n+4],[-3/4,-3*(n+1)*(3*n+1)*(3*n+2)^2]]
\end{verbatim}
$$\dfrac{\G(2/3)^2}{\G(1/3)}=3/4-\dfrac{3/4}{25-\dfrac{600}{88-\dfrac{4032}{184-\dfrac{14520}{316-\dfrac{38220}{484-\dfrac{83232}{688-\ddots}}}}}}$$
Convergence type $P^+$ with $P=1$ and $C=-4\pi^2/(27\G(1/3)^3)$, so that
$$\dfrac{\G(2/3)^2}{\G(1/3)}-\dfrac{p(n)}{q(n)}\sim-\dfrac{4\pi^2/(27\G(1/3)^3)}{n}$$
$$A=1-(10/9)/n+(284/243)/n^2-(2512/2187)/n^3+\cdots$$
Series:
\begin{align*}\dfrac{\G(2/3)^2}{\G(1/3)}&=\dfrac{3}{4}-\dfrac{3}{100}\sum_{n\ge0}\dfrac{(n+1)!(4/3)_n}{(8/3)_n^2}\\
\dfrac{\G(1/3)}{\G(2/3)^2}&=1+\dfrac{1}{3}\sum_{n\ge0}\dfrac{(2/3)_n^2}{(n+1)!(4/3)_n}\end{align*}
Parametric family for $k\ge0$:
\begin{verbatim}
[()->gamma(2/3)^2/gamma(1/3),18*n^2+6*n+4+9*k*(k+1),
                                 -3*(n+1)*(3*n+1)*(3*n+2)^2]
\end{verbatim}
Convergence type $P^+$ with $P=2k+1$.
\end{cf}

\smallskip

\begin{cf}\label{5.11.132}{\ }
\begin{verbatim}
[()->gamma(2/3)^2/gamma(1/3),[3/4,22,21*(2*n-1)],
                                 [-3/2,(3*n-1)^2*(3*n+1)^2]]
\end{verbatim}
$$\dfrac{\G(2/3)^2}{\G(1/3)}=3/4-\dfrac{3/2}{22+\dfrac{64}{63+\dfrac{1225}{105+\dfrac{6400}{147+\dfrac{20449}{189+\dfrac{50176}{231+\ddots}}}}}}$$
Convergence type $P^-$ with $P=14/3$ and $C=-\G(2/3)^4/(9\pi^2)$, so that
$$\dfrac{\G(2/3)^2}{\G(1/3)}-\dfrac{p(n)}{q(n)}\sim(-1)^{n+1}\dfrac{\G(2/3)^4/(9\pi^2)}{n^{14/3}}$$
$$A=1-(7/3)/n+(7/81)/n^2+(1715/243)/n^3+\cdots$$
Series:
$$\dfrac{\G(2/3)^2}{\G(1/3)}=\dfrac{3}{4}-6\sum_{n\ge0}(-1)^n\dfrac{(4/3)_n^2}{(9n^2+9n+4)(9n^2+27n+22)(5/3)_n^2}$$
Parametric family for $k\ge0$:
\begin{verbatim}
[()->gamma(2/3)^2/gamma(1/3),3*(6*k+7)*(2*n-1),(3*n-1)^2*(3*n+1)^2]
\end{verbatim}
Convergence type $P^-$ with $P=4k+14/3$.
\end{cf}

\smallskip

\begin{cf}\label{5.11.133}{\ }
\begin{verbatim}
[()->gamma(2/3)^2/gamma(1/3),[1/2,6*n-1],[2,9*n^2*(3*n+2)^2]]
\end{verbatim}
$$\dfrac{\G(2/3)^2}{\G(1/3)}=1/2+\dfrac{2}{5+\dfrac{225}{11+\dfrac{2304}{17+\dfrac{9801}{23+\dfrac{28224}{29+\dfrac{65025}{35+\ddots}}}}}}$$
Convergence type $P^-$ with $P=2/3$ and $C=3\G(2/3)^4/(4\pi^2)$, so that
$$\dfrac{\G(2/3)^2}{\G(1/3)}-\dfrac{p(n)}{q(n)}\sim(-1)^n\dfrac{3\G(2/3)^4/4\pi^2}{n^{2/3}}$$
$$A=1-(5/9)/n+(2/9)/n^2+(170/2187)/n^3-(3169/19683)/n^4+\cdots$$
Series:
$$\dfrac{\G(1/3)}{\G(2/3)^2}=2\sum_{n\ge0}(-1)^n\dfrac{(2/3)_n^2}{n!^2}$$
Parametric family for $k\ge0$:
\begin{verbatim}
[()->gamma(2/3)^2/gamma(1/3),(6*k+1)*(6*n-1),9*n^2*(3*n+2)^2]
\end{verbatim}
Convergence type $P^-$ with $P=4k+2/3$.
\end{cf}

\smallskip

\begin{cf}\label{5.11.134}{\ }
\begin{verbatim}
[()->gamma(2/3)^2/gamma(1/3),[1/2,6,18*n^2-30*n+20],[1,-3*n*(3*n-2)*(3*n-1)^2]]
\end{verbatim}
$$\dfrac{\G(2/3)^2}{\G(1/3)}=1/2+\dfrac{1}{6-\dfrac{12}{32-\dfrac{600}{92-\dfrac{4032}{188-\dfrac{14520}{320-\dfrac{38220}{488-\ddots}}}}}}$$
Convergence type $P^+$ with $P=5/3$ and $C=1/(15\G(1/3))$, so that
$$\dfrac{\G(2/3)^2}{\G(1/3)}-\dfrac{p(n)}{q(n)}\sim\dfrac{1/(15\G(1/3))}{n^{5/3}}$$
$$A=1-(5/18)/n-(20/297)/n^2+(145/2187)/n^3+\cdots$$
Series:
$$\dfrac{\G(2/3)^2}{\G(1/3)}=\dfrac{1}{2}+\dfrac{1}{3}\sum_{n\ge0}\dfrac{(1/3)_n}{(3n+2)(n+1)!}$$
Parametric family for $k\ge0$:
\begin{verbatim}
[()->gamma(2/3)^2/gamma(1/3),18*n^2-30*n+20+3*k*(3*k+5),
                                 -3*n*(3*n-2)*(3*n-1)^2]
\end{verbatim}
Convergence type $P^+$ with $P=2k+5/3$.
\end{cf}

\smallskip

\begin{cf}\label{5.11.135}{\ }
\begin{verbatim}
[()->gamma(2/3)^2/gamma(1/3),[5,7*n],[-14,-2*n*(6*n+7)]]
\end{verbatim}
$$\dfrac{\G(2/3)^2}{\G(1/3)}=5-\dfrac{14}{7-\dfrac{26}{14-\dfrac{76}{21-\dfrac{150}{28-\dfrac{248}{35-\dfrac{370}{42-\ddots}}}}}}$$
Convergence type $E$ with $E=4/3$, $P=-7/6$, and
$C=-(3\pi)^{3/2}/(2^{10/3}\G(1/3)^2)$, so that
$$\dfrac{\G(2/3)^2}{\G(1/3)}-\dfrac{p(n)}{q(n)}\sim-\dfrac{(3\pi)^{3/2}/(2^{10/3}\G(1/3)^2)}{(4/3)^nn^{-7/6}}$$
$$A=1-(35/72)/n-(169351/10368)/n^2+\cdots$$
Parametric families for $k\ge0$ and $u\ge0$:
\begin{verbatim}
[()->gamma(2/3)^2/gamma(1/3),7*n+k,-2*n*(6*n+7-6*u)]
[()->gamma(2/3)^2/gamma(1/3),7*n+k,-2*(2*n-1)*(3*n+4-3*u)]
\end{verbatim}
Convergence type $E$ with $E=4/3$ and $P=2k+7u-7/6$ or $P=2k+7u+7/6$
respectively.
\end{cf}

\smallskip

\begin{cf}\label{5.11.136}{\ }
\begin{verbatim}
[()->gamma(2/3)^2/gamma(1/3),[1/2,2*n],[1/2,n*(3*n+1)]]
\end{verbatim}
$$\dfrac{\G(2/3)^2}{\G(1/3)}=1/2+\dfrac{1/2}{2+\dfrac{4}{4+\dfrac{14}{6+\dfrac{30}{8+\dfrac{52}{10+\dfrac{80}{12+\ddots}}}}}}$$
Convergence type $E$ with $E=-3$, $P=1/3$, and $C=2^{5/3}\pi^3/(3\G(1/3)^5)$, so that
$$\dfrac{\G(2/3)^2}{\G(1/3)}-\dfrac{p(n)}{q(n)}\sim(-1)^n\dfrac{2^{5/3}\pi^3/\G(1/3)^5}{3^{n+1}n^{1/3}}$$
$$A=1-(25/72)/n+(1751/10368)/n^2+\cdots$$
Parametric families for $k\ge0$ and $u\ge0$:
\begin{verbatim}
[()->gamma(2/3)^2/gamma(1/3),2*n+4*k+u,n*(3*n-3*u+1)]
[()->gamma(2/3)^2/gamma(1/3),4*n+8*k+2*u+2,(2*n-1)*(6*n-6*u-1)]
\end{verbatim}
Convergence type $E$ with $E=-3$ and $P=2k+u+1/3$ or $P=2k+u+4/3$ respectively.
\end{cf}

\smallskip

\begin{cf}\label{5.11.137}{\ }
\begin{verbatim}
[()->gamma(2/3)^2/gamma(1/3),[3/2,5*n],[-3,-2*n*(2*n+3)]]
\end{verbatim}
$$\dfrac{\G(2/3)^2}{\G(1/3)}=3/2-\dfrac{3}{5-\dfrac{10}{10-\dfrac{28}{15-\dfrac{54}{20-\dfrac{88}{25-\dfrac{130}{30-\ddots}}}}}}$$
Convergence type $E$ with $E=4$, $P=-5/6$, and $C=-3^{4/3}\pi^{5/2}/(4\G(1/3)^4)$, so that
$$\dfrac{\G(2/3)^2}{\G(1/3)}-\dfrac{p(n)}{q(n)}\sim-\dfrac{3^{4/3}\pi^{5/2}/\G(1/3)^4}{2^{2n+2}n^{-5/6}}$$
$$A=1+(275/216)/n-(50615/93312)/n^2+\cdots$$
Parametric family for $k$ and $u\ge0$:
\begin{verbatim}
[()->gamma(2/3)^2/gamma(1/3),5*n-u+3*k,-2*n*(2*n+3-2*u)]
\end{verbatim}
Convergence type $E$ with $E=4$ and $P=u+2k-5/6$.
\end{cf}

\smallskip

\begin{cf}\label{5.11.138}{\ }
\begin{verbatim}
[()->gamma(2/3)^2/gamma(1/3),[2/3,15*n+5],[1/3,-2*(3*n+1)*(6*n-1)]]
\end{verbatim}
$$\dfrac{\G(2/3)^2}{\G(1/3)}=2/3+\dfrac{1/3}{20-\dfrac{40}{35-\dfrac{154}{50-\dfrac{340}{65-\dfrac{598}{80-\dfrac{928}{95-\ddots}}}}}}$$
Convergence type $E$ with $E=4$, $P=5/2$, and $C=\G(2/3)^3/(3\cdot2^{5/3}\pi^{3/2})$, so that
$$\dfrac{\G(2/3)^2}{\G(1/3)}-\dfrac{p(n)}{q(n)}\sim\dfrac{\G(2/3)^3/(3\pi^{3/2})}{2^{2n+5/3}n^{5/2}}$$
$$A=1-(35/8)/n+(5995/384)/n^2+\cdots$$
Series:
$$\dfrac{\G(2/3)^2}{\G(1/3)}=\dfrac{2}{3}+\dfrac{1}{30}\sum_{n\ge0}\dfrac{(4/3)_n}{(n+1)(n+2)(11/6)_n}2^{-2n}$$
Parametric family for $k\ge0$:
\begin{verbatim}
[()->gamma(2/3)^2/gamma(1/3),15*n+5+9*k,-2*(3*n+1)*(6*n-1)]
\end{verbatim}
Convergence type $E$ with $E=4$ and $P=2k+5/2$.

There exist two-variable families as above which I have been too lazy to
write down explicitly.
\end{cf}

\smallskip

\begin{cf}\label{5.11.139}{\ }
\begin{verbatim}
[()->gamma(2/3)^2/gamma(1/3),
[2/3,125,270*n^3-216*n^2+90*n-17],[2,-2*(3*n+1)^3*(6*n-1)^3]]
\end{verbatim}
$$\dfrac{\G(2/3)^2}{\G(1/3)}=2/3+\dfrac{2}{125-\dfrac{16000}{1459-\dfrac{913066}{5599-\dfrac{9826000}{14167-\dfrac{53461798}{28783-\dfrac{199794688}{51067-\ddots}}}}}}$$
Convergence type $E$ with $E=4$, $P=3/2$, and $C=16\pi^{9/2}/(3^{5/2}\G(1/3)^9)$, so that
$$\dfrac{\G(2/3)^2}{\G(1/3)}-\dfrac{p(n)}{q(n)}\sim\dfrac{\pi^{9/2}/(3^{5/2}\G(1/3)^9)}{2^{2n-4}n^{3/2}}$$
$$A=1-(17/8)/n+(1747/384)/n^2+\cdots$$
Series:
$$\dfrac{\G(2/3)^2}{\G(1/3)}=\dfrac{2}{3}+\dfrac{2}{125}\sum_{n\ge0}\dfrac{(4/3)_n^3}{(11/6)_n^3}2^{-2n}$$
\end{cf}

\smallskip

\begin{cf}\label{5.11.140}{\ }
\begin{verbatim}
[()->gamma(2/3)^2/gamma(1/3),[2/3,63*n^2+9*n-2],[4/3,72*n^2*(3*n+1)*(3*n+2)]]
\end{verbatim}
$$\dfrac{\G(2/3)^2}{\G(1/3)}=2/3+\dfrac{4/3}{70+\dfrac{1440}{268+\dfrac{16128}{592+\dfrac{71280}{1042+\dfrac{209664}{1618+\dfrac{489600}{2320+\ddots}}}}}}$$
Convergence type $E$ with $E=-8$, $P=1$, and $C=(\G(2/3)/\G(1/3))^3/6$, so that
$$\dfrac{\G(2/3)^2}{\G(1/3)}-\dfrac{p(n)}{q(n)}\sim(-1)^n\dfrac{(\G(2/3)/\G(1/3))^3/3}{2^{3n+1}n}$$
$$A=1-(10/9)/n+(10/9)/n^2-(1930/2187)/n^3+\cdots$$
Series:
$$\dfrac{\G(1/3)}{\G(2/3)^2}=\dfrac{3}{2}\sum_{n\ge0}(-1)^n\dfrac{(1/3)_n(2/3)_n}{n!^2}2^{-3n}$$
\end{cf}

\smallskip

\begin{cf}\label{5.11.141}{\ }
\begin{verbatim}
[()->gamma(2/3)^2/gamma(1/3),[1/2,6*n-1],[1,n*(3*n+2)]]
\end{verbatim}
$$\dfrac{\G(2/3)^2}{\G(1/3)}=1/2+\dfrac{1}{5+\dfrac{5}{11+\dfrac{16}{17+\dfrac{33}{23+\dfrac{56}{29+\dfrac{85}{35+\ddots}}}}}}$$
Convergence type $E$ with $E=-(2+\sqrt{3})^2$, $P=0$, and
$C=3(\G(2/3)^2/\G(1/3))/(2+\sqrt{3})^{5/3}$, so that
$$\dfrac{\G(2/3)^2}{\G(1/3)}-\dfrac{p(n)}{q(n)}\sim(-1)^n\dfrac{3(\G(2/3)^2/\G(1/3))}{(2+\sqrt{3})^{2n+5/3}}$$
$$A=1-(5d/72)/n+(25d/432+25/3456)/n^2+\cdots$$
\end{cf}

\smallskip

\begin{cf}\label{5.11.142}{\ }
\begin{verbatim}
[()->gamma(2/3)^2/gamma(1/3),
[2/3,338/45,495*n^3-132*n^2-113*n+10],
[91/675,208,3*n*(3*n+1)*(3*n+2)^2*(5*n-3)*(5*n+7)]]
\end{verbatim}
$$\dfrac{\G(2/3)^2}{\G(1/3)}=2/3+\dfrac{91/675}{338/45+\dfrac{208}{3216+\dfrac{319872}{11848+\dfrac{2874960}{29126+\dfrac{14034384}{58020+\dfrac{48829440}{101500+\ddots}}}}}}$$
Convergence type $E$ with $E=-((1+\sqrt{5})/2)^{10}$, $P=0$, and $C=3(\G(2/3)^2/\G(1/3))/((1+\sqrt{5})/2)^{29/3}$, so that
$$\dfrac{\G(2/3)^2}{\G(1/3)}-\dfrac{p(n)}{q(n)}\sim(-1)^n\dfrac{3\G(2/3)^2/\G(1/3)}{((1+\sqrt{5})/2)^{10n+29/3}}$$
$$A=1-(2*d/45)/n+(32*d/675+2/405)/n^2+\cdots$$
\end{cf}

\smallskip

For the next CFs, note that
$$2^{-2/3}\dfrac{\G(2/3)^2}{\G(1/3)}=\dfrac{\G(1/2)\G(2/3)}{\G(1/6)}$$

\smallskip

\begin{cf}\label{5.11.143}{\ }
\begin{verbatim}
[()->2^(-2/3)*gamma(2/3)^2/gamma(1/3),[1,2,1],[-2,9*n^2-1]]
\end{verbatim}
$$2^{-2/3}\dfrac{\G(2/3)^2}{\G(1/3)}=1-\dfrac{2}{2+\dfrac{8}{1+\dfrac{35}{1+\dfrac{80}{1+\dfrac{143}{1+\dfrac{224}{1+\ddots}}}}}}$$
Convergence type $P^-$ with $P=1/3$ and $C=-\G(2/3)/\G(1/3)$, so that
$$2^{-2/3}\dfrac{\G(2/3)^2}{\G(1/3)}-\dfrac{p(n)}{q(n)}\sim(-1)^{n+1}\dfrac{\G(2/3)/\G(1/3)}{n^{1/3}}$$
$$A=1-(1/6)/n-(1/81)/n^2+(14/243)/n^3+(14/6561)/n^4+\cdots$$
Series:
$$2^{-2/3}\dfrac{\G(2/3)^2}{\G(1/3)}=1-\sum_{n\ge0}(-1)^n\dfrac{(4/3)_n}{(5/3)_n}$$
Parametric family for $k\ge0$:
\begin{verbatim}
[()->2^(-2/3)*gamma(2/3)^2/gamma(1/3),6*k+1,9*n^2-1]
\end{verbatim}
Convergence type $P^-$ with $P=2k+1/3$.
\end{cf}

\smallskip

\begin{cf}\label{5.11.144}{\ }
\begin{verbatim}
[()->2^(-2/3)*gamma(2/3)^2/gamma(1/3),[1,3],[-2,3*n*(3*n-2)]]
\end{verbatim}
$$2^{-2/3}\dfrac{\G(2/3)^2}{\G(1/3)}=1-\dfrac{2}{3+\dfrac{3}{3+\dfrac{24}{3+\dfrac{63}{3+\dfrac{120}{3+\dfrac{195}{3+\ddots}}}}}}$$
Convergence type $P^-$ with $P=1$ and $C=-(2^{-2/3}/3)\G(2/3)^2/\G(1/3)$, so that
$$2^{-2/3}\dfrac{\G(2/3)^2}{\G(1/3)}-\dfrac{p(n)}{q(n)}\sim(-1)^{n+1}\dfrac{(2^{-2/3}/3)\G(2/3)^2/\G(1/3)}{n}$$
$$A=1-(1/6)/n-(7/36)/n^2+(23/216)/n^3+(11/48)/n^4+\cdots$$
Series:
\begin{align*}
2^{-2/3}\dfrac{\G(2/3)^2}{\G(1/3)}&=1-\dfrac{1}{2}\sum_{n\ge0}\dfrac{n!(1/6)_n}{(3/2)_n(5/3)_n}\\
2^{2/3}\dfrac{\G(1/3)}{\G(2/3)^2}&=3\sum_{n\ge0}\dfrac{(-1/3)_n(1/2)_n}{n!(7/6)_n}\end{align*}
Parametric family for $k\ge0$:
\begin{verbatim}
[()->2^(-2/3)*gamma(2/3)^2/gamma(1/3),6*k+3,3*n*(3*n-2)]
\end{verbatim}
Convergence type $P^-$ with $P=2k+1$.
\end{cf}

\smallskip

\begin{cf}\label{5.11.145}{\ }
\begin{verbatim}
[()->2^(-2/3)*gamma(2/3)^2/gamma(1/3),
[1/2,10,9*(2*n-1)],[-1,(9*n^2-1)*(36*n^2-1)]]
\end{verbatim}
$$2^{-2/3}\dfrac{\G(2/3)^2}{\G(1/3)}=1/2-\dfrac{1}{10+\dfrac{280}{27+\dfrac{5005}{45+\dfrac{25840}{63+\dfrac{82225}{81+\dfrac{201376}{99+\ddots}}}}}}$$
Convergence type $P^-$ with $P=1$ and $C=-2^{-4/3}(\G(2/3)/\G(1/3))^3$, so that
$$2^{-2/3}\dfrac{\G(2/3)^2}{\G(1/3)}-\dfrac{p(n)}{q(n)}\sim(-1)^{n+1}\dfrac{2^{-4/3}(\G(2/3)/\G(1/3))^3}{n}$$
$$A=1-(1/2)/n-(1/36)/n^2+(7/24)/n^3+(1/216)/n^4+\cdots$$
Series:
$$2^{-2/3}\dfrac{\G(2/3)^2}{\G(1/3)}=-\dfrac{1}{2}+\sum_{n\ge0}(-1)^n\dfrac{(1/3)_n(1/6)_n}{(2/3)_n(5/6)_n}$$
Parametric family for $k\ge0$:
\begin{verbatim}
[()->2^(-2/3)*gamma(2/3)^2/gamma(1/3),
9*(4*k+1)*(2*n-1),(9*n^2-1)*(36*n^2-1)]
\end{verbatim}
Convergence type $P^-$ with $P=4k+1$.
\end{cf}

\smallskip

\begin{cf}\label{5.11.146}{\ }
\begin{verbatim}
[()->2^(-2/3)*gamma(2/3)^2/gamma(1/3),
[3/8,25,27*(2*n-1)],[3/2,(9*n^2-4)*(36*n^2-1)]]
\end{verbatim}
$$2^{-2/3}\dfrac{\G(2/3)^2}{\G(1/3)}=3/8+\dfrac{3/2}{25+\dfrac{175}{81+\dfrac{4576}{135+\dfrac{24871}{189+\dfrac{80500}{243+\dfrac{198679}{297+\ddots}}}}}}$$
Convergence type $P^-$ with $P=3$ and $C=(2^{-4/3}/3)(\G(2/3)/\G(1/3))^3$, so
that
$$2^{-2/3}\dfrac{\G(2/3)^2}{\G(1/3)}-\dfrac{p(n)}{q(n)}\sim(-1)^n\dfrac{(2^{-4/3}/3)(\G(2/3)/\G(1/3))^3}{n^3}$$
$$A=1-(3/2)/n+(5/12)/n^2+(35/24)/n^3+(41/216)/n^4+\cdots$$
Series:
$$2^{-2/3}\dfrac{\G(2/3)^2}{\G(1/3)}=-\dfrac{3}{8}+\dfrac{3}{4}\sum_{n\ge0}(-1)^n\dfrac{(-2/3)_n(1/6)_n}{(5/3)_n(5/6)_n}$$
Parametric family for $k\ge0$:
\begin{verbatim}
[()->2^(-2/3)*gamma(2/3)^2/gamma(1/3),
9*(4*k+3)*(2*n-1),(9*n^2-4)*(36*n^2-1)]
\end{verbatim}
Convergence type $P^-$ with $P=4k+3$.
\end{cf}

\smallskip

\begin{cf}\label{5.11.147}{\ }
\begin{verbatim}
[()->2^(-2/3)*gamma(2/3)^2/gamma(1/3),
[3/8,36*n^2-18*n+7],[3/4,-9*n*(2*n+1)*(3*n+2)*(6*n-1)]]
\end{verbatim}
$$2^{-2/3}\dfrac{\G(2/3)^2}{\G(1/3)}=3/8+\dfrac{3/4}{25-\dfrac{675}{115-\dfrac{7920}{277-\dfrac{35343}{511-\dfrac{104328}{817-\dfrac{244035}{1195-\ddots}}}}}}$$
Convergence type $P^+$ with $P=1$ and $C=(2^{-2/3}/9)\G(2/3)^2/\G(1/3)$,
so that
$$2^{-2/3}\dfrac{\G(2/3)^2}{\G(1/3)}-\dfrac{p(n)}{q(n)}\sim\dfrac{(2^{-2/3}/9)\G(2/3)^2/\G(1/3)}{n}$$
Series:
\begin{align*}2^{-2/3}\dfrac{\G(2/3)^2}{\G(1/3)}&=\dfrac{3}{8}+\dfrac{3}{100}\sum_{n\ge0}\dfrac{n!(3/2)_n}{(8/3)_n(11/6)_n}\\
  2^{2/3}\dfrac{\G(1/3)}{\G(2/3)^2}&=\dfrac{8}{3}\sum_{n\ge0}\dfrac{(2/3)_n(-1/6)_n}{n!(3/2)_n}\end{align*}
Parametric family for $k\ge0$:
\begin{verbatim}
[()->2^(-2/3)*gamma(2/3)^2/gamma(1/3),
36*n^2-18*n+7+18*k*(k+1),-9*n*(2*n+1)*(3*n+2)*(6*n-1)]
\end{verbatim}
Convergence type $P^+$ with $P=2k+1$.
\end{cf}

\smallskip

\begin{cf}\label{5.11.148}{\ }
\begin{verbatim}
[()->2^(-2/3)*gamma(2/3)^2/gamma(1/3),
[1/2,36*n^2-54*n+28],[-1/2,-9*n*(2*n-1)*(3*n-1)*(6*n-1)]]
\end{verbatim}
$$2^{-2/3}\dfrac{\G(2/3)^2}{\G(1/3)}=1/2-\dfrac{1/2}{10-\dfrac{90}{64-\dfrac{2970}{190-\dfrac{18360}{388-\dfrac{63756}{658-\dfrac{164430}{1000-\ddots}}}}}}$$
Convergence type $P^+$ with $P=1$ and $C=-(2^{-5/3}/9)\G(2/3)^2/\G(1/3)$,
so that
$$2^{-2/3}\dfrac{\G(2/3)^2}{\G(1/3)}-\dfrac{p(n)}{q(n)}\sim-\dfrac{(2^{-5/3}/9)\G(2/3)^2/\G(1/3)}{n}$$
Series:
\begin{align*}2^{-2/3}\dfrac{\G(2/3)^2}{\G(1/3)}&=\dfrac{1}{2}-\dfrac{1}{20}\sum_{n\ge0}\dfrac{n!(1/2)_n}{(5/3)_n(11/6)_n}\\
  2^{2/3}\dfrac{\G(1/3)}{\G(2/3)^2}&=2\sum_{n\ge0}\dfrac{(-1/3)_n(-1/6)_n}{n!(1/2)_n}\end{align*}
Parametric family for $k\ge0$:
\begin{verbatim}
[()->2^(-2/3)*gamma(2/3)^2/gamma(1/3),
36*n^2-54*n+28+18*k*(k+1),-9*n*(2*n-1)*(3*n-1)*(6*n-1)]
\end{verbatim}
Convergence type $P^+$ with $P=2k+1$.
\end{cf}

\smallskip

\begin{cf}\label{5.11.149}{\ }
\begin{verbatim}
[()->2^(-2/3)*gamma(2/3)^2/gamma(1/3),
[32/75,121,72*n^2+50],[32/75,-(6*n+1)^2*(6*n+5)^2]]
\end{verbatim}
$$2^{-2/3}\dfrac{\G(2/3)^2}{\G(1/3)}=32/75+\dfrac{32/75}{121-\dfrac{5929}{338-\dfrac{48841}{698-\dfrac{190969}{1202-\dfrac{525625}{1850-\dfrac{1177225}{2642-\ddots}}}}}}$$
Convergence type $P^+$ with $P=7/3$ and $C=(2^{7/3}/63)(\G(2/3)/\G(1/3))^4$,
so that
$$2^{-2/3}\dfrac{\G(2/3)^2}{\G(1/3)}-\dfrac{p(n)}{q(n)}\sim\dfrac{(2^{7/3}/63)(\G(2/3)/\G(1/3))^4}{n^{7/3}}$$
$$A=1-(7/3)/n+(15575/4212)/n^2-(4655/972)/n^3+\cdots$$
Series:
$$2^{-2/3}\dfrac{\G(2/3)^2}{\G(1/3)}=\dfrac{32}{75}+\dfrac{32}{9075}\sum_{n\ge0}\dfrac{(7/6)_n^2}{(17/6)_n^2}$$
Parametric family for $k\ge0$:
\begin{verbatim}
[()->2^(-2/3)*gamma(2/3)^2/gamma(1/3),
72*n^2+50+12*k*(3*k+7),-(6*n+1)^2*(6*n+5)^2]
\end{verbatim}
Convergence type $P^+$ with $P=2k+7/3$.
\end{cf}

\smallskip

\begin{cf}\label{5.11.150}{\ }
\begin{verbatim}
[()->2^(-2/3)*gamma(2/3)^2/gamma(1/3),
[32/75,6*(6*n+1)],[16/75,3*n*(3*n+2)^2*(3*n+4)]]
\end{verbatim}
$$2^{-2/3}\dfrac{\G(2/3)^2}{\G(1/3)}=32/75+\dfrac{16/75}{42+\dfrac{525}{78+\dfrac{3840}{114+\dfrac{14157}{150+\dfrac{37632}{186+\dfrac{82365}{222+\ddots}}}}}}$$
Convergence type $P^-$ with $P=4$ and $C=(25/81)2^{-11/3}\G(2/3)^2/\G(1/3)$,
so that
$$2^{-2/3}\dfrac{\G(2/3)^2}{\G(1/3)}-\dfrac{p(n)}{q(n)}\sim(-1)^n\dfrac{(25/81)2^{-11/3}\G(2/3)^2/\G(1/3)}{n^4}$$
$$A=1-(14/3)/n+(104/9)/n^2-(469/27)/n^3+\cdots$$
Series:
\begin{align*}
  2^{-2/3}\dfrac{\G(2/3)^2}{\G(1/3)}&=\dfrac{32}{75}+\dfrac{32}{21}\sum_{n\ge0}\dfrac{(12n+13)n!(4/3)_n^2(5/3)}{(72n^2+84n+25)(72n^2+228n+181)(3/2)_n(11/6)_n^2(13/6)_n}\\
  2^{2/3}\dfrac{\G(1/3)}{\G(2/3)^2}&=\dfrac{225}{8}\sum_{n\ge0}\dfrac{(12n+7)(1/2)_n(5/6)_n^2(7/6)_n}{(72n^2+12n+1)(72n^2+156n+85)n!(4/3)_n^2(5/3)_n}\end{align*}
Parametric family for $k\ge0$:
\begin{verbatim}
[()->2^(-2/3)*gamma(2/3)^2/gamma(1/3),
6*(k+1)*(6*n+1),3*n*(3*n+2)^2*(3*n+4)]
\end{verbatim}
Convergence type $P^-$ with $P=4k+4$.
\end{cf}
    
\smallskip

\begin{cf}\label{5.11.151}{\ }
\begin{verbatim}
[()->2^(-2/3)*gamma(2/3)^2/gamma(1/3),[3/2,17*n],[-21/2,-12*n*(6*n+7)]]
\end{verbatim}
$$2^{-2/3}\dfrac{\G(2/3)^2}{\G(1/3)}=3/2-\dfrac{21/2}{17-\dfrac{156}{34-\dfrac{456}{51-\dfrac{900}{68-\dfrac{1488}{85-\dfrac{2220}{102-\ddots}}}}}}$$
Convergence type $E$ with $E=9/8$, $P=-17/6$, and $C=-\G(2/3)^4/(2^{2/3}3^{9/2}\pi^{3/2})$, so that
$$2^{-2/3}\dfrac{\G(2/3)^2}{\G(1/3)}-\dfrac{p(n)}{q(n)}\sim-\dfrac{\G(2/3)^4/(2^{2/3}3^{9/2}\pi^{3/2})}{(9/8)^nn^{-17/6}}$$
$$A=1-(2125/72)/n-(48773/384)/n^2+\cdots$$
Parametric family for $k\ge0$:
\begin{verbatim}
[()->2^(-2/3)*gamma(2/3)^2/gamma(1/3),17*n+k,-12*n*(6*n+7)]
\end{verbatim}
Convergence type $E$ with $E=9/8$ and $P=2k-17/6$.

There exist two-variable families as above which I have been too lazy to
write down explicitly.
\end{cf}

\smallskip

\begin{cf}\label{5.11.152}{\ }
\begin{verbatim}
[()->2^(-2/3)*gamma(2/3)^2/gamma(1/3),[2/3,14*n^2-5*n+1],
                          [-4/3,-8*n^2*(2*n+1)*(3*n+2)]]
\end{verbatim}
$$2^{-2/3}\dfrac{\G(2/3)^2}{\G(1/3)}=2/3-\dfrac{4/3}{10-\dfrac{120}{47-\dfrac{1280}{112-\dfrac{5544}{205-\dfrac{16128}{326-\dfrac{37400}{475-\ddots}}}}}}$$
Convergence type $E$ with $E=4/3$, $P=5/6$, and
$C=-2^{2/3}3^{1/2}\pi^{5/2}/\G(1/3)^5$, so that
$$2^{-2/3}\dfrac{\G(2/3)^2}{\G(1/3)}-\dfrac{p(n)}{q(n)}\sim-\dfrac{2^{2/3}3^{1/2}\pi^{5/2}/\G(1/3)^5}{(4/3)^nn^{5/6}}$$
$$A=1-(257/72)/n+(240065/10368)/n^2+\cdots$$
Series:
$$2^{2/3}\dfrac{\G(1/3)}{\G(2/3)^2}=\dfrac{3}{2}\sum_{n\ge0}\dfrac{(1/2)_n(2/3)_n}{n!^2}(3/4)^n$$
\end{cf}

\smallskip

\begin{cf}\label{5.11.153}{\ }
\begin{verbatim}
[()->2^(-2/3)*gamma(2/3)^2/gamma(1/3),
[-1/4,-24,56*n^3-52*n^2-22*n+1],[-7/4,-8*n*(2*n-3)*(2*n+1)^3*(6*n+7)]]
\end{verbatim}
$$2^{-2/3}\dfrac{\G(2/3)^2}{\G(1/3)}=-1/4+\dfrac{7/4}{24-\dfrac{2808}{197-\dfrac{38000}{979-\dfrac{617400}{2665-\dfrac{3615840}{5591-\dfrac{13789160}{10093-\ddots}}}}}}$$
Convergence type $E$ with $E=4/3$, $P=-1/6$, and $C=2^{-2/3}3\sqrt{3\pi}/\G(1/3)^2$, so that
$$2^{-2/3}\dfrac{\G(2/3)^2}{\G(1/3)}-\dfrac{p(n)}{q(n)}\sim\dfrac{2^{-2/3}3\sqrt{3\pi}/\G(1/3)^2}{(4/3)^nn^{-1/6}}$$
$$A=1-(17/72)/n+(15041/10368)/n^2+\cdots$$
Series:
$$2^{-2/3}\dfrac{\G(2/3)^2}{\G(1/3)}=\dfrac{1}{4}\sum_{n\ge0}\dfrac{(2n-1)(7/6)_n}{(2n+1)n!}(3/4)^n$$
\end{cf}

\smallskip

\begin{cf}\label{5.11.154}{\ }
\begin{verbatim}
[()->2^(-2/3)*gamma(2/3)^2/gamma(1/3),[1/2,9*n],[-1/2,-6*n*(3*n+1)]]
\end{verbatim}
$$2^{-2/3}\dfrac{\G(2/3)^2}{\G(1/3)}=1/2-\dfrac{1/2}{9-\dfrac{24}{18-\dfrac{84}{27-\dfrac{180}{36-\dfrac{312}{45-\dfrac{480}{54-\ddots}}}}}}$$
Convergence type $E$ with $E=2$, $P=2$, and $C=-(2^{5/3}/9)\G(2/3)^2/\G(1/3)$,
so that
$$2^{-2/3}\dfrac{\G(2/3)^2}{\G(1/3)}-\dfrac{p(n)}{q(n)}\sim-\dfrac{(1/9)\G(2/3)^2/\G(1/3)}{2^{n-5/3}n^2}$$
$$A=1-5/n+(223/9)/n^2-(11905/81)/n^3+\cdots$$
Parametric family for $k\ge0$:
\begin{verbatim}
[()->2^(-2/3)*gamma(2/3)^2/gamma(1/3),9*n+3*k,-6*n*(3*n+1)]
\end{verbatim}
Convergence type $E$ with $E=2$ and $P=2k+2$.

There exist two-variable families as above which I have been too lazy to
write down explicitly.
\end{cf}

\smallskip

\begin{cf}\label{5.11.155}{\ }
\begin{verbatim}
[()->2^(-2/3)*gamma(2/3)^2/gamma(1/3),[1/2,7*n],[-1/2,4*n*(2*n-1)]]
\end{verbatim}
$$2^{-2/3}\dfrac{\G(2/3)^2}{\G(1/3)}=1/2-\dfrac{1/2}{7+\dfrac{4}{14+\dfrac{24}{21+\dfrac{60}{28+\dfrac{112}{35+\dfrac{180}{42+\ddots}}}}}}$$
Convergence type $E$ with $E=-8$, $P=7/6$, and
$C=-\G(2/3)^4/(2^{1/3}3^{7/3}\pi^{3/2})$, so that
$$2^{-2/3}\dfrac{\G(2/3)^2}{\G(1/3)}-\dfrac{p(n)}{q(n)}\sim(-1)^{n+1}\dfrac{\G(2/3)^4/(3^{7/3}\pi^{3/2})}{2^{3n+1/3}n^{7/6}}$$
$$A=1-(455/648)/n+(313201/839808)/n^2+\cdots$$
Parametric family for $k\ge0$:
\begin{verbatim}
[()->2^(-2/3)*gamma(2/3)^2/gamma(1/3),7*n+9*k,4*n*(2*n-1)]
\end{verbatim}
Convergence type $E$ with $E=-8$ and $P=2k+7/6$.

There exist two-variable families as above which I have been too lazy to
write down explicitly.
\end{cf}

\smallskip

\begin{cf}\label{5.11.156}{\ }
\begin{verbatim}
[()->2^(-2/3)*gamma(2/3)^2/gamma(1/3),[1/2,20,21*n-2],[-3/2,4*(3*n+1)*(6*n-1)]]
\end{verbatim}
$$2^{-2/3}\dfrac{\G(2/3)^2}{\G(1/3)}=1/2-\dfrac{3/2}{20+\dfrac{80}{40+\dfrac{308}{61+\dfrac{680}{82+\dfrac{1196}{103+\dfrac{1856}{124+\ddots}}}}}}$$
Convergence type $E$ with $E=-8$, $P=1/2$, and $C=-2^{1/3}\pi^{3/2}/(3^{3/2}\G(1/3)^3)$, so that
$$2^{-2/3}\dfrac{\G(2/3)^2}{\G(1/3)}-\dfrac{p(n)}{q(n)}\sim(-1)^{n+1}\dfrac{\pi^{3/2}/(3^{3/2}\G(1/3)^3)}{2^{3n-1/3}}$$
$$A=1-(35/72)/n+(1043/3456)/n^2+\cdots$$
Series:
$$2^{-2/3}\dfrac{\G(2/3)^2}{\G(1/3)}=\dfrac{1}{2}-\dfrac{3}{40}\sum_{n\ge0}(-1)^n\dfrac{(4/3)_n}{(11/6)_n}2^{-3n}$$
Parametric family for $k\ge0$:
\begin{verbatim}
[()->2^(-2/3)*gamma(2/3)^2/gamma(1/3),21*n-2+27*k,4*(3*n+1)*(6*n-1)]
\end{verbatim}
Convergence type $E$ with $E=-8$ and $P=2k+1/2$.
\end{cf}

\smallskip

\begin{cf}\label{5.11.157}{\ }
\begin{verbatim}
[()->2^(-2/3)*gamma(2/3)^2/gamma(1/3),[3/4,10*n],[-21/8,-3*n*(3*n+7)]]
\end{verbatim}
$$2^{-2/3}\dfrac{\G(2/3)^2}{\G(1/3)}=3/4-\dfrac{21/8}{10-\dfrac{30}{20-\dfrac{78}{30-\dfrac{144}{40-\dfrac{228}{50-\dfrac{330}{60-\ddots}}}}}}$$
Convergence type $E$ with $E=9$, $P=-5/3$, and $C=-2^{29/3}\pi^3/(3^7\G(1/3)^5)$, so that
$$2^{-2/3}\dfrac{\G(2/3)^2}{\G(1/3)}-\dfrac{p(n)}{q(n)}\sim-\dfrac{2^{29/3}\pi^3/\G(1/3)^5}{3^{2n+7}n^{-5/3}}$$
$$A=1+(475/144)/n+(74195/41472)/n^2+\cdots$$
Parametric families for $k\ge0$:
\begin{verbatim}
[()->2^(-2/3)*gamma(2/3)^2/gamma(1/3),10*n+8*k,-3*n*(3*n+7)]
[()->2^(-2/3)*gamma(2/3)^2/gamma(1/3),20*n+16*k,-3*(2*n+5)*(6*n-1)]
\end{verbatim}
Convergence type $E$ with $E=9$ and $P=2k-5/3$.

There exist two-variable families as above which I have been too lazy to
write down explicitly.
\end{cf}

\smallskip

\begin{cf}\label{5.11.158}{\ }
\begin{verbatim}
[()->2^(-2/3)*gamma(2/3)^2/gamma(1/3),[1/2,15,20*n-6],[-1,-3*(2*n-1)*(6*n-1)]]
\end{verbatim}
$$2^{-2/3}\dfrac{\G(2/3)^2}{\G(1/3)}=\dfrac{1}{2}-\dfrac{1}{15-\dfrac{15}{34-\dfrac{99}{54-\dfrac{255}{74-\dfrac{483}{94-\dfrac{783}{114-\ddots}}}}}}$$
Convergence type $E$ with $E=9$, $P=4/3$, and $C=-2^{-2/3}\G(2/3)/(8\G(1/3))$, so that
$$2^{-2/3}\dfrac{\G(2/3)^2}{\G(1/3)}-\dfrac{p(n)}{q(n)}\sim-\dfrac{2^{-2/3}\G(2/3)/(8\G(1/3))}{3^{2n}n^{4/3}}$$
$$A=1-(19/18)/n+1603/1296/n^2+\cdots$$
Series:
$$2^{-2/3}\dfrac{\G(2/3)^2}{\G(1/3)}=-\dfrac{1}{2}+\sum_{n\ge0}\dfrac{(-1/2)_n}{(5/6)_n}3^{-2n}$$
Parametric family for $k\ge0$:
\begin{verbatim}
[()->2^(-2/3)*gamma(2/3)^2/gamma(1/3),20*n-6+16*k,-3*(2*n-1)*(6*n-1)]
\end{verbatim}
Convergence type $E$ with $E=9$ and $P=2k+4/3$.
\end{cf}

\smallskip

\begin{cf}\label{5.11.159}{\ }
\begin{verbatim}
[()->2^(-2/3)*gamma(2/3)^2/gamma(1/3),[1/2,18*n-3],[-1,-3*n*(3*n+2)]]
\end{verbatim}
$$2^{-2/3}\dfrac{\G(2/3)^2}{\G(1/3)}=1/2-\dfrac{1}{15-\dfrac{15}{33-\dfrac{48}{51-\dfrac{99}{69-\dfrac{168}{87-\dfrac{255}{105-\ddots}}}}}}$$
Convergence type $E$ with $E=(1+\sqrt{2})^4$, $P=0$, and
$C=-3\G(2/3)^3/(2^{2/3}\pi(1+\sqrt{2})^{10/3})$, so that
$$2^{-2/3}\dfrac{\G(2/3)^2}{\G(1/3)}-\dfrac{p(n)}{q(n)}\sim-\dfrac{3\G(2/3)^3/(2^{2/3}\pi)}{(1+\sqrt{2})^{4n+10/3}}$$
$$A=1-(5d/48)/n+(25d/288+25/2304)/n^2+\cdots$$
\end{cf}

\smallskip

\begin{cf}\label{5.11.160}{\ }
\begin{verbatim}
[()->2^(-2/3)*gamma(2/3)^2/gamma(1/3),[1,18*(6*n^2-8*n+1)],
                                      [10,-3*n*(3*n-1)*(6*n-7)*(6*n+5)]]
\end{verbatim}
$$2^{-2/3}\dfrac{\G(2/3)^2}{\G(1/3)}=1+\dfrac{10}{-18+\dfrac{66}{162-\dfrac{2550}{558-\dfrac{18216}{1170-\dfrac{65076}{1998-\dfrac{169050}{3042-\ddots}}}}}}$$
Convergence type $E$ with $E=(1+\sqrt{2})^4$, $P=0$, and $C=-2^{1/3}3^{1/2}(\G(2/3)^2/\G(1/3))/(1+\sqrt{2})^{4/3}$, so that
$$2^{-2/3}\dfrac{\G(2/3)^2}{\G(1/3)}-\dfrac{p(n)}{q(n)}\sim-\dfrac{2^{1/3}3^{1/2}\G(2/3)^2/\G(1/3)}{(1+\sqrt{2})^{4n+4/3}}$$
$$A=1-(d/48)/n+(d/144+1/2304)/n^2+\cdots$$
\end{cf}

\smallskip

For the next CFs, note that
$$2^{-1/3}\dfrac{\G(2/3)^2}{\G(1/3)}=\dfrac{\G(1/6)\G(2/3)^3}{2\G(1/2)\G(1/3)^2}$$

\smallskip

\begin{cf}\label{5.11.161}{\ }
\begin{verbatim}
[()->2^(-1/3)*gamma(2/3)^2/gamma(1/3),
[1/2,4*(3*n-1)],[1/2,9*n^2*(3*n+1)^2]]
\end{verbatim}
$$2^{-1/3}\dfrac{\G(2/3)^2}{\G(1/3)}=1/2+\dfrac{1/2}{8+\dfrac{144}{20+\dfrac{1764}{32+\dfrac{8100}{44+\dfrac{24336}{56+\dfrac{57600}{68+\ddots}}}}}}$$
Convergence type $P^-$ with $P=4/3$ and $C=2^{-2/3}(\G(2/3)/\G(1/3))^4$,
so that
$$2^{-1/3}\dfrac{\G(2/3)^2}{\G(1/3)}-\dfrac{p(n)}{q(n)}\sim(-1)^n\dfrac{2^{-2/3}(\G(2/3)/\G(1/3))^4}{n^{4/3}}$$
$$A=1-(8/9)/n+(22/81)/n^2+(920/2187)/n^3+\cdots$$
Series:
$$2^{1/3}\dfrac{\G(1/3)}{\G(2/3)^2}=2\sum_{n\ge0}(-1)^n\dfrac{(1/3)_n^2}{n!^2}$$
Parametric family for $k\ge0$ and $u\ge0$:
\begin{verbatim}
[()->2^(-1/3)*gamma(2/3)^2/gamma(1/3),
4*(3*k+1)*(3*n-1+3*u),9*n^2*(3*n+1+6*u)^2]
\end{verbatim}
Convergence type $P^-$ with $P=4k+4/3$.
\end{cf}

\smallskip

\begin{cf}\label{5.11.162}{\ }
\begin{verbatim}
[()->2^(-1/3)*gamma(2/3)^2/gamma(1/3),
[3/5,10*(3*n-1)],[-6/5,9*n*(n+1)*(3*n-2)*(3*n+1)]]
\end{verbatim}
$$2^{-1/3}\dfrac{\G(2/3)^2}{\G(1/3)}=3/5-\dfrac{6/5}{20+\dfrac{72}{50+\dfrac{1512}{80+\dfrac{7560}{110+\dfrac{23400}{140+\dfrac{56160}{170+\ddots}}}}}}$$
Convergence type $P^-$ with $P=10/3$ and $C=-(5/9)2^{-2/3}(\G(2/3)/\G(1/3))^4$,
so that
$$2^{-1/3}\dfrac{\G(2/3)^2}{\G(1/3)}-\dfrac{p(n)}{q(n)}\sim(-1)^{n+1}\dfrac{(5/9)2^{-2/3}(\G(2/3)/\G(1/3))^4}{n^{10/3}}$$
$$A=1-(20/9)/n+(335/162)/n^2+(560/2187)/n^3+\cdots$$
Series:
$$2^{1/3}\dfrac{\G(1/3)}{\G(2/3)^2}=\dfrac{5}{3}\sum_{n\ge0}(-1)^n\dfrac{(-2/3)_n(1/3)_n}{n!(n+1)!}$$
Parametric family for $k\ge0$ and $u\ge0$:
\begin{verbatim}
[()->2^(-1/3)*gamma(2/3)^2/gamma(1/3),
2*(6*k+5)*(3*n-1+3*u),9*n*(n+1)*(3*n-2+6*u)*(3*n+1+6*u)]
\end{verbatim}
Convergence type $P^-$ with $P=4k+10/3$.
\end{cf}

\smallskip

\begin{cf}\label{5.11.163}{\ }
\begin{verbatim}
[()->2^(-1/3)*gamma(2/3)^2/gamma(1/3),
[1/2,15,36*n^2-48*n+28],[1/2,-3*n*(3*n-1)*(6*n-1)^2]]
\end{verbatim}
$$2^{-1/3}\dfrac{\G(2/3)^2}{\G(1/3)}=1/2+\dfrac{1/2}{15-\dfrac{150}{76-\dfrac{3630}{208-\dfrac{20808}{412-\dfrac{69828}{688-\dfrac{176610}{1036-\ddots}}}}}}$$
Convergence type $P^+$ with $P=4/3$ and $C=1/(48\G(2/3))$, so that
$$2^{-1/3}\dfrac{\G(2/3)^2}{\G(1/3)}-\dfrac{p(n)}{q(n)}\sim\dfrac{1/(48\G(2/3))}{n^{4/3}}$$
$$A=1-(4/9)/n+(49/810)/n^2+(133/2187)/n^3+\cdots$$
Series:
$$2^{-1/3}\dfrac{\G(2/3)^2}{\G(1/3)}=\dfrac{7}{2}\sum_{n\ge0}\dfrac{(-1/3)_n}{(6n+5)n!}$$
Parametric family for $k\ge0$:
\begin{verbatim}
[()->2^(-1/3)*gamma(2/3)^2/gamma(1/3),
36*n^2-48*n+28+6*k*(3*k+4),-3*n*(3*n-1)*(6*n-1)^2]
\end{verbatim}
Convergence type $P^+$ with $P=2k+4/3$
\end{cf}

\smallskip

\begin{cf}\label{5.11.164}{\ }
\begin{verbatim}
[()->2^(-1/3)*gamma(2/3)^2/gamma(1/3),
[3/5,33,36*n^2-12*n+7],[-6/5,-3*n*(3*n+2)*(6*n-1)*(6*n+5)]]
\end{verbatim}
$$2^{-1/3}\dfrac{\G(2/3)^2}{\G(1/3)}=3/5-\dfrac{6/5}{33-\dfrac{825}{127-\dfrac{8976}{295-\dfrac{38709}{535-\dfrac{112056}{847-\dfrac{258825}{1231-\ddots}}}}}}$$
Convergence type $P^+$ with $P=4/3$ and $C=-1/(16\G(2/3))$, so that
$$2^{-1/3}\dfrac{\G(2/3)^2}{\G(1/3)}-\dfrac{p(n)}{q(n)}\sim-\dfrac{1/(16\G(2/3))}{n^{4/3}}$$
$$A=1-(10/9)/n+(847/810)/n^2-(3955/4374)/n^3+\cdots$$
Series:
$$2^{-1/3}\dfrac{\G(2/3)^2}{\G(1/3)}=-3\sum_{n\ge0}\dfrac{(2/3)_n}{(6n-1)(6n+5)n!}$$
Parametric family for $k\ge0$:
\begin{verbatim}
[()->2^(-1/3)*gamma(2/3)^2/gamma(1/3),
36*n^2-12*n+7+6*k*(3*k+4),-3*n*(3*n+2)*(6*n-1)*(6*n+5)]
\end{verbatim}
Convergence type $P^+$ with $P=2k+4/3$.
\end{cf}

\smallskip

\begin{cf}\label{5.11.165}{\ }
\begin{verbatim}
[()->2^(-1/3)*gamma(2/3)^2/gamma(1/3),[2/3,14*n^2-7*n+2],
                          [-2/3,-8*n^2*(2*n+1)*(3*n+1)]]
\end{verbatim}
$$2^{-1/3}\dfrac{\G(2/3)^2}{\G(1/3)}=2/3-\dfrac{2/3}{9-\dfrac{96}{44-\dfrac{1120}{107-\dfrac{5040}{198-\dfrac{14976}{317-\dfrac{35200}{464-\ddots}}}}}}$$
Convergence type $E$ with $E=4/3$, $P=7/6$, and $C=-2^{7/3}\pi^{7/2}/\G(1/3)^7$,
so that
$$2^{-1/3}\dfrac{\G(2/3)^2}{\G(1/3)}-\dfrac{p(n)}{q(n)}\sim-\dfrac{2^{7/3}\pi^{7/2}/\G(1/3)^7}{(4/3)^nn^{7/6}}$$
$$A=1-(353/72)/n+(129451/3456)/n^2+\cdots$$
Series:
$$2^{1/3}\dfrac{\G(1/3)}{\G(2/3)^2}=\dfrac{3}{2}\sum_{n\ge0}\dfrac{(1/2)_n(1/3)_n}{n!^2}(3/4)^n$$
\end{cf}

\medskip

For the next CFs, note that $\G(1/3)/\G(2/3)=\G(1/3)^2/(2\pi/\sqrt{3})$ and
that
$$\left(\dfrac{\G(1/3)}{\G(2/3)}\right)^3=\dfrac{1}{4}\left(\dfrac{\G(1/6)}{\G(1/2)}\right)^3$$

In addition, note that $(\G(1/3)/\G(2/3))^3=\CS(-3)$, see Section \ref{sec:CS}
for CFs related to general Chowla--Selberg gamma quotients.

\smallskip

\begin{cf}\label{5.11.166}{\ }
\begin{verbatim}
[()->(gamma(1/3)/gamma(2/3))^3,[9/2,17,3*(2*n-1)*(18*n^2-18*n+7)],
                                 [36,-4*(9*n^2-4)*(9*n^2-1)^2]]
\end{verbatim}
$$\left(\dfrac{\G(1/3)}{\G(2/3)}\right)^3=9/2+\dfrac{36}{17-\dfrac{1280}{387-\dfrac{156800}{1725-\dfrac{1971200}{4683-\dfrac{11451440}{9909-\ddots}}}}}$$
Convergence type $P^+$ with $P=4/3$ and
$C=3^{7/2}\G(1/3)^{10}/(2^{31/3}\pi^5)$, so that
$$\left(\dfrac{\G(1/3)}{\G(2/3)}\right)^3-\dfrac{p(n)}{q(n)}\sim\dfrac{3^{7/2}\G(1/3)^{10}/(2^{31/3}\pi^5)}{n^{4/3}}\;.$$
$$A=1-(2/3)/n+(113/324)/n^2-(145/972)/n^3+(1554539/26873856)/n^4+\cdots$$
Parametric family for $k\ge0$:
\begin{verbatim}
[()->(gamma(1/3)/gamma(2/3))^3,
3*(2*n-1)*(18*n^2-18*n+7+12*k*(3*k+2)),-4*(9*n^2-4)*(9*n^2-1)^2]
\end{verbatim}
Convergence type $P^+$ with $P=4k+4/3$.
\end{cf}

\smallskip

\begin{cf}\label{5.11.167}{\ }
\begin{verbatim}
[()->(gamma(1/3)/gamma(2/3))^3,[6,16,9*(2*n-1)*(3*n^2-3*n+2)],
                              [24,-(9*n^2-4)*(9*n^2-1)^2]]
\end{verbatim}
$$\left(\dfrac{\G(1/3)}{\G(2/3)}\right)^3=6+\dfrac{24}{16-\dfrac{320}{216-\dfrac{39200}{900-\dfrac{492800}{2394-\dfrac{2862860}{5022-\dfrac{11088896}{9108-\ddots}}}}}}$$
Convergence type $P^+$ with $P=2$ and $C=2(\G(1/3)/\G(2/3))^3/27$,
so that
$$\left(\dfrac{\G(1/3)}{\G(2/3)}\right)^3-\dfrac{p(n)}{q(n)}\sim\dfrac{2(\G(1/3)/\G(2/3))^3/27}{n^2}\;.$$
$$A=1-1/n+(17/27)/n^2-(7/27)/n^3+(137/2187)/n^4+\cdots$$
Series:
\begin{align*}\left(\dfrac{\G(1/3)}{\G(2/3)}\right)^3&=6+12\sum_{n\ge0}\dfrac{(3n+2)(2/3)_n^3}{(3n+4)^2(4/3)_n^3}\\
\left(\dfrac{\G(2/3)}{\G(1/3)}\right)^3&=-\dfrac{1}{6}-\dfrac{4}{3}\sum_{n\ge0}\dfrac{(1/3)_n^3}{(3n-2)(3n+2)(2/3)_n^3}\end{align*}
Parametric family for $k\ge0$:
\begin{verbatim}
[()->(gamma(1/3)/gamma(2/3))^3,9*(2*n-1)*(3*n^2-3*n+2+6*k*(k+1)),
                               -(9*n^2-4)*(9*n^2-1)^2]
\end{verbatim}
Convergence type $P^+$ with $P=4k+2$.
\end{cf}

\smallskip

\begin{cf}\label{5.11.168}{\ }
\begin{verbatim}
[()->(gamma(1/3)/gamma(2/3))^3,[8,4*(27*n^4+9*n^2-4)],
                               [-24,-(2*n-1)*(2*n+3)*(3*n+1)^3*(3*n+2)^3]]
\end{verbatim}
$$\left(\dfrac{\G(1/3)}{\G(2/3)}\right)^3=8-\dfrac{24}{128-\dfrac{40000}{1856-\dfrac{3687936}{9056-\dfrac{59895000}{28208-\dfrac{464199736}{68384-\dfrac{2354466816}{141248-\ddots}}}}}}$$
Convergence type $P^+$ with $P=2$ and $C=-(1/27)(\G(1/3)/\G(2/3))^3$, so that
$$\left(\dfrac{\G(1/3)}{\G(2/3)}\right)^3-\dfrac{p(n)}{q(n)}\sim-\dfrac{(1/27)(\G(1/3)/\G(2/3))^3}{n^2}$$
$$A=1-2/n+\cdots$$
Series:
\begin{align*}\left(\dfrac{\G(1/3)}{\G(2/3)}\right)^3&=8-\dfrac{1}{8}\sum_{n\ge0}\dfrac{(2n+3)(5/3)_n^3}{(n+1)(n+2)(7/3)_n^3}\\
  \left(\dfrac{\G(2/3)}{\G(1/3)}\right)^3&=\dfrac{1}{8}\sum_{n\ge0}\dfrac{(2n+1)(1/3)_n^3}{(5/3)_n^3}\end{align*}
\end{cf}

\smallskip

\begin{cf}\label{5.11.169}{\ }
\begin{verbatim}
[()->(gamma(1/3)/gamma(2/3))^3,[8,64,9*(2*n-1)*(3*n^2-3*n+7)],
                               [-16,-(9*n^2-1)^3]]
\end{verbatim}
$$\left(\dfrac{\G(1/3)}{\G(2/3)}\right)^3=8-\dfrac{16}{64-\dfrac{512}{351-\dfrac{42875}{1125-\dfrac{512000}{2709-\dfrac{2924207}{5427-\dfrac{11239424}{9603-\ddots}}}}}}$$
Convergence type $P^+$ with $P=4$ and $C=-(4/729)(\G(1/3)/\G(2/3))^3$, so that
$$\left(\dfrac{\G(1/3)}{\G(2/3)}\right)^3-\dfrac{p(n)}{q(n)}\sim-\dfrac{(4/729)(\G(1/3)/\G(2/3))^3}{n^4}$$
$$A=1-2/n+(155/81)/n^2-(20/27)/n^3+\cdots$$
Series:
\begin{align*}\left(\dfrac{\G(1/3)}{\G(2/3)}\right)^3&=8-\dfrac{1}{4}\sum_{n\ge0}\dfrac{(2/3)_n^3}{(7/3)_n^3}\\
  \left(\dfrac{\G(2/3)}{\G(1/3)}\right)^3&=-\dfrac{1}{8}+16\sum_{n\ge0}\dfrac{(1/3)_n^3}{(27n^2-27n+8)(27n^2+27n+8)(2/3)_n^3}\end{align*}
Parametric family for $k\ge0$:
\begin{verbatim}
[()->(gamma(1/3)/gamma(2/3))^3,
9*(2*n-1)*(3*n^2-3*n+7+6*k*(k+2)),-(9*n^2-1)^3]
\end{verbatim}
Convergence type $P^+$ with $P=4k+4$.
\end{cf}

\smallskip

\begin{cf}\label{5.11.170}{\ }
\begin{verbatim}
[()->(gamma(1/3)/gamma(2/3))^3,[6,36,(6*n-1)*(9*n^2-3*n+2)],
                               [24,-9*n^2*(3*n+1)^2*(3*n+2)^2]]
\end{verbatim}
$$\left(\dfrac{\G(1/3)}{\G(2/3)}\right)^3=6+\dfrac{24}{36-\dfrac{3600}{352-\dfrac{112896}{1258-\dfrac{980100}{3082-\dfrac{4769856}{6148-\dfrac{16646400}{10780-\ddots}}}}}}$$
Convergence type $P^+$ with $P=2/3$ and $C=9\G(1/3)^2/(4\pi^2)$, so that
$$\left(\dfrac{\G(1/3)}{\G(2/3)}\right)^3-\dfrac{p(n)}{q(n)}\sim\dfrac{9\G(1/3)^2/(4\pi^2)}{n^{2/3}}$$
$$A=1-(5/9)/n+(1/3)/n^2-(370/2187)/n^3+\cdots$$
Series:
$$\left(\dfrac{\G(1/3)}{\G(2/3)}\right)^3=6+\dfrac{8}{3}\sum_{n\ge0}\dfrac{(5/3)_n^2}{(3n+4)(n+1)!^2}$$
Parametric family for $k\ge0$:
\begin{verbatim}
[()->(gamma(1/3)/gamma(2/3))^3,(6*n-1)*(9*n^2-3*n+2+6*k*(3*k+1)),
                               -9*n^2*(3*n+1)^2*(3*n+2)^2]
\end{verbatim}
Convergence type $P^+$ with $P=4k+2/3$.
\end{cf}

\smallskip

\begin{cf}\label{5.11.171}{\ }
\begin{verbatim}
[()->(gamma(1/3)/gamma(2/3))^3,[9/2,49,4*(3*n-2)*(36*n^2-48*n+31)],
                    [144,-27*(2*n-1)^2*(2*n+1)*(6*n-5)*(6*n+1)^2]]
\end{verbatim}
$$\left(\dfrac{\G(1/3)}{\G(2/3)}\right)^3=9/2+\dfrac{144}{49-\dfrac{3969}{1264-\dfrac{1437345}{5908-\dfrac{22174425}{16600-\dfrac{141395625}{35932-\ddots}}}}}$$
Convergence type $P^+$ with $P=2$ and $C=2(\G(1/3)/\G(2/3))^3/27$, so that
$$\left(\dfrac{\G(1/3)}{\G(2/3)}\right)^3-\dfrac{p(n)}{q(n)}\sim\dfrac{2(\G(1/3)/\G(2/3))^3/27}{n^2}\;.$$
$$A=1-(2/3)/n+(23/108)/n^2+\cdots$$
Series:
\begin{align*}\left(\dfrac{\G(1/3)}{\G(2/3)}\right)^3&=\dfrac{9}{2}+144\sum_{n\ge0}\dfrac{(2n+1)(1/2)_n^3}{(6n+7)^2(7/6)_n^3}\\
\left(\dfrac{\G(2/3)}{\G(1/3)}\right)^3&=-\dfrac{6}{5}-\dfrac{64}{9}\sum_{n\ge0}\dfrac{(1/6)_n^3}{(2n+1)(6n-5)(1/2)_n^3}\end{align*}
Parametric family for $k\ge0$:
\begin{verbatim}
[()->(gamma(1/3)/gamma(2/3))^3,4*(3*n-2)*(36*n^2-48*n+31+72*k*(k+1)),
                             -27*(2*n-1)^2*(2*n+1)*(6*n-5)*(6*n+1)^2]
\end{verbatim}
Convergence type $P^+$ with $P=4k+2$.
\end{cf}

\smallskip

\begin{cf}\label{5.11.172}{\ }
\begin{verbatim}
[()->(gamma(1/3)/gamma(2/3))^3,[6,(6*n-5)*(9*n^2-15*n+20)],
                               [24,-9*n^2*(3*n-1)^2*(3*n-2)^2]]
\end{verbatim}
$$\left(\dfrac{\G(1/3)}{\G(2/3)}\right)^3=6+\dfrac{24}{14-\dfrac{36}{182-\dfrac{14400}{728-\dfrac{254016}{1976-\dfrac{1742400}{4250-\dfrac{7452900}{7874-\ddots}}}}}}$$
Convergence type $P^+$ with $P=10/3$ and $C=\G(1/3)^{10}/(320\pi^6)$, so that
$$\left(\dfrac{\G(1/3)}{\G(2/3)}\right)^3-\dfrac{p(n)}{q(n)}\sim\dfrac{\G(1/3)^{10}/(320\pi^6)}{n^{10/3}}$$
$$A=1-(5/9)/n-(20/81)/n^2+(740/2187)/n^3+\cdots$$
Series:
$$\left(\dfrac{\G(2/3)}{\G(1/3)}\right)^3=-\dfrac{1}{6}\sum_{n\ge0}\dfrac{(-2/3)_n^2}{(3n-1)n!^2}$$
Parametric family for $k\ge0$:
\begin{verbatim}
[()->(gamma(1/3)/gamma(2/3))^3,(6*n-5)*(9*n^2-15*n+20+6*k*(3*k+5)),
                               -9*n^2*(3*n-1)^2*(3*n-2)^2]
\end{verbatim}
Convergence type $P^+$ with $P=4k+10/3$.
\end{cf}

\smallskip

\begin{cf}\label{5.11.173}{\ }
\begin{verbatim}
[()->(gamma(1/3)/gamma(2/3))^3,[15,14*n+7],[-96,-8*(2*n+1)*(3*n+4)]]
\end{verbatim}
$$\left(\dfrac{\G(1/3)}{\G(2/3)}\right)^3=15-\dfrac{96}{21-\dfrac{168}{35-\dfrac{400}{49-\dfrac{728}{63-\dfrac{1152}{77-\dfrac{1672}{91-\ddots}}}}}}$$
Convergence type $E$ with $E=4/3$, $P=7/6$, and $C=-2^{-19/6}3^{11/4}(\G(1/3)/\G(2/3))^{5/2}$, so that
$$\left(\dfrac{\G(1/3)}{\G(2/3)}\right)^3-\dfrac{p(n)}{q(n)}\sim-\dfrac{2^{-19/6}3^{11/4}(\G(1/3)/\G(2/3))^{5/2}}{(4/3)^nn^{7/6}}$$
$$A=1-(329/72)/n+(123515/3456)/n^2+\cdots$$
Parametric family for $k\ge0$:
\begin{verbatim}
[()->(gamma(1/3)/gamma(2/3))^3,14*n+2*k+1,-8*(2*n+1)*(3*n+1)]
\end{verbatim}
Convergence type $E$ with $E=4/3$ and $P=2k+13/6$.
\end{cf}

\smallskip

\begin{cf}\label{5.11.174}{\ }
\begin{verbatim}
[()->(gamma(1/3)/gamma(2/3))^3,[15,16,14*n-3],[-96,-8*(2*n+1)*(3*n-1)]]
\end{verbatim}
$$\left(\dfrac{\G(1/3)}{\G(2/3)}\right)^3=15-\dfrac{96}{16-\dfrac{48}{25-\dfrac{200}{39-\dfrac{448}{53-\dfrac{792}{67-\dfrac{1232}{81-\ddots}}}}}}$$
Convergence type $E$ with $E=4/3$, $P=17/6$, and $C=-2^{13/6}3^{3/4}(\G(1/3)/\G(2/3))^{7/2}$, so that
$$\left(\dfrac{\G(1/3)}{\G(2/3)}\right)^3-\dfrac{p(n)}{q(n)}\sim-\dfrac{2^{13/6}3^{3/4}(\G(1/3)/\G(2/3))^{7/2}}{(4/3)^nn^{17/6}}$$
$$A=1-(1169/72)/n+(2558945/10368)/n^2+\cdots$$
Parametric family for $k\ge0$:
\begin{verbatim}
[()->(gamma(1/3)/gamma(2/3))^3,14*n+2*k-3,-8*(2*n+1)*(3*n-1)]
\end{verbatim}
Convergence type $E$ with $E=4/3$ and $P=2k+17/6$.
\end{cf}

\smallskip

\begin{cf}\label{5.11.175}{\ }
\begin{verbatim}
[()->(gamma(1/3)/gamma(2/3))^3,[15,49,14*n+3],[-480,32*(3*n+2)*(6*n+1)]]
\end{verbatim}
$$\left(\dfrac{\G(1/3)}{\G(2/3)}\right)^3=15-\dfrac{480}{49+\dfrac{1120}{31+\dfrac{3328}{45+\dfrac{6688}{59+\dfrac{11200}{73+\dfrac{16864}{87+\ddots}}}}}}$$
Convergence type $E$ with $E=-16/9$, $P=1/6$, and $C=-2^{-7/2}3^{11/4}5^{-1/3}(\G(1/3)/\G(2/3))^{5/2}$, so that
$$\left(\dfrac{\G(1/3)}{\G(2/3)}\right)^3-\dfrac{p(n)}{q(n)}\sim(-1)^{n+1}\dfrac{2^{-7/2}3^{11/4}5^{-1/3}(\G(1/3)/\G(2/3))^{5/2}}{(4/3)^{2n}n^{1/6}}$$
$$A=1-(373/1800)/n+(236771/2160000)/n^2+\cdots$$
Parametric family for $k\ge0$:
\begin{verbatim}
[()->(gamma(1/3)/gamma(2/3))^3,14*n+3+50*k,32*(3*n+2)*(6*n+1)]
\end{verbatim}
Convergence type $E$ with $E=-16/9$ and $P=2k+1/6$.
\end{cf}

\smallskip

\begin{cf}\label{5.11.176}{\ }
\begin{verbatim}
[()->(gamma(1/3)/gamma(2/3))^3,[15,64,14*n+33],[-480,32*(3*n-2)*(6*n-1)]]
\end{verbatim}
$$\left(\dfrac{\G(1/3)}{\G(2/3)}\right)^3=15-\dfrac{480}{64+\dfrac{160}{61+\dfrac{1408}{75+\dfrac{3808}{89+\dfrac{7360}{103+\dfrac{12064}{117+\ddots}}}}}}$$
Convergence type $E$ with $E=-16/9$, $P=11/6$, and
$C=-2^{7/2}3^{3/4}5^{-11/3}(\G(1/3)/\G(2/3))^{7/2}$, so that
$$\left(\dfrac{\G(1/3)}{\G(2/3)}\right)^3-\dfrac{p(n)}{q(n)}\sim(-1)^{n+1}\dfrac{2^{7/2}3^{3/4}5^{-11/3}(\G(1/3)/\G(2/3))^{7/2}}{(4/3)^{2n}n^{11/6}}$$
$$A=1-(793/1800)/n-(2071087/6580000)/n^2+\cdots$$
Parametric family for $k\ge0$:
\begin{verbatim}
[()->(gamma(1/3)/gamma(2/3))^3,14*n+33+50*k,32*(3*n-2)*(6*n-1)]
\end{verbatim}
Convergence type $E$ with $E=-16/9$ and $P=2k+11/6$.
\end{cf}

\smallskip

\begin{cf}\label{5.11.177}{\ }
\begin{verbatim}
[()->(gamma(1/3)/gamma(2/3))^3,[9/2,4*n+3],[36,2*(2*n+3)*(3*n+4)]]
\end{verbatim}
$$\left(\dfrac{\G(1/3)}{\G(2/3)}\right)^3=9/2+\dfrac{36}{7+\dfrac{70}{11+\dfrac{140}{15+\dfrac{234}{19+\dfrac{352}{23+\dfrac{494}{27+\ddots}}}}}}$$
Convergence type $E$ with $E=-3$, $P=-1/6$ and
$C=2^{-1/6}3^{-3/4}(\G(1/3)/\G(2/3))^{7/2}$, so that
$$\left(\dfrac{\G(1/3)}{\G(2/3)}\right)^3-\dfrac{p(n)}{q(n)}\sim(-1)^n\dfrac{2^{-1/6}(\G(1/3)/\G(2/3))^{7/2}}{3^{n+3/4}n^{-1/6}}$$
$$A=1+(7/36)/n-(235/10368)/n^2+\cdots$$
Parametric family for $k\ge0$:
\begin{verbatim}
[()->(gamma(1/3)/gamma(2/3))^3,4*n+3+8*k,2*(2*n+3)*(3*n+4)]
\end{verbatim}
Convergence type $E$ with $E=-3$ and $P=2k-1/6$.
\end{cf}

\smallskip

\begin{cf}\label{5.11.178}{\ }
\begin{verbatim}
[()->(gamma(1/3)/gamma(2/3))^3,[0,2,5*(6*n-7)],[9,-16*(3*n-2)^2]]
\end{verbatim}
$$\left(\dfrac{\G(1/3)}{\G(2/3)}\right)^3=\dfrac{9}{2-\dfrac{16}{25-\dfrac{256}{55-\dfrac{784}{85-\dfrac{1600}{115-\dfrac{2704}{145-\ddots}}}}}}$$
Convergence type $E$ with $E=4$, $P=0$ and
$C=2^{-2/3}3^{3/2}(\G(1/3)/\G(2/3))^3$, so that
$$\left(\dfrac{\G(1/3)}{\G(2/3)}\right)^3-\dfrac{p(n)}{q(n)}\sim\dfrac{3^{3/2}(\G(1/3)/\G(2/3))^3}{2^{2n+2/3}}$$
$$A=1-(5/12)/n+(5/288)/n^2-(1255/3456)/n^3+\cdots$$
Parametric family for $k\ge0$:
\begin{verbatim}
[()->(gamma(1/3)/gamma(2/3))^3,30*n-35+18*k,-16*(3*n-2)^2]
\end{verbatim}
Convergence type $E$ with $E=4$ and $P=2k$.
\end{cf}

\smallskip

\begin{cf}\label{5.11.179}{\ }
\begin{verbatim}
[()->(gamma(1/3)/gamma(2/3))^3,[0,1,10*n-11],[3,-4*(2*n-1)^2]]
\end{verbatim}
$$\left(\dfrac{\G(1/3)}{\G(2/3)}\right)^3=\dfrac{3}{1-\dfrac{4}{9-\dfrac{36}{19-\dfrac{100}{29-\dfrac{196}{39-\dfrac{324}{49-\ddots}}}}}}$$
Convergence type $E$ with $E=4$, $P=-1/3$ and
$C=2^{-1/3}3^{4/3}(\G(1/3)/\G(2/3))^4$, so that
$$\left(\dfrac{\G(1/3)}{\G(2/3)}\right)^3-\dfrac{p(n)}{q(n)}\sim\dfrac{3^{4/3}(\G(1/3)/\G(2/3))^4}{2^{2n+1/3}}$$
$$A=1-(25/54)/n-(325/2916)/n^2+\cdots$$
Parametric family for $k\ge0$:
\begin{verbatim}
[()->(gamma(1/3)/gamma(2/3))^3,10*n-11+6*k,-4*(2*n-1)^2]
\end{verbatim}
Convergence type $E$ with $E=4$ and $P=2k-1/3$.
\end{cf}

\smallskip

\begin{cf}\label{5.11.180}{\ }
\begin{verbatim}
[()->(gamma(1/3)/gamma(2/3))^3,[9,26,21*n+1],[-36,8*(3*n+1)^2]]
\end{verbatim}
$$\left(\dfrac{\G(1/3)}{\G(2/3)}\right)^3=9-\dfrac{36}{26+\dfrac{128}{43+\dfrac{392}{64+\dfrac{800}{85+\dfrac{1352}{106+\dfrac{2048}{127+\ddots}}}}}}$$
Convergence type $E$ with $E=-8$, $P=1/3$, and $C=-3^{1/3}(\G(1/3)/\G(2/3))^2/4$, so that
$$\left(\dfrac{\G(1/3)}{\G(2/3)}\right)^3-\dfrac{p(n)}{q(n)}\sim(-1)^{n+1}\dfrac{3^{1/3}(\G(1/3)/\G(2/3))^2}{2^{3n+2}n^{1/3}}$$
$$A=1-(40/81)/n+(2657/6561)/n^2-(470750/1594323)/n^3+\cdots$$
Parametric family for $k\ge0$:
\begin{verbatim}
[()->(gamma(1/3)/gamma(2/3))^3,21*n+27*k+1,8*(3*n+1)^2]
\end{verbatim}
Convergence type $E$ with $E=-8$ and $P=2k+1/3$.
\end{cf}

\smallskip

\begin{cf}\label{5.11.181}{\ }
\begin{verbatim}
[()->(gamma(1/3)/gamma(2/3))^3,[9,28,21*n+5],[-36,8*(3*n-1)^2]]
\end{verbatim}
$$\left(\dfrac{\G(1/3)}{\G(2/3)}\right)^3=9-\dfrac{36}{28+\dfrac{32}{47+\dfrac{200}{68+\dfrac{512}{89+\dfrac{968}{110+\dfrac{1568}{131+\ddots}}}}}}$$
Convergence type $E$ with $E=-8$, $P=5/3$, and
$C=-4\cdot3^{-13/3}(\G(1/3)/\G(2/3))^4$, so that
$$\left(\dfrac{\G(1/3)}{\G(2/3)}\right)^3-\dfrac{p(n)}{q(n)}\sim(-1)^{n+1}\dfrac{3^{-13/3}(\G(1/3)/\G(2/3))^4}{2^{3n-2}n^{5/3}}$$
$$A=1-(82/81)/n+(4394/6561)/n^2-(115988/1594323)/n^3+\cdots$$
Parametric family for $k\ge0$:
\begin{verbatim}
[()->(gamma(1/3)/gamma(2/3))^3,21*n+27*k+5,8*(3*n-1)^2]
\end{verbatim}
Convergence type $E$ with $E=-8$ and $P=2k+5/3$.
\end{cf}

\smallskip

\begin{cf}\label{5.11.182}{\ }
\begin{verbatim}
[()->(gamma(1/3)/gamma(2/3))^3,[6,8,10*n],[12,-(3*n+1)*(3*n+2)]]
\end{verbatim}
$$\left(\dfrac{\G(1/3)}{\G(2/3)}\right)^3=6+\dfrac{12}{8-\dfrac{20}{20-\dfrac{56}{30-\dfrac{110}{40-\dfrac{182}{50-\dfrac{272}{60-\ddots}}}}}}$$
Convergence type $E$ with $E=9$, $P=0$, and $C=(\G(1/3)/\G(2/3))^3/3$, so that
$$\left(\dfrac{\G(1/3)}{\G(2/3)}\right)^3-\dfrac{p(n)}{q(n)}\sim\dfrac{(\G(1/3)/\G(2/3))^3}{3^{2n+1}}$$
$$A=1-(5/18)/n+(205/648)/n^2-(3755/8748)/n^3+\cdots$$
Parametric family for $k\ge0$:
\begin{verbatim}
[()->(gamma(1/3)/gamma(2/3))^3,10*n+8*k,-(3*n+1)*(3*n+2)]
\end{verbatim}
Convergence type $E$ with $E=9$ and $P=2k$.
\end{cf}

\smallskip

\begin{cf}\label{5.11.183}{\ }
\begin{verbatim}
[()->(gamma(1/3)/gamma(2/3))^3,[6,70,(2*n-1)*(82*n^2-82*n-1)],
                               [120,-560,-(n^2-1)*(9*n^2-1)*(36*n^2-1)]]
\end{verbatim}
$$\left(\dfrac{\G(1/3)}{\G(2/3)}\right)^3=6+\dfrac{120}{70-\dfrac{560}{489-\dfrac{15015}{2455-\dfrac{206720}{6881-\dfrac{1233375}{14751-\dfrac{4833024}{27049-\ddots}}}}}}$$
Convergence type $E$ with $E=81$, $P=0$, and $C=(\G(1/3)/\G(2/3))^3/3$, so that
$$\left(\dfrac{\G(1/3)}{\G(2/3)}\right)^3-\dfrac{p(n)}{q(n)}\sim\dfrac{(\G(1/3)/\G(2/3))^3}{3^{4n+1}}$$
$$A=1-(5/36)/n+(205/2592)/n^2-(3755/69984)/n^3+\cdots$$
\end{cf}

This is the contraction of the previous CF.

\smallskip

\begin{cf}\label{5.11.184}{\ }
\begin{verbatim}
[()->(gamma(1/3)/gamma(2/3))^3,[6,124,(2*n-1)*(153*n^2-153*n-2)],
                               [216,-80,-(n^2-1)*(9*n^2-1)*(9*n^2-4)]]
\end{verbatim}
$$\left(\dfrac{\G(1/3)}{\G(2/3)}\right)^3=6+\dfrac{216}{124-\dfrac{80}{912-\dfrac{3360}{4580-\dfrac{49280}{12838-\dfrac{300300}{27522-\dfrac{1188096}{50468-\ddots}}}}}}$$
Convergence type $E$ with $E=(1+\sqrt{2})^8$, $P=0$, and $C=3^{3/2}(\G(1/3)/\G(2/3))^3/((1+\sqrt{2})^4)$, so that
$$\left(\dfrac{\G(1/3)}{\G(2/3)}\right)^3-\dfrac{p(n)}{q(n)}\sim\dfrac{3^{3/2}(\G(1/3)/\G(2/3))^3}{(1+\sqrt{2})^{8n+4}}$$
$$A=1+(19d/96)/n+(-19d/192+361/9216)/n^2+\cdots$$
\end{cf}

Note that this CF converges quite fast, but the next one even faster.

\smallskip

\begin{cf}\label{5.11.185}{\ }
\begin{verbatim}
[()->(gamma(1/3)/gamma(2/3))^3,[0,31,1012*(n-1)],[240,-(6*n-1)*(6*n-5)]]
\end{verbatim}
$$\left(\dfrac{\G(1/3)}{\G(2/3)}\right)^3=\dfrac{240}{31-\dfrac{5}{1012-\dfrac{77}{2024-\dfrac{221}{3036-\dfrac{437}{4048-\dfrac{725}{5060-\ddots}}}}}}$$
Convergence type $E$ with $E=(16+5\sqrt{10})^4/36$, $P=0$, and
$C=3^{3/2}(\G(1/3)/\G(2/3))^3$, so that
$$\left(\dfrac{\G(1/3)}{\G(2/3)}\right)^3-\dfrac{p(n)}{q(n)}\sim\dfrac{3^{3/2}(\G(1/3)/\G(2/3))^3}{(16+5\sqrt{10})^{4n}6^{-2n}}$$
$$A=1-(253d/5760)/n+(64009/6635520)/n^2+\cdots$$
Parametric family for $3\nmid k$:
\begin{verbatim}
[()->(gamma(1/3)/gamma(2/3))^3,1012*(n-1),-((6*n-3)^2-4*k^2)]
\end{verbatim}
Convergence type $E$ with $E=(16+5\sqrt{10})^4/36$ and $P=0$.
\end{cf}

\smallskip

\begin{cf}\label{5.11.186}{\ }
\begin{verbatim}
[()->2^(1/3)*(gamma(1/3)/gamma(2/3))^3,[0,-1,8*n-10],[12,3*(4*n-1)^2]]
\end{verbatim}
$$2^{1/3}\left(\dfrac{\G(1/3)}{\G(2/3)}\right)^3=\dfrac{12}{-1+\dfrac{27}{6+\dfrac{147}{14+\dfrac{363}{22+\dfrac{675}{30+\dfrac{1083}{38+\ddots}}}}}}$$
Convergence type $E$ with $E=-3$, $P=-1/2$ and
$C=2^{7/3}3^{-1/2}(\G(1/3)\G(1/4)/\G(2/3)^2)^2$, so that
$$2^{1/3}\left(\dfrac{\G(1/3)}{\G(2/3)}\right)^3-\dfrac{p(n)}{q(n)}\sim(-1)^n\dfrac{2^{7/3}(\G(1/3)\G(1/4)/\G(2/3)^2)^2}{3^{n+1/2}n^{-1/2}}$$
$$A=1-(1/32)/n+(205/2048)/n^2+(2417/65536)/n^3+\cdots$$
Parametric family for $k\ge0$:
\begin{verbatim}
[()->2^(1/3)*(gamma(1/3)/gamma(2/3))^3,8*n-10+16*k,3*(4*n-1)^2]
\end{verbatim}
Convergence type $E$ with $E=-3$ and $P=2k-1/2$.
\end{cf}

\smallskip

\begin{cf}\label{5.11.187}{\ }
\begin{verbatim}
[()->2^(1/3)*(gamma(1/3)/gamma(2/3))^3,[0,1,8*n-6],[12,3*(4*n-3)^2]]
\end{verbatim}
$$2^{1/3}\left(\dfrac{\G(1/3)}{\G(2/3)}\right)^3=\dfrac{12}{1+\dfrac{3}{10+\dfrac{75}{18+\dfrac{243}{26+\dfrac{507}{34+\dfrac{867}{42+\ddots}}}}}}$$
Convergence type $E$ with $E=-3$, $P=1/2$ and
$C=2^{-2/3}3^{3/2}(\G(1/3)^2/(\G(2/3)\G(1/4)))^2$, so that
$$2^{1/3}\left(\dfrac{\G(1/3)}{\G(2/3)}\right)^3-\dfrac{p(n)}{q(n)}\sim(-1)^n\dfrac{2^{-2/3}(\G(1/3)^2/(\G(2/3)\G(1/4)))^2}{3^{n-3/2}n^{1/2}}$$
$$A=1-(1/32)/n-(203/2048)/n^2+(2825/65536)/n^3+\cdots$$
Parametric family for $k\ge0$:
\begin{verbatim}
[()->2^(1/3)*(gamma(1/3)/gamma(2/3))^3,8*n-6+16*k,3*(4*n-3)^2]
\end{verbatim}
Convergence type $E$ with $E=-3$ and $P=2k+1/2$.
\end{cf}

\smallskip

\begin{cf}\label{5.11.188}{\ }
\begin{verbatim}
[()->2^(1/3)*(gamma(1/3)/gamma(2/3))^3,[0,1,10*(n-1)],[6,-(4*n-1)*(4*n-3)]]
\end{verbatim}
$$2^{1/3}\left(\dfrac{\G(1/3)}{\G(2/3)}\right)^3=\dfrac{6}{1-\dfrac{3}{10-\dfrac{35}{20-\dfrac{99}{30-\dfrac{195}{40-\dfrac{323}{50-\ddots}}}}}}$$
Convergence type $E$ with $E=4$, $P=0$, and
$C=2^{5/6}3^{1/2}(\G(1/3)/\G(2/3))^3$, so that
$$2^{1/3}\left(\dfrac{\G(1/3)}{\G(2/3)}\right)^3-\dfrac{p(n)}{q(n)}\sim\dfrac{3^{1/2}(\G(1/3)/\G(2/3))^3}{2^{2n-5/6}}$$
$$A=1-(5/16)/n+(25/512)/n^2-(6925/24576)/n^3+\cdots$$
Parametric family for $k\ge0$:
\begin{verbatim}
[()->2^(1/3)*(gamma(1/3)/gamma(2/3))^3,10*n-10+6*k,-(4*n-1)*(4*n-3)]
\end{verbatim}
Convergence type $E$ with $E=4$ and $P=2k$.
\end{cf}

\smallskip

\begin{cf}\label{5.11.189}{\ }
\begin{verbatim}
[()->2^(1/3)*(gamma(1/3)/gamma(2/3))^3,[0,1,12*(n-1)],[12,3*(2*n-1)^2]]
\end{verbatim}
$$2^{1/3}\left(\dfrac{\G(1/3)}{\G(2/3)}\right)^3=\dfrac{12}{1+\dfrac{3}{12+\dfrac{27}{24+\dfrac{75}{36+\dfrac{147}{48+\dfrac{243}{60+\ddots}}}}}}$$
Convergence type $E$ with $E=-(2+\sqrt{3})^2$, $P=0$, and
$C=2^{4/3}3^{1/2}(\G(1/3)/\G(2/3))^3$, so that
$$2^{1/3}\left(\dfrac{\G(1/3)}{\G(2/3)}\right)^3-\dfrac{p(n)}{q(n)}\sim(-1)^n\dfrac{2^{4/3}3^{1/2}(\G(1/3)/\G(2/3))^3}{(2+\sqrt{3})^{2n}}$$
$$A=1-(d/8)/n+(3/128)/n^2+(21d/1024)/n^3+\cdots$$
\end{cf}

\smallskip

\begin{cf}\label{5.11.190}{\ }
\begin{verbatim}
[()->2^(1/3)*(gamma(1/3)/gamma(2/3))^3,[0,3,66*(n-1)],[30,(6*n-1)*(6*n-5)]]
\end{verbatim}
$$2^{1/3}\left(\dfrac{\G(1/3)}{\G(2/3)}\right)^3=\dfrac{30}{3+\dfrac{5}{66+\dfrac{77}{132+\dfrac{221}{198+\dfrac{437}{264+\dfrac{725}{330+\ddots}}}}}}$$
Convergence type $E$ with $E=-((1+\sqrt{5})/2)^{10}$, $P=0$, and
$C=2^{4/3}3^{1/2}(\G(1/3)/\G(2/3))^3$, so that
$$2^{1/3}\left(\dfrac{\G(1/3)}{\G(2/3)}\right)^3-\dfrac{p(n)}{q(n)}\sim(-1)^n\dfrac{2^{4/3}3^{1/2}(\G(1/3)/\G(2/3))^3}{((1+\sqrt{5})/2)^{10n}}$$
$$A=1-(11d/180)/n+(121/12960)/n^2+\cdots$$
\end{cf}

\smallskip

\begin{cf}\label{5.11.191}{\ }
\begin{verbatim}
[()->2^(1/3)*3^(1/2)*(gamma(1/3)/gamma(2/3))^3,
[0,2,30*(n-1)],[36,2*(3*n-1)*(3*n-2)]]
\end{verbatim}
$$2^{1/3}3^{1/2}\left(\dfrac{\G(1/3)}{\G(2/3)}\right)^3=\dfrac{36}{2+\dfrac{4}{30+\dfrac{40}{60+\dfrac{112}{90+\dfrac{220}{120+\dfrac{364}{150+\ddots}}}}}}$$
Convergence type $E$ with $E=-(2+\sqrt{3})^3$, $P=0$, and
$C=2^{7/3}3^{1/2}(\G(1/3)/\G(2/3))^3$,
$$2^{1/3}3^{1/2}\left(\dfrac{\G(1/3)}{\G(2/3)}\right)^3-\dfrac{p(n)}{q(n)}\sim(-1)^n\dfrac{2^{7/3}3^{1/2}(\G(1/3)/\G(2/3))^3}{(2+\sqrt{3})^{3n}}$$
$$A=1-(10d/81)/n+(50/2187)/n^2+(4450d/531441)/n^3+\cdots$$
\end{cf}

\smallskip

\begin{cf}\label{5.11.192}{\ }
\begin{verbatim}
[()->5^(1/6)*(gamma(1/3)/gamma(2/3))^3,[0,3,54*(n-1)],[30,-(3*n-1)*(3*n-2)]]
\end{verbatim}
$$5^{1/6}\left(\dfrac{\G(1/3)}{\G(2/3)}\right)^3=\dfrac{30}{3-\dfrac{2}{54-\dfrac{20}{108-\dfrac{56}{162-\dfrac{110}{216-\dfrac{182}{270-\ddots}}}}}}$$
Convergence type $E$ with $E=((1+\sqrt{5})/2)^{12}$, $P=0$,
and $C=5^{7/6}(\G(1/3)/\G(2/3))^3$, so that
$$5^{1/6}\left(\dfrac{\G(1/3)}{\G(2/3)}\right)^3-\dfrac{p(n)}{q(n)}\sim\dfrac{5^{7/6}(\G(1/3)/\G(2/3))^3}{((1+\sqrt{5})/2)^{12n}}$$
$$A=1-(d/10)/n+(1/40)/n^2+(d/450)/n^3+\cdots$$
\end{cf}

\smallskip

\begin{cf}\label{5.11.193}{\ }
\begin{verbatim}
[()->3^(1/2)*7^(-1/6)*(gamma(1/3)/gamma(2/3))^3,
[0,13,330*(n-1)],[126,-(3*n-1)*(3*n-2)]]
\end{verbatim}
$$3^{1/2}7^{-1/6}\left(\dfrac{\G(1/3)}{\G(2/3)}\right)^3=\dfrac{126}{13-\dfrac{2}{330-\dfrac{20}{660-\dfrac{56}{990-\dfrac{110}{1320-\dfrac{182}{1650-\ddots}}}}}}$$
Convergence type $E$ with $E=((5+\sqrt{21})/2)^6$, $P=0$, and
$C=3^{1/2}7^{5/6}(\G(1/3)/\G(2/3))^3$, so that
$$3^{1/2}7^{-1/6}\left(\dfrac{\G(1/3)}{\G(2/3)}\right)^3-\dfrac{p(n)}{q(n)}\sim\dfrac{3^{1/2}7^{5/6}(\G(1/3)/\G(2/3))^3}{((5+\sqrt{21})/2)^{6n}}$$
$$A=1-(55d/1134)/n+(3025/122472)/n^2+\cdots$$
\end{cf}

\smallskip

\begin{cf}\label{5.11.194}{\ }
\begin{verbatim}
[()->2^(1/3)*(2+sqrt(3))*(gamma(1/3)/gamma(2/3))^3,
[30,7*n+5],[90,(2*n+3)*(4*n+5)]]
\end{verbatim}
$$2^{1/3}(2+\sqrt{3})\left(\dfrac{\G(1/3)}{\G(2/3)}\right)^3=30+\dfrac{90}{12+\dfrac{45}{19+\dfrac{91}{26+\dfrac{153}{33+\dfrac{231}{40+\dfrac{325}{47+\ddots}}}}}}$$
Convergence type $E$ with $E=-8$, $P=-1/4$ and
$C=2^{-25/6}3^{3/4}(2+\sqrt{3})\G(1/3)^{5/2}\G(1/4)/\G(2/3)^{7/2}$, so that
$$2^{1/3}(2+\sqrt{3})\left(\dfrac{\G(1/3)}{\G(2/3)}\right)^3-\dfrac{p(n)}{q(n)}\sim(-1)^n\dfrac{3^{3/4}(2+\sqrt{3})\G(1/3)^{5/2}\G(1/4)/\G(2/3)^{7/2}}{2^{3n+25/6}n^{-1/4}}$$
$$A=1+(79/288)/n-(1541/55296)/n^2+\cdots$$
Parametric family for $k\ge0$:
\begin{verbatim}
[()->2^(1/3)*(2+sqrt(3))*(gamma(1/3)/gamma(2/3))^3,
7*n+5+9*k,(2*n+3)*(4*n+5)]
\end{verbatim}
Convergence type $E$ with $E=-8$ and $P=2k-1/4$.
\end{cf}

\smallskip

\begin{cf}\label{5.11.195}{\ }
\begin{verbatim}
[()->2^(1/3)*(2+sqrt(3))*(gamma(1/3)/gamma(2/3))^3,
[0,1,10*(4*n-5)],[24,-9*(4*n-3)^2]]
\end{verbatim}
$$2^{1/3}(2+\sqrt{3})\left(\dfrac{\G(1/3)}{\G(2/3)}\right)^3=\dfrac{24}{1-\dfrac{9}{30-\dfrac{225}{70-\dfrac{729}{110-\dfrac{1521}{150-\dfrac{2601}{190-\ddots}}}}}}$$
Convergence type $E$ with $E=9$, $P=0$ and
$C=3\cdot2^{4/3}(2+\sqrt{3})(\G(1/3)/\G(2/3))^3$, so that
$$2^{1/3}(2+\sqrt{3})\left(\dfrac{\G(1/3)}{\G(2/3)}\right)^3-\dfrac{p(n)}{q(n)}\sim\dfrac{2^{4/3}(2+\sqrt{3})(\G(1/3)/\G(2/3))^3}{3^{2n-1}}$$
$$A=1-(5/16)/n-(15/512)/n^2-(645/8192)/n^3+\cdots$$
Parametric family for $k\ge0$:
\begin{verbatim}
[()->2^(1/3)*(2+sqrt(3))*(gamma(1/3)/gamma(2/3))^3,
40*n-50+32*k,-9*(4*n-3)^2]
\end{verbatim}
Convergence type $E$ with $E=9$ and $P=2k$.
\end{cf}

\smallskip

\begin{cf}\label{5.11.196}{\ }
\begin{verbatim}
[()->2^(1/3)*(2+sqrt(3))*(gamma(1/3)/gamma(2/3))^3,
[0,1,4*(5*n-6)],[12,-9*(2*n-1)^2]]
\end{verbatim}
$$2^{1/3}(2+\sqrt{3})\left(\dfrac{\G(1/3)}{\G(2/3)}\right)^3=\dfrac{12}{1-\dfrac{9}{16-\dfrac{81}{36-\dfrac{225}{56-\dfrac{441}{76-\dfrac{729}{96-\ddots}}}}}}$$
Convergence type $E$ with $E=9$, $P=-1/2$ and
$C=2^{23/6}3^{-1/2}(2+\sqrt{3})(\G(1/3)\G(1/4)/\G(2/3)^2)^2$, so that
$$2^{1/3}(2+\sqrt{3})\left(\dfrac{\G(1/3)}{\G(2/3)}\right)^3-\dfrac{p(n)}{q(n)}\sim\dfrac{2^{23/6}(2+\sqrt{3})(\G(1/3)\G(1/4)/\G(2/3)^2)^2}{3^{2n+1/2}n^{-1/2}}$$
$$A=1-(25/64)/n-(515/8192)/n^2-(61255/524288)/n^3+\cdots$$
Parametric family for $k\ge0$:
\begin{verbatim}
[()->2^(1/3)*(2+sqrt(3))*(gamma(1/3)/gamma(2/3))^3,
20*n-24+16*k,-9*(2*n-1)^2]
\end{verbatim}
Convergence type $E$ with $E=9$ and $P=2k-1/2$.
\end{cf}

\smallskip

\begin{cf}\label{5.11.197}{\ }
\begin{verbatim}
[()->gamma(1/3)^4/gamma(2/3)^2,
[27,2*(3*n-2)*(9*n^2-12*n+16)],[27,-27*n^3*(3*n-1)^3]]
\end{verbatim}
$$\dfrac{\G(1/3)^4}{\G(2/3)^2}=27+\dfrac{27}{26-\dfrac{216}{224-\dfrac{27000}{854-\dfrac{373248}{2240-\dfrac{2299968}{4706-\dfrac{9261000}{8576-\ddots}}}}}}$$
Convergence type $P^+$ with $P=3$ and $C=\G(1/3)^8/(3\G(2/3))^7$, so that
$$\dfrac{\G(1/3)^4}{\G(2/3)^2}-\dfrac{p(n)}{q(n)}\sim\dfrac{\G(1/3)^8/(3\G(2/3))^7}{n^3}$$
$$A=1-1/n+(11/45)/n^2+\cdots$$
Series:
$$\dfrac{\G(2/3)^2}{\G(1/3)^4}=\dfrac{1}{27}\sum_{n\ge0}\dfrac{(-1/3)_n^3}{n!^3}$$
Parametric family for $k\ge0$:
\begin{verbatim}
[()->gamma(1/3)^4/gamma(2/3)^2,
2*(3*n-2)*(9*n^2-12*n+16+9*k*(2*k+3)),-27*n^3*(3*n-1)^3]
\end{verbatim}
Convergence type $P^+$ with $P=4k+3$.
\end{cf}

\smallskip

\begin{cf}\label{5.11.198}{\ }
\begin{verbatim}
[()->gamma(1/3)^4/gamma(2/3)^2,[18,18,14*n^2+7*n+2],
                               [90,-2*(n+1)*(2*n+1)^2*(6*n+5)]]
\end{verbatim}
$$\dfrac{\G(1/3)^4}{\G(2/3)^2}=18+\dfrac{90}{18-\dfrac{396}{72-\dfrac{2550}{149-\dfrac{9016}{254-\dfrac{23490}{387-\dfrac{50820}{548-\ddots}}}}}}$$
Convergence type $E$ with $E=4/3$, $P=7/6$, and $C=2^{-7/3}3^{7/2}\G(1/3)^2/\sqrt{\pi}$, so that
$$\dfrac{\G(1/3)^4}{\G(2/3)^2}-\dfrac{p(n)}{q(n)}\sim\dfrac{2^{-7/3}3^{7/2}\G(1/3)^2/\sqrt{\pi}}{(4/3)^nn^{7/6}}$$
$$A=1-(395/72)/n+(49633/1152)/n^2+\cdots$$
Series:
$$\dfrac{\G(1/3)^4}{\G(2/3)^2}=18\sum_{n\ge0}\dfrac{n!(5/6)_n}{(3/2)_n^2}(3/4)^n$$
\end{cf}
    
\smallskip

\begin{cf}\label{5.11.199}{\ }
\begin{verbatim}
[()->gamma(1/3)^4/gamma(2/3)^2,
[6561/196,260,45*n^3+144*n^2+61*n-30],
[-65610/49,-2*n*(n+4)*(3*n-2)*(3*n+5)*(3*n+7)*(6*n+7)]]
\end{verbatim}
$$\dfrac{\G(1/3)^4}{\G(2/3)^2}=6561/196-\dfrac{65610/49}{260-\dfrac{10400}{1028-\dfrac{260832}{2664-\dfrac{1646400}{5398-\dfrac{6408320}{9500-\dfrac{19047600}{15240-\ddots}}}}}}$$
Convergence type $E$ with $E=4$, $P=3/2$, and $C=-3^{1/2}\pi^{3/2}/(2^{4/3}\G(2/3)^3)$, so that
$$\dfrac{\G(1/3)^4}{\G(2/3)^2}-\dfrac{p(n)}{q(n)}\sim-\dfrac{3^{1/2}\pi^{3/2}/\G(2/3)^3}{2^{2n+4/3}n^{3/2}}$$
$$A=1-(1/24)/n-(3535/1152)/n^2+(14735/1024)/n^3+\cdots$$
Series:
$$\dfrac{\G(1/3)^4}{\G(2/3)^2}=\dfrac{6561}{196}-\dfrac{10935}{182}\sum_{n\ge0}\dfrac{(n+2)(n+3)(n+4)(8/3)_n}{(3n+1)(3n+4)(3n+7)(3n+10)(19/6)_n}2^{-2n}$$
\end{cf}

\smallskip

\begin{cf}\label{5.11.200}{\ }
\begin{verbatim}
[()->gamma(2/3)^4/gamma(1/3)^2,[1/2,2*(3*n-1)*(9*n^2-6*n+4)],
                               [-1/2,-27*n^3*(3*n+1)^3]]
\end{verbatim}
$$\dfrac{\G(2/3)^4}{\G(1/3)^2}=1/2-\dfrac{1/2}{28-\dfrac{1728}{280-\dfrac{74088}{1072-\dfrac{729000}{2728-\dfrac{3796416}{5572-\dfrac{13824000}{9928-\ddots}}}}}}$$
Convergence type $P^+$ with $P=1$ and $C=-2\G(2/3)^8/\G(1/3)^7$, so that
$$\dfrac{\G(2/3)^4}{\G(1/3)^2}-\dfrac{p(n)}{q(n)}\sim-\dfrac{2\G(2/3)^8/\G(1/3)^7}{n}$$
Series:
$$\dfrac{\G(1/3)^2}{\G(2/3)^4}=2\sum_{n\ge0}\dfrac{(1/3)_n^3}{n!^3}$$
Parametric family for $k\ge0$:
\begin{verbatim}
[()->gamma(2/3)^4/gamma(1/3)^2,
2*(3*n-1)*(9*n^2-6*n+4+9*k*(2*k+1)),-27*n^3*(3*n+1)^3]
\end{verbatim}
Convergence type $P^+$ with $P=4k+1$.
\end{cf}

\smallskip

\begin{cf}\label{5.11.201}{\ }
\begin{verbatim}
[()->gamma(2/3)^4/gamma(1/3)^2,[0,14*n^2-11*n+3],
                               [2,-2*n*(2*n-1)*(2*n+1)*(6*n+1)]]
\end{verbatim}
$$\dfrac{\G(2/3)^4}{\G(1/3)^2}=\dfrac{2}{6-\dfrac{42}{37-\dfrac{780}{96-\dfrac{3990}{183-\dfrac{12600}{298-\dfrac{30690}{441-\ddots}}}}}}$$
Convergence type $E$ with $E=4/3$, $P=11/6$, and $C=2^{-5/3}3^{3/2}\G(2/3)^2/\sqrt{\pi}$, so that
$$\dfrac{\G(2/3)^4}{\G(1/3)^2}-\dfrac{p(n)}{q(n)}\sim\dfrac{2^{-5/3}3^{3/2}\G(2/3)^2/\sqrt{\pi}}{(4/3)^nn^{11/6}}$$
$$A=1-(551/72)/n+(793265/10368)/n^2+\cdots$$
Series:
$$\dfrac{\G(2/3)^4}{\G(1/3)^2}=\dfrac{1}{3}\sum_{n\ge0}\dfrac{n!(7/6)_n}{(3/2)_n(5/2)_n}(3/4)^n$$
\end{cf}

\smallskip

\begin{cf}\label{5.11.202}{\ }
\begin{verbatim}
[()->gamma(2/3)^4/gamma(1/3)^2,[4/9,14*n^3+27*n^2-13*n+2],
                               [8/9,32*n^3*(2*n+1)*(3*n+1)*(3*n+2)]]
\end{verbatim}
$$\dfrac{\G(2/3)^4}{\G(1/3)^2}=4/9+\dfrac{8/9}{30+\dfrac{1920}{196+\dfrac{71680}{584+\dfrac{665280}{1278+\dfrac{3354624}{2362+\dfrac{11968000}{3920+\ddots}}}}}}$$
Convergence type $E$ with $E=-16/9$, $P=3/2$, and $C=32\cdot3^{1/2}\pi^{13/2}/(25\G(1/3)^{12})$, so that
$$\dfrac{\G(2/3)^4}{\G(1/3)^2}-\dfrac{p(n)}{q(n)}\sim(-1)^n\dfrac{32\cdot3^{1/2}\pi^{13/2}/(25\G(1/3)^{12})}{(4/3)^{2n}n^{3/2}}$$
$$A=1-(2353/1800)/n+(1233581/1296000)/n^2+\cdots$$
Series:
$$\dfrac{\G(1/3)^2}{\G(2/3)^4}=\dfrac{9}{4}\sum_{n\ge0}(-1)^n\dfrac{(1/2)_n(1/3)_n(2/3)_n}{n!^3}(3/4)^{2n}$$
\end{cf}

\smallskip

For the next CFs, note that
$$2^{-2/3}\dfrac{\G(1/3)^4}{\G(2/3)^2}=\dfrac{\G(1/2)\G(1/3)^5}{\G(1/6)\G(2/3)^3}$$

\smallskip

\begin{cf}\label{5.11.203}{\ }
\begin{verbatim}
[()->2^(-2/3)*gamma(1/3)^4/gamma(2/3)^2,
[12,12,12*n^2+2],[24,-(2*n+1)^2*(3*n+1)*(3*n+2)]]
\end{verbatim}
$$2^{-2/3}\dfrac{\G(1/3)^4}{\G(2/3)^2}=12+\dfrac{24}{12-\dfrac{180}{50-\dfrac{1400}{110-\dfrac{5390}{194-\dfrac{14742}{302-\dfrac{32912}{434-\ddots}}}}}}$$
Convergence type $P^+$ with $P=2/3$ and $C=3\G(1/3)/\G(2/3)$, so that
$$2^{-2/3}\dfrac{\G(1/3)^4}{\G(2/3)^2}-\dfrac{p(n)}{q(n)}\sim\dfrac{3\G(1/3)/\G(2/3)}{n^{2/3}}$$
$$A=1-(2/3)/n+(655/1296)/n^2-(175/486)/n^3+\cdots$$
Series:
$$2^{-2/3}\dfrac{\G(1/3)^4}{\G(2/3)^2}=12+6\sum_{n\ge0}\dfrac{(5/3)_n}{(2n+3)(7/3)_n}$$
Parametric families for $k\ge0$:
\begin{verbatim}
[()->2^(-2/3)*gamma(1/3)^4/gamma(2/3)^2,
12*n^2+2+2*k*(3*k+2),-(2*n+1)^2*(3*n+1)*(3*n+2)]
[()->2^(1/3)*gamma(1/3)^4/gamma(2/3)^2,
24*n^2+2+4*k*(3*k+1),-(2*n+1)^2*(6*n+1)*(6*n+5)]
\end{verbatim}
Convergence type $P^+$ with $P=2k+2/3$ or $P=2k+1/3$ respectively.
\end{cf}

\smallskip

\begin{cf}\label{5.11.204}{\ }
\begin{verbatim}
[()->2^(-2/3)*gamma(1/3)^4/gamma(2/3)^2,[12,12*n-8],[24,n^2*(3*n-1)^2]]
\end{verbatim}
$$2^{-2/3}\dfrac{\G(1/3)^4}{\G(2/3)^2}=12+\dfrac{24}{4+\dfrac{4}{16+\dfrac{100}{28+\dfrac{576}{40+\dfrac{1936}{52+\dfrac{4900}{64+\ddots}}}}}}$$
Convergence type $P^-$ with $P=4$ and $C=2^{-2/3}\G(1/3)^4/(9\G(2/3)^2)$, so that
$$2^{-2/3}\dfrac{\G(1/3)^4}{\G(2/3)^2}-\dfrac{p(n)}{q(n)}\sim(-1)^n\dfrac{2^{-2/3}\G(1/3)^4/(9\G(2/3)^2)}{n^4}$$
$$A=1-(4/3)/n-(4/3)/n^2+(112/27)/n^3+(853/162)/n^4+\cdots$$
Series:
\begin{align*}
2^{-2/3}\dfrac{\G(1/3)^4}{\G(2/3)^2}&=12+48\sum_{n\ge0}\dfrac{(3n+2)n!^2(5/6)_n^2}{(12n^2+4n+1)(12n^2+28n+17)(3/2)_n^2(4/3)_n^2}\\
2^{2/3}\dfrac{\G(2/3)^2}{\G(1/3)^4}&=\dfrac{1}{6}\sum_{n\ge0}\dfrac{(6n+1)(1/2)_n^2(1/3)_n^2}{(6n^2-4n+1)(6n^2+8n+3)n!^2(5/6)_n^2}\end{align*}
Parametric families for $k\ge0$:
\begin{verbatim}
[()->2^(-2/3)*gamma(1/3)^4/gamma(2/3)^2,4*(k+1)*(3*n-2),n^2*(3*n-1)^2]
[()->2^(-2/3)*gamma(1/3)^4/gamma(2/3)^2,2*(k+1)*(6*n-5),n^2*(3*n-2)^2]
\end{verbatim}
Convergence type $P^-$ with $P=4k+4$.
\end{cf}

\smallskip

\begin{cf}\label{5.11.205}{\ }
\begin{verbatim}
[()->2^(-2/3)*gamma(1/3)^4/gamma(2/3)^2,[18,112,63*n^2+45*n+2],
                            [-36,4*(3*n+1)^2*(3*n+2)*(6*n+1)]]
\end{verbatim}
$$2^{-2/3}\dfrac{\G(1/3)^4}{\G(2/3)^2}=18-\dfrac{36}{112+\dfrac{2240}{344+\dfrac{20384}{704+\dfrac{83600}{1190+\dfrac{236600}{1802+\dfrac{539648}{2540+\ddots}}}}}}$$
Convergence type $E$ with $E=-8$, $P=3/2$, and $C=-\G(1/3)^3/(3\cdot2^{4/3}\pi^{3/2})$, so that
$$2^{-2/3}\dfrac{\G(1/3)^4}{\G(2/3)^2}-\dfrac{p(n)}{q(n)}\sim(-1)^{n+1}\dfrac{\G(1/3)^3/(3\pi^{3/2})}{2^{3n+4/3}n^{3/2}}$$
$$A=1-(15/8)/n+(3145/1152)/n^2-\cdots$$
Series:
$$2^{-2/3}\dfrac{\G(1/3)^4}{\G(2/3)^2}=18\sum_{n\ge0}(-1)^n\dfrac{(2/3)_n}{(3n+1)(7/6)_n}2^{-3n}$$
\end{cf}

\smallskip

\begin{cf}\label{5.11.206}{\ }
\begin{verbatim}
[()->2^(-2/3)*3^(1/2)*gamma(1/3)^4/gamma(2/3)^2,
[32,40*n^3+68*n^2+22*n+3],[-160,-32*n*(n+1)^2*(2*n+1)^2*(2*n+5)]]
\end{verbatim}
$$2^{-2/3}3^{1/2}\dfrac{\G(1/3)^4}{\G(2/3)^2}=32-\dfrac{160}{133-\dfrac{8064}{639-\dfrac{129600}{1761-\dfrac{827904}{3739-\dfrac{3369600}{6813-\dfrac{10454400}{11223-\ddots}}}}}}$$
Convergence type $E$ with $E=4$, $P=3/2$, and
$C=2^{11/3}\pi^{13/2}/(243\G(2/3)^{12})$, so that
$$2^{1/3}3^{1/2}\dfrac{\G(1/3)^4}{\G(2/3)^2}-\dfrac{p(n)}{q(n)}\sim-\dfrac{\pi^{13/2}/(243\G(2/3)^{12})}{2^{2n-11/3}n^{3/2}}$$
$$A=1-(19/8)/n+(2011/384)/n^2+\cdots$$
Series:
$$2^{2/3}3^{-1/2}\dfrac{\G(2/3)^2}{\G(1/3)^4}=\dfrac{1}{32}\sum_{n\ge0}\dfrac{(1/2)_n^2(5/2)_n}{n!(n+1)!^2}2^{-2n}$$
\end{cf}

\smallskip

For the next CFs, note that
$$2^{-1/3}\dfrac{\G(2/3)^4}{\G(1/3)^2}=\dfrac{\G(1/6)\G(2/3)^5}{2\G(1/2)\G(1/3)^3}$$

\smallskip

\begin{cf}\label{5.11.207}{\ }
\begin{verbatim}
[()->2^(-1/3)*gamma(2/3)^4/gamma(1/3)^2,
[1/3,15,12*n^2+4],[1/3,-(2*n+1)^2*(3*n+1)*(3*n+2)]]
\end{verbatim}
$$2^{-1/3}\dfrac{\G(2/3)^4}{\G(1/3)^2}=1/3+\dfrac{1/3}{15-\dfrac{180}{52-\dfrac{1400}{112-\dfrac{5390}{196-\dfrac{14742}{304-\dfrac{32912}{436-\ddots}}}}}}$$
Convergence type $P^+$ with $P=4/3$ and $C=\G(2/3)/(12\G(1/3))$, so that
$$2^{2/3}\dfrac{\G(2/3)^4}{\G(1/3)^2}-\dfrac{p(n)}{q(n)}\sim\dfrac{\G(2/3)/(12\G(1/3))}{n^{4/3}}$$
$$A=1-(4/3)/n+(1169/810)/n^2-(329/243)/n^3+\cdots$$
Series:
$$2^{-1/3}\dfrac{\G(2/3)^4}{\G(1/3)^2}=\dfrac{1}{3}+\dfrac{1}{15}\sum_{n\ge0}\dfrac{(4/3)_n}{(2n+3)(8/3)_n}$$
Parametric families for $k\ge0$:
\begin{verbatim}
[()->2^(-1/3)*gamma(2/3)^4/gamma(1/3)^2,
12*n^2+4+2*k*(3*k+4),-(2*n+1)^2*(3*n+1)*(3*n+2)]
[()->2^(-1/3)*gamma(2/3)^4/gamma(1/3)^2,
24*n^2+10+4*k*(3*k+5),-(2*n+1)^2*(6*n+1)*(6*n+5)]
\end{verbatim}
Convergence type $P^+$ with $P=2k+4/3$ or $P=2k+5/3$ respectively.
\end{cf}

\smallskip

\begin{cf}\label{5.11.208}{\ }
\begin{verbatim}
[()->2^(-1/3)*gamma(2/3)^4/gamma(1/3)^2,[1/3,12*n-4],[1/3,n^2*(3*n+1)^2]]
\end{verbatim}
$$2^{-1/3}\dfrac{\G(2/3)^4}{\G(1/3)^2}=1/3+\dfrac{1/3}{8+\dfrac{16}{20+\dfrac{196}{32+\dfrac{900}{44+\dfrac{2704}{56+\dfrac{6400}{68+\ddots}}}}}}$$
Convergence type $P^-$ with $P=4$ and $C=2^{-1/3}\G(2/3)^4/(9\G(1/3)^2)$, so that
$$2^{-1/3}\dfrac{\G(2/3)^4}{\G(1/3)^2}-\dfrac{p(n)}{q(n)}\sim(-1)^n\dfrac{2^{-1/3}\G(2/3)^4/(9\G(1/3)^2)}{n^4}$$
$$A=1-(8/3)/n+2/n^2+(104/27)/n^3-(869/162)/n^4+\cdots$$
Series:
\begin{align*}
2^{-1/3}\dfrac{\G(2/3)^4}{\G(1/3)^2}&=\dfrac{1}{3}+\dfrac{1}{12}\sum_{n\ge0}\dfrac{(6n+5)n!^2(7/6)_n^2}{(6n^2+4n+1)(6n^2+16n+11)(3/2)_n^2(5/3)_n^2}\\
2^{1/3}\dfrac{\G(1/3)^2}{\G(2/3)^4}&=24\sum_{n\ge0}\dfrac{(3n+1)(1/2)_n^2(2/3)_n^2}{(12n^2-4n+1)(12n^2+20n+9)n!^2(7/6)_n^2}\end{align*}
Parametric families for $k\ge0$:
\begin{verbatim}
[()->2^(-1/3)*gamma(2/3)^4/gamma(1/3)^2,4*(k+1)*(3*n-1),n^2*(3*n+1)^2]
[()->2^(-1/3)*gamma(2/3)^4/gamma(1/3)^2,2*(k+1)*(6*n-1),n^2*(3*n+2)^2]
\end{verbatim}
Convergence type $P^-$ with $P=4k+4$.
\end{cf}

\smallskip

\begin{cf}\label{5.11.209}{\ }
\begin{verbatim}
[()->2^(-1/3)*gamma(2/3)^4/gamma(1/3)^2,[1/3,100,63*n^2+45*n-4],
                             [4,4*(3*n-1)*(3*n+2)*(3*n+4)*(6*n-1)]]
\end{verbatim}
$$2^{-1/3}\dfrac{\G(2/3)^4}{\G(1/3)^2}=1/3+\dfrac{4}{100+\dfrac{1400}{338+\dfrac{17600}{698+\dfrac{77792}{1184+\dfrac{226688}{1796+\dfrac{524552}{2534+\ddots}}}}}}$$
Convergence type $E$ with $E=-8$, $P=3/2$, and
$C=2^{-2/3}\G(2/3)^3/(9\pi^{3/2})$, so that
$$2^{-1/3}\dfrac{\G(2/3)^4}{\G(1/3)^2}-\dfrac{p(n)}{q(n)}\sim(-1)^n\dfrac{\G(2/3)^3/(9\pi^{3/2})}{2^{3n+2/3}n^{3/2}}$$
$$A=1-(11/8)/n+(1841/1152)/n^2-\cdots$$
Series:
$$2^{-1/3}\dfrac{\G(2/3)^4}{\G(1/3)^2}=-\dfrac{2}{3}-2\sum_{n\ge0}(-1)^n\dfrac{(4/3)_n}{(3n-1)(3n+2)(5/6)_n}2^{-3n}$$
\end{cf}

\smallskip

\begin{cf}\label{5.11.210}{\ }
\begin{verbatim}
[()->2^(-1/3)*3^(1/2)*gamma(2/3)^4/gamma(1/3)^2,
[2/3,40*n^3-12*n^2+6*n-1],[-2/3,-32*n^3*(2*n+1)^3]]
\end{verbatim}
$$2^{-1/3}3^{1/2}\dfrac{\G(2/3)^4}{\G(1/3)^2}=2/3-\dfrac{2/3}{33-\dfrac{864}{283-\dfrac{32000}{989-\dfrac{296352}{2391-\dfrac{1492992}{4729-\dfrac{5324000}{8243-\ddots}}}}}}$$
Convergence type $E$ with $E=4$, $P=3/2$, and
$C=-2^{19/3}\pi^{13/2}/(27\G(1/3)^{12})$, so that
$$2^{-1/3}3^{1/2}\dfrac{\G(2/3)^4}{\G(1/3)^2}-\dfrac{p(n)}{q(n)}\sim-\dfrac{\pi^{13/2}/(27\G(1/3)^{12})}{2^{2n-19/3}n^{3/2}}$$
$$A=1-(19/8)/n+(2107/384)/n^2+\cdots$$
Series:
$$2^{1/3}3^{-1/2}\dfrac{\G(1/3)^2}{\G(2/3)^4}=\dfrac{3}{2}\sum_{n\ge0}\dfrac{(1/2)_n^3}{n!^3}2^{-2n}$$
\end{cf}

\smallskip

\begin{cf}\label{5.11.211}{\ }
\begin{verbatim}
[()->gamma(1/3)^5/gamma(2/3)^4,[18,16,18*n^2+2],[72,-(3*n+1)^2*(3*n+2)^2]]
\end{verbatim}
$$\dfrac{\G(1/3)^5}{\G(2/3)^4}=18+\dfrac{72}{16-\dfrac{400}{74-\dfrac{3136}{164-\dfrac{12100}{290-\dfrac{33124}{452-\dfrac{73984}{650-\ddots}}}}}}$$
Convergence type $P^+$ with $P=1/3$ and $C=8\pi^2/\G(2/3)^4$, so that
$$\dfrac{\G(1/3)^5}{\G(2/3)^4}-\dfrac{p(n)}{q(n)}\sim\dfrac{8\pi^2/\G(2/3)^4}{n^{1/3}}$$
$$A=1-(1/3)/n+(113/567)/n^2-(29/243)/n^3+\cdots$$
Series:
$$\dfrac{\G(1/3)^5}{\G(2/3)^4}=18+\dfrac{9}{2}\sum_{n\ge0}\dfrac{(5/3)_n^2}{(7/3)_n^2}$$
Parametric family for $k\ge0$:
\begin{verbatim}
[()->gamma(1/3)^5/gamma(2/3)^4,18*n^2+2+3*k*(3*k+1),-(3*n+1)^2*(3*n+2)^2]
\end{verbatim}
Convergence type $P^+$ with $P=2k+1/3$.
\end{cf}

\smallskip

\begin{cf}\label{5.11.212}{\ }
\begin{verbatim}
[()->gamma(1/3)^5/gamma(2/3)^4,[18,6*(6*n-5)],[144,3*n*(3*n-1)^2*(3*n-2)]]
\end{verbatim}
$$\dfrac{\G(1/3)^5}{\G(2/3)^4}=18+\dfrac{144}{6+\dfrac{12}{42+\dfrac{600}{78+\dfrac{4032}{114+\dfrac{14520}{150+\dfrac{38220}{186+\ddots}}}}}}$$
Convergence type $P^-$ with $P=4$ and $C=\G(1/3)^9/(18\pi^4)$, so that
$$\dfrac{\G(1/3)^5}{\G(2/3)^4}-\dfrac{p(n)}{q(n)}\sim(-1)^n\dfrac{\G(1/3)^9/(18\pi^4)}{n^4}$$
$$A=1-(2/3)/n-(19/9)/n^2+(62/27)/n^3+\cdots$$
Series:
\begin{align*}
  \dfrac{\G(1/3)^5}{\G(2/3)^4}&=18+72\sum_{n\ge0}\dfrac{(12n+7)n!(2/3)_n(5/6)_n^2}{(18n^2+3n+1)(18n^2+39n+22)(3/2)_n(4/3)_n^2(7/6)_n}\\
  \dfrac{\G(2/3)^4}{\G(1/3)^5}&=\dfrac{2}{3}\sum_{n\ge0}\dfrac{(12n+1)(1/6)_n(1/2)_n(1/3)_n^2}{(18n^2-15n+4)(18n^2+21n+7)n!(2/3)_n(5/6)_n^2}\end{align*}
Parametric family for $k\ge0$:
\begin{verbatim}
[()->gamma(1/3)^5/gamma(2/3)^4,6*(k+1)*(6*n-5),3*n*(3*n-1)^2*(3*n-2)]
\end{verbatim}
Convergence type $P^-$ with $P=4(k+1)$.
\end{cf}

\smallskip

\begin{cf}\label{5.11.213}{\ }
\begin{verbatim}
[()->gamma(2/3)^5/gamma(1/3)^4,[1/12,25,18*n^2+8],[1/12,-(3*n+1)^2*(3*n+2)^2]]
\end{verbatim}
$$\dfrac{\G(2/3)^5}{\G(1/3)^4}=1/12+\dfrac{1/12}{25-\dfrac{400}{80-\dfrac{3136}{170-\dfrac{12100}{296-\dfrac{33124}{458-\dfrac{73984}{656-\ddots}}}}}}$$
Convergence type $P^+$ with $P=5/3$ and $C=\G(2/3)^4/(60\pi^2)$, so that
$$\dfrac{\G(2/3)^5}{\G(1/3)^4}-\dfrac{p(n)}{q(n)}\sim\dfrac{\G(2/3)^4/(60\pi^2)}{n^{1/3}}$$
$$A=1-(5/3)/n+(1850/891)/n^2-(530/243)/n^3+\cdots$$
Series:
$$\dfrac{\G(2/3)^5}{\G(1/3)^4}=\dfrac{1}{12}+\dfrac{1}{300}\sum_{n\ge0}\dfrac{(4/3)_n^2}{(8/3)_n^2}$$
Parametric family for $k\ge0$:
\begin{verbatim}
[()->gamma(2/3)^5/gamma(1/3)^4,18*n^2+8+3*k*(3*k+5),-(3*n+1)^2*(3*n+2)^2]
\end{verbatim}
Convergence type $P^+$ with $P=2k+5/3$.
\end{cf}

\smallskip

\begin{cf}\label{5.11.214}{\ }
\begin{verbatim}
[()->gamma(2/3)^5/gamma(1/3)^4,[1/12,6*(6*n-1)],[1/6,3*n*(3*n+1)^2*(3*n+2)]]
\end{verbatim}
$$\dfrac{\G(2/3)^5}{\G(1/3)^4}=1/12+\dfrac{1/6}{30+\dfrac{240}{66+\dfrac{2352}{102+\dfrac{9900}{138+\dfrac{28392}{174+\dfrac{65280}{210+\ddots}}}}}}$$
Convergence type $P^-$ with $P=4$ and $C=\G(2/3)^9/(18\pi^4)$, so that
$$\dfrac{\G(2/3)^5}{\G(1/3)^4}-\dfrac{p(n)}{q(n)}\sim(-1)^n\dfrac{\G(2/3)^9/(18\pi^4)}{n^4}$$
$$A=1-(10/3)/n+(41/9)/n^2+(10/27)/n^3+\cdots$$
Series:
\begin{align*}
  \dfrac{\G(2/3)^5}{\G(1/3)^4}&=\dfrac{1}{12}+\dfrac{1}{15}\sum_{n\ge0}\dfrac{(12n+11)n!(4/3)_n(7/6)_n^2}{(18n^2+15n+4)(18n^2+51n+37)(3/2)_n(5/3)_n^2(11/6)_n}\\
  \dfrac{\G(1/3)^4}{\G(2/3)^5}&=36\sum_{n\ge0}\dfrac{(12n+5)(5/6)_n(1/2)_n(2/3)_n^2}{(18n^2-3n+1)(18n^2+33n+16)n!(4/3)_n(7/6)_n^2}\end{align*}
Parametric family for $k\ge0$:
\begin{verbatim}
[()->gamma(2/3)^5/gamma(1/3)^4,6*(k+1)*(6*n-1),3*n*(3*n+1)^2*(3*n+2)]
\end{verbatim}
Convergence type $P^-$ with $P=4(k+1)$.
\end{cf}

\smallskip

\begin{cf}\label{5.11.215}{\ }
\begin{verbatim}
[()->(gamma(1/3)/gamma(2/3))^6,[[0,11],[6*n-3,11]],
[[2750/3,32],[2*(3*n-2)*(3*n+2)^2,2*(3*n-2)^2*(3*n+2)]]]
\end{verbatim}
$$\left(\dfrac{\G(1/3)}{\G(2/3)}\right)^6=\dfrac{2750/3}{11+\dfrac{32}{3+\dfrac{50}{11+\dfrac{10}{9+\dfrac{512}{11+\dfrac{256}{15+\ddots}}}}}}$$
Convergence type $P^-$ with $P=4$ and $C=(1375/1728)\G(1/3)^7/\G(2/3)^5$,
so that
$$\left(\dfrac{\G(1/3)}{\G(2/3)}\right)^6-\dfrac{p(n)}{q(n)}\sim(-1)^n\dfrac{(1375/1728)\G(1/3)^7/\G(2/3)^5}{n^4}$$
\end{cf}

\smallskip

There exist infinitely many CFs of the same type, with $a(n)=(6n-3,a_0)$
for instance for $a_0=11,17,27,35,45,65,77,81,\dotsc$.

\smallskip

\begin{cf}\label{5.11.216}{\ }
\begin{verbatim}
[()->(gamma(1/3)/gamma(2/3))^6,[54,43,9*(2*n-1)*(6*n^2-6*n+5)],
                               [216,-4*(9*n^2-1)^3]]
\end{verbatim}
$$\left(\dfrac{\G(1/3)}{\G(2/3)}\right)^6=54+\dfrac{216}{43-\dfrac{2048}{459-\dfrac{171500}{1845-\dfrac{2048000}{4851-\dfrac{11696828}{10125-\dfrac{44957696}{18315-\ddots}}}}}}$$
Convergence type $P^+$ with $P=2$ and $C=\G(1/3)^7/(96\G(2/3)^5)$, so that
$$\left(\dfrac{\G(1/3)}{\G(2/3)}\right)^6-\dfrac{p(n)}{q(n)}\sim\dfrac{\G(1/3)^7/(96\G(2/3)^5)}{n^2}$$
Parametric family for $k\ge0$:
\begin{verbatim}
[()->(gamma(1/3)/gamma(2/3))^6,9*(2*n-1)*(6*n^2-6*n+5+12*k*(k+1)),
                               -4*(9*n^2-1)^3]
\end{verbatim}
Convergence type $P^+$ with $P=4k+2$.
\end{cf}

\smallskip

\begin{cf}\label{5.11.217}{\ }
\begin{verbatim}
[()->(gamma(1/3)/gamma(2/3))^6,[20626/343,2197,36*n*(12*n^2+49)],
                               [-131072/343,-(6*n-1)^3*(6*n+7)^3]]
\end{verbatim}
$$\left(\dfrac{\G(1/3)}{\G(2/3)}\right)^6=20626/343-\dfrac{131072/343}{2197-\dfrac{274625}{6984-\dfrac{9129329}{16956-\dfrac{76765625}{34704-\dfrac{362467097}{62820-\ddots}}}}}$$
Convergence type $P^+$ with $P=6$ and $C=-(256/3^{10})(\G(1/3)/\G(2/3))^6$,
so that
$$\left(\dfrac{\G(1/3)}{\G(2/3)}\right)^6-\dfrac{p(n)}{q(n)}\sim-\dfrac{(256/3^{10})(\G(1/3)/\G(2/3))^6}{n^6}$$
$$A=1-6/n+(182/9)/n^2-(448/9)/n^3+\cdots$$
Series:
\begin{align*}\left(\dfrac{\G(1/3)}{\G(2/3)}\right)^6&=-322+\dfrac{131072}{343}\sum_{n\ge0}\dfrac{(-1/6)_n^3}{(13/6)_n^3}\\
  \left(\dfrac{\G(2/3)}{\G(1/3)}\right)^6&=\dfrac{5}{54}-\dfrac{256}{3375}\sum_{n\ge0}\dfrac{(1/6)_n^3}{(11/6)_n^3}\end{align*}
Parametric family for $k\ge0$:
\begin{verbatim}
[()->(gamma(1/3)/gamma(2/3))^6,
36*n*(12*n^2+49+24*k*(k+3)),-(6*n-1)^3*(6*n+7)^3]
\end{verbatim}
Convergence type $P^+$ with $P=4k+6$.
\end{cf}

\smallskip

For the next CFs, note that
$$2^{-1/3}\left(\dfrac{\G(1/3)}{\G(2/3)}\right)^6=\dfrac{\G(1/6)\G(1/3)^5}{2\G(1/2)\G(2/3)^5}$$

\smallskip

\begin{cf}\label{5.11.218}{\ }
\begin{verbatim}
[()->2^(-1/3)*(gamma(1/3)/gamma(2/3))^6,
[18,308,531*n^2-177*n+62],[1944,-648*(2*n+1)*(3*n+1)^2*(6*n+1)]]
\end{verbatim}
$$2^{-1/3}\left(\dfrac{\G(1/3)}{\G(2/3)}\right)^6=18+\dfrac{1944}{308-\dfrac{217728}{1832-\dfrac{2063880}{4310-\dfrac{8618400}{7850-\dfrac{24640200}{12452-\dfrac{56567808}{18116-\ddots}}}}}}$$
Convergence type $E$ with $E=32/27$, $P=0$, and $C=(3^{5/2}/20)(\G(1/3)/\G(2/3))^6$,
so that
$$2^{-1/3}\left(\dfrac{\G(1/3)}{\G(2/3)}\right)^6-\dfrac{p(n)}{q(n)}\sim\dfrac{(3^{5/2}/20)(\G(1/3)/\G(2/3))^6}{(32/27)^n}$$
$$A=1-(125/36)/n+(23125/2592)/n^2+\cdots$$
\end{cf}

\smallskip

\begin{cf}\label{5.11.219}{\ }
\begin{verbatim}
[()->2^(-1/3)*(gamma(1/3)/gamma(2/3))^6,
[18,21,28*n^2+2],[162,-3*(2*n+1)^2*(4*n+1)*(4*n+3)]]
\end{verbatim}
$$2^{-1/3}\left(\dfrac{\G(1/3)}{\G(2/3)}\right)^6=18+\dfrac{162}{21-\dfrac{945}{114-\dfrac{7425}{254-\dfrac{28665}{450-\dfrac{78489}{702-\dfrac{175329}{1010-\ddots}}}}}}$$
Convergence type $E$ with $E=4/3$, $P=0$, and $C=3\cdot2^{-11/6}\CS(-3)^2$,
so that
$$2^{-1/3}\left(\dfrac{\G(1/3)}{\G(2/3)}\right)^6-\dfrac{p(n)}{q(n)}\sim\dfrac{3\cdot2^{-11/6}\CS(-3)^2}{(4/3)^n}$$
$$A=1-(33/16)/n+(2145/512)/n^2-(486939/8192)/n^3+\cdots$$
\end{cf}

\smallskip

\begin{cf}\label{5.11.220}{\ }
\begin{verbatim}
[()->2^(-1/3)*(gamma(1/3)/gamma(2/3))^6,
[90,245,108*n^2+2],[-12150,9*(4*n+1)*(4*n+3)*(6*n+1)*(6*n+5)]]
\end{verbatim}
$$2^{-1/3}\left(\dfrac{\G(1/3)}{\G(2/3)}\right)^6=90-\dfrac{12150}{245+\dfrac{24255}{434+\dfrac{196911}{974+\dfrac{766935}{1730+\dfrac{2107575}{2702+\dfrac{4716495}{3890+\ddots}}}}}}$$
Convergence type $E$ with $E=-4$, $P=0$, and $C=-2^{-5/6}(\G(1/3)/\G(2/3))^6$,
so that
$$2^{-1/3}\left(\dfrac{\G(1/3)}{\G(2/3)}\right)^6-\dfrac{p(n)}{q(n)}\sim(-1)^{n+1}\dfrac{(\G(1/3)/\G(2/3))^6}{2^{2n+5/6}}$$
$$A=1-(25/144)/n+(7825/41472)/n^2+\cdots$$
\end{cf}

\smallskip

\begin{cf}\label{5.11.221}{\ }
\begin{verbatim}
[()->2^(-1/3)*(gamma(1/3)/gamma(2/3))^6,
[36,33,40*n^2+2],[324,-9*(2*n+1)^4]]
\end{verbatim}
$$2^{-1/3}\left(\dfrac{\G(1/3)}{\G(2/3)}\right)^6=36+\dfrac{324}{33-\dfrac{729}{162-\dfrac{5625}{362-\dfrac{21609}{642-\dfrac{59049}{1002-\dfrac{131769}{1442-\ddots}}}}}}$$
Convergence type $E$ with $E=9$, $P=0$, and
$C=(2^{5/3}/9)(\G(1/3)/\G(2/3))^6$, so that
$$2^{-1/3}\left(\dfrac{\G(1/3)}{\G(2/3)}\right)^6-\dfrac{p(n)}{q(n)}\sim\dfrac{2^{5/3}(\G(1/3)/\G(2/3))^6}{3^{2n+2}}$$
$$A=1-(1/2)/n+(5/8)/n^2-(7/8)/n^3+(167/128)/n^4+\cdots$$
Parametric families for all $u$:
\begin{verbatim}
[()->2^(-1/3)*(gamma(1/3)/gamma(2/3))^6,
40*n^2-24*u*n+2,-9*(2*n+1)^2*(2*n+2*u+1)^2]
[()->2^(-1/3)*(gamma(1/3)/gamma(2/3))^6,
40*n^2+104*u*n+64*u^2+2,-9*(2*n+1)^2*(2*n+2*u+1)^2]
[()->2^(-1/3)*(gamma(1/3)/gamma(2/3))^6,
40*n^2+48*u*n+2,-9*(2*n+1)^2*(2*n+4*u+1)^2]
[()->2^(-1/3)*(gamma(1/3)/gamma(2/3))^6,
40*n^2+144*u*n+108*u^2+2,-9*(2*n+1)^2*(2*n+2*u+1)*(2*n+6*u+1)]
\end{verbatim}
Convergence type $E$ with $E=9$ and $P=-4u$, $4u$, $-2u$, and $4u$
respectively.
\end{cf}

\smallskip

\begin{cf}\label{5.11.222}{\ }
\begin{verbatim}
[()->2^(-1/3)*(gamma(1/3)/gamma(2/3))^6,
[[99/2,38],[45*n^2+45*n+11,45*n^2+81*n+38]],
[[-81,400],[-9*(2*n+1)^2*(3*n+1)*(3*n+2),(3*n+4)^2*(3*n+5)^2]]]
\end{verbatim}
$$2^{-1/3}\left(\dfrac{\G(1/3)}{\G(2/3)}\right)^6=99/2-\dfrac{81}{38+\dfrac{400}{101-\dfrac{1620}{164+\dfrac{3136}{281-\dfrac{12600}{380+\dfrac{12100}{551-\ddots}}}}}}$$
Convergence type $E$ with $E=-((1+\sqrt{5})/2)^5$, $P=0$, and $C=-2^{5/3}3^{1/2}(\G(1/3)/\G(2/3))^6/((1+\sqrt{5})/2)^{10}$, so that
$$2^{-1/3}\left(\dfrac{\G(1/3)}{\G(2/3)}\right)^6-\dfrac{p(n)}{q(n)}\sim(-1)^{n+1}\dfrac{2^{5/3}3^{1/2}(\G(1/3)/\G(2/3))^6}{((1+\sqrt{5})/2)^{5n+10}}$$
$$A=1-(31d/90)/n+(31d/45+961/3240)/n^2+\cdots$$
\end{cf}

\smallskip

\begin{cf}\label{5.11.223}{\ }
\begin{verbatim}
[()->2^(-1/3)*(gamma(1/3)/gamma(2/3))^6,
[45,159,164*n^2+4],[405,-9*(2*n+1)^2*(3*n+1)*(3*n+2)]]
\end{verbatim}
$$2^{-1/3}\left(\dfrac{\G(1/3)}{\G(2/3)}\right)^6=45+\dfrac{405}{159-\dfrac{1620}{660-\dfrac{12600}{1480-\dfrac{48510}{2628-\dfrac{132678}{4104-\dfrac{296208}{5908-\ddots}}}}}}$$
Convergence type $E$ with $E=81$, $P=0$, and $C=2^{5/3}3^{-7/2}(\G(1/3)/\G(2/3))^6$, so that
$$2^{-1/3}\left(\dfrac{\G(1/3)}{\G(2/3)}\right)^6-\dfrac{p(n)}{q(n)}\sim\dfrac{2^{5/3}(\G(1/3)/\G(2/3))^6}{3^{4n+7/2}}$$
$$A=1-(125/288)/n+(87625/165888)/n^2+\cdots$$
\end{cf}

\smallskip

\begin{cf}\label{5.11.224}{\ }
\begin{verbatim}
[()->(gamma(1/3)/gamma(2/3))^9,
[216,24,24*n^3+6*n],[1296,-(2*n+1)^4*(3*n+1)*(3*n+2)]]
\end{verbatim}
$$\left(\dfrac{\G(1/3)}{\G(2/3)}\right)^9=216+\dfrac{1296}{24-\dfrac{1620}{204-\dfrac{35000}{666-\dfrac{264110}{1560-\dfrac{1194102}{3030-\dfrac{3982352}{5220-\ddots}}}}}}$$
Convergence type $P^+$ with $P=1/3$ and $C=9\cdot2^{-1/3}(\G(1/3)/\G(2/3))^5$,
so that
$$\left(\dfrac{\G(1/3)}{\G(2/3)}\right)^9-\dfrac{p(n)}{q(n)}\sim\dfrac{9\cdot2^{-1/3}(\G(1/3)/\G(2/3))^5}{n^{1/3}}$$
$$A=1-(1/3)/n+(2137/11340)/n^2-(457/4860)/n^3+\cdots$$
Parametric family for $k\ge0$:
\begin{verbatim}
[()->(gamma(1/3)/gamma(2/3))^9,
24*n^3+(8*k*(6*k+1)+6)*n,-(2*n+1)^4*(3*n+1)*(3*n+2)]
\end{verbatim}
Convergence type $P^+$ with $P=4k+1/3$.
\end{cf}

\medskip

Note that in the above we have used the formulas
$$\G(1/6)=\dfrac{3^{1/2}\G(1/3)^2}{2^{1/3}\pi^{1/2}}\text{\quad and\quad}\G(5/6)=\dfrac{3^{1/2}\G(2/3)^2}{2^{2/3}\pi^{1/2}}\;.$$

\smallskip

\subsection{Chowla--Selberg Gamma Quotients}\label{sec:CS}

\smallskip

Let $D$ be a negative fundamental discriminant and denote by $w(D)$ and
$h(D)$ the number of roots of unity and the class number of $\Q(\sqrt{D})$.

We define the \emph{Chowla--Selberg gamma quotient} by
  $$\CS(D)=\Bigg(\prod_{j=1}^{|D|}\G\bigg(\frac j{|D|}\bigg)^{\left(\frac{D}{j}\right)}\Bigg)^{w(D)/(2h(D))}\;.$$

In particular we have $\CS(-3)=(\G(1/3)/\G(2/3))^3$ and
$\CS(-4)=(\G(1/4)/\G(3/4))^2$, and the CFs involving these quantities have
been given in the preceding two subsections. We give here CFs for larger
values of $|D|$.

\medskip

\smallskip

\begin{cf}\label{5.11.225}{\ }
\begin{verbatim}
[()->CS(-7),[0,12,189*(n-1)],[105,-84*(6*n-1)*(6*n-5)]]
\end{verbatim}
$$\CS(-7)=\dfrac{105}{12-\dfrac{420}{189-\dfrac{6468}{378-\dfrac{18564}{567-\dfrac{36708}{756-\dfrac{60900}{945-\ddots}}}}}}$$
Convergence type $E$ with $E=(157+15\sqrt{105})/32$, $P=0$, and
$C=7^{1/2}\CS(-7)$, so that
$$\CS(-7)-\dfrac{p(n)}{q(n)}\sim\dfrac{7^{1/2}\CS(-7)}{((157+15\sqrt{105})/32)^n}\;.$$
$$A=1-(d/60)/n+(7/480)/n^2+\cdots$$
\end{cf}

\smallskip

\begin{cf}\label{5.11.226}{\ }
\begin{verbatim}
[()->CS(-7),[0,13,248*n-250],[168,63*(4*n-1)^2]]
\end{verbatim}
$$\CS(-7)=\dfrac{168}{13+\dfrac{567}{246+\dfrac{3087}{494+\dfrac{7623}{742+\dfrac{14175}{990+\dfrac{22743}{1238+\ddots}}}}}}$$
Convergence type $E$ with $E=-63$, $P=-1/2$ and
$C=2^{7/2}3^{-1}\CS(-7)\G(1/4)/\G(3/4)$, so that
$$\CS(-7)-\dfrac{p(n)}{q(n)}\sim(-1)^n\dfrac{2^{7/2}3^{-1}\CS(-7)\G(1/4)/\G(3/4)}{63^nn^{-1/2}}$$
$$A=1-(91/512)/n+(20485/524288)/n^2+\cdots$$
Parametric family for $k\ge0$:
\begin{verbatim}
[()->CS(-7),248*n-250+256*k,63*(4*n-1)^2]
\end{verbatim}
Convergence type $E$ with $E=-63$ and $P=2k-1/2$.
\end{cf}

\smallskip

\begin{cf}\label{5.11.227}{\ }
\begin{verbatim}
[()->CS(-7),[0,15,248*n-246],[168,63*(4*n-3)^2]]
\end{verbatim}
$$\CS(-7)=\dfrac{168}{15+\dfrac{63}{250+\dfrac{1575}{498+\dfrac{5103}{746+\dfrac{10647}{994+\dfrac{18207}{1242+\ddots}}}}}}$$
Convergence type $E$ with $E=-63$, $P=1/2$ and
$C=21\cdot2^{-5/2}\CS(-7)\G(3/4)/\G(1/4)$, so that
$$\CS(-7)-\dfrac{p(n)}{q(n)}\sim(-1)^n\dfrac{21\cdot2^{-5/2}\CS(-7)\G(3/4)/\G(1/4)}{63^nn^{1/2}}$$
$$A=1-(91/512)/n-(3923/524288)/n^2+\cdots$$
Parametric family for $k\ge0$:
\begin{verbatim}
[()->CS(-7),248*n-246+256*k,63*(4*n-3)^2]
\end{verbatim}
Convergence type $E$ with $E=-63$ and $P=2k+1/2$.
\end{cf}

\smallskip

\begin{cf}\label{5.11.228}{\ }
\begin{verbatim}
[()->CS(-7),[0,4,65*(n-1)],[42,-4*(4*n-1)*(4*n-3)]]
\end{verbatim}
$$\CS(-7)=\dfrac{42}{4-\dfrac{12}{65-\dfrac{140}{130-\dfrac{396}{195-\dfrac{780}{260-\dfrac{1292}{325-\ddots}}}}}}$$
Convergence type $E$ with $E=64$, $P=0$, and $C=2^{1/2}7^{1/2}\CS(-7)$, so that
$$\CS(-7)-\dfrac{p(n)}{q(n)}\sim\dfrac{7^{1/2}\CS(-7)}{2^{6n-1/2}}$$
$$A=1-(65/336)/n+(4225/225792)/n^2+\cdots$$
Parametric family for $k\ge0$:
\begin{verbatim}
[()->CS(-7),65*n-65+63*k,-4*(4*n-1)*(4*n-3)]
\end{verbatim}
Convergence type $E$ with $E=64$ and $P=2k$.
\end{cf}

\smallskip

\begin{cf}\label{5.11.229}{\ }
\begin{verbatim}
[()->CS(-7),[0,5,84*(n-1)],[56,7*(2*n-1)^2]]
\end{verbatim}
$$\CS(-7)=\dfrac{56}{5+\dfrac{7}{84+\dfrac{63}{168+\dfrac{175}{252+\dfrac{343}{336+\dfrac{567}{420+\ddots}}}}}}$$
Convergence type $E$ with $E=-(8+3\sqrt{7})^2$, $P=0$, and
$C=2\cdot7^{1/2}\CS(-7)$, so that
$$\CS(-7)-\dfrac{p(n)}{q(n)}\sim(-1)^n\dfrac{2\cdot7^{1/2}\CS(-7)}{(8+3\sqrt{7})^{2n}}$$
$$A=1-(3d/32)/n+(63/2048)/n^2+\cdots$$
\end{cf}

\smallskip

\begin{cf}\label{5.11.230}{\ }
\begin{verbatim}
[()->CS(-7),[0,324,10773*(n-1)],[3570,84*(6*n-1)*(6*n-5)]]
\end{verbatim}
$$\CS(-7)=\dfrac{3570}{324+\dfrac{420}{10773+\dfrac{6468}{21546+\dfrac{18564}{32319+\dfrac{36708}{43092+\dfrac{60900}{53865+\ddots}}}}}}$$
Convergence type $E$ with $E=-(614093+14535\sqrt{1785})/32$, $P=0$, and
$C=2\cdot7^{1/2}\CS(-7)$, so that
$$\CS(-7)-\dfrac{p(n)}{q(n)}\sim(-1)^n\dfrac{2\cdot 7^{1/2}\CS(-7)}{((614093+14535\sqrt{1785})/32)^n}$$
$$A=1-(19d/5780)/n+(7581/786080)/n^2+\cdots$$
\end{cf}

\smallskip

\begin{cf}\label{5.11.231}{\ }
\begin{verbatim}
[()->(4+sqrt(7))*CS(-7),[0,33,130*(4*n-5)],[1008,-63^2*(4*n-3)^2]]
\end{verbatim}
$$(4+\sqrt{7})\CS(-7)=\dfrac{1008}{33-\dfrac{3969}{390-\dfrac{99225}{910-\dfrac{321489}{1430-\dfrac{670761}{1950-\dfrac{1147041}{2470-\ddots}}}}}}$$
Convergence type $E$ with $E=81/49$, $P=0$ and $C=2(4+\sqrt{7})\CS(-7)$,
so that
$$(4+\sqrt{7})\CS(-7)-\dfrac{p(n)}{q(n)}\sim\dfrac{2(4+\sqrt{7})\CS(-7)}{(9/7)^{2n}}$$
$$A=1-(65/64)/n+(2145/8192)/n^2+\cdots$$
Parametric family for $k\ge0$:
\begin{verbatim}
[()->(4+sqrt(7))*CS(-7),520*n-650+128*k,-63^2*(4*n-3)^2]
\end{verbatim}
Convergence type $E$ with $E=81/49$ and $P=2k$.
\end{cf}

\smallskip

\begin{cf}\label{5.11.232}{\ }
\begin{verbatim}
[()->(4+sqrt(7))*CS(-7),[0,41,4*(65*n-69)],[504,-63^2*(2*n-1)^2]]
\end{verbatim}
$$(4+\sqrt{7})\CS(-7)=\dfrac{504}{41-\dfrac{3969}{244-\dfrac{35721}{504-\dfrac{99225}{764-\dfrac{194481}{1024-\dfrac{321489}{1284-\ddots}}}}}}$$
Convergence type $E$ with $E=81/49$, $P=-1/2$ and
$C=(16/9)(4+\sqrt{7})\CS(-7)\G(1/4)/\G(3/4)$, so that
$$(4+\sqrt{7})\CS(-7)-\dfrac{p(n)}{q(n)}\sim\dfrac{(16/9)(4+\sqrt{7})\CS(-7)\G(1/4)/\G(3/4)}{(9/7)^{2n}n^{-1/2}}$$
$$A=1-(325/256)/n-(233915/131072)/n^2+\cdots$$
Parametric family for $k\ge0$:
\begin{verbatim}
[()->(4+sqrt(7))*CS(-7),260*n-276+64*k,-63^2*(2*n-1)^2]
\end{verbatim}
Convergence type $E$ with $E=81/49$ and $P=2k-1/2$.
\end{cf}

\smallskip

\begin{cf}\label{5.11.233}{\ }
\begin{verbatim}
[()->CS(-7)/Pi,[7/2,512,(10*n-1)*(52*n^2+4*n+1)],[7/2,-512*n^3*(2*n+1)^3]]
\end{verbatim}
$$\dfrac{\CS(-7)}{\pi}=7/2+\dfrac{7/2}{512-\dfrac{13824}{4123-\dfrac{512000}{13949-\dfrac{4741632}{33111-\dfrac{23887872}{64729-\ddots}}}}}$$
Convergence type $E$ with $E=64$, $P=3/2$, and $C=1/(18\pi^{3/2})$, so that
$$\dfrac{\CS(-7)}{\pi}-\dfrac{p(n)}{q(n)}\sim\dfrac{1/(18\pi^{3/2})}{2^{6n}n^{3/2}}$$
$$A=1-(319/168)/n+(505987/169344)/n^2+\cdots$$
Series:
$$\dfrac{\CS(-7)}{\pi}=\dfrac{7}{2}\sum_{n\ge0}\dfrac{(1/2)_n^3}{n!^3}2^{-6n}$$
\end{cf}

\smallskip

\begin{cf}\label{5.11.234}{\ }
\begin{verbatim}
[()->CS(-8),[0,1,8*(n-1)],[8,-2*(2*n-1)^2]]
\end{verbatim}
$$\CS(-8)=\dfrac{8}{1-\dfrac{2}{8-\dfrac{18}{16-\dfrac{50}{24-\dfrac{98}{32-\dfrac{162}{40-\ddots}}}}}}$$
Convergence type $E$ with $E=(1+\sqrt{2})^2$, $P=0$ and $C=2^{3/2}\CS(-8)$,
so that
$$\CS(-8)-\dfrac{p(n)}{q(n)}\sim\dfrac{2^{3/2}\CS(-8)}{(1+\sqrt{2})^{2n}}$$
$$A=1-(d/4)/n+(1/16)/n^2-(d/8)/n^3+\cdots$$
\end{cf}

\smallskip

\begin{cf}\label{5.11.235}{\ }
\begin{verbatim}
[()->CS(-8),[0,3,56*(n-1)],[40,6*(6*n-1)*(6*n-5)]]
\end{verbatim}
$$\CS(-8)=\dfrac{40}{3+\dfrac{30}{56+\dfrac{462}{112+\dfrac{1326}{168+\dfrac{2622}{224+\dfrac{4350}{280+\ddots}}}}}}$$
Convergence type $E$ with $E=-(223+70\sqrt{10})/27$, $P=0$ and
$C=2^{3/2}\CS(-8)$, so that
$$\CS(-8)-\dfrac{p(n)}{q(n)}\sim(-1)^n\dfrac{2^{3/2}\CS(-8)}{((223+70\sqrt{10})/27)^n}$$
$$A=1-(7d/180)/n+(49/6480)/n^2+\cdots$$
\end{cf}

\smallskip

\begin{cf}\label{5.11.236}{\ }
\begin{verbatim}
[()->6^(1/2)*CS(-8),[0,36,960*(n-1)],[1008,6*(4*n-1)*(4*n-3)]]
\end{verbatim}
$$6^{1/2}\CS(-8)=\dfrac{1008}{36+\dfrac{18}{960+\dfrac{210}{1920+\dfrac{594}{2880+\dfrac{1170}{3840+\dfrac{1938}{4800+\ddots}}}}}}$$
Convergence type $E$ with $E=-(5+2\sqrt{6})^4$, $P=0$, and
$C=6^{3/2}\CS(-8)$, so that
$$6^{1/2}\CS(-8)-\dfrac{p(n)}{q(n)}\sim(-1)^n\dfrac{6^{3/2}\CS(-8)}{(5+2\sqrt{6})^{4n}}$$
$$A=1-(15d/196)/n+(675/38416)/n^2+\cdots$$
\end{cf}

\smallskip

\begin{cf}\label{5.11.237}{\ }
\begin{verbatim}
[()->CS(-8)^2,
[320/7,1287/7,236*n^2+21],[163840/49,-384*(2*n+1)^2*(3*n+1)*(3*n+2)]]
\end{verbatim}
$$\CS(-8)^2=320/7+\dfrac{163840/49}{1287/7-\dfrac{69120}{965-\dfrac{537600}{2145-\dfrac{2069760}{3797-\dfrac{5660928}{5921-\dfrac{12638208}{8517-\ddots}}}}}}$$
Convergence type $E$ with $E=32/27$, $P=0$, and $C=(1/16)3^{5/2}\CS(-8)^2$,
so that
$$\CS(-8)^2-\dfrac{p(n)}{q(n)}\sim\dfrac{(1/16)3^{5/2}\CS(-8)^2}{(32/27)^n}$$
$$A=1-(125/36)/n+(24625/2592)/n^2+\cdots$$
\end{cf}

\smallskip

\begin{cf}\label{5.11.238}{\ }
\begin{verbatim}
[()->CS(-8)^2,[64,9,12*n^2+1],[256,-2*(2*n+1)^4]]
\end{verbatim}
$$\CS(-8)^2=64+\dfrac{256}{9-\dfrac{162}{49-\dfrac{1250}{109-\dfrac{4802}{193-\dfrac{13122}{301-\dfrac{29282}{433-\ddots}}}}}}$$
Convergence type $E$ with $E=2$, $P=0$, and $C=\CS(-8)^2$, so that
$$\CS(-8)^2-\dfrac{p(n)}{q(n)}\sim\dfrac{\CS(-8)^2}{2^n}$$
$$A=1-1/n+(3/2)/n^2-(25/4)/n^3+(155/8)/n^4+\cdots$$
Parametric family for $k\ge0$:
\begin{verbatim}
[()->CS(-8)^2,12*n^2+1,-2*(2*n+1)^3*(2*n+1-2*k)]
\end{verbatim}
Convergence type $E$ with $E=2$ and $P=3k$.
\end{cf}
        
\smallskip

\begin{cf}\label{5.11.239}{\ }
\begin{verbatim}
[()->CS(-8)^2,[192,45,28*n^2+1],[-3072,8*(2*n+1)^4]]
\end{verbatim}
$$\CS(-8)^2=192-\dfrac{3072}{45+\dfrac{648}{113+\dfrac{5000}{253+\dfrac{19208}{449+\dfrac{52488}{701+\dfrac{117128}{1009+\ddots}}}}}}$$
Convergence type $E$ with $E=-8$, $P=0$, and $C=-\CS(-8)^2/2$, so that
$$\CS(-8)^2-\dfrac{p(n)}{q(n)}\sim(-1)^{n+1}\dfrac{\CS(-8)^2}{2^{3n+1}}$$
$$A=1-(1/3)/n+(7/18)/n^2-(115/324)/n^3+\cdots$$
\end{cf}

\smallskip

\begin{cf}\label{5.11.240}{\ }
\begin{verbatim}
[()->CS(-8)^2,[320/3,99,116*n^2+3],[6400/3,-6*(2*n+1)^2*(6*n+1)*(6*n+5)]]
\end{verbatim}
$$\CS(-8)^2=320/3+\dfrac{6400/3}{99-\dfrac{4158}{467-\dfrac{33150}{1047-\dfrac{128478}{1859-\dfrac{352350}{2903-\dfrac{787710}{4179-\ddots}}}}}}$$
Convergence type $E$ with $E=27/2$, $P=0$, and $C=(8/27)\CS(-8)^2$, so that
$$\CS(-8)^2-\dfrac{p(n)}{q(n)}\sim\dfrac{4\CS(-8)^2}{(27/2)^{n+1}}$$
$$A=1-(88/225)/n+(23672/50625)/n^2+\cdots$$
\end{cf}

\smallskip

\begin{cf}\label{5.11.241}{\ }
\begin{verbatim}
[()->CS(-11),[0,15,308*(n-1)],[176,-33*(6*n-1)*(6*n-5)]]
\end{verbatim}
$$\CS(-11)=\dfrac{176}{15-\dfrac{165}{308-\dfrac{2541}{616-\dfrac{7293}{924-\dfrac{14421}{1232-\dfrac{23925}{1540-\ddots}}}}}}$$
Convergence type $E$ with $E=(1051+224\sqrt{22})/27$, $P=0$ and $C=11^{1/2}\CS(-11)$, so that
$$\CS(-11)-\dfrac{p(n)}{q(n)}\sim\dfrac{11^{1/2}\CS(-11)}{((1051+224\sqrt{22})/27)^n}$$
$$A=1-(35d/1152)/n+(13475/1327104)/n^2+\cdots$$
\end{cf}

\smallskip

\begin{cf}\label{5.11.242}{\ }
\begin{verbatim}
[()->CS(-15),[0,2,15*(n-1)],[15,-5*(3*n-1)*(3*n-2)]]
\end{verbatim}
$$\CS(-15)=\dfrac{15}{2-\dfrac{10}{15-\dfrac{100}{30-\dfrac{280}{45-\dfrac{550}{60-\dfrac{910}{75-\ddots}}}}}}$$
Convergence type $E$ with $E=((1+\sqrt{5})/2)^2$, $P=0$, and
$C=5^{1/2}\CS(-15)$, so that
$$\CS(-15)-\dfrac{p(n)}{q(n)}\sim\dfrac{5^{1/2}\CS(-15)}{((1+\sqrt{5})/2)^{2n}}$$
$$A=1-(2d/9)/n+(10/81)/n^2-(1010d/2187)/n^3+\cdots$$
\end{cf}

\smallskip

\begin{cf}\label{5.11.243}{\ }
\begin{verbatim}
[()->CS(-15),[0,2,33*(n-1)],[30,(3*n-1)*(3*n-2)]]
\end{verbatim}
$$\CS(-15)=\dfrac{30}{2+\dfrac{2}{33+\dfrac{20}{66+\dfrac{56}{99+\dfrac{110}{132+\dfrac{182}{165+\ddots}}}}}}$$
Convergence type $E$ with $E=-((1+\sqrt(5))/2)^{10}$, $P=0$, and
$C=2\cdot 5^{1/2}\CS(-15)$, so that
$$\CS(-15)-\dfrac{p(n)}{q(n)}\sim(-1)^n\dfrac{2\cdot 5^{1/2}\CS(-15)}{((1+\sqrt{5})/2)^{10n}}$$
$$A=1-(22d/225)/n+(242/10125)/n^2+\cdots$$
\end{cf}

\smallskip

\begin{cf}\label{5.11.244}{\ }
\begin{verbatim}
[()->CS(-19),[0,75,2052*(n-1)],[912,-57*(6*n-1)*(6*n-5)]]
\end{verbatim}
$$\CS(-19)=\dfrac{912}{75-\dfrac{285}{2052-\dfrac{4389}{4104-\dfrac{12597}{6156-\dfrac{24909}{8208-\dfrac{41325}{10260-\ddots}}}}}}$$
Convergence type $E$ with $E=1025+96\sqrt{114}$, $P=0$ and
$C=19^{1/2}\CS(-19)$, so that
$$\CS(-19)-\dfrac{p(n)}{q(n)}\sim\dfrac{19^{1/2}\CS(-19)}{(1025+96\sqrt{114})^n}$$
$$A=1-(5d/384)/n+(475/49152)/n^2+\cdots$$
\end{cf}

\smallskip

\begin{cf}\label{5.11.245}{\ }
\begin{verbatim}
[()->CS(-20),[0,3,40*(n-1)],[40,-5*(4*n-1)*(4*n-3)]]
\end{verbatim}
$$\CS(-20)=\dfrac{40}{3-\dfrac{15}{40-\dfrac{175}{80-\dfrac{495}{120-\dfrac{975}{160-\dfrac{1615}{200-\ddots}}}}}}$$
Convergence type $E$ with $E=((1+\sqrt{5})/2)^6$, $P=0$ and
$C=10^{1/2}\CS(-20)$, so that
$$\CS(-20)-\dfrac{p(n)}{q(n)}\sim\dfrac{10^{1/2}\CS(-20)}{((1+\sqrt{5})/2)^{6n}}$$
$$A=1-(3d/32)/n+(45/2048)/n^2+\cdots$$
\end{cf}

\smallskip

\begin{cf}\label{5.11.246}{\ }
\begin{verbatim}
[()->CS(-24),[0,1/2,6*(n-1)],[12,(3*n-1)*(3*n-2)]]
\end{verbatim}
$$\CS(-24)=\dfrac{12}{1/2+\dfrac{2}{6+\dfrac{20}{12+\dfrac{56}{18+\dfrac{110}{24+\dfrac{182}{30+\ddots}}}}}}$$
Convergence type $E$ with $E=-(1+\sqrt{2})^2$, $P=0$,
and $C=2^{3/2}\CS(-24)$, so that
$$\CS(-24)-\dfrac{p(n)}{q(n)}\sim(-1)^n\dfrac{2^{3/2}\CS(-24)}{(1+\sqrt{2})^{2n}}$$
$$A=1-(d/9)/n+(1/81)/n^2+(71d/2187)/n^3+\cdots$$
\end{cf}

\smallskip

\begin{cf}\label{5.11.247}{\ }
\begin{verbatim}
[()->2^(1/2)*CS(-24),[0,2,32*(n-1)],[48,2*(4*n-1)*(4*n-3)]]
\end{verbatim}
$$2^{1/2}\CS(-24)=\dfrac{48}{2+\dfrac{6}{32+\dfrac{70}{64+\dfrac{198}{96+\dfrac{390}{128+\dfrac{646}{160+\ddots}}}}}}$$
Convergence type $E$ with $E=-(1+\sqrt{2})^4$, $P=0$, and
$C=24^{1/2}\CS(-24)$, so that
$$2^{1/2}\CS(-24)-\dfrac{p(n)}{q(n)}\sim(-1)^n\dfrac{24^{1/2}\CS(-24)}{(1+\sqrt{2})^{4n}}$$
$$A=1-(d/8)/n+(1/64)/n^2+(d/96)/n^3+\cdots$$
\end{cf}

\smallskip

\begin{cf}\label{5.11.248}{\ }
\begin{verbatim}
[()->CS(-40),[0,8,160*(n-1)],[120,5*(4*n-1)*(4*n-3)]]
\end{verbatim}
$$\CS(-40)=\dfrac{120}{8+\dfrac{15}{160+\dfrac{175}{320+\dfrac{495}{480+\dfrac{975}{640+\dfrac{1615}{800+\ddots}}}}}}$$
Convergence type $E$ with $E=-((1+\sqrt{5})/2)^{12}$, $P=0$ and
$C=2\cdot5^{1/2}\CS(-40)$, so that
$$\CS(-40)-\dfrac{p(n)}{q(n)}\sim(-1)^n\dfrac{2\cdot5^{1/2}\CS(-40)}{((1+\sqrt{5})/2)^{12n}}$$
$$A=1-(d/12)/n+(5/288)/n^2+(55d/20736)/n^3+\cdots$$
\end{cf}

\smallskip

\begin{cf}\label{5.11.249}{\ }
\begin{verbatim}
[()->CS(-43),[0,2367,97524*(n-1)],
             [20640,-129*(6*n-1)*(6*n-5)]]
\end{verbatim}
$$\CS(-43)=\dfrac{20640}{2367-\dfrac{645}{97524-\dfrac{9933}{195048-\dfrac{28509}{292572-\dfrac{56373}{390096-\ddots}}}}}$$
Convergence type $E$ with $E=1024001+40320\sqrt{645}$, $P=0$ and
$C=43^{1/2}\CS(-43)$, so that
$$\CS(-43)-\dfrac{p(n)}{q(n)}\sim\dfrac{43^{1/2}\CS(-43)}{(1024001+40320\sqrt{645})^n}$$
$$A=1-(7d/1280)/n+(6321/655360)/n^2+\cdots$$
\end{cf}

\smallskip

\begin{cf}\label{5.11.250}{\ }
\begin{verbatim}
[()->CS(-51),[0,7,102*(n-1)],[102,-17*(3*n-1)*(3*n-2)]]
\end{verbatim}
$$\CS(-51)=\dfrac{102}{7-\dfrac{34}{102-\dfrac{340}{204-\dfrac{952}{306-\dfrac{1870}{408-\dfrac{3094}{510-\ddots}}}}}}$$
Convergence type $E$ with $E=(4+\sqrt{17})^2$, $P=0$ and
$C=17^{1/2}\CS(-51)$, so that
$$\CS(-51)-\dfrac{p(n)}{q(n)}\sim\dfrac{17^{1/2}\CS(-51)}{(4+\sqrt{17})^{2n}}$$
$$A=1-(d/18)/n+(17/648)/n^2-(d/4374)/n^3+\cdots$$
\end{cf}

\smallskip

\begin{cf}\label{5.11.251}{\ }
\begin{verbatim}
[()->CS(-52),[0,23,520*(n-1)],[312,-13*(4*n-1)*(4*n-3)]]
\end{verbatim}
$$\CS(-52)=\dfrac{312}{23-\dfrac{39}{520-\dfrac{455}{1040-\dfrac{1287}{1560-\dfrac{2535}{2080-\dfrac{4199}{2600-\ddots}}}}}}$$
Convergence type $E$ with $E=((3+\sqrt{13})/2)^6$, $P=0$ and
$C=26^{1/2}\CS(-52)$, so that
$$\CS(-52)-\dfrac{p(n)}{q(n)}\sim\dfrac{26^{1/2}\CS(-52)}{((3+\sqrt{13})/2)^{6n}}$$
$$A=1-(5d/96)/n+(325/18432)/n^2+\cdots$$
\end{cf}

\smallskip

\begin{cf}\label{5.11.252}{\ }
\begin{verbatim}
[()->CS(-67),[0,30531,1570212*(n-1)],[176880,-201*(6*n-1)*(6*n-5)]]
\end{verbatim}
$$\CS(-67)=\dfrac{176880}{30531-\dfrac{1005}{1570212-\dfrac{15477}{3140424-\dfrac{44421}{4710636-\dfrac{87837}{6280848-\dfrac{145725}{7851060-\ddots}}}}}}$$
Convergence type $E$ with $E=170368001+1145760\sqrt{22110}$, $P=0$ and
$C=67^{1/2}\CS(-67)$, so that
$$\CS(-67)-\dfrac{p(n)}{q(n)}\sim\dfrac{67^{1/2}\CS(-67)}{(170368001+1145760\sqrt{22110})^n}$$
$$A=1-(217d/232320)/n+(3154963/327106560)/n^2+\cdots$$
\end{cf}

\smallskip

\begin{cf}\label{5.11.253}{\ }
\begin{verbatim}
[()->2^(1/2)*CS(-88),[0,38,1120*(n-1)],[528,2*(4*n-1)*(4*n-3)]]
\end{verbatim}
$$2^{1/2}\CS(-88)=\dfrac{528}{38+\dfrac{6}{1120+\dfrac{70}{2240+\dfrac{198}{3360+\dfrac{390}{4480+\dfrac{646}{5600+\ddots}}}}}}$$
Convergence type $E$ with $E=-(1+\sqrt{2})^{12}$, $P=0$, and
$C=88^{1/2}\CS(-88)$, so that
$$2^{1/2}\CS(-88)-\dfrac{p(n)}{q(n)}\sim(-1)^n\dfrac{88^{1/2}\CS(-88)}{(1+\sqrt{2})^{12n}}$$
$$A=1-(35*d/264)/n+(1225/69696)/n^2+\cdots$$
\end{cf}

\smallskip

\begin{cf}\label{5.11.254}{\ }
\begin{verbatim}
[()->CS(-123),[0,53,1230*(n-1)],[492,-41*(3*n-1)*(3*n-2)]]
\end{verbatim}
$$\CS(-123)=\dfrac{492}{53-\dfrac{82}{1230-\dfrac{820}{2460-\dfrac{2296}{3690-\dfrac{4510}{4920-\dfrac{7462}{6150-\ddots}}}}}}$$
Convergence type $E$ with $E=(32+5\sqrt{41})^2$, $P=0$ and
$C=41^{1/2}\CS(-123)$, so that
$$\CS(-123)-\dfrac{p(n)}{q(n)}\sim\dfrac{41^{1/2}\CS(-123)}{(32+5\sqrt{41})^{2n}}$$
$$A=1-(5d/144)/n+(1025/41472)/n^2+(275d/279936)/n^3+\cdots$$
\end{cf}

\smallskip

\begin{cf}\label{5.11.255}{\ }
\begin{verbatim}
[()->CS(-148),[0,1123,42920*(n-1)],[6216,-37*(4*n-1)*(4*n-3)]]
\end{verbatim}
$$\CS(-148)=\dfrac{6216}{1123-\dfrac{111}{42920-\dfrac{1295}{85840-\dfrac{3663}{128760-\dfrac{7215}{171680-\dfrac{11951}{214600-\ddots}}}}}}$$
Convergence type $E$ with $E=(6+\sqrt{37})^6$, $P=0$, and
$C=74^{1/2}\CS(-148)$, so that
$$\CS(-148)-\dfrac{p(n)}{q(n)}\sim\dfrac{74^{1/2}\CS(-148)}{(6+\sqrt{37})^{6n}}$$
$$A=1-(145d/4704)/n+(777925/44255232)/n^2+\cdots$$
\end{cf}

\smallskip

\begin{cf}\label{5.11.256}{\ }
\begin{verbatim}
[()->CS(-163),[0,40774227,3270840804*(n-1)],
              [52186080,-489*(6*n-1)*(6*n-5)]]
\end{verbatim}
$$\CS(-163)=\dfrac{52186080}{40774227-\dfrac{2445}{3270840804-\dfrac{37653}{6541681608-\dfrac{108069}{9812522412-\ddots}}}}$$
Convergence type $E$ with $E=303862746112001+237944192640\sqrt{1630815}$,
$P=0$, and $C=163^{1/2}\CS(-163)$, so that
$$\CS(-163)-\dfrac{p(n)}{q(n)}\sim\dfrac{163^{1/2}\CS(-163)}{(303862746112001+237944192640\sqrt{1630815})^n}$$
$$A=1-(185801d/1708373760)/n+\cdots$$
\end{cf}

\smallskip

\begin{cf}\label{5.11.257}{\ }
\begin{verbatim}
[()->CS(-232),[0,8824,422240*(n-1)],[22968,29*(4*n-1)*(4*n-3)]]
\end{verbatim}
$$\CS(-232)=\dfrac{22968}{8824+\dfrac{87}{422240+\dfrac{1015}{844480+\dfrac{2871}{1266720+\dfrac{5655}{1688960+\dfrac{9367}{2111200+\ddots}}}}}}$$
Convergence type $E$ with $E=-((5+\sqrt{29})/2)^{12}$, $P=0$ and
$C=116^{1/2}\CS(-232)$, so that
$$\CS(-232)-\dfrac{p(n)}{q(n)}\sim(-1)^n\dfrac{116^{1/2}\CS(-232)}{((5+\sqrt{29})/2)^{12n}}$$
$$A=1-(455d/13068)/n+(6003725/341545248)/n^2+\cdots$$
\end{cf}

\smallskip

\begin{cf}\label{5.11.258}{\ }
\begin{verbatim}
[()->CS(-267),[0,827,28302*(n-1)],[2670,-89*(3*n-1)*(3*n-2)]]
\end{verbatim}
$$\CS(-267)=\dfrac{2670}{827-\dfrac{178}{28302-\dfrac{1780}{56604-\dfrac{4984}{84906-\dfrac{9790}{113208-\dfrac{16198}{141510-\ddots}}}}}}$$
Convergence type $E$ with $E=(500+53\sqrt{89})^2$, $P=0$,
and $C=89^{1/2}\CS(-267)$, so that
$$\CS(-267)-\dfrac{p(n)}{q(n)}\sim\dfrac{89^{1/2}\CS(-267)}{(500+53\sqrt{89})^{2n}}$$
$$A=1-(53d/2250)/n+(250001/10125000)/n^2+\cdots$$
\end{cf}

\section{Constants Coming from Bessel Functions}

\smallskip

\begin{cf}\label{5.12.1}{\ }
\begin{verbatim}
[()->besselj(1,2)/besselj(0,2),[n],[1,-1]]
\end{verbatim}
$$\dfrac{J_1(2)}{J_0(2)}=\dfrac{1}{1-\dfrac{1}{2-\dfrac{1}{3-\dfrac{1}{4-\dfrac{1}{5-\dfrac{1}{6-\ddots}}}}}}$$
Convergence type $F^2$ with $E=1$, $P=1$, and $C=1/J_0(2)^2$, so that
$$\dfrac{J_1(2)}{J_0(2)}-\dfrac{p(n)}{q(n)}\sim\dfrac{1/J_0(2)^2}{n!^2n}\;.$$
$$A=1-3/n+7/n^2-(52/3)/n^3+45/n^4-(1924/15)/n^5+\cdots$$
\end{cf}

\smallskip

\begin{cf}\label{5.12.2}{\ }
\begin{verbatim}
[()->besselj(1,2)/besselj(0,2),[2,7,(2*n+1)*(2*n^2+2*n-1)],[4,-n*(n+2)]]
\end{verbatim}
$$\dfrac{J_1(2)}{J_0(2)}=2+\dfrac{4}{7-\dfrac{3}{55-\dfrac{8}{161-\dfrac{15}{351-\dfrac{24}{649-\dfrac{35}{1079-\ddots}}}}}}$$
Convergence type $F^4$ with $E=16$, $P=4$, and $C=\pi/(32J_0(2)^2)$, so that
$$\dfrac{J_1(2)}{J_0(2)}-\dfrac{p(n)}{q(n)}\sim\dfrac{\pi/(32J_0(2)^2)}{n!^42^{4n}n^4}\;.$$
$$A=1-(21/4)/n+(565/32)/n^2-(18805/384)/n^3+(252115/2048)/n^4+\cdots$$
\end{cf}

This is the contraction of the previous CF.

\smallskip

\begin{cf}\label{5.12.3}{\ }
\begin{verbatim}
[()->besseli(1,2)/besseli(0,2),[n],[1]]
\end{verbatim}
$$\dfrac{I_1(2)}{I_0(2)}=\dfrac{1}{1+\dfrac{1}{2+\dfrac{1}{3+\dfrac{1}{4+\dfrac{1}{5+\dfrac{1}{6+\ddots}}}}}}$$
Convergence type $F^2$ with $E=-1$, $P=1$, and $C=1/I_0(2)^2$, so that
$$\dfrac{I_1(2)}{I_0(2)}-\dfrac{p(n)}{q(n)}\sim(-1)^n\dfrac{1/I_0(2)^2}{n!^2n}\;.$$
$$A=1+1/n-1/n^2+(4/3)/n^3-(5/3)/n^4-(32/5)/n^5+\cdots$$
\end{cf}

The contraction of this CF is also rather simple and analogous to
\ref{5.12.2}.

\smallskip

\begin{cf}\label{5.12.4}{\ }
\begin{verbatim}
[()->besselj(0,2),[0,4,n*(n+2)],[(n+1)^2]]
[()->besselj(0,2),[0,4,n+2],[1,(n+1)/n]]
\end{verbatim}
$$J_0(2)=\dfrac{1}{4+\dfrac{4}{8+\dfrac{9}{15+\dfrac{16}{24+\dfrac{25}{35+\dfrac{36}{48+\ddots}}}}}}=\dfrac{1}{4+\dfrac{2}{4+\dfrac{3/2}{5+\dfrac{4/3}{6+\dfrac{5/4}{7+\dfrac{6/5}{8+\ddots}}}}}}$$
Convergence type $F^2$ with $E=-1$, $P=4$, and $C=1$, so that
$$J_0(2)-\dfrac{p(n)}{q(n)}\sim(-1)^n\dfrac{1}{n!^2n^4}\;.$$
$$A=1-6/n+22/n^2-60/n^3+116/n^4-62/n^5+\cdots$$
Series:
$$J_0(2)=\sum_{n\ge0}(-1)^n\dfrac{1}{(n+2)!^2}$$
\end{cf}

\smallskip

\begin{cf}\label{5.12.5}{\ }
\begin{verbatim}
[()->besseli(0,2),[1,1,n^2+1],[1,-n^2]]
[()->besseli(0,2),[2,4,n^2+2*n+2],[1,-(n+1)^2]]
\end{verbatim}
$$I_0(2)=1+\dfrac{1}{1-\dfrac{1}{5-\dfrac{4}{10-\dfrac{9}{17-\dfrac{16}{26-\ddots}}}}}=2+\dfrac{1}{4-\dfrac{4}{10-\dfrac{9}{17-\dfrac{16}{26-\dfrac{25}{37-\ddots}}}}}$$
Convergence type $F^2$ with $E=1$, $P=2$ or $4$, and $C=1$, so that
$$I_0(2)-\dfrac{p(n)}{q(n)}\sim\dfrac{1}{n!^2n^{2\text{ or }4}}\;.$$
\begin{align*}
  A&=1-2/n+4/n^2-10/n^3+29/n^4-90/n^5+\cdots\text{\quad or}\\
  A&=1-6/n+24/n^2-84/n^3+288/n^4-1022/n^5+\cdots
\end{align*}
Series:
$$I_0(2)=\sum_{n\ge0}\dfrac{1}{n!^2}$$
\end{cf}

\smallskip

\begin{cf}\label{5.12.6}{\ }
\begin{verbatim}
[()->exp(-1)*besseli(0,1),[0,4,n^2],[(2*n+3)*(n+1)^2]]
\end{verbatim}
$$\dfrac{I_0(1)}{e}=\dfrac{3}{4+\dfrac{20}{4+\dfrac{63}{9+\dfrac{144}{16+\dfrac{275}{25+\dfrac{468}{36+\ddots}}}}}}$$
Convergence type $F^1$ with $E=-1/2$, $P=5/2$, and $C=4/\sqrt{\pi}$, so that
$$\dfrac{I_0(1)}{e}-\dfrac{p(n)}{q(n)}\sim(-1)^n\dfrac{1/\sqrt{\pi}}{n!(1/2)^{n+2}n^{5/2}}$$
$$A=1-(49/8)/n+(4033/128)/n^2-(169867/1024)/n^3+\cdots$$
Series:
$$\dfrac{I_0(1)}{e}=\sum_{n\ge0}(-1)^n\dfrac{(1/2)_n}{n!^2}2^n$$
\end{cf}

\smallskip

\begin{cf}\label{5.12.7}{\ }
\begin{verbatim}
[()->exp(1)*besseli(0,1),[1,1,n^2+2*n-1],[1,-n^2*(2*n+1)]]
\end{verbatim}
$$e\cdot I_0(1)=1+\dfrac{1}{1-\dfrac{3}{7-\dfrac{20}{14-\dfrac{63}{23-\dfrac{144}{34-\dfrac{275}{47-\ddots}}}}}}$$
Convergence type $F^1$ with $E=1/2$, $P=3/2$, and $C=2/\sqrt{\pi}$, so that
$$e\cdot I_0(1)-\dfrac{p(n)}{q(n)}\sim\dfrac{1/\sqrt{\pi}}{n!(1/2)^{n+1}n^{3/2}}$$
$$A=1+(3/8)/n-(263/128)/n^2-(759/1024)/n^3+\cdots$$
Series:
$$e\cdot I_0(1)=\sum_{n\ge0}\dfrac{(1/2)_n}{n!^2}2^n$$
\end{cf}

\smallskip

\begin{cf}\label{5.12.8}{\ }
\begin{verbatim}
[()->cos(1)*besselj(0,1),[1,4,16*n^4-16*n^3-12*n^2+16*n-3],
                         [-3,4*n^2*(2*n-1)^2*(4*n+1)*(4*n+3)]]
\end{verbatim}
$$\cos(1)J_0(1)=1-\dfrac{3}{4+\dfrac{140}{109+\dfrac{14256}{801+\dfrac{175500}{2941+\dfrac{1012928}{7777+\dfrac{3912300}{16941+\ddots}}}}}}$$
Convergence type $F^2$ with $E=-1$, $P=2$, and $C=-1/\sqrt{2}$, so that
$$\cos(1)J_0(1)-\dfrac{p(n)}{q(n)}\sim(-1)^{n+1}\dfrac{1/\sqrt{2}}{n!^2n^2}$$
$$A=1-(31/16)/n+(929/512)/n^2+(18843/8192)/n^3+\cdots$$
Series:
$$\cos(1)J_0(1)=\sum_{n\ge0}(-1)^n\dfrac{(1/4)_n(3/4)_n}{n!^2(1/2)_n^2}$$
\end{cf}

\smallskip

\begin{cf}\label{5.12.9}{\ }
\begin{verbatim}
[()->sin(1)*besselj(0,1),[1,36,16*n^4+16*n^3-12*n^2+1],
                         [-15,4*n^2*(2*n+1)^2*(4*n+3)*(4*n+5)]]
\end{verbatim}
$$\sin(1)J_0(1)=1-\dfrac{15}{36+\dfrac{2268}{337+\dfrac{57200}{1621+\dfrac{449820}{4929+\dfrac{2068416}{11701+\dfrac{6957500}{23761+\ddots}}}}}}$$
Convergence type $F^2$ with $E=-1$, $P=3$, and $C=-1/\sqrt{2}$, so that
$$\sin(1)J_0(1)-\dfrac{p(n)}{q(n)}\sim(-1)^{n+1}\dfrac{1/\sqrt{2}}{n!^2n^3}$$
$$A=1-(59/16)/n+(4201/512)/n^2-(88065/8192)/n^3+\cdots$$
Series:
$$\sin(1)J_0(1)=\sum_{n\ge0}(-1)^n\dfrac{(3/4)_n(5/4)_n}{n!^2(3/2)_n^2}$$
\end{cf}

\smallskip

\begin{cf}\label{5.12.10}{\ }
\begin{verbatim}
[()->besseli(1,1)/besseli(0,1),[0,3,n+3],[1,-(2*n+1)]]
\end{verbatim}
$$\dfrac{I_1(1)}{I_0(1)}=\dfrac{1}{3-\dfrac{3}{5-\dfrac{5}{6-\dfrac{7}{7-\dfrac{9}{8-\dfrac{11}{9-\ddots}}}}}}$$
Convergence type $F^1$ with $E=1/2$, $P=5/2$, and $C=2/(I_0(1)^2e^2\sqrt{\pi})$,
so that
$$\dfrac{I_1(1)}{I_0(1)}-\dfrac{p(n)}{q(n)}\sim\dfrac{2/(I_0(1)^2e^2\sqrt{\pi})}{n!(1/2)^nn^{5/2}}\;.$$
$$A=1+(3/8)/n-(391/128)/n^2+(1929/1024)/n^3+(861499/32768)/n^4+\cdots$$
\end{cf}

\smallskip

\begin{cf}\label{5.12.11}{\ }
\begin{verbatim}
[()->besseli(1,2)/besseli(0,2),[0,4,n+5],[2,-(4*n+2)]]
\end{verbatim}
$$\dfrac{I_1(2)}{I_0(2)}=\dfrac{2}{4-\dfrac{6}{7-\dfrac{10}{8-\dfrac{14}{9-\dfrac{18}{10-\dfrac{22}{11-\ddots}}}}}}$$
Convergence type $F^1$ with $E=1/4$, $P=5/2$, and $C=4/(I_0(2)^2e^4\sqrt{\pi})$,
so that
$$\dfrac{I_1(2)}{I_0(2)}-\dfrac{p(n)}{q(n)}\sim \dfrac{4/(I_0(2)^2e^4\sqrt{\pi})}{n!(1/4)^nn^{5/2}}\;.$$
$$A=1+(35/8)/n+(1081/128)/n^2-(3671/1024)/n^3+(209851/32768)/n^4+\cdots$$
\end{cf}

\smallskip

\begin{cf}\label{5.12.12}{\ }
\begin{verbatim}
[()->besselk(1,1/2)/besselk(0,1/2),[4*n+2],[-(2*n+1)^2]]
[()->besselk(1,1/2)/besselk(0,1/2),[2],[-(2*n+1)/(2*n+3)]]
\end{verbatim}
$$\dfrac{K_1(1/2)}{K_0(1/2)}=2-\dfrac{1}{6-\dfrac{9}{10-\dfrac{25}{14-\dfrac{49}{18-\dfrac{81}{22-\ddots}}}}}=2-\dfrac{1/3}{2-\dfrac{3/5}{2-\dfrac{5/7}{2-\dfrac{7/9}{2-\dfrac{9/11}{2-\ddots}}}}}$$
Convergence type $D^+$ with $D=16$ and $C=-4\pi/K_0(1/2)^2$, so that
$$\dfrac{K_1(1/2)}{K_0(1/2)}-\dfrac{p(n)}{q(n)}\sim-\dfrac{4\pi/K_0(1/2)^2}{e^{4\sqrt{n}}}\;.$$
$$A=1-(55/24)/n^{1/2}+(3025/1152)/n-(600593/414720)/n^{3/2}+\cdots$$
\end{cf}

\smallskip

\begin{cf}\label{5.12.13}{\ }
\begin{verbatim}
[()->besselk(1,1/2)/besselk(0,1/2),[5/3,7,4*(n+1)],
                                  [2/3,-(2*n+1)*(2*n+3)]]
\end{verbatim}
$$\dfrac{K_1(1/2)}{K_0(1/2)}=\dfrac{5}{3}+\dfrac{2/3}{7-\dfrac{15}{12-\dfrac{35}{16-\dfrac{63}{20-\dfrac{99}{24-\dfrac{143}{28-\ddots}}}}}}$$
Convergence type $D^+$ with $D=16$ and $C=4\pi/K_0(1/2)^2$, so that
$$\dfrac{K_1(1/2)}{K_0(1/2)}-\dfrac{p(n)}{q(n)}\sim\dfrac{4\pi/K_0(1/2)^2}{e^{4\sqrt{n}}}\;.$$
$$A=1-(67/24)/n^{1/2}+(4489/1152)/n-(1088933/414720)/n^{3/2}+\cdots$$
Parametric family:
\begin{verbatim}
[()->besselk(1,1/2)/besselk(0,1/2),
         4*n+2*u+2*v+2,-(2*n+2*u+1)*(2*n+2*v+1)]
\end{verbatim}
Convergence type $D^+$ with $D=16$.
\end{cf}

\smallskip

\begin{cf}\label{5.12.14}{\ }
\begin{verbatim}
[()->besselk(3/4,1)/besselk(1/4,1),[1,4],[1,4*n+2]]
\end{verbatim}
$$\dfrac{K_{3/4}(1)}{K_{1/4}(1)}=1+\dfrac{1}{4+\dfrac{6}{4+\dfrac{10}{4+\dfrac{14}{4+\dfrac{18}{4+\dfrac{22}{4+\ddots}}}}}}$$
Convergence type $D^-$ with $D=16$ and $C=\pi\sqrt{2}/K_{1/4}(1)^2$, so that
$$\dfrac{K_{3/4}(1)}{K_{1/4}(1)}-\dfrac{p(n)}{q(n)}\sim(-1)^n\dfrac{\pi\sqrt{2}/K_{1/4}(1)^2}{e^{4\sqrt{n}}}\;.$$
$$A=1-(8/3)/n^{1/2}+(32/9)/n-(4013/1620)/n^{3/2}+(346/1215)/n^2+\cdots$$
\end{cf}

\smallskip

\begin{cf}\label{5.12.15}{\ }
\begin{verbatim}
[()->besselk(3/4,1)/besselk(1/4,1),[1,11,4*(2*n+1)],[2,-(16*n^2-1)]]
\end{verbatim}
$$\dfrac{K_{3/4}(1)}{K_{1/4}(1)}=1+\dfrac{2}{11-\dfrac{15}{20-\dfrac{63}{28-\dfrac{143}{36-\dfrac{255}{44-\dfrac{399}{52-\ddots}}}}}}$$
Convergence type $D^+$ with $D=32$ and $C=\pi\sqrt{2}/K_{1/4}(1)^2$, so that
$$\dfrac{K_{3/4}(1)}{K_{1/4}(1)}-\dfrac{p(n)}{q(n)}\sim\dfrac{\pi\sqrt{2}/K_{1/4}(1)^2}{e^{4\sqrt{2n}}}\;.$$
$$A=1-(4d/3)/n^{1/2}+(16/9)/n-(4013d/6480)/n^{3/2}+(173/2430)/n^2+\cdots$$
\end{cf}

This is simply the contraction of the previous CF.

\smallskip

\begin{cf}\label{5.12.16}{\ }
\begin{verbatim}
[()->besselk(3/4,1)/besselk(1/4,1),[5/4,8*(n+1)],[-3/4,-(4*n+1)*(4*n+3)]]
\end{verbatim}
$$\dfrac{K_{3/4}(1)}{K_{1/4}(1)}=5/4-\dfrac{3/4}{16-\dfrac{35}{24-\dfrac{99}{32-\dfrac{195}{40-\dfrac{323}{48-\dfrac{483}{56-\ddots}}}}}}$$
Convergence type $D^+$ with $D=32$ and $C=-\pi\sqrt{2}/K_{1/4}(1)^2$, so that
$$\dfrac{K_{3/4}(1)}{K_{1/4}(1)}-\dfrac{p(n)}{q(n)}\sim-\dfrac{\pi\sqrt{2}/K_{1/4}(1)^2}{e^{4\sqrt{2n}}}\;.$$
$$A=1-(7d/3)/n^{1/2}+(49/9)/n-(23363d/6480)/n^{3/2}+(19481/9720)/n^2+\cdots$$
Parametric family:
\begin{verbatim}
[()->besselk(3/4,1)/besselk(1/4,1),4*(2*n+u+v),-(4*n+4*u-1)*(4*n+4*v-3)]
\end{verbatim}
Convergence type $D^+$ with $D=32$.
\end{cf}

\smallskip

\begin{cf}\label{5.12.17}{\ }
\begin{verbatim}
[()->besselk(3/4,1)/besselk(1/4,1),
[[1,8],[6*n+7,6*n+9]],[[2,24],[-(36*n^2-1),48*n+24]]]
\end{verbatim}
$$\dfrac{K_{3/4}(1)}{K_{1/4}(1)}=1+\dfrac{2}{8+\dfrac{24}{13-\dfrac{35}{15+\dfrac{72}{19-\dfrac{143}{21+\dfrac{120}{25-\ddots}}}}}}$$
Convergence type $D^-$ with $D=24$ and $C=\pi\sqrt{2}/K_{1/4}(1)^2$, so that
$$\dfrac{K_{3/4}(1)}{K_{1/4}(1)}-\dfrac{p(n)}{q(n)}\sim(-1)^n\dfrac{\pi\sqrt{2}/K_{1/4}(1)^2}{e^{2\sqrt{6n}}}\;.$$
$$A=1-(8d/9)/n^{1/2}+(64/27)/n-(4013d/7290)/n^{3/2}+(1384/10935)/n^2+\cdots$$
\end{cf}

\smallskip

Note that this is the \emph{unique} CF among those given in this dictionary
for which $b(2n)$ and $b(2n+1)$ do not have the same degrees, although it
is of course possible to construct others (\ref{5.3.77} has different
degrees for $a(2n)$ and $a(2n+1)$).

\medskip

\section{Other Constants}

The reader will also find some more CFs for constants in the next chapter,
obtained by specializing the CFs for \emph{functions} that we will
give. We have included only those who seemed
the simplest/nicest, but of course one can give infinitely many.

There are evidently also infinitely many variations on the CFs that I have
given in this chapter (for instance, I have tried to give some interesting
parametric solutions). I would be very interested in knowing genuinely
new CFs for interesting constants, i.e., not variations of the above nor
specializations of the CFs for functions that we give in the next chapter.

\chapter{Encyclopedic Dictionary: Functions}

\section{Trivial Functions with an Arbitrary Sequence}

Let $(f(n))_{n\ge0}$ be an (a priori arbitrary) sequence. The CFs in this
section converge to a value \emph{independent of $f$}, and are valid under
suitable assumptions on $f$, which will not be specified since it is difficult
to determine the exact domain of validity, but usually for most reasonable $f$,
including for instance non-constant polynomials.

All the continued fractions are of the form
\begin{verbatim}
[()->U(2)/U(1),[0,U(n)*f(n)+V(n)],
              [U(n+2)*f(n+1)+(U(n+2)+V(n)*U(n+1))/U(n)]]
\end{verbatim}
for suitable functions $U$ and $V$. The game consists not so much in giving
random examples, which have little interest, but in choosing $U$ and $V$ so
that the coefficients of the CF are on the one hand polynomials, and second
of degree as small as possible. In particular, we can choose $\deg(V)<\deg(U)$.

Note that the speed of convergence depends
on $f$ so cannot be given, but it is almost always of factorial type.
We will thus only give examples, directly taken from
\cite{Bow-Lau}. Many more examples are of course possible.

\bigskip

\begin{cf}\label{6.1.1}{\ }
\begin{verbatim}
[()->1,[0,f(n)],[f(n+1)+1]]
\end{verbatim}
$$1=\dfrac{f(1)+1}{f(1)+\dfrac{f(2)+1}{f(2)+\dfrac{f(3)+1}{f(3)+\dfrac{f(4)+1}{f(4)+\dfrac{f(5)+1}{f(5)+\dfrac{f(6)+1}{f(6)+\ddots}}}}}}$$
\end{cf}
This is the specialization of \ref{6.1.2} to $z=0$, $k=1$, and $f(n)$
replaced by $f(n)+2$, and is essentially the only example with $\deg(U)=0$.

\smallskip

\begin{cf}\label{6.1.2}{\ }
\begin{verbatim}
[(k,z)->k/(k-z),[0,((n-2)*z+k)*f(n)-2],[(n*z+k)*f(n+1)-1]]
\end{verbatim}
$$\dfrac{k}{k-z}=\dfrac{f(1)k-1}{-f(1)z+f(1)k-2+\dfrac{f(2)z+f(2)k-1}{f(2)k-2+\dfrac{2f(3)z+f(3)k-1}{f(3)z+f(3)k-2+\dfrac{3f(4)z+f(4)k-1}{2f(4)z+f(4)k-2+\ddots}}}}$$
\end{cf}
This is the general example with $\deg(U)=1$.

\smallskip

\begin{cf}\label{6.1.3}{\ }
\begin{verbatim}
[(z)->1,[0,((n-1)*(n-2)*z+1)*f(n)+(n-2)*z-2],
      [(n*(n+1)*z+1)*f(n+1)+(n+1)*z-1]]
\end{verbatim}
$$1=\dfrac{z+f(1)-1}{-z+f(1)-2+\dfrac{(2f(2)+2)z+f(2)-1}{f(2)-2+\dfrac{(6f(3)+3)z+f(3)-1}{(2f(3)+1)z+f(3)-2+\dfrac{(12f(4)+4)z+f(4)-1}{(6f(4)+2)z+f(4)-2+\ddots}}}}$$
\end{cf}

\smallskip

\begin{cf}\label{6.1.4}{\ }
\begin{verbatim}
[()->4,[0,n^2*f(n)+4*n-4],[(n+2)^2*f(n+1)+4*n+9]]
\end{verbatim}
$$4=\dfrac{4f(1)+9}{f(1)+\dfrac{9f(2)+13}{4f(2)+4+\dfrac{16f(3)+17}{9f(3)+8+\dfrac{25f(4)+21}{16f(4)+12+\dfrac{36f(5)+25}{25f(5)+16+\dfrac{49f(6)+29}{36f(6)+20+\ddots}}}}}}$$
\end{cf}

\smallskip

\begin{cf}\label{6.1.5}{\ }
\begin{verbatim}
[(z)->2*z^2-2*z+1,[0,(n^2*z^2-n*z*(z+2)+2*z+1)*f(n)+(n-1)*z-3],
                [((n+2)^2*z^2-(n+2)*z*(z+2)+2*z+1)*f(n+1)+(n+2)*z-2]]
\end{verbatim}
$$2z^2-2z+1= \dfrac{2f(1)z^2+(-2f(1)+2)z+f(1)-2}{f(1)-3+\dfrac{6f(2)z^2+(-4f(2)+3)z+f(2)-2}{2f(2)z^2+(-2f(2)+1)z+f(2)-3+\ddots}}$$
\end{cf}

\smallskip

\begin{cf}\label{6.1.6}{\ }
\begin{verbatim}
[(z)->2*z+1,[0,(n^2*z^2-n*z*(3*z-2)+2*z^2-2*z+1)*f(n)+(n-2)*z-1],
          [((n+2)^2*z^2-(n+2)*z*(3*z-2)+2*z^2-2*z+1)*f(n+1)+(n+1)*z]]
\end{verbatim}
$$2z+1=\dfrac{(2f(1)+1)z+f(1)}{-z+f(1)-1+\dfrac{2f(2)z^2+(4f(2)+2)z+f(2)}{2f(2)z+f(2)-1+\dfrac{6f(3)z^2+(6f(3)+3)z+f(3)}{2f(3)z^2+(4f(3)+1)z+f(3)-1+\ddots}}}$$
\end{cf}

\smallskip

\begin{cf}\label{6.1.7}{\ }
\begin{verbatim}
[(z)->2*z+1,[0,(n*(n-1)*z+1)*f(n)+(n-1)*z-2],
          [((n+1)*(n+2)*z+1)*f(n+1)+(n+2)*z-1]]
\end{verbatim}
$$2z+1=\dfrac{(2f(1)+2)z+f(1)-1}{f(1)-2+\dfrac{(6f(2)+3)z+f(2)-1}{(2f(2)+1)z+f(2)-2+\dfrac{(12f(3)+4)z+f(3)-1}{(6f(3)+2)z+f(3)-2+\ddots}}}$$
\end{cf}

\smallskip

\begin{cf}\label{6.1.8}{\ }
\begin{verbatim}
[(z)->2*z+3,[0,(n*(n-1)*z+2*n-1)*f(n)-((n-1)*z^2+z+2)],
          [((n+1)*(n+2)*z+2*n+3)*f(n+1)-((n+2)*z^2+z+1)]]
\end{verbatim}
$$2z+3=-\dfrac{2z^2+(-2f(1)+1)z-3f(1)+1}{-z+f(1)-2-\dfrac{3z^2+(-6f(2)+1)z-5f(2)+1}{-z^2+(2f(2)-1)z+3f(2)-2-\ddots}}$$
\end{cf}

\smallskip

\begin{cf}\label{6.1.9}{\ }
\begin{verbatim}
[(z)->-1,[0,((n-2)^2*z+(n-2)*(z-2)-1)*f(n)-((n-2)*z^2-z+2)],
       [(n^2*z+n*(z-2)-1)*f(n+1)-((n+1)*z^2-z+1)]]
\end{verbatim}
$$-1=-\dfrac{z^2-z+f(1)+1}{z^2+z+f(1)-2-\dfrac{2z^2+(-2f(2)-1)z+3f(2)+1}{z-f(2)-2-\dfrac{3z^2+(-6f(3)-1)z+5f(3)+1}{-z^2+(2f(3)+1)z-3f(3)-2-\ddots}}}$$
\end{cf}

The above CFs are examples with $\deg(U)=2$ (more general examples can be
given), and the following two are samples with $\deg(U)=3$.

\smallskip

\begin{cf}\label{6.1.10}{\ }
\begin{verbatim}
[(z)->1,[0,(n*(n-1)*(n-2)*z+1)*f(n)+2*n*(n-2)*z-2],
      [((n+2)*(n+1)*n*z+1)*f(n+1)+2*(n+2)*(n+1)*z-1]]
\end{verbatim}
\end{cf}
$$1=\dfrac{f(1)+4z-1}{f(1)-2z-2+\dfrac{(6z+1)f(2)+12z-1}{f(2)-2+\dfrac{(24z+1)f(3)+24z-1}{(6z+1)f(3)+6z-2+\ddots}}}$$

\smallskip

\begin{cf}\label{6.1.11}{\ }
\begin{verbatim}
[(z)->6*z+1,[0,((n+1)*n*(n-1)*z+1)*f(n)+2*(n+1)*(n-1)*z-2],
          [((n+3)*(n+2)*(n+1)*z+1)*f(n+1)+2*(n+3)*(n+2)*z-1]]
\end{verbatim}
\end{cf}
$$6z+1=\dfrac{(6z+1)f(1)+12z-1}{f(1)-2+\dfrac{(24z+1)f(2)+24z-1}{(6z+1)f(2)+6z-2+\dfrac{(60z+1)f(3)+40z-1}{(24z+1)f(3)+16z-2+\ddots}}}$$

Note that this CF is the same as the previous one with $n$ replaced by
$n+1$ and $f(n)$ by $f(n-1)$.

\smallskip

For instance, let us consider the CF \ref{6.1.4}, which does not
depend on any variable. Choosing for instance $f(n)=1$ and $f(n)=n$
gives the following continued fractions:

\smallskip

\begin{cf}\label{6.1.12}{\ }
\begin{verbatim}
[()->4,[0,n^2+4*n-4],[n^2+8*n+13]]
\end{verbatim}
$$4=\dfrac{13}{1+\dfrac{22}{8+\dfrac{33}{17+\dfrac{46}{28+\dfrac{61}{41+\dfrac{78}{56+\ddots}}}}}}$$
Convergence type $F^2$ with $E=-1$, $P=2$, and $C=...$, so that
$$4-\dfrac{p(n)}{q(n)}\sim(-1)^n\dfrac{C}{n!^2n^2}\;.$$
$$A=1-9/n+61/n^2-364/n^3+(4055/2)/n^4-(108261/10)/n^5+\cdots$$
\end{cf}

\smallskip

\begin{cf}\label{6.1.13}{\ }
\begin{verbatim}
[()->4,[0,n^3+4*n-4],[n^3+5*n^2+12*n+13]]
\end{verbatim}
$$4=\dfrac{13}{1+\dfrac{31}{12+\dfrac{65}{35+\dfrac{121}{76+\dfrac{205}{141+\dfrac{323}{236+\ddots}}}}}}$$
Convergence type $F^3$ with $E=-1$, $P=-2$, and $C=...$, so that
$$4-\dfrac{p(n)}{q(n)}\sim(-1)^n\dfrac{C}{n!^3n^{-2}}\;.$$
$$A=1+8/n+(49/2)/n^2+(191/6)/n^3+(235/24)/n^4+(45/4)/n^5+\cdots$$
\end{cf}

\medskip

\section{Other Trivial Functions}

\medskip

Note that some of these CFs for ``other'' trivial functions are special cases
of the CFs for arbitrary sequences seen above.

\smallskip

\begin{cf}\label{6.2.1}{\ }
\begin{verbatim}
[(k,z)->0,[n-k+z],[-z*(n+1)]]
\end{verbatim}
$$0=z-k-\dfrac{z}{z-k+1-\dfrac{2z}{z-k+2-\dfrac{3z}{z-k+3-\dfrac{4z}{z-k+4-\dfrac{5z}{z-k+5-\ddots}}}}}$$
Only for $k\in\Z_{\ge0}$. Convergence type $F^1$ with $E=1/z$, $P=-2k$,
and $C=-z^{1-2k}$, so that
$$0-\dfrac{p(n)}{q(n)}\sim-\dfrac{z^{n-2k+1}n^{2k}}{n!}\;.$$
\begin{align*}A&=1+(-k^2+(-2z+2)k-z)/n\\
  &\phantom{=}+(k^4/2+(2z-7/3)k^3+(2z^2-5z+3)k^2+(z^2-z-7/6)k)/n^2+\cdots\end{align*}
\end{cf}

\smallskip

\begin{cf}\label{6.2.2}{\ }
\begin{verbatim}
[(z)->1,[0,n+z-1],[n+z+1]]
\end{verbatim}
$$1=\dfrac{z+1}{z+\dfrac{z+2}{z+1+\dfrac{z+3}{z+2+\dfrac{z+4}{z+3+\dfrac{z+5}{z+4+\dfrac{z+6}{z+5+\ddots}}}}}}$$
Convergence type $F^1$ with $E=-1$, $P=z$, and $C=\G(z+1)$, so that
$$1-\dfrac{p(n)}{q(n)}\sim(-1)^n\dfrac{\G(z+1)}{n!z^n}\;.$$
$$A=1+(-z^2/2-z/2+1)/n+(z^4/8+5z^3/12-z^2/8-17z/12)/n^2+\cdots$$
\end{cf}

\smallskip

\begin{cf}\label{6.2.3}{\ }
\begin{verbatim}
[(a,z)->1,[0,a],[z+a,(z+n*a)^2-a^2]]
\end{verbatim}
$$1=\dfrac{z+a}{a+\dfrac{z^2+2az}{a+\dfrac{z^2+4az+3a^2}{a+\dfrac{z^2+6az+8a^2}{a+\dfrac{z^2+8az+15a^2}{a+\dfrac{z^2+10az+24a^2}{a+\ddots}}}}}}$$
Convergence type $P^-$ with $P=1$ and $C=2z/a$, so that
$$1-\dfrac{p(n)}{q(n)}\sim(-1)^n\dfrac{2z/a}{n}\;.$$
$$A=1-(1/2+z/a)/n+(1/2+z/a)^2/n^2-(1/2+z/a)^3/n^3+\cdots$$
Series:
$$1=\dfrac{z+a}{a}-\dfrac{z}{a}\sum_{n\ge1}\dfrac{1}{n(n+1)}$$
\end{cf}

\smallskip

\begin{cf}\label{6.2.4}{\ }
\begin{verbatim}
[(z)->1,[0,2*n-1],[z,z^2-n^2]]
\end{verbatim}
$$1=\dfrac{z}{1+\dfrac{z^2-1}{3+\dfrac{z^2-4}{5+\dfrac{z^2-9}{7+\dfrac{z^2-16}{9+\dfrac{z^2-25}{11+\ddots}}}}}}$$
Convergence type $P^+$ with $P=2z$ and $C=-2\G(z)/\G(-z)$, so that
$$1-\dfrac{p(n)}{q(n)}\sim-\dfrac{2\G(z)/\G(-z)}{n^{2z}}\;.$$
$$A=1-z/n+(z(z+1)(2z+1)/6)/n^2-(z^2(z+1)(2z+1)/6)/n^3+\cdots$$
\end{cf}

\smallskip

\begin{cf}\label{6.2.5}{\ }
\begin{verbatim}
[(a,b,c)->min(a,b),[0,a+b+(2*n-1)*c],[a*b,-(a+n*c)*(b+n*c)]]
\end{verbatim}
$$a=\dfrac{ba}{a+b+c-\dfrac{(b+c)a+cb+c^2}{a+b+3c-\dfrac{(b+2c)a+2cb+4c^2}{a+b+5c-\dfrac{(b+3c)a+3cb+9c^2}{a+b+7c-\ddots}}}}$$
Valid for $\Re((b-a)/c)>0$. Convergence type $P^+$ with
$P=(b-a)/c$ and $C=(b-a)\G(b/c)/\G(a/c)$, so that
$$a-\dfrac{p(n)}{q(n)}\sim\dfrac{(b-a)\G(b/c)\G(a/c)}{n^{(b-a)/c}}\;.$$
\end{cf}

\smallskip

\begin{cf}\label{6.2.6}{\ }
\begin{verbatim}
[(z)->z,[2*z+2*n+1],[-(n+z+1)^2]]
\end{verbatim}
$$z=2z+1-\dfrac{z^2+2z+1}{2z+3-\dfrac{z^2+4z+4}{2z+5-\dfrac{z^2+6z+9}{2z+7-\dfrac{z^2+8z+16}{2z+9-\dfrac{z^2+10z+25}{2z+11-\dfrac{z^2+12z+36}{2z+13-\ddots}}}}}}$$
Convergence type $L$ with $C=-1$, so that
$$z-\dfrac{p(n)}{q(n)}\sim-\dfrac{1}{\log(n)}\;.$$
\end{cf}

\smallskip

\begin{cf}\label{6.2.7}{\ }
\begin{verbatim}
[(z)->1-1/z,[1,z+2,2*n+1],[-(n+1)^2]]
\end{verbatim}
$$1-\dfrac{1}{z}=1-\dfrac{1}{z+2-\dfrac{4}{5-\dfrac{9}{7-\dfrac{16}{9-\dfrac{25}{11-\ddots}}}}}\;.$$
Convergence type $L$ with $C=-1/z^2$, so that
$$1-\dfrac{1}{z}-\dfrac{p(n)}{q(n)}\sim-\dfrac{1/z^2}{\log(n)}\;.$$
\end{cf}

These two CFs are among the rare CFs of this dictionary with convergence
type $L$, others being \ref{5.3.92}, \ref{5.4.15}, \ref{5.4.27} and
\ref{5.5.10}.

\smallskip

\begin{cf}\label{6.2.8}{\ }
\begin{verbatim}
[(a,b,z)->a*z,[0,n+b-1-(n+a)*z],[(n+a)*(n+b)*z]]
\end{verbatim}
$$az=\dfrac{zba}{-za+b-z+\dfrac{(zb+z)a+zb+z}{-za+b-2z+1+\dfrac{(zb+2z)a+2zb+4z}{-za+b-3z+2+\ddots}}}$$
Convergence type $E$ with $E=-1/z$, $P=b-a-1$, and $C=z(z+1)\G(b)/\G(a)$, so
that
$$az-\dfrac{p(n)}{q(n)}\sim(-1)^n\dfrac{(z+1)\G(b)/\G(a)}{(1/z)^{n+1}n^{b-a-1}}\;.$$
\end{cf}

\smallskip

\begin{cf}\label{6.2.9}{\ }
\begin{verbatim}
[(a,z)->1+a/(z+1),[0,z+(n-1)*a-1],[z+(n+1)*a]]
\end{verbatim}
$$\dfrac{z+a+1}{z+1}=\dfrac{a+z}{z-1+\dfrac{2a+z}{a+z-1+\dfrac{3a+z}{2a+z-1+\dfrac{4a+z}{3a+z-1+\ddots}}}}$$
Convergence type $F^1$ with $E=-a$, $P=z/a-2$, $C=azG(z/a)/(z+1)^2$, so that
$$\dfrac{z+a+1}{z+1}-\dfrac{p(n)}{q(n)}\sim(-1)^n\dfrac{z\G(z/a)/(z+1)^2}{n!a^{n-1}n^{z/a-2}}\;.$$
$$A=1-((z^2-3az-2a^2-6a)/(2a^2))/n+\cdots$$
\end{cf}

\smallskip

\begin{cf}\label{6.2.10}{\ }
\begin{verbatim}
[(z)->(z^2+z+1)/(z^2-z+1),[0,n+z-4],[n+z]]
\end{verbatim}
$$\dfrac{z^2+z+1}{z^2-z+1}=\dfrac{z}{z-3+\dfrac{z+1}{z-2+\dfrac{z+2}{z-1+\dfrac{z+3}{z+\dfrac{z+4}{z+1+\dfrac{z+5}{z+2+\ddots}}}}}}$$
Convergence type $F^1$ with $E=-1$, $P=z-5$, and $C=\G(z)/(z^2-z+1)^2$, so that
$$\dfrac{z^2+z+1}{z^2-z+1}-\dfrac{p(n)}{q(n)}\sim(-1)^n\dfrac{\G(z)/(z^2-z+1)^2}{n!n^{z-5}}$$
$$A=1-((z^2-9z-2)/2)/n+((3z^4-50z^3+177z^2+86z+144)/24)/n^2+\cdots$$
\end{cf}

\smallskip

\begin{cf}\label{6.2.11}{\ }
\begin{verbatim}
[(z)->(z^3+2*z+1)/((z-1)^3+2*z-1),[0,n+z-5],[n+z]]
\end{verbatim}
$$\dfrac{z^3+2z+1}{z^3-3z^2+5z-2}=\dfrac{z}{z-4+\dfrac{z+1}{z-3+\dfrac{z+2}{z-2+\dfrac{z+3}{z-1+\dfrac{z+4}{z+\dfrac{z+5}{z+1+\ddots}}}}}}$$
Convergence type $F^1$ with $E=-1$, $P=z-7$, and $C=\G(z)/(z^3-3z^2+5z-2)^2$, so that
$$\dfrac{z^3+2z+1}{z^3-3z^2+5z-2}-\dfrac{p(n)}{q(n)}\sim(-1)^n\dfrac{\G(z)/(z^3-3z^2+5z-2)^2}{n!n^{z-7}}$$
$$A=1-((z^2-13z+4)/2)/n+((3z^4-74z^3+453z^2-262z+264)/24)/n^2+\cdots$$
\end{cf}

\medskip

\section{Exponential Functions}\label{sec:exp}

\medskip

Preliminary remark: from the point of view of CFs of polynomial type with
rational coefficients, which are the only ones considered in this dictionary,
the functions $f(z)$ and $f(\pi z)$ (such as $\exp(z)$ and $\exp(\pi z)$,
$\sinh(z)$ and $\sinh(\pi z)$, $\tanh(z)$ and $\tanh(\pi z)$ etc...) are
quite different in nature. For instance, usually the CFs for $f(z)$ converge
factorially, while those for $f(\pi z)$ converge only exponentially or even
polynomially.

\medskip

\begin{cf}\label{6.3.1}{\ }
\begin{verbatim}
[(z)->exp(z),[1,1,n+z],[z,-n*z]]
\end{verbatim}
$$e^z=1+\dfrac{z}{1-\dfrac{z}{z+2-\dfrac{2z}{z+3-\dfrac{3z}{z+4-\dfrac{4z}{z+5-\dfrac{5z}{z+6-\ddots}}}}}}$$
Convergence type $F^1$ with $E=1/z$, $P=1$, and $C=z$, so that
$$e^z-\dfrac{p(n)}{q(n)}\sim\dfrac{z}{n!(1/z)^nn}$$
$$A=1+(z-1)/n+(z^2-3z+1)/n^2+(z^3-6z^2+7z-1)/n^3+\cdots$$
Series:
$$e^z=1+z\sum_{n\ge0}\dfrac{z^n}{(n+1)!}$$
\end{cf}
This is the term-by-term continued fraction corresponding to the Taylor
expansion of $e^z$.

\smallskip

\begin{cf}\label{6.3.2}{\ }
\begin{verbatim}
[(z)->exp(z),[1,n-z],[z,n*z]]
\end{verbatim}
$$e^z=1+\dfrac{z}{-z+1+\dfrac{z}{-z+2+\dfrac{2 z}{-z+3+\dfrac{3 z}{-z+4+\dfrac{4 z}{-z+5+\dfrac{5 z}{-z+6+\ddots}}}}}}$$
Convergence type $F^1$ with $E=-1/z$, $P=1$, and $C=ze^{2z}$, so that
$$e^z-\dfrac{p(n)}{q(n)}\sim(-1)^n\dfrac{ze^{2z}}{n!(1/z)^nn}\;.$$
$$A=1-(z+1)/n+(z^2+3z+1)/n^2-(z^3+6z^2+7z+1)/n^3+\cdots$$
Series:
$$e^{-z}=\sum_{n\ge0}\dfrac{(-z)^n}{n!}$$
\end{cf}
This is trivially equivalent to the previous one using $e^z=1/e^{-z}$.

\smallskip

\begin{cf}\label{6.3.3}{\ }
\begin{verbatim}
[(k,z)->exp(z)-sum(j=0,k-1,z^j/j!),[0,n-z+k-1],[z^k/(k-1)!,z*n]]
\end{verbatim}
$$e^z-\sum_{j=0}^{k-1}\dfrac{z^j}{j!}=\dfrac{z^k/(k-1)!}{-z+k+\dfrac{z}{-z+k+1+\dfrac{2z}{-z+k+2+\dfrac{3z}{-z+k+3+\ddots}}}}$$
Convergence type $F^1$ with $E=-1/z$, $P=2k-1$, and $C=(k-1)!z^ke^{2z}$,
so that
$$e^z-\sum_{j=0}^{k-1}\dfrac{z^j}{j!}-\dfrac{p(n)}{q(n)}\sim(-1)^n\dfrac{(k-1)!e^{2z}}{n!(1/z)^{n+k}n^{2k-1}}\;.$$
$$A=1-((2k-1)z+k^2)/n+((2k^2-k)z^2+(2k^3+k^2-k+1)z+(3k^4+2k^3+k)/6)/n^2+\cdots$$
\end{cf}

\smallskip

\begin{cf}\label{6.3.4}{\ }
\begin{verbatim}
[(k,z)->exp(z)-sum(j=0,k-1,z^j/j!),[0,k,n+z+k-1],
                                   [z^k/(k-1)!,-(n+k-1)*z]]
\end{verbatim}
$$e^z-\sum_{j=0}^{k-1}\dfrac{z^j}{j!}=\dfrac{z^k/(k-1)!}{k-\dfrac{kz}{z+(k+1)-\dfrac{(k+1)z}{z+(k+2)-\dfrac{(k+2)z}{z+(k+3)-\ddots}}}}\;.$$
Convergence type $F^1$ with $E=1/z$, $P=k$, and $C=z^k$,
so that
$$e^z-\sum_{j=0}^{k-1}\dfrac{z^j}{j!}-\dfrac{p(n)}{q(n)}\sim\dfrac{1}{n!(1/z)^{n+k}n^{k}}\;.$$
$$A=1+(z-(k^2+k)/2)/n+(z^2-((k^2+3k+2)/2)z+(3k^4+10k^3+9k^2+2k)/24)/n^2+\cdots$$
\end{cf}

\smallskip

\begin{cf}\label{6.3.5}{\ }
\begin{verbatim}
[(z)->exp(z),[1,2-z,4*n-2],[2*z,z^2]]
\end{verbatim}
$$e^z=1+\dfrac{2z}{-z+2+\dfrac{z^2}{6+\dfrac{z^2}{10+\dfrac{z^2}{14+\dfrac{z^2}{18+\dfrac{z^2}{22+\ddots}}}}}}$$
Convergence type $F^2$ with $E=-16/z^2$, $P=0$, and $C=(\pi/2)ze^z$, so
that
$$e^z-\dfrac{p(n)}{q(n)}\sim(-1)^n\dfrac{(\pi/2)ze^z}{n!^2(16/z^2)^n}\;.$$
$$A=1+((z^2-2)/8)/n+((z^4-12z^2+20)/128)/n^2+\cdots$$
\end{cf}

\smallskip

\begin{cf}\label{6.3.6}{\ }
\begin{verbatim}
[(z)->exp(z),[0,z-1,2*z+2*(2*n-1)*(2*n-3)],[-1,(2*n+1)*(2*n-3)*z^2]]
\end{verbatim}
$$e^z=-\dfrac{1}{z-1-\dfrac{3z^2}{2z+6+\dfrac{5z^2}{2z+30+\dfrac{21z^2}{2z+70+\dfrac{45z^2}{2z+126+\dfrac{77z^2}{2z+198+\ddots}}}}}}$$
Convergence type $F^2$ with $E=-16/z^2$, $P=-1$, and $C=2\pi e^z$, so that
$$e^z-\dfrac{p(n)}{q(n)}\sim(-1)^n\dfrac{\pi e^z}{n!^2(16/z^2)^nn^{-1}}\;.$$
$$A=1+((z^2+4z+2)/8)/n+((z^2+4z+2)^2/128)/n^2+\cdots$$
\end{cf}

\smallskip

\begin{cf}\label{6.3.7}{\ }
\begin{verbatim}
[(z)->exp(z),[z^2/2+z+1,(2*n+1)*(2*n^2+2*n-z)],[2*z^3,n^2*(n+2)^2*z^2]]
\end{verbatim}
$$e^z=\dfrac{z^2}{2}+z+1+\dfrac{2z^3}{-3z+12+\dfrac{9z^2}{-5z+60+\dfrac{64z^2}{-7z+168+\dfrac{225z^2}{-9z+360+\dfrac{576z^2}{-11z+660+\ddots}}}}}$$
Convergence type $F^2$ with $E=-16/z^2$, $P=2$, and $C=\pi e^zz^3/32$, so that
$$e^z-\dfrac{p(n)}{q(n)}\sim(-1)^n\dfrac{\pi e^zz/2}{n!^2(16/z^2)^{n+1}n^2}\;.$$
$$A=1+((z^2-8z-10)/8)/n+((z^2-8z-2)(z^2-8z-42)/128)/n^2+\cdots$$
\end{cf}

\smallskip

\begin{cf}\label{6.3.8}{\ }
\begin{verbatim}
[(z)->exp(Pi*z/2),[1,1-z,2],[2*z,z^2+(2*n-1)^2]]
\end{verbatim}
$$e^{\pi z/2}=1+\dfrac{2z}{-z+1+\dfrac{z^2+1}{2+\dfrac{z^2+9}{2+\dfrac{z^2+25}{2+\dfrac{z^2+49}{2+\dfrac{z^2+81}{2+\ddots}}}}}}$$
Convergence type $P^-$ with $P=1$ and $C=(z/2)e^{\pi z/2}\cosh(\pi z/2)$, so that
$$e^{\pi z/2}-\dfrac{p(n)}{q(n)}\sim(-1)^n\dfrac{(z/2)e^{\pi z/2}\cosh(\pi z/2)}{n}\;.$$
$$A=1+(-z^2/16-1/4)/n^2+(z^4/128+7z^2/64+5/16)/n^4+\cdots$$
Parametric family for $k\ge0$:
\begin{verbatim}
[(z)->exp(Pi*z/2),4*k+2,z^2+(2*n-1)^2]
\end{verbatim}
Convergence type $P^-$ with $P=2k+1$.
\end{cf}

\smallskip

\begin{cf}\label{6.3.9}{\ }
\begin{verbatim}
[(z)->exp(Pi*z/2),[1,1-z,3,4],[2*z,z^2+1,z^2+4*(n-1)^2]]
\end{verbatim}
$$e^{\pi z/2}=1+\dfrac{2z}{-z+1+\dfrac{z^2+1}{3+\dfrac{z^2+4}{4+\dfrac{z^2+16}{4+\dfrac{z^2+36}{4+\dfrac{z^2+64}{4+\ddots}}}}}}$$
Convergence type $P^-$ with $P=2$ and $C=(z^2+1)e^{\pi z/2}\sinh(\pi z/2)/8$,
so that
$$e^{\pi z/2}-\dfrac{p(n)}{q(n)}\sim(-1)^n\dfrac{(z^2+1)e^{\pi z/2}\sinh(\pi z/2)/8}{n^2}\;.$$
$$A=1+1/n+(-z^2/8+1/8)/n^2+(-z^2/4-3/4)/n^3+\cdots$$
Parametric family for $k\ge0$:
\begin{verbatim}
[(z)->exp(Pi*z/2),4*(k+1),z^2+4*(n-1)^2]
\end{verbatim}
Convergence type $P^-$ with $P=2k+2$.
\end{cf}

\smallskip

\begin{cf}\label{6.3.10}{\ }
\begin{verbatim}
[(z)->exp(Pi*z/2),[1,1-z,2*n-1],[2*z,z^2+n^2]]
\end{verbatim}
$$e^{\pi z/2}=1+\dfrac{2z}{1-z+\dfrac{z^2+1}{3+\dfrac{z^2+4}{5+\dfrac{z^2+9}{7+\dfrac{z^2+16}{9+\dfrac{z^2+25}{11+\ddots}}}}}}$$
Convergence type $E$ with $E=-(1+\sqrt{2})^2$, $P=0$, and
$C=2e^{\pi z/2}\sinh(\pi z)/(1+\sqrt{2})$, so that
$$e^{\pi z/2}-\dfrac{p(n)}{q(n)}\sim(-1)^n\dfrac{2e^{\pi z/2}\sinh(\pi z)}{(1+\sqrt{2})^{2n+1}}\;.$$
$$A=1-((4z^2+1)/8)d/n+((4z^2+1)(4d+4z^2+1)/64)/n^2+\cdots$$
\end{cf}

\smallskip

\begin{cf}\label{6.3.11}{\ }
\begin{verbatim}
[(z)->exp(Pi*z/2),[1,-z,2*n-1],[[2*z,z^2+4],
                               [z^2+(2*n-1)^2,z^2+4*(n+1)^2]]]
\end{verbatim}
$$e^{\pi z/2}=1+\dfrac{2z}{-z+\dfrac{z^2+4}{3+\dfrac{z^2+1}{5+\dfrac{z^2+16}{7+\dfrac{z^2+9}{9+\dfrac{z^2+36}{11+\ddots}}}}}}$$
Convergence type $E$ with $E=-(1+\sqrt{2})^2$, $P=0$, and
$C=2e^{\pi z/2}\sinh(\pi z)/(1+\sqrt{2})$, so that
$$e^{\pi z/2}-\dfrac{p(n)}{q(n)}\sim(-1)^n\dfrac{2e^{\pi z/2}\sinh(\pi z)}{(1+\sqrt{2})^{2n+1}}\;.$$
$$A=1+((-z^2/2-1/8)d)/n+((z^2/4+1/16)d+(z^4/4+z^2/8+1/64))/n^2+\cdots$$
\end{cf}

\smallskip

\begin{cf}\label{6.3.12}{\ }
\begin{verbatim}
[(z)->exp(Pi*z/(3*sqrt(3))),[-1,-z,4*n-6],[-2*z,12*(n-1)^2+z^2]]
\end{verbatim}
$$e^{\pi z/(3\sqrt{3})}=-1-\dfrac{2z}{-z+\dfrac{z^2}{2+\dfrac{z^2+12}{6+\dfrac{z^2+48}{10+\dfrac{z^2+108}{14+\dfrac{z^2+192}{18+\ddots}}}}}}$$
Convergence type $E$ with $E=-3$, $P=0$, and $C=...$, so that
$$e^{\pi z/(3\sqrt{3})}-\dfrac{p(n)}{q(n)}\sim(-1)^n\dfrac{C}{3^n}\;.$$
$$A=1-((z^2+3)/24)/n+((z^2+3)(z^2-21)/1152)/n^2+\cdots$$
Parametric family for $k\ge0$:
\begin{verbatim}
[(z)->exp(Pi*z/(3*sqrt(3))),4*n-6+8*k,12*(n-1)^2+z^2]
\end{verbatim}
Convergence type $E$ with $E=-3$ and $P=2k$.
\end{cf}

\smallskip

\begin{cf}\label{6.3.13}{\ }
\begin{verbatim}
[(z)->exp(Pi*z/(3*sqrt(3))),[-1,-z-1,4*n-8],[-2*z,3*(2*n-1)^2+z^2]]
\end{verbatim}
$$e^{\pi z/(3\sqrt{3})}=-1-\dfrac{2z}{-z-1+\dfrac{z^2+3}{0+\dfrac{z^2+27}{4+\dfrac{z^2+75}{8+\dfrac{z^2+147}{12+\dfrac{z^2+243}{16+\ddots}}}}}}$$
Convergence type $E$ with $E=-3$, $P=-1$, and $C=...$, so that
$$e^{\pi z/(3\sqrt{3})}-\dfrac{p(n)}{q(n)}\sim(-1)^n\dfrac{C}{3^nn^{-1}}\;.$$
$$A=1-((z^2+6)/24)/n+((z^2+12)(z^2+18)/1152)/n^2+\cdots$$
Parametric family for $k\ge0$:
\begin{verbatim}
[(z)->exp(Pi*z/(3*sqrt(3))),4*n-8+8*k,3*(2*n-1)^2+z^2]
\end{verbatim}
Convergence type $E$ with $E=-3$ and $P=2k-1$.
\end{cf}

\smallskip

\begin{cf}\label{6.3.14}{\ }
\begin{verbatim}
[(z)->exp(Pi*z/(3*sqrt(3))),[1,3-z,3*(2*n-1)],[2*z,z^2+3*n^2]]
\end{verbatim}
$$e^{\pi z/(3\sqrt{3})}=1+\dfrac{2z}{-z+3+\dfrac{z^2+3}{9+\dfrac{z^2+12}{15+\dfrac{z^2+27}{21+\dfrac{z^2+48}{27+\dfrac{z^2+75}{33+\ddots}}}}}}$$
Convergence type $E$ with $E=-(2+\sqrt{3})^2$, $P=0$, and $C=2e^{\pi z/(3\sqrt{3})}\sinh(\pi z/\sqrt{3})/(2+\sqrt{3})$, so
that
$$e^{\pi z/(3\sqrt{3})}-\dfrac{p(n)}{q(n)}\sim(-1)^n\dfrac{2e^{\pi z/(3\sqrt{3})}\sinh(\pi z/\sqrt{3})}{(2+\sqrt{3})^{2n+1}}\;.$$
$$A=1-((4z^2+3)d/24)/n+((4z^2+3)(8d+4z^2+3)/384)/n^2+\cdots$$
\end{cf}

More generally:

\smallskip

\begin{cf}\label{6.3.15}{\ }
\begin{verbatim}
[(a,z)->exp(2*z*atan(1/a)),[-1,-a*z-1,2*(a^2-1)*n-3*a^2+1],
                           [-2*a*z,a^2*((2*n-1)^2+z^2)]]
\end{verbatim}
$$e^{2z\atan(1/a)}=-1-\dfrac{2az}{-az-1+\dfrac{a^2z^2+a^2}{a^2-3+\dfrac{a^2z^2+9a^2}{3a^2-5+\dfrac{a^2z^2+25a^2}{5a^2-7+\ddots}}}}$$
Convergence type $E$ with $E=-a^2$, $P=-1$, and $C=...$, so that
$$e^{2z\atan(1/a)}-\dfrac{p(n)}{q(n)}\sim(-1)^n\dfrac{C}{a^{2n}n^{-1}}\;.$$
$$A=1+(((-a^2+1)/(4a^2+4))z^2+((-a^2+1)/(2a^2+2)))/n+\cdots$$
Parametric family for $k\ge0$:
\begin{verbatim}
[(a,z)->exp(2*z*atan(1/a)),2*(a^2-1)*n-3*a^2+1+2*k*(a^2+1),
                           a^2*((2*n-1)^2+z^2)]
\end{verbatim}
Convergence type $E$ with $E=-a^2$ and $P=2k-1$.
\end{cf}

\smallskip

\begin{cf}\label{6.3.16}{\ }
\begin{verbatim}
[(a,z)->exp(2*z*atan(1/a)),[-1,-a*z,(a^2-1)*(2*n-3)],
                           [-2*a*z,a^2*(4*(n-1)^2+z^2)]]
\end{verbatim}
$$e^{2z\atan(1/a)}=-1-\dfrac{2az}{-az+\dfrac{a^2z^2}{a^2-1+\dfrac{a^2z^2+4a^2}{3a^2-3+\dfrac{a^2z^2+16a^2}{5a^2-5+\ddots}}}}$$
Convergence type $E$ with $E=-a^2$, $P=0$, and $C=...$, so that
$$e^{2z\atan(1/a)}-\dfrac{p(n)}{q(n)}\sim(-1)^n\dfrac{C}{a^{2n}}\;.$$
$$A=1-(((a^2-1)/(4(a^2+1)))(z^2+1))/n+\cdots$$
Parametric family for $k\ge0$:
\begin{verbatim}
[(a,z)->exp(2*z*atan(1/a)),(a^2-1)*(2*n-3)+2*k*(a^2+1),
                           a^2*(4*(n-1)^2+z^2)]
\end{verbatim}
Convergence type $E$ with $E=-a^2$ and $P=2k$.
\end{cf}

\smallskip

\begin{cf}\label{6.3.17}{\ }
\begin{verbatim}
[(a,z)->exp(2*z*atan(1/a)),[1,a-z,(2*n-1)*a],[2*z,z^2+n^2]]
\end{verbatim}
$$e^{2z\atan(1/a)}=1+\dfrac{2z}{a-z+\dfrac{z^2+1}{3a+\dfrac{z^2+4}{5a+\dfrac{z^2+9}{7a+\dfrac{z^2+16}{9a+\dfrac{z^2+25}{11a+\ddots}}}}}}$$
Convergence type $E$ with $E=-(a+\sqrt{a^2+1})^2$, $P=0$, and $C=...$, so
that
$$e^{2z\atan(1/a)}-\dfrac{p(n)}{q(n)}\sim(-1)^n\dfrac{C}{(a+\sqrt{a^2+1})^{2n}}\;.$$
\end{cf}

\smallskip

\begin{cf}\label{6.3.18}{\ }
\begin{verbatim}
[(z)->2^z,[[1,1-z],[1-z,-1-z]],
        [[z,1-z],[-2*(n-z)*(2*n-z),-(2*n+1)*(2*n+1-z)]]]
\end{verbatim}
$$2^z=1+\dfrac{z}{-z+1-\dfrac{z-1}{-z+1-\dfrac{2z^2-6z+4}{-z-1+\dfrac{3z-9}{-z+1-\dfrac{2z^2-12z+16}{-z-1+\dfrac{5z-25}{-z+1-\ddots}}}}}}$$
Convergence type $P^i$ with $P=1$ and $C=z(1-z)2^{z-1}$, so that
$$2^z-\dfrac{p(n)}{q(n)}\sim(-1)^{\lfloor n/2\rfloor}\dfrac{z(1-z)2^{z-1}}{n}\;.$$
$$A=1+(z(z-1)/2)/n+(z(z-1)(z^2+3z-2)/8)/n^2+(z^2(z-1)^2(z^2+11z+22)/48)/n^3+\cdots$$
Series:
\begin{align*}2^z&=1-\dfrac{z}{z-2}\sum_{n\ge0}\dfrac{(1/2)_n(1-z)_n}{((3-z)/2)_n((4-z)/2)_n}\\
  2^{-z}&=\dfrac{2}{z+1}+\dfrac{z-1}{z+1}\sum_{n\ge0}\dfrac{(-z/2)_n(-(z+1)/2)_n}{(1/2)_n(1-z)_n}\end{align*}
\end{cf}

\smallskip

\begin{cf}\label{6.3.19}{\ }
\begin{verbatim}
[(z)->2^z,[1,2*z+2,8*n^2+(8*z-12)*n+(z-1)*(z-6)],
          [z*(z+1),-2*(2*n-1)*(n+z)*(2*n+z-1)*(2*n+z)]]
\end{verbatim}
$$2^z=1+\dfrac{z^2+z}{2z+2-\dfrac{2z^3+8z^2+10z+4}{z^2+9z+14-\dfrac{6z^3+54z^2+156z+144}{z^2+17z+42-\ddots}}}$$
Convergence type $P^+$ with $P=1$ and $C=z(z+1)2^{z-2}$, so that
$$2^z-\dfrac{p(n)}{q(n)}\sim\dfrac{z(z+1)2^{z-2}}{n}$$
$A=1+\cdots$

\noindent
Series:
\begin{align*}2^z&=1+\dfrac{z}{2}\sum_{n\ge0}\dfrac{((z+1)/2)_n(z/2+1)_n}{(3/2)_n(z+2)_n}\\
2^{-z}&=\dfrac{1-z}{2}+\dfrac{(1+z)}{2}\sum_{n\ge0}\dfrac{(z)_n(-1/2)_n}{((z+1)/2)_n(z/2+1)_n}\end{align*}
Parametric family for $k\ge0$:
\begin{verbatim}
[(z)->2^z,8*n^2+(8*z-12)*n+(z-1)*(z-6)+4*k*(k+1),
          -2*(2*n-1)*(n+z)*(2*n+z-1)*(2*n+z)]
\end{verbatim}
Convergence type $P^+$ with $P=2k+1$.
\end{cf}

\smallskip

\begin{cf}\label{6.3.20}{\ }
\begin{verbatim}
[(z)->2^z,[1,-z^2+5*z-2,(4*z-2)*(2*n-1)],
            [2*z*(z-1),(4*n^2-z^2)*(4*n^2-(z-1)^2)]]
\end{verbatim}
$$2^z=1+\dfrac{2z^2-2z}{-z^2+5z-2+\dfrac{z^4-2z^3-7z^2+8z+12}{12z-6+\dfrac{z^4-2z^3-31z^2+32z+240}{20z-10+\ddots}}}$$
Convergence type $P^-$ with $P=|1-2z|$ and $C=-\G(z)/(2\G(1-z))$, so that
$$2^z-\dfrac{p(n)}{q(n)}\sim(-1)^{n+1}\dfrac{\G(z)/(2\G(1-z))}{n^{|1-2z|}}\;.$$
$$A=1+(1/2-z)/n+(z^3/12+3z^2/8-17z/24+1/4)/n^2+\cdots$$
This is for $z<1/2$. When $z>1/2$, the limit is $2^{z-1}$ instead.
Parametric family for $k\ge0$:
\begin{verbatim}
[(z)->2^z,(4*z+8*k-2)*(2*n-1),(4*n^2-z^2)*(4*n^2-(z-1)^2)]
\end{verbatim}
Convergence type $P^-$ with $P=|4k+2z-1|$.
\end{cf}

\medskip

\section{Logarithm Functions}

\medskip

\begin{cf}\label{6.4.1}{\ }
\begin{verbatim}
[(z)->log(1+z),[0,n-(n-1)*z],[z,n^2*z]]
\end{verbatim}
$$\log(1+z)=\dfrac{z}{1+\dfrac{z}{-z+2+\dfrac{4z}{-2z+3+\dfrac{9z}{-3z+4+\dfrac{16z}{-4z+5+\dfrac{25z}{-5z+6+\ddots}}}}}}$$
Convergence type $E$ with $E=-1/z$, $P=1$, and $C=z/(z+1)$, so that
$$\log(1+z)-\dfrac{p(n)}{q(n)}\sim(-1)^n\dfrac{1/(z+1)}{(1/z)^{n+1}n}\;.$$
$$A=1-(1/(1+z))/n+((1-z)/(1+z)^2)/n^2-((z^2-4z+1)/(1+z)^3)/n^3+\cdots$$
Series:
$$\log(1+z)=z\sum_{n\ge0}(-1)^n\dfrac{z^n}{n+1}$$
Parametric family for $k\ge0$:
\begin{verbatim}
[(z)->log(1+z),n-(n-1)*z+(1+z)*k,n^2*z]
\end{verbatim}
Convergence type $E$ with $E=-1/z$ and $P=2k+1$.
\end{cf}

This is the CF corresponding to the term-by-term Taylor expansion of
$\log(1+z)$.

\smallskip

\begin{cf}\label{6.4.2}{\ }
\begin{verbatim}
[(z)->log(1+z),[0,n+(2*n-1)*z],[z,-n^2*z*(1+z)]]
\end{verbatim}
$$\log(1+z)=\dfrac{z}{z+1-\dfrac{z^2+z}{3z+2-\dfrac{4z^2+4z}{5z+3-\dfrac{9z^2+9z}{7z+4-\dfrac{16z^2+16z}{9z+5-\dfrac{25z^2+25z}{11z+6-\ddots}}}}}}$$
Convergence type $E$ with $E=(z+1)/z$, $P=1$, and $C=z$, so that
$$\log(1+z)-\dfrac{p(n)}{q(n)}\sim\dfrac{z}{((z+1)/z)^nn}\;.$$
$$A=1-(z+1)/n+(z+1)(2z+1)/n^2-(z+1)(6z^2+6z+1)/n^3+\cdots$$
Series:
$$\log(1+z)=\dfrac{z}{z+1}\sum_{n\ge0}\dfrac{(z/(z+1))^n}{n+1}$$
Parametric family for $k\ge0$:
\begin{verbatim}
[(z)->log(1+z),n+(2*n-1)*z+k,-n^2*z*(1+z)]
\end{verbatim}
Convergence type $E$ with $E=(z+1)/z$ and $P=2k+1$.
\end{cf}

\smallskip

\begin{cf}\label{6.4.3}
\begin{verbatim}
[(z)->log(1+z),[0,(-z^2+4*z+4)*n+z^2-2*z-2],
                        [z*(z+2),2*z^2*(z+1)*n*(2*n-1)]]
\end{verbatim}
$$\log(1+z)=\dfrac{z^2+2z}{4z+4+\dfrac{8z^3+8z^2}{-2z^2+12z+12+\dfrac{48z^3+48z^2}{-4z^2+20z+20+\dfrac{120z^3+120z^2}{-6z^2+28z+28+\ddots}}}}$$
Convergence type $E$ with $E=-4(z+1)/z^2$, $P=1/2$, and
$C=(z/(z+2))\sqrt{\pi}$, so that
$$\log(1+z)-\dfrac{p(n)}{q(n)}\sim(-1)^n\dfrac{(z/(z+2))\sqrt{\pi}}{((z+1)/z^2)^nn^{1/2}}$$
$$A=1+((z^2-12z-12)/(8(z+2)^2))/n+\cdots$$
Series:
$$\log(1+z)=\dfrac{z(z+2)}{2(z+1)}\sum_{n\ge0}\dfrac{n!}{(3/2)_n}(-z^2/(4(z+1)))^n$$
Parametric family for $k\ge0$:
\begin{verbatim}
[(z)->log(1+z),(-z^2+4*z+4)*n+(k+1)*z^2+(4*k-2)*(z+1),
               2*z^2*(z+1)*n*(2*n-1)]
\end{verbatim}
Convergence type $E$ with $E=-4(z+1)/z^2$ and $P=2k+1/2$.
\end{cf}

\smallskip

\begin{cf}\label{6.4.4}
\begin{verbatim}
[(z)->log(1+z),[0,(z^2+2*z+2)*(2*n-1)],
                        [z*(z+2),-z^2*(z+2)^2*n^2]]
\end{verbatim}
$$\log(1+z)=\dfrac{z^2+2z}{z^2+2z+2-\dfrac{z^4+4z^3+4z^2}{3z^2+6z+6-\dfrac{4z^4+16z^3+16z^2}{5z^2+10z+10-\dfrac{9z^4+36z^3+36z^2}{7z^2+14z+14-\ddots}}}}$$
Convergence type $E$ with $E=(z+2)^2/z^2$, $P=0$, and $C=(z/(z+2))\pi$, so that
$$\log(1+z)-\dfrac{p(n)}{q(n)}\sim\dfrac{\pi}{((z+2)/z)^{2n+1}}$$
$$A=1-((z^2+2z+2)/(8(z+1)))/n+\cdots$$
Parametric family for $k\ge0$:
\begin{verbatim}
[(z)->log(1+z),2*(z^2+2*z+2)*n-z^2+(4*k-2)*(z+1),-z^2*(z+2)^2*n^2]
\end{verbatim}
Convergence type $E$ with $E=(z+2)^2/z^2$ and $P=2k$.
\end{cf}

\smallskip

\begin{cf}\label{6.4.5}
\begin{verbatim}
[(z)->log(1+z),[0,(z+2)^2,4*((z^2+2*z+2)*n-(z^2+z+1))],
                        [2*z*(z+2),-z^2*(z+2)^2*(2*n-1)^2]]
\end{verbatim}
$$\log(1+z)=\dfrac{2z^2+4z}{z^2+4z+4-\dfrac{z^4+4z^3+4z^2}{4z^2+12z+12-\dfrac{9z^4+36z^3+36z^2}{8z^2+20z+20-\dfrac{25z^4+100z^3+100z^2}{12z^2+28z+28-\ddots}}}}$$
Convergence type $E$ with $E=(z+2)^2/z^2$, $P=1$, and $C=z(z+2)/(4(z+1))$, so
that
$$\log(1+z)-\dfrac{p(n)}{q(n)}\sim\dfrac{z(z+2)/(4(z+1))}{((z+2)/z)^{2n}n}$$
$$A=1-((z^2+2z+2)/(4(z+1)))/n+\cdots$$
Series:
$$\log(1+z)=\dfrac{2z}{z+2}\sum_{n\ge0}\dfrac{(z/(z+2))^{2n}}{2n+1}$$
Parametric family for $k\ge0$:
\begin{verbatim}
[(z)->log(1+z),4*(z^2+2*z+2)*n-4*z^2+4*(2*k-1)*(z+1),-z^2*(z+2)^2*(2*n-1)^2]
\end{verbatim}
Convergence type $E$ with $E=(z+2)^2/z^2$ and $P=2k+1$.
\end{cf}

\smallskip

\begin{cf}\label{6.4.6}{\ }
\begin{verbatim}
[(z)->log(1+z),[0,(2+z)*(2*n-1)],[2*z,-n^2*z^2]]
\end{verbatim}
$$\log(1+z)=\dfrac{2z}{z+2-\dfrac{z^2}{3z+6-\dfrac{4z^2}{5z+10-\dfrac{9z^2}{7z+14-\dfrac{16z^2}{9z+18-\dfrac{25z^2}{11z+22-\ddots}}}}}}$$
Convergence type $E$ with $E=(1+\sqrt{z+1})^4/z^2$, $P=0$, and
$C=2\pi z/(1+\sqrt{z+1})^2$, so that
$$\log(1+z)-\dfrac{p(n)}{q(n)}\sim\dfrac{2\pi}{((1+\sqrt{z+1})^2/z)^{2n+1}}\;.$$
Parametric family for $u\ge0$:
\begin{verbatim}
[(z)->log(1+z),(2+z)*(2*n+2*u-1),-z^2*n*(n+2*u)]
\end{verbatim}
Convergence type $E$ with $E=(1+\sqrt{z+1})^4/z^2$ and $P=0$.
\end{cf}
  
\smallskip

\begin{cf}\label{6.4.7}{\ }
\begin{verbatim}
[(z)->log(1+z),[z,(4*z+8)*n^2-2],[-3*z^2,-z^2*n*(n+1)*(2*n-1)*(2*n+3)]]
\end{verbatim}
$$\log(1+z)=z-\dfrac{3z^2}{4z+6-\dfrac{10z^2}{16z+30-\dfrac{126z^2}{36z+70-\dfrac{540z^2}{64z+126-\dfrac{1540z^2}{100z+198-\ddots}}}}}$$
Convergence type $E$ with $E=(1+\sqrt{z+1})^4/z^2$, $P=0$, and
$C=-2\pi z^2/(1+\sqrt{z+1})^4$, so that
$$\log(1+z)-\dfrac{p(n)}{q(n)}\sim-\dfrac{2\pi}{((1+\sqrt{z+1})^2/z)^{2n+2}}\;.$$
$A=...$.
\end{cf}

\smallskip

\begin{cf}\label{6.4.8}{\ }
\begin{verbatim}
[(z)->log(1+z),[z-z^2/2,(2*n+1)*((n^2+n+1)*z+2*n*(n+1))],
             [4*z^3,-z^2*n^3*(n+2)^3]]
\end{verbatim}
$$\log(1+z)=z-1/2z^2+\dfrac{4z^3}{9z+12-\dfrac{27z^2}{35z+60-\dfrac{512z^2}{91z+168-\dfrac{3375z^2}{189z+360-\ddots}}}}$$
Convergence type $E$ with $E=(1+\sqrt{z+1})^4/z^2$, $P=0$, and
$C=2\pi z^3/(1+\sqrt{z+1})^6$, so that
$$\log(1+z)-\dfrac{p(n)}{q(n)}\sim\dfrac{2\pi}{((1+\sqrt{z+1})^2/z)^{2n+3}}\;.$$
$$A=1+((15z+14)/(8d))/n+(((-360z-336)d+225z^2+420z+196)/(128(z+1)))/n^2+\cdots$$
with $d=\sqrt{z+1}$.
\end{cf}

\smallskip

Generalizing the above three CFs, one can obtain CFs of a similar
form for $\log(1+z)-\sum_{1\le j\le k}(-1)^{j-1}z^j/j$, all with
convergence type $E=(1+\sqrt{z+1})^4/z^2$, $P=0$, and
$C=(-1)^k2\pi/((1+\sqrt{z+1})^2/z)^{k-1}$.

\smallskip

\begin{cf}\label{6.4.9}{\ }
\begin{verbatim}
[(z)->1/(1+z)+log(1+z)/z,[0,n^2*(1-z)+z],[2,z*n^3*(n+2)]]
\end{verbatim}
$$\dfrac{1}{1+z}+\dfrac{\log(1+z)}{z}=\dfrac{2}{1+\dfrac{3z}{-3z+4+\dfrac{32z}{-8z+9+\dfrac{135z}{-15z+16+\dfrac{384z}{-24z+25+\ddots}}}}}$$
Convergence type $E$ with $E=-1/z$, $P=0$, and $C=1/(1+z)$, so that
$$\dfrac{1}{1+z}+\dfrac{\log(1+z)}{z}-\dfrac{p(n)}{q(n)}\sim(-1)^n\dfrac{1/(1+z)}{(1/z)^n}\;.$$
$$A=1+1/n-(1/(1+z))/n^2+((1-z)/(1+z)^2)/n^3-((z^2-4z+1)/(1+z)^3)/n^4+\cdots$$
Series:
$$\dfrac{1}{1+z}+\dfrac{\log(1+z)}{z}=\sum_{n\ge0}(-1)^n\dfrac{n+2}{n+1}z^n$$
\end{cf}

\smallskip

\begin{cf}\label{6.4.10}{\ }
\begin{verbatim}
[(z)->(1+z)*log(1+z)/z,[n+1],[[z,2*z],(n+1)*z*[n,n+2]]]
\end{verbatim}
$$\dfrac{1+z}{z}\log(1+z)=1+\dfrac{z}{2+\dfrac{2z}{3+\dfrac{2z}{4+\dfrac{6z}{5+\dfrac{6z}{6+\dfrac{12z}{7+\ddots}}}}}}$$
Convergence type $E$ with $E=-(1+\sqrt{z+1})^2/z$, $P=0$, and
$C=2\pi z(1+z)/(1+\sqrt{1+z})^4$, so that
$$\dfrac{1+z}{z}\log(1+z)-\dfrac{p(n)}{q(n)}\sim(-1)^n\dfrac{2\pi(1+z)/z}{((1+\sqrt{z+1})^2/z)^{n+2}}\;.$$
\end{cf}

\smallskip

\begin{cf}\label{6.4.11}{\ }
\begin{verbatim}
[(z)->(1+z)*log(1+z)/z,[1,(4*z+8)*n^2-2*z-2],
                     [3*z,-n*(n+1)*(2*n-1)*(2*n+3)*z^2]]
\end{verbatim}
$$\dfrac{1+z}{z}\log(1+z)=1+\dfrac{3z}{2z+6-\dfrac{10z^2}{14z+30-\dfrac{126z^2}{34z+70-\dfrac{540z^2}{62z+126-\ddots}}}}$$
Convergence type $E$ with $E=(1+\sqrt{1+z})^4/z^2$, $P=0$, and
$C=2\pi z(1+z)/(1+\sqrt{1+z})^4$, so that
$$\dfrac{1+z}{z}\log(1+z)-\dfrac{p(n)}{q(n)}\sim\dfrac{2\pi(1+z)/z}{((1+\sqrt{1+z})^2/z)^{2n+2}}\;.$$
\end{cf}

This is simply the contraction of the previous CF.

\smallskip

\begin{cf}\label{6.4.12}{\ }
\begin{verbatim}
[(z)->(1+z)*log(1+z)/z,[z/2+1,(z+2)*(2*n+1)],-z^2*[1,n*(n+2)]]
\end{verbatim}
$$\dfrac{1+z}{z}\log(1+z)=z/2+1-\dfrac{z^2}{3z+6-\dfrac{3z^2}{5z+10-\dfrac{8z^2}{7z+14-\dfrac{15z^2}{9z+18-\dfrac{24z^2}{11z+22-\ddots}}}}}$$
Convergence type $E$ with $E=(1+\sqrt{z+1})^4/z^2$, $P=0$, and
$C=-2\pi z^2(1+z)/(1+\sqrt{1+z})^6$, so that
$$\dfrac{1+z}{z}\log(1+z)-\dfrac{p(n)}{q(n)}\sim-\dfrac{2\pi(1+z)/z}{((1+\sqrt{z+1})^2/z)^{2n+3}}\;.$$
$A=...$
\end{cf}

Note that this is a special case of the parametric family given in
\ref{6.4.6}, so it of course also has the same parametric family.

\smallskip

Continued fractions for $\log((1+z)/(1-z))$ can be found below since
$$\log((1+z)/(1-z))=2\atanh(z)\;.$$

\smallskip

\begin{cf}\label{6.4.13}{\ }
\begin{verbatim}
[(z)->dilog(z),[0,n^2+z*(n-1)^2],[z,-z*n^4]]
\end{verbatim}
$$\Li_2(z)=\dfrac{z}{1-\dfrac{z}{z+4-\dfrac{16z}{4z+9-\dfrac{81z}{9z+16-\dfrac{256z}{16z+25-\dfrac{625z}{25z+36-\ddots}}}}}}$$
Convergence type $E$ with $E=1/z$, $P=2$, and $C=z/(1-z)$, so that
$$\Li_2(z)-\dfrac{p(n)}{q(n)}\sim\dfrac{1/(1-z)}{(1/z)^{n+1}n^2}\;.$$
$$A=1-(2/(1-z))/n+(3(1+z)/(1-z)^2)/n^2-(4(z^2+4z+1)/(1-z)^3)/n^3+\cdots$$
Series:
$$\Li_2(z)=z\sum_{n\ge0}\dfrac{z^n}{(n+1)^2}$$
\end{cf}

\smallskip

\begin{cf}\label{6.4.14}{\ }
\begin{verbatim}
[(z)->polylog(3,z),[0,n^3+z*(n-1)^3],[z,-z*n^6]]
\end{verbatim}
$$\Li_3(z)=\dfrac{z}{1-\dfrac{z}{z+8-\dfrac{64z}{8z+27-\dfrac{729z}{27z+64-\dfrac{4096z}{64z+125-\dfrac{15625z}{125z+216-\ddots}}}}}}$$
Convergence type $E$ with $E=1/z$, $P=3$, and $C=z/(1-z)$, so that
$$\Li_3(z)-\dfrac{p(n)}{q(n)}\sim\dfrac{1/(1-z)}{(1/z)^{n+1}n^3}\;.$$
$$A=1-(3/(1-z))/n+(6(1+z)/(1-z)^2)/n^2-(10(z^2+4z+1)/(1-z)^3)/n^3+\cdots$$
Series:
$$\Li_3(z)=z\sum_{n\ge0}\dfrac{z^n}{(n+1)^3}$$
\end{cf}

\smallskip

These two CFs are simply the CFs corresponding to the Taylor expansions of
$\Li_2(z)$ and $\Li_3(z)$.

\smallskip

\begin{cf}\label{6.4.15}{\ }
\begin{verbatim}
[(a,z)->suminf(k=0,z^k/(k+a)),[0,a,(z+1)*(n+a)-(2*z+1)],
                              [1,-z*(n+a-1)^2]]
\end{verbatim}
$$\sum_{k\ge0}\dfrac{z^k}{k+a}=\dfrac{1}{a-\dfrac{za^2}{(z+1)a+1-\dfrac{za^2+2za+z}{(z+1)a+(z+2)-\dfrac{za^2+4za+4z}{(z+1)a+(2z+3)-\ddots}}}}$$
Convergence type $E$ with $E=1/z$, $P=1$, and $C=1/(1-z)$, so that
$$\sum_{k\ge0}\dfrac{z^k}{k+a}-\dfrac{p(n)}{q(n)}\sim\dfrac{1/(1-z)}{(1/z)^nn}\;.$$
$$A=1-(((1-a)z+a)/(1-z))/n+(((1-a)^2z^2-(2a^2-2a-1)z+a^2)/(1-z)^2)/n^2+\cdots$$
Series:
$$\sum_{k\ge0}\dfrac{z^k}{k+a}=\sum_{n\ge0}\dfrac{z^n}{n+a}$$
Parametric family for $k\ge0$:
\begin{verbatim}
[(a,z)->suminf(k=0,z^k/(k+a)),(z+1)*(n+a)-(2*z+1)+k*(1-z),
                              -z*(n+a-1)^2]
\end{verbatim}
Convergence type $E$ with $E=1/z$ and $P=2k+1$.
\end{cf}

\smallskip

\begin{cf}\label{6.4.16}{\ }
\begin{verbatim}
[(a,z)->suminf(k=0,z^k/(k+a)),[0,(z+1)*n+(1-z)*a-1],[1,-z*n^2]]
\end{verbatim}
$$\sum_{k\ge0}\dfrac{z^k}{k+a}=\dfrac{1}{(-a+1)z+a-\dfrac{z}{(-a+2)z+(a+1)-\dfrac{4z}{(-a+3)z+(a+2)-\ddots}}}$$
Convergence type $E$ with $E=1/z$, $P=2a-1$, and $C=\G(a)^2/(1-z)^{2a-1}$, so
that
$$\sum_{k\ge0}\dfrac{z^k}{k+a}-\dfrac{p(n)}{q(n)}\sim\dfrac{\G(a)^2/(1-z)^{2a-1}}{(1/z)^nn^{2a-1}}\;.$$
$$A=1+((-(a-1)^2z-a^2)/(z-1))/n+\cdots$$
Parametric family for $k\ge0$:
\begin{verbatim}
[(a,z)->suminf(k=0,z^k/(k+a)),(z+1)*n+(1-z)*a-1+k*(1-z),-z*n^2]
\end{verbatim}
Convergence type $E$ with $E=1/z$ and $P=2k+2a-1$.
\end{cf}
      
\smallskip

\begin{cf}\label{6.4.17}{\ }
\begin{verbatim}
[(a,z)->suminf(k=0,z^k/(k+a)),[0,n+a-1],[[1,-a^2*z],-z*[n^2,(n+a)^2]]]
\end{verbatim}
$$\sum_{k\ge0}\dfrac{z^k}{k+a}=\dfrac{1}{a-\dfrac{a^2z}{a+1-\dfrac{z}{a+2-\dfrac{(a^2+2a+1)z}{a+3-\dfrac{4z}{a+4-\dfrac{(a^2+4a+4)z}{a+5-\ddots}}}}}}$$
Convergence type $E$ with $E=(1+\sqrt{1-z})^2/z$, $P=0$, and
$C=2\pi/(1+\sqrt{1-z})^{2a}$, so that
$$\sum_{k\ge0}\dfrac{z^k}{k+a}-\dfrac{p(n)}{q(n)}\sim\dfrac{2\pi}{(1+\sqrt{1-z})^{2n+2a}/z^n}\;.$$
\end{cf}

\smallskip

\begin{cf}\label{6.4.18}{\ }
\begin{verbatim}
[(a,z)->suminf(k=0,z^k/(k+a)),[0,(1-a)*z+a,n+a-1],
                              [[1,-z],-z*[(n+a-1)^2,(n+1)^2]]]
\end{verbatim}
$$\sum_{k\ge0}\dfrac{z^k}{k+a}=\dfrac{1}{(-a+1)z+a-\dfrac{z}{a+1-\dfrac{a^2z}{a+2-\dfrac{4z}{a+3-\dfrac{(a^2+2a+1)z}{a+4-\dfrac{9z}{a+5-\ddots}}}}}}$$
Convergence type $E$ with $E=(1+\sqrt{1-z})^2/z$, $P=0$, and
$C=2\pi/(1+\sqrt{1-z})^{2a}$, so that
$$\sum_{k\ge0}\dfrac{z^k}{k+a}-\dfrac{p(n)}{q(n)}\sim\dfrac{2\pi}{(1+\sqrt{1-z})^{2n+2a}/z^n}\;.$$
\end{cf}

\smallskip

Note that
$$\sum_{k\ge0}\dfrac{(-1)^k}{k+a}=\dfrac{1}{2}(\psi((a+1)/2)-\psi(a/2))\;,$$
so the above CFs give CFs for the RHS, and conversely the CFs
applicable to the RHS give CFs for the LHS, see \ref{6.10.3} and
\ref{6.10.21}.

\medskip

\section{Power Functions}

\medskip

Note that since $(1+z)^a=1/(1+z)^{-a}$, one can obtain more continued
fractions for $(1+z)^a$ by changing $a$ into $-a$ and inverting.

\smallskip

\begin{cf}\label{6.5.1}{\ }
\begin{verbatim}
[(a,z)->(1+z)^a,[1,z+2-a*z,(z+2)*(2*n-1)],[2*a*z,z^2*(a^2-n^2)]]
\end{verbatim}
$$(1+z)^a=1+\dfrac{2az}{(-a+1)z+2+\dfrac{(a^2-1)z^2}{3z+6+\dfrac{(a^2-4)z^2}{5z+10+\dfrac{(a^2-9)z^2}{7z+14+\dfrac{(a^2-16)z^2}{9z+18+\ddots}}}}}$$
Convergence type $E$ with $E=(1+\sqrt{1+z})^4/z^2$, $P=0$, and
$C=2z(1+z)^a\sin(\pi a)/(1+\sqrt{1+z})^2$, so that
$$(1+z)^a-\dfrac{p(n)}{q(n)}\sim\dfrac{2(1+z)^a\sin(\pi a)}{((1+\sqrt{1+z})^2/z)^{2n+1}}\;.$$
\end{cf}

\smallskip

\begin{cf}\label{6.5.2}{\ }
\begin{verbatim}
[(a,z)->(1+z)^a,[[0,1],[4*n-2,2*n+1]],
                [[1,-2*a*z],z*(2*n+1)*[n+a,n-a]]]
\end{verbatim}
$$(1+z)^a=\dfrac{1}{1-\dfrac{2za}{2+\dfrac{3za+3z}{3-\dfrac{3za-3z}{6+\dfrac{5za+10z}{5-\dfrac{5za-10z}{10+\ddots}}}}}}$$
Convergence type $E$ with $E=(1+\sqrt{1+z})^2/z$, $P=0$, and
$C=-2(1+z)^a\sin(\pi a)$, so that
$$(1+z)^a-\dfrac{p(n)}{q(n)}\sim-\dfrac{2(1+z)^a\sin(\pi a)}{((1+\sqrt{1+z})^2/z)^n}\;.$$
\end{cf}

\smallskip

\begin{cf}\label{6.5.3}{\ }
\begin{verbatim}
[(a,z)->(1+z)^a,[[0,z+1],[4*n-2,(2*n+1)*(z+1)]],
                [[z+1,-2*a*z],-z*(2*n+1)*[n-a,n+a]]]
\end{verbatim}
$$(1+z)^a=\dfrac{z+1}{z+1-\dfrac{2za}{2+\dfrac{3za-3z}{3z+3-\dfrac{3za+3z}{6+\dfrac{5za-10z}{5z+5-\dfrac{5za+10z}{10+\ddots}}}}}}$$
Convergence type $E$ with $E=(1+\sqrt{1+z})^2/z$, $P=0$, and
$C=2(1+z)^a\sin(\pi a)$, so that
$$(1+z)^a-\dfrac{p(n)}{q(n)}\sim\dfrac{2(1+z)^a\sin(\pi a)}{((1+\sqrt{1+z})^2/z)^n}\;.$$
\end{cf}

\smallskip

\begin{cf}\label{6.5.4}{\ }
\begin{verbatim}
[(a,z)->(1+z)^a,[0,-z/a+1/a^2,n],[[1/a^2,(1+1/a)*z^2],
                                z*[(n+1)*(n+a+1),n*(n-a)]]]
\end{verbatim}
$$(1+z)^a=\dfrac{1/a^2}{-z/a+1/a^2+\dfrac{(1+1/a)z^2}{2+\dfrac{(2a+4)z}{3-\dfrac{(a-1)z}{4+\dfrac{(3a+9)z}{5-\dfrac{(2a-4)z}{6+\ddots}}}}}}$$
Convergence type $E$ with $E=(1+\sqrt{1+z})^2/z$, $P=0$, and
$C=2z(1+z)^a\sin(\pi a)/(1+\sqrt{1+z})^2$, so that
$$(1+z)^a-\dfrac{p(n)}{q(n)}\sim\dfrac{2(1+z)^a\sin(\pi a)}{((1+\sqrt{1+z})^2/z)^{n+1}}\;.$$
\end{cf}

\smallskip

\begin{cf}\label{6.5.5}{\ }
\begin{verbatim}
[(a,z)->(1+z)^a,[0,(a-1)*z-1,4*(z+2)*n*(n-2)+2*(a+1)*z+6],
                [-z-1,-(n-a)*(n+a-1)*(2*n-3)*(2*n+1)*z^2]]
\end{verbatim}
$$(1+z)^a=-\dfrac{z+1}{(a-1)z-1-\dfrac{(3a^2-3a)z^2}{(2a+2)z+6+\dfrac{(5a^2-5a-10)z^2}{(2a+14)z+30+\dfrac{(21a^2-21a-126)z^2}{(2a+34)z+70+\ddots}}}}$$
Convergence type $E$ with $E=(1+\sqrt{1+z})^4/z^2$, $P=0$, and $C=...$, so that
$$(1+z)^a-\dfrac{p(n)}{q(n)}\sim\dfrac{C}{((1+\sqrt{1+z})^2/z)^{2n}}\;.$$
\end{cf}
      
Special cases:

\smallskip

\begin{cf}\label{6.5.6}{\ }
\begin{verbatim}
[(z)->(1+z)^(3/2),[1+3*z/2,2*n+2],[3*z^2/2,n*(n+3)*z]]
\end{verbatim}
$$(1+z)^{3/2}=3z/2+1+\dfrac{3z^2/2}{4+\dfrac{4z}{6+\dfrac{10z}{8+\dfrac{18z}{10+\dfrac{28z}{12+\dfrac{40z}{14+\ddots}}}}}}$$
Convergence type $E$ with $E=-(1+\sqrt{1+z})^2/z$, $P=0$, and
$C=2z^2(1+z)^{3/2}/(1+\sqrt{1+z})^4$, so that
$$(1+z)^{3/2}-\dfrac{p(n)}{q(n)}\sim(-1)^n\dfrac{2(1+z)^{3/2}}{((1+\sqrt{1+z})^2/z)^{n+2}}\;.$$
\end{cf}

\smallskip

\begin{cf}\label{6.5.7}{\ }
\begin{verbatim}
[(z)->(1+z)^(1/2),[0,z+2,2*z+4],[2*z+2,-z^2]]
\end{verbatim}
$$(1+z)^{1/2}=\dfrac{2z+2}{z+2-\dfrac{z^2}{2z+4-\dfrac{z^2}{2z+4-\dfrac{z^2}{2z+4-\dfrac{z^2}{2z+4-\dfrac{z^2}{2z+4-\ddots}}}}}}$$
Convergence type $E$ with $E=(1+\sqrt{1+z})^4/z^2$, $P=0$, and
$C=2\sqrt{1+z}$, so that
$$(1+z)^{1/2}-\dfrac{p(n)}{q(n)}\sim\dfrac{2\sqrt{1+z}}{((1+\sqrt{1+z})^2/z)^{2n}}\;.$$
$$A=1$$
\end{cf}

Of course this is simply a periodic CF for a quadratic irrationality.
        
\smallskip

\begin{cf}\label{6.5.8}{\ }
\begin{verbatim}
[(a,z)->(1+z)^a,[0,1,n-1+(2*n-3+a)*z],
                [1,-a*z,-(n-1)*(n-1+a)*z*(1+z)]]
\end{verbatim}
$$(1+z)^a=\dfrac{1}{1-\dfrac{za}{za+z+1-\dfrac{(z^2+z)a+z^2+z}{za+3z+2-\dfrac{(2z^2+2z)a+4z^2+4z}{za+5z+3-\dfrac{(3z^2+3z)a+9z^2+9z}{za+7z+4-\ddots}}}}}$$
Convergence type $E$ with $E=(z+1)/z$, $P=1-a$, and $C=(z+1)/\G(a)$, so that
$$(1+z)^a-\dfrac{p(n)}{q(n)}\sim\dfrac{(z+1)/\G(a)}{((z+1)/z)^nn^{1-a}}\;.$$
$$A=1+((a-1)(2z+a)/2)/n+((a-2)(a-1)(24z^2+12(a+1)z+3a^2-a)/24)/n^2+\cdots$$
Series:
$$(1+z)^a=\sum_{n\ge0}\dfrac{(a)_n}{n!}(z/(z+1))^n$$
Parametric family for $k\ge0$:
\begin{verbatim}
[(a,z)->(1+z)^a,n-1+(2*n-3+a)*z+k,-(n-1)*(n-1+a)*z*(1+z)]
\end{verbatim}
Convergence type $E$ with $E=(z+1)/z$ and $P=2k+1-a$.
\end{cf}

\smallskip

\begin{cf}\label{6.5.9}{\ }
\begin{verbatim}
[(a,z)->(1+z)^a,[0,1,(n-1)-(n-2-a)*z],[1,-a*z,(n-1)*(n-1-a)*z]]
\end{verbatim}
$$(1+z)^a=\dfrac{1}{1-\dfrac{za}{za+1-\dfrac{za-z}{za-z+2-\dfrac{2za-4z}{za-2z+3-\dfrac{3za-9z}{za-3z+4-\ddots}}}}}$$
Convergence type $E$ with $E=-1/z$, $P=a+1$, and $C=1/((z+1)\G(-a))$, so that
$$(1+z)^a-\dfrac{p(n)}{q(n)}\sim(-1)^n\dfrac{1/((z+1)\G(-a))}{(1/z)^nn^{a+1}}\;.$$
$$A=1+((a+1)((a+2)z+a)/(2(z+1)))/n+\cdots$$
Series:
$$(1+z)^a=\sum_{n\ge0}(-1)^n\dfrac{(-a)_n}{n!}z^n$$
Parametric family for $k\ge0$:
\begin{verbatim}
[(a,z)->(1+z)^a,(n-1)-(n-2-a)*z+k*(z+1),(n-1)*(n-1-a)*z]
\end{verbatim}
Convergence type $E$ with $E=-1/z$ and $P=2k+a+1$.
\end{cf}

This is the CF corresponding to the term-by-term Taylor expansion of
$(1+z)^a$.

\smallskip

For the next CFs for $((1+z)/(1-z))^a$ it is understood that $|z|<1$, but to
obtain CFs for $((z+1)/(z-1))^a$, simply change $z$ into $1/z$ everywhere.

\smallskip

\begin{cf}\label{6.5.10}{\ }
\begin{verbatim}
[(a,z)->((1+z)/(1-z))^a,[-1,z-a,2*(z^2+1)*n-(z^2+3)],
                        [-2*a,(a^2-1)*z,-((2*n-1)^2-a^2)*z^2]]
\end{verbatim}
$$((1+z)/(1-z))^a=-1-\dfrac{2a}{z-a+\dfrac{(a^2-1)z}{3z^2+1+\dfrac{(a^2-9)z^2}{5z^2+3+\dfrac{(a^2-25)z^2}{7z^2+5+\ddots}}}}$$
Convergence type $E$ with $E=1/z^2$, $P=-1$, and
$C=-4((1-z^2)/(az))\cos(\pi a/2)((1+z)/(1-z))^a$, so that
$$((1+z)/(1-z))^a-\dfrac{p(n)}{q(n)}\sim-\dfrac{4((1-z^2)/(az))\cos(\pi a/2)((1+z)/(1-z))^a}{(1/z)^{2n}n^{-1}}\;.$$
$$A=1+(((a^2-2)z^2+(a^2-2))/(4(1-z^2)))/n+\cdots$$
Parametric family for $k\ge0$:
\begin{verbatim}
[(a,z)->((1+z)/(1-z))^a,2*(z^2+1)*n-(z^2+3)+2*k*(1-z^2),
                        -((2*n-1)^2-a^2)*z^2]
\end{verbatim}
Convergence type $E$ with $E=1/z^2$ and $P=2k-1$.
\end{cf}

\smallskip

\begin{cf}\label{6.5.11}{\ }
\begin{verbatim}
[(a,z)->((1+z)/(1-z))^a,[1,z^2-a*z+1,(z^2+1)*(2*n-1)],
                        [2*a*z,-(4*n^2-a^2)*z^2]]
\end{verbatim}
$$((1+z)/(1-z))^a=1+\dfrac{2az}{z^2-az+1+\dfrac{(a^2-4)z^2}{3z^2+3+\dfrac{(a^2-16)z^2}{5z^2+5+\dfrac{(a^2-36)z^2}{7z^2+7+\ddots}}}}$$
Convergence type $E$ with $E=1/z^2$, $P=0$, and
$C=2z\sin(\pi a/2)((1+z)/(1-z))^a$, so that
$$((1+z)/(1-z))^a-\dfrac{p(n)}{q(n)}\sim\dfrac{2z\sin(\pi a/2)((1+z)/(1-z))^a}{(1/z)^{2n}}\;.$$
$$A=1+((z^2+1)/(4(1-z^2))(a^2-1))/n+\cdots$$
Parametric family for $k\ge0$:
\begin{verbatim}
[(a,z)->((1+z)/(1-z))^a,(z^2+1)*(2*n-1)-2*k*(1-z^2),-(4*n^2-a^2)*z^2]
\end{verbatim}
Convergence type $E$ with $E=1/z^2$ and $P=-2k$.
\end{cf}

\smallskip

\begin{cf}\label{6.5.12}{\ }
\begin{verbatim}
[(a,z)->((1+z)/(1-z))^a,[1,1-a*z,2*n-1],[2*a*z,(a^2-n^2)*z^2]]
\end{verbatim}
$$((1+z)/(1-z))^a=1+\dfrac{2az}{-az+1+\dfrac{(a^2-1)z^2}{3+\dfrac{(a^2-4)z^2}{5+\dfrac{(a^2-9)z^2}{7+\dfrac{(a^2-16)z^2}{9+\ddots}}}}}$$
Convergence type $E$ with $E=((1+\sqrt{1-z^2})/z)^2$, $P=0$, and
$C=2\sin(\pi a)((1+z)/(1-z))^a/((1+\sqrt{1-z^2})/z)$, so that
$$((1+z)/(1-z))^a-\dfrac{p(n)}{q(n)}\sim\dfrac{2\sin(\pi a)((1+z)/(1-z))^a}{((1+\sqrt{1-z^2})/z)^{2n+1}}\;.$$
\end{cf}

\smallskip

\begin{cf}\label{6.5.13}{\ }
\begin{verbatim}
[(z)->z^(1/z),[1,z^2+1,(2*n-1)*z*(z+1)],
              [2*(z-1),-(n^2*z^2-1)*(z-1)^2]]
\end{verbatim}
$$z^{1/z}=1+\dfrac{2z-2}{z^2+1-\dfrac{z^4-2z^3+2z-1}{3z^2+3z-\dfrac{4z^4-8z^3+3z^2+2z-1}{5z^2+5z-\dfrac{9z^4-18z^3+8z^2+2z-1}{7z^2+7z-\ddots}}}}$$
Convergence type $E$ with $E=(1+\sqrt{z})^4/(z-1)^2$, $P=0$, and $C=...$,
so that
$$z^{1/z}-\dfrac{p(n)}{q(n)}\sim\dfrac{C}{((1+\sqrt{z})^2/(z-1))^{2n}}\;.$$
$$A=1-((z-2)(z+1)(z+2)/(8z^{5/2}))/n+\cdots$$
\end{cf}

\smallskip

\begin{cf}\label{6.5.14}{\ }
\begin{verbatim}
[(z)->z^(2/z^2),[1,z^3+z^2-2*z+2,z^2*(z+1)*(2*n-1)],
                [4*z-4,-(z-1)^2*(z^4*n^2-4)]]
\end{verbatim}
$$z^{2/z^2}=1+\dfrac{4z-4}{z^3+z^2-2z+2-\dfrac{z^6-2z^5+z^4-4z^2+8z-4}{3z^3+3z^2-\dfrac{4z^6-8z^5+4z^4-4z^2+8z-4}{5z^3+5z^2-\ddots}}}$$
Convergence type $E$ with $E=(1+\sqrt{z})^4/(z-1)^2$, $P=0$, and $C=...$,
so that
$$z^{2/z^2}-\dfrac{p(n)}{q(n)}\sim\dfrac{C}{((1+\sqrt{z})^2/(z-1))^{2n}}\;.$$
$$A=1-((z-2)(z+1)(z+2)(z^2+4)/(8z^{9/2}))/n+\cdots$$
\end{cf}

\smallskip

\begin{cf}\label{6.5.15}{\ }
\begin{verbatim}
[(z)->(z/(z-1))^(z-1),[2,(z+1)*(n+1)],[2*z-4,-z*(n+2)*(n+2-z)]]
\end{verbatim}
$$(z/(z-1))^{z-1}=2+\dfrac{2z-4}{2z+2+\dfrac{3z^2-9z}{3z+3+\dfrac{4z^2-16z}{4z+4+\dfrac{5z^2-25z}{5z+5+\dfrac{6z^2-36z}{6z+6+\dfrac{7z^2-49z}{7z+7+\ddots}}}}}}$$
Convergence type $E$ with $E=z$, $P=z+1$, and $C=z^{2z-2}/((z-1)^{2z}\G(1-z))$,
so that
$$(z/(z-1))^{z-1}-\dfrac{p(n)}{q(n)}\sim\dfrac{1/((z-1)^{2z}\G(1-z))}{z^{n-2z+2}n^{z+1}}\;.$$
$$A=1-((z+1)(z^2-7z+4)/(2(z-1)))/n+\cdots$$
Series:
$$(z/(z-1))^{1-z}=\sum_{n\ge0}\dfrac{(2-z)_n}{(n+2)!}z^{-n}$$
Parametric family for $k\ge0$:
\begin{verbatim}
[(z)->(z/(z-1))^(z-1),(z+1)*(n+1)+k*(z-1),-z*(n+2)*(n+2-z)]
\end{verbatim}
Convergence type $E$ with $E=z$ and $P=2k+z+1$.
\end{cf}

\smallskip

\begin{cf}\label{6.5.16}{\ }
\begin{verbatim}
[(z)->(z/(z-1))^(z-1),[(z-1)/z,z,n*(z+1)+(z-1)],[z-1,-n*z*(n+z)]]
\end{verbatim}
$$(z/(z-1))^{z-1}=(z-1)/z+\dfrac{z-1}{z-\dfrac{z^2+z}{3z+1-\dfrac{2z^2+4z}{4z+2-\dfrac{3z^2+9z}{5z+3-\dfrac{4z^2+16z}{6z+4-\dfrac{5z^2+25z}{7z+5-\ddots}}}}}}$$
Convergence type $E$ with $E=z$, $P=1-z$, and $C=1/\G(z+1)$, so that
$$(z/(z-1))^{z-1}-\dfrac{p(n)}{q(n)}\sim\dfrac{1/\G(z+1)}{z^nn^{1-z}}\;.$$
$$A=1+(z(z+1)/2)/n+(z(z-2)(3z^3+5z^2+5z+11)/(24(z-1)))/n^2+\cdots$$
Series:
$$(z/(z-1))^{z-1}=\dfrac{z-1}{z}\left(1+\sum_{n\ge0}\dfrac{(z+1)_n}{(n+1)!}z^{-n}\right)$$
Parametric family for $k\ge0$:
\begin{verbatim}
[(z)->(z/(z-1))^(z-1),n*(z+1)+(k+1)*(z-1),-n*z*(n+z)]]
\end{verbatim}
Convergence type $E$ with $E=z$ and $P=1-z+2*k$.
\end{cf}

\smallskip

\begin{cf}\label{6.5.17}{\ }
\begin{verbatim}
[(z)->(z/(z-1))^(z-1),[2,z*(n+2)],
              [[2*z-4,-z^2-z],z*[-(n+2)*(n+2-z),-(n+1)*(n+1+z)]]]
\end{verbatim}
$$(z/(z-1))^{z-1}=2+\dfrac{2z-4}{3z-\dfrac{z^2+z}{4z+\dfrac{3z^2-9z}{5z-\dfrac{2z^2+4z}{6z+\dfrac{4z^2-16z}{7z-\dfrac{3z^2+9z}{8z+\ddots}}}}}}$$
Convergence type $E$ with $E=-(\sqrt{z}+\sqrt{z-1})^2$, $P=0$, and $C=...$,
so that
$$(z/(z-1))^{z-1}-\dfrac{p(n)}{q(n)}\sim(-1)^n\dfrac{C}{(\sqrt{z}+\sqrt{z-1})^{2n}}\;.$$
\end{cf}

\smallskip

\begin{cf}\label{6.5.18}{\ }
\begin{verbatim}
[(z)->(z/(z-1))^(z-1),[(z-1)/z,z-1,(2*z-1)*(2*n-1)],[2*z-2,-(n^2-z^2)]]
\end{verbatim}
$$(z/(z-1))^{z-1}=(z-1)/z+\dfrac{2z-2}{z-1+\dfrac{z^2-1}{6z-3+\dfrac{z^2-4}{10z-5+\dfrac{z^2-9}{14z-7+\dfrac{z^2-16}{18z-9+\ddots}}}}}$$
Convergence type $E$ with $E=(\sqrt{z}+\sqrt{z-1})^4$, $P=0$, and $C=...$,
so that
$$(z/(z-1))^{z-1}-\dfrac{p(n)}{q(n)}\sim\dfrac{C}{(\sqrt{z}+\sqrt{z-1})^{4n}}\;.$$
\end{cf}

\smallskip

\begin{cf}\label{6.5.19}{\ }
\begin{verbatim}
[(a,b,c,z)->((1+b*z)/(1+c*z))^a,
            [1,2+(b+c-a*(b-c))*z,(2*n-1)*(2+(b+c)*z)],
            [2*a*(b-c)*z,-(b-c)^2*(n^2-a^2)*z^2]]
\end{verbatim}
$$\left(\dfrac{1+bz}{1+cz}\right)^a=1+\dfrac{(2b-2c)za}{(-b+c)za+((b+c)z+2)+\dfrac{(b^2-2cb+c^2)z^2a^2+(-b^2+2cb-c^2)z^2}{(3b+3c)z+6+\ddots}}$$
Convergence type $E$ with
$$E=((4(b+c)z+2)\sqrt{(1+bz)(1+cz)}+(b^2+6bc+c^2)z^2+8(b+c)z+8)/((b-c)^2z^2)\;,$$
$P=0$, and $C=...$, so that
$$\left(\dfrac{1+bz}{1+cz}\right)^a-\dfrac{p(n)}{q(n)}\sim\dfrac{C}{E^n}\;.$$
\end{cf}

\smallskip

\begin{cf}\label{6.5.20}{\ }
\begin{verbatim}
[(z)->polrootsreal(x^3-z*x+z)[2],
[1,6*z,2*n*(2*n+1)*z+3*(3*n-1)*(3*n-2)],[6,-6*n*(2*n+1)*(3*n+1)*(3*n+2)*z]]
\end{verbatim}
Denote by $r$ the smallest positive root of $x^3-zx+z=0$, assuming $z$ real
and $|z|>27/4$:
$$r=1+\dfrac{6}{6z-\dfrac{360z}{20z+60-\dfrac{3360z}{42z+168-\dfrac{13860z}{72z+330-\dfrac{39312z}{110z+546-\dfrac{89760z}{156z+816-\ddots}}}}}}$$
Convergence type $E$ with $E=4z/27$, $P=3/2$, and $C=...$, so that
$$r-\dfrac{p(n)}{q(n)}\sim\dfrac{C}{(4z/27)^nn^{3/2}}$$
$$A=1-((604z-1161)/(72(4z-27)))/n+\cdots$$
Series:
$$r=1+\dfrac{1}{z}\sum_{n\ge0}\dfrac{(4/3)_n(5/3)_n}{(n+1)!(5/2)_n}(27/(4z))^n$$
\end{cf}

\smallskip

{\bf Remarks.}\begin{enumerate}\item Although the above is a special cubic
  polynomial and with the restriction $|z|>27/4$, it is not difficult to
  show that \emph{every} real root of a cubic polynomial with integer
  coefficients has a CF expansion of polynomial type as above (but not
  necessarily a convergent hypergeometric expansion), see \cite{Coh8}.
\item Similar CFs exist for roots of trinomials $x^k+ax+b=0$, and more
  generally for values of algebraic hypergeometric functions.
\end{enumerate}

\medskip

\section{Direct Trigonometric and Hyperbolic Functions}\label{sec:trighyp}

\medskip

Since $\sin(z)=\sinh(iz)/i$, $\cos(z)=\cosh(iz)$, and $\tan(z)=\tanh(iz)/i$,
all the continued fractions for trigonometric functions can be
deduced trivially from those of hyperbolic functions by replacing $z$ by $iz$,
so we only give the latter. Also, evidently CFs for $\tanh(z)$ can be
deduced from those of $\cotanh(z)$ and conversely.

\smallskip

As already mentioned for the exponential functions, there is a complete
difference in nature between CFs for $f(z)$ and CFs for $f(\pi z)$, where
$f$ is a hyperbolic or trigonometric function. But since we only give
CFs with rational coefficients, there is also a difference between CFs
for $f(\pi z)$ and CFs for $\pi f(\pi z)$: see for instance below the case
of $f(z)=\sinh(z)$ and $f(z)=\coth(z)$.

\smallskip

Also, note the following elementary transformation which is applicable
to many of the polynomially convergent CFs of this section:
$$1+\dfrac{z^2}{a(1)-z^2+\ddots}=\dfrac{1}{1-\dfrac{z^2}{a(1)+\ddots}}$$
(and similarly with $z^2$ replaced by $z^2+1$),
which can be written using the {\tt Pari/GP} syntax:

\begin{verbatim}
[[1,a(1)-z^2,a(n)],[z^2,b(n)]]=[[0,1,a(n-1)],[1,-z^2,b(n-1)]]
\end{verbatim}

Thus, every one of these CFs can be written in two (equivalent) ways; we only
give the first one.

\smallskip

Note finally that since
$$\dfrac{\sinh(\pi z)/(\pi z)-1}{z^2}=\dfrac{\pi^2}{6}+O(z^2)\text{\quad and\quad}\dfrac{\cosh(\pi z)-1}{z^2}=\dfrac{\pi^2}{2}+O(z^2)\;,$$
the CFs for $\sinh(\pi z)/\pi z$ and for $\cosh(\pi z)$ all specialize to
CFs for $\pi^2$. In particular, the CF \ref{6.6.12} is a generalization
to $\cosh(\pi z)$ of Ap\'ery's CF \ref{5.4.14} for $\z(2)$.

\medskip

\begin{cf}\label{6.6.1}{\ }
\begin{verbatim}
[(z)->sinh(Pi*z)/(Pi*z),[1,1,2*n^2-2*n+1+z^2],[z^2,-n^2*(n^2+z^2)]]
\end{verbatim}
$$\dfrac{\sinh(\pi z)}{\pi z}=1+\dfrac{z^2}{1-\dfrac{z^2+1}{z^2+5-\dfrac{4z^2+16}{z^2+13-\dfrac{9z^2+81}{z^2+25-\dfrac{16z^2+256}{z^2+41-\dfrac{25z^2+625}{z^2+61-\ddots}}}}}}$$
Convergence type $P^+$ with $P=1$ and $C=z\sinh(\pi z)/\pi$, so that
$$\dfrac{\sinh(\pi z)}{\pi z}-\dfrac{p(n)}{q(n)}\sim\dfrac{z\sinh(\pi z)/\pi}{n}\;.$$
$$A=1-((z^2+1)/2)/n+((z^2+1)^2/6)/n^2-(z^2(z^2+1)^2/24)/n^3+\cdots$$
Series:
\begin{align*}\dfrac{\sinh(\pi z)}{\pi z}&=1+z^2\sum_{n\ge0}\dfrac{(1-iz)_n(1+iz)_n}{(n+1)!^2}\\
\dfrac{\pi z}{\sinh(\pi z)}&=1-z^2\sum_{n\ge1}\dfrac{(n-1)!^2}{(1-iz)_n(1+iz)_n}\end{align*}
Parametric family for $k\ge0$:
\begin{verbatim}
[(z)->sinh(Pi*z)/(Pi*z),2*n^2-2*n+z^2+k^2+k+1,-n^2*(n^2+z^2)]
\end{verbatim}
Convergence type $P^+$ with $P=2k+1$.
\end{cf}

\smallskip

\begin{cf}\label{6.6.2}{\ }
\begin{verbatim}
[(z)->sinh(Pi*z)/(Pi*z),[1,1-z^2,2*n-1],[2*z^2,(n^2+z^2)^2]]
\end{verbatim}
$$\dfrac{\sinh(\pi z)}{\pi z}=1+\dfrac{2z^2}{-z^2+1+\dfrac{z^4+2z^2+1}{3+\dfrac{z^4+8z^2+16}{5+\dfrac{z^4+18z^2+81}{7+\dfrac{z^4+32z^2+256}{9+\dfrac{z^4+50z^2+625}{11+\ddots}}}}}}$$
Convergence type $P^-$ with $P=2$ and $C=\sinh^2(\pi z)/\pi^2$, so that
$$\dfrac{\sinh(\pi z)}{\pi z}-\dfrac{p(n)}{q(n)}\sim(-1)^n\dfrac{\sinh^2(\pi z)/\pi^2}{n^2}\;.$$
$$A=1-1/n-z^2/n^2+(2z^2+1)/n^3+z^4/n^4-(3z^4+5z^2+3)/n^5+\cdots$$
Series:
$$\dfrac{\pi z}{\sinh(\pi z)}=-1+2z^2\sum_{n\ge0}\dfrac{(-1)^n}{n^2+z^2}$$
Parametric family for $k\ge0$:
\begin{verbatim}
[(z)->sinh(Pi*z)/(Pi*z),(2*k+1)*(2*n-1),(n^2+z^2)^2]
\end{verbatim}
Convergence type $P^-$ with $P=4k+2$
\end{cf}
      
\smallskip

\begin{cf}\label{6.6.3}{\ }
\begin{verbatim}
[(z)->sinh(Pi*z)/(Pi*z),[[1,1-z^2],[5*n^2+z^2,5*n^2+4*n+1+z^2]],
[[2*z^2,(z^2+1)^2],[-4*n^2*(n^2+z^2),((n+1)^2+z^2)^2]]]
\end{verbatim}
$$\dfrac{\sinh(\pi z)}{\pi z}=1+\dfrac{2z^2}{-z^2+1+\dfrac{z^4+2z^2+1}{z^2+5-\dfrac{4z^2+4}{z^2+10+\dfrac{z^4+8z^2+16}{z^2+20-\dfrac{16z^2+64}{z^2+29+\ddots}}}}}$$
Convergence type $E$ with $E=((1+\sqrt{5})/2)^5i$, $P=0$, and
$C=4\sinh^3(\pi z)/(\pi z((1+\sqrt{5})/2)^5)$, so that
$$\dfrac{\sinh(\pi z)}{\pi z}-\dfrac{p(n)}{q(n)}\sim(-1)^{\lfloor n/2\rfloor}\dfrac{4\sinh^3(\pi z)/(\pi z)}{((1+\sqrt{5})/2)^{5n+5}}\;.$$
$$A=1-((2z^2+2/5)d)/n+(10z^4+(2d+4)z^2+(2/5)(d+1))/n^2+\cdots$$
\end{cf}

\smallskip

\begin{cf}\label{6.6.4}{\ }
\begin{verbatim}
[(z)->cosh(Pi*z),[1,1,4*n^2-5*n+2+2*z^2],[2*z^2,-2*n*(2*n-1)*(n^2+z^2)]]
\end{verbatim}
$$\cosh(\pi z)=1+\dfrac{2z^2}{1-\dfrac{2z^2+2}{2z^2+8-\dfrac{12z^2+48}{2z^2+23-\dfrac{30z^2+270}{2z^2+46-\dfrac{56z^2+896}{2z^2+77-\dfrac{90z^2+2250}{2z^2+116-\ddots}}}}}}$$
Convergence type $P^+$ with $P=1/2$ and $C=2z\sinh(\pi z)/\sqrt{\pi}$, so that
$$\cosh(\pi z)-\dfrac{p(n)}{q(n)}\sim\dfrac{2z\sinh(\pi z)/\sqrt{\pi}}{n^{1/2}}\;.$$
$$A=1-(z^2/3+5/24)/n+(z^4/10+z^2/8+21/640)/n^2+\cdots$$
Series:
$$\cosh(\pi z)=1+2z^2\sum_{n\ge0}\dfrac{(1-iz)_n(1+iz)_n}{(n+1)!(3/2)_n}$$
Parametric family for $k\ge0$:
\begin{verbatim}
[(z)->cosh(Pi*z),4*n^2-5*n+2*z^2+2*k^2+k+2,-2*n*(2*n-1)*(n^2+z^2)]
\end{verbatim}
Convergence type $P^+$ with $P=2k+1/2$.
\end{cf}

\smallskip

\begin{cf}\label{6.6.5}{\ }
\begin{verbatim}
[(z)->cosh(Pi*z),[1,1/2,4*n^2-8*n+5+2*z^2],
                 [2*z^2,-(2*n-1)^2*((n-1/2)^2+z^2)]]
\end{verbatim}
$$\cosh(\pi z)=1+\dfrac{2z^2}{1/2-\dfrac{z^2+1/4}{2z^2+5-\dfrac{9z^2+81/4}{2z^2+17-\dfrac{25z^2+625/4}{2z^2+37-\dfrac{49z^2+2401/4}{2z^2+65-\ddots}}}}}$$
Convergence type $P^+$ with $P=1$ and $C=z^2\cosh(\pi z)$, so that
$$\cosh(\pi z)-\dfrac{p(n)}{q(n)}\sim\dfrac{C}{n}\;.$$
$$A=1-(z^2/2)/n+((2z^4-2z^2-1)/12)/n^2-((z^4-4z^2-2)/24)/n^3+\cdots$$
Series:
\begin{align*}\cosh(\pi z)&=1+4z^2\sum_{n\ge0}\dfrac{(1/2-iz)_n(1/2+iz)_n}{(3/2)_n^2}\\
\dfrac{1}{\cosh(\pi z)}&=1-4z^2\sum_{n\ge1}\dfrac{(-1/2)_n^2}{(1/2-iz)_n(1/2+iz)_n}\end{align*}
Parametric family for $k\ge0$:
\begin{verbatim}
[(z)->cosh(Pi*z),4*n^2-8*n+2*z^2+2*k^2+2*k+5,-(2*n-1)^2*((n-1/2)^2+z^2)]
\end{verbatim}
Convergence type $P^+$ with $P=2k+1$.
\end{cf}

\smallskip

\begin{cf}\label{6.6.6}{\ }
\begin{verbatim}
[(z)->cosh(Pi*z),[4*z^2+1,2-z^2,4*n-2],[2*z^2*(4*z^2+1),(n^2+z^2)^2]]
\end{verbatim}
$$\cosh(\pi z)=4z^2+1+\dfrac{8z^4+2z^2}{-z^2+2+\dfrac{z^4+2z^2+1}{6+\dfrac{z^4+8z^2+16}{10+\dfrac{z^4+18z^2+81}{14+\ddots}}}}$$
Convergence type $P^-$ with $P=4$ and $C=(4z^2+1)\sinh(\pi z)^2/16$, so that
$$\cosh(\pi z)-\dfrac{p(n)}{q(n)}\sim(-1)^n\dfrac{(4z^2+1)\sinh(\pi z)^2/16}{n^4}\;.$$
$$A=1-2/n-2z^2/n^2+(6z^2+5)/n^3+(3z^4-z^2/2-1/8)/n^4+\cdots$$
Parametric family for $k\ge0$:
\begin{verbatim}
[(z)->cosh(Pi*z),(k+1)*(4*n-2),(n^2+z^2)^2]
\end{verbatim}
Convergence type $P^-$ with $P=4k+4$.
\end{cf}

\smallskip

\begin{cf}\label{6.6.7}{\ }
\begin{verbatim}
[(z)->cosh(Pi*z),[1,2,7*n^2-8*n+(3/4)*(9*z^2+4)],
                 [27*z^2/4,-(3/2)*n*(2*n-1)*(4*n^2+9*z^2)]]
\end{verbatim}
$$\cosh(\pi z)=1+\dfrac{(27/4)z^2}{2-\dfrac{(27/)2z^2+6}{(27/4)z^2+15-\dfrac{81z^2+144}{(27/4)z^2+42-\dfrac{(405/2)z^2+810}{(27/4)z^2+83-\ddots}}}}$$
Convergence type $E$ with $E=4/3$, $P=3/2$, and
$C=(9/2)z\sinh(3\pi z/2)/\sqrt{\pi}$, so that
$$\cosh(\pi z)-\dfrac{p(n)}{q(n)}\sim\dfrac{(9/2)z\sinh(3\pi z/2)/\sqrt{\pi}}{(4/3)^nn^{3/2}}\;.$$
$$A=1+(-9z^2/4-47/8)/n+(81z^4/32+675z^2/32+6561/128)/n^2+\cdots$$
Series:
$$\cosh(\pi z)=1+\dfrac{27}{8}z^2\sum_{n\ge0}\dfrac{(1-3iz/2)_n(1+3iz/2)_n}{(n+1)!(3/2)_n}(3/4)^n$$
\end{cf}

\smallskip

\begin{cf}\label{6.6.8}{\ }
\begin{verbatim}
[(z)->cosh(Pi*z),[1,1,3*n^2-3*n+4*z^2+1],[4*z^2,-n*(2*n-1)*(n^2+4*z^2)]]
\end{verbatim}
$$\cosh(\pi z)=1+\dfrac{4z^2}{1-\dfrac{4z^2+1}{4z^2+7-\dfrac{24z^2+24}{4z^2+19-\dfrac{60z^2+135}{4z^2+37-\ddots}}}}$$
Convergence type $E$ with $E=2$, $P=3/2$, and $C=2z\sinh(2\pi z)/\sqrt{\pi}$,
so that
$$\cosh(\pi z)-\dfrac{p(n)}{q(n)}\sim\dfrac{z\sinh(2\pi z)/\sqrt{\pi}}{2^{n-1}n^{3/2}}\;.$$
$$A=1-(4z^2+23/8)/n+(8z^4+35z^2/2+1361/128)/n^2+\cdots$$
Series:
$$\cosh(\pi z)=1+4z^2\sum_{n\ge0}\dfrac{(1-2iz)_n(1+2iz)_n}{(n+1)!(3/2)_n}2^{-n}$$
Parametric family for $k\ge0$:
\begin{verbatim}
[(z)->cosh(Pi*z),3*n^2-3*n+4*z^2+1,-n*(2*n-2*k-1)*(n^2+4*z^2)]
\end{verbatim}
Convergence type $E$ with $E=2$ and $P=3k+3/2$.
\end{cf}

\smallskip

\begin{cf}\label{6.6.9}{\ }
\begin{verbatim}
[(z)->cosh(Pi*z),[1,2,5*n^2-4*n+9*z^2+1],[9*z^2,-2*n*(2*n-1)*(n^2+9*z^2)]]
\end{verbatim}
$$\cosh(\pi z)=1+\dfrac{9z^2}{2-\dfrac{18z^2+2}{9z^2+13-\dfrac{108z^2+48}{9z^2+34-\dfrac{270z^2+270}{9z^2+65-\ddots}}}}$$
Convergence type $E$ with $E=4$, $P=3/2$, and $C=z\sinh(3\pi z)/\sqrt{\pi}$,
so that
$$\cosh(\pi z)-\dfrac{p(n)}{q(n)}\sim\dfrac{z\sinh(3\pi z)/\sqrt{\pi}}{2^{2n}n^{3/2}}\;.$$
$$A=1-(9z^2+15/8)/n+(81z^4/2+195z^2/8+481/128)/n^2+\cdots$$
Series:
$$\cosh(\pi z)=1+\dfrac{9}{2}z^2\sum_{n\ge0}\dfrac{(1-3iz)_n(1+3iz)_n}{(n+1)!(3/2)_n}2^{-2n}$$
Parametric family for $k\ge0$:
\begin{verbatim}
[(z)->cosh(Pi*z),5*n^2-4*n+9*z^2+1+k*(2*n-1),
               -(2*n-k)*(2*n-k-1)*(n^2+9*z^2)]
\end{verbatim}
Convergence type $E$ with $E=4$ and $P=3k+3/2$.
\end{cf}

\smallskip

\begin{cf}\label{6.6.10}{\ }
\begin{verbatim}
[(z)->cosh(Pi*z),[1,2-4*z^2,7*n^2-7*n+4*z^2+2],
                 [12*z^2,8*(n^2+z^2)*(n^2+4*z^2)]]
\end{verbatim}
$$\cosh(\pi z)=1+\dfrac{12z^2}{-4z^2+2+\dfrac{32z^4+40z^2+8}{4z^2+16+\dfrac{32z^4+160z^2+128}{4z^2+44+\dfrac{32z^4+360z^2+648}{4z^2+86+\ddots}}}}$$
Convergence type $E$ with $E=-8$, $P=0$, and $C=\cosh(\pi z)\sinh(\pi z)^2$,
so that
$$\cosh(\pi z)-\dfrac{p(n)}{q(n)}\sim(-1)^n\dfrac{\cosh(\pi z)\sinh(\pi z)^2}{2^{3n}}\;.$$
$$A=1-((9z^2+18z+10)/3)/n+((9z^2+18z+10)(9z^2+18z+13)/9)/n^2+\cdots$$
\end{cf}

\smallskip

\begin{cf}\label{6.6.11}{\ }
\begin{verbatim}
[(z)->cosh(Pi*z),[[1,-z^2],[5*n^2-3*n+z^2+1/2,5*n^2+n+z^2+1/2]],
[[2*z^2,(z^2+1)^2],[-(2*n-1)^2*((n-1/2)^2+z^2),((n+1)^2+z^2)^2]]]
\end{verbatim}
$$\cosh(\pi z)=1+\dfrac{2z^2}{-z^2+\dfrac{z^4+2z^2+1}{z^2+5/2-\dfrac{z^2+1/4}{z^2+13/2+\dfrac{z^4+8z^2+16}{z^2+29/2-\dfrac{9z^2+81/4}{z^2+45/2+\ddots}}}}}$$
Convergence type $E$ with $E=((1+\sqrt{5})/2)^5i$, $P=0$, and
$C=2\sinh(\pi z)\sinh(2\pi z)/((1+\sqrt{5})/2)^3$, so that
$$\cosh(\pi z)-\dfrac{p(n)}{q(n)}\sim(-1)^{\lfloor n/2\rfloor}\dfrac{2\sinh(\pi z)\sinh(2\pi z)}{((1+\sqrt{5})/2)^{5n+3}}\;.$$
$$A=1-((2z^2+1/5)d)/n+(10z^4+(8d/5+2)z^2+(4d/25+1/10))/n^2+\cdots$$
\end{cf}

\smallskip

\begin{cf}\label{6.6.12}{\ }
\begin{verbatim}
[(z)->cosh(Pi*z),[1,3*(1-z^2),11*n^2-11*n+3+9*z^2],
                 [15*z^2,(n^2+4*z^2)*(n^2+9*z^2)]]
\end{verbatim}
$$\cosh(\pi z)=1+\dfrac{15z^2}{-3z^2+3+\dfrac{36z^4+13z^2+1}{9z^2+25+\dfrac{36z^4+52z^2+16}{9z^2+69+\dfrac{36z^4+117z^2+81}{9z^2+135+\ddots}}}}$$
Convergence type $E$ with $E=-((1+\sqrt{5})/2)^{10}$, $P=0$, and $C=2\sinh(2\pi z)\sinh(3\pi z)/((1+\sqrt{5})/2)^5$, so
that
$$\cosh(\pi z)-\dfrac{p(n)}{q(n)}\sim(-1)^n\dfrac{2\sinh(2\pi z)\sinh(3\pi z)}{((1+\sqrt{5})/2)^{10n+5}}\;.$$
$$A=1+(-5dz^2-d/5)/n+(125z^4/2+(5d/2+5)z^2+d/10+1/10)/n^2+\cdots$$
\end{cf}

This CF converges quite fast and it is simply the generalization of Ap\'ery's
CF \ref{5.4.14} for $\z(2)=\pi^2/6$.

\smallskip

\begin{cf}\label{6.6.13}{\ }
\begin{verbatim}
[(z)->sinh(Pi*z/2)/z,[1,6,8*n^2-2*n+1+z^2],
                     [1+z^2,-2*n*(2*n+1)*((2*n+1)^2+z^2)]]
\end{verbatim}
$$\dfrac{\sinh(\pi z/2)}{z}=1+\dfrac{z^2+1}{6-\dfrac{6z^2+54}{z^2+29-\dfrac{20z^2+500}{z^2+67-\dfrac{42z^2+2058}{z^2+121-\dfrac{72z^2+5832}{z^2+191-\dfrac{110z^2+13310}{z^2+277-\ddots}}}}}}$$
Convergence type $P^+$ with $P=1/2$ and $C=\cosh(\pi z/2)/\sqrt{\pi}$, so that
$$\dfrac{\sinh(\pi z/2)}{z}-\dfrac{p(n)}{q(n)}\sim\dfrac{\cosh(\pi z/2)/\sqrt{\pi}}{n^{1/2}}\;.$$
$$A=1-(z^2/12+11/24)/n+(z^4/160+3z^2/32+181/640)/n^2+\cdots$$
Series:
$$\dfrac{\sinh(\pi z/2)}{z}=1+\dfrac{z^2+1}{6}\sum_{n\ge0}\dfrac{((3-iz)/2)_n((3+iz)/2)_n}{(n+1)!(5/2)_n}$$
Parametric family for $k\ge0$:
\begin{verbatim}
[(z)->sinh(Pi*z/2)/z,8*n^2-2*n+z^2+4*k^2+2*k+1,-2*n*(2*n+1)*((2*n+1)^2+z^2)]
\end{verbatim}
Convergence type $P^+$ with $P=2k+1/2$.
\end{cf}

\smallskip

\begin{cf}\label{6.6.14}{\ }
\begin{verbatim}
[(z)->sinh(Pi*z/2)/z,[3/2,6,5*n^2+n+(9*z^2+1)/4],
                     [(3/8)*(9*z^2+1),-n*(2*n+1)*((2*n+1)^2+9*z^2)/2]]
\end{verbatim}
$$\dfrac{\sinh(\pi z/2)}{z}=3/2+\dfrac{27/8z^2+3/8}{6-\dfrac{27/2z^2+27/2}{9/4z^2+89/4-\dfrac{45z^2+125}{9/4z^2+193/4-\ddots}}}$$
Convergence type $E$ with $E=4$, $P=3/2$, and
$C=\cosh(3\pi z/2)/(4\sqrt{\pi})$, so that
$$\dfrac{\sinh(\pi z/2)}{z}-\dfrac{p(n)}{q(n)}\sim\dfrac{\cosh(3\pi z/2)/\sqrt{\pi}}{2^{2n+2}n^{3/2}}\;.$$
$$A=1-(9z^2/4+21/8)/n+(81z^4/32+285z^2/32+841/128)/n^2+\cdots$$
Series:
$$\dfrac{\sinh(\pi z/2)}{z}=\dfrac{3}{2}+\dfrac{9z^2+1}{16}\sum_{n\ge0}\dfrac{((3-3iz)/2)_n((3+3iz)/2)_n}{(n+1)!(5/2)_n}2^{-2n}$$
Parametric family for $k\ge0$:
\begin{verbatim}
[(z)->sinh(Pi*z/2)/z,5*n^2+(2*k+1)*n+(9*z^2+1)/4,
                   -(2*n-k)*(2*n+1-k)*((2*n+1)^2+9*z^2)/4]
\end{verbatim}
Convergence type $E$ with $E=4$ and $P=3k+3/2$.
\end{cf}

\smallskip

\begin{cf}\label{6.6.15}{\ }
\begin{verbatim}
[(z)->sinh(Pi*z/2)/z,[[0,1],[10*n-2,10*n+2]],
                     [[3/2,-(9*z^2+1)/2],[36*(n^2+z^2),-((2*n+1)^2+9*z^2)]]]
\end{verbatim}
$$\dfrac{\sinh(\pi z/2)}{z}=\dfrac{3/2}{1-\dfrac{9/2z^2+1/2}{8+\dfrac{36z^2+36}{12-\dfrac{9z^2+9}{18+\dfrac{36z^2+144}{22-\dfrac{9z^2+25}{28+\ddots}}}}}}$$
Convergence type $E$ with $E=((1+\sqrt{5})/2)^5i$, $P=0$, and
$C=\sinh(\pi z)\cosh(3\pi z/2)/(z(1+\sqrt{5})/2)$, so that
$$\dfrac{\sinh(\pi z/2)}{z}-\dfrac{p(n)}{q(n)}\sim(-1)^{\lfloor n/2\rfloor}\dfrac{\sinh(\pi z)\cosh(3\pi z/2)/z}{((1+\sqrt{5})/2)^{5n+1}}\;.$$
$$A=1+(-5dz^2/2-d/10)/n+(125z^4/8+(d+5/4)z^2+d/25+1/40)/n^2+\cdots$$
\end{cf}

\smallskip

\begin{cf}\label{6.6.16}{\ }
\begin{verbatim}
[(z)->sinh(Pi*z/2)/z,
[[1,5-z^2],[40*n^3+2*(z^2-1)*n,2*((5*n+4)*(2*n+1)^2+z^2*(n+1))]],
[(z^2+1)*[3,2*(z^2+4)],
2*(2*n+1)*((2*n+1)^2+z^2)*[-n*(4*n-1)*(4*n+3),(n+1)*((2*n+2)^2+z^2)]]]
\end{verbatim}
$$\dfrac{\sinh(\pi z/2)}{z}=1+\dfrac{3z^2+3}{-z^2+5+\dfrac{2z^4+10z^2+8}{2z^2+38-\dfrac{126z^2+1134}{4z^2+162+\dfrac{12z^4+300z^2+1728}{4z^2+316-\ddots}}}}$$
Convergence type $E$ with $E=((1+\sqrt{5})/2)^5i$, $P=0$, and
$C=(2/z)\cosh(\pi z/2)\sinh(\pi z)/((1+\sqrt{5})/2)^5$, so that
$$\dfrac{\sinh(\pi z/2)}{z}-\dfrac{p(n)}{q(n)}\sim(-1)^{\lfloor n/2\rfloor}\dfrac{(2/z)\cosh(\pi z/2)\sinh(\pi z)}{((1+\sqrt{5})/2)^{5n+5}}\;.$$
$$A=1+((-z^2/2+1/10)d)/n+((5/8)z^4+(d/2-1/4)z^2-d/10+1/40)/n^2+\cdots$$
\end{cf}

\smallskip

\begin{cf}\label{6.6.17}{\ }
\begin{verbatim}
[(z)->sinh(Pi*z/2)/(z*sqrt(2)),[1,6,6*n^2+(4*z^2+1)/2],
                   [(4*z^2+1)/2,-n*(2*n+1)*((2*n+1)^2+4*z^2)]]
\end{verbatim}
$$\dfrac{\sinh(\pi z/2)}{z\sqrt{2}}=1+\dfrac{2z^2+1/2}{6-\dfrac{12z^2+27}{2z^2+49/2-\dfrac{40z^2+250}{2z^2+109/2-\dfrac{84z^2+1029}{2z^2+193/2-\ddots}}}}$$
Convergence type $E$ with $E=2$, $P=3/2$, and $C=\cosh(\pi z)/(2\sqrt{\pi})$,
so that
$$\dfrac{\sinh(\pi z/2)}{z\sqrt{2}}-\dfrac{p(n)}{q(n)}\sim\dfrac{\cosh(\pi z)/\sqrt{\pi}}{2^{n+1}n^{3/2}}\;.$$
$$A=1-(z^2+29/8)/n+(z^4/2+45z^2/8+1881/128)/n^2+\cdots$$
Series:
$$\dfrac{\sinh(\pi z/2)}{z\sqrt{2}}=1+\dfrac{4z^2+1}{12}\sum_{n\ge0}\dfrac{(3/2-iz)_n(3/2+iz)_n}{(n+1)!(5/2)_n}2^{-n}$$
Parametric family for $k\ge0$:
\begin{verbatim}
[(z)->sinh(Pi*z/2)/(z*sqrt(2)),6*n^2+(4*z^2+1)/2,
                  -(n-k)*(2*n+1)*((2*n+1)^2+4*z^2)]
\end{verbatim}
Convergence type $E$ with $E=2$ and $P=3k+3/2$.
\end{cf}

\smallskip

\begin{cf}\label{6.6.18}{\ }
\begin{verbatim}
[(z)->sinh(Pi*z/2)/(z*sqrt(2)),[[1,3],[6*n+1,6*n+3]],
[[z^2+1/4,-3*z^2-27/4],[4*n*(4*n^2+z^2),-(2*n+3)*((2*n+3)^2/4+z^2)]]]
\end{verbatim}
$$\dfrac{\sinh(\pi z/2)}{z\sqrt{2}}=1+\dfrac{z^2+1/4}{3-\dfrac{3z^2+27/4}{7+\dfrac{4z^2+16}{9-\dfrac{5z^2+125/4}{13+\dfrac{8z^2+128}{15-\dfrac{7z^2+343/4}{19+\ddots}}}}}}$$
Convergence type $E$ with $E=2\sqrt{2}i$, $P=1/2$, and $C=...$, so that
$$\dfrac{\sinh(\pi z/2)}{z\sqrt{2}}-\dfrac{p(n)}{q(n)}\sim(-1)^{\lfloor n/2\rfloor}\dfrac{C}{2^{3n/2}n^{1/2}}\;.$$
$$A=1+(-3z^2/2-47/36)/n+(9z^4/8+33z^2/8+1883/864)/n^2+\cdots$$
Series:
$$\dfrac{z\sqrt{2}}{\sinh(\pi z/2)}=\sum_{n\ge0}(-1)^n\dfrac{(6n+1)(1/2)_n(1/2-iz)_n(1/2+iz)_n}{n!(1-iz/2)_n(1+iz/2)_n}2^{-3n}$$
\end{cf}

\smallskip

\begin{cf}\label{6.6.19}{\ }
\begin{verbatim}
[(z)->sinh(Pi*z/3)/(z*sqrt(3)/2),[1,6,7*n^2-n+(3/4)*(z^2+1)],
[(3/4)*(z^2+1),-(3/2)*n*(2*n+1)*((2*n+1)^2+z^2)]]
\end{verbatim}
$$\dfrac{\sinh(\pi z/3)}{z\sqrt{3}/2}=1+\dfrac{(3/4)z^2+3/4}{6-\dfrac{(9/2)z^2+81/2}{(3/4)z^2+107/4-\dfrac{15z^2+375}{(3/4)z^2+243/4-\ddots}}}$$
Convergence type $E$ with $E=4/3$, $P=3/2$, and
$C=3\cosh(\pi z/2)/(2\sqrt{\pi})$, so that
$$\dfrac{\sinh(\pi z/3)}{z\sqrt{3}/2}-\dfrac{p(n)}{q(n)}\sim\dfrac{2\cosh(\pi z/2)/\sqrt{\pi}}{(4/3)^{n+1}n^{3/2}}\;.$$
$$A=1-(z^2/4+53/8)/n+(z^4/32+85z^2/32+7561/128)/n^2+\cdots$$
Series:
$$\dfrac{\sinh(\pi z/3)}{z\sqrt{3}/2}=1+\dfrac{z^2+1}{8}\sum_{n\ge0}\dfrac{((3-iz)/2)_n((3+iz)/2)_n}{(n+1)!(5/2)_n}(3/4)^n$$
\end{cf}
    
\smallskip

\begin{cf}\label{6.6.20}{\ }
\begin{verbatim}
[(z)->sinh(Pi*z/3)/(z*sqrt(3)/2),[1,6,5*n^2+2*n+z^2],
                                 [1+z^2,-2*n*(2*n+1)*((n+1)^2+z^2)]]
\end{verbatim}
$$\dfrac{\sinh(\pi z/3)}{z\sqrt{3}/2}=1+\dfrac{z^2+1}{6-\dfrac{6z^2+24}{z^2+24-\dfrac{20z^2+180}{z^2+51-\dfrac{42z^2+672}{z^2+88-\dfrac{72z^2+1800}{z^2+135-\dfrac{110z^2+3960}{z^2+192-\ddots}}}}}}$$
Convergence type $E$ with $E=4$, $P=1/2$, and $C=\sinh(\pi z)/(6z\sqrt{\pi})$,
so that
$$\dfrac{\sinh(\pi z/3)}{z\sqrt{3}/2}-\dfrac{p(n)}{q(n)}\sim\dfrac{\sinh(\pi z)/(6z\sqrt{\pi})}{2^{2n}n^{1/2}}\;.$$
$$A=1-(z^2+25/24)/n+(z^4/2+23z^2/8+683/384)/n^2+\cdots$$
Series:
$$\dfrac{\sinh(\pi z/3)}{z\sqrt{3}/2}=1+\dfrac{z^2+1}{6}\sum_{n\ge0}\dfrac{(2-iz)_n(2+iz)_n}{(n+1)!(5/2)_n}2^{-2n}$$
\end{cf}

\smallskip

The following CFs are for $\sin$ and $\cos$, and since they are not even
functions, they cannot give a CF for $\sinh$ and $\cosh$ without using complex
numbers.

\smallskip

\begin{cf}\label{6.6.21}{\ }
\begin{verbatim}
[(z)->sin(Pi*z)/(Pi*z),[1,2*n^2+(3*z-2)*n+(1-z)^2],[-z^2,-n*(n+z)^3]]
\end{verbatim}
$$\dfrac{\sin(\pi z)}{\pi z}=1-\dfrac{z^2}{z^2+z+1-\dfrac{z^3+3z^2+3z+1}{z^2+4z+5-\dfrac{2z^3+12z^2+24z+16}{z^2+7z+13-\dfrac{3z^3+27z^2+81z+81}{z^2+10z+25-\ddots}}}}$$
Convergence type $P^+$ with $P=|1-z|$ and $C=\sin^3(\pi z)\G(-z)/(\pi^3 (1-z))$, so that
$$\dfrac{\sin(\pi z)}{\pi z}-\dfrac{p(n)}{q(n)}\sim\dfrac{\sin^3(\pi z)\G(-z)/(\pi^3 (1-z))}{n^{|1-z|}}\;.$$
$$A=1-((z^3-3z^2+2)/(2(z-2)))/n+\cdots$$
Series:
$$\dfrac{\pi}{\sin(\pi z)}=\sum_{n\ge0}\dfrac{(z)_n}{(n+z)n!}$$
Parametric family for $k\ge0$:
\begin{verbatim}
[(z)->sin(Pi*z)/(Pi*z),2*n^2+(3*z-2)*n+(1-z)^2+k*(1-z)+k^2,-n*(n+z)^3]
\end{verbatim}
Convergence type $P^+$ with $P=2k+1-z$.
\end{cf}

\smallskip

\begin{cf}\label{6.6.22}{\ }
\begin{verbatim}
[(z)->cos(Pi*z)/(Pi/2),[1-2*z,8*n^2-(12*z+2)*n+(2*z+1)^2],
                       [(2*z-1)^3,-2*n*(2*n+1-2*z)^3]]
\end{verbatim}
$$\dfrac{\cos(\pi z)}{\pi/2}=-2z+1+\dfrac{8z^3-12z^2+6z-1}{4z^2-8z+7+\dfrac{16z^3-72z^2+108z-54}{4z^2-20z+29+\dfrac{32z^3-240z^2+600z-500}{4z^2-32z+67+\ddots}}}$$
Convergence type $P^+$ with $P=|z+1/2|$ and
$C=\cos^3(\pi z)\G(z+1/2)/(\pi^3(z+1/2))$, so that
$$\dfrac{\cos(\pi z)}{\pi/2}-\dfrac{p(n)}{q(n)}\sim\dfrac{\cos^3(\pi z)\G(z+1/2)/(\pi^3(z+1/2))}{n^{|z+1/2|}}\;.$$
$$A=1+((8z^3+12z^2-18z-11)/(16z+24))/n+\cdots$$
Series:
$$\dfrac{\pi/2}{\cos(\pi z)}=\sum_{n\ge0}\dfrac{(1/2-z)_n}{(2n+1-2z)n!}$$
Parametric family for $k\ge0$:
\begin{verbatim}
[(z)->cos(Pi*z)/(Pi/2),
8*n^2-(12*z+2)*n+(2*z+1)^2+2*k*(2*z+2*k+1),-2*n*(2*n+1-2*z)^3]
\end{verbatim}
Convergence type $P^+$ with $P=2k+z+1/2$.
\end{cf}

This CF is simply the previous one with $z$ replaced by $1/2-z$, but
we included it nonetheless.

\smallskip

\begin{cf}\label{6.6.23}{\ }
\begin{verbatim}
[(z)->sin(Pi*z),[1,1,3*n^2-3*n-4*z^2+4*z],
                [-(2*z-1)^2,-n*(2*n-1)*(n-2*z+1)*(n+2*z-1)]]
\end{verbatim}
$$\sin(\pi z)=1-\dfrac{4z^2-4z+1}{1+\dfrac{4z^2-4z}{-4z^2+4z+6+\dfrac{24z^2-24z-18}{-4z^2+4z+18+\dfrac{60z^2-60z-120}{-4z^2+4z+36+\ddots}}}}$$
Convergence type $E$ with $E=2$, $P=3/2$, and $C=(2z-1)\sin(2\pi z)/\sqrt{\pi}$, so that
$$\sin(\pi z)-\dfrac{p(n)}{q(n)}\sim\dfrac{(2z-1)\sin(2\pi z)/\sqrt{\pi}}{2^nn^{3/2}}$$
$$A=1+(4z^2-4z-15/8)/n+\cdots$$
Series:
$$\sin(\pi z)=1-(2z-1)^2\sum_{n\ge0}\dfrac{(2z)_n(2-2z)_n}{(n+1)!(3/2)_n}2^{-n}$$
Parametric family for $k\ge0$:
\begin{verbatim}
[(z)->sin(Pi*z),3*n^2-3*n-4*z^2+4*z,-n*(2*n-2*k-1)*(n-2*z+1)*(n+2*z-1)]
\end{verbatim}
Convergence type $E$ with $E=2$ and $P=3k+3/2$.
\end{cf}

\smallskip

\begin{cf}\label{6.6.24}{\ }
\begin{verbatim}
[(z)->sin(Pi*z),[1,4*z^2-4*z+3,7*n^2-7*n+2-(2*z-1)^2],
                [-3*(2*z-1)^2,2*(n^2-(2*z-1)^2)*(4*n^2-(2*z-1)^2)]]
\end{verbatim}
$$\sin(\pi z)=1-\dfrac{12z^2-12z+3}{4z^2-4z+3+\dfrac{32z^4-64z^3+8z^2+24z}{-4z^2+4z+15+\dfrac{32z^4-64z^3-112z^2+144z+90}{-4z^2+4z+43+\ddots}}}$$
Convergence type $E$ with $E=-8$, $P=0$, and $C=-\sin(\pi z)\cos^2(\pi z)$,
so that
$$\sin(\pi z)-\dfrac{p(n)}{q(n)}\sim(-1)^{n+1}\dfrac{\sin(\pi z)\cos^2(\pi z)}{2^{3n}}$$
$$A=1+(3z^2-3z+5/12)/n+\cdots$$
\end{cf}

\smallskip

\begin{cf}\label{6.6.25}{\ }
\begin{verbatim}
[(z)->sin(Pi*z),[1,12*z^2-12*z+15,44*n^2-44*n-36*z^2+36*z+3],
                [-15*(2*z-1)^2,4*(n^2-(2*z-1)^2)*(4*n^2-9*(2*z-1)^2)]]
\end{verbatim}
$$\sin(\pi z)=1-\dfrac{60z^2-60z+15}{12z^2-12z+15+\dfrac{576z^4-1152z^3+656z^2-80z}{-36z^2+36z+91+\dfrac{576z^4-1152z^3+32z^2+544z+84}{-36z^2+36z+267+\ddots}}}$$
Convergence type $E$ with $E=-((1+\sqrt{5})/2)^{10}$, $P=0$, and $C=...$, so
that
$$\sin(\pi z)-\dfrac{p(n)}{q(n)}\sim(-1)^n\dfrac{C}{((1+\sqrt{5})/2)^{10n}}$$
$$A=1+d(5z^2-5z+21/20)/n+\cdots$$
\end{cf}

\smallskip

\begin{cf}\label{6.6.26}{\ }
\begin{verbatim}
[(z)->cosh(z),[1,2,2*n*(2*n-1)+z^2],[z^2,-2*n*(2*n-1)*z^2]]
\end{verbatim}
$$\cosh(z)=1+\dfrac{z^2}{2-\dfrac{2z^2}{z^2+12-\dfrac{12z^2}{z^2+30-\dfrac{30z^2}{z^2+56-\dfrac{56z^2}{z^2+90-\dfrac{90z^2}{z^2+132-\ddots}}}}}}$$
Convergence type $F^2$ with $E=4/z^2$, $P=3/2$, and $C=(z^2/4)\sqrt{\pi}$,
so that
$$\cosh(z)-\dfrac{p(n)}{q(n)}\sim\dfrac{\sqrt{\pi}}{n!^2(2/z)^{2n+2}n^{3/2}}\;.$$
$$A=1-(11/8)/n+(z^2/4+201/128)/n^2-(39z^2/32+1713/1024)/n^3+\cdots$$
Series:
$$\cosh(z)=1+\dfrac{z^2}{2}\sum_{n\ge0}\dfrac{(z/2)^{2n}}{(n+1)!(3/2)_n}$$
\end{cf}

\smallskip

Since $\cosh^2(z)=(\cosh(2z)+1)/2$ and $\sinh^2(z)=(\cosh(2z)-1)/2$,
the above CF also gives CFs for these two functions.

\smallskip

\begin{cf}\label{6.6.27}{\ }
\begin{verbatim}
[(z)->sinh(z)/z,[1,6,2*n*(2*n+1)+z^2],[z^2,-2*n*(2*n+1)*z^2]]
\end{verbatim}
$$\dfrac{\sinh(z)}{z}=1+\dfrac{z^2}{6-\dfrac{6z^2}{z^2+20-\dfrac{20z^2}{z^2+42-\dfrac{42z^2}{z^2+72-\dfrac{72z^2}{z^2+110-\dfrac{110z^2}{z^2+156-\ddots}}}}}}$$
Convergence type $F^2$ with $E=4/z^2$, $P=5/2$, and $C=(z^2/8)\sqrt{\pi}$,
so that
$$\dfrac{\sinh(z)}{z}-\dfrac{p(n)}{q(n)}\sim\dfrac{\sqrt{\pi}/2}{n!^2(2/z)^{2n+2}n^{5/2}}\;.$$
$$A=1-(23/8)/n+(z^2/4+753/128)/n^2-(59z^2/32+10749/1024)/n^3+\cdots$$
Series:
$$\dfrac{\sinh(z)}{z}=1+\dfrac{z^2}{6}\sum_{n\ge0}\dfrac{(z/2)^{2n}}{(n+1)!(5/2)_n}$$
\end{cf}

The above two continued fractions for $\cosh(z)$ and $\sinh(z)/z$ are simply
the term-by-term expansion of the corresponding Taylor series.

\smallskip

\begin{cf}\label{6.6.28}{\ }
\begin{verbatim}
[()->cos(sqrt(2)),[0,6,n*(2*n+3)],[1,(n+1)*(2*n+1)]]
[()->cos(sqrt(2)),[0,6,2*n+3],[1,(2*n+1)/n]]
\end{verbatim}
$$\cos(\sqrt{2})=\dfrac{1}{6+\dfrac{6}{14+\dfrac{15}{27+\dfrac{28}{44+\dfrac{45}{65+\dfrac{66}{90+\ddots}}}}}}=\dfrac{1}{6+\dfrac{3}{7+\dfrac{5/2}{9+\dfrac{7/3}{11+\dfrac{9/4}{13+\dfrac{11/5}{15+\ddots}}}}}}$$
Convergence type $F^2$ with $E=-2$, $P=7/2$, and $C=\sqrt{\pi}/4$, so that
$$\cos(\sqrt{2})-\dfrac{p(n)}{q(n)}\sim(-1)^n\dfrac{\sqrt{\pi}/4}{n!^22^nn^{7/2}}\;.$$
$$A=1-(39/8)/n+(1937/128)/n^2-(37453/1024)/n^3+\cdots$$
Series:
$$\cos(\sqrt{2})=\dfrac{1}{3}\sum_{n\ge0}(-1)^n\dfrac{1}{(n+2)!(5/2)_n}2^{-n}$$
\end{cf}

\smallskip

\begin{cf}\label{6.6.29}{\ }
\begin{verbatim}
[()->sin(sqrt(2))/sqrt(2),[0,1,n*(2*n-3)],[1,1,(n-1)*(2*n-1)]]
[()->sin(sqrt(2))/sqrt(2),[n],[1,1,(n-1)/(2*n-3)]]
\end{verbatim}
$$\dfrac{\sin(\sqrt{2})}{\sqrt{2}}=\dfrac{1}{1+\dfrac{1}{2+\dfrac{3}{9+\dfrac{10}{20+\dfrac{21}{35+\dfrac{36}{54+\ddots}}}}}}=\dfrac{1}{1+\dfrac{1}{2+\dfrac{1}{3+\dfrac{2/3}{4+\dfrac{3/5}{5+\dfrac{4/7}{6+\ddots}}}}}}$$
Convergence type $F^2$ with $E=-2$, $P=1/2$, and $C=\sqrt{\pi}/2$, so that
$$\dfrac{\sin(\sqrt{2})}{\sqrt{2}}-\dfrac{p(n)}{q(n)}\sim(-1)^n\dfrac{\sqrt{\pi}/2}{n!^22^nn^{1/2}}\;.$$
$$A=1-(3/8)/n-(39/128)/n^2+(1367/1024)/n^3+\cdots$$
Series:
$$\dfrac{\sin(\sqrt{2})}{\sqrt{2}}=\sum_{n\ge0}(-1)^n\dfrac{1}{n!(3/2)_n}2^{-n}$$
\end{cf}

\smallskip

\begin{cf}\label{6.6.30}{\ }
\begin{verbatim}
[(z)->tanh(z),[0,2*n-1],[z,z^2]]
\end{verbatim}
$$\tanh(z)=\dfrac{z}{1+\dfrac{z^2}{3+\dfrac{z^2}{5+\dfrac{z^2}{7+\dfrac{z^2}{9+\dfrac{z^2}{11+\ddots}}}}}}$$
Convergence type $F^2$ with $E=-4/z^2$, $P=0$, and $C=(\pi/2)z/\cosh(z)^2$,
so that
$$\tanh(z)-\dfrac{p(n)}{q(n)}\sim(-1)^n\dfrac{\pi/\cosh(z)^2}{n!^2(2/z)^{2n+1}}\;.$$
$$A=1+(z^2/2-1/4)/n+(z^4/8-3z^2/8+5/32)/n^2+\cdots$$
\end{cf}

\smallskip

\begin{cf}\label{6.6.31}{\ }
\begin{verbatim}
[(z)->tan(z)+1/cos(z),[1,2-z,2*(2*n-1)],[2*z,-z^2]]
\end{verbatim}
$$\tan(z)+\dfrac{1}{\cos(z)}=1+\dfrac{2z}{-z+2-\dfrac{z^2}{6-\dfrac{z^2}{10-\dfrac{z^2}{14-\dfrac{z^2}{18-\dfrac{z^2}{22-\ddots}}}}}}$$
Convergence type $F^2$ with $E=16/z^2$, $P=0$, and
$C=\pi z(\tan(z)+1/cos(z))/(2\cos(z))$, so that
$$\tan(z)+\dfrac{1}{\cos(z)}-\dfrac{p(n)}{q(n)}\sim\dfrac{2\pi(\tan(z)+1/\cos(z))/\cos(z)}{n!^2(4/z)^{2n+1}}\;.$$
$$A=1-(z^2/8+1/4)/n+(z^4/128+3z^2/32+5/32)/n^2+\cdots$$
\end{cf}

\smallskip

\begin{cf}\label{6.6.32}{\ }
\begin{verbatim}
[(z)->tan(z)+1/cos(z),
[[0,z-1],[4*(2*n-1)*z+2*(4*n-1)*(4*n-3),-8*n*z+2*(4*n-1)*(4*n+1)]],
         [-1,(2*n-3)*(2*n+1)*z^2]]
\end{verbatim}
$$\tan(z)+\dfrac{1}{\cos(z)}=-\dfrac{1}{z-1-\dfrac{3z^2}{4z+6+\dfrac{5z^2}{-8z+30+\dfrac{21z^2}{12z+70+\dfrac{45z^2}{-16z+126+\ddots}}}}}$$
Convergence type $F^2$ with $E=-16/z^2$, $P=0$, and
$C=2\pi(\tan(z)+1/\cos(z))/\cos(z)$, so that
$$\tan(z)+\dfrac{1}{\cos(z)}-\dfrac{p(n)}{q(n)}\sim(-1)^n\dfrac{2\pi(\tan(z)+1/\cos(z))/\cos(z)}{n!^2(4/z)^{2n}}\;.$$
$$A=1+(-z^2/8+3/4)/n+(z^4/128-3z^2/32-z/4+9/32)/n^2+\cdots$$
\end{cf}

\smallskip
      
Note that the same CFs for the corresponding hyperbolic function would involve
$i=\sqrt{-1}$, and that to obtain CFs for $1/\cos(z)-\tan(z)$ one simply
replaces $z$ by $-z$.

Note that the CF \ref{6.6.30} for $\tanh(z)$ provides CFs for
$1/\sinh(z)-\cotanh(z)=\tanh(z/2)$ and for $1/\sinh(z)+\cotanh(z)=\cotanh(z/2)$
(and for their trigonometric counterparts).

\smallskip

\begin{cf}\label{6.6.33}{\ }
\begin{verbatim}
[(z)->z*cotanh(Pi*z/4),[1,2],[(2*n+1)^2+z^2]]
\end{verbatim}
$$z\coth(\pi z/4)=1+\dfrac{z^2+1}{2+\dfrac{z^2+9}{2+\dfrac{z^2+25}{2+\dfrac{z^2+49}{2+\dfrac{z^2+81}{2+\dfrac{z^2+121}{2+\ddots}}}}}}$$
Convergence type $P^-$ with $P=1$ and $C=z^2(1+\coth(\pi z/4)^2)/4$, so that
$$z\coth(\pi z/4)-\dfrac{p(n)}{q(n)}\sim(-1)^n\dfrac{(z^2/4)(1+\coth(\pi z/4)^2)}{n}\;.$$
$$A=1-1/n+(-z^2/16+3/4)/n^2+(3z^2/16-1/4)/n^3+\cdots$$
Parametric family for $k\ge0$:
\begin{verbatim}
[(z)->z*cotanh(Pi*z/4),4*k+2,(2*n+1)^2+z^2]
\end{verbatim}
Convergence type $P^-$ with $P=2k+1$.
\end{cf}

\smallskip

\begin{cf}\label{6.6.34}{\ }
\begin{verbatim}
[(z)->z*cotanh(Pi*z/4),[1,3,4],[z^2+1,4*n^2+z^2]]
\end{verbatim}
$$z\coth(\pi z/4)=1+\dfrac{z^2+1}{3+\dfrac{z^2+4}{4+\dfrac{z^2+16}{4+\dfrac{z^2+36}{4+\dfrac{z^2+64}{4+\dfrac{z^2+100}{4+\ddots}}}}}}$$
Convergence type $P^-$ with $P=2$ and $C=z(z^2+1)\coth(\pi z/4)/8$, so that
$$z\coth(\pi z/4)-\dfrac{p(n)}{q(n)}\sim(-1)^n\dfrac{z(z^2+1)\coth(\pi z/4)/8}{n^2}\;.$$
$$A=1-1/n+(-z^2/8+1/8)/n^2+(z^2/4+3/4)/n^3+\cdots$$
Series:
\begin{align*}z\coth(\pi z/4)&=1+4\dfrac{z^2+1}{z^2+16}\sum_{n\ge0}\dfrac{(1/2-iz/4)_n(1/2+iz/4)_n}{(2-iz/4)_n(2+iz/4)_n}\\
\dfrac{\tanh(\pi z/4)}{z}&=-\dfrac{1}{z^2}+4\dfrac{z^2+1}{z^2(z^2+4)}\sum_{n\ge0}\dfrac{(iz/4)_n(-iz/4)_n}{(3/2+iz/4)_n(3/2-iz/4)_n}\end{align*}
Parametric family for $k\ge0$:
\begin{verbatim}
[(z)->z*cotanh(Pi*z/4),4*k+4,4*n^2+z^2]
\end{verbatim}
Convergence type $P^-$ with $P=2k+2$.
\end{cf}

\smallskip

\begin{cf}\label{6.6.35}{\ }
\begin{verbatim}
[(z)->z*cotanh(Pi*z/4),[2*n+1],[z^2+(n+1)^2]]
\end{verbatim}
$$z\coth(\pi z/4)=1+\dfrac{z^2+1}{3+\dfrac{z^2+4}{5+\dfrac{z^2+9}{7+\dfrac{z^2+16}{9+\dfrac{z^2+25}{11+\dfrac{z^2+36}{13+\ddots}}}}}}$$
Convergence type $E$ with $E=-(1+\sqrt{2})^2$, $P=0$, and $C=...$, so that
$$z\coth(\pi z/4)-\dfrac{p(n)}{q(n)}\sim(-1)^n\dfrac{C}{(1+\sqrt{2})^{2n}}\;.$$
$$A=1-((z^2/2+1/8)d)/n+(z^4/4+(3d/4+1/8)z^2+3d/16+1/64)/n^2+\cdots$$
\end{cf}

\smallskip

\begin{cf}\label{6.6.36}{\ }
\begin{verbatim}
[(z)->z*cotanh(Pi*z/4),[0,2*n+1],[[4*(n+1)^2+z^2,(2*n+1)^2+z^2]]]
\end{verbatim}
$$z\coth(\pi z/4)=\dfrac{z^2+4}{3+\dfrac{z^2+1}{5+\dfrac{z^2+16}{7+\dfrac{z^2+9}{9+\dfrac{z^2+36}{11+\dfrac{z^2+25}{13+\ddots}}}}}}$$
Convergence type $E$ with $E=-(1+\sqrt{2})^2$, $P=0$, and $C=...$, so that
$$z\coth(\pi z/4)-\dfrac{p(n)}{q(n)}\sim(-1)^n\dfrac{C}{(1+\sqrt{2})^{2n}}\;.$$
$$A=1+((-z^2/2+15/8)d)/n+(z^4/4+(3d/4-15/8)z^2-45d/16+225/64)/n^2+\cdots$$
\end{cf}

\smallskip

\begin{cf}\label{6.6.37}{\ }
\begin{verbatim}
[(z)->(Pi*z/2)*cotanh(Pi*z/2),[1,2*n-1],[z^2,n^2*(n^2+z^2)]]
\end{verbatim}
$$(\pi z/2)\coth(\pi z/2)=1+\dfrac{z^2}{1+\dfrac{z^2+1}{3+\dfrac{4z^2+16}{5+\dfrac{9z^2+81}{7+\dfrac{16z^2+256}{9+\dfrac{25z^2+625}{11+\ddots}}}}}}$$
Convergence type $P^-$ with $P=2$ and $C=(\pi/4)z^3\cotanh(\pi z/2)$, so that
$$(\pi z/2)\coth(\pi z/2)-\dfrac{p(n)}{q(n)}\sim(-1)^n\dfrac{(\pi/4)z^3\cotanh(\pi z/2)}{n^2}\;.$$
$$A=1-1/n-(z^2/2)/n^2+(z^2+1)/n^3+(5z^4/16)/n^4+\cdots$$
Series:
\begin{align*}(\pi z/2)\coth(\pi z/2)&=1+\dfrac{z^2}{z^2+4}\sum_{n\ge0}\dfrac{n!^2(1/2-iz/2)_n(1/2+iz/2)_n(4n+3)}{(3/2)_n^2(2-iz/2)_n(2+iz/2)_n}\\
\dfrac{\tanh(\pi z/2)}{\pi z/2}&=\dfrac{1}{z^2+1}\sum_{n\ge0}\dfrac{(4n+1)(1/2)_n^2(iz/2)_n(-iz/2)_n}{(3/2-iz/2)_n(3/2+iz/2)_nn!^2}\end{align*}
Parametric family for $k\ge0$:
\begin{verbatim}
[(z)->(Pi*z/2)*cotanh(Pi*z/2),(2*k+1)*(2*n-1),n^2*(n^2+z^2)]
\end{verbatim}
Convergence type $P^-$ with $P=4k+2$.
\end{cf}

\smallskip

\begin{cf}\label{6.6.38}{\ }
\begin{verbatim}
[(z)->(Pi*z/2)*cotanh(Pi*z/2),
[1,z^2+4,8*n^2-8*n+2*z^2+4],[2*z^2,-(4*n^2+z^2)^2]]
\end{verbatim}
$$(\pi z/2)\coth(\pi z/2)=1+\dfrac{2z^2}{z^2+4-\dfrac{z^4+8z^2+16}{2z^2+20-\dfrac{z^4+32z^2+256}{2z^2+52-\dfrac{z^4+72z^2+1296}{2z^2+100-\ddots}}}}$$
Convergence type $P^+$ with $P=1$ and $C=z^2/2$, so that
$$(\pi z/2)\coth(\pi z/2)-\dfrac{p(n)}{q(n)}\sim\dfrac{z^2/2}{n}\;.$$
$$A=1-(1/2)/n+(-z^2/12+1/6)/n^2+(z^2/8)/n^3+\cdots$$
Series:
$$(\pi z/2)\coth(\pi z/2)=1+2z^2\sum_{n\ge1}\dfrac{1}{4n^2+z^2}$$
Parametric family for $k\ge0$:
\begin{verbatim}
[(z)->(Pi*z/2)*cotanh(Pi*z/2),8*n^2-8*n+2*z^2+4*(k^2+k+1),-(4*n^2+z^2)^2]
\end{verbatim}
Convergence type $P^-$ with $P=2k+1$.
\end{cf}
      
\smallskip

\begin{cf}\label{6.6.39}{\ }
\begin{verbatim}
[(z)->(Pi*z/2)*cotanh(Pi*z/2),
[[1,z^2/2+2],[z^2/2+5*n^2,z^2/2+5*n^2+6*n+2]],
[[z^2,-(z^2/2+2)^2],[n^2*(n^2+z^2),-(z^2/2+2*(n+1)^2)^2]]]
\end{verbatim}
$$(\pi z/2)\coth(\pi z/2)=1+\dfrac{z^2}{z^2/2+2-\dfrac{z^4/4+2z^2+4}{z^2/2+5+\dfrac{z^2+1}{z^2/2+13-\dfrac{z^4/4+8z^2+64}{z^2/2+20+\ddots}}}}$$
Convergence type $E$ with $E=-((1+\sqrt{5})/2)^5$, $P=0$, and $C=...$, so that
$$(\pi z/2)\coth(\pi z/2)-\dfrac{p(n)}{q(n)}\sim(-1)^n\dfrac{C}{((1+\sqrt{5})/2)^{5n}}\;.$$
$$A=1-((z^2+2/5)d)/n+((5/2)z^4+(d+2)z^2+(2/5)(d+1))/n^2+\cdots$$
\end{cf}

\smallskip

\begin{cf}\label{6.6.40}{\ }
\begin{verbatim}
[(z)->(Pi*z/2)*cotanh(Pi*z/2),
[[1,1],[5*n^2+z^2/2,5*n^2+4*n+1+z^2/2]],
[[z^2,z^2+1],[-(2*n^2+z^2/2)^2,(n+1)^2*((n+1)^2+z^2)]]]
\end{verbatim}
$$(\pi z/2)\coth(\pi z/2)=1+\dfrac{z^2}{1+\dfrac{z^2+1}{z^2/2+5-\dfrac{z^4/4+2z^2+4}{z^2/2+10+\dfrac{4z^2+16}{z^2/2+20-\dfrac{z^4/4+8z^2+64}{z^2/2+29+\ddots}}}}}$$
Convergence type $E$ with $E=-((1+\sqrt{5})/2)^5$, $P=0$, and $C=...$, so that
$$(\pi z/2)\coth(\pi z/2)-\dfrac{p(n)}{q(n)}\sim(-1)^n\dfrac{C}{((1+\sqrt{5})/2)^{5n}}\;.$$
$$A=1-((z^2+2/5)d)/n+((5/2)z^4+(d+2)z^2+(2/5)(d+1))/n^2+\cdots$$
\end{cf}

\smallskip

\begin{verbatim}
T(a,b)=(a*tanh(Pi*b/2)-b*tanh(Pi*a/2))/(a*tanh(Pi*a/2)-b*tanh(Pi*b/2));
\end{verbatim}

Note that
$$\dfrac{(a-b)}{(a+b)}\dfrac{\sinh(\pi(a+b)/2)}{\sinh(\pi(a-b)/2)}=-1-\dfrac{2}{-1+T(a,b)}$$
with $T(a,b)$ as above, so the next CFs give corresponding CFs for
$\sinh(\pi(a+b)/2)/\sinh(\pi(a-b)/2)$.

\smallskip

\begin{cf}\label{6.6.41}{\ }
\begin{verbatim}
[(a,b)->T(a,b),[0,2*n-1],[a*b,(a^2+n^2)*(b^2+n^2)]]
\end{verbatim}
$$\dfrac{a\tanh(\pi b/2)-b\tanh(\pi a/2)}{a\tanh(\pi a/2)-b\tanh(\pi b/2)}=\dfrac{ba}{1+\dfrac{(b^2+1)a^2+b^2+1}{3+\dfrac{(b^2+4)a^2+4b^2+16}{5+\dfrac{(b^2+9)a^2+9b^2+81}{7+\dfrac{(b^2+16)a^2+16b^2+256}{9+\dfrac{(b^2+25)a^2+25b^2+625}{11+\ddots}}}}}}$$
Convergence type $P^-$ with $P=2$ and $C=...$, so that
$$\dfrac{a\tanh(\pi b/2)-b\tanh(\pi a/2)}{a\tanh(\pi a/2)-b\tanh(\pi b/2)}-\dfrac{p(n)}{q(n)}\sim(-1)^n\dfrac{C}{n^2}\;.$$
$$A=1-1/n-((a^2+b^2)/2)/n^2+(a^2+b^2+1)/n^3+\cdots$$
Parametric family for $k\ge0$:
\begin{verbatim}
[(a,b)->T(a,b),(2*k+1)*(2*n-1),(a^2+n^2)*(b^2+n^2)]
\end{verbatim}
Convergence type $P^-$ with $P=4k+2$.
\end{cf}

\smallskip

\begin{cf}\label{6.6.42}{\ }
\begin{verbatim}
[(a,b)->T(a,b),[0,a^2+b^2+4,8*n^2-8*n+4+2*(a^2+b^2)],
               [2*a*b,-(4*n^2+(a-b)^2)*(4*n^2+(a+b)^2)]]
\end{verbatim}
$$\dfrac{a\tanh(\pi b/2)-b\tanh(\pi a/2)}{a\tanh(\pi a/2)-b\tanh(\pi b/2)}=\dfrac{2ab}{a^2+b^2+4-\dfrac{a^4+(-2b^2+8)a^2+b^4+8b^2+16}{2a^2+2b^2+20-\ddots}}$$
Convergence type $P^+$ with $P=1$ and $C=...$, so that
$$\dfrac{a\tanh(\pi b/2)-b\tanh(\pi a/2)}{a\tanh(\pi a/2)-b\tanh(\pi b/2)}-\dfrac{p(n)}{q(n)}\sim\dfrac{C}{n}\;.$$
Parametric family for $k\ge0$:
\begin{verbatim}
[(a,b)->T(a,b),8*n^2-8*n+2*(a^2+b^2)+4*(k^2+k+1),
               -(4*n^2+(a-b)^2)*(4*n^2+(a+b)^2)]]
\end{verbatim}
Convergence type $P^+$ with $P=2k+1$.
\end{cf}

\smallskip

\begin{cf}\label{6.6.43}{\ }
\begin{verbatim}
[(a,b)->T(a,b),[[0,a^2+b^2+4],[10*n^2+a^2+b^2,10*n^2+12*n+4+a^2+b^2]],
[[2*a*b,-((a-b)^2+4)*((a+b)^2+4)],
 [4*(n^2+a^2)*(n^2+b^2),-(4*(n+1)^2+(a-b)^2)*(4*(n+1)^2+(a+b)^2)]]]
\end{verbatim}
$$\dfrac{a\tanh(\pi b/2)-b\tanh(\pi a/2)}{a\tanh(\pi a/2)-b\tanh(\pi b/2)}=\dfrac{2ab}{a^2+b^2+4-\dfrac{a^4+(-2b^2+8)a^2+b^4+8b^2+16}{a^2+b^2+10+\ddots}}$$
Convergence type $E$ with $E=-((1+\sqrt{5})/2)^5$, $P=0$, and $C=...$, so that
$$\dfrac{a\tanh(\pi b/2)-b\tanh(\pi a/2)}{a\tanh(\pi a/2)-b\tanh(\pi b/2)}-\dfrac{p(n)}{q(n)}\sim(-1)^n\dfrac{C}{((1+\sqrt{5})/2)^{5n}}\;.$$
\end{cf}

\smallskip

\begin{cf}\label{6.6.44}{\ }
\begin{verbatim}
[(a,z)->tanh(a*atanh(z)),[0,2*n-1],[a*z,(a^2-n^2)*z^2]]
\end{verbatim}
$$\tanh(a\atanh(z))=\dfrac{az}{1+\dfrac{(a^2-1)z^2}{3+\dfrac{(a^2-4)z^2}{5+\dfrac{(a^2-9)z^2}{7+\dfrac{(a^2-16)z^2}{9+\dfrac{(a^2-25)z^2}{11+\ddots}}}}}}$$
Convergence type $E$ with $E=((1+\sqrt{1-z^2})/z)^2$, $P=0$, and $C=...$,
so that
$$\tanh(a\atanh(z))-\dfrac{p(n)}{q(n)}\sim\dfrac{C}{((1+\sqrt{1-z^2})/z)^{2n}}\;.$$

\end{cf}

\smallskip

\begin{cf}\label{6.6.45}{\ }
\begin{verbatim}
[(a,z)->sinh(a*asinh(z)),[a*z,6,2*n*(2*n+1)-z^2*((2*n-1)^2-a^2)],
                     [a*(a^2-1)*z^3,2*n*(2*n+1)*((2*n+1)^2-a^2)*z^2]]
\end{verbatim}
$$\sinh(a\asinh(z))=az+\dfrac{(a^3-a)z^3}{6-\dfrac{(6a^2-54)z^2}{(a^2-9)z^2+20-\dfrac{(20a^2-500)z^2}{(a^2-25)z^2+42-\ddots}}}$$
Convergence type $E$ with $E=-1/z^2$, $P=3/2$, and $C=...$, so that
$$\sinh(a\asinh(z))-\dfrac{p(n)}{q(n)}\sim(-1)^n\dfrac{C}{(1/z)^{2n}n^{3/2}}\;,$$
this for $|z|\le 1$. If $|z|>1$, change $E$ into $1/E$ and $P$ into $-P$.
$$A=1+(((5-2a^2)z^2+17-2a^2)/(8(z^2+1)))/n+\cdots$$
Series:
$$\sinh(a\asinh(z))=az+\dfrac{az^3(a^2-1)}{6}\sum_{n\ge0}(-1)^n\dfrac{(3/2-a/2)_n(3/2+a/2)_n}{(n+1)!(5/2)_n}z^{2n}$$
\end{cf}

\smallskip

\begin{cf}\label{6.6.46}{\ }
\begin{verbatim}
[(a,z)->cosh(a*asinh(z)),
[1,1,n*(2*n-1)-2*z^2*((n-1)^2-a^2/4)],z^2*[a^2/2,2*n*(2*n-1)*(n^2-a^2/4)]]
\end{verbatim}
$$\cosh(a\asinh(z))=1+\dfrac{(1/2)a^2z^2}{1-\dfrac{((1/2)a^2-2)z^2}{((1/2)a^2-2)z^2+6-\dfrac{(3a^2-48)z^2}{((1/2)a^2-8)z^2+15-\ddots}}}$$
Convergence type $E$ with $E=-1/z^2$, $P=3/2$, and $C=...$, so that
$$\cosh(a\asinh(z))-\dfrac{p(n)}{q(n)}\sim(-1)^n\dfrac{C}{(1/z)^{2n}n^{3/2}}\;,$$
this for $|z|\le1$. If $|z|>1$, change $E$ into $1/E$ and $P$ into $-P$;
$$A=1-(((2a^2+1)z^2+2a^2-11)/(8(z^2+1)))/n+\cdots$$
Series:
$$\cosh(a\asinh(z))=1+\dfrac{a^2z^2}{2}\sum_{n\ge0}(-1)^n\dfrac{(1-a/2)_n(1+a/2)_n}{(n+1)!(3/2)_n}z^{2n}$$
\end{cf}

Since $\cosh(2x)=2\sinh^2(x)+1=2\cosh^2(x)-1$, this also gives CFs for
$\sinh^2(a\asinh(z))$ and $\cosh^2(a\sinh(z))$.

\smallskip

\begin{verbatim}
S5(z)=(sinh(Pi*z)-sin(Pi*z))/(sinh(Pi*z)+sin(Pi*z));
\end{verbatim}

\smallskip

\begin{cf}\label{6.6.47}{\ }
\begin{verbatim}
[(z)->S5(z),[0,2*n-1],[2*z^2,4*z^4+n^4]]
\end{verbatim}
$$\dfrac{\sinh(\pi z)-\sin(\pi z)}{\sinh(\pi z)+\sin(\pi z)}=\dfrac{2z^2}{1+\dfrac{4z^4+1}{3+\dfrac{4z^4+16}{5+\dfrac{4z^4+81}{7+\dfrac{4z^4+256}{9+\dfrac{4z^4+625}{11+\ddots}}}}}}$$
Convergence type $P^-$ with $P=2$ and $C=...$, so that
$$\dfrac{\sinh(\pi z)-\sin(\pi z)}{\sinh(\pi z)+\sin(\pi z)}-\dfrac{p(n)}{q(n)}\sim(-1)^n\dfrac{C}{n^2}\;.$$
$$A=1-1/n+1/n^3-z^4/n^4+3(z^4-1)/n^5+z^4/n^6+\cdots$$
Parametric family for $k\ge0$:
\begin{verbatim}
[(z)->S5(z),(2*k+1)*(2*n-1),4*z^4+n^4]
\end{verbatim}
Convergence type $P^-$ with $P=4k+2$.
\end{cf}

\smallskip

\begin{cf}\label{6.6.48}{\ }
\begin{verbatim}
[(z)->S5(z),[0,2*n^2-2*n+1],[z^2,-n^4+z^4]]
\end{verbatim}
$$\dfrac{\sinh(\pi z)-\sin(\pi z)}{\sinh(\pi z)+\sin(\pi z)}=\dfrac{z^2}{1+\dfrac{z^4-1}{5+\dfrac{z^4-16}{13+\dfrac{z^4-81}{25+\dfrac{z^4-256}{41+\dfrac{z^4-625}{61+\ddots}}}}}}$$
Convergence type $P^+$ with $P=1$ and $C=...$, so that
$$\dfrac{\sinh(\pi z)-\sin(\pi z)}{\sinh(\pi z)+\sin(\pi z)}-\dfrac{p(n)}{q(n)}\sim\dfrac{C}{n}\;.$$
$A=1+\cdots$

\noindent
Parametric family for $k\ge0$:
\begin{verbatim}
[(z)->S5(z),2*n^2-2*n+k^2+k+1,-n^4+z^4]
\end{verbatim}
Convergence type $P^+$ with $P=2k+1$.
\end{cf}

\smallskip

\begin{cf}\label{6.6.49}{\ }
\begin{verbatim}
[(z)->S5(z),[[5*n^2,5*n^2+6*n+2]],
         [[2*z^2,4*z^4-4],[n^4+4*z^4,-4*((n+1)^4-z^4)]]]
\end{verbatim}
$$\dfrac{\sinh(\pi z)-\sin(\pi z)}{\sinh(\pi z)+\sin(\pi z)}=\dfrac{2z^2}{2+\dfrac{4z^4-4}{5+\dfrac{4z^4+1}{13+\dfrac{4z^4-64}{20+\dfrac{4z^4+16}{34+\dfrac{4z^4-324}{45+\ddots}}}}}}$$
Convergence type $E$ with $E=-((1+\sqrt{5})/2)^5$, $P=0$, and $C=...$, so that
$$\dfrac{\sinh(\pi z)-\sin(\pi z)}{\sinh(\pi z)+\sin(\pi z)}-\dfrac{p(n)}{q(n)}\sim(-1)^n\dfrac{C}{((1+\sqrt{5})/2)^{5n}}\;.$$
$$A=1-(2d/5)/n+(2d/5+2/5)/n^2+((-8d/3)z^4-(204d/625+4/5))/n^3+\cdots$$
\end{cf}

\medskip

\section{Inverse Trigonometric and Hyperbolic Functions}

\medskip

Since $\asin(z)=\asinh(iz)/i$, $\acos(z)=\acosh(z)/i$ (\emph{not} $iz$
here), and $\atan(z)=\atanh(iz)/i$, all the continued fractions for inverse
trigonometric functions can be deduced trivially from those of the inverse
hyperbolic functions, so we only give the latter.

\smallskip

\begin{cf}\label{6.7.1}{\ }
\begin{verbatim}
[(z)->asinh(z),[z,6,-(2*n-1)^2*z^2+2*n*(2*n+1)],[-z^3,2*n*(2*n+1)^3*z^2]]
\end{verbatim}
$$\asinh(z)=z-\dfrac{z^3}{6+\dfrac{54z^2}{-9z^2+20+\dfrac{500z^2}{-25z^2+42+\dfrac{2058z^2}{-49z^2+72+\dfrac{5832z^2}{-81z^2+110+\ddots}}}}}$$
Convergence type $E$ with $E=-1/z^2$, $P=3/2$, and 
$C=-z^3/(2(1+z^2)\sqrt{\pi})$, so that
$$\asinh(z)-\dfrac{p(n)}{q(n)}\sim(-1)^{n+1}\dfrac{1/(2(1+z^2)\sqrt{\pi})}{(1/z)^{2n+3}n^{3/2}}$$
$$A=1-((5z^2+17)/(8(z^2+1)))/n+\cdots$$
Series:
$$\asinh(z)=z\sum_{n\ge0}(-1)^n\dfrac{(1/2)_n}{(2n+1)n!}z^{2n}$$
\end{cf}

\smallskip

This is the CF corresponding to the Taylor expansion of $\asinh(z)$.

\smallskip

\begin{cf}\label{6.7.2}{\ }
\begin{verbatim}
[(z)->asinh(z)/sqrt(1+z^2),[z,3,1+2*n*(1-z^2)],[-2*z^3,2*z^2*(n+1)*(2*n+1)]]
\end{verbatim}
$$\dfrac{\asinh(z)}{\sqrt{1+z^2}}=z-\dfrac{2z^3}{3+\dfrac{12z^2}{-4z^2+5+\dfrac{30z^2}{-6z^2+7+\dfrac{56z^2}{-8z^2+9+\dfrac{90z^2}{-10z^2+11+\ddots}}}}}$$
Convergence type $E$ with $E=-1/z^2$, $P=1/2$, and
$C=-\sqrt{\pi}z^3/(2(1+z^2))$, so that
$$\dfrac{\asinh(z)}{\sqrt{1+z^2}}-\dfrac{p(n)}{q(n)}\sim(-1)^{n+1}\dfrac{\sqrt{\pi}/(2(1+z^2))}{(1/z)^{2n+3}n^{1/2}}\;.$$
$$A=1-((3z^2+7)/(8z^2+8))/n+\cdots$$
Series:
$$\dfrac{\asinh(z)}{\sqrt{1+z^2}}=z-\dfrac{2z^3}{3}\sum_{n\ge0}(-1)^n\dfrac{(n+1)!}{(5/2)_n}z^{2n}$$
Parametric family for $k\ge0$:
\begin{verbatim}
[(z)->asinh(z)/sqrt(1+z^2),1+2*n*(1-z^2)+2*k*(z^2+1),2*z^2*(n+1)*(2*n+1)]
\end{verbatim}
Convergence type $E$ with $E=-1/z^2$ and $P=2k+1/2$.
\end{cf}

\smallskip

\begin{cf}\label{6.7.3}{\ }
\begin{verbatim}
[(z)->asinh(z)/sqrt(1+z^2),[z,2*n*(1-z^2)+4*z^2+1],[-2*z^3,2*n*(2*n-1)*z^2]]
\end{verbatim}
$$\dfrac{\asinh(z)}{\sqrt{1+z^2}}=z-\dfrac{2z^3}{2z^2+3+\dfrac{2z^2}{5+\dfrac{12z^2}{-2z^2+7+\dfrac{30z^2}{-4z^2+9+\dfrac{56z^2}{-6z^2+11+\ddots}}}}}$$
Convergence type $E$ with $E=-1/z^2$, $P=5/2$, and
$C=-\sqrt{\pi}z^3/(4(1+z^2)^3)$, so that
$$\dfrac{\asinh(z)}{\sqrt{1+z^2}}-\dfrac{p(n)}{q(n)}\sim(-1)^{n+1}\dfrac{\sqrt{\pi}/(4(1+z^2)^3)}{(1/z)^{2n+3}n^{5/2}}\;.$$
$$A=1+((9z^2-19)/(8z^2+8))/n+\cdots$$
Series:
$$\dfrac{\asinh(z)}{\sqrt{1+z^2}}=z-2z^3\sum_{n\ge0}(-1)^n\dfrac{n!}{(2(z^2+1)n+1)(2(z^2+1)(n+1)+1)(3/2)_n}z^{2n}$$
Parametric family for $k\ge0$:
\begin{verbatim}
[(z)->asinh(z)/sqrt(1+z^2),2*n*(1-z^2)+4*z^2+1+2*k*(z^2+1),2*n*(2*n-1)*z^2]
\end{verbatim}
Convergence type $E$ with $E=-1/z^2$ and $P=2k+5/2$.
\end{cf}
    
\smallskip

\begin{cf}\label{6.7.4}{\ }
\begin{verbatim}
[(z)->asinh(z)/sqrt(1+z^2),[0,2*n-1],[z/(1+z^2),-z^2/(1+z^2)*n^2]]
\end{verbatim}
$$\dfrac{\asinh(z)}{\sqrt{1+z^2}}=\dfrac{z/(z^2+1)}{1-\dfrac{z^2/(z^2+1)}{3-\dfrac{4z^2/(z^2+1)}{5-\dfrac{9z^2/(z^2+1)}{7-\dfrac{16z^2/(z^2+1)}{9-\dfrac{25z^2/(z^2+1)}{11-\ddots}}}}}}$$
Convergence type $E$ with $E=(1+\sqrt{1+z^2})^2/z^2$, $P=0$, and $C=...$,
so that
$$\dfrac{\asinh(z)}{\sqrt{1+z^2}}-\dfrac{p(n)}{q(n)}\sim\dfrac{C}{((1+\sqrt{1+z^2})/z)^{2n}}\;.$$
\end{cf}

\smallskip

\begin{cf}\label{6.7.5}{\ }
\begin{verbatim}
[(z)->asinh(z)/sqrt(1+z^2),[0,2*n-1],
                         [[z,2*z^2],2*z^2*[n*(2*n-1),(n+1)*(2*n+1)]]]
\end{verbatim}
$$\dfrac{\asinh(z)}{\sqrt{1+z^2}}=\dfrac{z}{1+\dfrac{2z^2}{3+\dfrac{2z^2}{5+\dfrac{12z^2}{7+\dfrac{12z^2}{9+\dfrac{30z^2}{11+\ddots}}}}}}$$
Convergence type $E$ with $E=-(1+\sqrt{1+z^2})^2/z^2$, $P=0$, and $C=...$,
so that
$$\dfrac{\asinh(z)}{\sqrt{1+z^2}}-\dfrac{p(n)}{q(n)}\sim(-1)^n\dfrac{C}{((1+\sqrt{1+z^2})/z)^{2n}}\;.$$
\end{cf}

\smallskip

\begin{cf}\label{6.7.6}{\ }
\begin{verbatim}
[(z)->asinh(z)/sqrt(1+z^2),[z,2*z^2+3,2*n+1],
                         [[-2*z^3,2*z^2],2*z^2*(n+1)*(2*n+1)*[1,1]]]
\end{verbatim}
$$\dfrac{\asinh(z)}{\sqrt{1+z^2}}=z-\dfrac{2z^3}{2z^2+3+\dfrac{2z^2}{5+\dfrac{12z^2}{7+\dfrac{12z^2}{9+\dfrac{30z^2}{11+\dfrac{30z^2}{13+\ddots}}}}}}$$
Convergence type $E$ with $E=-(1+\sqrt{1+z^2})^2/z^2$, $P=0$, and $C=...$,
so that
$$\dfrac{\asinh(z)}{\sqrt{1+z^2}}-\dfrac{p(n)}{q(n)}\sim(-1)^n\dfrac{C}{((1+\sqrt{1+z^2})/z)^{2n}}\;.$$
\end{cf}

\smallskip

\begin{cf}\label{6.7.7}{\ }
\begin{verbatim}
[(z)->acosh(z)/sqrt(z^2-1),[z,3,-2*n*(z^2-2)+1],
                           (z^2-1)*[-2*z,2*(n+1)*(2*n+1)]]
\end{verbatim}
$$\dfrac{\acosh(z)}{\sqrt{z^2-1}}=z-\dfrac{2z^3-2z}{3+\dfrac{12z^2-12}{-4z^2+9+\dfrac{30z^2-30}{-6z^2+13+\dfrac{56z^2-56}{-8z^2+17+\ddots}}}}$$
Convergence type $E$ with $E=-1/(z^2-1)$, $P=1/2$, and
$C=-(z^2-1)\sqrt{\pi}/(2z)$, so that
$$\dfrac{\acosh(z)}{\sqrt{z^2-1}}-\dfrac{p(n)}{q(n)}\sim(-1)^{n+1}\dfrac{\sqrt{\pi}/(2z)}{(1/(z^2-1))^{n+1}n^{1/2}}$$
$$A=1-((3z^2+4)/(8z^2))/n+((25z^4+24z^2+96)/(128z^4))/n^2+\cdots$$
Series:
$$\dfrac{\acosh(z)}{\sqrt{z^2-1}}=z-\dfrac{2z(z^2-1)}{3}\sum_{n\ge0}(-1)^n\dfrac{(n+1)!}{(5/2)_n}(z^2-1)^n$$
Parametric family for $k\ge0$:
\begin{verbatim}
[(z)->acosh(z)/sqrt(z^2-1),-2*n*(z^2-2)+2*k*z^2+1,(z^2-1)*2*(n+1)*(2*n+1)]
\end{verbatim}
Convergence type $E$ with $E=-1/(z^2-1)$ and $P=2k+1/2$.
\end{cf}

\smallskip

\begin{cf}\label{6.7.8}{\ }
\begin{verbatim}
[(z)->acosh(z)/sqrt(z^2-1),[0,2*n-1],[[z,2*(z^2-1)],
                         2*(z^2-1)*[n*(2*n-1),(n+1)*(2*n+1)]]]
\end{verbatim}
$$\dfrac{\acosh(z)}{\sqrt{z^2-1}}=\dfrac{z}{1+\dfrac{2z^2-2}{3+\dfrac{2z^2-2}{5+\dfrac{12z^2-12}{7+\dfrac{12z^2-12}{9+\dfrac{30z^2-30}{11+\ddots}}}}}}$$
Convergence type $E$ with $E=-i(z+1)/(z-1)$, $P=0$, and $C=\pi/(z+1)$, so that
$$\dfrac{\acosh(z)}{\sqrt{z^2-1}}-\dfrac{p(n)}{q(n)}\sim(-1)^{\lfloor n/2\rfloor}\dfrac{\pi/(z+1)}{((z+1)/(z-1))^n}\;.$$
$$A=1-(z/4)/n+(z^2/32+z/8)/n^2+(-17z^3/128-z^2/32+5z/64)/n^3+\cdots$$
\end{cf}

\smallskip

\begin{cf}\label{6.7.9}{\ }
\begin{verbatim}
[(z)->acosh(z)/sqrt(z^2-1),[0,z*(2*n-1)],[1,-(z^2-1)*n^2]]
\end{verbatim}
$$\dfrac{\acosh(z)}{\sqrt{z^2-1}}=\dfrac{1}{z-\dfrac{z^2-1}{3z-\dfrac{4z^2-4}{5z-\dfrac{9z^2-9}{7z-\dfrac{16z^2-16}{9z-\dfrac{25z^2-25}{11z-\ddots}}}}}}$$
Convergence type $E$ with $E=(z+1)/(z-1)$, $P=0$, and $C=\pi/(z+1)$, so that
$$\dfrac{\acosh(z)}{\sqrt{z^2-1}}-\dfrac{p(n)}{q(n)}\sim\dfrac{\pi/(z+1)}{((z+1)/(z-1))^n}\;.$$
$$A=1-(z/4)/n+(z^2/32+z/8)/n^2-(17z^3/128+z^2/32-5z/64)/n^3+\cdots$$
Parametric family for $k\ge0$:
\begin{verbatim}
[(z)->acosh(z)/sqrt(z^2-1),z*(2*n-1)+2*k,-(z^2-1)*n^2]
\end{verbatim}
Convergence type $E$ with $E=(z+1)/(z-1)$ and $P=2k$.
\end{cf}

\smallskip

\begin{cf}\label{6.7.10}{\ }
\begin{verbatim}
[(z)->acosh(z)/sqrt(z^2-1),[0,z+1,4*z*(n-1)+2],[2,-(z^2-1)*(2*n-1)^2]]
\end{verbatim}
$$\dfrac{\acosh(z)}{\sqrt{z^2-1}}=\dfrac{2}{z+1-\dfrac{z^2-1}{4z+2-\dfrac{9z^2-9}{8z+2-\dfrac{25z^2-25}{12z+2-\dfrac{49z^2-49}{16z+2-\ddots}}}}}$$
Convergence type $E$ with $E=(z+1)/(z-1)$, $P=1$, and $C=1/2$, so that
$$\dfrac{\acosh(z)}{\sqrt{z^2-1}}-\dfrac{p(n)}{q(n)}\sim\dfrac{1/2}{((z+1)/(z-1))^n\cdot n}\;.$$
$$A=1-(z/2)/n+(z^2/2-1/4)/n^2+(-3z^3/4+5z/8)/n^3+\cdots$$
Series:
$$\dfrac{\acosh(z)}{\sqrt{z^2-1}}=\dfrac{2}{z+1}\sum_{n\ge0}\dfrac{((z-1)/(z+1))^n}{2n+1}$$
Parametric family for $k\ge0$:
\begin{verbatim}
[(z)->acosh(z)/sqrt(z^2-1),4*z*(n-1)+4*k+2,-(z^2-1)*(2*n-1)^2]
\end{verbatim}
Convergence type $E$ with $E=(z+1)/(z-1)$ and $P=2k+1$.
\end{cf}

\smallskip

\begin{cf}\label{6.7.11}{\ }
\begin{verbatim}
[(z)->acosh(z)/sqrt(z^2-1),[0,(z+1)*(2*n-1)],[2,-(z^2-1)*n^2]]
\end{verbatim}
$$\dfrac{\acosh(z)}{\sqrt{z^2-1}}=\dfrac{2}{z+1-\dfrac{z^2-1}{3z+3-\dfrac{4z^2-4}{5z+5-\dfrac{9z^2-9}{7z+7-\dfrac{16z^2-16}{9z+9-\dfrac{25z^2-25}{11z+11-\ddots}}}}}}$$
Convergence type $E$ with $E=(2+\sqrt{2z+2})^2/(2(z-1))$, $P=0$, and $C=...$,
so that
$$\dfrac{\acosh(z)}{\sqrt{z^2-1}}-\dfrac{p(n)}{q(n)}\sim\dfrac{C}{((2+\sqrt{2z+2})^2/(2(z-1)))^n}\;.$$
\end{cf}

\smallskip

\begin{cf}\label{6.7.12}{\ }
\begin{verbatim}
[(z)->atanh(z),
[z/(1-z^2),3*(1-z^2)^2,(2*z^4-12*z^2+2)*n^2+(z^4+6*z^2+1)*n-2*z^2],
[-2*z^3/(1-z^2),z^2*(1-z^2)^2*2*n*(2*n+1)^3]]
\end{verbatim}
$$\atanh(z)=z/(1-z^2)-\dfrac{2z^3/(1-z^2)}{3z^4-6z^2+3+\dfrac{54z^6-108z^4+54z^2}{10z^4-38z^2+10+\dfrac{500z^6-1000z^4+500z^2}{21z^4-92z^2+21+\ddots}}}$$
Convergence type $E$ with $E=-(1-z^2)^2/(4z^2)$, $P=3/2$, and $C=...$
if $|z|<\sqrt{2}-1$, $E=-4*z^2/(1-z^2)^2$, $P=-3/2$, and $C=...$ if
$|z|>\sqrt{2}-1$, so that for instance if $|z|<\sqrt{2}-1$:
$$\atanh(z)-\dfrac{p(n)}{q(n)}\sim(-1)^n\dfrac{C}{((1-z^2)/(2z))^{2n}n^{3/2}}$$
$$A=1-((17z^4-14z^2+17)/(8(1+z^2)^2))/n+\cdots$$
Series:
$$\atanh(z)=\dfrac{z}{1-z^2}-\dfrac{2z^3}{(1-z^2)^3}\sum_{n\ge0}(-1)^n\dfrac{(3/2)_n}{(2n+3)(n+1)!}(2z/(1-z^2))^{2n}$$
\end{cf}
      
\smallskip

\begin{cf}\label{6.7.13}{\ }
\begin{verbatim}
[(z)->atanh(z),[0,(2*n-1)-(4*n-3)*z^2],[z,2*n*(2*n-1)*z^2*(1-z^2)]]
\end{verbatim}
$$\atanh(z)=\dfrac{z}{-z^2+1-\dfrac{2z^4-2z^2}{-5z^2+3-\dfrac{12z^4-12z^2}{-9z^2+5-\dfrac{30z^4-30z^2}{-13z^2+7-\ddots}}}}$$
Convergence type $E$ with $E=1-1/z^2$, $P=1/2$, and $C=\sqrt{\pi}z/2$, so that
$$\atanh(z)-\dfrac{p(n)}{q(n)}\sim\dfrac{\sqrt{\pi}z/2}{(1-1/z^2)^nn^{1/2}}\;.$$
$$A=1+(z^2/2-3/8)/n+(3z^4/4-15z^2/16+25/128)/n^2+\cdots$$
Series:
$$\atanh(z)=\dfrac{z}{1-z^2}\sum_{n\ge0}\dfrac{n!}{(3/2)_n}(z^2/(z^2-1))^n$$
Parametric family for $k\ge0$:
\begin{verbatim}
[(z)->atanh(z),(2*n-1)-(4*n-3)*z^2+2*k,2*n*(2*n-1)*z^2*(1-z^2)]
\end{verbatim}
Convergence type $E$ with $E=1-1/z^2$ and $P=2k+1/2$.
\end{cf}

\smallskip

\begin{cf}\label{6.7.14}{\ }
\begin{verbatim}
[(z)->atanh(z),[0,(z^2+1)*(2*n-1)],[z,-4*z^2*n^2]]
\end{verbatim}
$$\atanh(z)=\dfrac{z}{z^2+1-\dfrac{4z^2}{3z^2+3-\dfrac{16z^2}{5z^2+5-\dfrac{36z^2}{7z^2+7-\dfrac{64z^2}{9z^2+9-\dfrac{100z^2}{11z^2+11-\ddots}}}}}}$$
Convergence type $E$ with $E=1/z^2$, $P=0$, and $C=\pi z/2$, so that
$$\atanh(z)-\dfrac{p(n)}{q(n)}\sim\dfrac{\pi/2}{(1/z)^{2n+1}}\;.$$
$$A=1-((1+z^2)/(4(1-z^2)))/n+((1+z^2)(5-3z^2)/(32(1-z^2)^2))/n^2+\cdots$$
Parametric family:
\begin{verbatim}
[(z)->atanh(z),(z^2+1)*(2*n-1)+2*k*(1-z^2),-4*z^2*n^2]
\end{verbatim}
Convergence type $E$ with $E=1/z^2$ and $P=2k$.
\end{cf}

\smallskip

\begin{cf}\label{6.7.15}{\ }
\begin{verbatim}
[(z)->atanh(z),[0,1,(2*n-3)*z^2+(2*n-1)],[z,-(2*n-1)^2*z^2]]
\end{verbatim}
$$\atanh(z)=\dfrac{z}{1-\dfrac{z^2}{z^2+3-\dfrac{9z^2}{3z^2+5-\dfrac{25z^2}{5z^2+7-\dfrac{49z^2}{7z^2+9-\dfrac{81z^2}{9z^2+11-\ddots}}}}}}$$
Convergence type $E$ with $E=1/z^2$, $P=1$, and $C=z/(2(1-z^2))$, so that
$$\atanh(z)-\dfrac{p(n)}{q(n)}\sim\dfrac{1/(2(1-z^2))}{(1/z)^{2n+1}n}\;.$$
$$A=1-((1+z^2)/(2(1-z^2)))/n+((z^4+6z^2+1)/(4(1-z^2)^2))/n^2+\cdots$$
Series:
$$\atanh(z)=z\sum_{n\ge0}\dfrac{z^{2n}}{2n+1}$$
Parametric family for $k\ge0$:
\begin{verbatim}
[(z)->atanh(z),(2*n-3)*z^2+(2*n-1)+2*k*(1-z^2),-(2*n-1)^2*z^2]
\end{verbatim}
Convergence type $E$ with $E=1/z^2$ and $P=2k+1$.
\end{cf}

This is the CF corresponding to the term-by-term Taylor expansion of
$\atanh(z)$.

\smallskip

\begin{cf}\label{6.7.16}{\ }
\begin{verbatim}
[(z)->atanh(z),[0,2*n-1],[z,-n^2*z^2]]
\end{verbatim}
$$\atanh(z)=\dfrac{z}{1-\dfrac{z^2}{3-\dfrac{4z^2}{5-\dfrac{9z^2}{7-\dfrac{16z^2}{9-\dfrac{25z^2}{11-\ddots}}}}}}$$
Convergence type $E$ with $E=(1+\sqrt{1-z^2})^2/z^2$, $P=0$, and
$C=\pi/\sqrt{E}$, so that
$$\atanh(z)-\dfrac{p(n)}{q(n)}\sim\dfrac{\pi}{((1+\sqrt{1-z^2})/z)^{2n+1}}\;.$$
\end{cf}

\smallskip

\begin{cf}\label{6.7.17}{\ }
\begin{verbatim}
[(z)->atanh(z),[0,z^2+1,2*n-1],[[z,-4*z^2],-z^2*[(2*n-1)^2,(2*n+2)^2]]]
\end{verbatim}
$$\atanh(z)=\dfrac{z}{z^2+1-\dfrac{4z^2}{3-\dfrac{z^2}{5-\dfrac{16z^2}{7-\dfrac{9z^2}{9-\dfrac{36z^2}{11-\ddots}}}}}}$$
Convergence type $E$ with $E=(1+\sqrt{1-z^2})^2/z^2$, $P=0$, and
$C=\pi/\sqrt{E}$, so that
$$\atanh(z)-\dfrac{p(n)}{q(n)}\sim\dfrac{\pi}{((1+\sqrt{1-z^2})/z)^{2n+1}}\;.$$
\end{cf}

\smallskip

\begin{cf}\label{6.7.18}{\ }
\begin{verbatim}
[(z)->atanh(z),[0,(1-z^2)*(2*n-1)],
       [[z,2*z^2*(1-z^2)],2*z^2*(1-z^2)*[n*(2*n-1),(n+1)*(2*n+1)]]]
\end{verbatim}
$$\atanh(z)=\dfrac{z}{-z^2+1-\dfrac{2z^4-2z^2}{-3z^2+3-\dfrac{2z^4-2z^2}{-5z^2+5-\dfrac{12z^4-12z^2}{-7z^2+7-\dfrac{12z^4-12z^2}{-9z^2+9-\ddots}}}}}$$
Convergence type $E$ with $E=-(1+\sqrt{1-z^2})^2/z^2$, $P=0$, and $C=...$, so
that
$$\atanh(z)-\dfrac{p(n)}{q(n)}\sim(-1)^n\dfrac{C}{((1+\sqrt{1-z^2})/z)^{2n}}\;.$$
\end{cf}

\smallskip

\begin{cf}\label{6.7.19}{\ }
\begin{verbatim}
[(z)->asinh(z)^2,[0,-2*z^2*(n-1)^2+n*(2*n-1)],z^2*[1,2*n^3*(2*n-1)]]
\end{verbatim}
$$\asinh^2(z)=\dfrac{z^2}{1+\dfrac{2z^2}{-2z^2+6+\dfrac{48z^2}{-8z^2+15+\dfrac{270z^2}{-18z^2+28+\dfrac{896z^2}{-32z^2+45+\ddots}}}}}$$
For $|z|\le1$: Convergence type $E$ with $E=-1/z^2$, $P=3/2$, and
$C=(z^2/(2(z^2+1)))\sqrt{\pi}$, so that
$$\asinh^2(z)-\dfrac{p(n)}{q(n)}\sim(-1)^n\dfrac{\sqrt{\pi}/(2(z^2+1))}{(1/z)^{2n+2}n^{3/2}}\;.$$
$$A=1-((11-z^2)/(8(1+z^2)))/n+((z^4-278z^2+201)/(128(1+z^2)^2))/n^2+\cdots$$
Series:
$$\asinh^2(z)=z^2\sum_{n\ge0}(-1)^n\dfrac{n!}{(n+1)(3/2)_n}z^{2n}$$
\end{cf}

\smallskip

\begin{cf}\label{6.7.20}{\ }
\begin{verbatim}
[(z)->acosh(z)^2,[0,-2*(z^2-1)*(n-1)^2+n*(2*n-1)],(z^2-1)*[1,2*n^3*(2*n-1)]]
\end{verbatim}
$$\acosh^2(z)=\dfrac{z^2-1}{1+\dfrac{2z^2-2}{-2z^2+8+\dfrac{48z^2-48}{-8z^2+23+\dfrac{270z^2-270}{-18z^2+46+\ddots}}}}$$
For $|z^2-1|\le1$: Convergence type $E$ with $E=-1/(z^2-1)$, $P=3/2$, and
$C=((z^2-1)/(2z^2))\sqrt{\pi}$, so that
$$\acosh^2(z)-\dfrac{p(n)}{q(n)}\sim(-1)^n\dfrac{\sqrt{\pi}/(2z^2)}{(1/(z^2-1))^{n+1}n^{3/2}}\;.$$
$$A=1-((12-z^2)/(8z^2))/n+((z^4-280z^2+480)/(128z^4))/n^2+\cdots$$
Series:
$$\acosh^2(z)=(z^2-1)\sum_{n\ge0}(-1)^n\dfrac{n!}{(n+1)(3/2)_n}(z^2-1)^n$$
\end{cf}

\smallskip

\begin{cf}\label{6.7.21}{\ }
\begin{verbatim}
[(z)->atanh(z)^2,[0,(2*z^4-12*z^2+2)*n^2-(z^4-18*z^2+1)*n-8*z^2],
               [z^2,8*z^2*(1-z^2)^2*n^3*(2*n-1)]]
\end{verbatim}
$$\atanh^2(z)=\dfrac{z^2}{z^4-2z^2+1+\dfrac{8z^6-16z^4+8z^2}{6z^4-20z^2+6+\dfrac{192z^6-384z^4+192z^2}{15z^4-62z^2+15+\ddots}}}$$
Convergence type $E$ with $E=-(1-z^2)^2/(4z^2)$, $P=3/2$, and
$C=(z^2/(2(z^2+1)^2))\sqrt{\pi}$, so that
$$\atanh^2(z)-\dfrac{p(n)}{q(n)}\sim(-1)^n\dfrac{(z^2/(2(z^2+1)^2))\sqrt{\pi}}{((1-z^2)/(2z))^{2n}n^{3/2}}\;.$$
\begin{align*}A&=1-((11z^4-26z^2+11)/(8(1+z^2)^2))/n\\
  &\phantom{=}+((201z^8-1916z^6+3446z^4-1916z^2+201)/(128(1+z^2)^4))/n^2+\cdots\end{align*}
  Series:
  $$\atanh^2(z)=\dfrac{z^2}{(z^2-1)^2}\sum_{n\ge0}(-1)^n\dfrac{n!}{(n+1)(3/2)_n}\left(\dfrac{2z}{1-z^2}\right)^{2n}$$
\end{cf}

\medskip

\section{Gamma Function}\label{sec:gamma}

\medskip

\subsection{Introduction}

\medskip

General remark: there exist CFs for quotients of products of $1$, $2$, $3$,
$4$, and $8$ gamma functions, those for $4$ and $8$ due to Ramanujan, and
those for $3$ can be found in \cite{CTZ}. We give most of them, and refer
to \cite{CTZ} for a few others.

Note also that there exist literally thousands of evaluations of hypergeometric
series in terms of linear combinations of quotients of gamma functions, which
can therefore be trivially converted into CFs for the latter using Euler's
transformation of series into CFs. We have decided not to include them,
except the simplest ones, and refer to \cite{Beu-Coh1} for a large number
of examples. See \ref{6.23.1} for a sample.

\smallskip

It is convenient to classify CFs related to $\G(z)$ by the number of
free variables. Specializing CFs for $k$ variables will of course give
CFs for $k'<k$ variables. We have included CFs for \emph{constants}
(i.e., $0$ variable) related to the gamma function in Section \ref{sec:gam0}.

\medskip

\subsection{One Free Variable}

\medskip

For the next CFs, note that
$$\sqrt{\pi}\dfrac{\G(z+1)}{\G(z+1/2)}=2^{2z}\dfrac{\G(z+1)^2}{\G(2z+1)}$$

\smallskip

\begin{cf}\label{6.8.1}{\ }
\begin{verbatim}
[(z)->sqrt(Pi)*gamma(z+1)/gamma(z+1/2),
[2,z+2,2*z],[2*(z-1),(n+z)*(n+1-z)]]
\end{verbatim}
$$\sqrt{\pi}\dfrac{\G(z+1)}{\G(z+1/2)}=2+\dfrac{2z-2}{z+1-\dfrac{z^2-z-2}{2z-\dfrac{z^2-z-6}{2z-\dfrac{z^2-z-12}{2z-\dfrac{z^2-z-20}{2z-\dfrac{z^2-z-30}{2z-\ddots}}}}}}$$
Convergence type $P^-$ with $P=2z$ and $C=-\G(1+z)/\G(1-z)$, so that
$$\sqrt{\pi}\dfrac{\G(z+1)}{\G(z+1/2)}-\dfrac{p(n)}{q(n)}\sim(-1)^{n+1}\dfrac{\G(1+z)/\G(1-z)}{n^{2z}}$$
$$A=1-2z/n+(z^3/3+3z^2/2+2z/3)/n^2+\cdots$$
Series:
$$\sqrt{\pi}\dfrac{\G(z+1)}{\G(z+1/2)}=2+\dfrac{2(z-1)}{z+1}\sum_{n\ge0}(-1)^n\dfrac{(2-z)_n}{(2+z)_n}$$
Parametric family for $k\ge0$:
\begin{verbatim}
[(z)->sqrt(Pi)*gamma(z+1)/gamma(z+1/2),2*z+2*k,(n+z)*(n+1-z)]
\end{verbatim}
Convergence type $P^-$ with $P=2z+2k$.
\end{cf}

\smallskip

\begin{cf}\label{6.8.2}{\ }
\begin{verbatim}
[(z)->sqrt(Pi)*gamma(z+1)/gamma(z+1/2),[2],[2*(z-1),n*(n+2*z-1)]]
\end{verbatim}
$$\sqrt{\pi}\dfrac{\G(z+1)}{\G(z+1/2)}=2+\dfrac{2z-2}{2+\dfrac{2z}{2+\dfrac{4z+2}{2+\dfrac{6z+6}{2+\dfrac{8z+12}{2+\dfrac{10z+20}{2+\ddots}}}}}}$$
Convergence type $P^-$ with $P=2$ and
$C=(z(z-1)/2)\sqrt{\pi}\G(z+1)/\G(z+1/2)$, so that
$$\sqrt{\pi}\dfrac{\G(z+1)}{\G(z+1/2)}-\dfrac{p(n)}{q(n)}\sim(-1)^n\dfrac{(z(z-1)/2)\sqrt{\pi}\G(z+1)/\G(z+1/2)}{n^2}$$
$$A=1-(2z)/n+(7z^2/2-z/2-1/2)/n^2+(-6z^3+2z^2+2z)/n^3+\cdots$$
Series:
\begin{align*}\sqrt{\pi}\dfrac{\G(z+1)}{\G(z+1/2)}&=2+2z(z-1)\sum_{n\ge0}\dfrac{n!(z+1/2)_n}{(2n+z)(2n+z+2)(3/2)_n(z+1)_n}\\
  \dfrac{\G(z+1/2)}{\sqrt{\pi}\G(z+1)}&=(z-1)\sum_{n\ge0}\dfrac{(1/2)_n(z)_n}{(2n+z-1)(2n+z+1)n!(z+1/2)_n}\end{align*}
Parametric family for $k\ge0$:
\begin{verbatim}
[(z)->sqrt(Pi)*gamma(z+1)/gamma(z+1/2),2*(k+1),n*(n+2*z-1)]
\end{verbatim}
Convergence type $P^-$ with $P=2(k+1)$.
\end{cf}

\smallskip

\begin{cf}\label{6.8.3}{\ }
\begin{verbatim}
[(z)->sqrt(Pi)*gamma(z+1)/gamma(z+1/2),
[1,2*(z+1),4*n^2+(4*z-3)*n-(z-1)],[z,-2*n*(2*n+1)*(n+z)^2]]
\end{verbatim}
$$\sqrt{\pi}\dfrac{\G(z+1)}{\G(z+1/2)}=1+\dfrac{z}{2z+2-\dfrac{6z^2+12z+6}{7z+11-\dfrac{20z^2+80z+80}{11z+28-\dfrac{42z^2+252z+378}{15z+53-\dfrac{72z^2+576z+1152}{19z+86-\ddots}}}}}$$
Convergence type $P^+$ with $P=1/2$ and $C=2z/\sqrt{\pi}$, so that
$$\sqrt{\pi}\dfrac{\G(z+1)}{\G(z+1/2)}-\dfrac{p(n)}{q(n)}\sim\dfrac{2z/\sqrt{\pi}}{n^{1/2}}\;.$$
$$A=1-(z/3+7/24)/n+(z^2/5+11z/40+61/640)/n^2+\cdots$$
Series:
$$\sqrt{\pi}\dfrac{\G(z+1)}{\G(z+1/2)}=1+\dfrac{z}{2}\sum_{n\ge0}\dfrac{(3/2)_n}{(n+z+1)(n+1)!}$$
Parametric family for $k\ge0$:
\begin{verbatim}
[(z)->sqrt(Pi)*gamma(z+1)/gamma(z+1/2),
4*n^2+(4*z-3)*n-(z-1)+k*(2*k+1),-2*n*(2*n+1)*(n+z)^2]
\end{verbatim}
Convergence type $P^+$ with $P=2k+1/2$.
\end{cf}
    
\smallskip

\begin{cf}\label{6.8.4}{\ }
\begin{verbatim}
[(z)->sqrt(Pi)*gamma(z+1)/gamma(z+1/2),[0,3*n-2],[2*z,-2*n*(n-z)]]
\end{verbatim}
$$\sqrt{\pi}\dfrac{\G(z+1)}{\G(z+1/2)}=\dfrac{2z}{1+\dfrac{2z-2}{4+\dfrac{4z-8}{7+\dfrac{6z-18}{10+\dfrac{8z-32}{13+\dfrac{10z-50}{16+\ddots}}}}}}$$
Convergence type $E$ with $E=2$, $P=3z-1$, and
$C=z\G^3(z)2^{4z-1}\sin(\pi z)/\pi$, so that
$$\sqrt{\pi}\dfrac{\G(z+1)}{\G(z+1/2)}-\dfrac{p(n)}{q(n)}\sim\dfrac{z\G^3(z)\sin(\pi z)/\pi}{2^{n+1-4z}n^{3z-1}}\;.$$
$$A=1-(9z^2/2-5z/2+1)/n+(81z^4/8-3z^3/4-17z^2/8+31z/4-2)/n^2+\cdots$$
Parametric family for $k\ge0$:
\begin{verbatim}
[(z)->sqrt(Pi)*gamma(z+1)/gamma(z+1/2),3*n-2+k,-2*n*(n-z)]
\end{verbatim}
Convergence type $E$ with $E=2$ and $P=3z+2k-1$.
\end{cf}

\smallskip

\begin{cf}\label{6.8.5}{\ }
\begin{verbatim}
[(z)->sqrt(Pi)*gamma(z+1)/gamma(z+1/2),
[0,2*z,3*n+4*z-4],[2*z,-2*(n+z-1)*(n+2*z-1)]]
\end{verbatim}
$$\sqrt{\pi}\dfrac{\G(z+1)}{\G(z+1/2)}=\dfrac{2z}{2z-\dfrac{4z^2}{4z+2-\dfrac{4z^2+6z+2}{4z+5-\dfrac{4z^2+12z+8}{4z+8-\dfrac{4z^2+18z+18}{4z+11-\ddots}}}}}$$
Convergence type $E$ with $E=2$, $P=1-z$, and
$C=z2^{2-2z}\sqrt{\pi}/\G(z+1/2)$, so that
$$\sqrt{\pi}\dfrac{\G(z+1)}{\G(z+1/2)}-\dfrac{p(n)}{q(n)}\sim\dfrac{z2^{2-2z}\sqrt{\pi}/\G(z+1/2)}{2^nn^{1-z}}\;.$$
$$A=1+(3z^2/2-z/2-1)/n+(9z^4/8-23z^3/12-5z^2/8-19z/12+3)/n^2+\cdots$$
Series:
$$\sqrt{\pi}\dfrac{\G(z+1)}{\G(z+1/2)}=\sum_{n\ge0}\dfrac{(2z)_n}{(z+1)_n}2^{-n}$$
Parametric family for $k\ge0$:
\begin{verbatim}
[(z)->sqrt(Pi)*gamma(z+1)/gamma(z+1/2),3*n+4*z-4+k,-2*(n+z-1)*(n+2*z-1)]
\end{verbatim}
Convergence type $E$ with $E=2$ and $P=2k+1-z$.
\end{cf}
    
\smallskip

\begin{cf}\label{6.8.6}{\ }
\begin{verbatim}
[(z)->sqrt(Pi)*gamma(z+1)/gamma(z+1/2),
[0,1,2*n+2*z-2],[[1,-2*z],[-2*n*(n-z),-2*(n+z)*(n+2*z)]]]
\end{verbatim}
$$\sqrt{\pi}\dfrac{\G(z+1)}{\G(z+1/2)}=\dfrac{1}{1-\dfrac{2z}{2z+2+\dfrac{2z-2}{2z+4-\dfrac{4z^2+6z+2}{2z+6+\dfrac{4z-8}{2z+8-\dfrac{4z^2+12z+8}{2z+10+\ddots}}}}}}$$
Convergence type $E$ with $E=-(1+\sqrt{2})^2$, $P=0$, and $C=4\sqrt{\pi}(\G(z+1)/\G(z+1/2))\sin(\pi z)/(1+\sqrt{2})^{2z}$, so that
$$\sqrt{\pi}\dfrac{\G(z+1)}{\G(z+1/2)}-\dfrac{p(n)}{q(n)}\sim(-1)^n\dfrac{4\sqrt{\pi}(\G(z+1)/\G(z+1/2))\sin(\pi z)}{(1+\sqrt{2})^{2n+2z}}\;.$$
$$A=1+((3z^2-2z+1/8)d)/n+(9z^4-(3d+12)z^3+(2d+19/4)z^2-(d/8+1/2)z+1/64)/n^2+\cdots$$
\end{cf}

\smallskip

\begin{cf}\label{6.8.7}{\ }
\begin{verbatim}
[(z)->sqrt(Pi)*gamma(z+1)/gamma(z+1/2),
[2,2,n+z],[[2*z-2,2*z],[(n+z)*(n+1-z),(n+1)*(n+2*z)]]]
\end{verbatim}
$$\sqrt{\pi}\dfrac{\G(z+1)}{\G(z+1/2)}=2+\dfrac{2z-2}{2+\dfrac{2z}{z+2-\dfrac{z^2-z-2}{z+3+\dfrac{4z+2}{z+4-\dfrac{z^2-z-6}{z+5+\dfrac{6z+6}{z+6-\ddots}}}}}}$$
Convergence type $E$ with $E=-(1+\sqrt{2})^2$, $P=0$, and $C=...$, so that
$$\sqrt{\pi}\dfrac{\G(z+1)}{\G(z+1/2)}-\dfrac{p(n)}{q(n)}\sim(-1)^n\dfrac{C}{(1+\sqrt{2})^{2n}}$$
$$A=1+(3z^2-4z+9/8)d/n+\cdots$$
\end{cf}

\smallskip

\begin{cf}\label{6.8.8}{\ }
\begin{verbatim}
[(z)->sqrt(Pi)*gamma(z)/gamma(z+1/2)+psi(z)+Euler(),
[2,2*(z+1),(2*z+1)*n+z],[2*(z-1),(n+1)^2*(n+z)*(n+1-z)]]
\end{verbatim}
$$\sqrt{\pi}\dfrac{\G(z)}{\G(z+1/2)}+\psi(z)+\gamma=2+\dfrac{2z-2}{2z+2-\dfrac{4z^2-4z-8}{5z+2-\dfrac{9z^2-9z-54}{7z+3-\dfrac{16z^2-16z-192}{9z+4-\dfrac{25z^2-25z-500}{11z+5-\ddots}}}}}$$
Convergence type $P^-$ with $P=2z+1$ and $C=-z\G(z)^2(\sin(\pi z)/\pi)$, so that
$$\sqrt{\pi}\dfrac{\G(z)}{\G(z+1/2)}+\psi(z)+\gamma-\dfrac{p(n)}{q(n)}\sim(-1)^{n+1}\dfrac{z\G(z)^2(\sin(\pi z)/\pi)}{n^{2z+1}}$$
$$A=1-(2z+3/2)/n+(z^3/3+3z^2/2+19z/6+2)/n^2+\cdots$$
Series:
$$\sqrt{\pi}\dfrac{\G(z)}{\G(z+1/2)}+\psi(z)+\gamma=2+2\dfrac{z-1}{z+1}\sum_{n\ge0}(-1)^n\dfrac{(2-z)_n}{(n+2)(2+z)_n}$$
\end{cf}

\smallskip

For the next CFs, note that
$$\dfrac{\G(z+1)^2}{\G(2z)}=z2^{1-2z}\sqrt{\pi}\dfrac{\G(z+1)}{\G(z+1/2)}\;.$$

\smallskip

\begin{cf}\label{6.8.9}{\ }
\begin{verbatim}
[(z)->gamma(z+1)^2/gamma(2*z),[z,z+1,2*n^2+(z-1)*(n-z)],
                              [-z^2*(z-1),-n*(n+z)^2*(n+1-z)]]
\end{verbatim}
$$\dfrac{\G(z+1)^2}{\G(2z)}=z-\dfrac{z^3-z^2}{z+1+\dfrac{z^3-3z-2}{-z^2+3z+6+\dfrac{2z^3+2z^2-16z-24}{-z^2+4z+15+\dfrac{3z^3+6z^2-45z-108}{-z^2+5z+28+\ddots}}}}$$
Convergence type $P^+$ with $P=z$ and $C=-1/\G(-z)$, so that
$$\dfrac{\G(z+1)^2}{\G(2z)}-\dfrac{p(n)}{q(n)}\sim-\dfrac{1/\G(-z)}{n^z}$$
$$A=1-((z^3-4z^2-2z)/(2(z+1)))/n+\cdots$$
Series:
$$\dfrac{\G(z+1)^2}{\G(2z)}=z-z^2(z-1)\sum_{n\ge0}\dfrac{(2-z)_n}{(n+z+1)(n+1)!}$$
Parametric family for $k\ge0$:
\begin{verbatim}
[(z)->gamma(z+1)^2/gamma(2*z),2*n^2+(z-1)*(n-z)+k*(z+k),-n*(n+z)^2*(n+1-z)]
\end{verbatim}
Convergence type $P^+$ with $P=2k+z$.
\end{cf}

\smallskip

\begin{cf}\label{6.8.10}{\ }
\begin{verbatim}
[(z)->gamma(z+1)^2/gamma(2*z),[1,2*n^2-(2*z-1)*n+z^2],
                              [-(z-1)^2,-n*(n+1)*(n+1-z)^2]]
\end{verbatim}
$$\dfrac{\G(z+1)^2}{\G(2z)}=1-\dfrac{z^2-2z+1}{z^2-2z+3-\dfrac{2z^2-8z+8}{z^2-4z+10-\dfrac{6z^2-36z+54}{z^2-6z+21-\dfrac{12z^2-96z+192}{z^2-8z+36-\ddots}}}}$$
Convergence type $P^+$ with $P=2z$ and $C=-2z^3\sin(\pi z)^2\G(z)^4/(\pi 2^{4z}\G(z+1/2)^2)$, so that
$$\dfrac{\G(z+1)^2}{\G(2z)}-\dfrac{p(n)}{q(n)}\sim-\dfrac{2z^3\sin(\pi z)^2\G(z)^4/(\pi 2^{4z}\G(z+1/2)^2)}{n^{2z}}\;.$$
$$A=1+((2z^3-4z^2-3z)/(2z+1))/n+\cdots$$
Series:
$$\dfrac{\G(2z)}{\G(z+1)^2}=\sum_{n\ge0}\dfrac{(1-z)_n^2}{(n+1)n!^2}$$
Parametric family for $k\ge0$:
\begin{verbatim}
[(z)->gamma(z+1)^2/gamma(2*z),2*n^2-(2*z-1)*n+(z+k)^2,-n*(n+1)*(n+1-z)^2]
\end{verbatim}
Convergence type $P^+$ with $P=2z+2k$.
\end{cf}

\smallskip

\begin{cf}\label{6.8.11}{\ }
\begin{verbatim}
[(z)->gamma(z+1)^2/gamma(2*z),[1,(z+1)^2,2*n^2+(2*z-1)*2*n+z^2],
                             [-(z-1)^2,-(n+1)*(n+2*z-1)*(n+z)^2]]
\end{verbatim}
$$\dfrac{\G(z+1)^2}{\G(2z)}=1-\dfrac{z^2-2z+1}{z^2+2z+1-\dfrac{4z^3+8z^2+4z}{z^2+8z+4-\dfrac{6z^3+27z^2+36z+12}{z^2+12z+12-\ddots}}}$$
Convergence type $P^+$ with $P=1$ and $C=-(z-1)^2\G(z+1)^2/\G(2z)$, so that
$$\dfrac{\G(z+1)^2}{\G(2z)}-\dfrac{p(n)}{q(n)}\sim-\dfrac{(z-1)^2\G(z+1)^2/\G(2z)}{n}\;.$$
$$A=1+(z^2/2-2z)/n+\cdots$$
Series:
\begin{align*}\dfrac{\G(z+1)^2}{\G(2z)}&=1-\dfrac{(z-1)^2}{(z+1)^2}\sum_{n\ge0}\dfrac{(n+1)!(2z)_n}{(z+2)_n^2}\\
  \dfrac{\G(2z)}{\G(z+1)^2}&=\dfrac{2z-1}{z^2}+\dfrac{(z-1)^2}{z^2}\sum_{n\ge0}\dfrac{(z)_n^2}{(n+1)!(2z)_n}\end{align*}
Parametric family for $k\ge0$:
\begin{verbatim}
[(z)->gamma(z+1)^2/gamma(2*z),2*n^2+(2*z-1)*2*n+z^2+k*(k+1),
                              -(n+1)*(n+2*z-1)*(n+z)^2]
\end{verbatim}
Convergence type $P^+$ with $P=2k+1$.
\end{cf}

\smallskip

\begin{cf}\label{6.8.12}{\ }
\begin{verbatim}
[(z)->gamma(z+1)^2/gamma(2*z),[1,z^2+2,(2*z+1)*(2*n-1)],
                             [-2*(z-1)^2,(n+z-1)^2*(n+1-z)^2]]
\end{verbatim}
$$\dfrac{\G(z+1)^2}{\G(2z)}=1-\dfrac{2z^2-4z+2}{z^2+2+\dfrac{z^4-4z^3+4z^2}{6z+3+\dfrac{z^4-4z^3-2z^2+12z+9}{10z+5+\dfrac{z^4-4z^3-12z^2+32z+64}{14z+7+\ddots}}}}$$
Convergence type $P^-$ with $P=4z+2$ and $C=...$, so that
$$\dfrac{\G(z+1)^2}{\G(2z)}-\dfrac{p(n)}{q(n)}\sim(-1)^n\dfrac{C}{n^{4z+2}}$$
$$A=1-(2z+1)/n+((2z+1)(z^2-5z+3)/3)/n^2+\cdots$$
Series:
$$\dfrac{\G(z+1)^2}{\G(2z)}=1-2z^2(z-1)^2\sum_{n\ge0}(-1)^n\dfrac{(2-z)_n^2}{(n^2+n+z^2)(n^2+3n+2+z^2)(1+z)_n^2}$$
Parametric family for $k\ge0$:
\begin{verbatim}
[(z)->gamma(z+1)^2/gamma(2*z),(2*z+2*k+1)*(2*n-1),(n+z-1)^2*(n+1-z)^2]
\end{verbatim}
Convergence type $P^-$ with $P=4z+4k+2$.
\end{cf}

\smallskip

\begin{cf}\label{6.8.13}{\ }
\begin{verbatim}
[(z)->(Pi/sin(Pi*z))*gamma(z+1)^2/gamma(2*z),
[2,(z+1)^2,(2*n+2*z-1)*(n^2+(2*z-1)*n+(z-1)^2)],
           [4*z^3,-n*(n+z)^4*(n+2*z)]]
\end{verbatim}
$$\dfrac{\pi}{\sin(\pi z)}\dfrac{\G(z+1)^2}{\G(2z)}=2+\dfrac{4z^3}{z^2+2z+1-\dfrac{2z^5+9z^4+16z^3+14z^2+6z+1}{2z^3+7z^2+12z+9-\ddots}}$$
Convergence type $P^+$ with $P=2z-2$ and $C=...$, so that
$$\dfrac{\pi}{\sin(\pi z)}\dfrac{\G(z+1)^2}{\G(2z)}-\dfrac{p(n)}{q(n)}\sim\dfrac{C}{n^{2z-2}}$$
$$A=1+(-2z^2+z+1)/n+((z-1)(12z^4-16z^3-17z^2+2z+6)/(6(z-2)))/n^2+\cdots$$
Series:
$$\dfrac{\pi}{\sin(\pi z)}\dfrac{\G(z+1)^2}{\G(2z)}=2+4z^3\sum_{n\ge0}\dfrac{(2z+1)_n}{(n+z+1)^2(n+1)!}$$
\end{cf}

\smallskip
    
\begin{cf}\label{6.8.14}{\ }
\begin{verbatim}
[(z)->z*(gamma(z)/gamma(z+1/2))^2,[0,4*z-1,8*z-2],[4*z,(2*n-1)^2]]
\end{verbatim}
$$z\left(\dfrac{\G(z)}{\G(z+1/2)}\right)^2=\dfrac{4z}{4z-1+\dfrac{1}{8z-2+\dfrac{9}{8z-2+\dfrac{25}{8z-2+\dfrac{49}{8z-2+\dfrac{81}{8z-2+\ddots}}}}}}$$
Convergence type $P^-$ with $P=4z-1$ and $C=z\G(z)^4/\pi^2$, so that
$$z\left(\dfrac{\G(z)}{\G(z+1/2)}\right)^2-\dfrac{p(n)}{q(n)}\sim(-1)^n\dfrac{z\G(z)^4/\pi^2}{n^{4z-1}}\;.$$
$$A=1-((4z-1)(4z^2-2z+3)/12)/n^2+\cdots$$
Parametric family for $k\ge0$:
\begin{verbatim}
[(z)->z*(gamma(z)/gamma(z+1/2))^2,8*z+4*k-2,(2*n-1)^2]
\end{verbatim}
Convergence type $P^-$ with $P=4z+2k-1$.
\end{cf}

\smallskip

\begin{cf}\label{6.8.15}{\ }
\begin{verbatim}
[(z)->z*(gamma(z)/gamma(z+1/2))^2,
[4*z/(4*z-1),8*z^2/(4*z-1),2],[-2*z/(4*z-1)^2,(n+2*z-1)^2]]
\end{verbatim}
$$z\left(\dfrac{\G(z)}{\G(z+1/2)}\right)^2=4z/(4z-1)-\dfrac{2z/(16z^2-8z+1)}{8z^2/(4z-1)+\dfrac{4z^2}{2+\dfrac{4z^2+4z+1}{2+\dfrac{4z^2+8z+4}{2+\dfrac{4z^2+12z+9}{2+\ddots}}}}}$$
Convergence type $P^-$ with $P=2$ and $C=-z(\G(z)/\G(z+1/2))^2/8$, so that
$$z\left(\dfrac{\G(z)}{\G(z+1/2)}\right)^2-\dfrac{p(n)}{q(n)}\sim(-1)^{n+1}\dfrac{z(\G(z)/\G(z+1/2))^2/8}{n^2}\;.$$
$$A=1+(1-4z)/n+(12z^2-6z+1/8)/n^2+(-32z^3+24z^2-z-3/4)/n^3+\cdots$$
Series:
\begin{align*}z\left(\dfrac{\G(z)}{\G(z+1/2)}\right)^2&=\dfrac{4z}{4z-1}-\dfrac{1}{16z}\sum_{n\ge0}\dfrac{(z+1/2)_n^2}{(n+z-1/4)(n+z+3/4)(z+1)_n^2}\\
\left(\dfrac{\G(z+1/2)}{\G(z)}\right)^2&=\dfrac{(2z-1)^2}{4z-3}-\sum_{n\ge0}\dfrac{(z)_n^2}{(4n+4z-3)(4n+4z+1)(z+1/2)_n^2}\end{align*}
Parametric family for $k\ge0$:
\begin{verbatim}
[(z)->z*(gamma(z)/gamma(z+1/2))^2,2*k+2,(n+2*z-1)^2]
\end{verbatim}
Convergence type $P^-$ with $P=2k+2$.
\end{cf}

\smallskip

\begin{cf}\label{6.8.16}{\ }
\begin{verbatim}
[(z)->z*(gamma(z)/gamma(z+1/2))^2,[1,8*z-1,8*(2*n-1)*z],
                                [2,(4*n^2-1)^2]]
\end{verbatim}
$$z\left(\dfrac{\G(z)}{\G(z+1/2)}\right)^2=1+\dfrac{2}{8z-1+\dfrac{9}{24z+\dfrac{225}{40z+\dfrac{1225}{56z+\dfrac{3969}{72z+\dfrac{9801}{88z+\ddots}}}}}}$$
Convergence type $P^-$ with $P=4z$ and $C=z^2\G(z)^4/\pi^2$, so that
$$z\left(\dfrac{\G(z)}{\G(z+1/2)}\right)^2-\dfrac{p(n)}{q(n)}\sim(-1)^n\dfrac{z^2\G(z)^4/\pi^2}{n^{4z}}\;.$$
$$A=1-2z/n-(z(8z^2-12z+1)/6)/n^2+(z(2z+1)(8z^2-4z+3))/n^3-\cdots$$
Parametric family for $k\ge0$:
\begin{verbatim}
[(z)->z*(gamma(z)/gamma(z+1/2))^2,8*(2*n-1)*(z+k),(4*n^2-1)^2]
\end{verbatim}
Convergence type $P^-$ with $P=4k+4z$.
\end{cf}

\smallskip

\begin{cf}\label{6.8.17}{\ }
\begin{verbatim}
[(z)->z*(gamma(z)/gamma(z+1/2))^2,[1,8*z-1,8*z],[2,4*n^2-1]]
\end{verbatim}
$$z\left(\dfrac{\G(z)}{\G(z+1/2)}\right)^2=1+\dfrac{2}{8z-1+\dfrac{3}{8z+\dfrac{15}{8z+\dfrac{35}{8z+\dfrac{63}{8z+\dfrac{99}{8z+\ddots}}}}}}$$
Convergence type $P^-$ with $P=4z$ and $C=z^2\G(z)^4/\pi^2$, so that
$$z\left(\dfrac{\G(z)}{\G(z+1/2)}\right)^2-\dfrac{p(n)}{q(n)}\sim(-1)^n\dfrac{z^2\G(z)^4/\pi^2}{n^{4z}}\;.$$
$$A=1-2z/n-(z(8z^2-12z+1)/6)/n^2+(z(2z+1)(8z^2-4z+3))/n^3-\cdots$$
Parametric family for $k\ge0$:
\begin{verbatim}
[(z)->z*(gamma(z)/gamma(z+1/2))^2,8*z+4*k,4*n^2-1]
\end{verbatim}
Convergence type $P^-$ with $P=4z+2k$.
\end{cf}

\smallskip

\begin{cf}\label{6.8.18}{\ }
\begin{verbatim}
[(z)->z*(gamma(z)/gamma(z+1/2))^2,[1,2*z,1],[1,(n+2*z-1)*(n+2*z)]]
\end{verbatim}
$$z\left(\dfrac{\G(z)}{\G(z+1/2)}\right)^2=1+\dfrac{1}{2z+\dfrac{4z^2+2z}{1+\dfrac{4z^2+6z+2}{1+\dfrac{4z^2+10z+6}{1+\dfrac{4z^2+14z+12}{1+\dfrac{4z^2+18z+20}{1+\ddots}}}}}}$$
Convergence type $P^-$ with $P=1$ and $C=z(\G(z)/\G(z+1/2))^2/2$, so that
$$z\left(\dfrac{\G(z)}{\G(z+1/2)}\right)^2-\dfrac{p(n)}{q(n)}\sim(-1)^n\dfrac{z(\G(z)/\G(z+1/2))^2/2}{n}\;.$$
$$A=1-2z/n+(4z^2-3/16)/n^2+(-8z^3+9z/8)/n^3+\cdots$$
Series:
\begin{align*}z\left(\dfrac{\G(z)}{\G(z+1/2)}\right)^2&=1+\dfrac{1}{4z}\sum_{n\ge0}\dfrac{(z+1/2)_n^2}{(n+z+1)(z+1)_n^2}\\
\left(\dfrac{\G(z+1/2)}{\G(z)}\right)^2&=\dfrac{2z-1}{2}+\dfrac{1}{2}\sum_{n\ge0}\dfrac{(z)_n^2}{(2n+2z+1)(z+1/2)_n^2}\end{align*}
Parametric family for $k\ge0$:
\begin{verbatim}
[(z)->z*(gamma(z)/gamma(z+1/2))^2,2*k+1,(n+2*z-1)*(n+2*z)]
\end{verbatim}
Convergence type $P^-$ with $P=2k+1$.
\end{cf}
    
\smallskip

\begin{cf}\label{6.8.19}{\ }
\begin{verbatim}
[(z)->z*(gamma(z)/gamma(z+1/2))^2,
[0,2*n+4*z-3],[[4*z,1],[4*(n+2*z-1)^2,(2*n+1)^2]]]
\end{verbatim}
$$z\left(\dfrac{\G(z)}{\G(z+1/2)}\right)^2=\dfrac{4z}{4z-1+\dfrac{1}{4z+1+\dfrac{16z^2}{4z+3+\dfrac{9}{4z+5+\dfrac{16z^2+16z+4}{4z+7+\dfrac{25}{4z+9+\ddots}}}}}}$$
Convergence type $E$ with $E=-(1+\sqrt{2})^2$, $P=0$, and $C=4z(\G(z)/\G(z+1/2))^2/(1+\sqrt{2})^{4z-1}$, so that
$$z\left(\dfrac{\G(z)}{\G(z+1/2)}\right)^2-\dfrac{p(n)}{q(n)}\sim(-1)^n\dfrac{4z(\G(z)/\G(z+1/2))^2}{(1+\sqrt{2})^{2n+4z-1}}\;.$$
\begin{align*}A&=1+((4z^2-2z-1/8)d)/n\\
  &\phantom{=}+(16z^4-(8d+16)z^3+(6d+3)z^2+(-3d/4+1/2)z-d/16+1/64)/n^2+\cdots\end{align*}
\end{cf}

\smallskip

\begin{cf}\label{6.8.20}{\ }
\begin{verbatim}
[(z)->z*(gamma(z)/gamma(z+1/2))^2,
[1,4*z,4*z+2*n-1],[[2,16*z^2+8*z],[4*n^2-1,4*(n+2*z)*(n+2*z+1)]]]
\end{verbatim}
$$z\left(\dfrac{\G(z)}{\G(z+1/2)}\right)^2=1+\dfrac{2}{4z+\dfrac{16z^2+8z}{4z+3+\dfrac{3}{4z+5+\dfrac{16z^2+24z+8}{4z+7+\dfrac{15}{4z+9+\dfrac{16z^2+40z+24}{4z+11+\ddots}}}}}}$$
Convergence type $E$ with $E=-(1+\sqrt{2})^2$, $P=0$, and $C=4z(\G(z)/\G(z+1/2))^2/(1+\sqrt{2})^{4z+1}$, so that
$$z\left(\dfrac{\G(z)}{\G(z+1/2)}\right)^2-\dfrac{p(n)}{q(n)}\sim(-1)^n\dfrac{4z(\G(z)/\G(z+1/2))^2}{(1+\sqrt{2})^{2n+4z+1}}\;.$$
$$A=1+((4z^2-2z+3/8)d)/n+(16z^4-(8d+16)z^3+(2d+7)z^2+(d/4-3/2)z-3d/16+9/64)/n^2+\cdots$$
\end{cf}

\smallskip

\begin{cf}\label{6.8.21}{\ }
\begin{verbatim}
[(z)->z*(gamma(z)/gamma(z+1/2))^2,
[[1,4*z*(z+1)],[20*n^2+(24*z-12)*n+8*z^2-8*z+1,
                20*n^2+(24*z+12)*n+8*z^2+8*z+1]],
[[1,-8*z*(z+1)*(2*z+1)^2],[(4*n^2-1)^2,
                           -16*(n+z)*(n+z+1)*(2*n+2*z+1)^2]]]
\end{verbatim}
$$z\left(\dfrac{\G(z)}{\G(z+1/2)}\right)^2=1+\dfrac{1}{4z^2+4z-\dfrac{32z^4+64z^3+40z^2+8z}{8z^2+16z+9+\dfrac{9}{8z^2+32z+33-\ddots}}}$$
Convergence type $E$ with $E=-((1+\sqrt{5})/2)^5$, $P=0$, and
$C=4z(\G(z)/\G(z+1/2))^2/((1+\sqrt{5})/2)^{4z+3}$, so that
$$z\left(\dfrac{\G(z)}{\G(z+1/2)}\right)^2-\dfrac{p(n)}{q(n)}\sim(-1)^n\dfrac{4z(\G(z)/\G(z+1/2))^2}{((1+\sqrt{5})/2)^{5n+4z+3}}\;.$$
$$A=1+((4z^2/5-4z/5+3/10)d)/n+\cdots$$
\end{cf}

\smallskip

{\bf Comment:} Note that \ref{6.8.17} is a simplification of \ref{6.8.16}
(divide by $2n-1$), both can be Ap\'ery accelerated, but the Ap\'ery
acceleration of the more complicated \ref{6.8.16} has a \emph{faster} speed of
convergence than that of \ref{6.8.17}.

\smallskip

\begin{cf}\label{6.8.22}{\ }
\begin{verbatim}
[(z)->Pi*z*(gamma(z)/gamma(z+1/2))^2,
[4,3*(z+1),4*n^2+2*z],[4*(1-z),-(2*n+1)^2*(n+z)*(n+1-z)]]
\end{verbatim}
$$\pi z\left(\dfrac{\G(z)}{\G(z+1/2)}\right)^2=4-\dfrac{4z-4}{3z+3+\dfrac{9z^2-9z-18}{2z+16+\dfrac{25z^2-25z-150}{2z+36+\dfrac{49z^2-49z-588}{2z+64+\ddots}}}}$$
Convergence type $P^+$ with $P=2z$ and $C=\G(z)^2\sin(\pi z)/\pi$, so that
$$\pi z\left(\dfrac{\G(z)}{\G(z+1/2)}\right)^2-\dfrac{p(n)}{q(n)}\sim\dfrac{\G(z)^2\sin(\pi z)/\pi}{n^{2z}}$$
$$A=1-2z/n+((4z^4+22z^3+32z^2+11z)/(12(z+1)))/n^2+\cdots$$
Series:
$$\pi z\left(\dfrac{\G(z)}{\G(z+1/2)}\right)^2=4+4\dfrac{1-z}{1+z}\sum_{n\ge0}\dfrac{(2-z)_n}{(2n+3)(2+z)_n}$$
Parametric family for $k\ge0$:
\begin{verbatim}
[(z)->Pi*z*(gamma(z)/gamma(z+1/2))^2,
4*n^2+(4*k+2)*z+2*k^2,-(2*n+1)^2*(n+z)*(n+1-z)]
\end{verbatim}
Convergence type $P^+$ with $P=2z+2k$.
\end{cf}

\smallskip

\begin{cf}\label{6.8.23}{\ }
\begin{verbatim}
[(z)->Pi*z*(gamma(z)/gamma(z+1/2))^2,[4,4*(n+z-1)],[4*(1-z),n^2*(n+2*z-1)^2]]
\end{verbatim}
$$\pi z\left(\dfrac{\G(z)}{\G(z+1/2)}\right)^2=4-\dfrac{4z-4}{4z+\dfrac{4z^2}{4z+4+\dfrac{16z^2+16z+4}{4z+8+\dfrac{36z^2+72z+36}{4z+12+\dfrac{64z^2+192z+144}{4z+16+\ddots}}}}}$$
Convergence type $P^-$ with $P=4$ and $C=(\pi/2)z^2(1-z)(\G(z)/\G(z+1/2))^2$, so that
$$\pi z\left(\dfrac{\G(z)}{\G(z+1/2)}\right)^2-\dfrac{p(n)}{q(n)}\sim(-1)^n\dfrac{(\pi/2)z^2(1-z)(\G(z)/\G(z+1/2))^2}{n^4}$$
$$A=1-4z/n+(12z^2-2z-2)/n^2+(-32z^3+12z^2+12z)/n^3+\cdots$$
Series:

There exist hypergeometric series for the CF and for its inverse, need
to use the square root of $z(z-1)$, to be done.

Parametric family for $k\ge0$:
\begin{verbatim}
[(z)->Pi*z*(gamma(z)/gamma(z+1/2))^2,4*(k+1)*(n+z-1),n^2*(n+2*z-1)^2]
\end{verbatim}
Convergence type $P^-$ with $P=4(k+1)$.
\end{cf}

\smallskip

\begin{cf}\label{6.8.24}{\ }
\begin{verbatim}
[(z)->Pi*z*(gamma(z)/gamma(z+1/2))^2,
[[4,4*z],[5*n^2+(4*z+1)*n+z,5*n^2+(4*z+5)*n+3*z+1]],
[-(2*n+1)^2*(n+z)*(n+1-z),(n+1)^2*(n+2*z)^2]]
\end{verbatim}
$$\pi z\left(\dfrac{\G(z)}{\G(z+1/2)}\right)^2=4-\dfrac{4z-4}{4z+\dfrac{4z^2}{5z+6+\dfrac{9z^2-9z-18}{7z+11+\dfrac{16z^2+16z+4}{9z+22+\dfrac{25z^2-25z-150}{11z+31+\ddots}}}}}$$
Convergence type $E$ with $E=-((1+\sqrt{5})/2)^5$, $P=0$, and $C=...$, so that
$$\pi z\left(\dfrac{\G(z)}{\G(z+1/2)}\right)^2-\dfrac{p(n)}{q(n)}\sim(-1)^n\dfrac{C}{((1+\sqrt{5})/2)^{5n}}$$
$$A=1+(14z^2/5-18z/5+9/10)d/n+\cdots$$
\end{cf}

\smallskip

\begin{cf}\label{6.8.25}{\ }
\begin{verbatim}
[(z)->(gamma((z+1)/6)/gamma((z+5)/6))^3,
[[0,z^2-1],[3*(2*n-1),z^2-1]],
[[36,8],2*[(3*n-2)*(3*n-1)^2,(3*n+2)*(3*n+1)^2]]]
\end{verbatim}
$$\left(\dfrac{\G((z+1)/6)}{\G((z+5)/6)}\right)^3=\dfrac{36}{z^2-1+\dfrac{8}{3+\dfrac{8}{z^2-1+\dfrac{160}{9+\dfrac{200}{z^2-1+\dfrac{784}{15+\ddots}}}}}}$$
Convergence type $P^-$ with $P=2z/3$ and $C=...$, so that
$$\left(\dfrac{\G((z+1)/6)}{\G((z+5)/6)}\right)^3-\dfrac{p(n)}{q(n)}\sim(-1)^n\dfrac{C}{n^{2z/3}}\;.$$
$$A=1-(z^3/162+z/6)/n^2+\cdots$$
\end{cf}

\smallskip

\begin{cf}\label{6.8.26}{\ }
\begin{verbatim}
[(z)->(gamma((z+1)/6)/gamma((z+5)/6))^3,
[[0,z^2-9],[3*(2*n-1),z^2-9]],[[36,32],2*[(3*n-2)^3,(3*n+2)^3]]]
\end{verbatim}
$$\left(\dfrac{\G((z+1)/6)}{\G((z+5)/6)}\right)^3=\dfrac{36}{z^2-9+\dfrac{32}{3+\dfrac{2}{z^2-9+\dfrac{250}{9+\dfrac{128}{z^2-9+\dfrac{1024}{15+\ddots}}}}}}$$
Convergence type $P^-$ with $P=2z/3$ and $C=...$, so that
$$\left(\dfrac{\G((z+1)/6)}{\G((z+5)/6)}\right)^3-\dfrac{p(n)}{q(n)}\sim(-1)^n\dfrac{C}{n^{2z/3}}\;.$$
$$A=1-(z^3/162+z/6)/n^2+\cdots$$
\end{cf}
Note that the fact that these two CFs are even functions of $z$ is an illusion
since the equality with the gamma quotient is only valid for $z>0$ when
$z$ is real.

\smallskip

\begin{cf}\label{6.8.27}{\ }
\begin{verbatim}
[(z)->(gamma((z+2)/6)/gamma((z+4)/6))^3,
z*[[0,3],[2*n-1,3]],[[18,4],2*[(3*n-1)^3,(3*n+1)^3]]]
\end{verbatim}
$$\left(\dfrac{\G((z+2)/6)}{\G((z+4)/6)}\right)^3=\dfrac{18}{3z+\dfrac{4}{z+\dfrac{16}{3z+\dfrac{128}{3z+\dfrac{250}{3z+\dfrac{686}{5z+\ddots}}}}}}$$
Convergence type $P^-$ with $P=2z/3$ and $C=...$, so that
$$\left(\dfrac{\G((z+2)/6)}{\G((z+4)/6)}\right)^3-\dfrac{p(n)}{q(n)}\sim(-1)^n\dfrac{C}{n^{2z/3}}\;.$$
$$A=1-((z^5+18z^3-252z)/(162(z^2-9)))/n+\cdots$$
\end{cf}

\smallskip

Replacing $z$ by $2z$ in the above gives a CF for
$(\G((z+1)/3)/\G((z+2)/3))^3$.

\smallskip

For the next CFs, note that
$$\dfrac{2\G((1+z)/4)\G((1-z)/4)\cos(\pi z/2)}{\pi^{3/2}}=\dfrac{\pi^{1/2}}{\G((3+z)/4)\G((3-z)/4)}\;.$$

\smallskip

\begin{cf}\label{6.8.28}{\ }
\begin{verbatim}
[(z)->gamma((3+z)/4)*gamma((3-z)/4)/sqrt(Pi),
[1,12*n^2-4*n+1-z^2],[1-z^2,-8*n^2*((2*n+1)^2-z^2)]]
\end{verbatim}
$$\dfrac{\G((3+z)/4)\G((3-z)/4)}{\sqrt{\pi}}=1+\dfrac{z^2-1}{-z^2+9+\dfrac{8z^2-72}{-z^2+41+\dfrac{32z^2-800}{-z^2+97+\dfrac{72z^2-3528}{-z^2+177+\ddots}}}}$$
Convergence type $E$ with $E=2$, $P=1$, and $C=-\cos(\pi z/2)\G((3+z)/4)^2\G((3-z)/4)^2/\pi^2$, so that
$$\dfrac{\G((3+z)/4)\G((3-z)/4)}{\sqrt{\pi}}-\dfrac{p(n)}{q(n)}\sim-\dfrac{\cos(\pi z/2)\G((3+z)/4)^2\G((3-z)/4)^2/\pi^2}{2^nn}\;.$$
$$A=1+(z^2/4-9/4)/n+(z^4/32-17z^2/16+225/32)/n^2+\cdots$$
Series:
$$\dfrac{\pi^{1/2}}{\G((3+z)/4)\G((3-z)/4)}=1-\dfrac{z^2-1}{8}\sum_{n\ge0}\dfrac{((3-z)/2)_n((3+z)/2)_n}{(n+1)!^2}2^{-n}$$
Parametric family for $k\ge0$:
\begin{verbatim}
[(z)->gamma((3+z)/4)*gamma((3-z)/4)/sqrt(Pi),
12*n^2-4*n+1-z^2,-8*(n-k)*n*((2*n+1)^2-z^2)]
\end{verbatim}
Convergence type $E$ with $E=2$ and $P=3k+1$.
\end{cf}

\smallskip

\begin{cf}\label{6.8.29}{\ }
\begin{verbatim}
[(z)->(9-z^2)*gamma((3+z)/4)*gamma((3-z)/4)/sqrt(Pi),
[[8,4],[6*n+2,6*n+6]],
[[2*(z^2-1),9-z^2],[-(n+1)*((2*n+3)^2-z^2)/2,(n+1)*((4*n+3)^2-z^2)]]]
\end{verbatim}
$$\dfrac{(9-z^2)\G((3+z)/4)\G((3-z)/4)}{\sqrt{\pi}}=8+\dfrac{2z^2-2}{4-\dfrac{z^2-9}{8+\dfrac{z^2-25}{12-\dfrac{2z^2-98}{14+\dfrac{(3/2)z^2-147/2}{18-\dfrac{3z^2-363}{20+\ddots}}}}}}$$
Convergence type $E$ with $E=-2\sqrt{2}$, $P=3/2$, and $C=...$, so that
$$\dfrac{(9-z^2)\G((3+z)/4)\G((3-z)/4)}{\sqrt{\pi}}-\dfrac{p(n)}{q(n)}\sim(-1)^n\dfrac{C}{2^{3n/2}n^{3/2}}\;.$$
Series:
\begin{align*}&\dfrac{\G((3+z)/4)\G((3-z)/4)}{\sqrt{\pi}}=\\&8(1-z^2)\sum_{n\ge0}(-1)^n\dfrac{(3n+1)((3-z)/2)_n((3+z)/2)_n}{(12n^2-4n+1-z^2)(12n^2+20n+9-z^2)((3-z)/4)_n((3+z)/4)_n}2^{-3n}\end{align*}
\end{cf}

\smallskip

For the next CFs, note that
$$z(z-1/2)\left(\dfrac{\G(z)}{\G(z+1/2)}\right)^4=\dfrac{\G(z+1)\G(z)^3}{\G(z-1/2)\G(z+1/2)^3}$$

\smallskip

\begin{cf}\label{6.8.30}{\ }
\begin{verbatim}
[(z)->z*(z-1/2)*(gamma(z)/gamma(z+1/2))^4,
[1,(4*z-1)^2+2,8*n^2-8*n+(4*z-1)^2+1],[-2,-(4*n^2-1)^2]]
\end{verbatim}
$$z(z-1/2)\left(\dfrac{\G(z)}{\G(z+1/2)}\right)^4=1-\dfrac{2}{16z^2-8z+3-\dfrac{9}{16z^2-8z+18-\dfrac{225}{16z^2-8z+50-\ddots}}}$$
Convergence type $P^+$ with $P=4z-1$ and $C=...$, so that
$$z(z-1/2)\left(\dfrac{\G(z)}{\G(z+1/2)}\right)^4-\dfrac{p(n)}{q(n)}\sim\dfrac{C}{n^{4z-1}}\;.$$
$$A=1-(2z-1/2)/n-(z(4z-1)(8z^3-56z^2+26z+9)/(6((4z-3)(4z+1))))/n^2+\cdots$$
Parametric family for $k\ge0$:
\begin{verbatim}
[(z)->z*(z-1/2)*(gamma(z)/gamma(z+1/2))^4,
8*n^2-8*n+1+(2*k+4*z-1)^2,-(4*n^2-1)^2]
\end{verbatim}
Convergence type $P^+$ with $P=2k+4z-1$.
\end{cf}

\smallskip

\begin{cf}\label{6.8.31}{\ }
\begin{verbatim}
[(z)->z*(z-1/2)*(gamma(z)/gamma(z+1/2))^4,
[1,4*z^2,2*n+4*z-3],[-1,(n+2*z-2)*(n+2*z-1)^2*(n+2*z)]]
\end{verbatim}
$$z(z-1/2)\left(\dfrac{\G(z)}{\G(z+1/2)}\right)^4=1-\dfrac{1}{4z^2+\dfrac{16z^4-4z^2}{4z+1+\dfrac{16z^4+32z^3+20z^2+4z}{4z+3+\dfrac{16z^4+64z^3+92z^2+56z+12}{4z+5+\ddots}}}}$$
Convergence type $P^-$ with $P=2$ and
$C=-(z/2)(z-1/2)(\G(z)/\G(z+1/2))^4$, so that
$$z(z-1/2)\left(\dfrac{\G(z)}{\G(z+1/2)}\right)^4-\dfrac{p(n)}{q(n)}\sim(-1)^{n+1}\dfrac{(z/2)(z-1/2)(\G(z)/\G(z+1/2))^4}{n^2}$$
$$A=1-(4z-1)/n+(12z^2-6z+1/2)/n^2-(32z^3-24z^2+4z)/n^3+\cdots$$
Series:
\begin{align*}z(z-1/2)\left(\dfrac{\G(z)}{\G(z+1/2)}\right)^4&=1-\dfrac{2z-1}{16z^3}\sum_{n\ge0}\dfrac{(4n+4z+1)(z+1/2)_n^4}{(2n+2z-1)(n+z+1)(z+1)_n^4}\\
\dfrac{1}{z(z-1/2)}\left(\dfrac{\G(z+1/2)}{\G(z)}\right)^4&=\dfrac{(2z-1)^2}{4z(z-1)}-\dfrac{1}{4z(2z-1)}\sum_{n\ge0}\dfrac{(4n+4z-1)(z)_n^4}{(n+z-1)(2n+2z+1)(z+1/2)_n^4}\end{align*}
Parametric family for $k\ge0$:
\begin{verbatim}
[(z)->z*(z-1/2)*(gamma(z)/gamma(z+1/2))^4,
(2*k+1)*(2*n+4*z-3),(n+2*z-2)*(n+2*z-1)^2*(n+2*z)]]
\end{verbatim}
Convergence type $P^-$ with $P=4k+2$.
\end{cf}
  
\smallskip

\begin{cf}\label{6.8.32}{\ }
\begin{verbatim}
[(z)->z*(z-1/2)*(gamma(z)/gamma(z+1/2))^4,
[[1,8*z^2],[10*n^2+(2*z-1)*8*n+8*z^2-8*z+1,10*n^2+16*z*n+8*z^2-1]],
[[-2,16*z^2*(4*z^2-1)],[-(4*n^2-1)^2,4*((n+2*z)^2-1)*(n+2*z)^2]]]
\end{verbatim}
$$z(z-1/2)\left(\dfrac{\G(z)}{\G(z+1/2)}\right)^4=1-\dfrac{2}{8z^2+\dfrac{64z^4-16z^2}{8z^2+8z+3-\dfrac{9}{8z^2+16z+9+\dfrac{64z^4+128z^3+80z^2+16z}{8z^2+24z+25-\ddots}}}}$$
Convergence type $E$ with $E=-((1+\sqrt{5})/2)^5$, $P=0$, and $C=-4z(2z-1)(\G(z)/\G(z+1/2))^4/((1+\sqrt{5})/2)^{12z-1}$, so that
$$z(2z-1)\left(\dfrac{\G(z)}{\G(z+1/2)}\right)^4-\dfrac{p(n)}{q(n)}\sim(-1)^{n+1}\dfrac{4z(2z-1)(\G(z)/\G(z+1/2))^4}{((1+\sqrt{5})/2)^{5n+12z-1}}\;.$$
$$A=1+((16z^2-16z+7)d/5)/n+\cdots$$
\end{cf}

\smallskip

\begin{cf}\label{6.8.33}{\ }
\begin{verbatim}
[(z)->z^2*(gamma(z)/gamma(2*z))^4*gamma(4*z),
[2*(4*z-1),(z+1)^2,2*(n^2+z^2)],[2*(4*z-1)*(z-1)^2,-(n+z)^2*(n+1-z)^2]]
\end{verbatim}
$$z^2\left(\dfrac{\G(z)}{\G(2z)}\right)^4\G(4z)=8z-2+\dfrac{8z^3-18z^2+12z-2}{z^2+2z+1-\dfrac{z^4-2z^3-3z^2+4z+4}{2z^2+8-\dfrac{z^4-2z^3-11z^2+12z+36}{2z^2+18-\ddots}}}$$
Convergence type $P^+$ with $P=|4z-1|$ and $C=...$, so that
$$z^2\left(\dfrac{\G(z)}{\G(2z)}\right)^4\G(4z)-\dfrac{p(n)}{q(n)}\sim\dfrac{C}{n^{|4z-1|}}
$$
$$A=1+(1-4z)/n+((8z^4+86z^3-2z^2-5z)/(3(4z+1)))/n^2+\cdots$$
Series:
$$z^2\left(\dfrac{\G(z)}{\G(2z)}\right)^4\G(4z)=8z-2+\dfrac{2(4z-1)(z-1)^2}{(z+1)^2}\sum_{n\ge0}\dfrac{(2-z)_n^2}{(2+z)_n^2}$$
Parametric family for $k\ge0$:
\begin{verbatim}
[(z)->z^2*(gamma(z)/gamma(2*z))^4*gamma(4*z),
2*n^2+2*z^2+4*k*z+k*(k-1),-(n+z)^2*(n+1-z)^2]
\end{verbatim}
Convergence type $P^+$ with $P=4z+2k-1$.
\end{cf}

\smallskip

\begin{cf}\label{6.8.34}{\ }
\begin{verbatim}
[(z)->z^2*(gamma(z)/gamma(2*z))^4*gamma(4*z),
[2*(4*z-1),4*n+8*z-6],[4*(4*z-1)*(z-1)^2,n*(n+4*z-2)*(n+2*z-1)^2]]
\end{verbatim}
$$z^2\left(\dfrac{\G(z)}{\G(2z)}\right)^4\G(4z)=8z-2+\dfrac{16z^3-36z^2+24z-4}{8z-2+\dfrac{16z^3-4z^2}{8z+2+\dfrac{32z^3+32z^2+8z}{8z+6+\dfrac{48z^3+108z^2+72z+12}{8z+10+\ddots}}}}$$
Convergence type $P^-$ with $P=4$ and $C=...$, so that
$$z^2\left(\dfrac{\G(z)}{\G(2z)}\right)^4\G(4z)-\dfrac{p(n)}{q(n)}\sim(-1)^n\dfrac{C}{n^4}$$
$$A=1+(2-8z)/n+(44z^2-24z+1)/n^2+\cdots$$

There exist hypergeometric series for the CF and for its inverse, need
to use the square root of $8z^2-8z+1$, to be done.

Parametric family for $k\ge0$:
\begin{verbatim}
[(z)->z^2*(gamma(z)/gamma(2*z))^4*gamma(4*z),
(k+1)*(4*n+8*z-6),n*(n+4*z-2)*(n+2*z-1)^2]
\end{verbatim}
Convergence type $P^-$ with $P=4(k+1)$.
\end{cf}
  
\medskip

\subsection{Two and Three Free Variables}

{\ }

\medskip

Recall the beta function:

\begin{verbatim}
B(a,b)=gamma(a)*gamma(b)/gamma(a+b)
\end{verbatim}

\smallskip

\begin{cf}\label{6.8.35}{\ }
\begin{verbatim}
[(a,b)->B(a,b),[1/a,a+1,2*n^2+(2*a-b-1)*n+b*(1-a)],
               [1-b,-n*(n+a)^2*(n+1-b)]]
\end{verbatim}
$$B(a,b)=1/a-\dfrac{b-1}{a+1+\dfrac{(b-2)a^2+(2b-4)a+b-2}{(-b+4)a-b+6+\dfrac{(2b-6)a^2+(8b-24)a+8b-24}{(-b+6)a-2b+15+\ddots}}}$$
Convergence type $P^+$ with $P=b$ and $C=...$, so that
$$B(a,b)-\dfrac{p(n)}{q(n)}\sim\dfrac{C}{n^b}\;.$$
$$A=1-(b(2a+2b+1-b^2)/(2(b+1)))/n+\cdots$$
Series:
$$B(a,b)=\dfrac{1}{a}+(1-b)\sum_{n\ge0}\dfrac{(2-b)_n}{(n+a+1)(n+1)!}$$
Parametric family for $k\ge0$:
\begin{verbatim}
[(a,b)->B(a,b),2*n^2+(2*a-b-1)*n-a*b+(k+1)*b+k^2,-n*(n+a)^2*(n+1-b)]
\end{verbatim}
Convergence type $P^+$ with $P=2k+b$.
\end{cf}

\smallskip

\begin{verbatim}
G1(a,z)=2*(gamma((z+a+1)/2)*gamma((z-a)/2+1))/
                           (gamma((z+a)/2)*gamma((z-a+1)/2));
\end{verbatim}

\smallskip

\begin{cf}\label{6.8.36}{\ }
\begin{verbatim}
[(a,z)->G1(a,z),[z,2*z],[(n+a)*(n+1-a)]]
\end{verbatim}
$$2\dfrac{\G((z+a+1)/2)\G((z-a)/2+1)}{\G((z+a)/2)\G((z-a+1)/2)}=z-\dfrac{a^2-a}{2z-\dfrac{a^2-a-2}{2z-\dfrac{a^2-a-6}{2z-\dfrac{a^2-a-12}{2z-\dfrac{a^2-a-20}{2z-\dfrac{a^2-a-30}{2z-\ddots}}}}}}$$
Convergence type $P^-$ with $P=2z$ and $C=...$, so that
$$2\dfrac{\G((z+a+1)/2)\G((z-a)/2+1)}{\G((z+a)/2)\G((z-a+1)/2)}-\dfrac{p(n)}{q(n)}\sim(-1)^n\dfrac{C}{n^{2z}}\;.$$
$$A=1-2z/n+(-z^3/6+2z^2+(a^2/2-a/2+2/3)z)/n^2+\cdots$$
Parametric family for $k\ge0$:
\begin{verbatim}
[(a,z)->G1(a,z),2*z+2*k,(n+a)*(n+1-a)]
\end{verbatim}
Convergence type $P^-$ with $P=2z+2k$.
\end{cf}

\smallskip

\begin{cf}\label{6.8.37}{\ }
\begin{verbatim}
[(a,z)->G1(a,z),[z,z+1,2],[-a^2+a,(n+z-a)*(n+z+a-1)]]
\end{verbatim}
$$2\dfrac{\G((z+a+1)/2)\G((z-a)/2+1)}{\G((z+a)/2)\G((z-a+1)/2)}=z-\dfrac{a^2-a}{z+1+\dfrac{z^2+z+(-a^2+a)}{2+\dfrac{z^2+3z+(-a^2+a+2)}{2+\dfrac{z^2+5z+(-a^2+a+6)}{2+\ddots}}}}$$
Convergence type $P^-$ with $P=2$ and $C=...$, so that
$$2\dfrac{\G((z+a+1)/2)\G((z-a)/2+1)}{\G((z+a)/2)\G((z-a+1)/2)}-\dfrac{p(n)}{q(n)}\sim(-1)^n\dfrac{C}{n^2}\;.$$
$$A=1-2z/n+(3z^2+a^2/2-a/2-1/2)/n^2+\cdots$$
Series:
\begin{align*}2&\dfrac{\G((z+a+1)/2)\G((z-a)/2+1)}{\G((z+a)/2)\G((z-a+1)/2)}\\
  &=z+\dfrac{2a(1-a)}{(z+1)(z+2)+a(1-a)}\sum_{n\ge0}\dfrac{((z+a)/2)_n((z+1-a)/2)_n}{((z+a+3)/2)_n((z-a)/2+2)_n}\\
  &\dfrac{\G((z+a)/2)\G((z-a+1)/2)}{2\G((z+a+1)/2)\G((z-a)/2+1)}\\
  &=\dfrac{z-1}{(z-a)(z+a-1)}+\dfrac{2a(1-a)}{(z^2-a^2)(z^2-(1-a)^2)}\sum_{n\ge0}\dfrac{((z+a-1)/2)_n((z-a)/2)_n}{((z+a)/2+1)_n((z-a+3)/2)_n}\end{align*}
Parametric family for $k\ge0$:
\begin{verbatim}
[(a,z)->G1(a,z),2*k+2,(n+z-a)*(n+z+a-1)]]
\end{verbatim}
Convergence type $P^-$ with $P=2k+2$.
\end{cf}

\smallskip

\begin{cf}\label{6.8.38}{\ }
\begin{verbatim}
[(a,z)->G1(a,z),
[z,n+z],[[-a*(a-1),(z+a)*(z+1-a)],[(n+a)*(n+1-a),(n+z+a)*(n+z+1-a)]]]
\end{verbatim}
$$2\dfrac{\G((z+a+1)/2)\G((z-a)/2+1)}{\G((z+a)/2)\G((z-a+1)/2)}=z-\dfrac{a^2-a}{z+1+\dfrac{z^2+z-a^2+a}{z+2-\dfrac{a^2-a-2}{z+3+\dfrac{z^2+3z-a^2+a+2}{z+4-\ddots}}}}$$
Convergence type $E$ with $E=-(1+\sqrt{2})^2$, $P=0$, and $C=...$, so that
$$2\dfrac{\G((z+a+1)/2)\G((z-a)/2+1)}{\G((z+a)/2)\G((z-a+1)/2)}-\dfrac{p(n)}{q(n)}\sim(-1)^n\dfrac{C}{(1+\sqrt{2})^{2n}}\;.$$
$$A=1+((z^2-2z+2a^2-2a+9/8)d)/n+\cdots$$
\end{cf}

\smallskip

\begin{verbatim}
G2(a,z)=(gamma((z+a+1)/4)*gamma((z-a+1)/4))/
                         (gamma((z+a+3)/4)*gamma((z-a+3)/4));
\end{verbatim}

\smallskip

\begin{cf}\label{6.8.39}{\ }
\begin{verbatim}
[(a,z)->G2(a,z),[0,z,2*z],[4,(2*n-1)^2-a^2]]
\end{verbatim}
$$\dfrac{\G((z+a+1)/4)\G((z-a+1)/4)}{\G((z+a+3)/4)\G((z-a+3)/4)}=\dfrac{4}{z-\dfrac{a^2-1}{2z-\dfrac{a^2-9}{2z-\dfrac{a^2-25}{2z-\dfrac{a^2-49}{2z-\dfrac{a^2-81}{2z-\ddots}}}}}}$$
Convergence type $P^-$ with $P=z$ and $C=...$, so that
$$\dfrac{\Gamma((z+a+1)/4)\Gamma((z-a+1)/4)}{\Gamma((z+a+3)/4)\Gamma((z-a+3)/4)}-\dfrac{p(n)}{q(n)}\sim(-1)^n\dfrac{C}{n^z}\;.$$
$$A=1+(-z^3/48+(a^2/16-11/48)z)/n^2+\cdots$$
Parametric family for $k\ge0$:
\begin{verbatim}
[(a,z)->G2(a,z),2*z+4*k,(2*n-1)^2-a^2]]
\end{verbatim}
Convergence type $P^-$ with $P=2k+z$.
\end{cf}

\smallskip

\begin{cf}\label{6.8.40}{\ }
\begin{verbatim}
[(a,z)->G2(a,z),[4/z,z+2+(1-a^2)/z,4],[4*(a^2-1)/z^2,(2*n+z-1)^2-a^2]]
\end{verbatim}
$$\dfrac{\G((z+a+1)/4)\G((z-a+1)/4)}{\G((z+a+3)/4)\G((z-a+3)/4)}=\dfrac{4}{z}+\dfrac{(4a^2-4)/z^2}{z+2+(1-a^2)/z+\dfrac{z^2+2z+(-a^2+1)}{4+\dfrac{z^2+6z+(-a^2+9)}{4+\ddots}}}$$
Convergence type $P^-$ with $P=2$ and $C=...$, so that
$$\dfrac{\G((z+a+1)/4)\G((z-a+1)/4)}{\G((z+a+3)/4)\G((z-a+3)/4)}-\dfrac{p(n)}{q(n)}\sim(-1)^n\dfrac{C}{n^2}\;.$$
$$A=1-z/n+(3z^2/4+a^2/8-5/8)/n^2+(-z^3/2+(-a^2/4+5/4)z)/n^3+\cdots$$
Series:
\begin{align*}&\dfrac{\G((z+a+1)/4)\G((z-a+1)/4)}{\G((z+a+3)/4)\G((z-a+3)/4)}\\
  &=\dfrac{4}{z}+\dfrac{a^2-1}{(z+1-a)(z+1+a)}\sum_{n\ge0}\dfrac{((z-a+3)/4)_n((z+a+3)/4)_n}{(n+z/4)(n+z/4+1)((z-a+5)/4)_n((z+a+5)/4)_n}\\
  &\dfrac{\G((z+a+3)/4)\G((z-a+3)/4)}{\G((z+a+1)/4)\G((z-a+1)/4)}\\
  &=\dfrac{(z+a-1)(z-a-1)}{4(z-2)}+(a^2-1)\sum_{n\ge0}\dfrac{((z-a+1)/4)_n((z+a+1)/4)_n}{(4n+z-2)(4n+z+2)((z-a+3)/4)_n((z+a+3)/4)_n}\end{align*}
Parametric family for $k\ge0$:
\begin{verbatim}
[(a,z)->G2(a,z),4*k+4,(2*n+z-1)^2-a^2]
\end{verbatim}
Convergence type $P^-$ with $P=2k+2$.
\end{cf}
  
\smallskip

\begin{cf}\label{6.8.41}{\ }
\begin{verbatim}
[(a,z)->G2(a,z),[0,z+2*n-2],[[4,1-a^2],[(2*n+z-1)^2-a^2,(2*n+1)^2-a^2]]]
\end{verbatim}
$$\dfrac{\G((z+a+1)/4)\G((z-a+1)/4)}{\G((z+a+3)/4)\G((z-a+3)/4)}=\dfrac{4}{z-\dfrac{a^2-1}{z+2+\dfrac{z^2+2z+1-a^2}{z+4-\dfrac{a^2-9}{z+6+\dfrac{z^2+6z+9-a^2}{z+8-\dfrac{a^2-25}{z+10+\ddots}}}}}}$$
Convergence type $E$ with $E=-(1+\sqrt{2})^2$, $P=0$, and $C=...$, so that
$$\dfrac{\G((z+a+1)/4)\G((z-a+1)/4)}{\G((z+a+3)/4)\G((z-a+3)/4)}-\dfrac{p(n)}{q(n)}\sim(-1)^n\dfrac{C}{(1+\sqrt{2})^{2n}}\;.$$
$$A=1+((z^2/5+a^2/2-3/8)d)/n+\cdots$$
\end{cf}

\smallskip

\begin{verbatim}
G3(a,b)=gamma(2*a+b+1)*gamma(a+1)*gamma(b-a+1)/
                (gamma(2*a+1)*gamma(b+1)*gamma(1-a)*gamma(b+a+1));
\end{verbatim}

\smallskip

\begin{cf}\label{6.8.42}{\ }
\begin{verbatim}
[(a,b)->G3(a,b),[1,2*a+b+1,(2*n+2*a-1)*(n^2+(2*a-1)*n-(2*a-1)*(b+1))],
                [-4*a^2*b,-n^2*(n+2*a)^2*(n-b)*(n+2*a+b)]]
\end{verbatim}
$$\dfrac{\G(2a+b+1)\G(a+1)\G(b-a+1)}{\G(2a+1)\G(b+1)\G(1-a)\G(b+a+1)}=1-\dfrac{4ba^2}{2a+b+1+\ddots}$$
Convergence type $P^+$ with $P=2(a-b-1)$ and $C=...$, so that
$$\dfrac{\G(2a+b+1)\G(a+1)\G(b-a+1)}{\G(2a+1)\G(b+1)\G(1-a)\G(b+a+1)}-\dfrac{p(n)}{q(n)}\sim\dfrac{C}{n^{2(a-b-1)}}$$
$$A=1+(-2a^2+2ab+a+b+1)/n+\cdots$$
Series:
$$\dfrac{\G(2a+b+1)\G(a+1)\G(b-a+1)}{\G(2a+1)\G(b+1)\G(1-a)\G(b+a+1)}=1-\dfrac{4a^2b}{2a+b+1}\sum_{n\ge0}\dfrac{(2a+1)_n^2(1-b)_n}{(2a+b+2)_n(n+1)!^2}$$
\end{cf}

\begin{verbatim}
G4(a,b)=gamma(2*a+b+1)*gamma(a+1)^2*gamma(b+1)/
                      (gamma(2*a+1)*gamma(a+b+1)^2);
\end{verbatim}

\smallskip

\begin{cf}\label{6.8.43}{\ }
\begin{verbatim}
[(a,b)->G4(a,b),
[1,(2*a+2)*(2*a+b+1),(2*n+2*a-1)*(2*n^2+2*(2*a-1)*n+2*(b+1)*(1-a))],
[-4*a^2*b,-4*n*(n-b)*(n+2*a)*(n+2*a+b)*(n+a)^2]]
\end{verbatim}
$$\dfrac{\G(2a+b+1)\G(a+1)^2\G(b+1)}{\G(2a+1)\G(a+b+1)^2}=1-\dfrac{4ba^2}{4a^2+(2b+6)a+2b+2+\ddots}$$
Convergence type $P^+$ with $P=2b+2$ and $C=...$, so that
$$\dfrac{\G(2a+b+1)\G(a+1)^2\G(b+1)}{\G(2a+1)\G(a+b+1)^2}-\dfrac{p(n)}{q(n)}\sim\dfrac{C}{n^{2b+2}}$$
$$A=1-(2ab+2a+b+1)/n+\cdots$$
Series:
$$\dfrac{\G(2a+b+1)\G(a+1)^2\G(b+1)}{\G(2a+1)\G(a+b+1)^2}=1-\dfrac{2a^2b}{2a+b+1}\sum_{n\ge0}\dfrac{(1-b)_n(2a+1)_n}{(n+a+1)(2a+b+2)_n(n+1)!}$$
\end{cf}

\smallskip

\begin{verbatim}
H1(a,b,z)=(gamma((z+a+b+1)/2)*gamma((z-a-b+1)/2))/
                            (gamma((z+a-b+1)/2)*gamma((z-a+b+1)/2));
\end{verbatim}

\smallskip

\begin{cf}\label{6.8.44}{\ }
\begin{verbatim}
[(a,b,z)->H1(a,b,z),[-1,1,(2*n-3)*z],[2,-a*b,(a^2-(n-1)^2)*(b^2-(n-1)^2)]]
\end{verbatim}
$$\dfrac{\G((z+a+b+1)/2)\G((z-a-b+1)/2)}{\G((z+a-b+1)/2)\G((z-a+b+1)/2)}=
 -1+\dfrac{2}{1+\dfrac{ba}{z+\dfrac{(b^4-2b^2+1)a^2+(-b^4+2b^2-1)}{3z+\ddots}}}$$
Convergence type $P^-$ with $P=2z$ and $C=...$, so that
$$\dfrac{\G((z+a+b+1)/2)\G((z-a-b+1)/2)}{\G((z+a-b+1)/2)\G((z-a+b+1)/2)}-\dfrac{p(n)}{q(n)}\sim(-1)^n\dfrac{C}{n^{2z}}\;.$$
$$A=1+z/n+(-z^3/6+z^2/2+(a^2/2+b^2/2-1/3)z)/n^2+\cdots$$
Parametric family for $k\ge0$:
\begin{verbatim}
[(a,b,z)->H1(a,b,z),(2*n-3)*(2*k+z),(a^2-(n-1)^2)*(b^2-(n-1)^2)]
\end{verbatim}
Convergence type $P^-$ with $P=4k+2z$.
\end{cf}

See also \ref{6.8.45} for an equivalent CF.

\smallskip

This CF can be Ap\'ery accelerated with convergence type
$E=-((1+\sqrt{5})/2)^5$, formulas too complicated to give here.

\medskip

\subsection{Ratios of Products of $3$, $4$, and $8$ Gamma Functions}

\medskip

Ramanujan considers the expressions
\begin{align*}
  P_k^+(x_0,x_1,\dotsc,x_k)&=\prod_{\substack{even\ number\\of\ -\ signs}}\G((1+x_0\pm x_1\pm\cdots\pm x_k)/2)\;,\\
  P_k^-(x_0,x_1,\dotsc,x_k)&=\prod_{\substack{odd\ number\\of\ -\ signs}}\G((1+x_0\pm x_1\pm\cdots\pm x_k)/2)\;,\text{\quad and}\\
  Q_k(x_0,x_1,\dotsc,x_k)&=(-1)^k\dfrac{P_k^+-P_k^-}{P_k^++P_k^-}\;.
\end{align*}
Note that $P^+$ and $P^-$ are products of $2^{k-1}$ gamma functions, that
\begin{align*}\lim_{x_1,x_2,\dotsc,x_k\to0}\dfrac{Q_k(x_0,x_1,\dotsc,x_k)}{x_1x_2\dotsc x_k}&=\dfrac{(k-1)!}{2}\z(k,(x_0+1)/2)\text{\quad and}\\
\lim_{x_k\to\infty}Q_k(x_0+x_k,x_1,\dotsc,x_k)&=Q_{k-1}(x_0,x_1,\dotsc,x_{k-1})\;,\end{align*}
and that the CFs that we will give specialize to the CFs \ref{6.11.5},
\ref{6.13.8}, and \ref{6.15.1}. These CFs can probably all be Ap\'ery
accelerated, but I have not done so.

\smallskip

We define $S_d=x_0^d+x_1^d+\cdots+x_k^d$.

\smallskip

\begin{cf}\label{6.8.45}{\ }
\begin{verbatim}
[(x_i)->Q2,[0,(2*n-1)*x0],[x1*x2,(x1^2-n^2)*(x2^2-n^2)]]
\end{verbatim}
$$Q_2=\dfrac{x_1x_2}{x_0+\dfrac{(x_1^2-1^2)(x_2^2-1^2)}{3x_0+\dfrac{(x_1^2-2^2)(x_2^2-2^2)}{5x_0+\ddots}}}\;.$$
Convergence type $P^-$ with $P=2x_0$ and $C=...$, so that
$$Q_2-\dfrac{p(n)}{q(n)}\sim(-1)^n\dfrac{C}{n^{2x_0}}\;.$$
$$A=1-x_0/n+(-x_0^3/6+x_0^2/2+(x_1^2/2+x_2^2/2-1/3)x_0)/n^2+\cdots$$
\end{cf}

This CF is equivalent to \ref{6.8.44}.

\smallskip

\begin{cf}\label{6.8.46}{\ }
\begin{verbatim}
[(x_i)->Q3,[0,(2*n-1)*(2*x0^2-S2+2*n^2-2*n+1)],
           [2*x1*x2*x3,4*(x1^2-n^2)*(x2^2-n^2)*(x3^2-n^2)]]
\end{verbatim}
$$Q_3=\dfrac{2x_1x_2x_3}{2x_0^2-S_2+1+\dfrac{4(x_1^2-1^2)(x_2^2-1^2)(x_3^2-1^2)}{3(2x_0^2-S_2+5)+\ddots}}\;.$$
Convergence type $P^-$ with $P=2x_0$ and $C=...$, so that
$$Q_3-\dfrac{p(n)}{q(n)}\sim(-1)^n\dfrac{C}{n^{2x_0}}\;.$$
$$A=1-x_0/n+(P(x)/(24(x_0^2-1)))/n^2+\cdots$$
with a complicated polynomial $P$.
\end{cf}

\smallskip

\begin{cf}\label{6.8.47}{\ }
\begin{verbatim}
[(x_i)->Q4,[0,(2*n-1)*(2*(S4+1)-(S2-2*n^2+2*n-1)^2-4*(n^2-n+1)^2)],
[8*x0*x1*x2*x3*x4,
       64*(x0^2-n^2)*(x1^2-n^2)*(x2^2-n^2)*(x3^2-n^2)*(x4^2-n^2)]]
\end{verbatim}
$$Q_4=\dfrac{8x_0x_1x_2x_3x_4}{2(S_4+1)-(S_2-1)^2-2^2+\dfrac{64(x_0^2-1^2)(x_1^2-1^2)\cdots(x_4^2-1^2)}{3(2(S_4+1)-(S_2-5)^2-6^2)+\ddots}}\;.$$
\end{cf}

Note that the domain of validity (not only of convergence) of this CF is quite
complicated and has been partly understood only comparatively recently
\cite{Mas}, so we do not give its speed of convergence.
For instance, it is known to be correct
if one of the $x_i$ for $i\ge1$ is an integer (in which case the CF is of
course finite). On the contrary, if for instance all the $x_i$ are equal to
$1/2$, the CF does converge, but not to the expected value.

\smallskip

In addition to Ramanujan's formulas for ratios of products of $2$, $4$, and
$8$ gamma functions, the paper \cite{CTZ} gives CFs for ratios of products of
three gamma functions, of which \ref{6.8.25}, \ref{6.8.26}, and \ref{6.8.27}
are special cases. Almost all the CFs have period 2, except two of them
which are similar to Ramanujan's. We give one of them (the other is similar
but more complicated), referring to \cite{CTZ} for all the others.

\smallskip

Here we set
\begin{align*}
  P^+(x_0,x_1,x_2)&=\G((1+x_0+x_1)/2)\G((1+x_0+x_2)/2)\G((1+x_0-x_1-x_2)/2)\;,\\
  P^-(x_0,x_1,x_2)&=\G((1+x_0-x_1)/2)\G((1+x_0-x_2)/2)\G((1+x_0+x_1+x_2)/2)\;,\\
  &\text{\quad and\quad}Q(x_0,x_1,x_2)=\dfrac{P^+-P^-}{P^++P^-}\;.
\end{align*}

\smallskip

\begin{cf}\label{6.8.48}{\ }
\begin{verbatim}
[(x_i)->Q,[0,(4*n-2)*(2*x0^2+4*n^2-4*n+2-(x1^2+x1*x2+x2^2))],
[x1*x2*(x1+x2),-(4*n^2-x1^2)*(4*n^2-x2^2)*(4*n^2-(x1+x2)^2)]]
\end{verbatim}
$$Q=\dfrac{x_2x_1^2+x_2^2x_1}{4x_0^2-2x_1^2-2x_2x_1-2x_2^2+4+\ddots}$$
Convergence type $P^+$ with $P=2x_0$ and $C=...$, so that
$$Q-\dfrac{p(n)}{q(n)}\sim\dfrac{C}{n^{2x_0}}\;.$$
$$A=1-x_0/n+(-x_0^3/24+x_0^2/2+(x_1^2/8+x_1x_2/8+x_2^2/8+1/24)x_0)/n^2+\cdots$$
\end{cf}

\medskip

\section{Function $\psi(z)=\G'(z)/\G(z)$}\label{sec:psi0}

\medskip

First, recall that the function $\psi$ and its derivatives take explicit
values at rational arguments, which involve values of $L$-functions of
Dirichlet characters at positive integers. For completeness, we give here
small tables of such values, and also for the functions $\beta_j$ that we
will define below.

\medskip

\centerline{
\begin{tabular}{|c||c|c|c||}
  \hline
  $r/m$ & $\psi(r/m)+\ga$ & $\psi'(r/m)$ & $\psi''(r/m)$ \\
  \hline\hline
  1 & $0$ & $\z(2)$ & $-2\z(3)$ \\
  1/2 & $-2\log(2)$ & $3\z(2)$ & $-14\z(3)$ \\
  1/3 & $-3\log(3)/2-\pi/(2\sqrt{3})$ & $4\z(2)+9G_3/2$ & $-26\z(3)-4\pi^3/(3\sqrt{3})$ \\
  2/3 & $-3\log(3)/2+\pi/(2\sqrt{3})$ & $4\z(2)-9G_3/2$ & $-26\z(3)+4\pi^3/(3\sqrt{3})$ \\
  1/4 & $-3\log(2)-\pi/2$ & $6\z(2)+8G$ & $-56\z(3)-2\pi^3$ \\
  3/4 & $-3\log(2)+\pi/2$ & $6\z(2)-8G$ & $-56\z(3)+2\pi^3$ \\
  1/6 & $-\log(432)/2-3\pi/(2\sqrt{3})$ & $12\z(2)+45G_3/2$ & $-182\z(3)-12\pi^3/\sqrt{3}$ \\
  5/6 & $-\log(432)/2+3\pi/(2\sqrt{3})$ & $12\z(2)-45G_3/2$ & $-182\z(3)+12\pi^3/\sqrt{3}$ \\
  \hline
\end{tabular}}

\bigskip

\centerline{
\begin{tabular}{|c||c|c|c||}
  \hline
  $r/m$ & $\beta(r/m)$ & $\beta_1(r/m)$ & $\beta_2(r/m)$ \\
  \hline\hline
  1 & $\log(2)$ & $\z(2)/2$ & $3\z(3)/2$ \\
  1/2 & $\pi/2$ & $4G$ & $\pi^3/2$ \\
  1/3 & $\log(2)+\pi/\sqrt{3}$ & $27G_3/4+2\z(2)$ & $39\z(3)/2+5\pi^3/(3\sqrt{3})$ \\
  2/3 & $-\log(2)+\pi/\sqrt{3}$ & $27G_3/4-2\z(2)$ & $-39\z(3)/2+5\pi^3/(3\sqrt{3})$ \\
  1/4 & $\pi/\sqrt{2}+\sqrt{2}\log(1+\sqrt{2})$ & $8G_8+\pi^2/\sqrt{2}$ & $3\pi^3/\sqrt{2}+64L(\chi_8,3)$ \\
  3/4 & $\pi/\sqrt{2}-\sqrt{2}\log(1+\sqrt{2})$ & $8G_8-\pi^2/\sqrt{2}$ & $3\pi^3/\sqrt{2}-64L(\chi_8,3)$ \\
  1/6 & $\pi+\sqrt{3}\log(2+\sqrt{3})$ & $20G+\sqrt{3}\pi^2$ & $7\pi^3+216L(\chi_{12},3)$ \\
  5/6 & $\pi-\sqrt{3}\log(2+\sqrt{3})$ & $20G-\sqrt{3}\pi^2$ & $7\pi^3-216L(\chi_{12},3)$ \\
  \hline
\end{tabular}}

\bigskip

In the above, we have set for simplicity $\chi_d(n)=\lgs{d}{n}$,
$G_d=L(\chi_{-d},2)$, and $G=G_4$ is Catalan's constant.

\smallskip

Second, since the functions $\psi^{(k)}(z)$ and $\be_j(z)$ satisfy
simple additive functional equations when $z$ is replaced by $z+u$ with
$u\in\Z$, in all the parametric CFs that we give below we can replace
$z$ by $z+u$, thus adding an extra parameter.

\bigskip

Note that
\begin{align*}&\int_0^\infty\dfrac{1-e^{-tz}}{\sinh(t)}\,dt=\psi((1+z)/2)+\ga+2\log(2)\;,\\
&\int_0^\infty(1-e^{-tz})(\coth(t)-1)\,dt=\psi(1+z/2)+\ga\;,\quad\text{and}\\
&\int_0^\infty\dfrac{e^{-tz}-1+tz}{\sinh^2(t)}\,dt=z(\psi(1+z/2)+\ga)\;,
\end{align*}
so all CFs for $\psi(1+z)+\ga$ give CFs for these integrals.

\smallskip

Note also that
$$\psi(1+z)+\ga=z\sum_{n\ge1}\dfrac{1}{n(n+z)}\;.$$

\smallskip

\begin{cf}\label{6.9.1}{\ }
\begin{verbatim}
[(z)->psi(1+z)+Euler(),[0,2*n^2+(z-1)*(2*n-1)],[z,-n^2*(n+z)^2]]
\end{verbatim}
$$\psi(1+z)+\ga=\dfrac{z}{z+1-\dfrac{z^2+2z+1}{3z+5-\dfrac{4z^2+16z+16}{5z+13-\dfrac{9z^2+54z+81}{7z+25-\dfrac{16z^2+128z+256}{9z+41-\ddots}}}}}$$
Convergence type $P^+$ with $P=1$ and $C=z$, so that
$$\psi(1+z)+\ga-\dfrac{p(n)}{q(n)}\sim\dfrac{z}{n}\;.$$
$$A=1-((z+1)/2)/n+((z+1)(2z+1)/6)/n^2-(z(z+1)^2/4)/n^3+\cdots$$
Series:
$$\psi(1+z)+\ga=z\sum_{n\ge1}\dfrac{1}{n(n+z)}$$
Parametric family for $k\ge0$:
\begin{verbatim}
[(z)->psi(1+z)+Euler(),2*n^2+(z-1)*(2*n-1)+k^2+k,-n^2*(n+z)^2]
\end{verbatim}
Convergence type $P^+$ with $P=2k+1$.
\end{cf}

\smallskip

\begin{cf}\label{6.9.2}{\ }
\begin{verbatim}
[(z)->psi(1+z)+Euler(),[0,2*n^2-2*n+1-(n-1)*z],[z,-n^3*(n-z)]]
\end{verbatim}
$$\psi(1+z)+\ga=\dfrac{z}{1+\dfrac{z-1}{-z+5+\dfrac{8z-16}{-2z+13+\dfrac{27z-81}{-3z+25+\dfrac{64z-256}{-4z+41+\dfrac{125z-625}{-5z+61+\ddots}}}}}}$$
Convergence type $P^+$ with $P=z+1$ and $C=-1/((z+1)\G(-z))$, so that
$$\psi(1+z)+\ga-\dfrac{p(n)}{q(n)}\sim-\dfrac{1/((z+1)\G(-z))}{n^{z+1}}\;.$$
$$A=1+((z+1)(z^2-2)/(2(z+2)))/n+\cdots$$
Series:
$$\psi(1+z)+\ga=z\sum_{n\ge0}\dfrac{(1-z)_n}{(n+1)(n+1)!}$$
Parametric family for $k\ge0$:
\begin{verbatim}
[(z)->psi(1+z)+Euler(),2*n^2-2*n+k^2+k+1-(n-k-1)*z,-n^3*(n-z)]
\end{verbatim}
Convergence type $P^+$ with $P=z+2k+1$.
\end{cf}

\smallskip

\begin{cf}\label{6.9.3}{\ }
\begin{verbatim}
[(z)->psi(1+z)+Euler(),[0,(z+1)*(2*n-1)],[2*z,n^2*(n^2-z^2)]]
\end{verbatim}
$$\psi(1+z)+\ga=\dfrac{2z}{z+1-\dfrac{z^2-1}{3z+3-\dfrac{4z^2-16}{5z+5-\dfrac{9z^2-81}{7z+7-\dfrac{16z^2-256}{9z+9-\dfrac{25z^2-625}{11z+11-\ddots}}}}}}$$
Convergence type $P^-$ with $P=2z+2$ and $C=\G(z+1)^2\sin(\pi z)/\pi$, so that
$$\psi(1+z)+\ga-\dfrac{p(n)}{q(n)}\sim\dfrac{\G(z+1)^2\sin(\pi z)/\pi}{n^{2z+2}}\;.$$
$$A=1-(z+1)/n+(z(z+1)(2z+1)/6)/n^2-((z+1)^2(z+2)(2z-3)/3)/n^3+\cdots$$
Series:
$$\psi(1+z)+\ga=\dfrac{2z}{z+1}\sum_{n\ge0}(-1)^n\dfrac{(1-z)_n}{(n+1)(2+z)_n}$$
Parametric family for $k\ge0$:
\begin{verbatim}
[(z)->psi(1+z)+Euler(),(z+2*k+1)*(2*n-1),n^2*(n^2-z^2)]
\end{verbatim}
Convergence type $P^-$ with $P=4k+2z+2$.
\end{cf}

\smallskip

\begin{cf}\label{6.9.4}{\ }
\begin{verbatim}
[(z)->psi(1+z)+Euler(),
[[0,1+z],[5*n^2+3*z*n,5*n^2+(3*z+4)*n+z+1]],
[[2*z,1-z^2],[-4*n^2*(n+z)^2,(n+1)^2*(n-z+1)*(n+z+1)]]]
\end{verbatim}
$$\psi(1+z)+\ga=\dfrac{2z}{z+1-\dfrac{z^2-1}{3z+5-\dfrac{4z^2+8z+4}{4z+10-\dfrac{4z^2-16}{6z+20-\dfrac{16z^2+64z+64}{7z+29-\dfrac{9z^2-81}{9z+45-\ddots}}}}}}$$
Convergence type $E$ with $E=-((1+\sqrt{5})/2)^5$, $P=0$, and $C=4\pi\sin(\pi z)/((1+\sqrt{5})/2)^{2z+5}$, so that
$$\psi(1+z)+\ga-\dfrac{p(n)}{q(n)}\sim(-1)^n\dfrac{4\pi\sin(\pi z)}{((1+\sqrt{5})/2)^{5n+2z+5}}\;.$$
\begin{align*}A&=1+((6z^2/5-2/5)d)/n\\
  &\phantom{=}+(18z^4/5-8z^3d/25-(6d/5+12/5)z^2+4zd/25+2d/5+2/5)/n^2+\cdots\end{align*}
\end{cf}

\smallskip

\begin{cf}\label{6.9.5}{\ }
\begin{verbatim}
[(z)->psi(1+z)+Euler(),
[[0,2+2*z],[5*n^2+3*z*n,5*n^2+(3*z+6)*n+2*z+2]],
[[2*z,-(2+2*z)^2],[n^2*(n^2-z^2),-4*(n+1)^2*(n+z+1)^2]]]
\end{verbatim}
$$\psi(1+z)+\ga=\dfrac{2z}{2z+2-\dfrac{4z^2+8z+4}{3z+5-\dfrac{z^2-1}{5z+13-\dfrac{16z^2+64z+64}{6z+20-\dfrac{4z^2-16}{8z+34-\dfrac{36z^2+216z+324}{9z+45-\ddots}}}}}}$$
Convergence type $E$ with $E=-((1+\sqrt{5})/2)^5$, $P=0$, and $C=4\pi\sin(\pi z)/((1+\sqrt{5})/2)^{2z+5}$, so that
$$\psi(1+z)+\ga-\dfrac{p(n)}{q(n)}\sim(-1)^n\dfrac{4\pi\sin(\pi z)}{((1+\sqrt{5})/2)^{5n+2z+5}}\;.$$
\begin{align*}A&=1+((6z^2/5-2/5)d)/n\\
  &\phantom{=}+(18z^4/5-8z^3d/25-(6d/5+12/5)z^2+4zd/25+2d/5+2/5)/n^2+\cdots\end{align*}
\end{cf}

\smallskip

\begin{cf}\label{6.9.6}{\ }
\begin{verbatim}
[(z)->(psi(1+z)+psi(1-z))/2+Euler,
[0,(2*n-1)*(n^2-n+1-z^2)],[-z^2,-n^2*(n^2-z^2)^2]]
\end{verbatim}
$$\dfrac{\psi(1+z)+\psi(1-z)}{2}+\ga=-\dfrac{z^2}{-z^2+1-\dfrac{z^4-2z^2+1}{-3z^2+9-\dfrac{4z^4-32z^2+64}{-5z^2+35-\dfrac{9z^4-162z^2+729}{-7z^2+91-\ddots}}}}$$
Convergence type $P^+$ with $P=2$ and $C=-z^2/2$, so that
$$\dfrac{\psi(1+z)+\psi(1-z)}{2}+\ga-\dfrac{p(n)}{q(n)}\sim-\dfrac{z^2/2}{n^2}\;.$$
$$A=1-1/n+(z^2/2+1/2)/n^2-z^2/n^3+(z^4/3+5z^2/6-1/6)/n^4+\cdots$$
Series:
$$\dfrac{\psi(1+z)+\psi(1-z)}{2}+\ga=-z^2\sum_{n\ge1}\dfrac{1}{n(n-z)(n+z)}$$
Parametric family for $k\ge0$:
\begin{verbatim}
[(z)->(psi(1+z)+psi(1-z))/2+Euler,
(2*n-1)*(n^2-n+1-z^2+2*k^2+2*k),-n^2*(n^2-z^2)^2]
\end{verbatim}
Convergence type $P^+$ with $P=4k+2$.
\end{cf}

\smallskip

\begin{cf}\label{6.9.7}{\ }
\begin{verbatim}
[(z)->(psi(1+z)+psi(1-z))/2+Euler,
[[0,1-z^2],[2*n*(3*n^2-z^2),(2*n+1)*(3*n^2+3*n+1-z^2)]],
[[-z^2,-(1-z^2)^2],[-(n^3-z^2*n)^2,-(n+1)^2*((n+1)^2-z^2)^2]]]
\end{verbatim}
$$\dfrac{\psi(1+z)+\psi(1-z)}{2}+\ga=-\dfrac{z^2}{-z^2+1-\dfrac{z^4-2z^2+1}{-2z^2+6-\dfrac{z^4-2z^2+1}{-3z^2+21-\dfrac{4z^4-32z^2+64}{-4z^2+48-\ddots}}}}$$
Convergence type $E$ with $E=-(1+\sqrt{2})^4$, $P=0$, and
$C=-4\pi\sin(\pi z)^2/(1+\sqrt{2})^4$, so that
$$\dfrac{\psi(1+z)+\psi(1-z)}{2}+\ga-\dfrac{p(n)}{q(n)}\sim(-1)^{n+1}\dfrac{4\pi\sin(\pi z)^2}{(1+\sqrt{2})^{4n+4}}\;.$$
$$A=1+((4z^2-15/16)d)/n+(16z^4-(4d+15/2)z^2+15d/16+225/256)/n^2+\cdots$$
\end{cf}

After dividing by $z^2$ and taking $z=0$, this gives a CF for $\z(3)$
which contracts to Ap\'ery's CF \ref{5.5.22}.

Also, note that
$$\dfrac{\psi(1+z)+\psi(1-z)}{2}=\psi(1+z)-\dfrac{1}{2z}+\dfrac{\pi}{2}\cotan(\pi z)\;.$$

\smallskip

\begin{cf}\label{6.9.8}{\ }
\begin{verbatim}
[(z)->(psi(1+z/Pi)-psi(1-z/Pi))/Pi,[0,2*n+1],[z,-z^2]]
\end{verbatim}
$$\dfrac{\psi(1+z/\pi)-\psi(1-z/\pi)}{\pi}=\dfrac{z}{3-\dfrac{z^2}{5-\dfrac{z^2}{7-\dfrac{z^2}{9-\dfrac{z^2}{11-\dfrac{z^2}{13-\ddots}}}}}}\;.$$
Convergence type $F^2$ with $E=4/z^2$, $P=2$, and $C=\pi z^3/(8\sin(z)^2)$,
so that
$$\dfrac{\psi(1+z/\pi)-\psi(1-z/\pi)}{\pi}-\dfrac{p(n)}{q(n)}\sim\dfrac{\pi/\sin(z)^2}{n!^2(2/z)^{2n+3}n^2}\;.$$
$$A=1-(z^2/2+9/4)/n+(z^4/8+15z^2/8+125/32)/n^2+\cdots$$
\end{cf}

Note that
$$\dfrac{\psi(1+z/\pi)-\psi(1-z/\pi)}{\pi}=\dfrac{1}{z}-\cotan(z)\;,$$
so this CF is the same as \ref{6.6.30}.

\smallskip

See Section \ref{sec:beta} for specific CFs for
\begin{align*}\be(z)&=\psi(1+z)-\psi(1+z/2)+1/z-\log(2)=\psi(z)-\psi(z/2)-\log(2)\\
  &=\dfrac{1}{2}(\psi((1+z)/2)-\psi(z/2))=\sum_{k\ge0}\dfrac{(-1)^k}{z+k}\;.\end{align*}

\smallskip

\begin{cf}\label{6.9.9}{\ }
\begin{verbatim}
[(z)->psi(1+z)-psi(1+z/3)+1/z-log(3),
[[0,z^2],[1,(2*n+1)*z^2]],[[2/3,1],(3*n+1)*(3*n+2)/2*[n,n+1]]]
\end{verbatim}
$$\psi(1+z)-\psi(1+z/3)+1/z-\log(3)=\dfrac{2/3}{z^2+\dfrac{1}{1+\dfrac{10}{3z^2+\dfrac{20}{1+\dfrac{56}{5z^2+\dfrac{84}{1+\ddots}}}}}}$$
Convergence type $P^-$ with $P=4z/3$ and $C=...$, so that
$$\psi(1+z)-\psi(1+z/3)+1/z-\log(3)-\dfrac{p(n)}{q(n)}\sim(-1)^n\dfrac{C}{n^{4z/3}}\;.$$
$$A=1-(4z/3)/n-4((4z^5-72z^4-48z^3+162z^2+81z)/(81(4z^2-9)))/n^2+\cdots$$
\end{cf}

\smallskip

\begin{cf}\label{6.9.10}{\ }
\begin{verbatim}
[(z)->psi(1+z)-psi(1+z/3)+1/z-log(3),[0,(2*n-1)*(9*n^2-9*n+2*z^2+2)],
                                   [4/3,-n^2*(9*n^2-1)*(9*n^2-4)]]
\end{verbatim}
$$\psi(1+z)-\psi(1+z/3)+1/z-\log(3)=\dfrac{4/3}{2z^2+2-\dfrac{40}{6z^2+60-\dfrac{4480}{10z^2+280-\dfrac{55440}{14z^2+770-\ddots}}}}$$
Convergence type $P^+$ with $P=4z/3$ and $C=...$, so that
$$\psi(1+z)-\psi(1+z/3)+1/z-\log(3)-\dfrac{p(n)}{q(n)}\sim\dfrac{C}{n^{4z/3}}\;.$$
$$A=1-(2z/3)/n-((4z^5-72z^4-48z^3+162z^2+81z)/(81(4z^2-9)))/n^2+\cdots$$
Parametric family for $k\ge0$:
\begin{verbatim}
[(z)->psi(1+z)-psi(1+z/3)+1/z-log(3),
(2*n-1)*(9*n^2-9*n+2+2*(3*k+z)^2),-n^2*(9*n^2-1)*(9*n^2-4)]
\end{verbatim}
Convergence type $P^+$ with $P=4k+4z/3$.
\end{cf}

This is simply the contraction of the previous one, but still nice.

It can be Ap\'ery accelerated with convergence type $E=-(1+\sqrt{2})^4$
and $P=0$, formula too complicated to give here.

\smallskip

\begin{cf}\label{6.9.11}{\ }
\begin{verbatim}
[(z)->psi(1+z)-psi(1+z/3),[0,6*n^2+(5*z-6)*n+(3-2*z)],
                           [2*z,-3*n^2*(n+z)*(3*n+2*z)]]
\end{verbatim}
$$\psi(1+z)-\psi(1+z/3)=\dfrac{2z}{3z+3-\dfrac{6z^2+15z+9}{8z+15-\dfrac{24z^2+120z+144}{13z+39-\dfrac{54z^2+405z+729}{18z+75-\ddots}}}}$$
Convergence type $P^+$ with $P=z/3+1$ and $C=...$, so that
$$\psi(1+z)-\psi(1+z/3)-\dfrac{p(n)}{q(n)}\sim\dfrac{C}{n^{z/3+1}}$$
$$A=1-((z+3)(5z^2+18z+18)/(18(z+6)))/n+\cdots$$
Series:
$$\psi(1+z)-\psi(1+z/3)=\dfrac{2z}{3(z+1)}\sum_{n\ge0}\dfrac{(1+2z/3)_n}{(n+1)(z+2)_n}$$
Parametric family for $k\ge0$:
\begin{verbatim}
[(z)->psi(1+z)-psi(1+z/3),6*n^2+(5*z-6)*n+(k-2)*z+3*(k^2+k+1),
                          -3*n^2*(n+z)*(3*n+2*z)]
\end{verbatim}
Convergence type $P^+$ with $P=z/3+2k+1$.
\end{cf}

\smallskip

\begin{cf}\label{6.9.12}{\ }
\begin{verbatim}
[(z)->psi(1+z)-psi(1+z/3),[0,(4*z+3)*(2*n-1)],[4*z,n^2*(9*n^2-4*z^2)]]
\end{verbatim}
$$\psi(1+z)-\psi(1+z/3)=\dfrac{4z}{4z+3-\dfrac{4z^2-9}{12z+9-\dfrac{16z^2-144}{20z+15-\dfrac{36z^2-729}{28z+21-\ddots}}}}$$
Convergence type $P^-$ with $P=8z/3+2$ and $C=...$, so that
$$\psi(1+z)-\psi(1+z/3)-\dfrac{p(n)}{q(n)}\sim(-1)^n\dfrac{C}{n^{8z/3+2}}\;.$$
$$A=1-(4z/3+1)/n-(2z(z-3)(4z+3)/81)/n^2+\cdots$$
Parametric family for $k\ge0$:
\begin{verbatim}
[(z)->psi(1+z)-psi(1+z/3),(2*n-1)*(4*z+3*(2*k+1)),n^2*(9*n^2-4*z^2)]
\end{verbatim}
Convergence type $P^-$ with $P=4k+8z/3+2$.
\end{cf}

\smallskip

\begin{cf}\label{6.9.13}{\ }
\begin{verbatim}
[(z)->psi(1+z)-psi(1+z/3),
[0,(z+1)*(z+3),6*n^2+(8*z-6)*n+2*z^2-4*z+3],[2*z,-(n+z)^2*(3*n+z)^2]]
\end{verbatim}
$$\psi(1+z)-\psi(1+z/3)=\dfrac{2z}{z^2+4z+3-\dfrac{z^4+8z^3+22z^2+24z+9}{2z^2+12z+15-\dfrac{z^4+16z^3+88z^2+192z+144}{2z^2+20z+39-\ddots}}}$$
Convergence type $P^+$ with $P=1$ and $C=2z/3$, so that
$$\psi(1+z)-\psi(1+z/3)-\dfrac{p(n)}{q(n)}\sim\dfrac{2z/3}{n}\;.$$
$$A=1-(2z/3+1/2)/n+(13z^2/27+2z/3+1/6)/n^2+\cdots$$
Series:
$$\psi(1+z)-\psi(1+z/3)=2z\sum_{n\ge1}\dfrac{1}{(n+z)(3n+z)}$$
Parametric family for $k\ge0$:
\begin{verbatim}
[(z)->psi(1+z)-psi(1+z/3),
6*n^2+(8*z-6)*n+2*z^2-4*z+3*(k^2+k+1),-(n+z)^2*(3*n+z)^2]
\end{verbatim}
Convergence type $P^+$ with $P=2k+1$.
\end{cf}
      
\smallskip

\begin{cf}\label{6.9.14}{\ }
\begin{verbatim}
[(a,b)->psi(b+1)-psi(a+1),[0,2*n^2+(2*a-b-2)*n+b-a+1],
                          [b-a,-n^2*(n+a)*(n+a-b)]]
\end{verbatim}
$$\psi(b+1)-\psi(a+1)=\dfrac{b-a}{a+1-\dfrac{a^2+(-b+2)a-b+1}{3a-b+5-\dfrac{4a^2+(-4b+16)a-8b+16}{5a-2b+13-\ddots}}}$$
Convergence type $P^+$ with $P=b+1$ and $C=...$, so that
$$\psi(b+1)-\psi(a+1)-\dfrac{p(n)}{q(n)}\sim\dfrac{C}{n^{b+1}}\;.$$
$$A=1+((b^2-(2b+2)a-2)(b+1)/(2(b+2)))/n+\cdots$$
Series:
$$\psi(b+1)-\psi(a+1)=\dfrac{b-a}{a+1}\sum_{n\ge0}\dfrac{(a-b+1)_n}{(n+1)(a+2)_n}$$
Parametric family for $k\ge0$:
\begin{verbatim}
[(a,b)->psi(b+1)-psi(a+1),2*n^2+(2*a-b-2)*n+(k+1)*b-a+k^2+k+1,
                          -n^2*(n+a)*(n+a-b)]
\end{verbatim}
Convergence type $P^+$ with $P=2k+b+1$.
\end{cf}

It can be Ap\'ery accelerated with convergence type $E=-((1+\sqrt{5})/2)^5$
and $P=0$, formula too complicated to give here.

\smallskip

\begin{cf}\label{6.9.15}{\ }
\begin{verbatim}
[(a,b)->psi(b+1)-psi(a+1),
[0,(b+1)*(a+1),2*n^2+2*(a+b-1)*n+2*a*b-(a+b)+1],[b-a,-(n+a)^2*(n+b)^2]]
\end{verbatim}
$$\psi(b+1)-\psi(a+1)=\dfrac{b-a}{(b+1)a+b+1-\dfrac{(b^2+2b+1)a^2+(2b^2+4b+2)a+(b^2+2b+1)}{(2b+3)a+3b+5-\ddots}}$$
Convergence type $P^+$ with $P=1$ and $C=b-a$, so that
$$\psi(b+1)-\psi(a+1)-\dfrac{p(n)}{q(n)}\sim\dfrac{b-a}{n}\;.$$
$$A=1-((a+b+1)/2)/n+((2(a^2+ab+b^2)+3(a+b)+1)/6)/n^2+\cdots$$
Series:
$$\psi(b+1)-\psi(a+1)=(b-a)\sum_{n\ge1}\dfrac{1}{(n+a)(n+b)}$$
Parametric family for $k\ge0$:
\begin{verbatim}
[(a,b)->psi(b+1)-psi(a+1),2*n^2+2*(a+b-1)*n+2*a*b-(a+b)+k^2+k+1,
                          -(n+a)^2*(n+b)^2]
\end{verbatim}
Convergence type $P^+$ with $P=2k+1$.
\end{cf}

\smallskip

\begin{cf}\label{6.9.16}{\ }
\begin{verbatim}
[(a,b)->psi(b+1)-psi(a+1),[0,(a+b+1)*(2*n-1)],
                          [2*(b-a),n^2*(n^2-(b-a)^2)]]
\end{verbatim}
$$\psi(b+1)-\psi(a+1)=\dfrac{2b-2a}{b+a+1-\dfrac{b^2-2ab+a^2-1}{3b+3a+3-\dfrac{4b^2-8ab+4a^2-16}{5b+5a+5-\ddots}}}$$
Convergence type $P^-$ with $P=2(a+b+1)$ and $C=...$, so that
$$\psi(b+1)-\psi(a+1)-\dfrac{p(n)}{q(n)}\sim(-1)^n\dfrac{C}{n^{2(a+b+1)}}\;.$$
$$A=1-(a+b+1)/n+(a^3/3-(b-1/2)a^2-(b^2+b-1/6)a+b^3/3+b^2/2+b/6)/n^2+\cdots$$
Parametric family for $k\ge0$:
\begin{verbatim}
[(a,b)->psi(b+1)-psi(a+1),(a+b+2*k+1)*(2*n-1),n^2*(n^2-(b-a)^2)]
\end{verbatim}
Convergence type $P^-$ with $P=2(a+b+2k+1)$.
\end{cf}

It can be Ap\'ery accelerated with convergence type $E=-((1+\sqrt{5})/2)^5$
and $P=0$, formula too complicated to give here, very similar to \ref{6.9.17}.

\smallskip

\begin{cf}\label{6.9.17}{\ }
\begin{verbatim}
[(a,b)->psi(b+1)-psi(a+1),
[[0,a+b+1],[5*n^2+3*(a+b)*n+2*a*b,5*n^2+(3*(a+b)+4)*n+2*a*b+a+b+1]],
[[2*(b-a),(a-b+1)*(b-a+1)],
           [-4*(n+a)^2*(n+b)^2,(n+1)^2*(n+a-b+1)*(n+b-a+1)]]]
\end{verbatim}
$$\psi(b+1)-\psi(a+1)=\dfrac{2b-2a}{a+b+1-\dfrac{a^2-2ba+b^2-1}{(2b+3)a+3b+5-\dfrac{4(a+1)^2(b+1)^2}{(2b+4)a+4b+10-\ddots}}}$$
Convergence type $E$ with $E=-((1+\sqrt{5})/2)^5$, $P=0$, and $C=...$, so that
$$\psi(b+1)-\psi(a+1)-\dfrac{p(n)}{q(n)}\sim(-1)^n\dfrac{C}{((1+\sqrt{5})/2)^{5n}}\;.$$
\end{cf}

\smallskip

\begin{cf}\label{6.9.18}{\ }
\begin{verbatim}
[(a,b)->psi(a+1)+psi(b+1)-psi(a+b+1)+Euler(),
[0,(2*n-1)*(n^2-n+(a+1)*(b+1))],[2*a*b,-n^2*(n^2-a^2)*(n^2-b^2)]]
\end{verbatim}
$$\psi(a+1)+\psi(b+1)-\psi(a+b+1)+\gamma=\dfrac{2ba}{(b+1)a+b+1-\dfrac{(b^2-1)a^2-b^2+1}{(3b+3)a+3b+9-\ddots}}$$
Convergence type $P^+$ with $P=2(a+b+1)$ and $C=...$, so that
$$\psi(a+1)+\psi(b+1)-\psi(a+b+1)+\gamma-\dfrac{p(n)}{q(n)}\sim\dfrac{C}{n^{2(a+b+1)}}$$
$$A=1-(a+b+1)/n+\cdots$$
Series:
$$\psi(a+1)+\psi(b+1)-\psi(a+b+1)+\gamma=2\sum_{n\ge1}\dfrac{(-a)_n(-b)_n}{n(a+1)_n(b+1)_n}$$
Parametric family for $k\ge0$:
\begin{verbatim}
[(a,b)->psi(a+1)+psi(b+1)-psi(a+b+1)+Euler(),
(2*n-1)*(n^2-n+(a+1)*(b+1)+2*k*(k+a+b+1)),-n^2*(n^2-a^2)*(n^2-b^2)]
\end{verbatim}
Convergence type $P^+$ with $P=2(a+b+2k+1)$.
\end{cf}

\smallskip

\begin{cf}\label{6.9.19}{\ }
\begin{verbatim}
[(a,b)->psi((a+b+3)/4)+psi((a-b+3)/4)-psi((a+b+1)/4)-psi((a-b+1)/4),
[0,a],[[4,1-b^2],[(2*n)^2,(2*n+1)^2-b^2]]]
\end{verbatim}
\begin{align*}\psi\left(\dfrac{a+b+3}{4}\right)+\psi\left(\dfrac{a-b+3}{4}\right)&-\psi\left(\dfrac{a+b+1}{4}\right)-\psi\left(\dfrac{a-b+1}{4}\right)=\\
  &=\dfrac{4}{a+\dfrac{1-b^2}{a+\dfrac{4}{a+\dfrac{9-b^2}{a+\dfrac{16}{a+\dfrac{25-b^2}{a+\ddots}}}}}}\end{align*}
Convergence type $P^-$ with $P=a$ and $C=...$, so that
$$\psi\left(\dfrac{a+b+3}{4}\right)+\psi\left(\dfrac{a-b+3}{4}\right)-\psi\left(\dfrac{a+b+1}{4}\right)-\psi\left(\dfrac{a-b+1}{4}\right)-\dfrac{p(n)}{q(n)}\sim(-1)^n\dfrac{C}{n^a}\;.$$
$$A=1-(a(a^2-b^2-1)/(2(a^2-1)))/n+\cdots$$
\end{cf}

\begin{cf}\label{6.9.20}{\ }
\begin{verbatim}
[(a,b)->psi((a+b+3)/4)+psi((a-b+3)/4)-psi((a+b+1)/4)-psi((a-b+1)/4),
[0,8*n^2-12*n+5+(a^2-b^2)],[4*a,4*n^2*(b^2-(2*n-1)^2)]]
\end{verbatim}
\begin{align*}&\psi\left(\dfrac{a+b+3}{4}\right)+\psi\left(\dfrac{a-b+3}{4}\right)-\psi\left(\dfrac{a+b+1}{4}\right)-\psi\left(\dfrac{a-b+1}{4}\right)=\\
&=\dfrac{4a}{a^2+(-b^2+1)+\dfrac{4b^2-4}{a^2+(-b^2+13)+\dfrac{16b^2-144}{a^2+(-b^2+41)+\dfrac{36b^2-900}{a^2+(-b^2+85)+\ddots}}}}\end{align*}
Convergence type $P^+$ with $P=a$ and $C=...$, so that
$$\psi\left(\dfrac{a+b+3}{4}\right)+\psi\left(\dfrac{a-b+3}{4}\right)-\psi\left(\dfrac{a+b+1}{4}\right)-\psi\left(\dfrac{a-b+1}{4}\right)-\dfrac{p(n)}{q(n)}\sim\dfrac{C}{n^a}\;.$$
$$A=1-(a(a^2-b^2-1)/(4(a^2-1)))/n+\cdots$$
Parametric family for $k\ge0$:
\begin{verbatim}
[(a,b)->psi((a+b+3)/4)+psi((a-b+3)/4)-psi((a+b+1)/4)-psi((a-b+1)/4),
a^2-b^2+8*n^2-12*n+5+4*k^2+4*k*a,4*n^2*(b^2-(2*n-1)^2)]
\end{verbatim}
Convergence type $P^+$ with $P=2k+a$.
\end{cf}

\smallskip

This is simply the contraction of the previous CF.

It can be Ap\'ery accelerated with convergence type $E=-((1+\sqrt{5})/2)^5$
and $P=0$, formula too complicated to give here.

\medskip

\section{Function $\be(z)=(\psi((1+z)/2)-\psi(z/2))/2$}\label{sec:beta}

\medskip

Since this function has so many incarnations and applications, it deserves a
separate section. Note the following formulas:
\begin{align*}
  \be(z)&=(\psi((1+z)/2)-\psi(z/2))/2=\psi(z)-\psi(z/2)-\log(2)\\
  &=\psi(1+z)-\psi(1+z/2)+1/z-\log(2)\,\\
  \be(z)&=\sum_{k\ge0}\dfrac{(-1)^k}{k+z}\;,\\
  \int_0^\infty\dfrac{e^{-zt}}{\cosh(t)}\,dt&=\be((z+1)/2)\;,\quad
  \int_0^\infty\dfrac{e^{-zt}}{\cosh^2(t)}\,dt=z\be(z/2)-1\;,\\
  \int_0^\infty e^{-zt}\tanh(t)\,dt&=\be(z/2)-1/z\;,\quad
  \int_0^\infty e^{-zt}\tanh^2(t)\,dt=1+1/z-z\be(z/2)\;.\end{align*}

Thus all CFs for $\be(z)$ give CFs for the integrals.

\smallskip

\begin{verbatim}
beta(z)=(psi((1+z)/2)-psi(z/2))/2;
\end{verbatim}

\smallskip

\begin{cf}\label{6.10.1}{\ }
\begin{verbatim}
[(z)->beta(z),[0,z,1],[1,(n+z-1)^2]]
\end{verbatim}
$$\be(z)=\dfrac{1}{z+\dfrac{z^2}{1+\dfrac{z^2+2z+1}{1+\dfrac{z^2+4z+4}{1+\dfrac{z^2+6z+9}{1+\dfrac{z^2+8z+16}{1+\ddots}}}}}}$$
Convergence type $P^-$ with $P=1$ and $C=1/2$, so that
$$\be(z)-\dfrac{p(n)}{q(n)}\sim(-1)^n\dfrac{1/2}{n}\;.$$
$$A=1+(1/2-z)/n+(z^2-z)/n^2-(z^3-3z^2/2+1/4)/n^3+\cdots$$
Series:
$$\beta(z)=\sum_{n\ge1}\dfrac{(-1)^{n+1}}{n+z-1}$$
Parametric family for $k\ge0$:
\begin{verbatim}
[(z)->beta(z),2*k+1,(n+z-1)^2]
\end{verbatim}
Convergence type $P^-$ with $P=2k+1$.
\end{cf}

\smallskip

\begin{cf}\label{6.10.2}{\ }
\begin{verbatim}
[(z)->beta(z),[0,z*(z+1),8*n^2+(8*z-20)*n+2*z^2-10*z+14],
            [1,-(2*n-1+z)^2*(2*n-2+z)^2]]
\end{verbatim}
$$\be(z)=\dfrac{1}{z^2+z-\dfrac{z^4+2z^3+z^2}{2z^2+6z+6-\dfrac{z^4+10z^3+37z^2+60z+36}{2z^2+14z+26-\ddots}}}$$
Convergence type $P^+$ with $P=1$ and $C=1/4$, so that
$$\be(z)-\dfrac{p(n)}{q(n)}\sim\dfrac{1/4}{n}\;.$$
$$A=1+(1/4-z/2)/n+(z^2/4-z/4)/n^2-(z^3/8-3z^2/16+1/32)/n^3+\cdots$$
Series:
$$\beta(z)=\sum_{n\ge1}\dfrac{1}{(2n+z-1)(2n+z-2)}$$
Parametric family for $k\ge0$:
\begin{verbatim}
[(z)->beta(z),8*n^2+(8*z-20)*n+2*z^2-10*z+14+4*k^2+4*k,
            -(2*n-1+z)^2*(2*n-2+z)^2]
\end{verbatim}
Convergence type $P^+$ with $P=2k+1$.
\end{cf}

This is the contraction of the previous CF.

\smallskip

\begin{cf}\label{6.10.3}{\ }
\begin{verbatim}
[(z)->beta(z),[1/z,z+1,1],[-1,(n+z)^2]]
\end{verbatim}
$$\be(z)=\dfrac{1}{z}-\dfrac{1}{z+1+\dfrac{z^2+2z+1}{1+\dfrac{z^2+4z+4}{1+\dfrac{z^2+6z+9}{1+\dfrac{z^2+8z+16}{1+\dfrac{z^2+10z+25}{1+\ddots}}}}}}$$
Convergence type $P^-$ with $P=1$ and $C=...$, so that
$$\be(z)-\dfrac{p(n)}{q(n)}\sim(-1)^n\dfrac{C}{n}\;.$$
$$A=1-(z+1/2)/n+(z^2+z)/n^2-(z^3+3z^2/2-1/4)/n^3+\cdots$$
Series:
$$\be(z)=\sum_{n\ge0}\dfrac{(-1)^n}{n+z}$$
Parametric family for $k\ge0$:
\begin{verbatim}
[(z)->beta(z),2*k+1,(n+z)^2]
\end{verbatim}
Convergence type $P^-$ with $P=2k+1$.
\end{cf}

\smallskip

\begin{cf}\label{6.10.4}{\ }
\begin{verbatim}
[(z)->beta(z),[0,(2*n-1)^2+z*n],[1/2,-n^2*(2*n+1)*(2*n+z-1)]]
\end{verbatim}
$$\be(z)=\dfrac{1/2}{z+1-\dfrac{3z+3}{2z+9-\dfrac{20z+60}{3z+25-\dfrac{63z+315}{4z+49-\dfrac{144z+1008}{5z+81-\ddots}}}}}$$
Convergence type $P^+$ with $P=z/2$ and $C=...$, so that
$$\be(z)-\dfrac{p(n)}{q(n)}\sim\dfrac{C}{n^{z/2}}\;.$$
$$A=1-((z^3+2z^2+4z)/(8(z+2)))/n+\cdots$$
Series:
$$\be(z)=\dfrac{1}{2z+2}\sum_{n\ge0}\dfrac{(3/2)_n}{(n+1)((z+3)/2)_n}$$
Parametric family for $k\ge0$:
\begin{verbatim}
[(z)->beta(z),(2*n-1)^2+z*(n+k)+2*k^2,-n^2*(2*n+1)*(2*n+z-1)]
\end{verbatim}
Convergence type $P^-$ with $P=2k+z/2$.
\end{cf}

\smallskip

\begin{cf}\label{6.10.5}{\ }
\begin{verbatim}
[(z)->beta(z),[0,2*z-1],[1,n^2]]
\end{verbatim}
$$\be(z)=\dfrac{1}{2z-1+\dfrac{1}{2z-1+\dfrac{4}{2z-1+\dfrac{9}{2z-1+\dfrac{16}{2z-1+\dfrac{25}{2z-1+\ddots}}}}}}$$
Convergence type $P^-$ with $P=2z-1$ and $C=\G(z)^2/2^{2z-1}$, so that
$$\be(z)-\dfrac{p(n)}{q(n)}\sim(-1)^n\dfrac{\G(z)^2/2^{2z-1}}{n^{2z-1}}\;.$$
$$A=1+(z-1/2)/n+((z(z+2)(2z-1))/12)/n^2+\cdots$$
Parametric family up to trivial change $n\mapsto n+j$:
\begin{verbatim}
[(z)->beta(z),2*z+u+2*k-1,n*(n+u)]
\end{verbatim}
Convergence type $P^-$ with $P=2z+2k+u-1$.
\end{cf}

\smallskip

When $z\in\Q$, an equivalent way of writing this CF is as follows:

\smallskip

\begin{cf}\label{6.10.6}{\ }
\begin{verbatim}
[(k,m)->beta((k+m)/(2*k)),[0,m],[k,k^2*n^2]]
\end{verbatim}
$$\be\left(\dfrac{m+k}{2k}\right)=\dfrac{k}{m+\dfrac{k^2}{m+\dfrac{4k^2}{m+\dfrac{9k^2}{m+\dfrac{16k^2}{m+\dfrac{25k^2}{m+\ddots}}}}}}$$
Convergence type $P^-$ with $P=m/k$ and $C=\G((k+m)/(2k))^2/2^{m/k}$, so that
$$\be\left(\dfrac{m+k}{2k}\right)-\dfrac{p(n)}{q(n)}\sim(-1)^n\dfrac{\G((k+m)/(2k))^2/2^{m/k}}{n^{m/k}}\;.$$
$$A=1-(m/(2k))/n-(m(m-k)(m-5k)/(48k^3))/n^2+\cdots$$
\end{cf}

\smallskip

Using explicit formulas for $\be(p/q)$ we can specialize to a few values of
$(k,m)$ with $(k,m)=1$ to obtain many of the given CFs for linear combinations
of zeta and $L$-values, and of course infinitely many others, for instance:

\smallskip

\begin{verbatim}
[()->log(2),[0,1],[1,n^2]]
[()->Pi/2,[2],[-1,n^2]]
[()->Pi/2,[4/3,4],[1,n^2]]
[()->Pi/sqrt(2)+sqrt(2)*log(1+sqrt(2)),[4,3],[-2,4*n^2]]
[()->Pi/sqrt(2)-sqrt(2)*log(1+sqrt(2)),[0,1],[2,4*n^2]]
[()->Pi/sqrt(3)+log(2),[3,5],[-3,9*n^2]]
[()->Pi/sqrt(3)-log(2),[0,1],[3,9*n^2]]
\end{verbatim}

\smallskip

The same game can of course be played with all the other CFs for $\be(z)$,
but this one is the simplest.

\smallskip

\begin{cf}\label{6.10.7}{\ }
\begin{verbatim}
[(z)->beta(z),[0,4*n^2-6*n+2*z^2-2*z+3],[z-1/2,-n^2*(2*n-1)^2]]
\end{verbatim}
$$\be(z)=\dfrac{z-1/2}{2z^2-2z+1-\dfrac{1}{2z^2-2z+7-\dfrac{36}{2z^2-2z+21-\dfrac{225}{2z^2-2z+43-\ddots}}}}$$
Convergence type $P^+$ with $P=2z-1$ and $C=...$, so that
$$\be(z)-\dfrac{p(n)}{q(n)}\sim\dfrac{C}{n^{2z-1}}\;.$$
$$A=1-(z/2-1/4)/n-((z-3)(z-1)(2z-1)/48)/n^2+\cdots$$
Parametric family for $k\ge0$:
\begin{verbatim}
[(z)->beta(z),4*n^2-6*n+2*z^2+(4*k-2)*z+2*k^2-2*k+3,-n^2*(2*n-1)^2]
\end{verbatim}
Convergence type $P^+$ with $P=2k+2z-1$.
\end{cf}

This is the contraction of \ref{6.10.5}.

\smallskip

Special case of the above:

\smallskip

\begin{cf}\label{6.10.8}{\ }
\begin{verbatim}
[(z)->intnum(t=0,[oo,1+z],exp(-t*z)/cosh(t)),
[0,8*n^2-12*n+5+z^2],[z,-4*n^2*(2*n-1)^2]]
\end{verbatim}
$$\int_0^\infty\dfrac{e^{-tz}}{\cosh(t)}\,dt=\dfrac{z}{z^2+1-\dfrac{4}{z^2+13-\dfrac{144}{z^2+41-\dfrac{900}{z^2+85-\dfrac{3136}{z^2+145-\dfrac{8100}{z^2+221-\ddots}}}}}}$$
Convergence type $P^+$ with $P=z$ and $C=...$, so that
$$\int_0^\infty\dfrac{e^{-tz}}{\cosh(t)}\,dt-\dfrac{p(n)}{q(n)}\sim\dfrac{C}{n^z}\;.$$
$$A=1-(z/4)/n+(-z^3/192+z^2/32-5z/192)/n^2+\cdots$$
Parametric family for $k\ge0$:
\begin{verbatim}
[(z)->intnum(t=0,[oo,1+z],exp(-t*z)/cosh(t)),
8*n^2-12*n+5+(2*k+z)^2,-4*n^2*(2*n-1)^2]
\end{verbatim}
Convergence type $P^+$ with $P=2k+z$.
\end{cf}

\smallskip

\begin{cf}\label{6.10.9}{\ }
\begin{verbatim}
[(z)->intnum(t=0,[oo,1+z],exp(-t*z)/cosh(t)),[0,z+1,2],[2,(2*n+z-1)^2]]
\end{verbatim}
$$\int_0^\infty\dfrac{e^{-tz}}{\cosh(t)}\,dt=\dfrac{2}{z+1+\dfrac{z^2+2z+1}{2+\dfrac{z^2+6z+9}{2+\dfrac{z^2+10z+25}{2+\dfrac{z^2+14z+49}{2+\dfrac{z^2+18z+81}{2+\ddots}}}}}}$$
Convergence type $P^-$ with $P=1$, and $C=...$, so that
$$\int_0^\infty\dfrac{e^{-tz}}{\cosh(t)}\,dt-\dfrac{p(n)}{q(n)}\sim(-1)^n\dfrac{C}{n}\;.$$
$$A=1-(z/2)/n+((z^2-1)/4)/n^2-((z^3-3)/8)/n^3+\cdots$$
Series:
$$\int_0^\infty\dfrac{e^{-tz}}{\cosh(t)}\,dt=2\sum_{n\ge1}\dfrac{(-1)^{n+1}}{2n+z-1}$$
Parametric family for $k\ge0$:
\begin{verbatim}
[(z)->intnum(t=0,[oo,1+z],exp(-t*z)/cosh(t)),4*k+2,(2*n+z-1)^2]
\end{verbatim}
Convergence type $P^-$ with $P=2k+1$.
\end{cf}

\smallskip

\begin{cf}\label{6.10.10}{\ }
\begin{verbatim}
[(z)->intnum(t=0,[oo,1+z],exp(-t*z)/cosh(t)),
[0,z,2*n+z-1],[[1,2],[(2*n+z-1)^2,4*(n+1)^2]]]
\end{verbatim}
$$\int_0^\infty\dfrac{e^{-tz}}{\cosh(t)}\,dt=\dfrac{1}{z+\dfrac{2}{z+3+\dfrac{z^2+2z+1}{z+5+\dfrac{16}{z+7+\dfrac{z^2+6z+9}{z+9+\dfrac{36}{z+11+\ddots}}}}}}$$
Convergence type $E$ with $E=-(1+\sqrt{2})^2$, $P=0$, and $C=...$, so that
$$\int_0^\infty\dfrac{e^{-tz}}{\cosh(t)}\,dt-\dfrac{p(n)}{q(n)}\sim(-1)^n\dfrac{C}{(1+\sqrt{2})^{2n}}\;.$$
$$A=1+((2z^2-4z-1)d/8)/n+(-(z+1)(2z^2-4z-1)d/16+(2z^2-4z-1)^2/64)/n^2+\cdots$$
\end{cf}
    
\smallskip

\begin{cf}\label{6.10.11}{\ }
\begin{verbatim}
[(z)->beta(z),[1/(2*z),2*z],[1/(2*z),n*(n+1)]]
\end{verbatim}
$$\be(z)=\dfrac{1}{2z}+\dfrac{1/(2z)}{2z+\dfrac{2}{2z+\dfrac{6}{2z+\dfrac{12}{2z+\dfrac{20}{2z+\dfrac{30}{2z+\ddots}}}}}}$$
Convergence type $P^-$ with $P=2z$ and $C=...$, so that
$$\be(z)-\dfrac{p(n)}{q(n)}\sim(-1)^n\dfrac{C}{n^{2z}}\;.$$
$$A=1-2z/n+(-z^3/6+2z^2+2z/3)/n^2+(z^4/3-z^3-4z^2/3)/n^3+\cdots$$
Parametric family for $k\ge0$:
\begin{verbatim}
[(z)->beta(z),2*k+2*z,n*(n+1)]
\end{verbatim}
Convergence type $P^-$ with $P=2k+2z$.
\end{cf}

This is the case $k=0$, $u=1$ of the parametric family of \ref{6.10.5}.

\smallskip

\begin{cf}\label{6.10.12}{\ }
\begin{verbatim}
[(z)->beta(z),[1/(2*z),(2*n-1)^2+2*z^2],[1/2,-n^2*(4*n^2-1)]]
\end{verbatim}
$$\be(z)=\dfrac{1}{2z}+\dfrac{1/2}{2z^2+1-\dfrac{3}{2z^2+9-\dfrac{60}{2z^2+25-\dfrac{315}{2z^2+49-\dfrac{1008}{2z^2+81-\ddots}}}}}$$
Convergence type $P^+$ with $P=2z$ and $C=...$, so that
$$\be(z)-\dfrac{p(n)}{q(n)}\sim\dfrac{C}{n^{2z}}\;.$$
$$A=1-z/n+(-z^3/24+z^2/2+z/6)/n^2+(z^4/24-z^3/8-z^2/6)/n^3+\cdots$$
Parametric family for $k\ge0$:
\begin{verbatim}
[(z)->beta(z),(2*n-1)^2+2*(z+k)^2,-n^2*(4*n^2-1)]
\end{verbatim}
Convergence type $P^+$ with $P=2z+2k$.
\end{cf}

This is the contraction of \ref{6.10.11}.

\smallskip

\begin{cf}\label{6.10.13}{\ }
\begin{verbatim}
[(z)->beta(z),[(2*z+1)/(4*z^2),4*n^2+2*z^2],
            [-1/(2*z)^2,-n*(n+1)*(2*n+1)^2]]
\end{verbatim}
$$\be(z)=\dfrac{2z+1}{4z^2}-\dfrac{1/(4z^2)}{2z^2+4-\dfrac{18}{2z^2+16-\dfrac{150}{2z^2+36-\dfrac{588}{2z^2+64-\dfrac{1620}{2z^2+100-\ddots}}}}}$$
Convergence type $P^+$ with $P=2z$ and $C=...$, so that
$$\be(z)-\dfrac{p(n)}{q(n)}\sim\dfrac{C}{n^{2z}}\;.$$
$$A=1-2z/n+(-z^3/24+2z^2+11z/12)/n^2+(z^4/12-5z^3/4-11z^2/6-z/2)/n^3+\cdots$$
Parametric family for $k\ge0$:
\begin{verbatim}
[(z)->beta(z),4*n^2+2*(z+k)^2,-n*(n+1)*(2*n+1)^2]
\end{verbatim}
Convergence type $P^+$ with $P=2k+2z$.
\end{cf}

\smallskip

\begin{cf}\label{6.10.14}{\ }
\begin{verbatim}
[(z)->beta(z),[1/(2*z),z+1,2],[1/(2*z),(n+z-1)*(n+z)]]
\end{verbatim}
$$\be(z)=1/(2z)+\dfrac{1/(2z)}{z+1+\dfrac{z^2+z}{2+\dfrac{z^2+3z+2}{2+\dfrac{z^2+5z+6}{2+\dfrac{z^2+7z+12}{2+\dfrac{z^2+9z+20}{2+\ddots}}}}}}$$
Convergence type $P^-$ with $P=2$ and $C=1/4$, so that
$$\be(z)-\dfrac{p(n)}{q(n)}\sim(-1)^n\dfrac{1/4}{n^2}\;.$$
$$A=1-(2z)/n+((6z^2-1)/2)n^2-(4z^3-2z)/n^3+\cdots$$
Series:
$$\be(z)=\dfrac{1}{2z}+\dfrac{1}{2}\sum_{n\ge1}\dfrac{(-1)^{n+1}}{(n+z-1)(n+z)}$$
Parametric family for $k\ge0$:
\begin{verbatim}
[(z)->beta(z),2*k+2,(n+z-1)*(n+z)]
\end{verbatim}
Convergence type $P^-$ with $P=2k+2$.
\end{cf}

\smallskip

\begin{cf}\label{6.10.15}{\ }
\begin{verbatim}
[(z)->beta(z),[1/(2*z),2*z+2,3*n+2*z],[1/(2*z),-2*(n+1)*(n+z)]]
\end{verbatim}
$$\be(z)=\dfrac{1}{2z}+\dfrac{1/(2z)}{2z+2-\dfrac{4z+4}{2z+6-\dfrac{6z+12}{2z+9-\dfrac{8z+24}{2z+12-\dfrac{10z+40}{2z+15-\dfrac{12z+60}{2z+18-\ddots}}}}}}$$
Convergence type $E$ with $E=2$, $P=z$, and $C=...$, so that
$$\be(z)-\dfrac{p(n)}{q(n)}\sim\dfrac{C}{2^nn^z}\;.$$
$$A=1-(z^2/2+5z/2)/n+(z^4/8+17z^3/12+43z^2/8+49z/12)/n^2+\cdots$$
Series:
$$\be(z)=\dfrac{1}{2z}+\dfrac{1}{4z(z+1)}\sum_{n\ge0}\dfrac{(n+1)!}{(z+2)_n}2^{-n}$$
Parametric family for $k\ge0$:
\begin{verbatim}
[(z)->beta(z),3*n+2*z+k,-2*(n+1)*(n+z)]
\end{verbatim}
Convergence type $E$ with $E=2$ and $P=2k+z$.
\end{cf}

\smallskip

\begin{cf}\label{6.10.16}{\ }
\begin{verbatim}
[(z)->beta(z),[1/(2*z),3*n+2*z-2],[1/(2*z),-2*n*(n+z-1)]]
\end{verbatim}
$$\be(z)=\dfrac{1}{2z}+\dfrac{1/(2z)}{2z+1-\dfrac{2z}{2z+4-\dfrac{4z+4}{2z+7-\dfrac{6z+12}{2z+10-\dfrac{8z+24}{2z+13-\ddots}}}}}$$
Convergence type $E$ with $E=2$, $P=z+2$, and $C=2\G(z+1)$, so that
$$\be(z)-\dfrac{p(n)}{q(n)}\sim\dfrac{2\G(z+1)}{2^nn^{z+2}}\;.$$
$$A=1+(-z^2/2-11z/2-3)/n+(z^4/8+35z^3/12+179z^2/8+343z/12+13)/n^2+\cdots$$
Series:
$$\be(z)=\dfrac{1}{2z}+\sum_{n\ge0}\dfrac{n!}{(n+2z)(n+2z+1)(2z)_n}2^{-n}$$
Parametric family for $k\ge0$:
\begin{verbatim}
[(z)->beta(z),3*n+2*z-2+k,-2*n*(n+z-1)]
\end{verbatim}
Convergence type $E$ with $E=2$ and $P=2k+z+2$.
\end{cf}

\smallskip

\begin{cf}\label{6.10.17}{\ }
\begin{verbatim}
[(z)->beta(z),[[0,2*z],[3*n+2*z-2,3*n+2*z]],
              [[1,-2*z],[2*n^2,-2*(2*n+1)*(n+z)]]]
\end{verbatim}
$$\be(z)=\dfrac{1}{2z-\dfrac{2z}{2z+1+\dfrac{2}{2z+3-\dfrac{6z+6}{2z+4+\dfrac{8}{2z+6-\dfrac{10z+20}{2z+7+\ddots}}}}}}$$
Convergence type $E$ with $E=2\sqrt{2}i$, $P=z-1/2$, and
$C=\sqrt{\pi}2^{z-1/2}3^{1-2z}\G(z)$, so that
$$\be(z)-\dfrac{p(n)}{q(n)}\sim(-1)^{\lfloor n/2\rfloor}\dfrac{\sqrt{\pi}3^{1-2z}\G(z)}{2^{(3n+1)/2-z}}$$
$$A=1+(7z^2/9-19z/9+25/36)/n+\cdots$$
\end{cf}

\smallskip

\begin{cf}\label{6.10.18}{\ }
\begin{verbatim}
[(z)->beta(z),[0,n+z-1],[[1,z^2],[n^2,(n+z)^2]]]
\end{verbatim}
$$\be(z)=\dfrac{1}{z+\dfrac{z^2}{z+1+\dfrac{1}{z+2+\dfrac{z^2+2z+1}{z+3+\dfrac{4}{z+4+\dfrac{z^2+4z+4}{z+5+\ddots}}}}}}$$
Convergence type $E$ with $E=-(1+\sqrt{2})^2$, $P=0$, and $C=2\pi/(1+\sqrt{2})^{2z}$, so that
$$\be(z)-\dfrac{p(n)}{q(n)}\sim(-1)^n\dfrac{2\pi}{(1+\sqrt{2})^{2n+2z}}\;.$$
$$A=1+((z^2-2z+5/8)d)/n+(z^4-(d+4)z^3+(2d+21/4)z^2-(5d/8+5/2)z+25/64)/n^2+\cdots$$
\end{cf}

\smallskip

\begin{cf}\label{6.10.19}{\ }
\begin{verbatim}
[(z)->beta(z),[0,2*z-1,n+z-1],[[1,1],[(n+z-1)^2,(n+1)^2]]]
\end{verbatim}
$$\be(z)=\dfrac{1}{2z-1+\dfrac{1}{z+1+\dfrac{z^2}{z+2+\dfrac{4}{z+3+\dfrac{z^2+2z+1}{z+4+\dfrac{9}{z+5+\ddots}}}}}}$$
Convergence type $E$ with $E=-(1+\sqrt{2})^2$, $P=0$, and $C=2\pi/(1+\sqrt{2})^{2z}$, so that
$$\be(z)-\dfrac{p(n)}{q(n)}\sim(-1)^n\dfrac{2\pi}{(1+\sqrt{2})^{2n+2z}}\;.$$
$$A=1+((z^2-2z+5/8)d)/n+(z^4-(d+4)z^3+(2d+21/4)z^2-(5d/8+5/2)z+25/64)/n^2+\cdots$$
\end{cf}

\smallskip

\begin{cf}\label{6.10.20}{\ }
\begin{verbatim}
[(z)->beta(z),[1/(2*z),2*z+1,2*(n+z)],
            [[1/(2*z),-2*z],-2*(n+1)*(n+z)*[1,1]]]
\end{verbatim}
$$\be(z)=\dfrac{1}{2z}+\dfrac{1/(2z)}{2z+1-\dfrac{2z}{2z+4-\dfrac{4z+4}{2z+6-\dfrac{4z+4}{2z+8-\dfrac{6z+12}{2z+10-\dfrac{6z+12}{2z+12-\ddots}}}}}}$$
Convergence type $E$ with $E=-(1+\sqrt{2})^2$, $P=0$, and $C=2\pi/(1+\sqrt{2})^{2z+2}$, so that
$$\be(z)-\dfrac{p(n)}{q(n)}\sim(-1)^n\dfrac{2\pi}{(1+\sqrt{2})^{2n+2z+2}}\;.$$
\begin{align*}A&=1+((z^2-2z+9/8)d)/n\\
  &\phantom{=}+(z^4-(d+4)z^3+(d+25/4)z^2+(7d/8-9/2)z-9d/8+81/64)/n^2+\cdots\end{align*}
\end{cf}

\smallskip

\begin{cf}\label{6.10.21}{\ }
\begin{verbatim}
[(z)->beta(z),[1/z,2*z+1,n+z],[[-1,1],[(n+z)^2,(n+1)^2]]]
\end{verbatim}
$$\be(z)=\dfrac{1}{z}-\dfrac{1}{2z+1+\dfrac{1}{z+2+\dfrac{z^2+2z+1}{z+3+\dfrac{4}{z+4+\dfrac{z^2+4z+4}{z+5+\dfrac{9}{z+6+\ddots}}}}}}$$
Convergence type $E$ with $E=-(1+\sqrt{2})^2$, $P=0$, and $C=-2\pi/(1+\sqrt{2})^{2z+2}$, so that
$$\be(z)-\dfrac{p(n)}{q(n)}\sim(-1)^{n+1}\dfrac{2\pi}{(1+\sqrt{2})^{2n+2z+2}}\;.$$
$$A=1+((z^2-3/8)d)/n+(z^4-dz^3-(d+3/4)z^2+(3d/8)z+3d/8+9/64)/n^2+\cdots$$
\end{cf}

\smallskip

\begin{cf}\label{6.10.22}{\ }
\begin{verbatim}
[(z)->beta(z),[1/(2*z),n+z],[[1/(2*z),z^2+z],[n*(n+1),(n+z)*(n+z+1)]]]
\end{verbatim}
$$\be(z)=\dfrac{1}{2z}+\dfrac{1/(2z)}{z+1+\dfrac{z^2+z}{z+2+\dfrac{2}{z+3+\dfrac{z^2+3z+2}{z+4+\dfrac{6}{z+5+\dfrac{z^2+5z+6}{z+6+\ddots}}}}}}$$
Convergence type $E$ with $E=-(1+\sqrt{2})^2$, $P=0$, and $C=2\pi/(1+\sqrt{2})^{2z+2}$, so that
$$\be(z)-\dfrac{p(n)}{q(n)}\sim(-1)^n\dfrac{2\pi}{(1+\sqrt{2})^{2n+2z+2}}\;.$$
\begin{align*}A&=1+((z^2-2z+9/8)d)/n\\
  &\phantom{=}+(z^4-(d+4)z^3+(d+25/4)z^2+(7d/8-9/2)z-9d/8+81/64)/n^2+\cdots\end{align*}
\end{cf}

\smallskip

\begin{cf}\label{6.10.23}{\ }
\begin{verbatim}
[(z)->beta(z),
[[0,2*z-1],[10*n^2+(6*z-9)*n+z^2-3*z+2,10*n^2+(6*z-1)*n+z^2-z+1]],
[[1,3],[-(2*n-1+z)^2*(2*n-2+z)^2,(n+1)^2*(2*n+1)*(2*n+3)]]]
\end{verbatim}
$$\be(z)=\dfrac{1}{2z-1+\dfrac{3}{z^2+3z+3-\dfrac{z^4+2z^3+z^2}{z^2+5z+10+\dfrac{60}{z^2+9z+24-\ddots}}}}$$
Convergence type $E$ with $E=-((1+\sqrt{5})/2)^5$, $P=0$, and $C=2\pi/((1+\sqrt{5})/2)^{2z+2}$, so that
$$\be(z)-\dfrac{p(n)}{q(n)}\sim(-1)^n\dfrac{2\pi}{((1+\sqrt{5})/2)^{5n+2z+2}}\;.$$
$$A=1+((z^2/5-3z/5+3/10)d)/n+\cdots$$
\end{cf}

\smallskip

\begin{cf}\label{6.10.24}{\ }
\begin{verbatim}
[(z)->beta(z),
[[1/(2*z),z^2+z],[5*n^2+(4*z-2)*n+z^2-z,5*n^2+(4*z+2)*n+z^2+z]],
[[1/2,z^2*(z+1)^2],[-n^2*(4*n^2-1),(n+z)^2*(n+z+1)^2]]]
\end{verbatim}
$$\be(z)=\dfrac{1}{2z}+\dfrac{1/2}{z^2+z+\dfrac{z^4+2z^3+z^2}{z^2+3z+3-\dfrac{3}{z^2+5z+7+\dfrac{z^4+6z^3+13z^2+12z+4}{z^2+7z+16-\ddots}}}}$$
Convergence type $E$ with $E=-((1+\sqrt{5})/2)^5$, $P=0$, and $C=2\pi/((1+\sqrt{5})/2)^{6z+2}$, so that
$$\be(z)-\dfrac{p(n)}{q(n)}\sim(-1)^n\dfrac{2\pi}{((1+\sqrt{5})/2)^{5n+6z+2}}\;.$$
$$A=1+((4z^2/5-4z/5+3/10)d)/n+\cdots$$
\end{cf}

\smallskip

\begin{cf}\label{6.10.25}{\ }
\begin{verbatim}
[(z)->beta(z),
[[(2*z+1)/(4*z^2),(z+1)^2],[5*n^2+(4*z+1)*n+z^2,5*n^2+(4*z+5)*n+(z+1)^2]],
[[-1/(2*z)^2,z^2*(z+1)^2],[-n*(n+1)*(2*n+1)^2,(n+z)^2*(n+z+1)^2]]]
\end{verbatim}
$$\be(z)=\dfrac{2z+1}{4z^2}-\dfrac{1/(4z^2)}{z^2+2z+1+\dfrac{z^4+2z^3+z^2}{z^2+4z+6-\dfrac{18}{z^2+6z+11+\dfrac{z^4+6z^3+13z^2+12z+4}{z^2+8z+22-\ddots}}}}$$
Convergence type $E$ with $E=-((1+\sqrt{5})/2)^5$, $P=0$, and $C=-2\pi/((1+\sqrt{5})/2)^{6z+4}$, so that
$$\be(z)-\dfrac{p(n)}{q(n)}\sim(-1)^{n+1}\dfrac{2\pi}{((1+\sqrt{5})/2)^{5n+6z+4}}\;.$$
$$A=1+((4z^2/5-8z/5+9/10)d)/n+\cdots$$
\end{cf}

\smallskip

Since $\be(z)$ is closely related to $\int_0^\infty e^{-zt}/\cosh(t)^k\,dt$
for $k=1$ and $k=2$, it is worth noting the more general CFs given starting
from \ref{6.21.1}.

\smallskip

\begin{cf}\label{6.10.26}{\ }
\begin{verbatim}
[(z)->beta(z)-log(2),[-1,4*n^2-6*n+3*z+1],[z+1,-n^2*(2*n+z-3)*(2*n-z+1)]]
\end{verbatim}
$$\be(z)-\log(2)=-1+\dfrac{z+1}{3z-1+\dfrac{z^2-4z+3}{3z+5+\dfrac{4z^2-16z-20}{3z+19+\dfrac{9z^2-36z-189}{3z+41+\dfrac{16z^2-64z-720}{3z+71+\ddots}}}}}$$
Convergence type $P^+$ with $P=z+1$ and $C=-\cos(\pi z/2)(\G(z)/\G(z/2+1))^2/2^{2z+1}$, so that
$$\be(z)-\log(2)-\dfrac{p(n)}{q(n)}\sim-\dfrac{\cos(\pi z/2)(\G(z)/\G(z/2+1))^2/2^{2z+1}}{n^{z+1}}$$
$$A=1-((3z^2+z-2)/(2z(z+2)))/n+((z^5-6z^4+11z^3-12z+6)/(24z^2))/n^2+\cdots$$
Series:
$$\be(z)-\log(2)=-1+(z^2-1)\sum_{n\ge0}\dfrac{((3-z)/2)_n}{(n+1)((2n+1)z-1)((2n+3)z-1)((z+1)/2)_n}$$
Parametric family for $k\ge0$:
\begin{verbatim}
[(z)->beta(z)-log(2),4*n^2-6*n+(2*k+3)*z+1+2*k*(k+1),-n^2*(2*n+z-3)*(2*n-z+1)]
\end{verbatim}
Convergence type $P^+$ with $P=2k+z+1$.
\end{cf}

\smallskip

\begin{cf}\label{6.10.27}{\ }
\begin{verbatim}
[(z)->beta(z)-log(2),[-1,2*z^2,2*n^2+(2*z-2)*n+(z-1)*(z-2)],
                     [z*(z+2),(n+z-1)^2*(n+1)^2*(2*n+z+2)*(2*n+z-2)]]
\end{verbatim}
$$\be(z)-\log(2)=-1+\dfrac{z^2+2z}{2z^2+\dfrac{4z^4+16z^3}{z^2+z+6+\dfrac{9z^4+90z^3+261z^2+288z+108}{z^2+3z+14+\ddots}}}$$
Convergence type $P^-$ with $P=1$ and $C=1$, so that
$$\be(z)-\log(2)-\dfrac{p(n)}{q(n)}\sim\dfrac{(-1)^n}{n}$$
$$A=1-((z+1)/2)/n+(z^2/2-z/2+1)/n^2-(z^3/2-3z^2/4+5/4)/n^3+\cdots$$
Series:
$$\be(z)-\log(2)=-1+\sum_{n\ge0}(-1)^n\dfrac{2n+z+2}{(n+2)(n+z)}$$
\end{cf}

\smallskip

\begin{cf}\label{6.10.28}{\ }
\begin{verbatim}
[(z)->beta(z)-log(2),[0,2*n+z-2],[-z+1,n^2*(n+z-1)^2]]
\end{verbatim}
$$\be(z)-\log(2)=\dfrac{1-z}{z+\dfrac{z^2}{z+2+\dfrac{4z^2+8z+4}{z+4+\dfrac{9z^2+36z+36}{z+6+\dfrac{16z^2+96z+144}{z+8+\ddots}}}}}$$
Convergence type $P^-$ with $P=2$ and $C=(1-z)/2$, so that
$$\be(z)-\log(2)-\dfrac{p(n)}{q(n)}\sim(-1)^n\dfrac{(1-z)/2}{n^2}$$
$$A=1-z/n+(z^2-z/2-1/2)/n^2-(z^3-z^2-z)/n^3+\cdots$$
Series:
$$\be(z)-\log(2)=(1-z)\sum_{n\ge0}\dfrac{(-1)^n}{(n+1)(n+z)}$$
Parametric family for $k\ge0$ and $u$:
\begin{verbatim}
[(z)->beta(z)-log(2),(2*k+1)*(2*n+2*u+z-2),n^2*(n+2*u+z-1)^2]
\end{verbatim}
Convergence type $P^-$ with $P=4k+2$.
\end{cf}

\smallskip

\begin{cf}\label{6.10.29}{\ }
\begin{verbatim}
[(z)->beta(z)-log(2),[0,4*n^2-4*n+z+1],[-z+1,-n^2*(4*n^2-(z-1)^2)]]
\end{verbatim}
$$\be(z)-\log(2)=\dfrac{1-z}{z+1+\dfrac{z^2-2z-3}{z+9+\dfrac{4z^2-8z-60}{z+25+\dfrac{9z^2-18z-315}{z+49+\dfrac{16z^2-32z-1008}{z+81+\ddots}}}}}$$
Convergence type $P^+$ with $P=z$ and $C=...$, so that
$$\be(z)-\log(2)-\dfrac{p(n)}{q(n)}\sim\dfrac{C}{n^z}$$
$$A=1-(z/2)/n+((z^4+2z^3+5z^2+4z)/(24(z+2)))/n^2+\cdots$$
Series:
$$\be(z)-\log(2)=\dfrac{1-z}{1+z}\sum_{n\ge0}\dfrac{((3-z)/2)_n}{(n+1)((3+z)/2)_n}$$
Parametric family for $k\ge0$ and $u$:
\begin{verbatim}
[(z)->beta(z)-log(2),4*n^2-4*n+(2*k+1)*(z+2*u)+2*k^2+1,
                      -n^2*(4*n^2-(z+2*u-1)^2)]
\end{verbatim}
Convergence type $P^+$ with $P=z+2u+2k$.
\end{cf}

\smallskip

\begin{cf}\label{6.10.30}{\ }
\begin{verbatim}
[(z)->beta(z)+log(2),[0,4*n^2-6*n+z+2],[z,-n^2*(2*n-2+z)*(2*n-z)]]
\end{verbatim}
$$\be(z)+\log(2)=\dfrac{z}{z+\dfrac{z^2-2z}{z+6+\dfrac{4z^2-8z-32}{z+20+\dfrac{9z^2-18z-216}{z+42+\dfrac{16z^2-32z-768}{z+72+\ddots}}}}}$$
Convergence type $P^+$ with $P=z$ and $C=\sin(\pi z/2)(\G(z)/\G((z+1)/2))^2/2^{2z-1}$, so that
$$\be(z)+\log(2)-\dfrac{p(n)}{q(n)}\sim\dfrac{\sin(\pi z/2)(\G(z)/\G((z+1)/2))^2/2^{2z-1}}{n^z}$$
$$A=1-(z/(2(z+1)))/n+(z^3/24-z^2/8+z/12)/n^2+\cdots$$
Series:
$$\be(z)+\log(2)=\sum_{n\ge0}\dfrac{(1-z/2)_n}{(n+1)(1+z/2)_n}$$
Parametric family for $k\ge0$:
\begin{verbatim}
[(z)->beta(z)+log(2),4*n^2-6*n+(2*k+1)*z+2*k^2+2,-n^2*(2*n-2+z)*(2*n-z)]
\end{verbatim}
Convergence type $P^+$ with $P=2k+z$.
\end{cf}

\smallskip

\begin{cf}\label{6.10.31}{\ }
\begin{verbatim}
[(z)->beta(z)+log(2),[0,2*n^2+(2*z-4)*n+(z-1)*(z-2)],
                     [z^2-1,(n+z-1)^2*n^2*(2*n+z+1)*(2*n+z-3)]]
\end{verbatim}
$$\be(z)+\log(2)=\dfrac{z^2-1}{z^2-z+\dfrac{z^4+2z^3-3z^2}{z^2+z+2+\dfrac{4z^4+32z^3+72z^2+64z+20}{z^2+3z+8+\ddots}}}$$
Convergence type $P^-$ with $P=1$ and $C=1$, so that
$$\be(z)+\log(2)-\dfrac{p(n)}{q(n)}\sim\dfrac{(-1)^n}{n}$$
$$A=1-(z/2)/n+(z^2/2-z/2)/n^2-(z^3/2-3z^2/4)/n^3+\cdots$$
Series:
$$\be(z)+\log(2)=\sum_{n\ge0}(-1)^n\dfrac{2n+z+1}{(n+1)(n+z)}$$
\end{cf}

\smallskip

\begin{cf}\label{6.10.32}{\ }
\begin{verbatim}
[(z)->beta(z)+log(2),[1/z,2*n+z-1],[z,n^2*(n+z)^2]]
\end{verbatim}
$$\be(z)+\log(2)=\dfrac{1}{z}+\dfrac{z}{z+1+\dfrac{z^2+2z+1}{z+3+\dfrac{4z^2+16z+16}{z+5+\dfrac{9z^2+54z+81}{z+7+\dfrac{16z^2+128z+256}{z+9+\ddots}}}}}$$
Convergence type $P^-$ with $P=2$ and $C=z/2$, so that
$$\be(z)+\log(2)-\dfrac{p(n)}{q(n)}\sim(-1)^n\dfrac{z/2}{n^2}$$
$$A=1-(z+1)/n+(z^2+3z/2)/n^2-(z^3+2z^2-1)/n^3+\cdots$$
Series:
$$\be(z)+\log(2)=\dfrac{1}{z}+z\sum_{n\ge0}\dfrac{(-1)^n}{(n+1)(n+z+1)}$$
Parametric family for $k\ge0$ and $u$:
\begin{verbatim}
[(z)->beta(z)+log(2),(2*k+1)*(2*n+2*u+z-1),n^2*(n+2*u+z)^2]
\end{verbatim}
Convergence type $P^-$ with $P=4k+2$.
\end{cf}

\smallskip

\begin{cf}\label{6.10.33}{\ }
\begin{verbatim}
[(z)->beta(z)+log(2),[1/z,4*n^2-4*n+z+2],[z,-n^2*(4*n^2-z^2)]]
\end{verbatim}
$$\be(z)+\log(2)=1/z+\dfrac{z}{z+2+\dfrac{z^2-4}{z+10+\dfrac{4z^2-64}{z+26+\dfrac{9z^2-324}{z+50+\dfrac{16z^2-1024}{z+82+\ddots}}}}}$$
Convergence type $P^+$ with $P=z+1$ and $C=...$, so that
$$\be(z)+\log(2)-\dfrac{p(n)}{q(n)}\sim\dfrac{C}{n^{z+1}}$$
$$A=1-((z+1)/2)/n+((z^4+6z^3+17z^2+24z+12)/(24(z+3)))/n^2+\cdots$$
Series:
$$\be(z)+\log(2)=\dfrac{1}{z}+\dfrac{z}{z+2}\sum_{n\ge0}\dfrac{(1-z/2)_n}{(n+1)(2+z/2)_n}$$
Parametric family for $k\ge0$ and $u$:
\begin{verbatim}
[(z)->beta(z)+log(2),4*n^2-4*n+2+(2*k+1)*(z+2*u)+2*k*(k+1),
                      -n^2*(4*n^2-(z+2*u)^2)]
\end{verbatim}
Convergence type $P^+$ with $P=z+1+2u+2k$.
\end{cf}

\smallskip

\begin{cf}\label{6.10.34}{\ }
\begin{verbatim}
[(z)->beta(z)-Pi/2,[-2,8*n^2-12*n+4*z+3],
                   [2*z+1,-n^2*(4*n-5+2*z)*(4*n+1-2*z)]]
\end{verbatim}
$$\be(z)-\dfrac{\pi}{2}=-2+\dfrac{2z+1}{4z-1+\dfrac{4z^2-12z+5}{4z+11+\dfrac{16z^2-48z-108}{4z+39+\dfrac{36z^2-108z-819}{4z+83+\ddots}}}}$$
Convergence type $P^+$ with $P=z+1/2$ and $C=...$, so that
$$\be(z)-\dfrac{\pi}{2}-\dfrac{p(n)}{q(n)}\sim\dfrac{C}{n^{z+1/2}}$$
$$A=1-((2z+1)(4z-3)/(2(2z-1)(2z+3)))/n+\cdots$$
Parametric family for $k\ge0$:
\begin{verbatim}
[(z)->beta(z)-Pi/2,8*n^2-12*n+4*(k+1)*z+2*k*(2*k+1)+3,
                    -n^2*(4*n-5+2*z)*(4*n+1-2*z)]
\end{verbatim}
Convergence type $P^+$ with $P=2k+z+1/2$.
\end{cf}

\smallskip

\begin{cf}\label{6.10.35}{\ }
\begin{verbatim}
[(z)->beta(z)-Pi/2,
[-2,3*z*(2*z-1),8*n^2+(8*z-12)*n+4*z^2-12*z+7],
[(2*z+3)*(2*z-1),(n+z-1)^2*(2*n+1)^2*(4*n+2*z+3)*(4*n+2*z-5)]]
\end{verbatim}
$$\be(z)-\dfrac{\pi}{2}=-2+\dfrac{4z^2+4z-3}{6z^2-3z+\dfrac{36z^4+108z^3-63z^2}{4z^2+4z+15+\dfrac{100z^4+900z^3+2325z^2+2350z+825}{4z^2+12z+43+\ddots}}}$$
Convergence type $P^-$ with $P=1$ and $C=1$, so that
$$\be(z)-\dfrac{\pi}{2}-\dfrac{p(n)}{q(n)}\sim\dfrac{(-1)^n}{n}$$
$$A=1-(z/2+1/4)/n+(z^2/2-z/2+3/8)/n^2+\cdots$$
Series:
$$\be(z)-\dfrac{\pi}{2}=-2+\sum_{n\ge0}(-1)^n\dfrac{4n+2z+3}{(n+z)(2n+3)}$$
\end{cf}

\smallskip

\begin{cf}\label{6.10.36}{\ }
\begin{verbatim}
[(z)->beta(z)-Pi/2,[0,z,4*n+2*z-5],[-2*z+1,(2*n-1)^2*(n+z-1)^2]]
\end{verbatim}
$$\be(z)-\dfrac{\pi}{2}=\dfrac{1-2z}{z+\dfrac{z^2}{2z+3+\dfrac{9z^2+18z+9}{2z+7+\dfrac{25z^2+100z+100}{2z+11+\ddots}}}}$$
Convergence type $P^-$ with $P=2$ and $C=(1-2z)/4$, so that
$$\be(z)-\dfrac{\pi}{2}\sim(-1)^n\dfrac{(1-2z)/4}{n^2}$$
$$A=1-(z-1/2)/n+(z^2-z-1/2)/n^2-(z^3-(3/2)z^2+(3/4)z-(5/8))/n^3+\cdots$$
Series:
$$\be(z)-\dfrac{\pi}{2}=(1-2z)\sum_{n\ge0}\dfrac{(-1)^n}{(2n+1)(n+z)}$$
Parametric family for $k\ge0$ and $u$:
\begin{verbatim}
[(z)->beta(z)-Pi/2,(2*k+1)*(4*n+4*u+2*z-5),(2*n-1)^2*(n+2*u+z-1)^2]
\end{verbatim}
Convergence type $P^-$ with $P=4k+2$.
\end{cf}

\smallskip

\begin{cf}\label{6.10.37}{\ }
\begin{verbatim}
[(z)->beta(z)-Pi/2,[0,8*n^2-8*n+3],[-2*z+1,-n^2*(16*n^2-(2*z-1)^2)]]
\end{verbatim}
$$\be(z)-\dfrac{\pi}{2}=\dfrac{1-2z}{3+\dfrac{4z^2-4z-15}{19+\dfrac{16z^2-16z-252}{51+\dfrac{36z^2-36z-1287}{99+\dfrac{64z^2-64z-4080}{163+\ddots}}}}}$$
Convergence type $P^+$ with $P=|z-1/2|$ and
$C=-\pi\sin((2z-1)\pi/4)\G(z)^2/((2z-1)\G((2z-1)/4)^22^{2z-4})$, so that
$$\be(z)-\dfrac{\pi}{2}-\dfrac{p(n)}{q(n)}\sim-\dfrac{\pi\sin((2z-1)\pi/4)\G(z)^2/((2z-1)\G((2z-1)/4)^22^{2z-4})}{n^{|z-1/2|}}$$
$$A=1-(z/2-1/4)/n+\cdots$$
Parametric family for $k\ge0$ and $u$:
\begin{verbatim}
[(z)->beta(z)-Pi/2,8*n^2-8*n+3+4*k*(z+2*u)+2*k*(2*k-1),
                    -n^2*(16*n^2-(2*z+4*u-1)^2)]
\end{verbatim}
Convergence type $P^+$ with $P=|z-1/2+2u+2k|$.
\end{cf}

\smallskip

\begin{cf}\label{6.10.38}{\ }
\begin{verbatim}
[(z)->beta(z)+Pi/2,[0,8*n^2-12*n+5],
                   [2*z-1,-n^2*(4*n+2*z-3)*(4*n-2*z-1)]]
\end{verbatim}
$$\be(z)+\dfrac{\pi}{2}=\dfrac{2z-1}{1+\dfrac{4z^2-4z-3}{13+\dfrac{16z^2-16z-140}{41+\dfrac{36z^2-36z-891}{85+\dfrac{64z^2-64z-3120}{145+\ddots}}}}}$$
Convergence type $P^+$ with $P=z-1/2$ and $C=...$, so that
$$\be(z)+\dfrac{\pi}{2}-\dfrac{p(n)}{q(n)}\sim\dfrac{C}{n^{z-1/2}}$$
$$A=1+((2z-1)/(2(2z-3)(2z+1)))/n+\cdots$$
Parametric family for $k\ge0$:
\begin{verbatim}
[(z)->beta(z)+Pi/2,8*n^2-12*n+5+4*k*z+2*k*(2*k-1),
                    -n^2*(4*n+2*z-3)*(4*n-2*z-1)]
\end{verbatim}
Convergence type $P^+$ with $P=z-1/2+2k$.
\end{cf}

\smallskip

\begin{cf}\label{6.10.39}{\ }
\begin{verbatim}
[(z)->beta(z)+Pi/2,
[0,2*z^2-3*z,8*n^2+(8*z-20)*n+4*z^2-12*z+11],
[(2*z+1)*(2*z-3),(n+z-1)^2*(2*n-1)^2*(4*n+2*z-7)*(4*n+2*z+1)]]
\end{verbatim}
$$\be(z)+\dfrac{\pi}{2}=\dfrac{4z^2-4z-3}{2z^2-3z+\dfrac{4z^4+4z^3-15z^2}{4z^2+4z+3+\dfrac{36z^4+252z^3+477z^2+342z+81}{4z^2+12z+23+\ddots}}}$$
Convergence type $P^-$ with $P=1$ and $C=1$, so that
$$\be(z)+\dfrac{\pi}{2}-\dfrac{p(n)}{q(n)}\sim(-1)^n\dfrac{1}{n}$$
$$A=1-(z/2-1/4)/n+(z^2/2-z/2-1/8)/n^2-(z^3/3-(3/4)z^2+1/8)/n^3+\cdots$$
Series:
$$\be(z)+\dfrac{\pi}{2}=\sum_{n\ge0}(-1)^n\dfrac{4n+2z+1}{(n+z)(2n+1)}$$
\end{cf}

\smallskip

\begin{cf}\label{6.10.40}{\ }
\begin{verbatim}
[(z)->beta(z)+Pi/2,[2,3*z,4*n+2*z-3],[-2*z+3,(2*n+1)^2*(n+z-1)^2]]
\end{verbatim}
$$\be(z)+\dfrac{\pi}{2}=2-\dfrac{2z-3}{3z+\dfrac{9z^2}{2z+5+\dfrac{25z^2+50z+25}{2z+9+\dfrac{49z^2+196z+196}{2z+13+\dfrac{81z^2+486z+729}{2z+17+\ddots}}}}}$$
Convergence type $P^-$ with $P=2$ and $C=(3-2z)/4$, so that
$$\be(z)+\dfrac{\pi}{2}\sim(-1)^n\dfrac{(3-2z)/4}{n^2}$$
$$A=1-(z+1/2)/n+z^2/n^2-(z^3-z^2/2-(3/4)z-1/8)/n^3+\cdots$$
Series:
$$\be(z)+\dfrac{\pi}{2}=2+(3-2z)\sum_{n\ge0}\dfrac{(-1)^n}{(n+z)(2n+3)}$$
Parametric family for $k\ge0$ and $u$:
\begin{verbatim}
[(z)->beta(z)+Pi/2,(2*k+1)*(4*n+4*u+2*z-3),(2*n+1)^2*(n+2*u+z-1)^2]
\end{verbatim}
Convergence type $P^-$ with $P=4k+2$.
\end{cf}

\smallskip

\begin{cf}\label{6.10.41}{\ }
\begin{verbatim}
[(z)->beta(z)+Pi/2,[2,8*n^2-8*n+4*z+1],[-2*z+3,-n^2*(16*n^2-(2*z-3)^2)]]
\end{verbatim}
$$\be(z)+\dfrac{\pi}{2}=2-\dfrac{2z-3}{4z+1+\dfrac{4z^2-12z-7}{4z+17+\dfrac{16z^2-48z-220}{4z+49+\dfrac{36z^2-108z-1215}{4z+97+\ddots}}}}$$
Convergence type $P^+$ with $P=|z+1/2|$ and
$C=\pi\sin((2z+1)\pi/4)\G(z)^2/((2z+1)\G((2z+1)/4)^22^{2z})$, so that
$$\be(z)+\dfrac{\pi}{2}-\dfrac{p(n)}{q(n)}\sim\dfrac{\pi\sin((2z+1)\pi/4)\G(z)^2/((2z+1)\G((2z+1)/4)^22^{2z})}{n^{|z+1/2|}}$$
$$A=1-(z/2+1/4)/n+\cdots$$
Parametric family for $k\ge0$ and $u$:
\begin{verbatim}
[(z)->beta(z)+Pi/2,8*n^2-8*n+1+4*(k+1)*(z+2*u)+2*k*(2*k+1),
                    -n^2*(16*n^2-(2*z+4*u-3)^2)]
\end{verbatim}
Convergence type $P^+$ with $P=|z+1/2+2u+2k|$.
\end{cf}

\smallskip

\begin{cf}\label{6.10.42}{\ }
\begin{verbatim}
[(a,b)->beta(a)-beta(b),[0,a*b,2*n+a+b-3],[b-a,(n+a-1)^2*(n+b-1)^2]]
\end{verbatim}
$$\be(a)-\be(b)=\dfrac{b-a}{ba+\dfrac{b^2a^2}{a+b+1+\dfrac{(b^2+2b+1)a^2+(2b^2+4b+2)a+(b^2+2b+1)}{a+b+3+\ddots}}}$$
Convergence type $P^-$ with $P=2$ and $C=(b-a)/2$, so that
$$\be(a)-\be(b)-\dfrac{p(n)}{q(n)}\sim(-1)^n\dfrac{(b-a)/2}{n^2}$$
$$A=1-(a+b-1)/n+(a^2+ab+b^2-(3/2)(a+b))/n^2+\cdots$$
Series:
$$\be(a)-\be(b)=(b-a)\sum_{n\ge0}\dfrac{(-1)^n}{(n+a)(n+b)}$$
Parametric family for $k\ge0$:
\begin{verbatim}
[(a,b)->beta(a)-beta(b),(2*k+1)*(2*n+a+b-3),(n+a-1)^2*(n+b-1)^2]
\end{verbatim}
Convergence type $P^-$ with $P=4k+2$.
\end{cf}

\smallskip

\begin{cf}\label{6.10.43}{\ }
\begin{verbatim}
[(a,b)->beta(a)-beta(b),[0,4*n^2-4*n+2*a*b-a-b+2],[b-a,-n^2*(4*n^2-(a-b)^2)]]
\end{verbatim}
$$\be(a)-\be(b)=\dfrac{b-a}{(2b-1)a-b+2+\dfrac{a^2-2ba+b^2-4}{(2b-1)a-b+10+\dfrac{4a^2-8ba+4b^2-64}{(2b-1)a-b+26+\ddots}}}$$
Convergence type $P^+$ with $P=a+b-1$ and $C=...$, so that
$$\be(a)-\be(b)-\dfrac{p(n)}{q(n)}\sim\dfrac{C}{n^{a+b-1}}$$
$$A=1-((a+b-1)/2)/n+\cdots$$
Parametric family for $k\ge0$:
\begin{verbatim}
[(a,b)->beta(a)-beta(b),4*n^2-4*n+2*a*b+(2*k-1)*(a+b)+2*k*(k-1)+2,
                         -n^2*(4*n^2-(a-b)^2)]
\end{verbatim}
Convergence type $P^+$ with $P=a+b+2k-1$.
\end{cf}

\smallskip

\begin{cf}\label{6.10.44}{\ }
\begin{verbatim}
[(a,b)->beta(a)-beta(b),
[[0,a*b],[5*n^2+2*(a+b-2)*n+(a-1)*(b-1),5*n^2+2*(a+b)*n+a*b]],
[[b-a,a^2*b^2],[-n^2*(4*n^2-(b-a)^2),(n+a)^2*(n+b)^2]]]
\end{verbatim}
$$\be(a)-\be(b)=\dfrac{b-a}{ba+\dfrac{b^2a^2}{(b+1)a+b+2+\dfrac{a^2-2ba+b^2-4}{(b+2)a+2b+5+\ddots}}}$$
Convergence type $E$ with $E=-((1+\sqrt{5})/2)^5$, $P=0$, and $C=...$, so that
$$\be(a)-\be(b)-\dfrac{p(n)}{q(n)}\sim(-1)^n\dfrac{C}{((1+\sqrt{5})/2)^{5n}}$$
\end{cf}

\smallskip

\begin{cf}\label{6.10.45}{\ }
\begin{verbatim}
[(a,b)->beta(a)+beta(b),[0,(2*a-1)*(2*b-1)+8*n^2-12*n+5],
                        [2*(a+b-1),-4*n^2*((2*n-1)^2-(a-b)^2)]]
\end{verbatim}
$$\be(a)+\be(b)=\dfrac{2a+(2b-2)}{(4b-2)a+(-2b+2)+\dfrac{4a^2-8ba+(4b^2-4)}{(4b-2)a+(-2b+14)+\dfrac{16a^2-32ba+(16b^2-144)}{(4b-2)a+(-2b+42)+\ddots}}}$$
Convergence type $P^+$ with $P=a+b-1$ and $C=...$, so that
$$\be(a)+\be(b)-\dfrac{p(n)}{q(n)}\sim\dfrac{C}{n^{a+b-1}}$$
$$A=1-(((2b-1)a-b)(a+b-1)/(2(a+b)(a+b-2)))/n+\cdots$$
Parametric family for $k\ge0$:
\begin{verbatim}
[(a,b)->beta(a)+beta(b),(2*a-1)*(2*b-1)+4*k*(a+b)+(2*k-1)^2+4*(n-1)*(2*n-1),
                        -4*n^2*((2*n-1)^2-(a-b)^2)]
\end{verbatim}
Convergence type $P^+$ with $P=2k+a+b-1$.
\end{cf}

\smallskip

\begin{cf}\label{6.10.46}{\ }
\begin{verbatim}
[(a,b)->beta(a)+beta(b),
[0,a*b*(a+b-2),2*n^2+2*(a+b-3)*n+a*(a-3)+b*(b-3)+4],
[(a+b)*(a+b-2),(n+a-1)^2*(n+b-1)^2*(2*n+a+b)*(2*n+a+b-4)]]
\end{verbatim}
$$\be(a)+\be(b)=\dfrac{a^2+(2b-2)a+b^2-2b}{ba^2+(b^2-2b)a+\dfrac{b^2a^4+2b^3a^3+(b^4-4b^2)a^2}{a^2+a+b^2+b+\ddots}}$$
Convergence type $P^-$ with $P=1$ and $C=1$, so that
$$\be(a)+\be(b)-\dfrac{p(n)}{q(n)}\sim\dfrac{(-1)^n}{n}$$
$$A=1-((a+b-1)/2)/n+((a^2+b^2-a-b)/2)/n^2+\cdots$$
Series:
$$\be(a)+\be(b)=\sum_{n\ge0}(-1)^n\dfrac{2n+a+b}{(n+a)(n+b)}$$
\end{cf}

\smallskip

\begin{cf}\label{6.10.47}{\ }
\begin{verbatim}
[(z)->sumalt(k=1,(-1)^(k-1)*k/(k^2+a)),[0,1,n^2-n-a],
[1/(a+1),2*a+2,(n^2+a)^2*(n^2-1)]]
\end{verbatim}
$$\sum_{k=1}^\infty(-1)^{k-1}\dfrac{k}{k^2+a}=\dfrac{1/(a+1)}{1+\dfrac{2a+2}{-a+2+\dfrac{3a^2+24a+48}{-a+6+\dfrac{8a^2+144a+648}{-a+12+\ddots}}}}$$
Convergence type $P^-$ with $P=1$ and $C=1/2$, so that
$$\sum_{k=1}^\infty(-1)^{k-1}\dfrac{k}{k^2+a}-\dfrac{p(n)}{q(n)}\sim(-1)^n\dfrac{1/2}{n}\;.$$
$$A=1-(1/2)/n-a/n^2+(3a/2+1/4)/n^3+a^2/n^4+\cdots$$
Series:
$$\sum_{k=1}^\infty(-1)^{k-1}\dfrac{k}{k^2+a}=\sum_{n\ge1}(-1)^{n+1}\dfrac{n}{n^2+a}$$
\end{cf}

\medskip

\section{Function $\psi'(z)$}\label{sec:psi1}

\medskip

Recall that
$$\psi'(z)=\sum_{k\ge0}\dfrac{1}{(k+z)^2}\;.$$

We also have the following integrals:
$$\int_0^\infty\dfrac{te^{-tz}}{\sinh(t)}\,dt=\psi'((z+1)/2)/2\text{\ and\ }\int_0^\infty te^{-tz}\cotanh(t)\,dt=1/z^2+\psi'(1+z/2)/2\;,$$
so the CFs in this section give CFs for the integrals.

\smallskip

\begin{cf}\label{6.11.1}{\ }
\begin{verbatim}
[(z)->psi'(z),[0,z^2,2*n^2+(4*z-6)*n+2*z^2-6*z+5],[1,-(n+z-1)^4]]
\end{verbatim}
$$\psi'(z)=\dfrac{1}{z^2-\dfrac{z^4}{2z^2+2z+1-\dfrac{z^4+4z^3+6z^2+4z+1}{2z^2+6z+5-\dfrac{z^4+8z^3+24z^2+32z+16}{2z^2+10z+13-\ddots}}}}$$
Convergence type $P^+$ with $P=1$ and $C=1$, so that
$$\psi'(z)-\dfrac{p(n)}{q(n)}\sim\dfrac{1}{n}\;.$$
$$A=1-(z-1/2)/n+(z^2-z+1/6)/n^2-(z^3-3z^2/2+z/2)/n^3+\cdots$$
Series:
$$\psi'(z)=\sum_{n\ge1}\dfrac{1}{(n+z-1)^2}$$
Parametric family for $k\ge0$:
\begin{verbatim}
[(z)->psi'(z),2*n^2+(4*z-6)*n+2*z^2-6*z+5+k^2+k,-(n+z-1)^4]
\end{verbatim}
Convergence type $P^+$ with $P=2k+1$.
\end{cf}

\smallskip

\begin{cf}\label{6.11.2}{\ }
\begin{verbatim}
[(z)->psi'(z),[0,2*n^2+(z-3)*n+1],[1,-n^4+(-z+1)*n^3]]
\end{verbatim}
$$\psi'(z)=\dfrac{1}{z-\dfrac{z}{2z+3-\dfrac{8z+8}{3z+10-\dfrac{27z+54}{4z+21-\dfrac{64z+192}{5z+36-\dfrac{125z+500}{6z+55-\ddots}}}}}}$$
Convergence type $P^+$ with $P=z$ and $C=\G(z)/z$, so that
$$\psi'(z)-\dfrac{p(n)}{q(n)}\sim\dfrac{\G(z)/z}{n^z}\;.$$
$$A=1-((z^3+z)/(2z+2))/n+\cdots$$
Series:
$$\psi'(z)=\dfrac{1}{z}\sum_{n\ge0}\dfrac{n!}{(n+1)(z+1)_n}$$
Parametric family for $k\ge0$:
\begin{verbatim}
[(z)->psi'(z),2*n^2+(z-3)*n+k*z+1+k^2,-n^4+(-z+1)*n^3]
\end{verbatim}
Convergence type $P^+$ with $P=2k+z$.
\end{cf}

\smallskip

\begin{cf}\label{6.11.3}{\ }
\begin{verbatim}
[(z)->psi'(z),[1/z+1/(2*z^2),z*(2*n+1)],[1/(2*z^2),n*(n+1)^2*(n+2)/4]]
\end{verbatim}
$$\psi'(z)=\dfrac{1}{z}+\dfrac{1}{2z^2}+\dfrac{1/(2z^2)}{3z+\dfrac{3}{5z+\dfrac{18}{7z+\dfrac{60}{9z+\dfrac{150}{11z+\dfrac{315}{13z+\ddots}}}}}}$$
Convergence type $P^-$ with $P=4z$ and $C=\G(z)^2\G(z+1)^2$, so that
$$\psi'(z)-\dfrac{p(n)}{q(n)}\sim(-1)^n\dfrac{\G(z)^2\G(z+1)^2}{n^{4z}}\;.$$
$$A=1-6z/n+(-4z^3/3+18z^2+13z/3)/n^2+(8z^4-32z^3-26z^2-4z)/n^3+\cdots$$
Parametric family for $k\ge0$:
\begin{verbatim}
[(z)->psi'(z),(z+k)*(2*n+1),n*(n+1)^2*(n+2)/4]
\end{verbatim}
Convergence type $P^-$ with $P=4k+4z$.
\end{cf}

\smallskip

\begin{cf}\label{6.11.4}{\ }
\begin{verbatim}
[(z)->psi'(z),[1/z+1/(2*z^2),(z+1)^2,2*(n+z-1)^2+2],
            [1/(2*z^2),-(n+z)^2*(n+z-1)^2]]
\end{verbatim}
$$\psi'(z)=\dfrac{1}{z}+\dfrac{1}{2z^2}+\dfrac{1/(2z^2)}{z^2+2z+1-\dfrac{z^4+2z^3+z^2}{2z^2+4z+4-\dfrac{z^4+6z^3+13z^2+12z+4}{2z^2+8z+10-\ddots}}}$$
Convergence type $P^+$ with $P=3$ and $C=1/6$, so that
$$\psi'(z)-\dfrac{p(n)}{q(n)}\sim\dfrac{1/6}{n^3}\;.$$
$$A=1-3z/n+(6z^2-1/5)/n^2+\cdots$$
Series:
$$\psi'(z)=\dfrac{2z+1}{2z^2}+\dfrac{1}{2}\sum_{n\ge1}\dfrac{1}{(n+z-1)^2(n+z)^2}$$
Parametric family for $k\ge0$:
\begin{verbatim}
[(z)->psi'(z),2*(n+z-1)^2+k^2+3*k+2,-(n+z)^2*(n+z-1)^2]
\end{verbatim}
Convergence type $P^+$ with $P=2k+3$.
\end{cf}

\smallskip

\begin{cf}\label{6.11.5}{\ }
\begin{verbatim}
[(z)->psi'(z),[0,(2*z-1)*(2*n-1)],[2,n^4]]
\end{verbatim}
$$\psi'(z)=\dfrac{2}{2z-1+\dfrac{1}{6z-3+\dfrac{16}{10z-5+\dfrac{81}{14z-7+\dfrac{256}{18z-9+\dfrac{625}{22z-11+\ddots}}}}}}$$
Convergence type $P^-$ with $P=|4z-2|$ and $C=\G(z)^4$, so that
$$\psi'(z)-\dfrac{p(n)}{q(n)}\sim(-1)^n\dfrac{\G(z)^4}{n^{|4z-2|}}\;.$$
$$A=1+(2z-1)/n+(4z^3/3-z/3)/n^2+(8z^4/3-16z^3/3+4z^2/3+z/3)/n^3+\cdots$$
This is true for $z>1/2$. The CF makes no sense for $z=1/2$, and converges
to $-\psi'(1-z)$ for $z<1/2$.

Parametric family for $k\ge0$:
\begin{verbatim}
[(z)->psi'(z),(2*z+2*k-1)*(2*n-1),n^4]
\end{verbatim}
Convergence type $P^-$ with $P=4k+4z-2$.
\end{cf}

\smallskip

When $z\in\Q$, an equivalent way of writing this CF is as follows:

\smallskip

\begin{cf}\label{6.11.6}{\ }
\begin{verbatim}
[(k,m)->psi'((m+k)/(2*k)),[0,m*(2*n-1)],[2*k,k^2*n^4]]
\end{verbatim}
$$\psi'\left(\dfrac{m+k}{2k}\right)=\dfrac{2k}{m+\dfrac{k^2}{3m+\dfrac{16k^2}{5m+\dfrac{81k^2}{7m+\dfrac{256k^2}{9m+\dfrac{625k^2}{11m+\ddots}}}}}}$$
Convergence type $P^-$ with $P=2m/k$ and $C=\G((m+k)/(2k))^4$, so that
$$\psi'\left(\dfrac{m+k}{2k}\right)-\dfrac{p(n)}{q(n)}\sim(-1)^n\dfrac{\G((m+k)/(2k))^4}{n^{2m/k}}\;.$$
$$A=1-(m/k)/n+(-m^3/(6k^3)+m^2/(2k^2)-m/(3k))/n^2+\cdots$$
\end{cf}

\smallskip

Using the formulas
$$\psi'(p/q)=\dfrac{q^2}{\phi(q)}\sum_{\chi\bmod q}\ov{\chi}(p)L(\chi,2)\;,$$
valid for $1\le p\le q$ with $\gcd(p,q)=1$, together with
$\psi'(z+1)=\psi'(z)-1/z^2$, we can specialize to a few values of
$(k,m)$ with $\gcd(k,m)=1$ to obtain many of the given CFs for linear
combinations of zeta and $L$-values, and of course infinitely many others,
for instance:

\smallskip

\begin{verbatim}
[()->Pi^2/6,[0,2*n-1],[2,n^4]]
[()->Pi^2/2,[4,2*(2*n-1)],[2,n^4]]
[()->Pi^2/6,[1,3*(2*n-1)],[2,n^4]]
[()->Pi^2/2,[40/9,4*(2*n-1)],[2,n^4]]
[()->Pi^2-8*Catalan,[0,2*n-1],[4,4*n^4]]
[()->Pi^2+8*Catalan,[16,3*(2*n-1)],[4,4*n^4]]
[()->Pi^2-8*Catalan,[16/9,5*(2*n-1)],[4,4*n^4]]
[()->Pi^2+8*Catalan,[416/25,7*(2*n-1)],[4,4*n^4]]
[()->4*Pi^2-27*lfun(-3,2),[0,2*n-1],[36,9*n^4]]
[()->4*Pi^2-45*lfun(-3,2),[0,2*(2*n-1)],[12,9*n^4]]
[()->4*Pi^2+45*lfun(-3,2),[72,4*(2*n-1)],[12,9*n^4]]
[()->4*Pi^2+27*lfun(-3,2),[54,5*(2*n-1)],[36,9*n^4]]
[()->Pi^2+8*Catalan-8*(lfun(-8,2)+lfun(8,2)),[0,2*n-1],[4,16*n^4]]
[()->Pi^2-8*Catalan-8*(lfun(-8,2)-lfun(8,2)),[0,3*(2*n-1)],[4,16*n^4]]
[()->Pi^2+8*Catalan+8*(lfun(-8,2)+lfun(8,2)),[32,5*(2*n-1)],[4,16*n^4]]
[()->Pi^2-8*Catalan+8*(lfun(-8,2)-lfun(8,2)),[32/9,7*(2*n-1)],[4,16*n^4]]
\end{verbatim}

Note also the following, coming from \ref{5.8.116} and \ref{5.8.120}, or
from \ref{6.12.3}:

\begin{verbatim}
[()->4*Pi^2-81*lfun(-3,2),[0,18*n^2-18*n+7],[-54,-81*n^4]]
[()->4*Pi^2+81*lfun(-3,2),[108,18*n^2-18*n+13],[-54,-81*n^4]]
\end{verbatim}

\smallskip

The same game can of course be played with all the other CFs for $\psi'(z)$,
but this one is the simplest.

\smallskip

\begin{cf}\label{6.11.7}{\ }
\begin{verbatim}
[(z)->psi'(z),[[0,z^2],[5*n^2+6*(z-1)*n+2*(z-1)^2,5*n^2+6*z*n+2*z^2]],
            [[1,-2*z^4],[n^4,-4*(n+z)^4]]]
\end{verbatim}
$$\psi'(z)=\dfrac{1}{z^2-\dfrac{2z^4}{2z^2+2z+1+\dfrac{1}{2z^2+6z+5-\dfrac{4z^4+16z^3+24z^2+16z+4}{2z^2+8z+10+\ddots}}}}$$
Convergence type $E$ with $E=-((1+\sqrt{5})/2)^5$, $P=0$, and
$C=4\pi^2/((1+\sqrt{5})/2)$, so that
$$\psi'(z)-\dfrac{p(n)}{q(n)}\sim(-1)^n\dfrac{4\pi^2}{((1+\sqrt{5})/2)^{5n+4z+1}}\;.$$
$$A=1+((4z^2/5-8z/5+2/5)d)/n+\cdots$$
\end{cf}

\smallskip

\begin{cf}\label{6.11.8}{\ }
\begin{verbatim}
[(z)->psi'(z),
[[0,2*z-1],[5*n^2+6*(z-1)*n+2*(z-1)^2,5*n^2+(6*z-2)*n+2*z^2-2*z+1]],
[[2,1],[-4*(n+z-1)^4,(n+1)^4]]]
\end{verbatim}
$$\psi'(z)=\dfrac{2}{2z-1+\dfrac{1}{2z^2+2z+1-\dfrac{4z^4}{2z^2+4z+4+\dfrac{16}{2z^2+8z+10-\ddots}}}}$$
Convergence type $E$ with $E=-((1+\sqrt{5})/2)^5$, $P=0$, and
$C=4\pi^2/((1+\sqrt{5})/2)^{4z+1}$, so that
$$\psi'(z)-\dfrac{p(n)}{q(n)}\sim(-1)^n\dfrac{4\pi^2}{((1+\sqrt{5})/2)^{5n+4z+1}}\;.$$
$$A=1+((4z^2/5-8z/5+2/5)d)/n+\cdots$$
\end{cf}

\smallskip

For the next CFs, recall that $\psi'(1)=\pi^2/6$.

\smallskip

\begin{cf}\label{6.11.9}{\ }
\begin{verbatim}
[(z)->psi'(z+1)-Pi^2/6,[0,(2*n-1)*(n^2-n+1+z)],[-2*z,-n^4*(n^2-z^2)]]
\end{verbatim}
$$\psi'(z+1)-\dfrac{\pi^2}{6}=-\dfrac{2z}{z+1+\dfrac{z^2-1}{3z+9+\dfrac{16z^2-64}{5z+35+\dfrac{81z^2-729}{7z+91+\dfrac{256z^2-4096}{9z+189+\ddots}}}}}$$
Convergence type $P^+$ with $P=2z+2$ and
$C=-(\G(2z+1)/\G(z+1/2))^2\sin(\pi z)/((z+1)2^{4z})$, so that
$$\psi'(z+1)-\dfrac{\pi^2}{6}-\dfrac{p(n)}{q(n)}\sim-\dfrac{(\G(2z+1)/\G(z+1/2))^2\sin(\pi z)/((z+1)2^{4z})}{n^{2z+2}}$$
$$A=1-(z+1)/n+((2z^4+7z^3+13z^2+14z+6)/(6(z+2)))/n^2+\cdots$$
Series:
$$\psi'(z+1)-\dfrac{\pi^2}{6}=-\dfrac{2z}{z+1}\sum_{n\ge0}\dfrac{(1-z)_n}{(n+1)^2(2+z)_n}$$
Parametric family for $k\ge0$:
\begin{verbatim}
[(z)->psi'(z+1)-Pi^2/6,(2*n-1)*(n^2-n+1+(2*k+1)*z+2*k*(k+1)),-n^4*(n^2-z^2)]
\end{verbatim}
Convergence type $P^+$ with $P=2z+4k+2$.
\end{cf}

\smallskip

\begin{cf}\label{6.11.10}{\ }
\begin{verbatim}
[a->sumnum(k=1,1/(k^2+a)),[0,a+1,2*n^2-2*n+2*a+1],[1,-(n^2+a)^2]]
\end{verbatim}
$$\sum_{k=1}^\infty\dfrac{1}{k^2+a}=\dfrac{1}{a+1-\dfrac{a^2+2a+1}{2a+5-\dfrac{a^2+8a+16}{2a+13-\dfrac{a^2+18a+81}{2a+25-\dfrac{a^2+32a+256}{2a+41-\dfrac{a^2+50a+625}{2a+61-\ddots}}}}}}$$
Convergence type $P^+$ with $P=1$ and $C=1$, so that
$$\sum_{k=1}^\infty\dfrac{1}{k^2+a}-\dfrac{p(n)}{q(n)}\sim\dfrac{1}{n}\;.$$
$$A=1-(1/2)/n+(-a/3+1/6)/n^2+(a/2)/n^3+(a^2/5-a/3-1/30)/n^4+\cdots$$
Series:
$$\sum_{k=1}^\infty\dfrac{1}{k^2+a}=\sum_{n\ge1}\dfrac{1}{n^2+a}$$
Parametric family for $k\ge0$:
\begin{verbatim}
[a->sumnum(k=1,1/(k^2+a)),2*n^2-2*n+2*a+1+k*(k+1),-(n^2+a)^2]
\end{verbatim}
Convergence type $P^+$ with $P=2k+1$.
\end{cf}

\smallskip

\begin{cf}\label{6.11.11}{\ }
\begin{verbatim}
[a->sumnum(k=1,1/(k^2+a)),[0,2*n-1],[2,n^2*(n^2+4*a)]]
\end{verbatim}
$$\sum_{k=1}^\infty\dfrac{1}{k^2+a}=\dfrac{2}{1+\dfrac{4a+1}{3+\dfrac{16a+16}{5+\dfrac{36a+81}{7+\dfrac{64a+256}{9+\dfrac{100a+625}{11+\ddots}}}}}}$$
Convergence type $P^-$ with $P=2$ and $C=\pi\sqrt{a}\cotanh(\pi\sqrt{a})$,
so that
$$\sum_{k=1}^\infty\dfrac{1}{k^2+a}-\dfrac{p(n)}{q(n)}\sim(-1)^n\dfrac{\pi\sqrt{a}\cotanh(\pi\sqrt{a})}{n^2}\;.$$
$$A=1-1/n-2a/n^2+(4a+1)/n^3+5a^2/n^4-(15a^2+10a+3)/n^5+\cdots$$
Series:
$$\sum_{k=1}^\infty\dfrac{1}{k^2+a}=\dfrac{1}{2(a+1)}\sum_{n\ge0}\dfrac{(4n+3)n!^2(1/2+i\sqrt{a})_n(1/2-i\sqrt{a})_n}{(2+i\sqrt{a})_n(2-i\sqrt{a})_n(3/2)_n^2}$$
Parametric family for $k\ge0$:
\begin{verbatim}
[a->sumnum(k=1,1/(k^2+a)),(2*k+1)*(2*n-1),n^2*(n^2+4*a)]
\end{verbatim}
Convergence type $P^-$ with $P=4k+2$.
\end{cf}

\smallskip

\begin{cf}\label{6.11.12}{\ }
\begin{verbatim}
[(a)->sumnum(k=1,1/(k^2+a)),[[2,2*a+2],[3*n+1,3*n+3]],
          [[-4*a-1,4*a+4],[(2*n+1)*((n+1)^2+4*a),8*(n+1)*((n+1)^2+a)]]]
\end{verbatim}
$$\sum_{k=1}^\infty\dfrac{1}{k^2+a}=2-\dfrac{4a+1}{2a+2+\dfrac{4a+4}{4+\dfrac{12a+12}{6+\dfrac{16a+64}{7+\dfrac{20a+45}{9+\dfrac{24a+216}{10+\ddots}}}}}}$$
Convergence type $E$ with $E=-2$, $P=5/2$, and $C=...$, so that
$$\sum_{k=1}^\infty\dfrac{1}{k^2+a}-\dfrac{p(n)}{q(n)}\sim(-1)^n\dfrac{C}{2^nn^{5/2}}\;.$$
$$A=1-(6a+35/4)/n+(18a^2+149a/2+2009/32)/n^2+\cdots$$
Series:
$$\sum_{k=1}^\infty\dfrac{1}{k^2+a}=2-\dfrac{4a+1}{6(a+1)}\sum_{n\ge0}\dfrac{(3n+4)n!(2+2i\sqrt{a})_n(2-2i\sqrt{a})_n}{(n+1)(n+2)(2+i\sqrt{a})_n(2-i\sqrt{a})_n(5/2)_n}2^{-2n}$$
\end{cf}

\smallskip
  
\begin{cf}\label{6.11.13}{\ }
\begin{verbatim}
[a->sumnum(k=1,1/(k^2+a)),[[0,1],[5*n^2+2*a,5*n^2+4*n+2*a+1]],
      [[2,4*a+1],[-4*(n^2+a)^2,(n+1)^2*((n+1)^2+4*a)]]]
\end{verbatim}
$$\sum_{k=1}^\infty\dfrac{1}{k^2+a}=\dfrac{2}{1+\dfrac{4a+1}{2a+5-\dfrac{4a^2+8a+4}{2a+10+\dfrac{16a+16}{2a+20-\dfrac{4a^2+32a+64}{2a+29+\dfrac{36a+81}{2a+45-\ddots}}}}}}$$
Convergence type $E$ with $E=-((1+\sqrt{5})/2)^5$, $P=0$, and
$C=2\pi\sinh(2\pi\sqrt{a})/\sqrt{a}/((1+\sqrt{5})/2)^5$, so that
$$\sum_{k=1}^\infty\dfrac{1}{k^2+a}-\dfrac{p(n)}{q(n)}\sim(-1)^n\dfrac{2\pi\sinh(2\pi\sqrt{a})/\sqrt{a}}{((1+\sqrt{5})/2)^{5n+5}}\;.$$
\begin{align*}A&=1-((24a/5+132/125)d)/n+(288a^2/5+(604d/125+3168/125)a\\
  &\phantom{=}+7936d/9375+8712/3125)/n^2+\cdots\end{align*}
\end{cf} 

\smallskip

\begin{cf}\label{6.11.14}{\ }
\begin{verbatim}
[a->sumalt(k=1,(-1)^(k-1)/(k^2+a)),[0,2*n^2-2*n+a+1],[1/2,-n^2*(n^2+a)]]
\end{verbatim}
$$\sum_{k=1}^\infty\dfrac{(-1)^{k-1}}{k^2+a}=\dfrac{1/2}{a+1-\dfrac{a+1}{a+5-\dfrac{4a+16}{a+13-\dfrac{9a+81}{a+25-\dfrac{16a+256}{a+41-\dfrac{25a+625}{a+61-\ddots}}}}}}$$
Convergence type $P^+$ with $P=1$ and $C=\pi\sqrt{a}/(2\sinh(\pi\sqrt{a}))$, so
that
$$\sum_{k=1}^\infty\dfrac{(-1)^{k-1}}{k^2+a}-\dfrac{p(n)}{q(n)}\sim\dfrac{\pi\sqrt{a}/(2\sinh(\pi\sqrt{a}))}{n}\;.$$
$$A=1-((a+1)/2)/n+((a+1)^2/6)/n^2-(a(a+1)^2/24)/n^3+\cdots$$
Series:
$$\sum_{k=1}^\infty\dfrac{(-1)^{k-1}}{k^2+a}=\dfrac{1}{2a+2}\sum_{n\ge0}\dfrac{n!^2}{(2+i\sqrt{a})_n(2-i\sqrt{a})_n}$$
Parametric family for $k\ge0$:
\begin{verbatim}
[a->sumalt(k=1,(-1)^(k-1)/(k^2+a)),2*n^2-2*n+a+1+k*(k+1),-n^2*(n^2+a)]
\end{verbatim}
Convergence type $P^+$ with $P=2k+1$.
\end{cf}
    
\smallskip

\begin{cf}\label{6.11.15}{\ }
\begin{verbatim}
[a->sumalt(k=1,(-1)^(k-1)/(k^2+a)),[0,a+1,2*n-1],[1,(n^2+a)^2]]
\end{verbatim}
$$\sum_{k=1}^\infty\dfrac{(-1)^{k-1}}{k^2+a}=\dfrac{1}{a+1+\dfrac{a^2+2a+1}{3+\dfrac{a^2+8a+16}{5+\dfrac{a^2+18a+81}{7+\dfrac{a^2+32a+256}{9+\dfrac{a^2+50a+625}{11+\ddots}}}}}}$$
Convergence type $P^-$ with $P=2$ and $C=1/2$, so that
$$\sum_{k=1}^\infty\dfrac{(-1)^{k-1}}{k^2+a}-\dfrac{p(n)}{q(n)}\sim(-1)^n\dfrac{1/2}{n^2}\;.$$
$$A=1-1/n-a/n^2+(2a+1)/n^3+a^2/n^4-(3a^2+5a+3)/n^5-\cdots$$
Series:
$$\sum_{k=1}^\infty\dfrac{(-1)^{k-1}}{k^2+a}=\sum_{n\ge1}\dfrac{(-1)^{n+1}}{n^2+a}$$
Parametric family for $k\ge0$:
\begin{verbatim}
[a->sumalt(k=1,(-1)^(k-1)/(k^2+a)),(2*k+1)*(2*n-1),(n^2+a)^2]
\end{verbatim}
Convergence type $P^-$ with $P=4k+2$.
\end{cf}

\smallskip

\begin{cf}\label{6.11.16}{\ }
\begin{verbatim}
[a->sumalt(k=1,(-1)^(k-1)/(k^2+a)),
[[0,2*a+2],[5*n^2+a,5*n^2+6*n+a+2]],
[[1,-4*(a+1)],[(n^2+a)^2,-4*(n+1)^2*((n+1)^2+a)]]]
\end{verbatim}
$$\sum_{k=1}^\infty\dfrac{(-1)^{k-1}}{k^2+a}=\dfrac{1}{2a+2-\dfrac{4a+4}{a+5+\dfrac{a^2+2a+1}{a+13-\dfrac{16a+64}{a+20+\dfrac{a^2+8a+16}{a+34-\dfrac{36a+324}{a+45+\ddots}}}}}}$$
Convergence type $E$ with $E=-((1+\sqrt{5})/2)^5$, $P=0$, and $C=2\pi\sinh(\pi\sqrt{a})/\sqrt{a}/((1+\sqrt{5})/2)^5$, so that
$$\sum_{k=1}^\infty\dfrac{(-1)^{k-1}}{k^2+a}-\dfrac{p(n)}{q(n)}\sim(-1)^n\dfrac{2\pi\sinh(\pi\sqrt{a})/\sqrt{a}}{((1+\sqrt{5})/2)^{5n+5}}\;.$$
$$A=1-((2a+2/5)d)/n+(10a^2+(2d+4)a+2d/5+2/5)/n^2+\cdots$$
\end{cf}

\smallskip

Note that
$$\sum_{k=1}^\infty\dfrac{1}{k^2+a}=\dfrac{(\pi/2)\cotanh(\pi\sqrt{a})}{\sqrt{a}}-\dfrac{1}{2a}\text{\quad and\quad}\sum_{k=1}^\infty\dfrac{(-1)^{k-1}}{k^2+a}=\dfrac{1}{2a}-\dfrac{\pi/2}{\sqrt{a}\sinh(\pi\sqrt{a})}\;,$$
thus linking the CFs for the LHS to those for $\sinh(z)$ and $\cotanh(z)$
in Section \ref{sec:trighyp}.

\section{Function $\be_1(z)=\psi'(z/2)/2-\psi'(z)$}\label{sec:beta1}
  
\smallskip

Analogously to the definition and properties of $\be(z)$ given above,
we define $\be_1(z)$ (my notation) by the following equivalent formulas:
\begin{align*}\be_1(z)&=-\be'(z)=\psi'(z/2)/2-\psi'(z)=\psi'(1+z/2)/2-\psi'(1+z)+1/z^2\\
  &=\dfrac{1}{4}(\psi'(z/2)-\psi'((1+z)/2))=\sum_{k\ge0}\dfrac{(-1)^k}{(z+k)^2}\;,\end{align*}
and note the following integral formulas, so that all CFs for $\be_1(z)$
give CFs for the integrals:

\begin{align*}
  \int_0^\infty e^{-zt}t\tanh(t)\,dt&=\be_1(z/2)/2-1/z^2\;,\\
  \int_0^\infty \dfrac{te^{-zt}}{\cosh(t)}\,dt&=\be_1((z+1)/2)/2\;,\\
  \int_0^\infty e^{-zt}(tz+1)\tanh^2(t)\,dt&=-z^2\be_1(z/2)/2+1+2/z\;,\\
  \int_0^\infty \dfrac{(tz+1)e^{-zt}}{\cosh^2(t)}\,dt&=z^2\be_1(z/2)/2-1\;.
\end{align*}

\smallskip

\begin{verbatim}
beta1(z)=psi'(z/2)/2-psi'(z);
\end{verbatim}

\smallskip

\begin{cf}\label{6.12.1}{\ }
\begin{verbatim}
[(z)->beta1(z),[1/z^2,(z+1)^2,2*(n+z)-1],[-1,(n+z)^4]]
\end{verbatim}
$$\be_1(z)=\dfrac{1}{z^2}-\dfrac{1}{z^2+2z+1+\dfrac{z^4+4z^3+6z^2+4z+1}{2z+3+\dfrac{z^4+8z^3+24z^2+32z+16}{2z+5+\ddots}}}$$
Convergence type $P^-$ with $P=2$ and $C=-1/2$, so that
$$\be_1(z)-\dfrac{p(n)}{q(n)}\sim(-1)^{n+1}\dfrac{1/2}{n^2}\;.$$
$$A=1-(2z+1)/n+(3z(z+1))/n^2-((2z+1)(2z^2+2z-1))/n^3+\cdots$$
Series:
$$\be_1(z)=\sum_{n\ge0}\dfrac{(-1)^n}{(n+z)^2}$$
Parametric family for $k\ge0$:
\begin{verbatim}
[(z)->beta1(z),(2*k+1)*(2*(n+z)-1),(n+z)^4]
\end{verbatim}
Convergence type $P^-$ with $P=4k+2$.
\end{cf}

\smallskip

\begin{cf}\label{6.12.2}{\ }
\begin{verbatim}
[(z)->beta1(z),[0,2*n*(n-1)+z^2-z+1],[1/2,-n^4]]
\end{verbatim}
$$\be_1(z)=\dfrac{1/2}{z^2-z+1-\dfrac{1}{z^2-z+5-\dfrac{16}{z^2-z+13-\dfrac{81}{z^2-z+25-\ddots}}}}$$
Convergence type $P^+$ with $P=2z-1$ and $C=(z-1/2)\G(z)^6/\G(2z)^2$, so that
$$\be_1(z)-\dfrac{p(n)}{q(n)}\sim\dfrac{(z-1/2)\G(z)^6/\G(2z)^2}{n^{2z-1}}\;.$$
$$A=1-(z-1/2)/n-((2z-1)(z^4-14z^3+16z^2+6z-3)/(12(2z-3)(2z+1)))/n^2+\cdots$$
Parametric family up to trivial change $n\mapsto n+j$:
\begin{verbatim}
[(z)->beta1(z),2*n^2+(u-2)*n+z^2+(u+2*k-1)*z+(k-1)*u+k^2-k+1,-n^3*(n+u)]
\end{verbatim}
Convergence type $P^+$ with $P=2z+2k+u-1$.
\end{cf}

When $z\in\Q$, an equivalent way of writing this CF is as follows:

\smallskip

\begin{cf}\label{6.12.3}{\ }
\begin{verbatim}
[(k,m)->beta1((m+k)/k),[0,2*k^2*n*(n-1)+m^2+k*m+k^2],[k^2/2,-k^4*n^4]]
\end{verbatim}
$$\be_1\left(\dfrac{m+k}{k}\right)=\dfrac{k^2/2}{m^2+km+k^2-\dfrac{k^4}{m^2+km+5k^2-\dfrac{16k^4}{m^2+km+13k^2-\ddots}}}$$
Convergence type $P^+$ with $P=2m/k+1$ and $C=((m+k/2)/k)\G((m+k)/k)^6/\G((2m+2k)/k)^2$, so that
$$\be_1(z)-\dfrac{p(n)}{q(n)}\sim\dfrac{((m+k/2)/k)\G((m+k)/k)^6/\G((2m+2k)/k)^2}{n^{2m/k+1}}\;.$$
$$A=1-(m/k+1/2)/n-((2m+k)(m^4-10km^3-20k^2m^2+6k^4)/(12k^3(2m-k)(2m+3k)))/n^2+\cdots$$
\end{cf}

\smallskip

Using explicit formulas for $\be_1(p/q)$ we can specialize to a few values of
$(k,m)$ with $(k,m)=1$ to obtain many of the given CFs for linear combinations
of zeta and $L$-values, and of course infinitely many others, for instance:

\begin{verbatim}
[()->Pi^2/6,[2,2*n^2-2*n+3],[-1,-n^4]]
[()->Catalan,[8/9,8*n^2-8*n+19],[1/2,-16*n^4]]
[()->4*Pi^2-81*lfun(-3,2),[-27,18*n^2-18*n+19],[54,-81*n^4]]
[()->4*Pi^2+81*lfun(-3,2),[108,18*n^2-18*n+13],[-54,-81*n^4]]
[()->lfun(-8,2)+lfun(8,2),[2,32*n^2-32*n+21],[-1,-256*n^4]]
[()->lfun(-8,2)-lfun(8,2),[0,32*n^2-32*n+13],[1,-256*n^4]]
\end{verbatim}

\smallskip

The same game can of course be played with all the other CFs for $\be_1(z)$,
but this one is the simplest.

\smallskip

Special case of the above:

\smallskip

\begin{cf}\label{6.12.4}{\ }
\begin{verbatim}
[(z)->intnum(t=0,[oo,1+z],t*exp(-t*z)/cosh(t)),
[0,8*n^2-8*n+3+z^2],[1,-16*n^4]]
\end{verbatim}
$$\int_0^\infty\dfrac{te^{-tz}}{\cosh(t)}\,dt=\dfrac{1}{z^2+3-\dfrac{16}{z^2+19-\dfrac{256}{z^2+51-\dfrac{1296}{z^2+99-\dfrac{4096}{z^2+163-\dfrac{10000}{z^2+243-\ddots}}}}}}$$
Convergence type $P^+$ with $P=z$ and $C=...$, so that
$$\int_0^\infty\dfrac{te^{-tz}}{\cosh(t)}\,dt-\dfrac{p(n)}{q(n)}\sim\dfrac{C}{n^z}\;.$$
$$A=1-(z/2)/n+((-z^5+24z^4+14z^3-96z^2-37z)/(192(z^2-4)))/n^2+\cdots$$
Parametric family for $k\ge0$:
\begin{verbatim}
[(z)->intnum(t=0,[oo,1+z],t*exp(-t*z)/cosh(t)),
8*n^2-8*n+z^2+4*k*z+4*k^2+3,-16*n^4]
\end{verbatim}
Convergence type $P^+$ with $P=2k+z$.
\end{cf}

\smallskip

\begin{cf}\label{6.12.5}{\ }
\begin{verbatim}
[(z)->intnum(t=0,[oo,1+z],t*exp(-t*z)/cosh(t)),
[0,(z+1)^2,4*(2*n+z-2)],[2,(2*n+z-1)^4]]
\end{verbatim}
$$\int_0^\infty\dfrac{te^{-tz}}{\cosh(t)}\,dt=\dfrac{2}{z^2+2z+1+\dfrac{z^4+4z^3+6z^2+4z+1}{4z+8+\dfrac{z^4+12z^3+54z^2+108z+81}{4z+16+\dfrac{z^4+20z^3+150z^2+500z+625}{4z+24+\ddots}}}}$$
Convergence type $P^-$ with $P=2$ and $C=1/4$, so that
$$\int_0^\infty\dfrac{te^{-tz}}{\cosh(t)}\,dt-\dfrac{p(n)}{q(n)}\sim(-1)^n\dfrac{1/4}{n^2}\;.$$
$$A=1-z/n+(3z^2/4-3/4)/n^2+(-z^3/2+3z/2)/n^3+\cdots$$
Series:
$$\int_0^\infty\dfrac{te^{-tz}}{\cosh(t)}\,dt=2\sum_{n\ge1}\dfrac{(-1)^{n+1}}{(2n+z-1)^2}$$
Parametric family for $k\ge0$:
\begin{verbatim}
[(z)->intnum(t=0,[oo,1+z],t*exp(-t*z)/cosh(t)),
4*(2*k+1)*(2*n+z-2),(2*n+z-1)^4]
\end{verbatim}
Convergence type $P^-$ with $P=4k+2$.
\end{cf}

\smallskip

\begin{cf}\label{6.12.6}{\ }
\begin{verbatim}
[(z)->beta1(z),[[1/z^2,2*z^2+2*z],[3*z*n+z^2,3*z*n+(z^2+2*z)]],
               [[-1,4*z*(z+1)],[z*(2*n+z)*(n+z)^2,4*z*(2*n+z+1)*(n+1)^2]]]
\end{verbatim}
$$\be_1(z)=\dfrac{1}{z^2}-\dfrac{1}{2z^2+2z+\dfrac{4z^2+4z}{z^2+3z+\dfrac{z^4+4z^3+5z^2+2z}{z^2+5z+\dfrac{16z^2+48z}{z^2+6z+\ddots}}}}$$
Convergence type $E$ with $E=-2$, $P=z+3/2$, and $C=-\sqrt{2\pi}\G(z+1)$,
so that
$$\be_1(z)-\dfrac{p(n)}{q(n)}\sim(-1)^{n+1}\dfrac{\sqrt{2\pi}\G(z+1)}{2^nn^{z+3/2}}\;.$$
$$A=1-(z^2/2+9z/2+15/4)/n+\cdots$$
Series:
$$\be_1(z)=\dfrac{1}{z^2}-\dfrac{3}{2(z+1)(z+2)}\sum_{n\ge0}\dfrac{(n+z/3+1)n!^2}{(n+z+1)(z/2+3/2)_n(z/2+2)_n}2^{-2n}$$
\end{cf}

\smallskip

\begin{cf}\label{6.12.7}{\ }
\begin{verbatim}
[(z)->beta1(z),
[[1/z^2,2*(z^2+z+1)],[5*n^2+4*z*n+z^2,5*n^2+(4*z+6)*n+z^2+2*z+2]],
                 [[-1,-4],[(n+z)^4,-4*(n+1)^4]]]
\end{verbatim}
$$\be_1(z)=\dfrac{1}{z^2}-\dfrac{1}{2z^2+2z+2-\dfrac{4}{z^2+4z+5+\dfrac{(z+1)^4}{z^2+6z+13-\dfrac{64}{z^2+8z+20+\ddots}}}}$$
Convergence type $E$ with $E=-((1+\sqrt{5})/2)^5$, $P=0$, and $C=-2\pi^2/((1+\sqrt{5})/2)^{6z+5}$, so
that
$$\be_1(z)-\dfrac{p(n)}{q(n)}\sim(-1)^{n+1}\dfrac{2\pi^2}{((1+\sqrt{5})/2)^{5n+6z+5}}\;.$$
\begin{align*}A&=1+((4z^2/5-2/5)d)/n+\\&\phantom{=}(8z^4/5-16dz^3/25-(4d/5+8/5)z^2+12dz/25+2d/5+2/5)/n^2+\cdots\end{align*}
\end{cf}

\medskip

\section{Function $\psi''(z)$}\label{sec:psi2}

\medskip

Recall that
$$\psi''(z)=-2\sum_{k\ge0}\dfrac{1}{(k+z)^3}\;.$$

Note also the following integrals:
\begin{align*}\int_0^\infty\dfrac{t^2e^{-tz}}{\sinh(t)}\,dt&=-\psi''((z+1)/2)/4\text{\quad and}\\
  \int_0^\infty t^2e^{-tz}\cotanh(t)\,dt&=2/z^3-\psi''(1+z/2)/4\;,\end{align*}
so the CFs in this section give CFs for the integrals.

\smallskip

\begin{cf}\label{6.13.1}{\ }
\begin{verbatim}
[(z)->psi''(z),[-2/z^3,(z+1)^3,(2*n+2*z-1)*(n^2+(2*z-1)*n+z^2-z+1)],
               [-2,-(n+z)^6]]
\end{verbatim}
$$\psi''(z)=-2/z^3-\dfrac{2}{z^3+3z^2+3z+1-\dfrac{z^6+6z^5+15z^4+20z^3+15z^2+6z+1}{2z^3+9z^2+15z+9-\ddots}}$$
Convergence type $P^+$ with $P=2$ and $C=-1$, so that
$$\psi''(z)-\dfrac{p(n)}{q(n)}\sim-\dfrac{1}{n^2}$$
$$A=1-(2z+1)/n+(3z^2+3z+1/2)/n^2-(4z^3+6z^2+2z)/n^3+\cdots$$
Series:
$$\psi''(z)=-2\sum_{n\ge0}\dfrac{1}{(z+n)^3}$$
Parametric family for $k\ge0$:
\begin{verbatim}
[(z)->psi''(z),(2*n+2*z-1)*(n^2+(2*z-1)*n+z^2-z+2*k*(k+1)+1),-(n+z)^6]
\end{verbatim}
Convergence type $P^+$ with $P=4k+2$.
\end{cf}

\smallskip

\begin{cf}\label{6.13.2}{\ }
\begin{verbatim}
[(z)->psi''(z),[[-1/z^2-1/z^3,3*z^2],[2*n+1,(2*n+3)*z^2]],
             [[-3/(2*z^2),3],(n+1)^2*(2*n+3)/2*[n,n+2]]]
\end{verbatim}
$$\psi''(z)=-\dfrac{1}{z^2}-\dfrac{1}{z^3}-\dfrac{3/(2z^2)}{3z^2+\dfrac{3}{3+\dfrac{10}{5z^2+\dfrac{30}{5+\dfrac{63}{7z^2+\dfrac{126}{7+\ddots}}}}}}$$
Convergence type $P^-$ with $P=4z$ and $C=...$, so that
$$\psi''(z)-\dfrac{p(n)}{q(n)}\sim(-1)^n\dfrac{C}{n^{4z}}\;.$$
$$A=1-8z/n-((16z^5-384z^4-104z^3+96z^2+28z)/(3(4z^2-1)))/n^2+\cdots$$
\end{cf}

\smallskip

\begin{cf}\label{6.13.3}{\ }
\begin{verbatim}
[(z)->psi''(z),[-1/z^2-1/z^3,(n+1)*z],[[-1/z^3,2],(n+1)^3*[n,n+2]]]
\end{verbatim}
$$\psi''(z)=-\dfrac{1}{z^2}-\dfrac{1}{z^3}-\dfrac{1/z^3}{2z+\dfrac{2}{3z+\dfrac{8}{4z+\dfrac{24}{5z+\dfrac{54}{6z+\dfrac{108}{7z+\ddots}}}}}}$$
Convergence type $P^-$ with $P=4z$ and $C=...$, so that
$$\psi''(z)-\dfrac{p(n)}{q(n)}\sim(-1)^n\dfrac{C}{n^{4z}}\;.$$
$$A=1-8z/n-((16z^5-384z^4-104z^3+96z^2+28z)/(3(4z^2-1)))/n^2+\cdots$$
\end{cf}

\smallskip

\begin{cf}\label{6.13.4}{\ }
\begin{verbatim}
[(z)->psi''(z),
[[-1/z^2-1/z^3,2*z],z*[2*n+1,2]],[[-1/z^3,2],(n+1)^2*[n,n+2]]]
\end{verbatim}
$$\psi''(z)=-\dfrac{1}{z^2}-\dfrac{1}{z^3}-\dfrac{1/z^3}{2z+\dfrac{2}{3z+\dfrac{4}{2z+\dfrac{12}{5z+\dfrac{18}{2z+\dfrac{36}{7z+\ddots}}}}}}$$
Convergence type $P^-$ with $P=4z$ and $C=...$, so that
$$\psi''(z)-\dfrac{p(n)}{q(n)}\sim(-1)^n\dfrac{C}{n^{4z}}\;.$$
$$A=1-8z/n-((16z^5-384z^4-104z^3+96z^2+28z)/(3(4z^2-1)))/n^2+\cdots$$
\end{cf}

\smallskip

\begin{cf}\label{6.13.5}{\ }
\begin{verbatim}
[(z)->psi''(z),[[-1/z^2-1/z^3,3*z],z*[2*n+1,2*n+3]],
             [[-3/(2*z^3),3],(n+1)^2*(2*n+3)/2*[n,n+2]]]
\end{verbatim}
$$\psi''(z)=-\dfrac{1}{z^2}-\dfrac{1}{z^3}-\dfrac{3/(2z^3)}{3z+\dfrac{3}{3z+\dfrac{10}{5z+\dfrac{30}{5z+\dfrac{63}{7z+\dfrac{126}{7z+\ddots}}}}}}$$
Convergence type $P^-$ with $P=4z$ and $C=...$, so that
$$\psi''(z)-\dfrac{p(n)}{q(n)}\sim(-1)^n\dfrac{C}{n^{4z}}\;.$$
$$A=1-8z/n-((16z^5-384z^4-104z^3+96z^2+28z)/(3(4z^2-1)))/n^2+\cdots$$
\end{cf}

\smallskip

\begin{cf}\label{6.13.6}{\ }
\begin{verbatim}
[(z)->psi''(z),[-1/z^2-1/z^3,4*n^4+(8*z^2-2)*n^2-2*z^2],
             [-3/z^2,-n^3*(n+1)^3*(2*n-1)*(2*n+3)]]
\end{verbatim}
$$\psi''(z)=-\dfrac{1}{z^2}-\dfrac{1}{z^3}-\dfrac{3/z^2}{6z^2+2-\dfrac{40}{30z^2+56-\dfrac{4536}{70z^2+306-\ddots}}}$$
Convergence type $P^+$ with $P=4z$ and $C=...$, so that
$$\psi''(z)-\dfrac{p(n)}{q(n)}\sim\dfrac{C}{n^{4z}}\;.$$
$$A=1-4z/n-((4z^5-96z^4-26z^3+24z^2+7z)/(3(4z^2-1)))/n^2+\cdots$$
\end{cf}

Note: this is the common even contraction of \ref{6.13.2}, \ref{6.13.3},
\ref{6.13.4}, and \ref{6.13.5}, and cannot be Ap\'ery accelerated.

\smallskip

\begin{cf}\label{6.13.7}{\ }
\begin{verbatim}
[(z)->psi''(z),[[0,z*(z-1)],
             [2*n,z*(z-1)*(2*n+1)]],[[-1,1],[n^4,(n+1)^4]]]
\end{verbatim}
$$\psi''(z)=-\dfrac{1}{z^2-z+\dfrac{1}{2+\dfrac{1}{3z^2-3z+\dfrac{16}{4+\dfrac{16}{5z^2-5z+\dfrac{81}{6+\ddots}}}}}}$$
Convergence type $P^-$ with $P=4z-2$ and $C=-2^{4z-4}\G(z)^8/(\G(2z)\G(2z-1))$, so that
$$\psi''(z)-\dfrac{p(n)}{q(n)}\sim(-1)^{n+1}\dfrac{2^{4z-4}\G(z)^8/(\G(2z)\G(2z-1))}{n^{4z-2}}\;.$$
$$A=1-(4z-2)/n-(2(2z-1)(z^2-7z+3)/3)/n^2+(8z(z-1)(z-2)(2z-1)/3)/n^3+\cdots$$
\end{cf}

\smallskip

\begin{cf}\label{6.13.8}{\ }
\begin{verbatim}
[(z)->psi''(z),[0,(2*n-1)*(n^2-n+1+2*z*(z-1))],[-2,-n^6]]
\end{verbatim}
$$\psi''(z)=-\dfrac{2}{2z^2-2z+1-\dfrac{1}{6z^2-6z+9-\dfrac{64}{10z^2-10z+35-\dfrac{729}{14z^2-14z+91-\ddots}}}}$$
Convergence type $P^+$ with $P=4z-2$ and $C=-\G(z)^8/(\G(2z)\G(2z-1))$, so that
$$\psi''(z)-\dfrac{p(n)}{q(n)}\sim-\dfrac{\G(z)^8/(\G(2z)\G(2z-1))}{n^{4z-2}}\;.$$
$$A=1-(2z-1)/n-((2z-1)(z^2-7z+3)/6)/n^2+((z(z-1)(z-2)(2z-1))/3)/n^3+\cdots$$
Parametric family for $k\ge0$:
\begin{verbatim}
[(z)->psi''(z),(2*n-1)*(n^2-n+2*z^2+(4*k-2)*z+2*k^2-2*k+1),-n^6]
\end{verbatim}
Convergence type $P^+$ with $P=4z+4k-2$.
\end{cf}

This is simply the contracted CF of the previous one, but I found
useful to give both.

\smallskip

As for $\psi'(z)$, when $z\in\Q$, an equivalent way of writing this CF is as
follows:

\begin{cf}\label{6.13.9}{\ }
\begin{verbatim}
[(k,m)->psi''((k+m)/k),
[0,(2*n-1)*(k^2*(n^2-n+1)+2*m*(k+m))],[-2*k^2,-k^4*n^6]]
\end{verbatim}
$$\psi''\left(\dfrac{k+m}{k}\right)=-\dfrac{2k^2}{2m^2+2km+k^2-\dfrac{k^4}{6m^2+6km+9k^2-\dfrac{64k^4}{10m^2+10km+35k^2-\ddots}}}$$
Convergence type $P^+$ with $P=4m/k+2$ and $C=-\G((k+m)/k)^8/(\G((2k+2m)/k)\G((k+2m)/k))$, so that
$$\psi''\left(\dfrac{k+m}{k}\right)-\dfrac{p(n)}{q(n)}\sim-\dfrac{\G((k+m)/k)^8/(\G((2k+2m)/k)\G((k+2m)/k))}{n^{4m/k+2}}\;.$$
$$A=1-(2m/k+1)/n-((2m+k)(m^2-5km-3k^2)/(6k^3))/n^2+\cdots$$
\end{cf}

\smallskip

Using the formulas
$$\psi''(p/q)=-\dfrac{2q^3}{\phi(q)}\sum_{\chi\bmod q}\ov{\chi}(p)L(\chi,3)\;,$$
valid for $1\le p\le q$ with $\gcd(p,q)=1$, together with
$\psi''(z+1)=\psi''(z)+2/z^3$, we can specialize to a few values of
$(k,m)$ with $\gcd(k,m)=1$ to obtain many of the given CFs for linear
combinations of zeta and $L$-values, and of course infinitely many others,
for instance:

\smallskip

\begin{verbatim}
[()->zeta(3),[1,(2*n-1)*(n^2-n+5)],[1,-n^6]]
[()->zeta(3),[9/8,(2*n-1)*(n^2-n+13)],[1,-n^6]]
[()->zeta(3),[8/7,(2*n-1)*(2*n^2-2*n+5)],[2/7,-4*n^6]]
[()->zeta(3),[32/27,(2*n-1)*(2*n^2-2*n+17)],[2/7,-4*n^6]]
[()->117*zeta(3)+2*Pi^3*sqrt(3),
                  [243,(2*n-1)*(9*n^2-9*n+17)],[81,-81*n^6]]
[()->117*zeta(3)-2*Pi^3*sqrt(3),
                  [243/8,(2*n-1)*(9*n^2-9*n+29)],[81,-81*n^6]]
[()->28*zeta(3)+Pi^3,[64,(2*n-1)*(8*n^2-8*n+13)],[8,-64*n^6]]
[()->28*zeta(3)-Pi^3,[64/27,(2*n-1)*(8*n^2-8*n+29)],[8,-64*n^6]]
[()->112*zeta(3)+4*Pi^3+3*Pi^3*sqrt(2)+128*lfun(8,3),
                  [512,(2*n-1)*(32*n^2-32*n+41)],[32,-1024*n^6]]
[()->112*zeta(3)-4*Pi^3+3*Pi^3*sqrt(2)-128*lfun(8,3),
                  [512/27,(2*n-1)*(32*n^2-32*n+65)],[32,-1024*n^6]]
[()->112*zeta(3)+4*Pi^3-3*Pi^3*sqrt(2)-128*lfun(8,3),
                  [512/125,(2*n-1)*(32*n^2-32*n+97)],[32,-1024*n^6]]
[()->112*zeta(3)-4*Pi^3-3*Pi^3*sqrt(2)+128*lfun(8,3),
                  [512/343,(2*n-1)*(32*n^2-32*n+137)],[32,-1024*n^6]]
\end{verbatim}

Special case of the above:

\smallskip

\begin{cf}\label{6.13.10}{\ }
\begin{verbatim}
[(z)->intnum(t=0,[oo,1+z],t^2*exp(-t*z)/sinh(t)),
[0,(2*n-1)*(2*n^2-2*n+1+z^2)],[1,-4*n^6]]
\end{verbatim}
$$\int_0^\infty\dfrac{t^2e^{-tz}}{\sinh(t)}\,dt=\dfrac{1}{z^2+1-\dfrac{4}{3z^2+15-\dfrac{256}{5z^2+65-\dfrac{2916}{7z^2+175-\dfrac{16384}{9z^2+369-\ddots}}}}}$$
Convergence type $P^+$ with $P=2z$ and $C=...$, so that
$$\int_0^\infty\dfrac{t^2e^{-tz}}{\sinh(t)}\,dt-\dfrac{p(n)}{q(n)}\sim\dfrac{C}{n^{2z}}\;.$$
$$A=1-z/n+(-z^3/24+z^2/2+z/24)/n^2+(z^4/24-z^3/8-z^2/24+z/8)/n^3+\cdots$$
Parametric family for $k\ge0$:
\begin{verbatim}
[(z)->intnum(t=0,[oo,1+z],t^2*exp(-t*z)/sinh(t)),
(2*n-1)*(2*n^2-2*n+1+(2*k+z)^2),-4*n^6]
\end{verbatim}
Convergence type $P^+$ with $P=4k+2z$.
\end{cf}

\smallskip

\begin{cf}\label{6.13.11}{\ }
\begin{verbatim}
[(z)->psi''(z),[[-1/z,z],[1,z]],
             [[1,(n^2-2*n+2)/(2*n-1),(n^2+1)/(2*n-1),
               n^3/(2*(n^2+1)),-n^3/(2*(n^2+1))]]]
\end{verbatim}
$$\psi''(z)=-\dfrac{1}{z}+\dfrac{1}{z+\dfrac{1}{1+\dfrac{2}{z+\dfrac{1/4}{1-\dfrac{1/4}{z+\dfrac{2/3}{1+\dfrac{5/3}{z+\dfrac{4/5}{1-\ddots}}}}}}}}$$
\end{cf}

\smallskip

\begin{cf}\label{6.13.12}{\ }
\begin{verbatim}
[(z)->psi''(z),[[0,2*z^2-2*z+1],[(2*n+z-1)*(3*n^2+3*(z-1)*n+(z-1)^2),
                               (2*n+z)*(3*n^2+3*z*n+z^2-z+1)]],
             [[-2,-1],[-(n+z-1)^6,-(n+1)^6]]]
\end{verbatim}
$$\psi''(z)=-\dfrac{2}{2z^2-2z+1-\dfrac{1}{z^3+2z^2+2z+1-\dfrac{z^6}{z^3+4z^2+8z+8-\ddots}}}$$
Convergence type $E$ with $E=(1+\sqrt{2})^4$, $P=0$, $C=-8\pi^3/(1+\sqrt{2})^{4z}$, so that
$$\psi''(z)-\dfrac{p(n)}{q(n)}\sim-\dfrac{8\pi^3}{(1+\sqrt{2})^{4n+4z}}\;.$$
\begin{align*}A&=1+((2z^2-4z+17/16)d)/n\\
  &\phantom{=}+(4z^4-(2d+16)z^3+(4d+81/4)z^2-(17d/16+17/2)z+289/256)/n^2+\cdots\end{align*}
\end{cf}

\smallskip

\begin{cf}\label{6.13.13}{\ }
\begin{verbatim}
[(z)->psi''(z),[[0,z^3],[(2*n+z-1)*(3*n^2+3*(z-1)*n+(z-1)^2),
                       (2*n+z)*(3*n^2+3*z*n+z^2)]],
             [[-2,-z^6],[-n^6,-(n+z)^6]]]
\end{verbatim}
$$\psi''(z)=-\dfrac{2}{z^3-\dfrac{z^6}{z^3+2z^2+2z+1-\dfrac{1}{z^3+5z^2+9z+6-\ddots}}}$$
Convergence type $E$ with $E=(1+\sqrt{2})^4$, $P=0$, $C=-8\pi^3/(1+\sqrt{2})^{4z}$, so that
$$\psi''(z)-\dfrac{p(n)}{q(n)}\sim-\dfrac{8\pi^3}{(1+\sqrt{2})^{4n+4z}}\;.$$
\begin{align*}A&=1+((2z^2-4z+17/16)d)/n\\
  &\phantom{=}+(4z^4-(2d+16)z^3+(4d+81/4)z^2-(17d/16+17/2)z+289/256)/n^2+\cdots\end{align*}
\end{cf}

\smallskip

\begin{cf}\label{6.13.14}{\ }
\begin{verbatim}
[a->sumnum(k=1,1/(k^3+a*k)),[0,(2*n-1)*(n^2-n+a+1)],
                            [1,-n^2*(n^2+a)^2]]
\end{verbatim}
$$\sum_{k=1}^\infty\dfrac{1}{k^3+ak}=\dfrac{1}{a+1-\dfrac{a^2+2a+1}{3a+9-\dfrac{4a^2+32a+64}{5a+35-\dfrac{9a^2+162a+729}{7a+91-\dfrac{16a^2+512a+4096}{9a+189-\ddots}}}}}$$
Convergence type $P^+$ with $P=2$ and $C=1/2$, so that
$$\sum_{k=1}^\infty\dfrac{1}{k^3+ak}-\dfrac{p(n)}{q(n)}\sim\dfrac{1/2}{n^2}\;.$$
$$A=1-1/n+(-a/2+1/2)/n^2+a/n^3+(a^2/3-5a/6-1/6)/n^4+\cdots$$
Series:
$$\sum_{k=1}^\infty\dfrac{1}{k^3+ak}=\sum_{n\ge1}\dfrac{1}{n(n^2+a)}$$
Parametric family for $k\ge0$:
\begin{verbatim}
[a->sumnum(k=1,1/(k^3+a*k)),(2*n-1)*(n^2-n+a+1+2*k^2+2*k),
                            -n^2*(n^2+a)^2]
\end{verbatim}
Convergence type $P^+$ with $P=4k+2$.
\end{cf}

\smallskip

\begin{cf}\label{6.13.15}{\ }
\begin{verbatim}
[a->sumnum(k=1,1/(k^3+a*k)),
[[0,a+1],[2*n*(a+3*n^2),(2*n+1)*(a+3*n^2+3*n+1)]],
[[1,-(a+1)^2],[-(n^3+a*n)^2,-((n+1)^3+a*(n+1))^2]]]
\end{verbatim}
$$\sum_{k=1}^\infty\dfrac{1}{k^3+ak}=\dfrac{1}{a+1-\dfrac{a^2+2a+1}{2a+6-\dfrac{a^2+2a+1}{3a+21-\dfrac{4a^2+32a+64}{4a+48-\dfrac{4a^2+32a+64}{5a+95-\dfrac{9a^2+162a+729}{6a+162+\ddots}}}}}}$$
Convergence type $E$ with $E=(1+\sqrt{2})^4$, $P=0$, and $C=...$, so that
$$\sum_{k=1}^\infty\dfrac{1}{k^3+ak}-\dfrac{p(n)}{q(n)}\sim\dfrac{C}{(1+\sqrt{2})^{4n}}\;.$$
$$A=1-((4a+15/16)d)/n+(16a^2+(4d+15/2)a+15d/16+225/256)/n^2+\cdots$$
\end{cf}

\medskip

\section{Function $\be_2(z)=\psi''(z)-\psi''(z/2)/4$}\label{sec:be2}

\medskip

In analogy to the CFs for the functions $\be(z)$ and $\be_1(z)$,
we could hope to find interesting CFs for
$$\be_2(z)=-\be_1'(z)=\psi''(z)-\psi''(z/2)/4=2\sum_{k\ge0}\dfrac{(-1)^k}{(z+k)^3}\;.$$
I have indeed found a family which is not really interesting, but which
I include anyway.

\smallskip

\begin{verbatim}
beta2(z)=psi''(z)-psi''(z/2)/4;
\end{verbatim}

\smallskip

\begin{cf}\label{6.14.1}{\ }
\begin{verbatim}
[(z)->beta2(z),[0,z^3,3*n^2+3*(2*z-3)*n+3*z^2-9*z+7],[2,(n+z-1)^6]]
\end{verbatim}
$$\be_2(z)=\dfrac{2}{z^3+\dfrac{z^6}{3z^2+3z+1+\dfrac{z^6+6z^5+15z^4+20z^3+15z^2+6z+1}{3z^2+9z+7+\ddots}}}$$
Convergence type $P^-$ with $P=3$ and $C=1$, so that
$$\be_2(z)-\dfrac{p(n)}{q(n)}\sim(-1)^n\dfrac{1}{n^3}\;.$$
$$A=1+(-3z+3/2)/n+(6z^2-6z)/n^2+(-10z^3+15z^2-5/2)/n^3+\cdots$$
Series:
$$\be_2(z)=2\sum_{n\ge1}\dfrac{(-1)^{n+1}}{(n+z-1)^3}$$
Parametric family of multiplied CF for $k\ge0$:
\begin{verbatim}
[(z)->beta2(z),(2*k+1)*(2*n+2*z-3)*(3*n^2+3*(2*z-3)*n+3*z^2-9*z+7),
             (n+z-1)^6*(4*(n+z-1)^2-(2*k+1)^2)]
\end{verbatim}
Convergence type $P^-$ with $P=6k+3$.
\end{cf}

Can be Ap\'ery accelerated with convergence type $E=8i$ and
$P=z-1$, formula too complicated to give here, dual also complicated.

\medskip

\section{Function $\psi'''(z)$}

\medskip

Recall that
$$\psi'''(z)=6\sum_{k\ge0}\dfrac{1}{(k+z)^4}\;.$$

Note also the following integrals:

\begin{align*}
\int_0^\infty\dfrac{t^3}{\sinh(t)}\,dt&=\psi'''((z+1)/2)/8\;,\\
\int_0^\infty t^3\cotanh(t)\,dt&=6/z^4+\psi'''(1+z/2)/8\;.\end{align*}

\medskip

\begin{cf}\label{6.15.1}{\ }
\begin{verbatim}
[(z)->psi'''(z),[0,(2*n-1)*(n^4-2*n^3-2*(z^2-z-1)*n^2+(2*z^2-2*z-1)*n
                -(z-1)*z*(z^2-z+1))],[4*z-2,-n^8*(n^2-(2*z-1)^2)]]
\end{verbatim}
$$\psi'''(z)\approx\dfrac{4z-2}{2z^4-4z^3+4z^2-2z+\dfrac{4z^2-4z}{6z^4-12z^3+24z^2-18z-18+\ddots}}$$
\end{cf}

This makes sense only as an \emph{asymptotic expansion} as $z\to\infty$. For
instance, if we denote by $w$ this continued fraction in {\tt GP} format,
\begin{verbatim}
? cftoser(subst(w,z,1/x))+O(x^20)
% = 2*x^3 + 3*x^4 + 2*x^5 - x^7 + 4/3*x^9 - 3*x^11 + 10*x^13
          - 691/15*x^15 + 280*x^17 - 10851/5*x^19 + O(x^20)
\end{verbatim}
shows the beginning of the asymptotic expansion of $\psi'''(1/x)$ for $x\to0$,
confirmed by the telltale presence of $691$ and of $10851=3\cdot 3617$.

The continued fraction itself is undefined (infinite) if $z$ is an
integer, and otherwise has a convergence type $P^+$ with $P=2$ to some
limit $f(z)$, so that
$$f(z)-\dfrac{p(n)}{q(n)}\sim\dfrac{C}{n^2}\;,$$
but I have no idea what $f(z)$ is.

\smallskip

{\bf Comment:} the above continued fraction comes from work done by G.~Rhin and
the author in 1980 on generalizing Ap\'ery's method to $\z(4)$, and is
a step in the proof of the continued fraction for $\z(4)=\pi^4/90$
given as \ref{5.8.1} found in the same paper. See Section \ref{sec:zeta4}
for its proof.

\medskip

\section{Ordinary Generating Functions for Bernoulli Type Sequences}

\smallskip

Note that here ``Ordinary'' is opposed to ``Exponential'', since exponential
generating functions for the sequences that we are going to mention are
exponential, hyperbolic, or trigonometric functions, whose CFs have been given
above in Sections \ref{sec:exp} and \ref{sec:trighyp}.

\smallskip

{\bf Definitions and Basic Properties}

\smallskip

For the convenience of the reader, and since some of the sequences are not
so well-known (and even possibly new), we recall their definitions and basic
properties.

\smallskip

{\bf Bernoulli Numbers $B_k$:}

\smallskip

\begin{enumerate}\item Exponential generating function:
  $$\dfrac{x}{e^x-1}=\sum_{k\ge0}\dfrac{B_k}{k!}x^k\;.$$
\item Vector of first few values:
  $$(B_k)_{k\ge0}=(1,-1/2,1/6,0,-1/30,0,1/42,0,-1/30,0,5/66,0,-691/2730,0,7/6,\dotsc)\;.$$
  \item L-value interpretation for $k\ge2$: $B_k=-k\z(1-k)$.
  \item {\tt Pari/GP} notation: {\tt bernfrac(k)}.
\end{enumerate}

\smallskip

{\bf Genocchi numbers $G_k$:}

\smallskip

\begin{enumerate}\item Exponential generating function:
  $$\dfrac{2x}{e^x+1}=\sum_{k\ge0}\dfrac{G_k}{k!}x^k\;.$$
  \item Vector of first few values:
$$(G_k)_{k\ge0}=(0,1,-1,0,1,0,-3,0,17,0,-155,0,2073,0,-38227,\dotsc)\;.$$
\item Relation with Bernoulli numbers: $G_k=-2(2^k-1)B_k$.
\item {\tt Pari/GP} notation:
\begin{verbatim}
genfrac(k)=-2*(2^k-1)*bernfrac(k);
\end{verbatim}
\end{enumerate}

\smallskip

{\bf Tangent numbers $T_k$:}

\smallskip

\begin{enumerate}\item Exponential generating function:
  $$\tanh(x)=\dfrac{e^{2x}-1}{e^{2x}+1}=\sum_{k\ge0}\dfrac{T_k}{k!}x^k\;.$$
\item Vector of first few values:
$$(T_k)_{k\ge0}=(0,1,0,-2,0,16,0,-272,0,7936,0,-353792,\dotsc)\;.$$
\item Relation to Bernoulli and Genocchi numbers for $k\ge1$:
$$T_k=2^{k+1}(2^{k+1}-1)\dfrac{B_{k+1}}{k+1}=-2^k\dfrac{G_{k+1}}{k+1}\;.$$
\item {\tt Pari/GP} notation:
\begin{verbatim}tanfrac(k)=if(k==0,0,2^(k+1)*(2^(k+1)-1)*bernfrac(k+1)/(k+1));
\end{verbatim}
\end{enumerate}

\smallskip

{\bf Euler numbers $E_k$:}

\smallskip

\begin{enumerate}\item Exponential generating function:
$$\dfrac{1}{\cosh(x)}=\dfrac{2}{e^x+e^{-x}}=\sum_{k\ge0}\dfrac{E_k}{k!}x^k\;.$$
\item Vector of first few values:
$$(E_k)_{k\ge0}=(1,0,-1,0,5,0,-61,0,1385,0,-50521,0,2702765,\dotsc)\;.$$
\item Relation to Bernoulli polynomials: $E_{2k+1}=0$ and
$$E_{2k}=-4^{2k+1}\dfrac{B_{2k+1}(1/4)}{2k+1}\;.$$
\item L-value interpretation: $E_k=2L(\chi_{-4},-k)$.
\item {\tt Pari/GP} notation: {\tt eulerfrac(k)}.
\end{enumerate}

The following are analogues modulo 3 of tangent and Euler numbers which are
modulo $4$. Although completely natural I have not found them in the literature
(not even in the OEIS), so I call them $T^{(3)}_k$ and $E^{(3)}_k$
respectively.

\smallskip

{\bf $3$-Tangent numbers $T^{(3)}_k$:}

\smallskip

\begin{enumerate}\item Exponential generating function:
  $$\dfrac{3\sinh(x)}{2\cosh(x)+1}=\dfrac{3}{2}\dfrac{e^{2x}-1}{e^{2x}+e^x+1}=\sum_{k\ge1}\dfrac{T^{(3)}_k}{k!}x^k\;.$$
\item Vector of first few values:
$$(T^{(3)}_k)_{k\ge0}=(0,1,0,-1,0,13/3,0,-41,0,671,0,-50443/3,\dotsc)\;.$$
\item Relation to Bernoulli numbers for $k\ge1$:
$$T^{(3)}_k=\dfrac{3}{2}(3^{k+1}-1)\dfrac{B_{k+1}}{k+1}\;.$$
\item {\tt Pari/GP} notation:
\begin{verbatim}tan3frac(k)=if(k==0,0,(3/2)*(3^(k+1)-1)*bernfrac(k+1)/(k+1));
\end{verbatim}
\end{enumerate}

It is possible that the ``correct'' normalization is without the facter $3/2$.

\smallskip

{\bf $3$-Euler numbers $E^{(3)}_k$:}

\smallskip

\begin{enumerate}\item Exponential generating function:
$$\dfrac{3}{2\cosh(x)+1}=\dfrac{3}{e^x+1+e^{-x}}=\sum_{k\ge0}\dfrac{E^{(3)}_k}{k!}x^k\;.$$
\item Vector of first few values:
$$(E^{(3)}_k)_{k\ge0}=(1,0,-2/3,0,2,0,-14,0,1618/9,0,-3694,\dotsc)\;.$$
\item Relation to Bernoulli polynomials: $E^{(3)}_{2k+1}=0$ and
$$E^{(3)}_{2k}=-2\cdot3^{2k+1}\dfrac{B_{2k+1}(1/3)}{2k+1}\;.$$
\item L-value interpretation: $E^{(3)}_k=3L(\chi_{-3},-k)$.
\item {\tt Pari/GP} notation: {\tt eul3frac(k)=3*lfun(-3,-k);}
\end{enumerate}

\smallskip

{\bf Springer numbers} $S_k$ (see \cite{Sok}):

\smallskip

\begin{enumerate}\item Exponential generating function:
$$\dfrac{1}{\cos(x)-\sin(x)}=\sum_{k\ge0}\dfrac{S_k}{k!}x^k\;.$$
\item Vector of first few values:
$$(S_k)_{k\ge0}=(1,1,3,11,57,361,2763,24611,250737,\dotsc)\;.$$
\item Bernoulli--Euler interpretation:
$$S_k=(-1)^{k(k-1)/2}\sum_{j=0}^k\binom{k}{j}2^jE_j\;.$$
\item {\tt Pari/GP} notation:
\begin{verbatim}
springfrac(k)=(-1)^(k*(k-1)/2)*
         sum(j=0,k\2,binomial(k,2*j)*2^(2*j)*eulerfrac(2*j));
\end{verbatim}
\end{enumerate}

\smallskip

{\bf Bernoulli--Hurwitz numbers} $H^{(w)}_k$ (see \cite{Coh6}):

\smallskip

I include these numbers in fond memory of Philippe Flajolet who
initiated the study of CFs associated to these numbers (I thank Wadim
Zudilin for reminding me of this). They require some
preliminaries. In what follows, $w=4$ or $w=6$, and $k\ge4$ is an even
integer.

For $\Im(\tau)>0$ we define
$G_k=\sum_{(m,n)\in\Z^2\setminus{(0,0)}}(m+n\tau)^{-k}$, and recall that
  if $q=e^{2\pi i\tau}$ we have
  $$G_{2k}(\tau)=(-1)^{k-1}\dfrac{(2\pi)^{2k}}{(2k)!}\left(B_{2k}-4k\sum_{n\ge1}\dfrac{n^{2k-1}q^n}{1-q^n}\right)\;.$$

  For $w=4$ and $w=6$ we define $\Om_w=B(1/w,2/w)$, where $B$ is the beta
  function, so that $\Om_4=\G(1/4)^2/(2\pi)^{1/2}$ and
  $\Om_6=3^{1/2}\G(1/3)^3/(2^{1/3}\pi)$. Finally, we define the
  Bernoulli--Hurwitz numbers $H^{(w)}_k$ by $H^{(w)}_k=0$ if $w\nmid k$,
  and otherwise
  $$H^{(w)}_k=(-1)^{k/w-1}\dfrac{k!}{\Om_w^k}G_k(\tau_w)\;,$$
  where $\tau_w=e^{2\pi i/w}$ (so $\tau_4=i$ and $\tau_6=(1+\sqrt{-3})/2$).

  \smallskip
  
  Vector of first few values:
\begin{align*}(H^{(4)}_{4k})_{k\ge1}&=(1/10,-3/10,567/130,-43659/170,392931/10,\dotsc)\;.\\
(H^{(6)}_{6k})_{k\ge1}&=(1/84,-25/1092,1375/1596,-257125/1092,\dotsc)\;.\end{align*}

We assume the existence of (easily written) {\tt GP} functions
{\tt h4frac(k)} and {\tt h6frac(k)} giving $H^{(4)}_k$ and $H^{(6)}_k$
respectively.

\bigskip

The following are formal power series expansions which converge only for
$z=0$, but whose corresponding CF expands to the same formal power series
and usually converge. They are included in this
section since they are closely linked to the expansions of $\psi(z)$ and
its derivatives. More precisely, changing $z$ into $1/z$ in each of the CF
given for those functions gives a CF related to generating functions of
the above numbers, so we only give samples.

An important warning, however. Consider for instance the very first CF
\ref{6.16.1} below, for $f(z)=\sum_{k\ge0}B_kz^k$. This CF means that if
we expand formally the CF as a power series in $z$, we will find the formal
expansion of $f(z)$, i.e. $f(z)=1-z/2+z^2/6-z^4/30+\cdots$, which has
radius of convergence $0$. Note that this is possible because all the
$b(n)$ are divisible by $z$ for $n\ge1$. Nonetheless the CF itself converges
for $z\ne1/2$, giving a value to $f(z)$, even though its series diverges
(according to \ref{6.11.5}, we have $f(z)=-\psi'(-1/z)/z$ when $-2<z<0$ and
$f(z)=\psi'(1+1/z)/z$ when $z>0$ or $z<-2$). Thus, subsequent CFs for
the same function, for instance \ref{6.16.2}, which do \emph{not} have
$b(n)$ divisible by $z$ so which cannot be expanded into formal power series
in $z$, are such that in their domain of convergence they converge to the same
function $f(z)$.

\smallskip

\begin{cf}\label{6.16.1}{\ }
\begin{verbatim}
[(z)->sum(k=0,oo,bernfrac(k)*z^k),[0,(2*n-1)*(z+2)],[2,n^4*z^2]]
\end{verbatim}
$$\sum_{k\ge0}B_kz^k=\dfrac{2}{z+2+\dfrac{z^2}{3z+6+\dfrac{16z^2}{5z+10+\dfrac{81z^2}{7z+14+\dfrac{256z^2}{9z+18+\dfrac{625z^2}{11z+22+\ddots}}}}}}$$
Convergence type $P^-$ with $P=|2+4/z|$ and $C=\G(1+1/z)^4/z$, so that
$$\sum_{k\ge0}B_kz^k-\dfrac{p(n)}{q(n)}\sim(-1)^n\dfrac{\G(1+1/z)^4/z}{n^{|2+4/z|}}\;.$$
$$A=1-((z+2)/z)/n+((z^2-4)/(3z^3))/n^2+\cdots$$
Assume here that $z\ge-1$.
Parametric family for $k\ge0$:
\begin{verbatim}
[(z)->sum(k=0,oo,bernfrac(k)*z^k),(2*n-1)*((2*k+1)*z+2),n^4*z^2]
\end{verbatim}
Convergence type $P^-$ with $P=4k+2+4/z$.
\end{cf}

\smallskip

\begin{cf}\label{6.16.2}{\ }
\begin{verbatim}
[(z)->sum(k=0,oo,bernfrac(k)*z^k),
[0,(z+1)^2,2*z^2*n^2-(2*z-4)*z*n+z^2-2*z+2],[z,-(z*n+1)^4]]
\end{verbatim}
$$\sum_{k\ge0}B_kz^k=\dfrac{z}{z^2+2z+1-\dfrac{z^4+4z^3+6z^2+4z+1}{5z^2+6z+2-\dfrac{16z^4+32z^3+24z^2+8z+1}{13z^2+10z+2-\ddots}}}$$
Convergence type $P^+$ with $P=1$ and $C=1/z$, so that
$$\sum_{k\ge0}B_kz^k-\dfrac{p(n)}{q(n)}\sim\dfrac{1/z}{n}\;.$$
$$A=1-((z+2)/(2z))/n+((z^2+6z+6)/(6z^2))/n^2+\cdots$$
Series:
$$\sum_{k\ge0}B_kz^k=z\sum_{n\ge1}\dfrac{1}{(nz+1)^2}$$
Parametric family for $k\ge0$:
\begin{verbatim}
[(z)->sum(k=0,oo,bernfrac(k)*z^k),
2*z^2*n^2-(2*z-4)*z*n+(k^2+k+1)*z^2-2*z+2,-(z*n+1)^4]
\end{verbatim}
Convergence type $P^+$ with $P=2k+1$.
\end{cf}
        
\smallskip

\begin{cf}\label{6.16.3}{\ }
\begin{verbatim}
[(z)->sum(k=0,oo,bernfrac(k)*z^k),
[[0,2*(z+1)^2],[5*z^2*n^2+6*z*n+2,5*z^2*n^2+6*z*(z+1)*n+2*(z+1)^2]],
[[2*z,-4*(z+1)^4],[z^4*n^4,-4*(z*n+z+1)^4]]]
\end{verbatim}
$$\sum_{k\ge0}B_kz^k=\dfrac{2z}{2z^2+4z+2-\dfrac{4z^4+16z^3+24z^2+16z+4}{5z^2+6z+2+\dfrac{z^4}{13z^2+10z+2-\dfrac{64z^4+128z^3+96z^2+32z+4}{20z^2+12z+2+\ddots}}}}$$
Convergence type $E$ with $E=-((1+\sqrt{5})/2)^5$, $P=0$, and
$C=4\pi^2/z/((1+\sqrt{5})/2)^{4/z+5}$, so that
$$\sum_{k\ge0}B_kz^k-\dfrac{p(n)}{q(n)}\sim(-1)^n\dfrac{4\pi^2/z}{((1+\sqrt{5})/2)^{5n+4/z+5}}\;.$$
$$A=1-((2z^2-4)d/(5z^2))/n+(((10d+10)z^4+8dz^3-(20d+40)z^2-24dz+40)/(25z^4))/n^2+\cdots$$
\end{cf}

\smallskip

\begin{cf}\label{6.16.4}{\ }
\begin{verbatim}
[(z)->sum(k=0,oo,bernfrac(k)*z^k),[1-z/2,4*n+2],z^2*[1,n*(n+1)^2*(n+2)]]
\end{verbatim}
$$\sum_{k\ge0}B_kz^k=1-z/2+\dfrac{z^2}{6+\dfrac{12z^2}{10+\dfrac{72z^2}{14+\dfrac{240z^2}{18+\dfrac{600z^2}{22+\dfrac{1260z^2}{26+\ddots}}}}}}$$
Convergence type $P^-$ with $P=4/|z|$ and $C=z\G(1+1/z)^4$, so that
$$\sum_{k\ge0}B_kz^k-\dfrac{p(n)}{q(n)}\sim(-1)^n\dfrac{z\G(1+1/z)^4}{n^{4/|z|}}\;.$$
$$A=1-(6/z)/n+((13z^2+54z-4)/(3z^3))/n^2+\cdots$$
Parametric family for $k\ge0$:
\begin{verbatim}
[(z)->sum(k=0,oo,bernfrac(k)*z^k),(k*z+1)*(4*n+2),z^2*n*(n+1)^2*(n+2)]
\end{verbatim}
Convergence type $P^-$ with $P=4k+4/z$.
\end{cf}

\smallskip

\begin{cf}\label{6.16.5}{\ }
\begin{verbatim}
[(z)->sum(k=0,oo,bernfrac(k)*z^k),
[1-z/2,(z+1)^2,2*(z^2*n^2+2*z*(1-z)*n+2*z^2-2*z+1)],
[z^3/2,-(z*n+1)^2*(z*n+1-z)^2]]
\end{verbatim}
$$\sum_{k\ge0}B_kz^k=-z/2+1+\dfrac{z^3/2}{z^2+2z+1-\dfrac{z^2+2z+1}{4z^2+4z+2-\dfrac{4z^4+12z^3+13z^2+6z+1}{10z^2+8z+2-\ddots}}}$$
Convergence type $P^+$ with $P=3$ and $C=1/(6z)$, so that
$$\sum_{k\ge0}B_kz^k-\dfrac{p(n)}{q(n)}\sim\dfrac{1/(6z)}{n^3}\;.$$
$$A=1-(3/z)/n+((-z^2+30)/(5z^2))/n^2+\cdots$$
Series:
$$\sum_{k\ge0}B_kz^k=-\dfrac{z}{2}+1+\dfrac{z^3}{2}\sum_{n\ge1}\dfrac{1}{(nz+1)^2((n-1)z+1)^2}$$
Parametric family for $k\ge0$:
\begin{verbatim}
[(z)->sum(k=0,oo,bernfrac(k)*z^k),
2*z^2*n^2-4*z*(z-1)*n+(k^2+3*k+4)*z^2-4*z+2,-(n*z+1)^2*((n-1)*z+1)^2]
\end{verbatim}
Convergence type $P^+$ with $P=2k+3$.
\end{cf}

\smallskip

\begin{cf}\label{6.16.6}{\ }
\begin{verbatim}
[(z)->sum(k=0,oo,bernfrac(k)*z^k),
[[1-z/2,2*(z+1)^2],[5*z^2*n^2+z*(z+6)*n+2,
                    5*z^2*n^2+z*(7*z+6)*n+2*(z+1)^2]],
[[z^3,-4*(z+1)^2],[z^4*n*(n+1)^2*(n+2),-4*(z*n+1)^2*(z*(n+1)+1)^2]]]
\end{verbatim}
$$\sum_{k\ge0}B_kz^k=1-z/2+\dfrac{z^3}{2z^2+4z+2-\dfrac{4z^2+8z+4}{6z^2+6z+2+\dfrac{12z^4}{14z^2+10z+2-\ddots}}}$$
Convergence type $E$ with $E=-((1+\sqrt{5})/2)^5$, $P=0$, and
$C=4\pi^2/z/((1+\sqrt{5})/2)^{4/z+9}$, so that
$$\sum_{k\ge0}B_kz^k-\dfrac{p(n)}{q(n)}\sim(-1)^n\dfrac{4\pi^2/z}{((1+\sqrt{5})/2)^{5n+4/z+9}}\;.$$
$$A=1+((12z^2-12z+4)d/(5z^2))/n+\cdots$$
\end{cf}

\smallskip

\begin{cf}\label{6.16.7}{\ }
\begin{verbatim}
[(z)->sum(k=0,oo,(bernfrac(k)-bernfrac(k+2))*z^k),
    [5/6-z/2,2*n+3],z^2*[1,n*(n+1)*(n+3)*(n+4)/4]]
\end{verbatim}
$$\sum_{k\ge0}(B_k-B_{k+2})z^k=5/6-\dfrac{z}{2}+\dfrac{z^2}{5+\dfrac{10z^2}{7+\dfrac{45z^2}{9+\dfrac{126z^2}{11+\dfrac{280z^2}{13+\dfrac{540z^2}{15+\ddots}}}}}}$$
Convergence type $P^-$ with $P=4/|z|$ and $C=\G(1/z)\G(1+1/z)\G(2+1/z)^2$,
so that
$$\sum_{k\ge0}(B_k-B_{k+2})z^k-\dfrac{p(n)}{q(n)}\sim(-1)^n\dfrac{\G(1/z)\G(1+1/z)\G(2+1/z)^2}{n^{4/|z|}}\;.$$
$$A=1-(10/z)/n+((49z^2+150z-4)/(3z^3))/n^2+\cdots$$
Parametric family for $k\ge0$:
\begin{verbatim}
[(z)->sum(k=0,oo,(bernfrac(k)-bernfrac(k+2))*z^k),
(k*z+1)*(2*n+3),z^2*n*(n+1)*(n+3)*(n+4)/4]
\end{verbatim}
Convergence type $P^-$ with $P=4k+4/z$.
\end{cf}

\smallskip

Can be Ap\'ery accelerated with convergence type $E=-((1+\sqrt{5})/2)^5$,
formula too complicated to give here.

\smallskip

\begin{cf}\label{6.16.8}{\ }
\begin{verbatim}
[(z)->sum(k=0,oo,(bernfrac(k)-bernfrac(k+2))*z^k),
[(5-3*z)/6,(z+1)*(2*z+1),2*z^2*n^2-4*z*(z-1)*n+6*z^2-4*z+2],
[z^3,-((n-2)*z+1)*((n-1)*z+1)*(n*z+1)*((n+1)*z+1)]]
\end{verbatim}
$$\sum_{k\ge0}(B_k-B_{k+2})z^k=5/6-z/2+\dfrac{z^3}{2z^2+3z+1+\dfrac{2z^3+z^2-2z-1}{6z^2+4z+2-\dfrac{6z^3+11z^2+6z+1}{12z^2+8z+2-\ddots}}}$$
Convergence type $P^+$ with $P=5$ and $C=(2z+1)(z+1)/z^3$, so that
$$\sum_{k\ge0}(B_k-B_{k+2})z^k-\dfrac{p(n)}{q(n)}\sim\dfrac{(2z+1)(z+1)/z^3}{n^5}\;.$$
$$A=1-(5/z)/n+((5z^2+105)/(7z^2))/n^2-((5z^2+35)/z^3)/n^3+\cdots$$
Series:
$$\sum_{k\ge0}(B_k-B_{k+2})z^k=-\dfrac{z}{2}+\dfrac{5}{6}+z^3(1-z^2)\sum_{n\ge0}\dfrac{1}{((n+2)z+1)((n+1)z+1)^2(nz+1)^2((n-1)z+1)}$$
Parametric family for $k\ge0$:
\begin{verbatim}
[(z)->sum(k=0,oo,(bernfrac(k)-bernfrac(k+2))*z^k),
2*z^2*n^2-4*z*(z-1)*n+(k^2+5*k+6)*z^2-4*z+2,
-((n-2)*z+1)*((n-1)*z+1)*(n*z+1)*((n+1)*z+1)]
\end{verbatim}
Convergence type $P^+$ with $P=2k+5$.
\end{cf}

\smallskip

\begin{cf}\label{6.16.9}{\ }
\begin{verbatim}
[(z)->sum(k=0,oo,(bernfrac(k)-bernfrac(k+2))*z^k),
[5/6,2*(z+4)*n*(n+2)+6*(z+1)],[-15*z,z^2*n^2*(n+3)^2*(2*n+1)*(2*n+5)]]
\end{verbatim}
$$\sum_{k\ge0}(B_k-B_{k+2})z^k=5/6-\dfrac{15z}{12z+30+\dfrac{336z^2}{22z+70+\dfrac{4500z^2}{36z+126+\dfrac{24948z^2}{54z+198+\ddots}}}}$$
Convergence type $P^-$ with $P=1+4/|z|$ and $C=-\G(1+1/z)^2\G(2+1/z)^2$,
so that
$$\sum_{k\ge0}(B_k-B_{k+2})z^k-\dfrac{p(n)}{q(n)}\sim(-1)^{n+1}\dfrac{\G(1+1/z)^2\G(2+1/z)^2}{n^{1+4/|z|}}\;.$$
$$A=1-((2z+8)/z)/n+((18z^3+82z^2+93z-4)/(3z^3))/n^2+\cdots$$
\end{cf}

\smallskip

By using the formula
$$\sum_{k\ge0}(B_k-B_{k+2})z^k=1/z^2-1/(2z)+(1-1/z^2)\sum_{k\ge0}B_kz^k$$
one can obtain from these CFs other CFs for $\sum_{k\ge0}B_kz^k$ which do not
seem worth recording.

\smallskip

\begin{cf}\label{6.16.10}{\ }
\begin{verbatim}
[(z)->sum(k=0,oo,(2^(k-1)-1)*bernfrac(k)*z^k),[0,2*n-1],[-1/2,z^2*n^4]]
\end{verbatim}
$$\sum_{k\ge0}(2^{k-1}-1)B_kz^k=-\dfrac{1/2}{1+\dfrac{z^2}{3+\dfrac{16z^2}{5+\dfrac{81z^2}{7+\dfrac{256z^2}{9+\dfrac{625z^2}{11+\ddots}}}}}}$$
Convergence type $P^-$ with $P=2/|z|$ and $C=-\G((1+1/z)/2)/(4z)$, so that
$$\sum_{k\ge0}(2^{k-1}-1)B_kz^k-\dfrac{p(n)}{q(n)}\sim(-1)^{n+1}\dfrac{\G((1+1/k)/2)/(4k)}{n^{2/|z|}}\;.$$
$$A=1-(1/z)/n-((z-1)(2z-1)/(6z^3))/n^2+((z+1)(3z^2-z+1)/(6z^4))/n^3+\cdots$$
Parametric family for $k\ge0$:
\begin{verbatim}
[(z)->sum(k=0,oo,(2^(k-1)-1)*bernfrac(k)*z^k),(2*k*z+1)*(2*n-1),z^2*n^4]
\end{verbatim}
Convergence type $P^-$ with $P=4k+2/z$.
\end{cf}

\smallskip

\begin{cf}\label{6.16.11}{\ }
\begin{verbatim}
[(z)->sum(k=0,oo,(2^(k-1)-1)*bernfrac(k)*z^k),
[0,(z+1)^2,8*z^2*n^2-8*z*(2*z-1)*n+10*z^2-8*z+2],[-z,-(z*(2*n-1)+1)^4]]
\end{verbatim}
$$\sum_{k\ge0}(2^{k-1}-1)B_kz^k=-\dfrac{z}{z^2+2z+1-\dfrac{z^4+4z^3+6z^2+4z+1}{10z^2+8z+2-\dfrac{81z^4+108z^3+54z^2+12z+1}{34z^2+16z+2-\ddots}}}$$
Convergence type $P^+$ with $P=1$ and $C=-1/(4z)$, so that
$$\sum_{k\ge0}(2^{k-1}-1)B_kz^k-\dfrac{p(n)}{q(n)}\sim-\dfrac{1/(4z)}{n}\;.$$
$$A=1-(1/(2z))/n+((-z^2+3)/(12z^2))/n^2+((z^2-1)/(8z^3))/n^3+\cdots$$
Series:
$$\sum_{k\ge0}(2^{k-1}-1)B_kz^k=-z\sum_{n\ge1}\dfrac{1}{((2n-1)z+1)^2}$$
Parametric family for $k\ge0$:
\begin{verbatim}
[(z)->sum(k=0,oo,(2^(k-1)-1)*bernfrac(k)*z^k),
8*z^2*n^2-8*z*(2*z-1)*n+(4*k^2+4*k+10)*z^2-8*z+2,-(z*(2*n-1)+1)^4]
\end{verbatim}
Convergence type $P^+$ with $P=2k+1$.
\end{cf}
      
\smallskip

\begin{cf}\label{6.16.12}{\ }
\begin{verbatim}
[(z)->sum(k=0,oo,(2^(k-1)-1)*bernfrac(k)*z^k),
  [[0,(z+1)^2],[10*z^2*n^2-6*z*(z-1)*n+(z-1)^2,
                10*z^2*n^2+6*z*(z+1)*n+(z+1)^2]],
  [[-z,-(z+1)^4],[4*z^4*n^4,-(z*(2*n+1)+1)^4]]]
\end{verbatim}
$$\sum_{k\ge0}(2^{k-1}-1)B_kz^k=-\dfrac{z}{z^2+2z+1-\dfrac{z^4+4z^3+6z^2+4z+1}{5z^2+4z+1+\dfrac{4z^4}{17z^2+8z+1-\ddots}}}$$
Convergence type $E$ with $E=-((1+\sqrt{5})/2)^5$, $P=0$, and
$C=-\pi^2/z/((1+\sqrt{5})/2)^{2/z+3}$, so that
$$\sum_{k\ge0}(2^{k-1}-1)B_kz^k-\dfrac{p(n)}{q(n)}\sim(-1)^{n+1}\dfrac{\pi^2/z}{((1+\sqrt{5})/2)^{5n+2/z+3}}\;.$$
$$A=1-((z^2+2z-1)d/(5z^2))/n+\cdots$$
\end{cf}

\smallskip

\begin{cf}\label{6.16.13}{\ }
\begin{verbatim}
[(z)->sum(k=0,oo,genfrac(k)*z^k),[z,(2*n^2-n)*z^2+1],
                               [-z^2,-z^4*n^3*(n+1)]]
\end{verbatim}
$$\sum_{k\ge0}G_kz^k=z-\dfrac{z^2}{z^2+1-\dfrac{2z^4}{6z^2+1-\dfrac{24z^4}{15z^2+1-\dfrac{108z^4}{28z^2+1-\dfrac{320z^4}{45z^2+1-\ddots}}}}}$$
Convergence type $P^+$ with $P=|2/z|$ and $C=-\G(1/z)^6/(2z^3\G(2/z)^2)$,
so that
$$\sum_{k\ge0}G_kz^k-\dfrac{p(n)}{q(n)}\sim-\dfrac{\G(1/z)^6/(2z^3\G(2/z)^2)}{n^{|2/z|}}\;.$$
$$A=1-((z^2-6)/(z^3-4z))/n+((z^4-z^3-15z^2+28z-1)/(6z^3(z-2)^2))/n^2+\cdots$$
Parametric family for $k\ge0$:
\begin{verbatim}
[(z)->sum(k=0,oo,genfrac(k)*z^k),(2*n^2-n)*z^2+(k*z+1)^2,-z^4*n^3*(n+1)]
\end{verbatim}
Convergence type $P^+$ with $P=2k+2/z$.
\end{cf}

\smallskip

\begin{cf}\label{6.16.14}{\ }
\begin{verbatim}
[(z)->sum(k=0,oo,genfrac(k)*z^k),[0,z^2*(2*n^2-2*n+1)+1+z],[z,-z^4*n^4]]
\end{verbatim}
$$\sum_{k\ge0}G_kz^k=\dfrac{z}{z^2+z+1-\dfrac{z^4}{5z^2+z+1-\dfrac{16z^4}{13z^2+z+1-\dfrac{81z^4}{25z^2+z+1-\ddots}}}}$$
Convergence type $P^+$ with $P=|2/z+1|$ and $C=\G(1/z)^6/(4z^4(z+2)\G(2/z)^2)$,
so that
$$\sum_{k\ge0}G_kz^k-\dfrac{p(n)}{q(n)}\sim\dfrac{\G(1/z)^6/(4z^4(z+2)\G(2/z)^2)}{n^{|2/z+1|}}\;.$$
$$A=1-((z+2)/(2z))/n+((z+2)(6z^4-20z^2-10z+1)/(12z^3(z-2)(3z+2)))/n^2+\cdots$$
Parametric family for $k\ge0$:
\begin{verbatim}
[(z)->sum(k=0,oo,genfrac(k)*z^k),
z^2*(2*n^2-2*n)+(k^2+k+1)*z^2+(2*k+1)*z+1,-z^4*n^4]
\end{verbatim}
Convergence type $P^+$ with $P=2k+1+2/z$.
\end{cf}

\begin{cf}\label{6.16.15}{\ }
\begin{verbatim}
[(z)->sum(k=0,oo,genfrac(k)*z^k),
[0,(z+1)^2,z^2*(2*n-1)+2*z],[2*z,(n*z+1)^4]]
\end{verbatim}
$$\sum_{k\ge0}G_kz^k=\dfrac{2z}{z^2+2z+1+\dfrac{z^4+4z^3+6z^2+4z+1}{3z^2+2z+\dfrac{16z^4+32z^3+24z^2+8z+1}{5z^2+2z+\ddots}}}$$
Convergence type $P^-$ with $P=2$ and $C=1/z$, so that
$$\sum_{k\ge0}G_kz^k-\dfrac{p(n)}{q(n)}\sim\dfrac{1/z}{n^2}\;.$$
$$A=1-((z+2)/z)/n+((3z+3)/z^2)/n^2+((z^3-6z-4)/z^3)/n^3+\cdots$$
Series:
$$\sum_{k\ge0}G_kz^k=2z\sum_{n\ge1}\dfrac{(-1)^{n+1}}{(nz+1)^2}$$
Parametric family for $k\ge0$:
\begin{verbatim}
[(z)->sum(k=0,oo,genfrac(k)*z^k),(2*k+1)*(z^2*(2*n-1)+2*z),(n*z+1)^4]
\end{verbatim}
Convergence type $P^-$ with $P=4k+2$.
\end{cf}
      
\smallskip

\begin{cf}\label{6.16.16}{\ }
\begin{verbatim}
[(z)->sum(k=0,oo,genfrac(k)*z^k),
[[2*z,1],[5*z^2*n^2-4*z*(z-1)*n+(z-1)^2,5*z^2*n^2+4*z*n+1]],
[[-2*z,1],[-4*z^4*n^4,(z*n+1)^4]]]
\end{verbatim}
$$\sum_{k\ge0}G_kz^k=2z-\dfrac{2z}{1+\dfrac{1}{2z^2+2z+1-\dfrac{4z^4}{5z^2+4z+1+\dfrac{z^4+4z^3+6z^2+4z+1}{13z^2+6z+1-\ddots}}}}$$
Convergence type $E$ with $E=-((1+\sqrt{5})/2)^5$, $P=0$, and
$C=-4\pi^2/z/((1+\sqrt{5})/2)^{6/z-1}$, so that
$$\sum_{k\ge0}G_kz^k-\dfrac{p(n)}{q(n)}\sim(-1)^{n+1}\dfrac{4\pi^2/z}{((1+\sqrt{5})/2)^{5n+6/z-1}}\;.$$
$$A=1+((2z^2-8z+4)d/(5z^2))/n+\cdots$$
\end{cf}

\smallskip

\begin{cf}\label{6.16.17}{\ }
\begin{verbatim}
[(z)->sum(k=0,oo,(genfrac(k)-genfrac(k+2))*z^k),
[[0,1-z],[n*(3*n+2)*z+2*n+1,-3*z*n^2-(4*z-2)*n+1-z]],
[[1,9*z^2],z^2*[n^2*(n+1)^2,(2*n+1)^2*(2*n+3)^2]]]
\end{verbatim}
$$\sum_{k\ge0}(G_k-G_{k+2})z^k=\dfrac{1}{-z+1+\dfrac{9z^2}{5z+3+\dfrac{4z^2}{-8z+3+\dfrac{225z^2}{16z+5+\dfrac{36z^2}{-21z+5+\ddots}}}}}$$
Convergence type $E$ with $E=-2$, $P=-1/z$, and $C=3\cdot2^{2/z+1}\G(1+1/z)^4\G(2+1/z)^2/((z+2)\G(1+2/z)^2)$, so that
$$\sum_{k\ge0}(G_k-G_{k+2})z^k-\dfrac{p(n)}{q(n)}\sim(-1)^n\dfrac{3\G(1+1/z)^4\G(2+1/z)^2/((z+2)\G(1+2/z)^2)}{2^{n-2/z-1}n^{-1/z}}\;.$$
$$A=1-(7(z+2)/(3z))/n+\cdots$$
\end{cf}

\medskip

I have not found any reasonable CF for $\exp(\pm\sum_{k\ge1}(B_k/k)z^k)$
or for $\exp(\pm\sum_{k\ge1}(G_k/k)z^k)$.

\smallskip

\begin{cf}\label{6.16.18}{\ }
\begin{verbatim}
[(z)->sum(k=0,oo,tanfrac(k)*z^k),[0,1],[z,n*(n+1)*z^2]]
\end{verbatim}
$$\sum_{k\ge0}T_kz^k=\dfrac{z}{1+\dfrac{2z^2}{1+\dfrac{6z^2}{1+\dfrac{12z^2}{1+\dfrac{20z^2}{1+\dfrac{30z^2}{1+\ddots}}}}}}$$
Convergence type $P^-$ with $P=1/|z|$ and $C=\G(1+2/z)^2/2^{1/z-1}$, so that
$$\sum_{k\ge0}T_kz^k-\dfrac{p(n)}{q(n)}\sim(-1)^n\dfrac{\G(1+2/z)^2/2^{1/z-1}}{n^{1/|z|}}\;.$$
$$A=1-(1/z)/n+((16z^2+24z-1)/(48z^3))/n^2+\cdots$$
Parametric family for $k\ge0$:
\begin{verbatim}
[(z)->sum(k=0,oo,tanfrac(k)*z^k),2*k*z+1,n*(n+1)*z^2]
\end{verbatim}
Convergence type $P^-$ with $P=2k+1/z$.
\end{cf}

\smallskip

\begin{cf}\label{6.16.19}{\ }
\begin{verbatim}
[(z)->sum(k=0,oo,tanfrac(k)*z^k),
[0,2*(2*n-1)^2*z^2+1],[z,-4*n^2*(4*n^2-1)*z^4]]
\end{verbatim}
$$\sum_{k\ge0}T_kz^k=\dfrac{z}{2z^2+1-\dfrac{12z^4}{18z^2+1-\dfrac{240z^4}{50z^2+1-\dfrac{1260z^4}{98z^2+1-\ddots}}}}$$
Convergence type $P^+$ with $P=1/|z|$ and $C=\G(1+2/z)^2/2^{2/z-1}$, so that
$$\sum_{k\ge0}T_kz^k-\dfrac{p(n)}{q(n)}\sim\dfrac{\G(1+2/z)^2/2^{2/z-1}}{n^{1/|z|}}\;.$$
$$A=1-(1/(2z))/n+((16z^2+24z-1)/(192z^3))/n^2+\cdots$$
Parametric family for $k\ge0$:
\begin{verbatim}
[(z)->sum(k=0,oo,tanfrac(k)*z^k),
8*z^2*(n^2-n)+(4*k^2+2)*z^2+4*k*z+1,-4*n^2*(4*n^2-1)*z^4]
\end{verbatim}
Convergence type $P^+$ with $P=2k+1/z$.
\end{cf}

This is simply the contraction of the previous CF.

\smallskip

\begin{cf}\label{6.16.20}{\ }
\begin{verbatim}
[(z)->sum(k=0,oo,tanfrac(k)*z^k),
[0,2*z+1,4*z],[2*z,(2*z*n+1)*(2*z*(n-1)+1)]]
\end{verbatim}
$$\sum_{k\ge0}T_kz^k=\dfrac{2z}{2z+1+\dfrac{2z+1}{4z+\dfrac{8z^2+6z+1}{4z+\dfrac{24z^2+10z+1}{4z+\dfrac{48z^2+14z+1}{4z+\ddots}}}}}$$
Convergence type $P^-$ with $P=2$ and $C=1/(4z)$, so that
$$\sum_{k\ge0}T_kz^k-\dfrac{p(n)}{q(n)}\sim(-1)^n\dfrac{1/(4z)}{n^2}\;.$$
$$A=1-(1/z)/n+((-2z^2+3)/(4z^2))/n^2+((2z^2-1)/(2z^3))/n^3+\cdots$$
Series:
$$\sum_{k\ge0}T_kz^k=2z\sum_{n\ge1}\dfrac{(-1)^{n+1}}{(2nz+1)((2n-2)z+1)}=1-2\sum_{n\ge1}\dfrac{(-1)^{n+1}}{2nz+1}$$
Parametric family for $k\ge0$:
\begin{verbatim}
[(z)->sum(k=0,oo,tanfrac(k)*z^k),4*(k+1)*z,(2*z*n+1)*(2*z*(n-1)+1)]
\end{verbatim}
Convergence type $P^-$ with $P=2k+2$.
\end{cf}

\smallskip

\begin{cf}\label{6.16.21}{\ }
\begin{verbatim}
[(z)->sum(k=0,oo,tanfrac(k)*z^k),
[1,2*z+2,2*n*z+1],[[-2,4*z^2],[(2*z*n+1)^2,4*(n+1)^2*z^2]]]
\end{verbatim}
$$\sum_{k\ge0}T_kz^k=1-\dfrac{2}{2z+2+\dfrac{4z^2}{4z+1+\dfrac{4z^2+4z+1}{6z+1+\dfrac{16z^2}{8z+1+\dfrac{16z^2+8z+1}{10z+1+\ddots}}}}}$$
Convergence type $E$ with $E=-(1+\sqrt{2})^2$, $P=0$, and
$C=-2\pi/z/(1+\sqrt{2})^{2+1/z}$, so that
$$\sum_{k\ge0}T_kz^k-\dfrac{p(n)}{q(n)}\sim(-1)^{n+1}\dfrac{2\pi/z}{(1+\sqrt{2})^{2n+1/z+2}}\;.$$
$$A=1-((3z^2-2)d/(8z^2))/n+\cdots$$
\end{cf}

\smallskip

\begin{cf}\label{6.16.22}{\ }
\begin{verbatim}
[(z)->sum(k=0,oo,tanfrac(k)*z^k),
[0,2*z*n+1],[[2*z,2*z+1],[4*z^2*n*(n+1),(2*z*n+1)*(2*z*(n+1)+1)]]]
\end{verbatim}
$$\sum_{k\ge0}T_kz^k=\dfrac{2z}{2z+1+\dfrac{2z+1}{4z+1+\dfrac{8z^2}{6z+1+\dfrac{8z^2+6z+1}{8z+1+\dfrac{24z^2}{10z+1+\dfrac{24z^2+10z+1}{12z+1+\ddots}}}}}}$$
Convergence type $E$ with $E=-(1+\sqrt{2})^2$, $P=0$, and
$C=2\pi/z/(1+\sqrt{2})^{2+1/z}$, so that
$$\sum_{k\ge0}T_kz^k-\dfrac{p(n)}{q(n)}\sim(-1)^n\dfrac{2\pi/z/(1+\sqrt{2})^{2+1/z}}{(1+\sqrt{2})^{2n}}\;.$$
$$A=1+((9z^2-8z+2)d/(8z^2))/n+\cdots$$
\end{cf}

\smallskip

\begin{cf}\label{6.16.23}{\ }
\begin{verbatim}
[(z)->sum(k=0,oo,(-1)^(k*(k-1)/2)*tanfrac(k)*z^k),
                                [0,1],[z,-n*(n+1)*z^2]]
\end{verbatim}
$$\sum_{k\ge0}(-1)^{k(k-1)/2}T_kz^k=\dfrac{z}{1-\dfrac{2z^2}{1-\dfrac{6z^2}{1-\dfrac{12z^2}{1-\dfrac{20z^2}{1-\dfrac{30z^2}{1-\ddots}}}}}}$$
Nonconvergent CF.
\end{cf}

\smallskip

\begin{cf}\label{6.16.24}{\ }
\begin{verbatim}
[(z)->sum(k=0,oo,(tanfrac(k)-tanfrac(k+2))*z^k),
[0,(8*n^2-2)*z^2+1],[3*z,-4*n^2*(2*n+1)*(2*n+3)*z^4]]
\end{verbatim}
$$\sum_{k\ge0}(T_k-T_{k+2})z^k=\dfrac{3z}{6z^2+1-\dfrac{60z^4}{30z^2+1-\dfrac{560z^4}{70z^2+1-\dfrac{2268z^4}{126z^2+1-\ddots}}}}$$
Convergence type $P^+$ with $P=1/|z|$ and $C=\G(1+1/(2z))^2(1+1/z)^2/2^{2/z-1}$,
so that
$$\sum_{k\ge0}(T_k-T_{k+2})z^k-\dfrac{p(n)}{q(n)}\sim\dfrac{\G(1+1/(2z))^2(1+1/z)^2/2^{2/z-1}}{n^{1/|z|}}\;.$$
$$A=1-((z^2-2)/(2(z^3-z)))/n+((16z^4-8z^3-93z^2+98z-1)/(192z^3(z-1)^2))/n^2+\cdots$$
Parametric family for $k\ge0$:
\begin{verbatim}
[(z)->sum(k=0,oo,(tanfrac(k)-tanfrac(k+2))*z^k),
(8*n^2+4*k^2-2)*z^2+4*k*z+1,-4*n^2*(2*n+1)*(2*n+3)*z^4]
\end{verbatim}
Convergence type $P^+$ with $P=2k+1/z$.
\end{cf}

\smallskip

\begin{cf}\label{6.16.25}{\ }
\begin{verbatim}
[(z)->sum(k=0,oo,(tanfrac(k)-tanfrac(k+2))*z^k),
[1,(8*n^2-4*n+2)*z^2+2*z+1],[z-1,-4*n^2*(2*n+1)^2*z^4]]
\end{verbatim}
$$\sum_{k\ge0}(T_k-T_{k+2})z^k=1+\dfrac{z-1}{6z^2+2z+1-\dfrac{36z^4}{26z^2+2z+1-\dfrac{400z^4}{62z^2+2z+1-\ddots}}}$$
Convergence type $P^+$ with $P=1+1/z$ and
$C=-\G(1+1/(2z))^2(1/z)(1/z^2-1)/2^{2/z+2}$, so that
$$\sum_{k\ge0}(T_k-T_{k+2})z^k-\dfrac{p(n)}{q(n)}\sim-\dfrac{\G(1+1/(2z))^2(1/z)(1/z^2-1)/2^{2/z+2}}{n^{1+1/z}}\;.$$
$$A=1-((3z-3)/(4z))/n+((z+1)(96z^2+52z-1)/(192z^3))/n^2+\cdots$$
Parametric family for $k\ge0$:
\begin{verbatim}
[(z)->sum(k=0,oo,(tanfrac(k)-tanfrac(k+2))*z^k),
(8*n^2-4*n+(2*k+1)^2+1)*z^2+(4*k+2)*z+1,-4*n^2*(2*n+1)^2*z^4]
\end{verbatim}
Convergence type $P^+$ with $P=2k+1+1/z$.
\end{cf}

\smallskip

\begin{cf}\label{6.16.26}{\ }
\begin{verbatim}
[(z)->sum(k=0,oo,(tanfrac(k)-tanfrac(k+2))*z^k),
[(z^2+z-1)/z^2,2*z+1,2*z],[-2*(1-1/z^2),(2*n*z+1)^2]]
\end{verbatim}
$$\sum_{k\ge0}(T_k-T_{k+2})z^k=(z^2+z-1)/z^2-\dfrac{(2z^2-2)/z^2}{2z+1+\dfrac{4z^2+4z+1}{2z+\dfrac{16z^2+8z+1}{2z+\dfrac{36z^2+12z+1}{2z+\ddots}}}}$$
Convergence type $P^-$ with $P=1$ and $C=(1-z^2)/(2z^3)$, so that
$$\sum_{k\ge0}(T_k-T_{k+2})z^k-\dfrac{p(n)}{q(n)}\sim(-1)^n\dfrac{(1-z^2)/(2z^3)}{n}\;.$$
$$A=1-((z+1)/2z)/n+((2z+1)/(4z^2))/n^2+((2z^3-3z-1)/(8z^3))/n^3+\cdots$$
Series:
$$\sum_{k\ge0}(T_k-T_{k+2})z^k=\dfrac{z^2+z-1}{z^2}+2\dfrac{z^2-1}{z^2}\sum_{n\ge1}\dfrac{(-1)^n}{2nz+1}$$
Parametric family for $k\ge0$:
\begin{verbatim}
[(z)->sum(k=0,oo,(tanfrac(k)-tanfrac(k+2))*z^k),(4*k+2)*z,(2*n*z+1)^2]
\end{verbatim}
Convergence type $P^-$ with $P=2k+1$.
\end{cf}

\smallskip

\begin{cf}\label{6.16.27}{\ }
\begin{verbatim}
[(z)->sum(k=0,oo,(tanfrac(k)-tanfrac(k+2))*z^k),
[(z^2+z-1)/z^2,z+1],[-(1-1/z^2),z^2*n^2]]
\end{verbatim}
$$\sum_{k\ge0}(T_k-T_{k+2})z^k=(z^2+z-1)/z^2-\dfrac{(z^2-1)/z^2}{z+1+\dfrac{z^2}{z+1+\dfrac{4z^2}{z+1+\dfrac{9z^2}{z+1+\dfrac{16z^2}{z+1+\ddots}}}}}$$
Convergence type $P^-$ with $P=1+1/z$ and $C=...$, so that
$$\sum_{k\ge0}(T_k-T_{k+2})z^k-\dfrac{p(n)}{q(n)}\sim(-1)^n\dfrac{C}{n^{1+1/z}}\;.$$
$$A=1-((z+1)/(2z))/n+((4z^2+3z-1)/(48z^3))/n^2+\cdots$$
Parametric family for $k\ge0$:
\begin{verbatim}
[(z)->sum(k=0,oo,(tanfrac(k)-tanfrac(k+2))*z^k),(2*k+1)*z+1,z^2*n^2]
\end{verbatim}
Convergence type $P^-$ with $P=2k+1+1/z$.
\end{cf}

\smallskip

\begin{cf}\label{6.16.28}{\ }
\begin{verbatim}
[(z)->sum(k=0,oo,(tanfrac(k)-tanfrac(k+2))*z^k),
[(z^2+z-1)/z^2,2*z*n+1],
[[-2*(1-1/z^2),(2*z+1)^2],[4*z^2*n^2,(2*n*z+2*z+1)^2]]]
\end{verbatim}
$$\sum_{k\ge0}(T_k-T_{k+2})z^k=(z^2+z-1)/z^2-\dfrac{(2z^2-2)/z^2}{2z+1+\dfrac{4z^2+4z+1}{4z+1+\dfrac{4z^2}{6z+1+\dfrac{16z^2+8z+1}{8z+1+\ddots}}}}$$
Convergence type $E$ with $E=-(1+\sqrt{2})^2$, $P=0$, and $C=...$, so that
$$\sum_{k\ge0}(T_k-T_{k+2})z^k-\dfrac{p(n)}{q(n)}\sim(-1)^n\dfrac{C}{(1+\sqrt{2})^{2n}}\;.$$
\end{cf}

\smallskip

By using the formula
$$\sum_{k\ge0}(T_k-T_{k+2})z^k=1/z+(1-1/z^2)\sum_{k\ge0}T_kz^k$$
one can obtain from these other CFs for $\sum_{k\ge0}T_kz^k$ which do not
seem worth recording.

\smallskip

\begin{cf}\label{6.16.29}{\ }
\begin{verbatim}
[(z)->sum(k=0,oo,(4*tanfrac(k)-tanfrac(k+2))*z^k),[0,1],[6*z,z^2*n*(n+3)]]
\end{verbatim}
$$\sum_{k\ge0}(4T_k-T_{k+2})z^k=\dfrac{6z}{1+\dfrac{4z^2}{1+\dfrac{10z^2}{1+\dfrac{18z^2}{1+\dfrac{28z^2}{1+\dfrac{40z^2}{1+\ddots}}}}}}$$
Convergence type $P^-$ with $P=1/z$ and $C=\G(2+1/(2z))^2/2^{1/z-3}$, so that
$$\sum_{k\ge0}(4T_k-T_{k+2})z^k-\dfrac{p(n)}{q(n)}\sim(-1)^n\dfrac{\G(2+1/(2z))^2/2^{1/z-3}}{n^{1/z}}\;.$$
$$A=1-(2/z)/n+((112z^2+96z-1)/(48z^3))/n^2+\cdots$$
Parametric family for $k\ge0$:
\begin{verbatim}
[(z)->sum(k=0,oo,(4*tanfrac(k)-tanfrac(k+2))*z^k),2*k*z+1,z^2*n*(n+3)]
\end{verbatim}
Convergence type $P^-$ with $P=2k+1/z$.
\end{cf}

\smallskip

\begin{cf}\label{6.16.30}{\ }
\begin{verbatim}
[(z)->sum(k=0,oo,(4*tanfrac(k)-tanfrac(k+2))*z^k),
[0,(8*n^2-4)*z^2+1],[6*z,-4*z^4*n*(n+1)*(2*n-1)*(2*n+3)]]
\end{verbatim}
$$\sum_{k\ge0}(4T_k-T_{k+2})z^k=\dfrac{6z}{4z^2+1-\dfrac{40z^4}{28z^2+1-\dfrac{504z^4}{68z^2+1-\dfrac{2160z^4}{124z^2+1-\dfrac{6160z^4}{196z^2+1-\ddots}}}}}$$
Convergence type $P^+$ with $P=1/z$ and $C=\G(2+1/(2z))^2/2^{2/z-3}$, so that
$$\sum_{k\ge0}(4T_k-T_{k+2})z^k-\dfrac{p(n)}{q(n)}\sim\dfrac{\G(2+1/(2z))^2/2^{2/z-3}}{n^{1/z}}\;.$$
$$A=1-(1/z)/n+((112z^2+96z-1)/(192z^3))/n^2+\cdots$$
Parametric family for $k\ge0$:
\begin{verbatim}
[(z)->sum(k=0,oo,(4*tanfrac(k)-tanfrac(k+2))*z^k),
8*n^2*z^2+(4*k^2-4)*z^2+4*k*z+1,-4*z^4*n*(n+1)*(2*n-1)*(2*n+3)]
\end{verbatim}
Convergence type $P^+$ with $P=2k+1/z$.
\end{cf}

This is the contraction of the previous CF. Both can be Ap\'ery accelerated
with convergence type $E=-(1+\sqrt{2})^2$, formula too complicated to give
here.

\smallskip

\begin{cf}\label{6.16.31}{\ }
\begin{verbatim}
[(z)->sum(k=0,oo,(4*tanfrac(k)-tanfrac(k+2))*z^k),
[0,4*z+1,8*z],[12*z,(2*n*z-4*z+1)*(2*n*z+2*z+1)]]
\end{verbatim}
$$\sum_{k\ge0}(4T_k-T_{k+2})z^k=\dfrac{12z}{4z+1-\dfrac{8z^2-2z-1}{8z+\dfrac{6z+1}{8z+\dfrac{16z^2+10z+1}{8z+\dfrac{40z^2+14z+1}{8z+\dfrac{72z^2+18z+1}{8z+\ddots}}}}}}$$
Convergence type $P^-$ with $P=4$ and $C=3(1-4z^2)/(8z^3)$, so that
$$\sum_{k\ge0}(4T_k-T_{k+2})z^k-\dfrac{p(n)}{q(n)}\sim(-1)^n\dfrac{3(1-4z^2)/(8z^3)}{n^4}\;.$$
$$A=1-(2/z)/n+(5/(2z^2))/n^2-(5/(2z^3))/n^3+\cdots$$
Series:
$$\sum_{k\ge0}(4T_k-T_{k+2})z^k=12z(4z^2-1)\sum_{n\ge1}\dfrac{(-1)^n}{((2n+2)z+1)(2nz+1)((2n-2)z+1)((2n-4)z+1)}$$
Parametric family for $k\ge0$:
\begin{verbatim}
[(z)->sum(k=0,oo,(4*tanfrac(k)-tanfrac(k+2))*z^k),
4*(k+2)*z,(2*n*z-4*z+1)*(2*n*z+2*z+1)]
\end{verbatim}
Convergence type $P^-$ with $P=2k+4$.
\end{cf}

\smallskip

\begin{cf}\label{6.16.32}{\ }
\begin{verbatim}
[(z)->exp(sum(k=1,oo,tanfrac(k)/k*z^k)),[z+1,2*z+2],[z^2*(2*n+1)^2]]
\end{verbatim}
$$\exp\left(\sum_{k\ge1}\dfrac{T_k}{k}z^k\right)=z+1+\dfrac{z^2}{2z+2+\dfrac{9z^2}{2z+2+\dfrac{25z^2}{2z+2+\dfrac{49z^2}{2z+2+\dfrac{81z^2}{2z+2+\ddots}}}}}$$
Convergence type $P^-$ with $P=1+1/z$ and
$C=z\G(1+1/(2z))^4/(\G(1/2+1/(4z))^42^{2/z-2})$, so that
$$\exp\left(\sum_{k\ge1}\dfrac{T_k}{k}z^k\right)-\dfrac{p(n)}{q(n)}\sim(-1)^n\dfrac{z\G(1+1/(2z))^4/(\G(1/2+1/(4z))^42^{2/z-2})}{n^{1+1/z}}\;.$$
$$A=1-((z+1)/z)/n+((z+1)(36z^2+22z-1)/(48z^3))/n^2+\cdots$$
Parametric family for $k\ge0$:
\begin{verbatim}
[(z)->exp(sum(k=1,oo,tanfrac(k)/k*z^k)),(4*k+2)*z+2,z^2*(2*n+1)^2]
\end{verbatim}
Convergence type $P^-$ with $P=2k+1+1/z$.
\end{cf}

\smallskip

\begin{cf}\label{6.16.33}{\ }
\begin{verbatim}
[(z)->exp(sum(k=1,oo,tanfrac(k)/k*z^k)),[z+1,3*z+1,4*z],[z^2,(2*n*z+1)^2]]
\end{verbatim}
$$\exp\left(\sum_{k\ge1}\dfrac{T_k}{k}z^k\right)=z+1+\dfrac{z^2}{3z+1+\dfrac{4z^2+4z+1}{4z+\dfrac{16z^2+8z+1}{4z+\dfrac{36z^2+12z+1}{4z+\ddots}}}}$$
Convergence type $P^-$ with $P=2$ and
$C=2^{1/z}\G(1/(4z))^4/(128\pi z\G(1/(2z))^2)$, so that
$$\exp\left(\sum_{k\ge1}\dfrac{T_k}{k}z^k\right)-\dfrac{p(n)}{q(n)}\sim(-1)^n\dfrac{2^{1/z}\G(1/(4z))^4/(128\pi z\G(1/(2z))^2)}{n^2}\;.$$
$$A=1-((z+1)/z)/n+((z^2+12z+6)/(8z^2))/n^2+((3z^3-z^2-6z-2)/(4z^3))/n^3+\cdots$$
Series:
\begin{align*}
\exp\left(\sum_{k\ge1}\dfrac{T_k}{k}z^k\right)&=z+1+\dfrac{4z^3}{(4z+1)^2}\sum_{n\ge0}\dfrac{((2z+1)/(4z))_n^2}{((8z+1)/(4z))_n^2}\\
\exp\left(-\sum_{k\ge1}\dfrac{T_k}{k}z^k\right)&=1-z+\dfrac{z}{4}\sum_{n\ge0}\dfrac{(1/(4z))_n^2}{((2z+1)/(4z))_{n+1}^2}\;.\end{align*}
Parametric family for $k\ge0$:
\begin{verbatim}
[(z)->exp(sum(k=1,oo,tanfrac(k)/k*z^k)),4*(k+1)*z,(2*n*z+1)^2]
\end{verbatim}
Convergence type $P^-$ with $P=2k+2$.
\end{cf}
   
\smallskip

\begin{cf}\label{6.16.34}{\ }
\begin{verbatim}
[(z)->exp(sum(k=1,oo,tanfrac(k)/k*z^k)),
[(2*n+1)*z+1],[[z^2,(2*z+1)^2],[z^2*(2*n+1)^2,((2*n+2)*z+1)^2]]]
\end{verbatim}
$$\exp\left(\sum_{k\ge1}\dfrac{T_k}{k}z^k\right)=z+1+\dfrac{z^2}{3z+1+\dfrac{4z^2+4z+1}{5z+1+\dfrac{9z^2}{7z+1+\dfrac{16z^2+8z+1}{9z+1+\ddots}}}}$$
Convergence type $E$ with $E=-(1+\sqrt{2})^2$, $P=0$, and $C=\pi z\G(1+1/(2z))^2/(\G(1/2+1/(4z))^42^{1/z-4}(1+\sqrt{2})^{1/z+3})$, so that
$$\exp\left(\sum_{k\ge1}\dfrac{T_k}{k}z^k\right)-\dfrac{p(n)}{q(n)}\sim(-1)^n\dfrac{\pi z\G(1+1/(2z))^2/(\G(1/2+1/(4z))^42^{1/z-4})}{(1+\sqrt{2})^{2n+1/z+3}}\;.$$
$$A=1-((z^2+4z-2)d/(8z^2))/n+\cdots$$
\end{cf}

\smallskip

\begin{cf}\label{6.16.35}{\ }
\begin{verbatim}
[(z)->sum(k=0,oo,tan3frac(k)*z^k),[0,(2*n-1)*(z^2*(9*n^2-9*n+2)+2)],
                                [2*z,-z^4*n^2*(9*n^2-1)*(9*n^2-4)]]
\end{verbatim}
$$\sum_{k\ge0}T^{(3)}_kz^k=\dfrac{2z}{2z^2+2-\dfrac{40z^4}{60z^2+6-\dfrac{4480z^4}{280z^2+10-\dfrac{55440z^4}{770z^2+14-\ddots}}}}$$
Convergence type $P^+$ with $P=4/(3|z|)$ and
$C=(\G(1+1/z)\G(1+1/(3z))/\G(1+2/(3z)))^2/3^{2/z-1}$, so that
$$\sum_{k\ge0}T^{(3)}_kz^k-\dfrac{p(n)}{q(n)}\sim\dfrac{(\G(1+1/z)\G(1+1/(3z))/\G(1+2/(3z)))^2/3^{2/z-1}}{n^{4/(3|z|)}}\;.$$
$$A=1-(2/(3z))/n+((81z^4+162z^3-48z^2-72z+4)/(81z^3(9z^2-4)))/n^2+\cdots$$
Parametric family for $k\ge0$:
\begin{verbatim}
[(z)->sum(k=0,oo,tan3frac(k)*z^k),
(2*n-1)*(z^2*(9*n^2-9*n+18*k^2+2)+12*k*z+2),-z^4*n^2*(9*n^2-1)*(9*n^2-4)]
\end{verbatim}
Convergence type $P^+$ with $P=4k+4/(3z)$.
\end{cf}

\smallskip

\begin{cf}\label{6.16.36}{\ }
\begin{verbatim}
[(z)->sum(k=0,oo,(k+1)*bernfrac(k)*z^k),
                     [1-z,n+1],z^2*[[1,2],(n+1)^3*[n,n+2]]]
\end{verbatim}
$$\sum_{k\ge0}(k+1)B_kz^k=1-z+\dfrac{z^2}{2+\dfrac{2z^2}{3+\dfrac{8z^2}{4+\dfrac{24z^2}{5+\dfrac{54z^2}{6+\dfrac{108z^2}{7+\ddots}}}}}}$$
Convergence type $P^-$ with $P=4/|z|$ and
$C=2^{1+4/z}z\G(1+1/z)^8/\G(1+2/z)^2$, so that
$$\sum_{k\ge0}(k+1)B_kz^k-\dfrac{p(n)}{q(n)}\sim(-1)^n\dfrac{2^{1+4/z}z\G(1+1/z)^8/\G(1+2/z)^2}{n^{4/|z|}}\;.$$
$$A=1-(8/z)/n+((28z^4+96z^3-104z^2-384z+16)/(3z^3(z^2-4)))/n^2+\cdots$$
\end{cf}

\smallskip

\begin{cf}\label{6.16.37}{\ }
\begin{verbatim}
[(z)->sum(k=0,oo,(k+1)*bernfrac(k)*z^k),
[0,(z+1)^3,((2*n-1)*z+2)*(n^2*z^2+z*n*(2-z)+z^2-z+1)],[2*z,-(n*z+1)^6]]
\end{verbatim}
$$\sum_{k\ge0}(k+1)B_kz^k=\dfrac{2z}{z^3+3z^2+3z+1-\dfrac{z^6+6z^5+15z^4+20z^3+15z^2+6z+1}{9z^3+15z^2+9z+2-\ddots}}$$
Convergence type $P^+$ with $P=2$ and $C=1/z^2$, so that
$$\sum_{k\ge0}(k+1)B_kz^k-\dfrac{p(n)}{q(n)}\sim\dfrac{C}{n^2}\;.$$
$$A=1-((z+2)/z)/n+((z^2+6z+6)/(2z^2))/n^2-((2z^2+6z+4)/z^3)/n^3+\cdots$$
Series:
$$\sum_{k\ge0}(k+1)B_kz^k=2z\sum_{n\ge1}\dfrac{1}{(nz+1)^3}$$
Parametric family for $k\ge0$:
\begin{verbatim}
[(z)->sum(k=0,oo,(k+1)*bernfrac(k)*z^k),
((2*n-1)*z+2)*(n^2*z^2+z*n*(2-z)+(2*k^2+2*k+1)*z^2-z+1),-(n*z+1)^6]
\end{verbatim}
Convergence type $P^+$ with $P=4k+2$.
\end{cf}

\smallskip

\begin{cf}\label{6.16.38}{\ }
\begin{verbatim}
[(z)->sum(k=0,oo,(k+1)*bernfrac(k)*z^k),
[0,(2*n-1)*(z^2*(n^2-n+1)+2*z+2)],[2,-z^4*n^6]]
\end{verbatim}
$$\sum_{k\ge0}(k+1)B_kz^k=\dfrac{2}{z^2+2z+2-\dfrac{z^4}{9z^2+6z+6-\dfrac{64z^4}{35z^2+10z+10-\ddots}}}$$
Convergence type $P^+$ with $P=|4/z+2|$ and
$C=\G(1+1/z)^8/(z^2\G(1+2/z)\G(2+2/z))$, so that
$$\sum_{k\ge0}(k+1)B_kz^k-\dfrac{p(n)}{q(n)}\sim\dfrac{\G(1+1/z)^8/(z^2\G(1+2/z)\G(2+2/z))}{n^{|4/z+2|}}\;.$$
$$A=1-((z+2)/z)/n+((z+2)(3z^2+5z-1)/(6z^3))/n^2+\cdots$$
Parametric family for $k\ge0$:
\begin{verbatim}
[(z)->sum(k=0,oo,(k+1)*bernfrac(k)*z^k),
(2*n-1)*(z^2*(n^2-n+2*k^2+2*k+1)+(4*k+2)*z+2),-z^4*n^6]
\end{verbatim}
Convergence type $P^+$ with $P=4k+2+4/z$.
\end{cf}

If desired, one can also obtain a CF for $\sum_{k\ge0}(k+2)(k+1)B_kz^k$
by changing $z$ into $1/z$ and multiplying by $z^4$ in \ref{6.15.1}.

\smallskip

Note that since $(k+1)T_k=-2^kG_{k+1}$, one can trivially obtain CFs for
$\sum_{k\ge0}(k+1)T_kz^k$ from those for $\sum_{k\ge0}G_kz^k$.

\smallskip

\begin{cf}\label{6.16.39}{\ }
\begin{verbatim}
[(z)->sum(k=0,oo,eulerfrac(k)*z^k),[0,1],[1,n^2*z^2]]
\end{verbatim}
$$\sum_{k\ge0}E_kz^k=\dfrac{1}{1+\dfrac{z^2}{1+\dfrac{4z^2}{1+\dfrac{9z^2}{1+\dfrac{16z^2}{1+\dfrac{25z^2}{1+\ddots}}}}}}$$
Convergence type $P^-$ with $P=1/|z|$ and $C=\G((1+1/z)/2)^2/(z\cdot2^{1/z})$,
so that
$$\sum_{k\ge0}E_kz^k-\dfrac{p(n)}{q(n)}\sim(-1)^n\dfrac{\G((1+1/z)/2)^2/(z\cdot2^{1/z})}{n^{1/|z|}}\;.$$
$$A=1-(1/(2z))/n-((z-1)(5z-1)/(48z^3))/n^2+\cdots$$
Parametric family for $k\ge0$:
\begin{verbatim}
[(z)->sum(k=0,oo,eulerfrac(k)*z^k),2*k*z+1,n^2*z^2]
\end{verbatim}
Convergence type $P^-$ with $P=2k+1/z$.
\end{cf}

\smallskip

\begin{cf}\label{6.16.40}{\ }
\begin{verbatim}
[(z)->sum(k=0,oo,eulerfrac(k)*z^k),[0,z+1,2*z],[2,(z*(2*n-1)+1)^2]]
\end{verbatim}
$$\sum_{k\ge0}E_kz^k=\dfrac{2}{z+1+\dfrac{z^2+2z+1}{2z+\dfrac{9z^2+6z+1}{2z+\dfrac{25z^2+10z+1}{2z+\dfrac{49z^2+14z+1}{2z+\ddots}}}}}$$
Convergence type $P^-$ with $P=1$ and $C=1/(2z)$, so that
$$\sum_{k\ge0}E_kz^k-\dfrac{p(n)}{q(n)}\sim(-1)^n\dfrac{1/(2z)}{n}\;.$$
$$A=1-(1/(2z))/n+((1-z^2)/(4z^2))/n^2+((3z^2-1)/(8z^3))/n^3+\cdots$$
Series:
$$\sum_{k\ge0}E_kz^k=2\sum_{n\ge1}\dfrac{(-1)^{n+1}}{(2n-1)z+1}$$
Parametric family for $k\ge0$:
\begin{verbatim}
[(z)->sum(k=0,oo,eulerfrac(k)*z^k),(4*k+2)*z,(z*(2*n-1)+1)^2]
\end{verbatim}
Convergence type $P^-$ with $P=2k+1$.
\end{cf}
        
\smallskip

\begin{cf}\label{6.16.41}{\ }
\begin{verbatim}
[(z)->sum(k=0,oo,eulerfrac(k)*z^k),
[0,(8*n^2-12*n+5)*z^2+1],[1,-4*n^2*(2*n-1)^2*z^4]]
\end{verbatim}
$$\sum_{k\ge0}E_kz^k=\dfrac{1}{z^2+1-\dfrac{4z^4}{13z^2+1-\dfrac{144z^4}{41z^2+1-\dfrac{900z^4}{85z^2+1-\ddots}}}}$$
Convergence type $P^+$ with $P=1/|z|$ and $C=\G((1+1/z)/2)^2/(z\cdot2^{2/z})$,
so that
$$\sum_{k\ge0}E_kz^k-\dfrac{p(n)}{q(n)}\sim\dfrac{\G((1+1/z)/2)^2/(z\cdot2^{2/z})}{n^{1/|z|}}\;.$$
$$A=1-(1/(4z))/n-((z-1)(5z-1)/(192z^3))/n^2+\cdots$$
Parametric family for $k\ge0$:
\begin{verbatim}
[(z)->sum(k=0,oo,eulerfrac(k)*z^k),
(8*n^2-12*n+4*k^2+5)*z^2+4*k*z+1,-4*n^2*(2*n-1)^2*z^4]
\end{verbatim}
Convergence type $P^+$ with $P=2k+1/z$.
\end{cf}

This is simply the contraction of \ref{6.16.39}.

\smallskip

\begin{cf}\label{6.16.42}{\ }
\begin{verbatim}
[(z)->sum(k=0,oo,eulerfrac(k)*z^k),
[0,(2*n-1)*z+1],[[2,(z+1)^2],[4*z^2*n^2,(z*(2*n+1)+1)^2]]]
\end{verbatim}
$$\sum_{k\ge0}E_kz^k=\dfrac{2}{z+1+\dfrac{z^2+2z+1}{3z+1+\dfrac{4z^2}{5z+1+\dfrac{9z^2+6z+1}{7z+1+\dfrac{16z^2}{9z+1+\dfrac{25z^2+10z+1}{11z+1+\ddots}}}}}}$$
Convergence type $E$ with $E=-(1+\sqrt{2})^2$, $P=0$, and $C=(2\pi/z)/(1+\sqrt{2})^{1+1/z}$, so that
$$\sum_{k\ge0}E_kz^k-\dfrac{p(n)}{q(n)}\sim(-1)^n\dfrac{2\pi/z}{(1+\sqrt{2})^{2n+1/z+1}}\;.$$
$$A=1-((z^2+4z-2)d/(8z^2))/n+\cdots$$
\end{cf}

\smallskip

\begin{cf}\label{6.16.43}{\ }
\begin{verbatim}
[(z)->sum(k=0,oo,(eulerfrac(k)-eulerfrac(k+2))*z^k),
               [0,1],[2,n*(n+2)*z^2]]
\end{verbatim}
$$\sum_{k\ge0}(E_k-E_{k+2})z^k=\dfrac{2}{1+\dfrac{3z^2}{1+\dfrac{8z^2}{1+\dfrac{15z^2}{1+\dfrac{24z^2}{1+\dfrac{35z^2}{1+\ddots}}}}}}$$
Convergence type $P^-$ with $P=1/|z|$ and
$C=\G((3+1/z)/2)^2/(z\cdot2^{1/z-2})$, so that
$$\sum_{k\ge0}(E_k-E_{k+2})z^k-\dfrac{p(n)}{q(n)}\sim(-1)^n\dfrac{\G((3+1/z)/2)^2/(z\cdot2^{1/z-2})}{n^{1/|z|}}\;.$$
$$A=1-(3/(2z))/n+((55z^2+54z-1)/(48z^3))/n^2+\cdots$$
Parametric family for $k\ge0$:
\begin{verbatim}
[(z)->sum(k=0,oo,(eulerfrac(k)-eulerfrac(k+2))*z^k),2*k*z+1,n*(n+2)*z^2]
\end{verbatim}
Convergence type $P^-$ with $P=2k+1/z$.
\end{cf}

\smallskip

\begin{cf}\label{6.16.44}{\ }
\begin{verbatim}
[(z)->sum(k=0,oo,(eulerfrac(k)-eulerfrac(k+2))*z^k),
[0,3*z+1,6*z],[4,((2*n+1)*z+1)*((2*n-3)*z+1)]]
\end{verbatim}
$$\sum_{k\ge0}(E_k-E_{k+2})z^k=\dfrac{4}{3z+1-\dfrac{3z^2-2z-1}{6z+\dfrac{5z^2+6z+1}{6z+\dfrac{21z^2+10z+1}{6z+\dfrac{45z^2+14z+1}{6z+\ddots}}}}}$$
Convergence type $P^-$ with $P=3$ and $C=(1-z^2)/(4*z^3)$, so that
$$\sum_{k\ge0}(E_k-E_{k+2})z^k-\dfrac{p(n)}{q(n)}\sim(-1)^n\dfrac{(1-z^2)/(4z^3)}{n^3}\;.$$
$$A=1-(3/(2z))/n+((3-z^2)/(2z^2))/n^2+((5z^2-5)/(4z^3))/n^3+\cdots$$
Series:
$$\sum_{k\ge0}(E_k-E_{k+2})z^k=4(1-z^2)\sum_{n\ge1}\dfrac{(-1)^{n+1}}{((2n+1)z+1)((2n-1)z+1)((2n-3)z+1)}$$
Parametric family for $k\ge0$:
\begin{verbatim}
[(z)->sum(k=0,oo,(eulerfrac(k)-eulerfrac(k+2))*z^k),
(4*k+6)*z,((2*n+1)*z+1)*((2*n-3)*z+1)]
\end{verbatim}
Convergence type $P^-$ with $P=2k+3$.
\end{cf}

\smallskip

\begin{cf}\label{6.16.45}{\ }
\begin{verbatim}
[(z)->sum(k=0,oo,(eulerfrac(k)-eulerfrac(k+2))*z^k),
[0,3*z+1,z*(2*n+1)+1],[[4,(1-z)*(1+3*z)],
[4*n*(n+2)*z^2,(z*(2*n-1)+1)*(z*(2*n+3)+1)]]]
\end{verbatim}
$$\sum_{k\ge0}(E_k-E_{k+2})z^k=\dfrac{4}{3z+1-\dfrac{3z^2-2z-1}{5z+1+\dfrac{12z^2}{7z+1+\dfrac{5z^2+6z+1}{9z+1+\dfrac{32z^2}{11z+1+\ddots}}}}}$$
Convergence type $E$ with $E=-(1+\sqrt{2})^2$, $P=0$, and $C=...$, so that
$$\sum_{k\ge0}(E_k-E_{k+2})z^k-\dfrac{p(n)}{q(n)}\sim(-1)^n\dfrac{C}{(1+\sqrt{2})^{2n}}\;.$$
$$A=1+((31z^2-12z+2)d/(8z^2))/n+\cdots$$
\end{cf}

\smallskip

By using the formula
$$\sum_{k\ge0}(E_k-E_{k+2})z^k=1/z^2+(1-1/z^2)\sum_{k\ge0}E_kz^k$$
one can obtain other CFs for $\sum_{k\ge0}E_kz^k$ which do not
seem worth recording.

\smallskip

\begin{cf}\label{6.16.46}{\ }
\begin{verbatim}
[(z)->sum(k=0,oo,(k+1)*eulerfrac(k)*z^k),
       [0,(8*n^2-8*n+3)*z^2+1],[1,-(2*n)^4*z^4]]
\end{verbatim}
$$\sum_{k\ge0}(k+1)E_kz^k=\dfrac{1}{3z^2+1-\dfrac{16z^4}{19z^2+1-\dfrac{256z^4}{51z^2+1-\dfrac{1296z^4}{99z^2+1-\dfrac{4096z^4}{163z^2+1-\ddots}}}}}$$
Convergence type $P^+$ with $P=1/|z|$ and $C=\G(1/2+1/(2z))^6/(4z\G(1/z)^2)$,
so that
$$\sum_{k\ge0}(k+1)E_kz^k-\dfrac{p(n)}{q(n)}\sim\dfrac{\G(1/2+1/(2z))^6/(4z\G(1/z)^2)}{n^{1/|z|}}\;.$$
$$A=1-(1/(2z))/n+((37z^4+96z^3-14z^2-24z+1)/(192z^3(4z^2-1)))/n^2+\cdots$$
Parametric family for $k\ge0$:
\begin{verbatim}
[(z)->sum(k=0,oo,(k+1)*eulerfrac(k)*z^k),
(8*n^2-8*n+4*k^2+3)*z^2+4*k*z+1,-(2*n)^4*z^4]
\end{verbatim}
Convergence type $P^+$ with $P=2k+1/z$.
\end{cf}

\smallskip

\begin{cf}\label{6.16.47}{\ }
\begin{verbatim}
[(z)->sum(k=0,oo,(k+1)*eulerfrac(k)*z^k),
[0,(z+1)^2,4*z*(2*z*(n-1)+1)],[2,(z*(2*n-1)+1)^4]]
\end{verbatim}
$$\sum_{k\ge0}(k+1)E_kz^k=\dfrac{2}{z^2+2z+1+\dfrac{z^4+4z^3+6z^2+4z+1}{12z^2+4z+\dfrac{81z^4+108z^3+54z^2+12z+1}{20z^2+4z+\ddots}}}$$
Convergence type $P^-$ with $P=2$ and $C=1/(4z^2)$, so that
$$\sum_{k\ge0}(k+1)E_kz^k-\dfrac{p(n)}{q(n)}\sim(-1)^n\dfrac{1/(4z^2)}{n^2}\;.$$
$$A=1-(1/z)/n+((-3z^2+3)/(4z^2))/n^2+((3z^2-1)/(2z^3))/n^3+\cdots$$
Series:
$$\sum_{k\ge0}(k+1)E_kz^k=2\sum_{n\ge1}\dfrac{(-1)^{n+1}}{((2n-1)z+1)^2}$$
Parametric family for $k\ge0$:
\begin{verbatim}
[(z)->sum(k=0,oo,(k+1)*eulerfrac(k)*z^k),
(2*k+1)*4*z*(2*z*(n-1)+1),(z*(2*n-1)+1)^4]
\end{verbatim}
Convergence type $P^-$ with $P=4k+2$.
\end{cf}

Can be Ap\'ery accelerated with convergence type $E=-((1+\sqrt{5})/2)^5$,
formula too complicated to give here.

\smallskip

\begin{cf}\label{6.16.48}{\ }
\begin{verbatim}
[(z)->exp(-sum(k=1,oo,eulerfrac(k)/k*z^k)),
                    [1,2],[z^2*(2*n+1)^2]]
\end{verbatim}
$$\exp\left(-\sum_{k\ge1}\dfrac{E_{k}}{k}z^{k}\right)=1+\dfrac{z^2}{2+\dfrac{9z^2}{2+\dfrac{25z^2}{2+\dfrac{49z^2}{2+\dfrac{81z^2}{2+\dfrac{121z^2}{2+\ddots}}}}}}$$
Convergence type $P^-$ with $P=1/z$ and
$C=z(\G(1/2+1/(2z))/\G(1/4+1/(4z)))^4/2^{2/z-4}$, so that
$$\exp\left(-\sum_{k\ge1}\dfrac{E_{k}}{k}z^{k}\right)-\dfrac{p(n)}{q(n)}\sim(-1)^n\dfrac{z(\G(1/2+1/(2z))/\G(1/4+1/(4z)))^4/2^{2/z-4}}{n^{1/z}}\;.$$
$$A=1-(1/z)/n+((13z^2+24z-1)/(48z^3))/n^2+\cdots$$
Parametric family for $k\ge0$:
\begin{verbatim}
[(z)->exp(-sum(k=1,oo,eulerfrac(k)/k*z^k)),4*k*z+2,z^2*(2*n+1)^2]
\end{verbatim}
Convergence type $P^-$ with $P=2k+1/z$.
\end{cf}

\smallskip

\begin{cf}\label{6.16.49}{\ }
\begin{verbatim}
[(z)->exp(-sum(k=1,oo,eulerfrac(k)/k*z^k)),
[1,2*z+1,4*z],[z^2,(z*(2*n-1)+1)^2]]
\end{verbatim}
$$\exp\left(-\sum_{k\ge1}\dfrac{E_{k}}{k}z^{k}\right)=1+\dfrac{z^2}{2z+1+\dfrac{z^2+2z+1}{4z+\dfrac{9z^2+6z+1}{4z+\dfrac{25z^2+10z+1}{4z+\dfrac{49z^2+14z+1}{4z+\ddots}}}}}$$
Convergence type $P^-$ with $P=2$ and $C=...$, so that
$$\exp\left(-\sum_{k\ge1}\dfrac{E_{k}}{k}z^{k}\right)-\dfrac{p(n)}{q(n)}\sim(-1)^n\dfrac{C}{n^2}\;.$$
$$A=1-(1/z)/n+((-5z^2+6)/(8z^2))/n^2+((5z^2-2)/(4z^3))/n^3+\cdots$$
Series:
\begin{align*}
  \exp\left(-\sum_{k\ge1}\dfrac{E_{k}}{k}z^{k}\right)&=1+\dfrac{4z^3}{(3z+1)^2}\sum_{n\ge0}\dfrac{((z+1)/(4z))_n^2}{((7z+1)/(4z))_n^2}\\
  \exp\left(\sum_{k\ge1}\dfrac{E_{k}}{k}z^{k}\right)&=\dfrac{1-2z}{(z-1)^2}+\dfrac{z}{4(z-1)^2}\sum_{n\ge0}\dfrac{((1-z)/(4z))_n^2}{((z+1)/(4z))_{n+1}^2}\;.\end{align*}
Parametric family for $k\ge0$:
\begin{verbatim}
[(z)->exp(-sum(k=1,oo,eulerfrac(k)/k*z^k)),4*(k+1)*z,(z*(2*n-1)+1)^2]
\end{verbatim}
Convergence type $P^-$ with $P=2k+2$.
\end{cf}
    
\smallskip

\begin{cf}\label{6.16.50}{\ }
\begin{verbatim}
[(z)->exp(-sum(k=1,oo,eulerfrac(k)/k*z^k)),
[1,2*n*z+1],[[z^2,(z+1)^2],[z^2*(2*n+1)^2,((2*n+1)*z+1)^2]]]
\end{verbatim}
$$\exp\left(-\sum_{k\ge1}\dfrac{E_{k}}{k}z^{k}\right)=1+\dfrac{z^2}{2z+1+\dfrac{z^2+2z+1}{4z+1+\dfrac{9z^2}{6z+1+\dfrac{9z^2+6z+1}{8z+1+\ddots}}}}$$
Convergence type $E$ with $E=-(1+\sqrt{2})^2$, $P=0$, and
$C=2\pi z\G(1/2+1/(2z))^2/(\G(1/4+1/(4z))^42^{1/z-4}(1+\sqrt{2})^{1/z+2})$,
so that
$$\exp\left(-\sum_{k\ge1}\dfrac{E_{k}}{k}z^{k}\right)-\dfrac{p(n)}{q(n)}\sim(-1)^n\dfrac{2\pi z\G(1/2+1/(2z))^2/(\G(1/4+1/(4z))^42^{1/z-4})}{(1+\sqrt{2})^{2n+1/z+2}}\;.$$
$$A=1+((5z^2-8z+2)d/(8z^2))/n+\cdots$$
\end{cf}

\smallskip

\begin{cf}\label{6.16.51}{\ }
\begin{verbatim}
[(z)->sum(k=0,oo,(eulerfrac(k)+tanfrac(k))*z^k),
               [1,(8*n^2-8*n+3)*z^2+z+1],[z,-z^4*n^2*(16*n^2-1)]]
\end{verbatim}
$$\sum_{k\ge0}(E_k+T_k)z^k=\dfrac{z}{3z^2+z+1-\dfrac{15z^4}{19z^2+z+1-\dfrac{252z^4}{51z^2+z+1-\dfrac{1287z^4}{99z^2+z+1-\ddots}}}}$$
Convergence type $P^+$ with $P=(2+z)/(2z)$ and
$C=\G((3+2/z)/4)^2\G(1+1/z)^42^{1/2+2/z}/(\pi z\G(1+2/z)\G(2+2/z))$, so that
$$\sum_{k\ge0}(E_k+T_k)z^k-\dfrac{p(n)}{q(n)}\sim\dfrac{\G((3+2/z)/4)^2\G(1+1/z)^42^{1/2+2/z}/(\pi z\G(1+2/z)\G(2+2/z))}{n^{(2+z)/(2z)}}\;.$$
\begin{align*}A&=1-((z+2)/(4z))/n\\
  &\phantom{=}+((z+2)(171z^4+126z^3-100z^2-44z+2)/(192z^3(3z-2)(5z+2)))/n^2+\cdots\end{align*}
Parametric family for $k\ge0$:
\begin{verbatim}
[(z)->sum(k=0,oo,(eulerfrac(k)+tanfrac(k))*z^k),
(8*n^2-8*n+4*k^2+2*k+3)*z^2+(4*k+1)*z+1,-z^4*n^2*(16*n^2-1)]
\end{verbatim}
Convergence type $P^+$ with $P=2k+(2+z)/(2z)$.
\end{cf}

\smallskip

\begin{cf}\label{6.16.52}{\ }
\begin{verbatim}
[(z)->sum(k=0,oo,(eulerfrac(k)+tanfrac(k))*z^k),
[1,(z+1)*(2*z+1),2*z^2*(4*n-3)+4*z],[2*z,(2*z*n+1)^2*(z*(2*n-1)+1)^2]]
\end{verbatim}
$$\sum_{k\ge0}(E_k+T_k)z^k=1+\dfrac{2z}{2z^2+3z+1+\dfrac{4z^4+12z^3+13z^2+6z+1}{10z^2+4z+\dfrac{144z^4+168z^3+73z^2+14z+1}{18z^2+4z+\ddots}}}$$
Convergence type $P^-$ with $P=2$ and $C=1/(4z)$, so that
$$\sum_{k\ge0}(E_k+T_k)z^k-\dfrac{p(n)}{q(n)}\sim(-1)^n\dfrac{1/(4z)}{n^2}\;.$$
$$A=1-((z+2)/(2z))/n-((2z^2-3z-3)/(4z^2))/n^2+\cdots$$
Series:
$$\sum_{k\ge0}(E_k+T_k)z^k=1+2z\sum_{n\ge1}\dfrac{(-1)^{n+1}}{(2nz+1)((2n-1)z+1)}$$
Parametric family for $k\ge0$:
\begin{verbatim}
[(z)->sum(k=0,oo,(eulerfrac(k)+tanfrac(k))*z^k),
(2*k+1)*(2*z^2*(4*n-3)+4*z),(2*z*n+1)^2*(z*(2*n-1)+1)^2]
\end{verbatim}
Convergence type $P^-$ with $P=4k+2$.
\end{cf}

Can be Ap\'ery accelerated with convergence type $E=-((1+\sqrt{5})/2)^5$,
formula too complicated to give here.

\smallskip

\begin{cf}\label{6.16.53}{\ }
\begin{verbatim}
[(z)->sum(k=0,oo,(-1)^(k*(k-1)/2)*(eulerfrac(k)+tanfrac(k))*z^k),
               [1-n*z],[z,-n*(n+1)*z^2/2]]
\end{verbatim}
$$\sum_{k\ge0}(-1)^{k(k-1)/2}(E_k+T_k)z^k=\dfrac{z}{-z+1-\dfrac{z^2}{-2z+1-\dfrac{3z^2}{-3z+1-\dfrac{6z^2}{-4z+1-\ddots}}}}$$
Non convergent CF.
\end{cf}

\smallskip

\begin{cf}\label{6.16.54}{\ }
\begin{verbatim}
[(z)->exp(-sum(k=1,oo,(eulerfrac(k)+tanfrac(k))/k*z^k)),
[1,2*z+2,z+2],[[-2*z,4*z^2],z^2*[(4*n-1)*(4*n+1),(4*n+2)^2]]]
\end{verbatim}
$$\exp\left(-\sum_{k\ge1}\dfrac{E_k+T_k}{k}z^k\right)=1-\dfrac{2z}{2z+2+\dfrac{4z^2}{z+2+\dfrac{15z^2}{z+2+\dfrac{36z^2}{z+2+\dfrac{63z^2}{z+2+\ddots}}}}}$$
Convergence type $P^-$ with $P=(z+2)/(2z)$ and $C=...$, so that
$$\exp\left(-\sum_{k\ge1}\dfrac{E_k+T_k}{k}z^k\right)-\dfrac{p(n)}{q(n)}\sim(-1)^n\dfrac{C}{n^{1/2+1/z}}\;.$$
$$A=1-((z^3-6z-4)/(6z^3-8z^2-8z))/n+\cdots$$
\end{cf}

\smallskip

\begin{cf}\label{6.16.55}{\ }
\begin{verbatim}
[(z)->exp(-sum(k=1,oo,(eulerfrac(k)+tanfrac(k))/k*z^k)),
[1,3*z^2+3*z+2,16*n^2*z^2-24*n*z^2+10*z^2+2*z+2],
[-z^2-2*z,-(2*n-1)^2*(16*n^2-1)*z^4]]
\end{verbatim}
$$\exp\left(-\sum_{k\ge1}\dfrac{E_k+T_k}{k}z^k\right)=1-\dfrac{z^2+2z}{3z^2+3z+2-\dfrac{15z^4}{26z^2+2z+2-\dfrac{567z^4}{82z^2+2z+2-\ddots}}}$$
Convergence type $P^+$ with $P=(z+2)/(2z)$ and $C=...$, so that
$$\exp\left(-\sum_{k\ge1}\dfrac{E_k+T_k}{k}z^k\right)-\dfrac{p(n)}{q(n)}\sim(-1)^n\dfrac{C}{n^{1/2+1/z}}\;.$$
$$A=1-((z^3-6z-4)/(12z^3-16z^2-16z))/n+\cdots$$
Parametric family for $k\ge0$:
\begin{verbatim}
[(z)->exp(-sum(k=1,oo,(eulerfrac(k)+tanfrac(k))/k*z^k)),
16*n^2*z^2-24*n*z^2+(8*k^2+4*k+10)*z^2+(8*k+2)*z+2,
-(2*n-1)^2*(16*n^2-1)*z^4]
\end{verbatim}
Convergence type $P^+$ with $P=2k+(z+2)/(2z)$.
\end{cf}

This is the contraction of the previous CF.

\smallskip

\begin{cf}\label{6.16.56}{\ }
\begin{verbatim}
[(z)->sum(k=0,oo,eul3frac(k)*z^k),[0,4*n-2],[2,z^2*n^2*(9*n^2-1)]]
\end{verbatim}
$$\sum_{k\ge0}E^{(3)}_kz^k=\dfrac{2}{2+\dfrac{8z^2}{6+\dfrac{140z^2}{10+\dfrac{720z^2}{14+\dfrac{2288z^2}{18+\dfrac{5600z^2}{22+\ddots}}}}}}$$
Convergence type $P^-$ with $P=4/(3|z|)$ and $C=...$, so that
$$\sum_{k\ge0}E^{(3)}_kz^k-\dfrac{p(n)}{q(n)}\sim(-1)^n\dfrac{C}{n^{4/(3|z|)}}\;.$$
$$A=1-(2/(3z))/n-((15z^2-18z+4)/(81z^3))/n^2+\cdots$$
Parametric family for $k\ge0$:
\begin{verbatim}
[(z)->sum(k=0,oo,eul3frac(k)*z^k),(3*k*z+1)*(4*n-2),z^2*n^2*(9*n^2-1)]
\end{verbatim}
Convergence type $P^-$ with $P=4k+4/(3z)$.
\end{cf}

\smallskip

\begin{cf}\label{6.16.57}{\ }
\begin{verbatim}
[(z)->sum(k=0,oo,(eul3frac(k)-eul3frac(k+2))*z^k),
               [0,4*n+2],[10,z^2*n*(n+2)*(3*n+1)*(3*n+5)]]
\end{verbatim}
$$\sum_{k\ge0}(E^{(3)}_k-E^{(3)}_{k+2})z^k=\dfrac{10}{6+\dfrac{96z^2}{10+\dfrac{616z^2}{14+\dfrac{2100z^2}{18+\dfrac{5304z^2}{22+\dfrac{11200z^2}{26+\ddots}}}}}}$$
Convergence type $P^-$ with $P=4/(3|z|)$ and $C=...$, so that
$$\sum_{k\ge0}(E^{(3)}_k-E^{(3)}_{k+2})z^k-\dfrac{p(n)}{q(n)}\sim(-1)^n\dfrac{C}{n^{4/(3|z|)}}\;.$$
$$A=1-(2/z)/n+((129z^2+162z-4)/(81z^3))/n^2+\cdots$$
Parametric family for $k\ge0$:
\begin{verbatim}
[(z)->sum(k=0,oo,(eul3frac(k)-eul3frac(k+2))*z^k),
               (3*k*z+1)*(4*n+2),z^2*n*(n+2)*(3*n+1)*(3*n+5)]
\end{verbatim}
Convergence type $P^-$ with $P=4k+4/(3z)$.
\end{cf}

\smallskip

\begin{cf}\label{6.16.58}{\ }
\begin{verbatim}
[(z)->sum(k=0,oo,(k+1)*eul3frac(k)*z^k),
[0,(2*n-1)*((9*n^2-9*n+4)*z^2+2)],[2,-9*n^4*(9*n^2-1)*z^4]]
\end{verbatim}
$$\sum_{k\ge0}(k+1)E^{(3)}_kz^k=\dfrac{2}{4z^2+2-\dfrac{72z^4}{66z^2+6-\dfrac{5040z^4}{290z^2+10-\dfrac{58320z^4}{784z^2+14-\ddots}}}}$$
Convergence type $P^+$ with $P=4/|3z|$ and $C=...$, so that
$$\sum_{k\ge0}(k+1)E^{(3)}_kz^k-\dfrac{p(n)}{q(n)}\sim\dfrac{C}{n^{4/|3z|}}\;.$$
$$A=1-(2/(3z))/n+((33z^4+162z^3-24z^2-72z+4)/(81z^3(9z^2-4)))/n^2+\cdots$$
Parametric family for $k\ge0$:
\begin{verbatim}
[(z)->sum(k=0,oo,(k+1)*eul3frac(k)*z^k),
(2*n-1)*((9*n^2-9*n+18*k^2+4)*z^2+12*k*z+2),-9*n^4*(9*n^2-1)*z^4]
\end{verbatim}
Convergence type $P^+$ with $P=4k+4/(3z)$.
\end{cf}

\smallskip

\begin{cf}\label{6.16.59}{\ }
\begin{verbatim}
[(z)->exp(sum(k=1,oo,eul3frac(k)/k*z^k)),6*[[0,1],[2*n-1,2*n+1]],
        6*[[1,2*z^2],z^2*(2*n+1)*[(3*n-1)^3,(3*n+1)^3]]]
\end{verbatim}
$$\exp\left(\sum_{k\ge1}\dfrac{E^{(3)}_k}{k}z^k\right)=\dfrac{6}{6+\dfrac{12z^2}{6+\dfrac{144z^2}{18+\dfrac{1152z^2}{18+\dfrac{3750z^2}{30+\dfrac{10290z^2}{30+\ddots}}}}}}$$
Convergence type $P^-$ with $P=4/(3z)$ and $C=...$, so that
$$\exp\left(\sum_{k\ge1}\dfrac{E^{(3)}_k}{k}z^k\right)-\dfrac{p(n)}{q(n)}\sim(-1)^n\dfrac{C}{n^{4/(3z)}}\;.$$
$$A=1-((252z^4-72z^2-16)/(81z^3(9z^2-4)))/n+\cdots$$
\end{cf}

\smallskip

\begin{cf}\label{6.16.60}{\ }
\begin{verbatim}
[(z)->sum(k=0,oo,springfrac(k)*z^k),[0,1-(2*n-1)*z],[1,-2*n^2*z^2]]
\end{verbatim}
$$\sum_{k\ge0}S_kz^k=\dfrac{1}{-z+1-\dfrac{2z^2}{-3z+1-\dfrac{8z^2}{-5z+1-\dfrac{18z^2}{-7z+1-\dfrac{32z^2}{-9z+1-\dfrac{50z^2}{-11z+1-\ddots}}}}}}$$
Nonconvergent CF.
\end{cf}

\smallskip

\begin{cf}\label{6.16.61}{\ }
\begin{verbatim}
[(z)->sum(k=0,oo,(-1)^k*springfrac(2*k)*z^(2*k)),
               [0,1],[[1,3*z^2],z^2*[(4*n)^2,(4*n+1)*(4*n+3)]]]
\end{verbatim}
$$\sum_{k\ge0}(-1)^kS_{2k}z^{2k}=\dfrac{1}{1+\dfrac{3z^2}{1+\dfrac{16z^2}{1+\dfrac{35z^2}{1+\dfrac{64z^2}{1+\dfrac{99z^2}{1+\ddots}}}}}}$$
Convergence type $P^-$ with $P=1/(2z)$ and $C=...$, so that
$$\sum_{k\ge0}(-1)^kS_{2k}z^{2k}-\dfrac{p(n)}{q(n)}\sim(-1)^n\dfrac{C}{n^{1/(2z)}}\;.$$
$$A=1+((-5z^2+1)/(16z^3-4z))/n+\cdots$$
\end{cf}

\smallskip

\begin{cf}\label{6.16.62}{\ }
\begin{verbatim}
[(z)->sum(k=0,oo,(-1)^k*springfrac(2*k)*z^(2*k)),
[0,(32*n^2-48*n+19)*z^2+1],[1,-(4*n-3)*(4*n)^2*(4*n-1)*z^4]]
\end{verbatim}
$$\sum_{k\ge0}(-1)^kS_{2k}z^{2k}=\dfrac{1}{3z^2+1-\dfrac{48z^4}{51z^2+1-\dfrac{2240z^4}{163z^2+1-\dfrac{14256z^4}{339z^2+1-\ddots}}}}$$
Convergence type $P^+$ with $P=1/(2|z|)$ and $C=...$, so that
$$\sum_{k\ge0}(-1)^kS_{2k}z^{2k}-\dfrac{p(n)}{q(n)}\sim\dfrac{C}{n^{1/(2|z|)}}\;.$$
$$A=1+((-5z^2+1)/(32z^3-8z))/n+\cdots$$
Parametric family for $k\ge0$:
\begin{verbatim}
[(z)->sum(k=0,oo,(-1)^k*springfrac(2*k)*z^(2*k)),
(32*n^2-48*n+16*k^2+19)*z^2+8*k*z+1,-(4*n-3)*(4*n)^2*(4*n-1)*z^4]
\end{verbatim}
Convergence type $P^+$ with $P=2k+1/(2z)$.
\end{cf}

This is the contraction of the previous CF.

\smallskip

\begin{cf}\label{6.16.63}{\ }
\begin{verbatim}
[(z)->sum(k=0,oo,(-1)^k*springfrac(2*k+1)*z^(2*k)),
[0,(32*n^2-32*n+11)*z^2+1],[1,-(4*n+1)*(4*n)^2*(4*n-1)*z^4]]
\end{verbatim}
$$\sum_{k\ge0}(-1)^kS_{2k+1}z^{2k}=\dfrac{1}{11z^2+1-\dfrac{240z^4}{75z^2+1-\dfrac{4032z^4}{203z^2+1-\dfrac{20592z^4}{395z^2+1-\ddots}}}}$$
Convergence type $P^+$ with $P=1/(2|z|)$ and $C=...$, so that
$$\sum_{k\ge0}(-1)^kS_{2k+1}z^{2k}-\dfrac{p(n)}{q(n)}\sim\dfrac{C}{n^{1/(2|z|)}}\;.$$
$$A=1-(1/(4z))/n+((709z^4+768z^3-62z^2-48z+1)/(1536z^3(16z^2-1)))/n^2+\cdots$$
Parametric family for $k\ge0$:
\begin{verbatim}
[(z)->sum(k=0,oo,(-1)^k*springfrac(2*k+1)*z^(2*k)),
(32*n^2-32*n+16*k^2+11)*z^2+8*k*z+1,-(4*n+1)*(4*n)^2*(4*n-1)*z^4]
\end{verbatim}
Convergence type $P^+$ with $P=2k+1/(2z)$.
\end{cf}

\smallskip

\begin{cf}\label{6.16.64}{\ }
\begin{verbatim}
[(z)->sum(k=0,oo,(-1)^(k*(k-1)/2)*springfrac(k)*z^k),
               [0,1-z],[1,4*n^2*z^2]]
\end{verbatim}
$$\sum_{k\ge0}(-1)^{k(k-1)/2}S_kz^k=\dfrac{1}{-z+1+\dfrac{4z^2}{-z+1+\dfrac{16z^2}{-z+1+\dfrac{36z^2}{-z+1+\dfrac{64z^2}{-z+1+\ddots}}}}}$$
Convergence type $P^-$ with $P=(1-z)/(2z)$ and $C=...$, so that
$$\sum_{k\ge0}(-1)^{k(k-1)/2}S_kz^k-\dfrac{p(n)}{q(n)}\sim(-1)^n\dfrac{C}{n^{(1-z)/(2z)}}\;.$$
$$A=1+((1-z)/(4z))/n-((9z+1)(z^2-1)/(384z^3))/n^2+\cdots$$
Parametric family for $k\ge0$:
\begin{verbatim}
[(z)->sum(k=0,oo,(-1)^(k*(k-1)/2)*springfrac(k)*z^k),
                                  (4*k-1)*z+1,4*n^2*z^2]
\end{verbatim}
Convergence type $P^-$ with $P=2k+(1-z)/(2z)$.
\end{cf}

\smallskip

\begin{cf}\label{6.16.65}{\ }
\begin{verbatim}
[(z)->sum(k=0,oo,(-1)^(k*(k-1)/2)*springfrac(k)*z^k),
               [0,z+1,4*z],[2,(z*(4*n-3)+1)^2]
\end{verbatim}
$$\sum_{k\ge0}(-1)^{k(k-1)/2}S_kz^k=\dfrac{2}{z+1+\dfrac{z^2+2z+1}{4z+\dfrac{25z^2+10z+1}{4z+\dfrac{81z^2+18z+1}{4z+\dfrac{169z^2+26z+1}{4z+\ddots}}}}}$$
Convergence type $P^-$ with $P=1$ and $C=1/(4z)$, so that
$$\sum_{k\ge0}(-1)^{k(k-1)/2}S_kz^k-\dfrac{p(n)}{q(n)}\sim(-1)^n\dfrac{1/(4z)}{n}\;.$$
$$A=1+((z-1)/(4z))/n+((-3z^2-2z+1)/(16z^2))/n^2+\cdots$$
Series:
$$\sum_{k\ge0}(-1)^{k(k-1)/2}S_kz^k=2\sum_{n\ge1}\dfrac{(-1)^{n+1}}{(4n-3)z+1}$$
Parametric family for $k\ge0$:
\begin{verbatim}
[(z)->sum(k=0,oo,(-1)^(k*(k-1)/2)*springfrac(k)*z^k),
               (8*k+4)*z,(z*(4*n-3)+1)^2]
\end{verbatim}
Convergence type $P^-$ with $P=2k+1$.
\end{cf}

\smallskip

\begin{cf}\label{6.16.66}{\ }
\begin{verbatim}
[(z)->sum(k=0,oo,(-1)^(k*(k-1)/2)*springfrac(k)*z^k),
    [0,(4*n-3)*z+1],[[2,(z+1)^2],[(4*z*n)^2,(4*z*n+z+1)^2]]]
\end{verbatim}
$$\sum_{k\ge0}(-1)^{k(k-1)/2}S_kz^k=\dfrac{2}{z+1+\dfrac{z^2+2z+1}{5z+1+\dfrac{16z^2}{9z+1+\dfrac{25z^2+10z+1}{13z+1+\ddots}}}}$$
Convergence type $E$ with $E=-(1+\sqrt{2})^2$, $P=0$, and $C=...$, so that
$$\sum_{k\ge0}(-1)^{k(k-1)/2}S_kz^k-\dfrac{p(n)}{q(n)}\sim(-1)^n\dfrac{C}{(1+\sqrt{2})^{2n}}\;.$$
$$A=1+((3z^2-6z+1)d/(16z^2))/n+\cdots$$
\end{cf}

\smallskip

\begin{cf}\label{6.16.67}{\ }
\begin{verbatim}
[(z)->sum(k=0,oo,(-1)^k*(springfrac(2*k)+springfrac(2*k+2))*z^(2*k)),
[0,(32*n^2-16*n-1)*z^2+1],[4,-16*n*(n+1)*(16*n^2-1)*z^4]]
\end{verbatim}
$$\sum_{k\ge0}(-1)^k(S_{2k}+S_{2k+2})z^{2k}=\dfrac{4}{15z^2+1-\dfrac{480z^4}{95z^2+1-\dfrac{6048z^4}{239z^2+1-\dfrac{27456z^4}{447z^2+1-\ddots}}}}$$
Convergence type $P^+$ with $P=1/(2|z|)$ and $C=...$, so that
$$\sum_{k\ge0}(-1)^k(S_{2k}+S_{2k+2})z^{2k}-\dfrac{p(n)}{q(n)}\sim\dfrac{C}{n^{1/(2|z|)}}\;.$$
\begin{align*}
  A&=1-((9z^2-3)/(8z(4z^2-1)))/n\\
  &\phantom{=}+((403z^4-160z^3-234z^2+112z-1)/(1536z^3(2z-1)^2))/n^2+\cdots
\end{align*}
Parametric family for $k\ge0$:
\begin{verbatim}
[(z)->sum(k=0,oo,(-1)^k*(springfrac(2*k)+springfrac(2*k+2))*z^(2*k)),
(32*n^2-16*n+16*k^2-1)*z^2+8*k*z+1,-16*n*(n+1)*(16*n^2-1)*z^4]
\end{verbatim}
Convergence type $P^+$ with $P=2k+1/(2z)$.
\end{cf}

There exist an Ap\'ery dual for the above analogous to \ref{6.16.69},
too complicated to give here.

\smallskip

\begin{cf}\label{6.16.68}{\ }
\begin{verbatim}
[(z)->sum(k=0,oo,(-1)^k*(springfrac(2*k+1)+springfrac(2*k+3))*z^(2*k)),
[0,(32*n^2-1)*z^2+1],[12,-16*n*(n+1)*(4*n+1)*(4*n+3)*z^4]]
\end{verbatim}
$$\sum_{k\ge0}(-1)^k(S_{2k+1}+S_{2k+3})z^{2k}=\dfrac{12}{31z^2+1-\dfrac{(280/3)z^4}{127z^2+1-\dfrac{9504z^4}{287z^2+1-\dfrac{37440z^4}{511z^2+1-\ddots}}}}$$
Convergence type $P^+$ with $P=1/(2|z|)$ and $C=...$, so that
$$\sum_{k\ge0}(-1)^k(S_{2k+1}+S_{2k+3})z^{2k}-\dfrac{p(n)}{q(n)}\sim\dfrac{C}{n^{1/(2|z|)}}\;.$$
$$A=1-(1/(2z))/n+((5773z^4+3072z^3-374z^2-192z+1)/(1536z^3(16z^2-1)))/n^2+\cdots$$
Parametric family for $k\ge0$:
\begin{verbatim}
[(z)->sum(k=0,oo,(-1)^k*(springfrac(2*k+1)+springfrac(2*k+3))*z^(2*k)),
(32*n^2+16*k^2-1)*z^2+8*k*z+1,-16*n*(n+1)*(4*n+1)*(4*n+3)*z^4]
\end{verbatim}
Convergence type $P^+$ with $P=2k+1/(2z)$.
\end{cf}

\smallskip

\begin{cf}\label{6.16.69}{\ }
\begin{verbatim}
[(z)->sum(k=0,oo,(-1)^k*(springfrac(2*k+1)+springfrac(2*k+3))*z^(2*k)),
[0,(3*z+1)*(5*z+1),16*z*(4*(n-1)*z+1)],
[24,((4*n-5)*z+1)*((4*n-3)*z+1)*((4*n-1)*z+1)*((4*n+1)*z+1)]]
\end{verbatim}
$$\sum_{k\ge0}(-1)^k(S_{2k+1}+S_{2k+3})z^{2k}=\dfrac{24}{15z^2+8z+1-\dfrac{15z^4+8z^3-14z^2-8z-1}{64z^2+16z-\ddots}}$$
Convergence type $P^-$ with $P=4$ and $C=3(1-z^2)/(4z^4)$, so that
$$\sum_{k\ge0}(-1)^k(S_{2k+1}+S_{2k+3})z^{2k}-\dfrac{p(n)}{q(n)}\sim(-1)^n\dfrac{3(1-z^2)/(4z^4)}{n^4}\;.$$
$$A=1-(1/z)/n+(5(1-z^2)/(8z^2))/n^2+((45z^2-5)/(16z^3))/n^3+\cdots$$
Series:
\begin{align*}&\sum_{k\ge0}(-1)^k(S_{2k+1}+S_{2k+3})z^{2k}\\&=24(1-z^2)\sum_{n\ge1}\dfrac{(-1)^{n+1}}{((4n+1)z+1)((4n-1)z+1)((4n-3)z+1)((4n-5)z+1)}\end{align*}
Parametric family for $k\ge0$:
\begin{verbatim}
[(z)->sum(k=0,oo,(-1)^k*(springfrac(2*k+1)+springfrac(2*k+3))*z^(2*k)),
(k+1)*16*z*(4*(n-1)*z+1),
((4*n-5)*z+1)*((4*n-3)*z+1)*((4*n-1)*z+1)*((4*n+1)*z+1)]
\end{verbatim}
Convergence type $P^-$ with $P=4k+4$.
\end{cf}

\smallskip

\begin{cf}\label{6.16.70}{\ }
\begin{verbatim}
[(z)->exp(sum(k=1,oo,(-1)^(k-1)*springfrac(2*k)/(2*k)*z^(2*k))),
                              [1,2],[(4*n+1)*(4*n+3)*z^2]]
\end{verbatim}
$$\exp\left(\sum_{k\ge1}(-1)^{k-1}\dfrac{S_{2k}}{2k}z^{2k}\right)=1+\dfrac{3z^2}{2+\dfrac{35z^2}{2+\dfrac{99z^2}{2+\dfrac{195z^2}{2+\dfrac{323z^2}{2+\dfrac{483z^2}{2+\ddots}}}}}}$$
Convergence type $P^-$ with $P=1/(2|z|)$ and $C=...$, so that
$$\exp\left(\sum_{k\ge1}(-1)^{k-1}\dfrac{S_{2k}}{2k}z^{2k}\right)-\dfrac{p(n)}{q(n)}\sim(-1)^n\dfrac{C}{n^{1/(2|z|)}}\;.$$
$$A=1-(1/(2z))/n+((55z^2+48z-1)/(384z^3))/n^2+\cdots$$
Parametric family for $k\ge0$:
\begin{verbatim}
[(z)->exp(sum(k=1,oo,(-1)^(k-1)*springfrac(2*k)/(2*k)*z^(2*k))),
                   8*k*z+2,(4*n+1)*(4*n+3)*z^2]
\end{verbatim}
Convergence type $P^-$ with $P=2k+1/(2z)$.
\end{cf}

\smallskip

\begin{cf}\label{6.16.71}{\ }
\begin{verbatim}
[(z)->exp(sum(k=1,oo,(-1)^(k-1)*springfrac(2*k)/(2*k)*z^(2*k))),
[1,4*z+1,8*z],[3*z^2,((4*n-3)*z+1)*((4*n-1)*z+1)]]
\end{verbatim}
$$\exp\left(\sum_{k\ge1}(-1)^{k-1}\dfrac{S_{2k}}{2k}z^{2k}\right)=1+\dfrac{3z^2}{4z+1+\dfrac{3z^2+4z+1}{8z+\dfrac{35z^2+12z+1}{8z+\dfrac{99z^2+20z+1}{8z+\dfrac{195z^2+28z+1}{8z+\ddots}}}}}$$
Convergence type $P^-$ with $P=2$ and $C=...$, so that
$$\exp\left(\sum_{k\ge1}(-1)^{k-1}\dfrac{S_{2k}}{2k}z^{2k}\right)-\dfrac{p(n)}{q(n)}\sim(-1)^n\dfrac{C}{n^2}\;.$$
$$A=1-(1/(2z))/n+((-19z^2+6)/(32z^2))/n^2+\cdots$$
Series:
\begin{align*}\exp\left(\sum_{k\ge1}(-1)^{k-1}\dfrac{S_{2k}}{2k}z^{2k}\right)&=1+\dfrac{24z^3}{(5z+1)(7z+1)}\cdot\\&\phantom{=}\cdot\sum_{n\ge0}\dfrac{((z+1)/(8z))_n((3z+1)/(8z))_n}{((13z+1)/(8z))_n((15z+1)/(8z))_n}\\
  \exp\left(-\sum_{k\ge1}(-1)^{k-1}\dfrac{S_{2k}}{2k}z^{2k}\right)&=-\dfrac{4z-1}{(z-1)(3z-1)}+\dfrac{24z^3}{(z+1)(3z+1)(z-1)(3z-1)}\cdot\\&\phantom{=}\cdot\sum_{n\ge0}\dfrac{((1-3z)/(8z))_n((1-z)/(8z))_n}{((1+9z)/(8z))_n((1+11z)/(8z))_n}\end{align*}
Parametric family for $k\ge0$:
\begin{verbatim}
[(z)->exp(sum(k=1,oo,(-1)^(k-1)*springfrac(2*k)/(2*k)*z^(2*k))),
(k+1)*8*z,((4*n-3)*z+1)*((4*n-1)*z+1)]
\end{verbatim}
Convergence type $P^-$ with $P=2k+2$.
\end{cf}

\smallskip

\begin{cf}\label{6.16.72}{\ }
\begin{verbatim}
[(z)->exp(sum(k=1,oo,(-1)^(k-1)*springfrac(2*k)/(2*k)*z^(2*k))),
[4*n*z+1],[[(4*n+1)*(4*n+3)*z^2,((4*n+1)*z+1)*((4*n+3)*z+1)]]]
\end{verbatim}
$$\exp\left(\sum_{k\ge1}(-1)^{k-1}\dfrac{S_{2k}}{2k}z^{2k}\right)=1+\dfrac{3z^2}{4z+1+\dfrac{3z^2+4z+1}{8z+1+\dfrac{35z^2}{12z+1+\dfrac{35z^2+12z+1}{16z+1+\ddots}}}}$$
Convergence type $E$ with $E=-(1+\sqrt{2})^2$, $P=0$, and $C=...$, so that
$$\exp\left(\sum_{k\ge1}(-1)^{k-1}\dfrac{S_{2k}}{2k}z^{2k}\right)-\dfrac{p(n)}{q(n)}\sim(-1)^n\dfrac{C}{(1+\sqrt{2})^{2n}}\;.$$
$$A=1+((12z^2-8z+1)d/(16z^2))/n+\cdots$$
\end{cf}

\smallskip

\begin{cf}\label{6.16.73}{\ }
\begin{verbatim}
[(z)->exp(sum(k=1,oo,(-1)^(k-1)*springfrac(2*k-1)/(2*k-1)*z^(2*k-1))),
[1,15*z^2-2*z+4,(128*n^2-192*n+78)*z^2+4],
[4*z,-(8*n-5)*(8*n-3)*(8*n-1)*(8*n+1)*z^4]]
\end{verbatim}
$$\exp\left(\sum_{k\ge1}(-1)^{k-1}\dfrac{S_{2k-1}}{2k-1}z^{2k-1}\right)=1+\dfrac{4z}{15z^2-2z+4-\dfrac{945z^4}{206z^2+4-\dfrac{36465z^4}{654z^2+4-\ddots}}}$$
Convergence type $P^+$ with $P=1/(2z)$ and $C=...$, so that
$$\exp\left(\sum_{k\ge1}(-1)^{k-1}\dfrac{S_{2k-1}}{2k-1}z^{2k-1}\right)-\dfrac{p(n)}{q(n)}\sim\dfrac{C}{n^{1/(2z)}}\;.$$
$$A=1-(1/(8z))/n-((17z^2-12z+1)/(1536z^3))/n^2+\cdots$$
Parametric family for $k\ge0$:
\begin{verbatim}
[(z)->exp(sum(k=1,oo,(-1)^(k-1)*springfrac(2*k-1)/(2*k-1)*z^(2*k-1))),
(128*n^2-192*n+64*k^2+78)*z^2+32*k*z+4,
-(8*n-5)*(8*n-3)*(8*n-1)*(8*n+1)*z^4]
\end{verbatim}
Convergence type $P^+$ with $P=2k+1/(2z)$.
\end{cf}

\smallskip

\begin{cf}\label{6.16.74}{\ }
\begin{verbatim}
[(z)->exp(sum(k=1,oo,(-1)^(k-1)*springfrac(2*k-1)/(2*k-1)*z^(2*k-1))),
[1,z+1,4*z],[2*z,((4*n-3)*z+1)*((4*n-1)*z+1)]]
\end{verbatim}
$$\exp\left(\sum_{k\ge1}(-1)^{k-1}\dfrac{S_{2k-1}}{2k-1}z^{2k-1}\right)=1+\dfrac{2z}{z+1+\dfrac{3z^2+4z+1}{4z+\dfrac{35z^2+12z+1}{4z+\dfrac{99z^2+20z+1}{4z+\ddots}}}}$$
Convergence type $P^-$ with $P=1$ and $C=...$, so that
$$\exp\left(\sum_{k\ge1}(-1)^{k-1}\dfrac{S_{2k-1}}{2k-1}z^{2k-1}\right)-\dfrac{p(n)}{q(n)}\sim(-1)^n\dfrac{C}{n}\;.$$
$$A=1-(1/(4z))/n+((-15z^2+4)/(64z^2))/n^2+\cdots$$
Series:
\begin{align*}\exp\left(\sum_{k\ge1}(-1)^{k-1}\dfrac{S_{2k-1}}{2k-1}z^{2k-1}\right)&=1+\dfrac{8z^2}{(z+1)(7z+1)}\cdot\\&\phantom{=}\cdot\sum_{n\ge0}\dfrac{((3z+1)/(8z))_n((5z+1)/(8z))_n}{((9z+1)/(8z))_n((15z+1)/(8z))_n}\\
 \exp\left(-\sum_{k\ge1}(-1)^{k-1}\dfrac{S_{2k-1}}{2k-1}z^{2k-1}\right)&=\dfrac{(3z-1)}{z-1}-\dfrac{8z^2}{(z-1)(3z+1)}\cdot\\&\phantom{=}\cdot\sum_{n\ge0}\dfrac{((1-z)/(8z))_n((1+z)/(8z))_n}{((5z+1)/(8z))_n((11z+1)/(8z))_n}\end{align*}
Parametric family for $k\ge0$:
\begin{verbatim}
[(z)->exp(sum(k=1,oo,(-1)^(k-1)*springfrac(2*k-1)/(2*k-1)*z^(2*k-1))),
                   (8*k+4)*z,((4*n-3)*z+1)*((4*n-1)*z+1)]
\end{verbatim}
Convergence type $P^-$ with $P=2k+1$.
\end{cf}

\smallskip

\begin{cf}\label{6.16.75}{\ }
\begin{verbatim}
[(z)->sum(k=1,oo,h4frac(4*k)*z^(2*k)),
[0,4*n+1],z^2*[[1/2,135],[2*n*(2*n+1)^3*(4*n+1)^2,(n+1)^2*(4*n+3)^3*(4*n+5)]]]
\end{verbatim}
$$\sum_{k\ge1}H^{(4)}_{4k}z^{2k}=\dfrac{1/2z^2}{5+\dfrac{135z^2}{9+\dfrac{1350z^2}{13+\dfrac{12348z^2}{17+\dfrac{40500z^2}{21+\dfrac{155727z^2}{25+\ddots}}}}}}$$
Convergence type $P^-$ with $P=3/2$ and $C=...$, so that
$$\sum_{k\ge1}H^{(4)}_{4k}z^{2k}-\dfrac{p(n)}{q(n)}\sim(-1)^n\dfrac{C}{n^{3/2}}$$
$$A=1-(15/8)/n+((2023z^2-768)/896z^2)/n^2+\cdots$$
\end{cf}

\smallskip

\begin{cf}\label{6.16.76}{\ }
\begin{verbatim}
[(z)->sum(k=1,oo,h6frac(6*k)*z^(2*k)),
[0,8*(6*n+1)],z^2*[2/3,n*(n+1)^2*(3*n+1)^2*(3*n+2)^2*(3*n+4)]]
\end{verbatim}
$$\sum_{k\ge1}H^{(6)}_{6k}z^{2k}=\dfrac{2/3z^2}{56+\dfrac{11200z^2}{104+\dfrac{564480z^2}{152+\dfrac{7550400z^2}{200+\dfrac{52998400z^2}{248+\ddots}}}}}$$
Convergence type $P^+$ with $P=2$ and $C=...$, so that
$$\sum_{k\ge1}H^{(6)}_{6k}z^{2k}-\dfrac{p(n)}{q(n)}\sim\dfrac{C}{n^2}$$
$A=1+\cdots$
\end{cf}

\medskip

\section{Incomplete Gamma Functions}

\medskip

We recall that
$$\ga(a,z)=\int_z^\infty t^ae^{-t}\,\dfrac{dt}{t}\text{\quad and\quad}
\gac(a,z)=\G(a)-\ga(a,z)=\int_0^z t^ae^{-t}\,\dfrac{dt}{t}\;.$$

\smallskip

\begin{cf}\label{6.17.1}{\ }
\begin{verbatim}
[(a,z)->incgam(a,z)/(z^a*exp(-z)),[0,z-a+2*n-1],[1,-n*(n-a)]]
\end{verbatim}
$$\dfrac{\ga(a,z)}{z^ae^{-z}}=\dfrac{1}{z-a+1+\dfrac{a-1}{z-a+3+\dfrac{2a-4}{z-a+5+\dfrac{3a-9}{z-a+7+\dfrac{4a-16}{z-a+9+\ddots}}}}}$$
Convergence type $D^+$ with $D=16z$ and $C=...$, so that
$$\dfrac{\ga(a,z)}{z^ae^{-z}}-\dfrac{p(n)}{q(n)}\sim\dfrac{C}{e^{4\sqrt{nz}}}\;.$$
$$A=1-((4z^2+24(1-a)z+3(1-4a^2))/(24z^{1/2}))/n^{1/2}+\cdots$$
\end{cf}

\smallskip

\begin{cf}\label{6.17.2}{\ }
\begin{verbatim}
[(a,z)->incgam(a,z)/(z^(a-1)*exp(-z)),[1,z-a+2*n],[a-1,-n*(n+1-a)]]
\end{verbatim}
$$\dfrac{\ga(a,z)}{z^{a-1}e^{-z}}= 1+\dfrac{a-1}{z-a+2+\dfrac{a-2}{z-a+4+\dfrac{2a-6}{z-a+6+\dfrac{3a-12}{z-a+8+\dfrac{4a-20}{z-a+10+\ddots}}}}}$$
Convergence type $D^+$ with $D=16z$ and $C=...$, so that
$$\dfrac{\ga(a,z)}{z^{a-1}e^{-z}}-\dfrac{p(n)}{q(n)}\sim\dfrac{C}{e^{4\sqrt{nz}}}\;.$$
$$A=1-((4z^2+24(2-a)z-3(2a-1)(2a-3))/(24z^{1/2}))/n^{1/2}+\cdots$$
\end{cf}

\smallskip

\begin{verbatim}
poch(a,k)=prod(j=0,k-1,a+j);
\end{verbatim}

\smallskip

\begin{cf}\label{6.17.3}{\ }
\begin{verbatim}
[(a,l,z)->incgam(a,z)/(z^a*exp(-z))+sum(k=0,l-1,z^k/poch(a,k+1)),
[0,2*n+z-a-l-1],[z^l/poch(a,l),-n*(n-a-l)]]
\end{verbatim}
$$\dfrac{\ga(a,z)}{z^ae^{-z}}+\sum_{k=0}^{l-1}\dfrac{z^k}{(a)_{k+1}}
=\dfrac{z^l/(a)_l}{-a-l+z+1+\dfrac{a+l-1}{-a-l+z+3+\dfrac{2a+2l-4}{-a-l+z+5+\ddots}}}$$
Convergence type $D^+$ with $D=16z$ and $C=...$, so that
$$\dfrac{\ga(a,z)}{z^ae^{-z}}+\sum_{k=0}^{l-1}\dfrac{z^k}{(a)_{k+1}}-\dfrac{p(n)}{q(n)}\sim\dfrac{C}{e^{4\sqrt{nz}}}\;.$$
$$A=1-((4z^2+24(1-(a+l))z+3(1-4(a+l)^2))/(24z^{1/2}))/n^{1/2}+\cdots$$
\end{cf}

\smallskip

\begin{cf}\label{6.17.4}{\ }
\begin{verbatim}
[(a,l,z)->incgam(a,z)/(z^a*exp(-z))+sum(k=1,l,poch(1-a,k-1)/(-z)^k),
[0,2*n+z-a+l-1],[poch(1-a,l)/(-z)^l,-n*(n-a+l)]]
\end{verbatim}
$$\dfrac{\ga(a,z)}{z^ae^{-z}}+\sum_{k=1}^l\dfrac{(1-a)_{k-1}}{(-z)^k}=
\dfrac{(1-a)_l/(-z)^l}{-a+l+z+1+\dfrac{a-l-1}{-a+l+z+3+\dfrac{2a-2l-4}{-a+l+z+5+\ddots}}}$$                            
Convergence type $D^+$ with $D=16z$ and $C=...$, so that
$$\dfrac{\ga(a,z)}{z^ae^{-z}}+\sum_{k=1}^l\dfrac{(1-a)_{k-1}}{(-z)^k}-\dfrac{p(n)}{q(n)}\sim\dfrac{C}{e^{4\sqrt{nz}}}\;.$$
$$A=1-((4z^2+24(1-a+l)z-3(2a-2l+1)(2a-2l-1))/(24z^{1/2}))/n^{1/2}+\cdots$$
\end{cf}

\smallskip

\begin{cf}\label{6.17.5}{\ }
\begin{verbatim}
[(a,z)->incgamc(a,z)/(z^a*exp(-z)),[0,n+a-1],[[1,-a*z],z*[n,-(n+a)]]]
\end{verbatim}
$$\dfrac{\gac(a,z)}{z^ae^{-z}}=\dfrac{1}{a-\dfrac{za}{a+1+\dfrac{z}{a+2-\dfrac{za+z}{a+3+\dfrac{2z}{a+4-\dfrac{za+2z}{a+5+\ddots}}}}}}$$
Convergence type $F^1$ with $E=-2/z$, $P=a-1/2$, and $C=...$, so that
$$\dfrac{\gac(a,z)}{z^ae^{-z}}-\dfrac{p(n)}{q(n)}\sim(-1)^n\dfrac{C}{n!(2/z)^nn^{a-1/2}}\;.$$
$$A=1+(z^2/4+(1-a)z+1/2-a)/n+\cdots$$
\end{cf}

\smallskip

\begin{cf}\label{6.17.6}{\ }
\begin{verbatim}
[(a,z)->incgamc(a,z)/(z^a*exp(-z)),[0,a-z+n-1],[1,n*z]]
\end{verbatim}
$$\dfrac{\gac(a,z)}{z^ae^{-z}}=\dfrac{1}{-z+a+\dfrac{z}{-z+a+1+\dfrac{2z}{-z+a+2+\dfrac{3z}{-z+a+3+\dfrac{4z}{-z+a+4+\ddots}}}}}$$
Convergence type $F^1$ with $E=-1/z$, $P=2a-1$, and $C=...$, so that
$$\dfrac{\gac(a,z)}{z^ae^{-z}}-\dfrac{p(n)}{q(n)}\sim(-1)^n\dfrac{C}{n!(1/z)^nn^{2a-1}}\;.$$
$$A=1-((2a-1)z+a^2)/n+((2a^2-a)z^2+(2a^3+a^2-a+1)z+(3a^4+2a^3+a)/6)/n^2+\cdots$$
\end{cf}

\smallskip

\begin{cf}\label{6.17.7}{\ }
\begin{verbatim}
[(a,z)->incgamc(a,z)/(z^a*exp(-z)),[0,a,n+z+a-1],[1,-z*(n+a-1)]]
\end{verbatim}
$$\dfrac{\gac(a,z)}{z^ae^{-z}}=\dfrac{1}{a-\dfrac{az}{z+a+1-\dfrac{(a+1)z}{z+a+2-\dfrac{(a+2)z}{z+a+3-\dfrac{(a+3)z}{z+a+4-\dfrac{(a+4)z}{z+a+5-\ddots}}}}}}$$
Convergence type $F^1$ with $E=1/z$, $P=a$, and $C=\G(a)$, so that
$$\dfrac{\gac(a,z)}{z^ae^{-z}}-\dfrac{p(n)}{q(n)}\sim\dfrac{\G(a)}{n!(1/z)^nn^a}\;.$$
$$A=1+(a(a+1)/2+z)/n+\cdots$$
Series:
$$\dfrac{\gac(a,z)}{z^ae^{-z}}=\sum_{n\ge0}\dfrac{z^n}{(a)_{n+1}}\;.$$
\end{cf}

\smallskip

\begin{cf}\label{6.17.8}{\ }
\begin{verbatim}
[(a,z)->incgamc(a,z)/z^a,[1/a,a+1,n^2+(a-z)*n+(1-a)*z],
                         [-z,n*z*(n+a)^2]]
\end{verbatim}
$$\dfrac{\gac(a,z)}{z^a}=\dfrac{1}{a}-\dfrac{z}{a+1+\dfrac{(a^2+2a+1)z}{(-a-1)z+2a+4+\dfrac{(2a^2+8a+8)z}{(-a-2)z+3a+9+\dfrac{(3a^2+18a+27)z}{(-a-3)z+4a+16+\ddots}}}}$$
Convergence type $F^1$ with $E=-1/z$, $P=2$, and $C=-z$, so that
$$\dfrac{\gac(a,z)}{z^a}-\dfrac{p(n)}{q(n)}\sim(-1)^{n+1}\dfrac{1}{n!(1/z)^{n+1}n^2}\;.$$
$$A=1-(a+2+z/2)/n+((a/2+7/4)z+a^2+3a+3)/n^2+\cdots$$
Series:
$$\dfrac{\gac(a,z)}{z^a}=\sum_{n\ge0}(-1)^n\dfrac{z^n}{n!(n+a)}\;.$$
\end{cf}

\smallskip

\begin{cf}\label{6.17.9}{\ }
\begin{verbatim}
[(z)->incgamc(z,z)/(z^(z-1)*exp(-z)),[n+1],[z*(n+2)]]
[(z)->incgamc(z,z)/(z^(z-1)*exp(-z)),[1],[z/(n+1)]]
\end{verbatim}
$$\dfrac{\gac(z,z)}{z^{z-1}e^{-z}}=1+\dfrac{2z}{2+\dfrac{3z}{3+\dfrac{4z}{4+\dfrac{5z}{5+\dfrac{6z}{6+\dfrac{7z}{7+\ddots}}}}}}=1+\dfrac{z}{1+\dfrac{z/2}{1+\dfrac{z/3}{1+\dfrac{z/4}{1+\dfrac{z/5}{1+\dfrac{z/6}{1+\ddots}}}}}}$$
Convergence type $F^1$ with $E=-1/z$, $P=2z+1$, and $C=...$, so that
$$\dfrac{\gac(z,z)}{z^{z-1}e^{-z}}-\dfrac{p(n)}{q(n)}\sim(-1)^n\dfrac{C}{n!(1/z)^nn^{2z+1}}\;.$$
$$A=1-(3z^2+3z+1)/n+((27z^4+74z^3+72z^2+37z+6)/6)/n^2+\cdots$$
\end{cf}

\smallskip

\begin{cf}\label{6.17.10}{\ }
\begin{verbatim}
[(z)->incgamc(z,z)/(z^(z-1)*exp(-z)),[1+z,n+2*z+1],[-z*(n+z)]]
\end{verbatim}
$$\dfrac{\gac(z,z)}{z^{z-1}e^{-z}}=z+1-\dfrac{z^2}{2z+2-\dfrac{z^2+z}{2z+3-\dfrac{z^2+2z}{2z+4-\dfrac{z^2+3z}{2z+5-\dfrac{z^2+4z}{2z+6-\dfrac{z^2+5z}{2z+7-\ddots}}}}}}$$
Convergence type $F^1$ with $E=1/z$, $P=z+3$, and $C=...$, so that
$$\dfrac{\gac(z,z)}{z^{z-1}e^{-z}}-\dfrac{p(n)}{q(n)}\sim\dfrac{C}{n!(1/z)^nn^{z+3}}\;.$$
$$A=1-((z^2+z+8)/2)/n+((3z^4+10z^3+57z^2-22z+264)/24)/n^2+\cdots$$
Series:
$$\dfrac{\gac(z,z)}{z^{z-1}e^{-z}}=z+1-\dfrac{z^2}{z+1}\sum_{n\ge0}\dfrac{1}{(n+1)(n+2)(z+2)_n}z^n$$
\end{cf}

\smallskip

Thanks to the identity
$$\ga(\nu,z)=z^{\nu-1}e^{-z}{}_2F_0(1,1-\nu;;-1/z)\;,$$
we can obtain more CFs for $\ga$ from those of ${}_2F_0$.

\smallskip

{\bf Comment:} Although $\ga$ and $\ga_c$ are trivially related by
$\ga_c(s,z)=\G(s)-\ga(s,z)$, the continued fractions for $\ga$ converge
\emph{subexponentially} because implicitly $z$ is large, while the ones
for $\ga_c$ converge factorially, because implicitly $z$ is small, so
the CFs for $\ga_c$ are not much help for the computation of $\ga$ for $z$
large.

\medskip

We recall that the generalized exponential integral is:
$$E_{\nu}(z)=z^{\nu-1}\int_z^\infty \dfrac{e^{-t}}{t^{\nu}}\,dt=z^{\nu-1}\ga(1-\nu,z)\;,$$
so the CFs that we give are restatements of CFs for $\ga(a,z)$.

\smallskip

\begin{verbatim}
eintn(nu,z)=z^(nu-1)*intnum(t=z,[oo,1],exp(-t)/t^nu);
\end{verbatim}

\smallskip

\begin{cf}\label{6.17.11}{\ }
\begin{verbatim}
[(nu,z)->exp(z)*eintn(nu,z),[0,nu+z+2*n-2],[1,-n*(nu+n-1)]]
\end{verbatim}
$$e^zE_{\nu}(z)=\dfrac{1}{\nu+z-\dfrac{\nu}{\nu+z+2-\dfrac{2\nu+2}{\nu+z+4-\dfrac{3\nu+6}{\nu+z+6-\dfrac{4\nu+12}{\nu+z+8-\ddots}}}}}$$
Convergence type $D^+$ with $D=16z$ and $C=...$, so that
$$e^zE_{\nu}(z)-\dfrac{p(n)}{q(n)}\sim\dfrac{C}{e^{4\sqrt{zn}}}\;.$$
$$A=1-((4z^2+24\nu z-12\nu^2+24\nu-9)/(24z^{1/2}))/n^{1/2}+\cdots$$
\end{cf}

\smallskip

\begin{cf}\label{6.17.12}{\ }
\begin{verbatim}
[(nu,z)->z*exp(z)*eintn(nu,z),[1,nu+z+2*n-1],[-nu,-n*(nu+n)]]
\end{verbatim}
$$ze^zE_{\nu}(z)=1-\dfrac{\nu}{\nu+z+1-\dfrac{\nu+1}{\nu+z+3-\dfrac{2\nu+4}{\nu+z+5-\dfrac{3\nu+9}{\nu+z+7-\dfrac{4\nu+16}{\nu+z+9-\ddots}}}}}$$
Convergence type $D^+$ with $D=16z$ and $C=...$, so that
$$ze^zE_{\nu}(z)-\dfrac{p(n)}{q(n)}\sim\dfrac{C}{e^{4\sqrt{zn}}}\;.$$
$$A=1-((4z^2+24(\nu+1)z-12\nu^2+3)/(24z^{1/2}))/n^{1/2}+\cdots$$
\end{cf}

\smallskip

Thanks to the identity
$$E_{\nu}(z)=\dfrac{e^{-z}}{z}{}_2F_0(1,\nu;;-1/z)\;,$$
we can obtain more CFs for $E_{\nu}$ from those for ${}_2F_0$.

\section{Error Function and Related}

\medskip

We recall that
$$\erfc(x)=\dfrac{2}{\sqrt{\pi}}\int_x^\infty e^{-t^2}\,dt\text{\quad and\quad}\erf(x)=1-\erfc(x)=\dfrac{2}{\sqrt{\pi}}\int_0^x e^{-t^2}\,dt\;.$$
Note that $\erfc(x)=\ga(1/2,x^2)/\sqrt{\pi}$ and
$\erf(x)=\gac(1/2,x^2)/\sqrt{\pi}$ for $x\ge0$, so most of the CFs that we
give are specializations of those for $\ga$ and $\gac$.

\smallskip

\begin{verbatim}
erf(x)=1-erfc(x);
\end{verbatim}

\smallskip

\begin{cf}\label{6.18.1}{\ }
\begin{verbatim}
[(z)->sqrt(Pi)*exp(z^2)*erf(z),
[0,2*n-1],[[2*z,-2*z^2],2*z^2*[2*n,-(2*n+1)]]]
\end{verbatim}
$$\sqrt{\pi}e^{z^2}\erf(z)=\dfrac{2z}{1-\dfrac{2z^2}{3+\dfrac{4z^2}{5-\dfrac{6z^2}{7+\dfrac{8z^2}{9-\dfrac{10z^2}{11+\ddots}}}}}}$$
Convergence type $F^1$ with $E=2i/z^2$, $P=0$, and $C=\pi ze^{z^2}$, so that
$$\sqrt{\pi}e^{z^2}\erf(z)-\dfrac{p(n)}{q(n)}\sim(-1)^{\lfloor n/2\rfloor}\dfrac{\pi ze^{z^2}}{n!(2/z^2)^n}\;.$$
$$A=1+(z^4/4+z^2/2)/n+((z^4+2z^2-2)^2/32)/n^2+\cdots$$
\end{cf}

\smallskip

\begin{cf}\label{6.18.2}{\ }
\begin{verbatim}
[(z)->exp(-z^2)*intnum(t=0,z,exp(t^2)),
[0,2*n-1],[[z,2*z^2],2*z^2*[-2*n,2*n+1]]]
\end{verbatim}
$$e^{-z^2}\int_0^z e^{t^2}\,dt=\dfrac{z}{1+\dfrac{2z^2}{3-\dfrac{4z^2}{5+\dfrac{6z^2}{7-\dfrac{8z^2}{9+\dfrac{10z^2}{11-\ddots}}}}}}$$
Convergence type $F^1$ with $E=-2/z^2$, $P=0$, and $C=(\pi/2)ze^{-z^2}$, so that
$$e^{-z^2}\int_0^z e^{t^2}\,dt-\dfrac{p(n)}{q(n)}\sim(-1)^n\dfrac{(\pi/2)ze^{-z^2}}{n!(2/z^2)^n}\;.$$
$$A=1+(z^4/4-z^2/2)/n+((z^4-2z^2-2)^2/32)/n^2+\cdots$$
\end{cf}

\smallskip

\begin{cf}\label{6.18.3}{\ }
\begin{verbatim}
[(z)->sqrt(Pi)*exp(z^2)*erf(z),[0,2*n-1-2*z^2],[2*z,4*n*z^2]]
\end{verbatim}
$$\sqrt{\pi}e^{z^2}\erf(z)=\dfrac{2z}{-2z^2+1+\dfrac{4z^2}{-2z^2+3+\dfrac{8z^2}{-2z^2+5+\dfrac{12z^2}{-2z^2+7+\dfrac{16z^2}{-2z^2+9+\ddots}}}}}$$
Convergence type $F^1$ with $E=-1/z^2$, $P=0$, and $C=\pi ze^{2z^2}$, so that
$$\sqrt{\pi}e^{z^2}\erf(z)-\dfrac{p(n)}{q(n)}\sim(-1)^n\dfrac{\pi ze^{2z^2}}{n!(1/z^2)^n}\;.$$
$$A=1-(1/4)/n+(z^2+5/32)/n^2-((224z^4+160z^2+11)/128)/n^3+\cdots$$
\end{cf}

\smallskip

\begin{cf}\label{6.18.4}{\ }
\begin{verbatim}
[(z)->exp(-z^2)*intnum(t=0,z,exp(t^2)),[0,2*n-1+2*z^2],[z,-4*n*z^2]]
\end{verbatim}
$$e^{-z^2}\int_0^z e^{t^2}\,dt=\dfrac{z}{2z^2+1-\dfrac{4z^2}{2z^2+3-\dfrac{8z^2}{2z^2+5-\dfrac{12z^2}{2z^2+7-\dfrac{16z^2}{2z^2+9-\ddots}}}}}$$
Convergence type $F^1$ with $E=1/z^2$, $P=0$, and $C=(\pi/2)ze^{-2z^2}$, so that
$$e^{-z^2}\int_0^z e^{t^2}\,dt-\dfrac{p(n)}{q(n)}\sim\dfrac{(\pi/2)ze^{-2z^2}}{n!(1/z^2)^n}\;.$$
$$A=1-(1/4)/n-(z^2-5/32)/n^2-((224z^4-160z^2+11)/128)/n^3+\cdots$$
\end{cf}

\smallskip

\begin{cf}\label{6.18.5}{\ }
\begin{verbatim}
[(z)->sqrt(Pi)*exp(z^2)*erf(z),[0,1,2*n+2*z^2-1],
                               [2*z,-2*z^2*(2*n-1)]]
\end{verbatim}
$$\sqrt{\pi}e^{z^2}\erf(z)=\dfrac{2z}{1-\dfrac{2z^2}{2z^2+3-\dfrac{6z^2}{2z^2+5-\dfrac{10z^2}{2z^2+7-\dfrac{14z^2}{2z^2+9-\dfrac{18z^2}{2z^2+11-\ddots}}}}}}$$
Convergence type $F^1$ with $E=1/z^2$, $P=1/2$, and $C=z\sqrt{\pi}$,
so that
$$\sqrt{\pi}e^{z^2}\erf(z)-\dfrac{p(n)}{q(n)}\sim\dfrac{\sqrt{\pi}}{n!(1/z)^{2n+1}n^{1/2}}\;.$$
$$A=1+(z^2-3/8)/n+\cdots$$
Series:
$$\sqrt{\pi}e^{z^2}\erf(z)=\sum_{n\ge0}\dfrac{z^{2n+1}}{(1/2)_{n+1}}\;.$$
\end{cf}

\begin{cf}\label{6.18.6}{\ }
\begin{verbatim}
[(z)->exp(-z^2)*intnum(t=0,z,exp(t^2)),
                [0,1,2*n-2*z^2-1],[z,2*z^2*(2*n-1)]]
\end{verbatim}
$$e^{-z^2}\int_0^ze^{t^2}\,dt=\dfrac{z}{1+\dfrac{2z^2}{-2z^2+3+\dfrac{6z^2}{-2z^2+5+\dfrac{10z^2}{-2z^2+7+\dfrac{14z^2}{-2z^2+9+\dfrac{18z^2}{-2z^2+11+\ddots}}}}}}$$
Convergence type $F^1$ with $E=-1/z^2$, $P=1/2$, and $C=z\sqrt{\pi}/2$,
so that
$$e^{-z^2}\int_0^ze^{t^2}\,dt-\dfrac{p(n)}{q(n)}\sim(-1)^n\dfrac{\sqrt{\pi}/2}{n!(1/z)^{2n+1}n^{1/2}}\;.$$
$$A=1+(-z^2-3/8)/n+\cdots$$
Series:
$$e^{-z^2}\int_0^ze^{t^2}\,dt=\dfrac{1}{2}\sum_{n\ge0}(-1)^n\dfrac{z^{2n+1}}{(1/2)_{n+1}}\;.$$
\end{cf}

\smallskip

\begin{cf}\label{6.18.7}{\ }
\begin{verbatim}
[(z)->sqrt(Pi)*erf(z),[2*z,3,2*n^2+(1-2*z^2)*n+z^2],
                      [-2*z^3,n*(2*n+1)^2*z^2]]
\end{verbatim}
$$\sqrt{\pi}\erf(z)=2z-\dfrac{2z^3}{3+\dfrac{9z^2}{-3z^2+10+\dfrac{50z^2}{-5z^2+21+\dfrac{147z^2}{-7z^2+36+\dfrac{324z^2}{-9z^2+55+\ddots}}}}}$$
Convergence type $F^1$ with $E=-1/z^2$, $P=2$, and $C=-z^3$, so that
$$\sqrt{\pi}\erf(z)-\dfrac{p(n)}{q(n)}\sim(-1)^{n+1}\dfrac{1}{n!(1/z)^{2n+3}n^2}\;.$$
$$A=1-(2z^2+5)/n+(4z^4+22z^2+19)/n^2+\cdots$$
Series:
$$\sqrt{\pi}\erf(z)=2\sum_{n\ge0}(-1)^n\dfrac{z^{2n+1}}{(2n+1)n!}\;.$$
\end{cf}

\smallskip

\begin{cf}\label{6.18.8}{\ }
\begin{verbatim}
[(z)->intnum(t=0,z,exp(t^2)),[2*z,3,2*n^2+(1+2*z^2)*n-z^2],
                      [z^3,-n*(2*n+1)^2*z^2]]
\end{verbatim}
$$\int_0^ze^{t^2}\,dt=2z+\dfrac{2z^3}{3-\dfrac{9z^2}{3z^2+10-\dfrac{50z^2}{5z^2+21-\dfrac{147z^2}{7z^2+36-\dfrac{324z^2}{9z^2+55-\ddots}}}}}$$
Convergence type $F^1$ with $E=1/z^2$, $P=2$, and $C=z^3$, so that
$$\int_0^ze^{t^2}\,dt-\dfrac{p(n)}{q(n)}\sim\dfrac{1}{n!(1/z)^{2n+3}n^2}\;.$$
$$A=1+(2z^2-5)/n+(4z^4-22z^2+19)/n^2+\cdots$$
Series:
$$\int_0^ze^{t^2}\,dt=\sum_{n\ge0}\dfrac{z^{2n+1}}{(2n+1)n!}\;.$$
\end{cf}

\smallskip

\begin{cf}\label{6.18.9}{\ }
\begin{verbatim}
[(z)->sqrt(Pi)*exp(z^2)*erfc(z),[0,2*z],[2,2*n]]
\end{verbatim}
$$\sqrt{\pi}e^{z^2}\erfc(z)=\dfrac{2}{2z+\dfrac{2}{2z+\dfrac{4}{2z+\dfrac{6}{2z+\dfrac{8}{2z+\dfrac{10}{2z+\ddots}}}}}}$$
Convergence type $D^-$ with $D=8z^2$ and $C=2\sqrt{\pi}e^{z^2}$, so that
$$\sqrt{\pi}e^{z^2}\erfc(z)-\dfrac{p(n)}{q(n)}\sim(-1)^n\dfrac{2\sqrt{\pi}e^{z^2}}{e^{2z\sqrt{2n}}}\;.$$
$$A=1-(dz(z^2+3)/6)/n^{1/2}+(z^2(z^2+3)^2/36)/n+\cdots$$
with $d=\sqrt{2}$.
\end{cf}

\smallskip

Another way of writing this CF is as follows:

\smallskip

\begin{cf}\label{6.18.10}{\ }
\begin{verbatim}
[(a,b)->sqrt(2*b/Pi)*exp(-a^2/(2*b))/erfc(a/sqrt(2*b)),[a],[b*(n+1)]]
\end{verbatim}
$$\dfrac{\sqrt{2b/\pi}\,e^{-a^2/(2b)}}{\erfc(a/\sqrt{2b})}=a+\dfrac{b}{a+\dfrac{2b}{a+\dfrac{3b}{a+\dfrac{4b}{a+\dfrac{5b}{a+\dfrac{6b}{a+\ddots}}}}}}$$
Convergence type $D^-$ with $D=4a^2/b$ and
$C=\sqrt{8b/\pi}e^{-a^2/(2b)}/\erfc(a/\sqrt{2b})^2$, so that
$$\dfrac{\sqrt{2b/\pi}\,e^{-a^2/(2b)}}{\erfc(a/\sqrt{2b})}-\dfrac{p(n)}{q(n)}\sim(-1)^n\dfrac{\sqrt{8b/\pi}e^{-a^2/(2b)}/\erfc(a/\sqrt{2b})^2}{e^{2a\sqrt{n/b}}}\;.$$
$$A=1-(a(a^2+18b)/(12b^{3/2}))/n^{1/2}+(a^2(a^2+18b)^2/(288b^3))/n+\cdots$$
\end{cf}

\smallskip

\begin{cf}\label{6.18.11}{\ }
\begin{verbatim}
[(z)->sqrt(Pi)*exp(z^2)*erfc(z),[0,4*n+2*z^2-3],[2*z,-2*n*(2*n-1)]]
\end{verbatim}
$$\sqrt{\pi}e^{z^2}\erfc(z)=\dfrac{2z}{2z^2+1-\dfrac{2}{2z^2+5-\dfrac{12}{2z^2+9-\dfrac{30}{2z^2+13-\dfrac{56}{2z^2+17-\ddots}}}}}$$
Convergence type $D^+$ with $D=16z^2$ and $C=2\sqrt{\pi}e^{z^2}$, so that
$$\sqrt{\pi}e^{z^2}\erfc(z)-\dfrac{p(n)}{q(n)}\sim\dfrac{2\sqrt{\pi}e^{z^2}}{e^{4z\sqrt{n}}}\;.$$
$$A=1-(z(z^2+3)/6)/n^{1/2}+(z^2(z^2+3)^2/72)/n+\cdots$$
\end{cf}

This CF is simply the contraction of \ref{6.18.9}, but I have found
useful to give both.

\smallskip

\begin{cf}\label{6.18.12}{\ }
\begin{verbatim}
[(z)->1/sqrt(Pi)*intnum(t=[-oo,1],[oo,1],exp(-t^2)/(z-t)),[0,2*z],[2,-2*n]]
\end{verbatim}
$$\dfrac{1}{\sqrt{\pi}}\int_{-\infty}^{\infty}\dfrac{e^{-t^2}}{z-t}\,dt=\dfrac{2}{2z-\dfrac{2}{2z-\dfrac{4}{2z-\dfrac{6}{2z-\dfrac{8}{2z-\dfrac{10}{2z-\ddots}}}}}}$$
For $\Im(z)>0$ convergence type $D^-$ with $D=-8z^2$, and
$C=2\sqrt{\pi}e^{-z^2}$, so that
$$\dfrac{1}{\sqrt{\pi}}\int_{-\infty}^{\infty}\dfrac{e^{-t^2}}{z-t}\,dt-\dfrac{p(n)}{q(n)}\sim(-1)^n\dfrac{2\sqrt{\pi}e^{-z^2}}{e^{\sqrt{-8z^2n}}}\;.$$
$$A=1-(dz(z^2-3)/6)/n^{1/2}+(z^2(z^2-3)^2/36)/n+\cdots$$
with $d=\sqrt{-2}$.
\end{cf}

Note in particular that this CF (and of course the integral) does not
converge when $z$ is real.

\smallskip

\begin{cf}\label{6.18.13}{\ }
\begin{verbatim}
[()->erf(1/sqrt(2))*sqrt(exp(1)*Pi/2),[1,2*(n+1)],[2,2*(n+2)]]
[()->erf(1/sqrt(2))*sqrt(exp(1)*Pi/2),[1],[1/(2*(n+1))]]
\end{verbatim}
$$\erf(1/\sqrt{2})\sqrt{\pi e/2}=1+\dfrac{2}{4+\dfrac{6}{6+\dfrac{8}{8+\dfrac{10}{10+\dfrac{12}{12+\ddots}}}}}=1+\dfrac{1/2}{1+\dfrac{1/4}{1+\dfrac{1/6}{1+\dfrac{1/8}{1+\dfrac{1/10}{1+\ddots}}}}}$$
Convergence type $F^1$ with $E=-2$, $P=0$, and $C=\pi e/8$, so that
$$\erf(1/\sqrt{2})\sqrt{\pi e/2}-\dfrac{p(n)}{q(n)}\sim(-1)^n\dfrac{\pi e}{n!2^{n+3}}\;.$$
$$A=1-(13/4)/n+(285/32)/n^2-(3199/128)/n^3+\cdots$$
\end{cf}

\smallskip

\begin{cf}\label{6.18.14}{\ }
\begin{verbatim}
[()->erf(1/sqrt(2))*sqrt(exp(1)*Pi/2),[0,1,2*n],[-(2*n-1)]]
\end{verbatim}
$$\erf(1/\sqrt{2})\sqrt{\pi e/2}=\dfrac{1}{1-\dfrac{1}{4-\dfrac{3}{6-\dfrac{5}{8-\dfrac{7}{10-\dfrac{9}{12-\ddots}}}}}}$$
Convergence type $F^1$ with $E=2$, $P=1/2$, and $C=\sqrt{\pi}/2$, so that
$$\erf(1/\sqrt{2})\sqrt{2\pi e}-\dfrac{p(n)}{q(n)}\sim\dfrac{\sqrt{\pi}}{n!2^{n+1}n^{1/2}}\;.$$
$$A=1+(1/8)/n-(63/128)/n^2+(443/1024)/n^3+\cdots$$
Series:
$$\erf(1/\sqrt{2})\sqrt{\pi e/2}=\sum_{n\ge1}\dfrac{2^{-n}}{(1/2)_n}\;.$$
\end{cf}
        
\smallskip

\begin{cf}\label{6.18.15}{\ }
\begin{verbatim}
[()->erfc(1/sqrt(2))*sqrt(exp(1)*Pi/2),[0,1],[1,n]]
\end{verbatim}
$$\erfc(1/\sqrt{2})\sqrt{\pi e/2}=\dfrac{1}{1+\dfrac{1}{1+\dfrac{2}{1+\dfrac{3}{1+\dfrac{4}{1+\dfrac{5}{1+\ddots}}}}}}$$
Convergence type $D^-$ with $D=4$ and $C=\sqrt{2\pi e}$, so that
$$\erfc(1/\sqrt{2})\sqrt{\pi e/2}-\dfrac{p(n)}{q(n)}\sim(-1)^n\dfrac{\sqrt{2\pi e}}{e^{2\sqrt{n}}}\;.$$
$$A=1-(7/12)/n^{1/2}+(49/288)/n-(3713/51840)/n^{3/2}+\cdots$$
\end{cf}

\smallskip

The contraction of this CF gives:

\smallskip

\begin{cf}\label{6.18.16}{\ }
\begin{verbatim}
[()->erfc(1/sqrt(2))*sqrt(exp(1)*Pi/2),[0,2*(2*n-1)],[1,-2*n*(2*n-1)]]
[()->erfc(1/sqrt(2))*sqrt(exp(1)*Pi/2),[0,2],[1,-2*n/(2*n+1)]]
\end{verbatim}
$$\erfc(1/\sqrt{2})\sqrt{\pi e/2}=\dfrac{1}{2-\dfrac{2}{6-\dfrac{12}{10-\dfrac{30}{14-\dfrac{56}{18-\ddots}}}}}=\dfrac{1}{2-\dfrac{2/3}{2-\dfrac{4/5}{2-\dfrac{6/7}{2-\dfrac{8/9}{2-\ddots}}}}}$$
Convergence type $D^+$ with $D=8$ and $C=\sqrt{2\pi e}$, so that
$$\erfc(1/\sqrt{2})\sqrt{\pi e/2}-\dfrac{p(n)}{q(n)}\sim\dfrac{\sqrt{2\pi e}}{e^{\sqrt{8n}}}\;.$$
$$A=1-(7d/24)/n^{1/2}+(49/576)/n-(3713d/207360)/n^{3/2}+\cdots$$
with $d=\sqrt{2}$.

Parametric family:
\begin{verbatim}
[()->erfc(1/sqrt(2))*sqrt(exp(1)*Pi/2),4*n+2*k,-(2*n+u)*(2*n+2*k+1-u)]
\end{verbatim}
Convergence type $D^+$ with $D=8$.
\end{cf}

\medskip

\section{Jacobi Elliptic Functions: Laplace Transforms}

\medskip

In what follows, we assume $k$ real with $|k|\le1$. Since $\sn(z,k)$ vanishes
at $z=0$, there are no formulas involving $\sn$ in the denominator. Since
$\cn(z,k)$ vanishes for $z=K(k)$ and all odd multiples, formulas involving
$\cn$ in the denominator are purely formal.

\smallskip

Note, however, that the formal (and even non-convergent) CF are still useful,
since they give explicitly the Mac-Laurin expansion at $0$ of the
corresponding elliptic function: if {\tt A} is the CF in the variable $z$
of the Laplace transform of some elliptic function, the command
\begin{verbatim}
serconvol(exp(z+O(z^N)),cftoser(cfsubst(A,z,1/z),N,z)/z)
\end{verbatim}
gives the Mac-Laurin expansion at $0$. See \ref{6.19.6} below for an explicit
example.

\smallskip

All of the CFs of this section are of the following two types:
either a period 2 CF with \emph{constant} $a(n)$ (in fact always $a(n)=z$),
and $b(2n)$, $b(2n+1)$ even polynomials of degree at most 2 in $k$ possibly
multiplied by $z^2$, (always followed by its even contraction) or a period 1
CF with $a(n)=z^2+u(n)k^2+v(n)$ for suitable $u(n)$ and $v(n)$, and $b(n)$
some even polynomial in $k$ of degree at most $4$.

\smallskip

In the following {\tt GP} scripts for CFs, it is assumed that all Jacobi
elliptic functions {\tt ellxx} are implemented.

\smallskip

\begin{cf}\label{6.19.1}{\ }
\begin{verbatim}
[(z,k)->intnum(t=0,[oo,z],ellsn(t,k)*exp(-z*t)),
[0,z^2+(1+k^2)*(2*n-1)^2],[1,-(2*n-1)*(2*n)^2*(2*n+1)*k^2]]
\end{verbatim}
$$\int_0^\infty\sn(t,k)e^{-zt}\,dt=\dfrac{1}{z^2+k^2+1-\dfrac{12k^2}{z^2+9k^2+9-\dfrac{240k^2}{z^2+25k^2+25-\ddots}}}$$
Convergence type $E$ with $E=1/k^2$, $P=0$, and $C=...$, so that
$$\int_0^\infty\sn(t,k)e^{-zt}\,dt-\dfrac{p(n)}{q(n)}\sim\dfrac{C}{(1/k)^{2n}}\;.$$
$$A=1-((2z^2-k^2-1)/(4(1-k^2)))/n+\cdots$$
\end{cf}

\smallskip

\begin{cf}\label{6.19.2}{\ }
\begin{verbatim}
[(z,k)->intnum(t=0,[oo,z],ellcn(t,k)*exp(-z*t)),
[0,z],[[1,1],[(2*n)^2*k^2,(2*n+1)^2]]]
\end{verbatim}
$$\int_0^\infty\cn(t,k)e^{-zt}\,dt=\dfrac{1}{z+\dfrac{1}{z+\dfrac{4k^2}{z+\dfrac{9}{z+\dfrac{16k^2}{z+\dfrac{25}{z+\ddots}}}}}}$$
Convergence type $E$ with $E=-1/k$, $P=1$, and $C=...$, so that
$$\int_0^\infty\cn(t,k)e^{-zt}\,dt-\dfrac{p(n)}{q(n)}\sim(-1)^n\dfrac{C}{(1/k)^nn}\;.$$
$$A=1-((2z^2-k^2-3)/(2(1-k^2)))/n+\cdots$$
\end{cf}

\smallskip

\begin{cf}\label{6.19.3}{\ }
\begin{verbatim}
[(z,k)->intnum(t=0,[oo,z],ellcn(t,k)*exp(-z*t)),
[0,z^2+(2*n-1)^2+(n-1)^2*4*k^2],[z,-(2*n-1)^2*(2*n)^2*k^2]]
\end{verbatim}
$$\int_0^\infty\cn(t,k)e^{-zt}\,dt=\dfrac{z}{z^2+1-\dfrac{4k^2}{z^2+4k^2+9-\dfrac{144k^2}{z^2+16k^2+25-\ddots}}}$$
Convergence type $E$ with $E=1/k^2$, $P=1$, and $C=...$, so that
$$\int_0^\infty\cn(t,k)e^{-zt}\,dt-\dfrac{p(n)}{q(n)}\sim\dfrac{C}{(1/k)^{2n}n}\;.$$
$$A=1-((2z^2-k^2-3)/(4(1-k^2)))/n+\cdots$$
\end{cf}

\smallskip

\begin{cf}\label{6.19.4}{\ }
\begin{verbatim}
[(z,k)->intnum(t=0,[oo,z],elldn(t,k)*exp(-z*t)),
[0,z],[[1,k^2],[(2*n)^2,(2*n+1)^2*k^2]]]
\end{verbatim}
$$\int_0^\infty\dn(t,k)e^{-zt}\,dt=\dfrac{1}{z+\dfrac{k^2}{z+\dfrac{4}{z+\dfrac{9k^2}{z+\dfrac{16}{z+\dfrac{25k^2}{z+\ddots}}}}}}$$
Convergence type $E$ with $E=-1/k$, $P=-1$, and $C=...$, so that
$$\int_0^\infty\dn(t,k)e^{-zt}\,dt-\dfrac{p(n)}{q(n)}\sim(-1)^n\dfrac{C}{(1/k)^nn^{-1}}\;.$$
$$A=1-((2z^2-3k^2-1)/(2(1-k^2)))/n+\cdots$$
\end{cf}

\smallskip

\begin{cf}\label{6.19.5}{\ }
\begin{verbatim}
[(z,k)->intnum(t=0,[oo,z],elldn(t,k)*exp(-z*t)),
[0,z^2+(2*n-1)^2*k^2+4*(n-1)^2],[z,-(2*n-1)^2*(2*n)^2*k^2]]
\end{verbatim}
$$\int_0^\infty\dn(t,k)e^{-zt}\,dt=\dfrac{z}{z^2+k^2-\dfrac{4k^2}{z^2+9k^2+4-\dfrac{144k^2}{z^2+25k^2+16-\ddots}}}$$
Convergence type $E$ with $E=1/k^2$, $P=-1$, and $C=...$, so that
$$\int_0^\infty\dn(t,k)e^{-zt}\,dt-\dfrac{p(n)}{q(n)}\sim\dfrac{C}{(1/k)^{2n}n^{-1}}\;.$$
$$A=1-((2z^2-3k^2-1)/(4(1-k^2)))/n+\cdots$$
\end{cf}

\smallskip

Set $k_1=2\sqrt{k}/(1+k)$. Since
\begin{align*}
  \dfrac{\sn(t,k)}{1+k\sn^2(t,k)}&=\dfrac{\sn((1+k)t,k_1)}{1+k}\\
  \dfrac{\cn(t,k)\dn(t,k)}{1+k\sn^2(t,k)}&=\cn((1+k)t,k_1)\\
  \dfrac{1-k\sn(t,k)^2}{1+k\sn^2(t,k)}&=\dn((1+k)t,k_1)\;,\end{align*}
it is immediate to obtain from the above CFs for the left-hand sides.
  
\smallskip

\begin{cf}\label{6.19.6}{\ }
\begin{verbatim}
[(z,k)->intnum(t=0,[oo,z],ellsc(t,k)*exp(-z*t)),
[0,z^2-(2*n-1)^2*(2-k^2)],[1,-(2*n-1)*(2*n)^2*(2*n+1)*(1-k^2)]]
\end{verbatim}
$$\int_0^\infty\ssc(t,k)e^{-zt}\,dt=\dfrac{1}{z^2+k^2-2+\dfrac{12k^2-12}{z^2+9k^2-18+\dfrac{240k^2-240}{z^2+25k^2-50+\ddots}}}$$
Convergence type $E$ with $E=1/(1-k^2)$, $P=0$, and $C=...$, so that
$$\int_0^\infty\ssc(t,k)e^{-zt}\,dt-\dfrac{p(n)}{q(n)}\sim\dfrac{C}{(1/(1-k^2))^n}\;.$$
$$A=1+((2z^2-k^2+2)/(4k^2))/n+\cdots$$
\end{cf}
Note: since $\ssc(t,k)$ has poles on the real $t$-axis, this CF is
purely formal.

\smallskip

As mentioned above, although formal this allows us to compute the Mac-Laurin
expansion of $\ssc$. Calling $C$ the above CF, the following {\tt Pari/GP}
command computes this expansion to a few terms:

\begin{verbatim}
? serconvol(exp(z),cftoser(cfsubst(C,z,1/z),2,z)/z)
% = z+(-1/6*k^2+1/3)*z^3+(1/120*k^4-2/15*k^2+2/15)*z^5+...
\end{verbatim}
The same expansion could of course be obtained directly with the command
{\tt ellsc(z,k)} assuming this function implemented for polynomial arguments.

\smallskip

\begin{cf}\label{6.19.7}{\ }
\begin{verbatim}
[(z,k)->intnum(t=0,[oo,z],ellsd(t,k)*exp(-z*t)),
[0,z^2-(2*n-1)^2*(2*k^2-1)],[1,(2*n-1)*(2*n)^2*(2*n+1)*k^2*(1-k^2)]]
\end{verbatim}
$$\int_0^\infty\sd(t,k)e^{-zt}\,dt=\dfrac{1}{z^2-2k^2+1-\dfrac{12k^4-12k^2}{z^2-18k^2+9-\dfrac{240k^4-240k^2}{z^2-50k^2+25-\ddots}}}$$
Convergence type $E$ with $E=-(1-k^2)/k^2$, $P=0$, and $C=...$, so that
$$\int_0^\infty\sd(t,k)e^{-zt}\,dt-\dfrac{p(n)}{q(n)}\sim(-1)^n\dfrac{C}{((1-k^2)/k^2)^n}\;.$$
$$A=1-((2z^2+2k^2-1)/4)/n+\cdots$$
Note: if $k^2>1/2$, replace $E$ by $1/E$.
\end{cf}

\smallskip

\begin{cf}\label{6.19.8}{\ }
\begin{verbatim}
[(z,k)->intnum(t=0,[oo,z],ellnc(t,k)*exp(-z*t)),
[0,z],[[1,-1],[-(2*n)^2*(1-k^2),-(2*n+1)^2]]]
\end{verbatim}
$$\int_0^\infty\nc(t,k)e^{-zt}\,dt=\dfrac{1}{z-\dfrac{1}{z+\dfrac{4k^2-4}{z-\dfrac{9}{z+\dfrac{16k^2-16}{z-\dfrac{25}{z+\ddots}}}}}}$$
Convergence type $E$ with $E=-1/\sqrt{1-k^2}$, $P=1$, and $C=...$, so that
$$\int_0^\infty\nc(t,k)e^{-zt}\,dt-\dfrac{p(n)}{q(n)}\sim(-1)^n\dfrac{C}{(1/(1-k^2))^{n/2}n}\;.$$
$$A=1+((2z^2-k^2+4)/(2k^2))/n+\cdots$$
\end{cf}

\smallskip

\begin{cf}\label{6.19.9}{\ }
\begin{verbatim}
[(z,k)->intnum(t=0,[oo,z],ellnc(t,k)*exp(-z*t)),
[0,z^2+(n-1)^2*4*k^2-8*n^2+12*n-5],[z,-(2*n-1)^2*(2*n)^2*(1-k^2)]]
\end{verbatim}
$$\int_0^\infty\nc(t,k)e^{-zt}\,dt=\dfrac{z}{z^2-1+\dfrac{4k^2-4}{z^2+4k^2-13+\dfrac{144k^2-144}{z^2+16k^2-41+\ddots}}}$$
Convergence type $E$ with $E=1/(1-k^2)$, $P=1$, and $C=...$, so that
$$\int_0^\infty\nc(t,k)e^{-zt}\,dt-\dfrac{p(n)}{q(n)}\sim\dfrac{C}{(1/(1-k^2))^nn}\;.$$
$$A=1+((2z^2-k^2+4)/(4k^2))/n+\cdots$$
\end{cf}

Since $\nc(t,k)$ has poles on the real $t$-axis, these two CFs are purely
formal.

\smallskip

\begin{cf}\label{6.19.10}{\ }
\begin{verbatim}
[(z,k)->intnum(t=0,[oo,z],ellcd(t,k)*exp(-z*t)),
[0,z],[[1,1-k^2],[-(2*n)^2*k^2,(2*n+1)^2*(1-k^2)]]]
\end{verbatim}
$$\int_0^\infty\cd(t,k)e^{-zt}\,dt=\dfrac{1}{z-\dfrac{k^2-1}{z-\dfrac{4k^2}{z-\dfrac{9k^2-9}{z-\dfrac{16k^2}{z-\dfrac{25k^2-25}{z-\ddots}}}}}}$$
Convergence type $E$ with $E=-\sqrt{(1-k^2)/k^2}$, $P=1$, and $C=...$, so that
$$\int_0^\infty\cd(t,k)e^{-zt}\,dt-\dfrac{p(n)}{q(n)}\sim(-1)^n\dfrac{C}{((1-k^2)/k^2)^{n/2}n}\;.$$
$$A=1-((2z^2+4k^2-3)/2)/n+\cdots$$
Note: if $k^2>1/2$, replace $E$ by $1/E$ and $P$ by $-P$.
\end{cf}

\smallskip

\begin{cf}\label{6.19.11}{\ }
\begin{verbatim}
[(z,k)->intnum(t=0,[oo,z],ellcd(t,k)*exp(-z*t)),
[0,z^2+(2*n-1)^2-k^2*(8*n^2-12*n+5)],[z,(2*n-1)^2*(2*n)^2*k^2*(1-k^2)]]
\end{verbatim}
$$\int_0^\infty\cd(t,k)e^{-zt}\,dt=\dfrac{z}{z^2-k^2+1-\dfrac{4k^4-4k^2}{z^2-13k^2+9-\dfrac{144k^4-144k^2}{z^2-41k^2+25-\ddots}}}$$
Convergence type $E$ with $E=-(1-k^2)/k^2$, $P=1$, and $C=...$, so that
$$\int_0^\infty\cd(t,k)e^{-zt}\,dt-\dfrac{p(n)}{q(n)}\sim(-1)^n\dfrac{C}{((1-k^2)/k^2)^nn}\;.$$
$$A=1-((2z^2+4k^2-3)/4)/n+\cdots$$
Note: if $k^2>1/2$, replace $E$ by $1/E$ and $P$ by $-P$.
\end{cf}

\smallskip

\begin{cf}\label{6.19.12}{\ }
\begin{verbatim}
[(z,k)->intnum(t=0,[oo,z],elldc(t,k)*exp(-z*t)),
[0,z],[[1,-(1-k^2)],[-(2*n)^2,-(2*n+1)^2*(1-k^2)]]]
\end{verbatim}
$$\int_0^\infty\dc(t,k)e^{-zt}\,dt=\dfrac{1}{z+\dfrac{k^2-1}{z-\dfrac{4}{z+\dfrac{9k^2-9}{z-\dfrac{16}{z+\dfrac{25k^2-25}{z-\ddots}}}}}}$$

Convergence type $E$ with $E=-1/\sqrt{1-k^2}$, $P=-1$, and $C=...$, so that
$$\int_0^\infty\dc(t,k)e^{-zt}\,dt-\dfrac{p(n)}{q(n)}\sim(-1)^n\dfrac{C}{(1/(1-k^2))^{n/2}n^{-1}}\;.$$
$$A=1+((2z^2-3k^2+4)/(2k^2))/n+\cdots$$
\end{cf}

\smallskip

\begin{cf}\label{6.19.13}{\ }
\begin{verbatim}
[(z,k)->intnum(t=0,[oo,z],elldc(t,k)*exp(-z*t)),
[0,z^2+(2*n-1)^2*k^2-(8*n^2-12*n+5)],[z,-(2*n-1)^2*(2*n)^2*(1-k^2)]]
\end{verbatim}
$$\int_0^\infty\dc(t,k)e^{-zt}\,dt=\dfrac{z}{z^2+k^2-1+\dfrac{4k^2-4}{z^2+9k^2-13+\dfrac{144k^2-144}{z^2+25k^2-41+\ddots}}}$$
Convergence type $E$ with $E=1/(1-k^2)$, $P=-1$, and $C=...$, so that
$$\int_0^\infty\dc(t,k)e^{-zt}\,dt-\dfrac{p(n)}{q(n)}\sim\dfrac{C}{(1/(1-k^2))^nn^{-1}}\;.$$
$$A=1+((2z^2-3k^2+4)/(4k^2))/n+\cdots$$
\end{cf}

Since $\dc(t,k)$ has poles on the real $t$-axis, these two CFs are purely
formal.

\smallskip

\begin{cf}\label{6.19.14}{\ }
\begin{verbatim}
[(z,k)->intnum(t=0,[oo,z],ellnd(t,k)*exp(-z*t)),
[0,z],[[1,-k^2],[(2*n)^2*(1-k^2),-(2*n+1)^2*k^2]]]
\end{verbatim}
$$\int_0^\infty\nd(t,k)e^{-zt}\,dt=\dfrac{1}{z-\dfrac{k^2}{z-\dfrac{4k^2-4}{z-\dfrac{9k^2}{z-\dfrac{16k^2-16}{z-\dfrac{25k^2}{z-\ddots}}}}}}$$
Convergence type $E$ with $E=-\sqrt{1-k^2}/k$, $P=-1$, and $C=...$, so that
$$\int_0^\infty\nd(t,k)e^{-zt}\,dt-\dfrac{p(n)}{q(n)}\sim(-1)^n\dfrac{C}{((1-k^2)/k^2)^{n/2}n^{-1}}\;.$$
$$A=1-((2z^2+4k^2-1)/2)/n+\cdots$$
Note: if $k^2>1/2$, replace $E$ by $1/E$ and $P$ by $-P$.
\end{cf}

\smallskip

\begin{cf}\label{6.19.15}{\ }
\begin{verbatim}
[(z,k)->intnum(t=0,[oo,z],ellnd(t,k)*exp(-z*t)),
[0,z^2+4*(n-1)^2-k^2*(8*n^2-12*n+5)],[z,(2*n-1)^2*(2*n)^2*k^2*(1-k^2)]]
\end{verbatim}
$$\int_0^\infty\nd(t,k)e^{-zt}\,dt=\dfrac{z}{z^2-k^2-\dfrac{4k^4-4k^2}{z^2-13k^2+4-\dfrac{144k^4-144k^2}{z^2-41k^2+16-\ddots}}}$$
Convergence type $E$ with $E=-(1-k^2)/k^2$, $P=-1$, and $C=...$, so that
$$\int_0^\infty\nd(t,k)e^{-zt}\,dt-\dfrac{p(n)}{q(n)}\sim(-1)^n\dfrac{C}{((1-k^2)/k^2)^nn^{-1}}\;.$$
$$A=1-((2z^2+4k^2-1)/4)/n+\cdots$$
Note: if $k^2>1/2$, replace $E$ by $1/E$ and $P$ by $-P$.
\end{cf}

\smallskip

The above CFs give all Laplace transforms of the Jacobi elliptic functions
which do not have $\sn$ in the denominator, since otherwise the integrals
would trivially diverge.

\smallskip

\begin{cf}\label{6.19.16}{\ }
\begin{verbatim}
[(z,k)->intnum(t=0,[oo,z],ellsn(t,k)*ellcn(t,k)*exp(-z*t)),
[0,z^2+(2*n-1)^2*k^2+(2*n)^2],[1,-4*k^2*n^2*(2*n+1)^2]]
\end{verbatim}
$$\int_0^\infty\sn(t,k)\cn(t,k)e^{-zt}\,dt=\dfrac{1}{z^2+k^2+4-\dfrac{36k^2}{z^2+9k^2+16-\dfrac{400k^2}{z^2+25k^2+36-\ddots}}}$$
Convergence type $E$ with $E=1/k^2$, $P=1$, and $C=...$, so that
$$\int_0^\infty\sn(t,k)\cn(t,k)e^{-zt}\,dt-\dfrac{p(n)}{q(n)}\sim\dfrac{C}{(1/k)^{2n}n}\;.$$
$$A=1+((2z^2+k^2-5)/(4(1-k^2)))/n+\cdots$$
\end{cf}

\smallskip

\begin{cf}\label{6.19.17}{\ }
\begin{verbatim}
[(z,k)->intnum(t=0,[oo,z],ellsn(t,k)*elldn(t,k)*exp(-z*t)),
[0,z^2+(2*n)^2*k^2+(2*n-1)^2],[1,-4*k^2*n^2*(2*n+1)^2]]
\end{verbatim}
$$\int_0^\infty\sn(t,k)\dn(t,k)e^{-zt}\,dt=\dfrac{1}{z^2+4k^2+1-\dfrac{36k^2}{z^2+16k^2+9-\dfrac{400k^2}{z^2+36k^2+25-\ddots}}}$$
Convergence type $E$ with $E=1/k^2$, $P=-1$, and $C=...$, so that
$$\int_0^\infty\sn(t,k)\dn(t,k)e^{-zt}\,dt-\dfrac{p(n)}{q(n)}\sim\dfrac{C}{(1/k)^{2n}n^{-1}}\;.$$
$$A=1+((2z^2+1-5k^2)/(4(1-k^2)))/n+\cdots$$
\end{cf}

\smallskip

\begin{cf}\label{6.19.18}{\ }
\begin{verbatim}
[(z,k)->intnum(t=0,[oo,z],ellcn(t,k)*elldn(t,k)*exp(-z*t)),
[0,z^2+(2*n-1)^2*(k^2+1)],[z,-4*k^2*n^2*(2*n-1)*(2*n+1)]]
\end{verbatim}
$$\int_0^\infty\cn(t,k)\dn(t,k)e^{-zt}\,dt=\dfrac{z}{z^2+k^2+1-\dfrac{12k^2}{z^2+9k^2+9-\dfrac{240k^2}{z^2+25k^2+25-\ddots}}}$$
Convergence type $E$ with $E=1/k^2$, $P=0$, and $C=...$, so that
$$\int_0^\infty\cn(t,k)\dn(t,k)e^{-zt}\,dt-\dfrac{p(n)}{q(n)}\sim\dfrac{C}{(1/k)^{2n}}\;.$$
$$A=1+((2z^2-1-k^2)/(4(1-k^2)))/n+\cdots$$
\end{cf}

\smallskip

\begin{cf}\label{6.19.19}{\ }
\begin{verbatim}
[(z,k)->intnum(t=0,[oo,z],ellsn(t,k)^2*exp(-z*t)),
[0,z^2+(2*n)^2*(1+k^2)],[2/z,-2*n*(2*n+1)^2*(2*n+2)*k^2]]
\end{verbatim}
$$\int_0^\infty\sn^2(t,k)e^{-zt}\,dt=\dfrac{2/z}{z^2+4k^2+4-\dfrac{72k^2}{z^2+16k^2+16-\dfrac{600k^2}{z^2+36k^2+36-\ddots}}}$$
Convergence type $E$ with $E=1/k^2$, $P=0$, and $C=...$, so that
$$\int_0^\infty\sn^2(t,k)e^{-zt}\,dt-\dfrac{p(n)}{q(n)}\sim\dfrac{C}{(1/k)^{2n}}\;.$$
$$A=1-((2z^2-k^2-1)/(4(1-k^2)))/n+\cdots$$
\end{cf}

Since $\cn^2=1-\sn^2$ and $\dn^2=1-k^2\sn^2$, this gives CFs for the Laplace
transforms of $\cn^2$ and $\dn^2$,

\smallskip

\begin{cf}\label{6.19.20}{\ }
\begin{verbatim}
[(z,k)->intnum(t=0,[oo,z],ellsc(t,k)^2*exp(-z*t)),
[0,z^2-4*n^2*(2-k^2)],[2/z,-2*n*(2*n+1)^2*(2*n+2)*(1-k^2)]]
\end{verbatim}
$$\int_0^\infty\ssc^2(t,k)e^{-zt}\,dt=\dfrac{2/z}{z^2+4k^2-8+\dfrac{72k^2-72}{z^2+16k^2-32+\dfrac{600k^2-600}{z^2+36k^2-72+\ddots}}}$$
Convergence type $E$ with $E=1/(1-k^2)$, $P=0$, and $C=...$, so that
$$\int_0^\infty\ssc^2(t,k)e^{-zt}\,dt-\dfrac{p(n)}{q(n)}\sim\dfrac{C}{(1/(1-k^2))^n}\;.$$
$$A=1+((2z^2-k^2+2)/(4k^2))/n+\cdots$$
\end{cf}
Note: since $\ssc(t,k)$ has poles on the real $t$-axis, this CF is
purely formal.

Since $\dc^2=1+(1-k^2)\ssc^2$ and $\nc^2=\ssc^2+1$, this gives formal CFs
for $\dc^2$ and $\nc^2$.

\smallskip

\begin{cf}\label{6.19.21}{\ }
\begin{verbatim}
[(z,k)->intnum(t=0,[oo,z],ellsd(t,k)^2*exp(-z*t)),
[0,z^2+4*n^2*(1-2*k^2)],[2/z,2*n*(2*n+1)^2*(2*n+2)*k^2*(1-k^2)]]
\end{verbatim}
$$\int_0^\infty\sd^2(t,k)e^{-zt}\,dt=\dfrac{2/z}{z^2-8k^2+4-\dfrac{72k^4-72k^2}{z^2-32k^2+16-\ddots}}$$
Convergence type $E$ with $E=-(1-k^2)/k^2$, $P=0$, and $C=...$, so that
$$\int_0^\infty\sd^2(t,k)e^{-zt}\,dt-\dfrac{p(n)}{q(n)}\sim(-1)^n\dfrac{C}{((1-k^2)/k^2)^n}\;.$$
$$A=1-((2z^2+2k^2-1)/4)/n+\cdots$$
Note: if $k^2>1/2$, replace $E$ by $1/E$.
\end{cf}

Since $\cd^2=1-(1-k^2)\sd^2$ and $\nd^2=1+k^2\sd^2$, this gives CFs for
$\cd^2$ and $\nd^2$.

\smallskip

\begin{cf}\label{6.19.22}{\ }
\begin{verbatim}
[(z,k)->intnum(t=0,[oo,z],ellsn(t,k)/elldn(t,k)^2*exp(-z*t)),
[0,z^2-(8*n^2-4*n+1)*k^2+(2*n-1)^2],[1,4*n^2*(2*n+1)^2*k^2*(1-k^2)]]
\end{verbatim}
$$\int_0^\infty\dfrac{\sn(t,k)}{\dn^2(t,k)}e^{-zt}\,dt=\dfrac{1}{z^2-5k^2+1-\dfrac{36k^4-36k^2}{z^2-25k^2+9-\dfrac{400k^4-400k^2}{z^2-61k^2+25-\ddots}}}$$
Convergence type $E$ with $E=-(1-k^2)/k^2$, $P=1$, and $C=...$, so that
$$\int_0^\infty\dfrac{\sn(t,k)}{\dn^2(t,k)}e^{-zt}\,dt-\dfrac{p(n)}{q(n)}\sim(-1)^n\dfrac{C}{((1-k^2)/k^2)^nn}\;.$$
$$A=1-((2z^2+4k^2+1)/4)/n+\cdots$$
Note: if $k^2>1/2$, replace $E$ by $1/E$.
\end{cf}

Note that $\sn(t,k)/\dn(t,k)^2=\sd(t,k)\nd(t,k)$.

\smallskip

\begin{cf}\label{6.19.23}{\ }
\begin{verbatim}
[(z,k)->intnum(t=0,[oo,z],ellcn(t,k)/elldn(t,k)^2*exp(-z*t)),
[0,z^2-(2*n-1)^2*(2*k^2-1)],[z,4*n^2*(4*n^2-1)*k^2*(1-k^2)]]
\end{verbatim}
$$\int_0^\infty\dfrac{\cn(t,k)}{\dn^2(t,k)}e^{-zt}\,dt=\dfrac{z}{z^2-2k^2+1-\dfrac{12k^4-12k^2}{z^2-18k^2+9-\dfrac{240k^4-240k^2}{z^2-50k^2+25-\ddots}}}$$
Convergence type $E$ with $E=-(1-k^2)/k^2$, $P=0$, and $C=...$, so that
$$\int_0^\infty\dfrac{\cn(t,k)}{\dn^2(t,k)}e^{-zt}\,dt-\dfrac{p(n)}{q(n)}\sim(-1)^n\dfrac{C}{((1-k^2)/k^2)^n}\;.$$
$$A=1-((2z^2+2k^2-1)/4)/n+\cdots$$
Note: if $k^2>1/2$, replace $E$ by $1/E$.
\end{cf}

Note that $\cn(t,k)/\dn(t,k)^2=\cd(t,k)\nd(t,k)$.

\smallskip

\begin{cf}\label{6.19.24}{\ }
\begin{verbatim}
[(z,k)->intnum(t=0,[oo,z],ellsn(t,k)*ellcn(t,k)/elldn(t,k)^2*exp(-z*t)),
[0,z^2-(8*n^2-4*n+1)*k^2+(2*n)^2],[1,4*n^2*(2*n+1)^2*k^2*(1-k^2)]]
\end{verbatim}
$$\int_0^\infty\dfrac{\sn(t,k)\cn(t,k)}{\dn^2(t,k)}e^{-zt}\,dt=\dfrac{1}{z^2-5k^2+4-\dfrac{36k^4-36k^2}{z^2-25k^2+16-\dfrac{400k^4-400k^2}{z^2-61k^2+36-\ddots}}}$$
Convergence type $E$ with $E=-(1-k^2)/k^2$, $P=-1$, and $C=...$, so that
$$\int_0^\infty\dfrac{\sn(t,k)\cn(t,k)}{\dn^2(t,k)}e^{-zt}\,dt-\dfrac{p(n)}{q(n)}\sim(-1)^n\dfrac{C}{((1-k^2)/k^2)^nn^{-1}}\;.$$
$$A=1-((2z^2+4k^2-5)/4)/n+\cdots$$
Note: if $k^2>1/2$, replace $E$ by $1/E$.
\end{cf}

Note that $\sn(t,k)\cn(t,k)/\dn(t,k)^2=\sd(t,k)\cd(t,k)$.

\smallskip

\begin{cf}\label{6.19.25}{\ }
\begin{verbatim}
[(z,k)->intnum(t=0,[oo,z],ellsn(t,k)/ellcn(t,k)^2*exp(-z*t)),
[0,z^2+(2*n-1)^2*k^2-(8*n^2-4*n+1)],[1,-4*n^2*(2*n+1)^2*(1-k^2)]]
\end{verbatim}
$$\int_0^\infty\dfrac{\sn(t,k)}{\cn^2(t,k)}e^{-zt}\,dt=\dfrac{1}{z^2+k^2-5+\dfrac{36k^2-36}{z^2+9k^2-25+\dfrac{400k^2-400}{z^2+25k^2-61+\ddots}}}$$
Convergence type $E$ with $E=1/(1-k^2)$, $P=1$, and $C=...$, so that
$$\int_0^\infty\dfrac{\sn(t,k)}{\cn^2(t,k)}e^{-zt}\,dt-\dfrac{p(n)}{q(n)}\sim\dfrac{C}{(1/(1-k^2))^nn}\;.$$
$$A=1+((2z^2+k^2+4)/(4k^2))/n+\cdots$$
\end{cf}

Note that $\sn(t,k)/\cn(t,k)^2=\ssc(t,k)\nc(t,k)$.

\smallskip

\begin{cf}\label{6.19.26}{\ }
\begin{verbatim}
[(z,k)->intnum(t=0,[oo,z],elldn(t,k)/ellcn(t,k)^2*exp(-z*t)),
[0,z^2+(2*n-1)^2*k^2-2*(2*n-1)^2],[z,-4*n^2*(4*n^2-1)*(1-k^2)]]
\end{verbatim}
$$\int_0^\infty\dfrac{\dn(t,k)}{\cn^2(t,k)}e^{-zt}\,dt=\dfrac{z}{z^2+k^2-2+\dfrac{12k^2-12}{z^2+9k^2-18+\dfrac{240k^2-240}{z^2+25k^2-50+\ddots}}}$$
Convergence type $E$ with $E=1/(1-k^2)$, $P=0$, and $C=...$, so that
$$\int_0^\infty\dfrac{\dn(t,k)}{\cn^2(t,k)}e^{-zt}\,dt-\dfrac{p(n)}{q(n)}\sim\dfrac{C}{(1/(1-k^2))^n}\;.$$
$$A=1+((2z^2-k^2+2)/(4k^2))/n+\cdots$$
\end{cf}

Note that $\dn(t,k)/\cn(t,k)^2=\dc(t,k)\nc(t,k)$.

\smallskip

\begin{cf}\label{6.19.27}{\ }
\begin{verbatim}
[(z,k)->intnum(t=0,[oo,z],ellsn(t,k)*elldn(t,k)/ellcn(t,k)^2*exp(-z*t)),
[0,z^2+(2*n)^2*k^2-(8*n^2-4*n+1)],[1,-4*n^2*(2*n+1)^2*(1-k^2)]]
\end{verbatim}
$$\int_0^\infty\dfrac{\sn(t,k)\dn(t,k)}{\cn^2(t,k)}e^{-zt}\,dt=\dfrac{1}{z^2+4k^2-5+\dfrac{36k^2-36}{z^2+16k^2-25+\dfrac{400k^2-400}{z^2+36k^2-61+\ddots}}}$$
Convergence type $E$ with $E=1/(1-k^2)$, $P=-1$, and $C=...$, so that
$$\int_0^\infty\dfrac{\sn(t,k)\dn(t,k)}{\cn^2(t,k)}e^{-zt}\,dt-\dfrac{p(n)}{q(n)}\sim\dfrac{C}{(1/(1-k^2))^nn^{-1}}\;.$$
$$A=1+((2z^2-5k^2+4)/(4k^2))/n+\cdots$$
\end{cf}

Note that $\sn(t,k)\dn(t,k)/\cn(t,k)^2=\ssc(t,k)\dc(t,k)$.

\smallskip

\begin{cf}\label{6.19.28}{\ }
\begin{verbatim}
[(z,k)->intnum(t=0,[oo,z],ellsn(t,k)*ellcd(t,k)*exp(-z*t)),
[0,z^2+2*(2*n-1)^2*(2-k^2)],[1,-(2*n-1)*(2*n)^2*(2*n+1)*k^4]]
\end{verbatim}
$$\int_0^\infty\sn(t,k)\cd(t,k)e^{-zt}\,dt=\dfrac{1}{z^2-2k^2+4-\dfrac{12k^4}{z^2-18k^2+36-\ddots}}$$
Convergence type $E$ with $E=((1+\sqrt{1-k^2})/k)^4$, $P=0$ and $C=...$, so
that
$$\int_0^\infty\sn(t,k)\cd(t,k)e^{-zt}\,dt-\dfrac{p(n)}{q(n)}\sim\dfrac{C}{((1+\sqrt{1-k^2})/k)^{4n}}\;.$$
\end{cf}

\smallskip

\begin{cf}\label{6.19.29}{\ }
\begin{verbatim}
[(z,k)->intnum(t=0,[oo,z],ellsc(t,k)*ellnd(t,k)*exp(-z*t)),
[0,z^2-2*(2*n-1)^2*(1+k^2)],[1,-(2*n-1)*(2*n)^2*(2*n+1)*(1-k^2)^2]]
\end{verbatim}
$$\int_0^\infty\ssc(t,k)\nd(t,k)e^{-zt}\,dt=\dfrac{1}{z^2-2k^2-2-\dfrac{12k^4-24k^2+12}{z^2-18k^2-18-\ddots}}$$
Convergence type $E$ with $E=((1+k)/(1-k))^2$, $P=0$, and $C=...$, so that
$$\int_0^\infty\ssc(t,k)\nd(t,k)e^{-zt}\,dt-\dfrac{p(n)}{q(n)}\sim\dfrac{C}{((1+k)/(1-k))^{2n}}\;.$$
$$A=1-((z^2+k^2+1)/(8k))/n+\cdots$$
\end{cf}
Note: since the integrand has poles on the positive real $t$-axis, this CF is
purely formal.

\smallskip

\begin{cf}\label{6.19.30}{\ }
\begin{verbatim}
[(z,k)->intnum(t=0,[oo,z],ellsn(t,k)*elldc(t,k)*exp(-z*t)),
[0,z^2+2*(2*n-1)^2*(2*k^2-1)],[1,-(2*n-1)*(2*n)^2*(2*n+1)]]
\end{verbatim}
$$\int_0^\infty\sn(t,k)\dc(t,k)e^{-zt}\,dt=\dfrac{1}{z^2+4k^2-2-\dfrac{12}{z^2+36k^2-18-\dfrac{240}{z^2+100k^2-50-\ddots}}}$$
Non convergent CF when $|k|<1$. It does converge for $|k|>1$, but we exclude
this.
\end{cf}
Note: since the integrand has poles on the positive real $t$-axis, this CF is
purely formal.

\smallskip

The above CFs give all Laplace transforms of products of two Jacobi elliptic
functions which do not have $\sn$ in the denominator, and which are not
themselves one of the twelve Jacobi elliptic functions, with the exception of
the following six functions:
\begin{align*}&\sn(t,k)\ssc(t,k),\quad \sn(t,k)\sd(t,k),\quad \cn(t,k)\cd(t,k),\\
  &\dn(t,k)\dc(t,k),\quad\nc(t,k)\nd(t,k),\quad\ssc(t,k)\sd(t,k),\end{align*}
for which no reasonably simple CF seems to exist.

\smallskip

\begin{cf}\label{6.19.31}{\ }
\begin{verbatim}
[(z,k)->intnum(t=0,[oo,z],ellsn(t,k)*ellcn(t,k)*elldn(t,k)*exp(-z*t)),
[0,z^2+4*n^2*(k^2+1)],[1,-4*k^2*n*(n+1)*(2*n+1)^2]]
\end{verbatim}
$$\int_0^\infty\sn(t,k)\cn(t,k)\dn(t,k)e^{-zt}\,dt=\dfrac{1}{z^2+4k^2+4-\dfrac{72k^2}{z^2+16k^2+16-\dfrac{600k^2}{z^2+36k^2+36-\ddots}}}$$
Convergence type $E$ with $E=1/k^2$, $P=0$, and $C=...$, so that
$$\int_0^\infty\sn(t,k)\cn(t,k)\dn(t,k)e^{-zt}\,dt-\dfrac{p(n)}{q(n)}\sim\dfrac{C}{(1/k)^{2n}}\;.$$
$$A=1+((2z^2-1-k^2)/(4(1-k^2)))/n+\cdots$$
\end{cf}

\smallskip

\begin{cf}\label{6.19.32}{\ }
\begin{verbatim}
[(z,k)->intnum(t=0,[oo,z],ellsn(t,k)*elldn(t,k)/ellcn(t,k)^3*exp(-z*t)),
[0,z^2-4*n^2*(2-k^2)],[1,-4*n*(n+1)*(2*n+1)^2*(1-k^2)]]
\end{verbatim}
$$\int_0^\infty\dfrac{\sn(t,k)\dn(t,k)}{\cn^3(t,k)}e^{-zt}\,dt=\dfrac{1}{z^2+4k^2-8+\dfrac{72k^2-72}{z^2+16k^2-32+\dfrac{600k^2-600}{z^2+36k^2-72+\ddots}}}$$
Convergence type $E$ with $E=1/(1-k^2)$, $P=0$, and $C=...$, so that
$$\int_0^\infty\dfrac{\sn(t,k)\dn(t,k)}{\cn^3(t,k)}e^{-zt}\,dt-\dfrac{p(n)}{q(n)}\sim\dfrac{C}{(1/(1-k^2))^n}\;.$$
$$A=1+((2z^2-k^2+2)/(4k^2))/n+\cdots$$
\end{cf}

\smallskip

\begin{cf}\label{6.19.33}{\ }
\begin{verbatim}
[(z,k)->intnum(t=0,[oo,z],ellsn(t,k)*ellcn(t,k)/elldn(t,k)^3*exp(-z*t)),
[0,z^2+4*n^2*(1-2*k^2)],[1,4*n*(n+1)*(2*n+1)^2*k^2*(1-k^2)]]
\end{verbatim}
$$\int_0^\infty\dfrac{\sn(t,k)\cn(t,k)}{\dn^3(t,k)}e^{-zt}\,dt=\dfrac{1}{z^2-8k^2+4-\dfrac{72k^4-72k^2}{z^2-32k^2+16-\dfrac{600k^4-600k^2}{z^2-72k^2+36-\ddots}}}$$
Convergence type $E$ with $E=-(1-k^2)/k^2$, $P=0$, and $C=...$, so that
$$\int_0^\infty\dfrac{\sn(t,k)\cn(t,k)}{\dn^3(t,k)}e^{-zt}\,dt-\dfrac{p(n)}{q(n)}\sim(-1)^n\dfrac{C}{((1-k^2)/k^2)^n}\;.$$
$$A=1-((2z^2+2k^2-1)/4)/n+\cdots$$
Note: if $k^2>1/2$, replace $E$ by $1/E$.
\end{cf}

\smallskip

\begin{cf}\label{6.19.34}{\ }
\begin{verbatim}
[(z,k)->intnum(t=0,[oo,z],ellsn(t,k)^2*ellcd(t,k)^2*exp(-z*t)),
[0,z^2+8*n^2*(2-k^2)],[2/z,-2*n*(2*n+1)^2*(2*n+2)*k^4]]
\end{verbatim}
$$\int_0^\infty\sn^2(t,k)\cd^2(t,k)e^{-zt}\,dt=\dfrac{2/z}{z^2-8k^2+16-\dfrac{72k^4}{z^2-32k^2+64-\ddots}}$$
Convergence type $E$ with $E=((1-\sqrt{1-k^2})/k)^4$, $P=0$, and $C=...$, so
that
$$\int_0^\infty\sn^2(t,k)\cd^2(t,k)e^{-zt}\,dt-\dfrac{p(n)}{q(n)}\sim\dfrac{C}{((1-\sqrt{1-k^2})/k)^{4n}}\;.$$
\end{cf}

Note that $$\sn^2(t,k)\cd^2(t,k)=\dfrac{1}{k^2}\dfrac{1-\dn(2t,k)}{1+\dn(2t,k)}\;,$$
so the above CF gives one for the Laplace transforms of $1/(1+\dn(t,k))$ and
$1/(1+\nd(t,k))$.

\smallskip

\begin{cf}\label{6.19.35}{\ }
\begin{verbatim}
[(z,k)->intnum(t=0,[oo,z],ellsc(t,k)^2*ellnd(t,k)^2*exp(-z*t)),
[0,z^2-8*n^2*(1+k^2)],[2/z,-2*n*(2*n+1)^2*(2*n+2)*(1-k^2)^2]]
\end{verbatim}
$$\int_0^\infty\ssc^2(t,k)\nd^2(t,k)e^{-zt}\,dt=\dfrac{2/z}{z^2-8k^2-8-\dfrac{72k^4-144k^2+72}{z^2-32k^2-32-\ddots}}$$
Convergence type $E$ with $E=((1+k)/(1-k))^2$, $P=0$, and $C=...$, so that
$$\int_0^\infty\ssc^2(t,k)\nd^2(t,k)e^{-zt}\,dt-\dfrac{p(n)}{q(n)}\sim\dfrac{C}{((1+k)/(1-k))^{2n}}\;.$$
$$A=1-((z^2+k^2+1)/(8k))/n+\cdots$$
\end{cf}
Note: since the integrand has poles on the positive real $t$-axis, this CF is
purely formal.

Note that
$$\ssc^2(t,k)\nd^2(t,k)=\dfrac{1}{1-k^2}\dfrac{1-\cd(2t,k)}{1+\cd(2t,k)}\;,$$
so the above CF gives one for the Laplace transforms of $1/(1+\cd(t,k))$
and $1/(1+\dc(t,k))$.

\smallskip

\begin{cf}\label{6.19.36}{\ }
\begin{verbatim}
[(z,k)->intnum(t=0,[oo,z],ellsn(t,k)^2*elldc(t,k)^2*exp(-z*t)),
[0,z^2+8*n^2*(2*k^2-1)],[2/z,-2*n*(2*n+1)^2*(2*n+2)]]
\end{verbatim}
$$\int_0^\infty\sn^2(t,k)\dc^2(t,k)e^{-zt}\,dt=\dfrac{2/z}{z^2+16k^2-8-\dfrac{72}{z^2+64k^2-32-\ddots}}$$
Non convergent CF when $|k|<1$. It does converge for $|k|>1$, but we exclude
this.
\end{cf}
Note: since the integrand has poles on the positive real $t$-axis, this CF is
purely formal.

Note that $$\sn^2(t,k)\dc^2(t,k)=\dfrac{1-\cn(2t,k)}{1+\cn(2t,k)}\;,$$
so the above CF gives one for the Laplace transforms of $1/(1+\cn(t,k))$
and $1/(1+\nc(t,k))$.

\smallskip

Since the Laplace transform of the cube of a Jacobi elliptic function
is a simple linear combination of $1$ and the Laplace transform of the function
itself, it is immediate to obtain CFs for these. For instance, we have
$$2k^2\int_0^\infty \sn^3(t,k)e^{-zt}\,dt=-1+(z^2+k^2+1)\int_0^\infty\sn(t,k)e^{-zt}\,dt\;,$$
so we can obtain a CF for the LHS from \ref{6.19.1}.

\smallskip

Note that since the Dixon functions $\sm$ and $\cm$ and the lemniscatic
functions $\ssl$ and $\cl$ essentially correspond to $k=(1-\sqrt{-3})/2$ and
$k=1/\sqrt{2}$ respectively, it is immediate to obtain CFs for the Laplace
transforms of these functions from the above.

\medskip

\section{Complete Elliptic Integrals}

\medskip

We recall their definitions:

\smallskip

\begin{align*}
  K(k)&=\int_0^{\pi/2}\dfrac{d\theta}{\sqrt{1-k^2\sin^2(\theta)}}=\int_0^1\dfrac{dt}{\sqrt{(1-t^2)(1-k^2t^2)}}\\
  E(k)&=\int_0^{\pi/2}\sqrt{1-k^2\sin^2(\theta)}\,d\theta=\int_0^1\dfrac{\sqrt{1-k^2t^2}}{\sqrt{1-t^2}}\,dt\\
  D(k)&=\int_0^{\pi/2}\dfrac{\sin^2(\theta)d\theta}{\sqrt{1-k^2\sin^2(\theta)}}=\int_0^1\dfrac{t^2dt}{\sqrt{(1-t^2)(1-k^2t^2)}}=\dfrac{K(k)-E(k)}{k^2}\;.\end{align*}

\smallskip

$K(k)$ and $E(k)$ are implemented in {\tt Pari/GP} under the names
{\tt ellK(k)} and {\tt ellE(k)} respectively, and so we also set:

\begin{verbatim}
ellD(k)=(ellK(k)-ellE(k))/k^2;
\end{verbatim}

\medskip

\begin{cf}\label{6.20.1}{\ }
\begin{verbatim}
[k->(2/Pi)*ellK(k),[1,4,(2*n-1)^2*k^2+4*n^2],k^2*[1,-4*n^2*(2*n+1)^2]]
\end{verbatim}
$$\dfrac{2}{\pi}K(k)=1+\dfrac{k^2}{4-\dfrac{36k^2}{9k^2+16-\dfrac{400k^2}{25k^2+36-\dfrac{1764k^2}{49k^2+64-\dfrac{5184k^2}{81k^2+100-\ddots}}}}}$$
Convergence type $E$ with $E=1/k^2$, $P=1$, and $C=k^2/(1-k^2)/\pi$, so that
$$\dfrac{2}{\pi}K(k)-\dfrac{p(n)}{q(n)}\sim\dfrac{1/(\pi(1-k^2))}{(1/k^2)^{n+1}n}\;.$$
$$A=1-((5-k^2)/(4(1-k^2)))/n+((k^2+7)^2/(32(1-k^2)^2))/n^2+\cdots$$
Series:
$$\dfrac{2}{\pi}K(k)=1+\dfrac{k^2}{4}\sum_{n\ge0}\dfrac{(3/2)_n^2}{(n+1)!^2}k^{2n}$$
\end{cf}

\smallskip

\begin{cf}\label{6.20.2}{\ }
\begin{verbatim}
[k->(2/Pi)*ellK(k),
[1,4,(4*n-3)^2*k^2*(1-k^2)+4*n^2],k^2*(1-k^2)*[1,-4*n^2*(4*n+1)^2]]
\end{verbatim}
$$\dfrac{2}{\pi}K(k)=1-\dfrac{k^4-k^2}{4+\dfrac{100k^4-100k^2}{-25k^4+25k^2+16+\dfrac{1296k^4-1296k^2}{-81k^4+81k^2+36+\ddots}}}$$
Convergence type $E$ with $E=1/(4k^2(1-k^2))$, $P=3/2$, and
$C=4k^2(1-k^2)/((1-2k^2)^2\G^2(1/4))$, so that
$$\dfrac{2}{\pi}K(k)-\dfrac{p(n)}{q(n)}\sim\dfrac{1/((1-2k^2)^2\G^2(1/4))}{(1/(4k^2(1-k^2)))^{n+1}n^{3/2}}\;.$$
$$A=1-(3(4k^4-4k^2+9)/(16(2k^2-1)^2))/n+\cdots$$
Series:
$$\dfrac{2}{\pi}K(k)=1+\dfrac{k^2(1-k^2)}{4}\sum_{n\ge0}\dfrac{(5/4)_n^2}{(n+1)!^2}(4k^2(1-k^2))^n$$
\end{cf}

\smallskip

\begin{cf}\label{6.20.3}{\ }
\begin{verbatim}
[k->(2/Pi)*ellE(k),[1-k^2,4,4*(k^2+1)*n^2-k^2],
                   [3*k^2*(1-k^2),-4*n^2*(2*n+1)*(2*n+3)*k^2]]
\end{verbatim}
$$\dfrac{2}{\pi}E(k)=1-k^2-\dfrac{3k^4-3k^2}{4-\dfrac{60k^2}{15k^2+16-\dfrac{560k^2}{35k^2+36-\dfrac{2268k^2}{63k^2+64-\ddots}}}}$$
Convergence type $E$ with $E=1/k^2$, $P=0$, and $C=2k^2/\pi$, so that
$$\dfrac{2}{\pi}E(k)-\dfrac{p(n)}{q(n)}\sim\dfrac{2/\pi}{(1/k^2)^{n+1}}\;.$$
$$A=1-(1/4)/n+((13-5k^2)/(32(1-k^2)))/n^2+\cdots$$
Series:
$$\dfrac{2}{\pi}E(k)=1-k^2+\dfrac{k^2(1-k^2)}{4}\sum_{n\ge0}\dfrac{(2n+3)(3/2)_n^2}{(n+1)!^2}k^{2n}$$
\end{cf}
      
\smallskip

\begin{cf}\label{6.20.4}{\ }
\begin{verbatim}
[k->(2/Pi)*ellE(k),[1,4,(2*n-1)*(2*n-3)*k^2+4*n^2],
                                       -k^2*[1,4*n^2*(4*n^2-1)]]
\end{verbatim}
$$\dfrac{2}{\pi}E(k)=1-\dfrac{k^2}{4-\dfrac{12k^2}{3k^2+16-\dfrac{240k^2}{15k^2+36-\dfrac{1260k^2}{35k^2+64-\dfrac{4032k^2}{63k^2+100-\ddots}}}}}$$
Convergence type $E$ with $E=1/k^2$, $P=2$, and $C=-k^2/(2\pi(1-k^2))$, so that
$$\dfrac{2}{\pi}E(k)-\dfrac{p(n)}{q(n)}\sim-\dfrac{1/(2\pi(1-k^2))}{(1/k^2)^{n+1}n^2}\;.$$
$$A=1-((k^2+7)/(4(1-k^2)))/n-((3k^4-46k^2+21)/(32(1-k^2)^2))/n^2+\cdots$$
Series:
$$\dfrac{2}{\pi}E(k)=1-\dfrac{k^2}{4}\sum_{n\ge0}\dfrac{(2n+1)(1/2)_n^2}{(n+1)!^2}k^{2n}$$
\end{cf}

\smallskip

\begin{cf}\label{6.20.5}{\ }
\begin{verbatim}
[k->(4/Pi)*ellD(k),
[1,8,(4*n^2-1)*k^2+4*n*(n+1)],k^2*[3,-4*n*(n+1)*(2*n+1)*(2*n+3)]]
\end{verbatim}
$$\dfrac{4}{\pi}D(k)=1+\dfrac{3k^2}{8-\dfrac{120k^2}{15k^2+24-\dfrac{840k^2}{35k^2+48-\dfrac{3024k^2}{63k^2+80-\dfrac{7920k^2}{99k^2+120-\ddots}}}}}$$
Convergence type $E$ with $E=1/k^2$, $P=1$, and $C=2k^2/(\pi(1-k^2))$, so that
$$\dfrac{4}{\pi}D(k)-\dfrac{p(n)}{q(n)}\sim\dfrac{2/(\pi(1-k^2))}{(1/k^2)^{n+1}n}\;.$$
$$A=1+((3k^2-7)/(4(1-k^2)))/n+(45/32)/n^2+\cdots$$
Series:
$$\dfrac{4}{\pi}D(k)=1+\dfrac{k^2}{4}\sum_{n\ge0}\dfrac{(2n+3)(3/2)_n^2}{(n+2)(n+1)!^2}k^{2n}$$
\end{cf}

\smallskip

\begin{cf}\label{6.20.6}{\ }
\begin{verbatim}
[k->(2/(3*Pi))*((4-2*k^2)*ellE(k)-(1-k^2)*ellK(k)),
[1,4,(2*n-1)*(2*n-5)*k^2+4*n^2],[-3*k^2,-4*n^2*(2*n+1)*(2*n-3)*k^2]]
\end{verbatim}
$$\dfrac{2}{3\pi}((4-2k^2)E(k)-(1-k^2)K(k))=1-\dfrac{3k^2}{4+\dfrac{12k^2}{-3k^2+16-\dfrac{80k^2}{5k^2+36-\dfrac{756k^2}{21k^2+64-\dfrac{2880k^2}{45k^2+100-\ddots}}}}}$$
Convergence type $E$ with $E=1/k^2$, $P=3$, and $C=3k^2/(4\pi(1-k^2))$, so that
$$\dfrac{2}{3\pi}((4-2k^2)E(k)-(1-k^2)K(k))-\dfrac{p(n)}{q(n)}\sim\dfrac{3/(4\pi(1-k^2))}{(1/k^2)^{n+1}n^3}\;.$$
$$A=1-((7k^2+5)/(4(1-k^2)))/n+((89k^4+238k^2+57)/(32(1-k^2)^2))/n^2+\cdots$$
Series:
$$\dfrac{2}{3\pi}((4-2k^2)E(k)-(1-k^2)K(k))=1+\dfrac{3k^2}{4}\sum_{n\ge0}\dfrac{(2n+1)(1/2)_n^2}{(2n-1)(n+1)!^2}k^{2n}$$
\end{cf}

There exist many CFs as above for suitable linear combinations of $E(k)$ and
$K(k)$.

\smallskip

\begin{cf}\label{6.20.7}{\ }
\begin{verbatim}
[k->(4/Pi^2)*ellK(k)^2,
[1,2,(2*n-1)^3*k^2*(1-k^2)+2*n^3],k^2*(1-k^2)*[1,-2*n^3*(2*n+1)^3]]
\end{verbatim}
$$\dfrac{4}{\pi^2}K(k)^2=1-\dfrac{k^4-k^2}{2+\dfrac{54k^4-54k^2}{-27k^4+27k^2+16+\dfrac{2000k^4-2000k^2}{-125k^4+125k^2+54+\ddots}}}$$
Convergence type $E$ with $E=1/(4k^2(1-k^2))$, $P=3/2$, and $C=4k^2(1-k^2)/(1-2k^2)^2/\pi^{3/2}$, so that
$$\dfrac{4}{\pi^2}K(k)^2-\dfrac{p(n)}{q(n)}\sim\dfrac{1/(\pi^{3/2}(1-2k^2)^2)}{(1/(4k^2(1-k^2)))^{n+1}n^{3/2}}\;.$$
$$A=1-(3(5k^4-4k^2+4)/(8(2-k^2)^2))/n+\cdots$$
Series:
$$\dfrac{4}{\pi^2}K(k)^2=1+\dfrac{k^2(1-k^2)}{2}\sum_{n\ge0}\dfrac{(3/2)_n^3}{(n+1)!^3}(4k^2(1-k^2))^n$$
\end{cf}

\smallskip

We now give a number of CFs for $K(k)/E(k)$.
Since $D(k)=(K(k)-E(k))/k^2$, it is immediate to obtain CFs for $K(k)/D(k)$
from CFs for $K(k)/E(k)$. We will thus simply indicate the other CF with no
extra information at the end.

\smallskip

\begin{cf}\label{6.20.8}{\ }
\begin{verbatim}
[k->ellK(k)/ellE(k),[1,2*k^2+(2-4*k^2)*n],[k^2,(4*n^2-1)*k^2*(1-k^2)]]
\end{verbatim}
$$\dfrac{K(k)}{E(k)}=1+\dfrac{k^2}{-2k^2+2-\dfrac{3k^4-3k^2}{-6k^2+4-\dfrac{15k^4-15k^2}{-10k^2+6-\dfrac{35k^4-35k^2}{-14k^2+8-\ddots}}}}$$
Convergence type $E$ with $E=-(1-k^2)/k^2$, $P=1$, and $C=(k^2\pi/4)/E^2(k)$,
so that
$$\dfrac{K(k)}{E(k)}-\dfrac{p(n)}{q(n)}\sim(-1)^n\dfrac{(k^2\pi/4)/E^2(k)}{((1-k^2)/k^2)^nn}\;.$$
Note: if $k^2>1/2$, replace $E$ by $1/E$ and $P$ by $-P$.
$$A=1+(k^2/2-3/4)/n+(11k^4/8-15k^2/8+21/32)/n^2+\cdots$$
Parametric family for $m\ge0$ and $u\in\Z$: (note: since $k$ is used, we
cannot use it):
\begin{verbatim}
[k->ellK(k)/ellE(k),(4*u+2)*k^2+(2-4*k^2)*n+2*m,
                    (2*n-1-2*u)*(2*n+1-2*u)*k^2*(1-k^2)]
\end{verbatim}
Convergence type $E$ with $E=-(1-k^2)/k^2$ and $P=2u+2m+1$.
\begin{verbatim}
[k->ellK(k)/ellD(k),[2-k^2,2-2*k^2+(2-4*k^2)*n],
                    [(2*n+1)*(2*n+3)*k^2*(1-k^2)]]
\end{verbatim}
\end{cf}

\smallskip

\begin{cf}\label{6.20.9}{\ }
\begin{verbatim}
[k->ellK(k)/ellE(k),[1,k^2+2,2*(k^2+1)*n],k^2*[1,-(2*n+1)^2]]
\end{verbatim}
$$\dfrac{K(k)}{E(k)}=1+\dfrac{k^2}{k^2+2-\dfrac{9k^2}{4k^2+4-\dfrac{25k^2}{6k^2+6-\dfrac{49k^2}{8k^2+8-\dfrac{81k^2}{10k^2+10-\ddots}}}}}$$
Convergence type $E$ with $E=1/k^2$, $P=0$, and $C=(k^2\pi/2)/E^2(k)$, so that
$$\dfrac{K(k)}{E(k)}-\dfrac{p(n)}{q(n)}\sim\dfrac{\pi/(2E^2(k))}{(1/k)^{2n+2}}\;.$$
$$A=1-((1+k^2)/(4(1-k^2)))/n+((-7k^4+2k^2+9)/(32(1-k^2)^2))/n^2+\cdots$$
Parametric family for all $m$:
\begin{verbatim}
[k->ellK(k)/ellE(k),2*(k^2+1)*n+2*m*(1-k^2),-(2*n+1)^2*k^2]
\end{verbatim}
Convergence type $E$ with $E=1/k^2$ and $P=2m$.
\begin{verbatim}
[k->ellK(k)/ellD(k),[2*(k^2+1)*(n+1)],[-(2*n+3)^2*k^2]]
\end{verbatim}
\end{cf}

\smallskip

\begin{cf}\label{6.20.10}{\ }
\begin{verbatim}
[k->ellK(k)/ellE(k),[1,k^2,2*(k^2+1)*n-4],k^2*[1,-(2*n-1)^2]]
\end{verbatim}
$$\dfrac{K(k)}{E(k)}=1+\dfrac{k^2}{k^2-\dfrac{k^2}{4k^2-\dfrac{9k^2}{6k^2+2-\dfrac{25k^2}{8k^2+4-\dfrac{49k^2}{10k^2+6-\dfrac{81k^2}{12k^2+8-\ddots}}}}}}$$
Convergence type $E$ with $E=1/k^2$, $P=-2$, and $C=2\pi(1-k^2)^2/(k^2E^2(k))$,
so that
$$\dfrac{K(k)}{E(k)}-\dfrac{p(n)}{q(n)}\sim\dfrac{2\pi(1-k^2)^2/E^2(k)}{(1/k)^{2n-2}n^{-2}}\;.$$
$$A=1-(5(k^2+1)/(4(1-k^2)))/n+((9k^4-62k^2+9)/(32(1-k^2)^2))/n^2+\cdots$$
Parametric family for all $m$ and $u$:
\begin{verbatim}
[k->ellK(k)/ellE(k),2*(k^2+1)*n+4*u+2*m*(1-k^2),-k^2*(2*n+2*u+1)^2]
\end{verbatim}
Convergence type $E$ with $E=1/k^2$ and $P=2u+2m$.
\begin{verbatim}
[k->ellK(k)/ellD(k),[2*(k^2+1)*(n+1)],[-k^2*(2*n+3)^2]]
\end{verbatim}
\end{cf}
      
\smallskip

\begin{cf}\label{6.20.11}{\ }
\begin{verbatim}
[k->ellK(k)/ellE(k),(2-k^2)*[0,1,4*(n-1)],[2,-(2*n-1)^2*k^4]]
\end{verbatim}
$$\dfrac{K(k)}{E(k)}=\dfrac{2}{-k^2+2-\dfrac{k^4}{-4k^2+8-\dfrac{9k^4}{-8k^2+16-\dfrac{25k^4}{-12k^2+24-\dfrac{49k^4}{-16k^2+32-\ddots}}}}}$$
Convergence type $E$ with $E=(1+\sqrt{1-k^2})^4/k^4$, $P=0$, and
$C=\pi/E^2(k)$, so that
$$\dfrac{K(k)}{E(k)}-\dfrac{p(n)}{q(n)}\sim\dfrac{\pi/E^2(k)}{((1+\sqrt{1-k^2})/k)^{4n}}\;.$$
\begin{verbatim}
[k->ellK(k)/ellD(k),[2,8-3*k^2,(8-4*k^2)*n],[-2*k^2,-(2*n+1)^2*k^4]]
\end{verbatim}
\end{cf}
      
\smallskip

\begin{cf}\label{6.20.12}{\ }
\begin{verbatim}
[k->ellK(k)/ellE(k),[1,8-5*k^2,4*(2-k^2)*n*(n-1)-k^2],
                    [4*k^2,-6*k^4,-(n^2-1)*(4*n^2-1)*k^4]]
\end{verbatim}
$$\dfrac{K(k)}{E(k)}=1+\dfrac{4k^2}{8-5k^2-\dfrac{6k^4}{16-9k^2-\dfrac{45k^4}{48-25k^2-\dfrac{280k^4}{96-49k^2-\dfrac{945k^4}{160-81k^2-\ddots}}}}}$$
Convergence type $E$ with $E=(1+\sqrt{1-k^2})^4/k^4$, $P=0$, and
$C=\pi/(((1+\sqrt(1-k^2))^2/k^2)E^2(k))$, so that
$$\dfrac{K(k)}{E(k)}-\dfrac{p(n)}{q(n)}\sim\dfrac{\pi/E^2(k)}{((1+\sqrt{1-k^2})/k)^{4n+2}}\;.$$
Parametric family for all $m$:
\begin{verbatim}
[k->ellK(k)/ellE(k),4*(2-k^2)*n*(n-1)+(2*m+1)*k^2,
                    -(n^2-1)*(4*n^2-(2*m+1)^2)*k^4]
\end{verbatim}
Convergence type $E$ with $E=(1+\sqrt{1-k^2})^4/k^4$ and $P=0$.
\begin{verbatim}
[k->ellK(k)/ellD(k),[2-k^2/4,4*(2-k^2)*n*(n+1)-k^2],
               -k^4*[3/2,n*(n+2)*(2*n+1)*(2*n+3)]]
\end{verbatim}
\end{cf}

\smallskip

\begin{cf}\label{6.20.13}{\ }
\begin{verbatim}
[k->ellK(2*k/(1+k^2))/ellE(2*k/(1+k^2)),[0,1,4*(2-k^4)*(n-1)],
                                        [(k^2+1)^2,-k^8*(2*n-1)^2]]
\end{verbatim}
$$\dfrac{K(2k/(1+k^2))}{E(2k/(1+k^2))}=\dfrac{k^4+2k^2+1}{1-\dfrac{k^8}{-4k^4+8-\dfrac{9k^8}{-8k^4+16-\dfrac{25k^8}{-12k^4+24-\dfrac{49k^8}{-16k^4+32-\dfrac{81k^8}{-20k^4+40-\ddots}}}}}}$$
Convergence type $E$ with $E=((1+\sqrt{1-k^4})/k^2)^4$, $P=0$, and $C=...$,
so that
$$\dfrac{K(2k/(1+k^2))}{E(2k/(1+k^2))}-\dfrac{p(n)}{q(n)}\sim\dfrac{C}{((1+\sqrt{1-k^4})/k^2)^{4n}}$$
\end{cf}

\smallskip

Note that this is simply a Landen transformation applied to \ref{6.20.11}.

\medskip

\section{Some More Integrals}

\medskip

In addition to the CFs given in the present section, note that Sections
\ref{sec:beta}, \ref{sec:beta1}, \ref{sec:psi0}, \ref{sec:psi1}, and
\ref{sec:psi2} give CFs for
\begin{align*}&\int_0^\infty\dfrac{e^{-zt}}{\cosh(t)}\,dt\;,\quad \int_0^\infty\dfrac{e^{-zt}}{\cosh^2(t)}\,dt\;,\quad \int_0^\infty e^{-zt}\tanh(t)\,dt\;,\\
  &\int_0^\infty e^{-zt}\tanh^2(t)\,dt\;,\quad \int_0^\infty\dfrac{te^{-zt}}{\cosh(t)}\,dt\;,\quad \int_0^\infty e^{-zt}t\tanh(t)\,dt\,\\
  &\int_0^\infty\dfrac{(tz+1)e^{-zt}}{\cosh^2(t)}\,dt\;,\quad \int_0^\infty e^{-zt}(tz+1)\tanh^2(t)\,dt\;,\\
  &\int_0^\infty\dfrac{1-e^{-tz}}{\sinh(t)}\,dt\;,\quad\int_0^\infty\dfrac{e^{-tz}-1+tz}{\sinh^2(t)}\,dt\;,\quad\int_0^\infty\dfrac{te^{-zt}}{\sinh(t)}\,dt\;,\\
  &\int_0^\infty(1-e^{-tz})(\cotanh(t)-1)\,dt\;,\quad\int_0^\infty\dfrac{t^2e^{-zt}}{\sinh(t)}\,dt\;,\\
  &\int_0^\infty te^{-zt}\cotanh(t)\,dt\text{\quad and\quad}\int_0^\infty t^2e^{-zt}\cotanh(t)\,dt\;.\end{align*}

In addition:

\medskip

\begin{cf}\label{6.21.1}{\ }
\begin{verbatim}
[(z,k)->intnum(t=0,[oo,k+z],exp(-z*t)/cosh(t)^k),[0,z],[1,n*(n+k-1)]]
\end{verbatim}
$$\int_0^\infty\dfrac{e^{-zt}}{\cosh^k(t)}\,dt=\dfrac{1}{z+\dfrac{k}{z+\dfrac{2k+2}{z+\dfrac{3k+6}{z+\dfrac{4k+12}{z+\dfrac{5k+20}{z+\ddots}}}}}}$$
Convergence type $P^-$ with $P=z$ and $C=2^{k-z-1}\G((k+z)/2)^2/\G(k)$, so that
$$\int_0^\infty\dfrac{e^{-zt}}{\cosh^k(t)}\,dt-\dfrac{p(n)}{q(n)}\sim(-1)^n\dfrac{2^{k-z-1}\G((k+z)/2)^2/\G(k)}{n^z}\;.$$
$$A=1-(kz/2)/n-((z^3-6k^2z^2-(9k^2-6k-8)z)/48)/n^2+\cdots$$
Parametric family for $m\ge0$:
\begin{verbatim}
[(z,k)->intnum(t=0,[oo,k+z],exp(-z*t)/cosh(t)^k),2*m+z,n*(n+k-1)]
\end{verbatim}
Convergence type $P^-$ with $P=2m+z$.
\end{cf}

\smallskip

\begin{cf}\label{6.21.2}{\ }
\begin{verbatim}
[(z,k)->intnum(t=0,[oo,k+z],exp(-z*t)/cosh(t)^k),
[0,z^2+8*n^2+(k-4)*4*n-(3*k-8)],[z,-2*n*(2*n-1)*(2*n+k-1)*(2*n+k-2)]]
\end{verbatim}
$$\int_0^\infty\dfrac{e^{-zt}}{\cosh^k(t)}\,dt=\dfrac{z}{z^2+k-\dfrac{2k^2+2k}{z^2+5k+8-\dfrac{12k^2+60k+72}{z^2+9k+32-\dfrac{30k^2+270k+600}{z^2+13k+72-\ddots}}}}$$
Convergence type $P^+$ with $P=z$ and $C=2^{k-2z-1}\G((k+z)/2)^2/\G(k)$,
so that
$$\int_0^\infty\dfrac{e^{-zt}}{\cosh^k(t)}\,dt-\dfrac{p(n)}{q(n)}\sim\dfrac{2^{k-2z-1}\G((k+z)/2)^2/\G(k)}{n^z}\;.$$
$$A=1-(kz/4)/n-((z^3-6k^2z^2-(9k^2-6k-8)z)/192)/n^2+\cdots$$
Parametric family for $m\ge0$:
\begin{verbatim}
[(z,k)->intnum(t=0,[oo,k+z],exp(-z*t)/cosh(t)^k),
z^2+8*n^2+(k-4)*4*n-(3*k-8)+4*m*(z+3),-2*n*(2*n-1)*(2*n+k-1)*(2*n+k-2)]
\end{verbatim}
Convergence type $P^+$ with $P=2m+z$.
\end{cf}

This is simply the contraction of the previous CF.

\smallskip

\begin{cf}\label{6.21.3}{\ }
\begin{verbatim}
[(z,k)->intnum(t=0,[oo,k+z],exp(-z*t)/cosh(t)^k),
[0,k+z,2*k],[2,(2*n+z-k)*(2*n+z+k-2)]]
\end{verbatim}
$$\int_0^\infty\dfrac{e^{-zt}}{\cosh^k(t)}\,dt=\dfrac{2}{z+k+\dfrac{z^2+2z-k^2+2k}{2k+\dfrac{z^2+6z-k^2+2k+8}{2k+\dfrac{z^2+10z-k^2+2k+24}{2k+\dfrac{z^2+14z-k^2+2k+48}{2k+\ddots}}}}}$$
Convergence type $P^-$ with $P=k$ and $C=\G((z+k)/2)/(2\G((z+2-k)/2))$, so that
$$\int_0^\infty\dfrac{e^{-zt}}{\cosh^k(t)}\,dt-\dfrac{p(n)}{q(n)}\sim(-1)^n\dfrac{\G((z+k)/2)/(2\G((z+2-k)/2))}{n^k}\;.$$
$$A=1-(zk/2)/n+(k^3/24+(z^2/8-1/8)k^2+(z^2/8-1/6)k)/n^2+\cdots$$
Series:
$$\int_0^\infty\dfrac{e^{-zt}}{\cosh^k(t)}\,dt=\dfrac{2}{z+k}\sum_{n\ge0}(-1)^n\dfrac{((z-k)/2+1)_n}{((z+k)/2+1)_n}$$
Parametric family for $m\ge0$:
\begin{verbatim}
[(z,k)->intnum(t=0,[oo,k+z],exp(-z*t)/cosh(t)^k),
2*k+4*m,(2*n+z-k)*(2*n+z+k-2)]
\end{verbatim}
Convergence type $P^-$ with $P=2m+k$.
\end{cf}
    
\smallskip

\begin{cf}\label{6.21.4}{\ }
\begin{verbatim}
[(z,k)->intnum(t=0,[oo,k+z],exp(-z*t)/cosh(t)^k),[0,2*n+z+k-2],
[[2,(z+k)*(2+z-k)],[4*n*(n+k-1),(2*n+z+k)*(2*n+2+z-k)]]]
\end{verbatim}
$$\int_0^\infty\dfrac{e^{-zt}}{\cosh^k(t)}\,dt=\dfrac{2}{2k+1+\dfrac{z^2+2z-k^2+2k}{z+k+2+\dfrac{4k}{z+k+4+\dfrac{z^2+6z-k^2+2k+8}{z+k+6+\ddots}}}}$$
Convergence type $E$ with $E=-(1+\sqrt{2})^2$, $P=0$, and
$C=\pi 2^k\G((z+k)/2)/(\G((z-k)/2+1)\G(k))/(1+\sqrt{2})^{z+k}$, so that
$$\int_0^\infty\dfrac{e^{-zt}}{\cosh^k(t)}\,dt-\dfrac{p(n)}{q(n)}\sim(-1)^n\dfrac{\pi 2^k(\G((z+k)/2)/(\G((z-k)/2+1)\G(k)))}{(1+\sqrt{2})^{2n+z+k}}\;.$$
\end{cf}

\smallskip

More generally:

\smallskip

\begin{cf}\label{6.21.5}{\ }
\begin{verbatim}
[(z,k,a)->intnum(t=0,[oo,k+z],exp(-z*t)/(cosh(t)+a*sinh(t))^k),
[0,z+a*(k+2*n-2)],[1,n*(n+k-1)*(1-a^2)]]
\end{verbatim}
$$\int_0^\infty\dfrac{e^{-zt}}{(\cosh(t)+a\sinh(t))^k}\,dt=\dfrac{1}{ka+z-\dfrac{ka^2-k}{(k+2)a+z-\dfrac{(2k+2)a^2-2k-2}{(k+4)a+z-\ddots}}}$$
Convergence type $E$ with $E=(a+1)/(a-1)$, $P=z$ and $C=...$, so that
$$\int_0^\infty\dfrac{e^{-zt}}{(\cosh(t)+a\sinh(t))^k}\,dt-\dfrac{p(n)}{q(n)}\sim\dfrac{C}{((a+1)/(a-1))^nn^z}\;.$$
$$A=1-((az^2+2kz-a(k^2-2k))/4)/n+\cdots$$
Parametric family for $m\ge0$:
\begin{verbatim}
[(z,k,a)->intnum(t=0,[oo,k+z],exp(-z*t)/(cosh(t)+a*sinh(t))^k),
z+a*(k+2*n-2)+2*m,n*(n+k-1)*(1-a^2)]
\end{verbatim}
Convergence type $E$ with $E=(a+1)/(a-1)$ and $P=2m+z$.
\end{cf}

\smallskip

\begin{cf}\label{6.21.6}{\ }
\begin{verbatim}
[(z,k,a)->intnum(t=0,[oo,k+z],exp(-z*t)/(cosh(t)+a*sinh(t))^k),
[0,(a+1)*(k+z),4*a*n+2*k+2*a*z-4*a],[2,(1-a^2)*(2*n+z-k)*(2*n+z+k-2)]]
\end{verbatim}
$$\int_0^\infty\dfrac{e^{-zt}}{(\cosh(t)+a\sinh(t))^k}\,dt=\dfrac{2}{(a+1)k+(a+1)z+\dfrac{N}{2k+2az+4a+\ddots}}$$
with $N=(a^2-1)k^2+(-2a^2+2)k+((-a^2+1)z^2+(-2a^2+2)z)$.

Convergence type $E$ with $E=(a+1)/(a-1)$, $P=k$, and
$C=\G((z+k)/2)/(2\G((z+2-k)/2))$, so that
$$\int_0^\infty\dfrac{e^{-zt}}{(\cosh(t)+a\sinh(t))^k}\,dt-\dfrac{p(n)}{q(n)}\sim\dfrac{\G((z+k)/2)/(2\G((z+2-k)/2))}{((a+1)/(a-1))^nn^k}\;.$$
$$A=1-(z/2+a/2)k/n+\cdots$$
Parametric family for $m\ge0$:
\begin{verbatim}
[(z,k,a)->intnum(t=0,[oo,k+z],exp(-z*t)/(cosh(t)+a*sinh(t))^k),
4*a*n+2*k+2*a*z-4*a+4*m,(1-a^2)*(2*n+z-k)*(2*n+z+k-2)]
\end{verbatim}
Convergence type $E$ with $E=(a+1)/(a-1)$ and $P=2m+k$.
\end{cf}
  
\smallskip

\begin{cf}\label{6.21.7}{\ }
\begin{verbatim}
[(z,k,a)->intnum(t=0,[oo,k+z],exp(-z*t)/(cosh(t)+a*sinh(t))^k),
[0,(a+1)*(2*n+z+k-2)],
[[2,(1-a^2)*(z+k)*(2+z-k)],(1-a^2)*[4*n*(n+k-1),(2*n+z+k)*(2*n+2+z-k)]]]
\end{verbatim}
$$\int_0^\infty\dfrac{e^{-zt}}{(\cosh(t)+a\sinh(t))^k}\,dt=\dfrac{2}{2za+2k+1-\dfrac{(z^2+2z-k^2+2k)a^2-z^2-2z+k^2-2k}{(z+k+2)a+z+k+2-\ddots}}$$
Convergence type $E$ with $E=(a+3+2\sqrt{2a+2})/(a-1)$, $P=0$, and $C=...$,
so that
$$\int_0^\infty\dfrac{e^{-zt}}{(\cosh(t)+a\sinh(t))^k}\,dt-\dfrac{p(n)}{q(n)}\sim\dfrac{C}{((a+3+2\sqrt{2a+2})/(a-1))^n}\;.$$
\end{cf}

\smallskip

\begin{cf}\label{6.21.8}{\ }
\begin{verbatim}
[(z,a)->intnum(t=0,[oo,z+1-a],exp(-z*t)*cosh(a*t)/cosh(t)),
           [0,z],[[1,1-a^2],[4*n^2,(2*n+1+a)*(2*n+1-a)]]]
\end{verbatim}
$$\int_0^\infty\dfrac{\cosh(at)}{\cosh(t)}e^{-zt}\,dt=\dfrac{1}{z-\dfrac{a^2-1}{z+\dfrac{4}{z-\dfrac{a^2-9}{z+\dfrac{16}{z-\dfrac{a^2-25}{z+\ddots}}}}}}$$
Convergence type $P^-$ with $P=z$ and $C=...$, so that
$$\int_0^\infty\dfrac{\cosh(at)}{\cosh(t)}e^{-zt}\,dt-\dfrac{p(n)}{q(n)}\sim(-1)^n\dfrac{C}{n^z}\;.$$
$$A=1-((z^3-(a^2+1)z)/(2(z^2-1)))/n+\cdots$$
\end{cf}

\smallskip

\begin{cf}\label{6.21.9}{\ }
\begin{verbatim}
[(z,a)->intnum(t=0,[oo,z+1-a],exp(-z*t)*cosh(a*t)/cosh(t)),
        [0,8*n^2-12*n+5+z^2-a^2],[z,-4*n^2*(2*n-1+a)*(2*n-1-a)]]
\end{verbatim}
$$\int_0^\infty\dfrac{\cosh(at)}{\cosh(t)}e^{-zt}\,dt=\dfrac{z}{z^2-a^2+1+\dfrac{4a^2-4}{z^2-a^2+13+\dfrac{16a^2-144}{z^2-a^2+41+\dfrac{36a^2-900}{z^2-a^2+85+\ddots}}}}$$
Convergence type $P^+$ with $P=z$ and $C=...$, so that
$$\int_0^\infty\dfrac{\cosh(at)}{\cosh(t)}e^{-zt}\,dt-\dfrac{p(n)}{q(n)}\sim\dfrac{C}{n^z}\;.$$
$$A=1-((z^3-(a^2+1)z)/(4(z^2-1)))/n+\cdots$$
Parametric family for $k\ge0$:
\begin{verbatim}
[(z,a)->intnum(t=0,[oo,z+1-a],exp(-z*t)*cosh(a*t)/cosh(t)),
8*n^2-12*n+5+(2*k+z)^2-a^2,-4*n^2*(2*n-1+a)*(2*n-1-a)]
\end{verbatim}
Convergence type $P^+$ with $P=2k+z$.
\end{cf}

This is simply the contraction of the previous CF.

\smallskip

\begin{cf}\label{6.21.10}{\ }
\begin{verbatim}
[(z,a)->intnum(t=0,[oo,z+1-a],exp(-z*t)*sinh(a*t)/cosh(t)),
[[0,z^2-1],[z^2,z^2-1]],[[a,z^2*(4-a^2)],z^2*[4*n^2,(2*n+2+a)*(2*n+2-a)]]]
\end{verbatim}
$$\int_0^\infty\dfrac{\sinh(at)}{\cosh(t)}e^{-zt}\,dt=\dfrac{a}{z^2-1-\dfrac{(a^2-4)z^2}{z^2+\dfrac{4z^2}{z^2-1-\dfrac{(a^2-16)z^2}{z^2+\dfrac{16z^2}{z^2-1-\dfrac{(a^2-36)z^2}{z^2+\ddots}}}}}}$$
Convergence type $P^-$ with $P=z$ and $C=...$, so that
$$\int_0^\infty\dfrac{\sinh(at)}{\cosh(t)}e^{-zt}\,dt-\dfrac{p(n)}{q(n)}\sim(-1)^n\dfrac{C}{n^z}\;.$$
$$A=1-z/n-((z^5-24z^4-(6a^2+14)z^3+96z^2+(3a^4-30a^2-37)z)/(48(z^2-4)))/n^2+\cdots$$
\end{cf}

\smallskip

\begin{cf}\label{6.21.11}{\ }
\begin{verbatim}
[(z,a)->intnum(t=0,[oo,z+1-a],exp(-z*t)*sinh(a*t)/cosh(t)),
           [0,8*n^2-8*n+3+z^2-a^2],[a,-4*n^2*(2*n+a)*(2*n-a)]]
\end{verbatim}
$$\int_0^\infty\dfrac{\sinh(at)}{\cosh(t)}e^{-zt}\,dt=\dfrac{a}{z^2-a^2+3+\dfrac{4a^2-16}{z^2-a^2+19+\dfrac{16a^2-256}{z^2-a^2+51+\dfrac{36a^2-1296}{z^2-a^2+99+\ddots}}}}$$
Convergence type $P^+$ with $P=z$ and $C=...$, so that
$$\int_0^\infty\dfrac{\sinh(at)}{\cosh(t)}e^{-zt}\,dt-\dfrac{p(n)}{q(n)}\sim\dfrac{C}{n^z}\;.$$
$$A=1-(z/2)/n-((z^5-24z^4-(6a^2+14)z^3+96z^2+(3a^4-30a^2-37)z)/(192(z^2-4)))/n^2+\cdots$$
Parametric family for $k\ge0$:
\begin{verbatim}
[(z,a)->intnum(t=0,[oo,z+1-a],exp(-z*t)*sinh(a*t)/cosh(t)),
8*n^2-8*n+3+(2*k+z)^2-a^2,-4*n^2*(2*n+a)*(2*n-a)]
\end{verbatim}
Convergence type $P^+$ with $P=2k+z$.
\end{cf}

This is simply the contraction of the previous CF.

\smallskip

\begin{cf}\label{6.21.12}{\ }
\begin{verbatim}
[(z,a)->intnum(t=0,[oo,z+1-a],exp(-z*t)*sinh(a*t)/cosh(t)),
[0,(z+1)^2-a^2,4*(2*n+z-2)],[2*a,(2*n+z-a-1)^2*(2*n+z+a-1)^2]]
\end{verbatim}
$$\int_0^\infty\dfrac{\sinh(at)}{\cosh(t)}e^{-zt}\,dt=\dfrac{2a}{z^2+2z+(-a^2+1)+\dfrac{N}{4z+8+\ddots}}$$
with $N=z^4+4z^3+(-2a^2+6)z^2+(-4a^2+4)z+a^4-2a^2+1$.

\noindent
Convergence type $P^-$ with $P=2$ and $C=...$, so that
$$\int_0^\infty\dfrac{\sinh(at)}{\cosh(t)}e^{-zt}\,dt-\dfrac{p(n)}{q(n)}\sim(-1)^n\dfrac{C}{n^2}\;.$$
$$A=1-z/n+(3z^2/4+a^2/4-3/4)/n^2+(-z^3/2-a^2z/2+3z/2)/n^3+\cdots$$
Series:
$$\int_0^\infty\dfrac{\sinh(at)}{\cosh(t)}e^{-zt}\,dt=2a\sum_{n\ge1}\dfrac{(-1)^{n+1}}{(2n+z+a-1)(2n+z-a-1)}$$
Parametric family for $k\ge0$:
\begin{verbatim}
[(z,a)->intnum(t=0,[oo,z+1-a],exp(-z*t)*sinh(a*t)/cosh(t)),
4*(2*k+1)*(2*n+z-2),(2*n+z-a-1)^2*(2*n+z+a-1)^2]
\end{verbatim}
Convergence type $P^-$ with $P=4k+2$.
\end{cf}

\medskip

\begin{cf}\label{6.21.13}{\ }
\begin{verbatim}
[(p,q,z)->intnum(t=0,z,t^p/(1+t^q)),
[0,p+1+(n-1)*q],[[z^(p+1),(p+1)^2*z^q],z^q*[n^2*q^2,(p+1+n*q)^2]]]
\end{verbatim}
$$\int_0^z\dfrac{t^p\,dt}{1+t^q}=\dfrac{z^{p+1}}{p+1+\dfrac{z^qp^2+2z^qp+z^q}{p+q+1+\dfrac{q^2z^q}{p+2q+1+\dfrac{z^qp^2+(2q+2)z^qp+q^2+2q+1)z^q}{p+3q+1+\dfrac{4q^2z^q}{p+4q+1+\ddots}}}}}$$
Convergence type $E$ with $E=-((1+\sqrt{1+z^q})^2/z^q)$, $P=0$, and $C=...$,
so that
$$\int_0^z\dfrac{t^p\,dt}{1+t^q}-\dfrac{p(n)}{q(n)}\sim(-1)^n\dfrac{C}{((1+\sqrt{1+z^q})^2/z^q)^n}\;.$$
\end{cf}

\smallskip

\begin{cf}\label{6.21.14}{\ }
\begin{verbatim}
[(p,q,z)->z^(p+1)/(p+1)-intnum(t=0,z,t^p/(1+t^q)),
[0,p+1+n*q],[[z^(p+1+q),(p+1+q)^2*z^q],z^q*[n^2*q^2,(p+1+(n+1)*q)^2]]]
\end{verbatim}
$$\dfrac{z^{p+1}}{p+1}-\int_0^z\dfrac{t^p\,dt}{1+t^q}=\dfrac{z^{p+1+q}}{p+q+1+\dfrac{z^qp^2+(2q+2)z^qp+(q^2+2q+1)z^q}{p+2q+1+\dfrac{q^2z^q}{p+3q+1+\ddots}}}$$
Convergence type $E$ with $E=-((1+\sqrt{1+z^q})^2/z^q)$, $P=0$, and $C=...$,
so that
$$\dfrac{z^{p+1}}{p+1}-\int_0^z\dfrac{t^p\,dt}{1+t^q}-\dfrac{p(n)}{q(n)}\sim(-1)^n\dfrac{C}{((1+\sqrt{1+z^q})^2/z^q)^n}\;.$$
\end{cf}

\smallskip

\begin{cf}\label{6.21.15}{\ }
\begin{verbatim}
[(k,m,p,q)->intnum(t=0,1,t^(k-1)/(1+t^m)^(p/q)),
[0,k,((2*n-3)*q-(n-2)*p)*m+(q-p)*k],
[1,p*k^2,(n-1)*q*(p+(n-1)*q)*(k+(n-1)*m)^2]]
\end{verbatim}
$$\int_0^1\dfrac{t^{k-1}}{(1+t^m)^{p/q}}\,dt=\dfrac{1}{k+\dfrac{pk^2}{(-p+q)k+qm+\dfrac{(qp+q^2)k^2+(2qp+2q^2)mk+(qp+q^2)m^2}{(-p+q)k+(-p+3q)m+\ddots}}}$$
Convergence type $P^-$ with $P=2-p/q$ and $C=...$, so that
$$\int_0^1\dfrac{t^{k-1}}{(1+t^m)^{p/q}}\,dt-\dfrac{p(n)}{q(n)}\sim(-1)^n\dfrac{C}{n^{2-p/q}}\;.$$
\end{cf}

\medskip

\section{Bessel and Related Functions}\label{sec:bessel}

\medskip

Since $J_{\nu}(z)/z^{\nu}=I_{\nu}(iz)/(iz)^{\nu}$, all the continued fractions
for $J$ can be deduced from those of $I$. The CFs for $I$ and $J$ all have
factorial convergence and those for $K$ or $H^{(i)}$ all have
subexponential convergence. There seems to be none involving $Y_{\nu}(z)$.

\smallskip

When we indicate a parametric solution, recall that in the case of
subexponential convergence Bauer--Muir does not accelerate, so the
speed of convergence is the same as that of the original CF, hence is
not indicated.

\medskip

\begin{cf}\label{6.22.1}{\ }
\begin{verbatim}
[(nu,z)->besseli(nu+1,z)/besseli(nu,z),[0,2*(nu+n)],[z,z^2]]
\end{verbatim}
$$\dfrac{I_{\nu+1}(z)}{I_{\nu}(z)}=\dfrac{z}{2\nu+2+\dfrac{z^2}{2\nu+4+\dfrac{z^2}{2\nu+6+\dfrac{z^2}{2\nu+8+\dfrac{z^2}{2\nu+10+\dfrac{z^2}{2\nu+12+\ddots}}}}}}$$
Convergence type $F^2$ with $E=-4/z^2$, $P=2\nu+1$, and
$C=(z/2)^{2\nu+1}/I_{\nu}(z)^2$, so that
$$\dfrac{I_{\nu+1}(z)}{I_{\nu}(z)}-\dfrac{p(n)}{q(n)}\sim(-1)^n\dfrac{1/I_{\nu}(z)^2}{n!^2(4/z^2)^{n-\nu-1/2}n^{2\nu+1}}\;.$$
$$A=1+(z^2/2-(\nu+1)^2)/n+\cdots$$
\end{cf}

\begin{cf}\label{6.22.2}{\ }
\begin{verbatim}
[(nu,z)->besseli(nu+2,z)/besseli(nu,z),[1,2*(nu+n)],[-2*(nu+1),z^2]]
\end{verbatim}
$$\dfrac{I_{\nu+2}(z)}{I_{\nu}(z)}=1-\dfrac{2\nu+2}{2\nu+2+\dfrac{z^2}{2\nu+4+\dfrac{z^2}{2\nu+6+\dfrac{z^2}{2\nu+8+\dfrac{z^2}{2\nu+10+\dfrac{z^2}{2\nu+12+\ddots}}}}}}$$
Convergence type $F^2$ with $E=-4/z^2$, $P=2\nu+1$, and
$C=-(\nu+1)(z/2)^{2\nu}/I_{\nu}(z)^2$, so that
$$\dfrac{I_{\nu+2}(z)}{I_{\nu}(z)}-\dfrac{p(n)}{q(n)}\sim(-1)^{n+1}\dfrac{(\nu+1)/I_{\nu}(z)^2}{n!^2(4/z^2)^{n-\nu}n^{2\nu+1}}\;.$$
$$A=1+(z^2/2-(\nu+1)^2)/n+\cdots$$
\end{cf}

\smallskip

\begin{cf}\label{6.22.3}{\ }
\begin{verbatim}
[(a,b,c)->sqrt(c)*besseli(b/a+1,2*sqrt(c)/a)/
                  besseli(b/a,2*sqrt(c)/a),[0,a*n+b],[c]]
\end{verbatim}
$$\sqrt{c}\dfrac{I_{b/a+1}(2\sqrt{c}/a)}{I_{b/a}(2\sqrt{c}/a)}=
\dfrac{c}{a+b+\dfrac{c}{2a+b+\dfrac{c}{3a+b+\dfrac{c}{4a+b+\dfrac{c}{5a+b+\dfrac{c}{6a+b+\ddots}}}}}}$$
Convergence type $F^2$ with $E=-a^2/c$, $P=1+2b/a$, and
$$C=c^{1+b/a}/(a^{1+2b/a}I_{b/a}(2\sqrt{c}/a)^2)\;,\quad\text{so that}$$
$$\sqrt{c}\dfrac{I_{b/a+1}(2\sqrt{c}/a)}{I_{b/a}(2\sqrt{c}/a)}-\dfrac{p(n)}{q(n)}\sim(-1)^n\dfrac{c^{1+b/a}/(a^{1+2b/a}I_{b/a}(2\sqrt{c}/a)^2)}{n!^2(a^2/c)^nn^{1+2b/a}}\;.$$
$$A=1-(((a+b)^2-2c)/a^2)/n+\cdots$$
\end{cf}

This is of course the same CF as \ref{6.22.1}, but presented differently.
See the same CF as \ref{6.23.4}.

\smallskip

\begin{cf}\label{6.22.4}{\ }
\begin{verbatim}
[(nu,z)->gamma(nu+1)*(2/z)^nu*besseli(nu,z),[1,4*nu+4,4*n*(n+nu)+z^2],
                                       [z^2,-4*z^2*n*(n+nu)]]
\end{verbatim}
$$\G(\nu+1)(2/z)^{\nu}I_{\nu}(z)=1+\dfrac{z^2}{4\nu+4-\dfrac{(4\nu+4)z^2}{z^2+8\nu+16-\dfrac{(8\nu+16)z^2}{z^2+12\nu+36-\ddots}}}$$
Convergence type $F^2$ with $E=4/z^2$, $P=\nu+2$, and $C=\G(\nu+1)z^2/4$,
so that
$$\G(\nu+1)(2/z)^{\nu}I_{\nu}(z)-\dfrac{p(n)}{q(n)}\sim\dfrac{\G(\nu+1)}{n!^2(2/z)^{2n+2}n^{\nu+2}}\;.$$
$$A=1-((\nu^2+3\nu+4)/2)/n+((6z^2+3\nu^4+22\nu^3+69\nu^2+98\nu+72)/24)/n^2+\cdots$$
Series:
$$\G(\nu+1)(2/z)^{\nu}I_{\nu}(z)=1+\dfrac{z^2}{4\nu+4}\sum_{n\ge0}\dfrac{(z/2)^{2n}}{(n+1)!(\nu+2)_n}$$
\end{cf}

This is simply the CF coming from the power series expansion of $I_{\nu}(z)$.

\smallskip

\begin{cf}\label{6.22.5}{\ }
\begin{verbatim}
[(nu,z)->gamma(nu+1)*(2/z)^nu*exp(-z)*besseli(nu,z),
[1,2*nu+1,n^2+(nu-z)*2*n+z*(1-2*nu)],
[-(2*nu+1)*z,z*n*(n+2*nu)*(2*n+2*nu+1)]]
\end{verbatim}
$$\G(\nu+1)(2/z)^{\nu}e^{-z}I_{\nu}(z)=1-\dfrac{(2\nu+1)z}{2\nu+1+\dfrac{(4\nu^2+8\nu+3)z}{(-2\nu-3)z+4\nu+4+\dfrac{(8\nu^2+28\nu+20)z}{(-2\nu-5)z+6\nu+9+\ddots}}}$$
Convergence type $F^1$ with $E=-1/(2z)$, $P=\nu+3/2$, and $C=...$, so that
$$\G(\nu+1)(2/z)^{\nu}e^{-z}I_{\nu}(z)-\dfrac{p(n)}{q(n)}\sim(-1)^n\dfrac{C}{n!(1/(2z))^nn^{\nu+3/2}}$$
$$A=1-(2z+3\nu^2/2+2\nu+13/8)/n+\cdots$$
Series:
$$\G(\nu+1)(2/z)^{\nu}e^{-z}I_{\nu}(z)=1-z\sum_{n\ge0}\dfrac{(\nu+3/2)_n}{(n+1)!(2\nu+2)_n}(-2z)^n$$
\end{cf}

\smallskip

\begin{cf}\label{6.22.6}{\ }
\begin{verbatim}
[(nu,z)->besseli(nu+1,z)/besseli(nu,z),[0,2*nu+2+z,2*nu+n+1+2*z],
                                       [z,-(2*nu+2*n+1)*z]]
\end{verbatim}
$$\dfrac{I_{\nu+1}(z)}{I_{\nu}(z)}=\dfrac{z}{z+2\nu+2-\dfrac{(2\nu+3)z}{2z+2\nu+3-\dfrac{(2\nu+5)z}{2z+2\nu+4-\dfrac{(2\nu+7)z}{2z+2\nu+5-\ddots}}}}$$
Convergence type $F^1$ with $E=1/(2z)$, $P=3\nu+5/2$, and
$C=2\G(2\nu+2)z^{2\nu+1}/(\G(\nu+1)I_{\nu}(z)^2e^{2z}\sqrt{\pi})$, so that
$$\dfrac{I_{\nu+1}(z)}{I_{\nu}(z)}-\dfrac{p(n)}{q(n)}\sim\dfrac{2\G(2\nu+2)z^{2\nu+1}/(\G(\nu+1)I_{\nu}(z)^2e^{2z}\sqrt{\pi})}{n!(1/(2z))^nn^{3\nu+5/2}}\;.$$
$$A=1+(4z(\nu+1)-(28\nu^2+56\nu+29)/8)/n+\cdots$$
\end{cf}

Note that there does not seem to exist a CF of this type for $J_{\nu+1}(z)/J_{\nu}(z)$ not involving complex numbers.

\smallskip

Since $I_{\nu}(z)=\dfrac{2^{-3\nu}z^{\nu}}{\G(\nu+1)}{}_0F_1(;\nu+1;z^2/4)$,
more CFs can be obtained from those for ${}_0F_1$.

\smallskip

\begin{cf}\label{6.22.7}{\ }
\begin{verbatim}
[(nu,z,a)->(z/2)^(1-nu)*gamma(nu)*(besseli(nu-1,z)+a*z*besseli(nu,z)),
[1,2*nu,2*n*(n+nu-1)*(2*a*(n-1)+1)+(2*a*n+1)*z^2/2],
z^2*[a+1/2,-n*(n+nu-1)*(2*a*(n-1)+1)*(2*a*(n+1)+1)]]
\end{verbatim}
$$(z/2)^{1-\nu}\G(\nu)(I_{\nu-1}(z)+azI_{\nu}(z))=1+\dfrac{(a+1/2)z^2}{2\nu-\dfrac{(4\nu a+\nu)z^2}{(2a+1/2)z^2+(8\nu+8)a+4\nu+4+\ddots}}$$
Convergence type $F^2$ with $E=4/z^2$, $P=\nu$, and $C=...$, so that
$$(z/2)^{1-\nu}\G(\nu)(I_{\nu-1}(z)+azI_{\nu}(z))-\dfrac{p(n)}{q(n)}\sim\dfrac{C}{n!^2(2/z)^{2n}n^{\nu}}\;.$$
$$A=1+(((-\nu^2-\nu)a+1)/(2a))/n+\cdots$$
\end{cf}

Special case:

\smallskip

\begin{cf}\label{6.22.8}{\ }
\begin{verbatim}
[(z)->besseli(0,z)-z*besseli(1,z),
[1,-4,4*n^2*(2*n-3)+z^2*(2*n-1)],z^2*[1,-4*n^2*(2*n-3)*(2*n+1)]]
\end{verbatim}
$$I_0(z)-zI_1(z)=1+\dfrac{z^2}{-4+\dfrac{12z^2}{3z^2+16-\dfrac{80z^2}{5z^2+108-\dfrac{756z^2}{7z^2+320-\dfrac{2880z^2}{9z^2+700-\ddots}}}}}$$
Convergence type $F^2$ with $E=4/z^2$, $P=1$, and $C=-z^2/2$, so that
$$I_0(z)-zI_1(z)-\dfrac{p(n)}{q(n)}\sim-\dfrac{2}{n!^2(2/z)^{2n+2}n}\;.$$
$$A=1-(3/2)/n+(z^2/4+2)/n^2-(9z^2/8+5/2)/n^3+\cdots$$
Series:
$$I_0(z)-zI_1(z)=1-\dfrac{z^2}{4}\sum_{n\ge0}\dfrac{(2n+1)}{(n+1)!^2}(z/2)^{2n}$$
\end{cf}

\smallskip

\begin{cf}\label{6.22.9}{\ }
\begin{verbatim}
[(nu,z)->(z/2)^(-2*nu)*gamma(nu+1)^2*besseli(nu,z)^2,
[1,2*(nu+1)*(2*nu+1),2*n*(n+nu)*(n+2*nu)+z^2*(2*n+2*nu-1)],
z^2*[2*nu+1,-2*n*(n+nu)*(n+2*nu)*(2*n+2*nu+1)]]
\end{verbatim}
$$(z/2)^{-2\nu}\G(\nu+1)^2I_{\nu}(z)^2=1+\dfrac{(2\nu+1)z^2}{4\nu^2+6\nu+2-\dfrac{(8\nu^3+24\nu^2+22\nu+6)z^2}{(2\nu+3)z^2+8\nu^2+24\nu+16+\ddots}}$$
Convergence type $F^2$ with $E=1/z^2$, $P=2\nu+5/2$, and $C=...$, so that
$$(z/2)^{-2\nu}\G(\nu+1)^2I_{\nu}(z)^2-\dfrac{p(n)}{q(n)}\sim\dfrac{C}{n!^2(1/z)^{2n}n^{2\nu+5/2}}\;.$$
$$A=1-(2\nu^2+7\nu/2+21/8)/n+\cdots$$
Series:
$$(z/2)^{-2\nu}\G(\nu+1)^2I_{\nu}(z)^2=1+\dfrac{z^2}{2(\nu+1)}\sum_{n\ge0}\dfrac{(\nu+3/2)_n}{(n+1)!(\nu+2)_n(2\nu+2)_n}z^{2n}$$
\end{cf}

Special case:

\smallskip

\begin{cf}\label{6.22.10}{\ }
\begin{verbatim}
[(z)->besseli(0,z)^2,[1,2,2*n^3+z^2*(2*n-1)],z^2*[1,-2*n^3*(2*n+1)]]
\end{verbatim}
$$I_0(z)^2=1+\dfrac{z^2}{2-\dfrac{6z^2}{3z^2+16-\dfrac{80z^2}{5z^2+54-\dfrac{378z^2}{7z^2+128-\dfrac{1152z^2}{9z^2+250-\dfrac{2750z^2}{11z^2+432-\ddots}}}}}}$$
Convergence type $F^2$ with $E=1/z^2$, $P=5/2$, and $C=z^2/\sqrt{\pi}$, so that
$$I_0(z)^2-\dfrac{p(n)}{q(n)}\sim\dfrac{1/\sqrt{\pi}}{n!^2(1/z)^{2n+2}n^{5/2}}\;.$$
$$A=1-(21/8)/n+(z^2+617/128)/n^2-(57z^2/8+7759/1024)/n^3+\cdots$$
Series:
$$I_0(z)^2=1+\dfrac{z^2}{2}\sum_{n\ge0}\dfrac{(3/2)_n}{(n+1)!^3}z^{2n}$$
\end{cf}

\smallskip

\begin{cf}\label{6.22.11}{\ }
\begin{verbatim}
[(nu,z)->(z/2)^(-2*nu)*gamma(nu+1)^2*(besseli(nu,z)^2-besseli(nu+1,z)^2),
[1,4*(nu+1)^2,2*n*(n+nu)*(n+2*nu+1)+z^2*(2*n+2*nu-1)],
z^2*[2*nu+1,-2*n*(n+nu)*(n+2*nu+1)*(2*n+2*nu+1)]]
\end{verbatim}
$$(z/2)^{-2\nu}\G(\nu+1)^2(I_{\nu}(z)^2-I_{\nu+1}(z)^2)=1+\dfrac{(2\nu+1)z^2}{4\nu^2+8\nu+4-\dfrac{(8\nu^3+28\nu^2+32\nu+12)z^2}{(2\nu+3)z^2+8\nu^2+28\nu+24+\ddots}}$$
Convergence type $F^2$ with $E=1/z^2$, $P=2\nu+7/2$, and $C=...$, so that
$$(z/2)^{-2\nu}\G(\nu+1)^2(I_{\nu}(z)^2-I_{\nu+1}(z)^2)-\dfrac{p(n)}{q(n)}\sim\dfrac{C}{n!^2(1/z)^{2n}n^{2\nu+7/2}}\;.$$
$$A=1-(2\nu^2+11\nu/2+37/8)/n+\cdots$$
Series:
$$(z/2)^{-2\nu}\G(\nu+1)^2(I_{\nu}(z)^2-I_{\nu+1}(z)^2)=1+\dfrac{2\nu+1}{4(\nu+1)^2}z^2\sum_{n\ge0}\dfrac{(\nu+3/2)_n}{(n+1)!(\nu+2)_n(2\nu+3)_n}z^{2n}$$
\end{cf}

Special case:

\smallskip

\begin{cf}\label{6.22.12}{\ }
\begin{verbatim}
[(z)->besseli(0,z)^2-besseli(1,z)^2,
[1,4,2*n^2*(n+1)+z^2*(2*n-1)],z^2*[1,-2*n^2*(n+1)*(2*n+1)]]
\end{verbatim}
$$I_0(z)^2-I_1(z)^2=1+\dfrac{z^2}{4-\dfrac{12z^2}{3z^2+24-\dfrac{120z^2}{5z^2+72-\dfrac{504z^2}{7z^2+160-\dfrac{1440z^2}{9z^2+300-\ddots}}}}}$$
Convergence type $F^2$ with $E=1/z^2$, $P=7/2$, and $C=z^2/\sqrt{\pi}$, so that
$$I_0(z)^2-I_1(z)^2-\dfrac{p(n)}{q(n)}\sim\dfrac{1/\sqrt{\pi}}{n!^2(1/z)^{2n+2}n^{7/2}}\;.$$
$$A=1-(37/8)/n+(z^2+1801/128)/n^2-(81z^2/8+36575/1024)/n^3+\cdots$$
Series:
$$I_0(z)^2-I_1(z)^2=1+\dfrac{z^2}{2}\sum_{n\ge0}\dfrac{(3/2)_n}{(n+2)(n+1)!^3}z^{2n}$$
\end{cf}

\smallskip

\begin{cf}\label{6.22.13}{\ }
\begin{verbatim}
[(nu,z)->(z/2)^(-2*nu)*gamma(nu+1)^2*(besseli(nu,z)^2+besseli(nu+1,z)^2),
[1,4*(nu+1)^2,2*n*(n+nu)*(n+2*nu+1)+z^2*(2*n+2*nu+1)],
z^2*[2*nu+3,-2*n*(n+nu)*(n+2*nu+1)*(2*n+2*nu+3)]]
\end{verbatim}
$$(z/2)^{-2\nu}\G(\nu+1)^2(I_{\nu}(z)^2+I_{\nu+1}(z)^2)=1+\dfrac{(2\nu+3)z^2}{4\nu^2+8\nu+4-\dfrac{(8\nu^3+36\nu^2+48\nu+20)z^2}{(2\nu+5)z^2+8\nu^2+28\nu+24+\ddots}}$$
Convergence type $F^2$ with $E=1/z^2$, $P=2\nu+5/2$, and $C=...$, so that
$$(z/2)^{-2\nu}\G(\nu+1)^2(I_{\nu}(z)^2+I_{\nu+1}(z)^2)-\dfrac{p(n)}{q(n)}\sim\dfrac{C}{n!^2(1/z)^{2n}n^{2\nu+5/2}}\;.$$
$$A=1-(2\nu^2+9\nu/2+25/8)/n+\cdots$$
Series:
$$(z/2)^{-2\nu}\G(\nu+1)^2(I_{\nu}(z)^2+I_{\nu+1}(z)^2)=1+\dfrac{2\nu+3}{4(\nu+1)^2}z^2\sum_{n\ge0}\dfrac{(\nu+5/2)_n}{(n+1)!(\nu+2)_n(2\nu+3)_n}z^{2n}$$
\end{cf}

Special case:

\smallskip

\begin{cf}\label{6.22.14}{\ }
\begin{verbatim}
[(z)->besseli(0,z)^2+besseli(1,z)^2,
[1,4,2*n^2*(n+1)+z^2*(2*n+1)],z^2*[3,-2*n^2*(n+1)*(2*n+3)]]
\end{verbatim}
$$I_0(z)^2+I_1(z)^2=1+\dfrac{3z^2}{4-\dfrac{20z^2}{5z^2+24-\dfrac{168z^2}{7z^2+72-\dfrac{648z^2}{9z^2+160-\dfrac{1760z^2}{11z^2+300-\ddots}}}}}$$
Convergence type $F^2$ with $E=1/z^2$, $P=5/2$, and $C=2z^2/\sqrt{\pi}$,
so that
$$I_0(z)^2+I_1(z)^2-\dfrac{p(n)}{q(n)}\sim\dfrac{2/\sqrt{\pi}}{n!^2(1/z)^{2n+2}n^{5/2}}\;.$$
$$A=1-(25/8)/n+(z^2+913/128)/n^2-(61z^2/8+14963/1024)/n^3+\cdots$$
Series:
$$I_0(z)^2+I_1(z)^2=1+\dfrac{3z^2}{2}\sum_{n\ge0}\dfrac{(5/2)_n}{(n+2)(n+1)!^3}z^{2n}$$
\end{cf}

\smallskip

\begin{cf}\label{6.22.15}{\ }
\begin{verbatim}
[(z)->besseli(0,z)*besseli(1,z),[z/2,8,2*n*(n+1)^2+z^2*(2*n+1)],
                                [3*z^3/2,-2*z^2*n*(n+1)^2*(2*n+3)]]
\end{verbatim}
$$I_0(z)I_1(z)=\dfrac{z}{2}+\dfrac{3z^3/2}{8-\dfrac{40z^2}{5z^2+36-\dfrac{252z^2}{7z^2+96-\dfrac{864z^2}{9z^2+200-\dfrac{2200z^2}{11z^2+360-\ddots}}}}}$$
Convergence type $F^2$ with $E=1/z^2$, $P=7/2$, and $C=z^3/\sqrt{\pi}$, so that
$$I_0(z)I_1(z)-\dfrac{p(n)}{q(n)}\sim\dfrac{1/\sqrt{\pi}}{n!^2(1/z)^{2n+3}n^{7/2}}\;.$$
$$A=1-(41/8)/n+(z^2+2225/128)/n^2-((85/8)z^2+50563/1024)/n^3+\cdots$$
Series:
$$I_0(z)I_1(z)=\dfrac{z}{2}+\dfrac{3z^3}{4}\sum_{n\ge0}\dfrac{(5/2)_n}{(n+1)!(n+2)!^2}z^{2n}$$
\end{cf}

\smallskip

There exist very similar CFs for many other linear or quadratic expressions
involving $I_{\nu}(z)$, for instance also for
$\nu I_{\nu}(z)^2+(\nu+1)I_{\nu+1}(z)^2$, $I_{\nu}(z)I_{\nu+1}(z)$,
$I_{\nu}(z)^2-I_{\nu-1}(z)I_{\nu+1}(z)$, etc...

\smallskip

\begin{cf}\label{6.22.16}{\ }
\begin{verbatim}
[(z)->besseli(0,z)*besselj(0,z),[1,32,32*n^3*(2*n-1)-z^4],
                                [-z^4,32*n^3*(2*n-1)*z^4]]
\end{verbatim}
$$I_0(z)J_0(z)=1-\dfrac{z^4}{32+\dfrac{32z^4}{-z^4+768+\dfrac{768z^4}{-z^4+4320+\dfrac{4320z^4}{-z^4+14336+\dfrac{14336z^4}{-z^4+36000+\ddots}}}}}$$
Convergence type $F^4$ with $E=-64/z^4$, $P=7/2$, and $C=-\sqrt{\pi}z^4/64$,
so that
$$I_0(z)J_0(z)-\dfrac{p(n)}{q(n)}\sim(-1)^{n+1}\dfrac{\sqrt{\pi}}{n!^4(z^2/8)^{2n+2}}$$
$$A=1-(27/8)/n+(937/128)/n^2-(13249/1024)/n^3+\cdots$$
Series:
$$I_0(z)J_0(z)=\sum_{n\ge0}\dfrac{(-z^4/64)^n}{(1/2)_nn!^3}$$
\end{cf}

\smallskip

This is the CF corresponding to the Taylor expansion of $I_0(z)J_0(z)$.
There exist similar CFs for general products $I_{\nu}(z)J_{\mu}(z)$.

\smallskip

\begin{cf}\label{6.22.17}{\ }
\begin{verbatim}
[(nu,z)->besselk(nu+1,z)/besselk(nu,z),
[0,z,2*z],[z*[1,-(2*nu+1)],z*[2*nu+1+2*n,-(2*nu+1-2*n)]]]
\end{verbatim}
$$\dfrac{K_{\nu+1}(z)}{K_{\nu}(z)}=\dfrac{z}{z-\dfrac{(2\nu+1)z}{2z+\dfrac{(2\nu+3)z}{2z-\dfrac{(2\nu-1)z}{2z+\dfrac{(2\nu+5)z}{2z-\dfrac{(2\nu-3)z}{2z+\ddots}}}}}}$$
Convergence type $D^-$ with $D=16z$ and $C=-2\pi\cos(\pi\nu)/(zK_{\nu}(z)^2)$,
so that
$$\dfrac{K_{\nu+1}(z)}{K_{\nu}(z)}-\dfrac{p(n)}{q(n)}\sim(-1)^{n+1}\dfrac{2\pi\cos(\pi\nu)/(zK_{\nu}(z)^2)}{e^{4\sqrt{zn}}}\;.$$
$$A=1-((2z^2/3-(2\nu^2+4\nu+15/8))/z^{1/2})/n^{1/2}+\cdots$$
\end{cf}

\smallskip

\begin{cf}\label{6.22.18}{\ }
\begin{verbatim}
[(nu,z)->besselk(nu+1,z)/besselk(nu,z),[0,2*z-2*nu-1,4*(n+z-1)],
                                       [2*z,4*(nu+1)^2-(2*n-1)^2]]
\end{verbatim}
$$\dfrac{K_{\nu+1}(z)}{K_{\nu}(z)}=\dfrac{2z}{2z-2\nu-1+\dfrac{4\nu^2+8\nu+3}{4z+4+\dfrac{4\nu^2+8\nu-5}{4z+8+\dfrac{4\nu^2+8\nu-21}{4z+12+\dfrac{4\nu^2+8\nu-45}{4z+16+\ddots}}}}}$$
Convergence type $D^+$ with $D=32z$ and $C=-2\pi\cos(\pi\nu)/(zK_{\nu}(z)^2)$, so that
$$\dfrac{K_{\nu+1}(z)}{K_{\nu}(z)}-\dfrac{p(n)}{q(n)}\sim-\dfrac{2\pi\cos(\pi\nu)/(zK_{\nu}(z)^2)}{e^{4\sqrt{2zn}}}\;.$$
$$A=1-((16z^2-(48\nu^2+96\nu+45))d/(48z^{1/2}))/n^{1/2}+\cdots$$
with $d=\sqrt{2}$.
Parametric family for $k\ge0$:
\begin{verbatim}
[(nu,z)->besselk(nu+1,z)/besselk(nu,z),
4*(n+z)+2*k-4,-(2*n+2*nu+1)*(2*n-2*nu+2*k-3)]
\end{verbatim}
Convergence type $D^+$ with $D=32z$.
\end{cf}

Note that this is simply the contracted continued fraction of the previous one,
and is only given for comparison with the next one.

\smallskip

\begin{cf}\label{6.22.19}{\ }
\begin{verbatim}
[(nu,z)->z*besselk(nu+1,z)/besselk(nu,z),
[z+nu+1/2,4*(n+z)],[2*nu^2-1/2,4*nu^2-(2*n+1)^2]]
\end{verbatim}
$$z\dfrac{K_{\nu+1}(z)}{K_{\nu}(z)}=z+\nu+1/2+\dfrac{2\nu^2-1/2}{4z+4+\dfrac{4\nu^2-9}{4z+8+\dfrac{4\nu^2-25}{4z+12+\dfrac{4\nu^2-49}{4z+16+\dfrac{4\nu^2-81}{4z+20+\ddots}}}}}$$
Convergence type $D^+$ with $D=32z$ and $C=-2\pi\cos(\pi\nu)/K_{\nu}(z)^2$, so that
$$z\dfrac{K_{\nu+1}(z)}{K_{\nu}(z)}-\dfrac{p(n)}{q(n)}\sim-\dfrac{2\pi\cos(\pi\nu)/K_{\nu}(z)^2}{e^{4\sqrt{2zn}}}\;.$$
$$A=1-((16z^2+96z-48\nu^2+3)d/(48z^{1/2}))/n^{1/2}+\cdots$$
with $d=\sqrt{2}$.
\end{cf}

\smallskip

Thanks to the identity
$$K_{\nu}(z)=(2z/\pi)^{-1/2}e^{-z}{}_2F_0(1/2+\nu,1/2-\nu;;-1/(2z))\;,$$
we can obtain more CFs for $K_{\nu}$ from those for ${}_2F_0$.

\smallskip

Note that
$$\dfrac{1}{i}\dfrac{H_{\nu+1}^{(2)}(z/i)}{H_{\nu}^{(2)}(z/i)}=i\dfrac{H_{\nu+1}^{(1)}(iz)}{H_{\nu}^{(1)}(iz)}\;,$$
so CFs for the left-hand side are trivially obtained from those for the
right-hand side.

\smallskip

\begin{cf}\label{6.22.20}{\ }
\begin{verbatim}
[(nu,z)->I*besselh1(nu+1,I*z)/besselh1(nu,I*z),
[0,2*z],[[2*z,-(4*nu+2)*z],z*[2*n+2*nu+1,2*n-2*nu-1]]]
\end{verbatim}
$$i\dfrac{H_{\nu+1}^{(1)}(iz)}{H_{\nu}^{(1)}(iz)}=\dfrac{2z}{2z-\dfrac{(4\nu+2)z}{2z+\dfrac{(2\nu+3)z}{2z-\dfrac{(2\nu-1)z}{2z+\dfrac{(2\nu+5)z}{2z-\dfrac{(2\nu-3)z}{2z+\ddots}}}}}}$$
Convergence type $D^-$ with $D=16z$ and
$C=8\cos(\pi\nu)e^{-\pi i\nu}/(\pi zH_{\nu}^{(1)}(iz)^2)$, so that
$$i\dfrac{H_{\nu+1}^{(1)}(iz)}{H_{\nu}^{(1)}(iz)}-\dfrac{p(n)}{q(n)}\sim(-1)^n\dfrac{8\cos(\pi\nu)e^{-\pi i\nu}/(\pi zH_{\nu}^{(1)}(iz)^2)}{e^{4\sqrt{zn}}}\;.$$
$$A=1+((-2z^2/3+2\nu^2+4\nu+15/8)/z^{1/2})/n^{1/2}+\cdots$$
\end{cf}

\smallskip

\begin{cf}\label{6.22.21}{\ }
\begin{verbatim}
[(nu,z)->z*besselh1(nu+1,I*z)/besselh1(nu,I*z),
[z+nu+1/2,4*n+4*z)],[2*nu^2-1/2,4*nu^2-(2*n+1)^2]]
\end{verbatim}
$$z\dfrac{H_{\nu+1}^{(1)}(iz)}{H_{\nu}^{(1)}(iz)}=z+\nu+1/2+\dfrac{2\nu^2-1/2}{4z+4+\dfrac{4\nu^2-9}{4z+8+\dfrac{4\nu^2-25}{4z+12+\dfrac{4\nu^2-49}{4z+16+\ddots}}}}$$
Convergence type $D^+$ with $D=32z$ and $C=8\cos(\pi\nu)e^{-\pi i\nu}/(\pi H_{\nu}^{(1)}(iz)^2)$, so that
$$z\dfrac{H_{\nu+1}^{(1)}(iz)}{H_{\nu}^{(1)}(iz)}-\dfrac{p(n)}{q(n)}\sim\dfrac{8\cos(\pi\nu)e^{-\pi i\nu}/(\pi H_{\nu}^{(1)}(iz)^2)}{e^{4\sqrt{2nz}}}\;.$$
$$A=1+((-16z^2-96z+48\nu^2-3)/(24(2z)^{1/2}))/n^{1/2}+\cdots$$
\end{cf}

\smallskip

{\tt Pari/GP}'s Airy function {\tt airy(z)} returns the 2-component vector
$[\Ai(z),\Bi(z)]$, so we define:

\begin{verbatim}
airyA(z)=airy(z)[1];
\end{verbatim}

\smallskip

\begin{cf}\label{6.22.22}{\ }
\begin{verbatim}
[(z)->-airyA'(z^2)/(z*airyA(z^2)),[0,2*z^2,4*z^2],
                              [[2*z^2,-2*z],2*z*[6*n+1,6*n-1]]]
[(z)->-airyA'(z^2)/(z*airyA(z^2)),[[0,z],[4*z^2,2*z]],[[z,-1],[6*n+1,6*n-1]]]
\end{verbatim}
\begin{align*}-\dfrac{1}{z}\dfrac{\Ai'}{\Ai}(z^2)&=\dfrac{2z^2}{2z^2-\dfrac{2z}{4z^2+\dfrac{14z}{4z^2+\dfrac{10z}{4z^2+\dfrac{26z}{4z^2+\dfrac{22z}{4z^2+\ddots}}}}}}\\
  &=\dfrac{z}{z-\dfrac{1}{4z^2+\dfrac{7}{2z+\dfrac{5}{4z^2+\dfrac{13}{2z+\dfrac{11}{4z^2+\ddots}}}}}}\end{align*}
Convergence type $D^-$ with $D=(32/3)z^3$ and $C=...$, so that
$$-\dfrac{1}{z}\dfrac{\Ai'}{\Ai}(z^2)-\dfrac{p(n)}{q(n)}\sim(-1)^n\dfrac{C}{e^{4z\sqrt{2z/3}}}\;.$$
$A=1+\cdots$
\end{cf}

\smallskip

\begin{cf}\label{6.22.23}{\ }
\begin{verbatim}
[(z)->-airyA'(z^2)/(z*airyA(z^2)),[0,4*z^3-1,12*n+8*z^3-12],
                              [4*z^3,-(6*n-7)*(6*n+1)]]
\end{verbatim}
$$-\dfrac{1}{z}\dfrac{\Ai'}{\Ai}(z^2)=\dfrac{4z^3}{4z^3-1+\dfrac{7}{8z^3+12-\dfrac{65}{8z^3+24-\dfrac{209}{8z^3+36-\dfrac{425}{8z^3+48-\dfrac{713}{8z^3+60-\ddots}}}}}}$$
Convergence type $D^+$ with $D=64z^3/3$ and $C=...$, so that
$$-\dfrac{1}{z}\dfrac{\Ai'}{\Ai}(z^2)-\dfrac{p(n)}{q(n)}\sim\dfrac{C}{e^{8z\sqrt{z/3}}}\;.$$
$A=1+\cdots$
\end{cf}

\smallskip

\begin{cf}\label{6.22.24}{\ }
\begin{verbatim}
[(z)->-airyA'(z^2)/(z*airyA(z^2)),[1,8*z^3+5,12*n+8*z^3-6],[2,-(36*n^2-1)]]
\end{verbatim}
$$-\dfrac{1}{z}\dfrac{\Ai'}{\Ai}(z^2)=1+\dfrac{2}{8z^3+5-\dfrac{35}{8z^3+18-\dfrac{143}{8z^3+30-\dfrac{323}{8z^3+42-\dfrac{575}{8z^3+54-\dfrac{899}{8z^3+66-\ddots}}}}}}$$
Convergence type $D^+$ with $D=64z^3/3$ and $C=...$, so that
$$-\dfrac{1}{z}\dfrac{\Ai'}{\Ai}(z^2)-\dfrac{p(n)}{q(n)}\sim\dfrac{C}{e^{8z\sqrt{z/3}}}\;.$$
$A=1+\cdots$

\noindent
Parametric family for all $k$:
\begin{verbatim}
[(z)->-airyA'(z^2)/(z*airyA(z^2)),12*n+8*z^3-6*k,-(6*n+1)*(6*n+5-6*k)]
\end{verbatim}
Convergence type $D^+$ with $D=64z^3/3$.
\end{cf}

\medskip

\section{General Hypergeometric Functions}

\medskip

\subsection{Preliminary Remark}

\smallskip

Thanks to Euler's transformation of hyper\-geometric-type series into
continued fractions (Proposition \ref{prop:euler} and Corollary
\ref{cor:euler}), any hypergeometric series can trivially be converted
into a termwise equal CF. Since there are myriads of specific evaluations
of hypergeometric functions (see for instance \cite{Beu-Coh2}), there exist
correspondingly myriad corresponding
CFs, and it would be as hopeless a task to give a reasonably complete list
as it would be to list all evaluations of hypergeometric functions. We have
therefore not even attempted to give even part of such a list. Note that
most of the evaluations are linear combinations of quotients of gamma
function values.

However, I cannot resist giving at least one example of such a (reasonably
complicated) CF. It is possible to show the following identity, see
\cite{Beu-Coh1}:
$${}_2F_1(2z,2z+1/3;z+5/6;-1/8)=(16/27)^z\dfrac{\G(z+5/6)/\G(5/6)}{\G(z+2/3)/\G(2/3)}\;.$$
Using Euler, this gives the following CF, which we have not included in
our dictionary:

\smallskip

\begin{cf}\label{6.23.1}{\ }
\begin{verbatim}
[(z)->(16/27)^z*gamma(z+5/6)*gamma(2/3)/(gamma(z+2/3)*gamma(5/6)),
[1,24*z+20,21*n^2+(12*z+1)*n-2*(2*z-1)*(3*z-1)],
[-2*z*(6*z+1),4*n*(2*z+n)*(6*z+3*n+1)*(6*z+6*n-1)]]
\end{verbatim}
$$(16/27)^z\dfrac{\G(z+5/6)/\G(5/6)}{\G(z+2/3)/\G(2/3)}=1-\dfrac{12z^2+2z}{24z+20+\dfrac{288z^3+576z^2+376z+80}{-12z^2+34z+84+\ddots}}$$
Convergence type $E$ with $E=-8$, $P=-3z+3/2$, and $C=...$, so that
$$(16/27)^z\dfrac{\G(z+5/6)/\G(5/6)}{\G(z+2/3)/\G(2/3)}-\dfrac{p(n)}{q(n)}\sim(-1)^n\dfrac{C}{2^{3n}n^{-3z+3/2}}\;.$$
$$A=1+(7z^2/2+z-11/8)/n+(49z^4/8+z^3-917z^2/144-7z/12+1585/1152)/n^2+\cdots$$
\end{cf}

\medskip

Note the two elementary cases ${}_0F_0(;;z)=e^z$ and
${}_1F_0(a;;z)=(1-z)^{-a}$, so we can trivially obtain CFs for these functions
from those given above for the exponential and power functions.

\medskip

\subsection{Function ${}_0F_1$}

Recall that $${}_0F_1(;b;z)=\sum_{n\ge0}\dfrac{z^n}{(b)_nn!}\;,$$
with infinite radius of convergence.

\medskip

\begin{verbatim}
F01(b,z)=hypergeom([],[b],z);
\end{verbatim}

\smallskip

\begin{cf}\label{6.23.2}{\ }
\begin{verbatim}
[(b,z)->F01(b,z),[1,b,n^2+(b-1)*n+z],[z,-n*(n+b-1)*z]]
\end{verbatim}
$${}_0F_1(;b;z)=1+\dfrac{z}{b-\dfrac{bz}{z+(2b+2)-\dfrac{(2b+2)z}{z+(3b+6)-\dfrac{(3b+6)z}{z+(4b+12)-\dfrac{(4b+12)z}{z+(5b+20)-\ddots}}}}}$$
Convergence type $F^2$ with $E=1/z$, $P=b+1$, and $C=\G(b)z$, so that
$${}_0F_1(;b;z)-\dfrac{p(n)}{q(n)}\sim\dfrac{\G(b)}{n!^2(1/z)^{n+1}n^{b+1}}\;.$$
$$A=1-((b^2+b+2)/2)/n+(z+(3b^4+10b^3+21b^2+14b+24)/24)/n^2+\cdots$$
Series:
$${}_0F_1(;b;z)=1+\dfrac{z}{b}\sum_{n\ge0}\dfrac{z^n}{(n+1)!(b+1)_n}$$
\end{cf}

This is the term-by-term continued fraction corresponding to the
Taylor expansion of ${}_0F_1(;b;z)$.

\smallskip

\begin{cf}\label{6.23.3}{\ }
\begin{verbatim}
[(b,z)->b*F01(b,z)/F01(b+1,z),[n+b],[z]]
\end{verbatim}
$$b\dfrac{{}_0F_1(;b;z)}{{}_0F_1(;b+1;z)}=b+\dfrac{z}{b+1+\dfrac{z}{b+2+\dfrac{z}{b+3+\dfrac{z}{b+4+\dfrac{z}{b+5+\ddots}}}}}$$
Convergence type $F^2$ with $E=-1/z$, $P=2b+1$, and
$C=z\G(b+1)^2/{}_0F_1(;b+1;z)^2$, so that
$$b\dfrac{{}_0F_1(;b;z)}{{}_0F_1(;b+1;z)}-\dfrac{p(n)}{q(n)}\sim(-1)^n\dfrac{\G(b+1)^2/{}_0F_1(;b+1;z)^2}{n!^2(1/z)^{n+1}n^{2b+1}}\;.$$
$$A=1+(2z-(b+1)^2)/n+(12z^2-12(b+1)(b+2)z+(b+1)(b+2)(3b^2+5b+3)/6)/n^2+\cdots$$
\end{cf}

\smallskip

\begin{cf}\label{6.23.4}{\ }
\begin{verbatim}
[(a,b,c)->b*F01(b/a,c/a^2)/F01(b/a+1,c/a^2),[a*n+b],[c]]
\end{verbatim}
$$b\dfrac{{}_0F_1(;b/a;c/a^2)}{{}_0F_1(;b/a+1;c/a^2)}=b+\dfrac{c}{a+b+\dfrac{c}{2a+b+\dfrac{c}{3a+b+\dfrac{c}{4a+b+\dfrac{c}{5a+b+\ddots}}}}}$$
Convergence type $F^2$ with $E=-a^2/c$, $P=2b/a+1$, and
$C=(c/a)\G(1+b/a)^2/{}_0F_1(;b/a+1,c/a^2)^2$, so that
$$b\dfrac{{}_0F_1(;b/a;c/a^2)}{{}_0F_1(;b/a+1;c/a^2)}-\dfrac{p(n)}{q(n)}\sim(-1)^n\dfrac{(c/a)\G(1+b/a)^2/{}_0F_1(;b/a+1,c/a^2)^2}{n!^2(a^2/c)^nn^{2b/a+1}}\;.$$
$$A=1-(((a+b)^2-2c)/a^2)/n+\cdots$$
\end{cf}

This is of course the same CF as the previous one, but presented differently.
See the same CF as \ref{6.22.3}.

\smallskip

\begin{cf}\label{6.23.5}{\ }
\begin{verbatim}
[(b,z)->2*b*F01(b,z^2)/F01(b+1,z^2),[2*(b+z),2*b+4*z+n],
                                    [-2*(2*b+2*n+1)*z]]
\end{verbatim}
$$2b\dfrac{{}_0F_1(;b;z^2)}{{}_0F_1(;b+1;z^2)}=2z+2b-\dfrac{(4b+2)z}{4z+2b+1-\dfrac{(4b+6)z}{4z+2b+2-\dfrac{(4b+10)z}{4z+2b+3-\ddots}}}$$
Convergence type $F^1$ with $E=1/(4z)$, $P=3b+2$, and
$C=-z\G(2b+1)\G(b+1)2^{2b+2}/(\sqrt{\pi}e^{4z}{}_0F_1(;b;z^2)^2)$, so that
$$2b\dfrac{{}_0F_1(;b;z^2)}{{}_0F_1(;b+1;z^2)}-\dfrac{p(n)}{q(n)}\sim-\dfrac{z\G(2b+1)\G(b+1)2^{2b+2}/(\sqrt{\pi}e^{4z}{}_0F_1(;b;z^2)^2)}{n!(1/(4z))^nn^{3b+2}}\;.$$
$$A=1+(8bz-(28b^2+24b+5)/8)/n+\cdots$$
\end{cf}

\smallskip

Since ${}_0F_1(;b;z^2)=2^{2b-2}\G(b)z^{1-b}I_{b-1}(2z)$, more CFs
can be obtained from those for $I_{\nu}$.

\medskip

\subsection{Function ${}_1F_1$}

Recall that $${}_1F_1(a;b;z)=\sum_{n\ge0}\dfrac{(a)_n}{(b)_n}\dfrac{z^n}{n!}\;,$$
with infinite radius of convergence.

\medskip

\begin{verbatim}
F11(a,b,z)=hypergeom([a],[b],z);
\end{verbatim}

\smallskip

\begin{cf}\label{6.23.6}{\ }
\begin{verbatim}
[(a,b,z)->F11(a,b,z),[1,b,n^2+(z+b-1)*n+(a-1)*z],[a*z,-n*(n+a)*(n+b-1)*z]]
\end{verbatim}
$${}_1F_1(a;b;z)=1+\dfrac{az}{b-\dfrac{(ba+b)z}{(a+1)z+(2b+2)-\dfrac{((2b+2)a+(4b+4))z}{(a+2)z+(3b+6)-\ddots}}}$$
Convergence type $F^1$ with $E=1/z$, $P=b-a+1$, and $C=z\G(b)/\G(a)$, so that
$${}_1F_1(a;b;z)-\dfrac{p(n)}{q(n)}\sim\dfrac{z\G(b)/\G(a)}{n!(1/z)^nn^{b-a+1}}\;.$$
$$A=1+(z+a(a+1)/2-b(b+1)/2-1)/n+\cdots$$
Series:
$${}_1F_1(a;b;z)=1+\dfrac{az}{b}\sum_{n\ge0}\dfrac{(a+1)_n}{(n+1)!(b+1)_n}z^n$$
\end{cf}

This is the term-by-term continued fraction corresponding to the
Taylor expansion of ${}_1F_1(;b;z)$.

\smallskip

\begin{cf}\label{6.23.7}{\ }
\begin{verbatim}
[(a,b,z)->b*F11(a,b,z)/F11(a+1,b+1,z),
[n+b],[[(a-b)*z,(a+1)*z],z*[a-b-n,n+a+1]]]
\end{verbatim}
$$b\dfrac{{}_1F_1(a;b;z)}{{}_1F_1(a+1;b+1;z)}=b+\dfrac{(a-b)z}{b+1+\dfrac{(a+1)z}{b+2+\dfrac{(a-b-1)z}{b+3+\dfrac{(a+2)z}{b+4+\dfrac{(a-b-2)z}{b+5+\ddots}}}}}$$
Convergence type $F^1$ with $E=-2/z$, $P=b+1/2$, and
$C=-\sqrt{\pi}ze^z\G(b+1)^2/({}_1F_1(a+1;b+1;z)^2\G(b-a)\G(a+1)2^{b+1/2})$, so that
\begin{align*}&b\dfrac{{}_1F_1(a;b;z)}{{}_1F_1(a+1;b+1;z)}-\dfrac{p(n)}{q(n)}\sim\\
  &\phantom{=}(-1)^{n+1}\dfrac{\sqrt{\pi}ze^z\G(b+1)^2/({}_1F_1(a+1;b+1;z)^2\G(b-a)\G(a+1)2^{b+1/2})}{n!(2/z)^nn^{b+1/2}}\;.\end{align*}
$$A=1+(z^2/4+(2a-b)z+2a^2-2ab-b-1/2)/n+\cdots$$
\end{cf}

\smallskip

\begin{cf}\label{6.23.8}{\ }
\begin{verbatim}
[(a,b,z)->b*F11(a,b,z)/F11(a+1,b+1,z),[n+b-z],[(n+a+1)*z]]
\end{verbatim}
$$b\dfrac{{}_1F_1(a;b;z)}{{}_1F_1(a+1;b+1;z)}=-z+b+\dfrac{(a+1)z}{-z+b+1+\dfrac{(a+2)z}{-z+b+2+\dfrac{(a+3)z}{-z+b+3+\ddots}}}$$
Convergence type $F^1$ with $E=-1/z$, $P=2b-a$, and
$C=\G(b+1)^2ze^{2z}/(\G(a+1){}_1F_1(a+1;b+1;z)^2)$, so that
$$b\dfrac{{}_1F_1(a;b;z)}{{}_1F_1(a+1;b+1;z)}-\dfrac{p(n)}{q(n)}\sim(-1)^n\dfrac{\G(b+1)^2ze^{2z}/(\G(a+1){}_1F_1(a+1;b+1;z)^2)}{n!(1/z)^nn^{2b-a}}\;.$$
$$A=1+((2a-2b+1)z+(a^2+3a-2b^2-4b)/2)/n+\cdots$$
\end{cf}

Another way of writing this CF is as follows:

\smallskip

\begin{cf}\label{6.23.9}{\ }
\begin{verbatim}
[(a,b,c,d)->(b+c/a)*F11(d/c-1,b/a+c/a^2,c/a^2)/
                    F11(d/c,b/a+c/a^2+1,c/a^2),[a*n+b],[c*n+d]]
\end{verbatim}
\begin{align*}(b+c/a)&\dfrac{{}_1F_1(d/c-1;b/a+c/a^2;c/a^2)}{{}_1F_1(d/c;b/a+c/a^2+1;c/a^2)}\\
  &=b+\dfrac{d}{a+b+\dfrac{c+d}{2a+b+\dfrac{2c+d}{3a+b+\dfrac{3c+d}{4a+b+\dfrac{4c+d}{5a+b+\dfrac{5c+d}{6a+b+\ddots}}}}}}\end{align*}
Convergence type $F^1$ with $E=-a^2/c$, $P=((c-d)a^2+2c^2+2abc)/(ca^2)$, and
$C=\G(b/a+c/a^2+1)^2(c/a)e^{2c/a^2}/(\G(d/c){}_1F_1(d/c;b/a+c/a^2+1;c/a^2)^2)$,
so that
\begin{align*}&(b+c/a)\dfrac{{}_1F_1(d/c-1;b/a+c/a^2;c/a^2)}{{}_1F_1(d/c;b/a+c/a^2+1;c/a^2)}-\dfrac{p(n)}{q(n)}\sim\\
  &\phantom{=}(-1)^n\dfrac{\G(b/a+c/a^2+1)^2(c/a)e^{2c/a^2}/(\G(d/c){}_1F_1(d/c;b/a+c/a^2+1;c/a^2)^2)}{n!(a^2/c)^nn^P}\;.\end{align*}
$$A=1+(((d^2+cd-2c^2)a^4-4bc^2a^3+(4c^2d-2c^2b^2-6c^3)a^2-8c^3ba-6c^4)/(2c^2a^4))/n+\cdots$$
\end{cf}

\smallskip

\begin{cf}\label{6.23.10}{\ }
\begin{verbatim}
[(b,z)->z*F11(1,b+1,z),[0,n+b-1],[b*z*[1,-1],z*[n,-(n+b)]]]
\end{verbatim}
$$z{}_1F_1(1,b+1,z)=\dfrac{bz}{b-\dfrac{bz}{b+1+\dfrac{z}{b+2-\dfrac{(b+1)z}{b+3+\dfrac{2z}{b+4-\dfrac{(b+2)z}{b+5+\ddots}}}}}}$$
Convergence type $F^1$ with $E=2i/z$, $P=b-1/2$, and
$C=\sqrt{\pi}ze^z\G(b+1)/2^{b-1/2}$, so that
$$z{}_1F_1(1,b+1,z)-\dfrac{p(n)}{q(n)}\sim(-1)^{\lfloor n/2\rfloor}\dfrac{\sqrt{\pi}e^z\G(b+1)/2^{b-3/2}}{n!(2/z)^{n+1}n^{b-1/2}}\;.$$
$$A=1+(z^2/4+(1-b)z-b+1/2)/n+\cdots$$
\end{cf}

\smallskip

\begin{cf}\label{6.23.11}{\ }
\begin{verbatim}
[(b,z)->F11(1,b+1,z),[0,b-1+n-z],[b,n*z]]
\end{verbatim}
$${}_1F_1(1;b+1;z)=\dfrac{b}{b-z+\dfrac{z}{b-z+1+\dfrac{2z}{b-z+2+\dfrac{3z}{b-z+3+\dfrac{4z}{b-z+4+\ddots}}}}}$$
Convergence type $F^1$ with $E=-1/z$, $P=2b-1$, and $C=b\G(b)^2e^{2z}$, so that
$${}_1F_1(1;b+1;z)-\dfrac{p(n)}{q(n)}\sim(-1)^n\dfrac{b\G(b)^2e^{2z}}{n!(1/z)^nn^{2b-1}}\;.$$
$$A=1-((2b-1)z+b^2)/n+(2(b^2-1)z^2+(2b^3+b^2-b+1)z+(3b^4+2b^3+b)/6)/n^2+\cdots$$
\end{cf}

\smallskip

\begin{cf}\label{6.23.12}{\ }
\begin{verbatim}
[(b,z)->F11(1,b+1,z),[1,b+1,n+b+z],[z,-z*(n+b)]]
\end{verbatim}
$${}_1F_1(1;b+1;z)=1+\dfrac{z}{b+1-\dfrac{(b+1)z}{z+(b+2)-\dfrac{(b+2)z}{z+(b+3)-\dfrac{(b+3)z}{z+(b+4)-\dfrac{(b+4)z}{z+(b+5)-\ddots}}}}}$$
Convergence type $F^1$ with $E=1/z$, $P=b+1$, and $C=z\G(b+1)$, so that
$${}_1F_1(1;b+1;z)-\dfrac{p(n)}{q(n)}\sim\dfrac{z\G(b+1)}{n!(1/z)^nn^{b+1}}\;.$$
$$A=1+(z-(b+1)(b+2)/2)/n+\cdots$$
Series:
$${}_1F_1(1;b+1;z)=1+\dfrac{z}{b+1}\sum_{n\ge0}\dfrac{z^n}{(b+2)_n}$$
\end{cf}

\smallskip

\begin{cf}\label{6.23.13}{\ }
\begin{verbatim}
[(z)->F11(1,z+1,z),[1+z,n+2*z+1],[-z*(n+z)]]
\end{verbatim}
$${}_1F_1(1;z+1;z)=z+1-\dfrac{z^2}{2z+2-\dfrac{z^2+z}{2z+3-\dfrac{z^2+2z}{2z+4-\dfrac{z^2+3z}{2z+5-\dfrac{z^2+4z}{2z+6-\ddots}}}}}$$
Convergence type $F^1$ with $E=1/z$, $P=z+3$, and $C=-z^3\G(z)$, so that
$${}_1F_1(1;z+1;z)-\dfrac{p(n)}{q(n)}\sim-\dfrac{\G(z)}{n!(1/z)^{n+3}n^{z+3}}\;.$$
$$A=1-((z^2+z+8)/2)/n+((3z^4+10z^3+57z^2-22z+264)/24)/n^2+\cdots$$
Series:
$${}_1F_1(1;z+1;z)=z+\dfrac{1-z^2}{z+1}\sum_{n\ge0}\dfrac{z^n}{(n+1)(n+2)(z+2)_n}$$
\end{cf}

\medskip

\subsection{Function ${}_2F_0$}

Recall that $${}_2F_0(a,b;;z)=\sum_{n\ge0}(a)_n(b)_n\dfrac{z^n}{n!}\;,$$
with radius of convergence $0$.

\medskip

\begin{verbatim}
F20(a,b,z)=hypergeom([a,b],[],z);
\end{verbatim}

\smallskip

\begin{cf}\label{6.23.14}{\ }
\begin{verbatim}
[(a,b,z)->F20(a,b,z)/F20(a,b+1,z),[1],[-z*[n+a,n+b+1]]]
\end{verbatim}
$$\dfrac{{}_2F_0(a,b;;z)}{{}_2F_0(a,b+1;;z)}=1-\dfrac{az}{1-\dfrac{(b+1)z}{1-\dfrac{(a+1)z}{1-\dfrac{(b+2)z}{1-\dfrac{(a+2)z}{1-\dfrac{(b+3)z}{1-\ddots}}}}}}$$
Convergence type $D^i$ with $D=-8/z$, and
$C=2\pi e^{-1/z}/({}_2F_0(a,b+1;;z)^2(-z)^{a+b}\G(a)\G(b+1))$, so that
$$\dfrac{{}_2F_0(a,b;;z)}{{}_2F_0(a,b+1;;z)}-\dfrac{p(n)}{q(n)}\sim(-1)^{\lfloor n/2\rfloor}\dfrac{2\pi e^{-1/z}/({}_2F_0(a,b+1;;z)^2(-z)^{a+b}\G(a)\G(b+1))}{e^{\sqrt{-8n/z}}}\;.$$
$$A=1+(((2(a-b)^2-1/2)z^2+4(a+b+1)z-2/3)/(-2z)^{3/2})/n^{1/2}+\cdots$$
\end{cf}
Special case:

\smallskip

\begin{cf}\label{6.23.15}{\ }
\begin{verbatim}
[(a,z)->2*F20(a,a-1/2,z)/F20(a,a+1/2,z),[2],[-2*z*(n+2*a)]]
\end{verbatim}
$$2\dfrac{{}_2F_0(a,a-1/2;;z)}{{}_2F_0(a,a+1/2;;z)}=2-\dfrac{4az}{2-\dfrac{(4a+2)z}{2-\dfrac{(4a+4)z}{2-\dfrac{(4a+6)z}{2-\dfrac{(4a+8)z}{2-\dfrac{(4a+10)z}{2-\ddots}}}}}}$$
Convergence type $D^-$ with $D=-8/z$, and
$C=4\pi e^{-1/z}/({}_2F_0(a,a+1/2;;z)^2(-z)^{2a-1/2}\G(a)\G(a+1/2))$, so that
$$2\dfrac{{}_2F_0(a,a-1/2;;z)}{{}_2F_0(a,a+1/2;;z)}-\dfrac{p(n)}{q(n)}\sim(-1)^n\dfrac{4\pi e^{-1/z}/({}_2F_0(a,a+1/2;;z)^2(-z)^{2a-1/2}\G(a)\G(a+1/2))}{e^{\sqrt{-8n/z}}}\;.$$
$$A=1+(6(4a+1)z-2)d/(12(-z)^{3/2})/n^{1/2}-((6(4a+1)z-2)^2/(144z^3))/n+\cdots$$
with $d=\sqrt{2}$.
\end{cf}

\smallskip

A different way of writing the above CF is as follows:

\smallskip

\begin{cf}\label{6.23.16}{\ }
\begin{verbatim}
[(a,b,c)->a*F20(c/(2*b),(c-b)/(2*b),-2*b/a^2)/
            F20(c/(2*b),(c+b)/(2*b),-2*b/a^2),[a],[b*n+c]]
\end{verbatim}
$$a\dfrac{{}_2F_0(c/(2b),(c-b)/(2b);;-2b/a^2)}{{}_2F_0(c/(2b),(c+b)/(2b);;-2b/a^2)}=a+\dfrac{c}{a+\dfrac{b+c}{a+\dfrac{2b+c}{a+\dfrac{3b+c}{a+\dfrac{4b+c}{a+\dfrac{5b+c}{a+\ddots}}}}}}$$
Convergence type $D^-$ with $D=4a^2/b$ and $C=2\pi ae^{a^2/(2b)}/({}_2F_0(c/(2b),(c+b)/(2b);;-2b/a^2)^2(2b/a^2)^{c/b-1/2}\G(c/(2b))\G((c+b)/(2b)))$, so that
\begin{align*}&a\dfrac{{}_2F_0(c/(2b),(c-b)/(2b);;-2b/a^2)}{{}_2F_0(c/(2b),(c+b)/(2b);;-2b/a^2)}-\dfrac{p(n)}{q(n)}\sim\\
  &\phantom{=}(-1)^n\dfrac{2\pi ae^{a^2/(2b)}/({}_2F_0(c/(2b),(c+b)/(2b);;-2b/a^2)^2(2b/a^2)^{c/b-1/2}\G(c/(2b))\G((c+b)/(2b)))}{e^{2|a|\sqrt{n/b}}}\;.\end{align*}
$$A=1-(a(a^2+6b+12c)/(12b^{3/2}))/n^{1/2}+\cdots$$
\end{cf}

\smallskip

\begin{cf}\label{6.23.17}{\ }
\begin{verbatim}
[(a,b,z)->F20(a,b,z)/F20(a,b+1,z),
[1,-(b+1)*z+1,1-(a+b+2*n-1)*z],[-a*z,-(n+a)*(n+b)*z^2]]
\end{verbatim}
$$\dfrac{{}_2F_0(a,b;;z)}{{}_2F_0(a,b+1;;z)}=1-\dfrac{az}{(-b-1)z+1-\dfrac{((b+1)a+(b+1))z^2}{(-a-b-3)z+1-\dfrac{((b+2)a+(2b+4))z^2}{(-a-b-5)z+1-\ddots}}}$$
Convergence type $D^+$ with $D=-16/z$, and
$C=2\pi e^{-1/z}/({}_2F_0(a,b+1;;z)^2(-z)^{a+b}\G(a)\G(b+1))$, so that
$$\dfrac{{}_2F_0(a,b;;z)}{{}_2F_0(a,b+1;;z)}-\dfrac{p(n)}{q(n)}\sim\dfrac{2\pi e^{-1/z}/({}_2F_0(a,b+1;;z)^2(-z)^{a+b}\G(a)\G(b+1))}{e^{\sqrt{-16n/z}}}\;.$$
$$A=1-(((12(a-b)^2-3)z^2-24(a+b+1)z+4)/(24(-z)^{3/2}))/n^{1/2}+\cdots$$
\end{cf}

This is simply the contraction of \ref{6.23.14}.

\smallskip

\begin{cf}\label{6.23.18}{\ }
\begin{verbatim}
[(a,b,z)->F20(a,b,z)/F20(a+1,b+1,z),
[1-(a+b+2*n+1)*z],[-(n+a+1)*(n+b+1)*z^2]]
\end{verbatim}
$$\dfrac{{}_2F_0(a,b;;z)}{{}_2F_0(a+1,b+1;;z)}=(-a-b-1)z+1-\dfrac{((b+1)a+(b+1))z^2}{(-a-b-3)z+1-\dfrac{((b+2)a+(2b+4))z^2}{(-a-b-5)z+1-\ddots}}$$
Convergence type $D^+$ with $D=-16/z$, and
$C=-2\pi e^{-1/z}/({}_2F_0(a+1,b+1;;z)^2(-z)^{a+b+1}\G(a+1)\G(b+1))$, so that
$$\dfrac{{}_2F_0(a,b;;z)}{{}_2F_0(a+1,b+1;;z)}-\dfrac{p(n)}{q(n)}\sim-\dfrac{2\pi e^{-1/z}/({}_2F_0(a+1,b+1;;z)^2(-z)^{a+b+1}\G(a+1)\G(b+1))}{e^{\sqrt{-16n/z}}}\;.$$
$$A=1-(((12(a-b)^2-3)z^2-24(a+b+3)z+4)/(24(-z)^{3/2}))/n^{1/2}+\cdots$$
\end{cf}

\smallskip

\begin{cf}\label{6.23.19}{\ }
\begin{verbatim}
[(a,b,z)->F20(a,b,z)/F20(a+1,b-1,z),
[1+(a+1-b)*z,1-(a+b+2*n-1)*z],-z^2*[(a+1)*(b-a-1),(n+a+1)*(n+b-1)]]
\end{verbatim}
$$\dfrac{{}_2F_0(a,b;;z)}{{}_2F_0(a+1,b-1;;z)}=(a+(-b+1))z+1+\dfrac{(a^2+(-b+2)a+(-b+1))z^2}{(-a+(-b-1))z+1-\dfrac{(ba+2b)z^2}{(-a+(-b-3))z+1+\ddots}}$$
Convergence type $D^+$ with $D=-16/z$, and $C=...$, so that
$$\dfrac{{}_2F_0(a,b;;z)}{{}_2F_0(a+1,b-1;;z)}-\dfrac{p(n)}{q(n)}\sim\dfrac{C}{e^{\sqrt{-16n/z}}}\;.$$
$$A=1+(3(2a-2b+3)(2a-2b+5)z^2+24(a+b+1)z-4)/(24(-z)^{3/2}))/n^{1/2}+\cdots$$
\end{cf}

\smallskip

\begin{cf}\label{6.23.20}{\ }
\begin{verbatim}
[(b,z)->F20(1,b,z),[0,1],[[1,-b*z],-z*[n,n+b]]]
\end{verbatim}
$${}_2F_0(1,b;;z)=\dfrac{1}{1-\dfrac{bz}{1-\dfrac{z}{1-\dfrac{(b+1)z}{1-\dfrac{2z}{1-\dfrac{(b+2)z}{1-\ddots}}}}}}$$
Convergence type $D^-$ with $D=-8/z$ and $C=2\pi e^{-1/z}/((-z)^b\G(b))$,
so that
$${}_2F_0(1,b;;z)-\dfrac{p(n)}{q(n)}\sim(-1)^n\dfrac{2\pi e^{-1/z}/((-z)^b\G(b))}{e^{\sqrt{-8n/z}}}\;.$$
$$A=1+(((2b^2-4b+3/2)z^2+4bz-2/3)/(-2z)^{3/2})/n^{1/2}+\cdots$$
\end{cf}

\smallskip

\begin{cf}\label{6.23.21}{\ }
\begin{verbatim}
[(b,z)->F20(1,b,z),[0,-(2*n+b-2)*z+1],[1,-z^2*n*(n+b-1)]]
\end{verbatim}
$${}_2F_0(1,b;;z)=\dfrac{1}{(-b)z+1-\dfrac{bz^2}{(-b-2)z+1-\dfrac{(2b+2)z^2}{(-b-4)z+1-\dfrac{(3b+6)z^2}{(-b-6)z+1-\ddots}}}}$$
Convergence type $D^+$ with $D=-16/z$ and $C=2\pi e^{-1/z}/((-z)^b\G(b))$, so that
$${}_2F_0(1,b;;z)-\dfrac{p(n)}{q(n)}\sim\dfrac{2\pi e^{-1/z}/((-z)^b\G(b))}{e^{\sqrt{-16n/z}}}\;.$$
$$A=1+((3(2b-1)(2b-3)z^2+24bz-4)/(24(-z)^{3/2}))/n^{1/2}+\cdots$$
\end{cf}

This is the contraction of the previous CF.

\smallskip

\smallskip

\begin{cf}\label{6.23.22}{\ }
\begin{verbatim}
[(b,z)->F20(1,b,z),[z*b+1,-(2*n+b)*z+1],z^2*[b*(b+1),-n*(n+b+1)]]
\end{verbatim}
$${}_2F_0(1,b;;z)=bz+1+\dfrac{(b^2+b)z^2}{(-b-2)z+1-\dfrac{(b+2)z^2}{(-b-4)z+1-\dfrac{(2b+6)z^2}{(-b-6)z+1-\ddots}}}$$
Convergence type $D^+$ with $D=-16/z$ and $C=2\pi e^{-1/z}/((-z)^b\G(b))$, so that
$${}_2F_0(1,b;;z)-\dfrac{p(n)}{q(n)}\sim\dfrac{2\pi e^{-1/z}/((-z)^b\G(b))}{e^{\sqrt{-16n/z}}}\;.$$
$$A=1-((3(2b+1)(2b+3)z^2+24(b+2)z-4)/(24(-z)^{3/2}))/n^{1/2}+\cdots$$
\end{cf}

Special cases: recall that
$(2k-1)!!=1\cdot3\cdots(2k-1)=(2k)!/(k!2^k)$, and that
$${}_2F_0(1,1/2;;2z)=\sum_{k\ge0}(2k-1)!!z^k$$

\smallskip

\begin{cf}\label{6.23.23}{\ }
\begin{verbatim}
[(z)->F20(1,1/2,2*z),[0,1],[1,-n*z]]
\end{verbatim}
$${}_2F_0(1,1/2;;2z)=\sum_{k\ge0}(2k-1)!!z^k=\dfrac{1}{1-\dfrac{z}{1-\dfrac{2z}{1-\dfrac{3z}{1-\dfrac{4z}{1-\dfrac{5z}{1-\ddots}}}}}}$$
Convergence type $D^-$ with $D=-4/z$ and $C=\sqrt{-2\pi/z}e^{-1/(2z)}$,
so that
$${}_2F_0(1,1/2;;2z)-\dfrac{p(n)}{q(n)}\sim(-1)^n\dfrac{\sqrt{-2\pi/z}e^{-1/(2z)}}{e^{\sqrt{-4n/z}}}\;.$$
$$A=1+((6z-1)/(12(-z)^{3/2}))/n^{1/2}-((6z-1)^2/(288z^3))/n+\cdots$$
\end{cf}

\smallskip

\begin{cf}\label{6.23.24}{\ }
\begin{verbatim}
[(z)->F20(1,1/2,2*z),[0,-z*(4*n-3)+1],[1,-z^2*2*n*(2*n-1)]]
\end{verbatim}
$${}_2F_0(1,1/2;;2z)=\dfrac{1}{-z+1-\dfrac{2z^2}{-5z+1-\dfrac{12z^2}{-9z+1-\dfrac{30z^2}{-13z+1-\dfrac{56z^2}{-17z+1-\ddots}}}}}$$
Convergence type $D^+$ with $D=-8/z$ and $C=\sqrt{-2\pi/z}e^{-1/(2z)}$, so that
$${}_2F_0(1,1/2;;2z)-\dfrac{p(n)}{q(n)}\sim\dfrac{\sqrt{-2\pi/z}e^{-1/(2z)}}{e^{\sqrt{-8n/z}}}\;.$$
$$A=1+((6z-1)d/(24(-z)^{3/2}))/n^{1/2}-((6z-1)^2/(576z^3))/n+\cdots$$
with $d=\sqrt{2}$.
\end{cf}

This is the contraction of the previous CF.

\smallskip

\begin{cf}\label{6.23.25}{\ }
\begin{verbatim}
[(z)->F20(1,1/2,2*z),[z+1,-z*(4*n+1)+1],z^2*[3,-2*n*(2*n+3)]]
\end{verbatim}
$${}_2F_0(1,1/2;;2z)=z+1+\dfrac{3z^2}{-5z+1-\dfrac{10z^2}{-9z+1-\dfrac{28z^2}{-13z+1-\dfrac{54z^2}{-17z+1-\ddots}}}}$$
Convergence type $D^+$ with $D=-8/z$ and $C=\sqrt{-2\pi/z}e^{-1/(2z)}$, so that
$${}_2F_0(1,1/2;;2z)-\dfrac{p(n)}{q(n)}\sim\dfrac{\sqrt{-2\pi/z}e^{-1/(2z)}}{e^{\sqrt{-8n/z}}}\;.$$
$$A=1-((24z^2+30z-1)d/(24(-z)^{3/2}))/n^{1/2}-((24z^2+30z-1)^2/(576z^3))/n+\cdots$$
with $d=\sqrt{2}$.
\end{cf}

\smallskip

Recall that
$${}_2F_0(1,1;;z)=\sum_{k\ge0}k!z^k$$

\smallskip

\begin{cf}\label{6.23.26}{\ }
\begin{verbatim}
[(z)->F20(1,1,z),[0,-(2*n-1)*z+1],[1,-z^2*n^2]]
\end{verbatim}
$${}_2F_0(1,1;;z)=\sum_{k\ge0}k!z^k=\dfrac{1}{-z+1-\dfrac{z^2}{-3z+1-\dfrac{4z^2}{-5z+1-\dfrac{9z^2}{-7z+1-\dfrac{16z^2}{-9z+1-\ddots}}}}}$$
Convergence type $D^+$ with $D=-16/z$ and $C=2\pi e^{-1/z}/(-z)$, so that
$${}_2F_0(1,1;;z)-\dfrac{p(n)}{q(n)}\sim\dfrac{2\pi e^{-1/z}/(-z)}{e^{\sqrt{-16n/z}}}\;.$$
$$A=1-((3z^2-24z+4)/(24(-z)^{3/2}))/n^{1/2}-((3z^2-24z+4)^2/(1152z^3))/n+\cdots$$
\end{cf}

Note that \ref{5.10.26} is the special case $z=-1$ of this CF.

\smallskip

\begin{cf}\label{6.23.27}{\ }
\begin{verbatim}
[(z)->F20(1,1,z),[z+1,-(2*n+1)*z+1],z^2*[2,-n*(n+2)]]
\end{verbatim}
$${}_2F_0(1,1;;z)=z+1+\dfrac{2z^2}{-3z+1-\dfrac{3z^2}{-5z+1-\dfrac{8z^2}{-7z+1-\dfrac{15z^2}{-9z+1-\dfrac{24z^2}{-11z+1-\ddots}}}}}$$
Convergence type $D^+$ with $D=-16/z$ and $C=2\pi e^{-1/z}/(-z)$, so that
$${}_2F_0(1,1;;z)-\dfrac{p(n)}{q(n)}\sim\dfrac{2\pi e^{-1/z}/(-z)}{e^{\sqrt{-16n/z}}}\;.$$
$$A=1-((45z^2+72z-4)/(24(-z)^{3/2}))/n^{1/2}+((45z^2+72z-4)^2/(1152z^3))/n+\cdots$$
\end{cf}

\smallskip

Note that the following trivial identity for nonnegative integers $b$ gives other
CFs for ${}_2F_0(1,1;;z)$:
$${}_2F_0(1,1;z)=\sum_{0\le k\le b-1}k!z^k+b!z^b{}_2F_0(1,b+1,z)\;.$$

\smallskip

Thanks to the identities
\begin{align*}
  {}_2F_0(1,\nu;;-1/z)&=ze^zE_{\nu}(z)\;,\\
  {}_2F_0(1,1/2-\nu;;-1/z)&=z^{1/2-\nu}e^z\ga(\nu+1/2,z)\;,\\
  {}_2F_0(1,1/2;;-1/z^2)&=\sqrt{\pi}ze^{z^2}\erfc(z)\;,\\
  {}_2F_0(1/2+\nu,1/2-\nu;;-1/(2z))&=(2z/\pi)^{1/2}e^zK_{\nu}(z)\;,\text{\quad and}\\
  {}_2F_0(a,a+1-b;;-1/z)&=z^aU(a,b,z)\;,\end{align*}
we can obtain more CFs for ${}_2F_0$ from those for $E_{\nu}$, $\ga$,
$\erfc$, $K_{\nu}$, and $U$.

\medskip

\subsection{Function $U$}

\medskip

\begin{verbatim}
U(a,b,z)=hyperu(a,b,z);
\end{verbatim}

\smallskip

\begin{cf}\label{6.23.28}{\ }
\begin{verbatim}
[(a,b,z)->-U(a,b,z)/U(a+1,b,z),[b-2*a-2-z-2*n],[(n+a+1)*(b-a-n-2)]]
\end{verbatim}
$$-\dfrac{U(a,b,z)}{U(a+1,b,z)}=-z-2a+b-2-\dfrac{a^2+(-b+3)a+(-b+2)}{-z-2a+b-4-\dfrac{a^2+(-b+5)a+(-2b+6)}{-z-2a+b-6-\ddots}}$$
Convergence type $D^+$ with $D=16z$ and
$C=2\pi e^z/(U(a+1,b,z)^2z^{b-1}\G(a+1)\G(a-b+2))$, so that
$$-\dfrac{U(a,b,z)}{U(a+1,b,z)}-\dfrac{p(n)}{q(n)}\sim\dfrac{2\pi e^z/(U(a+1,b,z)^2z^{b-1}\G(a+1)\G(a-b+2))}{e^{4\sqrt{nz}}}\;.$$
$$A=1-((4z^2-24(b-2a-4)z-3(2b-1)(2b-3))/(24z^{1/2}))/n^{1/2}+\cdots$$
\end{cf}

\smallskip

\begin{cf}\label{6.23.29}{\ }
\begin{verbatim}
[(a,b,z)->-z*derivnum(x=z,U(a,b,x))/U(a,b,z),
[a,b-2*a-z-2*n],[a*(1+a-b),(n+a)*(b-a-n-1)]]
\end{verbatim}
$$-z\dfrac{U'(a,b,z)}{U(a,b,z)}=a+\dfrac{a^2+(-b+1)a}{-z-2a+b-2-\dfrac{a^2+(-b+3)a+(-b+2)}{-z-2a+b-4-\dfrac{a^2+(-b+5)a+(-2b+6)}{-z-2a+b-6-\ddots}}}$$
Convergence type $D^+$ with $D=16z$ and $C=-2\pi e^z/(U(a,b,z)^2z^{b-1}\G(a)\G(a-b+1))$, so that
$$-z\dfrac{U'(a,b,z)}{U(a,b,z)}-\dfrac{p(n)}{q(n)}\sim-\dfrac{2\pi e^z/(U(a,b,z)^2z^{b-1}\G(a)\G(a-b+1))}{e^{4\sqrt{nz}}}\;.$$
$$A=1-((4z^2-24(b-2a-2)z-3(2b-1)(2b-3))/(24z^{1/2}))/n^{1/2}+\cdots$$
\end{cf}

\smallskip

Note that
$$U(a,b,z)=z^{-a}{}_2F_0(a,a+1-b;;-1/z)\;,$$
so we can obtain more CFs for $U$ from those for ${}_2F_0$.

\medskip

\subsection{Function ${}_2F_1$}

Recall that $${}_2F_1(a,b;c;z)=\sum_{n\ge0}\dfrac{(a)_n(b)_n}{(c)_n}\dfrac{z^n}{n!}\;,$$
with radius of convergence generically equal to $1$.

\medskip

\begin{verbatim}
F21(a,b,c,z)=hypergeom([a,b],[c],z);
\end{verbatim}

\smallskip

\begin{cf}\label{6.23.30}{\ }
\begin{verbatim}
[(a,b,c,z)->F21(a,b,c,z),[1,c,n*(n+c-1)+z*(n+a-1)*(n+b-1)],
                         [a*b*z,-n*(n+c-1)*(n+a)*(n+b)*z]]
\end{verbatim}
$$_2F_1(a,b;c;z)=1+\dfrac{abz}{c-\dfrac{(a+1)(b+1)cz}{2(c+1)+(a+1)(b+1)z-\dfrac{2(c+1)(a+2)(b+2)z}{3(c+2)+(a+2)(b+2)z-\ddots}}}$$
Convergence type $E$ with $E=1/z$, $P=c+1-a-b$, and $C=...$, so that
$$_2F_1(a,b;c;z)-\dfrac{p(n)}{q(n)}\sim\dfrac{C}{(1/z)^nn^{c+1-a-b}}\;.$$
$$A=1+(((c^2-a^2-b^2-c+a+b)z+a^2+b^2-c^2+a+b-c-2)/(2(z-1)))/n+\cdots$$
Series:
$$_2F_1(a,b;c;z)=1+\dfrac{abz}{c}\sum_{n\ge0}\dfrac{(a+1)_n(b+1)_n}{(n+1)!(c+1)_n}z^n$$
\end{cf}

This is simply the CF corresponding term by term to the Taylor expansion
of $_2F_1(a,b;c;z)$.

\smallskip

\begin{cf}\label{6.23.31}{\ }
\begin{verbatim}
[(a,b,c,z)->F21(a,b,c,z)/F21(a,b+1,c,z),
[1,c-(b+1)*z,c-(a+b+2*n-1)*z+n-1],[-a*z,(z-z^2)*(n+a)*(n+b)]]
\end{verbatim}
$$\dfrac{{}_2F_1(a,b;c;z)}{{}_2F_1(a,b+1,c,z)}=1-\dfrac{za}{c-zb-z-\dfrac{((z^2-z)a+(z^2-z))b+(z^2-z)a+z^2-z}{c-zb-za-3z+1-\ddots}}$$
Convergence type $E$ with $E=(z-1)/z$, $P=2c-(a+b+1)$, and $C=...$, so that
$$\dfrac{{}_2F_1(a,b;c;z)}{{}_2F_1(a,b+1,c,z)}-\dfrac{p(n)}{q(n)}\sim\dfrac{C}{((z-1)/z)^nn^{2c-a-b-1}}\;.$$
$A=1+\cdots$ explicit but complicated.

\noindent
Parametric family:
\begin{verbatim}
[(a,b,c,z)->F21(a,b,c,z)/F21(a,b+1,c,z),
c-(a+b+2*n-1)*z+n-k-1,(z-z^2)*(n+a)*(n+b)]
\end{verbatim}
Convergence type $E$ with $E=(z-1)/z$, $P=2c-(a+b+2k+1)$.
\end{cf}

\smallskip

\begin{cf}\label{6.23.32}{\ }
\begin{verbatim}
[(a,b,c,z)->F21(a,b,c,z)/F21(a,b+1,c,z),
[1,(b-c+1)*z,(a+b-2*c-1+2*n)*z-n+c+1],[a*z,(z-z^2)*(n+a-c)*(n+b-c)]]
\end{verbatim}
$$\dfrac{{}_2F_1(a,b;c;z)}{{}_2F_1(a,b+1,c,z)}=1+\dfrac{az}{(-c+(b+1))z-\dfrac{N}{(-2c+(a+(b+3)))z+(c-1)+\ddots}}\;,$$
\begin{align*}\text{where\quad}N&=(c^2+(-a+(-b-2))c+((b+1)a+(b+1)))z^2\\
  &\phantom{=}+(-c^2+(a+(b+2))c+((-b-1)a+(-b-1)))z\;.\end{align*}
Convergence type $E$ with $E=(z-1)/z$, $P=-(a+b+1)$, and $C=...$, so that
$$\dfrac{{}_2F_1(a,b;c;z)}{{}_2F_1(a,b+1,c,z)}-\dfrac{p(n)}{q(n)}\sim\dfrac{C}{((z-1)/z)^nn^{-a-b-1}}\;.$$
$A=1+\cdots$ explicit but complicated.

\noindent
Parametric family:
\begin{verbatim}
[(a,b,c,z)->F21(a,b,c,z)/F21(a,b+1,c,z),
(a+b-2*c-1+2*n)*z-n+c+1-k,(z-z^2)*(n+a-c)*(n+b-c)]
\end{verbatim}
Convergence type $E$ with $E=(z-1)/z$, $P=2k-(a+b+1)$.
\end{cf}

\smallskip

\begin{cf}\label{6.23.33}{\ }
\begin{verbatim}
[(a,b,c,z)->F21(a,b,c,z)/F21(a,b+1,c,z),
[1,n+c-1+(a-b+n-1)*z],[-a*z,-(n+a)*(n+c-b-1)*z]]
\end{verbatim}
$$\dfrac{{}_2F_1(a,b;c;z)}{{}_2F_1(a,b+1,c,z)}=1-\dfrac{za}{c-zb+za-\dfrac{(za+z)c+(-za-z)b}{c-zb+za+z+1-\ddots}}$$
Convergence type $E$ with $E=1/z$, $P=c+b-a$, and $C=...$, so that
$$\dfrac{{}_2F_1(a,b;c;z)}{{}_2F_1(a,b+1,c,z)}-\dfrac{p(n)}{q(n)}\sim\dfrac{C}{(1/z)^nn^{c+b-a}}\;.$$
$A=1+\cdots$ explicit but complicated.

\noindent
Parametric family for $k\ge0$:
\begin{verbatim}
[(a,b,c,z)->F21(a,b,c,z)/F21(a,b+1,c,z),
n+c-1+(a-b+n-1)*z+k*(1-z),-(n+a)*(n+c-b-1)*z]
\end{verbatim}
Convergence type $E$ with $E=1/z$, $P=c+b-a+2k$.
\end{cf}

\smallskip

\begin{cf}\label{6.23.34}{\ }
\begin{verbatim}
[(a,b,c,z)->F21(a,b,c,z)/F21(a,b+1,c,z),
[1,n+c-1],-z*[[a,-(b+1)*(c-a)*z],[(n+a)*(n+c-b-1),(n+b+1)*(n+c-a)]]]
\end{verbatim}
$$\dfrac{{}_2F_1(a,b;c;z)}{{}_2F_1(a,b+1,c,z)}=1-\dfrac{za}{c-\dfrac{(b+1)zc+(-b-1)za}{c+1-\dfrac{(a+1)zc-zba-zb}{c+2-\dfrac{(b+2)zc+((-b-2)a+(b+2))z}{c+3-\ddots}}}}$$
Convergence type $E$ with $E=(1+\sqrt{1-z})^2/z$, $P=0$, and $C=...$, so that
$$\dfrac{{}_2F_1(a,b;c;z)}{{}_2F_1(a,b+1,c,z)}-\dfrac{p(n)}{q(n)}\sim\dfrac{C}{((1+\sqrt{1-z})^2/z)^n}\;.$$
\end{cf}

\smallskip

\begin{cf}\label{6.23.35}{\ }
\begin{verbatim}
[(a,b,c,z)->F21(a,b,c,z)/F21(a,b,c+1,z),
[1,-b*z+c+1,(z+1)*n+(a-b-1)*z+c],[a*b*z/c,-(n+a)*(n+c-b)*z]]
\end{verbatim}
$$\dfrac{{}_2F_1(a,b;c;z)}{{}_2F_1(a,b,c+1,z)}=1+\dfrac{baz/c}{-bz+(c+1)+\dfrac{((b+(-c-1))a+(b+(-c-1)))z}{(a+(-b+1))z+(c+2)+\ddots}}$$
Convergence type $E$ with $E=1/z$, $P=b+c+1-a$, and $C=...$, so that
$$\dfrac{{}_2F_1(a,b;c;z)}{{}_2F_1(a,b,c+1,z)}-\dfrac{p(n)}{q(n)}\sim\dfrac{C}{(1/z)^nn^{b+c+1-a}}\;.$$
$A=1+\cdots$ explicit but complicated.

Parametric family for $k\ge0$:
\begin{verbatim}
[(a,b,c,z)->F21(a,b,c,z)/F21(a,b,c+1,z),
(z+1)*n+(a-b-1)*z+c+k*(1-z),-(n+a)*(n+c-b)*z]
\end{verbatim}
Convergence type $E$ with $E=1/z$ and $P=b+c+1-a-2*k$.
\end{cf}

\smallskip

\begin{cf}\label{6.23.36}{\ }
\begin{verbatim}
[(a,b,c,z)->F21(a,b,c,z)/F21(a,b,c+1,z),
[1,-a*z+c+1,n+c],[[a*b*z/c,(b+1)*(a-c-1)*z],
-z*[(n+a)*(n+c-b),(n+b+1)*(n+c-a+1)]]]
\end{verbatim}
$$\dfrac{{}_2F_1(a,b;c;z)}{{}_2F_1(a,b,c+1,z)}=1+\dfrac{1/cabz}{-az+(c+1)+\dfrac{((a+(-c-1))b+(a+(-c-1)))z}{c+2+\dfrac{((a+1)b+((-c-1)a+(-c-1)))z}{c+3+\ddots}}}$$
Convergence type $E$ with $E=-(1+\sqrt{1-z})^2/z$, $P=0$, and $C=...$, so that
$$\dfrac{{}_2F_1(a,b;c;z)}{{}_2F_1(a,b,c+1,z)}-\dfrac{p(n)}{q(n)}\sim(-1)^n\dfrac{C}{((1+\sqrt{1-z})^2/z)^n}\;.$$
\end{cf}

\smallskip

\begin{cf}\label{6.23.37}{\ }
\begin{verbatim}
[(a,b,c,z)->F21(a,b,c,z)/F21(a,b+1,c+1,z),
[1,n+c],[z*[a/c*(b-c),-(b+1)*(c-a+1)],
        -z*[(n+a)*(n+c-b),(n+b+1)*(n+1+c-a)]]]
\end{verbatim}
$$\dfrac{{}_2F_1(a,b;c;z)}{{}_2F_1(a,b+1,c+1,z)}=1-\dfrac{(azc-baz)/c}{c+1-\dfrac{(b+1)zc+((-b-1)a+(b+1))z}{c+2-\dfrac{(a+1)zc+((-b+1)a+(-b+1))z}{c+3-\dfrac{(b+2)zc+((-b-2)a+(2b+4))z}{c+4-\ddots}}}}$$
Convergence type $E$ with $E=(1+\sqrt{1-z})^2/z$, $P=0$, and $C=...$, so that
$$\dfrac{{}_2F_1(a,b;c;z)}{{}_2F_1(a,b+1,c+1,z)}-\dfrac{p(n)}{q(n)}\sim\dfrac{C}{((1+\sqrt{1-z})^2/z)^n}\;.$$
\end{cf}

\smallskip

\begin{verbatim}
C1(n)=-(c-a+n)*(n+b);
D1(n)=b-a+n+1;
E1(n)=n+c;
\end{verbatim}

\smallskip

\begin{cf}\label{6.23.38}{\ }
\begin{verbatim}
[(a,b,c,z)->c*F21(a,b,c,z)/F21(a,b+1,c+1,z),[E1(n)+D1(n)*z],[z*C1(n+1)]]
\end{verbatim}
$$c\dfrac{{}_2F_1(a,b;c;z)}{{}_2F_1(a,b+1;c+1;z)}=zb+(-a+1)z+c+\dfrac{(a-c-1)zb+(a-c-1)z}{zb+(-a+2)z+c+1+\ddots}$$
Convergence type $E$ with $E=1/z$, $P=c+a-b-1$, and $C=...$, so that
$$c\dfrac{{}_2F_1(a,b;c;z)}{{}_2F_1(a,b+1;c+1;z)}-\dfrac{p(n)}{q(n)}\sim\dfrac{C}{(1/z)^nn^{c+a-b-1}}\;.$$
$A=1+\cdots$ explicit but complicated.

\noindent
Parametric family for $k\ge0$:
\begin{verbatim}
[(a,b,c,z)->c*F21(a,b,c,z)/F21(a,b+1,c+1,z),E1(n)+D1(n)*z+k*(1-z),C1(n+1)]
\end{verbatim}
Convergence type $E$ with $E=1/z$ and $P=c+a-b-1+2k$.
\end{cf}

\smallskip

\begin{cf}\label{6.23.39}{\ }
\begin{verbatim}
[(a,b,c,z)->c*F21(a,b,c,z)/F21(a,b+1,c+1,z),[(b-a+1)*z+c,n+c],
z*[[(b+1)*(a-c-1),a*(b-c)],-[(n+b+1)*(n+c-a+1),(n+a)*(n+c-b)]]]
\end{verbatim}
$$c\dfrac{{}_2F_1(a,b;c;z)}{{}_2F_1(a,b+1;c+1;z)}=(-a+b+1)z+c-\dfrac{((b+1)c+((-b-1)a+(b+1)))z}{c+1-\dfrac{(ac-ba)z}{c+2-\dfrac{((b+2)c+((-b-2)a+(2b+4)))z}{c+3-\ddots}}}$$
Convergence type $E$ with $E=-(1+\sqrt{1-z})^2/z$, $P=0$, and $C=...$, so that
$$c\dfrac{{}_2F_1(a,b;c;z)}{{}_2F_1(a,b+1,c+1,z)}-\dfrac{p(n)}{q(n)}\sim(-1)^n\dfrac{C}{((1+\sqrt{1-z})^2/z)^n}\;.$$
\end{cf}

\smallskip

\begin{cf}\label{6.23.40}{\ }
\begin{verbatim}
[(a,b,c,z)->(b-a+1)*F21(b-c+1,b,b-a+1,z)/F21(b-c+1,b+1,b-a+2,z),
[E1(n)*z+D1(n)],[z*C1(n+1)]]
\end{verbatim}
$$(b-a+1)\dfrac{{}_2F_1(b-c+1,b;b-a+1;z)}{{}_2F_1(b-c+1,b+1;b-a+2;z)}=b+cz-a+1+\dfrac{(a-c-1)zb+(a-c-1)z}{b+(c+1)z-a+2+\ddots}$$
Convergence type $E$ with $E=1/z$, $P=b-a-c+1$, and $C=...$, so that
$$(b-a+1)\dfrac{{}_2F_1(b-c+1,b;b-a+1;z)}{{}_2F_1(b-c+1,b+1;b-a+2;z)}-\dfrac{p(n)}{q(n)}\sim\dfrac{C}{(1/z)^nn^{b-a-c+1}}\;.$$
$A=1+\cdots$ explicit but complicated.

\noindent
Parametric family for $k\ge0$:
\begin{verbatim}
[(a,b,c,z)->(b-a+1)*F21(b-c+1,b,b-a+1,z)/F21(b-c+1,b+1,b-a+2,z),
E1(n)*z+D1(n)+k*(1-z),z*C1(n+1)]
\end{verbatim}
Convergence type $E$ with $E=1/z$ and $P=b-a-c+1+2k$.
\end{cf}

\smallskip

\begin{cf}\label{6.23.41}{\ }
\begin{verbatim}
[(a,b,c,z)->(b-a+1)*F21(b-c+1,b,b-a+1,z)/F21(b-c+1,b+1,b-a+2,z),
[c*z+b-a+1,n+b-a+1],[z*[(b+1)*(a-c-1),(a-1)*(b+1-c)],
                    -z*[(n+b+1)*(n+c-a+1),(n+1-a)*(n+b-c+1)]]]
\end{verbatim}
\begin{align*}&(b-a+1)\dfrac{{}_2F_1(b-c+1,b;b-a+1;z)}{{}_2F_1(b-c+1,b+1;b-a+2;z)}=\\
  &\phantom{=}cz-a+b+1-\dfrac{((b+1)c+((-b-1)a+(b+1)))z}{-a+(b+2)-\dfrac{((a-1)c+((-b-1)a+(b+1)))z}{-a+(b+3)-\dfrac{((b+2)c+((-b-2)a+(2b+4)))z}{-a+(b+4)-\ddots}}}\end{align*}
Convergence type $E$ with $E=(1+\sqrt{1-z})^2/z$, $P=0$, and $C=...$, so that
$$(b-a+1)\dfrac{{}_2F_1(b-c+1,b;b-a+1;z)}{{}_2F_1(b-c+1,b+1;b-a+2;z)}-\dfrac{p(n)}{q(n)}\sim\dfrac{C}{((1+\sqrt{1-z})^2/z)^n}\;.$$
\end{cf}

\smallskip

\begin{verbatim}
C2(n)=(n+a)*(n+b);
D2(n)=-(a+b+2*n+1);
E2(n)=n+c;
\end{verbatim}

\smallskip

\begin{cf}\label{6.23.42}{\ }
\begin{verbatim}
[(a,b,c,z)->c*F21(a,b,c,z)/F21(a+1,b+1,c+1,z),
                           [E2(n)+D2(n)*z],[C2(n+1)*(z-z^2)]]
\end{verbatim}
\begin{align*}&c\dfrac{{}_2F_1(a,b;c;z)}{{}_2F_1(a+1,b+1;c+1;z)}=\\&\phantom{=}-zb+(-a-1)z+c-\dfrac{((a+1)z^2+(-a-1)z)b+((a+1)z^2+(-a-1)z)}{-zb+(-a-3)z+c+1-\ddots}\end{align*}
Convergence type $E$ with $E=(z-1)/z$, $P=2c-a-b-1$, and $C=...$, so that
$$c\dfrac{{}_2F_1(a,b;c;z)}{{}_2F_1(a+1,b+1;c+1;z)}-\dfrac{p(n)}{q(n)}\sim\dfrac{C}{((z-1)/z)^nn^{2c-a-b-1}}\;.$$
$A=1+\cdots$ explicit but complicated.
\end{cf}

\smallskip

\begin{cf}\label{6.23.43}{\ }
\begin{verbatim}
[(a,c,z)->c*F21(a,a-1/2,c,z)/F21(a,a+1/2,c+1,z),
[c,2*(n+c)],[-z*a*(2*(c-a)+1),-z*(n+2*a)*(n+2*(c-a)+1)]]
\end{verbatim}
$$c\dfrac{{}_2F_1(a,a-1/2;c;z)}{{}_2F_1(a,a+1/2;c+1;z)}=c-\dfrac{(2ac+(-2a^2+a))z}{2c+2-\dfrac{((4a+2)c+(-4a^2+2a+2))z}{2c+4-\dfrac{((4a+4)c+(-4a^2+2a+6))z}{2c+6-\dfrac{((4a+6)c+(-4a^2+2a+12))z}{2c+8-\ddots}}}}$$
Convergence type $E$ with $E=(1+\sqrt{1-z})^2$, $P=0$, and $C=...$, so that
$$c\dfrac{{}_2F_1(a,a-1/2;c;z)}{{}_2F_1(a,a+1/2;c+1;z)}-\dfrac{p(n)}{q(n)}\sim\dfrac{C}{(1+\sqrt{1-z})^{2n}}$$
\end{cf}

\smallskip

\begin{cf}\label{6.23.44}{\ }
\begin{verbatim}
[(a,z)->(1-z)*derivnum(x=z,F21(a,a+1/2,1,x))/F21(a,a+1/2,1,z),
[a,2*n],[a*(2*a-1),-z*(n+2*a)*(n-2*a+1)]]
\end{verbatim}
$$(1-z)\dfrac{{}_2F_1'(a,a+1/2;1;z)}{{}_2F_1(a,a+1/2;1;z)}=a+\dfrac{2a^2-a}{2+\dfrac{(4a^2-2a-2)z}{4+\dfrac{(4a^2-2a-6)z}{6+\dfrac{(4a^2-2a-12)z}{8+\dfrac{(4a^2-2a-20)z}{10+\dfrac{(4a^2-2a-30)z}{12+\ddots}}}}}}$$
Convergence type $E$ with $E=(1+\sqrt{1-z})^2$, $P=0$, and $C=...$, so that
$$(1-z)\dfrac{{}_2F_1'(a,a+1/2;1;z)}{{}_2F_1(a,a+1/2;1;z)}-\dfrac{p(n)}{q(n)}\sim\dfrac{C}{(1+\sqrt{1-z})^{2n}}$$
\end{cf}

\medskip

\subsection{Function ${}_2F_1$ with $b=1$}

\medskip

\begin{cf}\label{6.23.45}{\ }
\begin{verbatim}
[(a,c,z)->F21(a,1,c,z),[1,c,n+c-1+(n+a-1)*z],[a*z,-(n+c-1)*(n+a)*z]]
\end{verbatim}
$$_2F_1(a,1;c;z)=1+\dfrac{az}{c-\dfrac{(a+1)cz}{(a+1)z+c+1-\dfrac{((a+2)c+(a+2))z}{(a+2)z+c+2-\dfrac{((a+3)c+(2a+6))z}{(a+3)z+c+3-\ddots}}}}$$
Convergence type $E$ with $E=1/z$, $P=c-a$, and $C=...$, so that
$$_2F_1(a,1;c;z)-\dfrac{p(n)}{q(n)}\sim\dfrac{C}{(1/z)^nn^{c-a}}\;.$$
$$A=1+((c-a)((z-1)(c+a)-z-1)/(2z-2))/n+\cdots$$
Series:
$$_2F_1(a,1;c;z)=1+\dfrac{az}{c}\sum_{n\ge0}\dfrac{(a+1)_n}{(c+1)_n}z^n$$
Parametric family for $k\ge0$:
\begin{verbatim}
[(a,c,z)->F21(a,1,c,z),n+c-1+(n+a-1)*z+k*(1-z),-(n+c-1)*(n+a)*z]
\end{verbatim}
Convergence type $E$ with $E=1/z$ and $P=c-a+2k$.
\end{cf}

This is simply the CF corresponding term by term to the Taylor expansion
of $_2F_1(a,1;c;z)$.

\smallskip

\begin{cf}\label{6.23.46}{\ }
\begin{verbatim}
[(a,c,z)->F21(a,1,c,z),[1,c-1,n+c-2+(n+a-2)*z],[c-1,-(n+c-2)*(n+a-1)*z]]
\end{verbatim}
$$_2F_1(a,1;c;z)=\dfrac{c-1}{c-1-\dfrac{(ac-a)z}{az+c-\dfrac{(a+1)cz}{(a+1)z+(c+1)-\dfrac{((a+2)c+(a+2))z}{(a+2)z+(c+2)-\ddots}}}}$$
Convergence type $E$ with $E=1/z$, $P=c-a$, and $C=...$, so that
$$_2F_1(a,1;c;z)-\dfrac{p(n)}{q(n)}\sim\dfrac{C}{(1/z)^nn^{c-a}}\;.$$
$$A=1+((c-a)((z-1)(c+a)-3*z+1)/(2z-2))/n+\cdots$$

Parametric family for $k\ge0$:
\begin{verbatim}
[(a,c,z)->F21(a,1,c,z),n+c-2+(n+a-2)*z+k*(1-z),-(n+c-2)*(n+a-1)*z]
\end{verbatim}
Convergence type $E$ with $E=1/z$ and $P=c-a+2k$.
\end{cf}
      
\smallskip

\begin{cf}\label{6.23.47}{\ }
\begin{verbatim}
[(a,c,z)->F21(a,1,c,z),[1,n+c-1+(n-a-1)*z],[a*z,-n*(n+c-a-1)*z]]
\end{verbatim}
$$_2F_1(a,1;c;z)=1+\dfrac{az}{-az+c-\dfrac{(c-a)z}{(-a+1)z+c+1-\dfrac{(2c+(-2a+2))z}{(-a+2)z+c+2-\ddots}}}$$
Convergence type $E$ with $E=1/z$, $P=c+a$, and $C=...$, so that
$$_2F_1(a,1;c;z)-\dfrac{p(n)}{q(n)}\sim\dfrac{C}{(1/z)^nn^{c+a}}\;.$$
$$A=1+(((z-1)c^2-(2a+1)(z+1)c+(a-a^2)z+a^2+a)/(2z-2))/n+\cdots$$
Parametric family for $k\ge0$:
\begin{verbatim}
[(a,c,z)->F21(a,1,c,z),n+c-1+(n+a-1)*z+k*(1-z),-n*(n+c-a-1)*z]
\end{verbatim}
Convergence type $E$ with $E=1/z$ and $P=c+a+2k$.
\end{cf}

\smallskip

\begin{cf}\label{6.23.48}{\ }
\begin{verbatim}
[(a,c,z)->F21(a,1,c,z),[0,(n-a)*z+n+c-2],[c-1,-n*(n+c-a-1)*z]]
\end{verbatim}
$${}_2F_1(a,1;c;z)=\dfrac{c-1}{(-a+1)z+c-1+\dfrac{(a-c)z}{(-a+2)z+c+\dfrac{(2a-2c-2)z}{(-a+3)z+c+1+\ddots}}}$$
Convergence type $E$ with $E=1/z$, $P=a+c-2$, and $C=...$, so that
$${}_2F_1(a,1;c;z)-\dfrac{p(n)}{q(n)}\sim\dfrac{C}{(1/z)^nn^{a+c-2}}\;.$$
$$A=1-(((a^2+(2c-5)a-(c^2-c-2))z-a^2+(2c-1)a+(c-1)(c-2))/(2*(z-1)))/n+\cdots$$
Parametric family for $k\ge0$:
\begin{verbatim}
[(a,c,z)->F21(a,1,c,z),n+c-2+(n-a)*z+k*(1-z),-n*(n+c-a-1)*z]
\end{verbatim}
Convergence type $E$ with $E=1/z$ and $P=c+a-2+2k$.
\end{cf}

\smallskip

\begin{cf}\label{6.23.49}{\ }
\begin{verbatim}
[(a,c,z)->F21(a,1,c,z),[1,c-(a+2*n-1)*z+n-1],[a*z,n*(n+a)*(z-z^2)]]
\end{verbatim}
$${}_2F_1(a,1;c;z)=1+\dfrac{za}{c-za-z-\dfrac{(z^2-z)a+z^2-z}{c-za-3z+1-\dfrac{(2z^2-2z)a+(4z^2-4z)}{c-za-5z+2-\ddots}}}$$
Convergence type $E$ with $E=(z-1)/z$, $P=2c-a-1$, and $C=...$, so that
$${}_2F_1(a,1;c;z)-\dfrac{p(n)}{q(n)}\sim\dfrac{C}{((z-1)/z)^nn^{2c-a-1}}\;.$$
Replace $E$ by $1/E$ and $P$ by $-P$ if $z>1/2$.
$$A=1+((-2z+1)c^2+(2az+2z)c-a^2/2-a(z+1/2)-z)/n+\cdots$$
Parametric family:
\begin{verbatim}
[(a,c,z)->F21(a,1,c,z),c-(a+2*n-1)*z+n+k-1,(z-z^2)*n*(n+a)]
\end{verbatim}
Convergence type $E$ with $E=(z-1)/z$, $P=2c-a-1+2k$.
\end{cf}

\smallskip

\begin{cf}\label{6.23.50}{\ }
\begin{verbatim}
[(a,c,z)->F21(a,1,c,z),[0,n+c-2-(2*n+a-2)*z],[c-1,n*(n+a-1)*(z-z^2)]]
\end{verbatim}
$${}_2F_1(a,1;c;z)=\dfrac{c-1}{-az+c-1-\dfrac{az^2-az}{(-a-2)z+c-\dfrac{(2a+2)z^2+(-2a-2)z}{(-a-4)z+c+1-\ddots}}}$$
Convergence type $E$ with $E=(z-1)/z$, $P=2c-a-2$, and $C=...$, so that
$${}_2F_1(a,1;c;z)-\dfrac{p(n)}{q(n)}\sim\dfrac{C}{((z-1)/z)^nn^{2c-a-2}}\;.$$
$$A=1+(((2c-3)a-2(c-1)^2)z-(a^2-a-2(c-1)^2)/2)/n+\cdots$$
Parametric family:
\begin{verbatim}
[(a,c,z)->F21(a,1,c,z),n+c-2-(2*n+a-2)*z+k,n*(n+a-1)*(z-z^2)]
\end{verbatim}
Convergence type $E$ with $E=(z-1)/z$ and $P=2c-a-2+2k$.
\end{cf}

\smallskip

\begin{cf}\label{6.23.51}{\ }
\begin{verbatim}
[(a,c,z)->F21(a,1,c,z),[1,c-a*z,n+c-1],
[[a*z,(a-c)*z],-z*[(n+a)*(n+c-1),(n+1)*(n+c-a)]]]
\end{verbatim}
$$_2F_1(a,1;c;z)=1+\dfrac{az}{-az+c-\dfrac{(c-a)z}{c+1-\dfrac{(a+1)cz}{c+2-\dfrac{(2c+(-2a+2))z}{c+3-\dfrac{((a+2)c+(a+2))z}{c+4-\ddots}}}}}$$
Convergence type $E$ with $E=(1+\sqrt{1-z})^2/z$, $P=0$, and $C=...$,
so that
$$_2F_1(a,1;c;z)-\dfrac{p(n)}{q(n)}\sim\dfrac{C}{((1+\sqrt{1-z})^2/z)^n}\;.$$
\end{cf}

\smallskip

\begin{cf}\label{6.23.52}{\ }
\begin{verbatim}
[(a,c,z)->F21(a,1,c,z),
[0,n+c-2],[[c-1,-z*a*(c-1)],-z*[n*(c-a+n-1),(n+a)*(n+c-1)]]]
\end{verbatim}
$${}_2F_1(a,1;c;z)=\dfrac{c-1}{c-1-\dfrac{(ca-a)z}{c+\dfrac{(a-c)z}{c+1-\dfrac{(ca+c)z}{c+2+\dfrac{(2a-2c-2)z}{c+3-\dfrac{((c+1)a+(2c+2))z}{c+4-\ddots}}}}}}$$
Convergence type $E$ with $E=(1+\sqrt{1-z})^2/z$, $P=0$, and $C=...$, so that
$${}_2F_1(a,1;c;z)-\dfrac{p(n)}{q(n)}\sim\dfrac{C}{((1+\sqrt{1-z})^2/z)^n}\;.$$
\end{cf}

\smallskip

\begin{cf}\label{6.23.53}{\ }
\begin{verbatim}
[(c,z)->F21(1/2,1,c,z),[1,2*c-z,2*n+2*c-2],[z,-(n+1)*(n+2*c-2)*z]]
\end{verbatim}
$$_2F_1(1/2,1;c;z)=1+\dfrac{z}{-z+2c-\dfrac{(4c-2)z}{2c+2-\dfrac{6cz}{2c+4-\dfrac{(8c+4)z}{2c+6-\dfrac{(10c+10)z}{2c+8-\dfrac{(12c+18)z}{2c+10-\ddots}}}}}}$$
Convergence type $E$ with $E=(1+\sqrt{1-z})^2/z$, $P=0$, and $C=...$,
so that
$$_2F_1(1/2,1;c;z)-\dfrac{p(n)}{q(n)}\sim\dfrac{C}{((1+\sqrt{1-z})^2/z)^n}\;.$$
\end{cf}

\medskip

\subsection{Function ${}_3F_2$ with $z=1$}

Recall that $${}_3F_2(a,b,c;d,e;z)=\sum_{n\ge0}\dfrac{(a)_n(b)_n(c)_n}{(d)_n(e)_n}\dfrac{z^n}{n!}\;,$$
with radius of convergence generically equal to $1$.

\medskip

\begin{verbatim}
F32(a,b,c,d,e,z)=hypergeom([a,b,c],[d,e],z);
\end{verbatim}

\smallskip

\begin{cf}\label{6.23.54}{\ }
\begin{verbatim}
[(a,b,c,d,e)->F32(a,b,c,d,e,1)/F32(a+1,b,c,d,e,1),
[[1,d*(e-a-1)],[d+2*n-1,(d+2*n)*(e-a-1)]],
[[-b*c,(a+1)*(d-b)*(d-c)],[-(d-a+n-1)*(n+b)*(n+c),
                            (n+a+1)*(d-b+n)*(d-c+n)]]]
\end{verbatim}
$$\dfrac{{}_3F_2(a,b,c;d,e;1)}{{}_3F_2(a+1,b,c;d,e;1)}=1-\cdots$$
(too complicated to print).

Convergence type $P^-$ with $P=4e+d-2(a+b+c)-3$ and $C=...$, so that
$$\dfrac{{}_3F_2(a,b,c;d,e;1)}{{}_3F_2(a+1,b,c;d,e;1)}-\dfrac{p(n)}{q(n)}\sim(-1)^n\dfrac{C}{n^{4e+d-2(a+b+c)-3}}\;.$$
$A=1+\cdots$ explicit but complicated.
\end{cf}

\smallskip

\begin{cf}\label{6.23.55}{\ }
\begin{verbatim}
[(a,b,c,d,e)->F32(a,b,c,d,e,1)/F32(a,b,c,d+1,e,1),
[[1,(d+1)*(d+e-a-b-c)],[d+2*n,(d+2*n+1)*(d+e-a-b-c)]],
[[(n+a)*(n+b)*(n+c),-(d-a+n+1)*(d-b+n+1)*(d-c+n+1)]]]
\end{verbatim}
$$\dfrac{{}_3F_2(a,b,c;d,e;1)}{{}_3F_2(a,b,c;d+1,e;1)}=1+\cdots$$
(too complicated to print).

Convergence type $P^-$ with $P=4e+d-2(a+b+c)-2$ and $C=...$, so that
$$\dfrac{{}_3F_2(a,b,c;d,e;1)}{{}_3F_2(a,b,c;d+1,e;1)}-\dfrac{p(n)}{q(n)}\sim(-1)^n\dfrac{C}{n^{4e+d-2(a+b+c)-2}}\;.$$
$A=1+\cdots$ explicit but complicated.
\end{cf}

\smallskip

\begin{cf}\label{6.23.56}{\ }
\begin{verbatim}
[(a,b,c,d,e)->F32(a,b,c,d,e,1)/F32(a+1,b,c,d+1,e,1),
[[1,(d+1)*(e-a-1)],[d+2*n,(d+2*n+1)*(e-a-1)]],
[[-b*c,(a+1)*(d-b+1)*(d-c+1)],[-(d-a+n)*(n+b)*(n+c),
                                (n+a+1)*(d-b+n+1)*(d-c+n+1)]]]
\end{verbatim}
$$\dfrac{{}_3F_2(a,b,c;d,e;1)}{{}_3F_2(a+1,b,c;d+1,e;1)}=1+\cdots$$
(too complicated to print).

Convergence type $P^-$ with $P=4e+d-2(a+b+c)-2$ and $C=...$, so that
$$\dfrac{{}_3F_2(a,b,c;d,e;1)}{{}_3F_2(a+1,b,c;d+1,e;1)}-\dfrac{p(n)}{q(n)}\sim(-1)^n\dfrac{C}{n^{4e+d-2(a+b+c)-2}}\;.$$
$A=1+\cdots$ explicit but complicated.
\end{cf}

\smallskip

\begin{cf}\label{6.23.57}{\ }
\begin{verbatim}
[(a,b,c,d,e)->e*F32(a,b,c,d,e,1)/F32(a,b+1,c+1,d+1,e+1,1),
[[d*(e-a),d+1],[(d+2*n)*(e-a),d+2*n+1]],
[[(n+a)*(d-b+n)*(d-c+n),-(d-a+n+1)*(n+b+1)*(n+c+1)]]]
\end{verbatim}
$$\dfrac{{}_3F_2(a,b,c;d,e;1)}{{}_3F_2(a,b+1,c+1;d+1,e+1;1)}=d(e-a)+\cdots$$
(too complicated to print).

Convergence type $P^-$ with $P=4e+d-2(a+b+c)-2$ and $C=...$, so that
$$\dfrac{{}_3F_2(a,b,c;d,e;1)}{{}_3F_2(a,b+1,c+1;d+1,e+1;1)}-\dfrac{p(n)}{q(n)}\sim(-1)^n\dfrac{C}{n^{4e+d-2(a+b+c)-2}}\;.$$
$A=1+\cdots$ explicit but complicated.
\end{cf}

\chapter{Special Cases}

\section{Linear Polynomial Type Continued Fractions}

A polynomial type CF $(a(n),b(n))$ is \emph{linear} if $a(n)$ and $b(n)$
are linear (possibly constant) functions of $n$ for $n$ sufficiently large,
and we restrict to period $1$. We group here the four possible types,
depending on the ``bidegree'' $(d_a,d_b)$ of $(a(n),b(n))$ for $n$ large.
If $(d_a,d_b)=(0,0)$, this corresponds to quadratic irrationals, nothing
very interesting can be said, apart from determining which of the two
conjugates of the quadratic irrational $(a+\sqrt{a^2+4b})/2$ is the limit
of the CF {\tt [[a],[b]]}, assuming $a^2+4b\ge0$, otherwise the CF diverges (we
always assume $a$ and $b$ real).

\medskip

If $(d_a,d_b)=(1,0)$, this is essentially \ref{6.22.3} which we recall here:
\begin{cf}{\ }
\begin{verbatim}
[(a,b,c)->b+sqrt(c)*besseli(b/a+1,2*sqrt(c)/a)/
                    besseli(b/a,2*sqrt(c)/a),[a*n+b],[c]]
\end{verbatim}
$$b+\sqrt{c}\dfrac{I_{b/a+1}(2\sqrt{c}/a)}{I_{b/a}(2\sqrt{c}/a)}=
b+\dfrac{c}{a+b+\dfrac{c}{2a+b+\dfrac{c}{3a+b+\dfrac{c}{4a+b+\dfrac{c}{5a+b+\ddots}}}}}$$
\end{cf}

This is for $c>0$. If $c<0$, simply replace $I_{\nu}$ by $J_{\nu}$.

\smallskip

Note that we have
$$b+\sqrt{c}\dfrac{I_{b/a+1}(2\sqrt{c}/a)}{I_{b/a}(2\sqrt{c}/a)}=b\dfrac{{}_0F_1(;b/a;c/a^2)}{{}_0F_1(;b/a+1;c/a^2)}\;,$$
so the same CF is also given by \ref{6.23.4}.

\smallskip

Special cases:

\smallskip

\begin{verbatim}
[()->exp(1),[1,1,4*n-2],[2,1]]
[()->(exp(1)-1)/(exp(1)+1),[0,4*n-2],[1]]
[()->exp(2),[7,2*n+3],[2,1]]
[()->exp(3),[13,7,4*n+6],[54,9]]
[(k,l)->(exp(2*k/l)-1)/(exp(2*k/l)+1),[0,2*l*n-l],[k,k^2]]
[(k)->exp(2/k),[1,k-1,2*k*n-k],[2,1]]
[()->besselj(1,2)/besselj(0,2),[n],[1,-1]]
[()->besseli(1,2)/besseli(0,2),[n],[1]]
[(z)->exp(z),[1,-z+2,4*n-2],[2*z,z^2]]
[(z)->tanh(z),[0,2*n-1],[z,z^2]]
[(z)->tan(z)+1/cos(z),[1,-z+2,4*n-2],[2*z,-z^2]]
[(z)->(psi(1+z/Pi)-psi(1-z/Pi))/Pi,[0,2*n+1],[z,-z^2]]
[(nu,z)->besseli(nu+1,z)/besseli(nu,z),[0,2*n+2*nu],[z,z^2]]
[(nu,z)->besseli(nu+2,z)/besseli(nu,z),[1,2*n+2*nu],[-2*nu-2,z^2]]
[(b,z)->b*F01(b,z)/F01(b+1,z),[n+b],[z]]
[(a,b,c)->b*F01(b/a,c/a^2)/F01(1+b/a,c/a^2),[a*n+b],[c]]
\end{verbatim}

\medskip

If $(d_a,d_b)=(0,1)$, this is \ref{6.23.16} which we recall here:

\begin{cf}{\ }
\begin{verbatim}
[(a,b,c)->a*F20(c/(2*b),(c-b)/(2*b),-2*b/a^2)/
            F20(c/(2*b),(c+b)/(2*b),-2*b/a^2),[a],[b*n+c]]
\end{verbatim}
$$a\dfrac{{}_2F_0(c/(2b),(c-b)/(2b);;-2b/a^2)}{{}_2F_0(c/(2b),(c+b)/(2b);;-2b/a^2)}=a+\dfrac{c}{a+\dfrac{b+c}{a+\dfrac{2b+c}{a+\dfrac{3b+c}{a+\dfrac{4b+c}{a+\dfrac{5b+c}{a+\ddots}}}}}}$$
\end{cf}

\smallskip

Special case:

\smallskip

\begin{verbatim}
[()->besselk(3/4,1)/besselk(1/4,1),[1,4],[1,4*n+2]]
[(a,z)->2*F20(a,a-1/2,z)/F20(a,a+1/2,z),[2],[-2*z*(n+2*a)]]
\end{verbatim}

This is for $c\ne 0$ and $c\ne b$. Otherwise this is \ref{6.18.10}
(equivalently we could obtain it from \ref{6.23.23}), which we recall here:

\smallskip

\begin{cf}{\ }
\begin{verbatim}
[(a,b)->sqrt(2*b/Pi)*exp(-a^2/(2*b))/erfc(a/sqrt(2*b)),[[a],[b*(n+1)]]]
\end{verbatim}
$$\dfrac{a}{{}_2F_0(1,1/2;;-2b/a^2)}=\dfrac{\sqrt{2b/\pi}e^{-a^2/(2b)}}{\erfc(a/\sqrt{2b})}=a+\dfrac{b}{a+\dfrac{2b}{a+\dfrac{3b}{a+\dfrac{4b}{a+\dfrac{5b}{a+\dfrac{6b}{a+\ddots}}}}}}$$
\end{cf}

\smallskip

Special cases:

\smallskip

\begin{verbatim}
[z->sqrt(Pi)*exp(z^2)*erfc(z),[0,2*z],[2,2*n]]
[z->1/sqrt(Pi)*intnum(t=[-oo,1],[oo,1],exp(-t^2)/(z-t)),
                                       [0,2*z],[2,-2*n]]
[()->erfc(1/sqrt(2))*sqrt(exp(1)*Pi/2),[0,1],[1,n]]
[z->F20(1,1/2,2*z),[0,1],[1,-n*z]]
\end{verbatim}

\medskip

Finally, if $(d_a,d_b)=(1,1)$, this is \ref{6.23.9}, giving the following,
valid for $ac\ne0$:

\begin{cf}{\ }
\begin{verbatim}
[(a,b,c,d)->(b+c/a)*F11(d/c-1,b/a+c/a^2,c/a^2)/
                    F11(d/c,b/a+c/a^2+1,c/a^2),[a*n+b,c*n+d]]
\end{verbatim}
\begin{align*}(b+c/a)&\dfrac{{}_1F_1(d/c-1;b/a+c/a^2;c/a^2)}{{}_1F_1(d/c;b/a+c/a^2+1;c/a^2)}\\
  &=b+\dfrac{d}{a+b+\dfrac{c+d}{2a+b+\dfrac{2c+d}{3a+b+\dfrac{3c+d}{4a+b+\dfrac{4c+d}{5a+b+\dfrac{5c+d}{6a+b+\ddots}}}}}}\end{align*}
\end{cf}

\smallskip

Special cases:

\smallskip

\begin{verbatim}
[()->exp(1),[1,1,n+1],[1,-n]]
[()->exp(1),[5/2,n+1],[1/2,n]]
[()->exp(1),[2,n+1],[n+2]]
[()->exp(1),[n+3],[-n-1]]
[()->exp(2),[1,1,n+2],[2,-2*n]]
[()->exp(2),[5,n+1],[8,2*n+4]]
[()->exp(3),[1,1,n+3],[3,-3*n]]
[()->exp(3),[13,n+1],[27,3*n+6]]
[()->besseli(1,1)/besseli(0,1),[0,3,n+3],[1,-2*n-1]]
[()->besseli(1,2)/besseli(0,2),[0,4,n+5],[2,-4*n-2]]
[(k,z)->0,[n+(z-k)],[-z*n-z]]
[(z)->1,[0,n+(z-1)],[n+(z+1)]]
[(a,z)->1+a/(z+1),[0,a*n+(z+(-a-1))],[a*n+(z+a)]]
[(z)->(z^2+z+1)/(z^2-z+1),[0,n+(z-4)],[n+z]]
[(z)->(z^3+2*z+1)/((z-1)^3+2*z-1),[0,n+(z-5)],[n+z]]
[(z)->exp(z),[1,1,n+z],[z,-z*n]]
[(z)->exp(z),[1,n-z],[z,z*n]]
[(a,z)->incgamc(a,z)/(z^a*exp(-z)),[0,n+(-z+(a-1))],[1,z*n]]
[(a,z)->incgamc(a,z)/(z^a*exp(-z)),[0,a,n+(z+(a-1))],[1,-z*n+(-a+1)*z]]
[(z)->incgamc(z,z)/(z^(z-1)*exp(-z)),[n+1],[z*n+2*z]]
[(z)->incgamc(z,z)/(z^(z-1)*exp(-z)),[z+1,n+(2*z+1)],[-z*n-z^2]]
[(z)->sqrt(Pi)*exp(z^2)*erf(z),[0,2*n+(-2*z^2-1)],[2*z,4*z^2*n]]
[(z)->exp(-z^2)*intnum(t=0,z,exp(t^2)),[0,2*n+(2*z^2-1)],[z,-4*z^2*n]]
[(z)->sqrt(Pi)*exp(z^2)*erf(z),[0,1,2*n+(2*z^2-1)],[2*z,-4*z^2*n+2*z^2]]
[(z)->exp(-z^2)*intnum(t=0,z,exp(t^2)),[0,1,2*n+(-2*z^2-1)],[z,4*z^2*n-2*z^2]]
[()->erf(1/sqrt(2))*sqrt(exp(1)*Pi/2),[1,2*n+2],[2,2*n+4]]
[()->erf(1/sqrt(2))*sqrt(exp(1)*Pi/2),[0,1,2*n],[-2*n+1]]
[(nu,z)->besseli(nu+1,z)/besseli(nu,z),[0,z+(2*nu+2),n+(2*z+(2*nu+1))],
                                       [z,-2*z*n+(-2*nu-1)*z]]
[(b,z)->2*b*F01(b,z^2)/F01(b+1,z^2),[2*z+2*b,n+(4*z+2*b)],[-4*z*n+(-4*b-2)*z]]
[(a,b,z)->b*F11(a,b,z)/F11(a+1,b+1,z),[n+(-z+b)],[z*n+(a+1)*z]]
[(b,z)->F11(1,b+1,z),[0,n+(-z+(b-1))],[b,z*n]]
[(b,z)->F11(1,b+1,z),[1,b+1,n+(z+b)],[z,-z*n-b*z]]
[(z)->F11(1,z+1,z),[z+1,n+(2*z+1)],[-z*n-z^2]]
\end{verbatim}

\smallskip

It is interesting to note that the same CFs are sometimes obtained from
different functions. For instance, comparing CFs we find the identity
$$\dfrac{{}_0F_1(;b;z^2)}{{}_0F_1(;b+1,z^2)}=\dfrac{z}{b}+\dfrac{{}_1F_1(b-1/2;2b;4z)}{{}_1F_1(b+1/2;2b+1;4z)}\;.$$
which can easily be deduced from the known formulas
\begin{align*}
  {}_1F_1(a;2a;z)&=2^{2a-1}\G(a+1/2)z^{1/2-a}e^{z/2}I_{a-1/2}(z/2)\;,\\
  {}_1F_1(a;2a+1;z)&=2^{2a-1}\G(a+1/2)z^{1/2-a}e^{z/2}(I_{a-1/2}(z/2)-I_{a+1/2}(z/2))\;,\text{\quad and}\\
  {}_0F_1(;a,z^2)&=2^{2a-2}\G(a)z^{1-a}I_{a-1}(2z)\;.
\end{align*}

\appendix
\chapter{Asymptotic Expansions and Questions}\label{chap:asymp}

\section{Asymptotic Expansions}

The speed of convergence has been studied in Chapter \ref{chap:speed}, and in
{\tt Pari/GP} is obtained by the {\tt cftype} command up to a multiplicative
constant, which itself is obtained by the {\tt cfasymp()[3][5]} command. We
may want to obtain further terms in the asymptotic expansion. This is obtained
by adding a positive argument $L$ to the {\tt cftype} command asking for an
asymptotic expansion up to $1/n^L$ (for instance {\tt cftype(AB,10)}).
The purpose of this chapter is to explain how this is implemented.

\smallskip

To achieve this, we need four things:
first determine the asymptotics of $\prod_{0\le j\le n}b(j)=b(0)b!(n)$, then
those of $q(n)$, then of $S(n+1)-S(n)$, and finally of $S-S(n)$.

\smallskip

We have the following:

\begin{lemma}\label{lem:prod} Assume that a sequence $u(n)$ has an asymptotic
  expansion of the form
  $$u(n)=n^vu_0(1+u_1/n+u_2/n^2+\cdots)$$
  with $u_0\ne0$.
  Then $u!(n)=\prod_{1\le j\le n}u(n)$ has the asymptotic expansion
  $$u!(n)=n!^vu_0^nn^{u_1}C(1+e_1/n+e_2/n^2+...)\;,$$
  where $C$ is some nonzero constant, where the coefficients $e_i$ are obtained
  thanks to the following formulas. Write
  $$\log(1+u_1x+u_2x^2+\cdots)=\sum_{i\ge1}f_ix^i\text{\quad and\quad}
    g_m=\sum_{1\le i\le m+1}f_i\binom{m}{i-1}B_{m-i+1}\;,$$
  where the $B_k$ are the Bernoulli numbers. Then
  $$1+e_1x+e_2x^2+\cdots=\exp\left(-\sum_{m\ge1}(g_m/m)x^m\right)\;.$$
\end{lemma}

\begin{proof} From the Euler--MacLaurin summation formula, we obtain for
  any integer $k\ge1$ the asymptotic expansion:
  $$\sum_{1\le n\le N}\dfrac{1}{n^k}=C(N,k)-\sum_{j\ge0}\dfrac{(k+j-2)!}{(k-1)!j!}\dfrac{B_j}{N^{k+j-1}}\;,$$
  where $C(N,k)=\zeta(k)$ if $k\ge2$ and $C(N,1)=\log(N)+\ga$.
  The result now follows by taking logarithms and collecting the coefficients
  of $1/n^m$ for all $m\ge1$.\end{proof}

Thanks to this lemma we can compute the asymptotics of
$\prod_{0\le j\le n}b(j)$, and we will be able to compute that of $q(n)$
once we determine that of $v(n)=q(n)/q(n-1)$, since $q(n)=v!(n)$.

\smallskip

The analogous lemma when $u(n)$ has an asymptotic expansion involving
half-integral exponents is proved in exactly the same way, and is as
follows:

\begin{lemma}\label{lem:prod2} Assume that a sequence $u(n)$ has an asymptotic
  expansion of the form
  $$u(n)=n^vu_0(1+u_{1/2}/n^{1/2}+u_1/n+u_{3/2}/n^{3/2}+u_2/n^2+\cdots)$$
  with $u_0\ne0$.
  Then $u!(n)=\prod_{1\le j\le n}u(n)$ has the asymptotic expansion
  $$u!(n)=n!^vu_0^ne^{2u_{1/2}n^{1/2}}n^{u_1-u_{1/2}^2/2}C(1+e_{1/2}/n^{1/2}+e_1/n+e_{3/2}/n^{3/2}+e_2/n^2+...)\;,$$
  where $C$ is some nonzero constant, where the coefficients $e_i$ are obtained
  thanks to the following formulas. Write
  \begin{align*}\log(1+u_{1/2}y+u_1y^2+u_{3/2}y^3+\cdots)&=\sum_{i\ge1}f_{i/2}y^i\text{\quad and}\\
    g_{m/2}&=\sum_{\substack{1\le i\le m+2\\i\equiv m\pmod2}}f_{i/2}\binom{m/2}{i/2-1}B_{(m-i)/2+1}\;.\end{align*}
  $$\text{Then\quad}1+e_{1/2}y+e_1y^2+e_{3/2}y^3+\cdots=\exp\left(-\sum_{m\ge1}(g_{m/2}/(m/2))y^m\right)\;.$$
\end{lemma}

We now need to determine the asymptotics of $S-S(n)$ knowing those of
$S(n+1)-S(n)$. The following lemma is a strong refinement of Lemma
\ref{lem:sns} and the trivial proof is left to the reader:

\begin{lemma}\label{lem:sum} Keep the notation and assumptions of Lemma
  \ref{lem:sns}. Assume for now that $D=0$. Then if
  $$S(n+1)-S(n)=\dfrac{C}{n!^FE^nn^P}(1+u_1/n+u_2/n^2+\cdots)\;,$$
  we have
  $$S-S(n)=\dfrac{C'}{n!^FE^nn^P}n^u(1+e_1/n+e_2/n^2+\cdots)\;,$$
  where $u=0$ unless $F=0$ and $E=1$ (case $P^+$) where $u=1$, and
  where the coefficients $e_i$ are obtained thanks to the following formulas.
  Set $f(x)=x^u(1+u_1x+u_2x^2+\cdots)$ and $g(x)=1+e_1x+e_2x^2+\cdots$.
  Then for some explicit constant $C''$
  $$g(x)-(x^F/E)(1+x)^{-(P+F)}g(x/(x+1))=C''f(x)\;.$$
  More precisely, the $e_i$ are all given by linear equations, and
  \begin{enumerate}\item In case $F$ ($F>0$), $C''=1$ and the coefficient of
    $e_i$ is always equal to $1$.
  \item In cases $E$ and $P^-$ ($F=0$ and $E\ne1$), $C''=(E-1)/E$, and the
    coefficient of $e_i$ is also always equal to $C''$.
  \item In case $P^+$ ($F=0$ and $E=1$), $C''=P$, and the coefficient of
    $e_i$ is equal to $P+i$, so is nonzero.
  \end{enumerate}
\end{lemma}

The corresponding lemma for $D>0$ is as follows:

\begin{lemma}\label{lem:sum2} Assume that $D>0$, and let $E=\pm1$.
  Set $u=0$ if $E=-1$ and $u=1$ if $E=1$. Then if
  $$S(n+1)-S(n)=\dfrac{C}{E^ne^{(Dn)^{1/2}}n^{u/2}}(1+u_{1/2}/n^{1/2}+u_1/n+u_{3/2}/n^{3/2}+u_2/n^2+\cdots)\;,$$
  we have
  $$S-S(n)=\dfrac{C'}{E^ne^{(Dn)^{1/2}}}(1+e_{1/2}/n^{1/2}+e_1/n+e_{3/2}/n^{3/2}+e_2/n^2+\cdots)\;,$$
  where the coefficients $e_i$ are obtained thanks to the following formulas.
  Set $f(y)=y^u(1+u_{1/2}y+u_1y^2+\cdots)$ and $g(y)=1+e_{1/2}y+e_1y^2+\cdots$.
  Then for some explicit constant $C''$
  $$g(y)-Ee^{-\sqrt{D}(\sqrt{1+y^2}-1)/y}g(y/\sqrt{1+y^2})=C''f(y)\;.$$
  More precisely, the $e_i$ are all given by linear equations, and
  \begin{enumerate}\item In case $D^-$ ($E=-1$), $C''=2$ and the coefficient of
    $e_i$ is always equal to $2$.
  \item In case $D^+$ ($E=1$), $C''=\sqrt{D}/2$ and the coefficient of $e_i$
    is always equal to $C''$.
  \end{enumerate}
\end{lemma}

\medskip

As in the computation of the speed of convergence, we are thus reduced to
finding the asymptotic expansion of $v(n)=q(n)/q(n-1)$. In that proof
we reduced to the case $a(n)=1$ by dividing by $a(n)$. We could do that here,
but it is not necessary.
The recursion formula for $v(n)$ is $v(n+1)=a(n+1)+b(n)/v(n)$, and we have
seen in the proof that $v(n)=n^{\ga}(v_0+v_{1/2}/n^{1/2}+v_1/n+\cdots)$
for some $\ga$ and $v_i$.

We assume first that we are not in cases $D^{\pm}$, so that the expansion
is in integral powers of $1/n$, and consider the remaining cases later.
From the proof of the speed of convergence, we know that $\ga=\max(\al,\be/2)$.
With $x=1/n$, we set
\begin{align*}
  a(n)&=n^{\al}\sum_{i\ge0}a_i/n^i=x^{-\al}\sum_{i\ge0}a_ix^i=A(x)\;,\\
  b(n)&=n^{\be}\sum_{i\ge0}b_i/n^i=x^{-\be}\sum_{i\ge0}b_ix^i=B(x)\;,\\
  v(n)&=n^{\ga}\sum_{i\ge0}v_i/n^i=x^{-\ga}\sum_{i\ge0}v_ix^i=V(x)\;,
\end{align*}
(note that the $a_i$ and $b_i$ are not the same as before), so the recursion
is $V(x/(x+1))=A(x/(x+1))+B(x)/V(x)$, in other words $R(x)=0$, with
$$R(x)=V(x)(V(x/(x+1))-A(x/(x+1)))-B(x)\;.$$
We now consider the different cases of the proof, taking of course into
account that we have not divided by $a(n)$.

\begin{enumerate}
\item Case $\be\le2\al-1$: this is case $F$. We know
  that $\ga=\al$ and $v_0=a_0$. By expanding in indeterminate coefficients
  the recursion $R(x)=0$, it is easy to check that for $i\ge1$ each unknown
  $v_i$ occurs in a \emph{linear} equation with nonzero leading coefficient,
  so can be solved uniquely.
\item Case $\be=2\al+2$ and $b_0>0$: this is case $P^-$. We know that
  $\ga=\be/2$ and $v_0=b_0^{1/2}$. Once again we obtain only
  linear equations with nonzero leading coefficients.
\item Case $\be=2\al$ and $b_0>-a_0^2/4$: this is case $E$. We know that
  $\ga=\be/2=\al$ and $v_0=(|a_0|+\sqrt{a_0^2+4b_0})/2$. Once again we
  obtain only linear equations with nonzero leading coefficients.
\item Case $\be=2\al$, $b_0=-a_0^2/4$, $b_1/b_0=\al+2a_1/a_0$,
  and $B<0$ with $B=(2\al+2a_1/a_0-1)^2-2\al(\al-1)-4(b_2/b_0-2a_2/a_0)$:
  this is case $P^+$. Assume first that $B$ is not the square of an integer.
  We know that $\ga=\be/2=\al$, $v_0=|a_0|/2$, and
  $v_1=s_0((1+\sqrt{B})a_0+2a_1)/4$. Once again we obtain only linear
  equations, but now the leading coefficients are $2v_1-s_0(a_1+(k/2)a_0)$ for
  $k\ge2$ integral, and $$2v_1-s_0(a_1+(k/2)a_0)=|a_0|(1+\sqrt{B}-k)/2\;,$$
  and since we have assumed that $B$ is not the square of an integer, this
  is always nonzero.
\item Same, but now assume that $B=b^2$ is the square of an integer. This
  is more complicated. The above computation shows that we can at least
  find the coefficients of $1/n^k$ in the asymptotic expansion for $k\le b$.
  If we want more, note the following:
  \begin{enumerate}\item The expansion may involve $\log(n)$, as in case of
    logarithmic convergence: in fact, applying Bauer--Muir to a CF of type
    $L$ gives one of type $P^+$ with this property.
  \item When the expansion is in pure integral powers of $1/n$, it may be
    possible to compute the convergents explicitly. This is the case in
    particular when the CF comes from a summation, in which case we can apply
    the Euler--MacLaurin formula.
  \end{enumerate}
\end{enumerate}

Let us give an example of this last method. It is immediate to show that
$\sum_{n\ge1}1/(4n-3)^2=(\pi^2+G)/8$, where $G=L(\chi_{-4},2)$ is Catalan's
constant. Thus, applying Euler's transformation, we set
\begin{verbatim}
CF=[[0,1,(4*n-3)^2+(4*n-7)^2],[8,-(4*n-3)^4]]
\end{verbatim}
We know that the limit of {\tt CF} is $\pi^2+G$, and, forgetting for now
where it comes from, our program will immediately tell us that
we have convergence type $P^+$ with $P=1$ and $C=1/2$, so that
$\pi^2+G-p(n)/q(n)\sim1/(2n)$. However, our program will also tell us that
since we are in case $P^+$ with $P$ integral, it does not know how
to compute the asymptotics. Here this is very simple: using a toy
implementation of Euler--MacLaurin, we write
\begin{verbatim}
? E=eulermaclaurin(8/(4*n-3)^2)
% = -1/2*x - 1/8*x^2 + ...
/* Confirms convergence in 1/(2n) */
? E/(-x/2)
% = 1 + 1/4*x - 1/48*x^2 - 3/64*x^3 + 7/3840*x^4 + ...
\end{verbatim}
so the asymptotic is $A=1+(1/4)/n-(1/48)/n^2-...$.

\medskip

\medskip

We now consider the cases $D^{\pm}$. The expansions for $a(n)$ and $b(n)$
are unchanged, but now we set
$$v(n)=n^{\ga}\sum_{i\ge0}v_{i/2}/n^{i/2}=x^{-\ga}\sum_{i\ge0}v_{i/2}x^{i/2}=V(x)\;.$$
We now have two cases.
\begin{enumerate}
\item Case $\be=2\al+1$ and $b_0>0$: this is case $D^-$. We know that
  $\ga=\be/2=\al+1/2$ and $v_0=b_0^{1/2}$. Once again we obtain only
  linear equations with nonzero leading coefficients equal to $2b_0^{1/2}$.
\item Case $\be=2\al$, $b_0=-a_0^2/4<0$, and $b_1/b_0<\al+2a_1/a_0$: this is
  case $D^+$. We know that $\ga=\al$ and $v_0=a_0/2$. Here we choose
  $$v_1=(b_1-b_0(\al+2a_1/a_0))^{1/2}=(b_1+a_0/4(a_0\al+2a_1))^{1/2}\;,$$
  and after we obtain only linear equations with nonzero leading
  coefficients equal to $2v_1$.
\end{enumerate}

\smallskip

To summarize, apart from the logarithmic case, for which it is already
difficult to find the main term, the only remaining important case is the
case $P^+$ when $P=\sqrt{B}$ is an integer, and as already mentioned, I have no
general recipe to suggest.

\section{Questions}

In compiling and computing this huge list of continued fractions, a number
of sporadic questions have arisen. Here are a few:

\begin{enumerate}
\item As mentioned above, it is frustrating that in a single type of
  convergence (type $P^+$ with $P$ integral), I do not know of a systematic
  method to find the asymptotic expansion $A$. Closely related is the
  fact that such asymptotic expansions may involve $\log(n)$, and also
  have irrational coefficients not related to $d$, the square root
  of the discriminant.
\item The CF given as an example for $3^{1/2}$ (\ref{5.2.1}) has the
  curious property of being infinitely contractible: usually, computing
  the contraction of a CF of bidegree $(d_1,d_2)$ gives a CF of larger
  bidegree. For the present CF of bidegree $(2,4)$, the contractions
  become considerably more complicated, but can always be simplified
  (i.e., to a CF with the same partial quotients) to bidegree $(2,4)$.
  First, I do not know how to prove this, and second, I do not know how to
  characterize such CFs.
\item The parametric families, especially those for periods of degree $2$
  and $3$, have been obtained by searches. In many cases they are probably
  not complete, but in some cases there are so many that it would be
  desirable to classify them, which I have not always succeeded in doing.
  See \cite{Coh4} for detailed discussions of this.
\item In the case of logarithmic convergence (type $L$), I do not know of
  any reasonable method to compute the constant $C$ such that
  $S-p(n)/q(n)\sim C/\log(n)$, except when the CF is completely explicit,
  or when it is the specialization of a CF with more variables
  (for instance, this is how I found that the constant $C$ for \ref{5.4.27} is
  equal to $\pi^3/8$ and that of \ref{5.5.10} is $\pi^4/28$).
\item Since there are so many CFs for the functions $\be(z)$ and $\be_1(z)$,
  it is slightly surprising that there are no interesting ones for
  $\be_2(z)=-\be_1'(z)$; nonetheless, we have given a few (uninteresting)
  ones in Section \ref{sec:be2}.
\item Although $I_{\nu}(z)$ is an even function of $z$ (up to a power of $z$),
  the CF \ref{6.22.6} does not seem to define an even function. This implies
  that there does not seem to exist a corresponding CF for $J_{\nu}(z)$ (not
  using complex numbers), and consequently that there does not seem to
  exist analogs for $J$ of the CFs \ref{5.12.10} and \ref{5.12.11} for $I$.
\item There exist many CFs involving the Bessel functions $J_{\nu}(z)$,
  $I_{\nu}(z)$, and $K_{\nu}(z)$, as well as
  $H_{\nu}(z)=J_{\nu}(z)\pm i Y_{\nu}(z)$. Are there interesting ones involving
  $Y_{\nu}(z)$ alone ?
\item Using Ap\'ery, it is not difficult to find a parametric family of CFs
  for $\cosh(\pi z)$, see \ref{6.6.9}. On the other hand, the very similar
  CF \ref{6.6.20} for $\sinh(\pi z/3)/(z\sqrt{3}/2)$ does not accelerate
  nicely, so I have been unable to find a corresponding parametric family.
\end{enumerate}

\chapter{Implementation of Bauer--Muir--Ap\'ery}\label{chap:apimpl}

The different sections of this appendix are extremely technical and are
included only to help with possible implementations, but should be skipped
otherwise.

\section{Theory}

\medskip

Let $(a(n),b(n))$ be a CF with partial quotients $(p(n),q(n))$, so that
$(p(-1),q(-1))=(1,0)$, $(p(0),q(0))=(a(0),1)$, and
$$(p(n+1),q(n+1))=a(n+1)(p(n),q(n))+b(n)(p(n-1),q(n-1))\;.$$
The basic Bauer--Muir acceleration is to set
$$u'(n)=u(n+1)+r(n+1)u(n)$$
for some suitably chosen $r(n)$, where $u(n)=p(n)$ or $q(n)$.
Note that there are several possible variants, for instance we could
set $u'(n)=s(n+1)u(n+1)+r(n+1)u(n)$ for some function $s(n)$ close to $1$,
but it is not difficult to see that we can always reduce to the basic case
by suitable changes of functions.

\begin{proposition}
  Set $(p'(n),q'(n))=(p(n+1),q(n+1))+r(n+1)(p(n),q(n))$ for
  $n\ge0$, define
  $$d(n)=r(n)(r(n+1)+a(n+1))-b(n)\;,$$
  and assume that $d(n)\ne0$ for all $n$.
  Then if we set $(p'(-1),q'(-1))=(a(0)+r(0),1)$,
  $(p'(n),q'(n))$ are the $n+1$-st convergents of the CF $(a'(n),b'(n))$
  defined by
  \begin{align*}
    a'(-1)&=a(0)+r(0)\;,\quad b'(-1)=-d(0)\;,\quad a'(0)=r(1)+a(1)\;,\\
    a'(n)&=r(n+1)+a(n+1)-r(n-1)\dfrac{d(n)}{d(n-1)}\\
    b'(n)&=b(n)\dfrac{d(n+1)}{d(n)}\;.
  \end{align*}
\end{proposition}

\begin{proof} Here and elsewhere, let $u(n)=p(n)$ or $u(n)=q(n)$.
We have
\begin{align*}
  u'(n+1)&=u(n+2)+r(n+2)u(n+1)\\
  u'(n)&=u(n+1)+r(n+1)u(n)\\
  u'(n-1)&=u(n)+r(n)u(n-1)\\
  u(n+2)&=a(n+2)u(n+1)+b(n+1)u(n)\\
  u(n+1)&=a(n+1)u(n)+b(n)u(n-1)\;,\end{align*}
so keeping $u(n)$ and $u(n-1)$ as main variables $x$ and $y$, we have
\begin{align*}
  u(n+1)&=a(n+1)x+b(n)y\\
  u(n+2)&=(a(n+2)a(n+1)+b(n+1))x+a(n+2)b(n)y
\end{align*}
Thus,
\begin{align*}
  &u'(n+1)-a'(n+1)u'(n)-b'(n)u'(n-1)\\
  &=u(n+2)+(r(n+2)-a'(n+1))u(n+1)\\
  &\phantom{=}-(a'(n+1)r(n+1)+b'(n))u(n)-b'(n)r(n)u(n-1)\\
  &=x(a(n+2)a(n+1)+b(n+1)+r(n+2)a(n+1)\\
  &\phantom{=}-a'(n+1)(a(n+1)+r(n+1))-b'(n))\\
  &\phantom{=}+y((a(n+2)+r(n+2)-a'(n+1))b(n)-b'(n)r(n))\end{align*}
giving the system
$Aa'(n+1)+Bb'(n)=U$, $Ca'(n+1)+Db'(n)=V$ with
\begin{align*}
  A&=a(n+1)+r(n+1)\;,\ B=1\;,\ C=b(n)\;,\ D=r(n)\\
  U&=(a(n+2)+r(n+2))a(n+1)+b(n+1)\\
  V&=(a(n+2)+r(n+2))b(n)\end{align*}
We will thus set
$$d(n)=AD-BC=r(n)(r(n+1)+a(n+1))-b(n)\;,$$
so in particular
$$r(n+2)+a(n+2)=(d(n+1)+b(n+1))/r(n+1)\;,$$
so
\begin{align*}
  U&=(d(n+1)+b(n+1))a(n+1)/r(n+1)+b(n+1)\\
  V&=(d(n+1)+b(n+1))b(n)/r(n+1)\;,\end{align*}
so that
\begin{align*}
  &d(n)a'(n+1)=DU-BV=r(n)U-V\\
  &=(r(n+2)+a(n+2))(r(n)a(n+1)-b(n))+r(n)b(n+2))\\
  &=(r(n+2)+a(n+2))(d(n)-r(n)r(n+1))+r(n)b(n+2)\\
  &=(r(n+2)+a(n+2))d(n)-r(n)(r(n+1)(r(n+2)+a(n+2))-b(n+2))\\
  &=(r(n+2)+a(n+2))d(n)-r(n)d(n+1)\\
  &d(n)b'(n)=AV-CU=(a(n+1)+r(n+1))V-b(n)U\\
  &=b(n)/r(n+1)((a(n+1)+r(n+1))(d(n+1)+b(n+1))\\
  &\phantom{=}-a(n+1)(d(n+1)+b(n+1))-b(n+1)r(n+1))=b(n)d(n+1)\;.\end{align*}
Thus, with $a'(n)$ and $b'(n)$ as in the proposition, for $n\ge1$ we have
$u'(n+1)=a'(n+1)u'(n)+b'(n)u'(n-1)$. Finally, we have
\begin{align*}
  (p'(-1),q'(-1))&=(a(0),1)+r(0)(1,0)=(a(0)+r(0),1)\\
  (p'(0),q'(0))&=(a(0)a(1)+b(1),a(1))+r(1)(a(0),1)\\
  &=((r(1)+a(1))a(0)+b(1),r(1)+a(1))\end{align*}
We want $p'(-1)/q'(-1)=a'(-1)$ and $q'(-1)=1$, so we must set
$p'(-1)=a'(-1)=a(0)+r(0)$. We want
$$(p'(0),q'(0))=a'(0)(p'(-1),q'(-1))+b'(-1)(1,0)=(a'(0)p'(-1)+b'(-1),a'(0)q'(-1))\;,$$
so $a'(0)=q'(0)=q(1)+r(1)q(0)=a(1)+r(1)$, hence
\begin{align*}b'(-1)&=p'(0)-a'(0)p'(-1)=p(1)+r(1)p(0)-(a(1)+r(1))(a(0)+r(0))\\
  &=a(0)a(1)+b(0)+r(1)a(0)-a(1)a(0)\\
  &\phantom{=}-r(1)a(0)-(a(1)+r(1))r(0)\\
  &=b(0)-r(0)(r(1)+a(1))=-d(0)\end{align*}
\end{proof}

To summarize, we have proved that we have the formal identity
$$\dfrac{p'(n)}{q'(n)}=a'(-1)+\dfrac{b'(-1)}{a'(0)+\dfrac{b'(0)}{a'(1)+\dfrac{b'(1)}{a'(2)+\ddots+\dfrac{b'(n-1)}{a'(n)}}}}$$

\medskip

\section{Contraction}

Note that we have set $u'(n)=u(n+1)+r(n)u(n)$ instead of the
possibly more natural $u(n)+r(n-1)u(n-1)$: the reason for this is
that this is what we need for Ap\'ery acceleration, but has the effect
of creating a continued fraction of length $n+1$ from the first $n$ terms.
We thus need to contract it to a CF of length $n$, not by changing $n+1$ to
$n$, nor by contracting the last convergents (this could not be used in
Ap\'ery), but by contracting the \emph{first} convergents. This is done by
using the following trivial lemma:
\begin{lemma} We have (assuming $a(1)\ne0$)
  $$a(0)+\dfrac{b(0)}{a(1)+\dfrac{b(1)}{a(2)+x}}=a'(0)+\dfrac{b'(0)}{a'(1)+x}$$
  with $a'(0)=a(0)+b(0)/a(1)$, $b'(0)=-b(0)b(1)/a(1)^2$, and
  $a'(1)=a(2)+b(1)/a(1)$.\end{lemma}
\begin{corollary}
  With the notation of the proposition, set
  \begin{align*}
    a''(0)&=a'(-1)+\dfrac{b'(-1)}{a'(0)}=a(0)+\dfrac{b(0)}{r(1)+a(1)}\\
    a''(1)&=a'(1)+\dfrac{b'(0)}{a'(0)},\quad b''(0)=-\dfrac{b'(-1)b'(0)}{a'(0)^2}\\
    a''(n)&=a'(n)\text{\quad for $n\ge2$\;,\quad}b''(n)=b'(n)\text{\quad for $n\ge1$}\end{align*}
  Then for $n\ge2$, $(p'(n),q'(n))$ is the $n$th convergent of the continued
  fraction defined by $(a''(n),b''(n))$.\end{corollary}

In other words
$$\dfrac{p'(n)}{q'(n)}=a''(0)+\dfrac{b''(0)}{a''(1)+\dfrac{b''(1)}{a''(2)+\dfrac{b''(2)}{a''(3)+\ddots+\dfrac{b''(n-1)}{a''(n)}}}}$$

\section{Appendix: Gory Details II: Implementation of Bauer--Muir}

We must now consider the different convergence types to see how to choose
$r(n)$ so that $d(n)$ is as small as possible. We assume that $a(n)$ and
$b(n)$ are polynomials with rational coefficients for $n$ large enough, and
want to choose $r(n)$ also as a polynomial with rational coefficients.
For now, we will consider the simplest case $s=t=1$, and see how to modify
later. Recall that
$$a(n)=a_0n^{\al}(1+a_1/n+a_2/n^2+...)\text{\quad and\quad}b(n)=b_0n^{\be}(1+b_1/n+b_2/n^2+...)$$
with $a_0b_0\ne0$, and let $s_0=\pm1$ be the sign of $a_0$.

\medskip

\subsection{The Case $P^-$}

\medskip

This corresponds to $\be=2\al+2$ and $b_0>0$. We must choose
$$r(n)=s_0b_0^{1/2}n^{\be/2}(1+r_1/n+r_2/n^2+\cdots+r_{\be/2}/n^{\be/2})$$
for suitable $r_i$, and in fact we know that
$r_1=(b_1-\be/2-|a_0|/b_0^{1/2})/2$, but we can ignore this for now.
We will set $x=1/n$ and $A(x)=1+a_1x+a_2x^2+\cdots$, and similarly
$B(x)$ and $R(x)$. Note that changing $n=1/x$ into $n+1=(x+1)/x$ corresponds
to changing $x$ into $x/(x+1)$. Thus, translating the
identity $d(n)=r(n)(r(n+1)+a(n+1))-b(n)$ in terms of $x$:
\begin{align*}d(n)&=s_0b_0^{1/2}x^{-\be/2}R(x)(s_0b_0^{1/2}x^{-\be/2}(1+x)^{\be/2}R(x/(x+1))\\
  &\phantom{=}+a_0x^{1-\be/2}(1+x)^{\be/2-1}A(x/(x+1)))-b_0x^{-\be}B(x)\\
  &=b_0x^{-\be}(1+x)^{\be/2}Z(x)\text{\quad with}\\
  Z(x)&=R(x)(R(x/(x+1))+|a_0|/b_0^{1/2}(x/(1+x))A(x/(x+1)))-(1+x)^{-\be/2}B(x)\;.
\end{align*}
Clearly $Z(0)=0$, and the coefficient of $x$ in $Z$ is equal to
$2r_1+|a_0|/b_0^{1/2})x+\be/2-b_1$, giving the formula above for $r_1$, but
we do not need it. Indeed, we simply work with indeterminate coefficients:
We set $|a_0|/b_0^{1/2}(x/(1+x))A(x/(x+1))=\sum_{i\ge1}c_ix^i$,
$(1+x)^{-\be/2}B(x)=1+\sum_{i\ge1}d_ix^i$, and
$R(x/(x+1))=1+\sum_{i\ge1}s_ix^i$, where
$s_i=\sum_{1\le j\le i}(-1)^{i-j}\binom{i-1}{i-j}r_j$ depends \emph{linearly}
on $r_j$ for $1\le j\le i$, the coefficient of $r_i$ being equal to $1$.
Thus setting naturally $c_0=0$ and $r_0=s_0=d_0=1$, we have
\begin{align*}Z(x)&=(\sum_{i\ge0}r_ix^i)(\sum_{i\ge0}(s_i+c_i))-\sum_{m\ge0}d_mx^m\;,\text{\quad so expanding:}\\
Z(x)&=\sum_{m\ge0}(-d_m+\sum_{0\le i\le m}r_{m-i}(s_i+c_i))x^m\;.\end{align*}
The coefficient of $x^0$ is of course $0=-d_0+r_0(s_0+c_0)$. For $m\ge1$
the coefficient of $x^m$ is equal to
\begin{align*}
  &-d_m+r_m(s_0+c_0)+r_0(s_m+c_m)+\sum_{1\le i\le m-1}r_{m-i}(s_i+c_i)\\
  &=-d_m+c_m+r_m+s_m+\sum_{1\le i\le m-1}r_{m-i}(s_i+c_i)\;,\end{align*}
and since $$s_m=r_m+\sum_{1\le i\le m-1}(-1)^{m-i}\binom{m-1}{m-i}r_i$$
we obtain the recursion
\begin{align*}2r_m&=d_m-c_m-\sum_{1\le i\le m-1}(-1)^{m-i}\binom{m-1}{m-i}r_i\\
  &\phantom{=}-\sum_{1\le i\le m-1}r_{m-i}(c_i+s_i)
\end{align*}
which we solve for $1\le m\le\be/2$. This will ensure that
$Z(x)=O(x^{\be/2+1})$, hence that $d(n)=O(x^{1-\be/2})=O(n^{\al})$.
Of course, in some cases $d(n)$ will be smaller than that and even constant.

\medskip

\subsection{The Case $P^+$}

\medskip

This corresponds to $\be=2\al$, $a_0^2+4b_0=0$, $2a_1=b_1-\al$, and
$B_3>0$ with $B_3=(b_1+\al-1)^2-2\al(\al-1)-4(b_2-2a_2)$. Here we must choose
$$r(n)=-(a_0/2)n^{\al}(1+r_1/n+r_2/n^2+\cdots+r_{\al}/n^{\al})\;,$$
and in fact we know that $r_1=(b_1+1-\al\pm\sqrt{B_3})/2$ for a suitable
sign $\pm$. Here
\begin{align*}d(n)&=-(a_0/2)x^{-\al}R(x)(-(a_0/2)x^{-\al}(1+x)^{\al}R(x/(x+1))\\
  &\phantom{=}+a_0x^{-\al}(1+x)^{\al}A(x/(x+1)))+(a_0^2/4)x^{-2\al}B(x)\\
  &=(a_0^2/4)x^{-2\al}(1+x)^{\al}Z(x)\text{\quad with}\\
  Z(x)&=R(x)(R(x/(x+1))-2A(x/(x+1)))+(1+x)^{-\al}B(x)\;.\end{align*}
Here again $Z(0)=0$, the coefficient of $x$ in $Z$ is equal to
$b_1-2a_1-\al=0$, and that of $x^2$ is equal to
$$r_1^2+(\al-b_1-1)r_1+b_2-2a_2+(1-\al)b_1+\al(\al-1)/2\;.$$
We check that the discriminant of this second degree equation is equal to
$B_3$, so that $r_1=(b_1+1-\al\pm\sqrt{B_3})/2$. Once the sign chosen, we
proceed as before. We set $A(x/(x+1))=\sum_{i\ge0}c_ix^i$ (with $c_0=1$),
$(1+x)^{-\al}B(x)=\sum_{i\ge0}d_ix^i$ (with $d_0=1$), and
$R(x/(x+1))=\sum_{i\ge0}s_ix^i$ (with $s_0=1$), so
$$Z(x)=(\sum_{i\ge0}r_ix^i)(\sum_{i\ge0}(s_i-2c_i)x^i)+\sum_{m\ge0}d_mx^m\;.$$
The coefficient of $x^m$ for $m=0$ and $1$ vanish since we are in case $P^+$,
and that of $x^2$ vanishes by our choice of $r_1$. For $m\ge3$ the coefficient
$Z_m$ of $x^m$ is given by
\begin{align*}Z_m&=d_m+r_m(s_0-2c_0)+r_0(s_m-2c_m)+r_{m-1}(s_1-2c_1)\\
  &\phantom{=}+r_1(s_{m-1}-2c_{m-1})+S\,\text{\quad with\quad} S=\sum_{2\le i\le m-2}r_{m-i}(s_i-2c_i),\text{\quad so that}\\
Z_m&=d_m-2c_m+(s_m-r_m)+r_{m-1}(s_1-2c_1)+r_1(s_{m-1}-2c_{m-1})+S\;.
\end{align*}
Since
$$s_m=r_m-(m-1)r_{m-1}+\sum_{1\le i\le m-2}(-1)^{m-i}\binom{m-1}{m-i}r_i\;,$$
the coefficient of $r_m$ in $Z_m$ vanishes, and that of $r_{m-1}$ is equal to
$$-(m-1)+s_1-2c_1+r_1=1-m+2r_1-2a_1=2-m\pm\sqrt{B_3}\;.$$
Thus, if we choose the minus sign, this will always be nonzero. If we choose
the plus sign, it will always be nonzero if and only if $B_3$ is not the
square of an integer.

Thus, although we could complexify even more, we proceed as follows.
If $B_3$ is the square of an integer, we choose the minus sign, and
use the recursion giving $r_{m-1}$ by looking at the coefficient of $x^m$
for $m\ge3$. If $B_3$ is not the square of an integer, we try both
signs using the recursion in the same way, and after computing
$r_{\al}$ (by looking at the coefficient of $x^{\al+1}$) we compare the
degrees of $d(n)$ and choose the sign giving the smallest degree.

\medskip

\subsection{The Case $E$ with $E$ Rational}

\medskip

This corresponds to $\be=2\al$ and $\Delta=a_0^2+4b_0=d^2$, which is the
square of a rational number. We must choose
$$r(n)=((-a_0+s_0d)/2)n^{\al}(1+r_1/n+r_2/n^2+\cdots+r_{\al}/n^{\al})\;,$$
and in fact we know that $r_1=((b_1-\al)(|a_0|+d)/2-|a_0|a_1)/d$.
Here, setting $R_0=(-a_0+s_0d)/2$
\begin{align*}
  d(n)&=R_0x^{-\al}R(x)(R_0x^{-\al}(1+x)^{\al}R(x/(x+1))\\
  &\phantom{=}+a_0x^{-\al}(1+x)^{\al}A(x/(x+1)))-b_0x^{-2\al}B(x)\\
  &=x^{-2\al}(1+x)^{\al}Z(x)\text{\quad with}\\
  Z(x)&=R(x)(R_0^2R(x/(x+1))+a_0R_0A(x/(x+1)))-b_0(1+x)^{-\al}B(x)\;.\end{align*}
We have $Z(0)=R_0^2+a_0R_0-b_0=0$ by definition. Thus, setting as before
$r_0=1$ etc... (not to be confused with $R_0$) we have
$$Z(x)=(\sum_{i\ge0}r_ix^i)(\sum_{i\ge0}(R_0^2s_i+a_0R_0c_i))-b_0\sum_{m\ge0}d_mx^m\;,$$
so the coefficient of $x^m$ is equal to
\begin{align*}&-b_0d_m+r_m(R_0^2+a_0R_0)+R_0^2(s_m+a_0R_0c_m)\\
  &\phantom{=}+\sum_{1\le i\le m-1}r_{m-i}(R_0^2s_i+a_0R_0c_i)\\
  &=b_0(r_m-d_m)+R_0^2(s_m+a_0R_0c_m)+\sum_{1\le i\le m-1}r_{m-i}(R_0^2s_i+a_0R_0c_i)\;.\end{align*}
Since $s_m=r_m+\cdots$, the coefficient of $r_m$ is equal to
$$R_0^2+b_0=2b_0-a_0R_0=(4b_0+a_0^2-|a_0|d)/2=d(d-|a_0|)/2\;,$$
and since $d>0$ and $d\ne|a_0|$ (otherwise $b_0=0$), this is nonzero so
the recursion can proceed.

\medskip

\subsection{The Case $D^+$}

\medskip

This corresponds to $\be=2\al$, $a_0^2+4b_0=0$, and $2a_1>b_1-\al$.
We must choose $r(n)=(-a_0/2)n^{\al}(1+r_1/n+\cdots+r_{\al}/n^{\al})$,
The beginning of the computation is identical to that of $P^+$, we have
\begin{align*}d(n)&=(a_0^2/4)x^{-2\al}(1+x)^{\al}Z(x)\text{\quad with}\\
  Z(x)&=R(x)(R(x/(x+1))-2A(x/(x+1)))+(1+x)^{-\al}B(x)\;,\end{align*}
and the coefficient of $x$ in $Z$ is now equal to $b_1-2a_1-\al<0$, so
is nonzero. Thus, the coefficients $r_i$ will not improve the degree of $d(n)$,
which will always be equal to $2\al-1=\be-1$. However, in view of the
Bauer--Muir formulas, if we want polynomial coefficients (which of course
will rarely be possible), since $d(n)$ is coprime to $d(n+1)$, we must
have $d(n)\mid b(n)$ and $d(n)\mid r(n)$, which gives a condition. We thus
proceed as follows: if $b(n)$ does not have a degree $1$ divisor (over $\Q$),
nothing much can be done, choose for instance $r(n)=(-a_0/2)n^{\al}$.
If $b(n)$ does have a degree $1$ divisor, say $b(n)=(un+v)c(n)$ with
$\deg(c)=\be-1$ monic, we must have $d(n)=d_0c(n)$ for some constant $d_0$,
and hence $r(n)=(-a_0/2)c(n)$ since $c$ is monic, which implies that
$\al=\deg(r)=\be-1=2\al-1$, in other words this can happen only when
$\al=1$ and $\be=2$, and $b(n)$ a second degree polynomial which factors
into linear factors. We can redo the computation:
$r(n)=(-a_0/2)(n+r_1)$ with $-r_1$ one of the roots of $b(n)$,
$a(n)=a_0n+a_0a_1$, $b(n)=(-a_0^2/4)(n+r_1)(n+u)$, say, so
$$d(n)=(a_0^2/4)(n+r_1)(r_1+u-2a_1-1)\;.$$
This is always divisible by $r(n)$, but we also need it to be nonzero:
note that $r_1+u$ is the opposite of the sum of the roots of $b(n)$,
and since $b(n)=b_0(n^2+b_1n+b_2)$ we have $r_1+u=b_1$, so
$r_1+u-2a_1-1=b_1-2a_1-1<0$ never vanishes.

To summarize, in case $D^+$ we can find a new polynomial type continued
fraction if and only if $\al=1$, hence $\be=2$, and $b(n)$ factors into
two linear factors over $\Q$, in which case one must choose
$r(n)=(-a_0/2)(n+r_1)$, where $n+r_1$ is one of the linear factors of $b(n)$.

\section{Appendix: Gory Details III: Ap\'ery Acceleration}\label{sec:apgen}

The main idea is to iterate Bauer--Muir transforms, and to take a diagonal or
pseudo-diagonal. We must find the necessary formulas.

Let $(a(n),b(n),p(n),q(n))$ be an initial CF, which we will consider as our
$0$th iteration, so set $(a(n,0),b(n,0),p(n,0),q(n,0))=(a(n),b(n),p(n),q(n))$,
and assume by induction that at the $l$-th iteration we have the recursion
$$u(n+1,l)=a(n+1,l)u(n,l)+b(n,l)u(n-1,l)\;,$$
where here and below $u$ stands for $p$ or $q$; by assumption this is true for
$l=0$. In view of the Bauer--Muir formulas we set
\begin{align*}
  R(n,l)&=r(n+1,l)+a(n+1,l)\text{\quad and}\\
  d(n,l)&=r(n,l)R(n,l)-b(n,l)\;.\end{align*}

Now the Bauer--Muir formulas have denominators (unless, for instance,
$d(n,l)$ is independent of $n$), so we introduce a third function
$t(n,l)$, and define
\begin{align*}
  a(n,l+1)&=t(n,l)(R(n,l)-r(n-1,l)d(n,l)/d(n-1,l))\;,\\
  b(n,l+1)&=t(n,l)t(n+1,l)b(n,l)d(n+1,l)/d(n,l)\;,\\
  u(n,l+1)&=t!(n,l)(u(n+1,l)+r(n+1,l)u(n,l))\;,\end{align*}
the factorial being of course taken as a function of $n$.

Then an easy check shows that we indeed have the recursion
$$u(n+1,l+1)=a(n+1,l+1)u(n,l+1)+b(n,l+1)u(n-1,l+1)\;,$$
which is thus valid for all $l$.

We would now like to take the diagonal, i.e., find a recursion for
$(p(n,n),q(n,n))$. This is possible, but the formulas again involve
denominators, exactly the same as those occurring when \emph{contracting} a CF.
Thus, as we will see, it is much cleaner to use \emph{staircase steps} instead
of direct diagonal ones. In other words, we will find a 3-term recursion
between on the one hand $(u(n-1,n-1),u(n,n-1),u(n,n))$, and on the other hand
$(u(n,n-1),u(n,n),u(n+1,n))$, and putting both together we will obtain
a \emph{period 2} CF, in that the even and odd coefficients are given by
different formulas.

As for the Bauer--Muir formulas, this is simply done by solving trivial linear
systems. Before embarking on this, let us recall that our recursions
are as follows:
\begin{align*}
  u(n+1,l)&=a(n+1,l)u(n,l)+b(n,l)u(n-1,l)\text{\quad and}\\
  u(n,l+1)&=t!(n,l)(u(n+1,l)+r(n+1,l)u(n,l))\end{align*}

\begin{proposition}
  We have the recursions
  \begin{align*}
    u(n,l+1)&=t!(n,l)(R(n,l)u(n,l)+b(n,l)u(n-1,l))\;,\\
    u(n+1,l+1)&=t(n+1,l)R(n+1,l)u(n,l+1)-t!(n+1,l)d(n+1,l)u(n,l)\;,\\
    u(n,l+1)&=t!(n,l)(A(n,l)u(n,l)+B(n,l)u(n,l-1))\;,\text{\quad with}\\
    A(n,l)&=t(n+1,l-1)R(n+1,l-1)+r(n+1,l)\text{\quad and}\\
    B(n,l)&=-t!(n+1,l-1)d(n+1,l-1)\;.
  \end{align*}
\end{proposition}

\begin{proof} We have
  \begin{align*}
    u(n,l+1)&=t!(n,l)(u(n+1,l)+r(n+1,l)u(n,l))\\
    &=t!(n,l)((a(n+1,l)u(n,l)+b(n,l)u(n-1,l))+r(n+1,l)u(n,l))\\
    &=t!(n,l)R(n,l)u(n,l)+t!(n,l)b(n,l)u(n-1,l)\;,\end{align*}
  so the first recursion is immediate.

  For the second, setting $x=u(n+1,l)$ we have
  \begin{align*}
    u(n+1,l+1)&=t!(n+1,l)(u(n+2,l)+r(n+2,l)x)\\
    &=t!(n+1,l)((a(n+2,l)x+b(n+1,l)u(n,l))+r(n+2,l)x)\\
    &=t!(n+1,l)(R(n+1,l)x+b(n+1,l)u(n,l))\;,\end{align*}
  and since $u(n,l+1)=t!(n,l)(x+r(n+1,l)u(n,l))$ we have
  $$x=u(n,l+1)/t!(n,l)-r(n+1,l)u(n,l)\;,$$
  so
  \begin{align*}
    u(n+1,l+1)&=t!(n+1,l)(R(n+1,l)(u(n,l+1)/t!(n,l)-r(n+1,l)u(n,l))\\
    &\phantom{=}+b(n+1,l)u(n,l))\\
    &=t(n+1,l)R(n+1,l)u(n,l+1)\\
    &\phantom{=}-t!(n+1,l)(r(n+1,l)R(n+1,l)-b(n+1,l))u(n,l)\\
    &=t(n+1,l)R(n+1,l)u(n,l+1)-t!(n+1,l)d(n+1,l)u(n,l)\;,\end{align*}
  proving the second recursion.

  For the third, we have
  \begin{align*}
    u(n,l+1)&=t!(n,l)(u(n+1,l)+r(n+1,l)u(n,l))\\
    &=t!(n,l)(t(n+1,l-1)R(n+1,l-1)u(n,l)\\
    &\phantom{=}-t!(n+1,l-1)d(n+1,l-1)u(n,l-1))+r(n+1,l)u(n,l))\\
    &=A(n,l)u(n,l)+B(n,l)u(n,l-1)\;,\end{align*}
  proving the third recursion.\end{proof}

We can now play with different paths giving new CFs.

\begin{corollary} (Vertical walk). Let $m\ge0$ be fixed, and for $n\ge2$ set
  \begin{align*}
    a_V(n)&=t(m+1,n-2)R(m+1,n-2)+r(m+1,n-1)\text{\quad and}\\
    b_V(n)&=-t(m+1,n-1)d(m+1,n-1)\;.\end{align*}
  Then up to a choice of initial terms the CF $(a_V,b_V)$ converges to
  the same limit as the initial CF.
\end{corollary}

For instance, if $t=1$ this means that,
\emph{up to a M\"obius transformation}, for a fixed $m$ we have
$$S=R(m+1,0)+r(m+1,1)-\dfrac{d(m+1,1)}{R(m+1,1)+r(m+1,2)-\dfrac{d(m+1,2)}{R(m+1,2)+r(m+1,3)-\ddots}}\;.$$

\smallskip

In {\tt GP} notation, omitting the initial values which must be determined
case by case, this is written as:

\begin{verbatim}
[[t(m+1,n-2)*R(m+1,n-2)+r(m+1,n-1)],[-t(m+1,n-1)*d(m+1,n-1)]]
\end{verbatim}

\smallskip

\begin{corollary} (Staircase). Let $m\ge0$ be fixed. We have
  \begin{align*}
    u(n+m+1,n+1)&=A_0(n+1)u(n+m,n+1)+B_0(n)u(n+m,n)\text{\quad and}\\
    u(n+m,n+1)&=A_1(n+1)u(n+m,n)+B_1(n)u(n+m-1,n)\;,\text{\quad with}\\
    A_0(n)&=t(n+m,n-1)R(n+m,n-1)\\
    B_0(n)&=-t!(n+m+1,n)d(n+m+1,n)\\
    A_1(n)&=t!(n+m-1,n-1)R(n+m-1,n-1)\\
    B_1(n)&=t!(n+m,n)b(n+m,n)\;.
  \end{align*}
  In particular, if $t=1$:
  \begin{align*}
    A_0(n)&=R(n+m,n-1)\;,\quad B_0(n)=-d(n+m+1,n)\\
    A_1(n)&=R(n+m-1,n-1)\;,\quad B_1(n)=b(n+m,n)\;.\end{align*}
\end{corollary}
This means that, \emph{up to a M\"obius transformation} we have
$$S=R(m,0)-\dfrac{d(m+1,0)}{R(m+1,0)+\dfrac{b(m+1,1)}{R(m+1,1)-\dfrac{d(m+2,1)}{R(m+2,1)+\dfrac{b(m+2,2)}{R(m+2,2)-\ddots}}}}$$
(with evident modifications when $t$ is nontrivial).

\smallskip

In {\tt GP} notation, omitting the initial values which must be determined
case by case, this is written as:

\begin{verbatim}
[[[t(n+m,n-1)*R(n+m,n-1),tf(n+m,n)*R(n+m,n)]],
 [[tf(n+m,n)*b(n+m,n),-tf(n+m+1,n)*d(n+m+1,n)]]]
\end{verbatim}
where {\tt tf(n,l)} is the result of the computation of
$t!(n,l)=\prod_{1\le m\le n}t(m,l)$.

\chapter{Hypergeometric Series}\label{chap:hyper}

\section{Introduction}

\subsection{Interesting Series}

As already mentioned, in compiling the database of CFs given in this book
it was found partly experimentally that a substantial
proportion of these continued fractions correspond in fact to generalized
hypergeometric series through Euler's transformation of series into
continued fractions. The purpose of this chapter is to give a list of
the most interesting ones. Note that I do not claim to have proved all
of them.

The notion of ``interesting'' series is completely subjective, so here is a
non-exhaustive list of the series that I do \emph{not} include.

$\bullet$ Rational series $\sum_{n\ge0}P(n)/Q(n)$ where $Q(x)$ is a polynomial
having only rational roots. Indeed, by expanding into partial fractions the
sum of such a series can be computed explicitly, so is not interesting.
When $Q(x)$ has irreducible factors of degree at least $2$, it is also
easy to sum the series in terms of complex values of polygamma functions,
but the series are already more interesting so we keep them, except when they
sums to less interesting constants or functions, such as linear combinations
of zeta or L-values.

$\bullet$ Series comming by trivial specializations of classical power series.
For instance the series $\log(2)=\sum_{n\ge0}3^{-2n-1}/(n+1/2)$ is
uninteresting since it simply comes from the Taylor series expansion of
$-\log(1-x)$ for $x=1/9$.

$\bullet$ Series which are trivially equivalent to the definition, for
instance $e=1+\sum_{n\ge0}1/(n+1)!$ or $\psi'(z)=\sum_{n\ge1}1/(n+z-1)^2$,
or correspond to well-known expansions such as
$\atanh(z)=\sum_{n\ge0}z^{2n+1}/(2n+1)$.

\subsection{Use of the Pochhammer symbol}

Although all the series could be expressed as hypergeometric evaluations,
I have preferred to give the series more explicitly using the rising
Pochhammer symbol $(a)_n=a(a+1)\cdots(a+n-1)=\G(a+n)/\G(a)$, and replacing
expressions such as $(a+k)_n/(a)_n$ or (its inverse) with $k\in\Z$ by
rational functions. Several remarks:

\smallskip

\begin{enumerate}\item When $a$ is a positive integer, $(a)_n$ is (up to
  a multiplicative constant) a factorial, so we will write $n!$, $(n+1)!$,
  and $(n+2)!/2$ instead of $(1)_n$, $(2)_n$, or $(3)_n$ respectively.
\item When $a$ is half an odd integer, we will keep the notation as is,
  but of course note that $(1/2)_n=(2n)!/(2^{2n}n!)$,
  $(3/2)_n=(2n+1)!/(2^{2n}n!)$, and so on.
\end{enumerate}

\subsection{Numerical Check of the Formulas}

It is of course essential to be able to check numerically the given formulas.
Using {\tt Pari/GP} (and evidently other systems) this is very easy but
uses an essential trick.

All the series are of the form
$$S=u+\sum_{n\ge0}h(n)\prod_{1\le i\le g}(a_i)_n^{e_i}w^n\;,$$
where $h$ is a rational function, $e_i\in\Z$, and all the other quantities
do not depend on $n$. Note the trivial functional equation
$$S=u'+\sum_{n\ge0}h'(n)\prod_{1\le i\le g}(a_i+1)_n^{e_i}w^n\;,$$
with $u'=u+h(0)$ and $h'(n)=w\prod_{1\le i\le g}a_i^{e_i}h(n+1)$.

For convergence, we must have in particular
$E:=\sum_{1\le i\le g}e_i\le 0$, and if $E=0$ we must have $|w|\le1$.
 
\begin{enumerate}\item If $E<0$ or if $E=0$ and $|w|<1$, simply sum the series
  until the terms become negligible. In {\tt Pari/GP} this is the command
  {\tt suminf}.
\item If $E=0$ and $w=-1$, use the command {\tt sumalt}, very
  efficient for summing alternating series (note that there are no series
  in our list with $|w|=1$ and $w\ne\pm1$).
\item The more difficult case is when $E=0$ and $w=1$. Here, one
  must use the powerful discrete Euler--MacLaurin summation explained
  in \cite{Bel-Coh} and available in {\tt Pari/GP} under the name {\tt sumnum}.
  It works very well but must be used with care:
\end{enumerate}

The standard syntax is {\tt sumnum(n=a,f(n))}. This is OK if $f(n)$ is for
instance a rational function of $n$. But if $f$ involves Pochhammer symbols,
this cannot work for several reasons. The first essential reason is that
{\tt sumnum} uses \emph{nonintegral} values of $n$, so defining
{\tt poch(a,n)=prod(j=0,n-1,a+j)}, or using {\tt n!} would not work. Instead,
one can define {\tt poch(a,n)=gamma(a+n)/gamma(a)}, and replace {\tt n!}
by {\tt gamma(n+1)}. But even that is not sufficient.

Consider for example the following series:

$$S=\sum_{n\ge0}\dfrac{(3/2)_n}{(2n+3)(n+1)!}\;,$$

whose value is $\pi-2$. Writing

{\tt S=sumnum(n=0,gamma(n+3/2)/(gamma(3/2)*(2*n+3)*gamma(n+2)))}

results in an error message ``overflow in expo()'', because both the
numerator and denominator overflow for large values of $n$. Thus a more
reasonable command is

{\tt S=sumnum(n=0,exp(lngamma(n+3/2)-lngamma(n+2))/(2*n+3))/gamma(3/2)}

This seems to work but gives a completely wrong answer. The reason is because
internally $n$ is used as a real number with insufficient accuracy. The
solution if to write the following (assuming we are working at 38 decimals
of accuracy):

\begin{verbatim}
S=sumnum(n0=0,my(n=precision(n0,77));
     exp(lngamma(n+3/2)-lngamma(n+2))/(2*n+3))/gamma(3/2)
\end{verbatim}

\noindent
and now the result is correct to $31$ decimal digits, not perfect but
sufficient for numerically checking the formula.

\section{Hypergeometric Series: Constants}

I repeat that the following hypergeometric evaluations come from the
construction of a large database of continued fractions, so are not
necessarily the most interesting hypergeometric formulas. In addition,
although these formulas come from standard transformations of continued
fractions, in many cases I have simply extrapolated from the first few
terms, so I do not claim to have proved all of them rigorously, mainly for
lack of time.

\vfill\eject

\begin{align*}\root3\of2&=2\sum_{n\ge0}(-1)^n\dfrac{(2/3)_n}{n!}\\
\root3\of2&=\dfrac{5}{4}\sum_{n\ge0}\dfrac{(-1/3)_n}{n!}(-3/125)^n\\
\dfrac{1}{\root3\of2}&=2\sum_{n\ge0}\dfrac{(-2/3)_n}{n!}(3/4)^n\\
\dfrac{1}{\root3\of2}&=\sum_{n\ge0}\dfrac{(-1/3)_n}{n!}2^{-n}\\
\dfrac{1}{\root3\of2}&=\dfrac{4}{5}\sum_{n\ge0}\dfrac{(-1/3)_n}{n!}(3/128)^n\end{align*}
\begin{align*}\pi&=2\sum_{n\ge0}\dfrac{(1/2)_n}{(2n+1)n!}\\
  \pi&=4-\dfrac{4}{9}\sum_{n\ge0}\dfrac{n!(n+1)!}{(5/2)_n^2}\\
\pi&=4-4\sum_{n\ge0}\dfrac{n!^2}{(4n+1)(4n+5)(3/2)_n^2}\\
\pi&=3+\dfrac{3}{35}\sum_{n\ge0}\dfrac{n!^2}{(11/6)_n(13/6)_n}\\
\pi&=2\sum_{n\ge0}\dfrac{n!}{(3/2)_n}2^{-n}\\
\pi&=3\sum_{n\ge0}\dfrac{(1/2)_n}{(2n+1)n!}2^{-2n}\\
\dfrac{1}{\pi}&=\dfrac{1}{2}\sum_{n\ge0}(-1)^n\dfrac{(4n+1)(1/2)_n^3}{n!^3}\\
\dfrac{1}{\pi}&=\dfrac{1}{4}\sum_{n\ge0}\dfrac{(1/2)_n^2}{(n+1)n!^2}\\
\dfrac{1}{\pi}&=\dfrac{1}{4}+\dfrac{1}{16}\sum_{n\ge0}\dfrac{(1/2)_n^2}{(n+1)^2n!^2}\\
\dfrac{1}{\pi}&=-\dfrac{1}{2}\sum_{n\ge0}\dfrac{(1/2)_n^2}{(2n-1)n!^2}\\
\dfrac{1}{\pi}&=\dfrac{1}{4}\sum_{n\ge0}\dfrac{(1/2)_n^2}{(2n-1)^2n!^2}\\
\dfrac{1}{\pi}&=-\sum_{n\ge0}\dfrac{(1/2)_n^2}{(4n-1)(4n+3)n!^2}\\
\dfrac{1}{\pi}&=\dfrac{1}{3}\sum_{n\ge0}\dfrac{(-1/6)_n(1/6)_n}{n!^2}
\end{align*}
\begin{align*}\dfrac{\pi}{\sqrt{3}}&=\dfrac{1}{2}\sum_{n\ge0}\dfrac{(4/3)_n}{(n+1)(5/3)_n}\\
\dfrac{\pi}{\sqrt{3}}&=\dfrac{3}{2}+\dfrac{3}{16}\sum_{n\ge0}\dfrac{n!^2}{(5/3)_n(7/3)_n}\\
\dfrac{\pi}{\sqrt{3}}&=\dfrac{3}{4}\sum_{n\ge0}\dfrac{n!}{(3/2)_n}(3/4)^n\\
\dfrac{\pi}{\sqrt{3}}&=\dfrac{3}{2}\sum_{n\ge0}\dfrac{n!}{(3/2)_n}2^{-2n}\\
\dfrac{\sqrt{3}}{\pi}&=\dfrac{2}{3}\sum_{n\ge0}\dfrac{(-1/3)_n(1/3)_n}{n!^2}\end{align*}
\begin{align*}\dfrac{\pi}{\sqrt{2}}&=2+\dfrac{2}{15}\sum_{n\ge0}\dfrac{n!^2}{(7/4)_n(9/4)_n}\\
\dfrac{\pi}{\sqrt{2}}&=2\sum_{n\ge0}\dfrac{(1/2)_n}{(2n+1)n!}2^{-n}\\
\dfrac{\sqrt{2}}{\pi}&=\dfrac{1}{2}\sum_{n\ge0}\dfrac{(-1/4)_n(1/4)_n}{n!^2}\\
  \dfrac{\pi}{\sqrt{2}-1}&=8-\dfrac{8}{63}\sum_{n\ge0}\dfrac{(4n+3)n!^2(3/8)_n(5/8)_n}{(3/2)_n^2(15/8)_n(17/8)_n}\\
  \dfrac{\pi}{\sqrt{2}+1}&=\dfrac{8}{3}-\dfrac{24}{55}\sum_{n\ge0}\dfrac{(4n+3)n!^2(1/8)_n(7/8)_n}{(3/2)_n^2(13/8)_n(19/8)_n}\\
  \dfrac{\sqrt{2}-1}{\pi}&=\dfrac{2}{15}\sum_{n\ge0}\dfrac{(4n+1)(1/2)_n^2(-1/8)_n(1/8)_n}{n!^2(11/8)_n(13/8)_n}\\
  \dfrac{\sqrt{2}+1}{\pi}&=\dfrac{6}{7}\sum_{n\ge0}\dfrac{(4n+1)(1/2)_n^2(-3/8)_n(3/8)_n}{n!^2(9/8)_n(15/8)_n}\end{align*}
\begin{align*}\log(2)&=\dfrac{1}{2}\sum_{n\ge0}\dfrac{(7/6)_n}{(3n+2)(3/2)_n}\\
\log(2)&=\dfrac{1}{4}\sum_{n\ge0}\dfrac{(3/2)_n}{(n+1)(n+1)!}\\
\log(2)&=\dfrac{15}{16}\sum_{n\ge0}\dfrac{n!}{(3/2)_n}(-9/16)^n\\
\log(2)&=\dfrac{3}{4}\sum_{n\ge0}(-1)^n\dfrac{n!}{(3/2)_n}2^{-3n}\\
\log(2)&=\dfrac{1}{2}+\dfrac{1}{2}\sum_{n\ge0}(-1)^n\dfrac{(3n+4)n!}{(3n+5)(3n+2)(3/2)_n}2^{-3n}\end{align*}
\begin{align*}\pi^2&=4\sum_{n\ge0}\dfrac{n!}{(n+1)(3/2)_n}\\
\pi^2&=\dfrac{27}{4}\sum_{n\ge0}\dfrac{n!}{(n+1)(3/2)_n}(3/4)^n\\
\pi^2&=8\sum_{n\ge0}\dfrac{n!}{(n+1)(3/2)_n}2^{-n}\\
\pi^2&=9\sum_{n\ge0}\dfrac{n!}{(n+1)(3/2)_n}2^{-2n}\\
\pi^2&=10\sum_{n\ge0}(-1)^n\dfrac{(1/2)_n}{(2n+1)^2n!}2^{-2n}\\
\pi^2&=\dfrac{3}{4}\sum_{n\ge0}\dfrac{(21n+13)n!^3}{(3/2)_n^3}2^{-6n}\end{align*}

\begin{align*}L(\chi_{-3},2)&=\dfrac{8}{15}\sum_{n\ge0}\dfrac{n!(n+1)!}{(2n+1)(5/4)_n(7/4)_n}(3/4)^n\\
  \dfrac{L(\chi_{-3},2)}{\sqrt{3}}&=\dfrac{4}{9}\sum_{n\ge0}\dfrac{(1/2)_n}{(2n+1)^2n!}2^{-2n}\end{align*}

Recall that $G=L(\chi_{-4},2)$ is Catalan's constant.

\begin{align*}G&=\dfrac{1}{2}\sum_{n\ge0}\dfrac{n!}{(2n+1)(3/2)_n}\\
G&=1-\dfrac{1}{30}\sum_{n\ge0}\dfrac{(6n+7)n!^2}{(2n+3)(7/4)_n(9/4)_n}2^{-2n}\\
G&=\dfrac{1}{2}-3\sum_{n\ge0}(-1)^n\dfrac{n\cdot n!^3}{(3n^2-3n+1)(3n^2+3n+1)(1/2)_n^3}2^{-3n}\end{align*}
\begin{align*}\log^2(2)&=\dfrac{9}{16}\sum_{n\ge0}(-1)^n\dfrac{n!}{(n+1)(3/2)_n}(3/4)^{2n}\\
  \log^2(2)&=\dfrac{1}{2}\sum_{n\ge0}(-1)^n\dfrac{n!}{(n+1)(3/2)_n}2^{-3n}\\
  \log^2(3/2)&=\dfrac{25}{144}\sum_{n\ge0}(-1)^n\dfrac{n!}{(n+1)(3/2)_n}(5/12)^{2n}\\
  \log^2(3/2)&=\dfrac{1}{6}\sum_{n\ge0}(-1)^n\dfrac{n!}{(n+1)(3/2)_n}24^{-n}\\
\log^2((1+\sqrt{5})/2)&=\dfrac{1}{4}\sum_{n\ge0}(-1)^n\dfrac{n!}{(n+1)(3/2)_n}2^{-2n}\end{align*}
\begin{align*}
  \pi^3&=\dfrac{216}{7}\sum_{n\ge0}\dfrac{(1/2)_n}{(2n+1)^3n!}2^{-2n}\\
  \dfrac{1}{\pi^3}&=\dfrac{1}{32}\sum_{n\ge0}\dfrac{(6n+1)(28n^2+8n+1)(1/2)_n^7}{n!^7}2^{-6n}\\
  &=\dfrac{1}{32}\sum_{n\ge0}(6n+1)(28n^2+8n+1)\binom{2n}{n}^72^{-20n}\end{align*}
\begin{align*}\z(3)&=\dfrac{2}{21}\sum_{n\ge0}(-1)^n\dfrac{(20n^2+32n+13)n!^2}{(n+1)(2n+1)^2(5/4)_n(7/4)_n}2^{-2n}\\
\z(3)&=1+\dfrac{1}{12}\sum_{n\ge0}(-1)^n\dfrac{(5n^2+14n+10)n!}{(n+1)^2(n+2)^2(5/2)_n}2^{-2n}\\
\z(3)&=\dfrac{5}{4}\sum_{n\ge0}(-1)^n\dfrac{n!}{(n+1)^2(3/2)_n}2^{-2n}\\
\z(3)&=\dfrac{6}{7}+\dfrac{8}{21}\sum_{n\ge0}(-1)^n\dfrac{(28n^2+58n+31)n!(n+1)!^4}{(8n^2+7n+2)(8n^2+23n+17)(3/2)_n^3(5/3)_n(7/3)_n}3^{-3n}\end{align*}
$$\dfrac{\pi^4}{90}=\dfrac{18}{17}\sum_{n\ge0}\dfrac{n!}{(n+1)^3(3/2)_n}2^{-2n}$$
$$7\z(3)-2\pi G+3\pi=12+\dfrac{2}{9}\sum_{n\ge0}\dfrac{n!^2}{(2n+3)(5/2)_n^2}$$
\begin{align*}\pi-2\log(2)&=-2-6\sum_{n\ge0}\dfrac{(5/4)_n}{(n+1)(4n^2-1)(3/4)_n}\\
  \pi-2\log(2)&=\dfrac{2}{3}\sum_{n\ge0}\dfrac{(5/4)_n}{(n+1)(7/4)_n}\\
  \pi-2\log(2)&=\dfrac{5}{3}\sum_{n\ge0}\dfrac{n!(3/2)_n}{(4/3)_n(5/3)_n}(27/2)^{-n}\\
\pi+2\log(2)&=2\sum_{n\ge0}\dfrac{(3/4)_n}{(n+1)(5/4)_n}\\
\pi+2\log(2)&=4+\dfrac{2}{5}\sum_{n\ge0}\dfrac{(3/4)_n}{(n+1)(9/4)_n}\\
\pi-6\log(2)&=\dfrac{4}{3}\sum_{n\ge0}(-1)^{n+1}\dfrac{(5/4)_n}{(n+1)(7/4)_n}\\
\pi-6\log(2)&=-\dfrac{1}{2}\sum_{n\ge0}\dfrac{(5/4)_n}{(n+1)(n+1)!}\\
\pi+6\log(2)&=12\sum_{n\ge0}(-1)^n\dfrac{(7/4)_n}{(n+1)(5/4)_n}\\
\pi+6\log(2)&=\dfrac{3}{2}\sum_{n\ge0}\dfrac{(7/4)_n}{(n+1)(n+1)!}\end{align*}
\begin{align*}\dfrac{\pi}{\sqrt{3}}-2\log(2)&=-1-5\sum_{n\ge0}\dfrac{(7/6)_n}{(n+1)(4n-1)(4n+3)(5/6)_n}\\
\dfrac{\pi}{\sqrt{3}}-2\log(2)&=\dfrac{1}{5}\sum_{n\ge0}\dfrac{(7/6)_n}{(n+1)(11/6)_n}\\
\dfrac{\pi}{\sqrt{3}}+2\log(2)&=\sum_{n\ge0}\dfrac{(5/6)_n}{(n+1)(7/6)_n}\\
\dfrac{\pi}{\sqrt{3}}+2\log(2)&=3+\dfrac{1}{7}\sum_{n\ge0}\dfrac{(5/6)_n}{(n+1)(13/6)_n}\\
\dfrac{\pi}{\sqrt{3}}-\log(3)&=1-\dfrac{4}{3}\sum_{n\ge0}\dfrac{(3n+1)(1/3)_n^2}{(n+1)(3n+5)(5/3)_n^2}\\
\dfrac{\pi}{\sqrt{3}}-\log(3)&=\dfrac{2}{3}\sum_{n\ge0}\dfrac{n!}{(5/3)_n}3^{-2n}\\
\dfrac{\pi}{\sqrt{3}}+\log(3)&=2+\dfrac{4}{3}\sum_{n\ge0}\dfrac{(3n+2)(2/3)_n^2}{(n+1)(3n+4)(4/3)_n^2}\\
\dfrac{\pi}{\sqrt{3}}+\log(3)&=3-\dfrac{1}{2}\sum_{n\ge0}\dfrac{n!}{(7/3)_n}3^{-2n}\\
\dfrac{\pi}{\sqrt{3}}-3\log(3)&=2\sum_{n\ge0}(-1)^{n+1}\dfrac{(4/3)_n}{(n+1)(5/3)_n}\\
\dfrac{\pi}{\sqrt{3}}-3\log(3)&=-\dfrac{2}{3}\sum_{n\ge0}\dfrac{(4/3)_n}{(n+1)(n+1)!}\\
\dfrac{\pi}{\sqrt{3}}+3\log(3)&=8\sum_{n\ge0}(-1)^n\dfrac{(5/3)_n}{(n+1)(4/3)_n}\\
\dfrac{\pi}{\sqrt{3}}+3\log(3)&=\dfrac{4}{3}\sum_{n\ge0}\dfrac{(5/3)_n}{(n+1)(n+1)!}\\
\dfrac{\pi}{\sqrt{3}}-\dfrac{\log(432)}{3}&=\dfrac{4}{15}\sum_{n\ge0}(-1)^{n+1}\dfrac{(7/6)_n}{(n+1)(11/6)_n}\\
\dfrac{\pi}{\sqrt{3}}-\dfrac{\log(432)}{3}&=-\dfrac{1}{9}\sum_{n\ge0}\dfrac{(7/6)_n}{(n+1)(n+1)!}\\
\dfrac{\pi}{\sqrt{3}}+\dfrac{\log(432)}{3}&=\dfrac{20}{3}\sum_{n\ge0}(-1)^n\dfrac{(11/6)_n}{(n+1)(7/6)_n}\\
\dfrac{\pi}{\sqrt{3}}+\dfrac{\log(432)}{3}&=\dfrac{5}{9}\sum_{n\ge0}\dfrac{(11/6)_n}{(n+1)(n+1)!}\end{align*}
\begin{align*}\pi^2-9L(\chi_{-3},2)&=2\sum_{n\ge0}\dfrac{(4/3)_n}{(n+1)^2(5/3)_n}\\
\pi^2+9L(\chi_{-3},2)&=8\sum_{n\ge0}\dfrac{(5/3)_n}{(n+1)^2(4/3)_n}\\
11\pi^2-135L(\chi_{-3},2)&=\dfrac{12}{5}\sum_{n\ge0}\dfrac{(7/6)_n}{(n+1)^2(11/6)_n}\\
11\pi^2+135L(\chi_{-3},2)&=60\sum_{n\ge0}\dfrac{(11/6)_n}{(n+1)^2(7/6)_n}\\
\pi^2-16G&=-4\sum_{n\ge0}\dfrac{n!}{(n+1)(2n+1)(3/2)_n}\\
\pi^2+16G&=24+\dfrac{4}{3}\sum_{n\ge0}\dfrac{n!}{(n+1)(2n+3)(5/2)_n}\\
5\pi^2-48G&=4\sum_{n\ge0}\dfrac{(5/4)_n}{(n+1)^2(7/4)_n}\\
5\pi^2+48G&=36\sum_{n\ge0}\dfrac{(7/4)_n}{(n+1)^2(5/4)_n}\\
\pi^2-6\log(2)^2&=6\sum_{n\ge0}\dfrac{2^{-n}}{(n+1)^2}\\
\pi^2-12\log^2(2)&=4\sum_{n\ge0}\dfrac{n!(3/2)_n}{(n+1)(4/3)_n(5/3)_n}(27/2)^{-n}\\
\pi^2-3\log^2(3)&=\dfrac{16}{3}\sum_{n\ge0}\dfrac{n!(3/2)_n}{(n+1)(4/3)_n(5/3)_n}(81/32)^{-n}\\
\pi^2-18\log^2((1+\sqrt{5})/2)&=6\sum_{n\ge0}(-1)^n\dfrac{n!}{(2n+1)(3/2)_n}2^{-2n}\\
8G-\pi\log(2+\sqrt{3})&=3\sum_{n\ge0}\dfrac{n!}{(2n+1)(3/2)_n}2^{-2n}\\
\dfrac{4G+\pi\log(2)}{\sqrt{2}}&=4\sum_{n\ge0}\dfrac{(1/2)_n}{(2n+1)^2n!}2^{-n}\end{align*}
\begin{align*}3\z(3)-2\log(2)^3&=3\sum_{n\ge0}(-1)^n\dfrac{n!}{(n+1)^2(3/2)_n}2^{-3n}\\
\pi^3+12\pi\log(2)^2&=48\sum_{n\ge0}\dfrac{(1/2)_n}{(2n+1)^3n!}\\
\dfrac{G}{\pi}&=\dfrac{1}{4}\sum_{n\ge0}\dfrac{(1/2)_n^2}{(2n+1)n!^2}\\
\dfrac{1+2G}{\pi}&=-\sum_{n\ge0}\dfrac{(1/2)_n^2}{(2n-1)(2n+1)n!^2}\\
\log(2)-\dfrac{2G}{\pi}&=\dfrac{1}{16}\sum_{n\ge0}\dfrac{(3/2)_n^2}{(n+1)(n+1)!^2}\end{align*}
\begin{align*}e^{-1}&=2\sum_{n\ge0}\dfrac{(-1)^n}{(n^2+3n+1)(n^2+5n+5)n!}\\
e^{-1}&=-12\sum_{n\ge0}(-1)^n\dfrac{n+1}{(n^4+8n^3+17n^2+8n-1)(n^4+12n^3+47n^2+70n+33)n!}\\
e^{-1}&=1-\dfrac{1}{2}\sum_{n\ge0}\dfrac{1}{(1/2)_n(n+1)!}2^{-2n}=1-\dfrac{1}{2}\sum_{n\ge0}\dfrac{1}{(n+1)(2n)!}\\
e^{-1}&=\dfrac{1}{3}\sum_{n\ge0}\dfrac{1}{(5/2)_nn!}2^{-2n}=\sum_{n\ge0}\dfrac{1}{(2n+3)(2n+1)!}\\
e^{-3}&=162\sum_{n\ge0}\dfrac{1}{(n^2+11n+27)(n^2+13n+39)(n+2)!}(-3)^n\end{align*}
\begin{align*}\dfrac{\G(1/4)^3\G(3/4)}{\G(1/2)}&=\sqrt{2\pi}\G(1/4)^2=32\sum_{n\ge0}\dfrac{(1/2)_n}{(4n+1)^2n!}\\
\dfrac{\G(1/4)\G(3/4)^3}{\G(1/2)}&=\sqrt{2\pi}\G(3/4)^2=4\sum_{n\ge0}\dfrac{(-1/2)_n}{(4n-1)^2n!}\\
\dfrac{\G(1/4)^2}{\G(1/2)}&=\dfrac{\G(1/4)^2}{\sqrt{\pi}}=4\sum_{n\ge0}\dfrac{(3/4)_n}{(4n+1)n!}\\
\dfrac{\G(1/4)^2}{\G(1/2)}&=\dfrac{\G(1/4)^2}{\sqrt{\pi}}=8\sum_{n\ge0}(-1)^n\dfrac{(1/2)_n}{(4n+1)n!}\\
\dfrac{\G(1/4)^2}{\G(1/2)}&=\dfrac{\G(1/4)^2}{\sqrt{\pi}}=-144\sum_{n\ge0}(-1)^n\dfrac{(4n-1)(-1/2)_n}{(4n-3)^2(4n+1)^2n!}\\
\dfrac{\G(1/4)^2}{\G(1/2)}&=\dfrac{\G(1/4)^2}{\sqrt{\pi}}=-2+18\sum\dfrac{n!(-1/2)_n}{(5/4)_n^2}\\
\dfrac{\G(1/4)^2}{\G(1/2)}&=\dfrac{\G(1/4)^2}{\sqrt{\pi}}=-16+32\sum_{n\ge0}(-1)^n\dfrac{(3/4)_n^2}{(16n^2-16n+1)(16n^2+16n+1)(1/4)_n^2}\\
\dfrac{\G(1/4)^2}{\G(1/2)}&=\dfrac{\G(1/4)^2}{\sqrt{\pi}}=-32\sum_{n\ge0}(-1)^n\dfrac{(3n+1)(2n+1)(2n+3)(1/2)_n^2}{(12n^2-4n-3)(12n^2+20n+5)(1/4)_n^2}2^{-3n}\\
\dfrac{\G(1/2)}{\G(1/4)^2}&=\dfrac{\sqrt{\pi}}{\G(1/4)^2}=\dfrac{1}{8}\sum_{n\ge0}\dfrac{(-1/4)_n^2}{n!^2}\\
\dfrac{\G(1/2)}{\G(1/4)^2}&=\dfrac{\sqrt{\pi}}{\G(1/4)^2}=\dfrac{1}{16}\sum_{n\ge0}\dfrac{(3/4)_n^2}{(n+1)n!^2}\\
\dfrac{\G(1/2)}{\G(1/4)^2}&=\dfrac{\sqrt{\pi}}{\G(1/4)^2}=-\dfrac{1}{2}\sum_{n\ge0}\dfrac{(-3/4)_n^2}{n!(-1/2)_n}\\
\dfrac{\G(1/2)}{\G(1/4)^2}&=\dfrac{\sqrt{\pi}}{\G(1/4)^2}=\dfrac{1}{4}\sum_{n\ge0}\dfrac{(-1/2)_n(3/2)_n}{n!^2}2^{-n}\end{align*}
\begin{align*}\dfrac{\G(3/4)^2}{\G(1/2)}&=\dfrac{\G(3/4)^2}{\sqrt{\pi}}=2\sum_{n\ge0}\dfrac{(1/4)_n}{(4n+3)n!}\\
\dfrac{\G(3/4)^2}{\G(1/2)}&=\dfrac{\G(3/4)^2}{\sqrt{\pi}}=4\sum_{n\ge0}(-1)^n\dfrac{(3/2)_n}{(4n+3)n!}\\
\dfrac{\G(3/4)^2}{\G(1/2)}&=\dfrac{\G(3/4)^2}{\sqrt{\pi}}=8\sum_{n\ge0}(-1)^n\dfrac{(4n+1)(1/2)_n}{(4n-1)^2(4n+3)^2n!}\\
\dfrac{\G(3/4)^2}{\G(1/2)}&=\dfrac{\G(3/4)^2}{\sqrt{\pi}}=1-\dfrac{1}{9}\sum\dfrac{n!(1/2)_n}{(7/4)_n^2}\\
\dfrac{\G(3/4)^2}{\G(1/2)}&=\dfrac{\G(3/4)^2}{\sqrt{\pi}}=-\dfrac{8}{9}+144\sum_{n\ge0}(-1)^n\dfrac{(1/4)_n^2}{(16n^2-16n+9)(16n^2+16n+9)(3/4)_n^2}\\
\dfrac{\G(3/4)^2}{\G(1/2)}&=\dfrac{\G(3/4)^2}{\sqrt{\pi}}=8\sum_{n\ge0}(-1)^n\dfrac{(3n+1)(3/2)_n^2}{(12n^2-4n+1)(12n^2+20n+9)(3/4)_n^2}2^{-3n}\\
\dfrac{\G(1/2)}{\G(3/4)^2}&=\dfrac{\sqrt{\pi}}{\G(3/4)^2}\dfrac{3}{4}\sum_{n\ge0}\dfrac{(-3/4)_n^2}{n!^2}\\
\dfrac{\G(1/2)}{\G(3/4)^2}&=\dfrac{\sqrt{\pi}}{\G(3/4)^2}=\dfrac{9}{8}\sum_{n\ge0}\dfrac{(1/4)_n^2}{(n+1)n!^2}\\
\dfrac{\G(1/2)}{\G(3/4)^2}&=\dfrac{\sqrt{\pi}}{\G(3/4)^2}=\sum_{n\ge0}\dfrac{(-1/4)_n^2}{n!(1/2)_n}\\
\dfrac{\G(1/2)}{\G(3/4)^2}&=\dfrac{\sqrt{\pi}}{\G(3/4)^2}=\sum_{n\ge0}\dfrac{(1/2)_n^2}{n!^2}2^{-n}\end{align*}
\begin{align*}\dfrac{\G(1/4)\G(1/2)}{\G(3/4)}&=\dfrac{\G(1/4)^2}{\sqrt{2\pi}}=8\sum_{n\ge0}(-1)^n\dfrac{(3/4)_n}{(5/4)_n}\\
\dfrac{\G(1/4)\G(1/2)}{\G(3/4)}&=\dfrac{\G(1/4)^2}{\sqrt{2\pi}}=4+\dfrac{4}{5}\sum_{n\ge0}\dfrac{(3/4)_nn!}{(9/4)_n(3/2)_n}\\
\dfrac{\G(1/4)\G(1/2)}{\G(3/4)}&=\dfrac{\G(1/4)^2}{\sqrt{2\pi}}=-12\sum_{n\ge0}\dfrac{(-1/2)_n}{(4n-3)(4n+1)n!}\\
\dfrac{\G(1/4)\G(1/2)}{\G(3/4)}&=\dfrac{\G(1/4)^2}{\sqrt{2\pi}}=4\sum_{n\ge0}\dfrac{(1/2)_n}{(5/4)_n}2^{-n}\\
\dfrac{\G(3/4)}{\G(1/4)\G(1/2)}&=\dfrac{\sqrt{2\pi}}{\G(1/4)^2}=\dfrac{1}{6}\sum_{n\ge0}\dfrac{(1/4)_n(1/2)_n}{(7/4)_nn!}\\
\dfrac{\G(3/4)}{\G(1/4)\G(1/2)}&=\dfrac{\sqrt{2\pi}}{\G(1/4)^2}=\dfrac{1}{4}\sum_{n\ge0}(-1)^n\dfrac{(2n+1)(1/2)_n^2}{(n+1)n!^2}\\
\dfrac{\G(3/4)}{\G(1/4)\G(1/2)}&=\dfrac{\sqrt{2\pi}}{\G(1/4)^2}=\dfrac{3}{16}\sum_{n\ge0}\dfrac{(1/2)_n(1/4)_n}{n!(n+1)!}2^{-2n}\end{align*}
\begin{align*}\dfrac{\G(3/4)\G(1/2)}{\G(1/4)}&=\dfrac{\G(3/4)^2}{\sqrt{2\pi}}=-1+2\sum_{n\ge0}(-1)^n\dfrac{(1/4)_n}{(3/4)_n}\\
\dfrac{\G(3/4)\G(1/2)}{\G(1/4)}&=\dfrac{\G(3/4)^2}{\sqrt{2\pi}}=1-\dfrac{1}{3}\sum_{n\ge0}\dfrac{n!(1/4)_n}{(3/2)_n(7/4)_n}\\
\dfrac{\G(3/4)\G(1/2)}{\G(1/4)}&=\dfrac{\G(3/4)^2}{\sqrt{2\pi}}=-\dfrac{1}{2}\sum_{n\ge0}\dfrac{(-1/2)_n}{(4n-1)n!}\\
\dfrac{\G(3/4)\G(1/2)}{\G(1/4)}&=\dfrac{\G(3/4)^2}{\sqrt{2\pi}}=1-\dfrac{1}{3}\sum_{n\ge0}\dfrac{(1/2)_n}{(7/4)_n}2^{-n}\\
\dfrac{\G(1/4)}{\G(3/4)\G(1/2)}&=\dfrac{\sqrt{2\pi}}{\G(3/4)^2}=2\sum_{n\ge0}\dfrac{(1/2)_n(-1/4)_n}{n!(5/4)_n}\\
\dfrac{\G(1/4)}{\G(3/4)\G(1/2)}&=\dfrac{\sqrt{2\pi}}{\G(3/4)^2}=2\sum_{n\ge0}(-1)^n\dfrac{(1/2)_n^2}{n!^2}=2\sum_{n\ge0}(-1)^n\dfrac{\binom{2n}{n}^2}{2^{4n}}\\
\dfrac{\G(1/4)}{\G(3/4)\G(1/2)}&=\dfrac{\sqrt{2\pi}}{\G(3/4)^2}=\dfrac{3}{2}\sum_{n\ge0}\dfrac{(1/2)_n(3/4)_n}{n!^2}2^{-2n}\end{align*}
\begin{align*}\dfrac{1}{\G(1/4)^4}&=\dfrac{1}{512}\sum_{n\ge0}(-1)^n\dfrac{(2n+1)(4n+3)(1/2)_n^5}{n!(n+1)!^4}=\dfrac{1}{512}\sum_{n\ge0}(-1)^n\dfrac{(2n+1)(4n+3)}{(n+1)^4}\dfrac{\binom{2n}{n}^5}{2^{10n}}\\
\dfrac{1}{\G(3/4)^4}&=\dfrac{1}{2}\sum_{n\ge0}(-1)^n\dfrac{(4n+1)(1/2)_n^5}{n!^5}=\dfrac{1}{2}\sum_{n\ge0}(-1)^n(4n+1)\dfrac{\binom{2n}{n}^5}{2^{10n}}\\
\dfrac{\G(1/4)^4}{\G(1/2)^2}&=\dfrac{\G(1/4)^4}{\pi}=32\sum_{n\ge0}\dfrac{(3/4)_n}{(2n+1)(5/4)_n}\\
\dfrac{\G(1/4)^4}{\G(1/2)^2}&=\dfrac{\G(1/4)^4}{\pi}=32+96\sum_{n\ge0}\dfrac{(8n+5)n!^2(3/4)_n^2}{(16n^2+4n+1)(16n^2+36n+21)(3/2)_n^2(5/4)_n^2}\\
\dfrac{\G(1/2)^2}{\G(1/4)^4}&=\dfrac{\pi}{\G(1/4)^4}=\dfrac{3}{8}\sum_{n\ge0}\dfrac{(8n+1)(1/2)_n^2(1/4)_n^2}{(16n^2-12n+3)(16n^2+20n+7)n!^2(3/4)_n^2}\\
\dfrac{\G(1/2)^2}{\G(1/4)^4}&=\dfrac{\pi}{\G(1/4)^4}=\dfrac{1}{48}\sum_{n\ge0}\dfrac{(-1/2)_n^3}{n!^3}=-\dfrac{1}{48}\sum_{n\ge0}\dfrac{\binom{2n}{n}^3}{(2n-1)^32^{6n}}\\
\dfrac{\G(3/4)^4}{\G(1/2)^2}&=\dfrac{\G(3/4)^4}{\pi}=\dfrac{2}{3}\sum_{n\ge0}\dfrac{(1/4)_n}{(2n+1)(7/4)_n}\\
\dfrac{\G(3/4)^4}{\G(1/2)^2}&=\dfrac{\G(3/4)^4}{\pi}=\dfrac{2}{3}+\dfrac{2}{3}\sum_{n\ge0}\dfrac{(8n+7)n!^2(5/4)_n^2}{(16n^2+12n+3)(16n^2+44n+31)(3/2)_n^2(7/4)_n^2}\\
\dfrac{\G(1/2)^2}{\G(3/4)^4}&=\dfrac{\pi}{\G(3/4)^4}=\sum_{n\ge0}\dfrac{(1/2)_n^3}{n!^3}=\sum_{n\ge0}\dfrac{\binom{2n}{n}^3}{2^{6n}}\\
\dfrac{\G(1/2)^2}{\G(3/4)^4}&=\dfrac{\pi}{\G(3/4)^4}=6\sum_{n\ge0}\dfrac{(8n+3)(1/2)_n^2(3/4)_n^2}{(16n^2-4n+1)(16n^2+28n+13)n!^2(5/4)_n^2}\end{align*}
\begin{align*}\left(\dfrac{\G(1/4)}{\G(3/4)}\right)^2&=4+16\sum_{n\ge0}\dfrac{(3/4)_n^2}{(4n+5)(5/4)_n^2}\\
\left(\dfrac{\G(1/4)}{\G(3/4)}\right)^2&=-8+16\sum_{n\ge0}\dfrac{(-1/4)_n^2}{(5/4)_n^2}\\
\left(\dfrac{\G(1/4)}{\G(3/4)}\right)^2&=8\sum_{n\ge0}\dfrac{(1/2)_n^2}{(4n+1)n!^2}\\
\left(\dfrac{\G(3/4)}{\G(1/4)}\right)^2&=-\dfrac{1}{32}\sum_{n\ge0}\dfrac{(-1/2)_n^2}{(4n-1)n!^2}\\
\left(\dfrac{\G(3/4)}{\G(1/4)}\right)^2&=-\dfrac{1}{4}+\sum_{n\ge0}\dfrac{(1/4)_n^2}{(4n+3)(3/4)_n^2}\\
\left(\dfrac{\G(3/4)}{\G(1/4)}\right)^2&=\dfrac{1}{9}+\dfrac{1}{9}\sum_{n\ge0}\dfrac{(1/4)_n^2}{(7/4)_n^2}\end{align*}
\begin{align*}\left(\dfrac{\G(1/4)}{\G(3/4)}\right)^4&=48+1024\sum_{n\ge0}\dfrac{(n+1)(4n+3)(3/4)_n^4}{(4n+5)^3(5/4)_n^4}\\
\left(\dfrac{\G(1/4)}{\G(3/4)}\right)^4&=80-2048\sum_{n\ge0}\dfrac{n(-1/4)_n^4}{(5/4)_n^4}\\
\left(\dfrac{\G(3/4)}{\G(1/4)}\right)^4&=-\dfrac{1}{16}+\dfrac{2}{81}\sum_{n\ge0}\dfrac{(2n+1)(4n+1)(4n+3)(1/4)_n^4}{(7/4)_n^4}\\
\left(\dfrac{\G(3/4)}{\G(1/4)}\right)^4&=\dfrac{1}{16}-\dfrac{4}{81}\sum_{n\ge0}\dfrac{(2n+1)(1/4)_n^4}{(7/4)_n^4}\end{align*}
\begin{align*}\G(1/3)^3&=27\sum_{n\ge0}\dfrac{(2n+1)(2/3)_n}{(3n+1)(3n+2)(4/3)_n}\\
\G(1/3)^3&=\dfrac{81}{4}-\dfrac{81}{2}\sum_{n\ge0}\dfrac{(2/3)_n}{(3n+1)(3n+4)(3n+2)(3n+5)(7/3)_n}\\
\G(1/3)^3&=\dfrac{27}{2}+432\sum_{n\ge0}\dfrac{n!^3}{(27n^2+9n+2)(27n^2+63n+38)(4/3)_n^3}\\
\G(1/3)^3&=\dfrac{81}{4}-\dfrac{81}{16}\sum_{n\ge0}\dfrac{(n+1)^2n!^3}{(3n+7)(7/3)_n^3}\\
\G(1/3)^3&=18\sum_{n\ge0}\dfrac{(2/3)_n}{(3n+1)^2n!}\\
\G(1/3)^3&=\dfrac{77}{4}-\dfrac{8}{343}\sum_{n\ge0}\dfrac{n!^3}{(10/3)_n^3}\\
\dfrac{1}{\G(1/3)^3}&=\dfrac{2}{27}\sum_{n\ge0}\dfrac{(-2/3)_n^3}{n!^3}\\
\dfrac{1}{\G(1/3)^3}&=\dfrac{4}{81}+\dfrac{4}{81}\sum_{n\ge0}\dfrac{(3n+1)(1/3)_n^3}{(n+1)^2n!^3}\\
\G(2/3)^3&=-\dfrac{81}{4}\sum_{n\ge0}\dfrac{(2n+1)(1/3)_n}{(3n-1)(3n+2)(3n+1)(3n+4)(5/3)_n}\\
\G(2/3)^3&=\dfrac{9}{4}+\dfrac{9}{5}\sum_{n\ge0}\dfrac{(1/3)_n}{(3n+2)(3n+4)(8/3)_n}\\
\G(2/3)^3&=\dfrac{9}{4}+\dfrac{9}{2}\sum_{n\ge0}\dfrac{(n+1)n!^3}{(3n+5)^2(5/3)_n^3}\\
\G(2/3)^3&=3\sum_{n\ge0}\dfrac{(-2/3)_n}{(3n-2)^2n!}\\
\G(2/3)^3&=\dfrac{5}{2}-\dfrac{2}{125}\sum_{n\ge0}\dfrac{n!^3}{(8/3)_n^3}\\
\dfrac{1}{\G(2/3)^3}&=-\dfrac{4}{9}\sum_{n\ge0}\dfrac{(2/3)_n^3}{(3n-1)(n+1)n!^3}\\
\dfrac{1}{\G(2/3)^3}&=16\sum_{n\ge0}\dfrac{(2/3)_n^3}{(27n^2-9n+2)(27n^2+45n+20)n!^3}\\    
\dfrac{1}{\G(2/3)^3}&=-\dfrac{8}{27}\sum_{n\ge0}\dfrac{(-4/3)_n^3}{n!^3}\end{align*}
\begin{align*}\dfrac{\G(1/3)^2}{\G(2/3)}&=3\sum_{n\ge0}\dfrac{(2/3)_n}{(3n+1)n!}\\
\dfrac{\G(1/3)^2}{\G(2/3)}&=-3+12\sum_{n\ge0}\dfrac{n!(-1/3)_n}{(4/3)_n^2}\\
\dfrac{\G(1/3)^2}{\G(2/3)}&=-9+18\sum_{n\ge0}(-1)^n\dfrac{(2/3)_n^2}{(9n^2-9n+1)(9n^2+9n+1)(1/3)_n^2}\\
\dfrac{\G(1/3)^2}{\G(2/3)}&=-6\sum_{n\ge0}\dfrac{(-1/3)_n}{(3n-2)n!}\\
\dfrac{\G(1/3)^2}{\G(2/3)}&=24\sum_{n\ge0}\dfrac{(5/3)_n}{(n+1)(n+2)(n+3)(3/2)_n}(3/4)^n\\
\dfrac{\G(1/3)^2}{\G(2/3)}&=\dfrac{9}{2}\sum_{n\ge0}\dfrac{(2/3)_n}{(7/6)_n}2^{-2n}\\
\dfrac{\G(2/3)}{\G(1/3)^2}&=\dfrac{1}{9}\sum_{n\ge0}\dfrac{(2/3)_n^2}{(n+1)n!^2}\\
\dfrac{\G(2/3)}{\G(1/3)^2}&=-\dfrac{1}{3}\sum_{n\ge0}\dfrac{(-2/3)_n^2}{n!(-1/3)_n}\\
\dfrac{\G(2/3)}{\G(1/3)^2}&=\dfrac{1}{3}\sum_{n\ge0}(-1)^n\dfrac{(-2/3)_n^2}{n!^2}\\
\dfrac{\G(2/3)}{\G(1/3)^2}&=2\sum_{n\ge0}(-1)^n\dfrac{(6n+1)(1/3)_n^2}{(9n^2-6n+2)(9n^2+12n+5)n!^2}\end{align*}
\begin{align*}\dfrac{\G(2/3)^2}{\G(1/3)}&=\sum_{n\ge0}\dfrac{(1/3)_n}{(3n+2)n!}\\
\dfrac{\G(2/3)^2}{\G(1/3)}&=1-\dfrac{1}{4}\sum_{n\ge0}\dfrac{n!(1/3)_n}{(5/3)_n^2}\\
\dfrac{\G(2/3)^2}{\G(1/3)}&=-\dfrac{3}{4}+24\sum_{n\ge0}(-1)^n\dfrac{(1/3)_n^2}{(9n^2-9n+4)(9n^2+9n+4)(2/3)_n^2}\\
\dfrac{\G(2/3)^2}{\G(1/3)}&=-\dfrac{1}{2}\sum_{n\ge0}\dfrac{(-2/3)_n}{(3n-1)n!}\\
\dfrac{\G(2/3)^2}{\G(1/3)}&=\sum_{n\ge0}\dfrac{n(1/3)_n}{(3/2)_n}(3/4)^n\\
\dfrac{\G(2/3)^2}{\G(1/3)}&=\dfrac{2}{3}+\dfrac{1}{30}\sum_{n\ge0}\dfrac{(4/3)_n}{(n+1)(n+2)(11/6)_n}2^{-2n}\\
\dfrac{\G(2/3)^2}{\G(1/3)}&=-\dfrac{1}{3}+\sum_{n\ge0}\dfrac{(1/3)_n^3}{(5/6)_n^3}2^{-2n}\\
\dfrac{\G(1/3)}{\G(2/3)^2}&=\dfrac{4}{3}\sum_{n\ge0}\dfrac{(1/3)_n^2}{(n+1)n!^2}\\
  \dfrac{\G(1/3)}{\G(2/3)^2}&=\sum_{n\ge0}\dfrac{(-1/3)_n^2}{n!(1/3)_n}\\
\dfrac{\G(1/3)}{\G(2/3)^2}&=2\sum_{n\ge0}(-1)^n\dfrac{(2/3)_n^2}{n!^2}\\
\dfrac{\G(1/3)}{\G(2/3)^2}&=\dfrac{3}{2}\sum_{n\ge0}(-1)^n\dfrac{(1/3)_n(2/3)_n}{n!^2}2^{-3n}\end{align*}
\begin{align*}
  \dfrac{\G(1/6)\G(1/3)}{\G(1/2)}&=2^{2/3}\dfrac{\G(1/3)^2}{\G(2/3)}=12\sum_{n\ge0}(-1)^n\dfrac{(2/3)_n}{(4/3)_n}\\
  \dfrac{\G(1/6)\G(1/3)}{\G(1/2)}&=2^{2/3}\dfrac{\G(1/3)^2}{\G(2/3)}=6+\dfrac{3}{2}\sum_{n\ge0}\dfrac{n!(5/6)_n}{(3/2)_n(7/3)_n}\\
  \dfrac{\G(1/6)\G(1/3)}{\G(1/2)}&=2^{2/3}\dfrac{\G(1/3)^2}{\G(2/3)}=9\sum_{n\ge0}(-1)^n\dfrac{(2/3)_n}{(7/6)_n}2^{-3n}\\
  \dfrac{\G(1/6)\G(1/3)}{\G(1/2)}&=2^{2/3}\dfrac{\G(1/3)^2}{\G(2/3)}=-8+16\sum_{n\ge0}\dfrac{(-1/6)_n^2}{(7/6)_n^2}\\
  \dfrac{\G(1/6)\G(1/3)}{\G(1/2)}&=2^{2/3}\dfrac{\G(1/3)^2}{\G(2/3)}=8+\dfrac{8}{245}\sum_{n\ge0}\dfrac{(12n+11)n!(2/3)_n^2(4/3)_n}{(3/2)_n(11/6)_n(13/6)_n^2}\\
  \dfrac{\G(1/6)\G(1/3)}{\G(1/2)}&=2^{2/3}\dfrac{\G(1/3)^2}{\G(2/3)}=9-\dfrac{9}{28}\sum_{n\ge0}\dfrac{n!(3/2)_n}{(7/3)_n(13/6)_n}\\
  \dfrac{\G(1/6)\G(1/3)}{\G(1/2)}&=2^{2/3}\dfrac{\G(1/3)^2}{\G(2/3)}=6+\dfrac{12}{7}\sum_{n\ge0}\dfrac{n!(1/2)_n}{(4/3)_n(13/6)_n}\\
  \dfrac{\G(1/6)\G(1/3)}{\G(1/2)}&=2^{2/3}\dfrac{\G(1/3)^2}{\G(2/3)}=6+\dfrac{24}{7}\sum_{n\ge0}(-1)^n\dfrac{(5/3)_n(5/6)_n}{(4/3)_n(13/6)_n}\\
  \dfrac{\G(1/6)\G(1/3)}{\G(1/2)}&=2^{2/3}\dfrac{\G(1/3)^2}{\G(2/3)}=9-\dfrac{9}{14}\sum_{n\ge0}(-1)^n\dfrac{(2/3)_n(5/6)_
    n}{(7/3)_n(13/6)_n}\\
  \dfrac{\G(1/6)\G(1/3)}{\G(1/2)}&=2^{2/3}\dfrac{\G(1/3)^2}{\G(2/3)}=6+\dfrac{15}{4}\sum_{n\ge0}\dfrac{(5/6)_n}{(2n+1)n!}(3/4)^n\\
  \dfrac{\G(1/2)}{\G(1/6)\G(1/3)}&=2^{-2/3}\dfrac{\G(2/3)}{\G(1/3)^2}=\dfrac{1}{10}\sum_{n\ge0}\dfrac{(1/3)_n(1/2)_n}{n!(11/6)_n}\\
  \dfrac{\G(1/2)}{\G(1/6)\G(1/3)}&=2^{-2/3}\dfrac{\G(2/3)}{\G(1/3)^2}=\dfrac{3}{128}\sum_{n\ge0}\dfrac{(12n+5)(1/2)_n(1/6)_n^2(5/6)_n}{n!(4/3)_n(5/3)_n^2}\\
  \dfrac{\G(1/2)}{\G(1/6)\G(1/3)}&=2^{-2/3}\dfrac{\G(2/3)}{\G(1/3)^2}=\dfrac{1}{9}\sum_{n\ge0}\dfrac{(1/6)_n(1/3)_n}{n!(3/2)_n}\\
  \dfrac{\G(1/2)}{\G(1/6)\G(1/3)}&=2^{-2/3}\dfrac{\G(2/3)}{\G(1/3)^2}=\dfrac{1}{6}\sum_{n\ge0}\dfrac{(1/6)_n(-2/3)_n}{n!(1/2)_n}\\
  \dfrac{\G(1/2)}{\G(1/6)\G(1/3)}&=2^{-2/3}\dfrac{\G(2/3)}{\G(1/3)^2}=\dfrac{1}{32}\sum_{n\ge0}\dfrac{(3/2)_n(4/3)_n}{n!(n+1)!}(3/4)^n\end{align*}
  \begin{align*}
    \dfrac{\G(1/2)\G(1/3)^3}{\G(1/6)\G(2/3)^2}&=2^{-2/3}\dfrac{\G(1/3)^2}{\G(2/3)}=-\dfrac{15}{4}\sum_{n\ge0}\dfrac{(-2/3)_n}{(6n-5)n!}\\
    \dfrac{\G(1/2)\G(1/3)^3}{\G(1/6)\G(2/3)^2}&=2^{-2/3}\dfrac{\G(1/3)^2}{\G(2/3)}=3\sum_{n\ge0}\dfrac{(1/3)_n}{(6n+1)n!}\\
    \dfrac{\G(1/6)\G(2/3)^2}{\G(1/2)\G(1/3)^3}&=2^{2/3}\dfrac{\G(2/3)}{\G(1/3)^2}=\dfrac{1}{4}\sum_{n\ge0}\dfrac{(1/2)_n(2/3)_n}{n!(n+1)!}(3/4)^n\\
    \dfrac{\G(1/6)\G(2/3)^2}{\G(1/2)\G(1/3)^3}&=2^{2/3}\dfrac{\G(2/3)}{\G(1/3)^2}=\dfrac{4}{9}\sum_{n\ge0}(-1)^n\dfrac{(2/3)_n(5/3)_n}{n!(n+1)!}\\
   \dfrac{\G(1/6)\G(2/3)^2}{\G(1/2)\G(1/3)^3}&=2^{2/3}\dfrac{\G(2/3)}{\G(1/3)^2}=\dfrac{1}{3}\sum_{n\ge0}(-1)^n\dfrac{(-1/3)_n^2}{n!^2}\end{align*} 
\begin{align*}\dfrac{\G(1/2)\G(2/3)}{\G(1/6)}&=2^{-2/3}\dfrac{\G(2/3)^2}{\G(1/3)}=1-\sum_{n\ge0}(-1)^n\dfrac{(4/3)_n}{(5/3)_n}\\
  \dfrac{\G(1/2)\G(2/3)}{\G(1/6)}&=2^{-2/3}\dfrac{\G(2/3)^2}{\G(1/3)}=1-\dfrac{1}{2}\sum_{n\ge0}\dfrac{n!(1/6)_n}{(3/2)_n(5/3)_n}\\
  \dfrac{\G(1/2)\G(2/3)}{\G(1/6)}&=2^{-2/3}\dfrac{\G(2/3)^2}{\G(1/3)}=-\dfrac{1}{2}+\sum_{n\ge0}(-1)^n\dfrac{(1/3)_n(1/6)_n}{(2/3)_n(5/6)_n}\\
  \dfrac{\G(1/2)\G(2/3)}{\G(1/6)}&=2^{-2/3}\dfrac{\G(2/3)^2}{\G(1/3)}=-\dfrac{3}{8}+\dfrac{3}{4}\sum_{n\ge0}(-1)^n\dfrac{(-2/3)_n(1/6)_n}{(5/3)_n(5/6)_n}\\
  \dfrac{\G(1/2)\G(2/3)}{\G(1/6)}&=2^{-2/3}\dfrac{\G(2/3)^2}{\G(1/3)}=\dfrac{3}{8}+\dfrac{3}{100}\sum_{n\ge0}\dfrac{n!(3/2)_n}{(8/3)_n(11/6)_n}\\
  \dfrac{\G(1/2)\G(2/3)}{\G(1/6)}&=2^{-2/3}\dfrac{\G(2/3)^2}{\G(1/3)}=\dfrac{1}{2}-\dfrac{1}{20}\sum_{n\ge0}\dfrac{n!(1/2)_n}{(5/3)_n(11/6)_n}\\
\dfrac{\G(1/2)\G(2/3)}{\G(1/6)}&=2^{-2/3}\dfrac{\G(2/3)^2}{\G(1/3)}=-1+\dfrac{3}{2}\sum_{n\ge0}(-1)^n\dfrac{(1/3)_n}{(5/6)_n}2^{-3n}\\
\dfrac{\G(1/2)\G(2/3)}{\G(1/6)}&=2^{-2/3}\dfrac{\G(2/3)^2}{\G(1/3)}=\dfrac{32}{75}\sum_{n\ge0}\dfrac{(1/6)_n^2}{(11/6)_n^2}\\
\dfrac{\G(1/2)\G(2/3)}{\G(1/6)}&=2^{-2/3}\dfrac{\G(2/3)^2}{\G(1/3)}=\dfrac{32}{75}+\dfrac{32}{21}\sum_{n\ge0}\dfrac{(12n+13)n!(4/3)_n^2(5/3)}{P(n)P(n+1)(3/2)_n(11/6)_n^2(13/6)_n}\\
&\text{with\quad}P(n)=72n^2+84n+25\\
\dfrac{\G(1/6)}{\G(1/2)\G(2/3)}&=2^{2/3}\dfrac{\G(1/3)}{\G(2/3)^2}=3\sum_{n\ge0}\dfrac{(-1/3)_n(1/2)_n}{n!(7/6)_n}\\
\dfrac{\G(1/6)}{\G(1/2)\G(2/3)}&=2^{2/3}\dfrac{\G(1/3)}{\G(2/3)^2}=\dfrac{8}{3}\sum_{n\ge0}\dfrac{(2/3)_n(-1/6)_n}{n!(3/2)_n}\\
\dfrac{\G(1/6)}{\G(1/2)\G(2/3)}&=2^{2/3}\dfrac{\G(1/3)}{\G(2/3)^2}=2\sum_{n\ge0}\dfrac{(-1/3)_n(-1/6)_n}{n!(1/2)_n}\\
\dfrac{\G(1/6)}{\G(1/2)\G(2/3)}&=2^{2/3}\dfrac{\G(1/3)}{\G(2/3)^2}=\dfrac{225}{8}\sum_{n\ge0}\dfrac{(12n+7)(1/2)_n(5/6)_n^2(7/6)_n}{P(n)P(n+1)n!(4/3)_n^2(5/3)_n}\\
&\text{with\quad}P(n)=72n^2+12n+1\\
\dfrac{\G(1/6)}{\G(1/2)\G(2/3)}&=2^{2/3}\dfrac{\G(1/3)}{\G(2/3)^2}=\dfrac{3}{2}\sum_{n\ge0}\dfrac{(1/2)_n(2/3)_n}{n!^2}(3/4)^n\end{align*}
\begin{align*}\dfrac{\G(1/6)\G(2/3)^3}{\G(1/2)\G(1/3)^2}&=2^{2/3}\dfrac{\G(2/3)^2}{\G(1/3)}=7\sum_{n\ge0}\dfrac{(-1/3)_n}{(6n+5)n!}\\
  \dfrac{\G(1/6)\G(2/3)^3}{\G(1/2)\G(1/3)^2}&=2^{2/3}\dfrac{\G(2/3)^2}{\G(1/3)}=-6\sum_{n\ge0}\dfrac{(2/3)_n}{(6n-1)(6n+5)n!}\\
  \dfrac{\G(1/2)\G(1/3)^2}{\G(1/6)\G(2/3)^3}&=2^{-2/3}\dfrac{\G(1/3)}{\G(2/3)^2}=\sum_{n\ge0}(-1)^n\dfrac{(1/3)_n^2}{n!^2}\\
  \dfrac{\G(1/2)\G(1/3)^2}{\G(1/6)\G(2/3)^3}&=2^{-2/3}\dfrac{\G(1/3)}{\G(2/3)^2}=\dfrac{5}{6}\sum_{n\ge0}(-1)^n\dfrac{(-2/3)_n(1/3)_n}{n!(n+1)!}\\
\dfrac{\G(1/2)\G(1/3)^2}{\G(1/6)\G(2/3)^3}&=2^{-2/3}\dfrac{\G(1/3)}{\G(2/3)^2}=\dfrac{3}{4}\sum_{n\ge0}\dfrac{(1/2)_n(1/3)_n}{n!^2}(3/4)^n\end{align*}
\begin{align*}\left(\dfrac{\G(1/3)}{\G(2/3)}\right)^3&=\dfrac{1}{4}\left(\dfrac{\G(1/6)}{\G(1/2)}\right)^3=6+12\sum_{n\ge0}\dfrac{(3n+2)(2/3)_n^3}{(3n+4)^2(4/3)_n^3}\\
  \left(\dfrac{\G(1/3)}{\G(2/3)}\right)^3&=\dfrac{1}{4}\left(\dfrac{\G(1/6)}{\G(1/2)}\right)^3=8-\dfrac{1}{8}\sum_{n\ge0}\dfrac{(2n+3)(5/3)_n^3}{(n+1)(n+2)(7/3)_n^3}\\
  \left(\dfrac{\G(1/3)}{\G(2/3)}\right)^3&=\dfrac{1}{4}\left(\dfrac{\G(1/6)}{\G(1/2)}\right)^3=-8+16\sum_{n\ge0}\dfrac{(-1/3)_n^3}{(4/3)_n^3}\\
\left(\dfrac{\G(1/3)}{\G(2/3)}\right)^3&=\dfrac{1}{4}\left(\dfrac{\G(1/6)}{\G(1/2)}\right)^3=6\sum_{n\ge0}\dfrac{(2/3)_n^2}{(3n+1)n!^2}\\
\left(\dfrac{\G(1/3)}{\G(2/3)}\right)^3&=\dfrac{1}{4}\left(\dfrac{\G(1/6)}{\G(1/2)}\right)^3=\dfrac{9}{2}+144\sum_{n\ge0}\dfrac{(2n+1)(1/2)_n^3}{(6n+7)^2(7/6)_n^3}\\
\left(\dfrac{\G(2/3)}{\G(1/3)}\right)^3&=4\left(\dfrac{\G(1/2)}{\G(1/6)}\right)^3=-\dfrac{1}{6}-\dfrac{4}{3}\sum_{n\ge0}\dfrac{(1/3)_n^3}{(3n-2)(3n+2)(2/3)_n^3}\\
\left(\dfrac{\G(2/3)}{\G(1/3)}\right)^3&=4\left(\dfrac{\G(1/2)}{\G(1/6)}\right)^3=\dfrac{1}{8}\sum_{n\ge0}\dfrac{(2n+1)(1/3)_n^3}{(5/3)_n^3}\\
\left(\dfrac{\G(2/3)}{\G(1/3)}\right)^3&=4\left(\dfrac{\G(1/2)}{\G(1/6)}\right)^3=-\dfrac{1}{8}+16\sum_{n\ge0}\dfrac{(1/3)_n^3}{(27n^2-27n+8)(27n^2+27n+8)(2/3)_n^3}\\
\left(\dfrac{\G(2/3)}{\G(1/3)}\right)^3&=4\left(\dfrac{\G(1/2)}{\G(1/6)}\right)^3=-\dfrac{6}{5}-\dfrac{64}{9}\sum_{n\ge0}\dfrac{(1/6)_n^3}{(2n+1)(6n-5)(1/2)_n^3}\\
\left(\dfrac{\G(2/3)}{\G(1/3)}\right)^3&=4\left(\dfrac{\G(1/2)}{\G(1/6)}\right)^3=-\dfrac{1}{6}\sum_{n\ge0}\dfrac{(-2/3)_n^2}{(3n-1)n!^2}\end{align*}
\begin{align*}\dfrac{\G(1/3)^4}{\G(2/3)^2}&=-\dfrac{2187}{7}\sum_{n\ge0}\dfrac{(n+1)(n+2)(n+3)(5/3)_n}{(3n-2)(3n+1)(3n+4)(3n+7)(12/6)_n}2^{-2n}\\
  \dfrac{\G(1/3)^2}{\G(2/3)^4}&=\dfrac{9}{4}\sum_{n\ge0}(-1)^n\dfrac{(1/2)_n(1/3)_n(2/3)_n}{n!^3}(3/4)^{2n}\\
  \dfrac{\G(1/3)^2}{\G(2/3)^4}&=2\sum_{n\ge0}\dfrac{(1/3)_n^3}{n!^3}\\
  \dfrac{\G(2/3)^2}{\G(1/3)^4}&=\dfrac{1}{27}\sum_{n\ge0}\dfrac{(-1/3)_n^3}{n!^3}\\
\dfrac{\G(1/6)^2\G(1/3)^2}{\G(1/2)^2}&=2^{4/3}\dfrac{\G(1/3)^4}{\G(2/3)^2}=48\sum_{n\ge0}\dfrac{(2/3)_n}{(2n+1)(4/3)_n}\\
\dfrac{\G(1/6)^2\G(1/3)^2}{\G(1/2)^2}&=2^{4/3}\dfrac{\G(1/3)^4}{\G(2/3)^2}=48+192\sum_{n\ge0}\dfrac{(3n+2)n!^2(5/6)_n^2}{(12n^2+4n+1)(12n^2+28n+17)(3/2)_n^2(4/3)_n^2}\\
\dfrac{\G(1/6)^2\G(1/3)^2}{\G(1/2)^2}&=2^{4/3}\dfrac{\G(1/3)^4}{\G(2/3)^2}=72\sum_{n\ge0}(-1)^n\dfrac{(2/3)_n}{(3n+1)(7/6)_n}2^{-3n}\\
\dfrac{\G(1/2)^2}{\G(1/6)^2\G(1/3)^2}&=2^{-4/3}\dfrac{\G(2/3)^2}{\G(1/3)^4}=\dfrac{1}{24}\sum_{n\ge0}\dfrac{(6n+1)(1/2)_n^2(1/3)_n^2}{(6n^2-4n+1)(6n^2+8n+3)n!^2(5/6)_n^2}\\
\dfrac{\G(1/2)^2\G(2/3)^2}{\G(1/6)^2}&=2^{-10/3}\dfrac{\G(2/3)^4}{\G(1/3)^2}=\dfrac{1}{6}\sum_{n\ge0}\dfrac{(1/3)_n}{(2n+1)(5/3)_n}\\
\dfrac{\G(1/2)^2\G(2/3)^2}{\G(1/6)^2}&=2^{-10/3}\dfrac{\G(2/3)^4}{\G(1/3)^2}=\dfrac{1}{6}+\dfrac{1}{24}\sum_{n\ge0}\dfrac{(6n+5)n!^2(7/6)_n^2}{(6n^2+4n+1)(6n^2+16n+11)(3/2)_n^2(5/3)_n^2}\\
\dfrac{\G(1/2)^2\G(2/3)^2}{\G(1/6)^2}&=2^{-10/3}\dfrac{\G(2/3)^4}{\G(1/3)^2}=-\dfrac{1}{12}-\dfrac{1}{4}\sum_{n\ge0}(-1)^n\dfrac{(4/3)_n}{(3n-1)(3n+2)(5/6)_n}2^{-3n}\\
\dfrac{\G(1/6)^2}{\G(1/2)^2\G(2/3)^2}&=2^{10/3}\dfrac{\G(1/3)^2}{\G(2/3)^4}=48\sum_{n\ge0}\dfrac{(3n+1)(1/2)_n^2(2/3)_n^2}{(12n^2-4n+1)(12n^2+20n+9)n!^2(7/6)_n^2}\end{align*}
\begin{align*}\dfrac{\G(1/3)^5}{\G(2/3)^4}&=18\sum_{n\ge0}\dfrac{(2/3)_n^2}{(4/3)_n^2}\\
    \dfrac{\G(1/3)^5}{\G(2/3)^4}&=18+72\sum_{n\ge0}\dfrac{(12n+7)n!(2/3)_n(5/6)_n^2}{(18n^2+3n+1)(18n^2+39n+22)(3/2)_n(4/3)_n^2(7/6)_n}\\
  \dfrac{\G(2/3)^4}{\G(1/3)^5}&=\dfrac{2}{3}\sum_{n\ge0}\dfrac{(12n+1)(1/6)_n(1/2)_n(1/3)_n^2}{(18n^2-15n+4)(18n^2+21n+7)n!(2/3)_n(5/6)_n^2}\\
\dfrac{\G(2/3)^5}{\G(1/3)^4}&=\dfrac{1}{12}\sum_{n\ge0}\dfrac{(1/3)_n^2}{(5/3)_n^2}\\
    \dfrac{\G(2/3)^5}{\G(1/3)^4}&=\dfrac{1}{12}+\dfrac{1}{15}\sum_{n\ge0}\dfrac{(12n+11)n!(4/3)_n(7/6)_n^2}{(18n^2+15n+4)(18n^2+51n+37)(3/2)_n(5/3)_n^2(11/6)_n}\\
    \dfrac{\G(1/3)^4}{\G(2/3)^5}&=36\sum_{n\ge0}\dfrac{(12n+5)(5/6)_n(1/2)_n(2/3)_n^2}{(18n^2-3n+1)(18n^2+33n+16)n!(4/3)_n(7/6)_n^2}\\
\left(\dfrac{\G(1/3)}{\G(2/3)}\right)^6&=\dfrac{1}{2}\left(\dfrac{\G(1/6)}{\G(5/6)}\right)^3=-322+\dfrac{131072}{343}\sum_{n\ge0}\dfrac{(-1/6)_n^3}{(13/6)_n^3}\\
  \left(\dfrac{\G(2/3)}{\G(1/3)}\right)^6&=2\left(\dfrac{\G(5/6)}{\G(1/6)}\right)^3=\dfrac{5}{54}-\dfrac{256}{3375}\sum_{n\ge0}\dfrac{(1/6)_n^3}{(11/6)_n^3}
\end{align*}
\begin{align*}
  \dfrac{I_0(1)}{e}&=\sum_{n\ge0}(-1)^n\dfrac{(1/2)_n}{n!^2}2^n\\
  e\cdot I_0(1)&=\sum_{n\ge0}\dfrac{(1/2)_n}{n!^2}2^n\\
  \cos(1)J_0(1)&=\sum_{n\ge0}(-1)^n\dfrac{(1/4)_n(3/4)_n}{n!^2(1/2)_n^2}\\
  \sin(1)J_0(1)&=\sum_{n\ge0}(-1)^n\dfrac{(3/4)_n(5/4)_n}{n!^2(3/2)_n^2}
\end{align*}

\medskip

\section{Hypergeometric Series: Functions}

\begin{align*}2^z&=1-\dfrac{z}{z-2}\sum_{n\ge0}\dfrac{(1/2)_n(1-z)_n}{((3-z)/2)_n((4-z)/2)_n}\\
2^z&=1+\dfrac{z}{2}\sum_{n\ge0}\dfrac{((z+1)/2)_n(z/2+1)_n}{(3/2)_n(z+2)_n}\\
2^{-z}&=\dfrac{2}{z+1}+\dfrac{z-1}{z+1}\sum_{n\ge0}\dfrac{(-z/2)_n(-(z+1)/2)_n}{(1/2)_n(1-z)_n}\\
2^{-z}&=\dfrac{1-z}{2}+\dfrac{(1+z)}{2}\sum_{n\ge0}\dfrac{(z)_n(-1/2)_n}{((z+1)/2)_n(z/2+1)_n}\end{align*}
$$\log(1+z)=\dfrac{z(z+2)}{2(z+1)}\sum_{n\ge0}\dfrac{n!}{(3/2)_n}(-z^2/(4(z+1)))^n$$
\begin{align*}\dfrac{\sin(\pi z)}{\pi z}&=\sum_{n\ge0}\dfrac{(-z)_n(z)_n}{n!^2}\\
  \dfrac{\pi z}{\sin(\pi z)}&=1+z^2\sum_{n\ge1}\dfrac{(n-1)!^2}{(1-z)_n(1+z)_n}\\
\cos(\pi z)&=\sum_{n\ge0}\dfrac{(-z)_n(z)_n}{n!(1/2)_n}\\
\cos(\pi z)&=1-4z^2\sum_{n\ge0}\dfrac{(1/2-z)_n(1/2+z)_n}{(3/2)_n^2}\\
\cos(\pi z)&=\sum_{n\ge0}\dfrac{(-3z/2)_n(3z/2)_n}{n!(1/2)_n}(3/4)^n\\
\cos(\pi z)&=\sum_{n\ge0}\dfrac{(-2z)_n(2z)_n}{n!(1/2)_n}2^{-n}\\
\cos(\pi z)&=\sum_{n\ge0}\dfrac{(-3z)_n(3z)_n}{n!(1/2)_n}2^{-2n}\\
\dfrac{1}{\cos(\pi z)}&=1+4z^2\sum_{n\ge1}\dfrac{(-1/2)_n^2}{(1/2-z)_n(1/2+z)_n}\end{align*}
\begin{align*}\dfrac{\sin(\pi z/2)}{z}&=\sum_{n\ge0}\dfrac{((1-z)/2)_n((1+z)/2)_n}{n!(3/2)_n}\\
\dfrac{\sin(\pi z/2)}{z}&=\dfrac{1}{z}\sum_{n\ge0}\dfrac{(1-z)_n(z-1)_n}{n!(1/2)_n}2^{-n}\\
\dfrac{\sin(\pi z/2)}{z}&=\dfrac{3}{2}\sum_{n\ge0}\dfrac{((1-3z)/2)_n((1+3z)/2)_n}{n!(3/2)_n}2^{-2n}\\
\dfrac{\sin(\pi z/2)}{z\sqrt{2}}&=\sum_{n\ge0}\dfrac{(1/2-z)_n(1/2+z)_n}{n!(3/2)_n}2^{-n}\\
\dfrac{z\sqrt{2}}{\sin(\pi z/2)}&=\sum_{n\ge0}(-1)^n\dfrac{(6n+1)(1/2)_n(1/2-z)_n(1/2+z)_n}{n!(1-z/2)_n(1+z/2)_n}2^{-3n}\\
\dfrac{\sin(\pi z/3)}{z\sqrt{3}/2}&=\sum_{n\ge0}\dfrac{((1-z)/2)_n((1+z)/2)_n}{n!(3/2)_n}(3/4)^n\\
\dfrac{\sin(\pi z/3)}{z\sqrt{3}/2}&=\sum_{n\ge0}\dfrac{(1-z)_n(1+z)_n}{n!(3/2)_n}2^{-2n}\end{align*}
\begin{align*}\dfrac{\pi}{\sin(\pi z)}&=\sum_{n\ge0}\dfrac{(z)_n}{(n+z)n!}\\
\dfrac{\pi}{\cos(\pi z)}&=2\sum_{n\ge0}\dfrac{(1/2-z)_n}{(2n+1-2z)n!}\end{align*}
\begin{align*}\cos(\sqrt{2})&=\dfrac{1}{3}\sum_{n\ge0}(-1)^n\dfrac{1}{(n+2)!(5/2)_n}2^{-n}\\
\dfrac{\sin(\sqrt{2})}{\sqrt{2}}&=\sum_{n\ge0}(-1)^n\dfrac{1}{n!(3/2)_n}2^{-n}\end{align*}
\begin{align*}z\cot(\pi z)&=\dfrac{1}{4}+\dfrac{1-16z^2}{16(1-z^2)}\sum_{n\ge0}\dfrac{(1/2-z)_n(1/2+z)_n}{(2-z)_n(2+z)_n}\\
  \dfrac{\tan(\pi z)}{z}&=\dfrac{1}{4z^2}-\dfrac{1-16z^2}{4z^2(1-4z^2)}\sum_{n\ge0}\dfrac{(z)_n(-z)_n}{(3/2+z)_n(3/2-z)_n}\\
\pi z\cot(\pi z)&=1-\dfrac{z^2}{1-z^2}\sum_{n\ge0}\dfrac{n!^2(1/2-z)_n(1/2+z)_n(4n+3)}{(3/2)_n^2(2-z)_n(2+z)_n}\\
  \dfrac{\tan(\pi z)}{\pi z}&=\dfrac{1}{1-4z^2}\sum_{n\ge0}\dfrac{(4n+1)(1/2)_n^2(z)_n(-z)_n}{(3/2-z)_n(3/2+z)_nn!^2}\end{align*}
\begin{align*}\sinh(a\asinh(z))&=az\sum_{n\ge0}(-1)^n\dfrac{(1/2-a/2)_n(1/2+a/2)_n}{n!(3/2)_n}z^{2n}\\
\cosh(a\asinh(z))&=\sum_{n\ge0}(-1)^n\dfrac{(-a/2)_n(a/2)_n}{n!(1/2)_n}z^{2n}\end{align*}
\begin{align*}\dfrac{\asinh(z)}{\sqrt{1+z^2}}&=z-2z^3\sum_{n\ge0}(-1)^n\dfrac{n!}{((2(z^2+1)n+1)(2(z^2+1)(n+1)+1)(3/2)_n}z^{2n}\\
\atanh(z)&=\dfrac{z}{1-z^2}\sum_{n\ge0}\dfrac{n!}{(3/2)_n}(z^2/(z^2-1))^n\\
\atanh(z)&=\dfrac{z}{1-z^2}\sum_{n\ge0}(-1)^n\dfrac{(1/2)_n}{(2n+1)n!}(2z/(1-z^2))^{2n}\\
\asinh^2(z)&=z^2\sum_{n\ge0}(-1)^n\dfrac{n!}{(n+1)(3/2)_n}z^{2n}\\
\acosh^2(z)&=(z^2-1)\sum_{n\ge0}(-1)^n\dfrac{n!}{(n+1)(3/2)_n}(z^2-1)^n\\
\atanh^2(z)&=\dfrac{z^2}{(z^2-1)^2}\sum_{n\ge0}(-1)^n\dfrac{n!}{(n+1)(3/2)_n}\left(\dfrac{2z}{1-z^2}\right)^{2n}\end{align*}
\begin{align*}\sqrt{\pi}\dfrac{\G(z+1)}{\G(z+1/2)}&=2\sum_{n\ge0}(-1)^n\dfrac{(1-z)_n}{(1+z)_n}\\
\sqrt{\pi}\dfrac{\G(z+1)}{\G(z+1/2)}&=2+2z(z-1)\sum_{n\ge0}\dfrac{n!(z+1/2)_n}{(2n+z)(2n+z+2)(3/2)_n(z+1)_n}\\
\sqrt{\pi}\dfrac{\G(z+1)}{\G(z+1/2)}&=z\sum_{n\ge0}\dfrac{(1/2)_n}{(n+z)n!}\\
\sqrt{\pi}\dfrac{\G(z+1)}{\G(z+1/2)}&=\sum_{n\ge0}\dfrac{(2z)_n}{(z+1)_n}2^{-n}\\
\dfrac{\G(z+1/2)}{\sqrt{\pi}\G(z+1)}&=(z-1)\sum_{n\ge0}\dfrac{(1/2)_n(z)_n}{(2n+z-1)(2n+z+1)n!(z+1/2)_n}\end{align*}
$$\sqrt{\pi}\dfrac{\G(z)}{\G(z+1/2)}+\psi(z)+\gamma=2\sum_{n\ge0}(-1)^n\dfrac{(1-z)_n}{(n+1)(1+z)_n}$$
\begin{align*}\dfrac{\G(z+1)^2}{\G(2z)}&=z^2\sum_{n\ge0}\dfrac{(1-z)_n}{(n+z)n!}\\
\dfrac{\G(z+1)^2}{\G(2z)}&=1-\dfrac{(z-1)^2}{(z+1)^2}\sum_{n\ge0}\dfrac{(n+1)!(2z)_n}{(z+2)_n^2}\\
\dfrac{\G(z+1)^2}{\G(2z)}&=-1+2z^4\sum_{n\ge0}(-1)^n\dfrac{(1-z)_n^2}{(n^2-n+z^2)(n^2+n+z^2)(z)_n^2}\\
\dfrac{\G(2z)}{\G(z+1)^2}&=\dfrac{2z-1}{z^2}\sum_{n\ge0}\dfrac{(z-1)_n^2}{n!(2z-1)_n}\\
\dfrac{\G(2z)}{\G(z+1)^2}&=\sum_{n\ge0}\dfrac{(1-z)_n^2}{(n+1)n!^2}
\end{align*}
$$\dfrac{\pi}{\sin(\pi z)}\dfrac{\G(z+1)^2}{\G(2z)}=2z^2\sum_{n\ge0}\dfrac{(2z)_n}{(n+z)^2n!}$$
\begin{align*}z\left(\dfrac{\G(z)}{\G(z+1/2)}\right)^2&=\dfrac{4z}{4z-1}-\dfrac{1}{z}\sum_{n\ge0}\dfrac{(z+1/2)_n^2}{(4n+4z-1)(4n+4z+3)(z+1)_n^2}\\
z\left(\dfrac{\G(z)}{\G(z+1/2)}\right)^2&=1+\dfrac{1}{4z}\sum_{n\ge0}\dfrac{(z+1/2)_n^2}{(n+z+1)(z+1)_n^2}\\
\left(\dfrac{\G(z+1/2)}{\G(z)}\right)^2&=\dfrac{(2z-1)^2}{4z-3}-\sum_{n\ge0}\dfrac{(z)_n^2}{(4n+4z-3)(4n+4z+1)(z+1/2)_n^2}\\
\left(\dfrac{\G(z+1/2)}{\G(z)}\right)^2&=\dfrac{2z-1}{2}+\dfrac{1}{2}\sum_{n\ge0}\dfrac{(z)_n^2}{(2n+2z+1)(z+1/2)_n^2}\end{align*}
$$\pi z\left(\dfrac{\G(z)}{\G(z+1/2)}\right)^2=4\sum_{n\ge0}\dfrac{(1-z)_n}{(2n+1)(1+z)_n}$$
\begin{align*}&\dfrac{\sqrt{\pi}}{\G((3+z)/4)\G((3-z)/4)}=\sum_{n\ge0}\dfrac{((1-z)/2)_n((1+z)/2)_n}{n!^2}2^{-n}\\
&\dfrac{\G((3+z)/4)\G((3-z)/4)}{\sqrt{\pi}}=\\&8(1-z^2)\sum_{n\ge0}(-1)^n\dfrac{(3n+1)((3-z)/2)_n((3+z)/2)_n}{(12n^2-4n+1-z^2)(12n^2+20n+9-z^2)((3-z)/4)_n((3+z)/4)_n}2^{-3n}\end{align*}
\begin{align*}z(z-1/2)\left(\dfrac{\G(z)}{\G(z+1/2)}\right)^4&=1-\dfrac{2z-1}{16z^3}\sum_{n\ge0}\dfrac{(4n+4z+1)(z+1/2)_n^4}{(2n+2z-1)(n+z+1)(z+1)_n^4}\\
  \dfrac{1}{z(z-1/2)}\left(\dfrac{\G(z+1/2)}{\G(z)}\right)^4&=\dfrac{(2z-1)^2}{4z(z-1)}-\dfrac{1}{4z(2z-1)}\sum_{n\ge0}\dfrac{(4n+4z-1)(z)_n^4}{(n+z-1)(2n+2z+1)(z+1/2)_n^4}\end{align*}
$$z^2\left(\dfrac{\G(z)}{\G(2z)}\right)^4\G(4z)=(4z-1)\sum_{n\ge0}\dfrac{(-z)_n^2}{(z)_n^2}$$
$$\dfrac{\G(a)\G(b)}{\G(a+b)}=B(a,b)=\sum_{n\ge0}\dfrac{(1-b)_n}{(n+a)n!}$$
\begin{align*}2&\dfrac{\G((z+a+1)/2)\G((z-a)/2+1)}{\G((z+a)/2)\G((z-a+1)/2)}\\
    &=z+\dfrac{2a(1-a)}{(z+1)(z+2)+a(1-a)}\sum_{n\ge0}\dfrac{((z+a)/2)_n((z+1-a)/2)_n}{((z+a+3)/2)_n((z-a)/2+2)_n}\\
  &\dfrac{\G((z+a)/2)\G((z-a+1)/2)}{2\G((z+a+1)/2)\G((z-a)/2+1)}\\
  &=\dfrac{z-1}{(z-a)(z+a-1)}+\dfrac{2a(1-a)}{(z^2-a^2)(z^2-(1-a)^2)}\sum_{n\ge0}\dfrac{((z+a-1)/2)_n((z-a)/2)_n}{((z+a)/2+1)_n((z-a+3)/2)_n}\end{align*}
\begin{align*}&\dfrac{\G((z+a+1)/4)\G((z-a+1)/4)}{\G((z+a+3)/4)\G((z-a+3)/4)}\\
    &=\dfrac{4}{z}+\dfrac{a^2-1}{(z+1-a)(z+1+a)}\sum_{n\ge0}\dfrac{((z-a+3)/4)_n((z+a+3)/4)_n}{(n+z/4)(n+z/4+1)((z-a+5)/4)_n((z+a+5)/4)_n}\\
  &\dfrac{\G((z+a+3)/4)\G((z-a+3)/4)}{\G((z+a+1)/4)\G((z-a+1)/4)}\\
  &=\dfrac{(z+a-1)(z-a-1)}{4(z-2)}+(a^2-1)\sum_{n\ge0}\dfrac{((z-a+1)/4)_n((z+a+1)/4)_n}{(4n+z-2)(4n+z+2)((z-a+3)/4)_n((z+a+3)/4)_n}\end{align*}
$$\dfrac{\G(2a+b+1)\G(a+1)\G(b-a+1)}{\G(2a+1)\G(b+1)\G(1-a)\G(b+a+1)}=\sum_{n\ge0}\dfrac{(2a)_n^2(-b)_n}{(2a+b+1)_nn!^2}$$
$$\dfrac{\G(2a+b+1)\G(a+1)^2\G(b+1)}{\G(2a+1)\G(a+b+1)^2}=1-\dfrac{2a^2b}{2a+b+1}\sum_{n\ge0}\dfrac{(1-b)_n(2a+1)_n}{(n+a+1)(2a+b+2)_nn!}$$
\begin{align*}\psi(1+z)+\ga&=z\sum_{n\ge0}\dfrac{(1-z)_n}{(n+1)(n+1)!}\\
\psi(1+z)+\ga&=\dfrac{2z}{z+1}\sum_{n\ge0}(-1)^n\dfrac{(1-z)_n}{(n+1)(2+z)_n}\\
\psi(1+z)-\psi(1+z/3)&=\dfrac{2z}{3(z+1)}\sum_{n\ge0}\dfrac{(1+2z/3)_n}{(n+1)(z+2)_n}\\
\psi(b+1)-\psi(a+1)&=\dfrac{b-a}{a+1}\sum_{n\ge0}\dfrac{(a-b+1)_n}{(n+1)(a+2)_n}\\
\psi(a+1)+\psi(b+1)-\psi(a+b+1)+\gamma&=2\sum_{n\ge1}\dfrac{(-a)_n(-b)_n}{n(a+1)_n(b+1)_n}\end{align*}
Recall that $\be(z)=(\psi((1+z)/2)-\psi(z/2))/2$.
\begin{align*}\be(z)&=\dfrac{1}{2z+2}\sum_{n\ge0}\dfrac{(3/2)_n}{(n+1)((z+3)/2)_n}\\
\be(z)&=\dfrac{1}{2z}\sum_{n\ge0}\dfrac{n!}{(z+1)_n}2^{-n}\\
\be(z)&=\dfrac{1}{2z}+\sum_{n\ge0}\dfrac{n!}{(n+2z)(n+2z+1)(2z)_n}2^{-n}\\
\be(z)-\log(2)&=-1+(z^2-1)\sum_{n\ge0}\dfrac{((3-z)/2)_n}{(n+1)((2n+1)z-1)((2n+3)z-1)((z+1)/2)_n}\\
\be(z)-\log(2)&=\dfrac{1-z}{1+z}\sum_{n\ge0}\dfrac{((3-z)/2)_n}{(n+1)((3+z)/2)_n}\\
\be(z)+\log(2)&=\sum_{n\ge0}\dfrac{(1-z/2)_n}{(n+1)(1+z/2)_n}\\
\be(z)+\log(2)&=\dfrac{1}{z}+\dfrac{z}{z+2}\sum_{n\ge0}\dfrac{(1-z/2)_n}{(n+1)(2+z/2)_n}\end{align*}
\begin{align*}\psi'(z)&=\dfrac{1}{z}\sum_{n\ge0}\dfrac{n!}{(n+1)(z+1)_n}\\
\psi'(z+1)-\dfrac{\pi^2}{6}&=-\dfrac{2z}{z+1}\sum_{n\ge0}\dfrac{(1-z)_n}{(n+1)^2(2+z)_n}\end{align*}
\begin{align*}\sum_{k=1}^\infty\dfrac{1}{k^2-a}&=\dfrac{1}{2(1-a)}\sum_{n\ge0}\dfrac{(4n+3)n!^2(1/2+\sqrt{a})_n(1/2-\sqrt{a})_n}{(2+\sqrt{a})_n(2-\sqrt{a})_n(3/2)_n^2}\\
\sum_{k=1}^\infty\dfrac{1}{k^2-a}&=2-\dfrac{4a-1}{6(a-1)}\sum_{n\ge0}\dfrac{(3n+4)n!(2+2\sqrt{a})_n(2-2\sqrt{a})_n}{(n+1)(n+2)(2+\sqrt{a})_n(2-\sqrt{a})_n(5/2)_n}2^{-2n}\\
\sum_{k=1}^\infty\dfrac{(-1)^{k-1}}{k^2-a}&=\dfrac{1}{2-2a}\sum_{n\ge0}\dfrac{n!^2}{(2+\sqrt{a})_n(2-\sqrt{a})_n}\end{align*}
$$\psi'(z)-\dfrac{\psi'((z+1)/2)}{2}=\dfrac{1}{2z(z+1)}\sum_{n\ge0}\dfrac{(3n+z+2)n!^2}{(n+z)(z/2+1)_n(z/2+3/2)_n}2^{-2n}$$

\smallskip

For the next series involving Bernoulli and similar numbers, the formal series
on the LHS is to be interpreted as a CF, see \cite{Coh1} and \cite{Coh2} for
the precise meaning.

\smallskip

\begin{align*}\sum_{k\ge0}B_kz^k&=z\sum_{n\ge1}\dfrac{1}{(nz+1)^2}\\
\sum_{k\ge0}B_kz^k&=-\dfrac{z}{2}+1+\dfrac{z^3}{2}\sum_{n\ge1}\dfrac{1}{(nz+1)^2((n-1)z+1)^2}\\
\sum_{k\ge0}(B_k-B_{k+2})z^k&=-\dfrac{z}{2}+\dfrac{5}{6}+z^3(1-z^2)\sum_{n\ge0}\dfrac{1}{((n+2)z+1)((n+1)z+1)^2(nz+1)^2((n-1)z+1)}\\
\sum_{k\ge0}(2^{k-1}-1)B_kz^k&=-z\sum_{n\ge1}\dfrac{1}{((2n-1)z+1)^2}\\
\sum_{k\ge0}G_kz^k&=2z\sum_{n\ge1}\dfrac{(-1)^{n+1}}{(nz+1)^2}\\
\sum_{k\ge0}T_kz^k&=2z\sum_{n\ge1}\dfrac{(-1)^{n+1}}{(2nz+1)((2n-2)z+1)}=1-2\sum_{n\ge1}\dfrac{(-1)^{n+1}}{2nz+1}\\
\sum_{k\ge0}(T_k-T_{k+2})z^k&=\dfrac{z^2+z-1}{z^2}+2\dfrac{z^2-1}{z^2}\sum_{n\ge1}\dfrac{(-1)^n}{2nz+1}\\
\sum_{k\ge0}(4T_k-T_{k+2})z^k&=12z(4z^2-1)\sum_{n\ge1}\dfrac{(-1)^n}{((2n+2)z+1)(2nz+1)((2n-2)z+1)((2n-4)z+1)}\end{align*}
\begin{align*}
  \exp\left(\sum_{k\ge1}\dfrac{T_k}{k}z^k\right)&=z+1+\dfrac{4z^3}{(4z+1)^2}\sum_{n\ge0}\dfrac{((2z+1)/(4z))_n^2}{((8z+1)/(4z))_n^2}\\
\exp\left(-\sum_{k\ge1}\dfrac{T_k}{k}z^k\right)&=1-z+\dfrac{z}{4}\sum_{n\ge0}\dfrac{(1/(4z))_n^2}{((2z+1)/(4z))_{n+1}^2}\;.\end{align*}
\begin{align*}\sum_{k\ge0}(k+1)B_kz^k&=2z\sum_{n\ge1}\dfrac{1}{(nz+1)^3}\\
\sum_{k\ge0}E_kz^k&=2\sum_{n\ge1}\dfrac{(-1)^{n+1}}{(2n-1)z+1}\\
\sum_{k\ge0}(E_k-E_{k+2})z^k&=4(1-z^2)\sum_{n\ge1}\dfrac{(-1)^{n+1}}{((2n+1)z+1)((2n-1)z+1)((2n-3)z+1)}\\
\sum_{k\ge0}(k+1)E_kz^k&=2\sum_{n\ge1}\dfrac{(-1)^{n+1}}{((2n-1)z+1)^2}\end{align*}
\begin{align*}
    \exp\left(-\sum_{k\ge1}\dfrac{E_{k}}{k}z^{k}\right)&=1+\dfrac{4z^3}{(3z+1)^2}\sum_{n\ge0}\dfrac{((z+1)/(4z))_n^2}{((7z+1)/(4z))_n^2}\\
  \exp\left(\sum_{k\ge1}\dfrac{E_{k}}{k}z^{k}\right)&=\dfrac{1-2z}{(z-1)^2}+\dfrac{z}{4(z-1)^2}\sum_{n\ge0}\dfrac{((1-z)/(4z))_n^2}{((z+1)/(4z))_{n+1}^2}\;.\end{align*}
\begin{align*}\sum_{k\ge0}(E_k+T_k)z^k&=1+2z\sum_{n\ge1}\dfrac{(-1)^{n+1}}{(2nz+1)((2n-1)z+1)}\\
\sum_{k\ge0}(-1)^{k(k-1)/2}S_kz^k&=2\sum_{n\ge1}\dfrac{(-1)^{n+1}}{(4n-3)z+1}\end{align*}
\begin{align*}&\sum_{k\ge0}(-1)^k(S_{2k+1}+S_{2k+3})z^{2k}\\&=24(1-z^2)\sum_{n\ge1}\dfrac{(-1)^{n+1}}{((4n+1)z+1)((4n-1)z+1)((4n-3)z+1)((4n-5)z+1)}\end{align*}
\begin{align*}\exp\left(\sum_{k\ge1}(-1)^{k-1}\dfrac{S_{2k}}{2k}z^{2k}\right)&=1+\dfrac{24z^3}{(5z+1)(7z+1)}\cdot\\&\phantom{=}\cdot\sum_{n\ge0}\dfrac{((z+1)/(8z))_n((3z+1)/(8z))_n}{((13z+1)/(8z))_n((15z+1)/(8z))_n}\\
    \exp\left(-\sum_{k\ge1}(-1)^{k-1}\dfrac{S_{2k}}{2k}z^{2k}\right)&=-\dfrac{4z-1}{(z-1)(3z-1)}+\dfrac{24z^3}{(z+1)(3z+1)(z-1)(3z-1)}\cdot\\&\phantom{=}\cdot\sum_{n\ge0}\dfrac{((1-3z)/(8z))_n((1-z)/(8z))_n}{((1+9z)/(8z))_n((1+11z)/(8z))_n}\end{align*}
\begin{align*}\exp\left(\sum_{k\ge1}(-1)^{k-1}\dfrac{S_{2k-1}}{2k-1}z^{2k-1}\right)&=1+\dfrac{8z^2}{(z+1)(7z+1)}\cdot\\&\phantom{=}\cdot\sum_{n\ge0}\dfrac{((3z+1)/(8z))_n((5z+1)/(8z))_n}{((9z+1)/(8z))_n((15z+1)/(8z))_n}\\
   \exp\left(-\sum_{k\ge1}(-1)^{k-1}\dfrac{S_{2k-1}}{2k-1}z^{2k-1}\right)&=\dfrac{(3z-1)}{z-1}-\dfrac{8z^2}{(z-1)(3z+1)}\cdot\\&\phantom{=}\cdot\sum_{n\ge0}\dfrac{((1-z)/(8z))_n((1+z)/(8z))_n}{((5z+1)/(8z))_n((11z+1)/(8z))_n}\end{align*}
\begin{align*}\dfrac{\gac(a,z)}{z^ae^{-z}}&=\sum_{n\ge0}\dfrac{z^n}{(a)_{n+1}}\\
\dfrac{\gac(a,z)}{z^a}&=\sum_{n\ge0}(-1)^n\dfrac{z^n}{n!(n+a)}\\
\dfrac{\gac(z,z)}{z^{z-1}e^{-z}}&=z+1-\dfrac{z^2}{z+1}\sum_{n\ge0}\dfrac{1}{(n+1)(n+2)(z+2)_n}z^n\\
\sqrt{\pi}e^{z^2}\erf(z)&=\sum_{n\ge0}\dfrac{z^{2n+1}}{(1/2)_{n+1}}\\
e^{-z^2}\int_0^ze^{t^2}\,dt&=\dfrac{1}{2}\sum_{n\ge0}(-1)^n\dfrac{z^{2n+1}}{(1/2)_{n+1}}\\
\sqrt{\pi}\erf(z)&=2\sum_{n\ge0}(-1)^n\dfrac{z^{2n+1}}{(2n+1)n!}\\
\erf(1/\sqrt{2})\sqrt{\pi e/2}&=\sum_{n\ge1}\dfrac{2^{-n}}{(1/2)_n}\end{align*}
\begin{align*}\dfrac{2}{\pi}K(k)&=\sum_{n\ge0}\dfrac{(1/2)_n^2}{n!^2}k^{2n}\\
\dfrac{2}{\pi}K(k)&=\sum_{n\ge0}\dfrac{(1/4)_n^2}{n!^2}(4k^2(1-k^2))^n\\
\dfrac{2}{\pi}E(k)&=(1-k^2)\sum_{n\ge0}\dfrac{(2n+1)(1/2)_n^2}{n!^2}k^{2n}\\
\dfrac{2}{\pi}E(k)&=-\sum_{n\ge0}\dfrac{(2n-1)(-1/2)_n^2}{n!^2}k^{2n}\\
\dfrac{4}{\pi}D(k)&=\sum_{n\ge0}\dfrac{(2n+1)(1/2)_n^2}{(n+1)n!^2}k^{2n}\\
\dfrac{2}{3\pi}((4-2k^2)E(k)-(1-k^2)K(k))&=3\sum_{n\ge0}\dfrac{(2n-1)(-1/2)_n^2}{(2n-2)n!^2}k^{2n}\\
\dfrac{4}{\pi^2}K(k)^2&=\sum_{n\ge0}\dfrac{(1/2)_n^3}{n!^3}(4k^2(1-k^2))^n\\
\dfrac{4}{\pi^2}K(k)E(k)&=\sum_{n\ge0}\dfrac{(n+1-(2n+1)k^2)(1/2)_n^3}{n!^3}(4k^2(1-k^2))^n\\
\end{align*}
$$\int_0^\infty\dfrac{e^{-zt}}{\cosh^k(t)}\,dt=\dfrac{2}{z+k}\sum_{n\ge0}(-1)^n\dfrac{((z-k)/2+1)_n}{((z+k)/2+1)_n}$$
$$\int_0^\infty\dfrac{\sinh(at)}{\cosh(t)}e^{-zt}\,dt=2a\sum_{n\ge1}\dfrac{(-1)^{n+1}}{(2n+z+a-1)(2n+z-a-1)}$$
\begin{align*}\G(\nu+1)(2/z)^{\nu}I_{\nu}(z)&=1+\dfrac{z^2}{4\nu+4}\sum_{n\ge0}\dfrac{(z/2)^{2n}}{(n+1)!(\nu+2)_n}\\
\G(\nu+1)(2/z)^{\nu}e^{-z}I_{\nu}(z)&=1-z\sum_{n\ge0}\dfrac{(\nu+3/2)_n}{(n+1)!(2\nu+2)_n}(-2z)^n\\
I_0(z)-zI_1(z)&=1-\dfrac{z^2}{4}\sum_{n\ge0}\dfrac{(2n+1)}{(n+1)!^2}(z/2)^{2n}\\
(z/2)^{-2\nu}\G(\nu+1)^2I_{\nu}(z)^2&=1+\dfrac{z^2}{2(\nu+1)}\sum_{n\ge0}\dfrac{(\nu+3/2)_n}{(n+1)!(\nu+2)_n(2\nu+2)_n}z^{2n}\\
I_0(z)^2&=1+\dfrac{z^2}{2}\sum_{n\ge0}\dfrac{(3/2)_n}{(n+1)!^3}z^{2n}\\
(z/2)^{-2\nu}\G(\nu+1)^2(I_{\nu}(z)^2-I_{\nu+1}(z)^2)&=1+\dfrac{2\nu+1}{4(\nu+1)^2}z^2\sum_{n\ge0}\dfrac{(\nu+3/2)_n}{(n+1)!(\nu+2)_n(2\nu+3)_n}z^{2n}\\
I_0(z)^2-I_1(z)^2&=1+\dfrac{z^2}{2}\sum_{n\ge0}\dfrac{(3/2)_n}{(n+2)(n+1)!^3}z^{2n}\\
(z/2)^{-2\nu}\G(\nu+1)^2(I_{\nu}(z)^2+I_{\nu+1}(z)^2)&=1+\dfrac{2\nu+3}{4(\nu+1)^2}z^2\sum_{n\ge0}\dfrac{(\nu+5/2)_n}{(n+1)!(\nu+2)_n(2\nu+3)_n}z^{2n}\\
I_0(z)^2+I_1(z)^2&=1+\dfrac{3z^2}{2}\sum_{n\ge0}\dfrac{(5/2)_n}{(n+2)(n+1)!^3}z^{2n}\\
I_0(z)I_1(z)&=\dfrac{z}{2}+\dfrac{3z^3}{4}\sum_{n\ge0}\dfrac{(5/2)_n}{(n+1)!(n+2)!^2}z^{2n}\\
I_0(z)J_0(z)&=\sum_{n\ge0}\dfrac{(-z^4/64)^n}{(1/2)_nn!^3}\end{align*}

\section{An Observation on Gamma Quotients}

Looking at the list of hypergeometric formulas for quotients of values of
the gamma function, one can observe an interesting pattern for which I
have no explanation.

Let $$S=\prod_{1\le i\le g}\G(a_i)^{e_i}$$
be a gamma quotient, with $e_i\in\Z$. Adding if necessary $\G(1)$ factors,
we may always assume that $\sum_{1\le i\le g}e_i=0$. We observe that in a large
number of cases (given below), there exists a hypergeometric formula
$$S=\prod_{1\le i\le g}\G(a_i)^{e_i}=u+\sum_{n\ge0}h(n)\prod_{1\le i\le g}(a_i)_n^{-e_i}$$
for suitable $u$ and rational function $h(n)$. Note that adding an integer
to $a_i$ on the LHS simply multiplies it by a rational number, and on the RHS
simply changes the rational function $h(n)$, so the $a_i$ can be considered
modulo $1$.

Since $(a)_n=\G(a+n)/\G(a)$, note that the above identity can also be written
$$-u\prod_{1\le i\le g}\G(a_i)^{-e_i}=\sum_{n\ge0}h(n)\prod_{1\le i\le g}\G(n+a_i)^{-e_i}$$

I now give the list where this occurs, abbreviating $S$ above as
$\prod_{1\le i\le g}a_i^{e_i}$, and repeating as many times as it occurs.

\begin{align*}
  &\dfrac{(1/4)^2}{1^1(1/2)^1}\quad\dfrac{1^1(1/2)^1}{(1/4)^2}\quad\dfrac{(3/4)^2}{1^1(1/2)^1}\dfrac{1^1(1/2)^1}{(3/4)^2}\\
  &\dfrac{(1/4)^1(1/2)^1}{(3/4)^11^1}\quad\dfrac{(3/4)^11^1}{(1/4)^1(1/2)^1}\quad\dfrac{(3/4)^1(1/2)^1}{(1/4)^11^1}\quad\dfrac{(1/4)^11^1}{(3/4)^1(1/2)^1}\\
  &\dfrac{(1/4)^2(1/2)^2}{(3/4)^21^2}\quad\dfrac{(3/4)^21^2}{(1/4)^2(1/2)^2}\quad\dfrac{(3/4)^2(1/2)^2}{(1/4)^21^2}\quad\dfrac{(1/4)^21^2}{(2/4)^2(1/2)^2}\\
  &\dfrac{(1/4)^2}{(3/4)^2}\quad\dfrac{(1/4)^2}{(3/4)^2}\quad\dfrac{(3/4)^2}{(1/4)^2}\dfrac{(3/4)^2}{(1/4)^2}\quad
  \dfrac{(1/4)^4}{(3/4)^4}\quad\dfrac{(1/4)^4}{(3/4)^4}\quad\dfrac{(3/4)^4}{(1/4)^4}\dfrac{(3/4)^4}{(1/4)^4}\\
  &\dfrac{(1/3)^3}{1^3}\quad\dfrac{(1/3)^3}{1^3}\quad\dfrac{1^3}{(1/3)^3}\quad\dfrac{1^3}{(1/3)^3}\quad
  \dfrac{(2/3)^3}{1^3}\quad\dfrac{(2/3)^3}{1^3}\quad\dfrac{1^3}{(2/3)^3}\quad\dfrac{1^3}{(2/3)^3}\\
  &\dfrac{(1/3)^2}{(2/3)^1}\quad\dfrac{(2/3)^1}{(1/3)^2}\quad\dfrac{(2/3)^2}{(1/3)^1}\quad\dfrac{(1/3)^1}{(2/3)^2}\quad
  \dfrac{(1/3)^3}{(2/3)^3}\quad\dfrac{(2/3)^3}{(1/3)^3}\quad\dfrac{(1/6)^3}{(1/2)^3}\quad\dfrac{(1/2)^3}{(1/6)^3}\\
  &\dfrac{(1/2)^1(2/3)^2(5/6)^1}{1^1(1/6)^2(1/3)^1}\quad\dfrac{1^1(1/6)^2(1/3)^1}{(1/2)^1(2/3)^2(5/6)^1}\\
  &\dfrac{z^1(1/2)^1}{(z+1/2)^21^1}\quad\dfrac{(z+1/2)^11^1}{z^2(1/2)^1}\quad\dfrac{z^2}{(2z)^11^1}\quad\dfrac{(2z)^11^1}{z^2}\\
  &\dfrac{z^2}{(z+1/2)^2}\quad\dfrac{z^2}{(z+1/2)^2}\quad\dfrac{(z+1/2)^2}{z^2}\quad\dfrac{(z+1/2)^2}{z^2}\\
  &\dfrac{z^4}{(z+1/2)^4}\quad\dfrac{(z+1/2)^4}{z^4}\\
  &\dfrac{(a+1/2)^1b^1}{a^1(b+1/2)^1}\quad\dfrac{a^1(b+1/2)^1}{(a+1/2)^1b^1}\quad\dfrac{a^1b^1}{(a+1/2)^1(b+1/2)^1}\quad\dfrac{(a+1/2)^1(b+1/2)^1}{a^1b^1}
  \end{align*}

\chapter{The Pari/GP File of Continued Fractions}\label{chap:file}

This chapter is more or less the same as \cite{Coh5}, whose \emph{source
code} should be downloaded in any case to have the database.

It would be tedious to copy the 1900 or so continued fractions one by one
to be able to work on them. In addition, one may want to do searches for
CFs of a certain type, or other properties. For this, we provide a single
file, let us call it ``allcf.gp'', readable by {\tt Pari/GP} with the command

\begin{verbatim}
V=readvec("allcf.gp");
\end{verbatim}

\noindent
which allows to have a complete database of all the CFs given in this book,
together with their main properties. To obtain this file, download the
\TeX\ source of \cite{Coh5} and create the file ``allcf.gp'' by
extracting everything between the last two {\tt comments} of the source.

The aim of the present chapter is to describe in detail the format of this
database and how it is used. Once again, it is essential to have the
{\tt henri-ellCF3} branch of {\tt Pari/GP} installed,
although of course the format being transparent can be used with any other
software through an appropriate interface.

\medskip

First, a warning. We reserve the following variable names which must never
be used in working with the CFs of this file:

\begin{enumerate}\item ``n'', which is universally used to denote the CF index.
\item ``d'', which denotes the radicand of a quadratic extension, see below
  the explanation for {\tt v[3][1]}.
\item ``x'' which stands for $1/n$ and ``y'' which stands for $1/n^{1/2}$
  in the asymptotic expansions $A$, see the explanation for {\tt v[3][4]}.
\end{enumerate}

\smallskip

There are three kinds of inputs in the file. By far the most important ones are
of course the CFs themselves, which are also the most numerous.

\smallskip

$\bullet$ {\bf Format of the Continued Fractions Entries}

Such an entry, let's call it {\tt v}, has eight components, which we now
describe. To make the explanation clearer, we use the following example:

\begin{verbatim}
[(z)->cosh(Pi*z),[[1,1/2,4*n^2-8*n+(2*z^2+5)],
[2*z^2,-4*n^4+8*n^3+(-4*z^2-6)*n^2+(4*z^2+2)*n+(-z^2-1/4)]],
[1,[0,1,0,1],(z)->z^2*cosh(Pi*z),
1-1/2*z^2*x+(1/6*z^4-1/6*z^2-1/12)*x^2+(-1/24*z^4+1/6*z^2+1/12)*x^3+O(x^4)],
[[1,4*z^2,[1/2-I*z,1;1/2+I*z,1;3/2,-2],1],
 [1,-4*z^2,[-1/2,2;1/2-I*z,-1;1/2+I*z,-1],1]],
[],"3.2.4",["AP->3.2.5","APD->3.2.4.5"],""];
\end{verbatim}

\begin{enumerate}\item {\tt v[1]} is a closure giving the limit of the CF. Note
  that in CFs involving variables, the CF may converge to that limit only
  in a certain range of the variables, but we have not given that (we may
  in a future version).
\item {\tt v[2]} is the $[a(n),b(n)]$ of the CF in the format
  explained above, with one important difference: consider for instance the
  entry 1.5.10 beginning with:
\begin{verbatim}
  [k->zeta(k),[0,n^k+(n-1)^k],[1,-n^(2*k)]]
\end{verbatim}
Giving this to {\tt GP} would make it choke (give an error message) because
{\tt n\^{}k} is not a valid {\tt GP} input when {\tt k} is a variable.
  Thus we replace it in the file by
\begin{verbatim}
  [k->zeta(k),[k->[0,n^k+(n-1)^k],k->[1,-n^(2*k)]]]
\end{verbatim}
which is a valid input since the closure is not evaluated.
\item {\tt v[3]} is itself a 4-component vector {\tt [d,[F,E,D,P],C,A]}
  giving the exact speed of convergence of the CF, including a few terms of
  its asymptotic expansion. If $L$ denotes the value of the CF and $p(n)/q(n)$
  are the partial quotients, this notation means that
  $$L-\dfrac{p(n)}{q(n)}=\dfrac{C}{n!^FE^ne^{\sqrt{Dn}}n^P}A\;,$$
  with the following conventions:
  \begin{enumerate}
  \item {\tt v[3][1]} is a number or a polynomial (very often 1) with the
    following meaning: often, the speed of convergence and the asymptotic
    expansions involve elements of a quadratic extension (of $\Q$ or of
    $\Q[z]$ for instance). To avoid using square roots, we set {\tt v[3][1]}
    to the radicand, and we use the specific letter {\tt d} to indicate the
    square root of {\tt v[3][1]}. Thus, we must avoid the letter {\tt d} in
    variables of CFs. In the running example, we have {\tt v[3][1]=1} so
    no quadratic extension.
  \item {\tt v[3][2]} gives the speed of convergence the above FEDP format,
    using {\tt d} if necessary. The running example has
    {\tt v[3][2]=[0,1,0,1]}, meaning that we are in case $P^+$ with
    $P=1$. In some rare cases, the speed does not make any sense, in which
    case we set {\tt v[3][2]=0}. In some other rare
    cases, the CF does not converge for real values of the variables (which
    is implicitly what is assumed), in which case we set
    {\tt v[3][2]=[0,0,0,0]}.
  \item {\tt v[3][3]} gives the constant $C$ such that
    $S-p(n)/q(n)\sim C/[FEDP]$, or $0$ if unknown or has not been computed.
    Note that $C$ is always given as a closure, first because in some cases
    (as for the CF itself) it is not possible to do otherwise in {\tt GP},
    and second because writing {\tt C=()->Pi} is both much clearer to read
    than {\tt C=3.1415...} and also allows to have {\tt C} to any accuracy.
  \item {\tt v[3][4]} gives the asympotic expansion {\tt A=1+...}, such that
    $S-p(n)/q(n)=(C/[FEDP])A$, or $0$ if unknown or not computed.
    The reserved letter {\tt x} means $1/n$ and the reserved letter {\tt y}
    means $1/n^{1/2}$.

    Thus, the above example says that
    $$\cosh(\pi z)-\dfrac{p(n)}{q(n)}=\dfrac{z^2\cosh(\pi z)}{n}\left(1-\dfrac{z^2/2}{n}+\dfrac{z^4/6-z^2/2-1/12}{n^2}+\cdots\right)$$
  \end{enumerate}
\item In a large number of cases, the CF or its inverse corresponds to a more
  or less explicit hypergeometric series {\tt S} given in a special format
  explained below. Thus, {\tt v[4]} is set to {\tt 0} if neither CF nor
  its inverse has such a series, to {\tt S} if only the CF has a series,
  to the two-component vector {\tt [S,T]} if both have a series (with {\tt T}
  the series for the inverse), and finally the two-component vector
  {\tt [[],T]} if only the inverse has a series.

  These series are given in the following format: they are all of the form
  $$S=u+\sum_{n\ge0}h(n)\prod_{1\le i\le g}(a_i)_n^{e_i}w^n\;,$$
  where $h$ is a rational function, $e_i\in\Z$, and as usual $(a_i)_n$ denotes
  the rising Pochhammer symbol $\G(a_i+n)/\G(a_i)$ (note that the series
  always starts at $n=0$). This is represented by
  
\centerline{\tt S=[u,h(n),[a\_1,e\_1;a\_2,e\_2;...a\_g,e\_g],w]}

  \noindent
  (when $g=0$, the empty matrix is simply set to $1$).
  Note that for convergence, we always have
  $\sum_{1\le i\le g}e_i\le0$. Thus, the above example says that
  \begin{align*}\cosh(\pi z)&=1+4z^2\sum_{n\ge0}\dfrac{(1/2-iz)_n(1/2+iz)_n}{(3/2)_n^2}\text{\quad and}\\
    \dfrac{1}{\cosh(\pi z)}&=1-4z^2\sum_{n\ge1}\dfrac{(-1/2)_n^2}{(1/2-iz)_n(1/2+iz)_n}\end{align*}
  The fact that the series begins at $n=0$ has been chosen for uniformity,
  but often leads to slightly more cumbersome formulas, for instance
  instead of $\z(2)=\sum_{n\ge1}1/n^2$ we write instead
  $\z(2)=\sum_{n\ge0}1/(n+1)^2$, hence represented by
  {\tt [0,1/(n+1)\^{}2,1,1]}.
\item {\tt v[5]} is used for testing. When empty, it means that
  the default values for the variables used in the CF can be used
  (presently $2/3$, $3/4$, $1/2$, $2$ for closures with arity less or equal to
  $4$). If it is a vector of values, these are taken instead of the default
  values. For instance for a closure of arity $2$, {\tt v[5]=[2/3,-1]}
  means that one can take the default value $2/3$ for the first variable,
  but probably cannot use $3/4$ for the second, but $-1$ instead.
  Other values of {\tt v[5]} have specific meanings which may change:
  \begin{itemize}\item {\tt v[5]=1}: Bernoulli type sequence, test expansion.
  \item {\tt v[5]=2}: Jacobi type sequence, only test expansion.
  \item {\tt v[5]=3}: Jacobi type sequence, test expansion and numerics, (2,1/2).
  \item {\tt v[5]=4}: Change z into 1/z, multiply by z, test expansion.
  \item {\tt v[5]=5}: Choose $z=2i$, multiply by $-i$, test.
  \end{itemize}
  Note that in the present version, {\tt v[5]=4} and {\tt v[5]=5} occur
  each for a single CF.
\item {\tt v[6]} is the \TeX\ label used for the CF, so as to be
  able to find it easily. It is very useful to write a trivial little
  script to find conversely the given entry of the big vector {\tt V}
  corresponding to a given label. Note that the labels are rather arbitrary
  (but of course unique) since they were given after discovery of new CFs,
  and more or less in increasing order.
\item {\tt v[7]} is a 2-component vector [A1,AD1] whose entries are strings
  describing the behavior of the CF under Ap\'ery acceleration: A1 is
  for the Ap\'ery diagonal process, and AD1 for the Ap\'ery dual. The meanings
  of the strings are as follows: ``NO'' if Ap\'ery is not applicable,
  ``BAcomp'' if it is but the Bauer--Muir method introduces complicated
  denominators, ``AP\textrightarrow label'' (``APSI\textrightarrow label'' if followed by a simplification)
  and ``APD\textrightarrow label'' if Ap\'ery leads to the corresponding labeled CF,
  (``APD\textrightarrow SELF'' it is self-dual), ``APcomp'' and ``APDcomp'' if the result is
  too complicated to be included in the encyclopedia so with no label,
  and a possible additional ``sim'' if similar to.
\item {\tt v[8]} is a (usually empty) string giving any additional information.
\end{enumerate}

\smallskip

$\bullet$ {\bf Format of Parametrized Continued Fractions}

\smallskip

Many CFs are part of families of inequivalent CFs depending on one or more
parameters. These families are included in the file in a different format,
since the initial values of $a(n)$ and $b(n)$ are not given. We have
only given such families for CFs of period 1.

The corresponding vector {\tt v} has now only five components, and again
to make the explanation clearer we give an example:

\begin{verbatim}
[[()->Pi,2*k+2,(2*n-1)*(2*n+2*k-4*u-1)]
[1,[0,-1,0,k+1]],"1.2.1",[],"u<=k/2"];
\end{verbatim}

\begin{enumerate}\item The first component {\tt v[1]} is a three-component
  vector: {\tt v[1][1]} is the closure giving the limit $S$, and {\tt v[1][2]}
  and {\tt v[1][3]} are polynomials giving $a(n)$ and $b(n)$ for $n$
  sufficiently large. This means that for suitable initial values $a(0),a(1),\dotsc$, and $b(0),b(1),\dotsc$, and for given nonnegative integers $k$ and
  $u$, the corresponding CF converges to $S$. Equivalently, the CF
  {\tt [[v[1][2]],v[1][3]]} converges to a limit of the form $(aS+b)/(cS+d)$
  with $a$, $b$, $c$, $d$ integral and $ad-bc\ne0$. For instance, choosing
  $u=0$ and $k=0,1,2$, or $u=1$ and $k=2$, the above parametrized CF says that
  the limits of the CFs
  \begin{align*}
    &2+1/(2+1/(2+9/(2+25/(2+49/(2+81/(2+\cdots))))))\;,\\
    &4-1/(4+3/(4+15/(4+35/(4+63/(4+99/(4+\cdots))))))\;,\\
    &6-3/(6+5/(6+21/(6+45/(6+77/(6+117/(6+\cdots))))))\;,\text{\quad and}\\
    &6+1/(6+1/(6+9/(6+25/(6+49/(6+81/(6+\cdots))))))\end{align*}
  all belong to $\Q(\pi)$.
\item {\tt v[2]} contains the information on the speed of convergence
  in the shortened format {\tt [d,[F,E,D,P]]}, since the constant {\tt C}
  and the asymptotic expansion {\tt A} vary in a parametric family so
  are not given.
\item {\tt v[3]} is the label where the parametric CF occurs (same
  as {\tt v[5]} for ordinary CFs).
\item {\tt v[4]} is used for testing, but for now is the empty vector.
\item {\tt v[5]} is a string giving some additional information.
\end{enumerate}

Note that for integer $k$, the CF $(a(n+k),b(n+k))$ is essentially the
same as the CF $(a(n),b(n))$, so is of course not considered as a new CF,
so this trivial parametrization is not included.

\smallskip

$\bullet$ {\bf Format of Definitions}

\smallskip

In some of the CFs, it is necessary to add some function definitions,
either because the corresponding functions do not exist in {\tt Pari/GP}
(for example {\tt bernfrac(k)} and {\tt eulerfrac(k)} exist, but not
{\tt tanfrac(k)} for the tangent numbers), or simply because the existing
expression is cumbersome. The command {\tt V=readvec("file.gp");} will have two
effects: first, it defines the function so that it can be used in the CFs,
and second it is kept as an element of {\tt V} as a closure (but without
its name). Thus, if you do not know the meaning of a function,
typing {\tt ?functionname} will give you an explanation or a formula for it.

\medskip

{\bf More Details}

\medskip

Here is in more detail the different types of vectors $v$ that we can
encounter, with examples.

\begin{enumerate}
\item Definitions: as explained above. These are \emph{closures}. Example:
\begin{verbatim}
  R1(k,n)=(n^k*(2*n+1)+(n-1)^k*(2*n-3))/(2*n-1);
\end{verbatim}
\item Parametrized continued fractions: as explained above. These are
  \emph{vectors} with \emph{five} components. Example:
\begin{verbatim}
[[()->2^(1/3),7*n-5+k+v+4*u,-4*(n+u)*(3*n+v-2)],
            [1,[0,4/3,0,u+(5-v)/3+2*k]],"1.1.1",[],""];
\end{verbatim}
\item Ordinary CF with $a$ or $b$ expressed as a closure: These are
  \emph{vectors} $v$ with \emph{eight} components such that $v[2][1]$
  or $v[2][2]$ is a closure. Example:
\begin{verbatim}
[k->zeta(k),[k->[0,n^k+(n-1)^k],k->[1,-n^(2*k)]],[1,[0,1,0,k-1],
k->1/(k-1),1-((k-1)/2)*x+sum(m=1,4,binomial(k+2*m-2,2*m)*
bernfrac(2*m)*x^(2*m),O(x^10))],0,[3],"1.5.10",["NO","NO"],""];
\end{verbatim}
\item Ordinary CF with the last component of $a$ and $b$ a closure. Example:
\begin{verbatim}
  [()->1,[[0,n->f(n)],[n->f(n+1)+1]],[0,0,0,0],
                              0,[],"2.1.0.1",["NO","NO"],""];
\end{verbatim}
\item Ordinary CF of period 2 or more: These are vectors $v$ with eight
  components such that the first (hence any) component of $v[2][1]$ or
  $v[2][2]$ is a vector. Example:
\begin{verbatim}
[()->2^(1/3),[[[0,5],[2*n,12*n+6]],[[6,-1/3],[-(n+1)*(3*n+2),
-n*(3*n+1)]]],[2,[0,(1+d)^2,0,0],()->2^(1/3)*3^(1/2)/(1+sqrt(2))^2,
1+(9*d/8)*x+(-9*d/8+81/64)*x^2+(67105*d/41472-81/32)*x^3+O(x^4)],
0,[],"1.1.6",["NO","NO"],""];
\end{verbatim}
\item Ordinary CF of period 1 with no closures. Example:
\begin{verbatim}
[()->3^(1/2),[[3/2,(2*n+1)^2],[2,-n^2*(n+2)^2]],[3,[0,(2+d)^2,0,0],
()->2*sqrt(3)/(2+sqrt(3))^3,1+d*x+(3/2-3*d/2)*x^2+(3*d-9/2)*x^3
+(45/4-27*d/4)*x^4+O(x^5)],0,[],"1.1.0.5",
["NO","NO"],"Infinitely contractible"];
\end{verbatim}
A second example with Ap\'ery information and series for the inverse:
\begin{verbatim}
[()->2^(1/3),[[1/2,7*n-5],[1,-12*n^2+8*n]],[1,[0,4/3,0,5/3],
()->2^(8/3)/gamma(1/3),1-55/9*x+170/3*x^2-1579490/2187*x^3
+230596751/19683*x^4+O(x^5)],[[],[0,2,[-2/3,1;1,-1],3/4]],
[],"1.1.1",["APSI->1.1.8","APD->SELF"],""];
\end{verbatim}
where the string ``APSI\textrightarrow 1.1.8'' indicates that after
simplification, Ap\'ery leads to label 1.1.8, and ``APD\textrightarrow SELF''
that it is self-dual for Ap\'ery. The series entry says that
$2^{-1/3}=2\sum_{n\ge0}((-2/3)_n/n!)(3/4)^n$,
which is simply the special case $x=3/4$ of the expansion of $2(1-x)^{2/3}$.
\end{enumerate}

\medskip

Here is a little script which outputs a number from 1 to 6 giving the type
of a vector $v$:
\begin{verbatim}
/* Find the type number of an entry in the file */

typevec(v)=
{ my(tv=type(v),lv=#v,a,b);
  if(tv=="t_CLOSURE",return(1));
  if(tv!="t_VEC"||(lv!=5&&lv!=8),error("incorrect entry"));
  if(lv==5,return(2));
  [a,b]=v[2];
  if(type(a)=="t_CLOSURE"||type(b)=="t_CLOSURE",return(3));
  if(type(a)!="t_VEC"||type(b)!="t_VEC",error("incorrect entry 2"));
  if(type(a[#a])=="t_CLOSURE"||type(b[#b])=="t_CLOSURE",return(4));
  if(type(a[1])=="t_VEC"||type(b[1])=="t_VEC",return(5));
  return(6); }
\end{verbatim}

\smallskip

Here is a script which outputs the \emph{bidegree} for ordinary CFs of period
1 with no closures, and -1 otherwise.

\begin{verbatim}
bidegree(v)=
{ if(typevec(v)!=6,return(-1));
  return(apply(poldegree,cftopol(v))); }

findbidegree(da,db)=
{ my(V=readvec("allcf.gp"));
  for(i=1,#V,my(v=V[i]);if(bidegree(v)==[da,db],print1(i," "))); }
\end{verbatim}

\backmatter

\printindex

\end{document}